\documentclass{schwede}
\usepackage{amsmath} 
\usepackage{amsthm}
\usepackage{txfonts}

\usepackage{multind}
\makeindex{subject}
\makeindex{symbol}

\usepackage[curve,matrix,arrow,frame,tips]{xy}
\usepackage{mathrsfs}
\usepackage{enumerate}
\usepackage{stmaryrd}
\usepackage{manfnt}

\makeatletter
\providecommand*{\toclevel@schapter}{0} 
\makeatother
\usepackage[pdfborder={0 0 0}]{hyperref}

\DeclareMathOperator{\Wirth}{Wirth}
\DeclareMathOperator{\colim}{colim}
\DeclareMathOperator{\fc}{fc}
\DeclareMathOperator{\hocolim}{hocolim}

\DeclareMathOperator{\Ho}{Ho}
\DeclareMathOperator{\gl}{gl}
\DeclareMathOperator{\eig}{eig}

\DeclareMathOperator{\diag}{diag}
\DeclareMathOperator{\Hom}{Hom}

\DeclareMathOperator{\End}{End}
\DeclareMathOperator{\Nat}{Nat}
\DeclareMathOperator{\Mack}{-{\mathcal M}ack}
\DeclareMathOperator{\Pic}{Pic}
\DeclareMathOperator{\Prin}{Prin}
\DeclareMathOperator{\pic}{pic}
\DeclareMathOperator{\Tor}{Tor}

\DeclareMathOperator{\Cl}{Cl}
\DeclareMathOperator{\mCl}{{\mathbb C}l}
\DeclareMathOperator{\Out}{Out}

\DeclareMathOperator{\Sym}{Sym}
\DeclareMathOperator{\Spin}{Spin}
\DeclareMathOperator{\Pin}{Pin}

\DeclareMathOperator{\spn}{span}
\DeclareMathOperator{\coker}{coker}
\DeclareMathOperator{\map}{map}
\DeclareMathOperator{\Mod}{mod-}
\DeclareMathOperator{\Mon}{Mon}
\DeclareMathOperator{\mo}{-mod}

\DeclareMathOperator{\Id}{Id}
\DeclareMathOperator{\ind}{ind}
\DeclareMathOperator{\Ind}{Ind}

\DeclareMathOperator{\ad}{ad}
\DeclareMathOperator{\ev}{ev}
\DeclareMathOperator{\rel}{rel}
\DeclareMathOperator{\ort}{or}

\DeclareMathOperator{\odd}{odd}
\DeclareMathOperator{\op}{op}
\DeclareMathOperator{\str}{str}
\DeclareMathOperator{\Rep}{Rep}
\DeclareMathOperator{\sh}{sh}
\DeclareMathOperator{\sk}{sk}

\DeclareMathOperator{\tel}{tel}
\DeclareMathOperator{\tr}{tr}
\DeclareMathOperator{\gr}{gr}
\DeclareMathOperator{\Tr}{Tr}
\DeclareMathOperator{\Sing}{Sing}
\DeclareMathOperator{\Vect}{Vect}
\DeclareMathOperator{\res}{res}

\newcommand{\mA}{{\mathbb A}}
\newcommand{\mC}{{\mathbb C}}

\newcommand{\mF}{{\mathbb F}}
\newcommand{\mG}{{\mathbb G}}
\newcommand{\mH}{{\mathbb H}}
\newcommand{\mL}{{\mathbb L}}
\newcommand{\mN}{{\mathbb N}}
\newcommand{\mP}{{\mathbb P}}
\newcommand{\mQ}{{\mathbb Q}}
\newcommand{\mR}{{\mathbb R}}
\newcommand{\mS}{{\mathbb S}}
\newcommand{\RP}{{\mathbb R}\text{P}}
\newcommand{\mZ}{{\mathbb Z}}
\newcommand{\Ac}{{\mathcal A}}
\newcommand{\All}{{\mathcal A}ll}
\newcommand{\cyc}{cyc}
\newcommand{\Fin}{{{\mathcal F}in}}

\newcommand{\Ab}{{\mathcal A}b}

\newcommand{\Cc}{{\mathcal C}}

\newcommand{\Dc}{{\mathcal D}}
\newcommand{\Ec}{{\mathcal E}}
\newcommand{\Fc}{{\mathcal F}}
\newcommand{\Gc}{{\mathcal G}}
\newcommand{\Hc}{{\mathcal H}}

\newcommand{\Kc}{{\mathcal K}}
\newcommand{\Lc}{{\mathcal L}}
\newcommand{\M}{{\mathcal M}}

\newcommand{\Pc}{{\mathcal P}}

\newcommand{\Sc}{{\mathcal S}}
\newcommand{\Tc}{{\mathcal T}}
\newcommand{\Uc}{{\mathcal U}}
\newcommand{\Vc}{{\mathcal V}}
\newcommand{\Xc}{{\mathcal X}}
\newcommand{\Yc}{{\mathcal Y}}
\newcommand{\Zc}{{\mathcal Z}}
\newcommand{\bA}{{\mathbf A}}
\newcommand{\bB}{{\mathbf B}}
\newcommand{\bF}{\mathbf{F}}
\newcommand{\bGr}{\mathbf{Gr}}
\newcommand{\bMGr}{\mathbf{MGr}}
\newcommand{\bI}{{\mathbf I}}

\newcommand{\bK}{{\mathbf K}}
\newcommand{\bL}{{\mathbf L}}

\newcommand{\bO}{{\mathbf O}}

\newcommand{\bPin}{{\mathbf {Pin}}}
\newcommand{\bSpin}{{\mathbf {Spin}}}

\newcommand{\bP}{{\mathbf P}}
\newcommand{\bRO}{\mathbf{RO}}
\newcommand{\bIO}{\mathbf{IO}}
\newcommand{\bRU}{\mathbf{RU}}
\newcommand{\bRSp}{\mathbf{RSp}}
\newcommand{\bSO}{{\mathbf{SO}}}
\newcommand{\bSU}{{\mathbf{SU}}}

\newcommand{\bT}{{\mathbf T}}
\newcommand{\bU}{{\mathbf U}}
\newcommand{\bSp}{\mathbf{Sp}}

\newcommand{\bko}{{\mathbf{ko}}}
\newcommand{\bku}{{\mathbf{ku}}}
\newcommand{\bKO}{{\mathbf{KO}}}
\newcommand{\bKU}{{\mathbf{KU}}}
\newcommand{\bKSp}{{\mathbf{KSp}}}
\newcommand{\bmU}{{\mathbf{mU}}}
\newcommand{\bMU}{{\mathbf{MU}}}
\newcommand{\bmUP}{{\mathbf{mUP}}}
\newcommand{\bMUP}{{\mathbf{MUP}}}
\newcommand{\bBU}{{\mathbf{BU}}}
\newcommand{\bBUP}{{\mathbf{BUP}}}
\newcommand{\bBSpP}{{\mathbf{BSpP}}}

\newcommand{\bBO}{{\mathbf{BO}}}
\newcommand{\bbO}{{\mathbf{bO}}}
\newcommand{\bBSp}{{\mathbf{BSp}}}

\newcommand{\bbU}{{\mathbf{bU}}}
\newcommand{\bbSp}{{\mathbf{bSp}}}
\newcommand{\bbUP}{{\mathbf{bUP}}}
\newcommand{\bbSpP}{{\mathbf{bSpP}}}
\newcommand{\bbOP}{{\mathbf{bOP}}}
\newcommand{\bBOP}{{\mathbf{BOP}}}
\newcommand{\bMO}{{\mathbf{MO}}}
\newcommand{\bmO}{{\mathbf{mO}}}
\newcommand{\bmSO}{{\mathbf{mSO}}}

\newcommand{\bmOP}{{\mathbf{mOP}}}
\newcommand{\bMOP}{{\mathbf{MOP}}}
\newcommand{\bkuc}{{\mathbf{ku}^c}}
\newcommand{\bk}{{\mathscr{C}}}
\newcommand{\bDelta}{{\mathbf{\Delta}}}
\newcommand{\bGamma}{{\mathbf{\Gamma}}}

\newcommand{\bSpc}{\mathbf{Spc}}
\newcommand{\spec}{{\mathcal S}p}
\newcommand{\Spinf}{S\!p^\infty}
\newcommand{\spc}{spc}
\newcommand{\iso}{\cong}
\newcommand{\GF}{\mathcal{GF}}
\newcommand{\GlGre}{\mathcal{G}l\mathcal{G}re}
\newcommand{\GlPow}{\mathcal{G}l\mathcal{P}ow}
\newcommand{\sm}{\wedge}
\newcommand{\tensor}{\otimes}
\newcommand{\widebar}{\overline}
\newcommand{\xra}{\xrightarrow}
\newcommand{\xla}{\xleftarrow}
\newcommand{\bs}{\backslash}
\newcommand{\un}{\underline}
\newcommand{\Com}{{\mathcal{C}om}}
\newcommand{\GH}{{\mathcal{GH}}}
\newcommand{\SH}{{\mathcal{SH}}}
\newcommand{\upi}{{\underline \pi}}
\newcommand{\dslash}{\!\mathbin{/\mkern-11mu/}}
\newcommand{\td}[1]{\langle #1\rangle}
\newcommand{\gh}[1]{\llbracket #1\rrbracket}
\renewcommand{\to}{\longrightarrow}

\numberwithin{equation}{section}
\newtheorem{theorem}[equation]{Theorem}
\newtheorem{lemma}[equation]{Lemma}
\newtheorem{cor}[equation]{Corollary}
\newtheorem{prop}[equation]{Proposition}
 
\theoremstyle{definition}
\newtheorem{defn}[equation]{Definition}
\newtheorem{rk}[equation]{Remark}
\newtheorem{eg}[equation]{Example}
\newtheorem{construction}[equation]{Construction}

\newcommand{\Danger}{%
   \par\penalty-500\begingroup\parindent=0.0pt\clubpenalty=10000%
   \def\par{\endgraf\endgroup}\leavevmode
   \hangindent=2.4em\relax
   \hangafter=-2\relax
   \hbox to 0pt{\hss\relax
      \vbox to 0pt{\vskip-8pt\relax
        \hbox
          {\relax\dbend}
\vss}\kern.70em}}

\begin{document}

\title{Global homotopy theory}
\author{Stefan Schwede\\[25\baselineskip]
{Mathematisches Institut, Universit{\"a}t Bonn\\
{\tt schwede@math.uni-bonn.de}
}}

\bookkeywords{equivariant homotopy theory, compact Lie group, spectrum}

\frontmatter

\maketitle
\tableofcontents

\chapter*{Preface}

Equivariant stable homotopy theory has a long tradition,
starting from geometrically motivated questions about symmetries of manifolds. 
The homotopy theoretic foundations of the
subject were laid by tom Dieck, Segal and May and their students and
collaborators in the 70's, and during the last decades 
equivariant stable homotopy theory 
has been very useful to solve computational and conceptual problems 
in algebraic topology, geometric topology and algebraic $K$-theory. 
Various important equivariant theories naturally exist not just 
for a particular group, but in a uniform way for all groups in a specific class. 
Prominent examples of this are
equivariant stable homotopy, equivariant $K$-theory or equivariant bordism.
{\em Global} equivariant homotopy theory studies such uniform phenomena,
i.e., the adjective `global' refers to simultaneous 
and compatible actions of all compact Lie groups. 

This book introduces a new context for global homotopy theory.
Various ways to provide a home for global stable homotopy types have 
previously been
explored in \cite[Ch.\,II]{lms}, \cite[Sec.\,5]{greenlees-may-completion}, 
\cite{bohmann-thesis} and \cite{bohmann-orthogonal}.
We use a different approach: we work with the well-known category
of orthogonal spectra, but use a much finer notion of equivalence,
the {\em global equivalences}, than what is traditionally considered.
The basic underlying observation is that an orthogonal spectrum
gives rise to an orthogonal $G$-spectrum for every compact Lie group $G$,
and the fact that all these individual equivariant objects come from
one orthogonal spectrum implicitly encodes strong compatibility conditions
as the group $G$ varies. 
An orthogonal spectrum thus has $G$-equivariant homotopy groups 
for every compact Lie group,
and a global equivalence is a morphism of orthogonal spectra that induces
isomorphisms for all equivariant homotopy groups for all compact Lie groups
(based on `complete $G$-universes', compare Definition \ref{def:global equivalence}).

The structure on the equivariant homotopy groups of an orthogonal spectrum
gives an idea of the information encoded in a global homotopy type
in our sense: the equivariant homotopy groups $\pi_k^G(X)$ are contravariantly
functorial for continuous group homomorphisms (`restriction maps'), and they are
covariantly functorial for inclusions of closed subgroups (`transfer maps').
The restriction and transfer maps enjoy various transitivity properties
and interact via a double coset formula. This kind of algebraic structure has
been studied before under different names, 
e.g., `global Mackey functor', `inflation functor', \dots.
From a purely algebraic perspective, there are various parameters here than one can vary,
namely the class of groups to which a value is assigned and the classes
of homomorphisms to which restriction maps respectively transfer maps are assigned,
and lots of variations have been explored.
However, the decision to work with orthogonal spectra and equivariant
homotopy groups on complete universes dictates a canonical choice: we prove
in Theorem \ref{thm:Burnside category basis} 
that the algebra of natural operations between the equivariant homotopy groups 
of orthogonal spectra is freely generated by restriction maps 
along continuous group homomorphisms and transfer maps along closed subgroup inclusion, 
subject to explicitly understood relations.

We define the {\em global stable homotopy category} $\GH$
by localizing the category of orthogonal spectra at the class of
global equivalences.
Every global equivalence is in particular a non-equivariant stable equivalence,
so there is a `forgetful' functor $U:\GH\to \SH$ on localizations,
where $\SH$ denotes the traditional non-equivariant stable homotopy category.
By Theorem \ref{thm:change families} this forgetful functor
has a left adjoint $L$ and a right adjoint $R$, both fully faithful,
that participate in a recollement of triangulated categories:
\[ \xymatrix@C=15mm{ 
\GH^+ \ar@<-.3ex>[r]^-{i_*}  & 
\GH\ \ar@<-.3ex>[r]^-U  \ \ar@<.4ex>@/^1pc/[l]^-{i^!} \ar@<-.4ex>@/_1pc/[l]_-{i^*} & 
\SH \ \ar@<.4ex>@/^1pc/[l]^-R \ar@<-.4ex>@/_1pc/[l]_-L }
 \]
Here $\GH^+$ denotes the full subcategory of the global stable homotopy category
spanned by the orthogonal spectra that are stably contractible in the traditional,
non-equivariant sense.

The global sphere spectrum and suspension spectra are
in the image of the left adjoint 
(Example \ref{eg:suspensions are left induced}).
Global Borel cohomology theories are the image of the right adjoint
(Example \ref{eg:global Borel}).
The `natural' global versions of 
Eilenberg-Mac\,Lane spectra (Construction \ref{con:HM}),
Thom spectra (Section \ref{sec:global Thom}), 
or topological $K$-theory spectra 
(Sections \ref{sec:connective global K} and \ref{sec:global K})
are not in the image of either of the two adjoints.
Periodic global $K$-theory, however, is right induced from finite cyclic groups,
i.e., in the image of the analogous right adjoint from an intermediate 
global homotopy category $\GH_{\cyc}$ based on finite cyclic groups
(Example \ref{eg:Global K is right induced}).

Looking at orthogonal spectra through the eyes of global equivalences
is like using a prism: the latter breaks up white light into a spectrum
of colors, and  global equivalences split a traditional, non-equivariant homotopy
type into many different global homotopy types. 
The first example of this phenomenon that we will encounter
refines the classifying space of a compact Lie group $G$.
On the one hand, there is the constant orthogonal space with value
a non-equivariant model for $B G$; 
and there is the {\em global classifying space} $B_{\gl} G$ 
(see Definition \ref{def:global classifying}). 
The global classifying space is analogous to
the geometric classifying space of a linear algebraic group
in motivic homotopy theory \cite[4.2]{morel-voevodsky},
and it is the counterpart to the stack of $G$-principal bundles
in the world of stacks.

Another good example is the splitting up of the non-equivariant homotopy
type of the classifying space of the infinite orthogonal group.
Again there is the constant orthogonal space with value $B O$, 
the Grassmannian model $\bBO$ (Example \ref{eg:BOP}),
a different Grassmannian model $\bbO$ (Example \ref{eg:define bO}),
the bar construction model $\bB^\circ\bO$ (Example \ref{eg:bar construction BO}), 
and finally a certain `cofree' orthogonal space $R(B O)$.
The orthogonal space $\bbO$ is also a homotopy colimit, 
as $n$ goes to infinity, of the global classifying spaces $B_{\gl} O(n)$.
We discuss these different global forms of $B O$ 
is some detail in Section \ref{sec:forms of BO},
and the associated Thom spectra in Section \ref{sec:global Thom}.

In the stable global world, every non-equivariant homotopy type
has two extreme global refinements, the `left induced'
(the global analog of a constant orthogonal space,
see Example \ref{eg:reflexive example trivial})
and the `right induced' global homotopy type (representing Borel cohomology theories,
see Example \ref{eg:global Borel}).
Many important stable homotopy types have other natural global forms.
The non-equivariant Eilenberg-Mac\,Lane spectrum of the integers has a 
`free abelian group functor' model (Construction \ref{con:HM}),
and another incarnation as the Eilen\-berg-Mac\,Lane spectrum
of the constant global functor with value $\mZ$
(Remark \ref{rk:general Eilenberg Mac Lane}).
These two global refinements of the integral Eilenberg-Mac\,Lane spectrum
agree on finite groups, but differ for compact Lie groups of positive dimensions.

As already indicated, there is a great variety of orthogonal Thom spectra,
in real (or unoriented) flavors as $\bmO$ and $\bMO$,
as complex (or unitary) versions $\bmU$ and $\bMU$,
and there are periodic versions $\bmOP$, $\bMOP$, $\bmUP$ and $\bMUP$ of these;
we discuss these spectra in Section \ref{sec:global Thom}.
The theories represented by $\bmO$ and $\bmU$ have the closest ties to geometry;
for example, the equivariant homotopy groups of $\bmO$
receive Thom-Pontryagin maps from equi\-variant bordism rings,
and these are isomorphisms for products of finite groups and tori
(compare Theorem \ref{thm:TP is iso}).
The theories represented by $\bMO$ are tom Dieck's homotopical equivariant bordism,
isomorphic to `stable equivariant bordism'.

Connective topological $K$-theory also has two
fairly natural global refinements, in addition to the left and right induced ones.
The `orthogonal subspace' model $\bku$ 
(Construction \ref{con:connective global K-theory})
represents connective equivariant $K$-theory on the class of finite groups;
on the other hand, global connective $K$-theory $\bkuc$
(Construction \ref{con:global connective ku})
is the global synthesis of equivariant connective $K$-theory
in the sense of Greenlees \cite{greenlees-connective}.
The periodic global $K$-theory spectrum $\bKU$ is introduced in 
Construction \ref{con:global KU}; as the name suggests, $\bKU$
is Bott periodic and represents equivariant $K$-theory.

The global equivalences are part of a closed model structure
(see Theorem \ref{thm:All global spectra}), 
so the methods of homotopical algebra can be used to study the stable global homotopy
category. This works more generally relative to a class $\mathcal F$
of compact Lie groups, where we define $\mathcal F$-equivalences
by requiring that $\pi_k^G(f)$ is an isomorphism for all integers and all
groups in $\mathcal F$.
We call a class $\mathcal F$ of compact Lie groups a {\em global family}
if it is closed under isomorphisms, subgroups and quotients.
For global families we refine the $\mathcal F$-equivalences to a stable model structure,
the {\em $\Fc$-global model structure}, see Theorem \ref{thm:F-global spectra}.
Besides all compact Lie groups,
interesting global families are the classes of all finite groups,
or all abelian compact Lie groups. The class of trivial groups
is also admissible here, but then we just recover the `traditional' stable category.
If the family $\Fc$ is multiplicative, then the 
 $\Fc$-global model structure is monoidal with respect to
the smash product of orthogonal spectra and satisfies the monoid axiom
(Proposition \ref{prop:monoid axiom}).
Hence this model structure lifts to modules over an orthogonal ring
spectrum and to algebras over an ultra-commutative ring spectrum
(Corollary \ref{cor-lift modules}).

\smallskip

{\bf Ultra-commutativity.}
A recurring theme throughout this book is a phenomenon 
that I call {\em ultra-commutativity}.
I use this term in the unstable and stable context
for the homotopy theory of strictly commutative objects under global equivalences.
An ultra-commutative multiplication is significantly more structure
than just a coherently homotopy-commutative product
(usually called an $E_\infty$-multiplication). 
For example, the extra structure gives rise to power operations
that can be turned into transfer maps (in additive notation) 
respectively norm  maps (in multiplicative notation).
Another difference is that an unstable global $E_\infty$-structure 
would give rise to naive deloopings (i.e., by trivial representations). 
As I hope discuss elsewhere, a global ultra-commutative multiplication,
in contrast, gives rise to `genuine' deloopings (i.e., by non-trivial representations). 
As far as the objects are concerned, ultra-commutative monoids and
ultra-commutative ring spectra are not at all new and have been much studied before;
so one could dismiss the name 'ultra-commutativity' as a mere marketing maneuver.
However, the homotopy theory of ultra-commutative monoids and 
ultra-commutative ring spectra
with respect to global equivalences is new and, in the authors opinion,
important. And important concepts deserve catchy names.

\smallskip
\newpage

{\bf Global homotopy types as orbifold cohomology theories.} \index{subject}{orbifold|(}
I would like to briefly mention another reason for
why one might be interested in global stable homotopy theory.
In short, global stable homotopy types represent 
genuine cohomology theories on stacks, orbifolds, and orbispaces.
Stacks and orbifolds are concepts from algebraic geometry respectively
geometric topology that allow us to talk about objects
that locally look like the quotient of a smooth object by a group action,
in a way that remembers information about the isotropy groups of the action.
Such `stacky' objects can behave like smooth objects
even if the underlying spaces have singularities.
As for spaces, manifolds and schemes, cohomology theories
are important invariants also for stacks and orbifolds,
and examples such as ordinary cohomology or $K$-theory lend themselves
to generalization. Special cases of orbifolds are `global quotients',
often denoted $M\dslash G$, for example for a smooth action of a compact
Lie group $G$ on a smooth manifold $M$. In such examples, the orbifold
cohomology of $M\dslash G$ is supposed to be the $G$-equivariant cohomology of $M$.
This suggests a way to {\em define} orbifold cohomology theories by
means of equivariant stable homotopy theory, via suitable $G$-spectra $E_G$. 
However, since the group $G$ is not intrinsic and can vary, one needs
equivariant cohomology theories for all groups $G$, with some compatibility.

Part of the compatibility can be deduced from the fact that 
the same orbifold can be presented in different ways; for example, if $G$ is a closed
subgroup of $K$, then the global quotients
$M\dslash G$ and $(M\times_G K)\dslash K$ describe the same orbifold.
So if the orbifold cohomology theory is represented by equivariant spectra
$\{E_G\}_G$ as indicated above, then necessarily $E_G\simeq \res^K_G(E_K)$, i.e., 
the equivariant homotopy types are consistent under restriction.
This is the characteristic feature of {\em global} equivariant homotopy types, 
and it suggest that the latter ought to define orbifold cohomology theories.

The approach to global homotopy theory presented in this book
in particular provides a way to turn the above outline into rigorous mathematics.
There are different formal frameworks for 
stacks and orbifolds (algebro-geometric, smooth, topological), 
and these objects can be studied with respect to various notions of `equivalence'.
The approach that most easily feeds into our present context 
are the notions of {\em topological stacks}
respectively {\em orbispaces} as developed by Gepner and Henriques
in their paper \cite{gepner-henriques}. Their homotopy theory of topological stacks
is rigged up so that the derived mapping spaces out of 
the classifying stacks for principal $G$-bundles detect equivalences.
In our setup, the global classifying spaces of compact Lie groups
(see Definition \ref{def:global classifying})
play exactly the same role, and this is another hint of a deeper connection.
In fact, the global homotopy theory of orthogonal spaces
as developed in Chapter \ref{ch-unstable}
is a model for the homotopy theory of orbispaces in the sense of
Gepner and Henriques. 
For a formal comparison of the two models we refer the reader to
the author's paper \cite{schwede-orbispaces}.
The comparison proceeds through yet another model, the
global homotopy theory of `spaces with an action
of the universal compact Lie group'. 
Here the universal compact Lie group (which is neither compact nor a Lie group)
is the topological monoid $\Lc$ of linear isometric self-embeddings
of $\mR^\infty$, and in \cite{schwede-orbispaces} 
we establish a global model structure on the
category of $\Lc$-spaces.

If we now accept that one can pass between stacks, orbispaces and
orthogonal spaces in homotopically meaningful way, a consequence is
that every global stable homotopy type (i.e., every orthogonal spectrum)
gives rise to a cohomology theory on stacks and orbifolds.
Indeed, by taking the unreduced suspension spectrum,
every unstable global homotopy type is transferred into the triangulated
global stable homotopy category $\GH$. In particular, taking morphisms
in $\GH$ into an orthogonal spectrum $E$ defines  $\mZ$-graded 
$E$-cohomology groups. 
The counterpart of a global quotient $M\dslash G$ in the global homotopy 
theory of orthogonal spaces is the semifree orthogonal space $\bL_{G,V} M$
introduced in Construction \ref{con:semifree spc}.
By the adjunction relating the global and $G$-equivariant stable homotopy 
categories (see Theorem \ref{thm:change to groups}),
the morphisms $\gh{\Sigma^\infty_+ \bL_{G,V}M, E}$ in the global stable
homotopy category biject with the $G$-equi\-variant $E$-co\-homo\-logy groups of $M$.
In other words, when evaluated on a global quotient $M\dslash G$, our recipe
for generating an orbifold cohomology theory from a global stable homotopy type
precisely returns the $G$-equivariant cohomology of $M$,
which was the original design criterion.

The procedure sketched so far actually applies to more general objects 
than our global stable homotopy types: indeed, all that was needed to produce
the orbifold cohomology theory was a sufficiently exact functor from 
the homotopy theory of orbispaces to a triangulated category.
If we aim for a stable homotopy theory (as opposed to its triangulated
homotopy category), then there is a universal example, namely
the stabilization of the homotopy theory of orbispaces, obtained by formally inverting 
suspension. Our global theory is, however, richer than 
this `naive' stabilization. Indeed, the global stable homotopy category forgets to
the $G$-equivariant stable homotopy category based on a complete $G$-universe;
the equivariant cohomology theories represented by such objects are
usually called `genuine' (as opposed to `naive').
Genuine equivariant cohomology theories have much more structure
than naive ones; this structure manifests itself in different forms, 
for example as transfer maps, stability under `twisted suspension' 
(i.e., smash product with linear representation spheres), 
an extension of the $\mZ$-graded cohomology groups to an $R O(G)$-graded theory,
and an equivariant refinement 
of additivity (the so called {\em Wirthm{\"u}ller isomorphism}).
Hence global stable homotopy types in the sense of this book
represent {\em genuine} (as opposed to `naive') orbifold cohomology theories.
\index{subject}{orbifold|)}

\smallskip

{\bf Organization.}
In Chapter \ref{ch-unstable} we set up unstable global homotopy theory
using orthogonal spaces, i.e., continuous functors from 
the category of finite-dimensional inner product spaces 
and linear isometric embeddings to spaces.
We introduce global equivalences (Definition \ref{def:global equivalence spaces}),
discuss global classifying spaces of compact Lie groups
(Definition \ref{def:global classifying}),
and set up the global model structures on the category of orthogonal spaces
(Theorem \ref{thm:All global spaces}). 
In Section \ref{sec:unstable monoidal}
we investigate the box product of orthogonal spaces 
from a global equivariant perspective.
Section \ref{sec:global families unstable}
introduces a variant of unstable global homotopy
theory based on a {\em global family}, i.e., a class $\Fc$ of 
compact Lie groups with certain closure properties. 
We discuss the $\Fc$-global model structure 
and record that for multiplicative global families, it lifts to category of modules and
algebras (Corollary \ref{cor-lift to modules spaces}).
In Section \ref{sec:equivariant homotopy sets} we discuss 
the $G$-equivariant homotopy sets of orthogonal spaces
and identify the natural structure between them
(restriction maps along continuous group homomorphisms).
The study of natural operations on $\pi_0^G(Y)$ 
is a recurring theme throughout this book;
in the later chapters we return to it in the contexts of
ultra-commutative monoids, orthogonal spectra and ultra-commutative ring spectra.

Chapter \ref{ch-umon} is devoted to ultra-commutative monoids
(a.k.a. commutative mo\-noids with respect to the box product, 
or lax symmetric monoidal functors), which we want to advertise as 
a rigidified notion of `global $E_\infty$-space'.
In Section \ref{sec:global model monoid spaces}
we establish a global model structure for ultra-commutative monoids
(Theorem \ref{thm:global umon}). 
Section \ref{sec:power op spaces} introduces and studies
global power monoids, the algebraic structure that an ultra-commutative multiplication
gives rise to on the homotopy group $\Rep$-functor $\upi_0(R)$.
Section \ref{sec:umon examples} contains a large collection
of examples of ultra-commutative monoids and interesting morphisms between them.
In Section \ref{sec:forms of BO} 
we discuss and compare different global refinements of
the non-equivariant homotopy type $B O$, 
the classifying space for the infinite orthogonal group.
Section \ref{sec:group completion and units}
discusses `units' and `group completions' of ultra-commutative monoids.
As an application of this technology we formulate and prove 
a global, highly structured version of Bott periodicity,
see Theorem \ref{thm:global Bott periodicity}.

Chapter \ref{ch-equivariant} is a largely self-contained 
exposition of many basics about equivariant stable homotopy theory
for a fixed compact Lie group, modeled by orthogonal $G$-spectra.
In Section \ref{sec:equivariant homotopy groups} 
we recall orthogonal $G$-spectra and equivariant homotopy groups
and prove their basic properties, such as the suspension isomorphism 
and long exact sequences of mapping cones and homotopy fibers, 
and the additivity of equivariant homotopy groups on sums and products.
Section \ref{sec:Wirthmuller and transfer} discusses the Wirthm{\"u}ller isomorphism
and the closely related transfers.
In Section \ref{sec:geometric fixed points} we introduce and study
geometric fixed point homotopy groups, an alternative invariant
to characterize equivariant stable equivalences. 
Section \ref{sec:double coset} contains a proof of 
the double coset formula for the composite of a transfer followed by
the restriction to a closed subgroup.
We review Mackey functors for finite groups and show that
after inverting the group order, the category of $G$-Mackey functors
splits as a product, indexed by conjugacy classes of subgroups,
of module categories over the Weyl groups (Theorem \ref{thm:split Mackey functors}).
A topological consequence is that after inverting the group order,
equivariant homotopy groups and geometric fixed point homotopy groups
determine each other in a completely algebraic fashion,
compare Proposition \ref{prop:rational geometric fixed for GSp} 
and Corollary \ref{cor:equivariant from geometric}.
Section \ref{sec:products} is devoted to multiplicative aspects of
equivariant stable homotopy theory.

Chapter \ref{ch-stable} sets the stage for stable 
global homotopy theory, based on the notion of global equivalences
for orthogonal spectra (Definition \ref{def:global equivalence}).
We discuss semifree orthogonal spectra
and identify certain morphisms between semifree orthogonal spectra
as global equivalences (Theorem \ref{thm:faithful independence}).
In Section \ref{sec:global functors}
we investigate {\em global functors}, the natural algebraic structure
on the collection of equivariant homotopy groups of a global stable homotopy type.
Among other things, we explicitly calculate the 
algebra of natural operations on equivariant homotopy groups
(Theorem \ref{thm:Burnside category basis}).
In Section \ref{sec:global model structures} we complement 
the global equivalences of orthogonal spectra by a stable model structure
that we call the {\em global model structure}.
Its fibrant objects are the `global $\Omega$-spectra' 
(Definition \ref{def:global Omega}),
the natural concept of a `global infinite loop space' in our setting.
Here we work more generally relative to a global family $\Fc$ and consider the 
$\Fc$-equivalences (i.e., equivariant stable equivalences for all compact
Lie groups in the family $\Fc$).
We follow the familiar outline: a certain $\Fc$-level model structure
is Bousfield localized to an $\Fc$-global model structure
(see Theorem \ref{thm:F-global spectra}).
In Section \ref{sec:generators} we develop some basic theory 
around the global stable homotopy category;
since it comes from a stable model structure, this category is naturally
triangulated and we show that the suspension spectra of
global classifying spaces form a set of compact generators 
(Theorem \ref{thm:Bgl G weak generators}).
In Section \ref{sec:change family} we vary the global family:
we construct and study left and right adjoints
to the forgetful functors associated to a change of global family
(Theorem \ref{thm:change families}).
As an application of Morita theory for stable model categories
we show that rationally the global homotopy category for finite groups
has an algebraic model, namely the derived category of rational global functors
(Theorem \ref{thm:rational SH}).

Chapter \ref{ch:ultra} focuses on
{\em ultra-commutative ring spectra}, i.e.,  commutative orthogonal ring spectra
under multiplicative global equivalences. 
Section \ref{sec:power ops ring spectra} introduces `global power functors',
the algebraic structure on the equivariant homotopy
groups of ultra-commutative ring spectra.
Roughly speaking, global power functors are
global Green functors equipped with additional power operations,
satisfying various properties reminiscent
of those of the power maps $x\mapsto x^m$ in a commutative ring.
The power operations give rise to norm maps
(`multiplicative transfers') along finite index inclusions,
and in our global context, the norm maps conversely determine the power operations,
compare Remark \ref{rk:power vs Tambara}.
As we show  in Theorem \ref{thm:pi_0 is global power functor},
the 0-th equivariant homotopy groups of an ultra-commutative ring spectrum
form a global power functor.
In Section \ref{sec:algebra global power} 
we develop a description of the category of global power functors
via the comonad of `exponential sequences' 
(Theorem \ref{thm:power comonadic})
and discuss localization of global power functors
at a multiplicative subset of the underlying ring
(Theorem \ref{thm:localize global power functor}).
In Section \ref{sec:GPF examples} we give various examples of 
global power functors, such as the Burnside ring global power functor, 
the global functor represented by an abelian compact Lie group,
free global power functors, constant global power functors,
and the complex representation ring global functor.
In Section \ref{sec:ucom global model} we establish the global model structure
for ultra-\-commu\-tative ring spectra (Theorem \ref{thm:global ultra-commutative})
and show that every global power functor is realized 
by an ultra-commutative ring spectrum
(Theorem \ref{thm:realize global power functor}).

Chapter \ref{ch:examples} is devoted to interesting examples 
of ultra-commutative ring spectra.
Section \ref{sec:global Thom} discusses two orthogonal 
Thom spectra $\bmO$ and $\bMO$.
The spectrum $\bmO$ is globally connective
and closely related to equivariant bordism. 
The global functor $\upi_0(\bmO)$ admits a short and elegant algebraic presentation:
it is obtained from the Burnside ring global functor by imposing the 
single relation $\tr_e^{C_2}=0$, compare Theorem \ref{thm:pi_0 of mO}.
The Thom spectrum $\bMO$ was first considered by tom Dieck
and it represents `stable' equivariant bordism; it is periodic
for orthogonal representations of compact Lie groups, and admits Thom isomorphisms
for equivariant vector bundles.
The equivariant homology theory represented by $\bMO$ can be obtained
from the one represented by $\bmO$ in an algebraic fashion, by inverting
the collection of `inverse Thom classes', 
compare Corollary \ref{cor:MO is localized mO}.
Section \ref{sec:equivariant bordism} recalls the geometrically defined
equivariant bordism theories.
The Thom-Pontryagin construction 
maps the unoriented $G$-equivariant bordism ring $\mathcal N_*^G$
to the equivariant homotopy ring $\pi_*^G(\bmO)$,
and that map is an isomorphism when $G$ is a product of a finite group and a torus,
see Theorem \ref{thm:TP is iso}. 
We discuss global $K$-theory in Sections \ref{sec:connective global K} 
and \ref{sec:global K},
which comes in three interesting flavors as
connective global $K$-theory $\bku$,
global connective $K$-theory $\bkuc$ and 
periodic global $K$-theory $\bKU$ (and in the real versions
$\bko$, $\bko^c$ and $\bKO$).

We include three appendices where we collect material that is mostly
well-known, but that is either scattered through the literature
or where we found the existing expositions too sketchy.
Appendix \ref{app:CGWH} is a self-contained review of
compactly generated spaces, our basic category to work in.
Appendix \ref{app:equivariant spaces} deals with fundamental properties
of equivariant spaces, including the basic model structure 
in Proposition \ref{prop:proj model structures for G-spaces}.
We also provide an exposition of the
equivariant $\bGamma$-space machinery, culminating in a version
of the Segal-Shimakawa delooping machine.
In Appendix \ref{app:enriched functors} we review the basic definitions, 
properties and constructions involving categories of enriched functors.

While most of the material in the appendices is well-known, 
there are a few results I could not find in the literature.
These results include the fact that compactly generated spaces
are closed under geometric realization 
(Proposition \ref{prop:geometric realization in bT}~(iii)),
fixed points commute with geometric realization and latching objects
(Proposition \ref{prop:G-fix preserves pushouts}~(iv)),
and compactly generated spaces are closed under prolongation of $\bGamma$-spaces
(Proposition \ref{prop:prolongation in bT}).
Also apparently new are the results that prolongation of
$G$-cofibrant $\bGamma$-$G$-spaces to finite $G$-CW-complexes 
is homotopically meaningful (Proposition \ref{prop:cofibrant Gamma-G on G-CW}),
and that prolongation of $G$-cofibrant $\bGamma$-$G$-spaces 
to spheres gives rise to $G$-$\Omega$-spectra
(for very special $\bGamma$-$G$-spaces, see Theorem \ref{thm:prolonged delooping})
respectively positive $G$-$\Omega$-spectra (for special $\bGamma$-$G$-spaces,
see Theorem \ref{thm:special cofibrant Gamma positive Omega}).
Here the key ideas all go back to Segal \cite{segal-some_equivariant} 
and Shimakawa \cite{shimakawa}; however, we formulate our results for the
prolongation (i.e., categorical Kan extension), whereas Segal and Shimakawa
work with a bar construction (also known as a homotopy coend
or homotopy Kan extension) instead.
We also give a partial extension of the machinery to compact Lie groups,
whereas previous papers on the subject restrict attention to finite groups.
As we explain in Remark \ref{rk:Gamma limitations}, there is no hope
to obtain a $G$-$\Omega$-spectrum by evaluation on spheres
for compact Lie groups of positive dimension.
However, we do prove in Theorem \ref{thm:special cofibrant Gamma positive Omega} 
that evaluating a $G$-cofibrant special $\bGamma$-$G$-space
on spheres yields a `$G^\circ$-trivial positive $G$-$\Omega$-spectrum',
where $G^\circ$ is the identity component of $G$.
Our Appendix \ref{app:equivariant spaces} substantially overlaps
with the paper \cite{mmo} by May, Merling and Osorno that
provides comparisons of prolongation, 
bar construction and the operadic approach to equivariant deloopings.

\smallskip

{\bf Relation to other work.}
The idea of global equivariant homotopy theory is not at all new and has 
previously been explored in different contexts. For example,
in Chapter~II of \cite{lms}, Lewis and May define {\em coherent families}
of equivariant spectra; these consist of collections of equivariant 
coordinate free spectra in the sense Lewis, May and Steinberger,
equipped with comparison maps involving change of groups
and change of universe functors.

The approach closest to ours are the 
{\em global $\mathcal I_*$-functors} 
introduced by Greenlees and May in \cite[Sec.~5]{greenlees-may-completion}.
These objects are `global orthogonal spectra' in that they are indexed on pairs
$(G,V)$ consisting of a compact Lie group and a $G$-representation $V$.
The corresponding objects with commutative multiplication are called
{\em global $\mathcal I_*$-functors with 
smash products} in \cite[Sec.~5]{greenlees-may-completion} 
and it is for these that Greenlees and May define and study multiplicative
norm maps. Clearly, an orthogonal spectrum gives
rise to a  global $\mathcal I_*$-functors in the sense of Greenlees and May.
In the second chapter of her thesis \cite{bohmann-thesis},
Bohmann compares the approaches of Lewis-May and Greenlees-May;
in the paper \cite{bohmann-orthogonal} she also relates these to orthogonal spectra.

Symmetric spectra\index{subject}{symmetric spectrum}\index{subject}{spectrum!symmetric|see{symmetric spectrum}} 
in the sense of Hovey, Shipley and Smith \cite{HSS}
are another prominent model for the (non-equivariant) stable homotopy category.
Much of what we do here with orthogonal spectra can also be done 
with symmetric spectra, if one is willing to restrict to finite groups
(as opposed to general compact Lie groups). 
This restriction arises because only finite groups embed into symmetric groups,
while every compact Lie group embeds into an orthogonal group.
Hausmann \cite{hausmann-master thesis, hausmann-global} 
has established  a global model structure on the category of symmetric spectra, 
and he showed that the forgetful functor is a right Quillen equivalence
from the category of orthogonal spectra with the $\Fin$-global model structure
to the category of symmetric spectra with the global model structure.
While some parts of the symmetric and orthogonal theories are similar,
there are serious technical complications arising from
the fact that for symmetric spectra the naively defined equivariant homotopy groups 
are not `correct', a phenomenon that is already present non-equivariantly.

\smallskip

{\bf Prerequisites.} 
This book assumes a solid background in algebraic topology
and (non-equivariant) homotopy theory, including topics such as
singular homology and cohomology, CW-complexes, homotopy groups, 
mapping spaces, loop spaces, fibrations and fiber bundles,
Eilenberg-Mac Lane spaces, smooth manifolds, Grassmannian
and Stiefel manifolds. Two modern references that contain
all we need (and much more) are the textbooks by Hatcher \cite{hatcher}
and tom Dieck \cite{tomDieck-algebraic topology}.
Some knowledge of non-equivariant stable homotopy theory
is helpful to appreciate the equivariant and global features of the structures
and examples we discuss; from a strictly logical perspective, however,
the non-equivariant theory is a degenerate special case of the global theory for
the global family of trivial Lie groups.
In particular, by simply ignoring all group actions, 
the examples presented in this book give models for many interesting 
and prominent non-equivariant stable homotopy types.

Since actions of compact Lie groups are central to this book,
some familiarity with the structure and representation theory of
compact Lie groups is obviously helpful, but we give references
to the literature whenever we need any non-trivial facts.
Many of our objects of study organize themselves into model categories
in the sense of Quillen \cite{Q}, so some background on model
categories is necessary to understand the respective sections.
The article \cite{dwyer-spalinski} by Dwyer and Spalinski is a good introduction,
and Hovey's book \cite{hovey-book} is still the definitive reference.
Some acquaintance with unstable equivariant homotopy theory is useful
(but not logically necessary). In contrast, we do not assume any prior knowledge
of equivariant {\em stable} homotopy theory, and Chapter \ref{ch-equivariant}
is a self-contained introduction based on equivariant orthogonal spectra.
The last two sections of Chapter \ref{ch-stable} study the global stable
homotopy category, and here we freely use the language of triangulated categories.
The first chapter of Neeman's book \cite{neeman-triangulated categories}
is a possible reference for the necessary background.

Throughout the book we work in the category
of compactly generated spaces in the sense of McCord \cite{mccord},
i.e., a $k$-spaces (also called {\em Kelley spaces})
that satisfy the weak Hausdorff condition, see Definition \ref{def:cgwh}.
Since the various useful properties of
compactly generated spaces are scattered through the literature,
we include a detailed discussion in Appendix \ref{app:CGWH}.

{\bf Acknowledgments.} 
I would like to thank John Greenlees for being a reliable consultant
on matters of equivariant stable homotopy theory; Johannes Ebert and Michael Joachim
for tutorials on the $C^*$-algebraic approach to equivariant $K$-theory;
and Markus Hausmann and Irakli Patchkoria 
for careful reading of large parts of this book, along with
numerous suggestions for improvements. 
I am indebted to Benjamin B{\"o}hme, Jack Davies, Lars Hesselholt, Richard Hitzl,
Matthias Kreck, Jacob Lurie, Cary Malkiewich, Peter May,
Mona Merling, Luca Pol, Steffen Sagave, Albert Schulz, Neil Strickland, 
Karol Szumi{\l}o, Peter Teichner, Andreas Thom, Jorge Vit{\'o}ria,
Christian Wimmer and Mahmoud Zeinalian for various helpful discussions,
comments and suggestions.

\mainmatter

\chapter{Unstable global homotopy theory}\label{ch-unstable}

In this chapter we develop a framework for unstable global homotopy theory
via {\em orthogonal spaces}, i.e., continuous functors from the
linear isometries category $\bL$ to spaces. 
In Section \ref{sec:orthogonal spaces} we define
global equivalences of orthogonal spaces and establish many basic properties
of this class of morphisms. 
We also introduce global classifying spaces of compact Lie groups, 
the basic building blocks of global homotopy types.
In Section \ref{sec:global model structures spaces} we complement the
global equivalences by a global model structure 
on the category of orthogonal spaces.
The construction follows a familiar pattern, by Bousfield localization
of an auxiliary `strong level model structure'.
Section \ref{sec:global model structures spaces} also contains
a discussion of cofree orthogonal spaces, i.e., global homotopy types
that are `right induced' from non-equivariant homotopy types.
In Section \ref{sec:unstable monoidal} we recall the box product of orthogonal spaces, 
a Day convolution product based on orthogonal direct sum of inner product spaces.
The box product is a symmetric monoidal product,
fully homotopical under global equivalences, and globally equivalent to
the cartesian product.
Section \ref{sec:global families unstable}
introduces an important variation of our theme, where we discuss
unstable global homotopy theory for a `global family',
i.e., a class of compact Lie groups with certain closure properties.
In Section \ref{sec:equivariant homotopy sets} we introduce
the $G$-equivariant homotopy set $\pi_0^G(Y)$ of an orthogonal space
and identify the natural structure on these sets 
(restriction maps along continuous group homomorphisms).
The study of natural operations on the sets $\pi_0^G(Y)$ 
is a recurring theme throughout this book,
and we will revisit and extend the results in the later chapters
for ultra-commutative monoids, orthogonal spectra and ultra-commutative ring spectra.

Our main reason for working with orthogonal spaces
is that they are the direct unstable analog of orthogonal spectra,
and in this unstable model for global homotopy theory
the passage to the stable theory in Chapter \ref{ch-stable}
is especially simple. However, there are other models for 
unstable global homotopy theory, most notably the {\em topological stacks}
and the {\em orbispaces} as developed by Gepner and Henriques
in their paper \cite{gepner-henriques}. 
For a comparison of these two models to our orthogonal space model 
we refer to the author's paper \cite{schwede-orbispaces}.
The comparison proceeds through yet another model, the
global homotopy theory of `spaces with an action
of the universal compact Lie group'. 
Here the universal compact Lie group (which is neither compact nor a Lie group)
is the topological monoid $\Lc$ of linear isometric self-embeddings
of $\mR^\infty$, and in \cite{schwede-orbispaces} 
we establish a global model structure on the
category of $\Lc$-spaces.

\section{Orthogonal spaces and global equivalences}
\label{sec:orthogonal spaces}

In this section we introduce {\em orthogonal spaces}, 
along with the notion of {\em global equivalences},
our setup to rigorously formulate the idea of 
`compatible equivariant homotopy types for all compact Lie groups'.
We introduce various basic techniques to manipulate global equivalences
of orthogonal spaces, such as recognition criteria by 
homotopy or strict colimits over representations 
(Propositions \ref{prop:telescope criterion} and \ref{prop:global eq for closed}),
and a list of standard constructions that preserve global equivalences
(Proposition \ref{prop:global equiv basics}).  
Theorem \ref{thm:general shift of osp} is a cofinality result
for orthogonal spaces, showing that fairly general changes in the 
linear isometries indexing category do not affect the global homotopy type.
Definition \ref{def:global classifying} introduces
global classifying spaces of compact Lie groups,
the basic building blocks of global homotopy theory.
Proposition \ref{prop:classifying classifies} 
justifies the name by explaining in which sense the global classifying space $B_{\gl} G$
`globally classifies' principal $G$-bundles.

\medskip

Before we start, let us fix some notation and conventions.
By a `space' we mean a {\em compactly generated space} in the sense of \cite{mccord},
i.e., a $k$-space (also called {\em Kelley space})
that satisfies the weak Hausdorff condition,
see Definition \ref{def:cgwh}.
We denote the category of compactly generated spaces
by $\bT$ and review the basic properties in Appendix \ref{app:CGWH}.

An {\em inner product space}\index{subject}{inner product space}
is a finite-dimensional real vector space 
equipped with a scalar product, i.e., a positive definite symmetric bilinear form.
We denote by $\bL$\index{symbol}{$\bL$ - {category of finite-dimensional inner product
spaces}}  
the category with objects the 
inner product spaces and morphisms the linear isometric embeddings.
The category $\bL$ is a topological category in the sense
that the morphism spaces come with a preferred topology:
if $\varphi:V\to W$ is one linear isometric embedding,
then the action of the orthogonal group $O(W)$, by postcomposition,
induces a bijection\index{symbol}{$\bL(V,W)$ - {space of linear isometric embeddings}} 
\[ O(W)/O(\varphi^\perp) \ \iso \ \bL(V,W) \ ,\quad 
A \cdot O(\varphi^\perp) \ \longmapsto A\circ \varphi\ ,  \]
where $\varphi^\perp=W-\varphi(V)$ is the orthogonal complement of the image of $\varphi$.
We topologize $\bL(V,W)$ so that this bijection is a homeomorphism, and
this topology is independent of $\varphi$.
If $(v_1,\dots,v_k)$ is an orthonormal basis of $V$, then for every 
linear isometric embedding $\varphi:V\to W$ the tuple
$(\varphi(v_1),\dots,\varphi(v_k))$ is an orthonormal $k$-frame of $W$.
This assignment is a homeomorphism from $\bL(V,W)$ to the Stiefel manifold
of $k$-frames in $W$.

An example of an inner product spaces is the vector space $\mR^n$ 
with the standard scalar product
\[   \td{x,y} \ = \ x_1 y_1 \ + \ \dots \ +\   x_n y_n \ .  \]
In fact, every inner product space $V$ is isometrically isomorphic
to the inner product space $\mR^n$, for $n$ the dimension of $V$.
So the full topological subcategory with objects
the $\mR^n$ is a small skeleton of $\bL$.

\begin{defn}\label{def:orthogonal space} 
An {\em orthogonal space}\index{subject}{orthogonal space} 
is a continuous functor $Y:\bL\to\bT$ to the category of spaces.
A morphism of orthogonal spaces is a natural transformation.
We denote by $\spc$\index{symbol}{  $\spc$ - {category of orthogonal spaces}} 
the category of orthogonal spaces.\index{subject}{morphism!of orthogonal spaces} 
\end{defn}

The use of continuous functors from the category $\bL$ to spaces has a long history
in homotopy theory. 
The systematic use of inner product spaces (as opposed to numbers)
to index objects in stable homotopy theory seems to go back to Boardman's 
thesis \cite{boardman-thesis}.
The category $\bL$ (or its extension that
also contains countably infinite dimensional inner product spaces)
is denoted $\mathscr I$ by Boardman and Vogt \cite{boardman-vogt-homotopy everything},
and this notation is also used in \cite{may-quinn-ray};
other sources \cite{lind-diagram} use the symbol $\mathcal I$.
Accordingly, orthogonal spaces are sometimes referred to as $\mathscr I$-functors,
$\mathscr I$-spaces or $\mathcal I$-spaces.
Our justification for using yet another name is twofold: on the one hand,
our use of orthogonal spaces comes with a shift in emphasis, away from 
a focus on non-equivariant homotopy types, and towards viewing an
orthogonal space as representing compatible equivariant homotopy types for
all compact Lie groups. Secondly, we want to stress the analogy between
orthogonal spaces and orthogonal spectra, the former being an unstable global
world with the latter the corresponding stable global world.

Now we define our main new concept, the notion of `global equivalence'
between orthogonal spaces.
We let $G$ be a compact Lie group. By a {\em $G$-re\-presen\-tation}\index{subject}{G-representation@$G$-representation}
we mean a finite-dimensional orthogonal representation, 
i.e., a real inner product space
equipped with a continuous $G$-action by linear isometries.
In other words, a $G$-representation consists of an inner product space $V$
and a continuous homomorphism $\rho:G\to O(V)$.
In this context, and throughout the book,
we will often use without explicit mentioning that
continuous homomorphisms between Lie groups are automatically smooth,
see for example \cite[Prop.\,I.3.12]{broecker-tomDieck}.
For every orthogonal space $Y$ and every $G$-representation $V$,
the value $Y(V)$ inherits a $G$-action from the
$G$-action on $V$ and the functoriality of $Y$. 
For a $G$-equivariant linear isometric embedding $\varphi:V\to W$
the induced map $Y(\varphi):Y(V)\to Y(W)$ is $G$-equivariant.

We denote by
\[ D^k \ = \ \{x\in\mR^k \ : \ \td{x,x}\leq 1\} \text{\quad and\quad}
 \partial D^k \ = \ \{x\in\mR^k \ : \ \td{x,x}= 1\} \]
the unit disc in $\mR^k$ respectively
its boundary, a sphere of dimension $k-1$.
In particular, $D^0=\{0\}$ is a one-point space and $\partial D^0=\emptyset$
is empty.

\begin{defn}\label{def:global equivalence spaces}
A morphism $f:X\to Y$
of orthogonal spaces is a {\em global equivalence}\index{subject}{global equivalence!of orthogonal spaces} 
if the following condition holds: for every compact Lie group $G$,
every $G$-representation $V$, every $k\geq 0$ 
and all continuous maps $\alpha:\partial D^k\to X(V)^G$
and $\beta:D^k\to Y(V)^G$ such that $\beta|_{\partial D^k}=f(V)^G\circ\alpha$,
there is a $G$-representation $W$, a $G$-equivariant linear
isometric embedding $\varphi:V\to W$ and a continuous map $\lambda:D^k\to X(W)^G$
such that $\lambda|_{\partial D^k}=X(\varphi)^G\circ \alpha$ and
such that $f(W)^G\circ \lambda$
is homotopic, relative to $\partial D^k$, to $Y(\varphi)^G\circ \beta$.
\end{defn}

In other words, for every commutative square on the left
\[ \xymatrix@C=10mm{
\partial D^k\ar[r]^-\alpha \ar[d]_{\text{incl}} & X(V)^G \ar[d]^{f(V)^G} &
\partial D^k\ar[r]^-\alpha \ar[d]_{\text{incl}} & X(V)^G \ar[r]^-{X(\varphi)^G} & 
X(W)^G \ar[d]^{f(W)^G} 
\\
D^k\ar[r]_-\beta & Y(V)^G & 
D^k\ar[r]_-\beta \ar@{-->}[urr]^-(.4)\lambda & Y(V)^G \ar[r]_-{Y(\varphi)^G} & Y(W)^G  }\]
there exists the lift $\lambda$ on the right hand side that makes the upper left
triangle commute on the nose, and the lower right triangle up to homotopy
relative to $\partial D^k$. In such a situation we will often 
refer to the pair $(\alpha,\beta)$ as a `lifting problem' and we will say that
the pair $(\varphi,\lambda)$ {\em solves the lifting problem}.

\begin{eg} If $X=\underline{A}$ and $Y=\underline{B}$ 
are the constant orthogonal spaces with 
values the spaces $A$ respectively $B$, and $f=\underline{g}$ the constant 
morphism associated to a continuous map $g:A\to B$,
then $\underline{g}$ is a global equivalence if and only if
for every commutative square 
\[ \xymatrix@C=10mm{
\partial D^k\ar[r] \ar[d]_{\text{incl}} & A  \ar[d]^g \\
D^k\ar[r] \ar@{-->}[ur]^-(.4)\lambda & B }\]
there exists a lift $\lambda$ that makes the upper left
triangle commute, and the lower right triangle up to homotopy
relative to $\partial D^k$. But this is one of the equivalent ways of
characterizing weak equivalences of spaces,
compare \cite[Sec.\,9.6, Lemma]{may-concise}. 
So $\underline{g}$ is a global equivalence
if and only if $g$ is a weak equivalence.  
\end{eg}

\begin{rk}
The notion of global equivalence is meant to capture the idea that
for every compact Lie group $G$,
some induced morphism 
\[ \hocolim_V f(V)\ :\ \hocolim_V X(V)\ \to \ \hocolim_V Y(V) \]
is a $G$-weak equivalence, where `$\hocolim_V$' is a suitable homotopy
colimit over all $G$-representations $V$ along all 
equivariant linear isometric embeddings.
This is a useful way to think about global equivalences, 
and it could be made precise by letting $V$ run over the poset of
finite-dimensional subrepresentations of a complete $G$-universe and
using the Bousfield-Kan construction of a homotopy colimit over this poset.
Since the `poset of all $G$-representations' has a cofinal
subsequence, called an {\em exhaustive sequence} 
in Definition \ref{def:exhaustive sequence},
we can also model the `homotopy colimit over all $G$-representations'
as the mapping telescope over an exhaustive sequence.
However, the actual definition we work with has the advantage that
it does not refer to universes
and we do not have to define or manipulate homotopy colimits.

In many examples of interest, all the structure maps of an orthogonal space $Y$
are closed embeddings. When this is the case, the actual colimit
(over the subrepresentations of a complete universe) of the $G$-spaces $Y(V)$
serves the purpose of a `homotopy colimit over all representations',
and it can be used to detect global equivalences, 
compare Proposition \ref{prop:global eq for closed} below.
\end{rk}

We will now establish some useful criteria for detecting global equivalences.
We call a continuous map $f:A\to B$ 
an {\em h-cofibration}\index{subject}{h-cofibration} 
if it has the homotopy extension property, 
i.e., given a continuous map $\varphi:B\to X$ and a homotopy
$H:A\times [0,1]\to X$ starting with $\varphi f$,
there is a homotopy $\bar H:B\times [0,1]\to X$ starting with $\varphi$
such that $\bar H\circ(f\times [0,1])=H$.
Below we will write $H_t=H(-,t):A\to X$.
All h-cofibrations in the category of compactly generated spaces
are closed embeddings, compare Proposition \ref{prop:h-cof is closed embedding}.
The following somewhat technical lemma should be well known, but I was unable
to find a reference.

\begin{lemma} \label{lemma:equivalence of lifting problems}
Let $A$ be a subspace of a space $B$ such that the 
inclusion $A\to B$ is an h-cofibration.
Let $f:X\to Y$ be a continuous map and
\[ H\ : \ A\times [0,1]\ \to\ X \text{\qquad and\qquad}
K \ : \ B\times [0,1]\ \to\ Y\]
homotopies such that $K|_{A\times [0,1]}=f  H$.
Then the lifting problem $(H_0,K_0)$  has a solution if and only if
the lifting problem $(H_1,K_1)$ has a solution. 
\end{lemma}
\begin{proof}
The problem is symmetric, so we only show one direction.
We suppose that the lifting problem $(H_0,K_0)$  has a solution 
consisting of a continuous map $\lambda:B\to X$ 
such that $\lambda|_A=H_0$ and a homotopy $G:B\times[0,1]\to Y$
such that
\[ G_0 \ = \ f\circ\lambda \ , \quad G_1\ = \ K_0 \text{\qquad and\qquad}
(G_t)|_A \ = \ f\circ H_0  \]
for all $t\in[0,1]$. 
The homotopy extension property provides
a homotopy $H':B\times[0,1]\to X$ such that
\[ H'_0\ = \ \lambda \text{\qquad and\qquad} H'|_{A\times[0,1]} \ = \ H \ . \]
Then the map $\lambda'=H'_1:B\to X$ satisfies
\[ \lambda'|_A \ = \ ( H'_1)|_A \ = \ H_1\ . \]
We define a continuous map $J:B\times[0,3]\to Y$ by
\[ J_t \ = \
\begin{cases}
  f\circ H'_{1-t} & \text{ for $0\leq t\leq 1$,}\\
  \ G_{t-1} & \text{ for $1\leq t\leq 2$, and}\\
  \ K_{t-2} & \text{ for $2\leq t\leq 3$.}
\end{cases}\]
In particular,
\[ J_0 \ = \ f\circ \lambda'\text{\qquad and\qquad} 
J_3\ = \ K_1 \ ;\]
so $J$ almost witnesses the fact that $\lambda'$ solves the 
lifting problem $(H_1,K_1)$, except that $J$ is {\em not} a relative homotopy.

We improve $J$ to a relative homotopy from $f\circ\lambda'$ to $K_1$.
We define a continuous map $L:A\times[0,3]\times[0,1]\to Y$ by
\[ L(-,t,s) \ = \
\begin{cases}
f\circ H_{1-t} & \text{ for $0\leq t\leq s$, }\\
f\circ H_{1-s} & \text{ for $s\leq t\leq 3-s$, and}\\
f\circ H_{t-2} & \text{ for $3-s\leq t\leq 3$.}
\end{cases}\]
Then $L(-,-,0)$ is the constant homotopy at the map $f\circ H_1$, and  
\[ L(-,-,1) \ = \ J|_{A\times [0,3]}\ :\ A\times[0,3]\ \to\  Y\ .\]
Since the inclusion of $A$ into $B$ is an h-cofibration, 
the inclusion
of $B\times\{0\}\cup_{A\times\{0\}} A\times [0,1]$ into $B\times[0,1]$
has a continuous retraction; hence the inclusion
\[ B\times\{0\}\times [0,1] \cup_{A\times\{0\}\times [0,1]} A\times [0,1]\times [0,1]\ \to \ 
B\times[0,1]\times [0,1] \]
also has a continuous retraction.
We abbreviate $D=[0,3]\times \{1\}\cup \{0,3\}\times [0,1]$;
the pair of spaces $([0,3]\times [0,1],\, D)$
is pair homeomorphic to $([0,1]\times[0,1], \{0\}\times [0,1])$.
So the inclusion
\[ B\times D \cup_{A\times D} A\times [0,3]\times [0,1]\ \to \ B\times [0,3]\times [0,1] \]
has a continuous retraction.
The map $L$ and the map
\[  J \cup \text{const}_{f\lambda} \cup\text{const}_{K_1}\ : \ 
B\times D \ = \ B\times ( [0,3]\times \{1\}\cup  \{0,3\}\times [0,1])\ \to \ Y \]
agree on $A\times D$, so there is a continuous map
$\bar L:B\times[0,3]\times [0,1]\to Y$ such that
\[ \bar L(-,-,1)\ = \ J \ , \qquad
 \bar L|_{A\times[0,3]\times [0,1]} \ = \ L \ , \]
and
\[\bar L(-,0,s)\ = \ f\circ\lambda  \text{\qquad and\qquad}
\bar L(-,1,s)\ = \ K_1\] 
for all $s\in[0,1]$.
The map $\bar J=\bar L(-,-,0):B\times[0,3]\to Y$ then satisfies
\[ \bar J|_{A\times [0,3]}\ = \ \bar L(-,-,0)|_{A\times [0,3]}\ = \ 
L(-,-,0)\ , \]
which is the constant homotopy at the map $f\circ H_1$;
so $\bar J$ is a homotopy (parametrized by $[0,3]$ instead of $[0,1]$) 
{\em relative} to $A$.
Because
\[ \bar J_0 \ = \ \bar L(-,0,0)\ = \ f\circ\lambda \text{\qquad and\qquad}
 \bar J_3 \ = \ \bar L(-,3,0)\ = \ K_1\ , \]
the homotopy $\bar J$ witnesses that $\lambda'$ solves the lifting problem $(H_1,K_1)$.
\end{proof}

\begin{defn}\label{def:exhaustive sequence}
Let $G$ be a compact Lie group. 
An {\em exhaustive sequence}\index{subject}{exhaustive sequence!of representations}
is a nested sequence
\[ V_1 \ \subset \ V_2 \ \subset \ \dots \ \subset \ V_n \ \subset \ \dots \]
of finite-dimensional $G$-representations such that every finite-dimensional
$G$-representation admits a linear isometric $G$-embedding into some $V_n$.
\end{defn}

Given an exhaustive sequence $\{V_i\}_{i\geq 1}$ of $G$-representations 
and an orthogonal space $Y$, the values at the representations
and their inclusions form a sequence of $G$-spaces and $G$-equivariant continuous maps
\[ Y(V_1) \ \to \ Y(V_2) \ \to \ \cdots \ \to \ Y(V_i) \ \to \ \cdots \ .\]
We denote by
\[ \tel_i\, Y(V_i) \]
the mapping telescope of this sequence of $G$-spaces; this telescope inherits
a natural $G$-action.

We recall that a $G$-equivariant continuous map $f:A\to B$ between $G$-spaces
is a {\em $G$-weak equivalence}\index{subject}{G-weak equivalence@$G$-weak equivalence}
if for every closed subgroup $H$ of $G$
the map $f^H:A^H\to B^H$ of $H$-fixed points is a weak homotopy equivalence
(in the non-equivariant sense).

\begin{prop}\label{prop:telescope criterion} 
For every morphism of orthogonal spaces $f:X\to Y$, the following
three conditions are equivalent.
\begin{enumerate}[\em (i)]
\item The morphism $f$ is a global equivalence.
\item For every compact Lie group $G$,
every $G$-representation $V$, every finite $G$-CW-pair $(B,A)$
and all continuous $G$-maps $\alpha:A\to X(V)$
and $\beta:B\to Y(V)$ such that $\beta|_A=f(V)\circ\alpha$,
there is a $G$-representation $W$, a $G$-equivariant linear
isometric embedding $\varphi:V\to W$ and a continuous $G$-map $\lambda:B\to X(W)$
such that $\lambda|_A=X(\varphi)\circ \alpha$ and
such that $f(W)\circ \lambda$
is $G$-homotopic, relative to $A$, to $Y(\varphi)\circ \beta$.
\item For every compact Lie group $G$ and every exhaustive sequence $\{V_i\}_{i\geq 1}$
of $G$-representations the induced map
\[ \tel_i f(V_i)\ : \ \tel_i X(V_i)\ \to \ \tel_i Y(V_i) \]
is a $G$-weak equivalence.
\end{enumerate}
\end{prop}
\begin{proof}
At various places in the proof we use without explicit mentioning that
taking $G$-fixed points commutes with formation of the mapping telescopes;
this follows from the fact that taking $G$-fixed points
commutes with pushouts along closed embeddings and sequential colimits
along closed embeddings, compare Proposition \ref{prop:G-fix preserves pushouts}.

(i)$\Longrightarrow$(ii)  
We argue by induction over the number of the relative $G$-cells in $(B,A)$.
If $B=A$, then $\lambda=\alpha$ solves the lifting problem, and there is nothing to show.
Now we suppose that $A$ is a proper subcomplex of $B$. We choose a
$G$-CW-subcomplex $B'$ that contains $A$ and such that $(B,B')$ has
exactly one equivariant cell. Then $(B',A)$ has strictly fewer cells,
and the restricted equivariant lifting problem 
$(\alpha:A\to X(V),\beta'=\beta|_{B'}:B'\to Y(V))$ 
has a solution $(\varphi:V\to U, \lambda':B'\to X(U))$
by the inductive hypothesis.

We choose a characteristic map for the last cell, i.e., a pushout square of $G$-spaces
\[ \xymatrix{ 
G/H \times \partial D^k \ar[r]^-\chi \ar[d]_{\text{incl}} 
 & B' \ar[d]^{\text{incl}}  \\
G/H \times D^k \ar[r]_-\chi & B } \]
in which $H$ is a closed subgroup of $G$.
We arrive at the non-equivariant lifting problem on the left:
\[ \xymatrix@C=13mm{
\partial D^k\ar[r]^-{(\lambda')^H\circ\bar\chi} \ar[d]_{\text{incl}} & X(U)^H \ar[d]^{f(U)^H} &
\partial D^k\ar[r]^-{(\lambda')^H\circ\bar\chi}\ar[d]_{\text{incl}} & 
X(U)^H \ar[r]^-{X(\psi)^H} & 
X(W)^H \ar[d]^{f(W)^H} 
\\
D^k\ar[r]_-{Y(\varphi)^H\circ\beta^H\circ\bar\chi} & Y(U)^H & 
D^k\ar[r]_-{Y(\varphi)^H\circ\beta^H\circ\bar\chi} \ar@{-->}[urr]^-(.4)\lambda & 
Y(U)^H \ar[r]_-{Y(\psi)^H} & Y(W)^H  }\]
Here $\bar\chi=\chi(e H,-):D^k\to B^H$.
Since $f$ is a global equivalence, there is an $H$-equivariant
linear isometric embedding $\psi:U\to W$ and a continuous map 
$\lambda:D^k\to X(W)^H$ such that 
$\lambda|_{\partial D^k}=X(\psi)^H\circ(\lambda')^H\circ\bar\chi$
and $f(W)^H\circ\lambda$ is homotopic, relative $\partial D^k$, to
$Y(\psi)^H\circ Y(\varphi)^H\circ \beta^H\circ\bar\chi$, 
as illustrated by the diagram on the right above.
By enlarging $W$, if necessary, we can assume without loss of generality
that $W$ is underlying a $G$-representation and $\psi$ is even $G$-equivariant.

The $G$-equivariant extension of $\lambda$ 
\[ G/H\times D^k\to X(W) \ , \quad 
(g H,x)\ \longmapsto \ g\cdot \lambda(x) \]
and the map $X(\psi)\circ\lambda':B'\to X(W)$
then agree on $G/H\times\partial D^k$, so they glue to a $G$-map
$\tilde\lambda : B \to X(W)$.
The pair $(\psi\varphi:V\to W,\tilde\lambda:B\to X(W))$ then
solves the original lifting problem $(\alpha,\beta)$.

(ii)$\Longrightarrow$(iii)  
We suppose that $f$ satisfies~(ii), and we let $G$ be any compact Lie group
and $\{V_i\}_{i\geq 1}$ an exhaustive sequence of $G$-representations.
We consider an equivariant lifting problem, i.e., a finite $G$-CW-pair $(B,A)$ and
a commutative square:
\[ \xymatrix@C=15mm{
A\ar[r]^-\alpha \ar[d]_{\text{incl}} &  \tel_i  X(V_i)  \ar[d]^{\tel_i f(V_i)} \\
B\ar[r]_-\beta  & \tel_i  Y(V_i)  }\]
We show that every such lifting problem has an equivariant solution.
Since $B$ and $A$ are compact, there is an $n\geq 0$
such that $\alpha$ has image in the truncated telescope $\tel_{[0,n]} X(V_i)$
and $\beta$ has image in the truncated telescope $\tel_{[0,n]} Y(V_i)$
(see Proposition \ref{prop:filtered colim preserve weq}~(i)).
There is a natural equivariant homotopy from the identity of the truncated  
telescope $\tel_{[0,n]} X(V_i)$ to the composite
\[ \tel_{[0,n]} X(V_i) \ \xra{\text{proj}} \ X(V_n) \ \xra{\text{incl}}\
 \tel_{[0,n]} X(V_i) \ .\]
Naturality means that this homotopy is compatible with the 
same homotopy for the telescope of the $G$-spaces $Y(V_i)$.
Lemma \ref{lemma:equivalence of lifting problems}
(or rather its $G$-equivariant generalization)
applies to these homotopies, so instead of the original lifting problem
we may solve the homotopic lifting problem
\[ \xymatrix@C=15mm{
A\ar[r]^-{\alpha'} \ar[d]_{\text{incl}} &  X(V_n)\ar[r]^-{i_n}\ar[d]_{f(V_n)} &
\tel_i  X(V_i)  \ar[d]^{\tel_i f(V_i)} \\
B\ar[r]_-{\beta'}  & Y(V_n) \ar[r]_-{i_n} & \tel_i  Y(V_i)  }\]
where $\alpha'$ is the composite of the projection 
$\tel_{[0,n]} X(V_i) \to X(V_n)$
with $\alpha$, viewed as a map into the truncated telescope,
and similarly for $\beta'$.

Since $f$ satisfies~(ii), the lifting problem
$(\alpha':A\to X(V_n),\,\beta':B\to Y(V_n))$
has a solution after enlarging $V_n$ along some linear isometric $G$-em\-bedding.
Since the sequence $\{V_i\}_{i\geq 1}$ is exhaustive, we can take this embedding as
the inclusion $i:V_n\to V_m$ for some $m\geq n$, i.e., there is 
a continuous $G$-map $\lambda:B\to X(V_m)$ such that 
$\lambda|_A=X(i)^G\circ \alpha'$
and such that $f(V_m)^G\circ \lambda$ is $G$-homotopic, relative $A$,
to $Y(i)^G\circ \beta'$, compare the diagram:
\[ \xymatrix@C=15mm{
A\ar[r]^-{\alpha'} \ar[d]_{\text{incl}} &   
X(V_n)  \ar[r]^-{X(i)} & X(V_m)  \ar[d]^{f(V_m)} 
\\
B\ar[r]_-{\beta'} \ar@{-->}[urr]_-(.7)\lambda  & Y(V_n)  \ar[r]_-{Y(i)} & Y(V_m) }\]
The composite
\[ \xymatrix@C=10mm{
  X(V_n)\ar[r]^-{X(i)} & X(V_m) \ar[r]^-{i_m} & \tel_i  X(V_i)  }\]
does {\em not} agree with $i_n:X(V_n)\to \tel_i  X(V_i)$;
so the composite $i_n\circ \lambda:B\to \tel_i X(V_i)$ does not quite
solve the (modified) lifting problem $(i_n \circ\alpha', i_n\circ\beta')$.
But there is a $G$-equivariant homotopy $H:X(V_n)\times [0,1]\to \tel_i X(V_i)$
between $i_m\circ X(i)$ and $i_n$, and a similar homotopy 
$K:Y(V_n)\times [0,1]\to \tel_i Y(V_i)$ for $Y$ instead of $X$.
These homotopies satisfy
\[ K\circ (f(V_n)\times[0,1]) \ = \ (\tel_i f(V_i))\circ H \ ,\]
so Lemma \ref{lemma:equivalence of lifting problems} implies that the modified
lifting problem, and hence the original lifting problem,
has an equivariant solution.

(iii)$\Longrightarrow$(i)  
We let $G$ be a compact Lie group,
$V$ a $G$-representation, $k\geq 0$
and $(\alpha:\partial D^k\to X(V)^G, \beta:D^k\to Y(V)^G)$ 
a lifting problem, i.e., such that $\beta|_{\partial D^k}=f(V)^G\circ\alpha$.
We choose an exhaustive sequence $\{V_i\}$ of $G$-representations;
then we can embed $V$ into some $V_n$ by a linear isometric $G$-map 
and thereby assume without loss of generality that $V=V_n$.
 
We let $i_n:X(V_n)\to \tel_i X(V_i)$ and $i_n:Y(V_n)\to\tel_i Y(V_i)$
be the canonical maps. 
Since $\tel_i f(V_i):\tel_i X(V_i)\to \tel_i Y(V_i)$
is a $G$-weak equivalence, there is a continuous map
$\lambda:D^k\to (\tel_i X(V_i))^G$ such that $\lambda|_{\partial D^k}=i_n^G\circ \alpha$
and $(\tel_i f(V_i))^G\circ \lambda$ is homotopic, relative $\partial D^k$, to
$i_n^G\circ\beta$.
Since fixed points commute with mapping telescopes and since $D^k$ is compact,
there is an $m\geq n$ such that $\lambda$ and the relative homotopy
that witnesses the relation $(\tel_i f(V_i))^G\circ \lambda \simeq i_n^G\circ\beta$
both have image in
$\tel_{[0,m]} X(V_i)^G$, the truncated telescope of the $G$-fixed points.
The following diagram commutes
\[ \xymatrix@C=14mm{ 
X(V_n)^G \ar[r]_-{\text{can}} \ar[d]_{f(V_n)^G} \ar@/^1pc/[rrr]^(.3){X(\text{incl})^G}& 
\tel_{[0,n]} X(V_i)^G \ar[r]_-{\text{incl}} \ar[d]_{\tel f(V_i)^G} &
\tel_{[0,m]} X(V_i)^G \ar[r]_-{\text{proj}} \ar[d]^{\tel f(V_i)^G} &
X(V_m)^G \ar[d]^{f(V_m)^G} \\
Y(V_n)^G \ar[r]^-{\text{can}} \ar@/_1pc/[rrr]_(.3){Y(\text{incl})^G}& 
\tel_{[0,n]} Y(V_i)^G \ar[r]^-{\text{incl}} &
\tel_{[0,m]} Y(V_i)^G \ar[r]^-{\text{proj}} & Y(V_m)^G
}\]
where the right horizontal maps are the projections of 
the truncated telescope to the last term.
So projecting from $\tel_{[0,m]} X(V_i)^G$ to $X(V_m)^G$ 
and from $\tel_{[0,m]} Y(V_i)^G$ to $Y(V_m)^G$ 
produces the desired solution to the lifting problem. 
\end{proof}

We establish some more basic facts about the class of global equivalences.
A {\em homotopy} between two morphisms of orthogonal spaces $f,f':X\to Y$ is a morphism
\[ H \ : \ X\times [0,1] \ \to \ Y \]
such that $H(-,0)=f$ and $H(-,1)=f'$. 

\begin{defn}\label{def:strong level equivalence spaces}
  A morphism $f:X\to Y$ of orthogonal spaces is a {\em homotopy equivalence}
  if there is a morphism $g:Y\to X$ such that $g f$ and $f g$
  are homotopic to the respective identity morphisms.
  The morphism $f$ is a 
  {\em strong level equivalence}\index{subject}{strong level equivalence!of orthogonal spaces} 
  if for every compact Lie group $G$ and every $G$-representation $V$ the map
  $f(V)^G:X(V)^G\to Y(V)^G$ is a weak equivalence.
  The morphism $f$ is a {\em strong level fibration}\index{subject}{strong level fibration!of orthogonal spaces} 
  if for every compact Lie group $G$ and every $G$-representation $V$ the map
  $f(V)^G:X(V)^G\to Y(V)^G$ is a Serre fibration.
  \end{defn}

If $f,f':X\to Y$ are homotopic morphisms of orthogonal spaces, 
then the maps $f(V)^G, f'(V)^G:X(V)^G\to Y(V)^G$
are homotopic for every compact Lie group $G$ and every $G$-representation $V$.
So if $f$ is a homotopy equivalence of orthogonal spaces,
then the map $f(V)^G:X(V)^G\to Y(V)^G$
is a non-equivariant homotopy equivalence for every $G$-representation $V$.
So every homotopy equivalence is in particular a strong level equivalence.
By the following proposition, strong level equivalences are global equivalences.

A continuous map $\varphi:A\to B$ is a 
{\em closed embedding}\index{subject}{closed embedding}\index{subject}{embedding!closed|see{closed embedding}}
if it is injective and a closed map.
Such a map is then a homeomorphism of $A$ onto the closed subspace $\varphi(A)$ of $B$.
If a compact Lie group $G$ acts on two spaces $A$ and $B$ and
$\varphi:A\to B$ is a $G$-equivariant closed embedding,
then the restriction $\varphi^G:A^G\to B^G$ to $G$-fixed points 
is also a closed embedding.

We call a morphism $f:A\to B$ of orthogonal spaces
an {\em h-cofibration}\index{subject}{h-cofibration!of orthogonal spaces} 
if it has the homotopy extension property, 
i.e., given a morphism of orthogonal spaces $\varphi:B\to X$ and a homotopy
$H:A\times [0,1]\to X$ starting with $\varphi f$,
there is a homotopy $\bar H:B\times [0,1]\to X$ starting with $\varphi$
such that $\bar H\circ(f\times [0,1])=H$.

\begin{prop}\label{prop:global equiv basics} 
  \begin{enumerate}[\em (i)]
  \item 
    Every strong level equivalence is a global equivalence.
  \item The composite of two global equivalences is a global equivalence.
  \item 
    If $f, g$ and $h$ are composable morphisms of orthogonal spaces 
    such that $h g$ and $g f$ are global equivalences, then $f, g, h$ and $h g f$
    are also global equivalences.
  \item 
    Every retract of a global equivalence is a global equivalence.
  \item 
    A coproduct of any set of global equivalences is a global equivalence.
  \item 
    A finite product of global equivalences is a global equivalence.
  \item 
    For every space $K$ the functor $-\times K$ preserves global equivalences
    of orthogonal spaces.
  \item Let $e_n:X_n\to X_{n+1}$ and $f_n:Y_n\to Y_{n+1}$ be morphisms 
    of orthogonal spaces that are objectwise closed embeddings, for $n\geq 0$. 
    Let $\psi_n:X_n\to Y_n$ be global equivalences of orthogonal spaces
    that satisfy $\psi_{n+1}\circ e_n=f_n\circ\psi_n$ for all $n\geq 0$.
    Then the induced morphism $\psi_\infty:X_\infty\to Y_\infty$ 
    between the colimits of the sequences is a global equivalence.
  \item Let $f_n:Y_n\to Y_{n+1}$ be a global equivalence of orthogonal spaces
   that is objectwise a closed embedding, for $n\geq 0$. 
   Then the canonical morphism 
   $f_\infty:Y_0\to Y_\infty$ to the colimit of the sequence $\{f_n\}_{n\geq 0}$
   is a global equivalence.
 \item  Let 
   \[ \xymatrix{
     C  \ar[d]_\gamma & A \ar[l]_-g \ar[d]^\alpha \ar[r]^-f & B \ar[d]^\beta\\
     C'  & A' \ar[l]^-{g'}\ar[r]_-{f'} & B' } \]
   be a commutative diagram of orthogonal spaces such that $f$ and $f'$ are
   h-cofibrations.
   If the morphisms $\alpha,\beta$ and $\gamma$ are global equivalences,
   then so is the induced morphism of pushouts
    \[ \gamma\cup \beta\ : \ C\cup_A B \ \to \ C'\cup_{A'} B'\ . \]
  \item  Let 
    \[ \xymatrix{ A \ar[r]^-f \ar[d]_g & B \ar[d]^h\\
      C \ar[r]_-k & D } \]
    be a pushout square of orthogonal spaces such that $f$ is a global equivalence.
    If in addition $f$ or $g$ is an h-cofibration, 
    then the morphism $k$ is a global equivalence.
  \item Let 
    \[\xymatrix{ P\ar[r]^-k \ar[d]_g & X \ar[d]^f \\ Z \ar[r]_-h & Y }\]
    be a pullback square of orthogonal spaces in which $f$ is a global equivalence.
    If in addition $f$ or $h$ is a strong level fibration, 
    then the morphism $g$ is also a global equivalence.
 \end{enumerate}
\end{prop}
\begin{proof}
(i)  We let $f:X\to Y$ be a strong level equivalence, $G$ a compact Lie group,
$V$ a $G$-representation and $\alpha:\partial D^k\to X(V)^G$
and $\beta:D^k\to Y(V)^G$ continuous maps 
such that $f(V)^G\circ\alpha=\beta|_{\partial D^k}$.
Since $f$ is a strong level equivalence, the map $f(V)^G:X(V)^G\to Y(V)^G$
is a weak equivalence, so there is a continuous map $\lambda:D^k\to X(V)^G$ such that
$\lambda|_{\partial D^k}=\alpha$ and $f(V)^G\circ\lambda$ is homotopic to
$\beta$ relative $\partial D^k$. So the pair $(\Id_V,\lambda)$ solves the
lifting problem, and hence $f$ is a global equivalence.

(ii) We let $f:X\to Y$ and $g:Y\to Z$ be global equivalences, $G$ a compact Lie group,
$(B,A)$ a finite $G$-CW-pair, $V$ a $G$-representation and $\alpha:A\to X(V)$
and $\beta:B\to Z(V)$ continuous $G$-maps 
such that $(g f)(V)\circ\alpha=\beta|_A$.
Since $g$ is a global equivalence, the equivariant lifting problem 
$(f(V)\circ\alpha,\beta)$
has a solution $(\varphi:V\to W,\,\lambda:B\to Y(W))$ such that
\[ \lambda|_A\ =\ Y(\varphi)\circ f(V)\circ\alpha\ = \ 
f(W)\circ X(\varphi)\circ \alpha\ , \]
and $g(W)\circ\lambda$ is homotopic to $Z(\varphi)\circ\beta$ relative $A$. 
Since $f$ is a global equivalence,
the equivariant lifting problem $(X(\varphi)\circ\alpha,\lambda)$
has a solution $(\psi:W\to U,\, \lambda':B\to X(U))$ such that
\[ \lambda'|_A\ =\ X(\psi)\circ X(\varphi)\circ \alpha \]
and such that $f(U)\circ\lambda'$
is $G$-homotopic to $Y(\psi)\circ\lambda$ relative $A$.
Then $(g f)(U)\circ \lambda'$ is $G$-homotopic, relative $A$, to
\[ g(U)\circ Y(\psi)\circ\lambda  \ = \ Z(\psi)\circ g(W)\circ\lambda  \]
which in turn is $G$-homotopic to $Z(\psi\varphi)\circ\beta$,
also relative $A$. So the pair $(\psi\varphi,\lambda')$ solves
the original lifting problem for the morphism $g f:X\to Z$.

(iii) 
We let $f:X\to Y$, $g:Y\to Z$ and $h:Z\to Q$ be the three composable morphisms
such that $g f:X\to Z$ and $h g:Y\to Q$ are global equivalences.
We let $G$ be a compact Lie group and $\{V_i\}_{i\geq 1}$ an exhaustive sequence
of $G$-representations. Evaluating everything in sight on the
representations and forming mapping telescopes yields three composable
continuous $G$-maps
\[ 
\tel_i X(V_i)\ \xra{\tel_i f(V_i)} \
\tel_i Y(V_i)\ \xra{\tel_i g(V_i)} \
\tel_i Z(V_i)\ \xra{\tel_i h(V_i)} \
\tel_i Q(V_i)\ . \]
Proposition \ref{prop:telescope criterion} shows that the $G$-maps
\begin{align*}
   (\tel_i g(V_i))\circ(\tel_i f(V_i))\ &= \ \tel_i ( g f)(V_i)   \text{\quad and\quad}\\
(\tel_i h(V_i))\circ(\tel_i g(V_i))\ &= \ \tel_i (h g)(V_i) 
\end{align*}
are $G$-weak equivalences.
Since $G$-weak equivalences satisfy the 2-out-of-6-property,
we conclude that the $G$-maps
$\tel_i f(V_i)$, $\tel_i g(V_i)$, $\tel_i h(V_i)$ 
and
\[ (\tel_i h(V_i))\circ(\tel_i g(V_i))\circ(\tel_i f(V_i))\ = \ \tel_i (h g f)(V_i) \]
are $G$-weak equivalences.
Another application of Proposition \ref{prop:telescope criterion} then shows that 
$f$, $g$, $h$ and $h g f$ are global equivalences.

(iv) Let $g$ be a global equivalence and $f$ a retract of $g$. So there is a
commutative diagram
\[ \xymatrix{ 
X\ar[r]^-i \ar[d]_f & \bar X\ar[r]^-r\ar[d]^g & X\ar[d]^f \\
Y \ar[r]_-j & \bar Y \ar[r]_-s & Y} \]
such that $r i=\Id_X$ and $s j=\Id_Y$.
We let $G$ be a compact Lie group, $V$ a $G$-representation,
$(B,A)$ a finite $G$-CW-pair and $\alpha:A\to X(V)$
and $\beta:B\to Y(V)$ continuous $G$-maps 
such that $f(V)\circ\alpha=\beta|_A$.
Since $g$ is a global equivalence and
\[ g(V)\circ i(V)\circ \alpha\ = \ j(V)\circ f(V)\circ \alpha\ = \ 
 ( j(V)\circ \beta)|_A \ ,\]
there is a $G$-equivariant linear isometric embedding $\varphi:V\to W$ 
and a continuous $G$-map $\lambda:B\to \bar X(W)$ such that
$\lambda|_A=\bar X(\varphi)\circ i(V)\circ \alpha$ 
and $g(W)\circ\lambda$ is $G$-homotopic to 
$\bar Y(\varphi)\circ j(V)\circ\beta$ relative $A$. 
Then
\begin{align*}
 ( r(W)\circ\lambda )|_A\ &= \ 
 r(W)\circ \bar X(\varphi)\circ i(V)\circ \alpha\ 
= \  X(\varphi)\circ r(V)\circ i(V)\circ \alpha\ = \  
 X(\varphi)\circ \alpha   
\end{align*}
and
\[ f(W)\circ r(W)\circ\lambda \ = \  s(W)\circ g(W)\circ\lambda  \]
is $G$-homotopic to
\[  s(W)\circ  \bar Y(\varphi)\circ j(V)\circ\beta\ = \ 
Y(\varphi)\circ s(V)\circ j(V)\circ\beta\ = \ Y(\varphi)\circ\beta \]
relative $A$. 
So the pair $(\varphi, r(W)\circ\lambda)$ solves
the original lifting problem for the morphism $f:X\to Y$;
thus $f$ is a global equivalence.

Part~(v) holds because the disc $D^k$ is connected, so any lifting problem for
a coproduct of orthogonal spaces is located in one of the summands.

For part~(vi) it suffices to consider a product of two global equivalences $f:X\to Y$
and $f':X'\to Y'$.
Because global equivalences are closed under composition (part~(ii))
and $f\times f'=(f\times Y')\circ(X\times f')$, it suffices to show that for
every global equivalence $f:X\to Y$ and every orthogonal space $Z$ the morphism
$f\times Z:X\times Z\to Y\times Z$ is a global equivalence.
But this is straightforward:
we let $G$ be a compact Lie group, $V$ a $G$-representation,
$(B,A)$ a finite $G$-CW-pair and $\alpha:A\to (X\times Z)(V)$
and $\beta:B\to (Y\times Z)(V)$ continuous $G$-maps 
such that $(f\times Z)(V)\circ\alpha=\beta|_A$.
Because
\[ (X\times Z)(V) \ = \ X(V) \times Z(V)  \]
and similarly for $(Y\times Z)(V)$, we have
$\alpha=(\alpha_1,\alpha_2)$ and $\beta= (\beta_1,\beta_2)$
for continuous $G$-maps $\alpha_1:A\to X(V)$, $\alpha_2:A\to Z(V)$,
$\beta_1:B\to Y(V)$ and $\beta_2:B\to Z(V)$.
The relation $(f\times Z)(V)\circ(\alpha_1,\alpha_2)=(\beta_1,\beta_2)|_A$
shows that $\alpha_2=(\beta_2)|_A$.
Since $f$ is a global equivalence, 
the equivariant lifting problem $(\alpha_1,\beta_1)$
for $f(V)$ has a solution $(\varphi:V\to W,\,\lambda:B\to X(W))$ such that
$\lambda|_A=X(\varphi)\circ \alpha_1$ 
and $f(W)\circ\lambda$ is $G$-homotopic to 
$Y(\varphi)\circ\beta_1$ relative $A$. 
Then the pair $(\varphi,(\lambda,Z(\varphi)\circ\beta_2))$ solves
the original lifting problem, so $f\times Z$ is a global equivalence.

(vii) If $X$ is an orthogonal space and $K$ a space, then
$X\times K$ is the product of $X$ with the constant orthogonal space
with values $K$. So part~(vii) is a special case of~(vi).

(viii) We let $G$ be a compact Lie group,
$V$ a $G$-representation, $(B,A)$ a finite $G$-CW-pair
and $\alpha:A\to X_\infty(V)$
and $\beta:B\to Y_\infty(V)$ continuous $G$-maps 
such that $\psi_\infty(V)\circ\alpha=\beta|_A:A\to Y_\infty(V)$.
Since $A$ and $B$ are compact and $X_\infty(V)$ respectively $Y_\infty(V)$ 
are colimits of sequences of closed embeddings,
the maps $\alpha$ and $\beta$ factor through  maps
\[ \bar\alpha\ :\ A\ \to\ X_n(V)\text{\qquad respectively\qquad}
\bar\beta\ :\ B\ \to \ Y_n(V)\]
for some $n\geq 0$, see Proposition \ref{prop:filtered colim preserve weq}~(i).
Since the canonical maps $X_n(V)\to X_\infty(V)$ and $Y_n(V)\to Y_\infty(V)$ are injective,
$\bar\alpha$ and $\bar\beta$ are again $G$-equivariant.
Moreover, the relation $\psi_n(V)\circ\bar\alpha=\bar\beta|_A:A\to Y_n(V)$
holds because it holds after composition with the injective map $Y_n(V)\to Y_\infty(V)$. 

Since $\psi_n$ is a global equivalence, 
there is a $G$-equivariant linear isometric embedding $\varphi:V\to W$ 
and a continuous $G$-map  $\lambda:B\to X_n(W)$ such that
$\lambda|_A=X_n(\varphi)\circ \bar\alpha$ and 
$\psi_n(W)\circ\lambda$ is $G$-homotopic to
$Y_n(\varphi)\circ\bar\beta$ relative $A$. 
We let $\lambda':B\to X_\infty(W)$ be the composite of $\lambda$ and
the canonical map $X_n(W)\to X_\infty(W)$. 
Then the pair $(\varphi,\lambda')$ is a solution for the original lifting problem,
and hence $\psi_\infty:X_\infty\to Y_\infty$ is a global equivalence.

(ix) This is a special case of part~(viii) where we set $X_n=Y_0$,
$e_n=\Id_{Y_0}$ and $\psi_n=f_{n-1}\circ\dots\circ f_0:Y_0\to Y_n$.
The morphism $\psi_n$ is then a global equivalence by part~(ii), 
and $Y_0$ is a colimit of the constant first sequence.
Since the morphism $\psi_\infty$ induced on the colimits of the two sequences is
the canonical map $Y_0\to Y_\infty$, part~(viii) proves the claim.

(x)
Let $G$ be a compact Lie group. 
We consider an exhaustive sequence $\{V_i\}_{i\geq 1}$ 
of finite-dimensional $G$-representations.
Since $\alpha, \beta$ and $\gamma$ are global equivalences, the three
vertical maps in the following commutative diagram of $G$-spaces
are $G$-weak equivalences, by Proposition \ref{prop:telescope criterion}:
    \[ \xymatrix@C=15mm{
      \tel_i C(V_i)  \ar[d]_{\tel \gamma(V_i)} & 
      \tel_i A(V_i) \ar[l]_-{\tel g(V_i)} \ar[d]^{\tel \alpha(V_i)} \ar[r]^-{\tel f(V_i)} & 
      \tel_i B(V_i) \ar[d]^{\tel \beta(V_i)}\\
      \tel_i  C'(V_i)  & 
      \tel_i A'(V_i) \ar[l]^-{\tel g'(V_i)}\ar[r]_-{\tel f'(V_i)} & 
      \tel_i B'(V_i) } \]
Since mapping telescopes commute with product with $[0,1]$ and retracts,
the maps $\tel_i f(V_i)$ and $\tel_i f'(V_i)$ are h-cofibrations of $G$-spaces.
The induced map of the horizontal pushouts is thus a $G$-weak equivalence
by Proposition \ref{prop:gluing lemma G-spaces}.
Since formation of mapping telescopes commutes with pushouts, 
the map
    \[ \tel_i (\gamma\cup \beta)(V_i)\ : \ \tel_i (C\cup_A B)(V_i) \ \to \ 
    \tel_i (C'\cup_{A'} B')(V_i) \]
    is a $G$-weak equivalence. The claim thus
follows by another application of the telescope criterion for 
global equivalences, Proposition \ref{prop:telescope criterion}.

(xi)
In the commutative diagram
    \[ \xymatrix{
      C  \ar@{=}[d] & A \ar[l]_-g \ar@{=}[d] \ar@{=}[r] & A \ar[d]^f\\
      C  & A \ar[l]^-{g}\ar[r]_-{f} & B } \]
all vertical morphisms are global equivalences.
If $f$ is an h-cofibration, we apply part~(x) to this diagram
to get the desired conclusion.
If $g$ is an h-cofibration, we apply part~(x) after interchanging the roles
of left and right horizontal morphisms.

(xii) We let $G$ be a compact Lie group,
$V$ a $G$-representation, $(B,A)$ a finite $G$-CW-pair
and $\alpha:A\to P(V)$ and $\beta:B\to Z(V)$ continuous $G$-maps 
such that $g(V)\circ\alpha=\beta|_A$.
Since $f$ is a global equivalence, there is 
a $G$-equivariant linear isometric embedding $\varphi:V\to W$ 
and a continuous $G$-map $\lambda:B\to X(W)$
such that 
$\lambda|_A=X(\varphi)\circ k(V)\circ \alpha$ and
such that $f(W)\circ \lambda$
is $G$-homotopic, relative to $A$, to $Y(\varphi)\circ h(V)\circ \beta$.  
We let $H:B\times[0,1]\to Y(W)$ be a relative $G$-homotopy from 
$Y(\varphi)\circ h(V)\circ \beta =  h(W)\circ Z(\varphi)\circ \beta$
to $f(W)\circ \lambda$. 
Now we distinguish two cases.

Case~1: The morphism $h$ is a strong level fibration.
We can choose a lift $\bar H$ in the square
\[ \xymatrix@C25mm{ 
B\times 0\cup_{A\times 0} A\times [0,1]\ar[d]_\sim 
\ar[r]^-{(Z(\varphi)\circ\beta)\cup K} &
Z(W) \ar@{->>}[d]^{h(W)} \\
B\times [0,1] \ar[r]_-H \ar@{-->}[ur]_(.6){\bar H} & Y(W)
} \]
where $K:A\times [0,1]\to Z(W)$ is the constant homotopy
from $g(W)\circ P(\varphi)\circ\alpha$ to itself.
Since the square is a pullback and
$h(W)\circ \bar H(-,1)= H(-,1)=f(W)\circ\lambda$,
there is a unique continuous $G$-map $\bar\lambda : B\to P(W)$
that satisfies
\[ g(W) \circ \bar \lambda  \ = \ \bar H(-,1)
\text{\qquad and\qquad}  
k(W) \circ \bar\lambda \ = \ \lambda \ .\]
The restriction of $\bar\lambda$ to $A$ satisfies
\begin{align*}
   g(W) \circ \bar \lambda|_A \ &= \ \bar H(-,1)|_A \ = \ 
g(W)\circ P(\varphi)\circ\alpha
\text{\qquad and}\\ 
k(W) \circ \bar\lambda|_A  \ &= \ \lambda|_A \ = \ 
X(\varphi)\circ k(V)\circ \alpha \ = \ 
k(W)\circ P(\varphi)\circ \alpha \ .\end{align*}
The pullback property thus implies that
$\bar\lambda|_A =P(\varphi)\circ \alpha$. 

Finally, the composite $g(W) \circ \bar \lambda$ 
is homotopic, relative $A$ and via $\bar H$,
to $\bar H(-,0)= Z(\varphi)\circ\beta$. 
This is the required lifting data, and
we have verified the defining property of a global equivalence for the
morphism $g$.

Case~2: The morphism $f$ is a strong level fibration. 
The argument is similar as in the first case.
Now we can choose a lift $H'$ in the square
\[ \xymatrix@C30mm{ 
B\times 1\cup_{A\times 1} A\times [0,1]\ar[d]_\sim 
\ar[r]^-{\lambda\cup K'} &
X(W) \ar@{->>}[d]^{f(W)} \\
B\times [0,1] \ar[r]_-H \ar@{-->}[ur]_{H'} & Y(W)
} \]
where $K':A\times [0,1]\to X(W)$ is the constant homotopy
from $X(\varphi)\circ k(V)\circ\alpha$ to itself.
Since the square is a pullback and
$f(W)\circ H'(-,0)= H(-,0)=h(W)\circ Z(\varphi)\circ\beta$,
there is a unique continuous $G$-map $\bar\lambda : B\to P(W)$ that satisfies
\[ g(W) \circ \bar \lambda  \ = \ Z(\varphi)\circ \beta
\text{\qquad and\qquad}  
k(W) \circ \bar\lambda \ = \ H'(-,0) \ .\]
The restriction of $\bar\lambda$ to $A$ satisfies
\begin{align*}
   g(W) \circ \bar \lambda|_A \ &= \ 
Z(\varphi)\circ g(V)\circ\alpha \ = \ 
g(W)\circ P(\varphi)\circ\alpha
\text{\qquad and}\\ 
k(W) \circ \bar\lambda|_A  \ &= \ H'(-,0)|_A \ = \ 
X(\varphi)\circ k(V)\circ\alpha\ = \ 
k(W)\circ P(\varphi)\circ \alpha \ .\end{align*}
The pullback property thus implies that
$\bar\lambda|_A =P(\varphi)\circ \alpha$. 
Since $g(W) \circ \bar \lambda= Z(\varphi)\circ\beta$,
this is the required lifting data, and
we have verified the global equivalence criterion 
of Proposition \ref{prop:telescope criterion}~(ii) for the morphism $g$.
\end{proof}

\Danger The restriction to {\em finite} products is essential in part~(vi)
of the previous Proposition \ref{prop:global equiv basics}, i.e.,
an infinite product of global equivalences need not be a global equivalence.
The following simple example illustrates this. 
We let $Y_n$ denote the orthogonal space with
\[ Y_n(V) \ = \
\begin{cases}
\ \emptyset & \text{ if $\dim(V)< n$, and}\\
 \{\ast\} & \text{ if $\dim(V)\geq n$.}
\end{cases}\]
The unique morphism $Y_n\to Y_0$ 
is a global equivalence for every $n\geq 0$.
However, the product $\prod_{n\geq 0 } Y_n$ is
the empty orthogonal space, whereas the product
$\prod_{n\geq 0} Y_0 $ is a terminal orthogonal space.
The unique morphism from the initial (i.e., empty)
to a terminal orthogonal space is {\em not} a global equivalence.

\medskip

The following proposition provides a lot of flexibility
for changing an orthogonal space into a globally equivalent one by modifying
the input variable. We will use it multiple times in this book.

\begin{theorem}\label{thm:general shift of osp}
  Let $F:\bL\to\bL$ be a continuous endofunctor of the linear isometries
  category and $i:\Id\to F$ a natural transformation.
  Then for every orthogonal space $Y$ the morphism
  \[ Y\circ i\ : \ Y\ \to \ Y\circ F \]
  is a global equivalence of orthogonal spaces.
\end{theorem}
\begin{proof}
In a first step we show an auxiliary statement. We let $V$ be an inner product
space and $z\in F(V)$ an element that is orthogonal to the subspace $i_V(V)$,
the image of the linear isometric embedding $i_V:V\to F(V)$.
We claim that for every linear isometric embedding $\varphi:V\to W$
the element $F(\varphi)(z)$ of $F(W)$ is orthogonal to the subspace $i_W(W)$.
To prove the claim we write any given element of $W$ as $\varphi(v)+ y$ 
for some $v\in V$ and $y\in W$ orthogonal to $\varphi(V)$. Then 
\[ \td{ F(\varphi)(z), i_W(\varphi(v)) } \ = \ 
\td{ F(\varphi)(z), F(\varphi)(i_V(v)) } \ = \ \td{ z,  i_V(v) } \ = \ 0 \]
by the hypotheses on $z$. 
Now we define $A\in O(W)$ as the linear isometry that is the identity on $\varphi(V)$
and the negative of the identity on the orthogonal complement of $\varphi(V)$.
Then $A\circ\varphi=\varphi$ and
\begin{align*}
  \td{ F(\varphi)(z), i_W(y) } \ &= \ 
\td{ F(A)( F(\varphi)(z)) , F(A)(i_W(y)) } \\ 
&= \ \td{ F(A \varphi)(z) , i_W(A(y)) } \ = \ - \td{ F(\varphi)(z) , i_W(y) }  \ , 
\end{align*}
and hence $\td{ F(\varphi)(z), i_W(y)} =  0$.
Altogether this shows that $\td{ F(\varphi)(z), i_W(\varphi(v)+y)} =  0$,
which establishes the claim.

Now we consider a compact Lie group $G$, a $G$-representation $V$,
a finite $G$-CW-pair $(B,A)$ and a lifting problem $\alpha:A\to Y(V)$ 
and $\beta:B\to Y(F(V))$ for $(Y\circ i)(V)$. 
Then $\beta|_A=Y(i_V)\circ\alpha$ by hypothesis, 
and we claim that $Y(i_{F(V)})\circ\beta$
is $G$-homotopic to $Y(F(i_V))\circ\beta = (Y\circ F)(i_V)\circ\beta$,
relative $A$; granting this for the moment, we conclude that
the pair $(i_V:V\to F(V), \beta)$ solves the lifting problem.

It remains to construct the relative homotopy.
The two embeddings
  \[  F(i_V) \ , \ i_{ F(V)} \ : \ F(V)\ \to \  F( F(V))  \]
are homotopic, relative to $i_V:V\to F(V)$, 
through $G$-equivariant isometric embeddings, via
\begin{align*}
 H\  : \  F(V) \times [0,1]\ &\to \qquad  F(F(V)) \\
(v + w, t)\qquad &\longmapsto  F(i_V)(v) \ + \ t\cdot i_{F(V)}(w)\ + \ 
\sqrt{1-t^2}\cdot  F(i_V)(w) \ ,  
\end{align*}
where $v\in i_V(V)$ and $w$ is orthogonal to $i_V(V)$.
The verification that $H(-,t):F(V)\to F(F(V))$ is indeed
a linear isometric embedding for every $t\in[0,1]$ uses that
$i_{F(V)}=F(i_V)$ on the subspace  $i_V(V)$ of $F(V)$, and 
that $i_{F(V)}(w)$ is orthogonal to $F(i_V)(w)$, by the claim proved above.
The continuous functor $Y$ takes this homotopy of equivariant linear isometric
embeddings to a $G$-equivariant homotopy $Y(H(-,t))$
from $Y(F(i_V))$ to $Y(i_{F(V)})$,  relative to $Y(i_V)$.
Composing with $\beta$ gives the required relative $G$-homotopy
from $Y(F(i_V))\circ\beta$ to $Y(i_{F(V)})\circ\beta$.
\end{proof}

\begin{eg}[Additive and multiplicative shift]\label{eg:Additive and multiplicative shift}
Here are some typical examples where the previous theorem applies. 
Every inner product space $W$ defines an `additive shift functor'\index{subject}{shift!of an orthogonal space!additive}
and a `multiplicative shift functor'\index{subject}{shift!of an orthogonal space!multiplicative}
on the category of orthogonal spaces,
defined by precomposition with the continuous endofunctors
\[ -\oplus W \ : \ \bL \ \to \ \bL \text{\qquad respectively\qquad}
 -\tensor W \ : \ \bL \ \to \ \bL \ . \]
In other words, the additive respectively multiplicative $W$-shift 
of an orthogonal space $Y$ have values
\[ (\sh^W_\oplus\! Y)(V)\ = \ Y(V\oplus W)\text{\qquad respectively\qquad}
 (\sh^W_\tensor\! Y)(V)\ = \ Y(V\tensor W)\ . \]
Here, and in the rest of the book, we endow
the tensor product $V\tensor W$ of two inner product spaces $V$ and $W$ 
with the inner product characterized by
\[ \td{v\tensor w, \bar v\tensor \bar w}\ = \  \td{v, \bar v}\cdot  \td{w, \bar w} \]
for all $v,\bar v\in V$ and $w,\bar w\in W$.
Another way to say this is that for every orthonormal basis $\{b_i\}_{i\in I}$
of $V$ and every orthonormal basis $\{d_j\}_{j\in J}$ of $W$ the family
$\{b_i\tensor d_j\}_{(i,j)\in I\times J}$ forms an orthonormal basis of $V\tensor W$.
Theorem \ref{thm:general shift of osp}  then shows that the morphism
$Y\to \sh^W_\oplus\! Y$ given by applying $Y$ 
to the first summand embedding $V\to V\oplus W$
is a global equivalence.
To get a similar statement for the multiplicative shift we have to assume
that $W\ne 0$; then for every vector $w\in W$ of length~1 the map
\[ V \ \to \ V\tensor W \ ,\quad v\ \longmapsto \ v\tensor w\]
is a natural linear isometric embedding. 
So Theorem \ref{thm:general shift of osp} shows that the morphism
$Y(-\tensor w):Y\to \sh^W_\tensor\! Y$ is a global equivalence.
\end{eg}

For the following discussion of universes we recall that 
a finite group has finitely many isomorphism classes of
irreducible orthogonal representations, and a compact Lie group of
positive dimension has countably infinitely many such isomorphism classes.
I have not yet found an explicit reference with a proof of this well-known fact,
but one can argue as follows. 

We first consider {\em unitary} representations of $G$.
If $G$ is finite, then the characters of irreducible unitary representations form 
a basis of the $\mC$-vector space of conjugation invariant $\mC$-valued functions
on $G$. So the number of isomorphism classes of irreducible
unitary representations agrees with the number of conjugacy classes of elements
of $G$, and is thus finite.
If $G$ is of positive dimension, 
then the characters of irreducible unitary representations form 
an orthonormal basis of the complex Hilbert space of square integrable 
(with respect to the Haar measure), conjugation invariant functions on $G$, 
see for example \cite[\S 11, Thm.\,2]{kirillov}.
Since $G$ is compact and of positive dimension, 
this Hilbert space is infinite dimensional and separable, so there are
countably infinitely many isomorphism classes of irreducible unitary representations.

To treat the case of orthogonal representations of $G$,
we recall from \cite[II Prop.\,6.9]{broecker-tomDieck}
that complexification can be used to construct a map from
the set of isomorphism classes of irreducible orthogonal representations
to the set of isomorphism classes of irreducible unitary representations of $G$.
There is a caveat, however, namely that the complexification of
an irreducible orthogonal representation need not be irreducible.
More precisely, the reducibility behavior under complexification
depends on the `type' of the irreducible orthogonal representation $\lambda$.
By Schur's lemma, the endomorphism algebra $\Hom_\mR^G(\lambda,\lambda)$ 
is a finite-dimensional skew-field extension of $\mR$, hence isomorphic
to $\mR$, $\mC$ or $\mH$. 
\begin{itemize}
\item If $\Hom_\mR^G(\lambda,\lambda)$ is isomorphic to $\mR$, 
then $\lambda$ is of {\em real type}.
In this case the complexification $\lambda_\mC$
is irreducible as a unitary $G$-representation.
\item If $\Hom_\mR^G(\lambda,\lambda)$ is isomorphic to $\mC$ or $\mH$, 
then $\lambda$ is of {\em complex type} respectively of {\em quaternionic type}.
In this case there is an irreducible unitary $G$-representation $\rho$
such that $\lambda_\mC$ is isomorphic to the direct sum of $\rho$ and its
conjugate $\bar\rho$.
If $\lambda$ is of complex type, then $\rho$ is not isomorphic to its conjugate;
if $\lambda$ is of quaternionic type, then $\rho$ is self-conjugate,
i.e., isomorphic to its conjugate.
\end{itemize}
Since the underlying orthogonal representation 
of $\lambda_\mC$ is isomorphic to the direct sum of two copies of $\lambda$,
two non-isomorphic irreducible orthogonal representations
cannot become isomorphic after complexification. 
So the above construction gives an injective map from 
the set of irreducible orthogonal representations
to the set of irreducible unitary representations.
Altogether this shows that there at most countably many
isomorphism classes of irreducible orthogonal representations
of a compact Lie group.

\begin{defn} 
Let $G$ be a compact Lie group. 
A {\em $G$-universe}\index{subject}{G-universe@$G$-universe}\index{subject}{universe|see{$G$-universe}}
is an orthogonal $G$-representation $\Uc$ 
of countably infinite dimension with the following two properties:
\begin{itemize}
\item the representation $\Uc$ has non-zero $G$-fixed points,
\item if a finite-dimensional $G$-representation $V$ embeds into $\Uc$,
then a countable infinite direct sum
of copies of $V$ also embeds into $\Uc$.
\end{itemize}
A $G$-universe is {\em complete}\index{subject}{G-universe@$G$-universe!complete} 
if every finite-dimensional $G$-representation embeds into it.
\end{defn}

A $G$-universe is characterized, up to equivariant linear isometry,
by the set of irreducible $G$-represen\-ta\-tions that embed into it.
We let $\Lambda=\{\lambda\}$ be a complete set of 
pairwise non-isomorphic irreducible $G$-representations that embed into $\Uc$.
The first condition says that $\Lambda$ contains
a trivial 1-dimensional representation, and the second condition is equivalent
to the requirement that
\[ \Uc \ \iso \ \bigoplus_{\lambda\in\Lambda} \bigoplus_\mN\lambda \ .\]
Moreover, $\Uc$ is complete if and only if $\Lambda$ contains
(representatives of) {\em all} irreducible $G$-representations.
Since there are only countably many isomorphism classes of irreducible orthogonal 
$G$-representations, a complete $G$-universe exists.

\begin{rk}\label{rk:complete universes restrict} 
We let $H$ be a closed subgroup of a compact Lie group $G$.
We will frequently use that the underlying $H$-representation of a complete
$G$-universe $\Uc$ is a complete $H$-universe.
Indeed, if $U$ is an $H$-representation, 
then there is a $G$-representation $V$ 
and an $H$-equivariant linear isometric embedding $U\to V$,
see for example \cite[Prop.\,1.4.2]{palais-classification} 
or \cite[III Thm.\,4.5]{broecker-tomDieck}.
Since $V$ embeds $G$-equivariantly into $\Uc$,
the original representation $U$ embeds $H$-equivariantly into $\Uc$.
\end{rk}

In the following we fix, for every compact Lie group $G$,
a complete $G$-universe $\Uc_G$.\index{symbol}{$\Uc_G$ - {complete $G$-universe}}
We let $s(\Uc_G)$\index{symbol}{$s(\Uc_G)$ - {poset of $G$-subrepresentation of $\Uc_G$}}
denote the poset, under inclusion,
of finite-dimensio\-nal $G$-subrepresentations of $\Uc_G$.

\begin{defn}
For an orthogonal space $Y$ and a compact Lie group $G$ we define the 
{\em underlying $G$-space}\index{subject}{underlying $G$-space!of an orthogonal space} as
\[ Y(\Uc_G)\ = \ \colim_{V\in s(\Uc_G)}\, Y(V)\ , \]
the colimit of the $G$-spaces $Y(V)$. 
\end{defn}

\begin{rk}
The underlying $G$-space $Y(\Uc_G)$ can always be written as a sequential colimit
of values of $Y$. Indeed, we can choose a nested sequence of
finite-dimensional $G$-subrepresentations
\[ V_0 \ \subseteq V_1\ \subseteq\ V_2\ \subseteq\ \cdots \]
whose union is all of $\Uc_G$.
This is then in particular an exhaustive sequence in the sense of
Definition \ref{def:exhaustive sequence}.
Since the subposet
$\{V_n\}_{n\geq 0}$ is cofinal in $s(\Uc_G)$, the colimit of the functor
$V\mapsto Y(V)$ over $s(\Uc_G)$ is also a colimit over the 
subsequence $Y(V_n)$.

If the group $G$ is finite, then we can define a complete universe as 
\[ \Uc_G \ = \  \bigoplus_\mN\, \rho_G\ , \]
a countably infinite sum of copies of the regular representation 
$\rho_G=\mR[G]$, with $G$ as orthonormal basis.
Then $\Uc_G$ is filtered by the finite sums $n\cdot\rho_G$, and we get
\[ Y(\Uc_G)\ = \ \colim_n \, Y(n\cdot \rho_G) \ , \]
where the colimit is taken along the inclusions 
$n\cdot \rho_G\to(n+1)\cdot\rho_G$ that miss the last summand.
\end{rk}

\begin{defn}
An orthogonal space $Y$ is {\em closed}\index{subject}{orthogonal space!closed} 
if it takes every linear isometric embedding $\varphi:V\to W$ of inner product spaces
to a closed embedding $Y(\varphi):Y(V)\to Y(W)$.
\end{defn}

In particular, for every closed orthogonal space $Y$ 
and every $G$-equivariant linear isometric embedding $\varphi:V\to W$ of
$G$-representations, the induced map on $G$-fixed points
$Y(\varphi)^G:Y(V)^G\to Y(W)^G$ is also a closed embedding.

\begin{prop}\label{prop:global eq for closed} 
Let $f:X\to Y$ be a morphism between closed orthogonal spaces.
Then $f$ is a global equivalence if and only if for every compact Lie
group $G$ the map
\[ f(\Uc_G)^G\ : \ X(\Uc_G)^G \ \to \ Y(\Uc_G)^G\]
is a weak equivalence.
\end{prop}
\begin{proof}
The poset $s(\Uc_G)$ has a cofinal subsequence, so all colimits
over $s(\Uc_G)$ can be realized as sequential colimits.
The claim is then a straightforward consequence of the fact that
fixed points commute with sequential colimits along closed embeddings
(see Proposition \ref{prop:G-fix preserves pushouts}~(ii))
and continuous maps from compact spaces such as $D^k$ and $\partial D^k$ 
to sequential colimits along closed embeddings factor through a finite stage 
(see Proposition \ref{prop:filtered colim preserve weq}~(i)).
\end{proof}

Now we turn to semifree orthogonal spaces.
Important special cases of this construction are the global classifying spaces
of compact Lie groups (Definition \ref{def:global classifying}),
the basic building blocks of global homotopy theory.
Free and semifree orthogonal spaces are made from
spaces of linear isometric embeddings, so we start by recalling
various properties of certain spaces of linear isometric embeddings.
We consider two compact Lie groups $G$ and $K$,
a finite-dimensional $G$-representation $V$, 
and a $K$-representation $\Uc$, possibly of countably infinite dimension. 
If $\Uc$ is infinite dimensional, we topologize the space $\bL(V,\Uc)$ 
of linear isometric embeddings as the filtered colimit of the 
spaces $\bL(V,U)$, taken over the poset of finite-dimensional subspaces $U$ of $\Uc$.
The space $\bL(V,\Uc)$ 
inherits a continuous left $K$-action and a compatible continuous right $G$-action
from the actions on the target and source, respectively.
We turn these two actions into a single left action of the group $K\times G$
by defining
\begin{equation}\label{eq:turn_into_left}
 ((k,g)\cdot \varphi)(v) \ = \ k\cdot \varphi(g^{-1} \cdot v)  
\end{equation}
for $\varphi\in\bL(V,\Uc)$ and $(k,g)\in K\times G$.
We recall that a continuous $(K\times G)$-equivariant map is a
{\em $(K\times G)$-cofibration} if it has the right lifting
property with respect to all morphisms of $(K\times G)$-spaces
$f:X\to Y$ such that the map $f^\Gamma:X^\Gamma\to Y^\Gamma$ is a weak equivalence and
Serre fibration for every closed subgroup $\Gamma$ of $K\times G$.

\begin{prop}\label{prop:K G cofibration} 
Let $G$ and $K$ be compact Lie groups,
$V$ a finite-dimensional $G$-representation, and $\Uc$ a $K$-representation 
of finite or countably infinite dimension.
\begin{enumerate}[\em (i)]
\item 
For every finite-dimensional $K$-subrepresentation $U$ of $\Uc$, 
the inclusion induces a $(K\times G)$-cofibration
\[ \bL(V,U)\ \to \ \bL(V,\Uc) \]
and a $K$-cofibration of orbit spaces
\[ \bL(V,U)/G\ \to \ \bL(V,\Uc)/G \ .\]
\item The $(K\times G)$-space $\bL(V,\Uc)$ is $(K\times G)$-cofibrant.
The $K$-space $\bL(V,\Uc)/G$ is $K$-cofibrant.
\end{enumerate}
\end{prop}
\begin{proof}
(i) 
We consider two natural numbers $m,n\geq 0$. 
The space $\bL(V,\mR^{m+n})$ is homeomorphic to the Stiefel manifold of   
$\dim(V)$-frames in $\mR^{m+n}$, hence a compact smooth manifold, and the
action of $O(m)\times O(n)\times G$ is smooth.
Illman's theorem \cite[Cor.\,7.2]{illman} thus provides an 
$(O(m)\times O(n)\times G)$-CW-structure on $\bL(V,\mR^{m+n})$.
In particular, $\bL(V,\mR^{m+n})$ is cofibrant as an 
$(O(m)\times O(n)\times G)$-space.
The group $N= e\times O(n)\times e$ is a closed normal subgroup of 
$O(m)\times O(n)\times G$, so the inclusion of the $N$-fixed points
into $\bL(V,\mR^{m+n})$ is an $(O(m)\times O(n)\times G)$-cofibration
(compare Proposition \ref{prop:cofamily pushout property}).
The map
\begin{equation}\label{eq:L_2_L}
 \bL(V,\mR^m)\ \to \ \bL(V,\mR^{m+n})   \ ,
\end{equation}
induced by the embedding $\mR^m\to \mR^{m+n}$ 
as the first $m$ coordinates,
is a homeomorphism from $\bL(V,\mR^m)$ onto the $N$-fixed points
of $\bL(V,\mR^{m+n})$; so the map \eqref{eq:L_2_L}
is an $(O(m)\times O(n)\times G)$-cofibration.

Now we can prove the proposition when $\Uc$ is finite-dimensional. 
We can assume that $\Uc$ is $\mR^{m+n}$ with the standard scalar product, 
and that $U$ is the subspace where the last $n$ coordinates vanish.
The $K$-action on $\Uc$ is given by a continuous homomorphism $\psi:K\to O(m+n)$.
Since $U$ is a $K$-subrepresentation, the image of $\psi$ must be contained in the
subgroup $O(m)\times O(n)$.
The $(K\times G)$-action on the map \eqref{eq:L_2_L}
is then obtained by restriction of the $(O(m)\times O(n)\times G)$-action
along the homomorphism
\[ \psi\times\Id \ : \  K\times G \ \to \ O(m)\times O(n)\times G\ . \]
Restriction along any continuous homomorphism between compact Lie groups
preserves cofibrations by Proposition \ref{prop:cofibrancy preservers}~(i),
so the map \eqref{eq:L_2_L} is a $(K\times G)$-cofibration by the first part. 

Now we treat the case when the dimension of $\Uc$ is infinite.
We choose an exhausting nested sequence of $K$-subrepresentations
\[ U = U_0 \ \subset\ U_1 \ \subset\ U_2\ \subset\ \dots \ .\]
Then all the morphisms $\bL(V,U_{n-1})\to\bL(V,U_n)$ are $(K\times G)$-cofibrations
by the above. Since cofibrations are closed under sequential composites,
the morphism 
\[ \bL(V,U_0)\ \to\  \colim_n\, \bL(V,U_n) \ = \ \bL(V,\Uc) \]
is also a $(K\times G)$-cofibration.

Applying Proposition \ref{prop:cofibrancy preservers}~(iii)
to the normal subgroup $e\times G$ of $K\times G$ shows that
the functor
\[ (e\times G)\bs- \ : \ (K\times G)\bT \ \to \ K\bT \]
takes $(K\times G)$-cofibrations to $K$-cofibrations.
This proves the second claim.

(ii) This is the special case $U=\{0\}$. The space $\bL(V,0)$ is either empty
or consists of a single point; in any case $\bL(V,\{0\})$ 
is $(K\times G)$-cofibrant. Part~(i) then implies that
$\bL(V,\Uc)$ is $(K\times G)$-cofibrant
and $\bL(V,\Uc)/G$ is $K$-cofibrant.
\end{proof}

The following fundamental contractibility property goes back, 
at least, to Boardman and Vogt \cite{boardman-vogt-homotopy everything}.
The equivariant version that we need can be found in \cite[Lemma II 1.5]{lms}.

\begin{prop}\label{prop:isometry spaces contractible} 
Let $G$ be a compact Lie group, $V$ a $G$-representation and $\Uc$ a $G$-universe
such that $V$ embeds into $\Uc$. Then the space $\bL(V,\Uc)$, 
equipped with the conjugation action by $G$, is $G$-equivariantly contractible.
\end{prop}
\begin{proof}
We start by showing that the space $\bL^G(V,\Uc)$
of $G$-equivariant linear isometric embeddings is weakly contractible.
We let $\Uc$ be a $G$-representation of finite or countably infinite dimension. 
Then the map
\[ H \ : \ [0,1]\times \bL^G(V,\Uc)\ \to \ \bL^G(V,\Uc\oplus V) \]
defined by
\[ H(t,\varphi)(v)\ = \ ( t\cdot \varphi(v),\, \sqrt{1-t^2}\cdot v) \]
is a homotopy from the constant map with value $i_2:V\to \Uc\oplus V$
to the map $i_1\circ -$ (postcomposition with $i_1:\Uc\to \Uc\oplus V$).

Since $V$ embeds into $\Uc$ and $\Uc$ is a $G$-universe, it contains
infinitely many orthogonal copies of $V$. In other words, we can assume
that 
\[ \Uc \ = \ \Uc'\oplus V^\infty \]
for some $G$-representation $\Uc'$.
Then
\[ \bL^G(V,\Uc) \ = \ \bL^G(V, \Uc'\oplus V^\infty) \ = \ 
\colim_{n\geq 0} \bL^G(V, \Uc'\oplus V^n) \ ; \]
the colimit is formed along the postcomposition maps
with the direct sum embedding $\Uc'\oplus V^n\to \Uc'\oplus V^{n+1}$.
Every map in the colimit system is a closed embedding and homotopic
to a constant map, by the previous paragraph. 
So the colimit is weakly contractible.

Applying the previous paragraph to a closed subgroup $H$ of $G$
shows that the fixed point space $\bL^H(V,\Uc)$ is weakly contractible;
in other words, $\bL(V,\Uc)$ is $G$-weakly contractible.
The space $\bL(V,\Uc)$ comes with a $(G\times G)$-action
as in \eqref{eq:turn_into_left}, and it is 
$(G\times G)$-cofibrant by Proposition \ref{prop:K G cofibration}~(ii).  
Then $\bL(V,\Uc)$ is also cofibrant as a $G$-space for the
diagonal action, by Proposition \ref{prop:cofibrancy preservers}~(i). 
Since $\bL(V,\Uc)$ is $G$-cofibrant and weakly $G$-contractible,
it is actually equivariantly contractible.
\end{proof}

Now we turn to semifree orthogonal spaces.

\begin{construction}\label{con:semifree spc}
Given a compact Lie group $G$ and a $G$-representation $V$, the functor
\[ \ev_{G,V}\ : \ \spc \ \to \ G\bT \]
that sends an orthogonal space $Y$ to the $G$-space $Y(V)$
has a left adjoint
\begin{equation}\label{eq:define L_G,V}
 \bL_{G,V}\ : \ G\bT\ \to \ \spc \ .  
\end{equation}
To construct the left adjoint we note that $G$ acts from the right on $\bL(V,W)$ by
\[ (\varphi\cdot g)(v) \ = \ \varphi(g v)  \]
for $\varphi\in\bL(V,W)$, $g\in G$ and $v\in V$.
Given a $G$-space $A$, the value of $\bL_{G,V} A$ at an inner product space $W$ is
\[ (\bL_{G,V}A)(W) \ = \ \bL(V,W)\times_G A\ = \ 
\left( \bL(V,W)\times A \right) / (\varphi g,a)\sim (\varphi, g a)\ . \]
We refer to $\bL_{G,V} A$ as
the {\em semifree orthogonal space}\index{subject}{orthogonal space!semifree}
\index{subject}{semifree orthogonal space}\index{symbol}{$\bL_{G,V}$ - {semifree orthogonal space generated by $(G,V)$}}
generated by $A$ at $(G,V)$.
We also denote by $\bL_{G,V}$ the orthogonal space with
\[ \bL_{G,V}(W)\ = \ \bL(V,W)/G \ . \]
So $\bL_{G,V}$ is isomorphic to the semifree orthogonal space 
generated at $(G,V)$ by a one-point $G$-space.
\end{construction}

The `freeness' property of $\bL_{G,V} A$
is a consequence of the enriched Yoneda lemma,\index{subject}{Yoneda lemma!enriched}
see Remark \ref{rk:enriched Yoneda} or \cite[Sec.\,1.9]{kelly-enriched category};
it means explicitly that
for every orthogonal space $Y$ and every continuous $G$-map $f:A\to Y(V)$ 
there is a unique morphism 
$f^\flat:\bL_{G,V} A\to Y$ of orthogonal spaces such that the composite
\[  A \xra{\ [\Id, -]\ }\ 
\bL(V,V) \times_G A = (\bL_{G,V} A)(V) \ \xra{f^\flat(V)} \ Y(V) \]
is $f$.
Indeed, the map $f^\flat(W)$ is the composite
\[ \bL(V,W) \times_G A \ \xra{\Id\times_G f} \
\bL(V,W) \times_G Y(V) \ \xra{[\varphi,y]\mapsto Y(\varphi)(y)} \ Y(W)\ .\]

\begin{eg}\label{eg:free is closed}
For every compact Lie group $G$, every $G$-representation $V$
and every $G$-space $A$ the semifree orthogonal space $\bL_{G,V} A$ 
is closed.\index{subject}{orthogonal space!closed}
To see that we let $\varphi:U\to W$ be a linear isometric embedding. 
Since $\bL(V,U)$ is compact, the continuous injection
\[ \bL(V,\varphi)\ : \  \bL(V,U)\ \to \ \bL(V,W) \]
is a closed embedding. 
So the map $\bL(V,\varphi)\times A$ is a closed embedding as well.
The orbit map 
\[ \bL(V,\varphi)\times_G A\ : \  \bL(V,U)\times_G A\ \to \  \bL(V,W)\times_G A  \]
is then a closed embedding
by Proposition \ref{prop:G quotient properties} (iii).
\end{eg}

The next proposition identifies the fixed point spaces 
of a semifree orthogonal space $\bL_{G,V}$.
A certain family $\Fc(K;G)$ of subgroups of $K\times G$ arises naturally,
which we call `graph subgroups'.

\begin{defn}
Let $K$ and $G$ be compact Lie groups. The family
$\Fc(K;G)$\index{symbol}{$\Fc(K;G)$ - {family of graph subgroups of $K\times G$}}
of {\em graph subgroups}\index{subject}{graph subgroup}
consists of those closed subgroups $\Gamma$ of $K\times G$
that intersect $1\times G$ only in the neutral element $(1,1)$.
\end{defn}

The name `graph subgroup' stems from the fact that
$\Fc(K;G)$ consists precisely of the graphs of all `subhomomorphisms',
i.e., continuous homomorphisms $\alpha:L\to G$ from a closed subgroup $L$ of $K$.
Clearly, the graph $\Gamma(\alpha)=\{(l,\alpha(l))\ |\ l\in L\}$ 
of every such homomorphism belongs to $\Fc(K;G)$.
Conversely, for $\Gamma\in \Fc(K;G)$ we let $L\leq K$ be the
image of $\Gamma$ under the projection $K\times G\to K$.
Since $\Gamma\cap (1\times G)=\{(1,1)\}$, every element $l\in L$
then has a unique preimage $(l,\alpha(l))$ under the projection,
and the assignment $l\mapsto \alpha(l)$ is a continuous homomorphism
from $L$ to $G$ whose graph is $\Gamma$.

We recall that a 
{\em universal $G$-space}\index{subject}{universal space!for a set of subgroups}
for a family $\Fc$ of closed subgroups 
is a cofibrant $G$-space $E$ such that
\begin{itemize}
\item all isotropy groups of $E$ belong to the family $\Fc$, and
\item for every $H\in\Fc$ the fixed point space $E^H$ is weakly contractible.
\end{itemize}

If $V$ and $W$ are $G$-representations, then restriction of a linear isometry
from $V\oplus W$ to $V$ defines a $G$-equivariant morphism of orthogonal spaces\index{symbol}{$\rho_{G,V,W}$ - {fundamental global equivalence of orthogonal spaces}}
\[ \rho_{V,W}\ : \ \bL_{V\oplus W}\ \to \ \bL_V\ . \]
If $U$ is a $K$-representation, then we combine the left $K$-action
and the right $G$-action on $\bL(V,U)$ into a left action
of $K\times G$ as in \eqref{eq:turn_into_left}.

\begin{prop}\label{prop:free_orthogonal_space}
Let $G$ and $K$ be compact Lie groups and $V$ a faithful $G$-representation.
\begin{enumerate}[\em (i)]
\item 
The $(K\times G)$-space $\bL_{V}(\Uc_K)=\bL(V,\Uc_K)$ is a universal space 
for the family $\Fc(K;G)$ of graph subgroups.\index{subject}{universal space!for the family of graph subgroups}
\item If $W$ is another $G$-representation, then the restriction map
\[  \rho_{V,W}(\Uc_K)\ : \ \bL(V\oplus W, \Uc_K) \ \to \ \bL(V, \Uc_K) \]
is a $(K\times G)$-homotopy equivalence.
For every $G$-space $A$, the map
\[ (\rho_{V,W}\times_G A)(\Uc_K)\ : \ 
(\bL_{G,V\oplus W} A) (\Uc_K) \ \to \ ( \bL_{G,V} A)(\Uc_K)  \]
is a $K$-homotopy equivalence and the morphism of orthogonal spaces
\[ \rho_{V,W}\times_G A\ : \ \bL_{G,V\oplus W} A \ \to \  \bL_{G,V} A \]
is a global equivalence.
\end{enumerate}
\end{prop}
\begin{proof}
(i) We let $\Gamma$ be any closed subgroup of $K\times G$.
Since the $G$-action on $V$ is faithful, the induced right $G$-action
on  $\bL(V,\Uc_K)$ is free.
So if $\Gamma$ intersects $1\times G$ non-trivially, then
$\bL(V,\Uc_K)^\Gamma$ is empty.
On the other hand, if $\Gamma\cap (1\times G)=\{(1,1)\}$,
then $\Gamma$ is the graph of a unique continuous homomorphism $\alpha:L\to G$,
where $L$ is the projection of $\Gamma$ to $K$. Then
\[ \bL(V,\Uc_K)^\Gamma \ = \  \bL^L(\alpha^* V, \Uc_K) \]
is the space of $L$-equivariant linear isometric embeddings from
the $L$-represen\-tation $\alpha^* V$ to the underlying $L$-universe of $\Uc_K$.
Since $\Uc_K$ is a complete $K$-universe, the underlying $L$-universe is also
complete (Remark \ref{rk:complete universes restrict}), 
and so the space $\bL^L(\alpha^* V, \Uc_K)$  is contractible
by Proposition \ref{prop:isometry spaces contractible}.
The space $\bL(V,\Uc_K)$ is cofibrant as a
$(K\times G)$-space by Proposition \ref{prop:K G cofibration}~(ii).

(ii)
Since $G$ acts faithfully on $V$, and hence also on $V\oplus W$, the
$(K\times G)$-spaces $\bL(V\oplus W,\Uc_K)$ and $\bL(V,\Uc_K)$
are universal spaces for the same family $\Fc(K;G)$, by part~(i).
So the map $\rho_{V,W}(\Uc_K):\bL(V\oplus W,\Uc_K)\to\bL(V,\Uc_K)$
is a $(K\times G)$-equivariant homotopy equivalence, 
see Proposition \ref{prop:universal spaces}~(ii).
The functor $-\times_G A$ preserves homotopies, so 
the restriction map $(\rho_{V,W}\times_G A)(\Uc_K)$ is a $K$-homotopy equivalence.

The orthogonal spaces $\bL_{G,V\oplus W} A$ and $\bL_{G,V} A$ are closed
by Example \ref{eg:free is closed}, 
so Proposition \ref{prop:global eq for closed} applies 
and shows that $\rho_{V,W}\times_G A$ is a global equivalence.
\end{proof}

\index{subject}{global classifying space|(}

\begin{defn}\label{def:global classifying}
The {\em global classifying space}\index{subject}{global classifying space}\index{subject}{classifying space!global|see{global classifying space}}\index{symbol}{$B_{\gl} G$ - {global classifying space of the compact Lie group $G$}}
$B_{\gl} G$ of a compact Lie group $G$ is the semifree orthogonal space
\[ B_{\gl} G \ = \ \bL_{G,V}\ = \ \bL(V,-)/G\ , \]
where $V$ is any faithful $G$-representation.
\end{defn}

The global classifying space $B_{\gl}G$ is well-defined up to preferred 
zigzag of global equivalences of orthogonal spaces. Indeed, if $V$ and
$\bar V$ are two faithful $G$-representations, then $V\oplus \bar V$ is
yet another one, and the two restriction morphisms
\[ \bL_{G,V}\ \xleftarrow{\quad}\ 
 \bL_{G,V\oplus \bar V}\ \xra{\quad}\ \bL_{G,\bar V} \]
are global equivalences by 
Proposition \ref{prop:free_orthogonal_space}~(ii).

\begin{eg}\label{eg:B C_2}\index{subject}{global classifying space!of $C_2$}  
We make the global classifying space more explicit in 
the smallest non-trivial example, i.e., for the cyclic group
$C_2$ of order~2. The sign representation $\sigma$ 
of $C_2$ is faithful, so we can take $B_{\gl} C_2$ 
to be the semifree orthogonal space generated by
$(C_2,\sigma)$; its value at an inner product space $W$  is
\[ (B_{\gl}C_2)(W) \ = \ \bL_{C_2,\sigma}(W) \ = \ 
\bL(\sigma,W)/C_2\ .\]
Evaluation at any of the two unit vectors in $\sigma$ 
is a homeomorphism $\bL(\sigma,W)\iso S(W)$
to the unit sphere of $W$, and the $C_2$-action on the left
becomes the antipodal action on $S(W)$. So the map descends to a homeomorphism
between $\bL(\sigma,W)/C_2$ and $P(W)$, the projective space of $W$,
and hence
\[ (B_{\gl}C_2)(W) \ \iso \ P(W)\ .\]
So for a compact Lie group $K$, the underlying $K$-space of
$B_{\gl} C_2$ is $P(\Uc_K)$, the projective space of a complete $K$-universe. 
In particular, the underlying non-equivariant space is
homeomorphic to $\RP^\infty$.
\end{eg}

\begin{rk}[$B_{\gl}G$ globally classifies principal $G$-bundles]
  The term `glo\-bal classifying space' is justified by the fact that $B_{\gl} G$
  `globally classifies principal $G$-bundles'. 
  We recall that a {\em $(K,G)$-bundle},\index{subject}{K,G-bundle@$(K,G)$-bundle}
  also called a {\em $K$-equivariant $G$-principal bundle},
  is a principal $G$-bundle in the category of $K$-spaces, i.e., a $G$-principal bundle 
  $p:E\to B$ that is also a morphism of $K$-spaces and such that
  the actions of $G$ and $K$ on the total space $E$ commute (see for example \cite[Ch.\,I (8.7)]{tomDieck-transformation}).
  For every compact Lie group $K$, the quotient map
  \[ q\ : \ \bL(V,\Uc_K) \ \to \ \bL(V,\Uc_K)/G \ = \ \bL_{G,V}(\Uc_K) 
  \ = \ (B_{\gl}G)(\Uc_K) \]
  is a principal $(K,G)$-bundle. Indeed, the total space $\bL(V,\Uc_K)$
  is homeomorphic to the Stiefel manifold of $\dim(V)$-frames in $\mR^\infty$,
  and hence it admits a CW-structure. Every CW-complex is a normal Hausdorff space
  (see for example \cite[Prop.\,A.3]{hatcher} or \cite[Prop.\,1.2.1]{fritsch-piccinini}),
  hence completely regular. So $\bL(V,\Uc_K)$ is completely regular.
  Since the $G$-action on $\bL(V,\Uc_K)$ is free, the quotient map $q$
  is a $G$-principal bundle by \cite[Prop.\,1.7.35]{palais-classification}
  or \cite[II Thm.\,5.8]{bredon-intro}. 
  Moreover, this bundle is universal in the sense of the following proposition.
  Every $G$-space that admits a $G$-CW-structure is paracompact,
  see \cite[Thm.\,3.2]{murayama-GANRs} 
  (this reference is rather sketchy, but one can follow the non-equivariant
  argument spelled out in more detail in \cite[Thm.\,1.3.5]{fritsch-piccinini}).
  So the next proposition applies in particular to all $G$-CW-complexes.
\end{rk}

  \begin{prop}\label{prop:classifying classifies} 
    Let $V$ be a faithful representation of a compact Lie group $G$.
    Then for every paracompact $K$-space $A$ the map
    \[   [A, \bL(V,\Uc_K)/G]^K \ \to \ \Prin_{(K,G)}(A)\ ,\quad
    [f]\ \longmapsto \ [f^*(q)]  \]
    from the set of equivariant homotopy classes of $K$-maps to the
    set of isomorphism classes of $(K,G)$-bundles is bijective.
  \end{prop}
  \begin{proof}
    I do not know a reference for the result in precisely this form,
    so I sketch how to deduce it from various results in the
    literature about equivariant fiber bundles.
    A principal $(K,G)$-bundle $p:E\to B$~is {\em equivariantly trivializable} 
    if there is a closed subgroup $L$ of $K$, a continuous homomorphism
    $\alpha:L\to G$, an $L$-space $X$ and an isomorphism of $(K,G)$-bundles
    between $p$ and the projection
    \[ (K\times G)\times_L X \ \to \ K\times_L X \ ; \]
    here the source is the quotient space of
    $K\times G\times X$ by the equivalence relation
    $(k,g,l x)\sim (k l, g\alpha(l),x)$ 
    for all $(k,g,l,x)\in K\times G\times L\times X$.
    A principal $(K,G)$-bundle $p:E\to B$~is {\em numerable} if $B$ 
    has a trivializing (in the above sense) open cover
    by $K$-invariant subsets, such that moreover the cover admits 
    a subordinate partition of unity by $K$-invariant functions.

    A {\em universal $(K,G)$-principal bundle}
    is a numerable $(K,G)$-bundle $p^u:E^u\to B^u$ 
    such that for every $K$-space $A$ the map
    \[  
    [A, B^u]^K \ \to \ \Prin^{\text{num}}_{(K,G)}(A)\ , \quad
    [f]\ \longmapsto \ [f^*(p^u)]
    \]
    is bijective, where now the target is the set of isomorphism classes of
    {\em numerable} principal $(K,G)$-bundles.
    Over a paracompact base, every principal $(K,G)$-bundle is numerable
    by \cite[Cor.\,1.5]{lashof-equivariant bundles}. So we are done if we can show
    that $q:\bL(V,\Uc_K)\to \bL(V,\Uc_K)/G$ is a universal $(K,G)$-principal bundle 
    in the above sense.

    Universal $(K,G)$-principal bundles can be built in different ways; 
    the most common construction is a version of
    Milnor's infinite join \cite{milnor-universalII}, 
    see for example \cite[3.1 Satz]{tomDieck-faserbundel}
    or \cite[I Theorem (8.12)]{tomDieck-transformation}.
    Another method is via bar construction, 
    compare \cite{waner-equivariant classifying}.
    I do not know of a reference that explicitly identifies the projection
    $q:\bL(V,\Uc_K)\to \bL(V,\Uc_K)/G$ as a universal bundle in the present sense,
    so we appeal to Lashof's criterion \cite[Thm.\,2.14]{lashof-equivariant bundles}. 
    The base space $\bL(V,\Uc_K)/G$ is the union, along h-cofibrations,
    of the compact spaces $\bL(V,W_i)/G$, where $\{W_i\}_{i\geq 1}$
    is any exhausting sequence of subrepresentations of $\Uc_K$.
    Since compact spaces are paracompact, the union $\bL(V,\Uc_K)/G$ 
    is paracompact, see for example \cite[Prop.\,A.5.1\,(v)]{fritsch-piccinini}.
    Since $\bL(V,\Uc_K)/G$ is also normal, hence completely regular,
    the bundle $q$ is numerable by Corollaries~1.5 and~1.13 
    of \cite{lashof-equivariant bundles}.
    Moreover, the fixed points of $\bL(V,\Uc_K)$
    under any graph subgroup of $(K\times G)$ are contractible by
    Proposition \ref{prop:free_orthogonal_space}~(i),
    so Theorem~2.14 of \cite{lashof-equivariant bundles}
    applies and shows that $q:\bL(V,\Uc_K)\to \bL(V,\Uc_K)/G$ is strongly universal, 
    and hence a universal principal $(K,G)$-bundle.
    \end{proof}

    As another example we look at the case $G=O(n)$, the $n$-th orthogonal group. 
    The category of principal $O(n)$-bundles is equivalent to the category 
    of euclidean vector bundles of rank $n$, via the associated frame bundle.
    By the same construction, principal $(K,O(n))$-bundles can be identified
    with $K$-equivariant euclidean vector bundles of rank $n$ over $K$-spaces.
    The space $\bL(\mR^n,\Uc_K)/O(n)$ is homeomorphic to $G r_n(\Uc_K)$,
    the Grassmannian of $n$-planes in $\Uc_K$. When $K$ is a trivial group,
    the fact that $G r_n(\mR^\infty)$ is a classifying space for rank $n$
    vector vector bundles over paracompact spaces is proved in various textbooks.
    Since $O(1)$ is a cyclic group of order~2, this gives another perspective on
    Example~\ref{eg:B C_2}. 
\index{subject}{global classifying space|)}

\section{Global model structure for orthogonal spaces}
\label{sec:global model structures spaces}

In this section we establish the global model structure on the category of orthogonal spaces, see Theorem \ref{thm:All global spaces}.
Towards this aim we first discuss a `strong level model structure', 
which we then localize.
In Proposition \ref{prop:cofree properties}
we use the global model structure to relate unstable global homotopy
theory to the homotopy theory of $K$-spaces for a fixed compact Lie group $K$.
We also show that global classifying spaces of compact Lie groups 
with abelian identity component are `cofree', i.e.,
right induced from non-equivariant classifying spaces, see
Theorem \ref{thm:BA is right induced}. 
At the end of this section we briefly discuss the realization of
simplicial orthogonal spaces in Construction \ref{con:realize simplicial}; 
we show that the realization takes levelwise global equivalences 
to global equivalences, under a certain `Reedy flatness' condition
(Proposition \ref{prop:realization invariance in spc}).

\medskip

There is a functorial way to write an orthogonal space as a 
sequential colimit of orthogonal spaces which are made from the information
below a fixed dimension.
We refer to this as the {\em skeleton filtration} of an orthogonal
space. The word `filtration' should be used with caution
because the maps from the skeleta to the orthogonal space
need not be injective.

The skeleton filtration is in fact a special case of 
a more general skeleton filtration on certain enriched functor categories
that we discuss in Appendix \ref{app:enriched functors}.
Indeed, if we specialize the base category to $\Vc=\bT$,
the category of spaces under cartesian product, and the index category
to $\Dc=\bL$, then the functor category $\Dc^\ast$ becomes the
category $\spc$ of orthogonal spaces. The dimension function needed in
the construction and analysis of skeleta is the vector space dimension.

We denote by $\bL^{\leq m}$ the full topological subcategory of the linear isometries
category $\bL$ with objects all inner product spaces of dimension at most $m$.
We denote by $\spc^{\leq m}$ the category of continuous functors from $\bL^{\leq m}$
to $\bT$. The restriction functor 
\[  \spc \ \to \ \spc^{\leq m}\ , \quad Y \ \longmapsto \ Y^{\leq m}\ = \ Y|_{\bL^{\leq m}}\]
has a left adjoint 
\[ l_m \ : \ \ \spc^{\leq m} \ \to \ \spc\]
given by an enriched Kan extension as follows. 
The extension $l_m(Z)$ of a continuous functor $Z:\bL^{\leq m}\to\bT$ 
is a coequalizer of the two morphisms of orthogonal spaces
\begin{equation}  \label{eq:skeleton_coequalizer spc}
 \xymatrix@C=7mm{ 
\coprod_{0\leq j\leq k\leq m} \bL(\mR^k,-)\times\bL(\mR^j,\mR^k)\times Z(\mR^j) 
\ \ar@<-.4ex>[r]  \ar@<.4ex>[r]  &
\ \coprod_{0\leq i\leq m} \bL(\mR^i,-)\times Z(\mR^i) 
 }   
\end{equation}
One morphism arises from the composition morphisms
\[  \bL(\mR^k,-)\times\bL(\mR^j,\mR^k) \ \to \ \bL(\mR^j,-) \]
and the identity on $Z(\mR^j)$;
the other morphism arises from the action maps
\[  \bL(\mR^j,\mR^k)\times Z(\mR^j) \ \to \ Z(\mR^k) \]
and the identity on the free orthogonal space $ \bL(\mR^k,-)$.
Colimits in the category of orthogonal spaces are created objectwise, 
so the value $l_m(Z)(V)$ at an inner product space can be calculated 
by plugging $V$ into the variable slot 
in the coequalizer diagram \eqref{eq:skeleton_coequalizer spc}.

It is a general property of Kan extensions along a fully faithful functor 
(such as the inclusion $\bL^{\leq m}\to\bL$) that the values do not
change on the given subcategory,
see for example \cite[Prop.\,4.23]{kelly-enriched category}. 
More precisely, the adjunction unit
\[ Z \ \to \ (l_m(Z))^{\leq m} \]
is an isomorphism for every continuous functor $Z:\bL^{\leq m}\to\bT$. 

\begin{defn}
The {\em $m$-skeleton},\index{subject}{skeleton!of an orthogonal space}
for $m\geq 0$, of an orthogonal space $Y$ is the orthogonal space
\[ \sk^m Y\ = \ l_m(Y^{\leq m}) \ ,\]
the extension of the restriction of $Y$ to $\bL^{\leq m}$.
It comes with a natural morphism $i_m:\sk^m Y\to Y$, the counit of
the adjunction $(l_m,(-)^{\leq m})$.
The {\em $m$-th latching space}\index{subject}{latching space!of an orthogonal space}
of $Y$ is the $O(m)$-space
\[ L_m Y \ = \ (\sk^{m-1} Y)(\mR^m) \ ;\]
it comes with a natural $O(m)$-equivariant map
\[  \nu_m=i_{m-1}(\mR^m)\ :\ L_m Y\ \to \ Y(\mR^m) \ , \]
the {\em $m$-th latching map}.
\end{defn}

We agree to set $\sk^{-1} Y=\emptyset$, the empty orthogonal space,
and $L_0 Y=\emptyset$, the empty space. The value
\[ i_m(V)\ :\ (\sk^m Y)(V)\ \to \ Y(V) \]
of the morphism $i_m$ is an isomorphism for all inner product spaces $V$ of dimension
at most $m$.

The two morphisms $i_{m-1}:\sk^{m-1} Y\to Y$ and
$i_m:\sk^m Y\to Y$ both restrict to isomorphisms on $\bL^{\leq m-1}$, so there is a 
unique morphism $j_m:\sk^{m-1} Y\to \sk^m Y$ such that $i_m\circ j_m=i_{m-1}$.
The sequence of skeleta stabilizes to $Y$ in a 
very strong sense. For every inner product space $V$, the maps
$j_m(V)$ and $i_m(V)$ are homeomorphisms as soon as $m >\dim(V)$.
In particular, $Y(V)$ is a colimit, with respect to the morphisms $i_m(V)$, 
of the sequence of maps $j_m(V)$.
Since colimits in the category of orthogonal spaces are created objectwise,
we deduce that the orthogonal space $Y$ is a colimit, 
with respect to the morphisms $i_m$, of the sequence of morphisms $j_m$.

We denote by
\[ G_m \ : \ O(m)\bT \ \to \ \spc \]
the left adjoint to the functor $Y\mapsto Y(\mR^m)$.
So $G_m$ is a shorthand notation for $\bL_{O(m),\mR^m}$,
the semifree functor \eqref{eq:define L_G,V} indexed by the tautological
$O(m)$-representation.
Proposition \ref{prop:basic skeleton pushout square}
specializes to:

\begin{prop}
For every orthogonal space $Y$ and every $m\geq 0$ the commutative square
\begin{equation}\begin{aligned}\label{eq:filtration pushout spaces}
\xymatrix@C=15mm{ G_m L_m Y \ar[r]^{G_m\nu_m} \ar[d] & G_m Y(\mR^m) \ar[d]\\
\sk^{m-1} Y \ar[r]_-{j_m} & \sk^m Y}  
\end{aligned}\end{equation}
is a pushout of orthogonal spaces.
The two vertical morphisms are adjoint to the identity of $L_m Y$ 
respectively $Y(\mR^m)$.
\end{prop}

\begin{eg}
As an illustration of the definition, we describe the skeleta and  latching objects
for small values of $m$. We have 
\[ \sk^0 Y \ = \ \text{const}(Y(0)) \ ,\]
the constant orthogonal space with value $Y(0)$; the latching map
\[ \nu_1 \ : \ L_1 Y \ = \ (\sk^0 Y)(\mR) \ = \ Y(0) \ \xra{Y(u)} \ Y(\mR)\]
is the map induced by the unique linear isometric embedding $u:0\to\mR$.
Now we evaluate the pushout square \eqref{eq:filtration pushout spaces}
for $m=1$ at an inner product space $V$; the result is a pushout square
of $O(1)$-spaces
\[ 
\xymatrix@C=15mm{ P(V)\times Y(0) \ar[r] \ar[d]_{\text{proj}} & 
\bL(\mR,V)\times_{O(1)} Y(\mR) \ar[d]\\
Y(0)  \ar[r] & (\sk^1 Y)(V)}   \]
where $P(V)$ is the projective space of $V$.
Here we exploit that $O(1)$ acts trivially on $L_1 Y=Y(0)$
and we can thus identify
\[ (G_1 L_1 Y)(V)\ = \ \bL(\mR,V)\times_{O(1)} Y(0)\ \iso \ P(V)\times Y(0)\ , \quad
[\varphi, y]\ \longmapsto \ (\varphi(\mR),y) \ .\]
The upper horizontal map sends $(\varphi(\mR),y)$ to $[\varphi, Y(u)(y)]$.
\end{eg}

\begin{eg}[Latching objects of free orthogonal spaces]
  We let $V$ be an $n$-dimensional representation of compact Lie group $G$,
  and $A$ a $G$-space. 
  Then the semifree orthogonal space \eqref{eq:define L_G,V}
  generated by $A$ in level $V$ is `purely $n$-dimensional' in the following sense.
  The evaluation functor
  \[ \ev_{G,V}\ : \ \spc \ \to \ G\bT \]
  factors through the category $\bL^{\leq n}$ as the composite
  \[  \spc \ \to \  \spc^{\leq n} \ \xra{\ev_{G,V}} \ G\bT \ .\]
  So the left adjoint semifree functor $\bL_{G,V}$ can be chosen as the composite
  of the two individual left adjoints
  \[ \bL_{G,V}\ = \ l_n \circ l_{G,V} \ .\]
  Here $l_{G,V}:G\bT \to \spc^{\leq n}$ 
  is given at a $G$-space $A$ and an inner product
  space $W$ of dimension at most $n$ by
  \[ (l_{G,V} A)(W) \ = \
  \begin{cases}
    \bL(V,W)\times_G A & \text{ if $\dim(W)= n$,}\\
    \qquad \emptyset & \text{ if $\dim(W)< n$.}
  \end{cases}\]
  The space $(\bL_{G,V}A)_m$ is trivial for $m < n$, and hence the 
  latching space $L_m (\bL_{G,V}A)$ is trivial for $m\leq n$.
  For $m>n$ the latching map $\nu_m:L_m(\bL_{G,V}A) \to (\bL_{G,V}A)_m$ 
  is an isomorphism. So for $m<n$ the skeleton $\sk^m (\bL_{G,V}A)$ is trivial, 
  and for $m \geq n$ the skeleton $\sk^m (\bL_{G,V}A)=\bL_{G,V}A$ is the entire orthogonal space. 
\end{eg}

Now we work our way towards the strong level model structure of orthogonal spaces.
Proposition \ref{prop:general level model structure} 
is a fairly general recipe for constructing 
level model structures on a category such as orthogonal spaces.
We specialize the general construction to the situation at hand.
We recall from Definition \ref{def:strong level equivalence spaces}
that a morphism $f:X\to Y$ of orthogonal spaces is a
{\em strong level equivalence}\index{subject}{strong level equivalence!of orthogonal spaces} 
(respectively  {\em strong level fibration})\index{subject}{strong level fibration!of orthogonal spaces} 
if for every compact Lie group $G$ and every $G$-representation $V$ the map
 $f(V)^G:X(V)^G\to Y(V)^G$ is a weak equivalence (respectively a Serre fibration).

\begin{lemma}\label{lemma:characterize strong level equivalences} 
For every morphism $f:X\to Y$ of orthogonal spaces,
the following are equivalent.
  \begin{enumerate}[\em (i)]
\item The morphism $f$ is a strong level equivalence.
\item For every compact Lie group $G$ and every faithful $G$-representation $V$ the map
    $f(V):X(V)\to Y(V)$ is a $G$-weak equivalence.
\item The map $f(\mR^m):X(\mR^m)\to Y(\mR^m)$ is an $O(m)$-weak equivalence
  for every $m\geq 0$.
  \end{enumerate}
  \end{lemma}
  \begin{proof}
Clearly, condition~(i) implies condition~(ii), and that implies condition~(iii)
(because the tautological action of $O(m)$ on $\mR^m$ is faithful).
So we suppose that $f(\mR^m)$ is an $O(m)$-weak equivalence
for every $m\geq 0$, and we show that $f$ is a strong level equivalence.
Given a $G$-representation $V$ of dimension $m$, 
we choose a linear isometry $\varphi:V\iso\mR^m$; conjugation by $\varphi$ turns 
the $G$-action on $V$ into a homomorphism $\rho:G\to O(m)$, i.e.,
\[ \rho(g)\ = \ \varphi\circ(g\cdot-)\circ \varphi^{-1}\ . \]
The homeomorphism $X(\varphi):X(V)\to X(\mR^m)$ then restricts to a homeomorphism
\[ X(V)^G \ \iso \ X(\mR^m)^{\rho(G)} \ .\]
This homeomorphism is natural for morphisms of orthogonal spaces, so the hypothesis that
$ f(\mR^m)^{\rho(G)}: X(\mR^m)^{\rho(G)}\to Y(\mR^m)^{\rho(G)}$ is a weak equivalence
implies that also the map $f(V)^G:X(V)^G\to Y(V)^G$ is a weak equivalence.
\end{proof}

The same kind of reasoning as in Lemma \ref{lemma:characterize strong level equivalences}
shows:

\begin{lemma}\label{lemma:characterize strong level fibrations}  
For every morphism $f:X\to Y$ of orthogonal spaces,
the following are equivalent.
  \begin{enumerate}[\em (i)]
\item The morphism $f$ is a strong level fibration.
\item For every compact Lie group $G$ and every faithful $G$-representation $V$ the map
    $f(V):X(V)\to Y(V)$ is a fibration in the projective model structure of $G$-spaces.
\item The map $f(\mR^m):X(\mR^m)\to Y(\mR^m)$ is an $O(m)$-fibration
  for every $m\geq 0$.
  \end{enumerate}
  \end{lemma}

\begin{defn}
  A morphism of orthogonal spaces $i:A\to B$ 
  is a {\em flat cofibration}\index{subject}{flat cofibration!of orthogonal spaces} 
  if the latching morphism
  \[  \nu_m i\ = \ i(\mR^m)\cup \nu_m^B\ :\ A(\mR^m)\cup_{L_m A}L_m B\ \to\ B(\mR^m) \]
  is an $O(m)$-cofibration  for all $m\geq 0$.
  An orthogonal space $B$ is 
  {\em flat}\index{subject}{flat!orthogonal space}\index{subject}{orthogonal space!flat} 
  if the unique morphism from the empty orthogonal space
  to $B$ is a flat cofibration.
  Equivalently, for every $m\geq 0$ the latching map $\nu_m:L_m B \to B(\mR^m)$
  is an $O(m)$-cofibration.
\end{defn}

We are ready to establish the strong level model structure.

\begin{prop}\label{prop:strong level spaces}
    The strong level equivalences, strong level fibrations 
    and flat cofibrations form a topological cofibrantly generated model structure,
    the {\em strong level model structure},\index{subject}{strong level model structure!for orthogonal spaces}\index{subject}{model structure!strong level|see{strong level model structure}}
    on the category of orthogonal spaces.
\end{prop}
\begin{proof}
We apply Proposition \ref{prop:general level model structure} as follows.
We let $\Cc(m)$ be the projective
model structure on the category of $O(m)$-spaces
(with respect to the set of all closed subgroups of $O(m)$),
compare Proposition \ref{prop:proj model structures for G-spaces}.
The classes of level equivalences, level fibrations and cofibrations
in the sense of Proposition \ref{prop:general level model structure} 
then precisely become the strong level equivalences,
strong level fibrations and flat cofibrations.

In this situation the consistency condition 
(see Definition \ref{def:consistency condition})
is a consequence of a stronger property, namely that the functor
\[ \bL(\mR^m,\mR^{m+n})\times_{O(m)} - \ : \ O(m)\bT \ \to\ O(m+n)\bT \]
takes acyclic cofibrations to acyclic cofibrations 
(in the two relevant projective model structures).
Since the functor is a left adjoint, it suffices to prove the claim for the generating
acyclic cofibrations, i.e., the maps
\[    O(m)/H \times j_k \]
for all $k\geq 0$ and all closed subgroups $H$ of $O(m)$,
where $j_k:D^k\times\{0\}\to D^k\times [0,1]$ is the inclusion.
The functor under consideration takes this generator to the map
$\bL(\mR^m,\mR^{m+n})/ H\times j_k$,
which is an acyclic $O(m+n)$-cofibration because 
$\bL(\mR^m,\mR^{m+n})/H$
is cofibrant as an $O(m+n)$-space, by Proposition \ref{prop:K G cofibration} (iii). 

We describe explicit sets of generating cofibrations 
and generating acyclic cofibrations.
We let $I^{\str}$ be the set of all morphisms $G_m i$ for $m\geq 0$
and for $i$ in the set of generating cofibrations for the projective model
structure on the category of $O(m)$-spaces specified 
in \eqref{eq:I_for_F-proj_on_GT}. 
Then the set $I^{\str}$ detects the acyclic fibrations 
in the strong level model structure 
by Proposition \ref{prop:general level model structure}~(iii). 
Similarly, we let $J^{\str}$ be the set of all morphisms $G_m j$ 
for $m\geq 0$ and for $j$ in the set
of generating acyclic  cofibrations for the projective model
structure on the category of $O(m)$-spaces specified 
in \eqref{eq:J_for_F-proj_on_GT}.
Again by Proposition \ref{prop:general level model structure}~(iii),
$J^{\str}$ detects the fibrations in the strong level model structure. 

The model structure is topological by Proposition \ref{prop:topological criterion},
where we take $\Gc$ as the set of orthogonal spaces $L_{H,\mR^m}$
for all $m\geq 0$ and all closed subgroups $H$ of $O(m)$, and we take $\Zc=\emptyset$.
\end{proof}

For easier reference we make the generating (acyclic) cofibrations
of the strong level model structure even more explicit. Using the isomorphism  
\[ G_m ( O(m)/H ) \ = \ \bL(\mR^m,-)\times_{O(m)} ( O(m)/H ) \ \iso \ 
\bL(\mR^m,-)/H  \ = \bL_{H,\mR^m}  \ ,\]
we can identify $I^{\str}$ with the set of all morphisms
\[  \bL_{H,\mR^m}\times i_k \ : \ \bL_{H,\mR^m}\times \partial D^k 
\ \to \  \bL_{H,\mR^m}\times D^k  \]
for all $k,m\geq 0$ and all closed subgroups $H$ of $O(m)$.
The tautological action of $H$ on $\mR^m$ is faithful; 
conversely every pair $(G,V)$ consisting of a compact Lie group and a faithful
representation is isomorphic to a pair $(H,\mR^m)$ 
for some closed subgroup $H$ of $\mR^m$.
We conclude that $I^{\str}$ is a set of representatives of the 
isomorphism classes of morphisms
\[  \bL_{G,V}\times i_k \ : \  \bL_{G,V}\times  \partial D^k
  \ \to \  \bL_{G,V} \times D^k   \]
for $G$ a compact Lie group, $V$ a faithful $G$-representation and $k\geq 0$.
Similarly, $J^{\str}$ is a set of representatives of the 
isomorphism classes of morphisms
\[ \bL_{G,V}\times  j_k\ : \  \bL_{G,V}\times D^k\times\{0\} 
  \ \to \  \bL_{G,V}  \times D^k\times [0,1]   \]
for $G$ a compact Lie group, $V$ a faithful $G$-representation and $k\geq 0$.

\begin{prop}\label{prop:flat becomes K-cofibration}
  Let $K$ be a compact Lie group and $\varphi:W\to U$ a linear isometric embedding 
  of $K$-representations, where $W$ is finite-dimensional, and $U$ is finite-dimensional
  or countably infinite dimensional.
  \begin{enumerate}[\em (i)]
  \item    For every flat cofibration of orthogonal spaces $i:A\to B$ 
    the maps
    \begin{align*}
      i(U)\ &: \ A(U)\ \to \ B(U) \text{\quad and\quad}\\
    i(U)\cup B(\varphi)\ &: \ A(U)\cup_{A(W)} B(W)\ \to \ B(U)
    \end{align*}
    are $K$-cofibrations of $K$-spaces.
  \item For every flat orthogonal space $B$ the map 
    $B(\varphi): B(W)\to B(U)$ is a $K$-cofibration of $K$-spaces and the 
    $K$-space $B(U)$ is $K$-cofibrant. 
  \item Every flat orthogonal space is closed.\index{subject}{orthogonal space!closed}
  \end{enumerate}
\end{prop}
\begin{proof}
(i) The class of those morphisms of orthogonal spaces $i$
such that the map $i(U)\cup B(\varphi)$
is a $K$-cofibration of $K$-spaces is
closed under coproducts, cobase change, composition and
retracts. Similarly, the class of those morphisms of orthogonal spaces $i$
such that the map $i(U)$ is a $K$-cofibration of $K$-spaces is
closed under coproducts, cobase change, composition and
retracts. 
So it suffices to show each of the two claims for a set of generating cofibrations.
We do this for the  morphisms $\bL_{G,V}\times i_k$
for all $k\geq 0$, all compact Lie groups $G$ and all $G$-representations $V$,
where $i_k:\partial D^k\to D^k$ is the inclusion.
In this case the first map specializes to $\bL(V,U)/G\times i_k$.
The map $i_k$ is a cofibration and $\bL(V,U)/G$ is cofibrant as
a $K$-space by Proposition \ref{prop:K G cofibration}~(ii).
So $\bL(V,U)/G\times i_k$ is a $K$-cofibration of $K$-spaces.

The second map in question becomes the pushout product of
the sphere inclusion $i_k$ with the map
\[\bL(V,\varphi)/G \ : \ \bL(V,W)/G\ \to \  \bL(V,U)/G\ . \]
The map $i_k$ is a cofibration and $\bL(V,\varphi)/G$ is a $K$-cofibration
by Proposition \ref{prop:K G cofibration}~(i).
So their pushout product is again a $K$-cofibration.

Part~(ii) is the special case of part~(i) where $A=\emptyset$ 
is the empty orthogonal space.
Part~(iii) is the special case of~(ii) where $K$
is a trivial group, using that cofibrations of spaces are in particular 
h-cofibrations (Corollary \ref{cor-h-cofibration closures}) 
and hence closed embeddings (Proposition \ref{prop:h-cof is closed embedding}).
\end{proof}

Now we proceed towards the {\em global model structure} 
on the category of orthogonal spaces,
see Theorem \ref{thm:All global spaces}.
The weak equivalences in this model structure are the global equivalences
and the cofibrations are the flat cofibrations. The fibrations in the
global model structure are defined as follows.

\begin{defn}\index{subject}{global fibration!of orthogonal spaces}\index{subject}{fibration!global|see{global fibration}}
    A morphism $f:X\to Y$ of orthogonal spaces is a {\em global fibration}
    if it is a strong level fibration
    and for every compact Lie group $G$,
    every faithful $G$-representation $V$ 
    and every equivariant linear isometric embedding $\varphi:V\to W$ 
    of $G$-representations, the map
    \[ (f(V)^G, X(\varphi)^G)\ : \ X(V)^G \ \to \ 
    Y(V)^G \times_{Y(W)^G }  X(W)^G \]
    is a weak equivalence.
    
    An orthogonal space $X$ is {\em static}\index{subject}{static!orthogonal space}
    if for every compact Lie group $G$, every faithful $G$-representation $V$,
    and every $G$-equivariant linear isometric embedding $\varphi:V\to W$ 
    the structure map 
    \[ X(\varphi)\ :\  X(V)\ \to \ X(W) \]
    is a $G$-weak equivalence.
\end{defn}

Equivalently, a morphism $f$ is a global fibration if and only if $f$
is a strong level fibration and for every compact Lie group $G$,
every faithful $G$-representa\-tion $V$ and 
equivariant linear isometric embedding $\varphi:V\to W$ 
the square of $G$-fixed point spaces
    \begin{equation}  \begin{aligned}
        \label{eq:fibration characterization spaces}
        \xymatrix@C=18mm{ X(V)^G \ar[d]_{f(V)^G} \ar[r]^-{X(\varphi)^G} & 
           X(W)^G \ar[d]^{f(W)^G} \\
          Y(V)^G \ar[r]_-{Y(\varphi)^G} & Y(W)^G }
      \end{aligned}\end{equation}
    is homotopy cartesian.

Clearly, an orthogonal space $X$ is static if and only if the unique morphism
to a terminal orthogonal space is a global fibration;
the static orthogonal spaces will thus turn out to be the fibrant objects in the
global model structure.
The static orthogonal spaces are those which,
roughly speaking, don't change the equivariant homotopy type
once a faithful representation has been reached.

\begin{prop}\label{prop:global equiv preservation spc} 
  \begin{enumerate}[\em (i)]
   \item 
    Every global equivalence that is also a global fibration
    is a strong level equivalence.
  \item  
    Every global equivalence between static orthogonal spaces
    is a strong level equivalence.
\end{enumerate}
\end{prop}
\begin{proof}
(i) We let $f:X\to Y$ be a morphism of orthogonal spaces 
that is both a global fibration 
and a global equivalence. We consider a compact Lie group $G$, 
a faithful $G$-representation $V$, a finite $G$-CW-pair $(B,A)$ and a commutative square:
\[ \xymatrix{
A\ar[r]^-\alpha\ar[d]_{\text{incl}} & X(V) \ar[d]^{f(V)} \\
B\ar[r]_-\beta & Y(V) } \]
We will exhibit a continuous $G$-map $\mu:B\to X(V)$ such that
$\mu|_A=\alpha$ and such that $f(V)\circ\mu$
is homotopic, relative $A$, to $\beta$. 
This shows that the map $f(V)$ is a $G$-weak equivalence,
so $f$ is a strong level equivalence.

Since $f$ is a global equivalence, there is 
a $G$-equivariant linear isometric embedding $\varphi:V\to W$ 
and a continuous $G$-map $\lambda:B\to X(W)$
such that $\lambda|_A=X(\varphi)\circ \alpha:A\to X(W)$ and
such that $f(W)\circ \lambda:B\to Y(W)$
is $G$-homotopic, relative to $A$, to $Y(\varphi)\circ \beta$.
Since $f$ is a strong level fibration, 
we can improve $\lambda$ into a continuous $G$-map $\lambda':B\to X(W)$
such that $\lambda'|_A=\lambda|_A=X(\varphi)\circ \alpha$
and such that $f(W)\circ \lambda'$ is equal to $Y(\varphi)\circ \beta$.

Since $f$ is a global fibration the $G$-map
\[ (f(V), X(\varphi))\ : \ X(V) \ \to \ 
          Y(V) \times_{Y(W)}  X(W) \]
is a $G$-weak equivalence.
So we can find a continuous $G$-map $\mu:B\to X(V)$ such that
$\mu|_A=\alpha$ and $(f(V), X(\varphi))\circ\mu$
is $G$-homotopic, relative $A$, to 
$(\beta,\lambda'):B\to Y(V) \times_{Y(W) }  X(W)$:
\[ \xymatrix@C=15mm{
A\ar[r]^-\alpha\ar[d]_{\text{incl}} & X(V)  \ar[d]^{(f(V), X(\varphi))} \\
B\ar[r]_-{(\beta,\lambda')}\ar@{-->}[ur]^-(.4)\mu &   
Y(V) \times_{Y(W) }  X(W)  } \]
This is the desired map.

(ii) We let $f:X\to Y$ be a global equivalence between static orthogonal spaces.
We let $G$ be a compact Lie group, $V$ a faithful $G$-representation,
$(B,A)$ a finite $G$-CW-pair 
and $\alpha:A\to X(V)$ and $\beta:B\to Y(V)$ continuous $G$-maps 
such that $f(V)\circ\alpha=\beta|_A$.
Since $f$ is a global equivalence, there is 
a $G$-equivariant linear isometric embedding $\varphi:V\to W$ and
a continuous $G$-map $\lambda:B\to X(W)$ such that
$\lambda|_A=X(\varphi)\circ\alpha$ and $f(W)\circ\lambda$ is $G$-homotopic to
$Y(\varphi)\circ\beta$ relative $A$. 
Since $X$ is static, the map $X(\varphi):X(V)\to X(W)$
is a $G$-weak equivalence, so there is a continuous $G$-map $\bar\lambda:B\to X(V)$
such that $\bar\lambda|_A=\alpha$ and $X(\varphi)\circ\bar\lambda$ 
is $G$-homotopic to $\lambda$ relative $A$.
The two $G$-maps $f(V)\circ \bar\lambda$ and $\beta:B\to Y(V)$
then agree on $A$ and become $G$-homotopic, relative $A$,
after composition with $Y(\varphi):Y(V)\to Y(W)$.
Since $Y$ is static, the map $Y(\varphi)$ is a $G$-weak equivalence,
so $f(V)\circ \bar\lambda$ and $\beta:B\to Y(V)$
are already $G$-homotopic relative $A$.
This shows that $f(V):X(V)\to Y(V)$ is a $G$-weak equivalence,
and hence $f$ is a strong level equivalence.
\end{proof}

\begin{construction}\label{con:define Z(j)}
We let $j:A\to B$ be a morphism in a topological model category.
We factor $j$ through the mapping cylinder as the composite
\[ A\ \xra{\ c(j)\ } \ Z(j) =
(A\times [0,1])\cup_j B \ \xra{\ r(j)\ } \  B \ ,\]
where $c(j)$ is the `front' mapping  cylinder inclusion and $r(j)$ 
is the projection, which is a homotopy equivalence.
In our applications we will assume that both $A$ and $B$ are cofibrant,
and then the morphism $c(j)$ is a cofibration by the pushout product property.
We then define $\Zc(j)$\index{symbol}{$\Zc(j)$ - {set of pushout products of cylinder inclusions}}
 as the set of all pushout product maps
\[  c(j)\Box i_k \ : \   A\times D^k \cup_{A\times \partial D^k} Z(j)\times \partial D^k
\ \to \ Z(j)\times D^k \]
for $k\geq 0$, where $i_k:\partial D^k\to D^k$ is the inclusion.
\end{construction}

\begin{prop}\label{prop:hocartesian via RLP} 
Let $\Cc$ be a topological model category, $j:A\to B$
a morphism between cofibrant objects and $f:X\to Y$ a fibration.
Then the following two conditions are equivalent:
\begin{enumerate}[\em (i)]
\item The square of spaces
\begin{equation}\begin{aligned}\label{eq:hocart square}
 \xymatrix@C=15mm{ 
\map(B,X)\ar[r]^-{\map(j,X)} \ar[d]_{\map(B,f)} &
\map(A,X) \ar[d]^{\map(A,f)} \\
\map(B,Y)\ar[r]_-{\map(j,Y)} & \map(A,Y) }       
\end{aligned}\end{equation}
is homotopy cartesian.
\item The morphism $f$ has the right lifting property with respect to the
set $\Zc(j)$.
\end{enumerate}
\end{prop}
\begin{proof}
The square \eqref{eq:hocart square} maps to the square
\begin{equation}\begin{aligned}\label{eq:hocart square replaced}
 \xymatrix@C=20mm{ 
\map(Z(j),X)\ar[r]^-{\map(c(j),X)} \ar[d]_{\map(Z(j),f)} &
\map(A,X) \ar[d]^{\map(A,f)} \\
\map(Z(j),Y)\ar[r]_-{\map(c(j),Y)} & \map(A,Y) } 
\end{aligned}\end{equation}
via the map induced by $r(j):Z(j)\to B$ on the left part
and the identity on the right part. Since $r(j)$ is a homotopy equivalence,
the map of squares is a weak equivalence at all four corners.
So the square \eqref{eq:hocart square} is homotopy cartesian if and only if
the square \eqref{eq:hocart square replaced} is homotopy cartesian.

Since $A$ is cofibrant and $f$ a fibration, 
$\map(A,f)$ is a Serre fibration. So the 
square \eqref{eq:hocart square replaced} is homotopy cartesian
if and only if the map
\begin{align}  \label{eq:ppullback_map}
 (\map(Z(j),f),\map(c(j),X))\ &:  \\ 
\map(Z(j),X)\ \to \ &\map(Z(j),Y) \times_{\map(A,Y)}  \map(A,X)\nonumber
\end{align}
is a weak equivalence. Since $c(j)$ is a cofibration and $f$ is a fibration,
the map \eqref{eq:ppullback_map} is always a Serre fibration.
So \eqref{eq:ppullback_map} is a weak equivalence
if and only if it is an acyclic fibration,
which is equivalent to the right lifting property for
the inclusions $i_k:\partial D^k\to D^k$ for all $k\geq 0$.
By adjointness, the map \eqref{eq:ppullback_map} 
has the right lifting property with respect to the maps $i_k$
if and only if the morphism $f$ has the right lifting property with respect to the
set $\Zc(j)$.
\end{proof}

The set $J^{\str}$ was defined in 
the proof of Proposition \ref{prop:strong level spaces} as the set of morphisms $G_m j$ 
for $m\geq 0$ and for $j$ in the set of generating acyclic cofibrations 
for the projective model structure on the category of $O(m)$-spaces specified 
in \eqref{eq:J_for_F-proj_on_GT}.
The set $J^{\str}$ detects the fibrations 
in the strong level model structure. We add another set of morphisms $K$ 
that detects when the squares \eqref{eq:fibration characterization spaces}
are homotopy cartesian.
Given any compact Lie group $G$ and $G$-representations $V$ and $W$,
the restriction morphism
\[ \rho_{G,V,W} \ = \ \rho_{V,W}/G \ : \ \bL_{G,V\oplus W} \ \to \ \bL_{G,V}  \]
restricts (the $G$-orbit of) a linear isometric embedding
from $V\oplus W$ to $V$.
If the representation $V$ is faithful, 
then this morphism is a global equivalence 
by Proposition \ref{prop:free_orthogonal_space}~(ii).
We set
\[ K \ = \ \bigcup_{G,V,W} \Zc(\rho_{G,V,W}) \ ,\]
the set of all pushout products of boundary inclusions $\partial D^k\to D^k$
with the mapping cylinder inclusions of the morphisms $\rho_{G,V,W}$;
here the union is over a set of representatives
of the isomorphism classes of triples $(G,V,W)$ consisting of
a compact Lie group $G$, a faithful $G$-representation $V$ and
an arbitrary $G$-representation $W$.
The morphism $\rho_{G,V,W}$ represents the map of $G$-fixed point spaces
$X(i_{V,W})^G:X(V)^G\to X(V\oplus W)^G$;
every $G$-equivariant linear isometric embedding is isomorphic
to a direct summand inclusion $i_{V,W}$,
so by Proposition \ref{prop:hocartesian via RLP}, the right lifting property 
with respect to the union $J^{\str}\cup K$ 
characterizes the global fibrations, i.e., we have shown:

\begin{prop}\label{prop:global fibrations via RLP}
A morphism of orthogonal spaces is a global fibration if and only if it has
the right lifting property with respect to the set $J^{\str}\cup K$.  
\end{prop}

Now we are ready for the main result of this section.

\begin{theorem}[Global model structure]\label{thm:All global spaces} 
    The global equivalences, global fibrations and flat cofibrations 
    form a model structure,
    the {\em global model structure}\index{subject}{global model structure!for orthogonal spaces}\index{subject}{model structure!global|see{global model structure}}
    on the category of orthogonal spaces. 
    The fibrant objects in the global model structure are the static orthogonal spaces.
    The global model structure is proper, topological and cofibrantly generated.
\end{theorem}
\begin{proof}
We number the model category axioms as in \cite[3.3]{dwyer-spalinski}. 
The category of orthogonal spaces is complete and cocomplete,
so axiom MC1 holds.
Global equivalences satisfy the 2-out-of-6 property
by Proposition \ref{prop:global equiv basics}~(iii),
so they also satisfy the 2-out-of-3 property MC2.
Global equivalences are closed under retracts
by Proposition \ref{prop:global equiv basics}~(iv);
it is straightforward that cofibrations and global fibrations 
are closed under retracts, so axiom MC3 holds.

The strong level model structure shows that every morphism of orthogonal spaces
can be factored as $f\circ i$ for a flat cofibration $i$
followed by a strong level equivalence $f$ that is also a strong level fibration.
For every $G$-equivariant linear isometric embedding $\varphi:V\to W$
between faithful $G$-representations, 
 both vertical maps in the commutative square 
of fixed point spaces \eqref{eq:fibration characterization spaces}
are then weak equivalences, so the square is homotopy cartesian.
The morphism $f$ is thus a global fibration and a global equivalence,
so this provides one of the factorizations as required by MC5.
For the other half of the factorization axiom MC5
we apply the small object argument\index{subject}{small object argument} 
(see for example \cite[7.12]{dwyer-spalinski} or \cite[Thm.\,2.1.14]{hovey-book})
to the set $J^{\str}\cup K$.
All morphisms in $J^{\str}$ are flat cofibrations and strong level equivalences.
Since $\bL_{G,V\oplus W}$ and $\bL_{G,V}$ are flat, the morphisms in $K$
are also flat cofibrations, and they are global equivalences
because the morphisms $\rho_{G,V,W}$ are 
(Proposition \ref{prop:free_orthogonal_space}~(ii)).
The small object argument provides a functorial factorization
of every morphism $X\to Y$ of orthogonal spaces as a composite
\[ X \ \xra{\ i \ }\ W \ \xra{\ q \ }\ Y \]
where $i$ is a sequential composition of cobase changes of coproducts
of morphisms in $J^{\str}\cup K$, and $q$ has the right lifting property with respect 
to $J^{\str}\cup K$. 
Since all morphisms in $J^{\str}\cup K$ are flat cofibrations
and global equivalences, the morphism $i$ is a
flat cofibration and a global equivalence by the closure properties
of Proposition \ref{prop:global equiv basics}.
Moreover, $q$ is a global fibration by Proposition~\ref{prop:global fibrations via RLP}.

Now we show the lifting properties MC4. 
By Proposition \ref{prop:global equiv preservation spc}~(i)
a morphism that is both a global equivalence and a global fibration
is a strong level equivalence, and hence an acyclic fibration
in the strong level model structure. So every morphism that is
simultaneously a global equivalence and  a global fibration has the
right lifting property with respect to flat cofibrations.
Now we let $j:A\to B$ be a flat cofibration that is also a global equivalence and 
we show that it has the left lifting property with respect to all global fibrations.
We factor $j=q\circ i$, 
via the small object argument for $J^{\str}\cup K$,
where $i:A\to W$ is a $(J^{\str}\cup K)$-cell complex and $q:W\to B$ a global fibration.
Then $q$ is a global equivalence since $j$ and $i$ are,
and hence an acyclic fibration in the strong level model structure,
again by Proposition \ref{prop:global equiv preservation spc}~(i).
Since $j$ is a flat cofibration, a lifting in
\[\xymatrix{
A \ar[r]^-i \ar[d]_j & W \ar[d]^q_(.6)\sim \\
B \ar@{=}[r] \ar@{..>}[ur] & B }\]
exists. Thus $j$ is a retract of the morphism $i$ that has the left lifting property
with respect to global fibrations. But then $j$ itself has this lifting property.
This finishes the verification of the model category axioms.
Alongside we have also specified sets of generating flat cofibrations $I^{\str}$
and generating acyclic cofibrations $J^{\str}\cup K$.
Sources and targets of all morphisms in these sets are small with
respect to sequential colimits of flat cofibrations. So the global
model structure is cofibrantly generated.

Left properness of the global model structure follows from
Proposition \ref{prop:global equiv basics} (xi) and the fact that
flat cofibrations are h-cofibrations 
(Corollary \ref{cor-h-cofibration closures}~(iii)).
Right properness follows from
Proposition \ref{prop:global equiv basics}~(xii) because
global fibrations are in particular strong level fibrations.

The global model structure is topological by 
Proposition \ref{prop:topological criterion},
with $\Gc$ the set of semifree orthogonal spaces $\bL_{G,V}$
indexed by a set of representatives $(G,V)$ of the isomorphism classes
of pairs consisting of a compact Lie group $G$ and a faithful $G$-representation $V$,
and with $\Zc$ the set of mapping cylinder inclusions $c(\rho_{G,V,W})$
of the morphisms $\rho_{G,V,W}$.
\end{proof}

The global model structure of orthogonal spaces is also monoidal,
even with respect to two different monoidal structures.
Indeed, the categorical product of orthogonal spaces has the pushout product
property for flat cofibrations, 
by Proposition \ref{prop:ppp for cartesian product} below.
Moreover,  Proposition \ref{prop:ExF ppp spaces}~(iii)
(for the global family of all compact Lie groups)
shows that  global model structure of orthogonal spaces 
satisfies the pushout product property
with respect to the box product of orthogonal spaces.

\medskip

We also introduce a `positive' version of the global model structure for
orthogonal spaces; our main use of this variation is for the
global model structure of ultra-commutative monoids 
in Section \ref{sec:global model monoid spaces}.
As is well known from similar contexts (for example, the stable model
structure for commutative orthogonal ring spectra),
model structures cannot usually be lifted naively to 
multiplicative objects with strictly commutative products.
The solution is to lift a `positive'
version of the global model structure in which the values at the
trivial inner product space
are homotopically meaningless and where the fibrant objects
are the `positively static' orthogonal spaces.

\begin{defn} A morphism $f:A\to B$ of orthogonal spaces is a {\em positive cofibration}\index{subject}{cofibration!positive}\index{subject}{positive cofibration!of orthogonal spaces}
if it is a flat cofibration and the map $f(0):A(0)\to B(0)$ is a homeomorphism.  
An orthogonal space $Y$ 
is {\em positively static}\index{subject}{positively static}\index{subject}{static!positively}
if for every compact Lie group $G$, every faithful $G$-representation $V$
with $V\ne 0$ and every $G$-equivariant linear isometric embedding $\varphi:V\to W$  
the structure map 
\[ Y(\varphi)\ :\  Y(V)\ \to \ Y(W) \]
is a $G$-weak equivalence.
\end{defn}

If $G$ is a non-trivial compact Lie group, then any faithful $G$-representation
is automatically non-trivial. So a positively static orthogonal space
is static (in the absolute sense) if the structure map $Y(0)\to Y(\mR)$
is a non-equivariant weak equivalence.

\begin{prop}[Positive global model structure]\label{prop:positive global spaces} 
  The global equivalences and positive cofibrations are part of 
  a cofibrantly generated, proper, topological model structure,
  the {\em positive global model structure}\index{subject}{positive global model structure!for orthogonal spaces}\index{subject}{global model structure!positive}
  on the category of orthogonal spaces. 
  A morphism $f:X\to Y$ of orthogonal spaces is a fibration 
  in the positive global model structure if and only if
  for every compact Lie group $G$,
  every faithful $G$-representation $V$ with $V\ne 0$ 
  and every equivariant linear isometric embedding $\varphi:V\to W$  
  the map $f(V)^G:X(V)^G\to Y(V)^G$ is a Serre fibration and
  the square of $G$-fixed point spaces
    \[  
        \xymatrix@C=18mm{ X(V)^G \ar[d]_{f(V)^G} \ar[r]^-{X(\varphi)^G} & 
           X(W)^G \ar[d]^{f(W)^G} \\
          Y(V)^G \ar[r]_-{Y(\varphi)^G} & Y(W)^G }
    \]
    is homotopy cartesian.  
    The fibrant objects in the positive global model structure 
    are the positively static orthogonal spaces.
\end{prop}
\begin{proof}
  We start by establishing a {\em positive strong level model structure}.
  We call a morphism $f:X\to Y$ of orthogonal spaces a
  {\em positive strong level equivalence}\index{subject}{level equivalence!positive strong!of orthogonal spaces} 
  (respectively {\em positive strong level fibration})\index{subject}{level fibration!positive strong!of orthogonal spaces}
  if for every inner product space $V$ with $V\ne 0$
  the map $f(V):X(V)\to Y(V)$  is an $O(V)$-weak equivalence
  (respectively an $O(V)$-fibration).
  Then we claim that the positive strong level equivalences, 
  positive strong level fibrations and positive cofibrations 
  form a topological model structure on the category of orthogonal spaces.

  The proof is another application of the general construction method
  for level model structures in 
  Proposition \ref{prop:general level model structure}.
  Indeed, we let $\Cc(0)$ be the degenerate model structure on 
  the category $\bT$ of spaces in which every morphism
  is a weak equivalence and a fibration, but only the isomorphisms are cofibrations. 
  For $m\geq 1$ we let $\Cc(m)$ be the projective
  model structure (for the set of all closed subgroups)
  on the category of $O(m)$-spaces,
  compare Proposition \ref{prop:proj model structures for G-spaces}.
  With respect to these choices of model structures,
  the classes of level equivalences, level fibrations and cofibrations
  in the sense of Proposition \ref{prop:general level model structure} 
  become the positive strong level equivalences,
  positive strong level fibrations and positive cofibrations.
  The consistency condition 
  (Definition \ref{def:consistency condition})
  is now strictly weaker than for the strong level model structure,
  so it holds.
  
  We obtain the positive global model structure
  for orthogonal spaces by `mixing' the 
  positive strong level model structure 
  with the global model structure of Theorem \ref{thm:All global spaces}.
  Every positive strong level equivalence is a global equivalence and every 
  positive cofibration is a flat cofibration.
  The global equivalences and the positive cofibrations
  are part of a model structure by Cole's mixing theorem \cite[Thm.\,2.1]{cole-mixed},
  which is our first claim. 
  By \cite[Cor.\,3.7]{cole-mixed} (or rather its dual formulation),
  an orthogonal space is fibrant in the positive global model structure 
  if and only if it is weakly equivalent in the positive strong level model structure to 
  a static orthogonal space; this is equivalent to
  being positively static.

  Cofibrant generation, properness and topologicalness 
  of the positive global model structure are proved in much
  the same way as for the absolute global model structure
  in Theorem \ref{thm:All global spaces}.
\end{proof}

\begin{rk}
We can relate the unstable global homotopy theory
of orthogonal spaces to the homotopy theory of
$G$-spaces for a fixed compact Lie group $G$. 
Then evaluation at a faithful $G$-representation $V$ 
and the semifree functor at $(G,V)$ are a pair of adjoint functors
\[ \xymatrix{ \bL_{G,V}\ : \ G\bT \ \ar@<.4ex>[r] & 
  \ \spc \ : \ \ev_{G,V} \ar@<.4ex>[l] } \]
between the categories of $G$-spaces and orthogonal spaces.
This adjoint pair is a Quillen pair
with respect to the global model structure of orthogonal spaces and
the `genuine' model structure of $G$-spaces
(i.e., the projective model structure with respect to the family
of all subgroups, 
compare Proposition \ref{prop:proj model structures for G-spaces}).
The adjoint total derived functors
\[ \xymatrix{ 
  L(\bL_{G,V}) \ : \ \Ho( G\bT ) \ \ar@<.4ex>[r]   & 
  \ \Ho( \spc)\ : \ R(\ev_{G,V})  \ar@<.4ex>[l]
} \] 
are independent of the faithful representation $V$ 
up to preferred natural isomorphism, by 
Proposition \ref{prop:free_orthogonal_space}~(ii).

Every $G$-space is $G$-weakly equivalent to a $G$-CW-complex,
and these are built from the orbits $G/H$.
So the derived left adjoint $L(\bL_{G,V}):G\bT\to\spc$ 
is essentially determined by its values on the coset spaces $G/H$.
Since $\bL_{G,V}(G/H)$ is isomorphic to $\bL_{H,V}=B_{\gl}H$, 
the derived left adjoint takes the homogeneous space $G/H$
to a global classifying space of $H$.

The derived right adjoint also has a more explicit description,
at least for closed orthogonal spaces $Y$, as the underlying $G$-space $Y(\Uc_G)$.
Indeed, we can choose a fibrant replacement of $Y$ 
in the global model structure, i.e., a flat cofibration $j:Y\to Z$
that is also a global equivalence, and such that $Z$ is
globally fibrant (i.e., static).
Then $Z$ is also closed, and so the induced map
$j(\Uc_G):Y(\Uc_G)\to Z(\Uc_G)$ is a $G$-weak equivalence by
Proposition \ref{prop:global eq for closed}.
We may assume that $V$ is a subrepresentation of $\Uc_G$;
we choose a nested sequence 
\[ V = V_1 \ \subset \ V_2 \ \subset \ \dots \ \subset \ V_n \ \subset \ \dots \]
of finite-dimensional $G$-subrepresentations that exhaust $\Uc_G$.
Since $V$ is faithful and $Z$ is closed and static, the induced maps
\[ Z(V) = Z(V_1)\ \to \ Z(V_2) \ \to \ \dots \ \to \ Z(V_n) \ \to \ \dots \]
are all closed embeddings and $G$-weak equivalences.
So the canonical map
\[ Z(V) \ \to \ \colim_{n\geq 1}\, Z(V_n) \ = \ Z(\Uc_G)  \]
is also a $G$-weak equivalence. Since $Z$ is a globally fibrant
replacement of $Y$, the $G$-space $Z(V)$ calculates the right derived functor
of $\ev_{G,V}$ at $Y$. This exhibits a chain of two $G$-weak equivalences
\[ R(\ev_{G,V})(Y)\ = \ Z(V)\ \xra{\ \simeq\ } \ Z(\Uc_G) 
\ \xla{\ \simeq\ } \ Y(\Uc_G) \ . \]
\end{rk}

\begin{construction}[Cofree orthogonal spaces]
\label{con:cofree_orthogonal_space}\index{subject}{cofree orthogonal space}\index{subject}{orthogonal space!cofree|see{cofree orthogonal space}}
For every compact Lie group $K$, we will now define
a right adjoint to the functor that takes an orthogonal space $Y$
to the underlying $K$-space $Y(\Uc_K)$.
We refer to the right adjoint $R_K$ as the {\em cofree functor}.\index{symbol}{$R_K(A)$ - {cofree orthogonal space of a $K$-space $A$}}
We consider the continuous functor
\[ \bL(-,\Uc_K) \ : \ \bL \ \to \ (K\bT)^{\op} \ , \quad V \ \longmapsto \ 
\bL(V,\Uc_K)\ ,\]
with functoriality by precomposition with linear isometric embeddings.
The group $K$ acts on the values of this functor through the action on
the complete universe $\Uc_K$.
The cofree orthogonal space $R_K(A)$ associated to a $K$-space $A$ is then the composite
\[  \bL \ \xra{\ \bL(-,\Uc_K)} \ K\bT^{\op} \ \xra{\map^K(-,A)} \ \bT\ .\]
The unit of the adjunction is the morphism
\begin{equation}  \label{eq:cofree_unit}
 \eta_Y \ : \ Y \ \to \ R_K (Y(\Uc_K))   
\end{equation}
whose value at an inner product space $V$ is the adjoint of the action map
\[\bL(V,\Uc_K)\times  Y(V) \ \to \ Y(\Uc_K) \ , \quad
(\varphi,y)\ \longmapsto \ Y(\varphi)(y)\ .\]
The counit of the adjunction is the continuous $K$-map
\[  \epsilon_A \ : \ R_K(A)(\Uc_K)   \ \to \ A \]
assembled from the compatible $K$-maps
\[ R_K(A)(V)\ = \ \map^K(\bL(V,\Uc_K),A)\ \to \ A \ ,\quad
f\longmapsto f(i_V)\ , \]
for $V\in s(\Uc_K)$, where $i_V:V\to\Uc_K$ is the inclusion.
This data makes the functors
\[\xymatrix{ 
(-)(\Uc_K)\ : \ \spc \quad \ar@<.4ex>[r]  & 
\quad K\bT \ : \ R_K \ar@<.4ex>[l] }\]
into an adjoint pair.
\end{construction}

\begin{prop}\label{prop:cofree properties}
Let $K$ be a compact Lie group.  
\begin{enumerate}[\em (i)]
\item The adjoint functor pair $((-)(\Uc_K),R_K)$ is a Quillen pair for
the global model structure of orthogonal spaces and the projective model
structure of $K$-spaces.
\item For every $K$-space $A$ the orthogonal space $R_K(A)$ is static.
\item
For every closed orthogonal space $Y$ the map
\[ (\eta_Y)^*\circ R_K\ : \ \Ho(K\bT)(Y(\Uc_K),A) \ \to \  \Ho(\spc)(Y, R_K(A))  \]
is bijective.
\end{enumerate}
\end{prop}
\begin{proof}
(i)
We let $f:X\to Y$ be a fibration of $K$-spaces.
We let $G$ be another compact Lie group and $V$ a faithful $G$-representation.
Then the $K$-space $\bL(V,\Uc_K)/G$ is $K$-cofibrant
by Proposition \ref{prop:K G cofibration}~(ii).
The projective model structure on $K$-spaces is topological, 
so $\map^K(\bL(V,\Uc_K)/G,-)$ takes fibrations of $K$-spaces
to fibrations of spaces. Because
\[ \map^K(\bL(V,\Uc_K)/G,X) \ = \ (R_K(X)(V))^G\ ,\]
this means that $R_K$ takes fibrations of $K$-spaces to strong level fibrations
of orthogonal spaces. By the same argument, 
$R_K$ takes acyclic fibrations of $K$-spaces to acyclic fibrations
in the strong level model structure, which coincide with the
acyclic fibrations in the global model structure of orthogonal spaces.

Now we let $\varphi:V\to W$ be a $G$-equivariant linear isometric embedding.
Then the map
\[ \rho_{V,W}(\Uc_K)/G\ : \ \bL(V\oplus W,\Uc_K)/G \ \to \ \bL(V,\Uc_K)/G \]
is a $K$-homotopy equivalence by Proposition \ref{prop:free_orthogonal_space}~(ii).
So the induced map 
\[ (R_K(X)(\varphi))^G\ : \ (R_K(X)(V))^G\ \to \ (R_K(X)(V\oplus W))^G \]
is a homotopy equivalence. So in the commutative square
\[ \xymatrix@C=20mm{ 
(R_K(X)(V))^G\ar[r]^-{R_K(X)(\varphi)^G} \ar[d]_{(R_K(f)(V))^G} & 
(R_K(X)(V\oplus W))^G \ar[d]^{(R_K(f)(V\oplus W))^G} \\
(R_K(Y)(V))^G\ar[r]_-{R_K(Y)(\varphi)^G} & (R_K(Y)(V\oplus W))^G } \]
both vertical maps are Serre fibrations and both horizontal maps
are weak equivalences. 
The square is then homotopy cartesian, and so the morphism
$R_K(f):R_K(X)\to R_K(Y)$ is a global fibration of orthogonal spaces.
Altogether this shows that the right adjoint $R_K$ preserves fibrations and 
acyclic fibrations, so $((-)(\Uc_K),R_K)$ is a Quillen pair.

(ii) Every $K$-space $A$ is fibrant in the projective model structure.
So $R_K(A)$ is fibrant in the global model structure of orthogonal spaces;
by Theorem \ref{thm:All global spaces} these fibrant object
are precisely the static orthogonal spaces.

(iii)
We choose a global equivalence $f:Y^c\to Y$ with flat source.
Then $Y^c$ and $Y$ are both closed, the former by Proposition \ref{prop:flat becomes K-cofibration}~(iii). 
So the map $f(\Uc_K):Y^c(\Uc_K)\to Y(\Uc_K)$
is a $K$-weak equivalence by Proposition \ref{prop:global eq for closed}.
So the morphism $f$ induces bijections on both sides of the map in question,
and it suffices to prove the claim for $Y^c$ instead of $Y$.
But $Y^c$ is cofibrant and $A$ is fibrant, so in this case the claim
is just the derived adjunction isomorphism.
\end{proof}

For $K=e$ the trivial group we drop the subscript
and abbreviate the cofree functor $R_e$ to $R$. 

\begin{defn}
An orthogonal space $Y$ is {\em cofree}\index{subject}{cofree orthogonal space} 
if it is globally equivalent to an orthogonal space of the form $R A$
for some space $A$.
\end{defn}

We will now develop criteria for detecting cofree orthogonal spaces,
and then recall some non-tautological examples. 
One criterion involves the unit of the adjunction, 
the special case 
\[ \eta_Y \ : \ Y \ \to \ R (Y(\mR^\infty)) \]
of \eqref{eq:cofree_unit} for $\Uc_K=\mR^\infty$. The next proposition shows that the 
morphism $\eta_Y$ is always a non-equivariant weak equivalence,
provided $Y$ is closed.

\begin{prop}\label{prop:eta_Y is weak equivalence}
For every closed orthogonal space $Y$ the morphism   
$\eta_Y :Y \to R (Y(\mR^\infty))$ induces a weak equivalence
\[ \eta_Y(\mR^\infty) \ : Y(\mR^\infty) \ \to\  R (Y(\mR^\infty))(\mR^\infty) \]
on underlying non-equivariant spaces.
\end{prop}
\begin{proof}
We start with a general observation about cofree orthogonal spaces.
Since $R A$ is static (Proposition \ref{prop:cofree properties}~(ii)) and closed,
the map $(R A)(\varphi):(R A)(V)\to (R A)(W)$ 
induced by any linear isometric embedding $\varphi:V\to W$
is a weak equivalence and a closed embedding.
So the canonical map
\begin{align}  \label{eq:canonical_map_RA}
 A \ \iso \ \map(\bL(0,\mR^\infty),A)\ = \ &(R A)(0) \\ 
\to \ &\colim_{W\in s(\mR^\infty)}  (R A)(W) \ = \ (R A)(\mR^\infty)   \nonumber
\end{align}
is a weak equivalence as well.
The adjunction counit $\epsilon_A : (R A)(\mR^\infty) \to  A$
is a retraction to the map \eqref{eq:canonical_map_RA}, 
so $\epsilon_A$ is also a weak equivalence. 

Now we turn to the proof of the proposition.
Even though the map $\eta_Y(\mR^\infty)$ under consideration is {\em not} the same 
as the canonical map \eqref{eq:canonical_map_RA} for $A=Y(\mR^\infty)$,
the counit $\epsilon_{Y(\mR^\infty)}:R(Y(\mR^\infty))(\mR^\infty)\to Y(\mR^\infty)$ 
is also a retraction to $\eta_Y(\mR^\infty)$. 
Since $\epsilon_{Y(\mR^\infty)}$ is a weak equivalence, so is $\eta_Y(\mR^\infty)$. 
\end{proof}

While the morphism $\eta_Y:Y\to R(Y(\mR^\infty))$ 
tends to be a non-equivariant equivalence,
it is typically {\em not} a global equivalence. 
We will now see that for a closed orthogonal space $Y$ 
the morphism $\eta_Y$ is a global equivalence if and only if $Y$ is cofree.

We recall that a {\em universal free $K$-space},
for a compact Lie group $K$, is a $K$-cofibrant free $K$-space 
whose underlying space in non-equivariantly contractible.
Any two universal free $K$-spaces are $K$-homotopy equivalent,
see Proposition \ref{prop:universal spaces}.
We call a $K$-space $A$ {\em cofree}\index{subject}{cofree equivariant space} 
if for some (hence any) universal free $K$-space $E K$ the map
\[ \text{const}\ : \  A\ \to \ \map(E K, A) \]
that sends a point to the corresponding constant map is a $K$-weak equivalence.

\begin{prop}\label{prop:cofree criteria}
For a closed orthogonal space $Y$ the following three conditions are equivalent.
\begin{enumerate}[\em (i)]
\item The orthogonal space $Y$ is cofree.
\item For every compact Lie group $K$ the $K$-space
$Y(\Uc_K)$ is cofree.
\item The adjunction unit $\eta_Y:Y\to R(Y(\mR^\infty))$ is a global equivalence.
\end{enumerate}
\end{prop}
\begin{proof}
In a first step we show that for every space $A$ and every compact Lie group $K$,
the $K$-space $(R A)(\Uc_K)$ is cofree.
We choose a faithful $K$-representa\-tion $W$.
Then $\bL(W,\mR^\infty)$ is a universal 
free $K$-space by Proposition \ref{prop:free_orthogonal_space}~(i).
So the projection from $E K\times \bL(W,\mR^\infty)$ to the second factor
is a $K$-weak equivalence between cofibrant $K$-spaces, hence
a $K$-homotopy equivalence. So the induced map
\begin{align*}
  \text{const}\ : \ (R A)(W) \ = \ &\map(\bL(W,\mR^\infty), A) \\ 
\to \ &\map(E K\times \bL(W,\mR^\infty), A) \ \iso \ 
\map(E K,  (R A)(W)) 
\end{align*}
is a $K$-homotopy equivalence. Hence the $K$-space $(R A)(W)$
is cofree as soon as $K$ acts faithfully on $W$.
Since  $R A$ is static (by Proposition \ref{prop:cofree properties}~(ii))
and closed, the canonical map
\[ (R A)(W) \ \to \ (R A)(\Uc_K) \]
is a $K$-weak equivalence. So $(R A)(\Uc_K)$ is $K$-cofree.
Now we prove the equivalence of conditions~(i), (ii) and~(iii).

(i)$\Longrightarrow$(ii)
The global equivalences are part of the global model structure on
the category of orthogonal spaces, compare Theorem \ref{thm:All global spaces}. 
Moreover, the orthogonal space $R A$ is static, hence fibrant in the global model
structure.
So if $Y$ is globally equivalent to $R A$, then for some (hence any)
global equivalence $p:Y^c\to Y$ with cofibrant (i.e., flat)
source, there is a global equivalence $f:Y^c\to R A$.

Now we let $K$ be any compact Lie group.
The orthogonal space $Y^c$ is closed by 
Proposition \ref{prop:flat becomes K-cofibration}~(iii). 
Since $Y$ and $R A$ are also closed, the global equivalences induce
$K$-weak equivalences
\[ Y(\Uc_K)\ \xla[\sim]{\ p(\Uc_K)}\ Y^c(\Uc_K)\ \xra[\sim]{\ f(\Uc_K)}\ (R A)(\Uc_K)\]
by Proposition \ref{prop:global eq for closed}.
Since $(R A)(\Uc_K)$ is $K$-cofree by the introductory remark, so is $Y(\Uc_K)$.

(ii)$\Longrightarrow$(iii)
We start with a preliminary observation.
We let $Y$ and $Z$ be two closed orthogonal spaces
such that the $K$-spaces $Y(\Uc_K)$ and $Z(\Uc_K)$ are cofree
for all compact Lie groups $K$. We claim that every morphism $f:Y\to Z$
of orthogonal spaces such that $f(\mR^\infty):Y(\mR^\infty)\to Z(\mR^\infty)$
is a non-equivariant weak equivalence is already a global equivalence.
Indeed, for every compact Lie group $K$ the two vertical maps in 
the commutative square of $K$-spaces
\[ \xymatrix@C=25mm{ 
Y(\Uc_K)\ar[d]_{\text{const}}\ar[r]^-{f(\Uc_K)} &
Z(\Uc_K)\ar[d]^{\text{const}}\\
\map(E K, Y(\Uc_K))\ar[r]_-{\map(E K, f(\Uc_K))}^-\simeq &
\map(E K, Z(\Uc_K)) }  \]
are $K$-weak equivalences by hypothesis.
Since $\Uc_K$ is non-equivariantly isometrically isomorphic to $\mR^\infty$,
the $K$-map $f(\Uc_K):Y(\Uc_K)\to Z(\Uc_K)$ is a non-equivariant weak equivalence
by hypothesis. So the lower horizontal map is a $K$-weak equivalence.
We conclude that the upper  horizontal map is a $K$-weak equivalence.
Since $Y$ and $Z$ are closed, 
the criterion of Proposition \ref{prop:global eq for closed} 
shows that $f$ is a global equivalence.

Now we apply the criterion to the morphism $\eta_Y:Y\to R(Y(\mR^\infty))$. 
The map $\eta_Y(\mR^\infty)$ is a weak equivalence
by Proposition \ref{prop:eta_Y is weak equivalence}. 
Moreover, for every compact Lie group $K$
the space $Y(\Uc_K)$ is $K$-cofree by hypothesis~(ii),
and $R(Y(\mR^\infty))(\Uc_K)$ is $K$-cofree by the introductory remark.
The criterion of the previous paragraph thus applies and shows that the morphism
$\eta_Y:Y\to R(Y(\mR^\infty))$ is a global equivalence.

Condition~(i) is a special case of~(iii).
\end{proof}

We recall now that the global classifying spaces of certain compact Lie groups
are cofree, namely of those with abelian identity path component.
Said differently, the group must be an extension of a finite group by a torus.
The following theorem is a reinterpretation of the main result
of Rezk's paper \cite{rezk-1truncated},
who calls these groups `1-truncated' 
because the homotopy groups of the underlying spaces
vanish in dimensions larger than~1.
The special case of {\em abelian} compact Lie groups was proved
earlier by Lashof, May and Segal \cite{lashof-may-segal}.
The case of finite groups seems to be folklore,
going all the way back to Hurewicz \cite{hurewicz-beitraege IV} who proved 
that for finite groups $G$, homotopy classes of continuous maps $B K\to B G$ 
biject with conjugacy classes of continuous group homomorphisms
from $K$ to $G$.

\begin{theorem}\label{thm:BA is right induced} 
Let $G$ be a compact Lie group whose identity path component is abelian.
Then the global classifying space $B_{\gl} G$ is cofree.\index{subject}{global classifying space!of an abelian compact Lie group}
\end{theorem}
\begin{proof}
We let $V$ be any faithful $G$-representation, so that $B_{\gl}G=\bL_{G,V}$. 
Then $\bL_{G,V}(\Uc_K)=\bL(V,\Uc_K)/G$ is a classifying space
for principal $(K,G)$-bundles, by Proposition \ref{prop:free_orthogonal_space}.
Since $G$ has abelian identity component, it is 1-truncated
in the sense of Rezk, and so the $K$-space $\bL_{G,V}(\Uc_K)$ is cofree by 
\cite[Thm.\,1.4]{rezk-1truncated}.
So criterion~(ii) of Proposition \ref{prop:cofree criteria} is satisfied;
since the orthogonal space $\bL_{G,V}$ is closed, we have thus shown that it is cofree.
\end{proof}

\begin{rk}\index{subject}{global classifying space}
The global classifying space $B_{\gl} G$ is {\em not} cofree in general,
i.e., when the identity component of $G$ is not abelian.
Indeed, for another compact Lie group $K$, the homotopy set
$\pi_0^K(B_{\gl} G)$ (to be introduced in Definition \ref{def:homotopy set} below)
bijects with conjugacy classes of continuous homomorphisms from $K$ to $G$, 
by Proposition \ref{prop:fix of global classifying}~(ii).
On the other hand, the set $\pi_0^K(R(B G))$ bijects 
with homotopy classes of continuous maps from $B K$ to $B G$.
However, there are continuous maps $B K\to B G$ that are not homotopic to $B\alpha$
for any continuous homomorphism $\alpha:K\to G$.
Whenever this happens, the adjunction unit $\eta_{B_{\gl} G}:B_{\gl} G\to R(B G)$ 
is not surjective on $\pi_0^K$, and hence not a global equivalence 
(by Corollary \ref{cor:pi_0^G of closed} below).
The first examples of such `exotic' maps between classifying spaces
of compact Lie groups were constructed by Sullivan
and appeared in his widely circulated and highly influential
MIT lecture notes;
an edited version of Sullivan's notes was eventually published
in \cite{sullivan-localization reprint}.
Indeed, Corollary~5.11 of \cite{sullivan-localization reprint}
constructs `unstable Adams operations' $\psi^p:B U(n)\to B U(n)$
for a prime $p$ and all $n<p$;
for $n>1$ these maps are not induced by any continuous homomorphism.
\end{rk}

\begin{construction}[Realization of simplicial objects]\index{subject}{realization!of simplicial objects}\label{con:realize simplicial}
We will occasionally want to realize simplicial objects, so we quickly
recall the necessary background. We let $\Cc$ be a cocomplete category
tensored over the category $\bT$ of spaces. 
We let $\bDelta$\index{symbol}{$\bDelta$ - {simplicial index category}} denote
the simplicial indexing category, with objects the finite
totally ordered sets $[n]=\{0\leq 1\leq \dots\leq n\}$
for $n\geq 0$. Morphisms in $\bDelta$ are all weakly monotone maps.
We let\index{symbol}{$\Delta^n$ - {topological $n$-simplex}}
\[ \Delta^n \ = \ \{ (t_1,\dots,t_n)\in [0,1]^n \ |\ t_1\leq t_2\leq \dots\leq t_n \} \]
be the topological $n$-simplex.
As $n$ varies, these topological simplices assemble into a covariant functor
\[ \bDelta \ \to \bT \ , \quad [n]\ \longmapsto \ \Delta^n\ ; \]
the coface maps are given by
\[ (d_i)_*(t_1,\dots,t_n) \ = \left\lbrace \begin{array}{ll}
\quad (0,t_1,\dots,t_n) & \mbox{for $i=0$,} \\
(t_1,\dots,t_i,t_i,\dots,t_n)  & \mbox{for $0<i< n$,} \\
\quad (t_1,\dots,t_n,1) & \mbox{for $i=n$.}
\end{array} \right.  \]
For $0\leq i \leq n-1$, the codegeneracy map $(s_i)_*:\Delta^n\to\Delta^{n-1}$
drops the entry~$t_{i+1}$.

A {\em simplicial object}\index{subject}{simplicial object}
in $\Cc$ is functor $X:\bDelta^{\op}\to\Cc$, i.e., a contravariant functor from $\bDelta$.
We use the customary notation $X_n=X([n])$ for the value of a simplicial object at $[n]$. 
The {\em realization} of $X$ is the coend
\[ |X|\ = \ \int^{[n]\in\Delta} X_n\times \Delta^n \]
of the functor
\[ \bDelta^{\op}\times\bDelta\ \to \ \Cc\ , \quad 
([m],[n])\ \longmapsto \ X_m\times\Delta^n \ .\]
We also need to recall the {\em latching objects} of a simplicial object.
We let $\bDelta(n)$ denote category with objects the weakly 
monotone surjections $\sigma:[n]\to [k]$;
a morphism from $\sigma:[n]\to [k]$ to $\sigma':[n]\to [k']$
is a morphism $\alpha:[k]\to [k']$ in $\bDelta$ 
(necessarily surjective as well)
with $\alpha\circ\sigma=\sigma'$.
We let $\bDelta(n)_\circ$ denote the full subcategory of $\bDelta(n)$
with all objects {\em except} the identity of $[n]$.
A simplicial object $X$ can be restricted along the forgetful functor
\[ \bDelta(n)^{\op}_\circ\ \xra{\ u \ } \ \bDelta^{\op} \ ,\quad 
(\sigma:[n]\to [k])\ \longmapsto \ [k]\ .\]
The $n$-th {\em latching object} of $X$ is the colimit over $\bDelta(n)^{\op}_\circ$
of the restricted functor:
\begin{equation}  \label{eq:simplicial_latching}
 L_n^\bDelta(X) \ = \ \colim_{\bDelta(n)^{\op}_\circ} \, (X\circ u)\ .  
\end{equation}
The morphisms
\[ \sigma^* \ : \ X_k\ = \ (X\circ u)(\sigma:[n]\to [k])\ \to \ X_n \]
assemble into a 
{\em latching morphism}\index{subject}{latching morphism!of a simplicial object}
\[ l_n\ : \ L_n^\bDelta(X) \ \to \ X_n \]
from the $n$-th latching object of $X$ to the value at $[n]$.
For example, the category $\bDelta(0)_\circ$ is empty, so $L_0^\bDelta(X)$
is an initial object of $\Cc$. 
The category $\bDelta(1)_\circ$ has a unique object $s_0:[1]\to [0]$,
so $L_1^\bDelta(X) = X_0$ and the latching morphism is given by $s_0^*:X_0\to X_1$.
The category $\bDelta(2)_\circ$ has three objects and two non-identity morphisms,
and $L_2^\bDelta(X)$ is a pushout of the diagram
\[ X_1 \ \xla{\ s_0^* \ }\ X_0 \ \xra{\ s_0^*\ } X_1 \ .\]
\end{construction}

We specialize the above to realizations of simplicial orthogonal spaces,
i.e., simplicial objects in the category of orthogonal spaces.
For orthogonal spaces, coends and product with $\Delta^n$ 
are objectwise, and hence 
\[ |X|(V) \  = \ |X (V)|\ ,\]
i.e., the value of $|X|$ at an inner product space $V$ is 
the realization of the simplicial space $[n]\mapsto X_n(V)$,
as discussed in Construction \ref{con:realize simplicial space}.
By Proposition \ref{prop:geometric realization in bT}, the realization
can be formed in the ambient category of all topological spaces,
and the result is automatically compactly generated.

\begin{defn}\label{def:Reedy flat}
A simplicial orthogonal space $X$ is {\em Reedy flat}\index{subject}{Reedy flat!orthogonal space}
if the latching morphism $l_n:L_n^\bDelta(X)\to X_n$ is a flat cofibration
of orthogonal spaces for every $n\geq 0$.
\end{defn}

The terminology `Reedy flat' stems from Reedy's theorem \cite{reedy}
that the simplicial objects in any model category admit a certain model
structure, nowadays called the `Reedy model structure', 
in which the equivalences are the level equivalences of simplicial objects. 
Reedy's paper -- albeit highly influential -- remains unpublished, but an account of 
the Reedy model structure can 
for example be found in \cite[VII Prop.\,2.11]{goerss-jardine}.
If we form the Reedy model structure starting
with the global model structure of orthogonal spaces,
then the cofibrant objects are precisely the 
Reedy flat simplicial orthogonal spaces.

The most important consequence of Reedy flatness for our purposes is that for these
simplicial orthogonal spaces, realization is homotopical.
Indeed, the global model structure of orthogonal spaces is
a topological model structure by Theorem \ref{thm:All global spaces}.
We can turn this into a simplicial model structure by defining
the tensor of an orthogonal space $X$ with a simplicial set $A$ as
\[ X\tensor A \ = \ X\times |A| \ ,\]
the objectwise product with the geometric realization of $A$.
The following proposition is then a special case of the fact
that realization is a left Quillen functor for the Reedy model
structure on simplicial orthogonal spaces, see \cite[VII Prop.\,3.6]{goerss-jardine}.

\begin{prop}\label{prop:realization invariance in spc}\index{subject}{realization!of simplicial orthogonal spaces}
  \begin{enumerate}[\em (i)]
  \item 
    The realization of every Reedy flat simplicial orthogonal space is flat.
  \item Let $f:X\to Y$ be a morphism of Reedy flat
    simplicial orthogonal spaces. If $f_n:X_n\to Y_n$ is a global equivalence
    for every $n\geq 0$, then the morphism of realizations
    $|f|:|X|\to |Y|$ is a global equivalence.
  \end{enumerate}
\end{prop}

\section{Monoidal structures}
\label{sec:unstable monoidal}

This section is devoted to monoidal products on the category of orthogonal spaces, 
with emphasis on global homotopical features.
Our main focus is the box product of orthogonal spaces,
a special case of a Day type convolution product,
and the `good' monoidal structure for orthogonal spaces.
We prove in Theorem \ref{thm:box to times} that  the box product 
is fully homotopical with respect to global equivalences.
While the box product is the most useful monoidal
structure for orthogonal spaces, the cartesian product
is also relevant for our purposes.
We show in Proposition \ref{prop:ppp for cartesian product}
that the categorical product, too, satisfies the pushout product property
for flat cofibrations.
In particular, the product of two flat orthogonal spaces is again flat.

The last part of this section introduces complex analogues $\bL^\mC_{G,W}$ 
of the semi\-free orthogonal spaces,
indexed by unitary $G$-representations $W$,
see Construction \ref{con:complex-free spc}.
While these complex versions are not (semi)free in any categorical
sense, they are similar to the semifree orthogonal spaces
in many ways; for example, the orthogonal spaces $\bL^\mC_{G,W}$ 
are flat (Proposition \ref{prop:complex global classifying}~(ii)) 
and behave well multiplicatively under box and cartesian product
(Proposition \ref{prop:box complex almost-representables}).

\medskip

We define a {\em bimorphism} 
$b:(X,Y)\to Z$\index{subject}{bimorphism!of orthogonal spaces}
from a pair of orthogonal spaces $(X,Y)$
to another orthogonal space $Z$ as a collection of continuous maps 
\[ b_{V,W} \ : \ X(V) \times  Y(W) \ \to \ Z(V\oplus W),\ \]
for all inner product spaces $V$ and $W$, such that for all
linear isometric embeddings $\varphi:V\to V'$ and $\psi:W\to W'$ the following
square commutes:
\[ \xymatrix@C=20mm{ X(V) \times Y(W) \ar[r]^-{b_{V,W}} \ar[d]_{X(\varphi)\times Y(\psi)} &
Z(V\oplus W) \ar[d]^{Z(\varphi\oplus\psi)} \\
X(V') \times Y(W') \ar[r]_-{b_{{V'},{W'}}} &  Z(V'\oplus W') } \]
We define a box product\index{subject}{box product!of orthogonal spaces|(} 
of $X$ and $Y$ as a universal example
of an orthogonal space with a bimorphism from $X$ and $Y$.
More precisely, a box product is a pair $(X\boxtimes Y,i)$ 
consisting of an orthogonal space $X\boxtimes Y$
and a universal bimorphism $i:(X,Y)\to X\boxtimes Y$,
i.e., such that for every orthogonal space $Z$ the map\index{symbol}{$\boxtimes$ - {box product of orthogonal spaces}}
\[  \spc(X\boxtimes Y,Z) \ \to \ \text{Bimor}((X,Y),Z) \ , \quad
f\longmapsto f i = \{f(V\oplus W)\circ i_{V,W}\}_{V,W} \]
is bijective. 
We will often refer to this bijection as the
{\em universal property} of the box product of orthogonal spaces.
\index{subject}{universal property!of the box product}
Very often only the object $X\boxtimes Y$ will be referred to 
as the box product, but one should keep in mind that it comes equipped
with a specific, universal bimorphism.

The existence of a universal bimorphism out of any
pair of orthogonal spaces $X$ and $Y$, and thus of a box product $X\boxtimes Y$, 
is a special case of the existence of Day type convolution products on
certain functor categories; the construction is an enriched
Kan extension of the `pointwise' cartesian product of $X$ and $Y$
along the direct sum functor $\oplus:\bL\times\bL\to\bL$
(see Proposition \ref{prop:box exists}), or more explicitly an enriched
coend (see Remark \ref{rk:box as enriched coend}).

Also by the general theory of convolution products, 
the box product $X\boxtimes Y$ is a functor in both variables
(Construction \ref{con:box functoriality}) and it supports a preferred
symmetric monoidal structure (see Theorem \ref{thm:symmetric monoidal});
so there are specific natural associativity and symmetry isomorphisms
\[  (X\boxtimes Y)\boxtimes Z \ \to \ X\boxtimes(Y\boxtimes Z)
\text{\qquad respectively \qquad} 
 X\boxtimes Y \ \to \ Y\boxtimes X  \]
and a strict unit, the terminal orthogonal space $\mathbf 1$,
i.e., such that $\mathbf 1\boxtimes X= X= X\boxtimes\mathbf 1$.
The upshot is that the associativity and symmetry isomorphisms 
make the box product
of orthogonal spaces into a symmetric monoidal product with 
the terminal orthogonal space as unit object. 
The box product of orthogonal spaces is {\em closed} symmetric monoidal
in the sense that the box product is adjoint to an
internal Hom orthogonal space. We won't use the internal function
object, so we do not elaborate on it.

\medskip

The next result proves a key feature, namely that up to global equivalence,
the box product of orthogonal spaces coincides with the categorical product.
Given two orthogonal spaces $X$ and $Y$, the maps
\[ X(V)\times Y(W) \ \xra{X(i_1)\times Y(i_2)} \
X(V\oplus W)\times Y(V\oplus W) \ = \ (X\times Y)(V\oplus W) \]
form a bimorphism $(X,Y)\to X\times Y$ as $V$ and $W$ vary through 
all inner product spaces; here $i_1:V\to V\oplus W$ and
$i_2:W\to V\oplus W$ are the two direct summand embeddings. 
This bimorphism is represented by a morphism
\begin{equation}  \label{eq:define_box2times}
 \rho_{X,Y} \ : \ X\boxtimes Y\ \to \ X\times Y   
\end{equation}
of orthogonal spaces that is natural in both variables.

\begin{theorem}\label{thm:box to times}
Let $X$ and $Y$ be orthogonal spaces.
\begin{enumerate}[\em (i)]
\item 
The morphism 
$\rho_{X,Y}:X\boxtimes Y\to X\times Y$ is a global equivalence. 
\item The functor $X\boxtimes -$ preserves global equivalences.
\end{enumerate}
\end{theorem}
\begin{proof}
(i) For an orthogonal space $Z$ we denote by $\sh Z$ the orthogonal space defined by
\[ (\sh Z)(V) \ = \ Z( V\oplus V) \text{\qquad and\qquad}
(\sh Z)(\varphi)\ = \ Z(\varphi\oplus\varphi)\ ; \]
thus $\sh Z$ is isomorphic to $\sh^{\mR^2}_\tensor(Z)$, 
the multiplicative shift\index{subject}{shift!of an orthogonal space!multiplicative}
of $Z$ by $\mR^2$ as defined 
in Example \ref{eg:Additive and multiplicative shift}.
We define a morphism of orthogonal spaces
\[ \lambda\ : \ X\times Y \ \to \ \sh(X\boxtimes Y) \]
at an inner product space $V$ as the composite
\[ X(V)\times Y(V)\ \xra{i_{V,V}} \ 
(X\boxtimes Y)(V\oplus V) \ = \ (\sh (X\boxtimes Y))(V) \ . \]
Now we consider the two composites $\lambda\circ\rho_{X,Y}$ 
and $\sh(\rho_{X,Y})\circ\lambda$:
\[\xymatrix@C=12mm{  
X\boxtimes Y\ar[r]^-{\rho_{X,Y}}  &   X\times Y \ar[r]^-{\lambda} &
 \sh(X\boxtimes Y) \ar[r]^-{\sh(\rho_{X,Y})} & \sh( X\times Y) }\]
We claim that the composite $\lambda\circ\rho_{X,Y}:X\boxtimes Y\to\sh( X\boxtimes Y)$
is homotopic to the morphism $(X\boxtimes Y)\circ i_1$, where
$i_1$ is the natural linear isometric embedding $V\to V\oplus V$
as the first summand.
Indeed, for every $t\in [0,1]$ we define a natural linear isometric embedding
\[ j_t \ : \ V\oplus W \ \to \ V\oplus W\oplus V\oplus W\text{\quad by\quad}
j_t(v,w)\ = \ (v,\ t \cdot w,\ 0,\ \sqrt{1-t^2}\cdot  w)\ . \]
Then the maps
\begin{align*}
  X(V)\times Y(W)\ \xra{i_{V,W}}\  &(X\boxtimes Y)(V\oplus W)\\   
\xra{(X\boxtimes Y)(j_t)} \  &(X\boxtimes Y)(V\oplus W\oplus V\oplus W) 
\ = \ (\sh (X\boxtimes Y))(V\oplus W)    
\end{align*}
form a bimorphism as $V$ and $W$ vary; 
the universal property of the box product turns this into a morphism
\[ f_t \ : \ X\boxtimes Y \ \to \ \sh (X\boxtimes Y) \]
of orthogonal spaces.
The linear isometric embeddings $j_t$ vary continuously
with $t$, hence so do the morphisms $f_t$.
Moreover, $f_0=\lambda\circ\rho_{X,Y}$ and $f_1=(X\boxtimes Y)\circ i_1$, so this is
the desired homotopy.
The morphism $(X\boxtimes Y)\circ i_1$ is a global equivalence by
Theorem \ref{thm:general shift of osp}, hence so is the morphism $\lambda\circ\rho_{X,Y}$.

The shift functor preserves products, and under the canonical isomorphism
$\sh(X\times Y)\iso (\sh X)\times (\sh Y)$ the morphism
$\sh(\rho_{X,Y})\circ\lambda$ becomes the product of the two morphisms
\[ X\circ i_1\ : \ X\ \to\ \sh X \text{\qquad and\qquad}
Y\circ i_2\ : \ Y\ \to\ \sh Y \ .\]
The morphisms $X\circ i_1$ and $Y\circ i_2$ are global equivalences by
Theorem \ref{thm:general shift of osp}, hence so is their product
(by Proposition \ref{prop:global equiv basics}~(vi)).
The global equivalences satisfy the 2-out-of-6 property
by Proposition \ref{prop:global equiv basics}~(iii);
since $\lambda\circ\rho_{X,Y}$ and $\sh(\rho_{X,Y})\circ\lambda$ are
global equivalences, so is the morphism $\rho_{X,Y}$.

(ii) The cartesian product $X\times- $ preserves global equivalences 
by Proposition \ref{prop:global equiv basics}~(vi).
Together with part~(i) this implies part~(ii).
\end{proof}

\begin{eg}[Box product of semifree orthogonal spaces]\label{eg:box of free orthogonal}
We show that the box product of two semifree orthogonal spaces is another semifree
orthogonal space. This can be deduced from the general fact that a convolution
product of two representable functors is again representable 
(see Remark \ref{rk:box representable}); however, the argument
is simple enough that we make it explicit for orthogonal spaces.

We consider two compact Lie groups $G$ and $K$,
representations $V$ and $W$ of $G$ respectively $K$, a $G$-space $A$ and
a $K$-space $B$.
Then $V\oplus W$ is a $(G \times K)$-representation via 
\[ (g,k)\cdot (v,w)\ = \ (g v,\, k w)\ , \]
and $A\times B$ is a $(G \times K)$-space in much the same way.
The map
\begin{align*}
A\times B \ \xra{ (\Id_V\cdot G,-)\times (\Id_W\cdot K,- )} \  
 &( \bL_{G,V} A)(V)\times (\bL_{K,W} B) (W)\\ 
\xra{\quad i_{V,W}\quad}\quad &( ( \bL_{G,V} A)\boxtimes (\bL_{K,W} B) )(V\oplus W)
  \end{align*}
is $(G\times K)$-equivariant, so it extends freely to a morphism of orthogonal spaces
\begin{equation}  \label{eq:box_of_free}
\bL_{G\times K,V\oplus W}(A\times B) \ \to \ (\bL_{G,V} A)\boxtimes (\bL_{K,W} B) \ .
\end{equation}
The maps
\begin{align*}
 (\bL(V,U)\times_G A) \times (\bL(W,U')\times_K B) \ &\to \ 
\bL(V\oplus W,U\oplus U')\times_{G\times K}(A \times  B) \\
( [\varphi,a],\, [\psi,b]) \qquad\qquad &\longmapsto \qquad 
 [\varphi\oplus\psi,(a,b)] 
\end{align*}
form a bimorphism from $(\bL_{G,V}A,\bL_{K,W}B)$ to $\bL_{G\times K,V\oplus W}(A\times B)$
as the inner product spaces $U$ and $U'$ vary. The universal property of the box
product translates this into a morphism
\[  ( \bL_{G,V} A)\boxtimes (\bL_{K,W} B)\ \to \ \bL_{G\times K,V\oplus W}(A\times B)\ . \]
These two morphisms are mutually inverse isomorphisms, i.e., the box product
$ (\bL_{G,V} A)\boxtimes (\bL_{K,W} B)$ is isomorphic to $ \bL_{G\times K,V\times W}(A\times B)$.
A special case of this shows that the box product of two global classifying
spaces is another global classifying space.\index{subject}{global classifying space}
Indeed, if $G$ acts faithfully on $V$ and $K$ acts faithfully on $W$, then
the $(G\times K)$-action on $V\oplus W$ is also faithful, and hence
\begin{equation}  \label{eq:boxtimes_of_B_gl}
 ( B_{\gl}G) \boxtimes (B_{\gl} K) \ = \ 
 \bL_{G,V}\boxtimes \bL_{K,W} \ \iso \ \bL_{G\times K,V\oplus W}\ = \ B_{\gl}(G\times K)\ .
\end{equation}
If we compose the inverse of the isomorphism \eqref{eq:boxtimes_of_B_gl}
with the global equivalence $\rho_{\bL_{G,V}, \bL_{K, W}}$
from \eqref{eq:define_box2times}, we obtain a global 
equivalence of orthogonal spaces
\[  B_{\gl}(G\times K)\ \xra{\ \sim\ } \ ( B_{\gl}G) \times (B_{\gl} K) \ .\]
\end{eg}
\index{subject}{box product!of orthogonal spaces|)}

Now we show that the categorical product of orthogonal spaces has the pushout product
property for flat cofibrations,
see Proposition \ref{prop:ppp for cartesian product} below.
In particular, the product of two flat orthogonal spaces
is again flat, which is not completely obvious from the outset.
To this end we establish a useful sufficient condition for flatness;
the criterion is inspired by a flatness criterion for $\bI$-spaces
proved by Sagave and Schlichtkrull in \cite[Prop.\,3.11]{sagave-schlichtkrull-diagram}.
The conditions may seem technical at first sight,
but we give two examples where they are easily verified,
see Propositions \ref{prop:rho flat for free} 
and \ref{prop:complex global classifying}.
The linear isometries category $\bL$ does not have very many limits;
however, for property (b) of the next proposition we note that it does have
pullbacks. 

\begin{prop}\label{prop:spc flatness criterion}
Let $Y$ be an orthogonal space.
\begin{enumerate}[\em (i)]
\item 
Suppose that $Y$ satisfies the following conditions.
\begin{enumerate}[\em (a)]
\item For every inner product space $V$, the space $Y(V)$ is compact. 
\item As a functor from the category $\bL$ to sets, $Y$ preserves
pullbacks.
\end{enumerate}
Then for all $m\geq 0$ the canonical morphism  $i_m:\sk^m Y\to Y$ 
from the $m$-skeleton is objectwise a closed embedding.
\item Let $G$ be a compact Lie group acting continuously on $Y$
through automorphisms of orthogonal spaces. 
Suppose that in addition to the conditions {\em (a)} and {\em (b)},
the following also hold for all inner product spaces $V$:
\begin{enumerate}[\em (a)]
\setcounter{enumii}{2}
\item The $(G\times O(V))$-space $Y(V)$ admits a $(G\times O(V))$-CW-structure.
\item If an element of $Y(V)$ is fixed by a reflection in $O(V)$, then
it is in the image of $Y(U)$ for some proper subspace $U$ of $V$.
\end{enumerate}
Then the orthogonal orbit space $Y/G$ is flat.
\end{enumerate}
\end{prop}
\begin{proof}
We let $V$ be any inner product space and  consider the composite
\begin{equation}\label{eq:composite flat criterion}
 {\coprod}_{0\leq i \leq  m} \bL(\mR^i,V)\times Y(\mR^i) \ \to \
(\sk^m  Y)(V)\ \xra{i_m(V)} \ Y(V)  \ ,
\end{equation}
where the first map is from the definition of $\sk^m Y$
as a coequalizer \eqref{eq:skeleton_coequalizer spc}.
The first map in \eqref{eq:composite flat criterion}
is a quotient map, and its source is compact by hypothesis~(a)
and because the spaces $\bL(\mR^i,V)$ are all compact.
So the space $(\sk^m Y)(V)$ is quasi-compact. 
Since $Y(V)$ is Hausdorff, the map $i_m(V)$ is a closed map.
So all that remains is to show that $i_m(V)$ is injective.

We consider two pairs
\[ (\varphi, x) \ \in \  
 \bL(\mR^i,V)\times Y(\mR^i)    \text{\qquad and\qquad} 
(\varphi', x')\ \in \   \bL(\mR^j,V)\times Y(\mR^j)  \]
with $i, j\leq m$, such that
\[  Y(\varphi)(x) \ = \ Y(\varphi')(x')  \ \in \ Y(V) \ .  \]
We show that the two pairs represent the same element in the
coequalizer $(\sk^m Y)(V)$.
We choose a pullback square in the category $\bL$:
\[ \xymatrix{\mR^k \ar[r]^-{h'}\ar[d]_h & \mR^j\ar[d]^{\varphi'}\\
\mR^i\ar[r]_-{\varphi} & V} \]
Condition (b) provides
an element $z \in \ Y(\mR^k)$ such that $Y(h)(z)=x$ and $Y(h')(z)=x'$.
Then
\begin{align*}
  (\varphi,x) &= (\varphi, Y(h)(z)) \sim (\varphi h, z) =  (\varphi' h' ,z) \sim
 (\varphi', Y(h')(z)) = (\varphi', x')\ .
\end{align*}
So the pairs $(\varphi,x)$ and $(\varphi',x')$ 
describe the same element in $(\sk^m Y)(V)$.

Now we let $G$ act on $Y$ and also assume conditions~(c) and~(d).
By the first part, the latching morphism
\[ \nu_m^Y\ = \ i_{m-1}(\mR^m)\ : \ L_m(Y)\ = \ (\sk^{m-1}Y)(\mR^m)\ \to \ Y(\mR^m) \]
is a closed embedding for every $m\geq 0$,
and the space $Y(\mR^m)$ admits a $(G\times O(m))$-CW-structure by hypothesis~(c).
We claim that $\nu_m^Y$ is a $(G\times O(m))$-cofibration;
to this end we characterize its image by a stabilizer condition.
An element of $y\in Y(\mR^m)$ is in the image of the latching map 
if an only if it is in the image of the map $Y(U)\to Y(\mR^m)$
for some $(m-1)$-dimensional subspace $U$ of $\mR^m$. 
Then the orthogonal reflection in the hyperplane $U$ fixes $y$.
Conversely, if $y$ is fixed by a reflection,
then it is in the image of $Y(U)$ for some proper subspace $U$, by hypothesis~(d).
So the image of the latching morphism coincides with the subspace 
of all elements of $Y(\mR^m)$
whose stabilizer group contains the reflection in some hyperplane of $\mR^m$.

If $K$ is any compact Lie group and $X$ a $K$-CW-complex, then
the elements of $X$ whose stabilizer group contains a conjugate
of some fixed subgroup is automatically a $K$-CW-subcomplex.
In the situation at hand this means that the image of
the latching map $\nu_m^Y$ is a $(G\times O(m))$-subcomplex
in any $(G\times O(m))$-CW-structure on $Y(\mR^m)$.
So the latching map is a $(G\times O(m))$-cofibration.

The latching space construction commutes with colimits,
so the canonical map
\[ (L_m Y)/G\ \to \  L_m(Y/G) \]
is a homeomorphism, and the latching map for $Y/G$
is obtained from the latching map for $Y$ by passage to $G$-orbits.
Since the latching map for $Y$ 
is a $(G\times O(m))$-cofibration, 
the induced map on $G$-orbits is then an $O(m)$-cofibration by  
Proposition \ref{prop:cofibrancy preservers}~(iii).
This shows that the orthogonal space $Y/G$ is flat.
\end{proof}

We apply the flatness criterion to show that the
product of two semifree orthogonal spaces is flat.

\begin{prop}\label{prop:rho flat for free}
Let $G$ and $K$ be compact Lie groups and let $V$ and $W$ 
be faithful representations of $G$ respectively $K$. Then the
orthogonal space $\bL_{G,V}\times \bL_{K,W}$ is flat.
\end{prop}
\begin{proof}
We verify the conditions~(a)-(d) of Proposition \ref{prop:spc flatness criterion}
for the orthogonal space $\bL_V\times\bL_W$ with the action of $(G\times K)$
by precomposition on linear isometric embeddings.
The space $\bL(V,U)$ is compact and a Stiefel manifold,
so it comes with a `standard' smooth structure; 
moreover, the action of $G\times O(U)$ by pre- respectively postcomposition is smooth. 
By the same argument, $\bL(W,U)$ admits the structure of
a smooth closed $(K\times O(U))$-manifold.
Hence the product $\bL(V,U)\times\bL(W,U)$
underlies a closed smooth $(G\times K\times O(U))$-manifold.
So Illman's theorem \cite[Cor.\,7.2]{illman} provides
a finite $(G\times K\times O(U))$-CW-structure.
This verifies conditions~(a) and~(c).

Conditions~(b) is straightforward. Finally, if a pair
\[ (\varphi,\psi ) \ \in\ \bL(V,U)\times \bL(W,U) \]
is fixed by the reflection in some hyperplane $U'$ of $U$,
then the images of $\varphi$ and $\psi$ are contained in $U'$, 
and hence $(\varphi,\psi)$ lies in the image of
$\bL(V,U')\times \bL(W,U')$. This verifies condition~(d).
Proposition \ref{prop:spc flatness criterion} thus applies, and
shows that the orthogonal space $(\bL_V\times\bL_W)/(G\times K)$ is flat.
So the isomorphic orthogonal space $ \bL_{G,V}\times\bL_{K,W}$ is flat as well.
\end{proof}

Since the semifree orthogonal spaces $\bL_{G,V}$ `generate' the flat cofibrations,
the previous Proposition \ref{prop:rho flat for free}
is the key input for the pushout product property for
the categorical product of orthogonal spaces.

\begin{prop} \label{prop:ppp for cartesian product}
Let $f:A\to B$ and $g:X\to Y$ be flat cofibrations of orthogonal spaces.
Then the pushout product morphism
\[ f\Box g\ = \ (f\times Y)\cup(B\times g) \ : \ 
A\times Y\cup_{A\times X}B\times X \ \to \ B\times Y \]
is a flat cofibration.
In particular, the product of two flat orthogonal spaces is again flat.
\end{prop}
\begin{proof} 
Since the cartesian product preserves colimits in both variables,
it suffices to show the claim for two generating flat cofibrations
  \[   \bL_{G,V}\times i_k\ : \ \bL_{G,V}\times  \partial D^k 
  \ \to \  \bL_{G,V}  \times D^k  \]
  and $\bL_{K,W}\times i_m$, where $G$ and $K$ are compact Lie groups,
  $V$ and $W$ are faithful representations of $G$ respectively $K$,
  and $k,m\geq 0$.
  The pushout product $(\bL_{G,V}\times i_k)\Box(\bL_{K,W}\times i_m)$ 
  of two such generators is isomorphic to the morphism
  \[ \bL_{G,V} \times \bL_{K,W}\times i_{k+m} \ . \]
  This morphism is a flat cofibration 
  since the strong level model structure of orthogonal spaces is topological
  and $\bL_{G,V}\times \bL_{K,W}$ is flat (Proposition \ref{prop:rho flat for free}).
\end{proof}

In Definition \ref{def:global classifying}
we introduced global classifying spaces as the semifree orthogonal spaces
defined from faithful orthogonal representations of compact Lie groups. 
As we shall explain in Proposition \ref{prop:complex global classifying} below, 
we can also use faithful {\em unitary} representations instead. 
This extra flexibility will come in handy
when we study global objects with an intrinsic complex flavor,
such as global classifying spaces of unitary groups,
complex Grassmannians, or complex Bott periodicity.
The next construction introduces the orthogonal space
$\bL^\mC_{G,W}$ for a unitary $G$-representation $W$.
In contrast to their orthogonal cousins $\bL_{G,V}$,
the unitary analogs are not representable nor semifree in any sense.
However, the unitary versions also enjoy various useful properties,
for example that they are flat 
(see Proposition \ref{prop:complex global classifying}~(ii))
and behave well 
under box product (see Proposition \ref{prop:box complex almost-representables}).

\begin{construction}\label{con:complex-free spc}
To define the orthogonal space $\bL^\mC_{G,W}$ for a unitary $G$-representation $W$,
we introduce notation for going back and forth 
between euclidean inner product spaces over $\mR$ and hermitian inner
product spaces over $\mC$.
Throughout, we shall denote euclidean inner products on $\mR$-vector spaces
by pointy brackets $\td{-,-}$, and
hermitian inner products on $\mC$-vector spaces
by round parentheses $(-,-)$.
For a euclidean inner product space $V$ we let 
\[ V_\mC \ = \ \mC\tensor_\mR V \]
be the complexification;
\index{symbol}{$V_\mC$ - {complexification of the inner product space $V$}}
the euclidean inner product $\td{-,-}$ on $V$ induces a hermitian
inner product $(-,-)$ on the complexification $V_\mC$, defined
as the unique sesquilinear form that satisfies
\[ (1\tensor v, 1\tensor v') \ = \ \td{v,v'} \]
for all $v, v'\in V$.
For a hermitian inner product space $W$, we let $u W$ 
denote the underlying $\mR$-vector space,
equipped with the euclidean inner product
\[ \td{w, w'}\ = \ \text{Re} ( w ,w' ) \ ,\]
the real part of the given hermitian inner product. 
Every $\mC$-linear isometric embedding is in particular
an $\mR$-linear isometric embedding of underlying euclidean vector spaces;
so $U(W)\subseteq O(u W)$, i.e., the unitary group of $W$
is a subgroup of the orthogonal group of $u W$. We thus view a unitary
representation on $W$ as an orthogonal representation on $u W$.
If $V$ and $W$ are two finite-dimensional $\mC$-vector spaces equipped
with hermitian inner products, we denote by $\bL^\mC(V,W)$ 
the space of $\mC$-linear isometric embeddings.\index{symbol}{$\bL^\mC(V,W)$ - {space of $\mC$-linear isometric embeddings}}  
We topologize this as a complex Stiefel manifold, i.e., homeomorphic 
to the space of complex $\dim_\mC(V)$-frames in $W$.

Now we can define complex analogs of semifree orthogonal spaces.
We let $G$ be a compact Lie group and $W$ a finite-dimensional 
unitary $G$-representation. We define the orthogonal space $\bL^\mC_{G,W}$ by
\[ \bL^\mC_{G,W}(V) \ = \ \ \bL^\mC(W,V_\mC) / G\ .\]
We define a morphism of orthogonal spaces
\[ f_{G,W} \ : \ \bL_{G, u W} \ \to \ \bL^\mC_{G, W}\]
as follows.
The map
\[ j_W \ : \ W \ \to \  \mC\tensor_\mR W \ = \ ( u W)_\mC\ , \quad 
j_W(w) \ = \ 1/\sqrt{2}\cdot ( 1\tensor w - i\tensor(i w)  ) \]
is a $G$-equivariant $\mC$-linear isometric embedding. 
At a real inner product space $V$, we can thus define
\[ f_{G,W}(V)\ : \  \bL(u W,V)/ G  \ \to \ \bL^\mC(W,V_\mC)/G \text{\quad by\quad}
f_{G,W}(V)(\varphi G)\ = \ ( \varphi_\mC\circ j_W) G \ .\]
\end{construction}

\begin{prop}\label{prop:complex global classifying}
Let $G$ be a compact Lie group and $W$ a faithful unitary $G$-representation.
\begin{enumerate}[\em (i)]
\item 
The morphism $f_{G,W}:\bL_{G, u W} \to \bL^\mC_{G, W}$ is a global equivalence.
\item
The orthogonal space $\bL^\mC_{G,W}$ is flat.  
\end{enumerate}
\end{prop}
\begin{proof}
(i) Both source and target of $f$ are closed orthogonal spaces; 
so by Proposition \ref{prop:global eq for closed} 
we may show that for every compact Lie group $K$ the map
\[ f_{G,W}(\Uc_K)\ : \  \bL(u W,\Uc_K)/ G  \ \to \ 
\bL^\mC(W,\mC\tensor_\mR \Uc_K) / G \]
is a $K$-weak equivalence.
We consider the $(K\times G)$-equivariant continuous map
\[ \tilde f\ : \  \bL(u W,\Uc_K)  \ \to \ \bL^{\mC}(W ,\mC\tensor_\mR\Uc_K)
\ , \quad \varphi\  \longmapsto \  \varphi_\mC\circ j_W \]
that `covers' $f_{G,W}(\Uc_K)$.
The source of $\tilde f$ is a universal $(K\times G)$-space
for the family $\Fc(K;G)$ of graph subgroups, 
by Proposition \ref{prop:free_orthogonal_space}~(i).
Since $\mC\tensor_\mR \Uc_K$ is a complete complex $G$-universe, the complex analog
of Proposition \ref{prop:free_orthogonal_space}~(i), proved in much the same way,
shows that the target of $\tilde f$ is also such a universal space for the
same family of subgroups of $K\times G$.
So $\tilde f$ is a $(K\times G)$-equivariant homotopy equivalence,
compare Proposition \ref{prop:universal spaces}.
The map $f_{G,W}(\Uc_K)=\tilde f/G$ induced on $G$-orbit spaces
is thus a $K$-equivariant homotopy equivalence.

(ii)
We verify conditions~(a)-(d) of Proposition \ref{prop:spc flatness criterion}
for the orthogonal $G$-space $\bL^\mC_W$.
The space $\bL^\mC(W,V_\mC)$ is compact and a complex Stiefel manifold,
so it comes with a `standard' smooth structure; 
moreover, the action of $G\times O(V)$ by pre- respectively postcomposition is smooth. 
So Illman's theorem \cite[Cor.\,7.2]{illman} provides 
a $(G\times O(V))$-CW-structure on $\bL^\mC_W(V)$.
This verifies conditions~(a) and~(c).

For condition~(b) we observe that complexification
takes pullback squares in $\bL$ to pullback squares of hermitian 
inner product spaces and complex linear isometric embeddings.
So the functor $\bL^\mC(W,(-)_\mC)$ preserves pullbacks.
Finally, if $\varphi\in\bL^\mC(W,V_\mC)$
is fixed by the reflection in some hyperplane $U$ of $V$,
then the image of $\varphi$ is contained in $U_\mC$, and hence $\varphi$
lies in the image of $\bL^\mC(W,U_\mC)$. This verifies condition~(d).
Proposition \ref{prop:spc flatness criterion} thus applies, and
shows that the orthogonal space
$\bL^\mC_W/G = \bL^\mC_{G,W}$ is flat.
\end{proof}

Now we explain in which way the orthogonal spaces $\bL^\mC_{G,W}$
are multiplicative in the pair $(G,W)$.
We let $K$ be another compact Lie group and $U$ a unitary $K$-representation.
The maps
\begin{align*}
\bL^\mC(W,V_\mC)/G \times \bL^\mC(U,V'_\mC)/K \ &\to \ 
\bL^\mC(W\oplus U,(V\oplus V')_\mC)/ ( G\times K )\\
( \varphi\cdot G,\, \psi\cdot K ) \qquad \quad &\longmapsto \qquad 
 (\varphi\oplus\psi)\cdot(G\times K)
\end{align*}
form a bimorphism from $(\bL^\mC_{G,W},\bL^\mC_{K,U})$ to $\bL^\mC_{G\times K,W\oplus U}$
as the inner product spaces $V$ and $V'$ vary. The universal property of the box
product translates this into a morphism of orthogonal spaces
\[ \zeta_{G,K;W,U} \ : \  
  \bL^\mC_{G,W} \boxtimes \bL^\mC_{K,U} \ \to \ \bL^\mC_{G\times K,W\oplus U}\ . \]
The analogous morphism for semifree orthogonal spaces
(i.e., for orthogonal representations and without the superscript $^\mC$)
is an isomorphism, see Example \ref{eg:box of free orthogonal}.
For the complex analogs, an isomorphism would be too much to hope for, 
but the next best thing is true:

\begin{prop}\label{prop:box complex almost-representables}
Let $G$ and $K$ be compact Lie groups, $W$ a unitary $G$-representation
and $U$ a unitary $K$-representation.  
Then the morphism $\zeta_{G,K;W,U}$ is a global equivalence.
\end{prop}
\begin{proof}
The global equivalences discussed 
in Proposition \ref{prop:complex global classifying}~(i) 
make the following square commute:
\[ \xymatrix@C=20mm{ 
\bL_{G, u W}\boxtimes\bL_{K,u U}\ar[r]^-{\iso}_-{\eqref{eq:boxtimes_of_B_gl}}
\ar[d]_{f_{G,W}\boxtimes f_{K, U}}^{\simeq} &
\bL_{G\times K, u( W\oplus U)} \ar[d]_{\simeq}^{f_{G\times K, W\oplus U}}\\ 
 \bL^{\mC}_{G, W}\boxtimes\bL^{\mC}_{K,U}\ar[r]_-{\zeta_{G,K;W,U}} &
\bL^{\mC}_{G\times K, W\oplus U} }\]
The left vertical morphism is a global equivalence because these
are stable under box product (Theorem \ref{thm:box to times}~(ii)).
Since the vertical morphisms are global equivalences and the upper horizontal
one is an isomorphism, the lower horizontal morphism $\zeta_{G,K;W,U}$
is also a global equivalence.
\end{proof}

\section{Global families}
\label{sec:global families unstable}

In this section we explain a variant of unstable global homotopy
theory based on a {\em global family}, i.e., a class $\Fc$ of 
compact Lie groups with certain closure properties. 
We introduce $\Fc$-equivalences, a relative version of global equivalences,
and establish an $\Fc$-relative version of the global model
structure in Theorem \ref{thm:F-global spaces}. 
We also discuss the compatibility of the
$\Fc$-global model structure with the box product of
orthogonal spaces, see Proposition \ref{prop:ExF ppp spaces}.
Finally, we record that for multiplicative global families,
the $\Fc$-global model structure lifts to category of modules and
algebras, see Corollary \ref{cor-lift to modules spaces}.

\begin{defn}\label{def:global family} 
A {\em global family}\index{subject}{global family}\index{symbol}{$\Fc$ - {global family}} 
is a non-empty class of compact Lie groups
that is closed under isomorphisms, closed subgroups and quotient groups.
\end{defn}

Some relevant examples of global families are: 
all compact Lie groups;
all finite groups; 
all abelian compact Lie groups; 
all finite abelian groups;
all finite cyclic groups; 
all finite $p$-groups.
Another example is the global family $\td{G}$ generated by a compact Lie group $G$,
i.e., the class of all compact Lie groups isomorphic to a quotient 
of a closed subgroup of $G$.
A degenerate case is the global family $\td{e}$\index{symbol}{$\td{e}$ - {global family of trivial groups}} of all trivial groups.
In this case our theory specializes to the non-equivariant
homotopy theory of orthogonal spaces.
 
For a global family $\Fc$ and a compact Lie group $G$ we write $\Fc\cap G$\index{symbol}{$\Fc\cap G$ - {subgroups of $G$ belonging to $\Fc$}}
for the family of those closed subgroups of $G$
that belong to $\Fc$.  We also write $\Fc(m)$ for $\Fc\cap O(m)$,
the family of closed subgroups of $O(m)$ that belong to $\Fc$.\index{symbol}{$\Fc(m)$ - {subgroups of $O(m)$ belonging to $\Fc$}}
We recall that an equivariant continuous map of $O(m)$-spaces
is an {\em $\Fc(m)$-cofibration} if it has 
the right lifting property with respect to all morphisms $q:A\to B$
of $O(m)$-spaces such that the map $q^H:A^H\to B^H$ is a weak equivalence
and Serre fibration for all $H\in\Fc(m)$.

The following definitions of $\Fc$-level equivalences,
$\Fc$-level fibrations and $\Fc$-cofibrations 
are direct relativizations of the corresponding concepts in the
strong level model structure of orthogonal spaces.

\begin{defn} Let $\Fc$ be a global family.
  A morphism $f:X\to Y$ of orthogonal spaces is
  \begin{itemize}
  \item  an {\em $\Fc$-level equivalence}\index{subject}{F-level equivalence@$\Fc$-level equivalence!of orthogonal spaces} 
    if for every compact Lie group $G$ in $\Fc$ and every $G$-representation $V$ the map
    $f(V)^G:X(V)^G\to Y(V)^G$ is a weak equivalence;
  \item an {\em $\Fc$-level fibration}\index{subject}{F-level fibration@$\Fc$-level fibration!of orthogonal spaces} 
    if for every compact Lie group $G$ in $\Fc$ and every $G$-representation $V$ the map
    $f(V)^G:X(V)^G\to Y(V)^G$ is a Serre fibration;
  \item an $\Fc$-cofibration if the latching morphism
    $\nu_m f:X(\mR^m)\cup_{L_m X}L_m Y\to Y(\mR^m)$ is an $\Fc(m)$-cofibration 
    for all $m\geq 0$.
\end{itemize}
\end{defn}

Every inner product space $V$ is isometrically isomorphic to $\mR^m$
with the standard scalar product, where $m$ is the dimension of $V$.
So a morphism $f:X\to Y$ of orthogonal spaces is an
$\Fc$-level equivalence (respectively $\Fc$-level fibration)
precisely if for every $m\geq 0$
the map $f(\mR^m):X(\mR^m)\to Y(\mR^m)$  is an $\Fc(m)$-equivalence
(respectively $\Fc(m)$-projective fibration).
The formal argument is analogous to Lemmas \ref{lemma:characterize strong level equivalences} and \ref{lemma:characterize strong level fibrations} 
which treat the case $\Fc=\All$.
Clearly, the classes of $\Fc$-level equivalences,
$\Fc$-level fibrations and $\Fc$-cofibrations are closed under composition, retracts 
and coproducts.

Now we discuss the  $\Fc$-level model structures on orthogonal spaces.
When $\Fc=\All$ is the global family of all compact Lie groups, 
then $\All(m)$ is the family of all closed subgroups of $O(m)$.
For this maximal global family,
an $\All$-level equivalence is just a strong level equivalence
in the sense of Definition \ref{def:strong level equivalence spaces}.
Moreover, the $\All$-level fibrations coincide with the strong level fibrations.
The $\All$-cofibrations coincide with the flat cofibrations. 
So for the global family of all compact Lie groups the  
$\All$-level model structure on orthogonal spaces is the strong level model structure 
of Proposition \ref{prop:strong level spaces}.

\begin{prop}\label{prop:F-level spaces}
  Let $\Fc$ be a global family.
  The $\Fc$-level equivalences, $\Fc$-level fibrations 
  and $\Fc$-cofibrations form a topological and cofibrantly generated model structure,
  the {\em $\Fc$-level model structure},\index{subject}{F-level model structure@$\Fc$-level model structure!for orthogonal spaces}
  on the category of orthogonal spaces.
\end{prop}
\begin{proof}
We specialize Proposition \ref{prop:general level model structure} 
by letting $\Cc(m)$ be the $\Fc(m)$-projective
model structure on the category of $O(m)$-spaces,
compare Proposition \ref{prop:proj model structures for G-spaces}.
With respect to these choices of model structures $\Cc(m)$,
the classes of level equivalences, level fibrations and cofibrations
in the sense of Proposition \ref{prop:general level model structure} 
precisely become the $\Fc$-level equivalences,
 $\Fc$-fibrations and $\Fc$-cofibrations.
Every acyclic cofibration in the $\Fc(m)$-projective model structure of $O(m)$-spaces
is also an acyclic cofibration in the $\All$-projective model structure of $O(m)$-spaces.
So the consistency condition (see Definition \ref{def:consistency condition})
in the present situation is a special case of the consistency condition
for the strong level model structure that we established in the proof of
Proposition \ref{prop:strong level spaces}.

We describe explicit sets of generating cofibrations 
and generating acyclic cofibrations for the  $\Fc$-level model structure.
We let $I_{\Fc}$ be the set of all morphisms $G_mi$ for $m\geq 0$
and for $i$ in the set
of generating cofibrations for the $\Fc(m)$-projective model
structure on the category of $O(m)$-spaces specified in \eqref{eq:I_for_F-proj_on_GT}.
Then the set $I_{\Fc}$ detects the acyclic fibrations 
in the  $\Fc$-level model structure 
by Proposition \ref{prop:general level model structure}~(iii). 
Similarly, we let $J_{\Fc}$ be the set of all morphisms $G_m j$ 
for $m\geq 0$ and for $j$ in the set
of generating acyclic  cofibrations for the $\Fc(m)$-projective model
structure on the category of $O(m)$-spaces specified 
in \eqref{eq:J_for_F-proj_on_GT}.
Again by Proposition \ref{prop:general level model structure}~(iii),
$J_{\Fc}$ detects the fibrations in the  $\Fc$-level model structure. 
The $\Fc$-level model structure is topological 
by Proposition \ref{prop:topological criterion},
where we take $\Gc$ as the set of orthogonal spaces $L_{H,\mR^m}$
for all $m\geq 0$ and all $H\in\Fc(m)$.
\end{proof}

Now we proceed towards the construction of the 
{\em  $\Fc$-global model structure},
see Theorem \ref{thm:F-global spaces} below.
The weak equivalences in this model structures
are the $\Fc$-equivalences of the following definition,
the direct generalization of global equivalences in the presence of a global family.

\begin{defn}
  Let $\Fc$ be a global family. A morphism $f:X\to Y$ of orthogonal spaces
  is an {\em $\Fc$-equivalence}\index{subject}{F-equivalence@$\Fc$-equivalence!of orthogonal spaces}
  if the following condition holds: for every compact Lie group $G$ in $\Fc$,
  every $G$-representation $V$, every $k\geq 0$ 
  and all maps $\alpha:\partial D^k\to X(V)^G$
  and $\beta:D^k\to Y(V)^G$ such that $f(V)^G\circ\alpha=\beta|_{\partial D^k}$
  there is a $G$-representation $W$,
  a $G$-equivariant linear isometric embedding $\varphi:V\to W$
  and a continuous map $\lambda:D^k\to X(W)^G$
  such that $\lambda|_{\partial D^k}=X(\varphi)^G\circ \alpha$ and
  such that $f(W)^G\circ \lambda$
  is homotopic, relative to $\partial D^k$, to $Y(\varphi)^G\circ \beta$.
\end{defn}

When $\Fc=\All$ is the maximal global family of all compact Lie groups,
then $\All$-equivalences are precisely the global equivalences.
The following proposition generalizes Proposition \ref{prop:telescope criterion},
and it is proved in much the same way.

\begin{prop}\label{prop:F-telescope criterion} 
Let $\Fc$ be a global family.
For every morphism of orthogonal spaces $f:X\to Y$, the following
three conditions are equivalent.
\begin{enumerate}[\em (i)]
\item The morphism $f$ is an $\Fc$-equivalence.
\item Let $G$ be a compact Lie group,
$V$ a $G$-representation and $(B,A)$ a finite $G$-CW-pair 
all of whose isotropy groups belong to $\Fc$.
Then for all continuous $G$-maps $\alpha:A\to X(V)$
and $\beta:B\to Y(V)$ such that $\beta|_A=f(V)\circ\alpha$,
there is a $G$-representation $W$, a $G$-equivariant linear
isometric embedding $\varphi:V\to W$ and a continuous $G$-map $\lambda:B\to X(W)$
such that $\lambda|_A=X(\varphi)\circ \alpha$ and
such that $f(W)\circ \lambda$
is $G$-homotopic, relative to $A$, to $Y(\varphi)\circ \beta$.
\item For every compact Lie group $G$ in the family $\Fc$
and every exhaustive sequence $\{V_i\}_{i\geq 1}$
of $G$-represen\-tations the induced map
\[ \tel_i f(V_i)\ : \ \tel_i X(V_i)\ \to \ \tel_i Y(V_i) \]
is a $G$-weak equivalence.
\end{enumerate}
\end{prop}

\begin{defn}\index{subject}{F-global fibration@$\Fc$-global fibration!of orthogonal spaces}
    A morphism $f:X\to Y$ of orthogonal spaces is an {\em $\Fc$-global fibration}
    if it is an $\Fc$-level fibration
    and for every compact Lie group $G$ in the family $\Fc$,
    every faithful $G$-representation $V$ 
    and every equivariant linear isometric embedding $\varphi:V\to W$ 
    of $G$-representations, the map
    \[ (f(V)^G, X(\varphi)^G)\ : \ X(V)^G \ \to \ 
    Y(V)^G \times_{Y(W)^G }  X(W)^G \]
    is a weak equivalence.
\end{defn}

The next proposition contains various properties of $\Fc$-equivalences
that generalize Proposition \ref{prop:global equiv basics}
and Proposition \ref{prop:global equiv preservation spc} (i).

\begin{prop}\label{prop:F equiv preservation} 
  Let $\Fc$ be a global family.
  \begin{enumerate}[\em (i)]
  \item 
    Every $\Fc$-level equivalence is an $\Fc$-equivalence.
  \item The composite of two $\Fc$-equivalences is an $\Fc$-equivalence.
  \item 
    If $f, g$ and $h$ are composable morphisms of orthogonal spaces 
    such that $g f$ and $h g$ are $\Fc$-equivalences, then $f, g, h$ and $h g f$
    are also $\Fc$-equivalences.
  \item 
    Every retract of an $\Fc$-equivalence is an $\Fc$-equivalence.
  \item 
    A coproduct  of any set of $\Fc$-equivalences is an $\Fc$-equivalence.
  \item 
    A finite product of $\Fc$-equivalences is an $\Fc$-equivalence.
  \item 
    For every space $K$ the functor $-\times K$ preserves $\Fc$-equivalences
    of orthogonal spaces.
  \item Let $e_n:X_n\to X_{n+1}$ and $f_n:Y_n\to Y_{n+1}$ be morphisms 
    of orthogonal spaces that are objectwise closed embeddings, for $n\geq 0$. 
    Let $\psi_n:X_n\to Y_n$ be $\Fc$-equivalences of orthogonal spaces
    that satisfy $\psi_{n+1}\circ e_n=f_n\circ\psi_n$ for all $n\geq 0$.
    Then the induced morphism $\psi_\infty:X_\infty\to Y_\infty$ 
    between the colimits of the sequences is an $\Fc$-equivalence.
  \item Let $f_n:Y_n\to Y_{n+1}$ be an $\Fc$-equivalence
   and a closed embedding of orthogonal spaces, for $n\geq 0$. 
   Then the canonical morphism 
   $f_\infty:Y_0\to Y_\infty$ to the colimit of the sequence $\{f_n\}_{n\geq 0}$
   is an $\Fc$-equivalence.
\item  Let 
    \[ \xymatrix{
      C  \ar[d]_\gamma & A \ar[l]_-g \ar[d]^\alpha \ar[r]^-f & B \ar[d]^\beta\\
      C'  & A' \ar[l]^-{g'}\ar[r]_-{f'} & B' } \]
    be a commutative diagram of orthogonal spaces such that $f$ and $f'$ are
    h-cofibrations.
    If the morphisms $\alpha,\beta$ and $\gamma$ are $\Fc$-equivalences,
    then so is the induced morphism of pushouts
    \[ \gamma\cup \beta\ : \ C\cup_A B \ \to \ C'\cup_{A'} B'\ . \]
  \item  Let 
    \[ \xymatrix{ A \ar[r]^-f \ar[d]_g & B \ar[d]^h\\
      C \ar[r]_-k & D } \]
    be a pushout square of orthogonal spaces such that $f$ is an $\Fc$-equivalence.
    If in addition $f$ or $g$ is an h-cofibration, then the morphism 
    $k$ is an $\Fc$-equivalence.
 \item Let 
   \[\xymatrix{ P\ar[r]^-k \ar[d]_g & X \ar[d]^f \\ Z \ar[r]_-h & Y }\]
   be a pullback square of orthogonal spaces in which $f$ is an $\Fc$-equivalence.
   If in addition one of the morphisms $f$ or $h$ is
   an $\Fc$-level fibration, then the morphism $g$ is also an $\Fc$-equivalence.
   \item 
    Every $\Fc$-equivalence that is also an $\Fc$-global fibration 
    is an $\Fc$-level equivalence.
  \item The box product of two $\Fc$-equivalences is an $\Fc$-equivalence.
\end{enumerate}
\end{prop}
\begin{proof}
The proofs of~(i) through~(xii)
are almost verbatim the same as the corresponding parts
of Proposition \ref{prop:global equiv basics}, and we omit them.
Part~(xiii) is proved in the same way as 
Proposition \ref{prop:global equiv preservation spc}~(i). 

(xiv)
The product of orthogonal spaces preserves 
$\Fc$-equivalences in both variables by part~(vi).
The morphism $\rho_{X,Y}:X\boxtimes Y\to X\times Y$ is a global equivalence,
hence an $\Fc$-equivalence, for all orthogonal spaces $X$ and $Y$,
by Theorem \ref{thm:box to times}~(i); this implies the claim.
\end{proof}

Now we establish the  $\Fc$-global model structures on the
category of orthogonal spaces.
We spell out sets of generating cofibrations and
generating acyclic cofibrations for the $\Fc$-global model structures.
In Proposition \ref{prop:F-level spaces}
we introduced $I_{\Fc}$ as the set of all morphisms $G_m i$ for $m\geq 0$
and for $i$ in the set of generating cofibrations for the $\Fc(m)$-projective model
structure on the category of $O(m)$-spaces specified in \eqref{eq:I_for_F-proj_on_GT}.
The set $I_{\Fc}$ detects the acyclic fibrations 
in the  $\Fc$-level model structure, 
which coincide with the acyclic fibrations 
in the  $\Fc$-global model structure.

Also in Proposition \ref{prop:F-level spaces}
we defined  $J_{\Fc}$ as the set of all morphisms $G_m j$ 
for $m\geq 0$ and for $j$ in the set
of generating acyclic  cofibrations for the $\Fc(m)$-projective model
structure on the category of $O(m)$-spaces specified 
in \eqref{eq:J_for_F-proj_on_GT}.
The set $J_{\Fc}$ detects the fibrations in the  $\Fc$-level model structure. 
We add another set of morphisms $K_\Fc$ 
that detects when the squares \eqref{eq:fibration characterization spaces}
are homotopy cartesian for $G\in\Fc$.
We set
\[ K_\Fc \ = \ \bigcup_{G,V,W\ : \ G\in\Fc} \Zc(\rho_{G,V,W}) \ ,\]
the set of all pushout products of sphere inclusions $i_k:\partial D^k\to D^k$
with the mapping cylinder inclusions of the morphisms $\rho_{G,V,W}$,
compare Construction \ref{con:define Z(j)};
here the union is over a set of representatives
of the isomorphism classes of triples $(G,V,W)$ consisting of
a compact Lie group $G$ in $\Fc$, a faithful $G$-representation $V$ and
an arbitrary $G$-representation $W$.
By Proposition \ref{prop:hocartesian via RLP}, the right lifting property 
with respect to the union $J_\Fc \cup K_\Fc$ 
characterizes the  $\Fc$-global fibrations.

The proof of the following theorem
proceeds by mimicking the proof in the special case $\Fc=\All$,
and all arguments in the proof of Theorem \ref{thm:All global spaces} 
go through almost verbatim. Whenever the small object argument is used,
it now has to be taken with respect to the set
$J_\Fc\cup K_\Fc$ (as opposed to the set $J^{\str}\cup K$).

\begin{theorem}[$\Fc$-global model structure]\label{thm:F-global spaces} 
  Let $\Fc$ be a global family.
  \begin{enumerate}[\em (i)]
  \item 
    The $\Fc$-equivalences, $\Fc$-global fibrations and $\Fc$-cofibrations 
    form a model structure,
    the {\em $\Fc$-global model structure},\index{subject}{F-global model structure@$\Fc$-global model structure!for orthogonal spaces}\index{subject}{global model structure!$\Fc$-}
    on the category of orthogonal spaces. 
  \item
    The fibrant objects in the  $\Fc$-global model structure 
    are the {\em $\Fc$-static orthogonal spaces},\index{subject}{F-static@$\Fc$-static}
    i.e., those orthogonal spaces $X$ such that for every compact Lie group $G$ in $\Fc$,
    every faithful $G$-representation $V$
    and every $G$-equivariant linear isometric embedding $\varphi:V\to W$
    the map of $G$-fixed point spaces $X(\varphi)^G: X(V)^G\to X(W)^G$
    is a weak equivalence.
  \item
    A morphism of orthogonal spaces is:
    \begin{itemize}
    \item 
      an acyclic fibration in the  $\Fc$-global model structure
      if and only if it has the right lifting property 
      with respect to the set $I_{\Fc}$; 
    \item 
      a fibration in the  $\Fc$-global model structure
      if and only if it has the right lifting property 
      with respect to the set $J_\Fc \cup K_\Fc$.  
    \end{itemize}
  \item
    The  $\Fc$-global model structure is cofibrantly generated, 
    proper and topological. 
  \end{enumerate}
\end{theorem}

\begin{eg}
In the case $\Fc=\td{e}$ of the minimal global family of trivial groups,
the $\td{e}$-global homotopy theory of orthogonal spaces
just another model for the (non-equivariant) homotopy theory of spaces.
Indeed, the evaluation functor $\ev_0:\spc\to \bT$ is a right Quillen
equivalence with respect to the  $\td{e}$-global model structure.
So the total derived functor
\[ \Ho(\ev_0) \ :\ \Ho^{\td{e}}(\spc) \ \to \Ho(\bT) \]
is an equivalence of homotopy categories.

In fact, for the global family $\Fc=\td{e}$, most of what we do here has already been
studied before: The  $\td{e}$-global model structure 
and the fact that it is Quillen equivalent to the model category of spaces 
were established by Lind \cite[Thm.\,1.1]{lind-diagram};
in \cite{lind-diagram}, orthogonal spaces are called `$\mathcal I$-spaces' and
$\td{e}$-global equivalences 
are called `weak homotopy equivalences' and 
are defined as those morphisms that induce weak equivalences on homotopy colimits.
\end{eg}

\begin{cor}\label{cor-characterize F-equivalences spaces}
  Let $f:A\to B$ be a morphism of orthogonal spaces
  and $\Fc$ a global family.
  Then the following conditions are equivalent.
  \begin{enumerate}[\em (i)]
  \item The morphism $f$ is an $\Fc$-equivalence.\index{subject}{F-equivalence@$\Fc$-equivalence!of orthogonal spaces}
  \item
    The morphism can be written as $f=w_2\circ w_1$ for an
    $\Fc$-level equivalence $w_2$ and a global equivalence $w_1$.
  \item For some (hence any) $\Fc$-cofibrant approximation
    $f^c:A^c\to B^c$ in the  $\Fc$-level model structure
    and every $\Fc$-static orthogonal space $Y$ 
    the induced map
    \[ [f^c,Y] \ : \ [B^c,Y] \ \to \ [A^c,Y] \]
    on homotopy classes of morphisms is bijective.
\end{enumerate}
\end{cor}
\begin{proof}
(i)$\Longleftrightarrow$(ii) 
The $\Fc$-equivalences contain the global equivalences by definition,
and the $\Fc$-level equivalences by 
Proposition \ref{prop:F equiv preservation}~(i), and are closed under composition by
Proposition \ref{prop:F equiv preservation}~(ii);
so all composites $w_2\circ w_1$ as in the claim are $\Fc$-equivalences. 
On the other hand, every $\Fc$-equivalence $f$ can be factored 
in the global model structure of Theorem \ref{thm:All global spaces} as $f=q j$
where $j$ is a global equivalence and $q$ is a global fibration.
Since $f$ and $j$ are $\Fc$-equivalences, so is $q$ by 
Proposition \ref{prop:F equiv preservation}~(iii);
so $q$ is an $\Fc$-equivalence and a global fibration, 
hence an $\Fc$-level equivalence by
Proposition \ref{prop:F equiv preservation}~(xiii).

(i)$\Longleftrightarrow$(iii) The morphism $f$ is an $\Fc$-equivalence
if and only if the $\Fc$-co\-fibrant approximation $f^c:A^c\to B^c$ is
an $\Fc$-equivalence.
Since $A^c$ and $B^c$ are $\Fc$-cofibrant, they are
cofibrant in the  $\Fc$-global model structure.
So by general model model category theory
(see for example \cite[Cor.\,7.7.4]{hirschhorn}),
$f^c$ is an $\Fc$-equivalence if and only if the induced map $[f^c,X]$ 
is bijective for every fibrant object in the
 $\Fc$-global model structure.
By Theorem \ref{thm:F-global spaces}~(ii) 
these fibrant objects are precisely the $\Fc$-static orthogonal spaces.
\end{proof}

\begin{rk}[Mixed global model structures]
Cole's `mixing theorem' for model structures \cite[Thm.\,2.1]{cole-mixed}  
allows to construct many more `mixed' $\Fc$-global model structures 
on the category of orthogonal spaces.
We consider two global families such that $\Fc\subseteq \Ec$. 
Then every $\Ec$-equivalence is an $\Fc$-equivalence and every 
fibration in the  $\Ec$-global model structure
is a fibration in the  $\Fc$-global model structure.
By Cole's theorem \cite[Thm.\,2.1]{cole-mixed} 
the $\Fc$-equivalences and the fibrations of the
 $\Ec$-global model structure are part of a model structure,
the {\em $\Ec$-mixed $\Fc$-global model structure}\index{subject}{mixed global model structure}\index{subject}{global model structure!mixed}
on the category of orthogonal spaces.
By \cite[Prop.\,3.2]{cole-mixed} 
the cofibrations in the $\Ec$-mixed $\Fc$-global model structure 
are precisely the retracts of all composites
$h\circ g$ in which $g$ is an $\Fc$-cofibration and $h$ is simultaneously
an $\Ec$-equivalence and an $\Ec$-cofibration.
In particular, an orthogonal space is cofibrant in the
$\Ec$-mixed $\Fc$-global model structure if it is $\Ec$-cofibrant and
$\Ec$-equivalent to an $\Fc$-cofibrant 
orthogonal space \cite[Cor.\,3.7]{cole-mixed}.
The $\Ec$-mixed $\Fc$-global model structure is again proper 
(Propositions~4.1 and~4.2 of \cite{cole-mixed}).

When $\Fc=\td{e}$ is the minimal family of trivial groups,
this provides infinitely many $\Ec$-mixed model structures on
the category of orthogonal spaces that are all Quillen equivalent
to the model category of (non-equivariant) spaces, with respect to
weak equivalences.
\end{rk}

The next topic is the compatibility of the $\Fc$-global model structure 
with the box product of orthogonal spaces.
Given two morphisms $f:A\to B$ and $g:X\to Y$ of orthogonal spaces
we denote by 
\[ f\Box g = (f\boxtimes Y)\cup(B\boxtimes g) \ : \ 
A\boxtimes Y\cup_{A\boxtimes X}B\boxtimes X \ \to \ B\boxtimes Y\]
the {\em pushout product} morphism.\index{subject}{pushout product}
We recall that a  model structure on a symmetric monoidal category
satisfies the {\em pushout product property}
if the following two conditions hold:
\begin{itemize}
\item for every pair of cofibrations $f$ and $g$ the
  pushout product morphism $f\Box g$ is also a cofibration;
\item if in addition $f$ or $g$ is a weak equivalence, then so is the
  pushout product morphism  $f\Box g$.
\end{itemize}
We let $\Ec$ and $\Fc$ be two global families. We denote by
$\Ec\times\Fc$ the smallest global family that contains all
groups of the form $G\times K$ for $G\in\Ec$ and $K\in\Fc$.
So a compact Lie group $H$ belongs to $\Ec\times\Fc$ if and
only if $H$ is isomorphic to a closed subgroup of a 
group of the form $(G\times K)/N$ for some groups $G\in\Ec$ and $K\in\Fc$, 
and some closed normal subgroup $N$ of $G\times K$.

\begin{prop} \label{prop:ExF ppp spaces}
  Let $\Ec$ and $\Fc$ be two global families.
  \begin{enumerate}[\em (i)]
  \item The pushout product of an $\Ec$-cofibration with an $\Fc$-cofibration
    is an $(\Ec\times\Fc)$-cofibration.  
  \item The pushout product of a flat cofibration that is also
    an $\Fc$-equivalence with any morphism of orthogonal spaces is an $\Fc$-equivalence.
  \item 
    Let $\Fc$ be a multiplicative global family, i.e., $\Fc\times\Fc =\Fc$.
    Then the  $\Fc$-global model structure  
    satisfies the pushout product property
    with respect to the box product of orthogonal spaces.
  \item 
    The positive global model structure satisfies the pushout product property
    with respect to the box product of orthogonal spaces.
  \end{enumerate}
\end{prop}
\begin{proof} 
  (i) It suffices to show the claim for sets of generating cofibrations. 
  The $\Ec$-cofibrations are generated by the morphisms
  \[   \bL_{G,V}\times i_k \ : \  \bL_{G,V}\times \partial D^k
  \ \to \  \bL_{G,V}\times D^k   \]
  for $G\in\Ec$, $V$ a $G$-representation and $k\geq 0$.
  Similarly, the $\Fc$-cofibrations are generated by the morphisms
  $\bL_{K,W}\times i_m$ for $K\in\Fc$, $W$ a $K$-representation and $m\geq 0$.
  The pushout product of two such generators is isomorphic to the morphism
  \[
  \bL_{G\times K,V\oplus W}\times i_{k+m}  \ : \  \bL_{G\times K,V\oplus W}\times \partial D^{k+m}
  \  \to \ \bL_{G\times K,V\oplus W}\times D^{k+m} \ ,  \]
  compare Example \ref{eg:box of free orthogonal}.
  Since $G\times K$ belongs to the family $\Ec\times\Fc$,
  this pushout product morphism is an $(\Ec\times\Fc)$-cofibration.

  (ii) We let $i:A\to B$ and $j:K\to L$ be morphisms of orthogonal spaces
  such that $i$ is a flat cofibration and an $\Fc$-equivalence.
  Then $i\boxtimes K$ and $i\boxtimes L$ are $\Fc$-equivalences
  by Proposition \ref{prop:F equiv preservation}~(xiv).
  Moreover, $i$ is an h-cofibration 
  by Corollary \ref{cor-h-cofibration closures}~(iii),
  hence so is $i\boxtimes K:A\boxtimes K\to B\boxtimes K$.
  Thus its cobase change, the canonical morphism
  \[ A\boxtimes L \ \to \ A\boxtimes L \cup_{A\boxtimes K} B\boxtimes K \]
  is an $\Fc$-equivalence by Proposition \ref{prop:F equiv preservation}~(xi).
  Since $i\boxtimes L:A\boxtimes L\to B\boxtimes L$ is also an $\Fc$-equivalence,
  so is the pushout product map, by 2-out-of-6 
  (Proposition \ref{prop:F equiv preservation}~(iii)).

  (iii) The part of the pushout product property that refers only to cofibrations
  is true by part~(i) with $\Ec=\Fc$ and the hypothesis that $\Fc\times\Fc=\Fc$. 
  Every cofibration in the  $\Fc$-global model structure is
  in particular a flat cofibration,
  so the part of the pushout product property that also refers to 
  acyclic cofibrations in the  $\Fc$-global model structure
  is a special case of (ii).

  Part~(iv) is proved in the essentially the same way as (iii), 
  for the global family of all compact Lie groups.
\end{proof}

Finally, we will discuss another important relationship between
the $\Fc$-global model structures and the box product, 
namely the {\em monoid axiom} \cite[Def.\,3.3]{schwede-shipley-monoidal}.
We only discuss a slightly weaker form of the monoid axiom 
in the sense that we only cover sequential (as opposed to more general transfinite)
compositions.

\begin{prop}[Monoid axiom]\index{subject}{monoid axiom!for the box product of orthogonal spaces} 
\label{prop:monoid orthogonal spaces}
We let $\Fc$ be a global family.
For every flat cofibration $i:A\to B$ that is also an $\Fc$-equivalence
and every orthogonal space $Y$ the morphism
\[ i\boxtimes Y \ : \ A\boxtimes Y \ \to \ B\boxtimes Y \]
is an h-cofibration and an $\Fc$-equivalence.\index{subject}{h-cofibration}
Moreover, the class of h-cofibrations that are also $\Fc$-equivalences
is closed under cobase change, coproducts and sequential compositions.
\end{prop}
\begin{proof}
Every flat cofibration is an h-cofibration 
(Corollary \ref{cor-h-cofibration closures}~(iii)
applied to the strong level model structure),
and h-cofibrations are closed under box product with any orthogonal space
(Corollary \ref{cor-h-cofibration closures}~(ii)),
so $i\boxtimes Y$ is an h-cofibration.
Since $i$ is an $\Fc$-equivalence, so is $i\boxtimes Y$ 
by Proposition \ref{prop:F equiv preservation}~(xiv).

Proposition \ref{prop:F equiv preservation} shows that the class of h-cofibrations 
that are also $\Fc$-equivalences is closed under cobase change, 
coproducts and sequential compositions.
\end{proof}

We let $\Fc$ be a multiplicative global family, i.e., $\Fc\times\Fc=\Fc$.
The constant one-point orthogonal space $\mathbf 1$ 
is the unit object for the box product 
of orthogonal spaces, and it is `free', i.e., $\td{e}$-cofibrant.
So $\mathbf 1$ is cofibrant in the $\Fc$-global model structure. 
So with respect to the box product, the $\Fc$-global model structure is a
symmetric monoidal model category
in the sense of \cite[Def.\,4.2.6]{hovey-book}.
A corollary is that the unstable $\Fc$-global homotopy category, i.e., 
the localization of the category
of orthogonal spaces at the class of $\Fc$-equivalences,
inherits a closed symmetric monoidal structure, 
compare \cite[Thm.\,4.3.3]{hovey-book}.
This `derived box product'\index{subject}{box product!of orthogonal spaces} 
is nothing new, though: 
since the morphism $\rho_{X,Y}:X\boxtimes Y\to X\times Y$
is a global equivalence for all orthogonal spaces $X$ and $Y$, 
the derived box product is just a categorical product in $\Ho^\Fc(\spc)$.

\medskip

\begin{defn}\label{def:orthogonal monoid space}
\index{subject}{orthogonal monoid space} 
\index{subject}{monoid space!orthogonal|see{orthogonal monoid space}} 
An {\em orthogonal monoid space}
is an orthogonal space $R$ equipped with unit morphism $\eta:\mathbf 1\to R$
and a multiplication morphism $\mu:R\boxtimes R\to R$
that is unital and associative in the sense that the diagram
\[ \xymatrix@C=15mm{
(R\boxtimes R)\boxtimes R \ar[r]_-\iso^-{\text{associativity}} \ar[d]_{\mu\boxtimes R} &
R\boxtimes(R\boxtimes R)\ar[r]^-{R\boxtimes \mu} &
R\boxtimes R \ar[d]^{\mu} \\
R\boxtimes R \ar[rr]_-{\mu} && R }\]
commutes.
An orthogonal monoid space $R$ is 
{\em commutative}\index{subject}{orthogonal monoid space!commutative} if moreover
$\mu\circ\tau_{R,R}=\mu$, where $\tau_{R,R}:R\boxtimes R\to R\boxtimes R$
is the symmetry isomorphism of the box product.
A {\em morphism} of orthogonal monoid spaces is a morphism of orthogonal spaces
$f:R\to S$ such that $f\circ\mu^R= \mu^S\circ(f\boxtimes f)$ and
$f\circ\eta_R=\eta_S$.
\end{defn}

One can expand the data of an orthogonal monoid space into an `external' form
as follows. The unit morphism $\eta:\mathbf 1\to R$ is determined by
a unit element $0\in R(0)$, the image of the map $\eta(0):\mathbf 1(0)\to R(0)$.
The multiplication map corresponds to continuous maps
$\mu_{V,W}:R(V) \times R(W) \to R(V\oplus W)$ for all inner product spaces
$V$ and $W$ that form a bimorphism as $(V,W)$ varies and such that
\[ \mu_{V,0}(x,0)\ = \ x \text{\qquad and\qquad} \mu_{0,W}(0,y)\ = \ y\ . \]
Put another way, the data of an orthogonal monoid space in external form
is that of a lax monoidal functor.
The commutativity condition can be expressed in terms of the external multiplication
as the commutativity of the diagrams
 \[\xymatrix@C=15mm{  R(V) \times R(W) \ar[r]^-{\mu_{V,W}}
\ar[d]_{\text{twist}} &
R(V\oplus W) \ar[d]^{R(\tau_{V,W})} \\
R(W) \times R(V) \ar[r]_-{\mu_{W,V}} & R(W\oplus V) } \]
where $\tau_{V,W}:V\oplus W\to W\oplus V$ interchanges the summands.
So commutative orthogonal monoid spaces
in external form are lax symmetric monoidal functors.
We will later refer to commutative orthogonal monoid spaces as
{\em ultra-commutative monoids}.\index{subject}{ultra-commutative monoid}

Every $\Fc$-cofibration is in particular a flat cofibration,
so the monoid axiom in the $\Fc$-global model structure holds.
If the global family $\Fc$ is closed under products,
Theorem~4.1 of \cite{schwede-shipley-monoidal} applies
to the  $\Fc$-global model structure 
of Theorem \ref{thm:F-global spaces} and shows:

\begin{cor}\label{cor-lift to modules spaces} 
Let $R$ be an orthogonal monoid space and $\Fc$ a multiplicative 
  global family.
  \begin{enumerate}[\em (i)]
  \item The category of $R$-modules admits the {\em $\Fc$-global model structure}\index{subject}{F-global model structure@$\Fc$-global model structure!for $R$-modules}  
    in which a morphism is an equivalence (respectively fibration)
    if the underlying morphism of orthogonal spaces is an $\Fc$-equivalence
    (respectively $\Fc$-global fibration).
    This model structure is cofibrantly generated.
    Every cofibration in this $\Fc$-global model structure 
    is an h-cofibration of underlying orthogonal spaces. 
    If the underlying orthogonal space of $R$ is $\Fc$-cofibrant,
    then every cofibration of $R$-modules is an $\Fc$-cofibration of underlying
    orthogonal spaces.
  \item If $R$ is commutative, then with respect to $\boxtimes_R$ 
    the $\Fc$-global model structure of $R$-modules is a monoidal model category 
    that satisfies the monoid axiom.
  \item If $R$ is commutative, then the category of $R$-algebras 
    admits the {\em $\Fc$-global model structure}\index{subject}{F-global model structure@$\Fc$-global model structure!for $R$-algebras}  
    in which a morphism is an equivalence (respectively fibration)
    if the underlying morphism of orthogonal spaces is an $\Fc$-equivalence
    (respectively $\Fc$-global fibration).
    Every cofibrant $R$-algebra is also cofibrant as an $R$-module.
  \end{enumerate}
\end{cor}
\begin{proof}
Almost of the statements are in Theorem~4.1 of \cite{schwede-shipley-monoidal}.
The only additional claims that require justification are the
two statements in part~(i) that concern the behavior of
the forgetful functor on the cofibrations in the $\Fc$-global model structure.

Since the forgetful functor from $R$-modules to orthogonal spaces preserves
all colimits and the classes of h-cofibrations and of $\Fc$-cofibrations
of orthogonal spaces are both closed under coproducts, cobase change, 
sequential colimits and retracts, it suffices to show each claim
for the generating cofibrations in the $\Fc$-global model structure
on the category of $R$-modules. These are of the form 
\[ ( R\boxtimes \bL_{H,\mR^m}) \times i_k\]
for some $k,m\geq 0$ and $H$ a closed subgroup of $O(m)$ that belongs to the
global family $\Fc$; as usual $i_k:\partial D^k\to D^k$ is the inclusion.
Since $i_k$ is an h-cofibration of spaces, 
the morphisms $(R\boxtimes \bL_{H,\mR^m})\times i_k$
are h-cofibrations of orthogonal spaces.
This concludes the proof that every cofibration of $R$-modules
is an h-cofibration of underlying orthogonal spaces.

Now we suppose that the underlying orthogonal space of $R$ is $\Fc$-cofibrant.
Because $H$ belongs to $\Fc$, the orthogonal space $\bL_{H,\mR^m}$
is $\Fc$-cofibrant. Hence the orthogonal space $R\boxtimes \bL_{H,\mR^m}$
is $\Fc$-cofibrant by Proposition \ref{prop:ExF ppp spaces}~(iii).
So $(R\boxtimes \bL_{H,\mR^m})\times i_k$ is an $\Fc$-cofibration
of orthogonal spaces.
This concludes the proof that every cofibration of $R$-modules
is an $\Fc$-cofibration of underlying orthogonal spaces.
\end{proof}

Strictly speaking, Theorem~4.1 of \cite{schwede-shipley-monoidal}
does not apply verbatim to the $\Fc$-global model structure
because the hypothesis that every object is small (with respect to some regular cardinal)
is not satisfied, and our version of the monoid axiom 
in Proposition \ref{prop:monoid orthogonal spaces}
is weaker than Theorem~3.3 of \cite{schwede-shipley-monoidal}
in that we do not close under {\em transfinite} compositions.
However, in our situation the sources of the generating cofibrations
and generating acyclic cofibrations are small with respect to sequential
compositions of flat cofibrations, and this suffices to
run the countable small object argument 
(compare also Remark~2.4 of \cite{schwede-shipley-monoidal}).

\begin{prop}\label{prop:box cofibrant R-mod} 
Let $R$ be an orthogonal monoid space and $N$ a right $R$-module
that is cofibrant in the $\All$-global model structure 
of Corollary {\em \ref{cor-lift to modules spaces}~(i)}. 
Then for every global family $\Fc$,
the functor $N\boxtimes_R -$ takes $\Fc$-equivalences of left $R$-modules
to $\Fc$-equivalences of orthogonal spaces.  
\end{prop}
\begin{proof}
For the course of this proof we call an $R$-module $N$ {\em homotopical}
if the functor $N\boxtimes_R -$ takes $\Fc$-equivalences of left $R$-modules
to $\Fc$-equivalences of orthogonal spaces.  
Since the $\All$-global model structure on the category of right $R$-modules
is obtained by lifting the global model structure of orthogonal spaces
along the free and forgetful adjoint functor pair, every cofibrant right $R$-module
is a retract of an $R$-module that arises as the colimit of a sequence
\begin{equation}\label{eq:R-module sequence}
 \emptyset = M_0 \ \to \ M_1 \ \to \dots \ \to \ M_k \ \to \ \dots 
\end{equation}
in which each $M_k$ is obtained from $M_{k-1}$ as a pushout
\[ \xymatrix@C=15mm{ 
 A_k\boxtimes R \ar[r]^-{f_k\boxtimes R}\ar[d] & 
 B_k \boxtimes R \ar[d] & \\
M_{k-1}\ar[r] & M_k } \]
for some flat cofibration $f_k:A_k\to B_k$ of orthogonal spaces.
For example, $f_k$ can be chosen as a disjoint union of morphisms 
in the set $I^{\str}$ of generating flat cofibrations.
We show by induction on $k$ that each module $M_k$ is homotopical.
The induction starts with the empty $R$-module, where there is nothing to show.
Now we suppose that $M_{k-1}$ is homotopical, and we claim that then $M_k$
is homotopical as well.
To see this we consider an $\Fc$-equivalence of left $R$-modules
$\varphi:X\to Y$. Then $M_k\boxtimes_R \varphi$
is obtained by passing to horizontal pushouts in the following commutative diagram
of orthogonal spaces:
\[ \xymatrix@C=15mm{ 
M_{k-1}\boxtimes_R X \ar[d]_{M_{k-1}\boxtimes_R \varphi} & 
 A_k \boxtimes X\ar[d]^{A_k\boxtimes \varphi}\ar[l]\ar[r]^-{f_k\boxtimes X}
& B_k \boxtimes X\ar[d]^{B_k\boxtimes \varphi}\\
M_{k-1}\boxtimes_R Y &  A_k \boxtimes Y\ar[l]\ar[r]_-{f_k\boxtimes Y}& B_k \boxtimes Y
} \]
Here we exploit that $(A_k\boxtimes R)\boxtimes_R X$
is naturally isomorphic to $A_k\boxtimes X$.
In the diagram, the left vertical morphism is an $\Fc$-equivalence by hypothesis. 
The middle and right vertical morphisms are 
$\Fc$-equivalences because box product is homotopical for $\Fc$-equivalences
(Proposition \ref{prop:F equiv preservation}~(xiv)).
Moreover, since the morphism $f_k$ is a flat cofibration, it is an h-cofibration
(by Corollary \ref{cor-h-cofibration closures}~(iii)),
and so the morphisms $f_k\boxtimes X$ and $f_k\boxtimes Y$ are h-cofibrations.
Proposition \ref{prop:F equiv preservation}~(x)
then shows that the induced morphism on horizontal
pushouts $M_k\boxtimes_R \varphi$ is again an $\Fc$-equivalence.

Now we let $M$ be a colimit of the sequence \eqref{eq:R-module sequence}.
Then $M\boxtimes_R X$ is a colimit of the sequence $M_k \boxtimes_R X$.
Moreover, since $f_k:A_k\to B_k$ is an h-cofibration, so is the 
morphism $f_k\boxtimes R$, and hence also its cobase change $M_{k-1}\to M_k$.
So the sequence whose colimit is $M\boxtimes_R X$ consists of h-cofibrations,
which are objectwise closed embeddings 
by Proposition \ref{prop:h-cof is closed embedding}.
The same is true for $M\boxtimes_R Y$.
Since each $M_k$ is homotopical and colimits of orthogonal spaces
along closed embeddings are homotopical 
(by Proposition \ref{prop:F equiv preservation}~(viii)), 
we conclude that the morphism $M\boxtimes_R\varphi:M\boxtimes_R X\to M\boxtimes_R Y$ is an $\Fc$-equivalence, so that $M$ is homotopical.
Since $\Fc$-equivalences are closed under retracts, 
the class of homotopical $R$-modules is closed under retracts, and so every
cofibrant right $R$-module is homotopical.
\end{proof}

\section{Equivariant homotopy sets}
\label{sec:equivariant homotopy sets}

In this section we define the equivariant homotopy sets $\pi_0^G(Y)$
of orthogonal spaces and relate them by restriction maps defined from 
continuous homomorphisms between compact Lie groups.
As the Lie groups vary, the resulting structure is a `Rep-functor' $\upi_0(Y)$, 
i.e., a contravariant functor
from the category of compact Lie groups and conjugacy classes of 
continuous homomorphisms.
The Rep-functor $\upi_0(B_{\gl}G)$ associated to a global classifying space
is the Rep-functor represented by $G$, 
by Proposition \ref{prop:fix of global classifying}. 
We identify the category of all natural operations 
with the category Rep of conjugacy classes of continuous homomorphisms,
compare Corollary \ref{cor:operations spc}.
Construction \ref{con:pairing equivalence homotopy spaces}
introduces a pairing of equivariant homotopy sets
$\times:\pi_0^G(X)\times\pi_0^G(Y)\to\pi_0^G(X\boxtimes Y)$
for any pair of orthogonal spaces, and
Proposition \ref{prop:unstable product properties} summarizes its main properties.

\medskip

We recall that $\Uc_G$ is a chosen complete $G$-universe and $s(\Uc_G)$
denotes the poset, under inclusion,
of finite-dimensional $G$-subrepresentations of $\Uc_G$.

\begin{defn}\label{def:define [A,Y]^G}
  Let $Y$ be an orthogonal space, $G$ be a compact Lie group and $A$ a $G$-space.
  We define 
  \[ [A,Y]^G \ = \ \colim_{V\in s(\Uc_G)}\, [A, Y(V)]^G \ ,\]
  the colimit over the poset $s(\Uc_G)$ of the sets of $G$-homotopy 
  classes of $G$-maps from $A$ to $Y(V)$.
\end{defn}

The canonical $G$-maps $Y(V)\to Y(\Uc_G)$  
induce maps from $[A, Y(V)]^G$ to $[A, Y(\Uc_G)]^G$ and hence a canonical map
\begin{equation}  \label{eq:[A,Y]^G_to_colim}
 [A,Y]^G\ \to \ [A, Y(\Uc_G)]^G\ .  
\end{equation}
In general there is no reason for this map to be injective or surjective.
If $Y$ is closed and $A$ is compact, the situation improves:

\begin{prop}\label{prop:[A,Y]^G of closed} 
Let $G$ be a compact Lie group.
\begin{enumerate}[\em (i)]
\item For every closed orthogonal space $Y$ and every compact $G$-space $A$
the canonical map \eqref{eq:[A,Y]^G_to_colim} is bijective.
\item
Let $\Fc$ be a global family and $f:X\to Y$ an $\Fc$-equivalence of orthogonal spaces.
Then for every finite $G$-CW-complex $A$ 
all of whose isotropy groups belong to $\Fc$, the induced map
\[ [A,f]^G\ : \ [A,X]^G\ \to \ [A,Y]^G\]
is bijective.
\item For every pair of orthogonal spaces $X$ and $Y$ and
every $G$-space $A$, the canonical map
\[ ([A,p_X]^G,[A,p_Y]^G)\ : \ [A,X\times Y]^G \ \to \ [A,X]^G \times[A, Y]^G
 \]
is bijective, where $p_X$ and $p_Y$ are the projections.
\end{enumerate}
\end{prop}
\begin{proof}
(i) Since the poset $s(\Uc_G)$ contains a cofinal subsequence, 
$Y(\Uc_G)$ is a sequential colimit of values of $Y$ along closed embeddings. 
By Proposition \ref{prop:filtered colim preserve weq}~(i),
every continuous $G$-map $A\to Y(\Uc_G)$ thus factors through $Y(V)$
for some finite-dimensional $V\in s(\Uc_G)$,
which shows surjectivity. Injectivity follows by the same argument
applied to the compact $G$-space $A\times[0,1]$.

(ii) We let $\beta:A\to Y(V)$ be a continuous $G$-map, for some $V\in s(\Uc_G)$,
that represents an element of $[A,Y]^G$. 
Together with the unique map from the empty space
this specifies an equivariant lifting problem  on the left:
\[ \xymatrix@C=15mm{
\emptyset \ar[r] \ar[d] & X(V) \ar[d]^{f(V)} &
\emptyset \ar[r]\ar[d] & X(V) \ar[r]^-{X(\varphi)} &  X(W) \ar[d]^{f(W)} 
\\
A\ar[r]_-\beta & Y(V) & 
A\ar[r]_-\beta \ar@{-->}[urr]^-(.4)\lambda & Y(V) \ar[r]_-{Y(\varphi)} & Y(W)  }\]
Since $(A,\emptyset)$ is a finite $G$-CW-pair with isotropy in $\Fc$
and $f$ an $\Fc$-equivalence, Proposition \ref{prop:F-telescope criterion}~(ii) 
provides a $G$-equivariant linear isometric embedding $\varphi:V\to W$ 
and a continuous $G$-map $\lambda$ on the right hand side such that $f(W)\circ\lambda$ 
is $G$-homotopic to $Y(\varphi)\circ\beta$.
We choose a $G$-equivariant linear isometric embedding $j:W\to \Uc_G$
extending the inclusion of $V$. Then the class in $[A,X]^G$
represented by the $G$-map
\[ X(j)\circ\lambda\ : \ A\ \to \ X(j(W)) \]
is taken to $[\beta]$ by the map $[A,f]^G$.
This shows that $[A,f]^G$ is surjective.

For injectivity we consider two $G$-maps $g,g':A\to X(V)$, for some $V\in s(\Uc_G)$, 
such that $[A,f]^G[g]=[A,f]^G[g']$.
By enlarging $V$, if necessary, we can assume that the two composites
$f(V)\circ g$ and $f(V)\circ g'$ are $G$-homotopic.
A choice of such a homotopy specifies an equivariant lifting problem on the left:
\[ \xymatrix@C=12mm{
A\times\{0,1\}\ar[r]^-{g,g'} \ar[d] & X(V) \ar[d]^{f(V)} &
A\times \{0,1\} \ar[r]^-{g,g'}\ar[d] & X(V) \ar[r]^-{X(\varphi)} & 
X(W) \ar[d]^{f(W)} 
\\
A\times[0,1]\ar[r]_-\beta & Y(V) & 
A\times[0,1]\ar[r]_-\beta \ar@{-->}[urr]^-(.4)\lambda & 
Y(V) \ar[r]_-{Y(\varphi)} & Y(W)  }\]
Proposition \ref{prop:F-telescope criterion}~(ii) 
provides a $G$-equivariant linear isometric embedding $\varphi:V\to W$ 
and a lift $\lambda$ on the right hand side such that $\lambda(-,0)=g$, $\lambda(-,1)=g'$
and $f(W)\circ\lambda$ is $G$-homotopic, relative $A\times\{0,1\}$, 
to $Y(\varphi)\circ\beta$.
As in the first part, we use a $G$-equivariant 
linear isometric embedding $j:W\to \Uc_G$,
extending the inclusion of $V$, to transform $\lambda$ into the $G$-homotopy
\[ X(j)\circ\lambda\ : \ A\times[0,1] \ \to \ X(j(W))  \]
that connects the images of $g$ and $g'$ in $X(j(W))$.
This shows that $[g]=[g']$ in $[A,X]^G$, so $[A,f]^G$ is also injective.

(iii) Products of orthogonal spaces are formed objectwise, so the canonical map
\[  [A, (X\times Y)(V)]^G \ \to \ [A,X(V)]^G \times[A, Y(Y)]^G \]
is bijective for every $G$-representation $V$.
Filtered colimits commute with finite products, so the claim follows
by passage to colimits over the poset $s(\Uc_G)$.
\end{proof}

\begin{eg}\label{eg:[A, B_gl G]}\index{subject}{global classifying space}
We specialize to the case where $Y=B_{\gl}G$ is the global classifying space
of a compact Lie group $G$. Proposition \ref{prop:classifying classifies} 
above already gave an explanation for the name `global classifying space'
by exhibiting $(B_{\gl} G)(\Uc_K)$ as a classifying space for 
principal $(K,G)$-bundles over paracompact $K$-spaces.
We now reinterpret this result as follows.

We choose a faithful $G$-representation $V$, so that $B_{\gl}G=\bL_{G,V}$.
We let $A$ be a compact $K$-space and consider the composite
\[ [A, B_{\gl}G]^K \ \xra{\ \iso \ } \ 
[A, (B_{\gl}G)(\Uc_K) ]^K \ \xra{\ \iso\ } \  \Prin_{(K,G)}(A)\ ,\]
where the first map is the bijection of Proposition \ref{prop:[A,Y]^G of closed}~(i),
exploiting that the orthogonal space $\bL_{G,V}$ is closed.
The second map is the bijection provided 
by Proposition \ref{prop:classifying classifies}. 
The composite bijection
\[  [A,  B_{\gl} G]^K \ \iso \ \Prin_{(K,G)}(A) \]
sends the class represented by a continuous $K$-equivariant map $f:A\to \bL(V,W)/G$
to the class of the pullback principal $(K,G)$-bundle $f^*q$ over $A$,
where $q:\bL(V,W)\to\bL(V,W)/G$ is the projection.
\end{eg}

Now we specialize the equivariant homotopy sets $[A,Y]^G$ to
the case $A=\{\ast\}$ of a one-point $G$-space, and then give it
a new name.

\begin{defn}\label{def:homotopy set}
Let $G$ be a compact Lie group.
The {\em $G$-equivariant homotopy set} of an orthogonal 
space $Y$\index{subject}{equivariant homotopy set!of an orthogonal space} 
is the set
\begin{equation}\label{eq:define_pi_0_set}
   \pi_0^G(Y)  \ = \ \colim_{V\in s(\Uc_G)}\,  \pi_0( Y(V)^G)  \ .
\end{equation}  
\end{defn}

Specializing Proposition \ref{prop:[A,Y]^G of closed} 
to a one-point $G$-space yields:

\begin{cor}\label{cor:pi_0^G of closed} 
Let $G$ be a compact Lie group.
\begin{enumerate}[\em (i)]
\item For every closed orthogonal space $Y$ the canonical map
\[ \pi_0^G(Y)\ \to \ \pi_0 ( Y(\Uc_G)^G ) \]
is bijective.
\item
Let $\Fc$ be a global family and $f:X\to Y$ an $\Fc$-equivalence of orthogonal spaces.
Then for every compact Lie group $G$ in $\Fc$ the induced map
\[ \pi_0^G(f)\ : \ \pi_0^G(X)\ \to \ \pi_0^G(Y) \]
of equivariant homotopy sets is bijective.
\end{enumerate}
\end{cor}

As the group varies, the homotopy sets $\pi_0^G(Y)$ 
have contravariant functoriality in $G$:
every continuous group homomorphism $\alpha:K\to G$
between compact Lie groups 
induces a restriction map $\alpha^*:\pi_0^G(Y)\to\pi_0^K(Y)$,
as we shall now explain.
We denote by $\alpha^*$ the
restriction functor from $G$-spaces to $K$-spaces 
(or from $G$-representations to $K$-representations)
along $\alpha$, i.e., $\alpha^* Z$ (respectively $\alpha^* V$)
is the same topological space as $Z$ 
(respectively the same inner product space as $V$) endowed with
$K$-action via
\[  k\cdot z \ = \ \alpha(k)\cdot z \ .  \]
Given an orthogonal space $Y$,
we note that for every $G$-representation $V$, the $K$-spaces 
$\alpha^*(Y(V))$ and $Y(\alpha^*V)$ are equal (not just isomorphic).

The restriction $\alpha^*(\Uc_G)$ is a $K$-universe, 
but if $\alpha$ has a non-trivial kernel, then this $K$-universe is 
not complete. When $\alpha$ is injective,
then $\alpha^*(\Uc_G)$ is a complete $K$-universe
(by Remark \ref{rk:complete universes restrict}),  but typically
different from the chosen complete $K$-universe $\Uc_K$. 
To deal with this we explain how a $G$-fixed point $y\in Y(U)^G$,
for an arbitrary $G$-representation $U$,
gives rise to an unambiguously defined element $\td{y}$ in $\pi_0^G(Y)$.
The point here is that $U$ need not be a subrepresentation
of the chosen universe $\Uc_G$
and the resulting class does not depend on any additional choices.

To construct $\td{y}$ we choose a $G$-equivariant linear isometry $j:U\to V$
onto a $G$-subrepresentation $V$ of $\Uc_G$.
Then $Y(j)(y)$ is a $G$-fixed point of $Y(V)$, so we obtain an element 
\[  \td{y} \ = \ [Y(j)(y)] \ \in \ \pi_0^G (Y)\ . \]
It is crucial, but not completely obvious, 
that $\td{f}$ does not depend on the choice of isometry $j$.

\begin{prop}\label{prop:universal colimit spaces}
Let $Y$ be an orthogonal space, $G$ a compact Lie group,
$U$ a $G$-representation and $y\in Y(U)^G$ a $G$-fixed point.
\begin{enumerate}[\em (i)]
\item The class $\td{y}$ in $\pi_0^G(Y)$ is independent of the choice 
of linear isometry from $U$ to a subrepresentation of $\Uc_G$.
\item For every $G$-equivariant linear isometric embedding 
$\varphi:U\to W$ the relation
\[ \td{Y(\varphi)(y)} \ = \ \td{y} \text{\qquad holds in\quad $\pi_0^G(Y)$.}  \]
\end{enumerate}
\end{prop}
\begin{proof}
(i) We let $j:U\to V$ and $j':U\to V'$ be two $G$-equivariant linear
isometries, with $V,V'\in s(\Uc_G)$.
We choose a third $G$-equivariant linear isometry $j'':U\to V''$
such that $V''\in s(\Uc_G)$ and $V''$ is orthogonal to both $V$ and $V'$.
We let $W$ be the span of $V, V'$ and $V''$ inside $\Uc_G$. 
We can then view $j, j'$ and $j''$ as equivariant linear isometric embeddings
from $U$ to $W$. 

Since the images of $j$ and $j''$ are orthogonal,
they are homotopic through $G$-equivariant linear isometric embeddings into $W$, 
via the homotopy $H:U\times [0,1]\to W$ given by
\[ H(u,t)\ =\ \sqrt{1-t^2}\ \cdot j(u)  \ + \ t\cdot j''(u) \ .\]
By the same argument, $j'$ and $j''$ 
are homotopic through $G$-equivariant linear isometric embeddings.
In particular, $j$ and $j'$ are homotopic to each other. 
If $H(-,t):U\to W$  is a continuous 1-parameter family of 
$G$-equivariant linear isometric embeddings from $j$ to $j'$, then 
\[ t\ \longmapsto \ Y(H(-,t))(y) \]
is a path in $Y(W)^G$ from $Y(j)(y)$ to $Y(j')(y)$, so
$[Y(j)(y)]=[Y(j')(y)]$ in $\pi_0^G(Y)$.

(ii) If $j:W\to V$ is an equivariant linear isometry with $V\in s(\Uc_G)$,
we define $\bar V=j(\varphi(U))$ and we let $k:U\to\bar V$ be the equivariant
linear isometry that is defined by $k(u)=j(\varphi(u))$
(i.e., $k$ is essentially $j\circ\varphi$, but with range $\bar V$ instead of $V$).
Then
\[  \td{Y(\varphi)(y)} \ = \ [Y(j)(Y(\varphi)(y))] \ = \ [Y(j\varphi)(y)] 
\ = \ [Y(k)(y)] \ = \ \td{y} \ . \qedhere \]
\end{proof}

We can now define the {\em restriction map}\index{subject}{restriction map!of equivariant homotopy sets} associated to
a continuous group homomorphism $\alpha:K\to G$ by
\begin{equation}  \label{eq:define_alpha^*}
 \alpha^* \ : \ \pi^G_0(Y) \ \to \ \pi^K_0(Y) \ ,\quad
[y]\ \longmapsto \ \td{y}\ .   
\end{equation}
This makes sense because every $G$-fixed point of $Y(V)$ is also
a $K$-fixed point of $\alpha^*(Y(V))=Y(\alpha^*V)$.
For a second continuous group homomorphism $\beta:L\to K$ we have
\[ \beta^*\circ \alpha^* \ = \ (\alpha\beta)^* \ : \ \pi_0^G (Y) \ \to \ \pi_0^L(Y) \ .  \]
Clearly, restriction along the identity homomorphism is the identity,
so we have made the collection of equivariant homotopy sets $\pi_0^G(Y)$
into a contravariant functor in the group variable.

An important special case of the restriction homomorphisms
are conjugation maps.
Here we consider a closed subgroup $H$ of $G$, an element $g\in G$ and denote by 
\[ c_g \ : \ H \ \to \  H^g\ , \quad c_g(h)\ =\ g^{-1}h g\ \] 
the conjugation homomorphism, where 
$H^g=\{g^{-1} h g\ | \ h\in H\}$ is the conjugate subgroup.
As any group homomorphism, 
$c_g$ induces a restriction map\index{symbol}{$c_g$ - {conjugation by $g$}}
\[ g_\star \ = \ (c_g)^* \ : \ \pi^{H^g}_0 (Y) \ \to\  \pi^H_0(Y)      \]
of equivariant homotopy sets.
For $g,\bar g\in G$ we have $c_{g\bar g}=c_{\bar g}\circ c_g:H\to H^{g\bar g}$
and thus 
\[ (g\bar g)_\star\ = \  (c_{g\bar g})^*\ = \ (c_{\bar g}\circ c_g)^*\ = \
(c_g)^*\circ (c_{\bar g})^*\ = \ g_\star\circ \bar g_\star \ .\]
A key fact is that inner automorphisms act trivially, 
i.e., the restriction map $g_\star$ is the identity on $\pi_0^G (Y)$.
So the action, by the restriction maps, 
of the automorphism group of $G$ on $\pi_0^G(Y)$ 
factors through the outer automorphism group.

\begin{prop}
For every orthogonal space $Y$, every compact Lie group $G$, 
and every $g\in G$, the conjugation map\index{subject}{conjugation map!on equivariant homotopy sets}
$g_\star:\pi_0^G(Y) \to \pi_0^G(Y)$ is the identity.
\end{prop}
\begin{proof} 
We consider a finite-dimensional $G$-subrepresentation $V$ 
of $\Uc_G$ and a $G$-fixed point $y\in Y(V)^G$ that represents an element in $\pi_0^G(Y)$. 
Then the map $l_g:c_g^*(V)\to V$ given by left multiplication by $g$
is a $G$-equivariant linear isometry.
So 
\[ g_\star[y]\ = \ (c_g)^*[y]\ = \ [Y(l_g)(y)]\ = \ [g\cdot y]\ = \ [y] \ ,\]
by the very definition of the restriction map. The third equation is the definition
of the $G$-action on $Y(V)$ through the $G$-action on $V$.
The fourth equation is the hypothesis that $y$ is $G$-fixed.
\end{proof}

We note that because inner automorphisms act as the identity,
the restriction map $\alpha^*$ only depends on the homotopy class of $\alpha$.
More precisely, suppose that $\alpha,\alpha':K\to G$ 
are homotopic through continuous group homomorphisms. 
Then $\alpha$ and $\alpha'$ belong to the same path component of the space $\hom(K,G)$
of continuous homomorphisms, and so they are conjugate by an element of $G$,
compare Proposition \ref{prop:components of hom(K,G)}.

\medskip

We denote by Rep\index{symbol}{  $\Rep$ - {category of compact Lie groups and conjugacy classes of homomorphisms}} 
the category whose objects are the compact Lie groups
and whose morphisms are conjugacy classes of continuous group homomorphisms.
We can summarize the discussion thus far by saying that for every
orthogonal space $Y$ the restriction maps
make the equivariant homotopy sets
$\{\pi_0^G(Y)\}$ into a functor
\[ \upi_0(Y) \ : \ \Rep^{\op} \ \to \text{(sets)} \ . \]
We will refer to such a functor as a {\em Rep-functor}.\index{subject}{Rep-functor}

\medskip

Our next aim is to show that the homotopy $\Rep$-functor
$\upi_0(B_{\gl}G)$ of a global classifying space is represented by $G$.
For every $G$-representation $V$ 
we define the {\em tautological class}\index{subject}{tautological class!unstable|see{unstable tautological class}}\index{subject}{unstable tautological class}\index{symbol}{$u_{G,V}$ - {tautological class in $\pi_0^G(\bL_{G,V})$}}\index{subject}{global classifying space}
\begin{equation}\label{eq:tautological_class}
 u_{G,V}\ \in \ \pi_0^G(\bL_{G,V})   
\end{equation}
as the path component of the $G$-fixed point
\[  \Id_V\cdot G \ \in \ (\bL(V,V)/G)^G = (\bL_{G,V}(V))^G\ , \]
the $G$-orbit of the identity of $V$.

\begin{prop}\label{prop:fix of global classifying} 
Let $G$ and $K$ be compact Lie groups and $V$ a faithful $G$-representation.
\begin{enumerate}[\em (i)]
\item 
The $K$-fixed point space  $\left(\bL_{G,V}(\Uc_K)\right)^K$ is a 
disjoint union, indexed by conjugacy classes
of continuous group homomorphisms $\alpha:K\to G$,
of classifying spaces of the centralizer  of the image of $\alpha$.
\item
The map
\[  \Rep(K,G) \ \to \ \pi_0^K(\bL_{G,V})\ , \quad
[\alpha:K\to G] \ \longmapsto \ \alpha^*(u_{G,V})\]
is bijective. 
\end{enumerate}
\end{prop}
\begin{proof}
Part~(i) works for any universal $(K\times G)$-space $E$
for the family $\Fc(K;G)$ of graph subgroups, for example for $E=\bL(V,\Uc_K)$. 
The argument can be found in Theorem~2.17 of \cite{lashof-equivariant bundles}
or Proposition~5 of \cite{lewis-may-mcclure-classifying}.\index{subject}{universal space!for the family of graph subgroups}
We repeat the proof for the convenience of the reader.
For a continuous homomorphism $\alpha:K\to G$, we let $C(\alpha)$
denote the centralizer, in $G$, of the image of $\alpha$, and we set
\[ E^\alpha \ = \   \{x\in E\ |\ (k,\alpha(k))\cdot x = x\text{ for all $k\in K$}\}\ ,\]
the space of fixed points of the graph of $\alpha$.
Since the $G$-action on the universal space $E$ is free,
Proposition \ref{prop:fix of free cofibration} provides
a homeomorphism
\[ \coprod \alpha^\flat \ : \ 
{\coprod}_{\td{\alpha}}\   E^{\alpha}/C(\alpha)\ \to \  (E/G)^K \ , \]
where the disjoint union is indexed by conjugacy classes of continuous homomorphisms.
The graph of $\alpha$ belongs to the family $\Fc(K;G)$,
so $E^\alpha$ is a contractible space. 
The action of $C(\alpha)$ on $E^\alpha$ is a restriction of the $G$-action on $E$,
hence free. Since $E$ is $(K\times G)$-cofibrant, 
the fixed point space $E^\alpha$ is cofibrant for the action
of the normalizer (inside $K\times G$) of the graph of $\alpha$,
by Proposition \ref{prop:cofamily pushout property}.
Hence $E^\alpha$ is also cofibrant as a $C(\alpha)$-space, by 
Proposition \ref{prop:cofibrancy preservers}~(i). 
So for every homomorphism $\alpha$ the space 
$E^{\alpha}/C(\alpha)$ is a classifying space for the group $C(\alpha)$.
This shows part~(i).

(ii) Since the classifying space of a topological group is connected,
part~(i) identifies the path components of $\left(\bL_{G,V}(\Uc_K)\right)^K$
with the conjugacy classes of continuous homomorphisms $\alpha:K\to G$.
The bijection sends the class of $\alpha$ to $\alpha^* (u_{G,V})$.
The claim then follows by applying Corollary \ref{cor:pi_0^G of closed}~(i).
\end{proof}

Now we show that the restriction maps along continuous group
homomorphisms give all natural operations between equivariant homotopy sets
of orthogonal spaces. 
We will want to do similar calculations several other times
in this book, so we state the argument 
in the more general situation of a category $\Cc$ related to 
the category of orthogonal spaces by an adjoint functor pair.
Our present context is the degenerate case $\Cc=\spc$ and the identity functors.
Later we will also consider the cases of the categories
of ultra-commutative monoids, of orthogonal spectra and 
of ultra-commutative ring spectra.

We recall that the restriction morphism
\[ \rho_{G,V,W}\ = \ \rho_{V,W}/G \ : \ \bL_{G,V\oplus W}\ = \ \bL(V\oplus W,- )/G\ \to \ 
\bL(W,-)/G \ = \ \bL_{G,W} \]
restricts the orbit of a linear isometric embedding from $V\oplus W$
to the second summand $W$.
This morphism is a global equivalence of orthogonal spaces by 
Proposition \ref{prop:free_orthogonal_space}~(ii),
as long as $G$ acts faithfully on $W$.

\begin{prop}\label{prop:pi_0^G representability} 
Let $\Cc$ be a category and
\[ \xymatrix{ \Lambda \ : \ \spc \ \ar@<.4ex>[r] & \ \Cc \ : \ U \ar@<.4ex>[l] } \]
an adjoint functor pair such that the composite functor
$U\Lambda:\spc\to\spc$ is continuous.
Suppose in addition that
for every compact Lie group $G$, all $G$-representations $V$ 
and all non-zero faithful $G$-representations $W$, the morphism of $\Rep$-functors
\[ \upi_0( U\Lambda(\rho_{G,V,W})) \ : \ 
\upi_0( U\Lambda(\bL_{G,V\oplus W})) \ \to \ \upi_0( U\Lambda(\bL_{G,W})) \]
is an isomorphism.
Let $W$ be a non-zero faithful $G$-representation, and set
\[ u_{G,W}^\Cc\ = \ \eta_*(u_{G,W})\ \in \ \pi_0^G( U\Lambda(\bL_{G,W}))\ ,\]
where $\eta:\bL_{G,W}\to U\Lambda(\bL_{G,W})$ is the adjunction unit.
Then for every compact Lie group $K$,
evaluation at the class $u^\Cc_{G,W}$ is a bijection
\[  \Nat_{\Cc\rightarrow\emph{(sets)}}(\pi_0^G\circ U,\,\pi_0^K\circ U) 
\ \to \ \pi_0^K( U\Lambda( \bL_{G,W} ) ) \ , \quad
\tau\ \longmapsto \ \tau(u^\Cc_{G,W}) \]
between the set of natural transformations 
from the functor $\pi_0^G\circ U$ to $\pi_0^K\circ U$, 
and the set $\pi_0^K( U\Lambda(\bL_{G,W}))$.
\end{prop}
\begin{proof}
To show that the evaluation map is injective 
we show that every natural transformation $\tau:\pi_0^G\circ U\to \pi_0^K\circ U$ 
is determined by the element $\tau(u^\Cc_{G,W})$.
We let $X$ be any object of $\Cc$ and $[x]\in\pi_0^G(U X)$
a $G$-equivariant homotopy class, represented by a $G$-fixed point
\[ x \ \in \ (U X)(V)^G \ ,\]
for some $G$-representation $V$. We can stabilize with the representation $W$
and obtain another representative
\[ (U X)(i)(x) \ \in \ (U X)(V\oplus W)^G \ ,\]
where $i:V\to V\oplus W$ is the embedding of the first summand. 
This $G$-fixed point is adjoint to a morphism of orthogonal spaces
\[ \hat x\ : \ \bL_{G,V\oplus W} \ \to \ U X \]
and hence adjoint to a $\Cc$-morphism
\[  x^\flat \ : \ \Lambda(\bL_{G,V\oplus W}) \ \to \ X \]
that satisfies
\[ \pi_0^G(U x^\flat)(u^\Cc_{G,V\oplus W}) \ = \ [x] \quad \text{in } \ \pi_0^G(U X) \ .\]
The restriction morphism of orthogonal spaces 
$\rho_{G,V,W} :\bL_{G,V\oplus W}\to \bL_{G,W}$ sends $u_{G,V\oplus W}$ to $u_{G,W}$.
The morphism of orthogonal spaces
\[ U\Lambda(\rho_{G,V,W})\ : \ U\Lambda( \bL_{G,V\oplus W})\ \to \ U\Lambda( \bL_{G,W} )\]
then sends $u^\Cc_{G,V\oplus W}$ to $u^\Cc_{G,W}$.
The diagram
\[ \xymatrix@C=25mm{
\pi_0^G(U\Lambda( \bL_{G,W}))\ar[d]_\tau &
\pi_0^G(U\Lambda(\bL_{G,V\oplus W}))\ar[d]_\tau 
\ar[l]_-{\pi_0^G(U\Lambda(\rho_{G,V,W}))}^-\iso\ar[r]^-{\pi_0^G(U x^\flat)} &
\pi_0^G(U X)\ar[d]^\tau \\
\pi_0^K(U\Lambda(\bL_{G,W})) &
\pi_0^K(U\Lambda( \bL_{G,V\oplus W})) \ar[l]^-{\pi_0^K(U\Lambda(\rho_{G,V,W}))}_-\iso
\ar[r]_-{\pi_0^K(U x^\flat)} & \pi_0^K(U X)} \]
commutes and the two left horizontal maps are bijective by hypothesis.
Since 
\[ [x] \ = \ \pi_0^G(U x^\flat)(\pi_0^G(U\Lambda(\rho_{G,V,W}))^{-1}(u^\Cc_{G,W})) \ ,\]
naturality yields that
\begin{align*}
    \tau[x] \ &= \ \tau(\pi_0^G(U x^\flat)(\pi_0^G(U\Lambda(\rho_{G,V,W}))^{-1}(u^\Cc_{G,W})))\\ 
&= \ \pi_0^K(U x^\flat)(\pi_0^K(U\Lambda(\rho_{G,V,W}))^{-1}(\tau(u^\Cc_{G,W})))\ . 
\end{align*}
This shows that the transformation $\tau$ is determined by the value
$\tau(u^\Cc_{G,W})$.

It remains to construct, for every element 
$z\in \pi_0^K(U\Lambda(\bL_{G,W}))$, a natural transformation
$\tau:\pi_0^G\circ U\to \pi_0^K\circ U$ with $\tau(u^\Cc_{G,W})=z$.
The previous paragraph dictates what to do: 
we represent a given class in $\pi_0^G(U X)$ by
a $G$-fixed point $x\in(U X)(V\oplus W)^G$ and set
\[   \tau[x] \ = \ \pi_0^K(U x^\flat)(\pi_0^K(U\Lambda(\rho_{G,V,W}))^{-1}(z)) \ .\]
We must verify that the element $\tau[x]$ 
is independent of the representative.
If $y\in(U X)(V\oplus W)^G$ is in the same path component as $x$,
then any path adjoins to a homotopy of morphisms of orthogonal spaces
\[ H \ : \ \bL_{G,V\oplus W}\times [0,1]\  \to U X \]
from $\hat x$ to $\hat y$.
Since the functor $U\Lambda$ is continuous, the composite morphisms
\[ U\Lambda(\bL_{G,V\oplus W})\ \xra{U\Lambda(H(-,t))}\ U\Lambda(U X)\ 
\xra{U(\epsilon_X)} \ U X \]
form a continuous 1-parameter family of morphisms for $t\in[0,1]$, 
where $\epsilon_X:\Lambda(U X)\to X$ is the counit of the adjunction.
This witnesses that the morphism $U x^\flat$ is homotopic to the morphism $U y^\flat$.
So $\pi_0^K(U x^\flat)=\pi_0^K(U y^\flat)$, and the class $\tau[x]$
is independent of the representative in the given path component of $(U X)(V\oplus W)^G$.

Now we let $V'$ be another $G$-representation and $\varphi:V\to V'$ 
a $G$-equivariant linear isometric embedding. 
Then
\[ y \ = \ (U X)(\varphi\oplus W)(x)\text{\qquad in\quad} (U X)(V'\oplus W)^G \]
is another representative of the class $[x]$.
The restriction morphism 
\[ \varphi^\sharp\ = \ \bL(\varphi\oplus W,-)/G\ : \ 
\bL_{G,V'\oplus W} \ \to \ \bL_{G,V\oplus W} \]
makes the following diagram of orthogonal spaces commute:
\[ \xymatrix@C=15mm{ 
\bL_{G,W} & 
\bL_{G,V'\oplus W} \ar[d]_{\varphi^\sharp}\ar[r]^-{\hat y} \ar[l]_-{\rho_{G,V',W}}& U X\\
&\bL_{G,V\oplus W}\ar[ur]_{\hat x} \ar[ul]^{\rho_{G,V,W}}& } \]
Passing to adjoints and applying $U$ yields another commutative diagram:
\[ \xymatrix@C=15mm{ 
U\Lambda(\bL_{G,W}) & 
U\Lambda(\bL_{G,V'\oplus W}) \ar[d]_{U\Lambda(\varphi^\sharp)}
\ar[r]^-{U y^\flat} \ar[l]_-{U\Lambda(\rho_{G,V',W})}& U X\\
&U\Lambda(\bL_{G,V\oplus W})\ar[ur]_{U x^\flat} \ar[ul]^{U\Lambda(\rho_{G,V,W})}& } \]
So
\begin{align*}
  \pi_0^K(U x^\flat)\circ \pi_0^K(U\Lambda(\rho_{G,V,W}))^{-1}\ &= \
 \pi_0^K(U x^\flat)\circ \pi_0^K(U\Lambda(\varphi^\sharp))\circ\pi_0^K(U\Lambda(\rho_{G,V',W}))^{-1}\\
&= \  \pi_0^K(U y^\flat)\circ\pi_0^K(U\Lambda(\rho_{G,V',W}))^{-1} \ ,
\end{align*}
and hence the class $\tau[x]$ remains unchanged upon stabilization of $x$
along $\varphi$.
Altogether this shows that $\tau[x]$ is well-defined.

Naturality of $\tau$ is straightforward: if $\psi:X\to Y$ is a $\Cc$-morphism
and $x\in(U X)(V\oplus W)^G$ is as above,
then $(U \psi)(V\oplus W)(x)\in (U Y)(V\oplus W)^G$ 
represents the class $\pi_0^G(U\psi)[x]$.
Moreover, the adjoint of $(U\psi)(V\oplus W)(x)$ coincides with the composite
\[\Lambda(\bL_{G,V\oplus W}) \ \xra{\ x^\flat\ }\ X \ \xra{\ \psi\ }\ Y\ . \]
So naturality follows:
\begin{align*}
  \tau(\pi_0^G(U\psi)[x])\ &= \ 
\pi_0^K(U\psi\circ U x^\flat)(\pi_0^K(U\Lambda(\rho_{G,V,W}))^{-1}(z)) \\
&= \ \pi_0^K(U\psi)\left( \pi_0^K(U x^\flat)(\pi_0^K(U\Lambda(\rho_{G,V,W}))^{-1}(z))\right) \
= \ \pi_0^K(U\psi)(\tau[x]) \ .
\end{align*}
Finally, the class $u^\Cc_{G,W}$ 
is represented by the $G$-fixed point 
\[ \eta_{\bL_{G,W}}(\Id_W\cdot G)\ \in \ (U\Lambda(\bL_{G,W}))(W)^G\ , \]
which is adjoint to the identity of $\Lambda(\bL_{G,W})$. Hence
$\tau(u_{G,W}^\Cc)= \pi_0^K(\Id)(z) =z$.
\end{proof}

\begin{cor}\label{cor:operations spc}
Every natural transformation $\pi_0^G \to \pi_0^K$
of set valued functors on the category of orthogonal spaces is of the form
$\alpha^*$
for a unique conjugacy class of continuous group homomorphism $\alpha:K\to G$.  
\end{cor}
\begin{proof}
We let $W$ be any non-zero faithful $G$-representation.
The composite
\[ \Rep(K,G) \ \xra{[\alpha]\mapsto\alpha^*} \ \Nat(\pi_0^G,\pi_0^K) 
\ \xra{\ \ev\ } \ \pi_0^K(\bL_{G,W})\]
is bijective by Proposition \ref{prop:fix of global classifying}~(ii), 
where the second map is evaluation at the tautological class $u_{G,W}$.
The evaluation map is bijective by 
Proposition \ref{prop:pi_0^G representability}, applied to $\Cc=\spc$
and the identity adjoint functor pair.
So the first map is bijective as well.
\end{proof}

\begin{construction}\label{con:pairing equivalence homotopy spaces}
Given two orthogonal spaces $X$ and $Y$, we endow the equivariant homotopy sets
with a pairing\index{subject}{product!on equivariant homotopy sets}
\begin{equation}\label{eq:def_times}
\times \ : \
\pi_0^G(X) \ \times \ \pi_0^G(Y) \ \to \ \pi_0^G(X\boxtimes Y) \ ,
\end{equation}
where $G$ is any compact Lie group.
We suppose that $V$ and $W$ are $G$-represen\-tations
and $x\in X(V)^G$ and $y\in Y(W)^G$ are fixed points that
represent classes in $\pi^G_0(X)$ respectively $\pi^G_0(Y)$.
We denote by $x\times y$ the image of the $G$-fixed point $(x,y)$
under the $G$-map
\begin{align*}
i_{V,W}\ : \   X(V)\times Y(W)\ \to \ (X\boxtimes Y)(V\oplus W) 
\end{align*}
that is part of the universal bimorphism.
If $\varphi:V\to V'$ and $\psi:W\to W'$ are equivariant linear isometric embeddings,
then 
\begin{align*}
   X(\varphi)(x)\times Y(\psi)(y) \ &= \  i_{V',W'}(X(\varphi)(x),Y(\psi)(y))\\ 
&= \  (X\boxtimes Y)(\varphi\oplus\psi)(i_{V,W}(x,y)) 
\ = \  (X\boxtimes Y)(\varphi\oplus\psi)(x\times y) \ .
\end{align*}
So by Proposition \ref{prop:universal colimit spaces}~(ii)
the classes $\td{x\times y}$ and $\td{ X(\varphi)(x)\times Y(\psi)(y) }$
coincide in $\pi_0^G(X\boxtimes Y)$. The upshot is that the assignment
\[  [x]\times [y]\ = \ \td{x\times y}\ \in \ \pi_0^G(X\boxtimes Y) \]
is well-defined.
\end{construction}

The pairings of equivariant homotopy sets have several expected properties
that we summarize in the next proposition.

\begin{prop}\label{prop:unstable product properties} 
  Let $G$ be a compact Lie group and $X, Y$ and $Z$ orthogonal spaces.
  \begin{enumerate}[\em (i)]
  \item {\em (Unitality)}
    Let $1\in\pi_0^G(\mathbf 1)$ be the unique element.
    Then $1\times x=x=x\times 1$ for all $x\in\pi_0^G(X)$.
  \item {\em (Associativity)} For all classes $x\in \pi_0^G ( X )$, $y\in \pi_0^G( Y )$
    and $z\in\pi_0^G( Z )$ the relation
    \[ a_*( (x\times y)\times z) \ = \ x\times (y\times z) \]
    holds in $\pi_0^G(X\boxtimes (Y\boxtimes Z))$,
    where $a:(X\boxtimes Y)\boxtimes Z\iso X\boxtimes (Y\boxtimes Z)$ 
    is the associativity isomorphism.
  \item {\em (Commutativity)} For all classes $x\in \pi_0^G ( X )$ and $y\in \pi_0^G( Y )$
    the relation
    \[ \tau_*( x\times y ) \ = \  y\times x \]
    holds in $\pi_0^G( Y\boxtimes X)$,
    where $\tau:X\boxtimes Y\to Y\boxtimes X$ is the symmetry isomorphism.
  \item {\em (Restriction)} For all classes $x\in \pi_0^G ( X )$ and $y\in \pi_0^G( Y )$
    and all continuous homomorphisms $\alpha:K\to G$ the relation
    \[ \alpha^*(x)\times\alpha^*(y)\ = \ \alpha^*(x\times y) \]
    holds in $\pi_0^K(X\boxtimes Y)$.
      \end{enumerate}
\end{prop}
\begin{proof}
The unitality property~(i), the associativity property~(ii) 
and compatibility with restriction~(iv)
are straightforward from the definitions.
Part~(iii) exploits that the square
\[\xymatrix@C=18mm{ 
X(V)\times Y(W) \ar[r]^-{i_{V,W}}\ar[d]_{\text{twist}}&
(X\boxtimes Y)(V\oplus W)\ar[d]^{\tau(\chi_{V,W})} \\
Y(W)\times X(V) \ar[r]_{i_{W,V}} &
(Y\boxtimes X)(W\oplus V)}\]
commutes. The image of $(x,y)$ under the upper right composite represents
$\tau_*(x\times y)$, whereas the image of $(y,x)$ 
under the lower left composite represents $y\times x$, so
$\tau_*( x\times y ) =  y\times x$.
\end{proof}

\begin{rk}[External versus internal products]
If $G$ and $K$ are two compact Lie groups, we can define an `external' form
of the pairing \eqref{eq:def_times} as the composite
\begin{align}\label{eq:def_external_product}
 \pi_0^G(X)\times \pi_0^K(Y) \ \xra{\ p_1^*\times p_2^*\ } \
 \pi_0^{G\times K}(X)\times \pi_0^{G\times K}(Y)\  \xra{\ \times\ }\
\pi_0^{G\times K}(X\boxtimes Y)\ ,
\end{align}
where $p_1:G\times K\to G$ and $p_2:G\times K\to K$ are the two projections.
These external pairings also satisfy various naturality, unitality, associativity and
commutativity properties that we do not spell out.
On the other hand, the internal pairing \eqref{eq:def_times}
can be recovered from the external products \eqref{eq:def_external_product}
by taking $G=K$ and restricting along the diagonal
embedding $\Delta_G:G\to G\times G$.
Indeed, the $p_1\circ\Delta_G=p_2\circ\Delta_G=\Id_G$, and hence
\begin{align*}
    \Delta_G^* (p_1^*(x)\times p_2^*(y))\ = \ 
  \Delta_G^*( p_1^*(x))\times \Delta_G^*(p_2^*(y)) \ = \ x\times y\ .
\end{align*}
\end{rk}

Theorem \ref{thm:box to times}~(i) and the fact that
the functor $\pi_0^G$ commutes with finite products 
(Proposition \ref{prop:[A,Y]^G of closed}~(iii) for $A=\ast$)
imply:

\begin{cor}\label{cor-pi_0 of box}
For every compact Lie group $G$ and all orthogonal spaces $X$ and $Y$ the three maps
\[  \pi_0^G(X)\times \pi_0^G(Y)\ \xra{\ \times\ }\ 
\pi_0^G( X\boxtimes Y) \ \xra{(\rho_{X,Y})_*} \ 
\pi_0^G( X\times Y) \ \xra{((p_X)_*,(p_Y)_*)} \  \pi_0^G(X)\times \pi_0^G(Y)
 \]
are bijective, where $p_X:X\times Y\to X$ and $p_Y:X\times Y\to Y$
are the projections. Moreover, the composite is the identity.
\end{cor}

\begin{construction}[Infinite box products]\label{con:infinite box}
We end this section with a generalization of the previous corollary to
`infinite box products' of based orthogonal spaces, but we first have
to clarify what we mean by that. We let $I$ be an indexing set and $\{X_i\}_{i\in I}$
a family of based orthogonal spaces, i.e., each equipped with a 
distinguished point $x_i\in X_i(0)$.
If $K\subset J$ are two nested, {\em finite} subsets of $I$, then
the basepoints of $X_k$ for $k\in J-K$ provide a morphism
\begin{equation}  \label{eq:J2K}
 \boxtimes_{k\in K} X_k \ \to \ \boxtimes_{j\in J} X_j \ .  
\end{equation}
In terms of the universal property of the box product, this morphism arises
from the maps
\begin{align*}
  {\prod}_{k\in K}\, X_k(V_k) \ \to \ {\prod}_{k\in K} \, X_k(V_k) &\times 
{\prod}_{j\in J-K} \, X_j(0) \
\to \ ({\boxtimes}_{j\in J}\, X_j)\left({\oplus}_{k\in K}\, V_k\right) \ ,
\end{align*}
where the second map is part of the universal multi-morphism.
We can thus define the {\em infinite box product}\index{subject}{box product!of orthogonal spaces!infinite}
as the colimit of the finite box products over the filtered poset of
finite subsets of $I$:
\[ \boxtimes'_{i\in I} X_i\ = \ \colim_{J\subset I, |J|<\infty} \, \boxtimes_{j\in J} X_j \ .\]
If $I$ happens to be finite, then this recovers the iterated box product.
\end{construction}

The distinguished basepoint of $X_i$ represents a distinguished basepoint
in the equivariant homotopy set $\pi_0^G(X_i)$ for every compact Lie group $G$.
In fact, these points all arise from the basepoint in $\pi_0^e(X_i)$
by restriction along the unique homomorphism $G\to e$.
The {\em weak product} ${\prod}'_{i\in I} \pi_0^G(X_i)$
is the subset of the product consisting of all tuples $(y_i)_{i\in I}$
with the property that almost all $y_i$ are the distinguished basepoint.
Equivalently, the weak product is the filtered colimit,
over the poset of finite subsets of $I$, of the finite products.

If we iterate the pairing \eqref{eq:def_times}, it provides a multi-pairing
\[ {\prod}_{j\in J}\, \pi_0^G(X_j)\ \to \ \pi_0^G({\boxtimes}_{j\in J} \, X_j) \]
for every finite set $J$. Passing to colimits over finite subsets
of $I$ on both sides yields a map
\begin{equation}  \label{eq:pi_0^G_of_infinite_box}
   {\prod}'_{i\in I} \, \pi_0^G(X_i)\ \to \ \pi_0^G(\boxtimes'_{i\in I} X_i)\ .
\end{equation}

\begin{prop}\label{prop:pi_0 of infinite box}
Let $I$ be a set and $\{X_i\}_{i\in I}$ a family of based orthogonal spaces.
Then for every compact Lie group $G$ the map \eqref{eq:pi_0^G_of_infinite_box}
is bijective.  
\end{prop}
\begin{proof}
For every  $k\in I$ we define a `projection'
\[ \Pi_k \ : \ \boxtimes'_{i\in I} X_i \ \to \ X_k  \]
as follows. Since the infinite box product is defined as a colimit,
we must specify the `restriction' of $\Pi_k$ to $\boxtimes_{j\in J}X_j$
for every finite subset $J$ of $I$, compatibly as $J$ increases.
For $k\not\in J$ we define this restriction as the constant morphism
factoring through the basepoint of $X_k$.
For $k\in J$ we define the restriction
\[  \boxtimes_{j\in J} X_j \ \to \ X_k  \]
as the morphism corresponding, under the universal property of the box product,
to the multi-morphism with components
\[  {\prod}_{j\in J}\,  X_j(V_j) \ \xra{\text{proj}_k}\  X_k(V_k)\ 
\xra{X_k(\text{incl})}\ X_k({\oplus}_{j\in J}\, V_j)  \ .\]
Then the composite
\[ {\prod}'_{i\in I} \, \pi_0^G(X_i)\ \xra{\eqref{eq:pi_0^G_of_infinite_box}}\ 
\pi_0^G(\boxtimes'_{i\in I} X_i)\ 
\xra{\pi_0^G(\Pi_k)}\  \pi_0^G(X_k) \]
is the projection onto the $k$-th factor.
So if two tuples in the weak product have the same image
under the map \eqref{eq:pi_0^G_of_infinite_box}, they coincide.
This shows injectivity.

Now we show surjectivity. Every element of $\pi_0^G(\boxtimes'_{i\in I} X_i)$ 
is represented by a $G$-fixed point of $(\boxtimes'_{i\in I} X_i)(V)$
for some $G$-representation $V$.
Colimits of orthogonal spaces are formed objectwise, so
$(\boxtimes'_{i\in I} X_i)(V)$ is a colimit
of the spaces $(\boxtimes_{j\in J} X_j)(V)$,
formed over the filtered poset of finite subsets $J$ of $I$.
For every nested pair of finite subsets $K\subset J$
of $I$ the morphism \eqref{eq:J2K} has a retraction, by `projection'.
So at every inner product space $V$, the map
\[  (\boxtimes_{k\in K} X_k)(V) \ \to \ (\boxtimes_{j\in J} X_j)(V)  \]
is a closed embedding by Proposition \ref{prop:section is closed embedding}.
For fixed $V$, the colimit $(\boxtimes'_{i\in I} X_i)(V)$
in the category $\bT$ of compactly generated spaces
can thus be calculated in the ambient category of all topological spaces,
by Proposition \ref{prop:filtered colim in cg}~(ii).
In particular, every $G$-fixed point of $(\boxtimes'_{i\in I} X_i)(V)$
arises from a $G$-fixed point of $(\boxtimes_{j\in J} X_j)(V)$
for some finite subset $J$ of $I$.
In other words, the canonical map
\[ \colim_{J\subset I,|J|<\infty} \pi_0^G(\boxtimes_{j\in J} X_j) \ \to \ 
\pi_0^G(\boxtimes'_{i\in I} X_i)\]
is surjective. 
For finite sets $J$ the map $\prod_{j\in J}\pi_0^G(X_j)\to \pi_0^G(\boxtimes_{j\in J} X_j)$
is bijective by Corollary \ref{cor-pi_0 of box},
so this shows surjectivity.
\end{proof}

\chapter{Ultra-commutative monoids}
\label{ch-umon}

Orthogonal monoid spaces are the lax monoidal continuous functors
from the linear isometries category $\bL$ to the category of spaces.
Orthogonal monoid spaces with strictly commutative multiplication
(i.e., the lax {\em symmetric} monoi\-dal functors)
play a special role, and we honor this
by special terminology, referring to them 
as {\em ultra-commutative monoids}.
This chapter is devoted to the study of ultra-commutative monoids,
including a global model structure, an algebraic study of their
homotopy operations, and many examples.

I want to motivate the adjective `ultra-commutative'.
In various contexts of homotopy theory, 
highly structured multiplications that are
only associative or commutative up to higher coherence homotopies
can in fact be rigidified to multiplications that are strictly associative
or possibly strictly commutative.
One example is the fact that $E_\infty$-spaces can be rigidified to 
strictly commutative $\mathcal I$-space 
monoids \cite[Thm.\,1.3]{sagave-schlichtkrull-diagram};
another example is the fact that $E_\infty$-ring objects internal
to symmetric spectra can be rigidified to strictly commutative
symmetric ring spectra, see for example 
\cite[Thm.\ 1.4]{elmendorf-mandell}
and the paragraph immediately following thereafter.
More to the point of our present discussion, 
in \cite[Thm.\,1.3]{lind-diagram} Lind establishes
a Quillen equivalence between the non-equivariant homotopy theory of $E_\infty$-spaces
(i.e., spaces with an action of the linear isometries operad)
and the non-equivariant homotopy theory of commutative
orthogonal monoid spaces (there called `commutative $\mathcal I$-FCPs').

Our use of the word `ultra-commutative' is intended as a reminder
that the slogan `$E_\infty$=commutative' is no longer true in equivariant
or global contexts. 
More specifically, one can consider orthogonal monoid spaces
with an action of an $E_\infty$-operad;
up to non-equivariant equivalence, these objects model $E_\infty$-spaces,
and they can be replaced by equivalent strictly commutative orthogonal monoid spaces.
The analogous statement for global equivalences is false, i.e.,
$E_\infty$-orthogonal spaces cannot in general be replaced by
globally equivalent ultra-commutative monoids.
In fact, the definition of power operations and transfer maps requires
strict commutativity, 
and Remark \ref{rk:bO versus BO} illustrates 
how the lack of transfers obstructs ultra-commutativity.

The study of orthogonal monoid spaces goes back to 
Boardman and Vogt \cite{boardman-vogt-homotopy everything},
who introduce them as a delooping machine in a non-equivariant context. 
More precisely, they show that for every ultra-commutative monoid $R$
the space $R(\mR^\infty)$ has the structure of an `$E$-space'
(nowadays called an $E_\infty$-space), and if in addition $\pi_0(R(\mR^\infty))$
is a group, then  $R(\mR^\infty)$ is an infinite loop space.
Ultra-commutative monoids also appear, 
with an extra pointset topological hypothesis and
under the name ${\mathscr I}_*$-prefunctor, in \cite[IV Def.\,2.1]{may-quinn-ray};
in \cite{lind-diagram}, they are studied under the name
`commutative $\mathcal I$-FCPs'.

\smallskip
 
In Section \ref{sec:global model monoid spaces}
we formally define ultra-commutative monoids and establish the global model structure.
Section \ref{sec:power op spaces} is devoted to
the algebraic structure on the homotopy $\Rep$-functor $\upi_0(R)$
of an ultra-commutative monoid. 
We refer to this structure as a `global power monoid';
it consist of an abelian monoid structure
on the set $\pi_0^G(R)$  for every compact Lie group $G$, 
natural for restriction along continuous homomorphisms,
and additional structure that can equivalently be encoded
as power operations (see Definition \ref{def:power monoid})
or as transfer maps (see Construction \ref{con:transfer map}).
In this section we also show that these operations are the entire natural structure
(see Theorem \ref{thm:umon operation}).
Section \ref{sec:umon examples} collects various examples of ultra-commutative monoids:
among these are examples made from the infinite families of classical Lie groups
(orthogonal, special orthogonal, unitary, special unitary, symplectic,
spin and spin$^c$); examples consisting of Grassmannians under direct sum
of subspaces (in real, oriented, complex and quaternionic flavors);
examples made from Grassmannians under tensor product
of subspaces (in a real and complex version);
and ultra-commutative multiplicative models for global classifying spaces
of abelian compact Lie groups.

Section \ref{sec:forms of BO} is a case study of how
non-equivariant homotopy types can `fold up' into many different global
homotopy types. We define, discuss and compare different ultra-commutative
and $E_\infty$-orthogonal monoid spaces whose
underlying non-equivariant homotopy type is $B O$, 
a classifying space for the infinite orthogonal group; 
in all examples we also identify
the associated global power monoids and fixed point spaces.
Section \ref{sec:group completion and units}
discusses `units' and `group completion' of ultra-commutative monoids.
The two constructions are dual to each other, and they are
homotopically right adjoint respectively left adjoint to
the inclusion of group-like ultra-commutative monoids.
On the algebraic level of global power monoids, the topological constructions
pick out the invertible elements and perform the algebraic group completion 
respectively. A naturally occurring example of a global group
completion is the morphism from the additive Grassmannians to
the periodic global version of $B O$.
As an application we end the section 
with a global, highly structured version of Bott periodicity:
Theorem \ref{thm:global Bott periodicity}
shows that $\bBUP$ is globally equivalent, as an
ultra-commutative monoid, to $\Omega\bU$.

\section{Global model structure}
\label{sec:global model monoid spaces}

In this section we formally define ultra-commutative monoids
and establish various formal properties of the category $umon$ of
ultra-commutative monoids.
We introduce free ultra-commutative monoids in Example \ref{eg:free umon}.  
The main result is the model structure with global equivalences as the weak equivalences,
see Theorem \ref{thm:global umon}.

\medskip

\begin{defn}
An {\em ultra-commutative monoid} is a commutative orthogonal monoid space.  
We write $umon$ for the category of ultra-commutative monoids.\index{symbol}{  $umon$ - {category of ultra-commutative monoids}}
\end{defn}
\index{subject}{monoid!ultra-commutative|see{ultra-commutative monoid}}\index{subject}{ultra-commutative monoid|(}
\index{subject}{ultra-commutative monoid}

As we explained after Definition \ref{def:orthogonal monoid space},
the data of an ultra-commutative monoid is the same as that
of a lax symmetric monoidal continuous functor 
from the linear isometries category $\bL$ (under orthogonal direct sum)
to the category $\bT$ of spaces (under cartesian product).

\begin{rk}\label{rk:umon is global E_infty}
One can think of an ultra-commutative monoid 
as encoding a collection of $E_\infty$-$G$-spaces,
one for every compact Lie group $G$, compatible under restriction.
If $R$ is a closed orthogonal space and $G$ a compact Lie group,
then the $G$-equivariant homotopy type encoded in $R$
can be accessed as the underlying $G$-space\index{subject}{underlying $G$-space!of an orthogonal space} 
\[ R(\Uc_G) \ = \ \colim_{V\in s(\Uc_G)}\, R(V) \ .\]
The additional structure of an ultra-commutative monoid 
gives rise to an action of a specific $E_\infty$-$G$-operad on this $G$-space,
namely the linear isometries operad of the complete $G$-universe $\Uc_G$.
The $n$-th space of this operad is the space
$\bL(\Uc_G^n,\Uc_G)$\index{subject}{linear isometries operad}
of linear isometric embeddings (not necessarily equivariant)
of $\Uc_G^n$ into $\Uc_G$. 
The group $G$ acts on $\bL(\Uc_G^n,\Uc_G)$ by conjugation,
and the operad structure is by direct sum and composition of linear isometric embeddings.
The symmetric group $\Sigma_n$ permutes the summands in the source.
The space $\bL(\Uc_G^n,\Uc_G)$ 
is $G$-equivariantly contractible by \cite[II Lemma 1.5]{lms},
and the $\Sigma_n$-action is free; in fact,
$\bL(\Uc_G^n,\Uc_G)$ even has the $(G\times \Sigma_n)$-equivariant
homotopy type of a universal space for $(G,\Sigma_n)$-bundles.

By simultaneous passage to colimit over $s(\Uc_G)$ in all $n$ variables,
the iterated multiplication maps 
\[ R(V_1)\times \dots\times R(V_n) \ \to \ R(V_1\oplus\dots\oplus V_n) \]
give rise to a map $\mu_{(n)}:R(\Uc_G)^n\to R(\Uc_G^n)$.
A linear isometric embedding $\psi:\Uc\to\Uc'$ between 
countably infinite dimensional inner product spaces induces a map
$ R(\psi) :  R(\Uc) \to R(\Uc')$; the resulting `action map'
\[ \bL(\Uc,\Uc')\times R(\Uc)\ \to \ R(\Uc') \ , \quad 
(\psi,y)\ \longmapsto \ R(\psi)(y)\]
is continuous.
The operadic action map is then simply the composite
\[   \bL(\Uc_G^n,\Uc_G) \times R(\Uc_G)^n \
\xra{\bL(\Uc_G^n,\Uc_G) \times \mu_{(n)}}\   
  \bL(\Uc_G^n,\Uc_G) \times R(\Uc_G^n) \ \xra{\ \text{act}\  }\   
 R(\Uc_G) \ .\]
\end{rk}

Now we work towards the main result of this section,
the global model structure for ultra-commutative monoids.
Before we start with the homotopical considerations, we get some
of the necessary category theory out of the way.
For a moment we consider more generally 
any symmetric monoidal category $\Cc$ with monoidal
product $\boxtimes$ and unit object $I$. We can then study operads
in $\Cc$ and algebras over a fixed operad.
The following (co-)completeness and preservation results
can be found in \cite[Prop.\,2.3.5]{rezk-thesis} 
or \cite[Prop.\,3.3.1]{fresse-module functors}.

\begin{prop}\label{prop:filtered colim P-alg} 
Let $(\Cc,\boxtimes,I)$ be a complete and cocomplete symmetric monoidal category 
such that the monoidal product preserves colimits in each variable.
Let $\Pc$ be an operad in $\Cc$. 
Then the category of $\Pc$-algebras is complete and cocomplete.
Moreover, the forgetful functor from the category of $\Pc$-algebras
to the underlying category $\Cc$ creates all limits, 
all filtered colimits and those coequalizers that are reflexive in
the underlying category $\Cc$.
\end{prop}

We let $\Com$ denote the incarnation of the commutative operad
internal to the category of orthogonal spaces, under box product.
So for every $n\geq 0$ the orthogonal space $\Com(n)$
of $n$-ary operations is constant with value a one-point space.
Equivalently, $\Com$ is a terminal operad in orthogonal spaces.
Endowing an orthogonal space with an ultra-commutative multiplication
is the same as giving it an algebra structure over the commutative operad $\Com$. 
More formally, the category of ultra-commutative monoids is
isomorphic to the category of $\Com$-algebras.
So Proposition \ref{prop:filtered colim P-alg} has the following special case:

\begin{cor}\label{cor:co-limits in umon}  
The category of ultra-commutative monoids is complete and cocomplete.
The forgetful functor from the category of ultra-commutative mo\-noids
to the category of orthogonal spaces creates all limits, 
all filtered colimits and those coequalizers that are reflexive in
the category of orthogonal spaces.
\end{cor}

\begin{eg}[Free ultra-commutative monoids]\label{eg:free umon} 
We quickly recall that ultra-commutative monoids are monadic over the
category of orthogonal spa\-ces; this is not particular to our context,
and the analogous fact holds for 
commutative monoids in any cocomplete symmetric monoidal category.
For every orthogonal space $Y$ and $m\geq 0$ we denote by
\[ \mP^m(Y)\ = \  Y^{\boxtimes m}/\Sigma_m \]
the $m$-symmetric power, with respect to the box product, of $Y$.
In particular, $\mP^0(Y)$ is the terminal, constant one-point orthogonal space,
and $\mP^1(Y)=Y$.
Then the orthogonal space
\[ \mP(Y)\ = \ {\coprod}_{m\geq 0}\,  \mP^m(Y) \ = \ 
{\coprod}_{m\geq 0}\,  Y^{\boxtimes m}/\Sigma_m \]
is an ultra-commutative monoid under the concatenation product, 
and it is in fact the free ultra-commutative monoid generated by $Y$.\index{symbol}{$\mP(Y)$ - {free ultra-commutative monoid generated by $Y$}}\index{subject}{ultra-commutative monoid!free|see{free ultra-commutative monoid}}\index{subject}{free ultra-commutative monoid}
More precisely, the functor
\[ \mP \ : \ \spc \ \to \ umon \]
becomes a left adjoint to the forgetful functor with respect to
the morphism $\eta_Y:Y =\mP^1(Y) \to \mP Y$, the inclusion of the `linear' summand,
as the adjunction unit. 
In other words, the following composite is bijective
for every ultra-commutative monoid $R$:
\[ umon(\mP Y, R)\ \xra{\text{forget}} \ \spc(\mP Y, R)\ \xra{\ \eta_Y^*\ }\ 
\spc(Y, R)\ .\]
Moreover, this adjunction is monadic, i.e., the category of ultra-commutative monoids
is isomorphic to the category of algebras over the monad $\mP$.
\end{eg}

\begin{construction}\label{con:umon enrich}
We will also exploit that the category of ultra-commu\-tative monoids is
tensored and cotensored over the category $\bT$ of spaces,
so we spend a few words explaining this enrichment.
In fact, the constructions work more generally for algebras over 
continuous monads on any category enriched in spaces,
see for example \cite[Lemma 2.8]{mcclure-schwanzl-vogt}.  
We only spell out the case of ultra-commutative monoids,
which are the algebras over the free ultra-commutative monoid monad $\mP:\spc\to\spc$.

The mapping space of morphisms between two orthogonal spaces $X$ and $Y$ is
defined as follows. Since every inner product space is isometrically isomorphic
to $\mR^n$ for some $n$, the map
\[ \spc(X,Y)\ \to \ {\prod}_{n\geq 0}\, \map(X(\mR^n),Y(\mR^n)) \ , \quad
f\ \longmapsto \ \{ f(\mR^n)\}_{n\geq 0}\]
is injective with closed image.
So we endow $\spc(X,Y)$ with the subspace topology of the product 
(which is taken internal to the category $\bT$ of
compactly generated spaces, i.e., it is the Kelleyfied product topology).

If $R$ and $S$ are ultra-commutative monoids, then the set
$umon(R,S)$ of morphisms of ultra-commutative monoids is a closed subset
of the space $\spc(R,S)$, and we give it the subspace topology.
We omit the verification that composition is continuous in this topology,
so we have indeed defined an enrichment of the category of
ultra-commutative monoids in spaces.

The cotensors of ultra-commutative monoids are defined `pointwise'.
In more detail, we consider an ultra-commutative monoid $R$ and a space $A$.
Then the orthogonal space $\map(A,R)$ inherits an ultra-commutative multiplication
\[ \map(A,R)\boxtimes \map(A,R)\ \to \map(A,R) \]
from the bimorphism with $(V,W)$-component
\[ \map(A,R(V))\times \map(A,R(W))\ \to \ \map(A,R(V\oplus W)) \ , \
(f,g)\ \longmapsto \ \mu_{V,W}\circ(f,g)\ .\]
The multiplicative unit is the constant map $A\to R(0)$
with value the unit of $R$.

The tensor of an ultra-commutative monoid $R$ with a space $A$ is, 
however, {\em not} pointwise. To avoid confusion with
the objectwise product we denote this tensor by $R\tensor A$;
its defining property is that it represents the functor
\[  \map(A,umon(R,-))\ : \ umon \ \to \ \text{(sets)}\ .\]
So $R\tensor A$ comes equipped with a continuous map
$i:A\to umon(R,R\tensor A)$ such that the map
\[ umon(R\tensor A,S)\ \to \ \map(A,umon(R,S)) \ , \quad
f \ \longmapsto \ umon(R,f)\circ i\]
is bijective.
One way to construct a tensor $R\tensor A$ is as a coequalizer,
in the category of ultra-commutative monoids, 
of the two morphisms:
\[\xymatrix@C=18mm{ 
 \mP( (\mP R)\times A)\quad \ar@<.4ex>[r]^-{\mP(\alpha\times A)}
\ar@<-.4ex>[r]_-\nu & \quad \mP(R\times A)  }\]
Here $\alpha:\mP R\to R$ is the structure morphism
(i.e., the counit of the free-forgetful adjunction) and $\nu$
is adjoint to the morphism of orthogonal spaces
\[  (\mP R)\times A \ \to \ \mP(R\times A) \]
that in turn is adjoint to the composite
\[ A \ \xra{a\mapsto (-,a)}\ 
\map(R, R\times A) \ \xra{\ \mP \ } \ \map(\mP R, \mP(R\times A)) \ .\]
The above coequalizer defining $R\tensor A$ is reflexive in the underlying category
of orthogonal spaces, so it can be calculated in the underlying category,
by Proposition \ref{prop:filtered colim P-alg}. 
\end{construction}

In our discussion of global group completions 
in Section \ref{sec:group completion and units}
we will want to realize simplicial ultra-commutative monoids.
We refer to Construction \ref{con:realize simplicial} for generalities
about realization of simplicial objects.\index{subject}{realization!of simplicial ultra-commutative monoids|(}
For a simplicial ultra-commutative monoid $B:\bDelta^{\op}\to umon$, 
the term `geometric realization' actually has two potentially different interpretations,
and we take some time to clarify this issue.
On the one hand we can form the geometric realization $|B|_\text{un}$ 
in the underlying category of orthogonal spaces; this is, by definition, a coend, in the
category of orthogonal spaces, of the functor
\[ \bDelta^{\op}\times \bDelta \ \to \ \spc \ , \quad 
([m],[n])\ \longmapsto \ B_m\times\Delta^n \ .\]
We call this the {\em underlying realization} of $B$.
Coends of orthogonal spaces are calculated objectwise, so $|B|_\text{un}(V)$
is a realization of the simplicial space $[m]\mapsto B_m(V)$.
It is not a priori obvious, however, whether this realization inherits any multiplication.

On the other hand, we explained in Construction \ref{con:umon enrich}
that the category of ultra-commutative monoids 
is tensored over the category $\bT$ of spaces.
We continue to write $R\tensor A$ for the tensor of 
an ultra-commutative monoid $R$ with a space $A$, in order to distinguish it
from the (objectwise) product of the underlying orthogonal space of $R$ with $A$. 
We can also consider the realization $|B|_\text{in}$
{\em internal to ultra-commutative monoids},
i.e., a coend, in the category of ultra-commutative monoids, of the functor
\[ \bDelta^{\op}\times \bDelta \ \to \ umon \ , \quad 
([m],[n])\ \longmapsto \ B_m\tensor \Delta^n \ .\]
We call this the {\em internal realization}.
The internal realization is, by definition, an ultra-commutative monoid, but it is not 
immediately clear how it relates to the underlying realization of $|B|_\text{un}$
as an orthogonal space.
As we shall now show, the forgetful functor from a category of
ultra-commutative monoids to orthogonal spaces commutes with
realization of simplicial objects.
We do not claim any originality here, 
and many results of this kind can be found in the literature, see for example
\cite[Thm.\,12.2]{may-geometry iterated loop},
\cite[VII Prop.\,3.3]{EKMM},
\cite[Prop.\,4.5]{mcclure-schwanzl-vogt},
\cite[Prop.\,12.4]{mandell-algebraization}
or \cite[Thm.\,2.2]{Goerss-Hopkins-moduli}.

\begin{prop}\label{prop:U preserves realization} 
Let $B$ be a simplicial object in the category of ultra-commu\-ta\-tive monoids. 
Then the canonical morphism $|B|_\text{\em un} \to |B|_\text{\em in}$
from the underlying realization to the internal realization 
is an isomorphism of orthogonal spaces. 
\end{prop}
\begin{proof}
We adapt an argument given by Mandell
in an unpublished preprint \cite[Prop.\,12.4]{mandell-algebraization}.
We start by considering two simplicial orthogonal spaces
$X, Y:\bDelta^{\op}\to \spc$.
We denote by $X\boxtimes Y$ the diagonal of the external box product, 
i.e., the composite simplicial orthogonal space 
\[ \bDelta^{\op}\ \xra{\text{diagonal}}\ \bDelta^{\op}\times\bDelta^{\op}\ 
\xra{X\times Y}\ \spc\times\spc\ \xra{\ \boxtimes\ }\ \spc\ .\]
For every $n\geq 0$ we consider the composite
\begin{align*}
  (X_n\boxtimes Y_n)\times\Delta^n \ \xra{\Id\times\text{diag}}\ 
&(X_n\boxtimes Y_n)\times (\Delta^n\times\Delta^n)\\ 
\xra{\ \text{shuffle}\ }\ &(X_n\times \Delta^n)\boxtimes(Y_n\times \Delta^n)
\ \to\ |X|\boxtimes|Y| \ , 
\end{align*}
where the last morphism is the box product of the two canonical morphisms
$X_n\times\Delta^n\to |X|$ and $Y_n\times\Delta^n\to|Y|$.
These composites are compatible with the coend relations,
so they assemble into a morphism of orthogonal spaces
\[   | X\boxtimes Y| \ \to \  |X|\boxtimes |Y| \ .\]
We claim that this morphism is an isomorphism. 
Indeed, since $\boxtimes$ preserves colimits
in each variable, the right hand side is a coend of the functor
\[ (\bDelta^2)^{\op}\times \bDelta^2 \ \to \ \spc \ , \quad 
([k],[l],[m],[n])\ \longmapsto \ (X_k\boxtimes Y_l) \times \Delta^m\times\Delta^n\ .\]
Coends of orthogonal spaces are calculated objectwise, and
for bisimplicial spaces the bi-realization is homeomorphic to the
realization of the diagonal 
(see \cite[p.\,94, Lemma]{quillen-higher algebraic}
or Proposition \ref{prop:iterated geometric realization}~(iii)).

By iterating, we obtain a $\Sigma_m$-equivariant isomorphism of orthogonal spaces
\[ |X^{\boxtimes m}| \ \iso \  |X|^{\boxtimes m} \]
for every $m\geq 0$. Since coends commute with colimits, we can pass to
$\Sigma_m$-orbits and take coproduct over $m\geq 0$, resulting in an isomorphism 
\[ 
   | \mP (X) |_\text{un} \ = \  | \amalg_{m\geq 0} (X^{\boxtimes m})/\Sigma_m | \ \iso \ 
 \amalg_{m\geq 0}\, |X|^{\boxtimes m}/\Sigma_m  \ = \  \mP |X|\ .  
 \]
On the other hand, the ultra-commutative monoid 
$\mP |X|$ has the universal property of the internal realization
of the simplicial ultra-commutative monoid $\mP\circ X$. 
This shows the claim in the special case where $B$ is
freely generated by a simplicial orthogonal space.

Now we treat the general case. The diagram
\[
 \xymatrix{
 \mP(\mP B) \ar@<.4ex>[r]^-{\mP\alpha}\ar@<-.4ex>[r]_-{\mu} & 
\mP B\ar[r]^-\alpha & B } 
 \]
is a coequalizer diagram of simplicial ultra-commutative monoids.
Here $\alpha:\mP R\to R$ is the adjunction counit
and $\mu$ is the monad structure of the free functor.
Moreover, the coequalizer is reflexive in the underlying category of orthogonal spaces,
by the morphisms
\[\xymatrix{
 \mP(\mP B) & \mP B\ar[l]_-{\eta_{\mP B}} & B \ar[l]_-{\eta_B} } \]
where $\eta:R\to \mP R$ is the unit of the free-forget adjunction,
i.e., the inclusion as the `linear' summand $\mP^1(R)$.

Applying the two functors under consideration
gives a commutative diagram of orthogonal spaces
\begin{equation}\begin{aligned}\label{eq:two split coeq} 
\xymatrix{ |\mP(\mP B)|_\text{un} \ar@<.4ex>[r]\ar@<-.4ex>[r] \ar[d] 
& |\mP B|_\text{un}\ar[r] \ar[d]&  |B|_\text{un} \ar[d] \\
|\mP(\mP B)|_\text{in} \ar@<.4ex>[r]\ar@<-.4ex>[r] & 
|\mP B|_\text{in}\ar[r] & |B|_\text{in}  }
\end{aligned}\end{equation}
We claim that both rows are coequalizer diagrams of orthogonal spaces.
For the upper row we argue as follows. For every $n\geq 0$ the diagram
\[ \xymatrix{ \mP(\mP B_n) \ar@<.4ex>[r]\ar@<-.4ex>[r] & 
 \mP B_n\ar[r] & B_n } \] 
is a coequalizer in the category of ultra-commutative monoids.
Since the diagram splits in the underlying category
of orthogonal spaces, it is also a coequalizers diagram there,
by Proposition \ref{prop:filtered colim P-alg}  or \cite[IV.6, Lemma]{maclane-working}.
So the diagram
\[ \xymatrix{  \mP(\mP B) \ar@<.4ex>[r]\ar@<-.4ex>[r] & \mP B \ar[r] & B } \] 
is also a coequalizer diagram of simplicial orthogonal spaces.
Since product with $\Delta^n$ and coends commute with colimits
in the category of orthogonal spaces, the diagram stays
a coequalizer after (underlying) geometric realization.
Since coends commute with all colimits,
the bottom row of \eqref{eq:two split coeq}
is a coequalizer diagram of ultra-commutative monoids.
Again the diagram splits in the underlying category of orthogonal spaces,
so the lower diagram is also a coequalizer diagram of orthogonal spaces.
Since the two left vertical morphisms in \eqref{eq:two split coeq}
are isomorphisms of orthogonal spaces by the special case above, 
this proves that the morphism 
$|B|_\text{un}\to |B|_\text{in}$ is an isomorphism as well.
\end{proof}
\index{subject}{realization!of simplicial ultra-commutative monoids|)}

Ultra-commutative monoids form a pointed category: 
the constant one-point orthogonal monoid space is a zero object.
The enrichment, tensors and cotensors over spaces extend to
enrichment, tensors and cotensors over the category of {\em based} topological spaces. 
We shall write $R\rhd A$ for the tensor of an ultra-commutative monoid $R$
with a based space $(A,a_0)$, in order to distinguish it
from the (objectwise) smash product of 
the underlying based orthogonal space of $R$ with $A$.
Thus $R\rhd A$ is a pushout, in the category of ultra-commutative monoids,
of the diagram
\begin{equation}\label{eq:R rhd A}
 \ast \ \xla{\qquad}\  R\tensor \{a_0\}\ \xra{R\tensor\text{incl}} \ R\tensor A\ .  
\end{equation}
As may be familiar from similar contexts,
the bar construction $B(R)$ of an ultra-commutative monoid $R$
can be interpreted as $R\rhd S^1$, the based tensor of $R$ with the based space $S^1$,
see \eqref{eq:umon_bar_and_suspension} below.
Another way to say this is that 
the bar construction is the internal suspension, 
in the category of an ultra-commutative monoids.
We show a more general statement
and consider a based simplicial set $A$.
We define a simplicial object of ultra-commutative monoids by
\[ B_m(R,A)\ = \ R\rhd A_m   \ ,  \]
with simplicial structure induced by that of $A$. Since $A_m$ is a based set,
$R\rhd A_m$ is in fact a categorical coproduct (i.e., box product) 
of copies of $R$, indexed by the non-basepoint elements of $A_m$.

The next proposition constructs an isomorphism 
of ultra-commutative mo\-noids between $R\rhd |A|$
and $|B_\bullet(R,A) |_\text{in}$, the internal geometric realization.
By Proposition \ref{prop:U preserves realization},
we can (and will) confuse the internal realization with 
the underlying realization of $B_\bullet(R,A)$ 
in the category of orthogonal spaces.
Variations of the following proposition appear in various places in the literature,
and they go back, at least, to the interpretation,
by McClure, Schw{\"a}nzl and~Vogt \cite{mcclure-schwanzl-vogt},
of topological Hochschild homology of a commutative ring spectrum 
as the tensor with $S^1$.

\begin{prop}\label{prop:B bullet and A lhd} 
Let $R$ be an ultra-commutative monoid and $A$ a based simplicial set.
Then $R\rhd |A|$ is an internal realization of the simplicial ultra-commutative monoid
$B_\bullet(R,A)$.
\end{prop}
\begin{proof}
The geometric realization $|A|$ is a coend of the functor
\[ \bDelta^{\op}\times \bDelta \ \to \ \bT_* \ , \quad 
([m],[n])\ \longmapsto \ A_m\sm \Delta^n_+  \ . \]
Since the functor $R\rhd -$ preserves colimits, $R\rhd |A|$ is a coend, 
in the category of ultra-commutative monoids, of the functor
\[ 
 \bDelta^{\op}\times \bDelta \ \to \ umon \ ,\quad
([m],[n])\ \longmapsto \ R\rhd ( A_m\sm \Delta^n_+) \ .\]
The isomorphisms
\begin{align*}
R\rhd ( A_m\sm \Delta^n_+)\ \iso \  (R\rhd A_m)\rhd \Delta^n_+\ 
\iso \ (R\rhd A_m)\tensor \Delta^n\ = \ B_m(A,R)\tensor \Delta^n
\end{align*}
are natural in $([m],[n])\in \bDelta^{\op}\times\bDelta$,
and they show that $R\rhd|A|$ is an internal realization
of the simplicial ultra-commutative monoid $B_\bullet(R,A)$.
\end{proof}

Now we approach the global model structure on the category 
of ultra-commu\-ta\-tive monoids.
We will establish this model structure as a special case of a lifting theorem
for model structures to categories of commutative monoids that 
was formulated by White \cite[Thm.\,3.2]{white-commutative monoids}.
Like its predecessor for 
associative monoids \cite[Thm.\,4.1 (3)]{schwede-shipley-monoidal},
the input is a cofibrantly generated symmetric monoidal model category
that satisfies the monoid axiom.
However, lifting a model structure to {\em commutative} monoids is more subtle 
and needs extra hypotheses; the essence of the additional condition is
that, loosely speaking, symmetric powers must be `sufficiently homotopy invariant'.
Before White, Gorchinskiy and Guletskii \cite{gor-gul-symmetric}
also studied symmetric power constructions in a symmetric monoidal model category,
and there is a substantial overlap in the arguments of \cite{gor-gul-symmetric}
and \cite{white-commutative monoids}. 

\medskip

We let $\Cc$ be a symmetric monoidal category with monoidal product $\boxtimes$.
To simplify the exposition we follow the common abuse
to suppress the associativity and unit isomorphisms from the notation,
i.e., we pretend that the underlying monoidal structure is strict
(i.e., a {\em permutative} structure).
We let $i:A\to B$ be a $\Cc$-morphism and arrange the $n$-fold $\boxtimes$-power
of $i$ into an $n$-dimensional cube $K^n(i)$ in $\Cc$, i.e., a functor
\[ K^n(i)\ :\  \Pc(n) \ \to \ \Cc \]
from the poset category of subsets of $\{1,2,\dots,n\}$ and
inclusions to $\Cc$. 
More explicitly, if $S\subseteq \{1,2,\dots,n\}$ is a subset, 
then the vertex of the cube at $S$ is 
\[ K^n(i)(S) \ = C_1 \boxtimes C_2 \boxtimes \dots \boxtimes C_n \text{\qquad with\qquad}
 C_j \ = \left\{ \begin{array}{l@{\quad}l} A & \mbox{if } j \not\in  S \\
B & \mbox{if } j \in S. \end{array} \right. \]
All morphisms in the cube $K^n(i)$ are $\boxtimes$-products of identities and
copies of the morphism $i:A\to B$.
The initial vertex of the cube is $K^n(i)(\emptyset)=A^{\boxtimes n}$
and the terminal vertex is $K^n(i)(\{1,\dots,n\})=B^{\boxtimes n}$.

We denote by $Q^n(i)$ the colimit of the punctured cube,
i.e., the cube $K^n(i)$ with the terminal vertex removed, 
and by $i^{\Box n}:Q^n(i)\to K^n(i)(\{1,\dots,n\})=B^{\boxtimes n}$ 
the canonical morphism, an iterated  pushout product morphism.
Indeed, for $n=2$ the cube $K^2(i)$ is a square and looks like
\[ \xymatrix{
A\boxtimes A\ar[r]^-{A\boxtimes i}\ar[d]_{i\boxtimes A} & A\boxtimes B\ar[d]^{i\boxtimes B} \\ 
B\boxtimes A\ar[r]_-{B\boxtimes i}  & B\boxtimes B} \]
Hence 
\[ i^{\Box 2}\ =\ i\Box i \ = \ (B\boxtimes i)\cup(i\boxtimes B)\ : \ 
B\boxtimes A\cup_{A\boxtimes A} A\boxtimes B \ \to \ B\boxtimes B\ .\]
Similarly, $i^{\Box 3}$ is the morphism from the colimit of the punctured
cube to the terminal vertex of the following cube:
\[ 
\xymatrix@C=8mm@R=8mm{ A\boxtimes A\boxtimes A \ar[rr]^{A \boxtimes A \boxtimes i} 
\ar[dr]_{A \boxtimes i\boxtimes A} \ar[dd]_{i\boxtimes A \boxtimes A } && 
A\boxtimes A\boxtimes B \ar[dr]^{A \boxtimes i\boxtimes B} \ar'[d][dd]^{i\boxtimes A\boxtimes B}\\
& A\boxtimes B\boxtimes A \ar[rr]_(.35){A \boxtimes B\boxtimes i} \ar[dd]^(.7){i\boxtimes B\boxtimes A } 
&& A\boxtimes B\boxtimes B \ar[dd]^{i\boxtimes B\boxtimes B}\\
B\boxtimes A\boxtimes A \ar[dr]_{B \boxtimes i\boxtimes A} 
\ar'[r]^(.6){B \boxtimes A\boxtimes i}[rr] && B\boxtimes A\boxtimes B \ar[dr]^{B \boxtimes i\boxtimes B}\\
& B\boxtimes B\boxtimes A \ar[rr]_-{B\boxtimes B\boxtimes i}  && B\boxtimes B\boxtimes B } \]

We observe that the symmetric group $\Sigma_n$ acts on 
$Q^n(i)$ and $B^{\boxtimes n}$ by permuting the factors,
and the iterated pushout product morphism
$i^{\Box n}:Q^n(i)\to B^{\boxtimes n}$ is $\Sigma_n$-equivariant.
We recall from \cite{gor-gul-symmetric} 
the notions of {\em symmetrizable cofibration}
and {\em symmetrizable acyclic cofibration}.

\begin{defn}\cite[Def.\,3]{gor-gul-symmetric} \label{def:symmetrizable} Let $\Cc$ be a symmetric monoidal model category. 
A morphism $i:A\to B$ is a  {\em symmetrizable cofibration}\index{subject}{symmetrizable cofibration}\index{subject}{symmetrizable acyclic cofibration}
(respectively a {\em symmetrizable acyclic cofibration}) if the morphism
\[ i^{\Box n}/\Sigma_n \ :\ Q^n(i)/\Sigma_n \ \to \ B^{\boxtimes n}/\Sigma_n = \mP^n(B) \]
is a cofibration (respectively an acyclic cofibration) for every $n\geq 1$.  
\end{defn}

Since the morphism $i^{\Box 1}/\Sigma_1$ is the original morphism $i$, 
every symmetrizable cofibration
is in particular a cofibration and every symmetrizable acyclic cofibration
is in particular an acyclic cofibration.
We will now proceed to prove that in the category of orthogonal spaces,
all cofibrations and acyclic cofibrations in the positive global model structure
are symmetrizable with respect to the box product. 
The next proposition will be used to verify this 
for the generating acyclic cofibrations.
We recall from Construction \ref{con:define Z(j)}
that given a morphism $j:A\to B$, the set $\Zc(j)$ consists of all pushout product maps
\[  c(j)\Box i_k \ : \  A\times D^k \cup_{A\times \partial D^k} Z(j)\times \partial D^k\ \to \
Z(j)\times D^k \]
of the mapping cylinder inclusion $c(j):A\to Z(j)$
with the sphere inclusions for $k\geq 0$.

\begin{prop}\label{prop: Zc is symmetrizable}
Let $\Cc$ be a symmetric monoidal topological model category. 
  \begin{enumerate}[\em (i)]
  \item For every $n\geq 1$ the functor $\mP^n$ preserves the homotopy
    relation on morphisms and it preserves homotopy equivalences.
  \item Let $j:A\to B$ be a symmetrizable acyclic cofibration between cofibrant objects.
    Then for every $k\geq 0$, the pushout product map
    \[ j\Box i_k \ : \ A\times D^k\cup_{A\times \partial D^k} B\times \partial D^k\ \to \ 
    B\times D^k \]
    is a symmetrizable acyclic cofibration.
  \item 
    Let $j:A\to B$ be a morphism between cofibrant objects
    such that the morphism $\mP^n(j):\mP^n(A)\to\mP^n(B)$ is a weak equivalence
    for every $n\geq 1$. Then every morphism in the set $\Zc(j)$
    is a symmetrizable acyclic cofibration.
  \end{enumerate}
\end{prop}
\begin{proof}
(i) This is the topological version of \cite[Lemma~1]{gor-gul-symmetric}.
For every object $A$ of $\Cc$ and  every space $K$ the morphism
\[ A^{\boxtimes n}\times K \ \xra{A^{\boxtimes n}\times\Delta} \
A^{\boxtimes n}\times K^n \ \iso \ (A\times K)^{\boxtimes n} \]
is $\Sigma_n$-equivariant (with respect to the trivial $\Sigma_n$-action 
on $K$ in the source) and factors over a natural morphism
\[ \tilde\Delta \ : \ \mP^n(A)\times K \ = \ ( A^{\boxtimes n}\times K )/\Sigma_n \ \to \
(A\times K)^{\boxtimes n} /\Sigma_n \ = \ \mP^n(A\times K) \ . \]
If $H:A\times [0,1]\to B$ is a homotopy from a morphism $f=H(-,0)$
to another morphism $g=H(-,1)$, then the composite
\[ \mP^n(A)\times [0,1] \ \xra{\tilde\Delta} \ \mP^n(A\times [0,1]) \ \xra{\mP^n(H)} \
\mP^n(B) \]
is a homotopy from the morphism $\mP^n(f)$ to $\mP^n(g)$.
So $\mP^n$ preserves the homotopy relation, and hence also homotopy equivalences.

(ii) We argue by induction on $k$. 
For $k=0$ the pushout product map $j\Box i_0$
is isomorphic to $j$, hence a symmetrizable acyclic cofibration by hypothesis.
Now we assume the claim for some $k$, and deduce it for $k+1$.
Since $j$ is a symmetrizable acyclic cofibration between cofibrant objects,
the morphism $\mP^n(j)$ is a weak equivalence for every $n\geq 1$ 
by \cite[Cor.\,23]{gor-gul-symmetric}.  
Since the functors $\mP^n$ preserve the homotopy relation and 
the projections $A\times D^k\to A$ and
$B\times D^k\to B$ are homotopy equivalences,
the morphism $\mP^n(j\times D^k)$ is a weak equivalence for every $n\geq 1$.
So $j\times D^k:A\times D^k\to B\times D^k$ is a symmetrizable acyclic cofibration,
again by \cite[Cor.\,23]{gor-gul-symmetric}.  
We write $\partial D^{k+1}=D^k_+\cup_{\partial D^k}D^k_-$ as the union of the upper and lower
hemisphere along the equator. The upper morphism in the pushout square
\[ \xymatrix{D^k_+\times A\ar[r]^-{D^k_+\times j}\ar[d] & D^k_+\times B \ar[d]\\
\partial D^{k+1}\times A \ar[r] & \partial D^{k+1}\times A\cup_{D^k_+\times A} D^k_+\times B } \]
is a symmetrizable acyclic cofibration by the previous paragraph.
The class of symmetrizable acyclic cofibrations is closed under cobase
change by \cite[Thm.\,7 (A)]{gor-gul-symmetric};
the lower morphism is thus a symmetrizable acyclic cofibration.

The square 
\[ \xymatrix@C=20mm{ 
A\times D^k_- \cup_{A\times \partial D^k} B\times \partial D^k \ar[d]\ar[r]^-{j\Box i_k} & 
B\times D^k_- \ar[d] \\
A\times \partial D^{k+1} \cup_{A\times D^k_+} B\times D^k_+ \ar[r]_-{j\Box\text{incl}} & 
B\times \partial D^{k+1} } \]
is a pushout. The upper morphism is a symmetrizable acyclic cofibration 
by the inductive hypothesis, 
hence so is the lower morphism,
again by stability under cobase change.
The morphism $j\times \partial D^{k+1}:A\times\partial D^{k+1}\to B\times \partial D^{k+1}$ 
is thus the composite of two symmetrizable acyclic cofibrations,
hence a symmetrizable acyclic cofibration itself, 
by \cite[Thm.\,7 (C)]{gor-gul-symmetric}.
As a cobase change, the morphism
\[ A\times D^{k+1}\ \to \ A\times D^{k+1}
\cup_{A\times \partial D^{k+1}} B\times \partial D^{k+1} \]
is then a symmetrizable acyclic cofibration
by \cite[Thm.\,7 (A)]{gor-gul-symmetric}.
The induced morphism
\[ \mP^n(A\times  D^{k+1})\ \to \ 
\mP^n( A\times D^{k+1}\cup_{A\times \partial D^{k+1}} B\times \partial D^{k+1}) \] 
is then a weak equivalence by \cite[Cor.\,23]{gor-gul-symmetric}. 
Since $\mP^n(j\times D^{k+1}):\mP^n( A\times D^{k+1})\to \mP^n(B\times D^{k+1})$ 
is a weak equivalence, so is the morphism
\[ \mP^n(j\Box i_{k+1})\ : \ 
\mP^n( A\times D^{k+1}\cup_{A\times \partial D^{k+1}} B\times \partial D^{k+1}) \ \to \
\mP^n( B\times D^{k+1})\ .\] 
One more time by \cite[Cor.\,23]{gor-gul-symmetric}, this shows that $j\Box i_{k+1}$ is
a symmetrizable acyclic cofibration. This completes the induction step.

(iii) Since $A$ and $B$ are cofibrant, the mapping cylinder inclusion  
\[ c(j) \ : \ A \ \to \ ( A\times [0,1] )\cup_j B = Z(j)\]
is a cofibration. Moreover, the projection $Z(j)\to B$ is a homotopy
equivalence, hence so is $\mP^n(Z(j))\to \mP^n(B)$ for every $n\geq 1$. 
Since $\mP^n(j)$ is a weak equivalence by hypothesis, the morphism 
$\mP^n(c(j)):\mP^n(A)\to \mP^n(Z(j))$ is a weak equivalence for every $n\geq 1$. 
So $c(j)$ is a symmetrizable acyclic cofibration by \cite[Cor.\,23]{gor-gul-symmetric}. 
Applying (ii) to the morphism $c(j)$ yields the claim.
\end{proof}

Now we can verify the symmetrizability of cofibrations and acyclic cofibrations
for the positive global model structure of orthogonal spaces.
The cofibration part~(i) is in fact slightly 
stronger in that it does not need any positivity hypothesis.

\begin{theorem}\label{thm:symmetrizable in orthogonal spaces}
  \begin{enumerate}[\em (i)]
  \item 
    Let $i:A\to B$ be a flat cofibration  of orthogonal spaces.
    Then for every $n\geq 1$ the morphism 
    \[ i^{\Box n}/\Sigma_n \ : \ Q^n(i)/\Sigma_n \ \to \ B^{\boxtimes n}/\Sigma_n \]
    is a flat cofibration. 
    In other words, all cofibrations in the global model structure
    of orthogonal spaces are symmetrizable.
  \item
    Let $i:A\to B$ be a positive flat cofibration of orthogonal spaces 
    that is also a global equivalence.
    Then for every $n\geq 1$ the morphism 
    \[ i^{\Box n}/\Sigma_n \ : \ Q^n(i)/\Sigma_n \ \to \ B^{\boxtimes n}/\Sigma_n \]
    is a global equivalence.
    In other words, all acyclic cofibrations in the positive global model structure
    of orthogonal spaces are symmetrizable.
  \end{enumerate}
\end{theorem}
\begin{proof}
(i) We recall from the proof of Proposition \ref{prop:strong level spaces}
the set 
\[ I^{\str} \ = \ \{\ G_m(  O(m)/H\times i_k )\ | \ m,k \geq 0, H\leq O(m)\} \]
of generating flat cofibrations of orthogonal spaces, where $i_k:\partial D^k\to D^k$
is the inclusion.
The set $I^{\str}$ detects the acyclic fibrations 
in the strong level model structure of orthogonal spaces.
In particular, every flat cofibration is a retract of an $I^{\str}$-cell complex.
By \cite[Cor.\,9]{gor-gul-symmetric} 
or \cite[Lemma A.1]{white-commutative monoids}, 
it suffices to show that the generating flat cofibrations in $I^{\str}$ are symmetrizable.

The orthogonal space $G_m( O(m)/H \times K)$ is isomorphic to $\bL_{H,\mR^m}\times K$,
so we show more generally that every morphism of the form
\[  \bL_{G,V}\times i_k  \ : \  \bL_{G,V}\times \partial D^k \ \to \ \bL_{G,V}\times D^k \]
is a symmetrizable cofibration, 
where $V$ is any representation of a compact Lie group $G$.
The symmetrized iterated pushout product
\begin{equation}  \label{eq:sym_iterated_ppp}
 (  \bL_{G,V}\times i_k )^{\Box n} / \Sigma_n \ : \ 
Q^n(\bL_{G,V}\times i_k ) /\Sigma_n\ \to \  (\bL_{G,V}\times i_k )^{\boxtimes n}/\Sigma_n 
\end{equation}
is isomorphic to 
\[ \bL_{\Sigma_n\wr G,V^n}( i_k^{\Box n}) \ : \ 
 \bL_{\Sigma_n\wr G,V^n} ( Q^n ( i_k) ) \ \to \  \bL_{\Sigma_n\wr G,V^n}( (D^k)^n) \ ,\]
where
\[ i_k^{\Box n}\ : \ Q^n ( i_k ) \ \to \ (D^k)^n \]
is the $n$-fold pushout product of the inclusion $i_k:\partial D^k\to D^k$, 
with respect to the cartesian product of spaces. 
Here the wreath product $\Sigma_n\wr G$ acts on $V^n$ by
\[ (\sigma;\, g_1,\dots,g_n)\cdot (v_1,\dots, v_n) \ = \ 
(g_{\sigma^{-1}(1)}v_{\sigma^{-1}(1)},\dots, g_{\sigma^{-1}(n)}v_{\sigma^{-1}(n)}) \ .\]
The map $i_k^{\Box n}$ is $\Sigma_n$-equivariant, 
and we claim that $i_k^{\Box n}$ is a cofibration of $\Sigma_n$-spaces.
One way to see this is to exploit that $i_k$ is homeomorphic to
the geometric realization of the inclusion $\iota_k:\partial \Delta[k]\to\Delta[k]$
of the boundary of the simplicial $k$-simplex.
So $i_k^{\Box n}$ is $\Sigma_n$-homeomorphic to
the geometric realization of the 
inclusion $\iota_k^{\Box n}:Q^n(\iota_k)\to \Delta[k]^n$
of $\Sigma_n$-simplicial sets. The geometric realization of an equivariant
embedding of simplicial sets is always an equivariant cofibration of spaces,
so altogether this shows that $i_k^{\Box n}$ is a cofibration of $\Sigma_n$-spaces.
Proposition \ref{prop:cofibrancy preservers}~(i)
then shows that $i_k^{\Box n}$ is also a cofibration of $(\Sigma_n\wr G)$-spaces
via restriction along the projection $\Sigma_n\wr G\to\Sigma_n$. 
So the morphism \eqref{eq:sym_iterated_ppp} is a flat cofibration.

(ii) Theorem \ref{thm:All global spaces} describes a set $J^{\str}\cup K$ 
of generating acyclic cofibrations for the global model structure
on the category of orthogonal spaces.
From this we obtain a set $J^+\cup K^+$ 
of generating acyclic cofibration for the {\em positive} 
global model structure of Proposition \ref{prop:positive global spaces} 
by restricting to those morphisms in  $J^{\str}\cup K$ 
that are positive cofibrations, i.e., homeomorphisms in level~0. 
So explicitly, we set
\[ J^+ \ = \ \{\ G_m(  O(m)/H \times j_k)\ | \ m\geq 1, k \geq 0, H\leq O(m)\} \ ,\]
where $j_k:D^k\times\{0\}\to D^k\times [0,1]$ is the inclusion, and
\[ K^+ \ = \ \bigcup_{G,V,W\ : \  V\ne 0} \Zc(\rho_{G,V,W}) \ ,\]
the set of all pushout products of sphere inclusions $i_k$
with the mapping cylinder inclusions of the global equivalences
$\rho_{G,V,W}:\bL_{G,V\oplus W}\to \bL_{G,V}$.
Here $(G,V,W)$ runs through a set of representatives
of the isomorphism classes of triples consisting of a compact Lie group $G$,
a non-zero faithful $G$-representation $V$ and an arbitrary $G$-representation $W$.
By \cite[Cor.\,9]{gor-gul-symmetric} or \cite[Lemma A.1]{white-commutative monoids} 
it suffices to show that all morphisms in $J^+\cup K^+$
are symmetrizable acyclic cofibrations. 

We start with a morphism $G_m( j_k\times O(m)/H )$ in $J^+$.
For every $n\geq 1$, the morphism
\[ (G_m(  O(m)/H\times j_k ))^{\Box n}/\Sigma_n \]
is a flat cofibration by part~(i), and a homeomorphism in level~0
because $m\geq 1$. Moreover, the map $j_k$ is a homotopy
equivalence of spaces, so  $G_m(  O(m)/H \times j_k)$ is a homotopy equivalence
of orthogonal spaces; the morphism $\mP^n(G_m(  O(m)/H\times j_k ))$ 
is then again a homotopy equivalence for every $n\geq 1$,
by Proposition \ref{prop: Zc is symmetrizable}~(i).
Then \cite[Cor.\,23]{gor-gul-symmetric} shows that $G_m( O(m)/H \times j_k)$
is a symmetrizable acyclic cofibration. 
This takes care of the set $J^+$.

Now we consider the morphisms in the set $K^+$.
Since $G$ acts faithfully on the non-zero inner product space $V$, the action
of the wreath product $\Sigma_n\wr G$ on $V^n$ is again faithful.
So the morphism
\[ \rho_{\Sigma_n\wr G,V^n, W^n}\ : \ \bL_{\Sigma_n\wr G, V^n\oplus W^n} \ \to \
\bL_{\Sigma_n\wr G, V^n} \]
is a global equivalence by Proposition \ref{prop:free_orthogonal_space}~(ii).
By the isomorphism
\[ \mP^n(\bL_{G,V}) \ = \ \bL_{G,V}^{\boxtimes n}/\Sigma_n \ \iso \ \bL_{\Sigma_n\wr G,V^n} \ ,\]
the morphism $\rho_{\Sigma_n\wr G,V^n, W^n}$ is isomorphic to
 $\mP^n(\rho_{G,V,W}) :\mP^n(\bL_{G,V\oplus W}) \to\mP^n(\bL_{G,V})$,
which is thus a global equivalence. 
Proposition \ref{prop: Zc is symmetrizable}~(iii)
then shows that all morphisms in $\Zc(\rho_{G,V,W})$ 
are symmetrizable acyclic cofibrations. 
\end{proof}

\Danger The hypothesis in Theorem \ref{thm:symmetrizable in orthogonal spaces}~(ii)
that $i$ is a {\em positive} flat cofibration is really necessary. 
Indeed, the unique morphism $\rho:\bL_{\mR}\to\ast$ to the terminal orthogonal space
is a global equivalence, and source and target of $\rho$ are flat, 
but only the source is positively flat.
Then the mapping cylinder inclusion $c(\rho):\bL_{\mR}\to C(\bL_{\mR})$
is a global equivalence between flat orthogonal spaces, but it is {\em not}
a homeomorphism at~0. And indeed, for no $n\geq 2$ is the morphism
$\mP^n(\bL_{\mR})\to\mP^n(C(\bL_{\mR}))$ a global equivalence,
because the source is isomorphic to $\bL_{\Sigma_n,\mR^n}=B_{\gl}\Sigma_n$,
whereas the target is homotopy equivalent to the terminal orthogonal space. 

\medskip

Now we put all the pieces together and prove the 
global model structure for ultra-commutative monoids.
We call a morphism of ultra-commutative monoids a
{\em global equivalence} (respectively {\em positive global fibration})
if the underlying morphism of orthogonal spaces is a
 global equivalence (respectively fibration in the positive global model structure).

\begin{theorem}[Global model structure for ultra-commutative monoids]\label{thm:global umon}\quad  
  \begin{enumerate}[\em (i)]
  \item 
    The global equivalences and positive global fibrations are part of a 
   cofibrantly generated, proper, topological model structure
    on the category of ultra-commutative monoids,
    the {\em global model structure}.\index{subject}{global model structure!for ultra-commutative monoids}
  \item Let $j:R\to S$ be a cofibration in the global model structure 
  of ultra-commutative monoids.
  \begin{enumerate}[\em (a)]
  \item The morphism of $R$-modules underlying $j$ is a cofibration
    in the global model structure of $R$-modules 
    of Corollary {\em \ref{cor-lift to modules spaces}~(i)}. 
  \item The morphism of orthogonal spaces underlying $j$ is an h-cofibration,
    and hence a closed embedding.
  \item If the underlying orthogonal space of $R$ is flat, then 
    $j$ is a flat cofibration of orthogonal spaces.
\end{enumerate} 
  \end{enumerate}
\end{theorem}
\begin{proof}
(i) The positive global model structure of orthogonal spaces
established in Proposition \ref{prop:positive global spaces} is
 cofibrantly generated and monoidal 
(by Proposition \ref{prop:ExF ppp spaces} (iv)).
The `unit axiom' also holds: we let $f:I\to\ast$ be any positive flat replacement
of the monoidal unit, the constant one-point orthogonal space.
Then for every orthogonal space $Y$ the induced morphism
$f\boxtimes X:I\boxtimes Y\to \ast\boxtimes Y$ is a global equivalence by
Theorem \ref{thm:box to times}~(ii).
The monoid axiom holds by Proposition \ref{prop:monoid orthogonal spaces}.
Cofibrations and acyclic cofibrations are symmetrizable by
Theorem \ref{thm:symmetrizable in orthogonal spaces},
so the model structure satisfies the `commutative monoid axiom'
in the sense of \cite[Def.\,3.1]{white-commutative monoids}.
The symmetric algebra functor $\mP$ commutes with filtered colimits
by Corollary \ref{cor:co-limits in umon}.  
Theorem~3.2 of \cite{white-commutative monoids} thus shows that the positive global
model structure of orthogonal spaces lifts to the category of 
ultra-commutative monoids.

The global model structure is topological 
by Proposition \ref{prop:topological criterion},
where we take $\Gc$ as the set of free ultra-commutative monoids $\mP(L_{H,\mR^m})$
for all $m\geq 1$ and all closed subgroups $H$ of $O(m)$,
and we take $\Zc$ as the set of acyclic cofibrations
$\mP(c(\rho_{G.V.W}))$ for the mapping cone inclusions $c(\rho_{G,V,W})$ of
the global equivalences $\rho_{G,V,W}:\bL_{G,V\oplus W}\to\bL_{G,V}$,
indexed by representatives as in the definition of the set $K^+$.
Since weak equivalences and fibrations of ultra-commutative monoids are
defined on underlying orthogonal spaces, and since pullbacks of
ultra-commutative monoids are created on underlying orthogonal spaces,
right properness is inherited from the positive global model structure
of orthogonal spaces (Proposition \ref{prop:positive global spaces}).
We defer the proof of left properness until after the proof of part~(ii).

  (ii) For (a) we recall that the global model structure on the category of $R$-modules
  is lifted, via the free and forgetful adjoint functor pair,
  from the absolute global model structure of Theorem \ref{thm:All global spaces}. 
  By Corollary \ref{cor-lift to modules spaces}~(i) and~(ii)
  this model structure of $R$-modules is a cofibrantly generated
  monoidal model category that satisfies the monoid axiom.
  Moreover, the unit object $R$ is cofibrant; for this it is relevant that
  we have lifted the {\em absolute} model structure (as opposed to the
  {\em positive model structure}).
  We claim that all cofibrations in this model structure
  are symmetrizable with respect to the box product of $R$-modules.
  By \cite[Cor.\,9]{gor-gul-symmetric}  or \cite[Lemma A.1]{white-commutative monoids} 
  it suffices to show this for a set of generating cofibrations, which can be taken
  of the form~ $R\boxtimes i$ for $i$~in a set of flat cofibrations
  of orthogonal spaces (for example the set $I^{\str}$ defined in
  the proof of the strong level model structure, 
  Proposition \ref{prop:strong level spaces}).
  A box product, over $R$, of free $R$-modules induced from orthogonal spaces
  is isomorphic to the free $R$-module generated by the box product of
  underlying orthogonal spaces:
  \[ (R\boxtimes X)\boxtimes_R (R\boxtimes Y)\ \iso \ R\boxtimes (X\boxtimes Y) \]
  Since $R\boxtimes-$ is a left adjoint, it commutes with
  pushouts and orbits by $\Sigma_n$-actions. Hence the analogous statement
  carries over to symmetrized iterated box products. In other words,
  for every morphism $i:A\to B$ of orthogonal spaces there is a natural
  isomorphism in the arrow category of $R$-modules between
  \[ (R\boxtimes i)^{\Box_R n}/\Sigma_n \ : \ 
  Q^n_R(R\boxtimes i)/\Sigma_n \ \to \ \mP^n_R(R\boxtimes B) \]
  and
  \[ R\boxtimes (i^{\Box n}/\Sigma_n)\ : \ R\boxtimes (Q^n(i)/\Sigma_n) \ \to \
  R\boxtimes  \mP^n(B) \ .\]
  If $i$ is a flat cofibration of orthogonal spaces, then so is
  the morphism $i^{\Box n}/\Sigma_n$,
  by Theorem \ref{thm:symmetrizable in orthogonal spaces}~(i).
  So the morphism $R\boxtimes (i^{\Box n}/\Sigma_n)$ is a cofibration of $R$-modules,
  hence so is the morphism $(R\boxtimes i)^{\Box_R n}/\Sigma_n$.
  This completes the proof that all cofibrations in the global
  model structure for $R$-modules of Corollary \ref{cor-lift to modules spaces}~(i)
  are symmetrizable with respect to $\boxtimes_R$. 

  Now we apply Corollary~3.6 of \cite{white-commutative monoids};
  there is a slight caveat here, because the hypotheses
  ask for the validity of the `strong commutative monoid axiom'
  (Definition~3.4 of \cite{white-commutative monoids}), which requires
  the symmetrizability of both the cofibrations and the acyclic cofibrations.
  Since the model structure on $R$-modules was lifted from 
  an absolute model structure, it is {\em not} the case that all acyclic 
  cofibrations are symmetrizable. However, \cite[Cor.\,3.6]{white-commutative monoids}
  and its proof are only about cofibrations, and don't involve the weak equivalences
  at all. So the proof of \cite[Cor.\,3.6]{white-commutative monoids}
  only needs the symmetrizability of the cofibrations, which we just
  established for the global model category of $R$-modules of 
  Corollary \ref{cor-lift to modules spaces}~(i).
  Since $R$ is cofibrant as an $R$-module, 
  \cite[Cor.\,3.6]{white-commutative monoids} shows that for every cofibrant
  commutative $R$-algebra $S$, the structure morphism $i:R\to S$
  is a cofibration of $R$-modules. Part~(i) now follows because
  commutative $R$-algebras are morphisms of ultra-commutative monoids
  with source $R$. More precisely, the category of commutative $R$-algebras
  is isomorphic to the category of ultra-commutative monoids under $R$.
  Moreover, a commutative $R$-algebra $S$ is cofibrant if and only if the
  structure morphism $i:R\to S$ is a cofibration of ultra-commutative monoids.

  (b) This is a combination of part~(a) and the fact, proved in 
  Corollary \ref{cor-lift to modules spaces}~(i),
  that all cofibrations of $R$-modules are h-cofibrations of orthogonal spaces.
  
  (c) This is a combination of part~(a) and the fact, also proved in 
  Corollary \ref{cor-lift to modules spaces}~(i),
  that if $R$ itself is flat, then all cofibrations of $R$-modules are 
  flat cofibrations of orthogonal spaces.

If remains to prove left properness of the model structure.
Pushouts in a category of commutative algebras are given
by the relative monoidal product. For ultra-commutative monoids this
means that a pushout square has the form
\[ \xymatrix@C=15mm{
R \ar[d]_j \ar[r]^-f_-\simeq & T\ar[d]^{j\boxtimes_R T}\\ 
S\ar[r]_-{S\boxtimes_R f} & S\boxtimes_R T} \]
where $S$ and $T$ are
considered as $R$-modules by restriction along $j$ respectively $f$.
For left properness we now suppose that $j$ is a cofibration and $f$ 
is a global equivalence.
By part~(a) of~(ii), the morphism $j$ is then a cofibration of $R$-modules 
in the global model structure of Corollary \ref{cor-lift to modules spaces}~(i). 
Since $R$ is cofibrant in that model structure, also $S$ is
cofibrant as an $R$-module.
Proposition \ref{prop:box cofibrant R-mod} then shows that the functor
$S\boxtimes_R -$ preserves global equivalences.
So the cobase change $S\boxtimes_R f$ of $f$ is a global equivalence.
This shows that the global model structure of ultra-commutative monoids is left proper. 
\end{proof}

\section{Global power monoids}
\label{sec:power op spaces}

In this section we investigate the algebraic structure that an
ultra-commutative multiplication produces on the $\Rep$-functor $\upi_0(R)$.
Besides an abelian monoid structure on $\pi_0^G(R)$ for every
compact Lie group $G$, this structure includes power operations and
transfer maps. We formalize this algebraic structure under the name
`global power monoid', see Definition \ref{def:power monoid}.
Theorem \ref{thm:umon operation} then says that global power monoids
are precisely the natural algebraic structure, i.e., they parametrize all
natural operations on $\upi_0(R)$ for ultra-commutative monoids.
In Construction \ref{con:transfer map} we introduce the transfer maps,
which are an equivalent way of packaging the power operations in a global
power monoid; the main properties of the transfers are summarized in
Proposition \ref{prop:transfer}. 

\medskip

Given an orthogonal monoid space $R$ (not necessarily commutative)
with multiplication morphism $\mu:R\boxtimes R\to R$ and a compact Lie group $G$,
we define a binary operation
\begin{equation}\label{eq:internal_product_monoid_space}
 +\ :\ \pi_0^G(R) \times \pi_0^G(R)\ \to\ \pi_0^G(R) \end{equation}
on the $G$-equivariant homotopy set of $R$ as the composite
\[  \pi_0^G(R) \times \pi_0^G(R) \ \xra{\ \times \ } \
\pi_0^G(R\boxtimes R) \ \xra{\ \mu_*\ }\ \pi_0^G(R)\ .  \]
The pairing $\times$ was defined 
in Construction \ref{con:pairing equivalence homotopy spaces}.
If we expand the definition, it boils down to the following explicit recipe:
if $V$ and $W$ are $G$-representations
and $x\in R(V)^G$ and $y\in R(W)^G$ are $G$-fixed points that
represent two classes in $\pi^G_0(R)$, 
then $[x]+[y]$ is represented by the $G$-fixed point
\begin{align*}
\mu_{V,W}(x,y)\ \in \ R(V\oplus W) \ .
\end{align*}

\Danger We write the pairing on the equivariant homotopy sets of $R$ 
additively because we will mostly be interested in {\em commutative}
orthogonal monoid spaces.
Obviously, the additive notation is slightly dangerous for 
non-commu\-tative orthogonal monoid spaces, because there the pairing
need not be commutative.

\medskip

The following properties of the operation `+` are direct consequences 
of the corresponding properties of the pairings `$\times$', 
compare Proposition \ref{prop:unstable product properties};
a direct proof from the explicit definition of the operation~`+' above
is also straightforward. 

\begin{cor}\label{cor-product properties monoid space} 
  Let $R$ be an orthogonal monoid space.
  \begin{enumerate}[\em (i)]
  \item 
    For every compact Lie group $G$ 
    the binary operation $+$ makes the set $\pi_0^G(R)$ into a monoid.
  \item 
    If the multiplication of $R$ is commutative, then so is the operation $+$.
  \item 
    The restriction map $\alpha^*:\pi_0^G(R)\to\pi_0^K(R)$ associated to
    a continuous homomorphism $\alpha:K\to G$ between compact Lie groups
    is a monoid homomorphism.
  \end{enumerate}
\end{cor}

Now we turn to special features that happen for ultra-commutative monoids.
If the multiplication on an orthogonal monoid space $R$ is commutative,
then this does not only imply commutativity of the monoids $\pi_0^G(R)$;
strict commutativity of the multiplication also gives rise to additional 
{\em power operations} that we discuss now.
An important special case will later be the 
multiplicative ultra-commutative monoid $\Omega^\bullet R$
arising from an ultra-commutative ring spectrum $R$.
In this situation the power operations satisfy further compatibility conditions
with respect to the addition and the transfer maps on $\upi_0(\Omega^\bullet R)=\upi_0(R)$;
altogether this structure makes altogether makes the 0-th equivariant homotopy groups
of an ultra-commutative ring spectrum into a {\em global power functor}.

\begin{construction}\label{con:wreath product morphisms}
We let $R$ be an ultra-commutative monoid, $G$ a compact Lie group
and $m\geq 1$.\index{subject}{power operation!in an ultra-commutative monoid}\index{symbol}{$[m]$ - {$m$-th power operation in an ultra-commutative monoid}} 
We construct a natural {\em power operation}
\begin{equation}  \label{eq:define_power_map}
 [m] \ : \ \pi_0^G ( R ) \ \to \ \pi^{\Sigma_m\wr G}_0 ( R )   
\end{equation}
that is an equivariant refinement of the map $x\longmapsto m\cdot x$.

We recall that the wreath product $\Sigma_m\wr G$\index{subject}{wreath product}\index{symbol}{$\Sigma_m\wr G$ - {wreath product}} of a symmetric group $\Sigma_m$ 
and a group $G$ is the semidirect product
\[ \Sigma_m\wr G \ = \ \Sigma_m \ltimes G^m\] 
formed with respect to the action of $\Sigma_m$ by permuting the
factors of $G^m$. So the multiplication in $\Sigma_m\wr G$ is given by
\[ (\sigma;\,g_1,\dots,g_m) \cdot (\tau;\,k_1,\dots, k_m) \ = \ 
(\sigma\tau;\,g_{\tau(1)}k_1,\dots,g_{\tau(m)}k_m) \ . \]
For every $G$-space $E$,
the wreath product $\Sigma_m\wr G$ acts on the space $E^m$ by
\[ (\sigma;\,g_1,\dots,g_m) \cdot (e_1,\dots, e_m) \ = \ 
(g_{\sigma^{-1}(1)}e_{\sigma^{-1}(1)},\dots,g_{\sigma^{-1}(m)} e_{\sigma^{-1}(m)}) \ . \]
For every $G$-representation $V$, this action even makes $V^m$ 
into a $(\Sigma_m\wr G)$-representation.
We let
\[ \mu_{V,\dots,V}\ : \
R(V)\times\dots\times R(V) \ \to \ R(V\oplus\dots\oplus V) \]
denote the $(V,\dots,V)$-component of the iterated multiplication map of $R$,
and we observe that this map is $(\Sigma_m\wr G)$-equivariant because the multiplication
of $R$ is commutative. 
If $x\in R(V)^G$ is a $G$-fixed point representing a class in $\pi^G_0 (R)$,
then $(x,\dots,x) \in R(V)^m$ is a $(\Sigma_m\wr G)$-fixed point.
So its image under the map $\mu_{V,\dots,V}$ 
is a $(\Sigma_m\wr G)$-fixed point of $R(V^m)$,
representing an element 
\[   [m]( [x] ) \ = \ \td{ \mu_{V,\dots,V}(x,\dots,x) } 
\ \in \ \pi_0^{\Sigma_m\wr G} ( R )\ .   \]
If we stabilize $x$ along a $G$-equivariant
linear isometric embedding $\varphi:V\to W$
to $R(\varphi)(x) \in R(W)^G$, 
then $\mu_{V,\dots,V}(x,\dots,x)$ changes into 
\[ \mu_{W,\dots,W} ( R(\varphi)(x),\dots, R(\varphi)(x))
\ = \ R(\varphi^m) (\mu_{V,\dots,V}(x,\dots,x)) \ 
\in \  R(W^m)^{\Sigma_m\wr G}\ . \]
Since $\varphi^m:V^m\to W^m$ is a $(\Sigma_m\wr G)$-equivariant
linear isometric embedding, 
this element represents the same class in $\pi_0^{\Sigma_m\wr G}(R)$
as $ \mu_{V,\dots,V}(x,\dots,x)$, so the class 
$ [m]([x])$ only depends on the class of $x$ in $\pi_0^G ( R )$.
We have thus constructed a well-defined power operation \eqref{eq:define_power_map}.

The power operations are clearly natural for homomorphisms $\varphi:R\to S$
of ultra-commutative monoids, i.e., for every compact Lie group $G$,
every $m\geq 0$ and all $x\in \pi_0^G(R)$ the relation
\[ [m]( \varphi_* (x))\ = \ \varphi_* ( [m](x) ) \]
holds in $\pi_0^{\Sigma_m\wr G}(S)$.
\end{construction}

The power operations $[m]$ satisfy various properties reminiscent 
of the map $x\mapsto m\cdot x$ in an abelian monoid.
We formalize these properties into the concept of a {\em global power monoid}.
In the definition we need certain homomorphisms between different wreath products, 
so we fix notation for these now.
We use the plus symbol for the `concatenation' group monomorphism
\[ +  \ : \ \Sigma_i\times \Sigma_j \ \to \ \Sigma_{i+j} \]
defined by
\[ (\sigma+\sigma')(k) \ = \ 
\begin{cases}
 \ \sigma(k) & \text{ for $1\leq k\leq i$, and}\\
  \sigma'(k-i)+i & \text{ for $i+1\leq k\leq i+j$.}
\end{cases}
\]
This operation is strictly associative, so we will leave out parentheses.
The operation + is {\em not} commutative, but the permutations
$\sigma+\sigma'$ and $\sigma'+\sigma$ differ by conjugation with the $(i,j)$-shuffle.
An embedding of a product of wreath products is now defined by
\begin{align}\label{eq:wreath_sum}  
\Phi_{i,j}\ : \ (\Sigma_i\wr G)\ \times\ (\Sigma_j\wr G)\hspace*{1.3cm} &\to\quad \Sigma_{i+j}\wr G \\
 ((\sigma;\, g_1,\dots,g_i),\, (\sigma';\, g_{i+1},\dots,g_{i+j})) \ &\longmapsto \
(\sigma+\sigma';\, g_1,\dots,g_{i+j})\ .\nonumber
\end{align}
Another group monomorphism
\[ \natural \ : \ \Sigma_k\wr \Sigma_m\ \to \ \Sigma_{k m}  \]
is defined by\index{symbol}{$\natural$ - {monomorphism $\Sigma_k\wr \Sigma_m\to\Sigma_{k m}$}}
\[ (\sigma\natural(\tau_1,\dots,\tau_k))((i-1)m +j)\ = \ (\sigma(i)-1)m + \tau_i(j) \ ,\]
for $1\leq i\leq k$ and $1\leq j\leq m$.
This yields an  embedding of an iterated wreath product 
\begin{align}\label{eq:wreath_iterate} 
\Psi_{k,m}\ : \ \Sigma_k\wr (\Sigma_m\wr G) \qquad
&\to\qquad \Sigma_{k m}\wr G \\
 (\sigma;\, (\tau_1;\, g^1),\dots,(\tau_k;\, g^k)) \ &\longmapsto \
(\sigma\natural(\tau_1,\dots,\tau_k);\, g^1+\dots+g^k)\ .\nonumber
\end{align}
Here each $g^i=(g^i_1,\dots,g^i_m)$ is an $m$-tuple of elements of $G$, and 
\[ g^1+\dots+g^k\ = \ 
(g^1_1,\dots,g^1_m,g^2_1,\dots,g^2_m,\dots,g^k_1,\dots,g^k_m) \]
denotes the concatenation of the tuples.

\begin{rk}
The formula for the homomorphism $\natural:\Sigma_k\wr \Sigma_m\to \Sigma_{k m}$
may seems slightly ad hoc, but it can be motivated 
in a more conceptual way as a composite 
\[ \Sigma_k\wr \Sigma_m\ \to \ \Sigma_{\{1,\dots,k\}\times\{1,\dots,m\}} \ \iso \ 
\Sigma_{k m} \ .\]
The first monomorphism sends $(\sigma;\, \tau_1,\dots,\tau_k)$ to the permutation
of the product set $\{1,\dots,k\}\times\{1,\dots,m\}$ defined by
\[ (i,j)\ \longmapsto \ (\sigma(i),\tau_i(j))\ . \]
The second isomorphism is conjugation by the lexicographic ordering
\[ \{1,\dots,k\}\times\{1,\dots,m\} \ \iso \ \{1,\dots,k m\} \ , \quad
(i,j)\ \longmapsto \ (i-1)m +j\ . \]
The use of the lexicographic ordering (and hence the precise formula
for the homomorphism $\natural$) is not essential here:
if we use a different bijection between the sets
$\{1,\dots,k\}\times\{1,\dots,m\}$ and $\{1,\dots,k m\}$, then
the homomorphisms $\natural$ and $\Psi_{k,m}$ change by inner automorphisms.
So the conjugacy classes of $\natural$ and $\Psi_{k,m}$
(but not the actual homomorphisms) are canonical.
Since we will always hit $\Psi_{k,m}$ with functors that are invariant
under conjugation, this should motivate that the construction is reasonably natural.
\end{rk}

\begin{defn}\label{def:power monoid} \index{subject}{global power monoid}
A {\em global power monoid} is a functor\index{subject}{abelian Rep-monoid}
\[ M \ : \ \Rep^{\op} \ \to \ \Ab Mon  \]
from the opposite of the category $\Rep$ of compact Lie groups 
and conjugacy classes of continuous homomorphisms to the category of abelian monoids,
equipped with monoid homomorphisms\index{symbol}{  $\Ab Mon$ - {category of abelian monoids}}
\[ [m] \ : \ M(G) \ \to \ M(\Sigma_m\wr G)\]
for all compact Lie groups $G$ and $m\geq 1$,
called {\em power operations}, that satisfy the following relations.
\begin{enumerate}[(i)]
\item (Identity) The operation $[1]$ is restriction along the 
preferred isomorphism $\Sigma_1\wr G\iso G$, $(1;g)\mapsto g$.
\item (Naturality) 
For every continuous homomorphism $\alpha:K\to G$ between compact Lie groups
and every $m\geq 1$ the relation
\[ [m]\circ \alpha^* \ = \ (\Sigma_m\wr\alpha)^*\circ [m]\]
holds as homomorphisms $M(G)\to M(\Sigma_m\wr K)$.
\item (Transitivity)
For all compact Lie groups $G$ and all $k,m\geq 1$ the relation
\[  \Psi^*_{k,m}\circ  [k m]\ = \ [k]\circ  [m] \]
holds as homomorphisms  $M(G)\to M(\Sigma_k\wr(\Sigma_m \wr G))$,
where $\Psi_{k,m}$ is the monomorphism \eqref{eq:wreath_iterate}. 
\item (Additivity)
For all compact Lie groups $G$, all $m>i>0$ and all $x\in M(G)$
the relation
\[  \Phi^*_{i,m-i}( [m](x) )\ = \ p_1^*([i](x)) \ + \ p_2^*([m-i](x)) \]
holds in $M((\Sigma_i\wr G)\times(\Sigma_{m-i}\wr G))$,
where $\Phi_{i,m-i}$ is the monomorphism \eqref{eq:wreath_sum} and
\[p_1:(\Sigma_i\wr G)\times(\Sigma_{m-i}\wr G)\to \Sigma_i\wr G  
\text{\quad and\quad}
p_2:(\Sigma_i\wr G)\times(\Sigma_{m-i}\wr G)\to \Sigma_{m-i}\wr G \]
are the two projections.  
\end{enumerate}
A {\em morphism} of global power monoids is a natural transformation
of abelian monoid valued functors that also commutes with the power operations $[m]$
for all $m\geq 1$.\index{subject}{morphism!of global power monoids}
\end{defn}

\begin{rk}\label{rk:power operation remarks}
In any abelian Rep-monoid $M$ we can define external pairings 
\[ 
 \oplus \ :\ M(G) \times M(K)\ \to\ M(G\times K) 
\text{\quad by\quad}
 x\oplus y \ = \ p_G^*(x)+p_K^*(y) \ ,
 \]
where $p_G:G\times K\to G$ and $p_K:G\times K\to K$ 
are the two projections.
In this notation, the additivity requirement in Definition \ref{def:power monoid}
becomes the relation 
\begin{equation} \label{eq:external_additivity}
  \Phi^*_{i,m-i}( [m](x) )\ = \ [i](x) \ \oplus [m-i](x)\ .   
\end{equation}
In a global power monoid, 
the power operations are also additive with respect to the external addition:
for all compact Lie groups $G$ and $K$ and all $m\geq 1$, and 
all classes $x\in M(G)$ and $y\in M(K)$ the relation
\[ [m](x\oplus y) \ = \ \Delta^*( [m](x)\oplus [m](y)) \]
holds in $M(\Sigma_m\wr(G\times K))$, where $\Delta$ is 
the `diagonal' monomorphism
\begin{align}\label{eq:wreath_diagonal}
\Delta\ :\ \Sigma_m\wr(G\times K) \qquad &\to\qquad (\Sigma_m\wr G)\times(\Sigma_m\wr K)\\
 (\sigma;\, (g_1,k_1),\dots,(g_m,k_m)) \ &\longmapsto \ 
((\sigma;\, g_1,\dots,g_m),\, (\sigma;\, k_1,\dots,k_m)) \ .\nonumber
\end{align}
Indeed, $\Delta$ factors as the composite
\begin{align*}
 \Sigma_m\wr(G\times K)\ \xra{\Delta_{\Sigma_m\wr(G\times K)}}\qquad  
&(\Sigma_m\wr(G\times K))\times (\Sigma_m\wr(G\times K))\\ 
\xra{(\Sigma_m\wr p_G)\times(\Sigma_m\wr p_K)}\  
&(\Sigma_m\wr G)\times(\Sigma_m\wr K)\ .  
\end{align*}
So
\begin{align*}
  [m](x\oplus y)\ &= \   [m](p_G^*(x) + p_K^*(y))\ = \ 
 [m](p_G^*(x)) + [m](p_K^*(y))\\ 
&= \ (\Sigma_m\wr p_G)^*([m](x)) + (\Sigma_m\wr p_K)^*([m](y))  \\
&= \ \Delta_{\Sigma_m\wr(G\times K)}^*( (\Sigma_m\wr p_G)^*([m](x)) \oplus (\Sigma_m\wr p_K)^*([m](y)))  \\
 &= \ \Delta_{\Sigma_m\wr(G\times K)}^*( ((\Sigma_m\wr p_G)\times(\Sigma_m\wr p_K))^*( [m](x) \oplus [m](y))) \\  
&= \ \Delta^*( [m](x)\oplus [m](y))\ .
\end{align*}
The class $[m](x)$ is an equivariant refinement of $m\cdot x=x+\ldots+x$
($m$ summands) in the following sense. 
Applying the relation \eqref{eq:external_additivity} repeatedly
shows that $[m](x)$ restricts to the external $m$-fold sum
\[ x\oplus\ldots\oplus x \ \in \ M(G^m) \]
on the normal subgroup $G^m$ of $\Sigma_m\wr G$. Restricting further to the
diagonal takes the $m$-fold external sum to $m\cdot x$ in $M(G)$.
  \end{rk}

We will soon discuss that the power operations of an ultra-commutative monoid
define a global power monoid. One aspect of this is the additivity 
of the power operations, which could be shown directly from the definition. 
However, we will use this opportunity to establish
a very general (and rather formal) additivity result 
that we will use several times in this book.
We let $\Cc$ be a category with a zero object and finite coproducts. 
We let $X\vee Y$ be a coproduct of two objects $X$ and $Y$
with universal morphisms $i:X\to X\vee Y$ and $j:Y\to X\vee Y$.
If $f:X\to A$ and $g:Y\to A$ are two morphisms with common target,
we denote by $f+g:X\vee Y\to A$ the unique morphism that satisfies
$(f+g)i=f$ and $(f+g)j=g$.

We call a functor from $\Cc$ to the category of
abelian monoids {\em reduced} if it takes every zero object in $\Cc$
to the trivial monoid. We call the functor $F$ {\em additive}
if for every pair of objects $X,Y$ of $\Cc$ the map
\[ ( F(\Id_X+0), F(0+\Id_Y)) \ : \ F(X\vee Y)\ \to \ F(X)\times F(Y) \]
is bijective (and hence an isomorphism of monoids).

\begin{prop}\label{prop:additivity prop}
Let $\Cc$ be a category with a zero object and finite coproducts and 
\[ F,G\ :\ \Cc\ \to\ \Ab Mon \]
two reduced functors to the category of
abelian monoids. Suppose that the functor $G$ is additive.
Then every natural transformation of set valued functors from $F$ to $G$
is automatically additive.
\end{prop}
\begin{proof}
We let $\tau:F\to G$ be a natural transformation of set valued functors.
We consider two classes $x, y\in F(X)$.
We let $i,j:X\to X\vee X$ be the two inclusions into the coproduct.
We claim that
\begin{equation}\label{eq:times y additive}
\tau_{X\vee X}(F(i)(x)+ F(j)(y))\  = \  G(i)(\tau_X(x)) + G(j)(\tau_X(y))
\end{equation}
in the abelian monoid $G(X\vee X)$.
Indeed,
\begin{align*}
 G(\Id_X+0) ( \tau_{X\vee X}(F(i)(x)+ F(j)(y)))\ &= \ 
\tau_X( F(\Id_X+0) (F(i)(x)+ F(j)(y)))\\ 
&= \ \tau_X( F(\Id_X)(x)+ F(0)(y))\ = \ \tau_X(x) \\ 
&= \ G(\Id_X)(\tau_X(x)) +  G(0)(\tau_X(x))\\
&= \  G(\Id_X+0)( G(i)(\tau_X(x)) + G(j)(\tau_X(y)))
\end{align*}
in $G(X)$. Similarly, 
\[  G(0+\Id_X) ( \tau_{X\vee X}(F(i)(x)+ F(j)(y)))\ = \ 
G(0+\Id_X)( G(i)(\tau_X(x)) + G(j)(\tau_X(y)))\ . \]
Since $G$ is additive, this shows the relation \eqref{eq:times y additive}.
We let $\nabla=(\Id+\Id):X\vee X\to X$ denote the fold morphism,
so that
\[ F(\nabla)(F(i)(x)+ F(j)(y)) \ = \ F(\nabla i)(x)+ F(\nabla j)(y)\ = \  x + y\ . \]
Then
\begin{align*}
\tau_X(x+y) \ &= \ \tau_X( F(\nabla)(F(i)(x)+ F(j)(y))) \\
&= \ G(\nabla)(\tau_{X\vee X}(F(i)(x)+ F(j)(y))) \\
_\eqref{eq:times y additive}& = \ G(\nabla)( G(i)(\tau_X(x)) + G(j)(\tau_X(y)))
\\
&= \ G(\nabla i)(\tau_X(x)) + G(\nabla j)(\tau_X(y)) \ = \ \tau_X(x)  +  \tau_X(y) \ . 
\qedhere\end{align*}
\end{proof}

\begin{prop}\label{prop:power op unstable} 
Let $R$ be an ultra-commutative monoid.
Then the binary operations \eqref{eq:internal_product_monoid_space}
and the power operations \eqref{eq:define_power_map} make the
Rep-functor $\upi_0(R)$ into a global power monoid.
\end{prop}
\begin{proof} 
Corollary \ref{cor-product properties monoid space} 
shows that the binary operations \eqref{eq:internal_product_monoid_space}
make the Rep-functor $\upi_0(R)$ into a functor to the category of
abelian monoids.
The coproduct of ultra-commutative monoids is given by the box product,
so the two reduced functors
\[ \pi_0^G\ , \ \pi_0^{\Sigma_m\wr G}\ : \ umon \ \to \ \Ab Mon  \]
are additive by Corollary \ref{cor-pi_0 of box}.
Since the power operation $[m]:\pi_0^G(R)\to\pi_0^{\Sigma_m\wr G}(R)$ 
is natural in $R$, 
Proposition \ref{prop:additivity prop}, applied to the category of
ultra-commutative monoids, shows that $[m]$ is additive.
The identity~(i), naturality~(ii), 
transitivity~(iii) and  additivity property~(iv)
in Definition \ref{def:power monoid} of global power monoids 
are straightforward, and we omit the proofs.
\end{proof}

\begin{eg} We let $M$ be a commutative topological monoid. Then the
constant orthogonal space $\underline{M}$ is 
naturally an ultra-commutative monoid.  
Moreover, the equivariant homotopy functor $\upi_0(\underline{M})$
is constant with value $\pi_0(M)$, and monoid structure induced
from the multiplication of $M$. The power operation
\[ [m]\ :\ \pi_0(M) = \pi_0^G(\underline{M}) \to \pi_0^{\Sigma_m\wr G}(\underline{M}) 
= \pi_0(M) \]
then sends an element $x$ to $m\cdot x$.
\end{eg}

\begin{eg}[Naive units of an orthogonal monoid space]\label{eg:R^naive times}\index{subject}{units!of an orthogonal monoid space!naive}
Every orthogonal monoid space $R$ contains an interesting  orthogonal monoid subspace
$R^{n \times}$, the {\em naive units} of $R$.\index{symbol}{$R^{n\times}$ - {naive units of the orthogonal monoid space $R$}} The value of $R^{n \times}$ 
at an inner product space $V$ is the union of those path components of $R(V)$ 
that are taken to invertible elements, with respect to the monoid structure 
on $\pi_0(R)$, under the map
\[ R(V)\ \to \ \pi_0(R(V)) \ \to \ \pi^e_0(R) \ . \]
In other words, a point $x\in R(V)$ belongs to $R^{n\times}(V)$ if and only
if there is an inner product space $W$ and a point $y\in R(W)$ such that
\[ \mu_{V,W}(x,y) \ \in  \ R(V\oplus W)
\text{\qquad and\qquad}
\mu_{W,V}(y,x) \ \in  \ R(W\oplus V) \]
are in the same path component as the respective unit elements.
We omit the verification that the subspaces $R^{n\times}(V)$ indeed form 
an orthogonal monoid subspace of $R$ as $V$ varies.
The induced map
\[ \upi_0(R^{n\times})\ \to \ \upi_0(R) \]
is also an inclusion, and the value $\pi_0^e(R^{n\times})$ at the trivial group is,
by construction, the set of invertible elements of $\pi_0^e(R)$.
For a general compact Lie group $G$,
\[ \pi_0^G(R^{n\times}) \ = \ 
\{ \, x\in \pi_0^G(R) \ |\ \text{$\res^G_e(x)$ is invertible in $\pi_0^e(R)$}\}\]
is the submonoid of $\pi_0^G(R)$ of elements that become invertible
when restricted to the trivial group. 
So contrary to what one could suspect at first sight, 
$\pi_0^G(R^{n\times})$ may contain non-invertible elements
and the orthogonal monoid space $R^{n\times}$ is not necessarily group-like;
this is why we use the adjective `naive'.
\end{eg}

\begin{eg}[Units of a global power monoid]\label{eg:M^times}\index{subject}{units!of a global power monoid}
Every global power mo\-noid $M$ has a global power submonoid $M^\times$ of {\em units}.
The value $M^\times(G)$ at a compact Lie group $G$ is the subgroup
of invertible elements of $M(G)$. Since the restriction maps
and the power operations are homomorphisms, 
the sets $M^\times(G)$ are closed under restriction maps
and power operations. So for varying $G$, the subgroups $M^\times(G)$ 
indeed form a global power submonoid of $M$.

We call a global power monoid $N$ {\em group-like}\index{subject}{group-like!global power monoid}
if the abelian monoid $N(G)$ is a group for every compact Lie group $G$.
If $f:N\to M$ is a homomorphism of global power monoids and $N$ is group-like, 
then the image of $f$ is contained in $M^\times$. 
So the functor $M\mapsto M^\times$ is right adjoint
to the inclusion of the full subcategory of group-like global power monoids.

If $R$ is an ultra-commutative monoid, then we introduce 
a global topological version $R^\times$ of the units 
in Construction \ref{con:R^times} below.\index{subject}{units!of an ultra-commutative monoid}
This construction comes with a homomorphism of ultra-commutative monoids
$R^\times\to R$ that realizes the inclusion of the units of $\upi_0(R)$,
compare Proposition \ref{prop:realize units in umon}.
\end{eg}

\begin{eg}[Group completion of a global power monoid]\label{eg:M^star}\index{subject}{group completion!of a global power monoid}
A morphism $j:M\to M^\star$ of global power monoids is a {\em group completion}
if for every group-like global power monoid $N$ the map
\[  j^*\ : \ \text{(global power monoids)}(M^\star,N)\ \to \ 
\text{(global power monoids)}(M,N)\]
is bijective. 
Since the pair $(M^\star,j)$ represents a functor,
it is unique up to preferred isomorphism under $M$.
Every global power monoid $M$ has a group completion, 
which can be constructed `objectwise'.
We define a global power monoid $M^\star$ at a compact Lie group $G$ by
letting $M^\star(G)$ be a group completion (Grothendieck construction) 
of the abelian monoid $M(G)$,
with $j(G):M(G)\to M^\star(G)$ the universal homomorphism.
Since the restriction maps $\alpha^*:M(G)\to M(K)$ and the power operations
$[m]:M(G)\to M(\Sigma_m\wr G)$ are monoid homomorphisms, the universal property
provides unique homomorphisms $\alpha^*:M^\star(G)\to M^\star(K)$ and 
$[m]:M^\star(G)\to M^\star(\Sigma_m\wr G)$ such that
\[ \alpha^*\circ j(G) = j(K)\circ \alpha^* \text{\qquad and\qquad}
 [m]\circ j(G) = j(\Sigma_m\wr G)\circ [m]\ . \]
The functoriality of the restriction maps $\alpha^*$ and the additional 
relations required of a global power monoid are relations between
monoid homomorphism; so they are inherited by $M^\star$ via the universal property
of group completion of abelian monoids.

If $R$ is an ultra-commutative monoid, then 
in Construction \ref{con:R^star} below we introduce 
a global topological version $R^\star$ of the group completion.\index{subject}{global group completion!of an ultra-commutative monoid}
This construction comes with a homomorphism of ultra-commutative monoids
$R\to R^\star$ that realizes the algebraic group completion,
compare Proposition \ref{prop:group completion completes pi_0}.
\end{eg}

\begin{eg}[Free ultra-commutative monoid of a global classifying space]\index{subject}{free ultra-commutative monoid!of a global classifying space}\label{eg:P B_gl G}
We look more closely at the free ultra-commutative monoid $\mP(B_{\gl}G)$ 
generated by the global classifying space $B_{\gl} G$ 
of a compact Lie group $G$.\index{subject}{global classifying space}
For every $G$-representation $V$, Example \ref{eg:box of free orthogonal} 
provides an isomorphism of orthogonal spaces
\[ \bL_{G,V}^{\boxtimes m}\ \iso \ \bL_{G^m,V^m}  \ . \]
At an inner product space $W$, 
the permutation  action of $\Sigma_m$ on the left hand side becomes
the action on
\[ \bL_{G^m,V^m} (W) \ = \ \bL(V^m,W)/G^m \]
by permuting the summands in $V^m$. 
Passage to $\Sigma_m$-orbits gives an isomorphism
  \[ \mP^m(\bL_{G,V}) \ = \ \bL_{G,V}^{\boxtimes m}/\Sigma_m \ 
\iso \ \bL_{G^m,V^m} /\Sigma_m \ \iso \ \bL_{\Sigma_m\wr G,V^m}  \ . \]
Thus
\[ \mP(\bL_{G,V}) \ = \ {\coprod}_{m\geq 0} \, \bL_{\Sigma_m\wr G,V^m} \ .  \]
If $G$ acts faithfully on $V$ and $V\ne 0$, then the action of
$\Sigma_m\wr G$ on $V^m$ is again faithful. So in terms of global classifying spaces
the free ultra-commutative monoid generated by $B_{\gl}G$ is given by
\begin{equation}\label{eq:P( B_gl G)}
  \mP( B_{\gl} G) \ = \ {\coprod}_{m\geq 0}\,  B_{\gl}(\Sigma_m\wr G) \ .
\end{equation}
 The tautological class $u_G\in \pi_0^G(B_{\gl} G)$ is represented by the orbit of the
identity of $V$ in
\[ (\bL_{G,V}(V))^G \ = \ (\bL(V,V)/G)^G\ ,\]
compare \eqref{eq:tautological_class}.
So the class $[m](u_G)\in \pi_0^{\Sigma_m\wr G}( \mP( B_{\gl} G)) $
is represented by the orbit of the identity of $V^m$ in
\[ (\bL_{\Sigma_m\wr G,V^m}(V^m))^G \ = \ (\bL(V^m,V^m)/{\Sigma_m\wr G})^{\Sigma_m\wr G}\ ;\]
so with respect to the identification \eqref{eq:P( B_gl G)}
we have
\begin{equation}  \label{eq:P^m_of_u_G}
 [m](u_G)\ = \ u_{\Sigma_m\wr G} \ .   
\end{equation}
\end{eg}

\begin{eg}[Coproducts of ultra-commutative monoids]
The category of ultra-commutative monoids is cocomplete;
in particular, every family $\{R_i\}_{i\in I}$
of ultra-commutative monoids has a coproduct that we denote
$ \boxtimes_{i\in I}'\, R_i$.
We claim that the functor
\[ \upi_0\ : \  umon\ \to \ \text{(global power monoids)} \]
preserves coproducts.
Indeed, if the indexing set $I$ is finite, then the underlying orthogonal space
of the coproduct $\boxtimes_{i\in I}'\, R_i$ is simply the iterated
box product of the underlying orthogonal spaces.
The functor $\upi_0$ takes box products 
of orthogonal spaces to the objectwise product of Rep-functors,
by Corollary \ref{cor-pi_0 of box}.
In the category of abelian monoids, finite products are also finite coproducts.
Coproducts of global power monoids are formed objectwise,
so a finite product of global power monoids is also a coproduct.
This proves the claim whenever the indexing set $I$ is finite.

In any category, an infinite coproduct is the filtered colimit of the finite coproducts.
Moreover, filtered colimits of ultra-commutative monoids are
formed on underlying orthogonal spaces,
compare Corollary \ref{cor:co-limits in umon}.
So if the set $I$ is infinite, then the underlying orthogonal space of
the coproduct is the `infinite box product'
in the sense of Construction \ref{con:infinite box},\index{subject}{box product!of orthogonal spaces!infinite} i.e., the filtered colimit, 
formed over the poset of finite subsets of $I$, of the finite coproducts,
\[ \boxtimes_{i\in I}'\, R_i\ \iso \ 
\colim_{J\subset I, J\text{ finite}} \, \boxtimes_{j\in J} R_j\ .\]
Proposition \ref{prop:pi_0 of infinite box}
thus provides a bijection
\[ {\prod}'_{i\in I} \, \pi_0^G(R_i)\ \to \ \pi_0^G(\boxtimes'_{i\in I} R_i) \]
from the weak product of the abelian monoids $\pi_0^G(R_i)$
to the abelian monoid $\pi_0^G(\boxtimes'_{i\in I} R_i)$.
For abelian monoids, the weak product is also the direct sum, i.e.,
the categorical coproduct.
Since colimits of global power monoids are calculated objectwise,
this proves the claim in general.
\end{eg}

Here is another family of global power monoids, with underlying Rep-func\-tor
represented by an abelian compact Lie group $A$.
In fact, the next proposition shows that $\Rep(-,A)$ is
a free global power monoid subject to a specific set of explicit `power relations'.
In Construction \ref{con:multiplicative B G} below we exhibit a 
multiplicative model for the global classifying space of $A$, 
i.e., an ultra-commutative monoid that realizes 
the global power monoid $\Rep(-,A)$ on $\upi_0$.

\begin{prop}\label{prop:unique power Rep(-,A)}
\index{subject}{compact Lie group!abelian}
For every abelian compact Lie group $A$, the Rep-functor $\Rep(-,A)$
has a unique structure of global power monoid. 
The monoid structure of $\Rep(G,A)$ is given by
pointwise multiplication of homomorphisms. The power operation
\[ [m]\ : \ \Rep(G,A)\ \to \ \Rep(\Sigma_m\wr G,A) \]
is given by
\[ [m](\alpha)\ = \ p_m\circ (\Sigma_m\wr\alpha)\ , \]
where $p_m:\Sigma_m\wr A\to A$ is the homomorphism defined by
$p_m(\sigma;\,a_1,\dots,a_m) = a_1\cdot\ldots\cdot a_m$.
Moreover, for every global power monoid $M$ the map  
\[ \text{\em (global power monoids)}(\Rep(-,A),M) \ \to \  M(A)\ ,
\quad f \longmapsto \ f(A)(\Id_A)\]
is injective with image those $x\in M(A)$ that satisfy $[m](x)=p_m^*(x)$ 
for all $m\geq 1$.
\end{prop}
\begin{proof}
Since $A$ is abelian, conjugate homomorphisms into $A$ are already equal,
i.e., we can ignore the difference between homomorphisms and conjugacy classes.
We establish the monoid structure and power opera\-tions first, 
which also shows the uniqueness.
Since $\Rep(e,A)$ has only one element, it is the additive unit.  
Since restriction maps are monoid homomorphisms,
the trivial homomorphism is the neutral element of $\Rep(G,A)$.
The sum
\[ q_1 + q_2 \ \in \ \Rep(A\times A,A) \]
of the two projections is a homomorphism from $A\times A$ to $A$
whose restriction along the two maps $(-,1),(1,-):A\to A\times A$
is the identity. The only such homomorphism is the multiplication
$\mu:A\times A\to A$, so we conclude that
\[ q_1 + q_2 \ = \ \mu\ .\]
Naturality now gives
\[ \alpha+\beta \ = \  (\alpha,\beta)^*(q_1) +(\alpha,\beta)^*(q_2) \ = \ 
(\alpha,\beta)^*(q_1+q_2) \ = \ (\alpha,\beta)^*(\mu) \ = \ \mu\circ(\alpha,\beta)  \ . \]
Since power operations refine power maps, the element
\[ [m](\Id_A)\  \in \ \Rep(\Sigma_m\wr A,A) \]
restricts to the sum of the $m$ projections on $A^m\leq \Sigma_m\wr A$. 
We let $1:e\to A$ be the unique homomorphism and claim that the composite
\[ \Sigma_m\wr e \ \xra{\Sigma_m\wr 1} \ \Sigma_m\wr A\  \xra{[m](\Id_A)}\ A \]
is trivial. Indeed,
\[ [m](\Id_A) \circ (\Sigma_m\wr 1)\ = \ 
 (\Sigma_m\wr 1)^*([m](\Id_A))\ = \  [m](1^*(\Id_A))\ = \ [m](1)\ = \ 1 \ ,\]
using that the operation $[m]$ is a monoid homomorphism. Thus 
\begin{align*}
 [m](\Id_A)(\sigma;\,a_1,\dots,a_m) \ &= \ 
[m](\Id_A)(\sigma;\,1,\dots,1) \cdot [m](\Id_A)(1;\,a_1,\dots,a_m) \\ 
&= \ a_1\cdot\ldots\cdot a_m\ .
\end{align*}
In other words, $[m](\Id_A)=p_m$. Naturality now gives
\begin{align*}
[m](\alpha) \ = \  [m](\alpha^*(\Id_A)) \ &= \ (\Sigma_m\wr\alpha)^*([m](\Id_A))\\ 
&= \ (\Sigma_m\wr\alpha)^*(p_m)\ =\ p_m\circ (\Sigma_m\wr\alpha) \ .
\end{align*}

It remains to show the existence of the global power monoid structure.
Clearly, pointwise multiplication of homomorphisms makes $\Rep(G,A)$ into
an abelian monoid (even an abelian group), and the monoid structure is
contravariantly functorial in $G$.
When we define $[m](\alpha)$ by the formula of the proposition,
then the remaining axioms of a global power monoid 
(compare Definition \ref{def:power monoid}) are similarly straightforward.
The identity property~(i) is clear, and naturality~(ii) follows from the relation
\[ [m](\alpha)\ = \ (\Sigma_m\wr \alpha)^*([m](\Id_A)) \ .\]
The transitivity relation~(iii) holds by inspection:
\begin{align*}
\Psi_{k,m}^*([k m](\alpha))(\sigma;(\tau_1;h^1),&\dots,(\tau_k,h^k)) \\ 
&= \  
([k m](\alpha))(\Psi_{k,m}(\sigma;(\tau_1;h^1),\dots,(\tau_k,h^k))) \\ 
&= \  \prod_{i=1}^k \left(\prod_{j=1}^m h_j^i\right) \
= \ \prod_{i=1}^k ([m](\alpha))(\tau_1;h^i) )\\ 
&= \ ([k]( [m](\alpha)))(\sigma;(\tau_1;h^1),\dots,(\tau_k,h^k)) \ ,
\end{align*}
and so does the additivity relation~(iv):
\begin{align*}
 \Phi^*_{i,m-i}(  [m](\alpha) ) &((\sigma;\,g_1,\dots,g_i),(\sigma';\,g_{i+1},\dots,g_m)) \\ 
&= \ ([m](\alpha)) (\sigma+\sigma';\,g_1,\dots,g_m) \\ 
&= \ \alpha(g_1)\cdot\dots\cdot\alpha(g_m) \\
&= \ ([i](\alpha)) (\sigma;\,g_1,\dots,g_i) \cdot
([m-i](\alpha)) (\sigma';\,g_i,\dots,g_m) \\  
&= \
([i](\alpha)\oplus [m-i](\alpha) )((\sigma;\,g_1,\dots,g_i),(\sigma';\,g_{i+1},\dots,g_m)) 
\ .
\end{align*}
It remains to identify the global power morphisms out of $\Rep(-,A)$.
Since the class $\Id_A$ generates $\Rep(-,A)$ as a $\Rep$-functor, 
the evaluation map in injective.
For surjectivity we consider a class $x\in M(A)$ such that $[m](x)=p_m^*(x)$
for all $m\geq 1$. The Yoneda lemma then provides a unique morphism of $\Rep$-functors
\[ f \ : \ \Rep(-,A)\ \to \ M \]
such that $f(A)(\Id_A)=x$,
and this morphism is given by $f(G)(\alpha)=\alpha^*(x)$, for $\alpha:G\to A$.
We need to show that $f$ is a morphism of global power monoids, 
i.e., additive and compatible with power operations. 
For additivity we recall that
\[   \alpha+\beta\ = \ \mu\circ (\alpha,\beta) \ = \ 
p_2\circ \text{incl}_{A^2}^{\Sigma_2\wr A}\circ (\alpha,\beta) \ . \]
Hence
\begin{align*}
  f(G)(\alpha+\beta)\ &= \ 
(p_2\circ \text{incl}_{A^2}^{\Sigma_2\wr A}\circ (\alpha,\beta))^*(x) \\ 
&= \ (\alpha,\beta)^*(\res_{A^2}^{\Sigma_2\wr A}(p_2^*(x))) \ 
= \ (\alpha,\beta)^*(\res_{A^2}^{\Sigma_2\wr A}([2](x))) \\ 
&= \ (\alpha,\beta)^*(x\oplus x) \
= \ \alpha^*(x)+ \beta^*(x) \ = \ f(G)(\alpha) + f(G)(\beta) \ .
\end{align*}
This shows that $f$ is a morphism of abelian Rep-monoids.
Finally, given a continuous homomorphism $\alpha:G\to A$, we have
\begin{align*}
  [m](f(G)(\alpha))\ &= \ [m](\alpha^*(x))\ = \ (\Sigma_m\wr\alpha)^*([m](x))\
  = \ (\Sigma_m\wr\alpha)^*(p_m^*(x))\\ 
&=\ (p_m\circ(\Sigma_m\wr\alpha))^*(x)\ 
=\ ( [m](\alpha) )^*(x)\ = \ f(\Sigma_m\wr G)([m](\alpha))\ .
\end{align*}
So the morphism $f$ is also compatible with power operations.
\end{proof}

Now we are going to show that the restriction maps along continuous group
homomorphisms and the power operations give all natural operations between
equivariant homotopy sets of ultra-commutative monoids. 
The strategy is the same as in the analogous situation for orthogonal spaces, 
see Corollary \ref{cor:operations spc}:
natural operations for ultra-commutative monoids from the functor $\pi_0^G$ 
to the functor $\pi_0^K$ biject with the $K$-equivariant homotopy set 
of $\mP(B_{\gl}G)$, the free ultra-commutative monoids
generated by a global classifying space of $G$. 
So ultimately we need to calculate $\pi_0^K(\mP(B_{\gl}G))$.

The tautological class $u_G$ in $\pi_0^G(B_{\gl}G)$
was defined in \eqref{eq:tautological_class}.
We set
\[ u_G^{umon}\ = \ \eta_*(u_G )\ \in \ \pi_0^G(\mP(B_{\gl} G))\ ,\]
where $\eta:B_{\gl} G\to \mP(B_{\gl} G)$ is the adjunction unit, i.e., the
inclusion of the homogeneous summand for $m=1$. 
The next theorem says in particular that 
the global power monoid $\upi_0(\mP (B_{\gl} G))$ is freely generated 
by the element $u_G^{umon}$.

\begin{theorem}\label{thm:umon operation}
  Let $G$ and $K$ be compact Lie groups.\index{subject}{global classifying space}
  \begin{enumerate}[\em (i)]
  \item 
    Every class in $\pi_0^K(\mP(B_{\gl}G))$ is of the form $\alpha^*([m](u_G^{umon}))$
    for a unique $m\geq 0$ and a unique conjugacy class of continuous homomorphisms
    $\alpha:K\to\Sigma_m\wr G$.
  \item 
    For every global power monoid $M$ and every $x\in M(G)$
    there is a unique morphism of global power monoids 
    $f:\upi_0(\mP(B_{\gl} G))\to M$ such that 
    \[ f(G)(u_G^{umon})\ = \ x \ .\]
  \item
  Every natural transformation $\pi_0^G \to \pi_0^K$
  of set valued functors on the category of ultra-commutative monoids 
  is of the form $\alpha^*\circ [m]$ for a unique $m\geq 0$ and
  a unique conjugacy class of continuous group homomorphisms $\alpha:K\to \Sigma_m\wr G$.
\end{enumerate}
\end{theorem}
\begin{proof}
  For the course of the proof we abbreviate $u=u_G^{umon}$. 

  (i) By \eqref{eq:P( B_gl G)} the underlying orthogonal space of $\mP( B_{\gl} G)$
  is the disjoint union of global classifying spaces for the wreath product
  groups $\Sigma_m\wr G$; moreover, the class $[m](u)$
  lies in the $m$-th summand of $\mP( B_{\gl} G)$
  and is a universal element for $B_{\gl}(\Sigma_m\wr G)$, by \eqref{eq:P^m_of_u_G}. 
  So part~(i) follows from Proposition \ref{prop:fix of global classifying}~(ii)
  and the fact that $\pi_0^K$ commutes with disjoint unions.

  (ii) By~(i) every element of $\pi_0^K(\mP(B_{\gl} G))$ is of the form $\alpha^*([m](u))$;
  every morphism of global power monoids
  $f:\upi_0(\mP(B_{\gl} G))\to M$ satisfies
  \[ f(K)(\alpha^*([m](u))) \ = \ \alpha^*([m](f(G)(u))) \ ,\]
  so $f$ is determined by its value on the class $u$.
  This shows uniqueness.

  Conversely, if $x\in M(G)$ is given we define $f(K):\pi_0^K(\mP (B_{\gl}G))\to M(K)$
  by the formula
  \[ f(K)(\alpha^*([m](u))) \ = \ \alpha^*([m](x)) \ .\]
  Then $f(G)(u)=x$.
  It remains to show that $f$ is indeed a morphism of global power monoids. 
  This is a routine  -- but somewhat lengthy -- calculation as follows.
  Given $\alpha:K\to\Sigma_m\wr G$ and  $\bar\alpha:K\to\Sigma_n\wr G$,
  we have
  \begin{align*}
    \alpha^*([m](u))\ + \ \bar\alpha^*([n](u)) \ 
    &= \ (\alpha,\bar\alpha)^*(p_1^*([m](u)) + p_2^*([n](u))) \\
    &= \ (\alpha,\bar\alpha)^*(\Phi_{m,n}^*([m+n](u))) \\ 
    &= \ (\Phi_{m,n}\circ(\alpha,\bar\alpha))^*([m+n](u)) \ .
  \end{align*}
  So $f(K)$ is additive because
  \begin{align*}
    f(K)(\alpha^*([m](u)) +  \bar\alpha^*([n](u))) \ &= \ 
    f(K)((\Phi_{m,n}\circ(\alpha,\bar\alpha))^*([m+n](u)) )\\
    &= \ (\Phi_{m,n}\circ(\alpha,\bar\alpha))^*([m+n](x))\\
    &= \ (\alpha,\bar\alpha)^*(\Phi_{m,n}^*([m+n](x)))\\
    &= \ (\alpha,\bar\alpha)^*( p_1^*([m](x)) + p_2^*([n](x)))\\
    &= \ \alpha^*([m](x))\ + \ \bar\alpha^*([n](x))\\
    &= \ f(K)(\alpha^*([m](u)))\ + \ f(K)(\bar\alpha^*([n](u)))\ .
  \end{align*}
  For a continuous homomorphism $\beta:L\to K$ we have
  \begin{align*}
    (\beta^*\circ  f(K))(\alpha^*([m](u))) \ &= \ 
    \beta^*(\alpha^*([m](x))) \ 
    = \ (\alpha\circ\beta)^*([m](x)) \\ 
    &= \ f(L)( (\alpha\circ\beta)^*([m](u))) \ 
    = \ (f(L)\circ\beta^*)( \alpha^*([m](u)))\ ;
  \end{align*}
  so the homomorphisms $f(K)$ form a natural transformation of $\Rep$-functors.
  For $k\geq 1$ we have
  \begin{align*}
    [k](\alpha^*([m](u)))\ &= \ 
    (\Sigma_k\wr\alpha^*)([k]([m](u)))\\ 
    &= \  (\Sigma_k\wr\alpha^*)(\Psi_{k, m}^*([k m](u)))\ 
    = \  (\Psi_{k,m}\circ(\Sigma_k\wr\alpha))^*([k m](u))\ ;
  \end{align*}
  hence
  \begin{align*}
    f(\Sigma_k\wr K)([k](\alpha^*&([m](u))))\ 
    = \ (\Psi_{k,m}\circ (\Sigma_k\wr\alpha))^*([k m](x)) \\ 
    &= \ (\Sigma_k\wr\alpha)^*(\Psi_{k,m}^*( [k m](x))) \ 
    = \ (\Sigma_k\wr\alpha)^*([k]([m](x))) \\ 
    &= \ [k](\alpha^*([m](x))) \ 
    =\ [k]( f(K)(\alpha^*([m](u)))) \ .
  \end{align*}
  So the homomorphisms $f(K)$ are compatible with power operations.

  (iii) 
  We apply the representability result 
  of Proposition \ref{prop:pi_0^G representability}
  to the category of ultra-commutative monoids and
  the free and forgetful adjoint functor pair: 
  \[ \xymatrix@C=10mm{ \mP \ : \ \spc \ \ar@<.4ex>[r]  & \ umon\ : \ U
    \ar@<.4ex>[l] } \]
  If $G$ is a compact Lie group, $V$ a $G$-representation and $W$
  a non-zero faithful $G$-representation, then the restriction morphism
  $\rho_{G,V,W}:\bL_{G,V\oplus W}\to \bL_{G,W}$ is a global equivalence
  between positive flat orthogonal spaces. We showed in the proof of
  Theorem \ref{thm:symmetrizable in orthogonal spaces}~(ii)
  that the induced morphism of free ultra-commutative monoids
  $\mP(\rho_{G,V,W}):\mP(\bL_{G,V\oplus W})\to \mP(\bL_{G,W})$ is a global equivalence;
  in particular, the morphism of $\Rep$-functors
  $\upi_0(\mP(\rho_{G,V,W}))$ is an isomorphism. So 
  Proposition \ref{prop:pi_0^G representability} applies and shows that
  evaluation at the tautological class is a bijection
  \[   \Nat^{umon}(\pi_0^G,\pi_0^K)  \ \to \ \pi_0^K(\mP(B_{\gl} G)) \ , \quad
  \tau\ \longmapsto \ \tau(u)\ .\]
  Part~(i) then completes the argument.
\end{proof}

\begin{rk}[Natural $n$-ary operations]
By similar arguments as in the previous theorem
we can also identify the natural $n$-ary operations on equivariant homotopy sets
of ultra-commutative monoids. For every $n$-tuple $G_1,\dots, G_n$ 
of compact Lie groups the functor
\[ \Ho(umon)\ \to \ \text{(sets)}\ ,\quad
X \ \longmapsto \ \pi_0^{G_1}(X)\times\dots\times\pi_0^{G_n}(X) \]
is represented by the free ultra-commutative monoid
$\mP(B_{\gl} G_1 \amalg \ldots \amalg  B_{\gl} G_n)$.
So the set of natural transformations from the functor
$\pi_0^{G_1}\times\dots\times\pi_0^{G_n}$ to the functor $\pi_0^K$,
for another compact Lie group $K$, bijects with the $K$-equivariant homotopy set
of this representing object.
Because
\begin{align*}
\mP(B_{\gl} G_1 \amalg \ldots \amalg  B_{\gl} G_n)\ &\iso \
 \mP(B_{\gl} G_1)\boxtimes \dots \boxtimes \mP( B_{\gl} G_n)\\ 
&\iso \ 
{\coprod}_{j_1,\dots,j_n\geq 0}\,
B_{\gl} (\Sigma_{j_1}\wr G_1)\boxtimes \dots \boxtimes B_{\gl} (\Sigma_{j_n}\wr G_n)\\
&\iso \ 
{\coprod}_{j_1,\dots,j_n\geq 0} \,
B_{\gl} ( (\Sigma_{j_1}\wr G_1)\times \dots \times  (\Sigma_{j_n}\wr G_n))\ ,
  \end{align*}
the group
$\pi_0^K( \mP(B_{\gl} G_1 \amalg \ldots \amalg  B_{\gl} G_n))$
bijects with the disjoint union of the sets
\[ \pi_0^K(B_{\gl} ( (\Sigma_{j_1}\wr G_1)\times \dots \times  (\Sigma_{j_n}\wr G_n))) \ \iso\
\Rep(K, (\Sigma_{j_1}\wr G_1)\times \dots \times  (\Sigma_{j_n}\wr G_n))\ .  \]
So every natural operation from $\pi_0^{G_1}\times\dots\times\pi_0^{G_n}$ to $\pi_0^K$
is of the form
\[ (x_1,\dots,x_n) \ \longmapsto \ \alpha^*( [j_1](x_1) \oplus \dots\oplus [j_n](x_n)) \]
for a unique tuple $(j_1,\dots, j_n)$ of non-negative integers
and a unique conjugacy class of continuous homomorphisms
$\alpha:K\to  (\Sigma_{j_1}\wr G_1)\times \dots \times  (\Sigma_{j_n}\wr G_n)$.
In particular, the $n$-ary operations are generated by unary operations
and external sum.
\end{rk}

\begin{construction}\label{con:define P and A^+}
We describe an alternative (but isomorphic) way 
to organize the book keeping of the natural operations between the
0-th equivariant homotopy sets of ultra-commutative monoids.
We denote by $\Nat^{umon}$ the category whose objects are all compact Lie groups
and where the morphism set $\Nat^{umon}(G,K)$ is the set of all natural
transformations, of functors from the ultra-commutative monoids to sets, 
from $\pi_0^G$ to $\pi_0^K$.
We define an isomorphic algebraic category $\mathbb A^+$, 
the {\em effective Burnside category}.\index{subject}{effective Burnside category}
Both $\Nat^{umon}$and $\mathbb A^+$ are
`pre-preadditive' in the sense that all morphism sets are abelian monoids 
and composition is biadditive. In $\Nat^{umon}$, the monoid structure is
objectwise addition of natural transformations.

The category $\mathbb A^+$ has the same objects as $\Nat^{umon}$, 
namely all compact Lie groups. 
In the effective Burnside category, the morphism set
$\mA^+(G,K)$\index{symbol}{$\mA^+(G,K)$ - {monoid of isomorphism classes of $G$-free $K$-$G$-sets}} 
is the set of isomorphism classes of those $K$-$G$-spaces
that are disjoint unions of finitely many free right $G$-orbits.
This set is an abelian monoid via disjoint union of $K$-$G$-spaces.
Composition
is induced by the balanced product over $K$:
\[ \circ \ : \ \mA^+(K,L) \times \mA^+(G,K) \ \to \ \mA^+(G,L)\ , \quad
 [T]\circ [S] \ = \  [ T\times_K S] \ .\] 
Here $T$ has a left $L$-action and a commuting free right $K$-action,
whereas $S$ has a left $K$-action and a commuting free right $G$-action.
The balanced product $T\times_K S$ then inherits a left $L$-action from $T$
and a free right $G$-action from $S$.
We define a functor
\[ B \ : \  \Nat^{umon}\  \to\  \mathbb A^+  \]
as the identity on objects; 
on morphisms, the functor is given by
\[ 
B \ : \ \Nat^{umon}(G,K) \ \to \ \mathbb A^+(G,K)\ , \quad 
B(\alpha^*\circ[m])\ = \ [ \alpha^* (\{1,\dots,m\}\times G)_G] \ .
 \]
In the definition we use the characterization of the natural operations
given by Theorem \ref{thm:umon operation}~(iii).
Also, we consider 
$\{1,\dots,m\}\times G$ as a right $G$-space by translation;
the wreath product $\Sigma_m\wr G$ acts from the left on 
$\{1,\dots,m\}\times G$ by 
\begin{equation}\label{eq:wreath_acts_on_mxG}
 (\sigma;\,g_1,\dots,g_m)\cdot (i,\gamma)\ = \ (\sigma(i),\, g_i\cdot\gamma) \ .  
\end{equation}
Then we let $K$ act by restriction of the $(\Sigma_m\wr G)$-action along $\alpha$.
\end{construction}

\begin{prop}\label{prop:B iso}
 The functor $B: \Nat^{umon}\to\mathbb A^+$ is additive and an isomorphisms of categories.
\end{prop}
\begin{proof}
We start by showing that the map $B:\Nat^{umon}(G,K)\to\mA^+(G,K)$ is additive.
The map
\[ (\{1,\dots,k\}\times G) \ \amalg\ (\{1,\dots,m\}\times G)   \ \to \ 
\Phi_{k,m}^*(\{1,\dots,k+m\}\times G) \]
that is the inclusion on the first summand and given by
$(j,g)\mapsto(k+j,g)$ on the second summand
is an isomorphism of $((\Sigma_k\wr G)\times (\Sigma_m\wr G))$-$G$-bispaces.
Restriction along the homomorphism 
$(\beta,\alpha):K\to (\Sigma_k\wr G)\times (\Sigma_m\wr G)$
provides an isomorphism of $K$-$G$-spaces between
\[ \beta^*(\{1,\dots,k\}\times G)_G \ \amalg \ 
\alpha^*(\{1,\dots,m\}\times G)_G \text{\ and\ }
(\Phi_{k,m}\circ (\beta,\alpha))^*(\{1,\dots,k+m\}\times G)_G  \ . \]
This shows that
\[ B((\beta^*\circ[k])+(\alpha^*\circ[m]))\ = \ 
 B((\Phi_{k,m}\circ(\beta,\alpha))^*\circ[k+m] )\ = \ 
 B(\beta^*\circ[k])\ +\ B(\alpha^*\circ[m])\ . \]
The identity operation in $\Nat^{umon}(G,G)$ can be written as $\Id_G^*\circ[1]$;
on the other hand, the $G$-bispace $\Id_G^*(\{1\}\times G)_G$ is isomorphic to $G$
under left and right translation. Since the isomorphism class of $_G G_G$ is the
identity of $G$ in $\mathbb A^+$, the construction $B$ preserves identities.

For the compatibility of $B$ with composition we consider another
operation $\beta^*\circ[k]\in\Nat^{umon}(K,M)$ and observe that 
\[ (\beta^*\circ[k])\circ(\alpha^*\circ[m])\ = \ 
(\Psi_{k,m}\circ(\Sigma_k\wr \alpha)\circ\beta)^*\circ[ k m] \ .\]
Moreover, an isomorphism of $M$-$G$-spaces
\[  \beta^* (\{1,\dots,k\}\times K) \times_K 
\alpha^* (\{1,\dots,m\}\times G) \ \iso \ 
(\Psi_{k,m}\circ (\Sigma_k\wr\alpha)\circ\beta)^*(\{1,\dots,k m\}\times G) \]
is given by
\[ [ ( i, \kappa), (j, \gamma) ] \ \longmapsto \ 
( (i-1)m + \sigma(j)-1,\ g_j\cdot \gamma) \ ,\]
where
\[  \alpha(\kappa)\ = \ (\sigma;\,g_1,\dots,g_m)\  \in \ \Sigma_m\wr G\ .\]
So $B$ is a functor, which is bijective on objects by definition.

To see that $B$ is full we let $S$ be any $K$-$G$-space that is 
a disjoint union of $m$ free right $G$-orbits.
We choose a $G$-equivariant homeomorphism
\[ \psi \ : \ S \ \to \ \{1,\dots,m\}\times G \ .\]
We transport the left $K$-action from $S$ to 
$\{1,\dots,m\}\times G$ along $\psi$, so that $\psi$ becomes an
isomorphism of $K$-$G$-spaces. The $(\Sigma_m\wr G)$-action 
on $\{1,\dots,m\}\times G$ defined in \eqref{eq:wreath_acts_on_mxG}
identifies the wreath product with the group of 
$G$-equivariant automorphisms of $\{1,\dots,m\}\times G$;
so the $K$-action on $\{1,\dots,m\}\times G$ corresponds to a continuous homomorphism
$\alpha:K\to\Sigma_m\wr G$. Altogether, $S$ is isomorphic 
to $\alpha^*(\{1,\dots,m\}\times G)_G$.

If the $K$-$G$-spaces constructed from $\alpha^*\circ[m]$ and $\beta^*\circ[n]$
are isomorphic, then we must have $m=n$. Moreover, a $K$-$G$-isomorphism
\[ \alpha^*(\{1,\dots,m\}\times G)_G\ \iso \ \beta^*(\{1,\dots,m\}\times G)_G \]
is given by the action of a unique element $\omega\in\Sigma_m\wr G$,
and then the homomorphisms $\alpha,\beta:K\to\Sigma_m\wr G$
are conjugate by $\omega$. So the functor $B$ is faithful.
\end{proof}

Now we define {\em transfer maps}
$\tr_H^G:M(H)\to M(G)$ in global power monoids,
for every subgroup $H$ of finite index in a compact Lie group $G$.
As we will see in Proposition \ref{prop:transfer} 
below, the set of operations from $\pi_0^G$ to $\pi_0^K$ 
is a free abelian monoid with an explicit basis involving transfers.

\begin{construction}[Transfer maps]\label{con:transfer map}
In the following we let $M$ be a global power monoid,
$G$ a compact Lie group and $H$ a closed subgroup of $G$
of finite index $m$.
We choose an `$H$-basis' of $G$, i.e., an ordered $m$-tuple 
$\bar g=(g_1,\dots,g_m)$ of elements in disjoint $H$-orbits such that
\[ G \ = \ {\bigcup}_{i=1}^m \ g_i H \ .\]
The wreath product $\Sigma_m\wr H$ acts freely and transitively 
from the right on the set of all such $H$-bases of $G$, by the formula
\[ (g_1,\dots,g_m)\cdot (\sigma;\,h_1,\dots,h_m)\ = \ 
(g_{\sigma(1)} h_1,\dots,g_{\sigma(m)}h_m)\ .  \]
We obtain a continuous homomorphism $\Psi_{\bar g}:G\to\Sigma_m\wr H$ by requiring that
\[ \gamma \cdot\bar g \ = \ \bar g\cdot\Psi_{\bar g}(\gamma)  \ .\]
We define the transfer $\tr_H^G:M(H)\to M(G)$\index{subject}{transfer!in global power monoids}\index{symbol}{$\tr_H^G$ - {transfer map}} as the composite
\[ M(H) \ \xra{\ [m]\ } \ M(\Sigma_m\wr H) \ \xra{\ \Psi_{\bar g}^*\ }\ M(G) \ .\]
Any other $H$-basis is of the form
$\bar g\omega$ for a unique $\omega\in\Sigma_m\wr H$. 
We have $\Psi_{\bar g\omega}=c_\omega\circ \Psi_{\bar g}$ as maps $G\to\Sigma_m\wr H$,
where $c_\omega(\gamma)=\omega^{-1}\gamma\omega$. Since inner automorphisms
induce the identity in any Rep-functor, we conclude that
\[ \Psi_{\bar g}^* = \Psi_{\bar g\omega}^* \ : \ M(\Sigma_m\wr H) \ \to \ M(G)\ . \]
So the transfer $\tr_H^G$ does not depend on the choice of basis $\bar g$. 
\end{construction}

The various properties of the power operations
imply certain properties of the transfer maps.
Moreover, the last item of the following proposition shows that
power operations in a global power monoid are determined by the transfer
and restriction maps.
 
\begin{prop}\label{prop:transfer} 
The transfer homomorphisms of a global power monoid $M$ satisfy the following relations,
where $H$ is any subgroup of finite index in a compact Lie group $G$.
\begin{enumerate}[\em (i)]
\item {\em (Transitivity)} 
We have $\tr_G^G=\Id_{M(G)}$ and for nested subgroups
$H\subseteq G\subseteq F$ of finite index the relation
\[ \tr_G^F \circ \tr_H^G \ = \  \tr_H^F  \]
holds as maps  $M(H)\to M(F)$.
\item {\em (Double coset formula)}
For every subgroup $K$ of $G$ (not necessarily of finite index)
the relation\index{subject}{double coset formula!for transfer map}
\[ \res^G_K\circ \tr_H^G  \ = \ \sum_{[g]\in K\backslash G/H}
\tr^K_{K\cap{^g H}}\circ (c_g)^* \circ \res^H_{K^g\cap H} \]
holds as homomorphisms $M(H)\to M(K)$.
Here $[g]$ runs over a set of representatives of the finite
set of $K$-$H$-double cosets.
\item {\em (Inflation)}
For every continuous epimorphism  $\alpha:K\to G$ of compact Lie groups
the relation
\[ \alpha^*\circ \tr_H^G \ =\  \tr_L^K\circ (\alpha|_L)^* \]
holds as maps from $M(H)$ to $M(K)$, where $L=\alpha^{-1}(H)$.
\item 
For every $m\geq 1$ the $m$-th power operation can be recovered as
\[ [m] \ =\  \tr_K^{\Sigma_m\wr G}\circ q^* \ ,\]
where $K$ is the subgroup of $\Sigma_m\wr G$
consisting of all $(\sigma;\,g_1,\dots,g_m)$ such that $\sigma(m)=m$
and $q:K\to G$ is defined by $q(\sigma;\,g_1,\dots,g_m)=g_m$.
\end{enumerate}
\end{prop}
\begin{proof}
(i)
For $G=H$ we can choose the unit~1 as the $G$-basis,
and with this choice $\Psi_1:G\to\Sigma_1\wr G$ 
is the preferred isomorphism that sends $g$~to $(1;g)$. 
The restriction of $[1](x)$ along this isomorphism is $x$,
so we get $\tr_G^G(x)=x$.

For the second claim we choose a $G$-basis $\bar f=(f_1,\dots,f_k)$ of $F$
and an $H$-basis $\bar g=(g_1,\dots,g_m)$ of $G$.
Then 
\[  
 \bar f\bar g\ = \ (f_1 g_1,\dots,f_1 g_m,\ f_2 g_1,\dots, f_2 g_m,
\ \dots,\ f_k g_1,\dots,f_k g_m) \]
is an $H$-basis of $F$. With respect to this basis, the homomorphism 
$\Psi_{\bar f\bar g}:F\to\Sigma_{k m}\wr H$ equals the composite
\[ 
F\ \xra{\ \Psi_{\bar f}\ } \Sigma_k\wr G
\ \xra{\Sigma_k\wr\Psi_{\bar g}} \
 \Sigma_k\wr(\Sigma_m\wr H)  \ \xra{\ \Psi_{k,m}\ } \
 \Sigma_{k m}\wr H \]
where the monomorphism $\Psi_{k,m}$ was defined in \eqref{eq:wreath_iterate}. 
So 
\begin{align*}
\tr_H^F \ &= \ \Psi_{\bar f \bar g}^*\circ  [k m] \ = \ 
\Psi_{\bar f}^*\circ (\Sigma_k\wr\Psi_{\bar g})^*\circ \Psi_{k,m}^*\circ[k m] \\ 
&= \ \Psi_{\bar f}^*\circ (\Sigma_k\wr\Psi_{\bar g})^*\circ [k] \circ[m] \
= \ \Psi_{\bar f}^*\circ [k]\circ \Psi_{\bar g}^* \circ[m] \
= \   \tr_G^F\circ \tr_H^G \ .
\end{align*}

(ii) We choose representatives $g_1,\dots,g_r$ for the $K$-$H$-double cosets
in $G$. Then we choose, for each $1\leq i\leq r$, a $(K\cap{^{g_i}H})$-basis
\[ \bar k^i \ = \ (k^i_1,\dots,k^i_{s_i}) \]
of $K$, where $s_i=[K:K\cap {^{g_i}H}]$. 
Then $s_1+\dots+s_r=m=[G:H]$ is the index of $H$ in $G$, 
and this data provides an $H$-basis of $G$, namely
\[ \bar g \ = \ (k^1_1 g_1,\,\dots,k^1_{s_i}g_1,\,
k^2_1 g_2,\,\dots,k^2_{s_2}g_2,\,\dots,
k^r_1 g_r,\,\dots,k^r_{s_r}g_r)\ .\]
The following diagram of group homomorphisms then commutes:
\[ \xymatrix@C=10mm{ 
K \ar[r]^-{\text{incl}} \ar[d]_{(\Psi_{\bar k^1},\dots,\Psi_{\bar k^r})} &
G \ar[r]^-{\Psi_{\bar g}} & \Sigma_m\wr H \\
\prod_{i=1}^r \,\Sigma_{s_i}\wr(K\cap{^{g_i}H})
\ar[rr]_-{\prod\Sigma_{s_i}\wr c_{g_i}} &&
\prod_{i=1}^r\, \Sigma_{s_1}\wr H\ar[u]_{\Phi_{s_1,\dots,s_r}}
} \]
The right vertical morphism is the generalization of
the embedding \eqref{eq:wreath_sum} to multiple factors.
From here the double coset formula is straightforward:
\begin{align*}
 \res^G_K\circ\tr_H^G  \ &= \  \res^G_K\circ\Psi_{\bar g}^*\circ[m] \\ 
&= \ ((\Sigma_{s_1}\wr c_{g_1})\circ \Psi_{\bar k^1},\dots,(\Sigma_{s_r}\wr c_{g_r})\circ \Psi_{\bar k^r})^*\circ \Phi_{s_1,\dots,s_r}^*\circ[m] \\ 
_\eqref{eq:external_additivity} &= \  ((\Sigma_{s_1}\wr c_{g_1})\circ \Psi_{\bar k^1},\dots,(\Sigma_{s_r}\wr c_{g_r})\circ \Psi_{\bar k^r})^*\circ( [s_1]\oplus\dots \oplus[s_r] )  \\ 
 &= \ \sum_{i=1}^r\ \Psi_{\bar k^i}^*\circ (\Sigma_{s_i}\wr c_{g_i})^*\circ [s_i]  \\ 
 &= \ \sum_{i=1}^r\  \Psi_{\bar k^i}^*\circ [s_i]\circ (c_{g_i})^* \ 
= \ 
\sum_{i=1}^r\, \tr^K_{K\cap{^{g_i} H}} \circ (c_{g_i})^* \circ \res^H_{K^{g_i}\cap H} \ .
\end{align*}

(iii) If $\bar k=(k_1,\dots,k_m)$ is an $L$-basis of $K$,
then $\alpha(\bar k)= (\alpha(k_1),\dots,\alpha(k_m))$ is an $H$-basis of $G$.
With respect to these bases we have
\[ \Psi_{\alpha(\bar k)}\circ \alpha \ = \ (\Sigma_m\wr (\alpha|_L)) \circ \Psi_{\bar k} 
\ : \ K \ \to \ \Sigma_m\wr H \ .  \]
So 
\begin{align*}
   \alpha^*\circ \tr_H^G \ = \ 
\alpha^*\circ \Psi^*_{\alpha(\bar k)} \circ [m]
\ &= \  \Psi_{\bar k}^*\circ(\Sigma_m\wr (\alpha|_L))^* \circ  [m] \\ 
&= \   \Psi_{\bar k}^*\circ  [m]\circ (\alpha|_L)^* \ = \ \tr_L^K \circ (\alpha|_L)^*  \ .
\end{align*}

(iv) The subgroup $K$ has index $m$ in $\Sigma_m\wr G$ and a $K$-basis
of $\Sigma_m\wr G$ is given by the elements $\tau_j=((j, m);\,1,\dots,1)$
for $j=1,\dots,m$. 
In order to determine the monomorphism $\Psi_{\bar\tau}:\Sigma_m\wr G\to \Sigma_m\wr K$
associated to this $K$-basis,
we consider any element $(\sigma;\,g_1,\dots,g_m)$ of $\Sigma_m\wr G$.
The permutation $(\sigma(j), m)\cdot \sigma\cdot (j, m)\in\Sigma_m$ 
fixes $m$, so the element
\[  l_j \ = \ \tau_{\sigma(j)}\cdot (\sigma;\,g_1,\dots,g_m)\cdot \tau_j 
\ \in \ \Sigma_m\wr G \]
in fact belongs to the subgroup $K$.
Then 
\[ (\sigma;\,g_1,\dots,g_m)\cdot \tau_j \ = \ \tau_{\sigma(j)}\cdot l_j  \]
in the group $\Sigma_m\wr G$, by definition. This means that
\[  (\sigma;\,g_1,\dots,g_m)\cdot (\tau_1,\dots,\tau_m) \ = \ 
(\tau_1,\dots,\tau_m) \cdot (\sigma;\,l_1,\dots,l_m) \ ,\]
and hence
\[ \Psi_{\bar\tau}(\sigma;\,g_1,\dots,g_m) \ = \ (\sigma;\,l_1,\dots,l_m) \ .\]
Because
\[ (\Sigma_m\wr q)(\Psi_{\bar\tau}(\sigma;\,g_1,\dots,g_m))\ = \ 
(\Sigma_m\wr q)(\sigma;\,l_1,\dots,l_m) \ =\
(\sigma;\,g_1,\dots,g_m) \ , \]
we conclude that the composite $(\Sigma_m\wr q)\circ\Psi_{\bar\tau}$
is the identity of $\Sigma_m\wr G$. So
\begin{align*}
 \tr_K^{\Sigma_m\wr G}\circ q^* \ = \ 
\Psi_{\bar\tau}^*\circ [m]\circ q^* \ &= \ 
\Psi_{\bar\tau}^*\circ (\Sigma_m\wr q)^*\circ [m] \\ 
&= \ ((\Sigma_m\wr q)\circ\Psi_{\bar\tau})^*\circ [m] \ = \ [m]\ .\qedhere
\end{align*}
\end{proof}

Theorem \ref{thm:umon operation}~(iii) gives a description of the
set of natural operations, on the category of ultra-commutative monoids, 
from $\pi_0^G$ to $\pi_0^K$. The next proposition gives an alternative
description that also captures the monoid structure given by pointwise
addition of operations.

\begin{prop}
  Let $G$ and $K$ be compact Lie groups.
  The monoid $\Nat^{umon}(\pi_0^G,\pi_0^K)$ is a free abelian monoid
  generated by the operations $\tr_L^K\circ\alpha^*$
  where $(L,\alpha)$ runs over all $(K\times G)$-conjugacy classes of pairs 
  consisting of
  \begin{itemize}
  \item a subgroup $L\leq K$ of finite index, and
  \item a continuous group homomorphism $\alpha:L\to G$.
  \end{itemize}    
\end{prop}
\begin{proof}
This a straightforward algebraic consequence of the calculation
of the category $\Nat^{umon}$ given in Proposition \ref{prop:B iso}.
Every $K$-$G$-space with finitely many free $G$-orbits is the disjoint union
of transitive $K$-$G$-spaces with the same property. So $\mathbb A^+(G,K)$
is a free abelian monoid with basis the isomorphism classes of the
transitive $K$-$G$-spaces. A transitive $K$-$G$-space 
with finitely many free $G$-orbits is isomorphic to one of the form
\[ K\times_{\alpha} G \ = \ (K\times G)/ (k l,g)\sim (k,\alpha(l) g) \]
for a pair $(L,\alpha:L\to G)$, with $L$ of finite index in $K$.
Moreover, $K\times_{\alpha} G$ is isomorphic to
$K\times_{\alpha'} G$ if and only if $(L,\alpha)$ is conjugate to
$(L',\alpha')$ by an element of $K\times G$.
So $\mathbb A^+(G,K)$ is freely generated by the classes of the
$K$-$G$-spaces $K\times_{\alpha} G$, where $(L,\alpha)$ 
runs through the $(K\times G)$-conjugacy classes of the relevant pairs.

The claim then follows from the verification that the isomorphism
of categories $B:\text{Nat}^{umon}\to \mA^+$ 
established in Proposition \ref{prop:B iso}
takes the operation $\tr_L^K\circ \alpha^*$ to the class of $K\times_{\alpha} G$. 
Indeed, if $\bar k=(k_1,\dots,k_m)$ is an $L$-basis of $K$ and
$\Psi_{\bar k}:K\to \Sigma_m\wr L$ the classifying homomorphism, then 
\[ B(\tr_L^K\circ\alpha^*) \ = \ B(\Psi^*_{\bar k}\circ[m]\circ\alpha^*)
\ = \  B(\Psi^*_{\bar k}\circ(\Sigma_m\wr\alpha)^*\circ[m])\ .\]
On the other hand, the map
\[ ((\Sigma_m\wr\alpha)\circ\Psi_{\bar k})^*(\{1,\dots,m\}\times G)_G \ \to \ 
K\times_{\alpha} G \ , \quad (i,\gamma)\ \longmapsto \ [k_i, \gamma]\]
is an isomorphism of $K$-$G$-bispaces, 
so $B(\tr_K^K\circ\alpha^*)=[K\times_{\alpha} G]$ in $\mathbb A^+(G,K)$.  
\end{proof}

\section{Examples}\label{sec:umon examples}

In this section we discuss various examples of ultra-commutative monoids,
mostly of a geometric nature, and several geometrically defined morphisms between them.
We start with the ultra-commutative 
monoids $\bO$ and $\bSO$ (Example \ref{eg:orthogonal group monoid space}), 
$\bU$ and $\bSU$ (Example \ref{eg:unitary group monoid space}), 
$\bSp$ (Example \ref{eg:symplectic group monoid space}) 
and $\bSpin$ and $\bSpin^c$ (Example \ref{eg:SPin monoid space}),
all made from the corresponding families of classical Lie groups.
All of these are examples of 
`symmetric monoid valued orthogonal spaces' 
in the sense of Definition \ref{def:symmetric monoid valued},
a more general source of examples of ultra-commutative monoids.
The additive Grassmannian $\bGr$ (Example \ref{eg:Gr additive}),
the oriented variant $\bGr^{\ort}$ (Example \ref{eg:Gr oriented})
and the complex and quaternionic analogs $\bGr^\mC$ and $\bGr^\mH$ 
(Example \ref{eg:Gr additive complex})
consist -- as the names suggest -- of Grassmannians with monoid
structure arising from direct sum of subspaces.
The multiplicative Grassmannian $\bGr_\tensor$ (Example \ref{eg:Gr multiplicative real})
is globally equivalent as an orthogonal space to $\bGr$, but the 
monoid structure arises from tensor product of subspaces;
the global projective space $\bP$ is the ultra-commutative submonoid
of $\bGr_\tensor$ consisting of lines (1-dimensional subspaces).
The global projective space $\bP$ is a multiplicative model of
a global classifying space for the cyclic group of order~2, and
Example \ref{con:multiplicative B G} describes multiplicative models
of global classifying spaces for all {\em abelian} compact Lie groups.
Example \ref{eg:unordered frames} introduces
the ultra-commutative monoid $\bF$ of unordered frames, with monoid structure
arising from disjoint union.
The ultra-commutative `multiplicative monoid of the sphere spectrum'
$\Omega^\bullet \mS$ is introduced in Example \ref{eg:Omega^bullet S};
this a special case of the multiplicative monoid of an ultra-commutative
ring spectrum, and we return to the more general construction 
in Example \ref{eg:suspension spectrum of orthogonal monoid space}.

\begin{construction}[Orthogonal monoid spaces from monoid valued
orthogonal spaces]\label{con:umon from group valued}
Our first series of examples involves orthogonal spaces made from
the infinite families of classical Lie groups, 
namely the orthogonal, unitary and symplectic groups,
the special orthogonal and unitary groups, and the pin, pin$^c$, spin and spin$^c$ groups. 
These orthogonal spaces have the special feature that they are {\em group valued};
we will now explain that a group valued (or even just a monoid valued)
orthogonal space automatically leads
to an orthogonal monoid space. In the cases of the orthogonal, special orthogonal,
unitary, special unitary, spin, spin$^c$ and symplectic groups, 
these multiplications are symmetric, so those examples yield ultra-commutative monoids.

\begin{defn}\label{def:monoid valued ospace} 
A {\em monoid valued orthogonal space}\index{subject}{orthogonal space!monoid valued} 
is a monoid object in the category of orthogonal spaces.
\end{defn}

Since orthogonal spaces are an enriched functor category, monoid valued orthogonal 
spaces are the same thing as continuous functors from the category $\bL$
to the category of topological monoids 
and continuous monoid homomorphisms
(i.e., monoid objects in the category $\bT$ of compactly generated spaces). 

Now we let $M$ be a monoid valued orthogonal space. 
In \eqref{eq:define_box2times}
we introduced a lax symmetric monoidal natural transformation
\[  \rho_{X,Y} \ : \ X\boxtimes Y\ \to \ X\times Y    \]
from the box product to the cartesian product of orthogonal spaces.
For $X=Y=M$ we can form the composite 
\begin{equation}  \label{eq:monoid_times2box}
  M\boxtimes M\ \xra{\ \rho_{M,M}\ } \ M\times M    \ \xra{\text{multiplication}} \ M   
\end{equation}
with the objectwise multiplication of $M$.
Since $\rho_{X,Y}$ is lax monoidal, this composite makes $M$ into
an orthogonal monoid space with unit the multiplicative unit $1\in M(0)$.
The morphism \eqref{eq:monoid_times2box} corresponds, 
via the universal property of $\boxtimes$, to the bimorphism
whose $(V,W)$-component is the composite
\[ M(V)\times M(W)\ \xra{M(i_V)\times M(i_W)}\
 M(V\oplus W)\times M(V\oplus W)\ \xra{\text{multiply}}\
 M(V\oplus W)\ , \]
where $i_V:V\to V\oplus W$ and $i_W:W\to V\oplus W$ are the direct summand embed\-dings. 
If $f:M\to M'$ is a morphism of monoid valued orthogonal spaces
(i.e., a morphism of orthogonal spaces that is objectwise a monoid homomorphism),
then $f$ is also a homomorphism of orthogonal monoid spaces
with respect to the multiplications \eqref{eq:monoid_times2box}.
\end{construction}

If $M$ is a commutative monoid valued orthogonal space, then the associated
$\boxtimes$-multiplication is also commutative, simply because the transformation
$\rho_{X,Y}$ is {\em symmetric} monoidal. However, there is a more general
condition that provides $M$ with the structure of an ultra-commutative monoid.

\begin{defn}\label{def:symmetric monoid valued} 
A monoid valued orthogonal space $M$ 
is~{\em symmetric}\index{subject}{orthogonal space!monoid valued!symmetric} 
if for all inner product spaces $V$ and $W$ the images of the two homomorphisms
\[ M(i_V)\ : \ M(V)\ \to \  M(V\oplus W)\text{\quad and\quad}
M(i_W)\ : \ M(W)\ \to \  M(V\oplus W)\]
commute. 
\end{defn}

We emphasize that the objectwise multiplications in a symmetric 
monoid valued orthogonal space need not be commutative -- 
we'll discuss many interesting examples of this kind below.
The proof of the following proposition is straightforward from the definitions,
and we omit it.

\begin{prop}
Let $M$ be a symmetric monoid valued orthogonal space.
Then the multiplication \eqref{eq:monoid_times2box}
makes $M$ into an ultra-commutative monoid.
\end{prop}

\begin{eg}[Orthogonal group ultra-commutative monoid]\label{eg:orthogonal group monoid space}
\index{subject}{orthogonal group ultra-commutative monoid}
We denote by $\bO$\index{symbol}{$\bO$ - {ultra-commutative monoid of orthogonal groups}} 
the orthogonal space that sends an inner product space $V$
to its orthogonal group $O(V)$.\index{subject}{orthogonal group}
A linear isometric embedding $\varphi:V\to W$ induces a continuous group
homomorphism $\bO(\varphi):\bO(V)\to\bO(W)$ by conjugation
and the identity on the orthogonal complement of the image of $\varphi$.
Construction \ref{con:umon from group valued} then gives $\bO$ the
structure of an orthogonal monoid space.
The $(V,W)$-component of the bimorphism 
\[ \mu_{V,W} \ : \ \bO(V)\times\bO(W) \to \ \bO(V\oplus W) \]
is direct sum of orthogonal transformations.
The unit element is the identity of the trivial vector space, 
the only element of $\bO(0)$.
So $\bO$ is a symmetric group valued orthogonal space,
and hence it becomes an ultra-commutative monoid.

If $G$ is a compact Lie group and $V$ a $G$-representation, 
then the $G$-action on the group $\bO(V)$ is by conjugation,
so the fixed points $\bO(V)^G$ are the group of $G$-equivariant orthogonal
automorphisms of $V$.
Moreover, $\bO(\Uc_G)$ is the orthogonal group of $\Uc_G$, 
i.e., $\mR$-linear isometries of $\Uc_G$
(not necessarily $G$-equivariant) that are the identity on the orthogonal 
complement of some finite-dimensional subspace;
the $G$-action is again by conjugation.
Any $G$-equivariant isometry preserves the decomposition of $\Uc_G$
into isotypical summands, and the restriction to
almost all of these isotypical summands must be the identity.
The $G$-fixed subgroup is thus given by
\[ \bO(\Uc_G)^G\ = \ O^G(\Uc_G) \ = \ {\prod}'_{[\lambda]}\  O^G(\Uc_\lambda)\ ,\]
where the weak product is indexed by the isomorphism classes of irreducible
orthogonal $G$-representations $\lambda$, and $\Uc_\lambda$ 
is the $\lambda$-isotypical summand. 
If the compact Lie group $G$ is finite, 
then there are only finitely many isomorphism classes 
of irreducible $G$-representations,
so in that case the weak product coincides with the product.

Irreducible orthogonal representations come in three different
flavors, and the group  $O^G(\Uc_\lambda)$ has one of three different forms.
If $\lambda$ is an irreducible orthogonal $G$-representation, then
the endomorphism ring $\Hom^G_{\mR}(\lambda,\lambda)$ is a finite-dimensional
skew field extension of $\mR$, so it is isomorphic to either $\mR$, $\mC$ 
or $\mH$; the representation $\lambda$ is accordingly called `real', `complex'
respectively `quaternionic'. We have
\[ O^G(\Uc_\lambda) \  \iso \ O^G(\lambda^\infty) \  \iso \
\begin{cases}
 \ O & \text{ if $\lambda$ is real,}\\
 \ U & \text{ if $\lambda$ is complex, and}\\
  Sp & \text{ if $\lambda$ is quaternionic.}
\end{cases}  \]
So we conclude that the $G$-fixed point space $\bO(\Uc_G)^G$
is a weak product of copies of the infinite orthogonal, unitary
and symplectic groups, indexed by the different types 
of irreducible orthogonal representations of $G$.
Since the infinite unitary and symplectic groups are connected, 
only the `real' factors contribute to $\pi_0(\bO(\Uc_G)^G)= \pi_0^G(\bO)$,
which is a weak product of copies of $\pi_0(O)=\mZ/2$ indexed
by the irreducible $G$-representations of real type.

There is a straightforward `special orthogonal' analog:
the property of having determinant~1 is preserved under conjugation by
linear isometric embeddings and under direct sum of linear isometries,
so the spaces $SO(V)$\index{subject}{special orthogonal group}\index{subject}{orthogonal group!special|see{special orthogonal group}}
form an ultra-commutative 
submonoid $\bSO$\index{symbol}{$\bSO$ - {ultra-commutative monoid of special orthogonal groups}} of $\bO$.
Here, $\bSO(\Uc_G)$ is the group of $\mR$-linear isometries of $\Uc_G$
(not necessarily $G$-equivariant) that have determinant~1 on some finite
dimensional subspace $V$ of $\Uc_G$ and are the identity on the orthogonal 
complement of $V$.
\end{eg}

\begin{eg}[Unitary group ultra-commutative monoid]\label{eg:unitary group monoid space}
\index{subject}{unitary group ultra-commutative monoid}
There is a straightforward unitary analog $\bU$\index{symbol}{$\bU$ - {ultra-commutative monoid of unitary groups}} 
of $\bO$, defined as follows. 
We recall that
\[ V_\mC \ = \ \mC\tensor_\mR V \]
denotes the complexification of an inner product space $V$.
The euclidean inner product $\td{-,-}$ on $V$ induces a hermitian
inner product $(-,-)$ on $V_\mC$, defined
as the unique sesquilinear form that satisfies
\[ (1\tensor v, 1\tensor w) \ = \ \td{v,w} \]
for all $v, w\in V$.
We now define an orthogonal space $\bU$ by
\[ \bU(V)\ = \ U(V_\mC) \ ,\]
the unitary group of the complexification of $V$.\index{subject}{unitary group}
The complexification of every $\mR$-linear isometric embedding $\varphi:V\to W$ 
preserves the hermitian inner products, 
so we can define a continuous group homomorphism
\[ \bU(\varphi)\ :\ \bU(V)\ \to\ \bU(W)  \]
by conjugation with $\varphi_\mC:V_\mC\to W_\mC$
and the identity on the orthogonal complement of the image of $\varphi_\mC$.
The $\boxtimes$-multiplication on $\bU$
produced by Construction \ref{con:umon from group valued}
is by direct sum of unitary transformations;
this multiplication is symmetric, and hence ultra-commutative.
There is a straightforward `special unitary' ultra-commutative
submonoid $\bSU$ of $\bU$; the value $\bSU(V)$ is the group of
unitary automorphisms of $V_\mC$ of determinant~1.\index{symbol}{$\bSU$ - {ultra-commutative monoid of special unitary groups}}\index{subject}{special unitary group}\index{subject}{unitary group!special|see{special unitary group}}

If $G$ is a compact Lie group, then the identification of
the $G$-fixed points of $\bU$ also works much like the orthogonal analog.
The outcome is an isomorphism between $\bU(\Uc_G)^G$ and $U^G(\Uc_G^\mC)$.
Here $\Uc_G^\mC=\mC\tensor_\mR \Uc_G$ is the complexified complete universe
for $G$, which happens to be a `complete complex $G$-universe' in
the sense that every finite-dimensional complex $G$-representation
embeds into it.
The complete complex $G$-universe breaks up into (unitary)
isotypical summands $\Uc_\lambda^\mC$,
indexed by the isomorphism classes of irreducible
unitary $G$-representations $\lambda$, and the group $U^G(\Uc_G^\mC)$
breaks up accordingly as a weak product.
In contrast to the orthogonal situation above, there is only one `type'
of irreducible unitary representation, and the group $U^G(\Uc^\mC_\lambda)$ 
is always isomorphic to the infinite unitary group $U$, independent of $\lambda$.
So in the unitary context, we get a decomposition
\[ \bU(\Uc_G)^G\ = \ U^G(\Uc_G^\mC) \ = \ {\prod}'_{[\lambda]}\ U^G(\Uc_\lambda^\mC)\
\ \iso \ {\prod}'_{[\lambda]}\  U  \ .\]
This weak product is indexed by the isomorphism classes of irreducible
unitary $G$-representations. Since the unitary group $U$ is connected,
the set $\pi_0^G(\bU)=\pi_0(\bU(\Uc_G)^G)$ has only one element,
and so $\bU$ is globally connected.

The orthogonal monoid space $\bU$ comes with an involution
\[ \psi \ : \ \bU \ \to \ \bU \]
by {\em complex conjugation}\index{subject}{complex conjugation!on $\bU$}
that is an automorphism of ultra-commutative monoids.
The value of $\psi$ at $V$ is the map
\[ \psi(V) \ : \ U(V_\mC) \ \to \ U(V_\mC) \ , \quad 
A \ \longmapsto \ \psi_V\circ A\circ\psi_V \ , \]
where 
\[ \psi_V\ : \ V_\mC \ \to \ V_\mC \ , \quad \lambda\tensor v\ \longmapsto \ 
\bar\lambda\tensor v \]
is the canonical $\mC$-semilinear conjugation map on $V_\mC$.

The {\em complexification morphism}\index{subject}{complexification morphism!from $\bO$ to $\bU$}
is the homomorphism of ultra-commutative monoids given by
\[ c \ : \ \bO \ \to \ \bU  \ , \quad
 c(V)\ : \ O(V)\ \to \ U(V_\mC) \ , \quad \varphi\ \longmapsto \ \varphi_\mC\ .\]
Complexification is an isomorphism onto the $\psi$-invariant ultra-commutative submonoid 
of $\bU$, and it takes $\bSO$ to $\bSU$.

Every hermitian inner product space $W$ has an underlying 
$\mR$-vector space equipped with the euclidean inner product
defined by
\[ \td{v, w}\ = \ \text{Re} ( v,w ) \ ,\]
the real part of the given hermitian inner product. 
Every $\mC$-linear isometric embedding is in particular
an $\mR$-linear isometric embedding of underlying euclidean vector spaces.
So the unitary group $U(W)$ is a subgroup of the special orthogonal group
of the underlying euclidean vector space of $W$. We can thus define 
the {\em realification morphism}\index{subject}{realification morphism!from $\bU$ to $\bO$}
\begin{equation}  \label{eq:realification_bu2shbO}
 r \ : \ \bU \ \to \ \sh^\tensor_\mC ( \bSO )   
\end{equation}
at $V$ as the inclusion
$r(V): U(V_\mC)\to O(\mC\tensor V)$.
Here $\sh^\tensor_\mC$ denotes the multiplicative shift by $\mC$
as defined in Example \ref{eg:Additive and multiplicative shift}.
\end{eg}

\begin{eg}[Symplectic group ultra-commutative monoid]\label{eg:symplectic group monoid space}
\index{subject}{symplectic group ultra-commutative monoid}\index{subject}{symplectic group}
There is also a quaternionic analog of $\bO$ and $\bU$,
the ultra-commutative monoid $\bSp$ made from symplectic groups.\index{symbol}{$\bSp$ - {ultra-commutative monoid of symplectic groups}} 
Given an $\mR$-inner product space $V$, we denote by
\[ V_\mH \ = \ \mH\tensor_\mR V \]
the extension of scalars to the skew field $\mH$ of quaternions.
The extension comes with a $\mH$-sesquilinear form 
\[ [-,-]\ : \ V_\mH\times V_\mH \ \to \ \mH \text{\quad characterized by\quad}
 [1\tensor v, 1\tensor w]\ = \ \td{v,w}  \]
for all $v,w\in V$. The symplectic group
\[ \bSp(V)\ = \ Sp(V_\mH) \]
is the compact Lie group of $\mH$-linear automorphisms $A:V_\mH\to V_\mH$
that leave the form invariant, i.e., such that
\[  [A x, A y] \ = \ [x,y]\]
for all $x,y\in V_\mH$.
The $\mH$-linear extension $\varphi_\mH=\mH\tensor_\mR \varphi:V_\mH\to W_\mH$ 
of an $\mR$-linear isometric embedding $\varphi:V\to W$ 
preserves the new inner products, so we can define a continuous group homomorphism
\[ \bSp(\varphi)\ :\ \bSp(V)\ \to\ \bSp(W)  \]
by conjugation with $\varphi_\mH$
and the identity on the orthogonal complement of the image of $\varphi_\mH$.
As for $\bO$ and $\bU$, the $\boxtimes$-multiplication on $\bSp$
produced by Construction \ref{con:umon from group valued}
is by direct sum of symplectic automorphisms;
so $\bSp$ is symmetric, hence ultra-commutative.

If $G$ is a compact Lie group, then the identification of
the $G$-fixed points of $\bSp$ also works much like the orthogonal case.
Quaternionic representations decompose canonically into isotypical
summands, and this results in a product decomposition for the 
$G$-fixed subgroup
\[ \bSp(\Uc_G)^G\ = \ ( Sp(\Uc_G^\mH))^G \ = \ 
{\prod}'_{[\lambda]}\, (Sp(\Uc_\lambda^\mH))^G\ ,\]
where the weak product is indexed by the isomorphism classes of irreducible
quaternionic $G$-representations $\lambda$, and $\Uc_\lambda^\mH$ 
is the $\lambda$-isotypical summand in $\Uc_G^\mH=\mH\tensor_\mR\Uc_G$. 
As in the real case, irreducible quaternionic representations $\lambda$ 
come in three different flavors, depending on whether 
$\Hom^G_{\mH}(\lambda,\lambda)$ -- again a finite-dimensional
skew field extension of $\mR$ -- is isomorphic to $\mR$, $\mC$ or $\mH$.
As for $\bO$, the $G$-fixed point space $\bSp(\Uc_G)^G$
is a weak product of copies of the infinite orthogonal, unitary
and symplectic groups, indexed by the different types 
of irreducible quaternionic representations of $G$.
\end{eg}

\begin{eg}[Pin and Spin group orthogonal spaces]\label{eg:SPin monoid space}
Given a real inner product space $V$ we denote by $\Cl(V)$\index{subject}{Clifford algebra|(}\index{symbol}{$\Cl(V)$ - {Clifford algebra}}
the Clifford algebra of the negative definite quadratic form on $V$, 
i.e., the quotient of
the $\mR$-tensor algebra of $V$ by the ideal generated by
$v\tensor v + |v|^2\cdot 1$ for all $v\in V$.
The Clifford algebra is $\mZ/2$-graded
with the even (respectively odd) part generated by an even
(respectively odd) number of vectors from $V$.   
The composite
\[ V \ \xra{\text{linear summand}}\ T V \ \xra{\text{proj}} \ \Cl(V) \]
is $\mR$-linear and injective, and we denote it by $v\mapsto[v]$.

We recall that orthogonal vectors of $V$ anti-commute in the Clifford algebra:
given $v, \bar v\in V$ with $\td{v,\bar v}=0$, then
\begin{align*}
[v][\bar v] + [\bar v][ v] \ = \ 
(|v|^2\cdot 1 + [v]^2)  +[v] [\bar v] &+ [\bar v] [v]  + ( |\bar v^2|\cdot 1 + [\bar v]^2)\\
&= \ |v+\bar v|^2\cdot 1  + [v+\bar v]^2 \  
 =\ 0 \ .
\end{align*}
In $\Cl(V)$ every unit vector $v\in S(V)$ satisfies $[v]^2=-1$,
so all unit vectors of $V$ become units in $\Cl(V)$.
The {\em pin group}\index{subject}{pin group} of $V$ is the subgroup
\[ \Pin(V) \ \subset \ \Cl(V)^\times \]
generated inside the multiplicative group of $\Cl(V)$
by $-1$ and all unit vectors of $V$.
A linear isometric embedding $\varphi:V\to W$ induces a morphism
of $\mZ/2$-graded $\mR$-algebras $\Cl(\varphi):\Cl(V)\to \Cl(W)$
that restricts to $\varphi$ on $V$.
So $\Cl(\varphi)$ restricts to a continuous homomorphism
\[ \Pin(\varphi)\ : \ \Cl(\varphi)|_{\Pin(V)} \ : \ \Pin(V)\ \to \ \Pin(W)\]
between the pin groups. The map $\Pin(\varphi)$ depends continuously
on $\varphi$ and satisfies $\Pin(\psi)\circ\Pin(\varphi)=\Pin(\psi\circ\varphi)$,
so we have defined a group-valued orthogonal space $\bPin$. 
Construction \ref{con:umon from group valued} then gives $\bPin$ the
structure of an orthogonal monoid space.
\index{subject}{pin group orthogonal monoid space}\index{symbol}{$\bPin$ - {orthogonal monoid space of pin groups}} 

Since the group $\Pin(V)$ is generated by homogeneous elements of
the Clifford algebra, all of its elements are homogeneous.
The $\mZ/2$-grading of $\Cl(V)$ provides a continuous homomorphism
$ \Pin(V) \to \mZ/2$ whose kernel 
\[ \Spin(V)\ = \  \Cl(V)_{\ev} \cap \Pin(V)\]
is the {\em spin group}\index{subject}{spin group} of $V$.
The map $\Pin(\varphi)$ induced by a linear isometric
embedding $\varphi:V\to W$ is homogeneous, so it restricts to a homomorphism
\[ \Spin(\varphi)\ = \ \Pin(\varphi)|_{\Spin(V)} \ : \ \Spin(V)\ \to \ \Spin(W)\]
between the spin groups. 
The spin groups thus form an orthogonal monoid subspace $\bSpin$
of $\bPin$.\index{subject}{spin group ultra-commutative monoid}\index{symbol}{$\bSpin$ - {ultra-commutative monoid of spin groups}}   
We claim that $\bSpin$ is symmetric, i.e.,
for all inner product spaces $V$ and $W$ the images of
the two group homomorphisms
\[ \Spin(V)\ \xra{\Spin(i_V)}\ \Spin(V\oplus W) \ \xla{\Spin(i_W)}\ \Spin(W) \]
commute.
Indeed, $\Spin(V)$ is generated by $-1$ and the elements $[v][v']$
for $v,v'\in S(V)$, and similarly for $W$. 
The elements $-1$ map to the central element $-1$ in $\Spin(V\oplus W)$,
and
\[ [v,0]\cdot[v',0]\cdot[0,w]\cdot [0,w']\ = \ 
[0,w]\cdot [0,w']\cdot [v,0]\cdot[v',0]  \]
because $[v,0]$ and $[v',0]$ anti-commute with $[0,w]$ and $[0,w']$.
So $\bSpin$ is an ultra-commutative monoid with respect to
the $\boxtimes$-multiplication of Construction \ref{con:umon from group valued}.

\medskip

\Danger In contrast to $\bSpin$, the group valued orthogonal space $\bPin$ is {\em not}
symmetric; equivalently, the continuous map
\[ \mu_{V,W}\ : \ \Pin(V)\times P(W) \to\ \Pin(V\oplus W) \]
is {\em not} a group homomorphism.
The issue is that for $v\in S(V)$ and $w\in S(W)$
the elements $[v,0]$ and $[0,w]$ {\em anti-commute} in $\Pin(V\oplus W)$.

\medskip

Now we turn to the groups pin$^c$ and spin$^c$,
the complex variations on the pin and spin groups.
These arise from the complexified  Clifford algebra
$ \mC\tensor_\mR \Cl(V)$, where $V$ is again a euclidean inner product space.
The complexified Clifford algebra is again $\mZ/2$-graded, and functorial 
for $\mR$-linear isometric embeddings.
The {\em pin$^c$ group}\index{subject}{pin$^c$ group} of $V$ is the subgroup
\[ \Pin^c(V) \ \subset \ (\mC\tensor_\mR\Cl(V))^\times \]
generated inside the multiplicative group
by the unit scalars $\lambda\tensor 1$ for all $\lambda\in U(1)$ 
and the elements $1\tensor [v]$ for all unit vectors $v\in S(V)$.
The pin$^c$-groups form a group valued orthogonal space $\bPin^c$,
and hence an orthogonal monoid space, in much the same way
as for the pin groups above. As for $\bPin$, the monoid valued orthogonal space $\bPin^c$
is {\em not} symmetric, 
so the associated $\boxtimes$-multiplication of $\bPin^c$ is {\em not} commutative.
\index{subject}{pin$^c$ group orthogonal monoid space}\index{symbol}{$\bPin^c$ - {orthogonal monoid space of pin$^c$ groups}} 

Since the group $\Pin^c(V)$ is generated by homogeneous elements of
the complexified Clifford algebra, all of its elements are homogeneous.
So the $\mZ/2$-grading of $\mC\tensor_\mR\Cl(V)$ provides a continuous homomorphism
$ \Pin^c(V)\to \mZ/2$ whose kernel 
\[ \Spin^c(V)\ = \  (\mC\tensor_\mR \Cl(V))_{\ev} \cap \Pin^c(V)\]
is the {\em spin$^c$ group}\index{subject}{spin$^c$ group} of $V$.
As $V$ varies, the spin$^c$ groups from a group valued orthogonal subspace $\bSpin^c$
of $\bPin^c$.\index{subject}{spin$^c$ group ultra-commutative monoid}\index{symbol}{$\bSpin^c$ - {ultra-commutative monoid of spin$^c$ groups}}   
As for $\bSpin$, the images of the homomorphisms
$\Spin^c(i_V)$ and $\Spin^c(i_W)$ commute, so
$\bSpin^c$ is even an ultra-commutative monoid. 

We embed the real Clifford algebra $\Cl(V)$
as an $\mR$-subalgebra of its complexification by
\[ \iota(V) = 1\tensor - \ : \  \Cl(V)\ \to \ \mC\tensor_\mR\Cl(V) \ .\quad \]
This homomorphism restricts to embeddings
\[ \iota(V)\ : \ \Pin(V)\ \to \ \Pin^c(V)\text{\quad and\quad}
\iota(V)\ : \ \Spin(V)\ \to \ \Spin^c(V)\ ; \]
as $V$ varies, these maps
form morphisms of group valued orthogonal spaces $\iota:\bPin\to\bPin^c$
and $\iota:\bSpin\to\bSpin^c$, and hence morphisms of
orthogonal monoid spaces.
These morphisms extend to isomorphisms 
\[ \bPin\times_{\{\pm 1\}} U(1)\ \iso \ \bPin^c \text{\qquad and\qquad}
 \bSpin\times_{\{\pm 1\}} U(1)\ \iso \ \bSpin^c   \]
of orthogonal monoid spaces.

We let $\alpha:\Cl(V)\to\Cl(V)$ denote the unique $\mR$-algebra automorphism
of the Clifford algebra such that $\alpha[v]=-[v]$ for all $v\in V$.
The map $\alpha$ is the grading involution, i.e., it is the
identity on the even part and $-1$ on the odd part of the Clifford algebra.
We also denote by $\alpha$ the automorphism  
of the complexification $\mC\tensor_\mR\Cl(V)$ obtained by complexifying the real version.
For every element $x\in\Pin^c(V)$ the twisted conjugation map
\[ c_x \ : \ \mC\tensor_\mR\Cl(V)\ \to\ \mC\tensor_\mR \Cl(V) \ , \quad 
c_x(y)\ = \ \alpha(x)y x^{-1} \]
is an automorphism of $\mZ/2$-graded $\mC$-algebras.
We let $v\in S(V)$ be a unit vector. Then twisted conjugation by $[v]$ 
takes $[v]$ to $-[v]$, and it fixes the elements $[w]$ for all
$w\in V$ that are orthogonal to $v$.
So twisted conjugation by $[v]$ `is' reflection
in the hyperplane orthogonal to $v$, hence a linear isometry of $V$
of determinant $-1$.
The complex scalars are central in $\mC\tensor_\mR\Cl(V)$, 
so conjugation by elements $\lambda\tensor 1$ with $\lambda\in U(1)$ is the identity.
So for all elements $x\in\Pin^c(V)$, the conjugation map $c_x$ 
restricts to a linear isometry on $V$, in the sense that 
there is a unique $\ad(x)\in O(V)$ satisfying
\[ \alpha(x)[v] x^{-1} \ =\ [\ad(x)(v)] \]
for all $v\in V$. We thus obtain a continuous group homomorphism
\[ \ad(V) \ : \ \Pin^c(V)\ \to \ O(V)\ , \quad x \ \longmapsto \ \ad(x) \ ,\]
the {\em twisted adjoint representation}.\index{subject}{twisted adjoint representation!of the pin group}
The kernel of the twisted adjoint representation
is the subgroup $U(1)\cdot 1$, compare \cite[Thm.\,3.17]{atiyah-bott-shapiro}.
This homomorphism takes the spin$^c$ group to the special orthogonal group,
and restricts to the {\em adjoint representation}\index{subject}{adjoint representation!of the spin$^c$ group}
\begin{equation} \label{eq:adjoint_rep_spin^c}
 \ad(V) \ : \ \Spin^c(V)\ \to \ SO(V)\ .   
\end{equation}
These homomorphisms form morphisms of group valued orthogonal spaces
\[ \ad \ : \ \bPin^c\ \to \ \bO   \text{\qquad and\qquad} 
   \ad \ : \ \bSpin^c\ \to \ \bSO . \]
Via Construction \ref{con:umon from group valued} we can interpret these as
morphisms of orthogonal monoid spaces.
\index{subject}{Clifford algebra|)} 

There is yet another interesting morphism of group valued orthogonal spaces 
\[ l \ : \ \bU \ \to \ \sh^\tensor_\mC(\bSpin^c)  \]
that lifts the forgetful realification morphism \eqref{eq:realification_bu2shbO}
through
\[ \sh^\tensor_\mC(\ad)\ :\ \sh^\tensor(\bSpin^c)\ \to\ \sh^\tensor_\mC(\bSO)\ . \]
The definition of $l(V_\mC)$ is a coordinate free description
of the homomorphism $U(n)\to\Spin^c(2n)$ that is defined 
for example in \cite[\S 3, p.10]{atiyah-bott-shapiro}.
Since we don't need this, we won't go into any details.
Since both $\bU$ and $\bSpin^c$ are {\em symmetric} group valued
orthogonal spaces, $l:\bU \to\sh^\tensor_\mC(\bSpin^c)$
is also a homomorphism of ultra-commutative monoids.
The morphism $l$ takes the special unitary group $S U(V_\mC)$ 
to the group $\Spin(\mC\tensor V)$, so it restricts to a morphism of ultra-commutative 
monoids $l:\bSU\to\sh^\tensor_\mC\bSpin$.

Most of the examples discussed so far 
can be summarized in the commutative diagram of orthogonal monoid spaces:
\[ \xymatrix{ 
\bSU \ar[r]^-{\text{incl}} \ar@{-->}[d]_l & \bU \ar@{-->}[d]_l \ar@{-->}[dr]^-{r}& \\
\bSpin \ar[d]_{\text{incl}} \ar[r]^-{\iota} & \bSpin^c \ar[d]^{\text{incl}}\ar[r]_-\ad 
 & \bSO \ar[d]^{\text{incl}}\\ 
 \bPin \ar[r]_-{\iota} &  \bPin^c \ar[r]_-\ad & \bO} \]
The two dotted arrows mean that the actual morphism goes to a multiplicative
shift of the target.
With the exception of $\bPin$ and $\bPin^c$, all the orthogonal monoid spaces
are ultra-commutative.
\end{eg}

\begin{eg}[Additive Grassmannian]\label{eg:Gr additive}
We define an ultra-commutative monoid $\bGr$,
the~{\em additive Grassmannian}.\index{subject}{additive Grassmannian}\index{subject}{Grassmannian!additive|see{additive Grassmannian}}
The value of $\bGr$\index{symbol}{$\bGr$ - {additive Grassmannian}}
at an inner product space $V$ is
\[ \bGr(V)\ = \ {\coprod}_{m\geq 0} \, Gr_m(V) \ ,\]
the disjoint union of all Grassmannians in $V$.
The structure map induced by a linear isometric embedding $\varphi:V\to W$
is given by  $\bGr(\varphi)(L) = \varphi(L)$.
A commutative multiplication on $\bGr$ is given by direct sum:
\begin{align*}
\mu_{V,W}\ : \  \bGr( V ) \times \bGr ( W ) \ \to \ \bGr( V\oplus W) \ ,\quad
(L,L')\ \longmapsto \ L\oplus L'
\end{align*}
The unit is the only point $\{0\}$ in $\bGr(0)$.
The orthogonal space $\bGr$ is naturally $\mN$-graded, with degree $m$
part given by $\bGr^{[m]}(V)=Gr_m(V)$. The multiplication is graded
in that it sends $\bGr^{[m]}( V ) \times \bGr^{[n]} ( W )$ to $\bGr^{[m+n]}( V\oplus W)$. 

As an orthogonal space, $\bGr$ is the disjoint union of global
classifying spaces of the orthogonal groups.\index{subject}{global classifying space!of $O(n)$} 
Indeed, the homeomorphisms
\[  \bL(\mR^m,V)/O(m) \ \iso \  \bGr^{[m]}(V) \ , \quad 
\varphi\cdot O(m)\ \longmapsto \ \varphi(\mR^m) \]
show that the summand $\bGr^{[m]}$ is isomorphic to the 
semifree orthogonal space $\bL_{O(m),\mR^m}$.
Since the tautological action of $O(m)$ on $\mR^m$ is faithful, 
this is a global classifying space for the orthogonal group.
So as orthogonal spaces,
\[  \bGr \ = \ {\coprod}_{m\geq 0} \, B_{\gl} O(m)\ . \]
Proposition \ref{prop:fix of global classifying}~(ii) identifies the
equivariant homotopy set $\pi_0^G ( B_{\gl} O(m))$
with the set of conjugacy classes of continuous homomorphisms from $G$ to $O(m)$;
by restricting the tautological $O(m)$-representation on $\mR^m$, this
set bijects with the set of isomorphism classes of 
$m$-dimensional $G$-representations.
In the union over all $m\geq 0$, this becomes a bijection between
$\pi_0^G(\bGr)$ and $\bRO^+(G)$, the set of isomorphism classes of 
orthogonal $G$-representations.

We make the bijection more explicit, showing at the same time that it is an 
isomorphism of monoids.
We let $V$ be a finite-dimensional orthogonal $G$-representation. 
The $G$-fixed points of $\bGr(V)$ are the $G$-invariant subspaces~of $V$, 
i.e., the $G$-subrepresentations.
We define a map
\[  \bGr(V)^G \ = \ {\coprod}_{m\geq 0} \,   \left( Gr_m(V) \right)^G 
\  \to \ \bRO^+(G) \]
from this fixed point space to the monoid of isomorphism classes
of $G$-re\-presentations
by sending $L \in \bGr(V)^G$ to its isomorphism class.
The isomorphism class of $L$ only depends on the path component of $L$ 
in $\bGr(V)^G$, and the resulting maps 
$\pi_0 ( \bGr(V)^G) \to \bRO^+(G)$
are compatible as $V$ runs through the finite-dimensional $G$-subrepresentations
of $\Uc_G$. So they assemble into a map
\begin{equation}\label{eq:pi^G Gr to RO^+ G}
\pi_0^G ( \bGr )\ = \ \colim_{V\in s(\Uc_G)} \, \pi_0( \bGr(V)^G ) \ \to \ \bRO^+(G)\ ,
\end{equation}
and this map is an isomorphism of monoids with respect to the direct sum
of representations on the target.
Moreover, the isomorphism is compatible with restriction maps,
and it takes the transfer maps induced by the commutative
multiplication of $\bGr$ to induction of representations on the right hand side;
so as $G$ varies, the maps \eqref{eq:pi^G Gr to RO^+ G} form an isomorphism
of global power monoids.

An interesting morphism of ultra-commutative monoids
\begin{equation}  \label{eq:define Gr2bO}
 \tau \ : \ \bGr \ \to \ \bO 
\end{equation}\index{subject}{orthogonal group ultra-commutative monoid}
is defined at an inner product space $V$ as the map  
\[ \tau(V) \ : \ \bGr(V) \ \to \ \bO(V) \ , \quad L \ \longmapsto \ 
p_{L^\perp}-p_L \ , \]
sending a subspace $L\subset V$ 
to the difference of the orthogonal projection onto $L^\perp=V-L$
and the orthogonal projection onto $L$.
Put differently, $\tau(V)(L)$ is the linear isometry that
is multiplication by $-1$ on $L$ and the identity on the orthogonal complement $L^\perp$.
We omit the straightforward verification that these maps
do define a morphism of ultra-commutative monoids.
The induced monoid homomorphism
\[ \pi_0^G(\tau)\ : \ \pi_0^G(\bGr)\ \to \ \pi_0^G(\bO) \]
is easily calculated. The isomorphism \eqref{eq:pi^G Gr to RO^+ G}
identifies the source with the monoid $\bRO^+(G)$ of isomorphism classes
of orthogonal $G$-representations, under direct sum.
By Example \ref{eg:orthogonal group monoid space}
the group $\pi_0^G(\bO)$ is a direct sum of copies of $\mZ/2$,
indexed by the isomorphism classes of irreducible
orthogonal $G$-rep\-resentations of real type. 
If $\lambda$ is any irreducible orthogonal $G$-representation,
then $\pi_0^G(\tau)$ sends its class to the automorphism $-\Id_\lambda$. 
The group $\bO(\lambda)^G$ is isomorphic to $O(1)$, $U(1)$
or $Sp(1)$ depending on whether $\lambda$ is of real, complex or
quaternionic type. In the real case, the map $-\Id_\lambda$ lies in the
non-identity path component; in the complex and quaternionic cases,
the group $\bO(\lambda)^G$ is path connected.
So under the previous isomorphisms, $\pi_0^G(\tau)$ becomes the homomorphism
\[ \bRO^+(G) \ \to \ \bigoplus_{[\lambda]\text{ real}} \mZ/2 \]
that sends the class of $\lambda$ to the generator of the $\lambda$-summand
if $\lambda$ is of real type, and to the trivial element 
if $\lambda$ is of complex or quaternionic type.
Since the classes of irreducible representations freely generate
$\bRO^+(G)$ as an abelian monoid, this determines the morphism $\pi_0^G(\tau)$.
Moreover, this also shows that $\pi_0^G(\bGr)\to\pi_0^G(\bO)$ is surjective.
\end{eg}

\begin{eg}[Oriented Grassmannian]\label{eg:Gr oriented}
A variation of the previous example is 
the orthogonal monoid space $\bGr^{\ort}$ of {\em oriented Grassmannians}.
\index{subject}{oriented Grassmannian}\index{subject}{Grassmannian!oriented|see{oriented Grassmannian}}\index{symbol}{$\bGr^{\ort}$ - {oriented Grassmannian}}
The value of $\bGr^{\ort}$ at an inner product space $V$ is
\[ \bGr^{\ort}(V)\ = \ {\coprod}_{m\geq 0} \, Gr_m^{\ort}(V) \ ,\]
the disjoint union of all oriented Grassmannians in $V$.
Here a point in $Gr_m^{\ort}(V)$ is a pair $(L,[b_1,\dots,b_m])$ consisting
of an $L\in Gr_m(V)$ and an orientation, 
i.e., a $G L^+(L)$-equivalence class  $[b_1,\dots,b_m]$ of bases of $L$.
The structure map induced by $\varphi:V\to W$ sends
$(L,[b_1,\dots,b_m])$ to $(\varphi(L),[\varphi(b_1),\dots\varphi(b_m)])$.
A multiplication on $\bGr^{\ort}$ is given by direct sum:
\begin{align*}
\mu_{V,W}\ : \  \bGr^{\ort}( V ) \times \bGr^{\ort} ( W ) \ &\to \ \bGr^{\ort}( V\oplus W) \\ 
((L,[b_1,\dots,b_m]),\, (L',[b_1',\dots,b_n']) )\ &\longmapsto \\ 
(L\oplus L',\ &[(b_1,0),\dots,(b_m,0),(0,b_1'),\dots,(0,b_n')]) \ .
\end{align*}
The unit is the only point $(\{0\},\emptyset)$ in $\bGr^{\ort}(0)$.

As an orthogonal space, $\bGr^{\ort}$ is the disjoint union of global
classifying spaces of the special orthogonal groups, via the homeomorphisms
\begin{align*}
   (B_{\gl} S O(m))(V)\ = \  \bL(\mR^m,V)/S O(m) \ &\iso \quad  Gr_m^{\ort}(V) \\
\varphi\cdot S O(m)\quad  &\longmapsto \ 
(\varphi(\mR^m),[\varphi(e_1),\dots\varphi(e_m)] )\ .
\end{align*}
\Danger The multiplication of $\bGr^{\ort}$ is {\em not} commutative. The issue is
that when pushing a pair around the two ways of the square
\[ \xymatrix@C=12mm{ 
Gr_m^{\ort}(V)\times Gr^{\ort}_n(W) \ar[r]^-{\mu_{V,W}}\ar[d]_{\text{twist}} & 
Gr_{m+n}^{\ort}(V\oplus W)\ar[d]^{Gr_{m+n}^{\ort}(\tau_{V,W})}\\
Gr_n^{\ort}(W)\times Gr_m^{\ort}(V) \ar[r]_-{\mu_{W,V}} & 
Gr_{n+m}^{\ort}(W\oplus V) } \]
then we end up with the same subspaces of $W\oplus V$, but they
come with {\em different} orientations if $m$ and $n$ are both odd.

We can arrange commutativity of the multiplication by passing
to the orthogonal submonoid $\bGr^{\ort,\ev}$ of {\em even-dimensional}
oriented Grassmannians, defined as
\[ \bGr^{\ort,\ev}(V)\ = \ {\coprod}_{n\geq 0} \, Gr_{2 n}^{\ort}(V) \ ;\]
the multiplication of $\bGr^{\ort,\ev}$ is then commutative.
Moreover, the forgetful map $\bGr^{\ort,\ev}\to\bGr$ 
to the additive Grassmannian is a homomorphism
of ultra-commutative monoids.
\end{eg}

\begin{eg}[Complex and quaternionic Grassmannians]\label{eg:Gr additive complex}
The~{\em complex additive Grassmannian} $\bGr^\mC$\index{subject}{additive Grassmannian!complex}\index{subject}{Grassmannian!complex additive|see{additive Grassmannian}}
\index{symbol}{$\bGr^\mC$ - {complex additive Grassmannian}}
and the~{\em quaternionic additive Grassmannian} $\bGr^\mH$\index{subject}{additive Grassmannian!quaternionic}\index{subject}{Grassmannian!quaternionic additive|see{additive Grassmannian}}
\index{symbol}{$\bGr^\mH$ - {quaternionic additive Grassmannian}}
are two more ultra-commutative monoids, the complex respectively quaternionic analogs 
of Example \ref{eg:Gr additive}.
The underlying orthogonal spaces send an inner product space $V$ to
\[ \bGr^\mC(V)\ = \ {\coprod}_{m\geq 0} \, Gr_m^\mC (V_\mC) 
\text{\quad respectively\quad}
   \bGr^\mH(V)\ = \ {\coprod}_{m\geq 0} \, Gr_m^\mH (V_\mH) \ ,\]
the disjoint union of all complex (respectively quaternionic) Grassmannians 
in the complexification $V_\mC=\mC\tensor_{\mR}V$ 
(respectively in $V_\mH=\mH\tensor_{\mR}V$).
As in the real analog, the structure maps are given by taking images 
under (complexified respectively quaternified) linear isometric embeddings;
direct sum of subspaces, plus identification 
along the isomorphism $V_\mC\oplus W_\mC\iso(V\oplus W)_\mC$
(respectively $V_\mH\oplus W_\mH\iso(V\oplus W)_\mH$),
provides an ultra-commutative multiplication on $\bGr^\mC$ and on $\bGr^\mH$.

The homogeneous summand $\bGr^{\mC,[m]}$ is a global classifying space
for the unitary group $U(m)$.   
Indeed, in Construction \ref{con:complex-free spc} we introduced the 
`complex semifree' orthogonal space $\bL^\mC_{U(m),\mC^m}$ 
associated to the tautological $U(m)$-representation on $\mC^m$.
This orthogonal space  is isomorphic to $\bGr^{\mC,[m]}$ via
\[ \bL^\mC(\mC^m, V_\mC) / U(m) \ \iso \ \bGr^{\mC,[m]}(V)\ ,\quad
\varphi\cdot U(m)\ \longmapsto \ \varphi(\mC^m)\ .\]
Proposition \ref{prop:complex global classifying}~(i) then exhibits a global equivalence
\[ B_{\gl} U(m)\ = \ \bL_{U(m),u(\mC^m)} \ \xra{\ \sim\ } \ \bL^\mC_{U(m),\mC^m}\ \iso \
\bGr^{\mC,[m]} \ . \]
Although Proposition \ref{prop:complex global classifying}~(i) 
does not literally apply to $\bGr^{\mH,[m]}$, it has an $\mH$-linear
analog for symplectic representations, 
showing that the homogeneous summand $\bGr^{\mH,[m]}$ is a global classifying space
for the symplectic group $S p(m)$.\index{subject}{global classifying space!of $U(n)$}    
\index{subject}{global classifying space!of $Sp(n)$} 

The ultra-commutative monoid $\bGr^\mC$ comes with an involutive automorphism
\[ \psi \ : \ \bGr^\mC \ \to \ \bGr^\mC \]
given by complex conjugation.\index{subject}{complex conjugation!on $\bGr^\mC$}
Here we exploit that the complexification of an $\mR$-vector space $V$
comes with a preferred $\mC$-semilinear morphism
\[ \psi_V \ : \ V_\mC \ \to \ V_\mC \ , \quad \lambda\tensor v \ \longmapsto \
\bar\lambda\tensor v \ . \]
The value of $\psi$ at $V$ takes a $\mC$-subspace $L\subset V_\mC$
to the conjugate subspace $\bar L=\psi_V(L)$.
Complexification of subspaces defines a morphism of ultra-commutative monoids\index{subject}{complexification morphism!from $\bGr$ to $\bGr^\mC$}
\[  c \ : \ \bGr \ \to \ \bGr^\mC  \ , \quad
\bGr(V)\ \to \  \bGr^\mC(V)\ , \quad L \ \longmapsto \ L_\mC \]
from the real to the complex additive Grassmannian.
A complex subspace of $V_\mC$ is invariant under $\psi_V$ if and
only if it is the complexification of an $\mR$-subspace of $V$
(namely the $\psi_V$-fixed subspace of $V$).
So the morphism $c$ is an isomorphism of $\bGr$ onto
the $\psi$-invariant ultra-commutative submonoid $(\bGr^\mC)^\psi$. 

Realification defines a morphism of ultra-commutative monoids
\index{subject}{realification morphism!from $\bGr^\mC$ to $\bGr^{\ort,\ev}$}
\[ r \ : \ \bGr^\mC \ \to \ \sh_\tensor^\mC (\bGr^{\ort,\ev} )  \]
to the multiplicative shift (see Example \ref{eg:Additive and multiplicative shift})
of the even part of the oriented Grassmannian of Example \ref{eg:Gr oriented}.
The value $r(V):\bGr^\mC(V)\to \bGr^{\ort,\ev}(V_\mC)$ 
takes a complex subspace of $V_\mC$ to the underlying real vector space,
endowed with the preferred orientation $[x_1, i x_1,\dots, x_n,i x _n]$,
where $(x_1,\dots,x_n)$ is any complex basis.

The isomorphism \eqref{eq:pi^G Gr to RO^+ G}
between $\pi_0^G ( \bGr)$ and $\bRO^+(G)$ has an obvious complex analog.
For every compact Lie group $G$ an isomorphism of commutative monoids
\begin{equation}\label{eq:pi^G Gr to RG}
\pi_0^G ( \bGr^\mC )\ = \ \colim_{V\in s(\Uc_G)} \, \pi_0( \bGr^\mC(V)^G ) 
\ \iso \ \bRU^+(G)
\end{equation}
is given by sending the class of a $G$-fixed point in $\bGr^\mC(V)^G$, 
i.e., a complex $G$-subrepresentation of $V_\mC$, to its isomorphism class.
The isomorphism is compatible with restriction maps,
takes the involution $\pi_0^G(\psi)$ of $\pi_0^G ( \bGr^\mC )$
to the complex conjugation involution of $\bRU^+(G)$,
and it takes the transfer maps induced by the commutative
multiplication of $\bGr^\mC$ to induction of representations.
So as $G$ varies, the maps form an isomorphism of global power monoids
$\upi_0( \bGr^\mC )\iso\bRU^+$. The isomorphisms are also compatible
with complexification and realification, in the sense of the commutative diagram:
\[ \xymatrix@C6mm{ 
\upi_0(\bGr) \ar[r]^-{\upi_0(c)} \ar[d]_{\eqref{eq:pi^G Gr to RO^+ G}}^\iso &
\upi_0(\bGr^\mC) \ar[r]^-{\upi_0(r)}  \ar[d]^{\eqref{eq:pi^G Gr to RG}}_\iso & 
\upi_0(\sh_\mC^\tensor( \bGr^{\ort,\ev} )) \ar[r]^-\iso & 
\upi_0( \bGr^{\ort,\ev}) \ar[r]^-{\text{forget}} & 
\upi_0( \bGr ) \ar[d]^{\eqref{eq:pi^G Gr to RO^+ G}}_\iso  \\
\bRO^+\ar[r]_-c & \bRU^+ \ar[rrr]_-r &&& \bRO^+ } \]
The lower horizontal maps are complexification 
and realification of representations.
The isomorphism in the upper row is inverse to the one induced by the homomorphism
$\bGr^{\ort,\ev}\circ i:\bGr^{\ort,\ev}\to \sh_\mC^\tensor( \bGr^{\ort,\ev} )$, i.e.,
precomposition with the natural linear isometric embedding
\[  i_V\ : \ V\ \to\ V_\mC\ , \quad v \ \longmapsto \ 1\tensor v \ ;\]
the morphism $\bGr^{\ort,\ev}\circ i$ is a global equivalence 
by Theorem \ref{thm:general shift of osp}.
\end{eg}

\begin{eg}[Multiplicative Grassmannians]\label{eg:Gr multiplicative real}
We define an ultra-commuta\-tive monoid $\bGr_\tensor$,
the {\em multiplicative Grassmannian}.\index{subject}{multiplicative Grassmannian}\index{subject}{Grassmannian!multiplicative|see{multiplicative Grassmannian}} 
We let 
\[ \Sym(V) \ = \ {\bigoplus}_{i\geq 0}\, \Sym^i(V)\ = \ 
{\bigoplus}_{i\geq 0}\, V^{\tensor i}/\Sigma_i \]
denote the symmetric algebra of an inner product space $V$.
If $W$ is another inner product space, then the two direct summand inclusions
induce algebra homomorphisms
\[  \Sym(V)\ \xra{\quad}\ \Sym(V\oplus W) \ \xla{\quad} \ \Sym(W)\ .\]
We use the commutative multiplication on $\Sym(V\oplus W)$
to combine these into an $\mR$-algebra isomorphism
\begin{equation}\label{eq:sym of sum and tensor}
 \Sym(V)\tensor \Sym(W)\ \iso \ \Sym(V\oplus W) \ .
\end{equation}
These isomorphisms are natural for linear isometric embeddings in $V$ and $W$.
The value of $\bGr_\tensor$\index{symbol}{$\bGr_\tensor$ - {multiplicative Grassmannian}} 
at an inner product space $V$ is then
\[ \bGr_\tensor(V)\ = \ {\coprod}_{n\geq 0} \, Gr_n(\Sym(V)) \ ,\]
the disjoint union of all Grassmannians 
in the symmetric algebra of $V$.
The structure map $\bGr_\tensor(\varphi):\bGr_\tensor(V)\to\bGr_\tensor(W)$
induced by a linear isometric embedding $\varphi:V\to W$
is given by  
\[ \bGr_\tensor(\varphi)(L) \ = \ \Sym(\varphi)(L)\ , \]
where $\Sym(\varphi):\Sym(V)\to \Sym(W)$ is the induced map of symmetric algebras.
A commutative multiplication on $\bGr_\tensor$ is given by tensor product, i.e.,
\begin{align*}
\mu_{V,W}\ : \  \bGr_\tensor( V ) \times \bGr_\tensor ( W ) \ \to \ \bGr_\tensor( V\oplus W) 
\end{align*}
sends $(L,L')\in\bGr_\tensor( V ) \times \bGr_\tensor ( W )$ to the image of
$L\otimes L'$ under the isomorphism \eqref{eq:sym of sum and tensor}. 
The multiplicative unit is the point $\mR$ in $\bGr_\tensor(0)=\mR$.
As the additive Grassmannian $\bGr$, the multiplicative Grassmannian $\bGr_\tensor$ 
is $\mN$-graded, with degree $n$
part given by $\bGr^{[n]}_\tensor(V)=Gr_n(\Sym(V))$. 
The multiplication sends $\bGr_\tensor^{[m]}( V ) \times \bGr_\tensor^{[n]} ( W )$ to 
$\bGr_\tensor^{[m\cdot n]}(V\oplus W)$. 

\medskip

\Danger As orthogonal spaces, the additive and multiplicative Grassmannians are
globally equivalent. 
For an inner product space $V$ we let
$i:V\to \Sym(V)$ be the embedding as the linear summand of the
symmetric algebra. Then as $V$ varies, the maps
\[ \bGr(V)\ = \ {\coprod}_{n\geq 0} Gr_n(V) \ \to \ 
{\coprod}_{n\geq 0} Gr_n(\Sym(V)) \ = \  \bGr_\tensor(V)\]
sending $L$ to $i(L)$ form a global equivalence $\bGr\to\bGr_\tensor$.
Indeed, for each $n\geq 0$, $\bGr_\tensor^{[n]}(V)$ is a sequential colimit, 
along closed embeddings, of a sequence of orthogonal spaces 
\[ \bGr^{[n]} \ \to \ \bGr^{[n]}_{\leq 1}\ \to \ \bGr^{[n]}_{\leq 2}\ \to \ 
\dots\ \to \ \bGr^{[n]}_{\leq k}\ \to \ \dots  \ ,\]
where $\bGr_{\leq k}^{[n]}(V)=G r_n(\bigoplus_{i=0}^k \Sym^i(V))$.
Each of the morphisms $\bGr^{[n]}\to \bGr_{\leq k}^{[n]}$ is a global equivalence 
by Theorem \ref{thm:general shift of osp}; so all morphisms in the
sequence are global equivalences, hence so is the map from $\bGr^{[n]}$
to the colimit $\bGr_{\tensor}^{[n]}$, by
Proposition \ref{prop:global equiv basics}~(ix).
This global equivalence induces a bijection
\[ \pi_0^G(\bGr) \ \iso \ \pi_0^G(\bGr_\tensor)\]
for every compact Lie group $G$, hence both are isomorphic
to the set $\bRO^+(G)$ of isomorphism classes of orthogonal $G$-representations. 
The commutative mo\-noid structures and transfer maps
induced by the products of $\bGr$ respectively $\bGr_\tensor$ are
quite different though: the monoid structure of $\pi_0^G(\bGr)$ 
corresponds to direct sum of representations, and the transfer maps
are additive transfers;
the monoid structure of $\pi_0^G(\bGr_\tensor)$ 
corresponds to tensor product of representations, and the transfer maps
are multiplicative transfers, also called norm maps.

The orthogonal subspace $\bP=\bGr_\tensor^{[1]}$ 
of the multiplicative Grassmannian $\bGr_\tensor$ 
is closed under the product and contains the multiplicative unit, 
hence $\bP$ is an ultra-commutative monoid 
in its own right. Because
\[ \bP(V)\ = \ \bGr_\tensor^{[1]}(V) \ = \ P(\Sym(V)) \]
is the projective space of the symmetric algebra of $V$,
we use the symbol $\bP$ and refer to it as the 
{\em global projective space}.\index{subject}{global projective space}\index{symbol}{$\bP$ - {global projective space}}
The multiplication is given by tensor product of lines,
and application of the isomorphism \eqref{eq:sym of sum and tensor}.
Since $\bP=\bGr_\tensor^{[1]}$ is globally equivalent
to the additive variant $\bGr^{[1]}$, it is a global classifying space
for the group $O(1)$, a cyclic group of order~2,
\[ \bP \ \simeq \ \bGr^{[1]} \ \simeq \ B_{\gl} O(1) \ = \ B_{\gl} C_2 \ . \]
In other words, $\bP$ is an ultra-commutative  
multiplicative model for $B_{\gl} C_2$.\index{subject}{global classifying space!of $C_2$} 

There is a straightforward complex analog of the multiplicative Grassmannian,
the ultra-commutative monoid $\bGr_\tensor^\mC$ with value at $V$ given by
\[ \bGr^\mC_\tensor(V)\ = \ {\coprod}_{n\geq 0} \, Gr_n^\mC(\Sym(V)_\mC) \ ,\]
the disjoint union of all Grassmannians 
in the complexified symmetric algebra of $V$.
The structure maps and multiplication (by tensor product) are as in the real case.
The orthogonal subspace $\bP^\mC=\bGr_\tensor^{\mC,[1]}$ 
consisting of 1-dimensional subspaces
is closed under the product and contains the multiplicative unit; 
hence $\bP^\mC$ is an ultra-commutative monoid 
in its own right, the {\em complex global projective space}.
As an orthogonal space, $\bP^\mC$ is globally equivalent
to the additive variant $\bGr^{\mC,[1]}$, and hence a global classifying space
for the group $U(1)$,\index{subject}{global classifying space!of $U(1)$}  
\begin{equation}\label{eq:complex_P}  
 \bP^\mC \ \simeq \ \bGr^{\mC,[1]} \ \simeq \ B_{\gl} U(1) \ . 
\end{equation}
\end{eg}

\Danger Since the multiplication in the skew field of quaternions $\mH$ 
is not commutative, there is no
tensor product of $\mH$-vector spaces; so there is no multiplicative
version of the quaternionic Grassmannian $\bGr^\mH$.

\begin{construction}[Bar construction]\label{con:bar construction}
For the next class of examples we quickly recall the {\em bar construction}
of a topological monoid $M$.\index{subject}{bar construction!of a topological monoid}
This is the geometric realization of the simplicial space $M^{\bullet}$ whose space 
of $n$-simplices is $M^n$, the $n$-fold cartesian power of $M$.
For $n\geq 1$ and $0\leq i \leq n$, the face map $d_i:M^n\to M^{n-1}$
is given by
\[ d_i(x_1,\dots x_n) \ = \left\lbrace \begin{array}{ll}
(x_2,\dots,x_n) & \mbox{for $i=0$,} \\
(x_1,\dots,x_{i-1},x_i\cdot x_{i+1},x_{i+2},\dots,x_n)  & \mbox{for $0<i< n$,} \\
(x_1,\dots,x_{n-1}) & \mbox{for $i=n$.}
\end{array} \right.  \]
For $n\geq 1$ and $0\leq i \leq n-1$ the degeneracy map  $s_i:M^{n-1}\to M^n$
is given by
\[ s_i(x_1,\dots, x_{n-1}) \ = \ (x_1,\dots,x_{i},1,x_{i+1},\dots,x_{n-1}) \ . \]
The bar construction is the geometric realization 
(see Construction \ref{con:realize simplicial space})
\[ B M \ = \ | M^\bullet| \ = \ |[n]\mapsto M^n| \]
of this simplicial space; the construction $M\mapsto B M$ 
is functorial in continuous monoid homomorphisms. 
The bar construction commutes with products in the sense that
for a pair of topological monoids $M$ and $N$, the canonical map
\begin{equation}\label{eq:B preserves product}
(B p_M, B p_N) \ : \ B(M\times N)\ \to \ B M\times B N  
\end{equation}
is a homeomorphism, where $p_M:M\times N\to M$ and $p_N:M\times N\to N$
are the projections. Indeed, this map factors as the composite of two maps
\[ | (M\times N)^{\bullet} | \ \xra{|(p_M^\bullet,p_N^\bullet)|} \ 
|  M^\bullet\times N^\bullet | \ \xra{(|p_{M^\bullet}|,|p_{N^\bullet}|)} \ 
|  M^\bullet|\times | N^\bullet | \ .\]
The first map is the realization of an isomorphism of simplicial spaces,
given levelwise by shuffling the factors.
The second map is a homeomorphism because realization commutes with products,
see Proposition \ref{prop:iterated geometric realization}~(ii).
\end{construction}

\begin{construction}[Multiplicative global classifying spaces]\label{con:multiplicative B G}
As we discuss now, all {\em abelian} compact Lie groups admit
multiplicative models of their global classifying spaces.
We use the bar construction, giving a non-equivariant classifying space, followed by the
cofree functor $R$ (see Construction \ref{con:cofree_orthogonal_space}).
The cofree functor takes a space $A$ to the orthogonal space $R A$ with
values
\[ (R A)(V)\ = \ \map(\bL(V,\mR^\infty),A) \ .\]
We endow the cofree functor with a lax symmetric monoidal transformation
\[ \mu_{A,B}\ : \ R A \boxtimes R B \ \to \ R(A\times B)\ . \]
To construct $\mu_{A,B}$ we start from the continuous maps
\begin{align*}
 \map(\bL(V,\mR^\infty), A) \times \map(\bL(W,\mR^\infty), B) 
\ \xra{\ \times \ } \
&\map(\bL(V,\mR^\infty)\times \bL(W,\mR^\infty), A\times B) \\
\xra{ (\res_{V,W})^*}\ 
&\map(\bL(V\oplus W,\mR^\infty), A\times B )\nonumber
\end{align*}
that constitute a bimorphism from $(R A,R B)$ to $R(A\times B)$.
Here 
\[ \res_{V,W}\ :\ \bL(V\oplus W,\mR^\infty)\ \to\  \bL(V,\mR^\infty)\times \bL(W,\mR^\infty) \]
restricts an embedding of $V\oplus W$ to the summands $V$ and $W$.
The morphism $\mu_{A,B}$ is associated to this bimorphism
via the universal property of the box product.

The bar construction preserves products in the sense that for every
pair of compact Lie groups $G$ and $K$ the canonical map
\[ (B p_G, B p_K)\ : \  B(G\times K)\ \to \ B G\times B K  \]
is a homeomorphism, compare \eqref{eq:B preserves product}.
So the composite $G\mapsto R(B G)$ is a lax symmetric
monoidal functor via the morphism of orthogonal spaces
\[ R(B G) \boxtimes R(B K) \ \xra{\ \mu_{B G,B K}\ } \ 
R(B G\times B K) \ \iso \ R(B(G\times K))\ .\]
The bar construction is functorial in continuous group homomorphisms, so for
an abelian compact Lie group $A$ the composite
\[ R(B A) \boxtimes R(B A) \ \to \  R(B(A\times A))\ \xra{R(B\mu_A)} \ R(B A) \]
is an ultra-commutative and associative multiplication on the orthogonal space
$R(B A)$, where $\mu_A:A\times A\to A$ is the multiplication of $A$.
Theorem \ref{thm:BA is right induced} shows that for abelian $A$
the cofree orthogonal space $R(B A)$ is a global classifying space for $A$.
\index{subject}{cofree orthogonal space}
In particular, the Rep-functor $\upi_0(R(B A))$ is representable by $A$.
We saw in Proposition \ref{prop:unique power Rep(-,A)}
that there is then a {\em unique} structure of global power monoid on
$\upi_0(R(B A))$, and the power operations are characterized by the relation
\[ [m](u_A)\ = \ p_m^*(u_A)\]
where $u_A\in\pi_0^A(R(B A))$ is a tautological class 
and $p_m:\Sigma_m\wr A\to A$ is the homomorphism defined by
\[ p_m(\sigma;\,a_1,\dots,a_m) \ = \ a_1\cdot\ldots\cdot a_m \ .\]
\end{construction}

\begin{eg}[Unordered frames]\label{eg:unordered frames}
The ultra-commutative monoid $\bF$ of {\em unordered frames}
\index{subject}{unordered frames}\index{symbol}{$\bF$ - {ultra-commutative monoid of unordered frames}}  
sends an inner product space $V$ to
\[ \bF(V)\ = \ \{ A \subset V\ |\ \text{$A$ is orthonormal}\}\ , \]
the space of all unordered frames in $V$, i.e., subsets of $V$ that consist of
pairwise orthogonal unit vectors. Since $V$ is finite-dimensional, such a subset
is necessarily finite. The topology on $\bF(V)$ is as the disjoint union,
over the cardinality of the sets, of quotient spaces of Stiefel manifolds.
The structure map induced by a linear isometric embedding $\varphi:V\to W$
is given by  $\bF(\varphi)(A) = \varphi(A)$.
A commutative multiplication on $\bF$ is given, essentially, by disjoint union:
\[ 
\mu_{V,W}\ : \  \bF( V ) \times \bF ( W ) \ \to \ \bF( V\oplus W) \ ,\quad
(A,A')\ \longmapsto \ i_V(A) \cup i_W(A')\ ;
 \]
here $i_V:V\to V\oplus W$ and $i_W:W\to V\oplus W$ are the direct summand embed\-dings.
The unit is the empty set, the only point in $\bF(0)$.
The orthogonal space $\bF$ is naturally $\mN$-graded, with degree $m$
part $\bF^{[m]}$ given by the unordered frames of cardinality $m$; 
the multiplication sends $\bF^{[m]}( V ) \times \bF^{[n]} ( W )$ 
to $\bF^{[m+n]}( V\oplus W)$. 

As an orthogonal space, $\bF$ is the disjoint union of global
classifying spaces of the symmetric groups. 
We let $\Sigma_m$ act on $\mR^m$ by permuting the coordinates,
which is also the permutation representation of the tautological
$\Sigma_m$-action on $\{1,\dots,m\}$.
This $\Sigma_m$-action is faithful, so the semifree orthogonal space $\bL_{\Sigma_m,\mR^m}$
is a global classifying space for the symmetric group.
The homeomorphisms
\[  \bL(\mR^m,V)/\Sigma_m \ \iso \  \bF^{[m]}(V) \ , \quad 
\varphi\cdot \Sigma_m\ \longmapsto \ \{\varphi(e_1),\dots,\varphi(e_m)\} \]
show that the homogeneous summand $\bF^{[m]}$ is isomorphic to
$\bL_{\Sigma_m\mR^m}=B_{\gl} \Sigma_m$; 
here $e_1,\dots,e_m$ is the canonical basis of $\mR^m$.
So as orthogonal spaces,\index{subject}{global classifying space!of $\Sigma_n$} 
\[ \bF \ = \ {\coprod}_{m\geq 0} \, B_{\gl} \Sigma_m\ . \]
Proposition \ref{prop:fix of global classifying}~(ii) identifies the
equivariant homotopy set $\pi_0^G ( B_{\gl} \Sigma_m)$
with the set of conjugacy classes of continuous homomorphisms from $G$ to $\Sigma_m$;
by restricting the tautological $\Sigma_m$-representation on $\{1,\dots,m\}$, 
this set bijects with the set of isomorphism classes of 
finite $G$-sets of cardinality $m$.

As $m$ varies, this gives an isomorphism of monoids from $\pi_0^G(\bF)$ 
to the set $\mA^+(G)$\index{symbol}{$\mA^+(G)$ - {monoid of isomorphism classes of finite $G$-sets}} of isomorphism classes of finite $G$-sets 
that we make explicit now.
We let $V$ be a $G$-representation. 
An unordered frame $A\in \bF(V)$ is a $G$-fixed point if and only if it is $G$-invariant.
So for such frames, the $G$-action restricts to an action on $A$,
making it a finite $G$-set.
We define a map
\[  \bF(V)^G \  \to \ \mA^+(G) \ , \quad A \ \longmapsto \ [A]\]
from this fixed point space to the monoid of isomorphism classes of finite $G$-sets.
The isomorphism class of $A$ as a $G$-set 
only depends on the path component of $A$ in $\bF(V)^G$, and the resulting maps 
$\pi_0 ( \bF(V)^G) \to \mA^+(G)$
are compatible as $V$ runs through the finite-dimensional $G$-subrepresentations
of $\Uc_G$. So they assemble into a map
\begin{equation}\label{eq:pi^G F to A^+ G}
\pi_0^G ( \bF )\ = \ \colim_{V\in s(\Uc_G)} \, \pi_0( \bF(V)^G ) \ \to \ \mA^+(G)\ ,
\end{equation}
and this map is a monoid isomorphism with respect to the disjoint union
of $G$-sets on the target.
Moreover, the isomorphisms are compatible with restriction maps,
and they take the transfer maps induced by the commutative
multiplication of $\bF$ to induction of equivariant sets on the right hand side.

\medskip

\Danger It goes without saying that actions of compact Lie groups
are required to be {\em continuous} and that the use of the term `set'
(as opposed to `space') implies the discrete topology on the set; 
so the identity path component $G^\circ$ acts trivially
on every $G$-set.
Hence the monoids $\pi_0^G ( \bF )$ and $\mA^+(G)$
only see the finite group $\pi_0(G)=G/G^\circ=\bar G$ of path components, i.e.,
for every compact Lie group $G$, the inflation maps
\[ p^* \ : \ \pi_0^{\bar G} ( \bF ) \ \to \ \pi_0^G ( \bF )
\text{\qquad and\qquad}
p^* \ : \ \mA^+(\bar G) \ \to \ \mA^+(G) \]
along the projection $p:G\to\bar G$ are isomorphisms. 
So if $G$ has positive dimension, then 
the group completion of the monoid $\mA^+(G)$ need not be isomorphic 
to what is sometimes called the Burnside ring of $G$ 
(the 0-th $G$-equivariant stable stem).

\medskip

A morphism of $\mN$-graded ultra-commutative monoids
\[ \spn \ : \ \bF \ \to \ \bGr \text{\qquad is defined by\qquad}
 \spn(V)(A)\ = \ \spn(A)\ ,\]
i.e., a frame is sent to its linear span.
The induced morphism of global power monoids
is linearization: the square of monoid homomorphisms
\[ 
 \xymatrix@C=18mm{  
\pi_0^G(\bF)\ar[d]_{\eqref{eq:pi^G F to A^+ G}}^\iso  \ar[r]^-{\pi_0^G(\spn)}&
\pi_0^G(\bGr)\ar[d]^{\eqref{eq:pi^G Gr to RO^+ G}}_\iso \\
\mA^+(G)\ar[r]_-{[S]\,\longmapsto\, [\mR S]} & \bRO^+(G)}    
 \]
commutes, where the lower map sends the class of a $G$-set to the class of its
permutation representation.
\end{eg}

\begin{eg}[Multiplicative monoid of the sphere spectrum]\label{eg:Omega^bullet S}
We define an ultra-commutative monoid $\Omega^\bullet \mS$,
the `multiplicative monoid of the sphere spectrum'.\index{subject}{sphere spectrum} 
The notation and terminology indicate that this is a special case
of a more general construction that associates to an ultra-commutative
ring spectrum $R$ its multiplicative monoid $\Omega^\bullet R$,
see Example \ref{eg:suspension spectrum of orthogonal monoid space} below.

The values of the orthogonal space $\Omega^\bullet \mS$ are given by
\[ ( \Omega^\bullet \mS)(V) \ = \ \map_*(S^V,S^V)\ , \]
the space of continuous based self-maps of the sphere $S^V$.
A linear isometric embedding $\varphi:V\to W$ acts by conjugation
and extension by the identity, i.e., the map 
\[ (\Omega^\bullet \mS)(\varphi)\ : \ \map_*(S^V,S^V)\ \to \ \map_*(S^W,S^W)\]
sends a continuous based map $f:S^V\to S^V$ to the composite
\[ S^W \iso \ S^V\sm S^{W-\varphi(V)}
\ \xra{f\sm S^{W-\varphi(V)} } \ S^V\sm S^{W-\varphi(V)}
\ \iso\  S^W\ . \]
The two unnamed homeomorphisms between $S^V\sm S^{W-\varphi(V)}$ and $S^W$
use the map $\varphi$ on the factor $S^V$.
In particular, the orthogonal group $O(V)$ acts on $\map_*(S^V,S^V)$
by conjugation.

The multiplication of $\Omega^\bullet \mS$ is by smash product, i.e., the map
\[ \mu_{V,W}\ : \  ( \Omega^\bullet \mS)(V) \times ( \Omega^\bullet \mS)(W) \ \to \ 
( \Omega^\bullet \mS)(V\oplus W) \]
smashes a self-map of $S^V$ with a self-map of $S^W$ and
conjugates with the canonical homeomorphism between $S^V\sm S^W$ and $S^{V\oplus W}$.
The unit is the identity of $S^V$.

The equivariant homotopy set
$\pi_0^G(\Omega^\bullet \mS)$ is equal 
to the stable $G$-equivariant 0-stem $\pi_0^G(\mS)$,
compare Construction \ref{con:Omega^bullet} below.
The monoid structure on $\pi_0^G(\Omega^\bullet \mS)$ 
arising from the multiplication on $\Omega^\bullet \mS$
is the {\em multiplicative} (rather than the additive) 
monoid structure of $\pi_0^G(\mS)$.
The set $\pi_0^G(\Omega^\bullet \mS)$
thus bijects with the underlying set of the Burnside ring $\mA(G)$ of the group $G$
(compare Example \ref{eg:sphere spectrum}), which is
additively a free abelian group with basis 
the conjugacy classes of closed subgroups of $G$ with finite Weyl group.
The multiplication on $\pi_0^G(\Omega^\bullet \mS)$ corresponds to the multiplication 
(not the addition !) in the Burnside ring $\mA(G)$.
When $G$ is finite, $\pi_0^G(\Omega^\bullet \mS)$ bijects with the
underlying set of the Grothendieck group of finite $G$-sets,
and the multiplication corresponds to the product of $G$-sets.
The power operations in  $\upi_0(\Omega^\bullet \mS)$ 
are thus represented by `raising a $G$-set to the cartesian power',
and the transfer maps are known as `norm maps' or `multiplicative induction'.
\end{eg}

\begin{eg}[Exponential homomorphisms]
The classical $J$-homomor\-phism\index{subject}{J-homomorphism@$J$-homomorphism}
fits in nicely here, in the form of a global refinement 
to a morphism of ultra-commutative monoids\index{subject}{orthogonal group ultra-commutative monoid}
\[ J \ : \ \bO \ \to \  \Omega^\bullet\mS \]
defined at an inner product space $V$ as the map
\[ J(V) \ : \ \bO(V) \ \to \  \map_*(S^V,S^V) \]
sending a linear isometry $\varphi:V\to V$
to its one-point compactification $S^\varphi:S^V\to S^V$.
The fact that these maps are multiplicative and 
compatible with the structure maps is straightforward.
The induced map
\[ \pi_0^G(J)\ : \ \pi_0^G(\bO)\ \to \ \pi_0^G(\Omega^\bullet\mS) \ =\ \pi_0^G(\mS)\ ,\]
for $G$ a compact Lie group, can be described as follows.
By Example \ref{eg:orthogonal group monoid space}, 
the group $\pi_0^G(\bO)$ is a direct sum of copies of $\mZ/2$,
indexed by the isomorphism classes of irreducible
orthogonal $G$-representations of real type. 
If $\lambda$ is such an irreducible $G$-representation,
then the image of the $\lambda$-indexed copy of $\mZ/2$
is represented by the antipodal map of $S^\lambda$.

In \eqref{eq:define Gr2bO} we defined 
a morphism of ultra-commutative monoids $\tau:\bGr\to\bO$. 
The composite morphism
\begin{equation}\label{eq:sh_J_after_tau}
\bGr\ \xra{\ \tau\ }\  \bO \ \xra{\  J\ }\ \Omega^\bullet \mS 
\end{equation}
realizes an `exponential' homomorphism\index{subject}{exponential homomorphism!from $\bRO(G)$ to $\mA(G)^{\times}$}  
from the real representation ring to the  multiplicative group
of the Burnside ring of a compact Lie group $G$.
The exponential homomorphism was studied by tom Dieck 
and is sometimes called the `tom Dieck exponential map'.
Tom Dieck's definition of the exponential homomorphism 
in \cite[5.5.9]{tomDieck-and representation theory}
is completely algebraic: we start from the homomorphism
\[ s\ : \
\bRO^+(G)\ \to \ \text{Cl}(G,\mZ)^\times\ , \quad  s[V](H)\ = \ (-1)^{\dim (V^H)}\]
that sends an orthogonal representation $V$ to the class function $s(V)$ 
that records the parities of the dimensions of the fixed point spaces.
The Burnside ring embeds into the ring of class functions by `fixed point counting'
\[ \Phi\ : \ \mA(G)\ \to \ \Cl(G,\mZ)\ , \ \Phi[S](H)\ = \ |S^H|\ , \]
i.e., a (virtual) $G$-set $S$ is sent to the class function that counts the number
of fixed points. This map is injective, and for finite groups the image
can be characterized by an explicit system of congruences \cite[Prop.\,1.3.5]{tomDieck-and representation theory}.
The image of $s$ satisfies the congruences, so there is
a unique map $\exp:\bRO^+(G)\to\mA(G)^\times$
such that $\Phi\circ\exp=s$.
The map $s$ sends the direct sum of representations to the product
of the parity functions, so $\exp$ is a homomorphism of abelian monoids,
and thus extends to a homomorphism
from the orthogonal representation ring $\bRO(G)$.
The morphism of ultra-commutative monoids \eqref{eq:sh_J_after_tau} 
realizes the exponential morphism in the sense that the following diagram 
of monoid homomorphisms commutes:
\[ \xymatrix@C=6mm{  
\pi_0^G(\bGr)\ar[d]_{\eqref{eq:pi^G Gr to RO^+ G}}^\iso \ar[rr]^-{\pi_0^G(J\circ \tau)} && 
\pi_0^G(\Omega^\bullet \mS) \ar@{=}[r] & (\pi_0^G(\mS))^\times\\ 
\bRO^+(G)\ar[rrr]_{\exp} &&& \mA(G)^\times \ar[u]_\iso^\alpha}\]
To show the commutativity of this diagram we may compose with
the degree monomorphism 
\[ \deg\ :\ (\pi_0^G(\mS))^\times\to\text{Cl}(G,\mZ)\ , \quad
\deg[f](H) \ = \ \deg(f^H:S^{V^H}\to S^{V^H} )\]
that takes the class of an equivariant self-map $f:S^V\to S^V$
of a representation sphere to the class function that records the
fixed point dimensions. The composite $\deg\circ\alpha:\mA(G)\to \text{Cl}(G,\mZ)$
coincides with the fixed point counting map $\Phi$, so 
\[ \deg\circ\alpha\circ\exp\ = \ \Phi\circ\exp \ = \ s \ : 
\bRO^+(G)\ \to \ \text{Cl}(G,\mZ)\ . \]
On the other hand, the map $\pi_0^G(J\circ\tau)$ 
sends the class of a $G$-representation $V$ to
the involution $S^{-\Id_V}:S^V\to S^V$. The diagram thus commutes because
\[ \deg\left( (S^{-\Id_V})^H \right)\ = \ (-1)^{\dim(V^H)}\ = \ s[V](H)\ . \]
\end{eg}

\section{Global forms of \texorpdfstring{$B O$}{BO}}
\label{sec:forms of BO}

In this section we discuss different orthogonal monoid spaces whose
underlying non-equivariant homotopy type is $B O$, 
a classifying space for the infinite orthogonal group.
Each example is interesting in its own right, 
and as a whole, the global forms of $B O$ are a great illustration of how
non-equivariant homotopy types `fold up' into many different global
homotopy types.
The different forms of $B O$ have associated orthogonal Thom spectra
with underlying non-equivariant stable homotopy type $M O$;
we will return to these Thom spectra in Section \ref{sec:global Thom}.
The examples we discuss here all come with multiplications, some ultra-commutative, 
but some only $E_\infty$-commutative; 
so our case study also illustrates the different degrees of commutativity 
that arise `in nature'.

We can name five different global homotopy types that all have the same underlying
non-equivariant homotopy type, namely that of a classifying space of 
the infinite orthogonal group:
\begin{itemize}
\item the constant orthogonal space $\underline{B O}$ with value a classifying space
of the infinite orthogonal group;
\item the `full Grassmannian model' $\bBO$,
the degree~0 part of the periodic global 
Grassmannian $\bBOP$ (Example \ref{eg:BOP});
\item the bar construction model $\mathbf B^\circ\bO$ (Construction \ref{eg:bar construction BO});
\item the `restricted Grassmannian model' $\bbO$ that is also
a sequential homotopy colimit of the global classifying spaces $B_{\gl} O(n)$
(Example \ref{eg:define bO});
\item the cofree orthogonal space $R(BO)$ associated to a classifying space
of the infinite orthogonal group (Construction \ref{con:cofree_orthogonal_space}).
\end{itemize}
These global homotopy types are related by weak morphisms of orthogonal spaces:
\[ 
\xymatrix{
\underline{B O}\ar[r]\ar[d] & \bbO \ar[d] \\
\mathbf B^\circ\bO \ar[r] & \bBO \ar[r]&  R(BO) } 
 \]
The orthogonal spaces $\mathbf B^\circ\bO$ and $\bBO$ come with
ultra-commutative multiplications.
The global homotopy type of $R(BO)$ also admits an ultra-commutative multiplication;
we will not elaborate this point, but one way to see it is to
extend the cofree functor $R$ to a lax symmetric monoidal functor
on the category of orthogonal spaces, so that $R(\bBO)$
is an ultra-commutative monoid within this global homotopy type.
The orthogonal spaces $\underline{B O}$ and $\bbO$ admit $E_\infty$-multiplications; 
for $\underline{B O}$ this is a consequence of the non-equivariant $E_\infty$-structure of $BO$.
All weak morphisms above can be arranged to preserve the $E_\infty$-multiplications,
so they induce additive maps of abelian monoids on $\pi_0^G$
for every compact Lie group $G$. 

As we explain in Example \ref{eg:general bar construction},
the bar construction model makes sense more generally for
monoid valued orthogonal spaces; in particular, applying the bar construction
to the ultra-commutative monoids made from the families of classical Lie groups
discussed in the previous section provides ultra-commutative monoids 
$\mathbf B^\circ\bSO$, $\mathbf B^\circ\bU$, $\mathbf B^\circ\bSU$
$\mathbf B^\circ\bSp$, $\mathbf B^\circ\bSpin$ and $\mathbf B^\circ\bSpin^c$.
Example \ref{eg:BUP and BSpP} introduces the complex and quaternionic
analogs of $\bBO$ and $\bbO$, i.e., the ultra-commutative monoids $\bBU$ and $\bBSp$
and the $E_\infty$-orthogonal monoid spaces $\bbU$ and $\bbSp$.

\begin{eg}[Periodic Grassmannian]\label{eg:BOP} 
We define an ultra-commutative monoid $\bBOP$\index{symbol}{$\bBOP$ - {periodic global $BO$}}\index{subject}{periodic global $BO$} that is a global refinement
of the non-equivariant homotopy type $\mZ\times BO$,
and at the same time a global group completion\index{subject}{group completion}
of the additive Grassmannian $\bGr$ introduced in Example \ref{eg:Gr additive}. 
The orthogonal space $\bBOP$ comes with tautological vector bundles whose Thom spaces
form the periodic Thom spectrum $\bMOP$, 
discussed in Example \ref{eg:MO} below.

For an inner product space $V$ we set
\[  \bBOP(V) \ = \  {\coprod}_{m\geq 0} \, Gr_m(V^2)\ , \]
the disjoint union of the Grassmannians 
of $m$-dimensional subspaces in $V^2=V\oplus V$.
The structure map associated to a linear isometric embedding $\varphi:V\to W$
is given by 
\[ \bBOP(\varphi)(L) \ = \ \varphi^2(L) \ +\ ((W-\varphi(V))\oplus 0)\ ,\]
the internal orthogonal sum of the image of $L$ under $\varphi^2:V^2\to W^2$
and the orthogonal complement of the image of $\varphi:V\to W$, viewed as sitting
in the first summand of $W^2=W\oplus W$.
In particular, the orthogonal group $O(V)$ acts on $\bBOP(V)$
through its diagonal action on $V^2$.

We make the orthogonal space $\bBOP$ into an ultra-commutative monoid 
by endowing it with multiplication maps
\[ \mu_{V,W}\ : \ \bBOP(V) \times \bBOP(W) \ \to \ \bBOP(V\oplus W) \ ,\quad
(L,L') \ \longmapsto \  \kappa_{V,W}(L\oplus L')\ ,\]
where 
\[ \kappa_{V,W}\ : \ V^2\oplus W^2\ \iso\ (V\oplus W)^2 \text{\quad is defined by\quad}
\kappa_{V,W}(v,v',w,w')=(v,w,v',w')\ .\]
The unit is the unique element $\{0\}$ of $\bBOP(0)$.

The orthogonal space $\bBOP$ is naturally $\mZ$-graded: for $m\in\mZ$ we let
\[ \bBOP^{[m]}(V)\ \subset \ \bBOP(V) \]
be the path component consisting of all subspaces $L\subset V^2$
such that $\dim(L)=\dim(V) + m$. For fixed $m$ the spaces $\bBOP^{[m]}(V)$
form a subfunctor of $\bBOP$, i.e., $\bBOP^{[m]}$ is an orthogonal subspace
of $\bBOP$. The multiplication is graded in the sense that $\mu_{V,W}$
takes  $\bBOP^{[m]}(V) \times \bBOP^{[n]}(W)$ to $\bBOP^{[m+n]}(V\oplus W)$.
We write $\bBO=\bBOP^{[0]}$\index{symbol}{$\bBO$ - {global $BO$}} 
for the homogeneous summand of $\bBOP$ of degree~0,
which is thus an ultra-commutative monoid in its own right.
The underlying non-equivariant homotopy type of $\bBO$ is that of
a classifying space of the infinite orthogonal group.

\medskip

\Danger While $\bBOP$ and the additive Grassmannian $\bGr$
are both made from Grassmannians, one should beware of the different nature 
of their structure maps.
There is a variation $\bGr'$ of the additive Grassmannian with values
$\bGr'(V)=\coprod_{n\geq 0} Gr_n(V^2)$ and structure maps
$\bGr'(\varphi)(L)=\varphi^2(L)$. This orthogonal space
is a `multiplicative shift' of $\bGr$ 
in the sense of Example \ref{eg:Additive and multiplicative shift},
it admits a commutative multiplication in much the same way as $\bGr$, and the maps
\[ \bGr(V)\ \to \ \bGr'(V) \ , \quad L \ \longmapsto \ L\oplus 0\]
form a global equivalence of ultra-commutative monoids
by Theorem \ref{thm:general shift of osp}.
A source of possible confusion is the fact that
$\bGr'(V)$ and $\bBOP(V)$ are {\em equal} as spaces, but they come with
very different structure maps making them into two different global homotopy types.

\begin{eg}[$\bGr$ versus $\bBOP$]\label{eg:Gr to BOP}\index{subject}{additive Grassmannian}
In Example \ref{eg:Gr additive} we explained that 
for every compact Lie group $G$, the monoid $\pi_0^G(\bGr)$ 
is isomorphic to the monoid $\bRO^+(G)$, under direct sum, 
of isomorphism classes of orthogonal $G$-represen\-tations. 
In Theorem \ref{thm:upi_0 of BOP} we will identify
the monoid $\pi_0^G(\bBOP)$ with the orthogonal representation ring $\bRO(G)$.
The latter is the algebraic group completion of the former, and this group
completion is realized by a morphism of ultra-commutative monoids 
\begin{equation}\label{eq:Gr2BOP}
i \ : \  \bGr\ \to \ \bBOP \ . 
\end{equation}
The morphism $i$ is given at $V$ by the map
\begin{align*}
   \bGr(V)\ = \ {\coprod}_{m\geq 0}\,  Gr_m( V) \ &\to \ 
 {\coprod}_{n\geq 0}\, Gr_n (V^2) = \bBOP(V)   \ ,\quad L \ \longmapsto \ V\oplus L\ .
\end{align*}
The morphism is homogeneous in that it takes $\bGr^{[m]}$ to $\bBOP^{[m]}$.
\end{eg}

As we will now show, the morphism $i:\bGr\to\bBOP$
induces an algebraic group completion of abelian monoids
upon taking equivariant homotopy sets from any equivariant space.
This fact is the algebraic shadow of a more refined relationship: 
as we will show in Theorem \ref{thm:BOP is group completion} below,
the morphism $i:\bGr\to\bBOP$ is a group completion
in the world of ultra-commutative monoids, i.e.,
`homotopy initial', in the category of ultra-commutative monoids,
among morphisms from $\bGr$ to group-like ultra-commutative monoids.

We recall from Definition \ref{def:define [A,Y]^G}
the equivariant homotopy set
  \[ [A,R]^G \ = \ \colim_{V\in s(\Uc_G)}\, [A, R(V)]^G \ ,\]
where $R$ is  an orthogonal space, 
$G$ is a compact Lie group and $A$ a $G$-space.
If $R$ is an ultra-commutative monoid, then this set inherits an
abelian monoid structure defined as follows. We let
$\alpha:A\to R(V)$ and $\beta:A\to R(W)$ be two $G$-maps
that represent classes in $[A,R]^G$. Then their sum is defined as
\begin{equation}  \label{eq:addition_on_[A,R]^G}
 [\alpha] +[\beta] \ = \ [\mu_{V,W}(\alpha,\beta)] \ ,  
\end{equation}
where $\mu_{V,W}:R(V)\times R(W)\to R(V\oplus W)$ is the
$(V,W)$-component of the  multiplication of $R$.
The monoid structure is contravariantly functorial for $G$-maps in $A$,
and covariantly functorial for morphisms of ultra-commutative monoids in $R$.
When $A=\ast$ is a one-point $G$-space, then $[A,R]^G$ becomes $\pi_0^G(R)$, 
and the addition \eqref{eq:addition_on_[A,R]^G} reduces to the
addition as previously defined in \eqref{eq:internal_product_monoid_space}. 

\begin{prop}\label{prop:i is group completion}
For every compact Lie group $G$ and every $G$-space $A$, the homomorphism
\[ [A,i]^G \ : \ [A,\bGr]^G \ \to \ [A,\bBOP]^G \]
is a group completion of abelian monoids.
\end{prop}
\begin{proof}
We start by showing that the abelian monoid $[A,\bBOP]^G$ is a group.
We consider a $G$-representation $V$.
For a linear subspace $L\subseteq V^2$ 
we consider the 1-parameter family of linear isometric embeddings
\[ H_L \ : \ [0,1]\times L^\perp \ \to \ V^2\oplus V^2 \ , \quad (t,x)\ \longmapsto \
( t\cdot x,  \sqrt{1-t^2}\cdot x ) \ .\]
For every $t\in[0,1]$, the image of $H_L(t,-)$ is isomorphic to $L^\perp$
and orthogonal to the space $L\oplus 0\oplus 0$.
We can thus define a $G$-equivariant homotopy
\begin{align*}
   K \ : \ [0,1]\times G r(V^2) \ \to \ G r(V^2\oplus V^2) \text{\quad by\quad} 
K(t,L) \ =\  (L\oplus 0\oplus 0) +  H_L(t,L^\perp)\ .
\end{align*}
Then
\[ K(0,L)\ =\ (L\oplus 0\oplus 0) + H_L(0,L^\perp) \ = \ 
(L\oplus 0) +(0\oplus L^\perp) \ = \ L\oplus L^\perp \]
and
\begin{align*}
 K(1,L)\ &=\ (L\oplus 0\oplus 0) + H_L(1,L^\perp)\\ 
&=\ (L\oplus 0\oplus 0) + (L^\perp\oplus 0\oplus 0) \ = \ V\oplus V\oplus 0\oplus 0 \ .   
\end{align*}
We recall that the multiplication of $\bBOP$ is given by
\[ \mu_{V,V}^{\bBOP} \ : \ \bBOP(V)\times \bBOP(V)\ \to \bBOP(V\oplus V) \ , \quad
\mu_{V,V}^{\bBOP}(L,L')\ = \ \kappa_{V,V}(L\oplus L')\ , \]
where $\kappa_{V,V}(v,v',w,w')=(v,w,v',w')$.
So the equivariant homotopy $G r(\kappa_{V,V})\circ K$ interpolates between the composite
\[ \bBOP(V)\ \xra{(\Id,(-)^\perp)} \ \bBOP(V)\times \bBOP(V)\ \xra{\mu_{V,V}^{\bBOP}} \
\bBOP(V\oplus V) \]
and the constant map with value
\[ \kappa_{V,V}(V\oplus V\oplus 0\oplus 0)\ = \ V\oplus 0\oplus V\oplus 0 \ . \]
The subspace $V\oplus 0\oplus V\oplus 0$ lies in the same path component
of $\bBOP(V^2)^G= (G r(V^2\oplus V^2))^G$ as the 
subspace $V\oplus V\oplus 0\oplus 0$.
So altogether this shows that the composite $\mu_{V,V}^{\bBOP}\circ (\Id,(-)^\perp)$
is $G$-equivariantly homotopic to the constant map with value
$V\oplus V\oplus 0\oplus 0$.

Now we let $\alpha:A\to \bBOP(V)$ be a $G$-map, representing a class in $[A,\bBOP]^G$.
The composite of $\alpha$ and the orthogonal complement map
$(-)^\perp:\bBOP(V)\to \bBOP(V)$ represents another class in $[A,\bBOP]^G$,
and
\begin{align*}
  [\alpha]+ [(-)^\perp\circ \alpha]\ &= \ 
[\mu^{\bBOP}_{V,V}\circ(\Id, (-)^\perp)\circ \alpha]\ = \ 
[c_{V\oplus V\oplus 0\oplus 0}\circ \alpha]\ = \ 0 
\end{align*}
in the monoid structure of $[A,\bBOP]^G$,
because the subspace $V\oplus V\oplus 0\oplus 0$ 
is the neutral element in $\bBOP(V\oplus V)$.
So the class $[\alpha]$ has an additive inverse, 
and this concludes the proof that the abelian monoid $[A,\bBOP]^G$ is a group.

To show that the homomorphism $[A,i]^G$ is a group completion we 
show two separate statements that amount to the surjectivity, respectively injectivity,
of the extension of $[A,i]^G$ to a homomorphism on the Grothendieck group
of the monoid $[A,\bGr]^G$.

(a) We show that every class in $[A,\bBOP]^G$ is the difference of two
classes in the image of $i_*=[A,i]^G:[A,\bGr]^G\to[A,\bBOP]^G$. 
To see this, we represent a given class $x\in[A,\bBOP]^G$
by a $G$-map $\alpha:A\to\bBOP(V)$, for some $G$-representation $V$.
Because $\bBOP(V)=\bGr(V\oplus V)$, the same map $\alpha$ also represents 
a class in $[A,\bGr]^G$; to emphasize the different role, we write this map as
$\alpha^\sharp:A\to\bGr(V\oplus V)$. 
We let $c_V:A\to \bGr(V)$ denote the constant map with value $V$
and $\chi:V^4\to V^4$ the linear isometry defined by
\[ \chi(v_1,v_2,v_3,v_4)\ = \ (v_2,v_3,v_1,v_4)\ . \]
We observe that the following diagram commutes:
\[ \xymatrix@C=10mm{ 
\bBOP(V)\ar@{=}[d]  \ar[rr]^-{(i(V)\circ c_V,\Id)} &&  
\bBOP(V)\times \bBOP(V)\ar[d]^{\mu_{V,V}^{\bBOP}}\\
\bGr(V\oplus V)\ar[r]_-{i(V\oplus V)} &\bBOP(V\oplus V)\ar[r]_-{G r(\chi)} &
\bBOP(V\oplus V)} \]
Since $\chi$ is equivariantly homotopic, through linear isometries,
to the identity, this shows that the composite
$\mu_{V,V}^{\bBOP}\circ (i(V)\circ c_V,\Id)$ is $G$-equivariantly homotopic
to $i(V\oplus V)$. Thus
\begin{align*}
i_*[c_V] + x \ &= \ [\mu_{V,V}^{\bBOP}\circ(i(V)\circ c_V,\alpha)]\\ &= \ 
[\mu_{V,V}^{\bBOP}\circ(i(V)\circ c_V,\Id)\circ \alpha]\
= \ [i(V\oplus V)\circ \alpha^\sharp]  \ = \ i_*[\alpha^\sharp] \ .
\end{align*}
Thus $x = i_*[\alpha^\sharp] - i_*[c_V]$, which shows the claim.

(b) Now we consider two classes $a,b \in [A,\bGr]^G$ such that
$i_*(a)=i_*(b)$ in $[A,\bBOP]^G$. We show that there exist another
class $c\in [A,\bGr]^G$ such that $c+a=c+b$.
We can represent $a$ and $b$ by two $G$-maps
$\alpha:A\to \bGr(V)$ and $\beta:A\to\bGr(V)$ such that the two composites
\[ i(V)\circ \alpha\ , \quad \ i(V)\circ \beta \ : \ A\ \to \ \bBOP(V) \]
are equivariantly homotopic.
As before we let $c_V:A\to \bGr(V)$ be the constant map with value $V$.
The map $i(V):\bGr(V)\to \bBOP(V)=\bGr(V\oplus V)$
factors as the composite
\[ \bGr(V)\ \xra{(c_V,\Id)}\ \bGr(V)\times\bGr(V)\ \xra{\mu_{V,V}^{\bGr}}\ 
\bGr(V\oplus V) \ ,\]
so
\begin{align*}
[c_V]\ + a \ &= \ [c_V]+[\alpha]\ = \ [\mu_{V,V}^{\bGr}\circ (c_V,\alpha)] \\ 
&= \  [\mu_{V,V}\circ (c_V,\Id)\circ \alpha] \ = \ [i(V)\circ \alpha] \ .
  \end{align*}
Similarly, $[c_V]+ b = [i(V)\circ \beta]$. So $[c_V] + a =[c_V] + b$
in $[A,\bGr]^G$, as claimed.
\end{proof}

Our next aim is to show that the ultra-commutative monoid $\bGr$
represents equivariant vector bundles, and $\bBOP$ 
represents equivariant $K$-theory, at least for compact $G$-spaces.

\begin{construction}
We let $G$ be a compact Lie group and $A$ a $G$-space.
We recall that a {\em $G$-vector bundle} over $A$\index{subject}{equivariant vector bundle}
consists of a vector bundle $\xi:E\to A$ 
and a continuous $G$-action on the total space $E$ such that
\begin{itemize}
\item the bundle projection $\xi:E\to A$ is a $G$-map,
\item for every $g\in G$ and $a\in A$ the map 
$g\cdot-:E_a\to E_{g a}$ is $\mR$-linear.
\end{itemize}
We let $\Vect_G(A)$\index{symbol}{$\Vect_G(A)$ - {monoid of isomorphism classes
of $G$-vector bundles over $A$}} 
be the commutative monoid, under Whitney sum, 
of isomorphism classes of $G$-vector bundles over $A$.
We define a homomorphism of monoids
\begin{equation}\label{eq:Gr_to_Vect}
 \td{-} \ : \ [A,\bGr]^G \ = \ 
\colim_{V\in s(\Uc_G)}\, [A,\, \bGr(V) ]^G  \ \to \  \Vect_G(A)   
\end{equation}
that will turn out to be an isomorphism 
for compact $A$ and that specializes 
to the isomorphism \eqref{eq:pi^G Gr to RO^+ G} from $\pi_0^G(\bGr)$ to $\bRO^+(G)$
when $A$ is a one-point $G$-space.
We let $f:A\to \bGr(V)$ be a continuous $G$-map, for some $G$-representation $V$.
We pull back the tautological $G$-vector bundle $\gamma_V$ over $\bGr(V)$ and obtain
a $G$-vector bundle $f^\star(\gamma_V):E\to A$ over $A$ with total space
\[ E\ = \ \{ (v,a)\in V\times A\ | \ v\in f(a)\} \ . \]
The $G$-action and bundle structure are as 
a $G$-subbundle of the trivial bundle $V\times A$.
Since the base $\bGr(V)$ of the tautological bundle is a disjoint union of
compact spaces, the isomorphism class of the bundle $f^\star(\gamma_V)$ 
depends only on the $G$-homotopy class of $f$, 
see for example \cite[Prop.\,1.3]{segal-equivariant K-theory}. 
So the construction yields a well-defined map
\[ [A,\, \bGr(V) ]^G \ \to \  \Vect_G(A) \ , \quad
[f]\ \longmapsto \ [f^\star(\gamma_V)]\ .\]
If $\varphi:V\to W$ is a linear isometric embedding of $G$-representations, 
then the restriction along $\bGr(\varphi):\bGr(V)\to\bGr(W)$
of the tautological $G$-vector bundle over $\bGr(W)$
is isomorphic to the tautological $G$-vector bundle over $\bGr(V)$.
So the two $G$-vector bundles $f^\star$ and $(\bGr(\varphi)\circ f)^\star$ 
over $A$ are isomorphic. We can thus pass to the colimit over the poset $s(\Uc_G)$
and get a well-defined map \eqref{eq:Gr_to_Vect}.
The map \eqref{eq:Gr_to_Vect} is a monoid homomorphism
because all additions in sight arise from direct sum of inner product spaces.

Now we `group complete' the picture.
We denote by $\bKO_G(A)$\index{symbol}{$\bKO_G(A)$ - {equivariant $K$-group of $A$}}\index{subject}{equivariant $K$-theory}\index{subject}{K-theory@$K$-theory!equivariant|see{equivariant $K$-theory}}
the $G$-equi\-variant $K$-group of $A$, i.e.,
the group completion (Grothendieck group) of the abelian monoid $\Vect_G(A)$.
The composite
\[ [A,\bGr]^G \ \xra{\ \td{-}\ }\ \Vect_G(A)\ \to \ \ \bKO_G(A) \]
of \eqref{eq:Gr_to_Vect} and the group completion map
is a monoid homomorphism into an abelian group. 
The morphism $[A,i]^G: [A,\bGr]^G \to [A,\bBOP]^G$
is a group completion of abelian monoids by Proposition \ref{prop:i is group completion}, 
where $i: \bGr\to\bBOP$ was defined in \eqref{eq:Gr2BOP}.
So there is a unique homomorphism of abelian groups 
\begin{equation}\label{eq:BOP_to_KO_G}
 \td{-} \ : \ [A,\bBOP]^G  \ = \ 
\colim_{V\in s(\Uc_G)}\, [A,\, \bBOP(V) ]^G \ \to \  \bKO_G(A)   
\end{equation}
such that the following square commutes:
\begin{equation}\begin{aligned}\label{eq:bGr to bBOP square}
 \xymatrix@C=20mm{ 
[A,\bGr]^G \ar[d]_{[A,i]^G} \ar[r]^-{\td{-}} &
\Vect_G(A)\ar[d]\\
[A,\bBOP]^G \ar[r]_-{\td{-}} & \bKO_G(A) }     
\end{aligned}\end{equation}
We can make the homomorphism \eqref{eq:BOP_to_KO_G} more explicit as follows.
We let $f:A\to \bBOP(V)$ be a $G$-map for some $G$-representation $V$.
We pull back the tautological $G$-vector bundle over $\bBOP(V)=Gr(V^2)$ and obtain
a $G$-vector bundle $f^\star(\gamma_{V^2}):E\to A$ over $A$ with total space
\[ E\ = \ \{ (v,a)\in V^2 \times A\ | \ v\in f(a)\} \ . \]
Again the $G$-action and bundle structure are as 
a $G$-subbundle of the trivial bundle $V^2\times A$.
The homomorphism \eqref{eq:BOP_to_KO_G} then sends the $G$-homotopy class of $f$
to the class of the virtual $G$-vector bundle
\[ [f^\star(\gamma_{V^2})] - [ V\times A]\ \in \, \bKO_G(A)\ .\]
In contrast to the definition 
of the earlier map \eqref{eq:Gr_to_Vect}, we now subtract the class
of the trivial $G$-vector bundle $V\times A$ over $A$.
The class in $\bKO_G(A)$ only depends on the class of $f$ in $[A,\bBOP]^G$,
and this recipe defines the map \eqref{eq:BOP_to_KO_G}.
\end{construction}

\begin{theorem}\label{thm:BOP_to_KO_G}
For every compact Lie group $G$ and every compact $G$-space $A$ 
the homomorphisms
\[  \td{-} \ : \ [A,\bGr]^G  \ \to \  \Vect_G(A)   \text{\quad and\quad} 
 \td{-} \ : \ [A,\bBOP]^G  \ \to \  \bKO_G(A)    \]
defined in \eqref{eq:Gr_to_Vect} respectively \eqref{eq:BOP_to_KO_G} are isomorphisms.
As $G$ varies the isomorphisms are compatible with
restriction along continuous homomorphism. 
\end{theorem}
\begin{proof}
The Grassmannian $\bGr$ is the disjoint union of the homogeneous pieces $\bGr^{[n]}$,
and the latter is isomorphic to the semifree orthogonal space $\bL_{O(n),\mR^n}$, via
\[ \bL(\mR^n,V)/ O(n) \ \to \ \bGr^{[n]}(V)\ , \quad \varphi\cdot O(n)\ \longmapsto\
\varphi(\mR^n)\ . \]
Since the tautological action of $O(n)$ on $\mR^n$ is faithful,
$\bL_{O(n),\mR^n}$ is a global classifying space for $O(n)$;
Example \ref{eg:[A, B_gl G]} thus provides a bijection
\[ [A, \bGr^{[n]} ]^G \ \to \ \Prin_{(G,O(n))}(A) \]
to the set of isomorphism classes of $G$-equivariant principal $O(n)$-bundles over $A$, 
by pulling back the $(G,O(n))$-principal bundle $\bL(\mR^n,V)\to \bGr^{[n]}(V)$,
the frame bundle of the tautological vector bundle over $\bGr^{[n]}(V)$.
On the other hand, we can consider the map
\[ \Prin_{(G,O(n))}(A) \ \to \ \Vect^{[n]}_G(A)\]
to the set of isomorphism classes of $G$-vector bundles of rank $n$ over $A$,
sending a $(G,O(n))$-bundle $\gamma:E\to A$ to the associated 
$G$-vector bundle with total space $E\times_{O(n)}\mR^n$.
Since $A$ is compact, every $G$-vector bundle admits a $G$-invariant
euclidean inner product, so it arises from a 
$(G,O(n))$-bundle; hence the latter map is bijective as well.
Altogether this shows that map
\begin{equation}  \label{eq:Gr^n_Vect^n}
 [A, \bGr^{[n]}]^G \ \to \ \Vect^{[n]}_G(A)   
\end{equation}
given by pulling back the tautological vector bundles is bijective.

A general $G$-vector bundle need not have constant rank, 
so it remains to assemble the results for varying $n$.
We let $\xi$ be any $G$-vector bundle over $A$, not necessarily
of constant rank.
Then the subset
\[ A_{(n)}\ = \ \{a\in A\ | \ \dim(\xi_a)= n\} \]
of points over which $\xi$ is $n$-dimensional is open
by local triviality of vector bundles.
So $A$ is the disjoint union of the sets $A_{(n)}$ for $n\geq 0$,
and each subset $A_{(n)}$ is also closed and hence compact
in the subspace topology.
Moreover, $A_{(n)}$ is $G$-invariant,
so the restriction $\xi_{(n)}$ of the bundle to $A_{(n)}$ is classified
by a $G$-map $f_{(n)}:A_{(n)}\to \bGr^{[n]}(V_n)$ for some
finite-dimensional $G$-representation $V_n$.
Since $A$ is compact, almost all $A_{(n)}$ are empty, so by increasing
the representations, if necessary, we can assume that the
classifying maps have target in $\bGr^{[n]}(V)$ for a fixed finite-dimensional
$G$-representation $V$, independent of $n$. Then
\[ {\amalg}_{n\geq 0}\ f_{(n)}\ : \ {\amalg}_{n\geq 0} \ A_{(n)}\ =  \ A \ 
\to \  {\amalg}_{n\geq 0}\ \bGr^{[n]}(V) \ = \ \bGr(V)\]
is a classifying $G$-map for the original bundle $\xi$.
This shows that the map \eqref{eq:Gr_to_Vect} is surjective. 

The argument for injectivity is similar.
Any pair of classes in $[A, \bGr]^G$ 
can be represented by $G$-maps $f,\bar f:A\to \bGr(V)$
for some finite-dimensional $G$-representation $V$.
Since $\bGr(V)$ is the disjoint union of the subspaces $\bGr^{[n]}(V)$
for $n\geq 0$, their inverse images under $f$ and $\bar f$
provide disjoint union decompositions of $A$ by fiber dimension.
If the bundles $f^\star(\gamma_V)$ and $\bar f^\star(\gamma_V)$
are isomorphic, the decompositions of $A$ induced by $f$ and $\bar f$
must be the same.
The rank $n$ summands $f_{(n)},\bar f_{(n)}:A_{(n)}\to\bGr^{[n]}(V)$
become equivariantly homotopic after increasing the representation $V$,
because the map \eqref{eq:Gr^n_Vect^n} is injective.
Moreover, almost all summands are empty, one more time by compactness.  
So there is a single finite-dimensional representation $W$
and a $G$-equivariant linear isometric embedding $\varphi:V\to W$
such that $f,\bar f:A\to\bGr(V)$ 
become equivariantly homotopic after composition with
$\bGr(\varphi):\bGr(V)\to \bGr(W)$. Hence $f$ and $\bar f$
represent the same class in $[A,\bGr]^G$, so the map \eqref{eq:Gr_to_Vect} is injective.
This completes the proof that the map is an isomorphism for compact $A$.

The left vertical map in the commutative square \eqref{eq:bGr to bBOP square}
is a group completion by Proposition \ref{prop:i is group completion}, 
and the right vertical map is a group completion by definition.
So the lower horizontal map \eqref{eq:BOP_to_KO_G} is also an isomorphism.
\end{proof}

We take the time to specialize Theorem \ref{thm:BOP_to_KO_G}
to the one-point $G$-space.
This special case identifies the global power monoid
$\upi_0(\bBOP)$ with the global power 
monoid $\bRO$\index{symbol}{$\bRO$ - {orthogonal representation ring global functor}}\index{subject}{representation ring!orthogonal}
of orthogonal representation rings.
For every compact Lie group $G$ the abelian monoid $\bRO(G)$ 
is the Grothendieck group, under direct sum,
of finite-dimensional $G$-representations.
The restriction maps are induced by restriction of representations,
and the power operation $[m]:\bRO(G)\to\bRO(\Sigma_m\wr G)$
takes the class of a $G$-representation $V$ 
to the class of the $(\Sigma_m\wr G)$-representation $V^m$.
The resulting transfer $\tr_H^G:\bRO(H)\to\bRO(G)$
of Construction \ref{con:transfer map}, for $H$ of finite index in $G$,
is then the transfer (or induction), sending
the class of an $H$-representation $V$ 
to the class of the induced $G$-representation $\map^H(G,V)$.
A $G$-vector bundle over a one-point space `is' a $G$-representation and the map
\[ \bRO(G)\ \to \ \bKO_G(\ast) \]
that considers a (virtual) representation as a (virtual) vector bundle is an 
isomorphism of groups and compatible with restriction along continuous homomorphisms
of compact Lie groups.

For easier reference we spell out the isomorphism \eqref{eq:BOP_to_KO_G}
in the special case $A=\ast$ more explicitly.
We let $V$ be a finite-dimensional orthogonal $G$-representation. 
The $G$-fixed points of $\bBOP(V)$ are the $G$-invariant subspaces 
of $V^2$, i.e., the $G$-subrepresentations $W$ of $V^2$.
Representations of compact Lie groups are discrete
(compare the example after \cite[Prop.\,1.3]{segal-equivariant K-theory}), 
so two fixed points in the same path component of $\bBOP(V)^G$
are isomorphic as $G$-representations. Hence we obtain a well-defined map
\[ \pi_0( \bBOP(V)^G) \  \to \ \bRO(G) \]
by sending $W \in \bBOP(V)^G$ to $[W]-[V]$, 
the formal difference in $\bRO(G)$ of the classes of $W$ and $V$. 
These maps are compatible as $V$ runs through 
the finite-dimensional $G$-subrepresentations
of $\Uc_G$, so they assemble into a map
\begin{equation}\label{eq:pi^G BOP to ROG}
\pi_0^G ( \bBOP )\ = \ \colim_{V\in s(\Uc_G)} \, \pi_0( \bBOP(V)^G ) \ \to \ \bRO(G)\ .
\end{equation}
\end{eg}

\begin{theorem}\label{thm:upi_0 of BOP} 
For every compact Lie group $G$
the map \eqref{eq:pi^G BOP to ROG} is an isomorphism of groups. 
As $G$ varies, these isomorphisms form an isomorphism of global power monoids
\[ \upi_0(\bBOP) \ \iso \ \bRO \ .\]
\end{theorem}
\begin{proof}
The special case $A=\ast$ of Theorem \ref{thm:BOP_to_KO_G}
shows that the map \eqref{eq:pi^G BOP to ROG} is an isomorphism
and compatible with restriction along continuous homomorphisms.
We have to argue that in addition, the maps \eqref{eq:pi^G BOP to ROG} 
are also compatible with transfers (or equivalently, with power operations).
The compatibility with transfers can either be deduced directly from
the definitions; equivalently it can be formally deduced from the 
compatibility of the isomorphisms $\upi_0(\bGr) \iso\bRO^+$
with transfers by the universal property of a group completion.
\end{proof}

The bijection \eqref{eq:pi^G BOP to ROG} sends elements of 
$\pi_0^G( \bBOP^{[k]})$ to virtual representations of dimension $k$,
so we can also identify the global power monoid of the homogeneous degree~0
part $\bBO=\bBOP^{[0]}$. Indeed, the map \eqref{eq:pi^G BOP to ROG} restricts
to an isomorphism of abelian groups
\[ \pi_0^G(\bBO) \ \iso \  \bIO(G) \]
to the augmentation ideal $\bIO(G)\subset \bRO(G)$ of the orthogonal representation ring,
compatible with restriction maps, power operations
and transfer maps.\index{subject}{augmentation ideal!of the orthogonal representation ring}

\begin{eg}[Bar construction model $\mathbf B^\circ\bO$]\label{eg:bar construction BO}
Using the functorial bar construction we define another global refinement
$\mathbf B^\circ\bO$ of the classifying space of the infinite orthogonal group.
This ultra-commutative monoid is globally connected, and it
comes with a weak homomorphism to $\bBO$ that `picks out' the path components
of the neutral element in the $G$-fixed point spaces $\bBO(\Uc_G)^G$. 

We define $\mathbf B^\circ\bO$ by applying the bar construction 
(see Construction \ref{con:bar construction})\index{subject}{bar construction!of a topological monoid} 
objectwise to the monoid valued orthogonal 
space $\bO$ of Example \ref{eg:orthogonal group monoid space}.\index{subject}{orthogonal group ultra-commutative monoid}
So the value at an inner product space $V$ is
\[ (\mathbf B^\circ\bO)(V) \ = \  B ( O(V) )\ ,\]
the bar construction of the orthogonal group of $V$.
The structure map of a linear isometric embedding $\varphi:V\to W$ is 
obtained by applying the bar construction to the continuous homomorphism
$\bO(\varphi):\bO(V)\to\bO(W)$.
We make $\mathbf B^\circ\bO$ into an ultra-commutative monoid by endowing 
it with multiplication maps
\[ \mu_{V,W}\ : \ (\mathbf B^\circ\bO)(V) \times (\mathbf B^\circ \bO)(W) 
\ \to \ (\mathbf B^\circ\bO)(V\oplus W)\]
defined as the composite
\[   B ( O(V) ) \times B ( O(W) ) \ \xra{\ \iso \ } \ 
  B(O(V) \times O(W) ) \ \xra{\ B\oplus\ } \ B ( O( V\oplus W ) ) \ , \]
where the first map is inverse to the homeomorphism \eqref{eq:B preserves product}.

Now we let $G$ be a compact Lie group and $V$ an orthogonal $G$-representation.
Taking fixed points commutes with geometric realization
(see Proposition \ref{prop:G-fix preserves pushouts} (iv)) 
and with products, so
\[  ( (\mathbf B^\circ\bO)(V) )^G \ = \ |O(V)^\bullet|^G \ \iso \ 
| (O(V)^G)^\bullet| \ = \  B( O^G(V) ) \ . \]
Taking colimit over the poset $s(\Uc_G)$ gives
\[ ((\mathbf B^\circ\bO)(\Uc_G))^G \ \iso \ 
\colim_{V\in s(\Uc_G)} \,  B( O^G(V) ) \  \iso \  B( O^G(\Uc_G) ) \ 
\iso \ {\prod_{[\lambda]}}'\, B(O^G(\Uc_\lambda)) \ . \]
Here the last weak product is indexed by isomorphism classes of irreducible 
$G$-representations, and each of the groups $O^G(\Uc_\lambda)$ 
is either an infinite orthogonal, unitary or symplectic group, depending
on the type of the irreducible representation,
compare Example \ref{eg:orthogonal group monoid space}.
In particular, the space $ ((\mathbf B^\circ\bO)(\Uc_G))^G $ is connected,
so the equivariant homotopy set $\pi_0^G (\mathbf B^\circ\bO)$ has one
element for every compact Lie group $G$;
the global power monoid structure is then necessarily trivial.
In particular, $\mathbf B^\circ\bO$ is {\em not} globally equivalent to $\bBO$.

However, the difference seen by $\upi_0$ is the only difference between
$\mathbf B^\circ\bO$ and $\bBO$, as shall now explain. We construct a weak morphism of
ultra-commu\-tative monoids
that exhibits $\mathbf B^\circ\bO$ as the `globally connected component' of $\bBO$.
We define an ultra-commutative monoid $\mathbf B'\bO$ by combining the
constructions of $\mathbf B^\circ\bO$ (bar construction) and $\bBO$ (Grassmannians)
into one definition. 
The value of $\mathbf B'\bO$ at an inner product space $V$ is
\[ (\mathbf B'\bO)(V)\ = \ |B_\bullet(\bL(V,V^2),\bO(V),\ast)| \ , \]
the two-sided bar construction (homotopy orbit construction) of the right $O(V)$-action
on the space $\bL(V,V^2)$ by precomposition.
Here $B_\bullet(\bL(V,V^2),\bO(V),\ast)$ is the simplicial space
whose space of $n$-simplices is $\bL(V,V^2)\times \bO(V)^n$.
For $n\geq 1$ and $0\leq i \leq n$, the face map $d_i$
is given by
\[ d_i(\varphi,A_1,\dots A_n) \ = \left\lbrace \begin{array}{ll}
(\varphi\circ A_1, A_2,\dots,A_n) & \mbox{for $i=0$,} \\
(\varphi, A_1,\dots,A_{i-1},A_i\circ A_{i+1},A_{i+2},\dots,A_n)  & \mbox{for $0<i< n$,} \\
(\varphi, A_1,\dots,A_{n-1}) & \mbox{for $i=n$.}
\end{array} \right.  \]
For $n\geq 1$ and $0\leq i \leq n-1$ the degeneracy map  $s_i$
is given by
\[ s_i(\varphi,A_1,\dots, A_{n-1}) \ = \ (\varphi,A_1,\dots,A_i,\Id,A_{i+1},\dots,A_{n-1}) \ . \]
Then $(\mathbf B'\bO)(V)$ is the realization of the
simplicial space $B_\bullet(\bL(V,V^2),\bO(V),\ast)$.

To define the structure map 
associated to a linear isometric embedding $\varphi:V\to W$ we recall that
the structure map $\bO(\varphi):\bO(V)\to\bO(W)$  of the orthogonal space $\bO$
is given by conjugation by $\varphi$ and direct sum with the identity on $W-\varphi(V)$.
We define a continuous map
\[ \varphi_\sharp \ : \ \bL(V,V^2)\ \to \ \bL(W,W^2) \]
by
\[ (\varphi_\sharp \psi)(\varphi(v)+ w)\ = \  \varphi^2( \psi(v) + (w,0) )\ ;\]
here $v\in V$ and $w\in W-\varphi(V)$ is orthogonal to $\varphi(V)$.
The map $\varphi_\sharp$ is compatible with the actions of the orthogonal groups,
i.e., the following square commutes:
\[\xymatrix@C=18mm{ 
\bL(V,V^2)\times \bO(V) \ar[r]^-{\varphi_\sharp\times\bO(\varphi)}\ar[d]_{\circ}&
\bL(W,W^2)\times\bO(W) \ar[d]^{\circ}\\
\bL(V,V^2)\ar[r]_-{\varphi_\sharp} &\bL(W,W^2)} \]
This equivariance property of $\varphi_\sharp$ ensures that it passes to
the two-sided bar construction, i.e., we can define the structure map 
$(\mathbf B'\bO)(\varphi): (\mathbf B'\bO)(V) \to (\mathbf B'\bO)(W)$
as the geometric realization of the morphism of simplicial spaces
\[  B_\bullet(\varphi_\sharp,\bO(\varphi),\ast) \ : \ 
B_\bullet(\bL(V,V^2),\bO(V),\ast) \ \to \ B_\bullet(\bL(W,W^2),\bO(W),\ast) \ . \]
A commutative multiplication
\[ \mu_{V,W}\ : \ (\mathbf B'\bO)(V)\times (\mathbf B'\bO)(W)\ \to \ 
(\mathbf B'\bO)(V\oplus W) \]
is obtained by combining the multiplications of $\mathbf B^\circ\bO$ and $\bBO$.
The construction comes with two collections of continuous maps:
\begin{align*}
  (\mathbf B^\circ\bO)(V)\  \xla{\alpha(V)} \  B(\bL(V,V^2),\bO(V),\ast) &\
= (\mathbf B'\bO)(V)\\ &\xra{\beta(V)}\ \bL(V,V^2)/\bO(V) = \bBO(V) 
\end{align*}
The left map $\alpha(V)$ is defined by applying the bar construction
to the unique map from $\bL(V,V^2)$ to the one-point space.
The right map $\beta(V)$ is the canonical map from homotopy orbits to
strict orbits. As $V$ varies, the $\alpha$ and $\beta$ maps
form morphisms of ultra-commutative monoids
\[
\mathbf B^\circ\bO \ \xleftarrow{\ \alpha\ } \ \mathbf B'\bO
\ \xra{\ \beta\ }\ \bBO \ , 
\]
essentially by construction.
As we shall now see, the morphism $\alpha$ is a global equivalence;
so we can view the chain as a weak morphism of ultra-commuta\-tive
monoids from $\mathbf B^\circ\bO$ to $\bBO$. 
The ultra-commutative monoid $\mathbf B^\circ\bO$ 
is globally connected, whereas $\bBO$ is not, so the morphism $\beta$ cannot
be a global equivalence. However, the second part of the next proposition shows that
it is as close to a global equivalence as it can be.
\end{eg}

\begin{prop}\label{prop:compare BOs}
  \begin{enumerate}[\em (i)]
  \item   The morphism $\alpha:\mathbf B'\bO\to \mathbf B^\circ\bO$ is a global equivalence 
    of ultra-commutative monoids.
  \item
    For every compact Lie group $G$, the morphism $\beta:\mathbf B'\bO\to \bBO$ 
    induces a weak equivalence 
    from the $G$-fixed point space $(\mathbf B'\bO(\Uc_G))^G$
    to the path component of the unit element in 
    the $G$-fixed point space $(\bBO(\Uc_G))^G$.
  \end{enumerate}
\end{prop}
\begin{proof}
(i) We let $V$ be a representation of a compact Lie group $G$ and 
compare the two-sided bar constructions for the $O(V)$-equivariant map
from $\bL(V,V^2)$ to the one-point space
\[  \tilde\alpha(V)\ : \ |B_\bullet (\bL(V,V^2),O(V), O(V))| \ 
\to \ | B_\bullet(\ast,O(V), O(V))| = E O(V) \ .\]
The group $G$ acts by conjugation on $\bL(V,V^2)$ and on $O(V)$. 
The group $O(V)$ acts freely from the right on the
last factor in the bar construction,
and this right $O(V)$-action then commutes with
the $G$-action. The map $\tilde\alpha(V)$ is $(G\times O(V))$-equivariant.
The map $\alpha(V):\mathbf B'\bO(V)\to\bBO(V)$
is obtained from $\tilde\alpha(V)$ by passage to $O(V)$-orbits.
So Proposition \ref{prop:fix of free cofibration}
allows us to analyze and compare the $G$-fixed points of $\alpha(V)$.
Indeed, the proposition shows that the $G$-fixed points of 
$\mathbf B'\bO(V)=|B_\bullet(\bL(V,V^2),O(V), O(V))|/O(V)$
are a disjoint union, indexed by conjugacy classes of continuous homomorphisms
$\gamma:G\to O(V)$ of the spaces
\[   |B_\bullet(\bL(V,V^2),O(V), O(V))| ^{\Gamma(\gamma)} / C(\gamma) \ ,\]
where $\Gamma(\gamma)$ is the graph of $\gamma$.
Since fixed points commute with geometric realization
(see Proposition \ref{prop:G-fix preserves pushouts}~(iv))
and with products, we have
\[  | B(\bL(V,V^2),O(V), O(V)) |^{\Gamma(\gamma)} \ = \ 
 | B_\bullet(\bL^G(V,V^2),O^G(V), \bL^G(\gamma^*(V),V))|\  . \]
If $\gamma^*(V)$ and $V$ are not isomorphic as $G$-representations, then
$\bL^G(\gamma^*(V),V)$ and hence also the bar construction is empty. 
So there is in fact only one summand in the disjoint union
decomposition of $( \mathbf B'\bO(V))^G$, namely the one indexed
by the representation homomorphism $G\to O(V)$
that specifies the given $G$-action on $V$. We conclude that the 
inclusions of $G$-fixed points
\[ \bL^G(V,V^2)\ \to \ \bL(V,V^2)\text{\qquad and\qquad} O^G(V)\ \to \  O(V)  \]
induce a homeomorphism
\[  
 |B_\bullet(\bL^G(V,V^2),O^G(V),\ast)|\ \xra{\ \iso\ } \ 
|B_\bullet(\bL(V,V^2),O(V),\ast)|^G\ = \ ( \mathbf B'\bO(V))^G \ .\]
The same argument identifies the $G$-fixed points of
the bar construction $B( O(V))= \mathbf B^\circ\bO(V)$, and we arrive at a commutative square
\[ \xymatrix{ 
 |B_\bullet(\bL^G(V,V^2),O^G(V),\ast)|\ar[d]_-\iso \ar[r]&
 B(O^G(V))\ar[d]^\iso \\
( \mathbf B'\bO(V))^G  \ar[r]_-{\alpha(V)^G} & (\mathbf B^\circ\bO(V))^G } \]
in which both vertical maps are homeomorphisms.
The space $\bL^G(V,V^2)$ becomes arbitrarily highly connected as $V$ exhausts 
the complete $G$-universe $\Uc_G$.
So the upper horizontal quotient map also becomes 
arbitrarily highly connected as $V$ grows.
Hence the map $\alpha(V)^G$ becomes an equivalence
\[ \tel_i\, \alpha(V_i)^G\ : \ 
\tel_i\, (\mathbf B'\bO(V_i))^G \ \to \ \tel_i\, (\mathbf B^\circ\bO(V_i))^G \]
on the mapping telescopes over an exhaustive sequence $\{V_i\}_{i\geq 1}$
of $G$-represen\-tations.
Fixed points commute with mapping telescopes, so we conclude that the map 
\[ \tel_i\, \alpha(V_i)\ : \ 
\tel_i\, \mathbf B'\bO(V_i) \ \to \ \tel_i\, \mathbf B^\circ\bO(V_i) \]
is a $G$-weak equivalence.
The mapping telescope criterion of Proposition \ref{prop:telescope criterion} 
thus shows that
the morphism $\alpha:\mathbf B'\bO\to \mathbf B^\circ\bO$ is a global equivalence. 

(ii) We let $V$ be a $G$-representation, specified by a continuous
homomorphism $\rho:G\to O(V)$.
We show that the map
\[ \beta(V)^G \ : \  (\mathbf B'\bO(V))^G\ = \  
| B( \bL(V,V^2), O(V),\ast) |^G \ \to \ \bBO(V)^G \]
is a weak equivalence of the source onto the path component of $\bBO(V)^G$
that contains the neutral element of the addition.
The claim then follows by passing to colimits over $V$ in $s(\Uc_G)$.

We showed in part~(i) that the 
inclusions of $G$-fixed points $\bL^G(V,V^2)\to \bL(V,V^2)$
and $O^G(V)\to O(V)$ induce a homeomorphism
\[  
| B_\bullet(\bL^G(V,V^2),O^G(V),\ast)|\ \xra{\ \iso\ } \ 
|B_\bullet(\bL(V,V^2),O(V),\ast)|^G\ = \ ( \mathbf B'\bO(V))^G \ .\]
Since the precomposition action of $O^G(V)$ on $\bL^G(V,V^2)$ is free, 
and $\bL^G(V,V^2)$ is cofibrant as an $O^G(V)$-space,
the homotopy orbits map by a weak equivalence to the strict orbits, 
\[  (\mathbf B'\bO(V))^G \ = \ 
 | B( \bL(V,V^2), O(V),\ast) |^G \ \xra{\ \simeq\ }\ 
 \bL^G(V,V^2) /  O^G(V) \ . \]
The map $\beta(V)^G$ factors as the composite
\[   (\mathbf B'\bO(V))^G \ \xra{\ \simeq\ } \ 
 \bL^G(V,V^2) /  O^G(V) \ \to \ 
 \left( \bL(V,V^2) / O(V) \right)^G \ = \ \bBO(V)^G  \ ,\]
where the second map is induced by the inclusion $\bL^G(V,V^2) \to \bL(V,V^2)$.
The space $\bBO(V)^G$ is the Grassmannian of $G$-invariant subspaces of $V^2$
of the same dimension as $V$, 
and the space $\bL^G(V,V^2) /  O^G(V)$ consists of 
those subspaces that are $G$-isomorphic to $V$.
This is precisely the path component of $\bBO(V)^G$ containing
the neutral element.
\end{proof}

We have now identified the $G$-equivariant path components of the three ultra-commutative
monoids $\mathbf B^\circ\bO$, $\bBO$ and $\bBOP$,
and they are isomorphic to the trivial group, the augmentation ideal $\bIO(G)$
respectively the real representation ring $\bRO(G)$.
Now we determine the entire homotopy types of the
$G$-fixed point spaces of the three ultra-commutative
monoids $\mathbf B^\circ\bO$, $\bBO$ and $\bBOP$.

\begin{cor}\label{cor-components of the BOs} 
Let $G$ be a compact Lie group.
  \begin{enumerate}[\em (i)]
  \item The $G$-fixed point space of $\mathbf B^\circ\bO$ is a
  classifying space of the group $O^G(\Uc_G)$ 
  of $G$-equivariant orthogonal isometries of the complete $G$-universe:
  \[  (\mathbf B^\circ\bO(\Uc_G))^G \ \simeq \  B( O^G(\Uc_G) ) \]
  \item 
    The $G$-fixed point space of $\bBO$ is a disjoint union,
    indexed by the augmentation ideal $\bIO(G)$, \index{subject}{augmentation ideal!of the orthogonal representation ring} 
    of classifying spaces of the group $O^G(\Uc_G)$: 
        \[  (\bBO(\Uc_G))^G \ \simeq \ \bIO(G)\times B( O^G(\Uc_G) ) \]
  \item 
    The $G$-fixed point space of $\bBOP$ is a disjoint union,
    indexed by $\bRO(G)$, of classifying spaces of the group $O^G(\Uc_G)$: 
    \[  (\bBOP(\Uc_G))^G \ \simeq \ \bRO(G)\times B( O^G(\Uc_G) ) \]
  \end{enumerate}
\end{cor}
\begin{proof}
(i) This is almost a tautology.
Since $G$-fixed points commute with geometric realization
(see Proposition \ref{prop:G-fix preserves pushouts}~(iv))
and with products, they commute with the bar construction.
So the space $(\mathbf B^\circ\bO(V))^G$ is homeomorphic to $B (O^G(V))$ 
for every finite-dimensional $G$-representation $V$.
Since $G$-fixed points also commute with the filtered colimit at hand
(see Proposition \ref{prop:G-fix preserves pushouts}~(ii)), we have
\begin{align*}
  (\mathbf B^\circ\bO(\Uc_G))^G \ &= \ \left( \colim_{V\in s(\Uc_G)} \mathbf B^\circ\bO(V) \right)^G
  \\ 
  &\iso \ \colim_{V\in s(\Uc_G)} \left( \mathbf B^\circ\bO(V) \right)^G\\
  &\iso \ \colim_{V\in s(\Uc_G)}  B \left( O^G(V)\right) \\
  &\iso \  B \left(\colim_{V\in s(\Uc_G)}  O^G(V)\right) \ 
  =\ B  \left(  O^G(\Uc_G) \right)  \ .
\end{align*}

(ii) 
As we explained in Remark \ref{rk:umon is global E_infty},
the commutative multiplication of $\bBO$ makes the $G$-space $\bBO(\Uc_G)$ 
into an $E_\infty$-$G$-space; so the fixed points $\bBO(\Uc_G)^G$
come with the structure of a non-equivariant $E_\infty$-space.
The abelian monoid of path components $\pi_0(\bBO(\Uc_G)^G)$
is isomorphic to $\pi_0^G(\bBO)\iso \bIO(G)$, hence an abelian group.
So all path components of the space $\bBO(\Uc_G)^G$
are homotopy equivalent. 
Proposition \ref{prop:compare BOs} identifies the zero path component of $\bBO(\Uc_G)^G$
with $\mathbf B^\circ\bO(\Uc_G)^G$, so part~(i) finishes the proof.

The proof of part~(iii) is the same as for part~(ii), the only difference being
that $\upi_0(\bBOP(\Uc_G)^G)$ is isomorphic to the abelian group $\bRO(G)$.
\end{proof}

The fixed point spaces described in Corollary \ref{cor-components of the BOs} 
can be decomposed even further.
As we explained in Example \ref{eg:orthogonal group monoid space},
the group $O^G(\Uc_G)$ is a weak product of infinite 
orthogonal, unitary and symplectic groups, indexed by the
isomorphism classes of irreducible $G$-representations $\lambda$.
The classifying space construction commutes with weak products, 
which gives a weak equivalence
\[ B (O^G(\Uc_G)) \ \simeq \ {\prod}' B( O^G(\Uc_\lambda))  \ . \]
Moreover, the group $ O^G(\Uc_\lambda)$ is isomorphic to an infinite orthogonal, unitary
or symplectic group, depending on the type of the irreducible representation $\lambda$.

\begin{eg}[More bar construction models]\label{eg:general bar construction}
Since the bar construction is functorial and continuous for 
continuous homomorphisms between topological monoids, we can apply it 
objectwise to every monoid valued orthogonal space $M$ 
in the sense of Definition \ref{def:monoid valued ospace};
the result is an orthogonal space $\mathbf B^\circ M$.\index{subject}{bar construction!of a monoid valued orthogonal space}
The bar construction is symmetric monoidal, so if $M$ is symmetric
(and hence an ultra-commutative monoid),
then $\mathbf B^\circ M$ inherits an ultra-commutative multiplication.
By the same argument as for $\mathbf B^\circ\bO$, the orthogonal space $\mathbf B^\circ M$ 
is globally connected.

So the ultra-commutative monoid $\mathbf B^\circ\bO$
has variations with $\bO$ replaced by $\bSO$, $\bU$, $\bSU$, $\bSp$,
$\bPin$, $\bSpin$, $\bPin^c$ and $\bSpin^c$.\index{subject}{unitary group ultra-commutative monoid}\index{subject}{symplectic group ultra-commutative monoid}\index{subject}{spin group ultra-commutative monoid}\index{subject}{spin$^c$ group ultra-commutative monoid}\index{subject}{pin group orthogonal monoid space}\index{subject}{pin$^c$ group orthogonal monoid space}
We can also apply the bar construction to morphisms of
monoid valued orthogonal spaces, i.e., morphism of orthogonal 
spaces that are objectwise monoid homomorphisms. 
So hitting all the previous examples with the bar construction
yields a commutative diagram of globally connected orthogonal spaces:
\[ \xymatrix@C=15mm{ 
\mathbf B^\circ \bSU \ar[r]^-{\mathbf B^\circ\text{incl}} \ar@{-->}[d]_{\mathbf B^\circ l} & 
\mathbf B^\circ \bU \ar@{-->}[d]_{\mathbf B^\circ l} \ar@{-->}[dr]^-{\mathbf B^\circ r}& \\
\mathbf B^\circ \bSpin \ar[d]_-{\mathbf B^\circ\text{incl}} \ar[r]^-{\mathbf B^\circ\iota} & 
\mathbf B^\circ \bSpin^c \ar[d]_{\mathbf B^\circ\text{incl}}\ar[r]_-{\mathbf B^\circ \ad} 
 & \mathbf B^\circ \bSO \ar[d]^{\mathbf B^\circ\text{incl}}\\ 
 \mathbf B^\circ \bPin \ar[r]_-{\mathbf B^\circ\iota} &  
\mathbf B^\circ \bPin^c \ar[r]_-{\mathbf B^\circ \ad} & 
\mathbf B^\circ \bO} \]
As before, the two dotted arrows mean that the actual morphism goes 
to a multiplicative shift of the target.
With the exception of $\mathbf B^\circ\bPin$ and $\mathbf B^\circ\bPin^c$, 
all these orthogonal spaces inherit ultra-commutative multiplications.
\end{eg}

\begin{eg}\label{eg:define bO}
We define an orthogonal space $\bbO$\index{symbol}{$\bbO$ - {global $BO$}};
its values come with tautological vector bundles whose Thom spaces
form the global Thom spectrum $\bmO$, 
compare Example \ref{eg:geometric global bordism} below.
For finite and abelian compact Lie groups $G$, the equivariant
homotopy groups of $\bmO$ are isomorphic to the bordism groups
of smooth closed $G$-manifolds, compare Theorem \ref{thm:TP is iso};
so $\bbO$ and $\bmO$ are also geometrically relevant.
In Remark \ref{rk:bO versus BO} we define an $E_\infty$-multiplication
on $\bbO$ and show, using power operations, that this $E_\infty$-multiplication
cannot be refined to an ultra-commutative multiplication.

For an inner product space $V$ of dimension $n$ we set
\[  \bbO(V) \ = \  Gr_n(V\oplus\mR^\infty)\ , \]
the Grassmannian of $n$-dimensional subspaces in $V\oplus\mR^\infty$.
The structure map $\bbO(\varphi):\bbO(V)\to \bbO(W)$ is given by
\[ \bbO(\varphi)(L) \ = \ (\varphi\oplus\mR^\infty)(L) + ( (W-\varphi(V))\oplus 0)\ ,\]
the internal orthogonal sum of the image of $L$ under
$\varphi\oplus\mR^\infty:V\oplus\mR^\infty\to W\oplus\mR^\infty$
and the orthogonal complement of the image of $\varphi:V\to W$, viewed as sitting
in the first summand of $W\oplus\mR^\infty$.
\end{eg}

We want to describe the equivariant homotopy sets $\pi_0^G(\bbO)$
and the homotopy types of the fixed point spaces $\bbO(\Uc_G)^G$,
for every compact Lie group $G$.
We denote by $\bRO^\sharp(G)$ the abelian submonoid
of $\bRO^+(G)$ consisting of the isomorphism classes of $G$-representations
with trivial $G$-fixed points.
We let $V$ be a $G$-representation. 
The $G$-fixed points of $\bbO(V)$ are the $G$-subrepresentations $L$
of $V\oplus\mR^\infty$ of the same dimension as $V$.
Since $G$ acts trivially on $\mR^\infty$, the `non-trivial summand' 
$L^\perp=L-L^G$ is contained in $V^\perp=V-V^G$. 
So $V^\perp-L^\perp$ is a $G$-representation with trivial fixed points.
We can thus define a map
\[  (\bbO(V))^G \ = \   \left( Gr_{|V|}(V\oplus\mR^\infty) \right)^G \  \to \ \bRO^\sharp(G) 
\]
from this fixed point space by sending $L \in \bbO(V)^G$ to $[V^\perp-L^\perp]$.
As before, the isomorphism type of $L$ only depends on the path component of $L$ 
in $\bbO(V)^G$. Moreover, for every linear isometric embedding $\varphi:V\to W$
the relation
\begin{align*}
  \left(\bbO(\varphi)(L)\right)^\perp \ &= \  
\left( (\varphi\oplus\mR^\infty)(L) +(W-\varphi(V))\oplus 0 \right)^\perp \\
&= \  
\left(\varphi(L^\perp) + (W^\perp-\varphi(V^\perp))\right)\oplus 0  \
= \  \left( W^\perp- \varphi(V^\perp-L^\perp) \right)\oplus 0  
\end{align*}
shows that
\[ [W^\perp  - \left(\bbO(\varphi)(L)\right)^\perp] \ = \  [\varphi(V^\perp-L^\perp)]
\ = \  [V^\perp-L^\perp] \ .\]
So the class in $\bRO^\sharp(G)$ depends only on the class of $L$
in $\pi_0^G(\bbO)$, and the assignments assemble into a well-defined map
\begin{equation}\label{eq:pi^G bO to RO^sharp}
\pi_0^G ( \bbO )\ = \ \colim_{V\in s(\Uc_G)} \, \pi_0\left( \bbO(V)^G \right) 
\ \to \ \bRO^\sharp(G)\ .
\end{equation}
For the description of the $G$-fixed points of $\bbO$ we introduce the abbreviation
\[ G r_j^{G,\perp}\ = \ \left( G r_j(\Uc_G^\perp) \right)^G \]
for the space of $j$-dimensional $G$-invariant subspaces of $\Uc_G^\perp=\Uc_G-(\Uc_G)^G$.
The space $G r_j^{G,\perp}$ can be decomposed further:
before taking $G$-fixed points, $G r_j(\Uc_G^\perp)$ 
is $G$-equivariantly homeomorphic to $\bL(\mR^j,\Uc_G^\perp)/O(j)$.
So Proposition \ref{prop:fix of free cofibration} provides a decomposition of
$G r_j^{G,\perp}$ as the disjoint union, indexed over conjugacy classes of
continuous homomorphisms $\alpha:G\to O(j)$, of the spaces
\[ \bL^G(\alpha^*(\mR^j),\Uc_G^\perp)/ C(\alpha) \ , \]
where $C(\alpha)$ is the centralizer, in $O(j)$, of the image of $\alpha$.
Conjugacy classes of homomorphism from $G$ to $O(j)$ biject -- by restriction
of the tautological $O(j)$-representation --
with isomorphism classes of $j$-dimensional $G$-representations. 
If $V=\alpha^*(\mR^j)$ is such a $G$-representation, then the space
$\bL^G(V,\Uc_G^\perp)$ is empty if $V$ has non-trivial $G$-fixed points,
and contractible otherwise. Moreover, the centralizer $C(\alpha)$
is precisely the group of $G$-equivariant linear self-isometries of $V$,
which acts freely on $\bL^G(V,\Uc_G^\perp)$.
So if $V^G=0$, then the orbit 
space $\bL^G(V,\Uc_G^\perp)/C(\alpha)$ is 
a classifying space for the group $\bL^G(V,V)=O^G(V)$.
So altogether,
\begin{equation}\label{eq:decompose_Gr_j^G}
 G r_j^{G,\perp}\ \simeq \ \coprod_{[V]\in\bRO^\sharp(G),\ |V|=j}\ B( O^G(V))\ .  
\end{equation}
Every $G$-representation $V$ is the direct sum of its isotypical components $V_\lambda$, 
indexed by the isomorphism classes of irreducible orthogonal $G$-representations.
If $V^G=0$, then only the non-trivial irreducibles occur,
and the group $O^G(V)$ decomposes accordingly as a product
\[ O^G(V)\ \iso \ {\prod}_{[\lambda]} \, O^G(V_\lambda) \ , \]
indexed by non-trivial irreducible $G$-representations. 
The irreducibles come in three flavors (real, complex or quaternionic),
and so the group $O^G(V_\lambda)$ is isomorphic to
one of the groups $O(m)$, $U(m)$, and $S p(m)$, where $m$ 
is the multiplicity of $\lambda$ in $V$.
So altogether, $G r_j^{G,\perp}$ is a disjoint union
of products of classifying spaces of orthogonal, unitary and symplectic groups.

\begin{prop}\label{prop:pi_0 of bbO}
Let $G$ be a compact Lie group.
\begin{enumerate}[\em (i)]
\item 
The $G$-fixed point space $\bbO(\Uc_G)^G$ is weakly equivalent to the space
\[ {\coprod}_{j\geq 0} \,  G r_j^{G,\perp} \times  B O \ .\]
\item
The map \eqref{eq:pi^G bO to RO^sharp} is a bijection from
$\pi_0^G ( \bbO )$ to $\bRO^\sharp(G)$.
\item
If $U$ is a $G$-representation with trivial fixed points,
then the path component of $\bbO(\Uc_G)^G$ indexed by $U$ 
is a classifying space for the group $O^G(U)\times O$.
\end{enumerate}
\end{prop}
\begin{proof}
(i)
We let $V$ be a $G$-representation with $V^G=0$ and $W$ a trivial $G$-representation.
Then every $G$-invariant subspace $L$ of $V\oplus W\oplus \mR^\infty$ is the 
internal direct sum of the fixed part $L^G$ 
(which is contained in $W \oplus \mR^\infty$)
and its orthogonal complement $L^\perp=L-L^G$ (which is contained in the summand $V$).
This canonical decomposition provides a homeomorphism
\[
 \bbO(V\oplus W)^G \ = \ (G r_{|V|+|W|}(V\oplus W\oplus\mR^\infty))^G \ \iso  
\coprod_{j=0,\dots,|V|} \left( G r_j(V) \right)^G \times  G r_{j+|W|}(W\oplus \mR^\infty)  \]
sending $L$ to the pair $(V- L^\perp,\ L^G)$.
Every $G$-invariant subspace of $\Uc_G$ is the direct sum of its fixed
points and their orthogonal complement, so the poset $s(\Uc_G)$
is the product of the two posets $s(\Uc_G^\perp)$ and $s((\Uc_G)^G)$.
We can thus calculate the colimit over $s(\Uc_G)$ in two steps.
For fixed $W$, passing to the colimit over $s(\Uc_G^\perp)$ gives a homeomorphism
\[  \colim_{V\in s(\Uc_G^\perp)} \bbO(V\oplus W)^G \ \iso \ 
{\coprod}_{j\geq 0} \  G r_j^{G,\perp} \times  G r_{j+|W|}(W\oplus\mR^\infty) \ . \]
The factor $ G r_{j+|W|}(W\oplus\mR^\infty)$ is a classifying space
for the group $O(j+|W|)$. Passing to the colimit over $s((\Uc_G)^G)$
then provides a weak equivalence
\begin{align*}
 \bbO(\Uc_G)^G \ &= \ 
\colim_{W\in s( (\Uc_G)^G)} \colim_{V\in s(\Uc_G^\perp)} \bbO(V\oplus W)^G \\ 
&\simeq \ 
\colim_{W\in s( (\Uc_G)^G)} \left( {\coprod}_{j\geq 0} \  G r_j^{G,\perp} \times  B O(j+|W|) \right) \\  
&\simeq \  {\coprod}_{j\geq 0} \  G r_j^{G,\perp} \times  B O \ .
\end{align*}

(ii) Since the orthogonal space $\bbO$ is closed, we can calculate $\pi_0^G(\bbO)$ 
as the set of path components of the space $\bbO(\Uc_G)^G$, by
Corollary \ref{cor:pi_0^G of closed}~(i). 
Since the space $B O$ is path connected, part~(i) 
allows us to identify $\pi_0(\bbO(\Uc_G)^G)$ with
the disjoint union, over $j\geq 0$, of the path components of $G r_j^{G,\perp}$.
As we explained in the weak equivalence \eqref{eq:decompose_Gr_j^G},
the set $\pi_0(G r_j^{G,\perp} )$ bijects with the set of isomorphism classes of
$j$-dimensional $G$-representations with trivial fixed points.
Altogether this identifies $\pi_0(\bbO(\Uc_G)^G)$ with $\bRO^\sharp(G)$,
and unraveling all definitions shows 
that the combined bijection between $\pi_0^G(\bbO)$ and $\bRO^\sharp(G)$
is the map \eqref{eq:pi^G bO to RO^sharp}.

(iii) This is a direct consequence of part~(i) and the
description \eqref{eq:decompose_Gr_j^G} of the path components of the space
$G r_j^{G,\perp} $.
\end{proof}

\Danger If $H$ is a closed subgroup of a compact Lie group $G$ 
and $V$ a $G$-represen\-tation with $V^G=0$,
then $V$ may have non-zero $H$-fixed points.
So the restriction homomorphism $\res^G_H:\bRO^+(G)\to\bRO^+(H)$
does not in general take $\bRO^\sharp(G)$ to $\bRO^\sharp(H)$.
So the monoids $\bRO^\sharp(G)$ do {\em not} form a sub Rep-functor
of $\bRO$, and Proposition \ref{prop:pi_0 of bbO} does not describe
$\upi_0(\bbO)$ as a Rep functor.
We will give a description of $\upi_0(\bbO)$ as sub-Rep monoid
of the augmentation ideal global power monoid $\bIO$ in 
Proposition \ref{prop:pi bbO inside IO} below.\medskip

As we shall now explain, the global homotopy type of the orthogonal space $\bbO$
is that of a sequential homotopy colimit, in the category of orthogonal spaces, 
of global classifying spaces of the orthogonal groups $O(m)$:
\[ \bbO \ \simeq \ \hocolim_{m\geq 1}\, B_{\gl} O(m)  \]
The homotopy colimit is taken over morphisms 
$B_{\gl} O(m)\to B_{\gl} O(m+1)$ that classify the homomorphisms
$O(m)\to O(m+1)$ given by $A\mapsto A\oplus\mR$.
To make this relation rigorous, we define a filtration 
\begin{equation}  \label{eq:bbO_filtration}
 \ast \ \iso \ \bbO_{(0)}\ \subset \ \bbO_{(1)}\ \subset \ \dots\ \subset \ \bbO_{(m)}\ \subset \ \dots   
\end{equation}
of $\bbO$ by orthogonal subspaces. At an inner product space $V$ we define
\begin{equation}  \label{eq:define_bmO_m}
  \bbO_{(m)}(V) \ = \  G r_{|V|}(V\oplus\mR^m)\ ;   
\end{equation}
here we consider $\mR^m$ as the subspace 
of vectors of the form $(x_1,\dots,x_m,0,0,\dots)$
in $\mR^\infty$.
The inclusion of $\bbO_{(m)}$ into $\bbO_{(m+1)}$ is a closed embedding,
so the global invariance property of Proposition \ref{prop:global equiv basics}~(viii)
entitles us to view the union $\bbO$ as a global homotopy colimit of
the filtration.

The tautological action of $O(m)$ on $\mR^m$ is faithful, 
so the semifree orthogonal space $ B_{\gl} O(m) = \bL_{O(m),\mR^m}$
is a global classifying space for $O(m)$. We define a morphism 
$\gamma_m:B_{\gl} O(m)\to \bbO_{(m)}$ by
\begin{align*}
 \gamma_m(V)\ : \ \bL(\mR^m,V)/ O(m) \ &\to \quad G r_{|V|}(V\oplus\mR^m)\ = \ \bbO_{(m)}(V)
\ , \\
 \varphi \cdot O(m)\quad &\longmapsto \ (V-\varphi(\mR^m))\oplus\mR^m \ .
\end{align*}
We omit the straightforward verification that these maps indeed form
a morphism of orthogonal spaces.
The semifree orthogonal space $B_{\gl}O(m)$ comes with a
tautological class $u_{O(m),\mR^m}$, defined in \eqref{eq:tautological_class},
which freely generates the Rep-functor $\upi_0(B_{\gl} O(m))$.
We denote by
\[ u_m \ = \ (\gamma_m)_*( u_{O(m),\mR^m} )\ \in \ \pi_0^{O(m)}\left( \bbO_{(m)} \right) \]
the image  in $\bbO_{(m)}$ of the tautological class.
The following proposition justifies the claim 
that $\bbO$ is a homotopy colimit of the orthogonal spaces $B_{\gl} O(m)$.

\begin{prop}\label{prop:bbO as hocolim}
For every $m\geq 0$ the morphism\index{subject}{global classifying space!of $O(n)$} 
\[ \gamma_m \  : \ B_{\gl} O(m)\ \to \ \bbO_{(m)} \]
is a global equivalence of orthogonal spaces.
The inclusion $\bbO_{(m)} \to \bbO_{(m+1)}$
takes the class $u_m$ to the class
\[ \res^{O(m+1)}_{O(m)}(u_{m+1})\ \in \ \pi_0^{O(m)}\left(\bbO_{(m+1)} \right) \ .\]
\end{prop}
\begin{proof}
  The morphism $\gamma_m$ factors as the composite of two morphisms of
  orthogonal spaces
  \[ \bL_{O(m),\mR^m} \ \xra{\bL_{O(m),\mR^m}\circ i} \ 
  \sh_\oplus^{\mR^m} ( \bL_{O(m),\mR^m})  \ \xra[\iso]{\ (-)^\perp\ } \ \bbO_{(m)}\ .\]
  For the first morphism we let $i:V\to V\oplus\mR^m$ denote the embedding 
  of the first summand, and $\sh_\oplus^{\mR^m}$ is the additive shift by $\mR^m$ 
  as defined in Example \ref{eg:Additive and multiplicative shift};
  the first morphism is a global equivalence 
  by Theorem \ref{thm:general shift of osp}.
  At an inner product space $V$, the second morphism is the map
  \begin{align*}
 (\sh^{\mR^m}_\oplus \bL_{O(m),\mR^m})(V)\ = \ 
  \bL(\mR^m,V\oplus\mR^m)/O(m) \ &\to \ G r_{|V|}(V\oplus\mR^m)\ , \\
  \varphi\cdot O(m)\quad &\longmapsto\quad \varphi(\mR^m)^\perp \ ,    
  \end{align*}
  the orthogonal complement of the image.
  This is a homeomorphism, so the second morphism is an isomorphism.
  Altogether this shows that $\gamma_m$ is a global equivalence.
  
  The second claim is also reasonably straightforward from the definitions,
  but it needs one homotopy. 
  We let $j:\mR^m\to\mR^{m+1}$ denote the linear isometric embedding defined by
  $j(x_1,\dots,x_m)=(x_1,\dots,x_m,0)$.
  The tautological class $u_{O(m),\mR^m}$ is defined as the path component of the point
  \[  \Id_{\mR^m}\cdot O(m) \ \in \ (\bL(\mR^m,\mR^m)/O(m))^{O(m)} \ , \]
  the $O(m)$-orbit of the identity of $\mR^m$.
  So the classes
  \[ \text{incl}_*( u_m ) \text{\qquad and\qquad}
  \res^{O(m)}_{O(m+1)}(u_{m+1}) \text{\quad in\ }
  \pi_0^{O(m)}(\bbO_{(m+1)})\]
  are represented by the two subspaces
  \[ (\mR^{m+1}-j(\mR^m))\oplus j(\mR^m) \text{\qquad respectively\qquad} 
  0\oplus \mR^{m+1} \]
  of $\mR^{m+1}\oplus\mR^{m+1}$. These two representatives are {\em not} the same.
  However, the $O(m+1)$-action on $\bbO_{(m+1)}(\mR^{m+1})$,
  and hence also the restricted $O(m)$-action, 
  is through the first copy of $\mR^{m+1}$
  (and {\em not} diagonally!). So there is a path of
  $O(m)$-invariant subspaces of $\mR^{m+1}\oplus\mR^{m+1}$ connecting the
  two representatives.
  The two points thus represent the same class in $\pi_0^{O(m)}( \bbO_{(m+1)})$,
  and this proves the second claim.
\end{proof}

\begin{rk}[Commutative versus $E_\infty$-orthogonal monoid spaces]\label{rk:bO versus BO}
Non-equivariantly, every $E_\infty$-multiplication on an 
orthogonal monoid space can be rigidified to a strictly commutative multiplication.
We will now see that this is not the case globally,
with power operations being an obstruction.

To illustrate the difference between a strictly commutative multiplication 
and an $E_\infty$-multiplication,
we take a closer look at the orthogonal space $\bbO$.
If we try to define a multiplication on $\bbO$ in a similar way as for $\bBO$, 
we run into the problem that $\mR^\infty\oplus\mR^\infty$
is different from $\mR^\infty$; even worse, although
$\mR^\infty\oplus\mR^\infty$ and $\mR^\infty$ are isometrically isomorphic, there is no
preferred isomorphism. The standard way out is to use {\em all} isomorphisms
at once, i.e., to parametrize the multiplications by the $E_\infty$-operad
of linear isometric self-embeddings of $\mR^\infty$. We recall that the $n$-th
space of the linear isometries operad\index{subject}{linear isometries operad}
is 
\[ \Lc(n) \ = \ \bL((\mR^\infty)^n,\mR^\infty) \ ,\]
with operad structure by direct sum and composition of linear isometric embeddings
(see for example \cite[Def.\,1.2]{may-quinn-ray} for details).
For all $n\geq 0$ and all inner product spaces $V_1,\dots,V_n$
we define a linear isometry
\[ \kappa \ : \ (V_1\oplus\mR^\infty)\oplus\dots\oplus(V_n\oplus\mR^\infty)\ \iso \ 
V_1\oplus\dots\oplus V_n\oplus(\mR^\infty)^n \]
by shuffling the summands, i.e.,
\[ \kappa(v_1,x_1,\dots,v_n,x_n)\ = \ (v_1,\dots,v_n,x_1,\dots,x_n)\ . \]
We can then define a continuous map
\[ \mu_n\ : \ \Lc(n)\times \bbO(V_1)\times\dots\times\bbO(V_n)\ \to \ 
\bbO(V_1\oplus\dots\oplus V_n) \]
by
\[ \mu_n(\varphi,L_1,\dots,L_n)\ = \ 
((V_1\oplus\dots\oplus V_n)\oplus\varphi)\circ \kappa)(L_1\oplus\dots\oplus L_n) \ .\]
For fixed $\varphi$ these maps form a multi-morphism, so the universal
property of the box product produces a morphism of orthogonal spaces 
\[ \mu_n(\varphi,-)\ : \ \bbO\boxtimes\dots\boxtimes\bbO\ \to \ \bbO\ . \]
For varying $\varphi$, these maps define a morphism 
of orthogonal spaces 
\[ \mu_n\ : \ \Lc(n)\times (\bbO\boxtimes\dots\boxtimes\bbO)\ \to \ \bbO\ . \]
As $n$ varies, all these morphism together make the orthogonal space $\bbO$
into an algebra (with respect to the box product)
over the linear isometries operad $\Lc$. Since the linear isometries operad
is an $E_\infty$-operad, we call $\bbO$, endowed with this $\Lc$-action,
an $E_\infty$-orthogonal monoid space.\index{subject}{E infty@$E_\infty$-orthogonal monoid space}

Now we explain why the $E_\infty$-structure on $\bbO$ {\em cannot}
be refined to an ultra-commutative multiplication.
An $E_\infty$-structure gives rise to abelian monoid
structures on the equivariant homotopy sets. In more detail,
we let $R$ be any $E_\infty$-orthogonal monoid space, such as for example $\bbO$.
We obtain binary pairings
\[ \pi_0^G(R)\times \pi_0^G(R)\ \xra{\ \times\ } \
 \pi_0^G(R\boxtimes R)\ \xra{\pi_0^G(\mu_2(\varphi,-))} \ \pi_0^G(R)\ , \]
where $\varphi\in\Lc(2)$ is any linear isometric embedding of
$(\mR^\infty)^2$ into $\mR^\infty$.
The second map (and hence the composite) is independent of $\varphi$
because the space $\Lc(2)$ is contractible.
In the same way as for strict multiplications in \eqref{eq:internal_product_monoid_space},
this binary operation makes $\pi_0^G(R)$ 
into an abelian monoid for every compact Lie group $G$,
such that all restriction maps are homomorphisms.
In other words, the $E_\infty$-structure provides a lift of the
Rep-functor $\upi_0(R)$ to an abelian $\Rep$-monoid, i.e., a functor
\[ \upi_0(R)\ : \ \Rep^{\op} \ \to \ \Ab Mon \ . \]
This structure is natural for homomorphisms of
$E_\infty$-orthogonal monoid spaces.\index{subject}{abelian Rep-monoid}\index{subject}{Rep-monoid!abelian|see{abelian Rep-monoid}}

An ultra-commutative monoid can be viewed as
an $E_\infty$-orthogonal monoid space by letting every element of $\Lc(n)$
act as the iterated multiplication. Equivalently: we let the 
linear isometries operad act along the unique homomorphism to the terminal operad
(whose algebras, with respect to the box product, are the ultra-commutative monoids).
For $E_\infty$-orthogonal monoid spaces
arising in this way from  ultra-commutative monoids, 
the products on $\upi_0$ defined here coincide with those originally defined 
in \eqref{eq:internal_product_monoid_space}.
For ultra-commutative monoids $R$, the abelian $\Rep$-monoid $\upi_0(R)$
is underlying a global power monoid, i.e., it comes with power operations and
transfer maps that satisfy various relations.
We show in Proposition \ref{prop:pi bbO inside IO} below 
that the abelian $\Rep$-monoid $\upi_0(\bbO)$
cannot be extended to a global power monoid whatsoever;
hence $\bbO$ is not globally equivalent, as an $E_\infty$-orthogonal monoid space,
to any ultra-commutative monoid.
A curious fact, however, is that {\em after global group completion}
the $E_\infty$-multiplication of $\bbO$ can be refined to an
ultra-commutative multiplication, compare Remark \ref{rk:completion makes ucom}.
\end{rk}

We compare the $E_\infty$-orthogonal monoid space $\bbO$ to the
ultra-commutative monoid $\bBO$ in the most highly structured way
possible. Every ultra-commuta\-tive monoid can be viewed as
an $E_\infty$-orthogonal monoid space, 
and we now define a `weak $E_\infty$-morphism' from $\bbO$ to $\bBO$.
The zigzag of morphisms passes through the orthogonal space $\bBO'$ with values 
\begin{equation}  \label{eq:define_BO'}
  \bBO'(V) \ = \  Gr_{|V|}(V^2\oplus\mR^\infty)\ .   
\end{equation}
The structure maps of $\bBO'$ are a mixture of those for $\bbO$ and $\bBO$, i.e.,
\[ \bBO'(\varphi)(L) \ = \ (\varphi^2\oplus\mR^\infty)(L) \oplus 
((W-\varphi(V))\oplus 0\oplus 0)\ ,\]
where now the orthogonal complement of the image of $\varphi$
is viewed as sitting in the first summand of $W\oplus W\oplus\mR^\infty$.
The linear isometries operad acts on $\bBO'$ in much the same
way as for $\bbO$, making it an  $E_\infty$-orthogonal monoid space. 
Postcomposition with the direct summand embeddings
\[ V\oplus\mR^\infty \ \xra{(v,x)\mapsto(v,0,x)} \ 
V^2 \oplus\mR^\infty \ \xla{(v,v',0)\mapsfrom(v,v')} V^2  \]
induces maps of Grassmannians
\[ \bbO(V) \ \xra{\ a(V)\ } \ \bBO'(V) \ \xla{\ b(V)\ } \ \bBO(V)  \]
that form morphisms of $E_\infty$-orthogonal monoid spaces
\begin{equation}  \label{eq:a_and_b_morphisms}
 \bbO \ \xra{\ a\ } \ \bBO' \ \xla{\ b\ } \ \bBO   \ .
\end{equation}

\begin{prop}\label{prop:bBO to bBO'} 
The morphism $b:\bBO\to\bBO'$ is a global equivalence of
orthogonal spaces.
\end{prop}
\begin{proof}
We define an exhaustive filtration 
\[  \bBO \ = \ \bBO'_{(0)}\ \subset \ \bBO'_{(1)}\ \subset \ \dots\ \subset \ \bBO'_{(m)}\ \subset \ \dots    \]
of $\bBO'$ by orthogonal subspaces by setting
\[  \bBO'_{(m)}(V) \ = \  G r_{|V|}(V^2\oplus\mR^m)\ . \] 
We denote by $\sh=\sh^\oplus_{\mR}$ the additive shift functor 
defined in Example \ref{eg:Additive and multiplicative shift},
and by $ i_X:X \to \sh X$ the morphism of orthogonal spaces 
given by applying $X$ to the direct summand embeddings $V\to V\oplus\mR$. 
The morphism $i_X$ is a global equivalence 
for every orthogonal space $X$, by Theorem \ref{thm:general shift of osp}.
We define a morphism
\[ j \ : \  \bBO'_{(m+1)} \ \to \ \sh (\bBO'_{(m)}) \]
at an inner product space $V$ by
\begin{align*}
  j(V)\ : \  \bBO'_{(m+1)}(V)\ &= \ G r_{|V|}(V\oplus V\oplus\mR^{m+1})\\ 
&\to \ G r_{|V|+1}(V\oplus \mR\oplus V\oplus\mR \oplus\mR^m)= \sh (\bBO'_{(m)})(V) 
\end{align*}
by applying the linear isometric embedding
\begin{align*}
  V\oplus V\oplus\mR^{m+1}\quad &\to \quad 
V\oplus \mR\oplus V\oplus\mR \oplus\mR^m\\
(v,v',(x_1,\dots, x_{m+1}))\ &\longmapsto \ (v,0,v',x_1,(x_2,\dots, x_{m+1})) 
\end{align*}
and adding the first copy of $\mR$ 
(the orthogonal complement of this last embedding).
Then the following diagram commutes:
\[ \xymatrix@C=15mm{ 
\bBO'_{(m)} \ar[r]^-{\text{incl}} \ar[dr]_{i_{\bBO'_{(m)}}} & 
\bBO'_{(m+1)} \ar[d]_j \ar[dr]^-{i_{\bBO'_{(m+1)}}} \\ 
& \sh ( \bBO'_{(m)}) \ar[r]_-{\sh(\text{incl})} &
\sh ( \bBO'_{(m+1)}) } \]
The two diagonal morphisms are global equivalences, hence so is
the inclusion $\bBO'_{(m)}\to \bBO'_{(m+1)}$, by the 2-out-of-6 property 
for global equivalences (Proposition \ref{prop:global equiv basics}~(iii)).
The inclusion of $\bBO'_{(m)}$ into $\bBO'_{(m+1)}$ is also objectwise a closed embedding,
so the inclusion 
\[ \bBO\ =\ \bBO'_{(0)}\ \to\ {\bigcup}_{m\geq 0} \, \bBO'_{(m)}\ = \ \bBO' \]
is a global equivalence, 
by Proposition \ref{prop:global equiv basics}~(ix).
\end{proof}

We recall that $\bIO(G)$ denotes the augmentation ideal in the real representation
ring $\bRO(G)$ of a compact Lie group $G$.\index{subject}{augmentation ideal!of the orthogonal representation ring}
In the same way as for $\bBO$ we define a monoid homomorphism 
\[ \gamma \ : \ \pi_0^G(\bBO') \ \to \ \bIO(G)\]
by sending the path component of $W\in \bBO'(V)^G$ to the
class $[W]-[V]$.
Then the triangle of monoid homomorphisms
on the right of following diagram commutes:
\[ \xymatrix{
 \pi_0^G(\bbO) \ar[r]^-{a_*} \ar[dr] & \pi_0^G(\bBO') \ar[d]_\gamma ^\iso& 
\pi_0^G(\bBO)  \ar[l]^-{\iso}_-{b_*} \ar[dl]_(.6){\iso}^{\eqref{eq:pi^G BOP to ROG}} \\
&  \bIO(G)  } \]
The two maps on the right are isomorphisms, hence so is the map $\gamma$.

The following proposition says that $\pi_0^G(\bbO)$ is isomorphic to the
free abelian submonoid of $\bIO(G)$ generated by
$\dim(\lambda)\cdot 1-[\lambda]$ as $\lambda$ runs over all
isomorphism classes of non-trivial irreducible $G$-representations.
We emphasize that for non-trivial groups $G$,
the monoid $\pi_0^G(\bbO)$ does not have inverses,
so $\bbO$ is {\em not} group-like.

\begin{prop}\label{prop:pi bbO inside IO}
The composite morphism of abelian $\Rep$-monoids
\[   \upi_0(\bbO) \ \xra{\ a_*\ } \ \upi_0(\bBO') \ \xra{\ \gamma\ } \bIO \]
is a monomorphism. For every compact Lie group $G$ the image of
$\pi_0^G(\bbO)$ in the augmentation ideal $\bIO(G)$ consists of
the submonoid of elements of the form $\dim(U)\cdot 1- [U]$, for $G$-representations $U$.
The abelian $\Rep$-monoid $\upi_0(\bbO)$ 
cannot be extended to a global power monoid.\index{subject}{augmentation ideal!of the orthogonal representation ring}
\end{prop}
\begin{proof}
If $L\subset V\oplus\mR^\infty$ is a $G$-invariant subspace of the same dimension
as $V$, then
\begin{align*}
 [L] - [V]  \ &= \  
 (\dim(L)-\dim(L^\perp))\cdot 1  + [L^\perp]  - 
 (\dim(V)-\dim(V^\perp))\cdot 1  - [V^\perp]\\
 &= \  
 \dim(V^\perp-L^\perp)\cdot 1 - [V^\perp-L^\perp]  
\end{align*}
in the group $\bIO(G)$.
This show that the following square commutes:
\[ \xymatrix@C=20mm{ 
  \pi_0^G(\bbO) \ar[r]^-{a_*} \ar[d]_{\eqref{eq:pi^G bO to RO^sharp}} & 
\pi_0^G(\bBO') \ar[d]^\gamma\\
\bRO^\sharp(G)\ar[r]_{[U]\ \mapsto\ \dim(U)\cdot 1- [U]} & \bIO(G) } \]
The left vertical map is bijective by Proposition \ref{prop:pi_0 of bbO}~(ii).
The right vertical map is an isomorphism as explained above. 
This shows the first two claims because the lower horizontal map is injective and
has the desired image.

Now  we argue, by contradiction, that the abelian $\Rep$-monoid $\upi_0(\bbO)$ 
cannot be extended to a global power monoid.
The additional structure would in particular specify a transfer map
\begin{equation}\label{eq:tr A_3 to Sigma_3}
\tr_{A_3}^{\Sigma_3}\ : \  \pi_0^{A_3}(\bbO)\ \to \ \pi_0^{\Sigma_3}(\bbO) 
\end{equation}
from the alternating group $A_3$ to the symmetric group $\Sigma_3$.
Since the monoid $\pi_0^e(\bbO)$ has only one element,
the double coset formula shows that
\[ \res^{\Sigma_3}_{(1 2)}\circ \tr_{A_3}^{\Sigma_3}\ =\ \tr_e^{(1 2)}\circ \res^{A_3}_e\ = \ 0 \
 : \  \pi_0^{A_3}(\bbO)\ \to \ \pi_0^{(1 2)}(\bbO) \  .\]
The group $\Sigma_3$ has two non-trivial irreducible orthogonal representations,
the 1-dimension sign representation $\sigma$ and
the 2-dimensional reduced natural representation $\nu$.
So $\pi_0^{\Sigma_3}(\bbO)$ `is' (via $\gamma\circ a_*$)
the free abelian submonoid of $\bIO(\Sigma_3)$ generated by
\[ 1 -\sigma \text{\qquad and\qquad} 2 - \nu \ .\]
We abuse notation and also write $\sigma$ for the 1-dimensional
sign representation of the cyclic subgroup of $\Sigma_3$ generated by
the transposition $(1 2)$.
Then
\[ \res^{\Sigma_3}_{(12)}(1 -\sigma) \ = \ 1 -\sigma \text{\quad and\quad}
 \res^{\Sigma_3}_{(12)}(2 -\nu) \ = \ 2 - ( 1 +\sigma)\ = \ 1 - \sigma\ ;\]
so the only element of $\pi_0^{\Sigma_3}(\bbO)$ that restricts to~0 in
$\pi_0^{(1 2)}(\bbO)$ is the zero element. 
Hence the transfer map \eqref{eq:tr A_3 to Sigma_3} must be the zero map.
However, another instance of the double coset formula is
\[ \res^{\Sigma_3}_{A_3}\circ \tr_{A_3}^{\Sigma_3}\ =\ \Id\, +\, (c_{(1 2)})^* \
= \ 2\cdot \Id  \  : \  \pi_0^{A_3}(\bbO)\ \to \ \pi_0^{A_3}(\bbO) \  .\]
The second equality uses that conjugation by the transposition $(1 2)$ 
is the non-trivial automorphism of $A_3$,
which acts trivially on $\bRO(A_3)$. 
Since $\bbO(A_3)$ is a non-trivial free abelian monoid, 
this contradicts the relation $\tr_{A_3}^{\Sigma_3}=0$.
So the abelian $\Rep$-monoid $\upi_0(\bbO)$ cannot be endowed with transfer maps
that satisfy the double coset formulas.
\end{proof}

\begin{eg}[Periodic $\bbO$]\label{eg:define bOP}
The $E_\infty$-orthogonal monoid space $\bbO$
also has a straightforward `periodic' variant $\bbOP$ that we briefly discuss.
\index{symbol}{$\bbOP$ - {periodic global $BO$}}
For an inner product space $V$ we set
\[  \bbOP(V) \ = \  {\coprod}_{n\geq 0}\, Gr_n(V\oplus\mR^\infty)\ , \]
the disjoint union of all the Grassmannians in $V\oplus\mR^\infty$.
The structure maps of $\bbOP$ are defined in exactly the same way as for $\bbO$.
The orthogonal space $\bbOP$ is naturally $\mZ$-graded: for $m\in\mZ$ we let
\[ \bbOP^{[m]}(V)\ \subset \ \bbOP(V) \]
be the path component consisting of all subspaces $L\subset V\oplus\mR^\infty$
such that $\dim(L)=\dim(V) + m$. For fixed $m$ these spaces
form an orthogonal subspace $\bbOP^{[m]}$ of $\bbOP$. 
The $E_\infty$-multiplication of $\bbO$ extends naturally
to a $\mZ$-graded $E_\infty$-multiplication on $\bbOP$, 
taking  $\bbOP^{[m]} \boxtimes \bbOP^{[n]}$ to $\bbOP^{[m+n]}$.
Moreover, $\bbO=\bbOP^{[0]}$, the homogeneous summand of $\bbOP$ of degree~0.
\end{eg}

We offer two descriptions of the global homotopy type of the orthogonal space $\bbOP$
in terms of other global homotopy types previously discussed.
As we show in Proposition \ref{prop:decompose bOP} below, 
each of the homogeneous summands $\bbOP^{[m]}$
is in fact globally equivalent to the degree~0 summand $\bbO$, and hence
\[ \bbOP \ \simeq \ \mZ\times \bbO \]
globally as orthogonal spaces.
When combined with Proposition \ref{prop:pi_0 of bbO}, this
yields a description of the homotopy types of the fixed point spaces $\bbOP(\Uc_G)^G$.

To compare the different summands of $\bbOP$
we choose a linear isometric embedding $\psi:\mR^\infty\oplus\mR\to\mR^\infty$
and define an endomorphism $\psi_\sharp:\bbOP\to\bbOP$ of the orthogonal space $\bbOP$
at an inner product space $V$ as the map
\[ \psi_\sharp(V)\ : \  \bbOP(V) \ \to \bbOP(V) \ , \quad L \ \longmapsto \
(V\oplus\psi)( L\oplus\mR )\ . \]
The morphism $\psi_\sharp$ increases the dimension of subspaces by~1, so
it takes the summand $\bbOP^{[k]}$ to the summand $\bbOP^{[k+1]}$.
Any two linear isometric embeddings from 
$\mR^\infty\oplus\mR$ to $\mR^\infty$ are homotopic through
linear isometric embeddings, so the homotopy class of $\psi_\sharp$
is independent of the choice of $\psi$.

\begin{prop}\label{prop:decompose bOP}
For every linear isometric embedding $\psi:\mR^\infty\oplus\mR\to\mR^\infty$
the morphism of orthogonal spaces  
$\psi_\sharp:\bbOP \to \bbOP$ is a global equivalence.
Hence for every $m\in\mZ$ the restriction is a global equivalence
\[ \psi_\sharp \ : \ \bbOP^{[m]} \ \to \bbOP^{[m+1]} \ .\]
\end{prop}
\begin{proof}
We let $\sh=\sh^\mR_\oplus$ denote the additive shift of an orthogonal space
as defined in Example \ref{eg:Additive and multiplicative shift},\index{subject}{shift!of an orthogonal space!additive}
and $i:\Id\to-\oplus\mR$ is the natural transformation given by the embedding
of the first summand.
The morphism $\psi_\sharp$ factors as the composite
\[ \bbOP \ \xra{\bbOP\circ i} \ \sh \bbOP \ \xra{\ \psi_! \ }\ \bbOP\ . \]
The second morphism is defined at $V$ as the map
\[ (\sh\bbOP)(V)\ = \ \bbOP(V\oplus\mR)\ \to \ \bbOP(V) \]
that sends a subspace of $V\oplus\mR\oplus\mR^\infty$ to its image under 
the linear isometric embedding
\[  V\oplus\mR\oplus\mR^\infty \ \to \ V\oplus\mR^\infty\ , \quad
(v,x,y)\ \longmapsto \ (v,\psi(y,x))\ . \]

The morphism $\bbOP\circ i: \bbOP \to\sh \bbOP$
is a global equivalence by Theorem \ref{thm:general shift of osp}.  
Any two linear isometric embeddings from $\mR^\infty\oplus\mR$ to $\mR^\infty$
are homotopic through linear isometric embeddings;
in particular, the linear isometric embedding $\psi$
is homotopic to a linear isometric isomorphism.
Thus $\psi_!$ is homotopic to an isomorphism, hence a global equivalence.
Since both $\bbOP\circ i$ and $\psi_!$ are global equivalences, 
so is the composite $\psi_\sharp$.
\end{proof}

Another description of the global homotopy type of $\bbOP$ 
is as a global homotopy colimit of a sequence of self-maps of $\bGr$:
\[ \bbOP \ \simeq \ \hocolim_{m\geq 1}\, \bGr\ .     \]
The homotopy colimit is taken over iterated instances of a morphism 
$\bGr\to\bGr$ that classifies the map `add a trivial summand $\mR$'.
To make this relation rigorous, we use the filtration 
\[ 
 \ast \ \iso \ \bbOP_{(0)}\ \subset \ \bbOP_{(1)}\ \subset \ \dots\ \subset \ \bbOP_{(m)}\ \subset \ \dots   
 \]
of $\bbOP$ by orthogonal subspaces, 
analogous to the filtration \eqref{eq:bbO_filtration} for $\bbO$.
In other words,
\[  \bbOP_{(m)}(V) \ = \  {\coprod}_{n\geq 0} \, G r_n(V\oplus\mR^m)\ ; \] 
as before we consider $\mR^m$ as the subspace of $\mR^\infty$
of all vectors of the form $(x_1,\dots,x_m,0,0,\dots)$.

We write $\sh^m$ for $\sh^{\mR^m}_\oplus$, the additive shift by $\mR^m$
in the sense of Example \ref{eg:Additive and multiplicative shift}.
We define a morphism 
$\gamma_m:\sh^m \bGr\to \bbOP_{(m)}$ by
\begin{align*}
  \gamma_m(V)\ : \ (\sh^m\bGr)(V)\ = \ \bGr(V\oplus\mR^m) \ &\to \  \bbOP_{(m)}(V)\\
  L \qquad &\longmapsto \ L^\perp  = (V\oplus\mR^m) - L\ ,
\end{align*}
the orthogonal complement of $L$ inside $V\oplus\mR^m$.
We omit the straightforward verification that these maps indeed form
a morphism of orthogonal spaces. Since each map $\gamma_m(V)$ is a homeomorphism,
the morphism $\gamma_m$ is in fact an isomorphism of orthogonal spaces.

The morphism
\[ i_m \ : \ \bGr \ \to \ \sh^m \bGr \]
induced by the direct summand embedding $V\to V\oplus\mR^m$
is a global equivalence by Theorem \ref{thm:general shift of osp}.
The inclusion of $\bbOP_{(m)}$ into $\bbOP_{(m+1)}$ is a closed embedding,
so the global invariance property of Proposition \ref{prop:global equiv basics}~(ix)
entitles us to view the union $\bbOP$ as a global homotopy colimit of
the filtration.
This justifies the interpretation of $\bbOP$ as a global homotopy colimit
of a sequence of copies of $\bGr$.
For this description to be useful we should
identify the global homotopy classes of the morphisms in the sequence, i.e.,
the weak morphisms
\[ \bGr \ \xra[\sim]{\ \gamma_m\circ i_m\ } \ \bbOP_{(m)}\ \xra{\ \text{incl}\ } \
\bbOP_{(m+1)}\ \xla[\sim]{\ \gamma_{m+1}\circ i_{m+1}\ } \ \bGr \ . \]
As is straightforward from the definition, this weak morphism models
`adding a summand $\mR$ with trivial action'.

Finally, we compare the $E_\infty$-orthogonal monoid space $\bbOP$ to the
ultra-commutative monoid $\bBOP$.
In analogy with the non-periodic version in \eqref{eq:define_BO'},
we introduce the orthogonal space $\bBOP'$ with values 
\[ \bBOP'(V) \ = \ {\coprod}_{n\geq 0} \, Gr_n(V^2\oplus\mR^\infty)\ .    \]
The structure maps of $\bBOP'$ are defined in the same way as for $\bBO'$,
and they mix the structure maps for $\bbOP$ and $\bBOP$.
Postcomposition with the direct summand embeddings
\[ V\oplus\mR^\infty \ \xra{(v,x)\mapsto(v,0,x)} \ 
V^2 \oplus\mR^\infty \ \xla{(v,v',0)\mapsfrom(v,v')} V^2  \]
induces morphisms of $E_\infty$-orthogonal monoid spaces
\[ \bbOP \ \xra{\ a\ } \ \bBOP' \ \xla{\ b\ } \ \bBOP \ . \]
These morphisms preserve the $\mZ$-grading;
the restrictions to the homogeneous degree~0 summand are precisely
the morphisms with the same names introduced in \eqref{eq:a_and_b_morphisms}.
The same argument as in Proposition \ref{prop:bBO to bBO'}
also shows that the morphism $b:\bBOP\to\bBOP'$ is a global equivalence of
orthogonal spaces.

We define a monoid homomorphism 
\[ \gamma \ : \ \pi_0^G(\bBOP') \ \to \ \bRO(G)\]
by sending the path component of $W\in \bBOP'(V)^G$ to the
class $[W]-[V]$.
Then the triangle of monoid homomorphisms
on the right of following diagram commutes:
\[ \xymatrix{
 \pi_0^G(\bbOP) \ar[r]^-{a_*} \ar[dr] & \pi_0^G(\bBOP') \ar[d]_\gamma ^\iso& 
\pi_0^G(\bBOP)  \ar[l]^-{\iso}_-{b_*} \ar[dl]_(.6){\iso}^{\eqref{eq:pi^G BOP to ROG}} \\
&  \bRO(G)  } \]
The two maps on the right are isomorphisms, hence so is the map $\gamma$.
The same argument as in Proposition \ref{prop:pi bbO inside IO}
shows that the composite morphism of abelian Rep-monoids
$\gamma\circ a_*:\upi_0(\bbOP) \to \bRO$ is a monomorphism, and
for every compact Lie group $G$ the image of
$\pi_0^G(\bbOP)$ in the representation ring $\bRO(G)$ consists of
the submonoid of elements of the form $n\cdot 1- [U]$, for $n\in\mZ$ and $U$ 
any $G$-representation.

\begin{eg}[Complex and quaternionic periodic Grassmannians]\label{eg:BUP and BSpP}
The ultra-commutative monoids $\bBO$ and $\bBOP$ and the $E_\infty$-orthogonal monoid
spa\-ces $\bbO$ and $\bbOP$ have straightforward complex and quaternionic analogs;
we quickly give the relevant definitions for the sake of completeness.
We define the periodic Grassmannians $\bBUP$ and $\bBSpP$ by\index{symbol}{$\bBUP$ - {periodic global $BU$}}\index{subject}{periodic global $BU$}
\index{symbol}{$\bBSpP$ - {periodic global $BSp$}}\index{subject}{periodic global $BSp$}\index{symbol}{$\bBU$ - {global $BU$}}\index{symbol}{$\bBSp$ - {global $BSp$}}
\[  \bBUP(V) \ = \  {\coprod}_{m\geq 0} \, Gr_m^\mC(V_\mC^2)\text{\quad respectively\quad}
 \bBSpP(V) \ = \  {\coprod}_{m\geq 0} \, Gr_m^\mH(V_\mH^2)\ , \]
the disjoint union of the respective Grassmannians,
with structure maps  as for $\bBOP$. External direct sum
of subspaces defines a $\mZ$-graded ultra-commu\-tative multiplication
on $\bBUP$ and on $\bBSpP$, again in much the same way as for $\bBOP$.
The homogeneous summands of degree zero are closed under the multiplication
and form ultra-commutative monoids $\bBU=\bBUP^{[0]}$ 
respectively $\bBSp=\bBSpP^{[0]}$. More explicitly,
\[  \bBU(V) \ = \   Gr_{|V|}^\mC(V_\mC^2)\text{\qquad respectively\qquad}
 \bBSp(V) \ = \  Gr_{|V|}^\mH(V_\mH^2)\ . \]
The complex and quaternionic analogues of Theorem \ref{thm:upi_0 of BOP} 
provide isomorphisms of global power monoids
\[ \upi_0(\bBUP) \ \iso \ \bRU \text{\qquad and\qquad}
\upi_0(\bBSpP) \ \iso \ \bRSp \ ; \]
here $\bRU(G)$ and $\bRSp(G)$ are the Grothendieck groups,
under direct sum, of isomorphism classes of unitary respectively symplectic
$G$-representations. The isomorphisms above match the $\mZ$-grading of
$\bBUP$ and $\bBSpP$ with the grading by virtual dimension of representations,
so they restrict to isomorphisms of global power monoids from
$\upi_0(\bBU)$ respectively $\upi_0(\bBSp)$ to the augmentation ideal
global power monoids inside $\bRU$ respectively $\bRSp$.

Theorem \ref{thm:BOP_to_KO_G} also generalizes to natural group isomorphisms,
compatible with restrictions,
\[  \td{-} \ : \ \bBUP_G(A)  \ \to \  \bKU_G(A)  \text{\quad and\quad}
  \td{-} \ : \ \bBSpP_G(A)  \ \to \  \bKSp_G(A)    \]
to the equivariant unitary respectively symplectic $K$-groups,
where $G$ is any compact Lie group and $A$ a compact $G$-space.
\index{subject}{equivariant $K$-theory}

Example \ref{eg:define bO} can be modified to define $E_\infty$-orthogonal monoid spaces
$\bbU$\index{symbol}{$\bbU$ - {global $B U$}} and $\bbSp$\index{symbol}{$\bbSp$ - {global $B S p$}} with values
\[  \bbU(V) \ = \  Gr_{|V|}^\mC(V_\mC\oplus\mC^\infty)
\text{\quad respectively\quad}
 \bbSp(V) \ = \  Gr_{|V|}^\mH(V_\mH\oplus\mH^\infty)\ .\]
The structure maps and $E_\infty$-multiplication are defined as for $\bbO$.
As orthogonal spaces, $\bbU$ and $\bbSp$ are global homotopy colimits
of the sequence of global classifying spaces $B_{\gl}U(m)$ respectively $B_{\gl}S\!p(m)$.
Periodic versions $\bbUP$\index{symbol}{$\bbUP$ - {periodic global $B U$}} 
and $\bbSpP$\index{symbol}{$\bbSpP$ - {periodic global $B S p$}}
are defined by taking the full Grassmannian
inside $V_\mC\oplus\mC^\infty$ respectively $V_\mH\oplus\mH^\infty$,
as in the real case in Example \ref{eg:define bOP}.
As orthogonal spaces, the periodic versions $\bbUP$ and $\bbSpP$ 
are global homotopy colimits
of iterated instances of the self-maps of $\bGr^\mC$ respectively $\bGr^\mH$
that represent `adding a trivial 1-dimensional representation'.
\end{eg}

\section{Global group completion and units}
\label{sec:group completion and units}

For every orthogonal monoid space $R$ and every compact Lie group $G$,
the operation \eqref{eq:internal_product_monoid_space}
makes the equivariant homotopy set $\pi_0^G(R)$ into a monoid,
and this multiplication is natural with respect to restriction maps in $G$.
If the multiplication of $R$ is commutative,
then so is the multiplication of $\pi_0^G(R)$.
In this section we look more closely at the {\em group-like} ultra-commutative 
monoids, i.e., the ones where all these monoid structures have inverses.
There are two universal ways to make an ultra-commutative monoid group-like:
the `global units' (Construction \ref{con:R^times})
are universal from the left and
the `global group completion' (Construction \ref{con:R^star}
and Corollary \ref{cor:umon Omega Sigma is group completion})
is universal from the right.
In the homotopy category of ultra-commutative monoids, these constructions
are right adjoint respectively left adjoint to the inclusion of group-like objects.
On $\pi_0^G$ both constructions have the expected effect:
the global units pick out the invertible elements of $\pi_0^G$
(see Proposition \ref{prop:realize units in umon}),
and the effect of global group completion is group completion of
the abelian monoids $\pi_0^G$
(see Proposition \ref{prop:group completion completes pi_0}).
Naturally occurring examples of global group completions 
are the morphism $i:\bGr\to\bBOP$ from the additive Grassmannian to
the periodic global version of $B O$, and its complex and quaternionic
versions, see Theorem \ref{thm:BOP is group completion}.
At the end of this section we use global group completion to prove
a global, highly structured version of Bott periodicity:
Theorem \ref{thm:global Bott periodicity}
shows that $\bBUP$ is globally equivalent, as an
ultra-commutative monoid, to $\Omega\bU$.

\medskip

The category of ultra-commutative monoids is pointed,
and product respectively box product are the
categorical product respectively coproduct 
in the category of ultra-commutative monoids.
These descend to product and respectively coproduct in the 
homotopy category $\Ho(umon)$ of ultra-commutative monoids,
with respect to the global model structure of Theorem \ref{thm:global umon}. 
The morphism $\rho_{R,S}:R\boxtimes S\to R\times S$ is a global
equivalence by Theorem \ref{thm:box to times}~(i),
so in $\Ho(umon)$ the canonical morphism from a coproduct to a product is an isomorphism.
Various features of units and group completions only depend on these
formal properties, and work just as well in any pointed model
category in which coproducts and products coincide up to weak equivalence.
So we develop large parts of the theory in this generality.

\begin{construction}\label{con:hom addition} 
Let $\Dc$ be a category which has finite products and a zero object.
We write $A\times B$ for any product of the objects $A$ and $B$ and
leave the projections $A\times B\to A$ and $A\times B\to B$ implicit.
Given morphisms $f:T\to A$ and $g:T\to B$ we write
$(f,g):T\to A\times B$ for the unique morphism that projects to $f$
respectively $g$. We write~0 for any morphism that factors through a zero object.

We call the category $\Dc$ {\em pre-additive}\index{subject}{pre-additive category}
if `finite products are coproducts';
more precisely, we require that every product $A\times B$
of two objects $A$ and $B$ is also a {\em co-}product, with respect to 
the morphisms 
\[ i_1\ =\ (\Id_A,0)\ :\ A\ \to\  A\times B \text{\qquad and\qquad}
i_2\ = \ (0,\Id_B)\ :\ B\ \to \ A\times B\ . \]
In other words, we demand that for every object $X$ the map
\[ \Dc(A\times B,X) \ \to \ \Dc(A,X) \times \Dc(B,X) \ , \quad
f\mapsto (f i_1,f i_2)  \]
is bijective. The main example we care about is $\Dc=\Ho(umon)$,
the homotopy category of ultra-commutative monoids.

In this situation we can define a binary operation on the
morphism set $\Dc(A,X)$ for every pair of objects $A$ and $X$.
Given morphisms $a,b:A\to X$ we let $a\bot b:A\times A\to X$
be the unique morphism such that $(a\bot b)i_1=a$ and $(a\bot b)i_2=b$.
Then we define 
\[  a+b\ = \ (a\bot b)\Delta \ :\ A\ \to \ X \ , \] 
where $\Delta=(\Id_A,\Id_A):A\to A\times A$ is the diagonal morphism.
\end{construction}

The next proposition is well known, but I do not know a convenient reference.

\begin{prop}\label{prop:additive category}
  Let $\Dc$ be a pre-additive category.
    For every pair of objects $A$ and $X$ of $\Dc$
    the binary operation $+$ makes the set
    $\Dc(A,X)$ of  morphisms into an abelian monoid
    with the zero morphism as neutral element.
    Moreover, the monoid structure is natural 
    for all morphisms in both variables, or, equivalently, 
    composition is biadditive.
\end{prop}
\begin{proof} The proof is lengthy, but completely formal.
For the associativity of `+' we consider three morphisms $
a,b,c:A\to X$. Then $a+(b+c)$ respectively $(a+b)+c$ are the two 
outer composites around the diagram:
\[ \xymatrix@R=5mm{ & A \ar[dl]_\Delta \ar[dr]^\Delta \\
A\times A \ar[d]_{\Id\times\Delta} && A\times A\ar[d]^{\Delta\times \Id}\\
A\times(A\times A) \ar[dr]_{a\bot(b\bot c)} \ar[rr]^-\alpha&& 
(A\times A)\times A\ar[dl]^{(a\bot b)\bot c}\\
& X} \]
Here $\alpha$ is the associativity isomorphism.
The upper part of the diagram commutes because the diagonal morphism
is coassociative. The lower triangle then commutes since
the two morphisms 
\[ a\bot(b\bot c)\ ,\, ((a\bot b)\bot c)\circ\alpha\ : A\times(A\times A)\ \to\ X \]
have the same `restrictions', namely $a$, $b$ respectively $c$.

The commutativity is a consequence of two elementary facts:
first, $b\bot a=(a\bot b)\tau$
where $\tau:A\times A\to A\times A$
is the automorphism that interchanges the two factors; this follows from
$\tau i_1=i_2$ and  $\tau i_2=i_1$. Second, the diagonal morphism
is cocommutative, i.e., $\tau\Delta=\Delta:A\to A\times A$.
Altogether we get 
\[ a+b\ =\ (a\bot b)\Delta\ =\ (a\bot b)\tau\Delta\ =\ (b\bot a)\Delta \ =\ b+a \ .\]
As before we denote by $0\in\Dc(A,X)$ the unique morphism that factors
through a zero object.
We let $p_1:A\times A\to A$ be the projection to the first factor. Then
\begin{align*}
  (a\bot 0)i_1 \ &= \ a \ = \ a p_1 (\Id,0) \ = \ a p_1 i_1
\text{\quad and}\\
 (a\bot 0)i_2 \ &= \ 0 \ = \ a p_1 (0,\Id) \ = \ a p_1 i_2\ .
\end{align*}
So we have $a\bot 0=a p_1$ in $\Dc(A\times A,X)$.
Hence $a+0=(a\bot 0)\Delta=a p_1\Delta=a$; by commutativity we
also have $0+a=a$.

Now we verify naturality of the addition.
To check $(a+b)c=ac+b c$ for $a,b:A\to X$ and $c:A'\to A$ we consider
the commutative diagram
\[
\xymatrix@C=12mm{
A'\ar[d]_{\Delta} \ar[r]^-{c} & 
A\ar[r]^-{\Delta} \ar[d]^\Delta & A\times A\ar[d]^{a\bot b}\\
A'\times A' \ar@/_1pc/[rr]_-{ac\bot b c} \ar[r]^-{c\times c} & 
A\times A \ar[r]^-{a\bot b} & X}
  \]
in which the composite through the upper right corner is $(a+b)c$.
We have 
\[ (a\bot b)(c\times c)i_1\ =\ (a\bot b)(c,0)\ = \ a c\ =\ (ac\bot b c)i_1  \]
and similarly for $i_2$ instead of $i_1$. 
So $(a\bot b)(c\times c)=ac\bot b c$ since both sides
have the same `restrictions' to the two factors of $A'\times A'$.
Since the composite through the lower left corner is
$ac+b c$, we have shown $(a+b)c=ac+b c$. Naturality in $X$ is even easier.
For a morphism $d:X\to Y$ we have
$d(a\bot b)=d a\bot db:A\times A\to Y$ since both sides 
have the same `restrictions' $d a$ respectively
$db$ to the two factors of $A\times A$. Thus 
$d(a+ b)=d a+ d b$ by the definition of `+'.
\end{proof}

Now we introduce the group-like objects in a pre-additive category.

\begin{prop}\label{prop:characterize group-like}  
  Let $\Dc$ be a pre-additive category.
    For every object $A$ of $\Dc$ the following two conditions are equivalent:
    \begin{enumerate}[\em (a)]
    \item The shearing morphism 
      $\Delta\bot i_2=(\Delta p_1)+i_2 p_2:A\times A\to A\times A$ 
      is an isomorphism.
    \item The identity of $A$ has an inverse in the abelian monoid $\Dc(A,A)$. 
  \end{enumerate}
We call $A$ {\em group-like}\index{subject}{group-like}
if it satisfies {\em (a)} and {\em (b)}.
If $A$ is group-like, then moreover for every object $X$ of $\Dc$ 
the abelian monoids $\Dc(A,X)$ and $\Dc(X,A)$ have inverses. 
\end{prop}
\begin{proof}
(a)$\Longrightarrow$(b) 
Since the shearing map is an isomorphism, there is a morphism $(k,j):A\to A\times A$
such that
\[ (\Id_A, 0)\ = \ (\Delta\bot i_2)\circ (k,j) \ = \ (k,k+j) \ . \]
So $k=\Id_A$ and $\Id_A+j=0$, i.e., $j$ is an additive inverse of the identity of $A$.

(b)$\Longrightarrow$(a) 
If $j\in \Dc(A,A)$ is an inverse of the identity of $A$, 
then the morphism
\[ (p_1,j\bot \Id_A)\ = \ (\Id_A,j)\bot i_2 \ : \ A\times A \ \to \ A\times A \]
is a two-sided inverse to the shearing morphism, which is thus an isomorphism.

If $j\in\Dc(A,A)$ is an additive inverse to the identity of $A$, 
then for all $f\in\Dc(X,A)$
\[ f +(j f)\ = \ (\Id_A\circ f)+ (j\circ f)\ = \ 
(\Id_A+j)\circ f\ = \ 0\circ f \ = \ 0\ ; \]
so $j\circ f$ is inverse to $f$. Similarly, for every $g\in\Dc(A,X)$ 
the morphism $g\circ j$ is additively inverse to $g$.
\end{proof}

In the next definition and in what follows, we denote by $M^\times$
the subgroup of invertible elements in an abelian monoid $M$.

\begin{defn}\label{def:abstract group completion and unit}
Let $\Dc$ be a pre-additive category.
A morphism $u:R^\times \to R$ is a 
{\em unit morphism}\index{subject}{unit morphism!in a pre-additive category}
if for every object $T$ the map
\[ \Dc(T,u)\ : \ \Dc(T,R^\times)\ \to \ \Dc(T,R)  \]
is injective with image the subgroup $\Dc(T,R)^\times$ of invertible elements.
A morphism $i:R\to R^\star$ in $\Dc$ is a 
{\em group completion}\index{subject}{group completion!in a pre-additive category}  
if for every object $T$ the map
\[ \Dc(i,T)\ : \ \Dc(R^\star,T )\ \to \ \Dc(R,T)  \]
is injective with image the subgroup $\Dc(R,T)^\times$ of invertible elements.
\end{defn}

\begin{rk}
If $u:R^\times\to R$ is a unit morphism 
then the abelian monoid $\Dc(R^\times,R^\times)$ is a group, by the defining property;
so the object $R^\times$ is in particular group-like.
Since the pair $(R^\times,u)$ represents the functor
\[ \Dc\ \to \ \text{(sets)} \ , \quad
T \ \longmapsto \ \Dc(T,R)^\times \ ,\]
it is unique up to preferred isomorphism.
A formal consequence is that if we choose  
a unit morphism $u_R:R^\times\to R$  for every object $R$, 
then this extends canonically to a functor
\[ (-)^\times\ : \  \Dc\ \to \  \Dc \]
and a natural transformation $u:(-)^\times\to\Id$.
Since the functor $(-)^\times$ takes values in group-like objects, 
it is effectively a right adjoint to the inclusion of
the full subcategory of group-like objects.

A category $\Dc$ is pre-additive if and only if its opposite
category $\Dc^{\op}$ is pre-additive. 
Moreover, in that situation $\Dc(A,X)$ and $\Dc^{\op}(X,A)$
are not only the same set (by definition), but they also have the
same monoid structure. Thus the concepts of unit morphism and
group completion are `dual' (or `opposite') to each other:
a morphism is a unit morphism in $\Dc$ if and only
if it is a group completion in $\Dc^{\op}$. 
This is why many properties of unit maps have a corresponding `dual'
property for group completions, and why most proofs for unit maps 
have `dual' proofs for group completions.
Since the identity of any object of $\Dc$ is also the identity of
the same object in $\Dc^{\op}$, 
part~(b) of Proposition \ref{prop:characterize group-like}  
shows that `group-like' is a self-dual property:
an object is group-like in $\Dc$ if and only if it is group-like in $\Dc^{\op}$.

So the above properties of unit morphisms dualize:
if $i:R\to R^\star$ is a group completion, then $R^\star$
is in particular group-like.
The pair $(R^\star,i)$ is unique up to preferred isomorphism,
and if we choose a group completion $i_R:R\to R^\star$ 
for every object $R$, then this extends canonically to a functor
\[ (-)^\star\ : \  \Dc\ \to \  \Dc  \]
and a natural transformation $i:\Id\to (-)^\star$,
producing a left adjoint to the inclusion of group-like objects.
\end{rk}

\begin{eg}[Unit morphisms and group completion for abelian monoids]\label{eg:Grothendieck group} 
The category of abelian monoids is the prototypical example of
a pre-additive category, and the general theory of units and group completions
is an abstraction of this special case. 
So we take the time to convince ourselves that the concepts of `unit morphism'
and `group completion' have their familiar meanings in the motivating example.

A basic observation is that in the category of abelian monoids, 
the abstract addition of morphism as in Proposition \ref{prop:additive category}
is simply pointwise addition of homomorphisms. So an abelian monoid
is group-like in the abstract sense of Proposition \ref{prop:characterize group-like} 
if and only if every element has an inverse; so the group-like objects
are precisely the abelian groups.

A given homomorphism $f:M\to N$ of abelian monoids is invertible 
if and only if it is pointwise invertible in $N$, 
which is the case if any only if the image of $f$ lies in
the subgroup $N^\times$ of invertible elements.
So the inclusion $u:N^\times\to N$ of the subgroup of 
invertible elements is a unit morphism in the sense of Definition \ref{def:abstract group completion and unit}.

We recall the {\em Grothendieck group}\index{subject}{Grothendieck group!of an abelian monoid} of an abelian monoid $M$.
An equivalence relation $\sim$ on $M^2$ is defined by
declaring $(x,y)$ equivalent to $(x',y')$ if and only if there is an element $z\in M$
with 
\[ x + y' + z \ = \ x' + y + z \ . \]
The componentwise addition on $M^2$ is well-defined on equivalence classes,
so the set of equivalence classes
\[ M^\star \ = \ M^2/\sim \]
inherits an abelian monoid structure. We write $[x,y]$
for the equivalence class in $M^\star$ of a pair $(x,y)$.
The pair $(x+y,y+x)$ is equivalent to $(0,0)$, so 
\[ [x,y] + [y,x] \ = \  [x+y,y+x]\ = \ 0\]
in the monoid $M^\star$. This shows that every element of $M^\star$ has an inverse,
and $M^\star$ is an abelian group.
We claim that the monoid homomorphism\index{subject}{group completion!of an abelian monoid} 
\[ i \ : \ M \ \to \ M^\star \ , \quad  i(x)\ = \ [x,0] \]
is a group completion in the sense of 
Definition \ref{def:abstract group completion and unit}. 
Indeed, given a monoid homomorphism $h:M\to N$ that is pointwise invertible,
then we can define $f:M^\star\to N$ by
\[ f[x,y]\ = \ h(x) - h(y)\ . \]
A routine verification shows that $f$ is indeed a well-defined homomorphism
and that sending $h$ to $f$ is inverse to the restriction map
\[ \Ab Mon(i,N) \ : \ \Ab Mon(M^\star,N) \ \to \ \Ab Mon(M,N)^\times\ .\]
A slightly different way to summarize the construction of the Grothendieck group
is to say that the group completion of an abelian monoid $M$ is a cokernel,
in the category of commutative monoids, of the diagonal morphism
$\Delta:M \to M\times M$. 

We observe that
\[ [x,y]\ = \ [x,0] + [0,y]\ = \ i(x) - i(y) \ ,\]
so every element in $M^\star$ is the difference of two elements in the
image of $i:M\to M^\star$. Moreover, if $x,x'\in M$ satisfy $i(x)=i(x')$,
then the pairs $(x,0)$ and $(x',0)$ are equivalent, which happens if and only if
there is an element $z\in M$ such that $x+z=x'+z$.
Conversely, these properties of the Grothendieck construction
characterize group completions of abelian monoids:
a homomorphism $j:M\to N$ of abelian monoids is a group completion
if and only if the following three conditions are satisfied:
\begin{itemize}
\item the monoid $N$ is a group;
\item every element in $N$ is the difference of two elements in the
image of $j$;
\item if $x,x'\in M$ satisfy $j(x)=j(x')$, then 
there is an element $z\in M$ such that $x+z=x'+z$.
\end{itemize}
Indeed, the first condition guarantees that $j$ extends (necessarily uniquely)
to a homomorphism $M^\star\to N$; the second and third conditions guarantee that
the extension is surjective respectively injective, hence an isomorphism.
\end{eg}

We mostly care about the situation where $\Dc=\Ho(\Cc)$ is the homotopy category
of a pointed model category $\Cc$, such as the category of ultra-commutative monoids.
As we shall now proceed to prove, in this situation
units and group completions always exist. 

We consider two composable morphisms
$f:A\to B$ and $g:B\to C$ in a pointed category $\Dc$.
We recall that $f$ is a {\em kernel} of $g$ if $g f=0$ and
for every morphism $\beta:T\to B$ such that $g\beta=0$,
there is a unique morphism $\alpha:T\to A$  such that $f\alpha=\beta$.
Dually, $g$ is a {\em cokernel} of $f$ if $g f=0$ and
for every morphism $\beta:B\to Y$ such that $\beta f=0$,
there is a unique morphism $\gamma:C\to Y$  such that $\gamma g=\beta$.

\begin{prop}\label{prop:kernel of codiagonal}
Let $R$ be an object of a pre-additive category $\Dc$.
\begin{enumerate}[\em (i)]
\item Let $e:R^\times\to R\times R$ be a kernel 
of the codiagonal morphism $\Id\bot\Id:R\times R\to R$.
Then the composite
\[ u \ = \ (\Id\bot 0)\circ e \ : \ R^\times \ \to \ R \]
is a unit morphism and 
\[ e \ = \ (u, - u)\ :\ R^\times \ \to \ R\times R \ . \]
Conversely, if $u:R^\times\to R$ is a unit morphism,
then the  morphism $(u,-u): R^\times \to  R\times R$
is a kernel of the codiagonal morphism $\Id\bot\Id:R\times R\to R$.
\item Let $d:R\times R\to R^\star$ be a cokernel 
of the diagonal morphism $(\Id,\Id):R\to R\times R$.
Then the composite
\[ i \ = \ d\circ(\Id,0)\ : \ R \  \to \ R^\star\]
is a group completion and 
\[ d \ = \ i\bot (- i )\ :\ R\times R\ \to \ R^\star\ . \]
Conversely, if $i:R\to R^\star$ is a group completion, 
then the  morphism $i\bot (- i) : R\times R \to  R^\star$
is a cokernel of the diagonal morphism $(\Id,\Id):R\to R\times R$.
\end{enumerate}
\end{prop}
\begin{proof}
We prove part~(ii). Part~(i) is dual, i.e., equivalent to
part~(ii) in the opposite category $\Dc^{\op}$.
We let $T$ be any object $\Dc$. Then the map
\[ \{ (f,g)\in \Dc(R,T)^2\ | \  f +g = 0 \} \ \to \  \Dc(R,T)^\times  \ ,\quad 
(f,g) \ \longmapsto \ f\]
is bijective because inverses in abelian monoids, if they exist, are unique.
A cokernel of the diagonal morphism is a morphism $d:R\times R\to R^\star$ 
that represents the left hand side of this bijection;
a group completion is a morphism $i:R\to R^\star$ 
that represents the right hand side of this bijection.
Hence $d:R\times R\to R^\star$ is a cokernel of the diagonal 
if and only if $d\circ(\Id,0)$ is a group completion.

The relation $d\circ (\Id,0)=i$ holds by definition.
The relation
\[ 
 ( d\circ (0,\Id) ) + i   \ = \  
 d\circ  ( (0,\Id)  +   (\Id,0)) \ = \  d\circ   (\Id,\Id)   \ = \ 0 
 \]
holds in the monoid $\Dc(R,R^\star)$, and thus 
$d\circ (0,\Id) = - i$. This shows that $d= i\bot (- i)$.
\end{proof}

The previous characterization of unit morphisms and group completions
as certain kernels respectively cokernels formally implies the following corollary.

\begin{cor}\label{cor-preserves localization}
Let $F:\Dc\to \Ec$ be a functor between pre-additive categories that
preserves products.
\begin{enumerate}[\em (i)]
\item If $F$ preserves kernels of splitable epimorphisms, then for every
unit morphism $u:R^\times \to R$ in $\Dc$, the morphism
$F u:F( R^\times)\to F R$ is a unit morphism.
\item If $F$ preserves cokernels of splitable monomorphisms, then for every
group completion $i:R \to R^\star$ in $\Dc$, the morphism
$F i:F R\to F( R^\star)$ is a group completion.
\end{enumerate}
\end{cor}

Now we consider a pointed model category $\Cc$ whose homotopy category is
pre-additive. The main example we have in mind is $\Cc=umon$,
the category of ultra-commutative monoids with the global model structure 
of Theorem \ref{thm:global umon}. 
The homotopy category $\Ho(\Cc)$ then comes with 
an adjoint functor pair $(\Sigma,\Omega)$ of suspension and loop,
compare \cite[I.2]{Q}.

\begin{prop}\label{prop:loop and suspension make group-like} 
Let $\Cc$ be a pointed model category whose homotopy category is
pre-additive. 
\begin{enumerate}[\em (i)]
\item 
For every object $R$ of $\Cc$,
the loop object $\Omega R$  and the suspension $\Sigma R$ 
are group-like in $\Ho(\Cc)$.
\item If $u:R^\times\to R$ is a unit morphism, then its loop
$\Omega u:\Omega (R^\times)\to \Omega R$ is an isomorphism in $\Ho(\Cc)$.
\item If $i:R\to R^\star$ is a group completion, then its suspension
$\Sigma i:\Sigma R\to \Sigma( R^\star)$ is an isomorphism in $\Ho(\Cc)$.
\end{enumerate}
\end{prop}
\begin{proof}
(i) This is a version of the Eckmann-Hilton argument. 
For every object $T$ of $\Cc$, the set $[T,\Omega R]$
has one abelian monoid structure 
via Construction \ref{con:hom addition}, coming from the fact that $\Ho(\Cc)$
is pre-additive.
A second binary operation on the set $[T,\Omega R]$
arises from the fact that $\Omega R$ is a group object in $\Ho(\Cc)$,
compare \cite[I.2]{Q}. This operation makes $[T,\Omega R]$ into a group. 
The monoid structure of Construction \ref{con:hom addition}
is natural for morphisms in the second variable $\Omega R$, in particular
for the group structure morphism  $\Omega R\times \Omega R\to\Omega R$.
This means that the two binary operations satisfy the interchange law.
Since they also share the same neutral element, they coincide.
Since one of the two operations has inverses, so does the other.
 
The argument that $\Sigma R$ is group-like is dual, using 
that $\Sigma R$ is the loop object of $R$ in $\Ho(\Cc)^{\op}=\Ho(\Cc^{\op})$,
and that `group-like' is a self-dual property.

(ii) 
Since the functor $\Omega:\Ho(\Cc)\to\Ho(\Cc)$ is right adjoint to $\Sigma$,
it preserves products and kernels.
So $\Omega u:\Omega(R^\times)\to\Omega R$ 
is a unit morphism by Corollary \ref{cor-preserves localization}.
Since $\Omega R$ is already group-like by part~(i), $\Omega u$ is an isomorphism.
Part~(iii) is dual to part~(ii); so it admits the dual proof,
or can be obtained by applying part~(ii) to the opposite model category.
\end{proof}

\begin{prop}\label{prop:kernel and cokernel}
Consider a commutative square 
\[  \xymatrix@C=10mm{ 
A \ar[r]^-f \ar[d] & B \ar[d]^g \\
C \ar[r] & D } \]
in a pointed model category $\Cc$ such that the object $C$ is weakly contractible.
\begin{enumerate}[\em (i)]
\item If the square is homotopy cartesian and $g$ admits a section in $\Ho(\Cc)$,
then the morphism $f$ is a kernel of $g$ in $\Ho(\Cc)$. 
\item If the square is homotopy cocartesian and $f$ admits a retraction in $\Ho(\Cc)$,
then the morphism $g$ is a cokernel of $f$ in $\Ho(\Cc)$.
\end{enumerate}
\end{prop}
\begin{proof}
We prove part~(i). Part~(ii) can be proved by dualizing the argument or by
applying part~(i) to the opposite category with the opposite model structure.
Since the square is homotopy cartesian and $C$ is weakly contractible,
the object $A$ is weakly equivalent to the homotopy fiber, 
in the abstract model category sense, of the morphism $g$. 
As Quillen explains in Section~I.3 of \cite{Q},
there is an action map (up to homotopy)
\[ A \times (\Omega D) \ \to \ A \ ,\]
by an abstract version of `fiber transport'.
For every other object $T$ of $\Cc$, Proposition~4 of \cite[I.3]{Q}
provides a sequence of based sets
\[  [T, \Omega B] \ \xra{[T,\Omega g]}\ 
[T, \Omega D] \ \xra{[T,\partial]}\  [T, A] \ \xra{[T,f]}\ 
 [T,B] \ \xra{[T,g]}\  [T,D] \]
that is exact in the sense explained in \cite[I.3 Prop.\,4]{Q},
where $[-,-]$ denotes the morphism sets in the homotopy category of $\Cc$.
In particular, the image of $[T,f]$ is equal to the preimage of
the zero morphism under the map $[T,g]$.

So in order to show that $f$ is a kernel of $g$ it remains to 
check that the map $[T,f]$ is injective. So we consider 
two morphisms $\alpha_1,\alpha_2\in[T,A]$ such that $f\circ\alpha_1=f\circ\alpha_2$.
Then by Proposition~4~(ii) of \cite[I.3]{Q}, 
there is an element $\lambda\in [A,\Omega D]$
such that $\alpha_2=\alpha_1\cdot\lambda$. 
Since the morphism $g:B\to D$ has a section, 
so does the morphism $\Omega g:\Omega B\to\Omega D$. So there is a morphism
$\bar\lambda\in [T,\Omega B]$ 
such that $\lambda=(\Omega g)\circ\bar\lambda$.
But all elements in the image of $[T,\Omega g]$ act trivially on $[T,A]$,
so 
then
\[ \alpha_2\ = \ 
\alpha_1\cdot\lambda \ = \  \alpha_1\cdot( (\Omega g)\circ\bar\lambda) \ = \ \alpha_1\ .
\qedhere \]
\end{proof}

\begin{theorem}\label{thm:unit completion existence}
Let $\Cc$ be a pointed model category whose homotopy category is
pre-additive. 
\begin{enumerate}[\em (i)]
\item Every object of $\Cc$ has a unit morphism and a group completion in $\Ho(\Cc)$. 
\item If $\Cc$ is right proper, then every object $R$ 
admits a $\Cc$-morphism $u:R^\times\to R$ 
that becomes a unit morphism  in the homotopy category $\Ho(\Cc)$. 
\item If $\Cc$ is left proper, then every object $R$ 
admits a $\Cc$-morphism $i:R\to R^\star$ 
that becomes a group completion in the homotopy category $\Ho(\Cc)$. 
\end{enumerate}
\end{theorem}
\begin{proof}
(i) We let $R$ be any object of $\Cc$. 
It suffices to show, by Proposition \ref{prop:kernel of codiagonal},
that the codiagonal morphism $\Id\bot\Id:R\times R\to R$ has a kernel in $\Ho(\Cc)$
and the diagonal morphism $(\Id,\Id):R\to R\times R$ has a cokernel in $\Ho(\Cc)$.
The arguments are again dual to each other, so we only show the first one.

We can assume without loss of generality that $R$ is cofibrant and fibrant.
Then the fold map $\nabla:R\amalg R\to R$ in the model category $\Cc$
becomes the codiagonal morphism of $R$ in $\Ho(\Cc)$.
We factor $\nabla=q\circ j$ for some weak equivalence $j:R\amalg R\xra{\sim} Q$
followed by a fibration $q:Q\to R$. Then we choose a pullback,
so that we arrive at the homotopy cartesian square:
\[ \xymatrix{ P \ar[r]^f\ar[d] & Q\ar[d]^q \\
\ast \ar[r] & R} \]
The morphism $q$ still becomes a codiagonal morphism in $\Ho(\Cc)$,
and so it has a section. By
Proposition \ref{prop:kernel and cokernel}~(i) the morphism $f$
becomes a kernel of $q$ in $\Ho(\Cc)$. So the codiagonal morphism
of $R$ has a kernel.  

(ii) We choose a weak equivalence $q:R\to \bar R$
to a fibrant object. A unit morphism $\bar R^\times\to\bar R$ exists in $\Ho(\Cc)$ 
by part~(i). By replacing the source $\bar R^\times$
by a weakly equivalent object, if necessary, we can assume that it
is cofibrant as an object in the model category $\Cc$.
Every morphism in $\Ho(\Cc)$ from a cofibrant to a fibrant object 
is the image of some $\Cc$-morphism under the localization functor, i.e.,
there is a $\Cc$-morphism $\bar u:\bar R^\times \to\bar R$ 
that becomes a unit morphism in $\Ho(\Cc)$. 
By factoring $\bar u$ as a weak equivalence
followed by a fibration we can moreover assume without loss of generality
that $\bar u$ is a fibration. We form a pullback
\[ \xymatrix{ 
R^\times \ar[r]^-u \ar[d]_p & R\ar[d]^q_\sim\\ 
\bar R^\times \ar[r]_-{\bar u} & \bar R  } \]
Since $q$ is a weak equivalence, $\bar u$ a fibration and $\Cc$ is right proper,
the base change $p$ of $q$ is also a weak equivalence.
So $u$ is isomorphic to $\bar u$ in the arrow category in $\Ho(\Cc)$,
hence $u$ is also a unit morphism of $\Ho(\Cc)$.
Part~(iii) is dual to part~(ii).
\end{proof}

\begin{rk}
We claim that unit morphisms and group completions also behave nicely
on derived mapping spaces. We explain this in detail for unit morphisms,
the other case being dual, one more time.
Model categories have derived mapping spaces (i.e., simplicial sets) $\map^h(-,-)$, 
giving well-defined homotopy types such that
\begin{equation}\label{eq:pi_0 of Rmap}
 \pi_0 (\map^h(T,R)) \ \iso \ \Ho(\Cc)(T,R)\ ,  
\end{equation}
compare \cite[Sec.\,5.4]{hovey-book} or \cite[Ch.\,18]{hirschhorn}.
We let $u:R^\times\to R$ be a $\Cc$-morphism that becomes a unit morphism in $\Ho(\Cc)$,
and $T$ any other object of $\Cc$. Because of the bijection \eqref{eq:pi_0 of Rmap}
the map
\[ u_*\ : \ \map^h(T,R^\times)\ \to \ \map^h(T,R)  \]
lands in the subspace $\map^{h,\times}(T,R)$, defined as the union of those 
path components that represent invertible elements in the monoid $\Ho(\Cc)(T,R)$.
We claim that $u_*$ is a weak equivalence onto the subspace $\map^{h,\times}(T,R)$.
To see this we exploit that both 
$\map^h(T,R^\times)$ and $\map^{h,\times}(T,R)$ are group-like H-spaces,
the multiplication arising from the fact $T$ is a comonoid object up to homotopy.
Moreover, the map $u_*$ is an H-map and bijection on path components
(by the universal property of unit morphisms and the bijection \eqref{eq:pi_0 of Rmap}).
So it suffices to show that the restriction of $u$ to the identity path components
is a weak equivalence. For this it suffices in turn to show that
the looped map~
\[ \Omega (u_*)\ : \ \Omega( \map^h(T,R^\times))\ \to \ \Omega(\map^{h,\times}(T,R))  \]
is a weak equivalence.
But this map is weakly equivalent to 
\[ ( \Omega u)_*\ : \ \map^h(T,\Omega( R^\times))\ \to \ \map^h(T,\Omega R) \ . \]
Since $\Omega u$ is an isomorphism in $\Ho(\Cc)$ 
(by Proposition \ref{prop:loop and suspension make group-like}~(ii))
it is a weak equivalence in $\Cc$,
hence so is the induced map on derived mapping spaces. 
\end{rk}

The next proposition will be used to show that loops on the bar construction 
provide functorial global group completions of ultra-commutative monoids.

\begin{prop}\label{prop:Omega Sigma group completion} 
Let $\Cc$ be a pointed model category whose homotopy category is
pre-additive. Suppose that for every group-like object $R$ of $\Cc$ the
adjunction unit $\eta:R\to \Omega(\Sigma R)$ is an isomorphism in $\Ho(\Cc)$.
Then for every object $R$ of $\Cc$ the
adjunction unit $\eta:R\to \Omega(\Sigma R)$ is a group completion.
\end{prop}
\begin{proof}
We let $i:R\to R^\star$ be a group completion, which exists by 
Theorem \ref{thm:unit completion existence}.
In the commutative square in $\Ho(\Cc)$
\[ \xymatrix@C=12mm{ 
R \ar[r]^-{\eta_R}\ar[d]_i &  \Omega(\Sigma R) \ar[d]^{\Omega(\Sigma i)}_\iso \\
R^\star \ar[r]_-{\eta_{R^\star}}^-\iso &  \Omega(\Sigma( R^\star)) } \]
the lower horizontal morphism is an isomorphism by hypothesis because $R^\star$
is group-like.
The morphism $\Sigma i:\Sigma R\to \Sigma(R^\star)$ is an isomorphism by   
Proposition \ref{prop:loop and suspension make group-like}~(iii),
hence the right vertical morphism
$\Omega(\Sigma i)$ is also an isomorphism.
So $\eta_R$ is isomorphic, as an object in the comma category $R\downarrow\Ho(\Cc)$,
to $i$, and hence also a group completion.
\end{proof}

The previous proposition also has a dual statement
(with the dual proof):  
if for every group-like object $R$ of $\Cc$ the
adjunction counit $\epsilon:\Sigma(\Omega R)\to R$ is an isomorphism in $\Ho(\Cc)$,
then $\epsilon$ is a unit morphism.
In practice, however, this dual formulation is less useful:
in the important examples that arise `in nature',
for example for ultra-commutative monoids,
the adjunction unit $\eta:R\to \Omega(\Sigma R)$ is an isomorphism 
for all group-like objects $R$, whereas
the adjunction counit $\epsilon:\Sigma(\Omega R)\to R$ is {\em not}
always an isomorphism. 

\medskip

Now we specialize the theory of units and group completions to
ultra-commu\-tative monoids. We recall that the category of ultra-commutative
monoids has the trivial monoid as zero object, and
the canonical morphism $\rho_{R,S}:R\boxtimes S\to R\times S$
from the coproduct to the product of two ultra-commutative monoids
is a global equivalence by Theorem \ref{thm:box to times}~(i).
So the homotopy category $\Ho(umon)$ is pre-additive.

\begin{defn}
An ultra-commutative monoid $R$ is 
{\em group-like}\index{subject}{group-like!ultra-commutative monoid} 
if it is group-like as an object of the pre-additive category $\Ho(umon)$.
A morphism $u:R^\times \to R$ of ultra-commutative monoids is a 
{\em global unit morphism}\index{subject}{global unit morphism}\index{subject}{units!of an ultra-commutative monoid}  
if it is a unit morphism in the pre-additive category $\Ho(umon)$.
A morphism $i:R\to R^\star$ of ultra-commutative monoids is a 
{\em global group completion}\index{subject}{global group completion}\index{subject}{group completion!global|see{global group completion}}  
if it is a group completion in the pre-additive category $\Ho(umon)$.
\end{defn}

The global model structure on the category of ultra-commutative monoids
is proper, see Theorem \ref{thm:global umon}.
Theorem \ref{thm:unit completion existence} thus guarantees that every ultra-commutative
monoid $R$ admits 
a global unit morphism $u:R^\times\to R$ and
a global group completion $i:R\to R^\star$.

As a reality check we show that for ultra-commutative monoids $R$, 
the abstract definition of `group-like' is equivalent to the requirement
that all the abelian monoids $\pi_0^G(R)$ are groups.
This part works more generally for all orthogonal monoid spaces,
not necessarily commutative. 
A monoid $M$ (not necessarily abelian) is a group if and only if the shearing map
\[ \chi \ : \ M^2\ \to \ M^2 \ , \quad (x,y) \ \longmapsto \ (x, x y)  \]
is bijective. Indeed, if $M$ is a group, then the map
$(x,z)\mapsto (x, x^{-1} z )$ is inverse to $\chi$.
Conversely, if $\chi$ is bijective, then for every $x\in M$
there is a $y\in M$ such that $\chi(x,y)=(x,1)$, i.e., with $x y=1$. 
Then $\chi(x,y x)=(x,x y x)=(x,x)=\chi(x,1)$, so $y x=1$ 
by injectivity of $\chi$. Thus $y$ is a two-sided inverse for $x$.

For orthogonal monoid spaces $R$ (not necessarily commutative), 
the group-like condition has a
similar characterization as follows. 
The {\em shearing morphism}\index{subject}{shearing morphism!of an  orthogonal monoid space}
is the morphism of orthogonal spaces
\[ \chi = (\rho_1,\mu)\ : \  R\boxtimes R \ \to \ R\times R\]
whose first component is the projection $\rho_1$ to the first factor
and whose second component is the multiplication morphism $\mu:R\boxtimes R\to R$.

\medskip

\Danger
The multiplication morphism $\mu:R\boxtimes R\to R$, and hence the
shearing morphism $\chi$, is a homomorphism of orthogonal monoid spaces
only if $R$ is commutative.

\begin{prop}
Let $R$ be an orthogonal monoid space. Then the following two conditions are equivalent:
\begin{enumerate}[\em (i)]
\item The shearing morphism $\chi:R\boxtimes R\to R\times R$ is a global equivalence
of orthogonal spaces.
\item For every compact Lie group $G$ the monoid $\pi_0^G(R)$ is a group.
\end{enumerate}
For commutative orthogonal monoid spaces, conditions {\em (i)} and {\em (ii)} 
are moreover equivalent to being group-like in the pre-additive 
homotopy category of ultra-commutative monoids.
\end{prop}
\begin{proof}
(i)$\Longrightarrow$(ii)
The vertical maps in the commutative diagram
\begin{equation} \begin{aligned}\label{eq:pi_0 umon square}
 \xymatrix@C=25mm{ 
\pi_0^G(R\boxtimes R) \ar[d]_{(\pi_0^G(\rho_1),\pi_0^G(\rho_2))}^\iso \ar[r]^-{\pi_0^G(\chi)} &
\pi_0^G( R\times R) \ar[d]^{(\pi_0^G(p_1),\pi_0^G(p_2))}_\iso  \\
\pi_0^G( R ) \times \pi_0^G( R ) \ar[r]_-{(x,y)\ \longmapsto \ (x, x y)} &
\pi_0^G( R ) \times \pi_0^G( R ) }    
  \end{aligned}\end{equation}
are bijective by Corollary \ref{cor-pi_0 of box}.
If the shearing morphism is a global equivalence,
then the map $\pi_0^G(\chi)$ is bijective, hence so is the algebraic shearing map
of the monoid $\pi_0^G(R)$. This monoid is thus a group.

(ii)$\Longrightarrow$(i)
Now we assume that all the monoids $\pi_0^G(R)$ are groups.
We assume first that $R$ is flat as an orthogonal space;
then $R\boxtimes R$ is also flat, 
by Proposition \ref{prop:ExF ppp spaces}~(i) for the global family of
all compact Lie groups.
The product $R\times R$ is also flat, 
by Proposition \ref{prop:ppp for cartesian product}.
Since $R\boxtimes R$ and $R\times R$ are flat, they are also closed as
orthogonal spaces by Proposition \ref{prop:flat becomes K-cofibration}~(iii).
We may thus show that for every compact Lie group $G$ the continuous map
\[ \chi^G \ = \ \chi(\Uc_G)^G \ : \ 
(R\boxtimes R)(\Uc_G)^G \ \to \ (R\times R)(\Uc_G)^G\ = \ R(\Uc_G)^G \times R(\Uc_G)^G \]
is a weak equivalence, compare Proposition \ref{prop:global eq for closed}. 
Since the monoid $\pi_0^G(R)$ has inverses, the shearing morphism
$\chi:R\boxtimes R\to R\times R$ induces a bijection on $\pi_0^G$,
by the commutative diagram \eqref{eq:pi_0 umon square} with vertical bijections.
On path components we have
\[ \pi_0^G(R\boxtimes R) \ \iso \ \pi_0((R\boxtimes R)(\Uc_G)^G)
\text{\quad and\quad}
\pi_0^G(R\times R) \ \iso \  \pi_0((R\times R)(\Uc_G)^G)\ ,\]
compare Corollary \ref{cor:pi_0^G of closed}. 
So we conclude that the map $\chi^G$ induces a bijection on path components. 

Now we show that $\chi^G$ also induces bijections on homotopy groups
in positive dimensions. We consider a point
$x\in ((R\boxtimes R)(\Uc_G))^G$ and $k\geq 1$. 
We let $\varphi:\Uc_G^2\to\Uc_G$ be any $G$-equivariant linear isometric embedding.
As we explained in Remark \ref{rk:umon is global E_infty}, 
this map induces an H-space structure
on $(R\boxtimes R)(\Uc_G)^G$, and hence a continuous map 
\[ \varphi_*(-,x)^G\ : \ (R\boxtimes R)(\Uc_G)^G \ \to \ 
(R\boxtimes R)(\Uc_G)^G  \ .\]
Since the unit element~1 is a homotopy unit for the H-space structure, 
the element $x'=\varphi_*(1,x)$ belongs to the same path component as $x$. 
Since the monoid $\pi_0(R(\Uc_G)^G)$ is isomorphic to the group $\pi_0^G(R)$,
the H-space structure has inverses.
So the map $\varphi_*(-,x)^G$ is a homotopy equivalence.
The same argument applies to $R\times R$ instead of $R\boxtimes R$,
and we obtain a commutative diagram
\[ \xymatrix@C=20mm{ 
\pi_k^G((R\boxtimes R)(\Uc_G)^G, 1)
\ar[r]^-{\pi_k(\chi^G,1)} \ar[d]^\iso_{\pi_k(\varphi_*(-,x)^G)} &
\pi_k^G( (R\times R)(\Uc_G)^G, 1) \ar[d]^{\pi_k(\varphi_*(-,\chi^G(x))^G)}_\iso \\
\pi_k^G((R\boxtimes R)(\Uc_G)^G, x')
\ar[r]_-{\pi_k(\chi^G,x')} &
\pi_k^G( (R\times R)(\Uc_G)^G, \chi^G(x'))  } \]
in which both vertical maps are bijective.
So to show that the map $\chi^G$ induces bijections of homotopy
groups based at $x$, it suffices to show this for the special case $x=1$
of the unit element.

Now the map $\pi_k(\chi^G,1)$ is a group homomorphism such that the composite
\begin{align*}
 (\pi_k(R(\Uc_G)^G,1))^2 \ \xra{(x,y)\mapsto x\times y} \
&\pi_k((R\boxtimes R)(\Uc_G)^G,1)  \\
 \xra{\pi_k(\chi^G,1)}\ 
 &\pi_k((R\times R)(\Uc_G)^G,1) \xra{((\rho_1)_*,(\rho_2)_*)} \ 
(\pi_k(R(\Uc_G)^G,1))^2   
\end{align*}
sends $(x,y)$ to $(x,\mu_*(x\times y))$, 
where $\mu:R\boxtimes R\to R$ is the multiplication
map. By the Eckmann-Hilton argument, $\mu_*(x\times y)=x y$, the product with respect
to the group structure of $\pi_k(R(\Uc_G)^G,1)$.
The first and third maps are bijective, and so is the composite 
(because $\pi_k(R(\Uc_G)^G,1)$ is a group). So the middle map is bijective.
Altogether this shows that the map $\chi^G$ is a weak equivalence.

For general $R$ we choose a global equivalence $f:R'\to R$ of
orthogonal monoid spaces such that $R'$ is flat as an orthogonal space.
One way to arrange this is by cofibrant replacement in the
global model structure of orthogonal monoid spaces
(Corollary \ref{cor-lift to modules spaces}~(ii) with $R=\ast$ and $\Fc=\All$).
Then $f\boxtimes f$ is a global equivalence by Theorem \ref{thm:box to times}
and $f\times f$ is a global equivalence 
by Proposition \ref{prop:global equiv basics}~(vi).
Since $\chi^{R'}$ is a global equivalence by the previous paragraph and
$\chi^R\circ (f\boxtimes f)=(f\times f)\circ\chi^{R'}$,
the morphism $\chi^R$ is also a global equivalence.

Finally, if $R$ is ultra-commutative, then the point-set level shearing morphism $\chi$
becomes the shearing morphism in the sense of 
Proposition \ref{prop:characterize group-like}~(a) in the pre-additive homotopy
category $\Ho(umon)$. So $\chi$ is a global equivalence if and only if
the shearing morphism in $\Ho(umon)$ is an isomorphism, i.e.,
precisely when $R$ is group-like.
\end{proof}

Now we look more closely at global unit morphisms, and we give an explicit,
functorial point set level construction.
For elements in an abelian monoid $M$, left inverses are automatically
right inverses, and they are unique (if they exist). 
So the subgroup of invertible elements of an abelian monoid
is isomorphic to the kernel of the multiplication map, by
\[ M^\times \ \iso \ \ker(+:M^2\to M) \ , \quad x \ \longmapsto \ (x,-x)\ .\]
Proposition \ref{prop:kernel of codiagonal}~(i) gives an abstract
formulation of this and explains how an abstract kernel of the multiplication map
gives rise to a unit morphism. 
The proof of Theorem \ref{thm:unit completion existence}
then shows that the homotopy fiber of the multiplication map,
formed at the model category level, constructs such a kernel.
If we make all this explicit for the model category of
ultra-commutative monoids, we arrive at the following construction.

\begin{construction}[Units of an ultra-commutative monoid]\label{con:R^times}\index{subject}{units!of an ultra-commutative monoid}
We introduce a functorial pointset level construction
of the global units of an ultra-commutative monoid $R$.
We define $R^\times$ as the homotopy fiber, over the additive unit element~0,
of the multiplication morphism $\mu:R\boxtimes R\to R$, i.e.,
\[ R^\times\ = \ F(\mu) \ = \ (R\boxtimes R)\times_{\mu} R^{[0,1]}\times _R \{0\} \ .\]
So at an inner product space $V$, we have
\[ (R^\times)(V)\ = \ (R\boxtimes R)(V)\times_{\mu(V)} R(V)^{[0,1]}\times _{R(V)} \{0\} \ ,\]
the space of pairs $(x,\omega)$ consisting of a point $x\in (R\boxtimes R)(V)$
and a path $\omega:[0,1]\to R(V)$ such that $\mu(V)(x)=\omega(0)$
and $\omega(1)=0$, the unit element in $R(V)$.
Since limits and cotensors with topological spaces of ultra-commutative mo\-noids 
are formed on underlying orthogonal spaces, this homotopy fiber 
inherits a preferred structure of ultra-commutative monoid.

We claim that the  composite
\[ u \ : \ R^\times\ \xra{\ p\ }\ R\boxtimes R \ \xra{\ \rho_1\ } \ R \]
is a global unit morphism, where $p$ is the projection
onto the first factor. 
Indeed, the commutative square
\[ \xymatrix{ R^\times \ar[r]^-p \ar[d]_q  & R\boxtimes R \ar[d]^\mu \\
R^{[0,1]}\times _R \{0\}  \ar[r]_-{\ev_0} & R} \]
is a pullback of ultra-commutative monoids,
by definition, and both horizontal morphisms are strong level fibrations,
where $q$ denotes the projection to the second factor.
So the square is homotopy cartesian.
The multiplication morphism $\mu$ has a section, so
Proposition \ref{prop:kernel and cokernel}~(i) shows that 
the morphism $p:R^\times\to R\boxtimes R$ becomes a kernel of the multiplication
morphism in the homotopy category $\Ho(umon)$.
So $u$ is a unit morphism by Proposition \ref{prop:kernel of codiagonal}~(i).
\end{construction}

We recall from Example \ref{eg:M^times}\index{subject}{units!of a global power monoid}
that every global power monoid $M$ has a global power submonoid $M^\times$ of {\em units};
the value $M^\times(G)$ at a compact Lie group $G$ is
the group of invertible elements of $M(G)$.
The next proposition verifies that global unit morphisms have the
expected behavior on $\upi_0$.

\begin{prop}\label{prop:realize units in umon} 
Let $u:R^\times \to R$ be a unit morphism of ultra-commutative monoids.
  Then the morphism of global power monoids  
  \[ \upi_0(u) \ : \ \upi_0(R^\times)\ \to \ \upi_0(R)\]
  is an isomorphism onto the global power submonoid $(\upi_0(R))^\times$
  of units of $\upi_0(R)$.
\end{prop}
\begin{proof}
This is a formal consequence of the fact that the functor $\pi_0^G$
from the homotopy category of ultra-commutative monoids is representable.
We let $G$ be a compact Lie group and $V$ a non-zero faithful $G$-representation.
Then the global classifying space $B_{\gl}G=\bL_{G,V}$ 
supports the tautological class $u_{G,V}\in \pi_0^G(B_{\gl} G)$,
compare \eqref{eq:tautological_class}.
We recall that
\[ u_G^{umon}\ =\ \eta_*(u_{G,V})\ \in \pi_0^G(\mP(B_{\gl} G)) \]
where $\eta:B_{\gl}G\to\mP(B_{\gl} G)$ is the inclusion of the linear summand.
We claim that evaluation at $u_G^{umon}$ is an isomorphism of abelian monoids
\[ \Ho(umon)( \mP(B_{\gl} G), T ) \ \iso \ \pi_0^G(T) \ , \quad
[f]\ \longmapsto \ f_*(u_G^{umon}) \]
for every ultra-commutative monoid $T$.
Indeed, both sides take global equivalence in $T$ to isomorphisms,
so we may assume that $T$ is fibrant in the global model structure
of ultra-commutative monoids, hence positively static.
Now we consider the composite
\begin{align*}
\pi_0( T(V)^G)\ &\iso \ 
\pi_0(\map^{umon}(\mP(B_{\gl}G),T)) \ \xra{\ \gamma_*\ } \\
&\Ho(umon)( \mP(B_{\gl} G), T ) \ 
\xra{[f]\mapsto f_*(u_G^{umon})} \ \pi_0^G(T)
\end{align*}
with the adjunction bijection and
the map induced by the localization 
functor $\gamma:umon\to\Ho(umon)$.
Since $V$ is non-zero and faithful, $B_{\gl}G=\bL_{G,V}$
is positively flat, so $\mP(B_{\gl}G)$ is cofibrant in the global
model structure of ultra-commutative monoids.
Since $\mP(B_{\gl}G)$ is cofibrant and $T$ is fibrant,
the middle map is bijective.
The composite is the canonical map $\pi_0(T(V)^G)\to\pi_0^G(T)$
to the colimit, which is bijective because $T$ is positively static.

The evaluation isomorphism is natural in the second variable, so we arrive at a 
commutative square of abelian monoids
\[ \xymatrix{
\Ho(umon)(\mP(B_{\gl} G),R^\times) \ar[r]^-\iso \ar[d]_{u_*} & 
\pi_0^G(R^\times) \ar[d]^{u_*}\\
\Ho(umon)(\mP(B_{\gl} G),R) \ar[r]_-\iso  & \pi_0^G(R)}  \]
in which both horizontal maps are bijective.
The left vertical map is injective with image the subgroup of
invertible elements; hence the same is true for the right vertical map.
\end{proof}

\index{subject}{global group completion!of an ultra-commutative monoid|(}
Now we turn to global group completions of ultra-commutative monoids.
Every ultra-commutative monoid has a group completion in the homotopy category,
by Theorem \ref{thm:unit completion existence}~(i).
Even better: since the model category of ultra-commutative monoids is left proper,
every ultra-commutative monoid is the source of a global group completion 
in the model category of ultra-commutative monoids,
by  Theorem \ref{thm:unit completion existence}~(iii).
Now we discuss two functorial pointset level constructions of global group
completions. The first one is dual to Construction \ref{con:R^times} of global units.

\begin{construction}[Global group completion of an ultra-commutative mo\-noid]\label{con:R^star}
We let $R$ be an ultra-commutative monoid
that is cofibrant in the global model structure of Theorem \ref{thm:global umon}. 
We define the cone of $R$ as a pushout in the category of ultra-commutative monoids:
\[ \xymatrix{ 
R\tensor \{0\} \ar[r] \ar[d]_{R\tensor\text{incl}} & \ast \ar[d]\\ 
R\tensor [0,1]\ar[r] & C R } \]
Here $\tensor$ is the tensor of an  ultra-commutative monoid
with a topological space, as explained in Construction \ref{con:umon enrich}
(not to be confused with the objectwise product of an orthogonal space
with a space).
So the cone is $R\rhd [0,1]$,
the tensor of $R$ with the based space $([0,1], 0)$, 
as defined more generally in \eqref{eq:R rhd A}.
Since $R$ is cofibrant and the global model structure is topological, 
the left vertical morphism
is an acyclic cofibration, and so the cone $C R$ 
is globally equivalent to the zero monoid. 

We can then construct a global group completion as a homotopy cofiber
of the diagonal morphism $\Delta:R\to  R\times R$, i.e.,
a pushout in the category of ultra-commutative monoids:
\[ \xymatrix{ 
R \ar[r]^-{\Delta} \ar[d] & R\times R \ar[d]^d\\ 
R\rhd [0,1]\ar[r] & R^\star } \]
The left vertical morphism is induced by $1\in [0,1]$,
and it is a cofibration since $R$ is cofibrant.
We claim that the  composite
\[ i \ : \ R\ \xra{(\Id,0)}\ R\times R \ \xra{\ d\ } \ R^\star \]
is a global group completion.
Indeed, the square above is homotopy cocartesian by construction
and the diagonal morphism has a retraction. 
So Proposition \ref{prop:kernel and cokernel}~(ii) shows that 
the morphism $d:R\times R\to R^\star$ becomes a cokernel of the diagonal
morphism in the homotopy category $\Ho(umon)$.
So $i$ is a global group completion by Proposition \ref{prop:kernel of codiagonal}~(ii).

The previous construction of the global group completion can be
rewritten as two-sided bar construction as follows.
The unit interval $[0,1]=\Delta^1$ is also the topological 1-simplex, 
and hence canonically homeomorphic to the geometric
realization of the simplicial 1-simplex $\Delta[1]$,\index{symbol}{$\Delta[1]$ - {simplicial $1$-simplex}} by
\[ [0,1]\ \to \ |\Delta[1]|\ , \quad t \ \longmapsto \ [ \Id_{[1]}, t] \ .   \]
So Proposition \ref{prop:B bullet and A lhd} 
shows that $C R=R\rhd [0,1]$ is effectively a bar construction
of $R$ with respect to the box product, i.e., the internal realization
of the simplicial object of ultra-commutative monoids 
$B(R,\Delta[1])$, where $\Delta[1]$ is pointed by the vertex~0.
Since $\Delta[1]$ has $(m+1)$ non-basepoint simplices of dimension $m$, 
expanding this yields 
\[ B_m(R,\Delta[1]) \ = \ R\rhd \Delta[1]_m \ \iso \ R^{\boxtimes (m+1)} \ ,  \]
with simplicial structure induced by that of $\Delta[1]$.

Since the internal realization of simplicial ultra-commutative monoids
is a coend, it commutes with pushouts. 
So the defining pushout for $R^\star$ can be rewritten
as the realization of a simplicial ultra-commutative monoid,
the two-sided bar construction
with respect to the box product:
\[ R^\star \ = \ (R\rhd [0,1])\boxtimes_R (R\times R)\ \iso \ 
B(R,\Delta[1])\boxtimes_R (R\times R)\ \iso \ B^\boxtimes(\ast, R, R\times R)\ .\]
In simplicial dimension $m$,
this bar construction is given by
\[ B_m^\boxtimes(\ast, R, R\times R)\ = \ 
\ R^{\boxtimes m }\boxtimes\, (R\times R)\ ;  \]
the simplicial face morphisms are given by
projections away from the first factor (for $d_0$),
the multiplication of $R$ on two adjacent factors (for $d_1,\dots,d_{m-1}$)
respectively the action of $R$ on $R\times R$ through the diagonal.
The simplicial degeneracy morphisms are inserting the unit of $R$.
In the context of topological monoids, this bar construction
of a group completion for `sufficiently homotopy commutative' monoids
is sketched by Segal on p.\,305 of \cite{segal-cat coho}.
\end{construction}

The next proposition is a reality check, showing that global group completion
has the expected effect on equivariant homotopy sets.

\begin{prop}\label{prop:group completion completes pi_0}\index{subject}{global group completion!of an ultra-commutative monoid}
Let $i:R\to R^\star$ be a global group completion of ultra-commutative monoids.
Then for every compact Lie group $G$ the map
\[ \pi_0^G(i)\ : \ \pi_0^G(R)\ \to\ \pi_0^G( R^\star) \]
is an algebraic group completion of abelian monoids and 
\[ \upi_0(i)\ : \ \upi_0(R)\ \to\ \upi_0( R^\star) \]
is a group completion of global power monoids.
\end{prop}
\begin{proof}
We may assume that $R$ is cofibrant in the global model
structure of ultra-commutative monoids of Theorem \ref{thm:global umon}.
As we explained in Construction \ref{con:R^star}, a group completion $R^\star$
can then be constructed as the realization of a certain simplicial 
ultra-commutative monoid,
the two-sided bar construction $B^\boxtimes(\ast, R, R\times R)$, 
where $R$ acts on $R\times R$ through the diagonal morphism.
By Proposition \ref{prop:U preserves realization}, 
the realization can equivalently be taken internal to the category
of ultra-commutative monoids, or in the underlying category of
orthogonal spaces. The `underlying' realization
is the sequential colimit of partial realizations $B^{[n]}$,
i.e., `skeleta' in the simplicial direction, defined as the coend
\[ B^{[n]}\ = \ \int^{[m]\in\bDelta_\leq n} B^\boxtimes_m(\ast, R, R\times R)
\times \Delta^m \]
of the restriction to the full subcategory $\bDelta_{\leq n}$
of $\bDelta$ with objects all $[m]$ with $m\leq n$.
Since $\bDelta_{\leq n}$ is contained in $\bDelta_{\leq n+1}$,
there is a canonical morphism 
$B^{[n]}\to B^{[n+1]}$ and the realization $B^\boxtimes(\ast, R, R\times R)$ 
is the colimit of the sequence of orthogonal spaces
\[ R\times R = B^{[0]} \ \to \ B^{[1]}\ \to \ \cdots \ 
\to \ B^{[n]}\ \to \ \cdots \ .  \]
The 1-skeleton $B^{[1]}$ is the pushout of orthogonal spaces:
\[\xymatrix@C=12mm{ 
 R\boxtimes(R\times R)\times\{0,1\} \ar[r]^-{\text{incl}}\ar[d] &
 R\boxtimes(R\times R)\times [0,1]\ar[d] \\
R\times R\ar[r] &B^{[1]}
} \]
The orthogonal space $R\boxtimes(R\times R)\times\{0,1\}$
is the disjoint union of two copies of
$R\boxtimes(R\times R)$, and the left vertical map is projection to $R\times R$
on one copy; on the other copy the morphism has the two components
\[ \mu\circ(R\boxtimes p_1)\ ,\ \mu\circ(R\boxtimes p_2)\ : 
R\boxtimes(R\times R)\ \to \ R\ ,\]
where $\mu:R\boxtimes R\to R$ is the multiplication
and $p_1,p_2:R\times R\to R$ are the two projections.
The functor $\pi_0^G$ takes both the box product and the product
of orthogonal spaces to products of sets, by Corollary \ref{cor-pi_0 of box}.
So for a compact Lie group $G$, the set $\pi_0^G(B^{[1]})$ is the coequalizer 
of the two maps
\[ \xymatrix{
\pi_0^G(R)\times\pi_0^G(R)\times\pi_0^G(R)\ \ar@<.4ex>[r]^-a \ar@<-.4ex>[r]_-b &
\ \pi_0^G(R)\times\pi_0^G(R)} \]
given by
\[ a(x,y,z)\ = \ (y,z)\text{\qquad respectively\qquad}
 b(x,y,z)\ = \ (x+y,x+z)\ .  \]
The equivalence relation on the set $\pi_0^G(R)\times \pi_0^G(R)$
generated by declaring $a(x,y,z)$ equivalent to $b(x,y,z)$ for all
$(x,y,z)\in\pi_0^G(R)$ is precisely the equivalence relation
that constructs the Grothendieck group of $\pi_0^G(R)$,
compare Example \ref{eg:Grothendieck group}. 
So $\pi_0^G(B^{[1]})$ bijects with the algebraic group completion of $\pi_0^G(R)$.
For $n\geq 1$ the passage from $B^{[n]}$ to $B^{[n+1]}$ involves
attaching simplices of dimension at least~2 along their boundary,
and this process does not change the path components of $G$-fixed points
of the values at any $G$-representation. So the morphism
$B^{[1]}\to B^\boxtimes(\ast, R, R\times R)$ induces a bijection on $\pi_0^G$.
This proves that the morphism $R\to B^\boxtimes(\ast, R, R\times R)=R^\star$
is a group completion of abelian monoids.

The second claim then follows because group completions 
of global power monoids are calculated `group-wise', compare Example \ref{eg:M^star}.
\end{proof}

For topological monoids, the loop space of the bar construction
(see Construction \ref{con:bar construction})\index{subject}{bar construction!of a topological monoid}
provides a functorial group completion. We will now explain that a similar construction
also provides global group completion for ultra-commutative monoids,
before passing to the homotopy category.
Part of this works for arbitrary orthogonal monoid spaces, not necessarily commutative.
We refer to Construction \ref{con:realize simplicial} for generalities
about realization of simplicial objects, in particular simplicial orthogonal spaces.

If $R$ is an orthogonal monoid space, 
then the {\em bar construction}\index{subject}{bar construction!of an orthogonal monoid space}
is the simplicial object of orthogonal spaces
\[
 B_\bullet(R)\ = \ \left([n]\mapsto  \ R^{\boxtimes n} \right) \ .  
\]
The simplicial face morphisms are induced by the multiplication in $R$, and
the degeneracy morphisms are induced by the unit morphism of $R$,
much like for the bar construction with respect to cartesian product
(as opposed to box product) in Construction \ref{con:bar construction}.
The geometric realization in the category of orthogonal spaces
is then the orthogonal space
\begin{equation}\label{eq:bar for umon}
 B(R) \ = \ | B_\bullet(R) | \ .
\end{equation}
Geometric realization of orthogonal spaces is `objectwise',
i.e., for an inner product space $V$ we have
\[ B(R)(V) \  = \ |B_\bullet(R) (V) |\ ,\]
the realization of the simplicial space $[n]\mapsto R^{\boxtimes n}(V)$.

The next proposition shows that the bar construction of orthogonal
monoid spaces preserves global equivalences under a mild non-degeneracy
condition on the unit.

\begin{defn}
  An orthogonal monoid space $R$ has a {\em flat unit}\index{subject}{flat unit!of an orthogonal monoid space}
  if the unit morphism $\ast\to R$ is a flat cofibration of orthogonal spaces.
\end{defn}

The condition that the unit morphism $\ast\to R$ is a flat cofibration
is equivalent to the requirements that the underlying orthogonal space of $R$
is flat and the unit map $\ast\to R(0)$ is a cofibration of spaces.

We recall from Definition \ref{def:Reedy flat}
that a simplicial orthogonal space $X$ is {\em Reedy flat}\index{subject}{Reedy flat!orthogonal space}
if the latching morphism $l_n^\Delta:L_n^\Delta(X)\to X_n$ in the simplicial
direction is a flat cofibration of orthogonal spaces for every $n\geq 0$.

\begin{prop}\label{prop:bar global invariance}
  \begin{enumerate}[\em (i)]
  \item 
    For every orthogonal monoid space $R$ with flat unit
    the simplicial orthogonal space $B_\bullet(R)$ is Reedy flat.
  \item 
    Let $f:R\to S$ be a morphism of orthogonal monoid spaces with flat units.
    If $f$ is a global equivalence, so is the morphism $B(f):B(R)\to B(S)$.  
  \end{enumerate}\end{prop}
\begin{proof}
  (i) The $n$-th latching morphism $L_n^\Delta(B_\bullet(R))\to B_n(R)$ in the simplicial
  direction is the iterated pushout product
  \[ i^{\Box n} \ : \ Q^n(i)\ \to \ R^{\boxtimes n}\]
  with respect to the unit morphism $\ast\to R$. Since this unit morphism
  is a flat cofibration, the pushout product property of the flat cofibrations
  shows that $i^{\Box n}$, and hence the latching morphism,
  is a flat cofibration for all $n\geq 0$.

  (ii) Since $R$ and $S$ have flat units, the simplicial orthogonal spaces
  $B_\bullet(R)$ and $B_\bullet(S)$ are Reedy flat by part~(i).
  Moreover, the morphism $B_n(f)=f^{\boxtimes n}:R^{\boxtimes n}\to S^{\boxtimes n}$
  is a global equivalence since $f$ is and because
  the box product is homotopical for global equivalences 
  (by Theorem \ref{thm:box to times}).
  So the claim follows from the global invariance of realizations
  between Reedy flat simplicial orthogonal spaces 
  (Proposition \ref{prop:realization invariance in spc}~(ii)). 
\end{proof}

\begin{rk}[Comparing bar constructions]
We let $M$ be a monoid valued orthogonal space 
in the sense of Definition \ref{def:monoid valued ospace}.\index{subject}{orthogonal space!monoid valued} 
Then we have two bar constructions available: on the one hand we can
take the bar construction 
objectwise as in Example \ref{eg:general bar construction},
resulting in the orthogonal space $\mathbf B^\circ M$.\index{subject}{bar construction!of a monoid valued orthogonal space}
On the other hand, we can first pass to the associated orthogonal monoid space
as in \eqref{eq:monoid_times2box}, and then perform the bar construction
with respect to the $\boxtimes$-multiplication as in \eqref{eq:bar for umon}. 
There is a natural comparison map: the symmetric monoidal transformation
$\rho_{X,Y} : X\boxtimes Y\to X\times Y$   
defined in \eqref{eq:define_box2times} has an analog 
for any finite number of factors, and for varying $n$, the morphisms
\[ \rho_{M,\dots,M}\ : \ M^{\boxtimes n}\ \to \ M^n  \]
form a morphism of simplicial orthogonal spaces
\[ \rho_\bullet\ : \ B_\bullet^\boxtimes(M)\ \to \ B^\times_\bullet(M) \]
from the $\boxtimes$-bar construction to the $\times$-bar construction.
If $M$ has a flat unit, then the simplicial orthogonal spaces
$B_\bullet^{\boxtimes}(M)$ and $B_\bullet^{\times}(M)$ are both Reedy flat. 
Indeed, for the former, this is 
Proposition \ref{prop:bar global invariance}~(i),
and for the latter the same proof works because the categorical
product of orthogonal spaces also satisfies the pushout product property
with respect to flat cofibrations 
(Proposition \ref{prop:ppp for cartesian product}).
So $\rho_\bullet$ is a morphism between Reedy flat simplicial orthogonal spaces that 
is a global equivalence
in every simplicial degree (by Theorem \ref{thm:box to times}~(i)).
The induced morphism $|\rho_\bullet|:B(M)\to\mathbf B^\circ M$ 
between the realizations is then a global equivalence 
by Proposition \ref{prop:realization invariance in spc}~(ii).
\end{rk}

The canonical morphism
\[ R\times [0,1] \ = \  B_1(R)\times\Delta^1 \ \to \ 
 | B_\bullet(R) |\ = \ B(R) \]
takes $R\times \{0,1\}$ to the basepoint, so it factors over a morphism
of orthogonal spaces
\[  R\sm ([0,1]/\{0,1\}) \ \to \  B(R)\ .\]
We let $u :S^1 \to[0,1]/\{0,1\}$ be the composite homeomorphism
\begin{equation}  \label{eq:simplex u}
   S^1 \ \xra{\ c \ }\  U(1)\ \xra{\log}\ [0,1]/\{0,1\}\ ,
\end{equation}
of the Cayley transform\index{subject}{Cayley transform} 
\[ c \ : \ S^1 \ \to \ U(1) \ , \quad x \ \longmapsto \ (x+i)(x-i)^{-1} \ ,\]
and the logarithm, i.e., the inverse of the exponential homeomorphism
\[ [0,1]/\{0,1\}\ \iso \ U(1)\ , \quad t \ \longmapsto \ e^{2\pi i t}\ . \]
This yields a composite morphism
\[  R\sm S^1 \ \xra[\iso]{R\sm u} \  R\sm ([0,1]/\{0,1\}) \ \to \  B(R) \]
which is adjoint to a morphism of orthogonal spaces
\begin{equation}\label{eq:define_bar_eta}
 \eta_R \ : \  R \ \to \  \Omega B(R)\ .  
\end{equation}

\begin{prop}\label{prop:Omega B deloops}
  Let $R$ be an orthogonal monoid space with flat unit.
  If for every compact Lie group $G$ the monoid $\pi_0^G(R)$ is a group, 
  then the morphism $\eta_R : R\to\Omega B(R)$ is a global equivalence.
\end{prop}
\begin{proof}
  Since $R$ has a flat unit, the simplicial orthogonal space $B_\bullet(R)$
  is Reedy flat by Proposition \ref{prop:bar global invariance}~(i),
  so the underlying orthogonal space of $B(R)$ is flat by 
  Proposition \ref{prop:realization invariance in spc}~(i). 
  As a flat orthogonal space, $B(R)$ is in particular closed.
  The loop space functor preserves closed inclusions 
  by \cite[Prop.\,7.7]{lewis-thesis},
  so the pointwise loop space $\Omega B(R)$ is also closed as an orthogonal space.
  Since $R$ and $\Omega B(R)$ are both closed orthogonal spaces,
  we can detect global equivalences on $G$-fixed points,
  see Proposition \ref{prop:global eq for closed}.

  So we let $G$ be a compact Lie group. Then
  $((\Omega B(R))(\Uc_G))^G$ is homeomorphic to
  $ \Omega\left( ( B(R)(\Uc_G))^G \right)$, which is in turn homeomorphic to
  the loop space of the geometric realization of the simplicial space
  \begin{equation} \label{eq:simplicial_space_for_Segal}
     [n]\ \longmapsto\  (R^{\boxtimes n}(\Uc_G))^G \ .     
  \end{equation}
  We define $i_k:[1]\to [n]$ by $i_k(0)=k-1$ and $i_k(1)=k$.
  Then the morphism 
  \[ (i_1^*,\dots,i_n^*)\ : \ R^{\boxtimes n}\ = \ B_n(R)\ \to \
  (B_1(R))^n \ = \  R^n \]
  is precisely the morphism $\rho_{R,\dots,R}:R^{\boxtimes n}\to R^n$,
  and hence a global equivalence by Theorem \ref{thm:box to times}~(i). 
  Since $\pi_0^G(R)$ is a group, 
  the simplicial space \eqref{eq:simplicial_space_for_Segal}
  satisfies the hypotheses of Segal's theorem \cite[Prop.\,1.5]{segal-cat coho};
  so the adjoint of the canonical map
  \[  R(\Uc_G)^G\sm S^1 \ \to \ 
  |[n]\mapsto ((R^{\boxtimes n})(\Uc_G))^G| \ \iso \   ( B(R)(\Uc_G) )^G \]
  is a weak equivalence. Here we have used again that fixed points commute
  with geometric realization, 
  see Proposition \ref{prop:G-fix preserves pushouts}~(iv).
  This adjoint is precisely the underlying map 
  of $G$-fixed points of the morphism $\eta_R:R\to \Omega B(R)$.
\end{proof}

The previous proposition works for general orthogonal monoid spaces,
not necessarily commutative; in that generality the bar construction
$B(R)$ is an orthogonal space, but it does not have any natural multiplication. 
When we apply the bar construction to ultra-commutative monoids,
then something special happens: since the multiplication morphism
$\mu:R\boxtimes R\to R$ is then a homomorphism of ultra-commutative monoids,
$B_\bullet(R)$ is a simplicial object in the category of
ultra-commutative monoids, i.e., a simplicial ultra-commutative monoid.
The geometric realization $B(R)$ is then canonically an
ultra-commu\-tative monoid, and it coincides with the realization of 
$B_\bullet(R)$ internal to the category of ultra-commutative monoids,
compare Proposition \ref{prop:U preserves realization}.

Moreover, we claim that for ultra-commutative monoids,
the bar construction $B(R)$ is naturally isomorphic to $R\rhd S^1$,
the `suspension' of $R$ internal to the category of
ultra-commutative monoids.
To see this we consider the `simplicial circle' $\mathbf S^1$,\index{subject}{simplicial circle}
the simplicial set given by
\[ (\mathbf S^1)_n \ = \ \{0,1,\dots,n\} \ ,  \]
with face maps $d_i:(\mathbf S^1)_n\to (\mathbf S^1)_{n-1}$ given by
\[ d_i(j) \ = \
\begin{cases}
j-1 & \text{ for $i <j$, and}\\  
\quad j & \text{ for $i \geq j$ and $j\ne n$,}\\
\quad 0 & \text{ for $i=j=n$,}
\end{cases}\]
and degeneracy maps $s_i:(\mathbf S^1)_n\to(\mathbf S^1)_{n+1}$ given by
\[
 s_i(j) \ = \
\begin{cases}
j+1 & \text{ for $i  <j$, and}\\  
\quad j & \text{ for $i\geq j$.}
\end{cases}\]
The simplicial set $\mathbf S^1$ is based by~0; it is isomorphic to 
the simplicial 1-simplex modulo its boundary,
and its realization is homeomorphic to a circle, whence the name.

The  `obvious' isomorphisms
\[ p_n \ : \ R^{\boxtimes n}\ \xra{\ \iso \ }\  
 R\rhd \{0,1,\dots,n\} \ = \ R\rhd (\mathbf S^1)_n \ = \ B_n(R,\mathbf S^1)\ , \]
are compatible with the simplicial structure maps as $n$ varies,
so they define an isomorphism of simplicial ultra-commutative monoids
\[ p_\bullet \ : \ B_\bullet(R)\ \iso \ B_\bullet(R,\mathbf S^1) \ .\]
When we specialize Proposition \ref{prop:B bullet and A lhd} 
to $A=\mathbf S^1$, we obtain an isomorphism of ultra-commutative monoids
\[  R\rhd |\mathbf S^1|\iso \  B(R)  \ . \]
The homeomorphism $u:S^1\to[0,1]/\{0,1\}$ from \eqref{eq:simplex u} 
and the homeomorphism
\[ [0,1] / \{0,1\}\ \xra{\ \iso \ }\  |\mathbf S^1|  \ , \quad
t \ \longmapsto \ [1,t]\]
turn this into an isomorphism of ultra-commutative monoids
\begin{equation}  \label{eq:umon_bar_and_suspension}
  R\rhd S^1\ \xra[\iso]{R\rhd u} \ R\rhd ( [0,1] / \{0,1\} ) \ \iso \  
 \ R\rhd |\mathbf S^1| \ \iso \  B(R)    
\end{equation}
whose adjoint $R\to \Omega B(R)$ is the morphism $\eta_R$
of \eqref{eq:define_bar_eta}.

\begin{cor}\label{cor:umon Omega Sigma is group completion}
For every ultra-commutative monoid with flat unit $R$ the adjunction unit
\[ \eta_R \ : \ R\ \to \ \Omega(R\rhd S^1) \]
is a global group completion.\index{subject}{flat unit!of an ultra-commutative monoid}
\end{cor}
\begin{proof}
We let $R$ be a group-like cofibrant ultra-commutative monoid.
Then $R$ has a flat unit by Theorem \ref{thm:global umon}~(ii a).
Since ultra-commutative monoids form a topological model category,
$R\rhd S^1$ is an abstract suspension of $R$.
The isomorphism \eqref{eq:umon_bar_and_suspension}
transforms the adjunction unit $ R\to\Omega(R\rhd S^1)$
into the morphism $\eta_R:R\to \Omega B(R)$ defined in \eqref{eq:define_bar_eta}. 
Proposition \ref{prop:Omega B deloops} shows
that this adjunction unit is a global equivalence 
for every cofibrant group-like ultra-commutative monoid $R$.
In the homotopy category $\Ho(umon)$ this implies that for
every group-like ultra-commutative monoid $R$
the derived adjunction unit $\eta:R\to \Omega(\Sigma R)$ 
is an isomorphism.
So Proposition \ref{prop:Omega Sigma group completion}
shows that for every ultra-commutative monoid $R$
the derived adjunction unit $\eta:R\to \Omega(\Sigma R)$ 
is a group completion in the pre-additive category $\Ho(umon)$.
For cofibrant $R$ the pointset level adjunction unit
$R\to\Omega(R\rhd S^1)$ realizes the derived unit,
hence the claim follows for every ultra-commutative monoid that is
cofibrant in the global model structure of Theorem \ref{thm:global umon}.

In the general case we choose a cofibrant replacement $q:R^c\to R$ 
in the global model structure of Theorem \ref{thm:global umon},
i.e., a global equivalence of ultra-commutative monoids with cofibrant source.
Since $R^c$ and $R$ have flat units, the induced morphism of bar constructions
$B(q):B(R^c)\to B(R)$ is a global equivalence
by Proposition \ref{prop:bar global invariance}~(ii).
Hence the morphism $q\rhd S^1:R^c\rhd S^1\to R\rhd S^1$ is a global equivalence, and so 
is $\Omega(q\rhd S^1):\Omega(R^c\rhd S^1)\to \Omega(R\rhd S^1)$.
The morphism $\eta_{R^c}:R^c\to\Omega(R^c\rhd S^1)$ is a global group
completion by the previous paragraph. Since $\eta_R:R\to\Omega(R\rhd S^1)$ 
is isomorphic to the morphism $\eta_{R^c}$ in the homotopy category $\Ho(umon)$,
the morphism $\eta_R$ is also a global group completion.
\end{proof}

An example of a global group completion that comes up naturally 
is the morphism $i:\bGr\to\bBOP$ introduced in Example \ref{eg:Gr to BOP}. 
The verification of the group completion property
will be done through a homological criterion.
For that purpose we define the homology groups of an orthogonal space $Y$ as
\[ H_*( Y^G;\mZ) \ = \ \colim_{V\in s(\Uc_G)} \ H_*( Y(V)^G;\mZ) \ . \]
Every global equivalence induces isomorphisms on $H_*((-)^G;\mZ)$
for all compact Lie groups $G$. Indeed, the functor $H_*((-)^G;\mZ)$
takes strong level equivalences to isomorphisms, which reduces the claim
(by cofibrant approximation in the strong level model structure)
to global equivalences $f:X\to Y$ between flat orthogonal spaces. Flat orthogonal
spaces are closed, so the global equivalence induces weak equivalences
$f(\Uc_G)^G:X(\Uc_G)^G\to Y(\Uc_G)^G$ on $G$-fixed points.
The poset $s(\Uc_G)$ is filtered, so homology commutes with this colimit, i.e.,
\[ H_*( Y^G;\mZ) \ \iso \  H_*( Y(\Uc_G)^G;\mZ) \ . \]
Thus the morphism $f$ also induces an isomorphism on $H_*((-)^G;\mZ)$. 

The multiplication of an orthogonal monoid space $R$
induces a graded multiplication on the homology groups
$H_*(R^G;\mZ)$, by simultaneous passage to colimits in both variables
of the maps
\begin{align*}
    H_m( R(V)^G;\mZ) \tensor H_n( R(W)^G;\mZ) \ \xra{\ \times\ }\ 
 &H_{m+n}( R(V)^G\times  R(W)^G;\mZ)  \\ 
\xra{( (\mu_{V,W})^G)_*}\  &H_{m+n}( R(V\oplus W)^G;\mZ) \ . 
\end{align*}
Assigning to a path component its homology class is a map
\[ \pi_0(R(V)^G) \ \to \ H_0(R(V)^G;\mZ)\]
compatible with increasing $V$. On colimits over $s(\Uc_G)$
this provides a map
\[ \pi_0^G(R) \ \to \ H_0(R^G;\mZ)\ .\]
This map takes the addition in $R$ to the multiplication in $ H_0(R^G;\mZ)$,
so its image is a multiplicative subset of $H_0(R^G;\mZ)$.
If the multiplication of $R$ is commutative, then the  
product of $H_*( R^G;\mZ)$ is commutative in the graded sense.
In particular, the multiplicative subset of $\pi_0^G(R)$
is then automatically central.

\begin{prop}\label{prop:localization criterion}
A morphism $i:R\to R^\star$ of ultra-commutative monoids is a 
global group completion if and only if the following two conditions
are satisfied.
\begin{enumerate}[\em (i)]
\item The ultra-commutative monoid $R^\star$ is group-like, and
\item for every compact Lie group $G$ the map of graded commutative
rings
\[ H_*(i^G;\mZ)\ : \ H_*( R^G;\mZ) \ \to \ H_*( (R^\star)^G;\mZ)  \]
is a localization at the multiplicative subset $\pi_0^G(R)$ of
$H_0(R^G;\mZ)$.
\end{enumerate}
\end{prop}
\begin{proof}
We start by showing that a global group completion satisfies properties~(i)
and~(ii). Property~(i) holds by definition of `group completion'. 
We give two alternative proofs for why a global group completion
satisfies property~(ii), based on the two different bar construction models
in Construction \ref{con:R^star} respectively 
Corollary \ref{cor:umon Omega Sigma is group completion}.

The first argument uses the loop space of the bar construction $B(R)$,
which is isomorphic to internal suspension $R\rhd S^1$.
By Corollary \ref{cor:umon Omega Sigma is group completion} 
it suffices to show that for every cofibrant
ultra-commutative monoid $R$ the morphism $\eta_R:R\to\Omega(R\rhd S^1)=\Omega(B R)$
has property~(ii). 
Since $R$ is cofibrant and the global model structure
of ultra-commutative monoids is topological 
(Theorem \ref{thm:global umon}), 
the basepoint inclusion of $S^1$ induces a cofibration
\[ R\tensor \text{incl}\ : \  R\tensor \{\infty\} \ \to \ R\tensor S^1 .\]
The cobase change is the unique morphism $\ast\to R\rhd S^1$ from the terminal
ultra-commutative monoid to the reduced tensor, and this is thus a cofibration.
In other words, $B(R)=R\rhd S^1$ is again cofibrant as an ultra-commutative monoid.
Since $R$ and $B(R)$ are cofibrant as ultra-commutative monoids,
Theorem \ref{thm:global umon}~(ii) shows that their underlying orthogonal spaces
are flat, hence closed.
The loop space functor preserves closed inclusions by \cite[Prop.\,7.7]{lewis-thesis},
so the pointwise loop space $\Omega B (R)$ is also closed as an orthogonal space.

Since $R$ and $\Omega B(R)$ are closed as orthogonal spaces,
it suffices to show that for every compact Lie group $G$ the map
\[ H_*( R(\Uc_G)^G;\mZ) \ \to \ H_*( (\Omega B (R))(\Uc_G)^G;\mZ)  \]
is a localization at the multiplicative subset $\pi_0(R(\Uc_G)^G)$ of
the source. Since the H-space structure of $R(\Uc_G)^G$ comes from
the action of an $E_\infty$-operad, the graded ring 
$H_*( R(\Uc_G)^G;\mZ)$ is graded commutative.
We can thus apply Quillen's group completion theorem from the unpublished, 
but widely circulated preprint {\em `On the group completion of a simplicial monoid'}.
Quillen's manuscript was later published as
Appendix~Q of the Friedlander-Mazur paper \cite{friedlander-mazur},
where the relevant theorem appears on page 104 in Section~Q.9.

The second, alternative, argument first reduces to cofibrant ultra-commuta\-tive monoids
by cofibrant approximation in the global model of Theorem \ref{thm:global umon}.
As explained in Construction \ref{con:R^star}, 
a group completion $R^\star$ can then be constructed
as the homotopy cofiber of the diagonal $\Delta:R\to R\times R$,
which is concretely given by the geometric realization of 
the simplicial ultra-commutative monoid $B^\boxtimes(\ast, R, R\times R)$, 
the bar construction with respect to the box product.
Segal \cite[p.\,305 f]{segal-cat coho} sketches an argument why the 
resulting morphism $R\to B^\boxtimes(\ast, R, R\times R)$ is localization on
homology with field coefficients. The argument is reproduced in more
detail in the proof of \cite[Lemma 3.2.2.1]{dundas-goodwillie-mccarthy}.

Now we prove the reverse implication.
We let $i:R\to R^\star$ be a morphism of ultra-commutative monoids that 
satisfies properties~(i) and~(ii); we need to show that $i$ is a global group 
completion.
We assume first that both $R$ and $R^\star$ are cofibrant 
in the global model structure of ultra-commutative monoids 
of Theorem \ref{thm:global umon}.
Then the unit morphisms of $R$ and $R^\star$ are flat cofibrations
of underlying orthogonal spaces by Theorem \ref{thm:global umon}~(ii).
So the morphism $\eta_R:R\to \Omega(R\rhd S^1)=\Omega B(R)$
is a global group completion by Corollary \ref{cor:umon Omega Sigma is group completion}.
Since $R^\star$ is group-like, the morphism $\eta_{R^\star}$
is a global equivalence by Proposition \ref{prop:Omega B deloops}. 

We claim that the morphism $B (i):B(R)\to B( R^\star)$ is a global equivalence.
For every coefficient system $L$ on $(B (R^\star))(\Uc_G)^G$ 
we compare the spectral sequence
\[ E^2_{p,q}\ = \ \Tor_p^{H_*(R^G;k)}(k,L) \ \Longrightarrow \
H_*( (B (R))(\Uc_G)^G;L)\]
(obtained by filtering the bar construction by simplicial skeleta)
with the analogous one for the homology of $(B (R^\star))(\Uc_G)^G$.
The localization hypothesis implies that the map of Tor groups
\[\Tor_p^{H_*(R^G;k)}(k,L) \ \to  \ 
\Tor_p^{H_*((R^\star)^G;k)}(k,L) \]
is an isomorphism, 
see for example \cite[Prop.\,7.17]{rotman-homological algebra}
or \cite[Prop.\,3.2.9]{weibel-homological algebra}.
So we have a morphism of first quadrant spectral sequences
that is an isomorphism of $E^2$-terms; the map on abutments
is then an isomorphism as well.
This shows that $B(i):B (R)\to B (R^\star)$ is a global equivalence.
Since looping preserves global equivalences, the morphism 
$\Omega B (i)$ is a global equivalence.
Now we contemplate the commutative square
\[ \xymatrix@C=12mm{
R \ar[r]^-i\ar[d]_{\eta_R} & R^\star \ar[d]^{\eta_{R^\star}}\\
\Omega B (R) \ar[r]_-{\Omega B(i)} & \Omega B (R^\star) } \]
The left vertical morphism is a global group completion,
and the lower horizontal and right vertical morphisms are
global equivalences. So the upper horizontal morphism $i$ is global group completion.

Now we reduce the general case to the special case by cofibrant
replacement.
We choose a cofibrant replacement $q:R^c\to R$ in the global model structure of
ultra-commutative monoids of Theorem \ref{thm:global umon},
and then factor the morphism $i q:R^c\to R^\star$ as a cofibration
$i^c:R^c\to R^\dagger$ followed by a global equivalence $\varphi:R^\dagger\to R^\star$.
Properties~(i) and~(ii) are invariant under global equivalences of pairs,
so the morphism $i^c:R^c\to R^\dagger$ satisfies~(i) and~(ii).
Since $R^c$ and $R^\dagger$ are both cofibrant, the morphism $i^c$
is a global group completion by the special case above.
So the morphism $i$ is also a global group completion.
\end{proof}

We showed in Theorem \ref{thm:upi_0 of BOP} that the ultra-commutative 
monoid $\bBOP$ is group-like and that its equivariant homotopy sets
$\upi_0(\bBOP)$ realize the orthogonal representation rings additively.
In Example \ref{eg:Gr to BOP} we introduced a morphism $i:\bGr\to\bBOP$
of ultra-commutative monoids from the additive Grassmannian
and showed in Proposition \ref{prop:i is group completion}
that for every compact Lie group $G$ and every $G$-space $A$, the homomorphism
\[ [A,i]^G \ : \ [A,\bGr]^G \ \to \ [A,\bBOP]^G \]
is a group completion of abelian monoids.
In particular, the map $\pi_0^G(i):\pi_0^G(\bGr)\to\pi_0^G(\bBOP)$
is an algebraic group completion.
In much the same way we can define morphisms of ultra-commutative monoids
\[ i\ :\ \bGr^\mC\ \to\ \bBUP\text{\qquad and\qquad}i\ :\ \bGr^\mH\ \to\ \bBSpP\]
by replacing $\mR$-subspaces in $V$ by $\mC$-subspaces in $V_\mC$, 
respectively $\mH$-subspaces in $V_\mH$.

\begin{theorem}\label{thm:BOP is group completion}
The morphisms $i:\bGr\to\bBOP$, $i:\bGr^\mC\to\bBUP$
and $i:\bGr^\mH\to\bBSpP$ are global group completions  
of ultra-commutative monoids.
\end{theorem}
\begin{proof}
We prove the real case in detail and leave the complex and quaternionic cases
to the reader.
We verify the localization criterion of Proposition \ref{prop:localization criterion}.
To this end we define a bi-orthogonal space, i.e., a functor
\[ \bGr^\sharp \ : \ \bL\times \bL\ \to \ \bT\]
on objects by
\[ \bGr^\sharp(U,V)\ = \ \bGr(U\oplus V) \ .\]
For linear isometric embeddings $\varphi:U\to \bar U$ and $\psi:V\to \bar V$, 
the induced map is
\[ \bGr^\sharp(\varphi,\psi)\ : \ \bGr^\sharp(U,V)\ \to \ \bGr^\sharp(\bar U,\bar V) \ , \quad
L \ \longmapsto (\varphi\oplus\psi)(L) \ + \ ((\bar U-\varphi(U))\oplus 0)\ .\]
We emphasize that the behavior on objects is {\em not} symmetric in the
two variables, and in the first variable it is not just applying $\varphi$.

Now we fix a compact Lie group $G$ and consider the colimit
of the bi-orthogonal space $\bGr^\sharp$  over the poset $s(\Uc_G)\times s(\Uc_G)$.
Since the diagonal is cofinal in the poset $s(\Uc_G)\times s(\Uc_G)$,
this `double colimit' is also a colimit over the restriction to the diagonal.
But the diagonal of $\bGr^\sharp$ is precisely the orthogonal space $\bBOP$, and so
\[
 \colim_{(U,V)\in s(\Uc_G)^2} \bGr^\sharp(U,V)\ = \ \colim_{W\in s(\Uc_G)} \bBOP(W) \ = \ \bBOP(\Uc_G)\ .  
  \]
On the other hand, if we fix an inner product space $U$ as the
first variable, then $\bGr^\sharp(U,-)$ is isomorphic to
the additive $U$-shift 
(in the sense of Example \ref{eg:Additive and multiplicative shift})
of the Grassmannian $\bGr$. Hence for fixed $U$,
\[ \colim_{V\in s(\Uc_G)} \bGr^\sharp(U,V)\ = \ \bGr(U\oplus\Uc_G)\ .    \]
A colimit over $s(\Uc_G)\times s(\Uc_G)$ can be calculated in two steps,
first in one variable and then in the other, so we conclude that
\begin{equation}\label{eq:BOP as colimit}
 \bBOP(\Uc_G)\ = \ \colim^\sharp_{U\in s(\Uc_G)} \, \bGr(U\oplus\Uc_G)\ ;
\end{equation}
under this identification, the map $i(\Uc_G):\bGr(\Uc_G)\to\bBOP(\Uc_G)$ 
becomes the canonical morphism 
\[  i^\sharp \ : \ \bGr(\Uc_G) \ \to \ 
\colim^\sharp_{U\in s(\Uc_G)}\, \bGr(U\oplus\Uc_G) \]
to the colimit, for $U=0$.\medskip

\Danger The decoration~`$\sharp$' is meant to emphasize that the structure maps 
in this colimit system come from the functoriality of $\bGr^\sharp$ 
in the first variable, so they are {\em not}
the maps obtained by applying $\bGr(-\oplus\Uc_G)$ to an inclusion $U\subset\bar U$.
For example, the maps in the colimit \eqref{eq:BOP as colimit}
do {\em not} preserve the $\mN$-grading by dimension.
So one should not confuse the colimit \eqref{eq:BOP as colimit}
with the space $\bGr(\Uc_G\oplus\Uc_G)$, 
which is $G$-homeomorphic to $\bGr(\Uc_G)$ by a choice
is equivariant linear isometry $\Uc_G\oplus\Uc_G\iso\Uc_G$.

\medskip

We claim that the map
\begin{equation}\label{eq:homology comparison}
 H_*((i^\sharp)^G) \ : \ H_*(\bGr(\Uc_G)^G)\ \to \ 
H_*\big( \colim^\sharp_{U\in s(\Uc_G)}\bGr(U\oplus\Uc_G)^G \big)
\end{equation}
is a localization at the multiplicative subset $\pi_0^G(\bGr)$,
where homology stands for singular homology with integer coefficients.
To see this we observe that all the maps in the colimit system are
closed embeddings; so singular homology commutes with this particular colimit.

For $U\in s(\Uc_G)$ we denote by $j_U:\bGr(\Uc_G)^G\to\bGr(U\oplus\Uc_G)^G$
the map induced by applying the direct summand inclusion $\Uc_G\to U\oplus\Uc_G$.
The map $j_U$ is a homotopy equivalence because $\Uc_G$ is a complete $G$-universe.
For all $U\subset V$ in $s(\Uc_G)$ the following square commutes
\[ \xymatrix@C=15mm{
H_*(\bGr(\Uc_G)^G) \ar[r]^-{\cdot [V-U]} \ar[d]^\iso_{H_*(j_U)} &
H_*(\bGr(\Uc_G)^G)\ar[d]_\iso^{H_*(j_V)} \\
H_*(\bGr(U\oplus\Uc_G)^G) \ar[r]_-{H_*(i_{U,V}^\sharp)} & H_*(\bGr(V\oplus\Uc_G)^G)
} \]
and the vertical maps are isomorphisms.
So the target of \eqref{eq:homology comparison} is the colimit 
of the functor on $s(\Uc_G)$ that takes all objects to the
ring $H_*(\bGr(\Uc_G)^G;\mZ)$ and an inclusion $U\subset V$
to multiplication by the class $[V-U]$ in the multiplicative subset
under consideration. Hence the map \eqref{eq:homology comparison}
is indeed a localization as claimed. Since the ultra-commutative monoid $\bBOP$
is group-like, the criteria of Proposition \ref{prop:localization criterion}
are satisfied, and so the morphism $i:\bGr\to\bBOP$ is a global group
completion.
\end{proof}

\begin{rk}\label{rk:completion makes ucom}
We had earlier defined an $E_\infty$-orthogonal monoid space $\bBOP'$ as a mixture of
$\bbOP$ and $\bBOP$: the value at an inner product space $V$ is
\[  \bBOP'(V) \ = \  {\coprod}_{m\geq 0}\, Gr_m(V^2\oplus\mR^\infty)\ . \]
The structure maps and an $E_\infty$-multiplication can be defined 
in essentially the same way as for $\bBO'$,
which was defined in \eqref{eq:define_BO'}.
So $\bBOP'$ becomes the $\mZ$-graded periodic analog
of the orthogonal space $\bBO'$.
In the same way as for the homogeneous degree~0 summands in \eqref{eq:a_and_b_morphisms},
we defined two morphisms
of $E_\infty$-orthogonal monoid spaces
\[ \bbOP \ \xra{\ a\ }\ \bBOP' \ \xla[\simeq]{\ b\ } \ \bBOP\ . \]
The same arguments as in Proposition \ref{eq:bbO_filtration}
show that the morphism $b$ is a global equivalence.
A very similar argument as in Proposition \ref{thm:BOP is group completion}
shows that the morphism $a$
is a global group completion in the homotopy category of
$E_\infty$-orthogonal monoid spaces. Strictly speaking we would first have
to justify that the homotopy category is pre-additive (which we won't do),
so that the formalism of group completions applies.

As we argued in Proposition \ref{prop:pi bbO inside IO}, the $E_\infty$-structure
on $\bbO$ {\em cannot} be refined to an ultra-commutative multiplication.
The argument was based on an algebraic obstruction that exists 
in the same way in $\upi_0(\bbOP)$, so $\bbOP$ 
cannot be refined to an ultra-commutative monoid either.
The fact that $\bbOP$ has an ultra-commutative group completion 
can be interpreted as saying that in this particular case 
`global group completion kills to obstruction to ultra-commutativity'.
\end{rk}

\index{subject}{Bott periodicity!global|(}

Bott periodicity is traditionally seen as a homotopy equivalence between
the space $\mZ\times B U$ and the loop space of the infinite unitary group $U$.
Since a loop space only sees the basepoint component, 
and the loop space of $B U$ is weakly equivalent to $U$, 
the 2-fold periodicity then takes the form of  a chain of weak equivalences:
\[ \Omega^2(\mZ\times B U) \ \simeq \ \Omega(\Omega(B U)) \ \simeq \ 
\Omega U \ \simeq \ \mZ\times B U\ .\]
We are going to prove a highly structured version of complex
Bott periodicity, in the form of a global equivalence of ultra-commutative monoids
between $\bBUP$ and $\Omega\bU$.\index{subject}{unitary group ultra-commutative monoid} 
Bott periodicity has been elucidated from many different angles,
and before we start, we put our approach into perspective.
Since Bott's original geometric argument \cite{bott-periodicity}
a large number of different proofs became available,
see for example \cite{hohnhold-stolz-teichner} for an overview.
In essence, our proof of global Bott periodicity is an adaptation 
of Harris' proof \cite{harris-Bott} of complex Bott periodicity.
The reviewer for Math Reviews praises Harris proof
as `a beautiful well-motivated proof of the complex Bott periodicity theorem 
using only two essential properties of the complex numbers'.
Suslin \cite{suslin-beilinson} calls this the `trivial proof'
of Bott periodicity and extends it to a `Real' (i.e., $C_2$-equivariant)
context. I also think that for readers with a homotopy theory background, Harris'
proof may be particularly accessible and appealing.

Harris' argument uses two main ingredients.
On the one hand, the group completion theorem is used
to identify the loop space of the bar construction of $\amalg_{n\geq 0}\, G r_n$
(under the monoid structure induced by orthogonal direct sum)
with $\mZ\times B U$. On the other hand, Harris exhibits an explicit homeomorphism
between the bar construction of $\amalg_{n\geq 0}\, Gr _n$
and the infinite unitary group, essentially the inverse to the eigenspace decomposition
of a unitary matrix.
Together these two ingredients provide a chain of weak equivalences
\[ \mZ\times B U \ \simeq \ \Omega( B( \amalg_{n\geq 0}\, G r_n)) \ \simeq \
\Omega U \ . \]
Our global proof is analogous:
Theorem \ref{thm:BOP is group completion} above
shows that $\bBUP$ is a global group completion of $\bGr^\mC$,
essentially by applying the group completion theorem to
all fixed point spaces. This part of the argument works just
as well for the real and symplectic versions of $\bGr^\mC$ and $\bBUP$. 
Theorem \ref{thm:group completion Gr to Omega U} below
shows that $\Omega \bU$ is also a  global group completion of $\bGr^\mC$,
by globally identifying the bar construction of $\bGr^\mC$ 
(with respect to the box product of orthogonal spaces) with $\bU$,
using Harris' homeomorphism between 
the realization $|\bGr^\mC_{\td{\bullet}}(W)|$ and $U(W_\mC)$.
Two global group completions of the same ultra-commutative monoid
are necessarily globally equivalent, which yields the global
version of complex Bott periodicity of 
Theorem \ref{thm:global Bott periodicity}.

\begin{construction}[Global Bott periodicity]
After this outline, we now provide the necessary details.
The ultra-commutative monoid $\bU$ of unitary groups
was defined in Example \ref{eg:unitary group monoid space}.\index{subject}{unitary group ultra-commutative monoid}
The orthogonal space $\Omega \bU$ inherits an ultra-commutative multiplication by
pointwise multiplication of loops,
where $\Omega$ means objectwise continuous based maps from $S^1$.
We define a morphism of ultra-commutative monoids
\begin{equation}  \label{eq:define_beta}
 \beta \ : \ \bGr^\mC \ \to \ \Omega \bU    
\end{equation}
at an inner product space $V$ by
\[ \beta(V)(L)(x)\ = \ c(x)\cdot p_L\ +\ p_{L^\perp}\ .\]
Here $L$ is a complex subspace of $V_\mC$, $x\in S^1$,
\[ c \ : \ S^1 \ \to \ U(1) \ , \quad x \ \longmapsto \ \frac{x+i}{x-i} \ ,\]
is the Cayley transform,\index{subject}{Cayley transform}
and $p_L$ and $p_{L^\perp}$ denote the orthogonal projections to $L$ 
respectively its orthogonal complement.
In other words, $L$ and $L^\perp$ are the eigenspaces of $\beta(V)(L)(x)$,
for the eigenvalues $c(x)$ respectively~1.
Then 
\[ \beta(V)(0)(x)\ = \  p_{V_\mC} \ = \ \Id_{V_\mC} \ ;\]
so $\beta(V)(0)$ is the constant loop at the identity, which is the
unit element of $\Omega \bU(V)$.
Now we consider subspaces $L\in\bGr^\mC(V)$ and $L'\in\bGr^\mC(W)$. Then
\begin{align}\label{eq:beta additivity}
  \beta(V\oplus W)(L\oplus L')(x) \ &= \ 
 (c(x)\cdot p_{L\oplus L'})\ +\ p_{L^\perp\oplus(L')^\perp}\\
&=\ 
( (c(x)\cdot p_L)\ +\ p_{L^\perp}) \oplus
( (c(x)\cdot p_{L'})\ +\ p_{(L')^\perp}) \nonumber\\
&= \ \beta(V)(L)(x)\oplus \beta(W)(L')(x)  \ .\nonumber
\end{align}
In other words, the square
\[ \xymatrix@C=20mm{
\bGr^\mC(V)\times \bGr^\mC(W)\ar[r]^-{\beta(V)\times\beta(W)} \ar[d]_{\oplus} & 
 \Omega \bU(V) \times \Omega\bU(W) \ar[d]^{\mu^{\Omega \bU}_{V,W}} \\
\bGr^\mC(V\oplus W)\ar[r]_-{\beta(V\oplus W)} & \Omega \bU(V\oplus W) } \]
commutes, i.e., $\beta$ is compatible with the multiplications
on both side. Since $\beta$ respects multiplication and unit,
it also respects the structure maps. The upshot is that
$\beta$ is a morphism of ultra-commutative monoids.

The category of ultra-commutative monoids is tensored and cotensored
over based spaces, so the functor of taking objectwise loops
is right adjoint to the functor $-\rhd S^1$ defined in \eqref{eq:R rhd A}.
For an ultra-commutative monoid $R$, the based tensor $R\rhd S^1$
is isomorphic to the bar construction $B(R)$, compare \eqref{eq:umon_bar_and_suspension}.

\begin{theorem}\label{thm:group completion Gr to Omega U}
The adjoint $\beta^\flat:B(\bGr^\mC)=\bGr^\mC\rhd S^1\to\bU$ 
of the morphism $\beta:\bGr^\mC\to\Omega\bU$ is a global equivalence
of ultra-commutative monoids.
The morphism $\beta:\bGr^\mC\to\Omega\bU$ is a global group completion of 
ultra-commu\-tative monoids.
\end{theorem}
\begin{proof}
We factor $\beta^\flat$ as a composite
of two morphisms of ultra-commutative monoids
\[ B(\bGr^\mC)\ = \ |B_\bullet(\bGr^\mC)|\ \xra[\simeq]{\ \zeta\ }\ 
 | \bGr^\mC_{\td{\bullet}} | \ \xra[\iso]{\ \epsilon\ }\ \bU \ ;\]
then we show that the morphism $\zeta$ is a global equivalence and the
morphism $\epsilon$ is an isomorphism. Together this shows the first claim.

The middle object is the realization 
of a simplicial ultra-commutative monoid $\bGr^\mC_{\td{\bullet}}$, 
and the first morphism is the realization of a simplicial morphism.
The object of $n$-simplices $\bGr^\mC_{\td{n}}$ is the ultra-commutative monoid
of $n$-tuples of pairwise orthogonal complex subspaces, i.e.,
\[ \bGr^\mC_{\td{n}}(V)\ = \ \{ (L_1,\dots,L_n)\in ( Gr^\mC(V_\mC))^n\ : \ 
\text{$L_i$ is orthogonal to $L_j$ for $i\ne j$}\} \ .\]
For varying $n$, the ultra-commutative monoids $\bGr^\mC_{\td{n}}$
assemble into a simplicial ultra-commutative monoid: the face morphisms
\[ d_i \ : \ \bGr^\mC_{\td{n}}\ \to \ \bGr^\mC_{\td{n-1}} \]
are given by
\[ d_i(L_1,\dots,L_n) \ = \left\lbrace \begin{array}{ll}
(L_2,\dots,L_n) & \mbox{for $i=0$,} \\
(L_1,\dots,L_{i-1},L_i\oplus L_{i+1},L_{i+2},\dots,L_n)  & \mbox{for $0<i< n$,} \\
(L_1,\dots,L_{n-1}) & \mbox{for $i=n$.}
\end{array} \right.  \]
For $n\geq 1$ and $0\leq i \leq n-1$ the degeneracy morphisms are given by
\[ s_i(L_1,\dots, L_{n-1}) \ = \ (L_1,\dots,L_i,0,L_{i+1},\dots,L_{n-1}) \ . \]
The direct sum maps
\begin{align*}
 \bGr^\mC(V_1)\times\dots\times \bGr^\mC(V_n)\ &\to \ \bGr^\mC_{\td{n}}(V_1\oplus\dots\oplus V_n)   \\
(L_1,\dots,L_n)\qquad &\longmapsto \quad (i_1(L_1),\dots,i_n(L_n))
\end{align*}
form a multi-morphism, where $i_k:(V_k)_\mC\to (V_1\oplus\dots\oplus V_n)_\mC$
is the embedding as the $k$-th summand. 
The universal property of the box product
turns this multi-morphism into a morphism of orthogonal spaces
\[
\zeta_n \ : \  B_n(\bGr^\mC)\ = \  (\bGr^\mC)^{\boxtimes n} \ \to \ \bGr^\mC_{\td{n}} \ .  
\]
The morphisms $\zeta_n$ are compatible with the simplicial face and
degeneracy maps, since these are given by orthogonal direct sum respectively
insertion of~0 on both sides. So for varying $n$, they form a morphism of
simplicial ultra-commutative monoids
\[ \zeta_\bullet \ : \ B_\bullet(\bGr^\mC) \ \to \ \bGr^\mC_{\td{\bullet}}\ .\]

We claim that $\zeta_n$ is a global equivalence for every $n\geq 0$. 
Since $\bGr^\mC=\coprod_{j\geq 0}\bGr^{\mC,[j]}$ and 
the box product distributes over disjoint unions,
$(\bGr^\mC)^{\boxtimes n}$ is the disjoint union of the orthogonal spaces
\[ \bGr^{\mC,[j_1]}\boxtimes\dots\boxtimes \bGr^{\mC,[j_n]}  \]
indexed over all tuples $(j_1,\dots,j_n)\in\mN^n$.
The orthogonal space $\bGr^\mC_{\td{n}}$ has an analogous decomposition, where
$\bGr^{\mC,[j_1,\dots,j_n]}_{\td{n}}$ consists of those tuples
$(L_1,\dots,L_n)$ with $\dim(L_i)=j_i$. The morphism $\zeta_n$ respects
the decomposition, i.e., it matches the two summands indexed by the same tuple
$(j_1,\dots,j_n)$. A disjoint union of global equivalences is a global equivalence
(Proposition \ref{prop:global equiv basics}~(v)),
so we are reduced to showing that each of the  morphisms
\[ \zeta_{j_1,\dots,j_n} \ : \ \bGr^{\mC,[j_1]}\boxtimes\dots\boxtimes \bGr^{\mC,[j_n]}\ \to \
\bGr^{\mC,[j_1,\dots,j_n]}_{\td{n}} \]
is a global equivalence.
This is in fact a restatement of an earlier result about box products
of orthogonal spaces `represented' by unitary representations.
In Construction \ref{con:complex-free spc} we defined
an orthogonal space $\bL^\mC_{G,W}$ 
from a unitary representation $W$ of a compact Lie group $G$.
The value at a euclidean inner product space $V$ is
\[ \bL^\mC_{G,W}(V) \ = \ \ \bL^\mC(W,V_\mC) / G\ .\]
Here $\bL^\mC$ is the space of $\mC$-linear maps that preserve
the hermitian inner products.
In the special case of the tautological $U(n)$-representation on $\mC^n$,
the homeomorphisms
\[  \bL^\mC(\mC^n,V_\mC) / U(n) \ \to \ G r^\mC_n(V_\mC)\ = \ \bGr^{\mC,[n]}(V)\ ,\quad
\varphi\cdot U(n)\ \longmapsto \ \varphi(\mC^n) \]
form an isomorphism of orthogonal spaces $\bL^\mC_{U(n),\mC^n}\iso\bGr^{\mC,[n]}$.
Similarly, passage to images provides 
an isomorphism of orthogonal spaces
\[ \bL^\mC_{U(j_1)\times\dots\times U(j_n),\mC^{j_1}\oplus\dots\oplus\mC^{j_n}}\ \iso\ 
\bGr^{\mC,[j_1,\dots,j_n]}_{\td{n}}\ .  \]
Under these identifications, the morphism $\zeta_{j_1,\dots,j_n}$ becomes the morphism
\[ \zeta_{U(j_1),\dots,U(j_n);\mC^{j_1},\dots,\mC^{j_n}} \ : \ 
\bL^\mC_{U(j_1),\mC^{j_1}}\boxtimes\dots\boxtimes \bL^\mC_{U(j_n),\mC^{j_n}}\ \to \
\bL^\mC_{U(j_1)\times\dots\times U(j_n),\mC^{j_1}\oplus\dots\oplus\mC^{j_n}}\ , \]
the iterate of the morphism
discussed in Proposition \ref{prop:box complex almost-representables}.
Proposition \ref{prop:box complex almost-representables}
thus shows that the morphism $\zeta_{j_1,\dots,j_n}$ is a global equivalence.
This completes the proof that 
the morphism $\zeta_n:(\bGr^\mC)^{\boxtimes n}\to\bGr^\mC_{\td{n}}$ is a global equivalence. 

Now we observe that the underlying simplicial orthogonal spaces of source and target
of $\zeta_\bullet$ are Reedy flat in the sense of Definition \ref{def:Reedy flat}, 
i.e., all latching morphisms (in the simplicial direction) are
flat cofibrations of orthogonal spaces.
Indeed, the unit of the ultra-commutative monoid $\bGr^\mC$
is the inclusion of the summand $\bGr^{\mC,[0]}$
into the disjoint union of all $\bGr^{\mC,[n]}$.
Since the orthogonal space $\bGr^{\mC,[n]}$ is isomorphic to $\bL^\mC_{U(n),\mC^n}$,
it is flat by Proposition \ref{prop:complex global classifying}~(ii).
So $\bGr^\mC$ has a flat unit, and its bar construction is Reedy flat by 
Proposition \ref{prop:bar global invariance}~(i).
Since $\bGr^{\mC}_{\td{\bullet}}$ is not the bar construction of any
orthogonal monoid space, we must show Reedy flatness directly.
Each simplicial degeneracy morphism of $\bGr^{\mC}_{\td{\bullet}}$ 
inserts the zero vector space in one slot; so the degeneracy morphisms
are embeddings of summands in a disjoint union.
The latching morphism
\[ L_m^\bDelta(\bGr^{\mC}_{\td{\bullet}}) \ \to \ \bGr^{\mC}_{\td{m}} \]
is then also the inclusions of certain summands, namely those 
$\bGr^{\mC,[j_1,\dots,j_n]}_{\td{n}}$ for which $j_i=0$ for at least one $i$.
Since the summand $\bGr^{\mC,[j_1,\dots,j_n]}_{\td{n}}$ is isomorphic to
$\bL^\mC_{U(j_1)\times\dots\times U(j_n),\mC^{j_1}\oplus\dots\oplus\mC^{j_n}}$,
it is flat by Proposition \ref{prop:complex global classifying}~(ii).
This verifies the Reedy flatness condition for $\bGr^{\mC}_{\td{\bullet}}$. 
Since source and target of the morphism $\zeta_\bullet$
are Reedy flat as simplicial orthogonal spaces, and $\zeta_\bullet$
is a global equivalence in every simplicial dimension,
the induced morphism of realizations
\[ \zeta\ = \  |\zeta_\bullet|\ : \  B(\bGr^\mC)\ = \ 
| B_\bullet(\bGr^\mC)|\ \to \ |\bGr^\mC_{\td{\bullet}}| \]
is a global equivalence by Proposition \ref{prop:realization invariance in spc}~(ii).

The isomorphism of ultra-commutative monoids
$\epsilon:|\bGr^\mC_{\td{\bullet}}|\iso\bU$
is taken from Harris \cite[Sec.\,2, Thm.]{harris-Bott},
and we recall it in some detail for the convenience of the reader.
We let $V$ be an inner product space and consider the continuous map
\begin{align*}
 \epsilon_n\ : \ \bGr^\mC_{\td{n}}(V)\times \Delta^n\quad &\to \quad\qquad U(V_\mC) \\  
(L_1,\dots,L_n;\, t_1,\dots,t_n)\ &\longmapsto \ 
{\prod}_{j=1}^n \, \exp(2\pi i t_j\cdot p_{L_j})\ ,
\end{align*}
where $p_L:V_\mC\to V_\mC$ is the orthogonal projection onto $L$.
In other words, $\epsilon_n(L_1,\dots,L_n;\, t_1,\dots,t_n)$
is the unitary automorphism of $V_\mC$ that has $L_j$ as
eigen\-space with eigenvalue $\exp(2\pi i t_j)$,
and is the identity on the orthogonal complement of all $L_j$'s.
We have the relations
\[
 \epsilon_n(L_1,\dots,L_n;\, 0, t_1,\dots, t_{n-1})
 \ = \ \epsilon_{n-1}(L_2,\dots, L_n;\, t_1,\dots,t_{n-1})
\]
and
\[
\epsilon_n(L_1,\dots,L_n;\, t_1,\dots, t_{n-1},1)\ = \ 
 \epsilon_{n-1}(L_1,\dots, L_{n-1};\, t_1,\dots,t_{n-1})
 \]
because $\exp(0)=\exp(2\pi i \cdot p_L)=\Id_{V_\mC}$; moreover,
\[ 
\epsilon_n(L_1,\dots,L_n;\, t_1,\dots,t_i,t_i\dots,t_{n-1})
\ = \  \epsilon_{n-1}(L_1,\dots,L_i\oplus L_{i+1},\dots, L_n;\, t_1,\dots,t_{n-1})
\]
for all $0<i<n$, because $\exp(2\pi i t\cdot p_{L_i})\cdot
\exp(2\pi i t\cdot p_{L_{i+1}})=\exp(2\pi i t\cdot p_{L_i\oplus L_{i+1}})$;
and finally
\begin{align*}
   \epsilon_{n+1}(L_1,\dots,L_i,0,L_{i+1},\dots, &L_n;\, t_1,\dots,t_{n+1})\\
 &= \ \epsilon_n(L_1,\dots,L_n;\, t_1,\dots,\widehat{t_{i+1}},\dots,t_{n+1})
\end{align*}
for all $0\leq i\leq n$.
So the maps $\epsilon_n$ are compatible with the equivalence relation
defining geometric realization, and they induce a continuous map
\[ \epsilon(V)\ : \ |\bGr^\mC_{\td{n}}(V)|\ \to \ U(V_\mC) \ .\]
The map $\epsilon(V)$ is bijective because every unitary automorphism is diagonalizable
with pairwise orthogonal eigenspaces and eigenvalues in $U(1)$.
As a continuous bijection from a compact space
to a Hausdorff space, $\epsilon(V)$ is a homeomorphism.
The homeomorphisms $\epsilon(V)$ are compatible
with linear isometric embeddings in $V$ and the ultra-commutative multiplications
on both sides, i.e., they define an isomorphism of ultra-commutative monoids
$\epsilon:|\bGr^\mC_{\td{\bullet}}|\iso \bU$.

Now we can conclude the proof.
Unraveling all definitions shows that the composite
\[ \bGr^\mC\sm S^1\ \to \ B(\bGr^\mC) \ \xra{\ \zeta\ }\ |\bGr^\mC_{\td{\bullet}}|
 \ \xra{\ \epsilon\ }\ \bU \]
is given at an inner product space $V$ by the map
\begin{align*}
 G r^\mC(V_\mC)\sm &S^1\, \to \quad \bU(V_\mC) \\
L\sm x \ &\longmapsto \ \exp(2\pi i\cdot \log(c(x))\cdot p_L)\ =\
c(x)\cdot p_L \ + \ p_{L^\perp}\ = \ \beta(V)(L)(x)\ .  
\end{align*}
This means that the original morphism $\beta$ is the composite
\[ \bGr^\mC\ \xra{\ \eta_{\bGr^\mC}\ }\ \Omega B(\bGr^\mC)\ 
\xra{\Omega(\zeta)}\  \Omega|\bGr^\mC_{\bullet}(V)|\
\xra{\Omega(\epsilon)}\  \Omega\bU\ ,  \]
where $\eta_{\bGr^\mC}$ was defined in \eqref{eq:define_bar_eta}.
Thus the adjoint $\beta^\flat$ factors as
$\epsilon\circ\zeta$, a global equivalence followed by an isomorphism.
So $\beta^\flat$ is  a global equivalence of ultra-commutative monoids.

We showed above that $\bGr^\mC$ has a flat unit, so the 
adjunction unit $\eta_{\bGr^\mC}$ 
is a global group completion by Corollary \ref{cor:umon Omega Sigma is group completion}.
Since $\beta^\flat$ is a global equivalence, so is $\Omega(\beta^\flat)$.
Hence the composite $\beta:\bGr^\mC\to \Omega\bU$ is a global group completion.
\end{proof}

Theorem \ref{thm:BOP is group completion}
and Theorem \ref{thm:group completion Gr to Omega U}
show that the morphisms 
\[ i\ :\ \bGr^\mC\ \to\ \bBUP\text{\qquad respectively\qquad} 
\beta\ :\ \bGr^\mC\ \to\ \Omega\bU\]
are both global group completions.
The universal property of group completions 
already implies that $\bBUP$ is isomorphic to $\Omega\bU$
in the homotopy category of ultra-commutative monoids;
the two can thus be linked by a chain of global equivalences of ultra-commutative monoids.
In a sense we could stop here, and call this `complex global Bott periodicity'.
However, we elaborate a bit more and exhibit an explicit chain
of two global equivalences between $\bBUP$ and $\Omega\bU$,
see Theorem \ref{thm:global Bott periodicity} below.

We define a morphism of ultra-commutative 
monoids\index{subject}{additive Grassmannian}\index{symbol}{$\beta$@$\bar\beta$ - {Bott morphism from $\bBUP$ to $\Omega \bU$}}
\[ \bar\beta \ : \ \bBUP \ \to \ \Omega( \sh_\tensor \bU) \ . \]
Here $\sh_\tensor=\sh^{\mR^2}_\tensor$ is the multiplicative shift by $\mR^2$ 
defined in Example \ref{eg:Additive and multiplicative shift}.
The orthogonal space $\bU$ has a commutative multiplication by
direct sum of unitary automorphisms; 
thus $\Omega(\sh_\tensor \bU)$ inherits a commutative multiplication by
pointwise multiplication of loops.
The target of $\bar\beta$ is globally equivalent, as an ultra-commutative monoid,
to $\Omega\bU$, the objectwise loops of the unitary group monoid.
The definition of the map
\[ \bar\beta(V)\ :\ \bBUP(V)\ \to\  \Omega (\sh_\tensor \bU) (V) \ = \ 
\map(S^1, \bU( V_\mC^2))\ , \]
for an inner product space $V$, is similar to, but slightly more elaborate than
the definition of $\beta(V)$ in \eqref{eq:define_beta} above.
An element of $\bBUP(V)$ is a complex subspace $L$ of $V_\mC^2$;
as before we denote by $p_L$ and $p_{L^\perp}$ the orthogonal projections to $L$ 
respectively to its orthogonal complement.
We define the loop
\[ \bar\beta(V)(L)\ : \ S^1 \ \to \ U(V_\mC^2)  \]
by
\[ \bar\beta(V)(L)(x)\ = \ 
( (c(x)\cdot p_L)\ +\ p_{L^\perp}) \circ
( (c(-x)\cdot p_{V_\mC\oplus 0}) + p_{0\oplus V_\mC})\ .\]
As before $c:S^1 \to U(1)$ is the Cayley transform.
The map $\bar\beta(V)$ is continuous in $L$.

For every inner product space $V$ we have
\[ \bar\beta(V)(V_\mC\oplus 0)(x)\ = \  
( (c(x)\cdot p_{V_\mC\oplus 0})\ +\ p_{0\oplus V_\mC}) \circ
( (c(-x)\cdot p_{V_\mC\oplus 0}) + p_{0\oplus V_\mC})\ =\ \Id_{V_\mC} \ ;\]
so $\bar\beta(V)(V_\mC\oplus 0)$ is the constant loop at the identity, which is the
unit element of $\Omega (\sh_\tensor \bU)(V)$.
Now we consider subspaces $L\in\bBUP(V)$ and $L'\in\bBUP(W)$. 
We recall that $\kappa^{V,W}:V_\mC^2\oplus W_\mC^2\iso (V\oplus W)_\mC^2$
is the preferred natural isometry, which enters into the definition of
the multiplication of $\bBUP$.
The argument for the additivity relation
\begin{align*}
  \bar\beta(V\oplus W)(L\oplus L')(x) \
= \ \kappa^{V,W}_* \left( \bar\beta(V)(L)(x)\oplus \bar\beta(W)(L')(x)\right)
\end{align*}
is straightforward and similar to (but somewhat longer than)
the argument for the map $\beta$ in \eqref{eq:beta additivity}, and we omit it. 
Hence $\bar\beta$ is compatible with the multiplications
on both sides. Since $\bar\beta$ respects multiplication and unit,
it also respects the structure maps. The upshot is that
$\bar\beta$ is a morphism of ultra-commutative monoids.
We also have
\begin{align*}
   \det( \bar\beta(V)(L)(x) ) \ &= \ \det( (c(x)\cdot p_L)+p_{L^\perp})
\cdot \det( (c(-x)\cdot p_{V_\mC\oplus 0})+p_{0\oplus V_\mC})\\
  &= \ c(x)^{\dim(L)-\dim(V)} \ ,
\end{align*}
exploiting that $c(-x)=\widebar{c(x)}=c(x)^{-1}$.
So the map $\bar\beta(V):\bBUP(V) \to \Omega U(V_\mC^2)$
sends the subspace $\bBU(V)=\bBUP^{[0]}(V)$ to $\Omega (S U(V_\mC^2))$.
Hence the morphism $\bar\beta$ restricts to a morphism of 
ultra-commutative monoids
\[ \bar\beta^{[0]} \ : \ \bBU\ \to \ \Omega ( \sh_\tensor \bSU) \ .\]
\end{construction}

Now we can properly state our global version of complex Bott periodicity.
The embeddings $j:V_\mC\to V_\mC^2$ as the first summand
induce a morphism of ultra-commutative monoids $\bU\circ j:\bU\to \sh_\tensor \bU$.

\begin{theorem}[Global Bott periodicity] \label{thm:global Bott periodicity}
The morphisms of ultra-commutative monoids 
\[ \bBUP\ \xra{\ \bar\beta\ }\ \Omega (\sh_\tensor \bU)\ 
\xla{\Omega(\bU\circ j)}\ \Omega\bU  \]
are global equivalences.
The morphism  $\bar\beta^{[0]}:\bBU\to\Omega (\sh_\tensor \bSU)$ 
is a global equivalence.
\end{theorem}
\begin{proof}
The morphism $\bU\circ j$ is a global equivalence 
by Theorem \ref{thm:general shift of osp},
hence so is $\Omega(\bU\circ j)$.
The following diagram of homomorphisms of ultra-commutative monoids
commutes by direct inspection:
\[  \xymatrix@C=15mm{
\bGr^\mC\ar[r]^-{\beta}\ar[d]_i  & \Omega \bU\ar[d]^{\Omega(\bU\circ j)} \\
\bBUP\ar[r]_-{\bar\beta} & \Omega(\sh_\tensor \bU) }     
\]
The morphism $\beta:\bGr^\mC\to\Omega\bU$ is a global group completion
by Theorem \ref{thm:group completion Gr to Omega U}.
So the composite $\Omega(\bU\circ j)\circ\beta:\bGr^\mC\to\Omega(\sh_\tensor\bU)$ 
with a global equivalence is also a global group completion.
The morphism $i:\bGr^\mC\to\bBUP$ is a global group completion by
Theorem \ref{thm:BOP is group completion}.
The universal property of group completions then shows that 
$\bar\beta:\bBUP\to\Omega(\sh_\tensor \bU)$ becomes an
isomorphism in the homotopy category of ultra-commutative monoids.
Since the global equivalences are part of a model structure,
this implies that $\bar\beta:\bBUP\to\Omega(\sh_\tensor \bU)$ 
is a global equivalence.
Since the morphism $\bar\beta^{[0]}:\bBU\to\Omega(\sh_\tensor \bSU)$ is
a retract of the global equivalence $\bar\beta$, it is a global equivalence itself.
\end{proof}

\begin{cor}\label{cor:group completion Gr to Omega U}
For every compact Lie group $G$ and every finite $G$-CW-complex $A$, the map
\[ [A,\beta]^G\ : \ [A,\bGr^\mC]^G\ \to \ [A,\Omega\bU]^G \]
is a group completion of abelian monoids.
\end{cor}
\begin{proof}
We contemplate the following commutative square of abelian monoid homomorphisms:
\[  \xymatrix@C=15mm{
[A,\bGr^\mC]^G \ar[d]_{[A,i]^G} \ar[r]^-{[A,\beta]^G} &[A,\Omega \bU]^G 
\ar[d]^{[A,\Omega(\bU\circ j)]^G}_\iso \\
[A,\bBUP]^G \ar[r]^-\iso_-{[A,\bar\beta]^G}& [A,\Omega(\sh_\tensor \bU) ]^G  }     
\]
Since $\bar\beta$ and $\Omega(\bU\circ j)$
are global equivalences by Theorem \ref{thm:global Bott periodicity},
the two homomorphisms $[A,\bar\beta]^G$ and $[A,\Omega(\bU\circ j)]^G$ 
are isomorphisms, by Proposition \ref{prop:[A,Y]^G of closed}~(ii).
The morphism $[A,i]^G$ is a group completion of abelian monoids by
Proposition \ref{prop:i is group completion}
(or rather its complex analog, which is proved analogously).
So $[A,\beta]^G:[A,\bGr^\mC]^G\to[A,\Omega\bU]^G$ 
is also a group completion of abelian monoids.
\end{proof}

\index{subject}{Bott periodicity!global|)}
\index{subject}{global group completion!of an ultra-commutative monoid|)}
\index{subject}{ultra-commutative monoid|)}

\chapter{Equivariant stable homotopy theory}
\label{ch-equivariant}

In this chapter we give a largely self-contained 
exposition of many basics about equivariant stable homotopy theory
for a fixed compact Lie group; our model is the category of orthogonal $G$-spectra.
In Section \ref{sec:equivariant homotopy groups}
we review orthogonal spectra and orthogonal $G$-spectra;
we define equivariant stable homotopy groups and prove their basic properties,
such as the suspension isomorphism and long exact sequences of
mapping cones and homotopy fibers, and the additivity
of equivariant homotopy groups on sums and products.
Section \ref{sec:Wirthmuller and transfer} discusses the Wirthm{\"u}ller isomorphism
that relates the equivariant homotopy groups 
of a spectrum over a subgroup to the equivariant homotopy groups 
of the induced spectrum; intimately related to the 
Wirthm{\"u}ller isomorphism are various transfers that we also recall.
In Section \ref{sec:geometric fixed points} we introduce and study
geometric fixed point homotopy groups.
We establish the isotropy separation sequence that facilitates
inductive arguments, and show that equivariant equivalences can
also be detected by geometric fixed points.
We use geometric fixed points to derive a functorial description of the
0-th equivariant stable homotopy group of a $G$-space $Y$ in terms
of the path components of the fixed point spaces $Y^H$.
Section \ref{sec:double coset} gives a self-contained proof
of the double coset formula for the composite of a transfer 
followed by a restriction to a closed subgroup.
We also discuss various examples and end with a discussion of Mackey functors for
finite groups.
After inverting the group order, the category of $G$-Mackey functors
splits as a product, indexed by conjugacy classes of subgroups,
of module categories over the Weyl groups, 
see Theorem \ref{thm:split Mackey functors}.
We show that rationally and for finite groups, geometric fixed point homotopy groups
can be obtained from equivariant homotopy groups by dividing out transfers
from proper subgroups.
Section \ref{sec:products} is devoted to multiplicative aspects of
equivariant stable homotopy theory. In our model,
all multiplicative features can be phrased in terms of the
smash product of orthogonal spectra (or orthogonal $G$-spectra),
another example of a Day type convolution product.
The smash product gives rise to pairings of equivariant homotopy groups;
when specialized to equivariant ring spectra, these pairings turn 
the equivariant stable homotopy into graded rings.

We do not discuss model category structures for orthogonal $G$-spectra;
the interested reader can find different ones in the memoir of
Mandell and May \cite{mandell-may},
in the thesis of Stolz \cite{stolz-thesis},
the article by Brun, Dundas and Stolz \cite{brun-dundas-stolz} 
and (for finite groups) in the paper of Hill, Hopkins and Ravenel \cite{HHR-Kervaire}.

\section{Equivariant orthogonal spectra}
\label{sec:equivariant homotopy groups}

In this section we begin to develop some of the basic 
features of equivariant stable homotopy theory for compact Lie groups
in the context of equivariant orthogonal spectra. 
After introducing orthogonal $G$-spectra
and equivariant stable homotopy groups, 
we discuss shifts by a representation and show that they are
equivariantly equivalent to smashing with the representation sphere
(Proposition \ref{prop:lambda upi_* isos}).
We establish the loop and suspension isomorphisms 
(Proposition \ref{prop:loop and suspension isomorphisms})
and the long exact homotopy
group sequences of homotopy fibers and mapping cones
(Proposition \ref{prop:LES for homotopy of cone and fibre}).
We prove that equivariant homotopy groups take wedges to sums and
preserve finite products (Corollary \ref{cor-wedges and finite products}).
We end by showing that the equivariant homotopy group functor $\pi_0^H$,
for a closed subgroup $H$ of a compact Lie group $G$, is represented
by the unreduced suspension spectrum of the homogeneous space $G/H$
(Proposition \ref{prop:G/H represents pi_0^H}).  

\medskip

We recall {\em orthogonal spectra}. 
These objects are used, at least implicitly, already in \cite{may-quinn-ray};
the term `orthogonal spectrum' was introduced by
Mandell, May, Shipley and the author in \cite{mmss}, where the
(non-equivariant) stable model structure for orthogonal spectra
was constructed.
Orthogonal spectra are stable versions of orthogonal spaces, and before recalling
the formal definition we try to motivate it -- already with a view
towards the global perspective.
An orthogonal space $Y$ assigns values 
to all finite-dimensional inner product spaces.
The global homotopy type is encoded in the $G$-spaces $Y(\Uc_G)$, 
where $\Uc_G$ is a complete $G$-universe, which we can informally think of 
as `the homotopy colimit of $Y(V)$ over all $G$-representations $V$'.
So besides the values $Y(V)$, an orthogonal space uses the information
about the $O(V)$-action (which is turned into a $G$-action when $G$ acts on $V$)
and the information about inclusions of inner product spaces
(in order to be able to stabilize to the colimit $\Uc_G$).
The information about the $O(V)$-actions and how to stabilize are conveniently
encoded together as a continuous functor 
from the category $\bL$ of linear isometric embeddings.

An orthogonal spectrum $X$ is a stable analog of this: it assigns
a based space $X(V)$ to every inner product space, and it keeps
track of an $O(V)$-action on $X(V)$ (to get $G$-homotopy types when $G$ acts on $V$)
and of a way to stabilize by suspensions (needed when exhausting a complete universe
by its finite-dimensional subrepresentations). 
When doing this in a coordinate free way, the stabilization data
assigns to a linear isometric embedding $\varphi:V\to W$ a 
continuous based map
\[ \varphi_\star \ : \ S^{W-\varphi(V)}\sm  X(V)\ \to \ X(W) \]
where $W-\varphi(V)$ is the orthogonal complement of the image of $\varphi$.
This structure map should `vary continuously with $\varphi$',
but this phrase has no literal meaning because the source of $\varphi_\star$
depends on $\varphi$. The way to make the continuous dependence
rigorous is to exploit that the complements $W-\varphi(V)$ vary in 
a locally trivial way, i.e., they are the fibers of a distinguished vector bundle,
the `orthogonal complement bundle', 
over the space of $\bL(V,W)$ of linear isometric embeddings.
All the structure maps $\varphi_\star$ together define a map on
the smash product of $X(V)$ with the Thom space of this complement bundle,
and the continuity of the dependence on $\varphi$ is formalized by
requiring continuity of that map.
All these Thom spaces together form the morphism spaces of a based
topological category, and the data of an orthogonal spectrum 
can conveniently be packaged as a continuous based
functor on this category.

\begin{construction}
We let $V$ and $W$ be inner product spaces. Over the space
$\bL(V,W)$ of linear isometric embeddings sits a certain `orthogonal complement'
vector bundle with total space
\[ \xi(V,W) \ = \ 
\{\, (w,\varphi) \in W\times\bL(V,W) \ | \ w\perp\varphi(V)\,\} \ .\]
The structure map $\xi(V,W)\to\bL(V,W)$\index{symbol}{$\xi(V,W)$ - {orthogonal complement vector bundle}} 
is the projection to the second factor.
The vector bundle structure of $\xi(V,W)$ is 
as a vector subbundle of the trivial vector bundle $W\times\bL(V,W)$,
and the fiber over $\varphi:V\to W$ is the orthogonal complement $W-\varphi(V)$ 
of the image of $\varphi$.

We let $\bO(V,W)$\index{symbol}{$\bO(V,W)$ - {Thom space of orthogonal complement bundle}} be the Thom space of the bundle $\xi(V,W)$,
i.e., the one-point compactification of the total space of $\xi(V,W)$.
Up to non-canonical homeomorphism, we can describe the space $\bO(V,W)$ 
differently as follows.
If the dimension of $W$ is smaller than the dimension of $V$,
then the space $\bL(V,W)$ is empty and $\bO(V,W)$ consists
of a single point at infinity. Otherwise we can choose a linear isometric embedding
$\varphi:V\to W$, and then the maps
\begin{align*}
 O(W)/O(W-\varphi(V) ) \qquad &\to \ \bL(V,W)\ , \quad
A\cdot O(W-\varphi(V) ) \ \longmapsto \ A\varphi \text{\quad and}\\
O(W)\ltimes_{O(W-\varphi(V))} S^{W-\varphi(V) }\ &\to \ \bO(V,W)\ , \quad
[A,w] \ \longmapsto \ (A w, A \varphi)
\end{align*}
are homeomorphisms.
Here, and in the following, we write\index{symbol}{$G\ltimes_H A$ - {induced based $G$-space}}
\[ G\ltimes_H A \ = \ (G_+)\sm_H A \ = \ (G_+\sm A) / \sim\]
for a closed subgroup $H$ of $G$ and a based $G$-space $A$;
the equivalence relation is $g h\sm a\sim g \sm h a$
for all $(g,h,a)\in G\times H\times A$.
Put yet another way: if $\dim V=n$ and $\dim W=n+m$,
then $\bL(V,W)$ is homeomorphic to the homogeneous space $O(n+m)/O(m)$
and $\bO(V,W)$ is homeomorphic to $O(n+m)\ltimes_{O(m)}S^m$.
The vector bundle $\xi(V,W)$ becomes trivial upon product with
the trivial bundle $V$, via the trivialization
\[ \xi(V,W)\times V \ \iso \ W\times \bL(V,W)\ , \quad
 ((w,\varphi), v) \ \longmapsto \ (w+\varphi(v),\varphi) \ . \]
When we pass to Thom spaces on both sides this becomes the
{\em untwisting homeomorphism:}\index{subject}{untwisting homeomorphism}
\begin{equation}\label{eq:untwisting_homeo}
   \bO(V,W)\sm S^V \ \iso \ S^W\sm \bL(V,W)_+\ .
\end{equation}

The Thom spaces $\bO(V,W)$ are the morphism spaces of a based topological category.
Given a third inner product space $U$, the bundle map 
\[ \xi(V,W) \times \xi(U,V) \ \to \ \xi(U,W) \ , \quad
((w,\varphi),\,(v,\psi)) \ \longmapsto \ (w+\varphi(v),\,\varphi\psi)\]
covers the composition map $\bL(V,W)\times\bL(U,V)\to \bL(U,W)$.
Passage to Thom spaces gives a based map
\[ \circ \ : \ \bO(V,W) \sm \bO(U,V) \ \to \ \bO(U,W) \]
which is clearly associative, and is the composition in the category $\bO$. 
The identity of $V$ is $(0,\Id_V)$ in $\bO(V,V)$.
\end{construction}

\begin{defn}\label{def:orthogonal spectrum}
An {\em orthogonal spectrum}\index{subject}{orthogonal spectrum} 
\index{subject}{spectrum!orthogonal|see{orthogonal spectrum}} 
is a based continuous functor from $\bO$ to the category $\bT_*$ of based spaces.
A {\em morphism} of orthogonal spectra\index{subject}{morphism!of orthogonal spectra} 
is a natural transformation of functors.
We denote the category of orthogonal spectra by $\spec$.\index{symbol}{  $\spec$ - {category of orthogonal spectra}}
\end{defn}

Given two inner product spaces $V$ and $W$ we define a continuous based map
\[ i_V \ : \ S^V \ \to \ \bO(W,V\oplus W)\text{\qquad by\qquad} 
v \ \longmapsto \ ((v,0), (0,-)) \ ,\]
where $(0,-):W\to V\oplus W$ is the embedding of the second summand.
We define the {\em structure map}\index{subject}{structure map!of an orthogonal spectrum} 
$\sigma_{V,W}:S^V\sm X(W)\to X(V\oplus W)$ of the orthogonal spectrum $X$
as the composite\index{symbol}{$\sigma_{V,W}$ - {structure map of an orthogonal spectrum}}
\begin{equation}\label{eq:define_structure} 
 S^V \sm   X(W) \ \xra{i_V\sm X(W)}\ 
\bO(W,V\oplus W)\sm X(W)\ \xra{\ X\ }\  X(V\oplus W) \ .
\end{equation}
Often it will be convenient to use the {\em opposite structure map}\index{subject}{structure map!of an orthogonal spectrum!opposite}\index{symbol}{$\sigma_{V,W}^{\op}$ - {opposite structure map of an orthogonal spectrum}}
\begin{equation}\label{eq:define opposite structure}  
 \sigma^{\op}_{V,W}\ : \    X(V)\sm S^W \ \to\  X(V\oplus W) 
\end{equation}
which we define as the following composite:
\[ X(V)\sm S^W \ \xra{\text{twist}}\
S^W \sm   X(V) \ \xra{\sigma_{W,V}} \  X(W\oplus V) \ \xra{X(\tau_{V,W})}
\  X(V\oplus W) \ \]

\begin{rk}[Coordinatized orthogonal spectra]
Every inner product space is isometrically isomorphic to $\mR^n$
with standard inner product, for some $n\geq 0$.
So the topological category $\bO$ has a small skeleton, 
and the functor category of orthogonal spectra has `small' morphism sets.
This also leads to the following
more explicit coordinatized description of orthogonal spectra
in a way that resembles a presentation by generators and relations. 

Up to isomorphism, an orthogonal spectrum $X$ is determined by the values $X_n=X(\mR^n)$ 
and the following additional data relating these spaces:
\begin{itemize}
\item a based continuous left $O(n)$-action on $X_n$ for each $n\geq 0$, 
\item based maps $\sigma_n:S^1\sm X_n\to X_{1+n}$ for $n\geq 0$.
\end{itemize}
This data is subject to the following condition:
for all $m, n\geq 0$, the iterated structure map
$S^m\sm X_n \to  X_{m+n}$ defined as the composition
\[ \xymatrix@C=15mm{
S^m\sm X_n \ar^-{S^{m-1}\sm \sigma_n}[r] &
S^{m-1}\sm X_{1+n} \ar^-{S^{m-2}\sm \sigma_{1+n}}[r] & 
\quad \cdots \quad \ar^-{\sigma_{m-1+n}}[r] & \ X_{m+n} 
} \]
is $(O(m)\times O(n))$-equivariant. Here the group $O(m)\times O(n)$ 
acts on the target by restriction, along orthogonal sum, of the $O(m+n)$-action.
Indeed, the map 
\[ O(n)_+\ \to\ \bO(\mR^n,\mR^n) \ , \quad A\ \longmapsto \ (0,A)\]
is a homeomorphism, so $O(n)$ `is' the endomorphism monoid of $\mR^n$
as an object of the category $\bO$; via this map, $O(n)$ acts on the
value at $\mR^n$ of any functor on $\bO$.
The map $\sigma_n=\sigma_{\mR,\mR^n}$ 
is just one of the structure maps \eqref{eq:define_structure}.
\end{rk}

\begin{defn}
Let $G$ be a compact Lie group.
An {\em orthogonal $G$-spectrum}\index{subject}{orthogonal $G$-spectrum} 
is a based continuous functor from $\bO$ to the category $G\bT_*$ of based $G$-spaces.
A {\em morphism} 
of orthogonal $G$-spectra\index{subject}{morphism!of orthogonal $G$-spectra} 
is a natural transformation of functors.
We write $G\spec$\index{symbol}{  $G\spec$ - {category of orthogonal $G$-spectra}}
for the category of orthogonal $G$-spectra and $G$-equivariant morphisms.
\end{defn}

A continuous functor to based $G$-spaces is the same data as a $G$-object
of continuous functors. So orthogonal $G$-spectra could equivalently
be defined as orthogonal spectra equipped with a continuous $G$-action.
An orthogonal $G$-spectrum $X$ can be evaluated on a $G$-representation $V$,
and then $X(V)$ is a $(G\times G)$-space by the `external' $G$-action on $X$ 
and the `internal' $G$-action from the $G$-action on $V$
and the $O(V)$-functoriality of $X$.
We consider $X(V)$ as a $G$-space via the diagonal $G$-action.
If $V$ and $W$ are $G$-representations, then the 
structure map \eqref{eq:define_structure} 
and the opposite structure map \eqref{eq:define opposite structure}  
are $G$-equivariant where 
the group $G$ also acts on the representation spheres.

\begin{rk}
Our definition of orthogonal $G$-spectra is not the same as the one used
by Mandell and May \cite{mandell-may} and Hill, Hopkins and Ravenel \cite{HHR-Kervaire},
who define orthogonal $G$-spectra as $G$-functors on a $G$-enriched extension
of the category $\bO$ that contains all $G$-representations as objects.
However, our category of orthogonal $G$-spectra is equivalent to theirs by
\cite[V Thm.\,1.5]{mandell-may}. The substance of this equivalence is
the fact that for every orthogonal $G$-spectrum in the sense of Mandell and May,
the values at arbitrary $G$-representations are in fact determined
by the values at trivial representations.
\end{rk}

Next we recall the equivariant stable homotopy groups $\pi_*^G(X)$
(indexed by the complete $G$-universe) of an orthogonal $G$-spectrum $X$.
We introduce a convenient piece of notation.
If $\varphi:V\to W$ is a linear isometric embedding
and $f:S^V\to X(V)$ a continuous based map, 
we define $\varphi_*f:S^W\to X(W)$ as the composite
\begin{align}\label{eq:define varphi_* f}
 S^W \iso \ S^{W-\varphi(V)}\sm S^V\ &\xra{S^{W-\varphi(V)}\sm f} \ 
S^{W-\varphi(V)}\sm X(V) \\ 
&\xra{\sigma_{W-\varphi(V),V}} \ X((W-\varphi(V))\oplus V) \ \xra{\ \iso\ } X(W) \nonumber
\end{align}
where two unnamed homeomorphisms use the linear isometry
\[ (W-\varphi(V))\oplus V \ \iso \ W  \ , \quad (w,v)\ \longmapsto \ w+ \varphi(v) \ .\]
For example, if $\varphi$ is bijective (i.e., an equivariant isometry),
then $\varphi_* f$ becomes the $\varphi$-conjugate of $f$, i.e., the composite
\[ S^W \ \xra{\ S^{\varphi^{-1}}\ } \ S^V \ \xra{\ f\ } \ X(V) \
\xra{\ X(\varphi)\ } \ X(W) \ . \]
The construction is continuous in both variables, i.e., the map
\[ \bL(V,W)\times \map_*(S^V,X(V)) \ \to \ \map_*(S^W,X(W)) \ , \quad
(\varphi,f)\ \longmapsto \ \varphi_* f \]
is continuous.

As before we let $s(\Uc_G)$ denote the poset, under inclusion,
of finite-dimensional $G$-subrepresentations 
of the chosen complete $G$-universe $\Uc_G$.
We obtain a functor from $s(\Uc_G)$ to sets by sending $V\in s(\Uc_G)$ to 
\[  [S^V, X(V)]^G \ ,\]
the set of $G$-equivariant homotopy classes of based $G$-maps
from $S^V$ to $X(V)$.
For $V\subseteq W$ in $s(\Uc_G)$ the inclusion $i:V\to W$ is sent to the map
\[ i_* \ : \ [S^V, X(V)]^G \ \to \ [S^W,X(W)]^G \ ,\quad [f]\ \longmapsto [i_*f] \ .\]
The {\em 0-th equivariant homotopy group} $\pi_0^G(X)$ is then defined 
as\index{subject}{homotopy group!equivariant}\index{subject}{equivariant homotopy group!of an orthogonal $G$-spectrum}
\[  \pi_0^G(X) \ = \ \colim_{V\in s(\Uc_G)}\, [S^V, X(V)]^G  \ ,   \]
the colimit of this functor over the poset $s(\Uc_G)$.

The sets $\pi_0^G(X)$ have a lot of extra structure;
we start with the abelian group structure.
We consider a finite-dimensional $G$-subrepresentation $V$ of the 
universe $\Uc_G$ with non-zero fixed points.
We choose a $G$-fixed unit vector $v_0\in V$, 
and we let $V^\perp$ denote the orthogonal complement of $v_0$ in $V$.
This induces a decomposition
\[ \mR \oplus V^\perp \ \iso \ V \ , \quad (t,v)\ \longmapsto \ t v_0 +v \]
that extends to a $G$-equivariant homeomorphism
$S^1\sm S^{V^\perp}\iso S^V$ on one-point compactifications.
From this we obtain a bijection
\begin{equation}\label{eq:biject to pi_1} 
[S^V,X(V)]^G\ \iso \ [S^1,\map_*^G(S^{V^\perp},X(V))]_*   
\ =\ \pi_1(\map_*^G(S^{V^\perp},X(V)))\ ,
\end{equation}
natural in the orthogonal $G$-spectrum $X$. 
We use the bijection \eqref{eq:biject to pi_1}
to transfer the group structure on the fundamental group
into a group structure on the set $[S^V,X(V)]^G$.

Now we suppose that the dimension of the fixed point space $V^G$ is at least~2.
Then the space of $G$-fixed unit vectors in $V$ is connected and
similar arguments as for the commutativity of higher homotopy groups show:
\begin{itemize}
\item the group structure on the set $[S^V,X(V)]^G$ defined by
the bijection \eqref{eq:biject to pi_1}
is commutative and independent of the choice of $G$-fixed unit vector;
\item
if $W$ is another finite-dimensional $G$-subrepresentation of
$\Uc_G$ containing $V$, then the map
\[ i_* \ : \ [S^V,X(V)]^G\ \to \ [S^W,X(W)]^G \]
is a group homomorphism.
\end{itemize}
The $G$-subrepresentations $V$ of $\Uc_G$ with $\dim(V^G)\geq 2$
are cofinal in the poset $s(\Uc_G)$, so the two properties above
show that the abelian group structures on $[S^V,X(V)]^G $
for  $\dim(V^G)\geq 2$ assemble into a well-defined and
natural abelian group structure on the colimit $\pi_0^G(X)$.\medskip

We generalize the definition of $\pi_0^G(X)$ to
integer graded equivariant homotopy groups of an orthogonal $G$-spectrum $X$.
If $k$ is a positive integer, then we set
\begin{align}\label{eq:define_pi_k}
 \pi_k^G (X) \ &= \ \colim_{V\in s(\Uc_G)}\, [S^{V\oplus\mR^k}, X(V)]^G  
\text{\qquad and}\\
\pi_{-k}^G (X) \ &= \ \colim_{V\in s(\Uc_G)}\, [S^{V}, X(V\oplus \mR^k)]^G  \ .\nonumber
\end{align}
The colimits are taken over the analogous stabilization maps as for $\pi_0^G$,
and they come with abelian group structures by the same reasoning as for $\pi_0^G(X)$.

\begin{defn}
A morphism $f:X\to Y$ of orthogonal $G$-spectra 
is a {\em $\upi_*$-isomorphism}\index{subject}{pi star isomorphism@$\upi_*$-isomorphism!of orthogonal $G$-spectra}
if the induced map $\pi_k^H(f):\pi_k^H(X) \to \pi_k^H(Y)$ is an isomorphism
for all closed subgroups $H$ of $G$ and all integers $k$.
\end{defn}

\begin{construction}\label{con:general representative pi_0^G}
While the definition of $\pi_k^G(X)$ involves a case distinction
in positive and negative dimensions $k$, every class in $\pi_k^G( X)$
can be represented by a $G$-map
\[ f \ : \ S^{V\oplus \mR^{n+k}}\ \to \ X(V\oplus \mR^n)  \]
for suitable $n\in\mN$ such that $n+k\geq 0$.
Moreover, $V$ can be any finite-dimensional $G$-representation,
not necessarily a subrepresentation of the chosen complete $G$-universe.
Since we will frequently use this way to represent elements of $\pi_k^G(X)$,
we make the construction explicit here.

We start with the case $k\geq 0$. We choose a $G$-equivariant linear isometry 
$j:V\oplus\mR^n\to \bar V$ onto a $G$-subrepresentation $\bar V$ of $\Uc_G$.
Then the composite
\[ S^{\bar V\oplus\mR^k} \ \xra[\iso]{(S^{j\oplus\mR^k})^{-1}}\ S^{V\oplus\mR^{n+k}} \ \xra{\ f \ } \ 
X(V\oplus\mR^n) \ \xra[\iso]{X(j)} \ X(\bar V) \]
represents a class $\td{f}\in\pi_k^G(X)$.
For $k\leq 0$, we choose a $G$-equivariant linear isometry
$j:V\oplus\mR^{n+k}\to \bar V$ onto a $G$-subrepresentation $\bar V$ of $\Uc_G$.
Then the composite
\[ S^{\bar V} \ \xra[\iso]{(S^j)^{-1}} \ S^{V\oplus\mR^{n+k}} \ \xra{\ f \ } \ 
X(V\oplus\mR^n) \ \xra[\iso]{X(j\oplus\mR^{-k})} \ X(\bar V\oplus\mR^{-k} ) \]
represents a class $\td{f}\in\pi_k^G(X)$. 

We also need a way to recognize that `stabilization along a linear isometric embedding'
does not change the class in $\pi_k^G(X)$.
For this we let $\varphi:V\to W$ be a $G$-equivariant linear isometric embedding
and $f:S^{V\oplus\mR^{n+k}}\to X(V\oplus\mR^n)$ a continuous based $G$-map as above. 
We define $\varphi_*f:S^{W\oplus\mR^{n+k}}\to X(W\oplus\mR^n)$ as the composite
\begin{align*}
   S^{W\oplus\mR^{n+k}}\ \iso \ S^{W-\varphi(V)}\sm &S^{V\oplus\mR^{n+k}}
\ \xra{S^{W-\varphi(V)}\sm f} \ 
S^{W-\varphi(V)}\sm X(V\oplus\mR^n) \\ 
&\xra{\sigma_{W-\varphi(V),V\oplus\mR^n}}\  X((W-\varphi(V))\oplus V\oplus \mR^n)  
\xra{\ \iso \ } X(W\oplus\mR^n)  \ ;
\end{align*}
the two unnamed homeomorphisms use the linear isometry
\[ (W-\varphi(V))\oplus V\oplus\mR^m\ \iso \ W\oplus\mR^m  \ , \quad (w,v,x)\ \longmapsto \ (w+ \varphi(v), x) \]
for $m=n+k$ respectively $m=n$.
In the special case $n=k=0$, this construction reduces 
to \eqref{eq:define varphi_* f}.
\end{construction}

The same reasoning as in the unstable situation 
in Proposition \ref{prop:universal colimit spaces}
shows the following stable analog:

\begin{prop}\label{prop:invariant description stable}
Let $G$ be a compact Lie group and $X$ an orthogonal $G$-spectrum.
Let $V$ be a $G$-representation and $f:S^{V\oplus\mR^{n+k}}\to X(V\oplus\mR^n)$ a 
based continuous $G$-map, where $n\in\mN$ and $k\in\mZ$ are such that $n+k\geq 0$.
\begin{enumerate}[\em (i)]
\item The class $\td{f}$ in $\pi_k^G(X)$ is independent of the choice 
of linear isometry onto a subrepresentation of $\Uc_G$.
\item For every $G$-equivariant linear isometric embedding 
$\varphi:V\to W$ the relation
\[ \td{\varphi_*f} \ = \ \td{f} \text{\qquad holds in\quad $\pi_k^G(X)$.}  \]
\end{enumerate}
\end{prop}

Now we let $K$ and $G$ be two compact Lie groups.
Every continuous based functor $F:G\bT_\ast\to K\bT_\ast$
from based $G$-spaces to based $K$-spaces gives rise to a functor
\[ F\circ- \ : \ G\spec \ \to \ K\spec \]
from orthogonal $G$-spectra to orthogonal $K$-spectra
by postcomposition: if $X$ is a $G$-orthogonal spectrum,
then the composite
\[ \bO \ \xra{\ X \ }\ G\bT_*\ \xra{\ F\ } \ K\bT_\ast\ . \]
is an orthogonal $K$-spectrum.
The next construction is an example of this.

\begin{construction}[Restriction maps]\index{subject}{restriction homomorphism!for equivariant stable homotopy groups}\label{con:restriction map G-spectra}
We let $\alpha:K\to G$ be a continuous homomorphism between compact Lie groups.
Given an orthogonal $G$-spectrum $X$, we apply restriction of scalars levelwise
and obtain an orthogonal $K$-spectrum $\alpha^* X$.
We define the {\em restriction homomorphism}
\[ \alpha^* \ : \ \pi^G_0(X)\ \to \ \pi^K_0(\alpha^* X) \text{\qquad by\qquad}
 \alpha^*[f]\ = \ \td{\alpha^*(f)}\ . \]
In other words, the class represented by a based $G$-map 
$f:S^V\to X(V)$ is sent to the class represented by the $K$-map
\[ \alpha^*(f)\ : \ S^{\alpha^*(V)} \ = \
\alpha^*(S^V) \ \to \alpha^*(X(V))\ = \ (\alpha^* X)(\alpha^* V) \ , \]
appealing to Construction \ref{con:general representative pi_0^G}.
The restriction  maps $\alpha^*$ are clearly transitive 
(contravariantly functorial) for composition of group homomorphisms.
\end{construction}

For $g\in G$ the conjugation homomorphism is defined as
\[ c_g \ : \ G \ \to \ G \ , \quad c_g(h) \ = \ g^{-1} h g\ .\]
For every $G$-space $A$, left multiplication by $g$ is then a
$G$-equivariant homeomorphism $l^A_g:c_g^*(A)\to A$.
For an orthogonal $G$-spectrum $X$ the maps 
$l_g^{X(V)}:(c_g^* X)(V)=c_g^*(X(V))\to X(V)$
assemble into an isomorphism of orthogonal $G$-spectra $l_g^X:c_g^*X\to X$,
as $V$ runs over all inner product spaces (with trivial $G$-action).

\begin{prop}\label{prop:inner is identity}
Let $G$ be a compact Lie group, $X$ an orthogonal $G$-spectrum and $g\in G$.
Then the two isomorphisms
\[ c_g^* \ : \ \pi_0^G(X) \ \to \ \pi_0^G(c_g^* X) \text{\qquad and\qquad}
(l_g^X)_*\ : \  \pi_0^G(c_g^* X) \ \to \ \pi_0^G( X) 
\]
are inverse to each other.  
\end{prop}
\begin{proof}
We let $V$ be a $G$-representation,
and we recall that the $G$-action on $X(V)$ is diagonally,
from the external $G$-action on $X$ and the internal $G$-action on $V$.   
Hence the map $l_g^{X(V)}:c_g^*(X(V))\to X(V)$
is the composite of the map $l_g^X(c_g^* V):(c_g^*X)(c_g^* V)\to X(c_g^* V)$
and the map $X(l_g^V):X(c_g^* V)\to X(V)$.
Now we let $f:S^V\to X(V)$ be a $G$-map representing a class in $\pi_0^G(X)$.
The following diagram of $G$-spaces and $G$-maps commutes
because $f$ is $G$-equivariant:
\[ \xymatrix@C=20mm{ 
S^{c_g^*V} \ar[r]^-{c_g^* f}\ar[d]_{l_g^V} & 
c_g^*(X(V)) =(c_g^* X)(c_g^* V)\ar[r]^-{l_g^X(c_g^*V)} \ar@/^1pc/[dr]_(.6){l_g^{X(V)}}&
X(c_g^* V)\ar[d]^-{X(l_g^V)} \\
S^V \ar[rr]_-f && X(V)
}\]
The upper horizontal composite represents the class $(l_g^X)_*(c_g^*[f])$.
Since it differs from $f$ by conjugation with an equivariant isometry,
the upper composite represents the same class as $f$,
by Proposition \ref{prop:invariant description stable}~(ii).
Thus we conclude that $(l_g)_*(c_g^*[f])=[f]$.
\end{proof}

\begin{rk}[Weyl group action on equivariant homotopy groups]\label{rk:Weyl on pi_0^G}
We consider a closed subgroup $H$ of a compact Lie group $G$ 
and an orthogonal $G$-spectrum $X$.
Every $g\in G$ 
gives rise to a conjugation homomorphism
$c_g:H\to H^g$ by $c_g(h)=g^{-1}h g$,
where $H^g=\{g^{-1} h g\ |\ h\in H\}$ is the conjugate subgroup.
One should beware that while $c_g^*(\res^G_{H^g}(X))$ and $\res^G_H(X)$
have the same underlying orthogonal spectrum, they come with {\em different} $H$-actions. 
However, left translation by $g$ is an isomorphism of orthogonal $H$-spectra
$l_g^X:c_g^*X\to X$. So combining the restriction map along $c_g$
with the effect of $l_g^X$ 
gives an isomorphism\index{subject}{conjugation homomorphism!on equivariant homotopy groups}  
\begin{equation}\label{eq:define_G_star}
 g_\star\ : \ \pi_0^{H^g}(X)\ \xra{\ (c_g)^*\ }\ \pi_0^H(c_g^*X)\ \xra{\ (l_g^X)_*\ }\
\pi_0^H(X)\ .   
\end{equation}
Moreover, 
\begin{align*}
 g_\star\circ g'_\star\ &= \ (l_g^X)_*\circ (c_g)^* \circ (l_{g'}^X)_*\circ (c_{g'})^*  \ 
= \ (l_g^X)_*\circ ( (c_g)^*(l_{g'}^X))_*\circ (c_g)^* \circ  (c_{g'})^*  \\
&= (l_g^X\circ (c_g)^*(l_{g'}^X))_*\circ (c_{g'}\circ c_g)^*\
= \ (l_{g g'}^X)_*\circ (c_{g g'})^*\ = \ (g g')_\star   
\end{align*}
by naturality of $(c_g)^*$.
If $g$ normalizes $H$, then $g_\star$ is a self-map of the group $\pi_0^H(X)$.
If moreover $g$ belongs to $H$, then $g_\star$ is the identity
by Proposition \ref{prop:inner is identity};
so the maps $g_\star$ define an action of the Weyl group
$W_G H=N_G H/H$ on the equivariant homotopy group $\pi_0^H(X)$.

If $H$ has finite index in its normalizer, this is the end of the story
concerning Weyl group actions on $\pi_0^H(X)$.
In general, however the group $H$ need not have finite index in its normalizer $N_G H$,
and the Weyl group $W_G H$ may have positive dimension, and hence a non-trivial
identity path component $(W_G H)^\circ$. We will now show that the entire
identity path component acts trivially on $\pi_0^H(X)$
for any orthogonal $G$-spectrum $X$. This is a consequence of the
fact that every element of $(W_G H)^\circ$
has the form $z H$ for an element $z$ in $(C_G H)^\circ$,
the identity component of the centralizer of $H$ in $G$,
compare Proposition \ref{prop:1-component of W_G H}.
But then $c_z:H\to H$ is the identity because $z$ centralizes $H$.
On the other hand, any path from $z$ and~1 in $C_G H$
induces a homotopy of morphisms of orthogonal $H$-spectra
from $l_z:X\to X$ to the identity of $X$.
So
\[ z_\star \ = \ (l_z^X)_*\circ (c_z)^* \ = \ \Id_{\pi_0^H(X)}\ .\]
This shows that the identity component of the Weyl group $W_G H$
acts trivially on $\pi_0^H(X)$. So the Weyl group action factors over
an action of the discrete group
\[ \pi_0(W_G H)\ = \ (W_G H) / (W_G H)^\circ \ . \ \]
\end{rk}

\begin{construction}\label{con:postcompose spec}
If $A$ is a pointed $G$-space, then smashing with $A$ and taking based maps
out of $A$ are two continuous based endofunctors
on the category of based $G$-spaces. So for every orthogonal $G$-spectrum $X$,
we can define two new orthogonal $G$-spectra $X\sm A$ and $\map_*(A,X)$
by smashing with $A$ (and letting $G$ act diagonally)
or taking based maps from $A$ levelwise (and letting $G$ act by conjugation).
More explicitly, we have
\[  (X\sm A)(V)\ =\ X(V)\sm A  \text{\quad respectively\quad}
\map_*(A,X)(V) \ = \ \map_*(A,X(V)) \]
for an inner product space $V$. 
The structure maps and actions of the orthogonal groups 
do not interact with $A$: the group $O(V)$ acts through its action on $X(V)$, 
and the structure maps are given by the composite
\[ S^V\sm (X\sm A)(W)
 \ =\ S^V\sm X(W)\sm A \ \xra{\sigma_{V,W}\sm A}\
X(V\oplus W)\sm A\ = \ (X\sm A)(V\oplus W) \]
respectively by the composite
\[ S^V\sm \map_*(A,X(W)) \ \to\ \map_*(A, S^V\sm X(W)) \ \xra{\map_*(A,\sigma_{V,W})}\ 
\map_*(A, X(V\oplus W)) \]
where the first is an assembly map that sends $v\sm f$
to the map sending $a\in A$ to $v\sm f(a)$.

Just as the functors $-\sm A$ and $\map_*(A,-)$ are adjoint
on the level of based $G$-spaces,
the two functors just introduced are an adjoint pair on the level 
of orthogonal $G$-spectra. The adjunction
\begin{equation}\label{eq:adjunction box and map}
G\spec(X,\map_*(A,Y) ) \ \xra{\ \iso\ } \    G\spec(X\sm A,Y)
\end{equation}
takes a morphism $f:X\to \map_*(A,Y)$ to the morphism $f^\flat:X\sm A\to Y$ 
whose $V$-th level $f^\flat(V):X(V)\sm A\to Y(V)$ is $f^\flat(V)(x\sm a)=f(V)(x)(a)$.

An important special case of this construction is
when $A=S^W$ is a representation sphere, i.e., the one-point compactification
of an orthogonal $G$-representation.
The \emph{$W$-th suspension}\index{subject}{suspension!of an orthogonal spectrum}
$X\sm S^W$ is defined by
\[ (X\sm S^W)(V) \ = \ X(V)\sm S^W \ , \]
the smash product of the $V$-th level of $X$ with the sphere $S^W$.
The \emph{$W$-th loop spectrum}\index{subject}{loop spectrum}
$\Omega^W X=\map_*(S^W,X)$, defined by
\[ (\Omega^W X)(V) \ = \ \Omega^W X(V) \ = \ \map_*(S^W,X(V)) \ , \]
the based mapping space from  $S^W$ to the $V$-th level of $X$.
We obtain an adjunction between $-\sm S^W$ and $\Omega^W$
as the special case $A=S^W$ of \eqref{eq:adjunction box and map}.
\end{construction}

\begin{construction}[Shift of an orthogonal spectrum]\label{con:shift}
We introduce a spectrum analog of the additive shift
of orthogonal spaces defined 
in Example \ref{eg:Additive and multiplicative shift}. 
We let $V$ be an inner product space and denote by
\[ -\oplus V\ : \ \bO \ \to \ \bO \]
the continuous functor given on objects by
orthogonal direct sum with $V$, and on morphism spaces by
\[ \bO(U,W) \ \to \ \bO(U\oplus V, W\oplus V)\ , \quad
(w,\varphi)\ \longmapsto \ ((w,0),\varphi\oplus V) \ .\]
The {\em $V$-th shift}\index{subject}{shift!of an orthogonal spectrum}\index{symbol}{$\sh^V X$ - {$V$-th shift of the orthogonal spectrum $X$}} 
of an orthogonal spectrum $X$ is the composite
\begin{equation}  \label{eq:define_shift}
 \sh^V X \ = \ X\circ (-\oplus V)\ . 
\end{equation}
In other words, the value of $\sh^V X$ at an inner product space $U$ is
\[  (\sh^V X)(U) \ = \ X(U\oplus V) \ . \]
The orthogonal group $O(U)$ acts through the monomorphism
$-\oplus V:O(U)\to O(U\oplus V)$.
The structure map $\sigma_{U,W}^{\sh^V X}$ of $\sh^V X$ is 
the structure map $\sigma^X_{U,W\oplus V}$ of $X$.
\end{construction}

Since composition of functors is associative, the shift construction commutes 
on the nose with all constructions on orthogonal spectra that are given by
postcomposition with a continuous based functor as in 
Construction \ref{con:postcompose spec}. This applies in particular to smashing
with and taking mapping space from a based space $A$, i.e., 
\[ (\sh^V X)\sm A\ = \ \sh^V(X\sm A)\text{\quad and\quad}
\map_*(A, \sh^V X) \ = \ \sh^V(\map_*(A, X))\ .\]
So we can -- and will -- omit the parentheses in expressions such as $\sh^V X\sm A$.

The shift construction is also transitive in the following sense.
The values of $\sh^V(\sh^W X)$ and $\sh^{V\oplus W}X$ at an inner product space $U$
are given by
\[ (\sh^V(\sh^W X))(U)\ = \ X((U\oplus V)\oplus W)\]
respectively
\[(\sh^{V\oplus W}X)(U)\ = \ X(U\oplus (V\oplus W)) \ .\]
We use the effect of $X$ on the associativity isomorphism
\[  (U\oplus V)\oplus W\ \iso\ U\oplus (V\oplus W) \ , \quad
((u,v),w)\ \longmapsto\ (u,(v,w))\]
to identify these two spaces; then we abuse notation and write
\[ \sh^V(\sh^W X)\ = \ \sh^{V\oplus W}X\ .\]

The suspension and the shift of an orthogonal spectrum $X$ are related by
a natural morphism\index{symbol}{$\lambda^V_X$ - {natural $\upi_*$.isomorphism of orthogonal $G$-spectra $X\sm S^V\to \sh^V X$}}
\begin{equation}\label{eq:defn lambda_n} 
\lambda_X^V\ :\ X\sm S^V\ \to\ \sh^V X\ .
\end{equation}
In level $U$, this is defined as $\lambda_X^V(U)=\sigma_{U,V}^{\op}$,
the opposite structure map \eqref{eq:define opposite structure},
i.e., the composite
\[ X(U)\sm S^V\ \xra{\text{twist}}\ S^V\sm X(U)\ \xra{\ \sigma_{V,U}\ }
\ X(V\oplus U)\ \xra{X(\tau_{V,U})} \ X(U\oplus V) = (\sh^V X)(U) \ .  \]
In the special case $V=\mR$ we abbreviate $\lambda_X^{\mR}$
to $\lambda_X:X\sm S^1\to\sh X$.\index{symbol}{$\lambda_X$ - {natural global equivalence of orthogonal spectra $X\sm S^1\to \sh X$}}
The $\lambda$-maps are transitive in the sense that for another inner product space $W$, 
the morphism $\lambda^{V\oplus W}_X$ coincides with the two composites
in the commutative diagram: 
\[
\xymatrix{ 
X\sm S^{V\oplus W}\ar[r]^-{\iso} & 
X\sm S^V\sm S^W \ar[rr]^-{\lambda^W_{X\sm S^V}}\ar[d]_{\lambda^V_X\sm S^W} && 
 \sh^W X\sm S^V \ar[d]^{\lambda^V_{\sh^W X}} \\
& \sh^V X\sm S^W \ar[rr]_-{\sh^V(\lambda^W_X)}&& 
\sh^V(\sh^W X)\ar@{=}[r] & \sh^{V\oplus W}X }  \]

Now we let $G$ be a compact Lie group, $V$ a $G$-representation and $X$
an orthogonal $G$-spectrum. Then the orthogonal spectra $X\sm S^V$ and $\sh^V X$
become orthogonal $G$-spectra by letting $G$ act diagonally on $X$ and $V$.
With respect to these diagonal actions, the morphism $\lambda^V_X:X\sm S^V\to\sh^V X$
is a morphism of orthogonal $G$-spectra.
Our next aim is to show that $\lambda^V_X$ is in fact a $\upi_*$-isomorphism.
We define a homomorphism
\begin{equation}  \label{eq:shift2suspension}
 \psi^V_X \ : \ \pi_k^G(\sh^V X) \ \to \ \pi_k^G(X\sm S^V)   
\end{equation}
by sending the class represented by a $G$-map
\[ f\ :\ S^{U\oplus\mR^{n+k}}\ \to\ X(U\oplus \mR^n\oplus V)\ =\ (\sh^V X)(U\oplus\mR^n) \]
to the class represented by the composite
\begin{align*}
 S^{U\oplus V\oplus \mR^{n+k}}\   \xra{S^U\sm \tau_{V,\mR^{n+k}}}\ S^{U\oplus\mR^{n+k}\oplus V}\
\xra{f\sm S^V} \  &X(U\oplus\mR^n\oplus V) \sm S^V \\
\xra{X(U\oplus\tau_{\mR^n,V})\sm S^V}\  &X(U\oplus V\oplus\mR^n) \sm S^V \ .  
\end{align*}
We omit the straightforward verification that this assignment
is compatible with stabilization, and hence well-defined.
The map $\psi^V_X$ is natural for morphisms of orthogonal $G$-spectra in $X$.
Finally, we define
\[ \varepsilon_V\ : \ \pi_k^G(X\sm S^V) \ \to \ \pi_k^G(X\sm S^V) \]
as the effect of the involution
\[ X\sm S^{-\Id_V}\ : \ X\sm S^V \ \to \ X\sm S^V \]
induced by the antipodal map of $S^V$.

\begin{prop}\label{prop:lambda upi_* isos}
Let $G$ be a compact Lie group, $X$ an orthogonal $G$-spectrum 
and $V$ a $G$-representation.
\begin{enumerate}[\em (i)]
\item 
For every integer $k$, each of the three composites around the triangle  
\[ \xymatrix{ 
\pi_k^G(X\sm S^V) \ar[rr]^-{(\lambda^V_X)_*} &&
\pi_k^G(\sh^V X) \ar[dl]^-{\psi^V_X} \\
& \pi_k^G(X\sm S^V) \ar[ul]^-{\varepsilon_V}} \]
is the respective identity.
\item The morphism
\[ \lambda^V_X\ :\  X\sm S^V\ \to\ \sh^V X  \ , 
\text{\qquad its adjoint\qquad}
\tilde\lambda^V_X \ : \ X \ \to \ \Omega^V \sh^V X \ ,\]
the adjunction unit $\eta^V_X:X\to\Omega^V(X\sm S^V)$ 
and the adjunction counit $\epsilon^V_X:(\Omega^V X)\sm S^V\to X$ 
are $\upi_*$-isomorphisms of orthogonal $G$-spectra.
\end{enumerate}
\end{prop}
\begin{proof}
We introduce an auxiliary functor $\pi^G(A;-)$
from orthogonal $G$-spectra to abelian groups that generalizes 
equivariant homotopy groups and depends on a based $G$-space $A$.
We set
\[ \pi^G(A;X) \ = \ \colim_{U\in s(\Uc_G)}\, [S^U\sm A, X(U)]^G \ ,  \]
where the colimit is taken over the analogous stabilization maps as for $\pi_0^G$;
the set $\pi^G(A;X)$ comes with a natural abelian group structure
by the same reasoning as for $\pi_0^G(X)$.
Then $\pi^G(S^k;X)$ is naturally isomorphic to $\pi_k^G(X)$;
more generally, for a $G$-representation $W$,
the adjunction bijections
\[ [S^U\sm S^{\mR^k\oplus W}, X(U)]^G \ \iso \  [S^{U\oplus\mR^k}, \Omega^W X(U)]^G  \]
assemble into a natural isomorphism of abelian groups
between $\pi^G(S^{\mR^k\oplus W};X)$ and $\pi_k^G(\Omega^W X)$.

The definition of the map $\psi^V_X$ has a straightforward generalization
to a natural homomorphism
\[ \psi^V_X \ : \ \pi^G(A; \sh^V X) \ \to \ \pi^G(A;X\sm S^V) \]
by sending the class represented by a $G$-map
\[ f\ :\ S^U\sm A\ \to\ X(U\oplus V)\ =\ (\sh^V X)(U) \]
to the class represented by the composite
\[  S^{U\oplus V}\sm A\ \xra{S^U\sm \tau_{S^V,A}}\   S^U\sm A\sm S^V \
   \xra{f\sm S^V} \  X(U\oplus V) \sm S^V \ . \]
We omit the straightforward verification that this assignment
is compatible with stabilization, and hence well-defined.
We claim that each of the three composites around the triangle  
\[ \xymatrix{ 
\pi^G(A;X\sm S^V) \ar[rr]^-{(\lambda^V_X)_*} &&
\pi^G(A;\sh^V X) \ar[dl]^-{\psi^V_{X}} \\
& \pi^G(A;X\sm S^V) \ar[ul]^-{\varepsilon_V}} \]
is the respective identity. 
We consider a based continuous $G$-map $f:S^U\sm A\to X(U)\sm S^V$
that represents a class in $\pi^G(A;X\sm S^V)$.
Then the class $\varepsilon_V(\psi^V_{X}((\lambda^V_X)_*\td{f}))$
is represented by the composite
\[   S^{U\oplus V}\sm A\ \ \iso \  S^U\sm A\sm S^V\ 
\xra{f\sm S^V}\ X(U)\sm S^V\sm S^V \
\xra{\sigma^{\op}_{U,V}\sm S^{-\Id_V}}\ X(U\oplus V)\sm S^V \ . \]
The map $V\oplus(-\Id_V):V\oplus V\to V\oplus V$
is homotopic, through $G$-equivariant linear isometries, to the
twist map $\tau_{V,V}:V\oplus V\to V\oplus V$ that interchanges the two summands.
So $\varepsilon_V(\psi^V_{X}((\lambda^V_X)_*\td{f}))$ is also represented by
the left vertical composite in the following diagram of based continuous $G$-maps:
\[ \xymatrix@C=22mm@R=6mm{ 
S^U\sm S^V\sm A\ar[d]_{S^U\sm \tau_{A,S^V}} \ar[r]^-{\tau_{U,V}\sm A} & S^V\sm S^U\sm A\ar[ddd]^(.4){S^V\sm f} \\
S^U\sm A\sm S^V\ar[d]_{f\sm S^V} & \\
X(U)\sm S^V\sm S^V\ar[d]_{X(U)\sm \tau_{V,V}}\ar[dr]^-{\tau_{X(U)\sm S^V,S^V}} & \\
X(U)\sm S^V\sm S^V\ar[r]_-{\tau_{X(U),S^V}\sm S^V}\ar[d]_{\sigma^{\op}_{U,V}\sm S^V} &
S^V\sm X(U)\sm S^V\ar[d]^{\sigma_{V,U}\sm S^V} \\
X(U\oplus V)\sm S^V\ar[r]_-{X(\tau_{U,V})\sm S^V} &
X(V\oplus U)\sm S^V } \]
The right vertical composite $(\sigma_{V,U}\sm S^V)\circ(S^V\sm f)$
is the stabilization of $f$, so it
represents the same class in $\pi^G(A;X\sm S^V)$.
Since the left and right vertical composites differ by conjugation with
an equivariant isometry, they also represent the same class in $\pi^G(A;X\sm S^V)$,
by Proposition \ref{prop:invariant description stable}~(ii).
Altogether this shows that the composite
$\varepsilon_V\circ\psi^V_X\circ(\lambda^V_X)_*$ is the identity.
Since $\varepsilon_V^2$ is the identity, this also implies
that the composite $\psi^V_X\circ(\lambda^V_X)_*\circ\varepsilon_V$ is the identity.

The remaining case is similar.
We consider a based continuous $G$-map $g:S^U\sm A\to X(U\oplus V)=(\sh^V X)(U)$
that represents a class in $\pi^G(A;\sh^V X)$.
Then the class $(\lambda^V_X)_*(\varepsilon_V(\psi^V_{X}\td{g}))$
is represented by the composite
\begin{align*}
   S^U\sm S^V\sm A\  &\xra{S^U\sm \tau_{S^V,A}}\   S^U\sm A\sm S^V \
 \xra{g\sm S^V}\ X(U\oplus V)\sm S^V \\
&\xra{X(U\oplus V)\sm S^{-\Id_V}}\ X(U\oplus V)\sm S^V\\
&\xra{\sigma^{\op}_{U\oplus V,V}}\ X(U\oplus V\oplus V) \ = \ (\sh^V X)(U\oplus V)  \ .
\end{align*}
Since 
\[ \sigma^{\op}_{U\oplus V,V}\circ (X(U\oplus V)\sm S^{-\Id}) \ = \ 
X(U\oplus V\oplus (-\Id)) \circ  \sigma^{\op}_{U\oplus V,V}  \]
and $V\oplus(-\Id_V):V\oplus V\to V\oplus V$
is $G$-homotopic to the twist $\tau_{V,V}$,
the class $(\lambda^V_X)_*(\varepsilon_V(\psi^V_X\td{g}))$
is also represented by the left vertical composite in the following diagram:
\[ \xymatrix@C=25mm@R=6mm{ 
S^U\sm S^V\sm A\ar[d]_{S^U\sm\tau_{S^V,A}}\ar[r]^-{\tau_{U,V}\sm A} & 
S^V\sm S^U\sm A\ar[ddd]^{S^V\sm g} \\
S^U\sm A\sm S^V\ar[d]_{g\sm S^V}  & \\
X(U\oplus V)\sm S^V\ar[d]_{\sigma^{\op}_{U\oplus V,V}}\ar[dr]^-{\tau_{X(U)\sm S^V,S^V}} & \\
X(U\oplus V\oplus V)\ar[dr]^-{X(\tau_{U\oplus V,V})} \ar[d]_{X(U\oplus\tau_{V,V})} &
S^V\sm X(U\oplus V)\ar[d]^{\sigma_{V,U\oplus V}} \\
X(U\oplus V\oplus V)\ar[r]_-{X(\tau_{U,V}\oplus V)}\ar@{=}[d] &
X(U\oplus V\oplus V)  \ar@{=}[d] \\
(\sh^V X)(U\oplus V)\ar[r]_-{(\sh^V X)(\tau_{U,V})}  &
(\sh^V X)(U\oplus V) } \]
The right vertical composite $\sigma_{V,U\oplus V}\circ(S^V\sm g)$
is the stabilization of $g$, 
so it represents the same class in $\pi^G(A;\sh^V X)$.
Since the left and right vertical composites differ by conjugation with
an equivariant isometry, they represent the same class,
so the composite $(\lambda^V_X)_* \circ \varepsilon_V\circ\psi^V_X$ is the identity.

Now we prove claim~(i) of the proposition. 
For $k\geq 0$, it is the special case $A=S^k$ of the discussion above.
To deduce the claim for negative dimensional homotopy groups we use  
 the isomorphism of orthogonal $G$-spectra 
 \begin{equation}\label{eq:shift_interchange_iso}
 \tau_{k,V}\ :\ \sh^k(\sh^V X)\ \iso\  \sh^V(\sh^k X)    
 \end{equation}
whose value at an inner product space $U$ is the map
\[ X(U\oplus\tau_{\mR^k,V})\ : \ X(U\oplus\mR^k\oplus V) \ \iso\  
X(U\oplus V\oplus\mR^k)\ . \]
Then the following diagram commutes:
\[ \xymatrix@C=15mm@R=5mm{ 
\pi_{-k}^G(X\sm S^V)\ar[r]^-{(\lambda_X^V)_*}\ar@{=}[dd] &
\pi_{-k}^G(\sh^V X)\ar[r]^-{\varepsilon_V\circ \psi_X^V}\ar@{=}[d] &
\pi_{-k}^G(X\sm S^V)\ar@{=}[dd]\\
& \pi_0^G(\sh^k(\sh^V X))\ar[d]_{(\tau_{k,V})_*}^\iso & \\
\pi_0^G(\sh^k X\sm S^V)\ar[r]_-{(\lambda_{\sh^k X}^V)_*}&
\pi_0^G(\sh^V (\sh^k X))\ar[r]_-{\varepsilon_V\circ \psi_{\sh^k X}^V}&
\pi_0^G(\sh^k X\sm S^V)
} \]
So the claim in dimension $-k$ for the orthogonal $G$-spectrum $X$
is a consequence of the previously established claim
in dimension~0 for the orthogonal $G$-spectrum $\sh^k X$.

(ii) 
We start with the morphism $\tilde\lambda^V_X$, which can be treated fairly
directly. We discuss the case $k\geq 0$ and leave the analogous argument
for $k<0$ to the reader.
We define a map in the opposite direction
\[ \kappa \ : \ \pi_k^G(\Omega^V\sh^V X)\ \to \ \pi_k^G(X)\ . \]
We let $g:S^{U\oplus\mR^k}\to \Omega^V X(U\oplus V)= (\Omega^V\sh^V X)(U)$
represent a class of the left hand side. The map $\kappa$ sends $[g]$
to the class represented by the composite
\[ S^{U\oplus V\oplus\mR^k} \ \xra{S^U\sm \tau_{V,\mR^k}}\
S^{U\oplus \mR^k\oplus V} \ \xra{\ g^\flat\ }\ X(U\oplus V) \ ,\]
where $g^\flat$ is the adjoint of $g$.
This is compatible with stabilization.

We claim that the map $\kappa$ is injective.
Indeed, if $g:S^{U\oplus\mR^k}\to \Omega^V X(U\oplus V)$
represents an element in the kernel of $\kappa$,
then after increasing $U$, if necessary, the composite
$g^\flat\circ(S^U\sm \tau_{V,\mR^k})$ is $G$-equivariantly null-homotopic.
But then $g^\flat$, and hence also its adjoint $g$, are
equivariantly null-homotopic. So $\kappa$ is injective.

The composite $\kappa\circ(\tilde \lambda^V_X)_*$ sends the class of a $G$-map
$f:S^{U\oplus\mR^k}\to X(U)$ to the class of the composite
\begin{equation}\label{eq:another composite}
 S^{U\oplus V\oplus\mR^k} \ \xra{S^U\sm \tau_{V,\mR^k}}\
S^{U\oplus \mR^k\oplus V} \ \xra{\ (\eta_X^V(U)\circ f)^\flat\ }\ X(U\oplus V) \ .  
\end{equation}
The adjoint $(\eta_X^V(U)\circ f)^\flat$ coincides with the composite
\[ S^{U\oplus \mR^k\oplus V} \ \xra{\ f\sm S^V\ }\ X(U)\sm S^V\
\xra{\sigma^{\op}_{U,V}} \ X(U\oplus V)\ , \]
so the composite \eqref{eq:another composite}
represents the same class as $f$.
This proves that $\kappa\circ(\tilde \lambda^V_X)_*$ is the identity.
Since $\kappa$ is also injective, the map $(\tilde\lambda^V_X)_*$
is bijective. 

If $H$ is a closed subgroup of $G$ we apply the previous argument
to the underlying $H$-representation of $V$ to conclude that
$\tilde\lambda^V_X$ induces an isomorphism on $\pi_*^H$.
So  $(\tilde\lambda^V_X)_*$ is a $\upi_*$-isomorphism of orthogonal $G$-spectra.

Now we treat the morphism $\lambda^V_X$.
Again we show a more general statement, namely that
for every pair of $G$-representations $V$ and $W$ the morphism
\[  \Omega^W(\lambda^V_X)\ :\ \Omega^W(X\sm S^V) \ \to\  \Omega^W(\sh^V X)\]
is a $\upi_*$-isomorphism of orthogonal $G$-spectra.
We start with the effect on $G$-equivariant homotopy groups.
For $k\geq 0$ this follows by applying part~(i) with $A=S^{\mR^k\oplus W}$ 
and exploiting the natural isomorphism $\pi_k^G(\Omega^W Y)\iso \pi^G(S^{\mR^k\oplus W};Y)$.
To get the same conclusion
for negative dimensional homotopy groups we exploit that
$\pi_{-k}^G(Y)=\pi_0^G(\sh^k Y)$, by definition.
Moreover, the following diagram commutes
\[ \xymatrix@R=5mm@C=25mm{
\pi_{-k}^G(\Omega^W(X\sm S^V)) \ar@{=}[r] \ar[dd]_{\pi_{-k}^G(\Omega^W(\lambda^V_X))} &
\pi_0^G(\Omega^W (\sh^k X\sm S^V))  \ar[d]^{\pi_0^G(\Omega^W(\lambda^V_{\sh^k X}))} \\
&\pi_0^G(\Omega^W(\sh^V(\sh^k X)))  \\
\pi_{-k}^G(\Omega^W(\sh^V X)) \ar@{=}[r] &
\pi_0^G(\Omega^W(\sh^k(\sh^V X)))
\ar[u]^\iso_{\pi_0^G(\Omega^W \tau_{k,V})}  } \]
where the isomorphism $\tau_{k,V}$ was defined in \eqref{eq:shift_interchange_iso}.
So the previous argument applied to the spectrum $\sh^k X$ shows
that the morphism $\Omega^W(\lambda^V_X)$ also induces isomorphisms 
on $G$-equivariant homotopy groups in negative dimensions.
If $H$ is any closed subgroup of $G$, then we consider the underlying $H$-representations
of $V$ and $W$ and conclude that the morphism $\Omega^W(\lambda^V_X)$ 
induces isomorphisms on $\pi_*^H$.
This proves that $\Omega^W(\lambda^V_X)$ is a $\upi_*$-isomorphism 
of orthogonal $G$-spectra.
The special case $W=0$ proves that $\lambda^V_X$ is a $\upi_*$-isomorphism.

The morphism $\tilde\lambda^V_X$ factors as the composite
\[ X \ \xra{\ \eta^V_X \ }\ \Omega^V(X\sm S^V) \ \xra{\Omega^V(\lambda^V_X)}
\ \Omega^V\sh^V X\ .\]
Since both $\tilde\lambda^V_X$ and $\Omega^V( \lambda^V_X)$ are
$\upi_*$-isomorphisms of orthogonal $G$-spectra,
so is the adjunction unit $\eta^V_X$.

Finally, we treat the adjunction counit $\epsilon^V_X$.
The two homomorphisms of orthogonal $G$-spectra
\[ \Omega^V (\tilde\lambda^V_X)\ , \ \tilde\lambda^V_{\Omega^V X}\ : \ 
\Omega^V X \ \to\ \Omega^V(\Omega^V\sh^V X) \] 
are {\em not} the same; they differ by the involution on the target
that interchanges the two $V$-loop coordinates.
An equivariant homotopy of linear isometries from $\tau_V:V\oplus V\to V\oplus V$
to $(-\Id)\oplus\Id$ thus induces an equivariant homotopy between the morphism
$\Omega^V (\tilde\lambda^V_X)$ and the composite
\[ \Omega^V X \ \xra{\tilde\lambda^V_{\Omega^V X}} \
\Omega^V(\Omega^V\sh^V X) \ \xra{\Omega^V \map_*(S^{-\Id},\sh^V X)} \
\Omega^V(\Omega^V\sh^V X) \ .\]
Passing to adjoints shows that the square of morphisms of orthogonal $G$-spectra
\[ \xymatrix@C=25mm{ 
(\Omega^V X)\sm S^V\ar[r]^-{\epsilon_X^V} \ar[d]_{\lambda^V_{\Omega^V X}} & 
X\ar[d]^{\tilde\lambda^V_X} \\
\sh^V \Omega^V X\ar[r]^-\iso_-{\sh^V\map_*(S^{-\Id},X)}& 
\sh^V \Omega^V X} \]
commutes up to $G$-equivariant homotopy. Since the two vertical morphisms
are $\upi_*$-isomorphisms, and the lower horizontal one is even an isomorphism,
we conclude that $\epsilon^V_X$ is a $\upi_*$-isomorphism.
\end{proof}

Now we recall some important properties of equivariant homotopy groups,
such as stability under suspension
and looping, and the long exact sequences associated to
mapping cones and homotopy fibers.
We define the {\em loop isomorphism}\index{subject}{loop isomorphism}
\begin{equation}\label{eq:loop iso} 
\alpha\ : \  \pi^G_k(\Omega X)\ \to\ \pi^G_{k+1} (X)  \ .
\end{equation}
We represent a given class in $\pi_k^G(\Omega X)$ by a based $G$-map 
$f:S^{V\oplus \mR^{n+k}}\to \Omega X(V\oplus \mR^n)$ and 
let $f^\flat :S^{V\oplus \mR^{n+k+1}}\to X(V\oplus \mR^n)$
denote the adjoint of $f$, which represents an element of
$\pi_{k+1}^G(X)$. Then we set $\alpha[f]=[f^\flat]$.

Next we define the {\em suspension isomorphism}\index{subject}{suspension isomorphism}
\begin{equation}\label{eq:suspension iso}
 -\sm S^1 \ : \ \pi^G_k (X)\ \to \ \pi^G_{k+1}(X\sm S^1) \ .   
\end{equation}
We represent a given class in $\pi_k^G(X)$ by a based $G$-map 
$f:S^{V\oplus \mR^{n+k}}\to X(V\oplus\mR^n)$;
then $f\sm S^1 : S^{V\oplus \mR^{n+k+1}}\to  X(V\oplus\mR^n)\sm S^1$
represents a class in $\pi_{k+1}^G(X\sm S^1)$, and we set $[f]\sm S^1= [f\sm S^1]$.

\begin{prop}\label{prop:loop and suspension isomorphisms}
Let $G$ be a compact Lie group, $X$ an orthogonal $G$-spectrum
and $k$ an integer.
Then the loop isomorphism \eqref{eq:loop iso} 
and the suspension isomorphism \eqref{eq:suspension iso} 
are isomorphisms of abelian groups.
\end{prop}
\begin{proof}
The inverse to the loop map \eqref{eq:loop iso} is given by
sending the class of a $G$-map $S^{V\oplus \mR^{n+k+1}}\to X(V\oplus\mR^n)$
to the class of its adjoint $S^{V\oplus \mR^{n+k}}\to \Omega X(V\oplus\mR^n)$.
The suspension homomorphism is the composite of the two maps
\[ \pi_k^G(X)\ \xra{\ (\eta_X)_*\ } \ \pi_k^G(\Omega(X\sm S^1)) \ \xra[\iso]{\ \alpha\ } \
\pi_{k+1}^G(X\sm S^1)\ ;\]
since $(\eta_X)_*$ is an isomorphism by Proposition \ref{prop:lambda upi_* isos}~(ii),
so is the suspension isomorphism.
\end{proof}

A key feature that distinguishes stable from unstable
equivariant homotopy theory -- and at the same time an important calculational tool --
is the fact that mapping cones give rise to long exact sequences of equivariant
homotopy groups. Our next aim is to establish
this long exact sequence, see Proposition \ref{prop:LES for homotopy of cone and fibre}.

\begin{construction}[Mapping cone and homotopy fiber]\label{con:cone and fiber}
The \emph{reduced mapping cone}\index{subject}{mapping cone}
$C f$ of a morphism of based spaces $f:A\to B$ is defined by
\[ C f \ = \ (A\sm [0,1])\cup_f B \ . \]
Here the unit interval $[0,1]$ is pointed by $0$,
so that $A\sm [0,1]$ is the reduced cone of $A$. 
The mapping cone comes with an embedding $i:B\to C f$
and a projection $p:C f\to A\sm S^1$; the projection
sends $B$ to the basepoint and is given on $A\sm [0,1]$ 
by $p(a,x)=a\sm t(x)$, where
\begin{equation}  \label{eq:define_t}
 t\ :\ [0,1]\ \to\  S^1\text{\qquad is defined as\qquad}
t(x)\ =\ \frac{2x-1}{x(1-x)}\ .   
\end{equation}
What is relevant about the map $t$ is not the precise formula, 
but that it passes to a homeomorphism between the quotient space 
$[0,1]/\{0,1\}$ and $S^1=\mR\cup\{\infty\}$, and that it satisfies $t(1-x)=-t(x)$.

The \emph{mapping cone}\index{symbol}{$C f$ - {mapping cone of $f$}} 
$C f$ of a morphism $f:X\to Y$ of orthogonal $G$-spectra
is now defined levelwise by
\[ (C f)(V) \ = \ C(f(V)) \ = \ (X(V)\sm [0,1])\cup_{f(V)} Y(V) \ ,  \]
the reduced mapping cone of $f(V):X(V)\to Y(V)$. 
The groups $G$ and $O(V)$ act on $(C f)(V)$
through the given action on $X(V)$ and $Y(V)$ and trivially on the interval.
The embeddings $i(V):Y(V)\to C(f(V))$
and projections $p(V):C(f(V))\to X(V)\sm S^1$ assemble into
morphisms of orthogonal $G$-spectra 
\begin{equation}\label{eq:define i and p}
 i\ :\ Y\ \to\ C f \text{\qquad and\qquad} p\ :\ C f\ \to\  X\sm S^1\ .  
\end{equation}
We define a {\em connecting homomorphism}\index{subject}{connecting homomorphism!for equivariant homotopy groups} 
$\partial:\pi^G_{k+1}(C f)\to\pi^G_k(X)$ as the composite
\begin{equation}\label{def:boundary}
\pi^G_{k+1}(C f) \ \xrightarrow{\ \pi^G_{k+1}(p)\ } \
\pi^G_{k+1}(X\sm S^1)\ \xra{-\sm S^{-1}}\ \pi^G_k (X) \ , \end{equation}
the effect of the projection $p:C f\to X\sm S^1$ on homotopy groups, 
followed by the inverse of the suspension isomorphism \eqref{eq:suspension iso}.

The {\em homotopy fiber}\index{subject}{homotopy fiber}\index{symbol}{$F f$ - {homotopy fiber of $f$}}
is the construction `dual' to the mapping cone.
The homotopy fiber of a continuous map $f:A\to B$  of based spaces
is the fiber product
\[ F(f)\ = \ *\times_B B^{[0,1]}\times_B A 
\ = \ 
\{(\lambda,a)\in B^{[0,1]}\times A\ | \ \lambda(0)=*,\, \lambda(1)=f(a)\} \ ,\]
i.e., the space of paths in $B$ starting at the basepoint 
and equipped with a lift of the endpoint to $A$.
As basepoint of the homotopy fiber we take the pair consisting
of the constant path at the basepoint of $B$ and the basepoint of $A$.
The homotopy fiber comes with maps
\[ \Omega B\ \xra{\ i\ }\ F(f)\ \xra{\ p\ }\ A  \ .\]
The map $p$ is the projection to the second factor;
the value of the map $i$ on a based loop $\omega:S^1\to B$ is
\[ i(\omega) = (\omega\circ t,*)  \ , \]
where $t:[0,1]\to S^1$ was defined in \eqref{eq:define_t}.

The homotopy fiber $F(f)$ of a morphism $f:X\to Y$ of orthogonal $G$-spectra
is the orthogonal $G$-spectrum defined by
\[ F(f)(V) \ = \ F(f(V)) \ , \]
the homotopy fiber of $f(V):X(V)\to Y(V)$. 
The groups $G$ and $O(V)$ act on $F(f)(V)$
through the given action on $X(V)$ and $Y(V)$ and trivially on the interval.
Put another way, the homotopy fiber is the pullback in the cartesian square 
of orthogonal $G$-spectra:
\[\xymatrix@C=25mm{ F(f) \ar[d] \ar[r]^-p &
 X \ar[d]^-{(*,f)}\\
  Y^{[0,1]}  \ar[r]_-{\lambda\ \mapsto\ (\lambda(0),\lambda(1))} & Y\times Y }    
\]
The inclusions $i(V):\Omega Y(V)\to F(f(V))$
and projections $p(V):F(f(V))\to X(V)$ assemble into
morphisms of orthogonal $G$-spectra 
\[ i\ :\ \Omega Y\ \to\  F(f) \text{\qquad and\qquad} p\ :\ F(f)\ \to\  X\ .  \]
We define a {\em connecting homomorphism}
\index{subject}{connecting homomorphism} 
$\partial:\pi^G_{k+1}(Y)\to\pi^G_k (F(f))$ as the composite
\[ 
\pi^G_{k+1}(Y) \ \xra{\ \alpha^{-1}} \  \pi^G_k(\Omega Y)
\xrightarrow{\ \pi^G_k(i)\ } \ \pi^G_k (F(f))\ ,
\]
where the first map is the inverse of the loop isomorphism \eqref{eq:loop iso}.
\end{construction}

The proof of exactness for the mapping cone sequence will 
need some elementary homotopies that we spell out in the next proposition.

\begin{prop}\label{prop:homotopies} 
Let $G$ be a topological group.
\begin{enumerate}[\em (i)]
\item For every continuous based map $f:A\to B$ 
of based $G$-spaces the composites
\[ A \xra{\ f\ } B \xra{\ i\ } C f \text{\qquad and \qquad}
F(f) \xra{\ p\ } A \xra{\ f\ }  B\]
are naturally based $G$-null-homotopic. Moreover, the diagram
\[\xymatrix@C=5mm{ 
& CA\cup_f CB \ar[dl]_{p_A\cup *} \ar[dr]^{*\cup p_B} \\
A \sm S^1 \ar[rr]_-{f\sm \tau} && B\sm S^1 }\]
commutes up to natural, based $G$-homotopy, where $\tau$ 
is the sign involution of $S^1$ given by $x\mapsto -x$.
\item
For every based $G$-space $Z$ the map
$p_Z\cup *:CZ\cup_{Z\times 1}CZ\to Z\sm S^1$
which collapses the second cone is a based $G$-homotopy equivalence.
\item
Let $f:A\to B$ and $\beta:Z\to B$ be morphisms of based $G$-spaces 
such that the composite $i\beta:Z\to C f$ is equivariantly null-homotopic. 
Then there exists a based $G$-map $h:Z\sm S^1\to A\sm S^1$ such that 
$(f\sm S^1)\circ h:Z\sm S^1\to B\sm S^1$ 
is equivariantly homotopic to $\beta\sm S^1$.
\end{enumerate}
\end{prop}
\begin{proof}
(i)
We specify natural $G$-homotopies by explicit formulas.
The map $i f:A\to C f$ is null-homotopic by
$A\times [0,1]\to C f, (a,s)\mapsto(a,s)$, 
i.e., the composite of the canonical maps $A\times [0,1]\to A\sm [0,1]$
and $A\sm [0,1]\to C f$. The map $f p:F(f)\to B$ is null-homotopic by
$F(f)\times [0,1]\to B, (\lambda,a,s)\mapsto\lambda(s)$.

The homotopy for the triangle will be glued together from
two pieces. We define a based homotopy 
$H : C A\times [0,1] \to B\sm S^1$  by the formula
\[ H(a,s,u) \ = \ f(a)\sm t(2-s-u)  \]
which is to be interpreted as the basepoint if $2-s-u\geq 1$.
Another based homotopy $H':C B\times [0,1] \to B\sm S^1$ 
is given by the formula
\[ H'(b,s,u) \ = \ b\sm t(s-u) \ ,\]
where $t:[0,1]\to S^1$ was defined in \eqref{eq:define_t}.
This formula is to be interpreted as the basepoint if $s\leq u$.
The two homotopies are compatible in the sense that
\[ H(a,1,u) \ = \ f(a)\sm t(1-u) \ = \ H'(f(a),1,u) \ , \]
for all $a\in A$ and $u\in[0,1]$. So $H$ and $H'$ glue and yield a homotopy
\[  \left( C A\cup_f C B\right) \times [0,1]\iso
(C A\times [0,1])\cup_{f\times [0,1]} (C B\times [0,1]) \ \xra{H\cup H'} 
\ B\sm S^1 \ .\]
For $u=0$ this homotopy starts with the map $*\cup p_B$,
and it ends for $u=1$ with the map $(f\sm \tau)\circ(p_A\cup *)$.

(ii) 
Since the functor $Z\sm-$ is a left adjoint and $Z\sm \{0,1\}\iso Z\times 1$,
the space $C Z\cup_{1\times Z} C Z$ is homeomorphic to the smash product
of $Z$ and the pushout $[0,1]\cup_{\{0,1\}} [0,1]$.
This identification
\[ C Z\cup_{Z\times 1} C Z\ \iso \ Z\sm ([0,1]\cup_{\{0,1\}} [0,1] )\]
turns the map $p_Z$ into the map
\[  Z\sm (t\cup\ast)\ : \ Z\sm ([0,1]\cup_{\{0,1\}} [0,1] )\ \to \ Z\sm S^1 \ . \]
So the claim follows from the fact that the map
$t\cup\ast:[0,1]\cup_{\{0,1\}} [0,1] \to  S^1$ is a based homotopy equivalence.

(iii)
Let $H:C Z= Z\sm [0,1]\to C f$ be a based, equivariant null-homotopy of
the composite $i\beta:Z\to C f$, 
i.e., $H(z,1)=i(\beta(z))$ for all $z\in Z$. 
The composite $p_AH: C Z\to A\sm S^1$ then factors 
as $p_A H=h p_Z$ for a unique $G$-map $h:Z\sm S^1\to A\sm S^1$.
We claim that $h$ has the required property.

To prove the claim we need the $G$-homotopy equivalence
$p_Z\cup *:CZ\cup_{Z\times 1}CZ\to Z\sm S^1$
which collapses the second cone.
We obtain a sequence of equalities and $G$-homotopies
\begin{align*}
(f\sm S^1)\circ h\circ(p_Z\cup *)\ &= \ 
(f\sm S^1)\circ (p_A\cup *)\circ (H\cup C(\beta)) \\
&=\ (B\sm \tau)\circ(f\sm \tau)\circ (p_A\cup *)\circ (H\cup C(\beta)) \\
&\simeq \ (B\sm \tau)\circ(*\cup p_B)\circ (H\cup C(\beta)) \\
&= \ (B\sm\tau)\circ(\beta\sm S^1)\circ(*\cup p_Z) \\
&= \ (\beta\sm S^1)\circ(Z\sm\tau)\circ(*\cup p_Z) \ 
\simeq \ (\beta\sm S^1)\circ(p_Z\cup *) 
\end{align*}
Here $H\cup C(\beta):CZ\cup_{1\times Z}CZ\to 
C f\cup_B CB \iso CA\cup_f CB$
and $\tau$ is the sign involution of $S^1$.
The two homotopies result from part~(i)  
applied to $f$ respectively the
identity of $Z$, and we used the naturality of various constructions.
Since the map $p_Z\cup *$ is a $G$-homotopy equivalence by part~(ii), this proves
that $(f\sm S^1)\circ h$ is $G$-homotopic to $\beta\sm S^1$.
\end{proof}

Now we are ready to prove the long exact homotopy group sequences for
mapping cones and homotopy fibers.

\begin{prop}\label{prop:LES for homotopy of cone and fibre}
\index{subject}{long exact sequence!of equivariant homotopy groups}
For  every compact Lie group $G$ and
every morphism $f:X\to Y$ of orthogonal $G$-spectra
the long sequences of equivariant homotopy groups
\[ \cdots\ \to\  \pi^G_{k+1} (C f) \
\xrightarrow{\quad \partial\quad } \ 
\pi^G_k(X) \ \xrightarrow{\ \pi^G_k(f)\ } \ \pi^G_k(Y) \
\xrightarrow{\pi^G_k(i)} \ \pi^G_k (Cf) \ \to \ \cdots \]
and
\[ \cdots\ \to\  \pi^G_{k+1}(Y) \
\xrightarrow{\quad \partial\quad } \ \pi_k^G ( F(f)) \ 
\xrightarrow{\ \pi^G_k(p)\ } \
\pi^G_k(X) \ \xrightarrow{\ \pi^G_k(f)\ } \ \pi^G_k(Y) \
\ \to \ \cdots \]
are exact.
\end{prop}
\begin{proof}
We start with exactness of the first sequence at $\pi^G_k(Y)$. 
The composite of $f:X\to Y$ and the
inclusion $i:Y\to C f$ is equivariantly null-homotopic, 
so it induces the trivial map on $\pi^G_k$. 
It remains to show that every element in the
kernel of $\pi_k^G(i):\pi^G_k ( Y ) \to \pi^G_k (C f)$
is in the image of $\pi_k^G(f)$. 
We show this for $k\geq 0$; for $k<0$ we can either use a similar argument 
or exploit that $\pi_k^G(X)=\pi_0^G(\sh^{\mR^{-k}}X)$
and shifting commutes with the formation of mapping cones.
We let $\beta:S^{V\oplus \mR^k}\to Y(V)$ represent an element in the kernel
of $\pi_k^G(i)$.
By increasing $V$, if necessary, we can assume that  
the composite of $\beta$ with the inclusion
$i:Y(V)\to (C f)(V)=C(f(V))$ is equivariantly null-homotopic.
Proposition \ref{prop:homotopies}~(iii) provides a $G$-map
$h:S^{V\oplus \mR^k}\sm S^1\to X(V)\sm S^1$ such that 
$(f(V)\sm S^1)\circ h$ is $G$-homotopic to $\beta\sm S^1$. 
The composite
\[ S^{V\oplus\mR^{k+1}}\ \xra{\ h\ }\ X(V)\sm S^1 \ 
\xra{\sigma^{\op}_{V,\mR}}\  X(V\oplus \mR)\] 
represents an equivariant homotopy class in $\pi_k^G(X)$ and we have
\begin{align*}
\pi_k^G(f)\td{\sigma_{V,\mR}^{\op}\circ h} \ &= \ 
 \td{f(V\oplus\mR)\circ\sigma_{V,\mR}^{\op}\circ h}   \\ 
= \  &\td{\sigma^{\op}_{V,\mR}\circ(f(V)\sm S^1)\circ h}  \
= \td{\sigma_{V,\mR}^{\op}\circ(\beta\sm S^1)} \ = \ (-1)^k\cdot \td{\beta} \ .
\end{align*}
To justify the last equation we let $\varphi:V\to V\oplus\mR$ 
denote the embedding of the first summand. 
Then the maps
\[ \sigma_{V,\mR}^{\op}\circ(\beta\sm S^1)\ : \ S^{V\oplus\mR^k\oplus\mR}\ \to\ Y(V\oplus \mR)
\text{\quad and\quad} 
\varphi_*\beta\ :\ S^{V\oplus\mR\oplus\mR^k}\ \to\ Y(V\oplus\mR)\]
differ by the permutation of the source that moves one sphere coordinate
past $k$ other sphere coordinates; this permutation has degree $k$,
so Proposition \ref{prop:invariant description stable} (ii) establishes 
the last equation.
Altogether this shows that the class represented by $\beta$ 
is the image under $\pi_k^G(f)$ 
of the class $(-1)^k\cdot\td{\sigma^{\op}_{V,\mR}\circ h}$.

We now deduce the exactness at $\pi_k^G (C f)$ and $\pi^G_{k-1} ( X )$
by comparing the mapping cone sequence for $f:X\to Y$ 
to the mapping cone sequence for the morphism $i:Y\to Cf$ 
(shifted to the left).
We observe that the collapse map
\[ *\cup p\ : \ C i \iso C Y\cup_f C X \ \to \ X\sm S^1\]
is an equivariant homotopy equivalence, and thus induces an isomorphism 
of equivariant homotopy groups.
Indeed, a homotopy inverse
\[ r \ : \ X\sm S^1 \ \to \ C Y\cup_f C X \]
is defined by the formula
\[ r(x\sm s) \ = \ \begin{cases}
\qquad (x,2s) \quad \in C X & \text{\ for $0\leq s\leq 1/2$, and}\\
(f(x),2-2s) \in C Y & \text{\ for $1/2\leq s\leq 1$,}
\end{cases}\]
which is to be interpreted levelwise.
We omit the $G$-homotopies $r(*\cup p)\simeq \Id$ 
and $(*\cup p)r\simeq \Id$; explicit formulas can be found in the
proof of \cite[Note (4.6.1)]{tomDieck-algebraic topology}. 
Now we consider the diagram 
\[\xymatrix{ 
C f \ar[r]^-{i_i}\ar[dr]_p & Ci \ar[r]^-{p_i}\ar[d]_(.4){*\cup p} & Y\sm S^1\\
& X\sm S^1 \ar[ur]_-{f\sm S^1}}\]
whose upper row is part of the mapping cone sequence
for the morphism $i:Y\to C f$.
The left triangle commutes on the nose and the right triangle commutes
up to $G$-homotopy. We get a commutative diagram
\[\xymatrix@C=17mm{ 
\pi^G_k ( Y ) \ar@{=}[d] \ar[r]^-{\pi_k^G(i)} & 
\pi^G_k ( C f) \ar[r]^-{\pi_k^G(i_i)}\ar@{=}[d] & 
\pi^G_k ( C i) \ar[r]^-(.4){\partial}\ar[d]_{(-\sm S^{-1})\circ(*\cup p)_*}^\iso & 
\pi^G_{k-1}(Y) \ar@{=}[d] \\
\pi^G_k ( Y ) \ar[r]_-{\pi_k^G(i)}& \pi^G_k(C f)\ar[r]_-{\partial} & 
\pi^G_{k-1} (X) \ar[r]_-{\pi_k^G(f)} & 
\pi^G_{k-1}(Y)}\]
(using, for the right square, the naturality of the suspension isomorphism).
By the previous paragraph, applied to $i:Y\to Cf$ instead of $f$,
the upper row is exact at $\pi^G_k(C f)$. Since all vertical maps
are isomorphisms, the original lower row is exact at $\pi^G_k(C f)$.
But the morphism $f$ was arbitrary, so when applied to 
$i:Y\to Cf$ instead of $f$, we obtain that the upper row is
exact at $\pi^G_k( C i)$.  Since all vertical maps
are isomorphisms, the original lower row is exact at $\pi^G_{k-1} (X)$.
This finishes the proof of exactness of the first sequence.

Now we come to why the second sequence is exact. 
We show the claim for $k\geq 0$, the other case being similar.
For every $V\in s(\Uc_G)$
the sequence $F(f)(V)=F(f(V))\to X(V)\to Y(V)$ 
is an equivariant homotopy fiber sequence. 
So for every based $G$-CW-complex $A$, the long sequence of based sets
\begin{align*}
  \cdots \  \to \ [A,\Omega Y(V)]^G \ &\xra{[A,i(V)]} [A,F(f(V))]^G \\ 
&\xra{[A,p(V)]^G } \ [A,X(V)]^G  \ \xra{[A,f(V)]^G } \ [A,Y(V)]^G 
\end{align*}
is exact. We take $A=S^{V\oplus\mR^k}$ and form the colimit over the poset $s(\Uc_G)$.
Since filtered colimits are exact, the resulting sequence of colimits
is again exact, and that proves the second claim.
\end{proof}

\begin{cor}\label{cor-wedges and finite products}
\index{subject}{equivariant homotopy group!of a wedge}
Let $G$ be a compact Lie group and $k$ an integer.
\begin{enumerate}[\em (i)]
\item 
For every family of orthogonal $G$-spectra $\{X_i\}_{i\in I}$ the canonical map
\[ {\bigoplus}_{i\in I}\,  \pi^G_k ( X_i )\ \to \ 
\pi^G_k\left( {\bigvee}_{i\in I}\,  X_i \right)\]
is an isomorphism.
\item
For every finite index set $I$ and 
every family $\{X_i\}_{i\in I}$ of orthogonal $G$-spectra the canonical map
\[ \pi^G_k\left( {\prod}_{i\in I}\, X_i \right)  \ \to \
{\prod}_{i\in I}\,  \pi^G_k ( X_i )\]
is an isomorphism.\index{subject}{equivariant homotopy group!of a product}
\item
For every finite family of orthogonal $G$-spectra
the canonical morphism from the wedge to the product 
induces isomorphisms of $G$-equivariant homotopy groups.
\end{enumerate}
\end{cor}
\begin{proof}
(i) We first show the special case of two summands.
If $X$ and $Y$ are two orthogonal $G$-spectra, 
then the wedge inclusion $j:X\to X\vee Y$ has a retraction. 
So the associated long exact homotopy group sequence of
Proposition \ref{prop:LES for homotopy of cone and fibre}
splits into short exact sequences
\[ 0 \ \to \ \pi^G_k(X) \ \xra{\pi^G_k(j)} \ \pi^G_k(X\vee Y) \ 
\xra{\pi_k^G(i)} \ \pi^G_k(C j) \ \to \ 0 \ .\]
The mapping cone $C j$ is isomorphic to $(C X)\vee Y$
and thus $G$-homotopy equivalent to $Y$. So we can replace 
$\pi^G_k(C j)$ by $\pi^G_k(Y)$ and conclude that $\pi^G_k(X\vee Y)$
splits as the sum of $\pi^G_k(X)$ and $\pi^G_k(Y)$, via the canonical map.
The case of a finite index set $I$ now follows by induction.

In the general case we consider the composite
\[ {\bigoplus}_{i\in I}\,  \pi^G_k ( X_i )\ \to \ 
\pi^G_k\left( {\bigvee}_{i\in I}\, X_i \right) \ \to \ 
{\prod}_{i\in I}\, \pi^G_k (  X_i ) \ ,\]
where the second map is induced by the projections to the wedge summands.
This composite is the canonical map from a direct sum to a product
of abelian groups, hence injective.
So the first map is injective as well. For surjectivity we consider
a $G$-map $f:S^{V\oplus\mR^{n+k}}\to \bigvee_{i\in I} X_i(V\oplus\mR^n)$
that represents an element in the $k$-th $G$-homotopy group of
$ \bigvee_{i\in I} X_i$. 
Since the source of $f$ is compact, there is a finite subset $J$
of $I$ such that $f$ has image in $\bigvee_{j\in J} X_j(V\oplus\mR^n)$,
compare Proposition \ref{prop:compact to wedge}.
Then the given class is in the image of 
$\pi^G_k\left( \bigvee_{j\in J} X_j \right)$;
since $J$ is finite, the class is in the image of the canonical map,
by the previous paragraph.

(ii) The functor $X\mapsto [S^{V\oplus\mR^{n+k}}, X(V\oplus\mR^n)]^G$ commutes with products. 
Finite products commute with filtered colimits, such as the one defining $\pi_k^G$,
so passage to colimits gives the claim.

(iii) For a finite indexing set $I$ and every integer $k$ the composite map
\[ {\bigoplus}_{i\in I}\, \pi^G_k( X_i) \ \to\ \pi^G_k\left({\bigvee}_{i\in I}\, X_i \right) 
\ \to\  \pi^G_k\left({\prod}_{i\in I}\, X_i \right)\ \to\ {\prod}_{i\in I}\, \pi^G_k(X_i)\]
is an isomorphism, where the first and last maps are the canonical ones.
These canonical maps are isomorphisms by parts~(i) respectively~(ii),
hence so is the middle map.
\end{proof}

\Danger The equivariant homotopy group functor $\pi_k^G$
does not in general commute with infinite products. The issue is that
$\pi_k^G$ involves a filtered colimit, and these do not always commute
with infinite products.
However, this defect is cured when we pass to the equivariant or global stable
homotopy category, i.e., $\pi_k^G$  takes `derived' infinite products
to products. We refer to Remark \ref{rk:homotopy of infinite product}
below for more details.

\medskip

We recall that a morphism $f:A\to B$ of orthogonal $G$-spectra
is an {\em h-cofibration}\index{subject}{h-cofibration!of orthogonal spectra} 
if it has the homotopy extension property, 
i.e., given a morphism of orthogonal $G$-spectra $\varphi:B\to X$ and a homotopy
$H:A\sm [0,1]_+\to X$ starting with $\varphi f$,
there is a homotopy $\bar H:B\sm [0,1]_+\to X$ starting with $\varphi$
such that $\bar H\circ(f\sm [0,1]_+)=H$.
For every such h-cofibration $f:A\to B$
the collapse map $c:C f\to B/A$ from the mapping cone to the cokernel
of $f$ is a $G$-homotopy equivalence. Indeed, since the square
\[ \xymatrix{ A\ar[r]^-f \ar[d]_{- \sm 1} & B\ar[d] \\ A\sm [0,1] \ar[r]_-g & C f} \]
is a pushout, the cobase change $g:C A=A\sm [0,1]\to C f$
also has the homotopy extension property.
The cone $C A$ is contractible, so the claim follows from the following
more general statement: for every h-cofibration $i:C\to D$ such that $C$ is contractible,
the quotient map $D\to D/C$ is a based homotopy equivalence.
The standard proof of this fact in the category of topological spaces
(see for example \cite[Prop.\,0.17]{hatcher}
or \cite[Prop.\,5.1.10]{tomDieck-algebraic topology} with $B=\ast$) 
only uses formal properties of homotopies,
and carries over to the category of orthogonal $G$-spectra.

So for every h-cofibration $f:A\to B$ of orthogonal $G$-spectra,
the collapse map $c:C f\to B/A$ is an equivariant homotopy equivalence, 
hence it induces isomorphisms of all equivariant homotopy groups.
We can thus define a modified connecting homomorphism
$\partial:\pi^G_{k+1}(B/A)\to\pi^G_k(A)$ as the composite
\[
\pi^G_{k+1}(B/A) \ \xrightarrow{(\pi^G_{k+1}(c))^{-1}} \
\pi^G_{k+1}( C f)\ \xra{\ \partial\ }\ \pi^G_k (A) \ . 
\]

We call a morphism $f:X\to Y$ of orthogonal $G$-spectra 
a {\em strong level fibration}\index{subject}{strong level fibration!of orthogonal $G$-spectra} 
if for every closed subgroup $H$ of $G$ and every $H$-representation $V$ 
the map $f(V)^H:X(V)^H\to Y(V)^H$ is a Serre fibration.
For every such strong level fibration,
the embedding $j:F\to F(f)$ of the strict fiber into the homotopy fiber
then induces isomorphisms on $\pi_*^G$.
We can thus define a modified connecting homomorphism
$\partial:\pi^G_{k+1}(Y)\to\pi^G_k (F)$ as the composite
\[ \pi^G_{k+1}(Y) \ \xra{\ \partial\ } \ \pi^G_k (F(f))\ 
\ \xra{(\pi_k^G(j))^{-1}} \ \pi^G_k (F)\ . \]
So we deduce:

\begin{cor}\label{cor-long exact sequence h-cofibration}
Let $G$ be a compact Lie group.
\begin{enumerate}[\em (i)]
\item 
Let $f:A\to B$ be an h-cofibration of orthogonal $G$-spectra.
Then the long sequence of equivariant homotopy groups
\[ \cdots \ \to \ \pi_{k+1}^G (B/A) \ \xrightarrow{\ \partial\ } \ \pi_k^G (A) \
\xra{\ \pi_k^G(f)\ } \ \pi_k^G(B) \ \xra{\ \pi_k^G(q) \ } \ 
\pi_k^G (B/A) \ \to \ \cdots \]
is exact.
\item
Let $f:X\to Y$ be a strong level fibration of orthogonal $G$-spectra
and $j:F\to X$ the inclusion of the pointset level fiber of $f$.
Then the long sequence of equivariant homotopy groups
\[ \cdots \ \to \ 
\pi_{k+1}^G (Y) \ \xra{\ \partial\ } \  \ \pi_k^G (F)
\ \xrightarrow{\ \pi_{k-1}^G(j)\ } \
\pi_k^G(X) \ \xra{\ \pi_k^G(f)\ } \ \pi_k^G(Y) 
 \ \to \ \cdots \]
is exact.
\end{enumerate}
\end{cor}

\begin{cor}\label{cor-cobase change h-cofibration}
  Let $G$ be a compact Lie group and
  \[ \xymatrix{ A \ar[r]^-f \ar[d]_g & B \ar[d]^h\\
    C \ar[r]_-k & D } \]
  a commutative square of orthogonal $G$-spectra.
\begin{enumerate}[\em (i)]
 \item Suppose that the square is a pushout and
   the map $\pi_*^G(f):\pi_*^G(A)\to \pi_*^G(B)$ 
   of $G$-equivariant homotopy groups is an isomorphism.
   If in addition $f$ or $g$ is an h-cofibration, then 
   the map $\pi_*^G(k):\pi_*^G(C)\to \pi_*^G(D)$ is also an isomorphism.
 \item Suppose that the square is a pullback and
   the map $\pi_*^G(k):\pi_*^G(C)\to \pi_*^G(D)$ 
   of $G$-equivariant homotopy groups is an isomorphism.
   If in addition $k$ or $h$ is a strong level fibration, then 
   the map $\pi_*^G(f):\pi_*^G(A)\to \pi_*^G(B)$ is also an isomorphism.
 \end{enumerate}
\end{cor}
\begin{proof}
(i) If $f$ is an h-cofibration, then its long exact homotopy group
sequence (Corollary \ref{cor-long exact sequence h-cofibration})
shows that all $G$-equivariant homotopy groups of the cokernel $B/A$
are trivial. Since the square is a pushout, the induced morphism from
$B/A$ to any cokernel $D/C$ of $k$ is an isomorphism, so the groups
$\pi_*^G(D/C)$ are all trivial. As a cobase change of the
h-cofibration $f$, the morphism $k$ is again an h-cofibration,
so its long exact homotopy group sequence shows that $\pi_*^G(k)$
is an isomorphism.

If $g$ is an h-cofibration, then so is its cobase change $h$.
Moreover, any cokernel $C/A$ of $g$ maps by an isomorphism to 
any cokernel $D/B$ of $h$, since the square is a pushout.
The square induces compatible maps between 
the two long exact homotopy group sequences of $g$ and $h$,
and the five lemma then shows that $\pi_*^G(k)$ is an isomorphism.
The argument for part~(ii) is dual to~(i).
\end{proof}

Since the group $\pi_k^G(\Omega^m X)$ is naturally isomorphic
to $\pi_{k+m}^G(X)$, looping preserves equivariant stable equivalences. 
The next proposition generalizes this.

\begin{prop}\label{prop:map(A,-) preserves global}
Let $G$ be a compact Lie group and $A$ a finite based $G$-CW-complex.
\begin{enumerate}[\em (i)]
\item Let $f:X\to Y$ be a morphism of orthogonal $G$-spectra with the
following property: for every closed subgroup $H$ of $G$ that 
fixes some non-basepoint of $A$,
the map $\pi_*^H(f):\pi_*^H(X)\to\pi_*^H(Y)$ is an isomorphism.
Then the morphism $\map_*(A,f):\map_*(A,X)\to \map_*(A,Y)$ is a
$\upi_*$-isomorphism of orthogonal $G$-spectra.
\item
The functor $\map_*(A,-)$ preserves $\upi_*$-isomorphisms of orthogonal $G$-spec\-tra.
\end{enumerate}\end{prop}
\begin{proof}
(i) We start with a special case and let $X$ be an orthogonal $G$-spectrum
whose equivariant homotopy groups $\pi_*^H(X)$ vanish for all closed subgroups
$H$ of $G$ that fix some non-basepoint of $A$. 
We show first that then the $G$-equivariant homotopy groups
of the orthogonal $G$-spectrum $\map_*(A,X)$ vanish.

We argue by induction over the number of equivariant cells in a CW-structure of $A$.
The induction starts when $A$ consists only of the basepoint,
in which case $\map_*(A,X)$ is a trivial spectrum and there is nothing to show.
For the inductive step we assume that the groups $\pi_*^G(\map_*(B,X))$ 
vanish and $A$ is obtained from $B$ by attaching an equivariant $n$-cell 
$G/H\times D^n$ along its boundary $G/H\times \partial D^n$, 
for some closed subgroup $H$ of $G$.
Then the restriction map $\map_*(A,X)\to\map_*(B,X)$
is a strong level fibration of orthogonal $G$-spectra whose fiber is
isomorphic to
\[ \map_*(A/B,X)\ \iso \ \map_*(G/H_+\sm S^n,X)\ \iso \ \map_*(G/H_+,\Omega^n X)\ . \]
The $G$-equivariant stable homotopy groups of this spectrum
are isomorphic to the $H$-equivariant homotopy groups of $\Omega^n X$,
and hence to the shifted $H$-homotopy groups of $X$.
Since $H$ occurs as a stabilizer of a cell of $A$, 
the latter groups vanish by assumption.
The long exact sequence of Corollary \ref{cor-long exact sequence h-cofibration}~(ii)
and the inductive hypothesis for $B$ then show that the groups 
$\pi_*^G(\map_*(A,X))$ vanish.

When $H$ is a proper closed subgroup of $G$ we exploit
that the underlying $H$-space of $A$ is $H$-equivariantly homotopy equivalent
to a finite $H$-CW-complex, see \cite[Cor.\,B]{illman-restricting equivariance}.
We can thus apply the previous paragraph to the underlying $H$-spectrum of $X$
and the underlying $H$-space of $A$ and conclude that
the $H$-equivariant homotopy groups of $\map_*(A,X)$ vanish.
Altogether this proves the special case of the proposition.

The functor $\map_*(A,-)$ commutes with homotopy fibers;
so two applications of the long exact homotopy group sequence
of a homotopy fiber (Proposition \ref{prop:LES for homotopy of cone and fibre})
reduce the general case of the first claim to the special case.
Part~(ii) is a special case of~(i).
\end{proof}

\begin{prop}\label{prop:sequential colimit closed embeddings}
  Let $G$ be a compact Lie group.
  \begin{enumerate}[\em (i)]
  \item Let $e_n:X_n\to X_{n+1}$ be morphisms 
    of orthogonal $G$-spectra that are levelwise closed embeddings, for $n\geq 0$. 
    Let $X_\infty$ be a colimit of the sequence $\{e_n\}_{n\geq 0}$.
    Then for every integer $k$ the canonical map
    \[ \colim_{n\geq 0}\, \pi_k^G(X_n)\ \to \ \pi_k^G(X_\infty) \]
    is an isomorphism.
  \item Let $e_n:X_n\to X_{n+1}$ and $f_n:Y_n\to Y_{n+1}$ be morphisms 
    of orthogonal $G$-spectra that are levelwise closed embeddings, for $n\geq 0$. 
    Let $\psi_n:X_n\to Y_n$ be $\upi_*$-isomorphisms of orthogonal $G$-spectra 
    that satisfy $\psi_{n+1}\circ e_n=f_n\circ\psi_n$ for all $n\geq 0$.
    Then the induced morphism $\psi_\infty:X_\infty\to Y_\infty$ 
    between the colimits of the sequences is a $\upi_*$-isomorphism.
  \item Let $f_n:Y_n\to Y_{n+1}$ be $\upi_*$-isomorphisms of orthogonal $G$-spectra
   that are levelwise closed embeddings, for $n\geq 0$. 
   Then the canonical morphism 
   $f_\infty:Y_0\to Y_\infty$ to a colimit of the sequence $\{f_n\}_{n\geq 0}$
   is a $\upi_*$-isomorphism.
  \end{enumerate}
\end{prop}
\begin{proof}
 (i) We let $V$ be a $G$-representation, and $m\geq 0$ such that $m+k\geq 0$.
 Since the sphere $S^{V\oplus\mR^{m+k}}$ is compact and $X_\infty(V\oplus\mR^m)$ 
 is a colimit of the sequence of closed embeddings
 $X_n(V\oplus\mR^m)\to X_{n+1}(V\oplus\mR^m)$,
 any based $G$-map $f:S^{V\oplus\mR^{m+k}}\to X_\infty(V\oplus\mR^m)$
 factors through a continuous map
 \[ \bar f\ :\ S^{V\oplus\mR^{m+k}}\ \to\ X_n(V\oplus\mR^m) \]
 for some $n\geq 0$, see for example \cite[Prop.\,2.4.2]{hovey-book}
 or Proposition \ref{prop:filtered colim preserve weq}~(i).
 Since the canonical map $X_n(V\oplus\mR^m)\to X_\infty(V\oplus\mR^m)$ is injective,
 $\bar f$ is again based and $G$-equivariant. The same applies to homotopies,
 so the canonical map
 \[ \colim_{n\geq 0} \, [S^{V\oplus\mR^{m+k}},X_n(V\oplus\mR^m)]^G \ \to \ 
 [S^{V\oplus\mR^{m+k}},X_\infty(V\oplus\mR^m)]^G \]
 is bijective. Passing to colimits over $m$ and over the poset $s(\Uc_G)$ 
 proves the claim.

 Part~(ii) is a direct consequence of~(i), applied to $G$ and all its closed subgroups.
 Part~(iii) is a special case of part~(ii) where we set $X_n=Y_0$,
 $e_n=\Id_{Y_0}$ and $\psi_n=f_{n-1}\circ\dots\circ f_0:Y_0\to Y_n$.
 The morphism $\psi_n$ is then a $\upi_*$-isomorphism
 and $Y_0$ is a colimit of the constant first sequence.
 Since the morphism $\psi_\infty$ induced on the colimits 
 is the canonical morphism $Y_0\to Y_\infty$, part~(ii) specializes to claim~(iii).
\end{proof}

\begin{eg}[Equivariant suspension spectrum]
Every based $G$-space $A$ gives rise to a {\em suspension spectrum}\index{subject}{suspension spectrum!of a $G$-space}
$\Sigma^\infty A$. This is the orthogonal $G$-spectrum with $V$-term
\[ (\Sigma^\infty A)(V) \ = \ S^V\sm A\ ,\]
with $O(V)$-action only on $S^V$, with $G$-action only on $A$, 
and with structure map $\sigma_{U,V}$ given by the canonical homeomorphism
$S^U\sm S^V\sm A\iso S^{U\oplus V}\sm A$.
For an unbased $G$-space $Y$, we obtain an {\em unreduced suspension spectrum}
$\Sigma^\infty_+ Y$ by first adding a disjoint basepoint to $Y$
and then forming the suspension spectrum as above.
\end{eg}

The suspension spectrum functor is homotopical on a large class
of $G$-spaces; however, since the reduced suspension $S^1\sm -$ does not preserve
weak equivalences in complete generality, a non-degeneracy condition
on the basepoint is often needed.

\begin{defn}
A based $G$-space is {\em well-pointed}\index{subject}{well-pointed}
if the basepoint inclusion has the homotopy extension property
in the category of unbased $G$-spaces.  
\end{defn}

Equivalently, a based $G$-space $(A,a_0)$ is well-pointed
if and only if the inclusion of $A\times \{0\}\cup \{a_0\}\times[0,1]$
into $A\times[0,1]$ has a continuous $G$-equivariant retraction.
Cofibrant based $G$-spaces are always well-pointed.
Also, if we add a disjoint basepoint to any unbased $G$-space,
the result is always well-pointed.
The reduced suspension of a well-pointed space is again well-pointed.
Indeed, if $(A,a_0)$ is well-pointed, then the subspace inclusion
\[ A\times\{0,1\}\cup \{a_0\}\times [0,1]\ \to \ A\times[0,1] \]
is an h-cofibration, see for example
\cite[Satz 2]{puppe-bemerkungen} or \cite[Prop.\,(5.1.6)]{tomDieck-algebraic topology}.
Since  h-cofibrations are stable under cobase change
(see Corollary \ref{cor-h-cofibration closures}~(ii)), 
the inclusion of the basepoint into the quotient space
\[  A\times[0,1]/ (A\times \{0,1\}\cup \{a_0\}\times [0,1]) \]
is an h-cofibration.
Since $S^1\sm A$ is homeomorphic to this quotient space, this proves the claim.

We recall that a continuous map $f:X\to Y$ is {\em $m$-connected},
for some $m\geq 0$,\index{subject}{connectivity!of a continuous map}
if the following condition holds: for $0\leq k\leq m$ 
and all continuous maps $\alpha:\partial D^k\to X$
and $\beta:D^k\to Y$ such that $\beta|_{\partial D^k}=f\circ\alpha$,
there is a continuous map $\lambda:D^k\to X$
such that $\lambda|_{\partial D^k}=\alpha$ and
such that $f\circ \lambda$
is homotopic, relative to $\partial D^k$, to $\beta$.
An equivalent condition is that the map $\pi_0(f):\pi_0(X)\to \pi_0(Y)$ is surjective,
and for all $x\in X$ the map
$\pi_k(f):\pi_k(X,x)\to \pi_k(Y,f(x))$ is bijective for
all $1\leq k<m$ and surjective for $k=m$.
The map $f$ is a weak homotopy equivalence if and only if it is $m$-connected
for every $m\geq 0$.

\begin{prop}\label{prop:suspension spectrum homotopical}
Let $G$ be a compact Lie group.
\begin{enumerate}[\em (i)]
\item 
Let $f:X\to Y$ be a continuous based $G$-map 
between well-pointed based $G$-spaces, and $m\geq 0$.
Suppose that for every closed subgroup $H$ of $G$ 
the fixed point map $f^H:X^H\to Y^H$ is $(m-\dim(W_G H))$-connected.
Then the induced map of $G$-equivariant stable homotopy groups
\[ \pi_k^G(\Sigma^\infty f)\  :\ \pi_k^G(\Sigma^\infty X) \ \to 
\ \pi_k^G(\Sigma^\infty Y) \]
is bijective for every integer $k$ with $k<m$ and surjective for $k=m$.
\item For every well-pointed based $G$-space $X$ and every negative integer $k$,
the $G$-equivariant homotopy group $\pi_k^G(\Sigma^\infty X)$ is trivial. 
\end{enumerate}
In particular, the suspension spectrum functor $\Sigma^\infty$ takes
$G$-weak equivalences between well-pointed $G$-spaces to $\upi_*$-isomorphisms,
and the unreduced suspension spectrum functor $\Sigma^\infty_+$ takes
$G$-weak equivalences between arbitrary unbased $G$-spaces to $\upi_*$-isomorphisms.
\end{prop}
\begin{proof}
(i)
We let $V$ be a $G$-representation and we let $n\geq 0$ be such that $n+k\geq 0$. 
Reduced suspension increases the connectivity of continuous maps between
well-pointed spaces, compare \cite[Cor.\,6.7.10]{tomDieck-algebraic topology},
and it preserves well-pointedness.
So the natural homeomorphism 
\[  (S^{V\oplus\mR^n}\sm X)^H \ \iso\ S^{V^H\oplus\mR^n}\sm X^H \]
shows that the map
\[ (S^{V\oplus\mR^n}\sm f)^H \  :\ (S^{V\oplus\mR^n}\sm X)^H \ \to \ (S^{V\oplus\mR^n}\sm Y)^H \]
is $(\dim(V^H)+m+n-\dim(W_G H))$-connected.
On the other hand, we claim that the cellular dimension of $S^{V\oplus\mR^{n+k}}$
at $H$,
in the sense of \cite[II.2, p.\,106]{tomDieck-transformation},
is at most $\dim(V^H)+n+k-\dim(W_G H)$. 
Indeed, if any $G$-CW-structure for $S^{V\oplus\mR^{n+k}}$ contains
an equivariant cell of the form $G/H\times D^j$, then 
$(G/H)^H\times\mathring{D}^j$ embeds into $S^{V^H\oplus\mR^{n+k}}$, and hence
\[ \dim(W_G H)+ j \ = \  \dim((G/H)^H)+ j \ \leq \ \dim(V^H)+ n+ k \ . \]
The cellular dimension at $H$ is the maximal $j$ that occurs in this way, so 
this cellular dimension is at most $\dim(V^H)+n+k-\dim(W_G H)$. 
We conclude that for $k<m$ the cellular dimension of
$S^{V\oplus\mR^{n+k}}$ at $H$ is smaller than the connectivity of $(S^{V\oplus\mR^n}\sm f)^H$,
and for $k=m$ the former is less than or equal to the latter.
So the induced map
\[ [S^{V\oplus\mR^{n+k}},S^{V\oplus\mR^n}\sm f]^G \ : \ 
[S^{V\oplus\mR^{n+k}},S^{V\oplus\mR^n}\sm X]^G \ \to \
[S^{V\oplus\mR^{n+k}},S^{V\oplus\mR^n}\sm X]^G  \]
is bijective for $k<m$ and surjective for $k=m$, 
by \cite[II Prop.\,2.7]{tomDieck-transformation}.
Passage to the colimit over $V\in s(\Uc_G)$ and $n$ then proves the claim.  

(ii) The unique map $f:X\to \ast$ to a one-point $G$-space has the property 
that $f^H$ is 0-connected for every closed subgroup $H$ of $G$. 
So for every negative integer $k$ the induced map from
$\pi_k^G(\Sigma^\infty X)$ to $\pi_k^G(\Sigma^\infty \ast)$ is an isomorphism
by part~(i). Since the latter group is trivial, this proves the claim.
\end{proof}

We end this section with a useful representability result
for the functor $\pi_0^H$ on the category of orthogonal $G$-spectra,
where $H$ is any closed subgroup of a compact Lie group $G$.
We define a {\em tautological class}\index{subject}{tautological class!for $\pi_0^H$}
\begin{equation}\label{eq:tautological_e_H}
 e_H\ \in\ \pi_0^H(\Sigma^\infty_+ G/H)   
\end{equation}
as the class represented by the $H$-fixed point $e H$ of $G/H$. 
\index{symbol}{$e_H$ - {stable tautological class in $\pi_0^H(\Sigma^\infty_+ G/H)$}}

\begin{prop}\label{prop:G/H represents pi_0^H}  
Let $H$ be a closed subgroup of a compact Lie group $G$ and 
$\Phi:G\spec\to \text{\em (sets)}$ a set valued functor 
on the category of orthogonal $G$-spectra that takes $\upi_*$-isomorphisms
to bijections. 
Then evaluation at the tautological class is a bijection
\[   \Nat^{G\spec}(\pi_0^H,\Phi) 
\ \to \ \Phi(\Sigma^\infty_+ G/H) \ , \quad
\tau\ \longmapsto \ \tau(e_H)\ .\]
\end{prop}
\begin{proof}
To show that the evaluation map is injective 
we show that every natural transformation $\tau:\pi_0^H\to \Phi$ 
is determined by the element $\tau(e_H)$.
We let $X$ be any orthogonal $G$-spectrum and $x\in\pi_0^H(X)$
an $H$-equivariant homotopy class.
The class $x$ is represented by a continuous based $H$-map
\[ f \ : \ S^U\ \to \ X(U) \]
for some $H$-representation $U$. By increasing $U$, if necessary,
we can assume that $U$ is underlying a $G$-representation.
We can then view $f$ as an $H$-fixed point 
in $(\Omega^U\sh^U X)(0)=\Omega^U X(U)$.
There is thus a unique morphism of orthogonal $G$-spectra
\[ f^\sharp \ : \ \Sigma^\infty_+ G/H \ \to \ \Omega^U\sh^U X \]
such that $f^\sharp(0):G/H_+\to \Omega^U X(U)$ sends the distinguished coset $e H$ to $f$. 
This morphism then satisfies
\[ \pi_0^H(f^\sharp)(e_H) \ = \ \pi_0^H(\tilde\lambda^U_X)(x) 
\text{\qquad in\ } \ \pi_0^H(\Omega^U\sh^U X) \ ,\]
where $\tilde\lambda^U_X:X\to\Omega^U\sh^U X$ is the $\upi_*$-isomorphism
discussed in Proposition \ref{prop:lambda upi_* isos}~(ii).
The diagram
\[ \xymatrix@C=20mm{
\pi_0^H(\Sigma^\infty_+ G/H)\ar[d]_\tau \ar[r]^-{\pi_0^H(f^\sharp)}&
\pi_0^H(\Omega^U\sh^U X)\ar[d]_\tau 
 & \pi_0^H(X)\ar[d]^\tau\ar[l]^-\iso_-{\pi_0^H(\tilde\lambda^U_X)} \\
\Phi(\Sigma^\infty_+ G/H) \ar[r]_-{\Phi(f^\sharp)} &
\Phi(\Omega^U\sh^U X) 
 & \Phi(X)\ar[l]^-{\Phi(\tilde\lambda^U_X)}_-\iso } \]
commutes and the two right horizontal maps are bijective.
So
\begin{align*}
   \Phi(f^\sharp)(\tau(e_H)) \ = \   \tau(\pi_0^H(f^\sharp)(e_H)) \ = \ 
\tau(\pi_0^H(\tilde\lambda^U_X)(x)) 
\ = \ \Phi(\tilde\lambda^U_X)(\tau(x)) \ .
\end{align*}
Since $\Phi(\tilde\lambda^U_X)$ is bijective, this proves that $\tau$
is determined by its value on the tautological class $e_H$.

It remains to construct, for every element 
$y\in \Phi(\Sigma^\infty_+ G/H)$, a natural transformation
$\tau:\pi_0^H\to \Phi$ with $\tau(e_H)=y$.
The previous paragraph dictates what to do: 
we represent a given class $x\in\pi_0^H(X)$  by
a continuous based $H$-map $f:S^U\to X(U)$,
where $U$ is underlying a $G$-representation.
Then we  set
\[   \tau(x) \ = \ \Phi(\tilde\lambda^U_X)^{-1}(\Phi(f^\sharp)(y)) \ .\]
We verify that the element $\tau(x)$ only depends on the class $x$.
To this end we need to show that $\tau(x)$ does not change if we either
replace $f$ by a homotopic $H$-map, or increase it by stabilization along
a $G$-equivariant linear isometric embedding.
If $\bar f:S^U\to X(U)$ is $H$-equivariantly homotopic to $f$, then the morphism
$\bar f^\sharp$ is homotopic to $f^\sharp$
via a homotopy of morphisms of orthogonal $G$-spectra
\[ K \ : \ (\Sigma^\infty_+ G/H)\sm [0,1]_+\  \to \Omega^U\sh^U X\ . \]
The morphism $q:(\Sigma^\infty_+ G/H)\sm [0,1]_+\to \Sigma^\infty_+ G/H$
that maps $[0,1]$ to a single point is a homotopy equivalence, hence
a $\upi_*$-isomorphism, of orthogonal $G$-spectra.
So $\Phi(q)$ is a bijection. The two embeddings
$i_0,i_1:(\Sigma^\infty_+ G/H)\to \Sigma^\infty_+ G/H\sm [0,1]_+$
as the endpoints of the interval are right inverse to $q$,
so $\Phi(q)\circ\Phi(i_0)=\Phi(q)\circ\Phi(i_1)=\Id$.
Since $\Phi(q)$ is bijective, $\Phi(i_0)=\Phi(i_1)$.
Hence
\[ \Phi(\bar f^\sharp)\ = \ \Phi(K\circ i_0)\ = \ \Phi(K)\circ\Phi(i_0)\ = \ 
\Phi(K)\circ\Phi(i_1)\ = \ \Phi(K\circ i_1)\ = \ \Phi(f^\sharp)\ . \]
This shows that $\tau(x)$ does not change if we modify $f$ by an $H$-homotopy.

Now we let $V$ be another $G$-representation and $\varphi:U\to V$ 
a $G$-equivariant linear isometric embedding. Then $\varphi_*f:S^V\to X(V)$
is another representative of the class $x$.
A morphism of orthogonal $G$-spectra
\[ \varphi_\sharp\ : \ \Omega^U\sh^U X \ \to \ \Omega^V\sh^V X\]
is defined at an inner product space $W$ as the stabilization map
\begin{align*}
   \varphi_*\ : \  (\Omega^U\sh^U X)(W)\ &= \ \map_*(S^U,X(W\oplus U)) \\ 
&\to \ \map_*(S^V,X(W\oplus V)) \ = \ (\Omega^V\sh^V X)(W) \ .
\end{align*}
This morphism makes the following diagram commute:
\[ \xymatrix@C=15mm@R=7mm{ 
\Sigma^\infty_+ G/H\ar[r]^-{f^\sharp} \ar[dr]_(.4){ (\varphi_* f)^\sharp} & 
\Omega^U\sh^U X \ar[d]_{\varphi_\sharp} &
X\ar[l]_-{\tilde\lambda^U_X}\ar[dl]^(.4){\tilde\lambda^V_X}\\
&\Omega^V\sh^V X & } \]
So
\[ \Phi(\tilde\lambda^U_X)^{-1}\circ \Phi(f^\sharp)\ = \
 \Phi(\tilde\lambda^V_X)^{-1}\circ\Phi(\varphi_\sharp)\circ\Phi(f^\sharp) 
\ = \  \Phi(\tilde\lambda^V_X)^{-1}\circ \Phi( (\varphi_* f)^\sharp) \ ,\]
and hence the class $\tau(x)$ remains unchanged upon stabilization of $f$
along $\varphi$.

Now we know that $\tau(x)$ is independent of the choice of representative
for the class $x$, and it remains to show that $\tau$ is natural.
But this is straightforward: if $\psi:X\to Y$ is a morphism of orthogonal $G$-spectra
and $f:S^U\to X(U)$ a representative for $x\in\pi_0^H(X)$ as above,
then $\psi(U)\circ f:S^U\to Y(U)$ represents the class $\pi_0^H(\psi)(x)$.
Moreover, the following diagram of orthogonal $G$-spectra commutes:
\[ \xymatrix@C=15mm{ 
\Sigma^\infty_+ G/H\ar[r]^-{f^\sharp} \ar[dr]_(.4){ (\psi(U)\circ f)^\sharp} & 
\Omega^U\sh^U X \ar[d]^{\Omega^U\sh^U \psi} &
X\ar[l]_-{\tilde\lambda^U_X}\ar[d]^{\psi}\\
&\Omega^U\sh^U Y & Y \ar[l]^-{\tilde\lambda^U_Y}} \]
So naturality follows:
\begin{align*}
  \tau(\pi_0^H(\psi)(x))\ &= \ 
(\Phi(\tilde\lambda^U_Y)^{-1}\circ \Phi( (\psi(U)\circ f)^\sharp))(y) \\
 &= \ (\Phi(\tilde\lambda^U_Y)^{-1}\circ\Phi(\Omega^U\sh^U \psi)\circ\Phi(f^\sharp))(y) \\
 &= \ (\Phi(\psi)\circ\Phi(\tilde\lambda^U_X)^{-1}\circ\Phi(f^\sharp))(y) 
\ = \ \Phi(\psi)(\tau(x)) \ .
\end{align*}
Finally, the class $e_H$ is represented by the $H$-map $-\sm e H:S^0\to S^0\sm G/H_+$,
which is adjoint to the identity of $\Sigma^\infty_+ G/H$. Hence
$\tau(e_H)= \Phi(\Id)(y) =y$.
\end{proof}

\section{The Wirthm{\"u}ller isomorphism and transfers}
\label{sec:Wirthmuller and transfer}

This section establishes the Wirthm{\"u}ller isomorphism
that relates the equivariant homotopy groups 
of a spectrum over a subgroup to the equivariant homotopy groups 
of the induced spectrum, see Theorem \ref{thm:Wirth iso}.
Intimately related to the Wirthm{\"u}ller isomorphism are various transfers 
that we discuss in Constructions \ref{con:external transfer} and \ref{con:transfer}.
We show that transfers are transitive (Proposition \ref{prop:transfer transitive})
and commute with inflations (Proposition \ref{prop:transfer and epi}).
Theorem \ref{thm:double coset formula} below establishes the `double coset formula'
that expresses the composite of a transfer followed by a restriction 
as a linear combination of restrictions followed by transfers. 

\medskip

We let $H$ be a closed subgroup of a compact Lie group $G$. 
We write $i^*:G\bT_*\to H\bT_*$ for the restriction functor
from based $G$-spaces to based $H$-spaces.
We write
\[ G\ltimes_H - \ = \ (G_+)\sm_H - \ : \ H\bT_*\ \to \ G\bT_* \]
for the left adjoint induction functor.
For a based $H$-space $A$ and a based $G$-space $B$,
the {\em shearing isomorphism}\index{subject}{shearing isomorphism}
is the $G$-equivariant homeomorphism
\[ B\sm (G\ltimes_H A) \ \iso \ G\ltimes_H(i^*B\sm A) \ , \quad
b\sm [g, a]\ \longmapsto \ [g, (g^{-1} b)\sm a]\ . \]

\begin{construction}\label{con:Wirthmuller map}
As before we consider a closed subgroup $H$ of a compact Lie group $G$.
To define the Wirthm{\"u}ller morphism, 
we need a specific $H$-equivariant map
\[ l_A \ : \ G\ltimes_H A \ \to \ A\sm S^L \]
that is natural for continuous based $H$-maps in $A$, see \eqref{eq:define_l_A} below.
Here $L=T_{e H}(G/H)$ is the tangent $H$-representation,
i.e., the tangent space at the preferred coset $e H$ of the smooth
manifold $G/H$. Since the left translation by $H$ on $G/H$ fixes $e H$,
the tangent space $L$ inherits an $H$-action, and $S^L$ 
is its one-point compactification.
When $H$ has finite index in $G$, then $L$ is trivial
and the map $l_A:G\ltimes_H A \to A$
is simply the projection onto the wedge summand indexed by the preferred $H$-coset $e H$.
If the dimension of $G$ is bigger than the dimension of $H$,
then $L$ is non-zero and the definition of the map $l_A$ is substantially more involved.

We consider the left $H^2$-action on $G$ given by 
\[ H^2\times G \ \to \ G \ , \quad (h',h)\cdot g \ = \ h' g h^{-1} \ .\]
There is another action of the group $\Sigma_2$ on $G$ by inversion,
i.e., the non-identity element of $\Sigma_2$ acts by $\tau\cdot g= g^{-1}$.
These two actions combine into an action of the wreath product group
$\Sigma_2\wr H=\Sigma_2\ltimes H^2$ on $G$.
The subgroup $H$ is the $(\Sigma_2\ltimes H^2)$-orbit of $1\in G$, whose stabilizer
is the subgroup $\Sigma_2\times\Delta$, where $\Delta= \{ (h,h) \ | \ h\in H\}$
is the diagonal subgroup of $H^2$.
The differential of the projection $G\to G/H$ identifies
\[ \nu \ = \ ( T_1 G) / ( T_1 H) \ ,\]
the normal space at $1\in H$ of the inclusion $H\to G$, 
with the tangent representation $L =T_{e H} (G/H)$.
The representation $\nu$ is a representation 
of the stabilizer group $\Sigma_2\times\Delta$ of $1$, 
and if we identify $H$ with $\Delta$ via $h\mapsto(h,h)$, 
then this identification takes the $\Delta$-action on $\nu$  
to the tangent $H$-action on $L$. 
Moreover, the differential at~1 of the inversion map $g\mapsto g^{-1}$
is multiplication by $-1$ on the tangent space; 
so the above identification of $\nu$ with $L$ takes the $\Sigma_2$-action
to the sign action on $L$.

We choose an $H$-invariant inner product on the vector space $L$.
Since $H$ is the $H^2$-orbit of~1 inside $G$, there is a slice,
i.e., a $\Sigma_2\times\Delta$-equivariant smooth embedding
\[ s\ : \ D( L ) \ \to \ G\]
of the unit disc of $L$ with $s(0)=1$ and such that 
the differential at $0\in D(L)$ of the composite
\[ D(L)\ \xrightarrow{\ s\ }\ G \ \xrightarrow{\text{proj}}\ G/H \]
is the identity of $L$.
The property that $s$ is $(\Sigma_2\times \Delta)$-equivariant means in more concrete terms
that the relations
\[ s(h\cdot l)\ = \ h\cdot s(l)\cdot h^{-1} \text{\qquad and\qquad}
 s(- l)\ = \  s(l)^{-1} \]
hold in $G$ for all $(h,l)\in H\times D(L)$.
After scaling the slice, if necessary, the map
\[  D(L)\times H \ \to \ G \ , \quad (l,h)\ \longmapsto \ s(l)\cdot h \]
is an embedding whose image is a tubular neighborhood of $H$ in $G$.
For a proof, see for example \cite[II Thm.\,5.4]{bredon-intro}.
This embedding is equivariant for the action of $H^2$ on the source by
\[ (h_1,h_2)\cdot (l, h)\ = \ (h_1 l,\, h_1  h h_2^{-1} ) \]
and for the action of $\Sigma_2$ on the source by
\[ \tau\cdot (l, h)\ = \ (-h^{-1} l,\, h^{-1} )  \ .\]
The map 
\[ l_H^G \ : \ G \ \to \ S^L \sm H_+\]
is then defined as the $H^2$-equivariant collapse map 
with respect to the above tubular neighborhood.\index{subject}{Thom-Pontryagin construction}
So explicitly,
\[ l_H^G(g) \ = \ 
\begin{cases}
 ( l / (1-|l|) ) \sm h & \text{ if $g=s(l)\cdot h$ with $(l,h)\in D(L)\times H$, and }\\
\qquad \ast & \text{ if $g$ is not of this form.}
\end{cases} \]
Given any based $H$-space $A$, we can now form $l\sm_H A$,
where $-\sm_H -$ refers to the action of the second $H$-factor in $H^2$.
We obtain an $H$-equivariant based map
\begin{equation} \label{eq:define_l_A}
 l_A \ = \ l_H^G\sm_H A \ : \ G\ltimes_H A\ \to \ 
(S^L\sm H_+)\sm_H A \ \iso \ A\sm S^L
\end{equation}
that is natural in $A$. 
Here $H$ acts by left translation on the source, and diagonally on the target.
This map is thus given by
\[ l_A[g,a] \ = \ 
\begin{cases}
 h a\sm ( l / (1-|l|) )  & \text{ if $g=s(l)\cdot h$ with $(l,h)\in D(L)\times H$, and }\\
\quad \ast & \text{ if $g$ is not of this form.}
\end{cases} \]
\end{construction}

If $H$ has finite index in $G$, then $L=0$ and the triangle 
of the following proposition commutes on the nose, by direct inspection.
So the essential content of the next result is when the dimension of $G$
exceeds the dimension of $H$.  

\begin{prop}\label{prop:shear and l}
Let $H$ be a closed subgroup of a compact Lie group $G$, 
$A$ a based $H$-space and $B$ a based $G$-space.
Then the following triangle commutes up to $H$-equivariant based homotopy:
\[ \xymatrix@C=18mm{ B\sm (G\ltimes_H A) \ar[d]^\iso_{\text{\em shear}} \ar[dr]^-{B\sm l_A} &\\
 G\ltimes_H(i^*B\sm A) \ar[r]_-{l_{i^* B\sm A}} & B\sm A\sm S^L } \]
\end{prop}
\begin{proof}
We write down an explicit homotopy: we define the map
\[ K \ : \ ( B\sm (G\ltimes_H A) )\times [0,1]\ \to \ B\sm A\sm S^L  \]
by
\[ K(b\sm [g,a], t) \ = \ 
\begin{cases}
  s(t l)^{-1} b\sm h a\sm \frac{l}{1-|l|} & \text{ if $g=s(l)\cdot h$ for $(l,h)\in D(L)\times H$,}\\
\quad \ast & \text{ if $g$ is not of this form.}
\end{cases}
\]
Then for all $(l,h)\in D(L)\times H$ with $g=s(l)\cdot h$ we have
\[ K(b\sm [g, a],0)\ = \ 
   b\sm h a\sm (l/(1-|l|)) \ = \ (B\sm l_A)(b\sm [g, a])  \]
(because $s(0)=1$), and
\[ K(b\sm [g, a],1)\ = \ 
  h g^{-1} b\sm h a\sm (l/(1-|l|) \ = \ l_{i^*B\sm A}[g, g^{-1}b \sm a])  \]
(because $s(l)^{-1}=h g^{-1}$). So $K$ is the desired $H$-equivariant homotopy.
\end{proof}

The restriction functor from orthogonal $G$-spectra to orthogonal $H$-spectra 
has a left and a right adjoint, and both are given by
applying the space level adjoints $G\ltimes_H -$ respectively $\map^H(G,-)$ levelwise.
We will mostly be concerned with the left adjoint, and so we spell out
the construction in more detail.

\begin{construction}[Induced spectrum]
We let $H$ be a closed subgroup of a compact Lie group $G$ 
and $Y$ an orthogonal $H$-spectrum. 
We denote by $G\ltimes_H Y$\index{subject}{induced spectrum}\index{symbol}{$G\ltimes_H Y$ - {induced spectrum}} 
the {\em induced $G$-spectrum} whose value at an inner product space $V$ is
$G\ltimes_H Y(V)$.
When $G$-acts on $V$ by linear isometries, 
then this value has the diagonal $G$-action, 
through the action on the external $G$ and by functoriality in $V$.
With this diagonal $G$-action, $(G\ltimes_H Y)(V)$ 
is equivariantly homeomorphic to $G\ltimes_H Y(i^* V)$ 
(where $H$ acts diagonally on $Y(i^* V)$), via
\begin{equation}  \label{eq:level_of_induced}
 G\ltimes_H Y(i^* V)\ \iso \ (G\ltimes_H Y)(V)\ , \quad
[g,y] \ \longmapsto \ [g, Y(l_g)(y)]\ .  
\end{equation}
Under this isomorphism the structure map of the spectrum $G\ltimes_H Y$
becomes the combination of the shearing isomorphism
and the structure map of the $H$-spectrum $Y$, 
i.e., for every $G$-representation $W$ the following square commutes:
\[ 
\xymatrix@C13mm{  S^V \sm ( G\ltimes_H Y(i^* W)) 
\ar[r]^-{\eqref{eq:level_of_induced}}_-\iso\ar[d]_{\text{shear}}^\iso &
S^V\sm (G\ltimes_H Y)(W) \ar[dd]^{\sigma_{V,W}^{G\ltimes_H Y}}\\
G\ltimes_H ( S^{i^* V}\sm Y(i^* W)) \ar[d]_{G\ltimes_H (\sigma^Y_{V,W})} & \\
G\ltimes_H ( Y(i^* V\oplus i^* W )) \ar[r]_-{\eqref{eq:level_of_induced}}^-\iso & 
(G\ltimes_H Y)(V\oplus W) } \]
\end{construction}

If the group $H$ has finite index in $G$, then 
the tangent representation $L=T_{e H}(G/H)$ is trivial and
$l_A:G\ltimes_H A \to A$ is the projection onto the wedge summand indexed by $e H$;
the maps $l_{Y(V)}:G\ltimes_H Y(V)\to Y(V)$ then
form a morphism of orthogonal $H$-spectra $l_Y:G\ltimes_H Y\to Y$.
This morphism induces a map of $H$-equivariant homotopy groups.
In general, however, the diagram
\[ \xymatrix@C=18mm{ 
S^V\sm (G\ltimes_H Y(W))\ar[r]^-{S^V\sm l_{Y(W)}}\ar@<-8ex>@/_2pc/[dd]_(.7){\sigma_{V,W}^{G\ltimes_H Y}}
\ar[d]^{\text{shear}}_{\iso}& 
S^V\sm Y(W)\sm S^L\ar[dd]^{\sigma_{V,W}^{Y\sm S^L}}\\
G\ltimes_H (S^V\sm Y(W))\ar[ur]_{l_{S^V\sm Y(W)}}\ar[d]^{G\ltimes_H \sigma_{V,W}^Y}& \\
G\ltimes_H Y(V\oplus W)\ar[r]_-{l_{Y(V\oplus W)}}& Y(V\oplus W)\sm S^L} \]
does {\em not} commute on the nose (because the upper triangle does not commute);
so if the dimension of $G$ is bigger than the dimension of $H$ we do {\em not}
obtain a morphism of orthogonal $H$-spectra in the strict sense.
Still, Proposition \ref{prop:shear and l} shows 
that the above diagram does commute up to based $H$-equivariant homotopy,
and this is good enough to yield a well-defined homomorphism 
\[ (l_Y)_* \ : \ \pi^H_k(G\ltimes_H Y) \ \to \ \pi_k^H(Y\sm S^L) \ . \]
As we just explained, this is an abuse of notation, since $(l_Y)_*$
is in general {\em not} induced by a morphism of orthogonal $H$-spectra.

We can now consider the composite  
\begin{equation}  \label{eq:define_Wirthmuller}
 \Wirth_H^G \ : \ \pi^G_k(G\ltimes_H Y) \ \xra{\ \res^G_H\ } \
\pi^H_k(G\ltimes_H Y) \ \xra{\ (l_Y)_*\ } \ \pi_k^H(Y\sm S^L) \ ,  
\end{equation}
which we call the {\em Wirthm{\"u}ller map}.\index{subject}{Wirthm{\"u}ller map}
We will show in Theorem \ref{thm:Wirth iso} below that the Wirthm{\"u}ller map is
an isomorphism.

\begin{construction}[Transfers]\label{con:external transfer}
We let $H$ be a closed subgroup of a compact Lie group $G$,
and we let $Y$ be an orthogonal $H$-spectrum.
The {\em external transfer}\index{subject}{transfer!external} 
\begin{equation}  \label{eq:define external transfer}\index{symbol}{$G\ltimes_H -$ - {external transfer}}
 G\ltimes_H - \ : \ \pi_0^H(Y\sm S^L) \ \to \ \pi_0^G(G\ltimes_H Y) 
\end{equation}
is a map in the direction opposite to the Wirthm{\"u}ller map.
The construction involves another equivariant Thom-Pontryagin construction.\index{subject}{Thom-Pontryagin construction}
We choose a $G$-representation $V$ and a vector $v_0\in V$ 
whose stabilizer group is $H$;
this is possible, for example by \cite[Ch.\,0, Thm.\,5.2]{bredon-intro},
\cite[III Thm.\,4.6]{broecker-tomDieck}
or \cite[Prop.\,1.4.1]{palais-classification}.
This data determines a $G$-equivariant smooth embedding
\[ i\ : \ G/H \ \to \ V \ , \quad g H \ \longmapsto \ g v_0 \]
whose image is the orbit $G v_0$. We let
\[ W \ = \ V - (d i)_{e H}( L ) \ = \ V - T_{v_0}( G\cdot v_0) \]
denote the orthogonal complement of the image of the tangent space at $e H$;
this is an $H$-subrepresentation of $V$, and
\begin{equation}  \label{eq:tangent_plus_normal}
L\oplus W \ \iso \ V\ ,\quad (x,w)\ \longmapsto \ (d i)_{e H}(x) \ + \ w
\end{equation}
is an isomorphism of $H$-representations.

By multiplying the original vector $v_0$ with a sufficiently large scalar,
if necessary, we can assume that the embedding $i$ is `wide', 
in the sense that the exponential map
\[ 
j \ : \ G\times_H D(W) \ \to \ V \ , \quad 
[g , w] \ \longmapsto \ g\cdot ( v_0 \ +\ w )
\]
is an embedding, where $D(W)$ is the unit disc of the normal $H$-representation,
compare \cite[Ch.\,II, Cor.\,5.2]{bredon-intro}.
So $j$ defines an equivariant tubular neighborhood of
the orbit $G v_0$ inside $V$.
The associated collapse map 
\begin{equation}\label{eq:TPcollaps}
c \ : \  S^V \ \to  \ G\ltimes_H S^W
\end{equation}
then becomes the $G$-map defined by 
\[ c(v) \ = \
\begin{cases}
  \left[g, \frac{w}{1-|w|} \right] &\text{if $v= g\cdot(v_0 +w)$ for some $(g,w)\in G\times D(W)$, and}\\
\quad \infty & \text{else.}
\end{cases}\]
With the collapse map in place, we can now define 
the external transfer \eqref{eq:define external transfer}.
We let $U$ be an $H$-representation and $f:S^U \to Y(U)\sm S^L$ an $H$-equivariant 
based map that represents a class in $\pi_0^H(Y\sm S^L)$.
By enlarging $U$, if necessary, we can assume that it
is underlying a $G$-representation, 
see for example \cite[Prop.\,1.4.2]{palais-classification}
or \cite[III Thm.\,4.5]{broecker-tomDieck}.
We stabilize $f$ by $W$ from the right to obtain the $H$-map
\begin{align}  \label{eq:define_f_diamond_W}
 f\diamond W \ : \ 
 S^U\sm S^W \ \xra{f\sm S^W}\  &Y(U)\sm S^L\sm S^W\\ 
& _\eqref{eq:tangent_plus_normal} \iso \ 
 Y(U)\sm S^V\ \xra{\sigma^{\op}_{U,V}} \  Y(U\oplus V) \ .\nonumber
\end{align}
The composite $G$-map
\begin{align*}
  S^{U\oplus V} \ \xra{\ S^U\sm c\ }\  
S^U\sm (G\ltimes_H S^W) \ &\xra[\iso]{\text{shear}}\  G\ltimes_H(S^U\sm  S^W) \\ 
\xra{G\ltimes_H( f\diamond W)} \ 
&G\ltimes_H ( Y(U\oplus V))\ \iso_\eqref{eq:level_of_induced} \  (G\ltimes_H Y)(U\oplus V) 
\end{align*}
then represents the external transfer
\[ G\ltimes_H\td{f} 
\text{\qquad in \quad} \pi_0^G(G\ltimes_H Y)\ . \]
\end{construction}

The next proposition provides the main geometric input for
the Wirthm{\"u}ller isomorphism.

\begin{prop}\label{prop:collaps map interaction}
Let $H$ be a closed subgroup of a compact Lie group $G$.
\begin{enumerate}[\em (i)]
\item 
The composite
\[ S^V\ \xra{\ c\ }\ G\ltimes_H S^W \ \xra{l_H^G\sm_H S^W}\  S^L\sm S^W  \]
is $H$-equivariantly homotopic to the map induced by
the inverse of the isometry \eqref{eq:tangent_plus_normal}.
\item
Let $A$ be a based $H$-space and $f,f':B\to G\ltimes_H A$
two continuous based $G$-maps. If the composites
$l_A\circ f, l_A\circ f':B\to A\sm S^L$ are $H$-equivariantly homotopic,
then the maps $f\sm S^V,f'\sm S^V:B\sm S^V\to (G\ltimes_H A)\sm S^V$ 
are $G$-equivariantly homotopic.
\end{enumerate}
\end{prop}
\begin{proof}
(i)
The composite $(l_H^G\sm_H S^W)\circ c$
is the collapse map based on the composite closed embedding
\[ \zeta \ : \  D(L)\times D(W) \ 
\xra{(l,w)\mapsto [s(l), w]}\  G\times_H D(W)\ \xra{[g,w]\mapsto g\cdot(v_0+w)} \ V \ .\]
Then $\zeta(0,0)=v_0$ and we let 
\[ D\ = \ (d \zeta)_{(0,0)}\ : \ L\times W\ \to\ V \]
denote the differential of $\zeta$ at $(0,0)$. We observe that
$ D(0,w) =  w$ because $\zeta(0,w)=v_0+w$; 
on the other hand, $\zeta(l,0)=s(l)\cdot v_0$,
so the restriction of $D$ to $L$ is the differential at~0 of the composite
\[ D(L)\ \xra{\ s\ }\ G \ \xra{\text{proj}}\ G/H \ \xra{\ i \ } \ V \ . \]
Since the differential of the composite $\text{proj}\circ s$
is the identity, we obtain 
\[ D(l,0) \ = \ (d i)_{e H}(l) \ .\]
Since the differential is additive we conclude that
$D(l,w)=(d i)_{e H}(l)+w$, i.e., $D$ equals the isomorphism \eqref{eq:tangent_plus_normal}.

We consider the $H$-equivariant homotopy
\begin{align*}
  K \ : \ D(L)\times D(W)\times [0,1]\ &\to \quad V \\
K(l,w,t)\quad &=\
\begin{cases}
  \frac{\zeta(t\cdot (l,w))- v_0}{t} \ + \ t\cdot v_0 & \text{ for $t>0$, and}\\
  \quad  D(l,w) &  \text{ for $t=0$.}
\end{cases} 
\end{align*}
That this assignment is continuous (in fact smooth) when $t$ approaches~0
is the defining property of the differential.  
Moreover, for every $t\in [0,1]$ the map $K(-,-,t):D(L)\times D(W)\to V$
is a smooth equivariant embedding, 
so it gives rise to a collapse map $c_t:S^V\to S^L\sm S^W$ defined by
\[ c_t(v) \ = \
\begin{cases}
  \left( \frac{l}{1-|l|},\frac{w}{1-|w|} \right) &\text{if $v= K(l,w,t)$ for some $(l,w)\in D(L)\times D(W)$, and}\\
\qquad \infty & \text{else.}
\end{cases}
\]
The passage from the embedding to the collapse map is continuous,
so the 1-parameter family $c_t$ provides an $H$-equivariant
based homotopy from the collapse map $c_0$
to the collapse map $c_1$
associated to the embedding $\zeta=K(-,-,1):D(L)\times D(W)\to V$,
i.e., to the map $(l_H^G\sm_H S^W)\circ c$.
Since $D=(d \zeta)_{(0,0)}$ is the isometry \eqref{eq:tangent_plus_normal},
another scaling homotopy compares the collapse map $c_0$ to the one-point
compactification of the inverse of \eqref{eq:tangent_plus_normal}.

(ii) We define a based continuous $G$-map
\[ r \ : \ \map^H(G,A\sm S^L)\sm S^V\ \to\
 G\ltimes_H( A\sm S^L\sm S^W) \]
as the composite
\begin{align*}
  \map^H(G,A\sm S^L)\sm S^V\ &\xra{\ \Id\sm c\ }\ 
\map^H(G,A\sm S^L)\sm (G\ltimes_H S^W)\\ 
&\xra{\ \text{shear}\ }\ 
 G\ltimes_H\left( \map^H(G,A\sm S^L)\sm S^W\right)\\ 
&\xra{G\ltimes_H(\epsilon\sm S^W)}\ 
 G\ltimes_H( A\sm S^L\sm S^W) \ .
\end{align*}
Here $\epsilon:\map^H(G,A\sm S^L)\to A\sm S^L$ 
is the adjunction counit, i.e., evaluation at $1\in G$.
We expand this definition. We let $\psi:G\to A\sm S^L$ be an $H$-equivariant
based map and $v\in S^V$. If $v$ is not in the image of the tubular
neighborhood $j:G\times_H D(W)\to V$, then $r(\psi\sm v)$ is the basepoint.
Otherwise, $v=j[g,w]=g\cdot(v_0+w)$ for some $(g,w)\in G\times D(W)$, and then
\begin{align}\label{expand r}
r(\psi\sm v) \ 
&= \ ((G\ltimes_H (\epsilon\sm S^W))\circ\text{shear})( \psi\sm [g,w/(1-|w|)]) \\
&= \ [g, \epsilon (g^{-1}\cdot  \psi)\sm  w/(1-|w|)] \ 
= \ \left[g, \psi(g^{-1})\sm \frac{w}{1-|w|}\right] \ .\nonumber
\end{align}
We denote by $l_A^\sharp: A\to\map^H(G,A\sm S^L)$
the $H$-map defined by $l_A^\sharp(a)(g)=l_A[g,a]$,
and we let $\varphi:V\to L\oplus W$ be the inverse of the isometry \eqref{eq:tangent_plus_normal}.
Now we argue that the following square commutes up to $H$-equivariant based
homotopy:
\begin{equation}  \begin{aligned}\label{eq:Wirth-square}
 \xymatrix@C=10mm{ A\sm S^V \ar[r]^-{[1,-]} 
\ar[dd]_{l_A^\sharp\sm S^V} & 
 G\ltimes_H(A\sm S^V)\ar[d]_\iso^{G\ltimes_H S^\varphi} \\
  &G\ltimes_H (A\sm S^L\sm S^W) \ar[d]_{\iso}^{G\ltimes_H(S^{-\Id_L}\sm S^W)}\\
\map^H(G,A\sm S^L)\sm S^V\ar[r]_-r &G\ltimes_H (A\sm S^L\sm S^W) }  
  \end{aligned}\end{equation}
To see this we define an $H$-homotopy
\[ K\ : \  (A\sm S^V) \times [0,1]\ \to \ G\ltimes_H(A\sm S^L\sm S^W) \]
as follows. We exploit that the map
\[ \zeta \ : \  D(L)\times D(W) \ \to \ V \ ,\quad \zeta(l,w)\ = \ s(l)\cdot(v_0 + w )  \]
is a smooth embedding.
This map already featured in the proof of part~(i),
because the collapse map based on $\zeta$ is the composite $(l_H^G\sm_H S^W)\circ c$.
If a vector is of the form $v=\zeta(l,w)=s(l)\cdot (v_0+w)$
for some $(l,w)\in D(L)\times D(W)$ (necessarily unique), then we set
\[ K( a\sm v,t) \ =\ 
 K( a\sm (s(l)\cdot(v_0 +w)),\ t) \ =\ 
\left[ s(t\cdot l),\, a\sm \frac{-l}{1-|l|}\sm \frac{w}{1-|w|} \right] \ .\]
For $|l|=1$ or $|w|=1$ this formula yields the basepoint,
so we can extend this definition by sending all elements $a\sm v$
to the basepoint whenever $v$ is not in the image of the embedding $\zeta$.
We claim that for $t=0$ we obtain $K(-,1)=r\circ(l_A^\sharp\sm S^V)$.
Indeed, if $v=s(l)\cdot(v_0+w)$ for $(l,w)\in D(L)\times D(W)$
(necessarily unique), then 
\begin{align*}
  ( r\circ(l_A^\sharp\sm S^V))(a\sm v)\quad 
_\eqref{expand r}  &= \
  \left[s(l), l_A^\sharp(a)(s(l)^{-1})\sm \frac{w}{1-|w|}\right] \\
&= \ 
 \left[s(l),l_A[s(-l),a]\sm \frac{w}{1-|w|}\right]
\\
&= \ 
 \left[s(l), a\sm \frac{-l}{1-|l|}\sm \frac{w}{1-|w|}\right]
\ = \ K(a\sm v,1)\ .
\end{align*}
On the other hand, the map $K(-,0)$ agrees with the composite
\begin{align*}
 A\sm S^V \ &\xra{A\sm ((l_H^G\sm_H S^W)\circ c)}\  
A\sm S^L\sm S^W \\ 
&\xra{\ [1,-]\ }\ G\ltimes_H( A\sm S^L\sm S^W)  
\ \xra{G\ltimes_H(S^{-\Id_L}\sm S^W)}\ G\ltimes_H( A\sm S^L\sm S^W)  \ ;  
\end{align*}
part~(i) provides an $H$-homotopy between $(l_H^G\sm_H S^W)\circ c$ and
the homeomorphism $S^\varphi$; so this proves the claim.

We let $\Lambda:G\ltimes_H(A\sm S^V)\to \map^H(G,A\sm S^L)\sm S^V$
be the $G$-equivariant extension of the $H$-map $l_A^\sharp\sm S^V$.
Since the square \eqref{eq:Wirth-square} commutes up to $H$-equivariant homotopy,
the composite
\[ G\ltimes_H(A\sm S^V)\ \xra{\ \Lambda \ }\ \map^H(G,A\sm S^L)\sm S^V\ \xra{\ r\ }\ 
G\ltimes_H(A\sm S^L\sm S^W)  \]
is $G$-equivariantly homotopic to the $G$-homeomorphism 
$G\ltimes_H ((S^{-\Id_L}\sm S^W)\circ S^\varphi)$, by adjointness.
So the composite of $r \Lambda$ with the shearing homeomorphism
$ (G\ltimes_H A)\sm S^V\iso  G\ltimes_H(A\sm S^V)$
is also $G$-homotopic to a homeomorphism.
This composite equals
the composite
\[ ( G\ltimes_H A)\sm S^V\ \xra{\Psi_A\sm S^V}\ \map^H(G,A\sm S^L)\sm S^V\ \xra{\ r\ }\ 
G\ltimes_H(A\sm S^L\sm S^W)  \ ,\]
where $\Psi_A: G\ltimes_H A\to\map^H(G,A\sm S^L)$
is the adjoint of the $H$-map $l_A:G\ltimes_H A\to A\sm S^L$.

Now we consider based $G$-maps $f,f':B\to G\ltimes_H A$ such that
$l_A\circ f$ and $l_A\circ f':B\to A\sm S^L$ are $H$-equivariantly homotopic.
Then the two composites
\[ B\ \xra{\ f, f' \ }\ G\ltimes_H A\ \xra{\ \Psi_A} \ \map^H(G, A\sm S^L) \]
are $G$-equivariantly homotopic, by adjointness.
The composite $(\Psi_A\circ f)\sm S^V=(\Psi_A\sm S^V)\circ(f\sm S^V)$
is then $G$-equivariantly homotopic 
to $(\Psi_A\circ f')\sm S^V=(\Psi_A\sm S^V)\circ(f'\sm S^V)$.
The map $\Psi_A\sm S^V$ has a $G$-equivariant retraction,
up to $G$-homotopy, by the previous paragraph.
So already the maps $f\sm S^V$ and $f'\sm S^V$ are $G$-homotopic.
\end{proof}

Now we can establish the Wirthm{\"u}ller isomorphism,
which first appeared in \cite[Thm.\,2.1]{wirthmuller}.
Wirthm{\"u}ller attributes parts of the ideas to tom\,Dieck,
and his statement that $G$-spectra define a `complete $G$-homology theory'
amounts to Theorem \ref{thm:Wirth iso}.
Our proof is essentially Wirthm{\"u}ller's original argument,
adapted to the context of equivariant orthogonal spectra.
We recall that
\[ \varepsilon_L\ : \ \pi_0^H(Y\sm S^L) \ \to \ \pi_0^H(Y\sm S^L) \]
denotes the effect of the involution of $Y\sm S^L$
induced by the linear isometry $-\Id_L:L\to L$ given by multiplication by $-1$.

\begin{theorem}[Wirthm{\"u}ller isomorphism]\index{subject}{Wirthm{\"u}ller isomorphism}\label{thm:Wirth iso}
Let $H$ be a closed subgroup of a compact Lie group $G$
and $Y$ an orthogonal $H$-spectrum.
Let $L=T_{e H}(G/H)$ denote the tangent $H$-representation.
Then the maps
\[ \Wirth_H^G \ : \  \pi^G_0(G\ltimes_H Y) \ \to \  \pi^H_0(Y\sm S^L) \]
and
\[ (G\ltimes_H -)\circ \varepsilon_L \ : \ 
\pi^H_0(Y\sm S^L) \ \to \ \pi^G_0(G\ltimes_H Y)  \]
are independent of the choices made in their definitions, 
and they are mutually inverse isomorphisms.
\end{theorem}
\begin{proof}
We prove the various claims in a specific order.
In a first step we show that the Wirthm{\"u}ller map
is left inverse to the map $(G\ltimes_H -)\circ\varepsilon_L$,
independently of all the choices made in the definitions.
We let $U$ be a $G$-representation and $f:S^U \to Y(U)\sm S^L$ an $H$-equivariant 
based map that represents a class in $\pi_0^H(Y\sm S^L)$.
We also choose a wide $G$-equivariant embedding
$i: G/H \to V$ as in Construction \ref{con:external transfer}.
This provides a decomposition $L\oplus W\iso V$ 
of $H$-representations as in \eqref{eq:tangent_plus_normal}
and a $G$-equivariant Thom-Pontryagin collapse map $c:S^V\to G\ltimes_H S^W$.
We let $\varphi:V\iso L\oplus W$ denote the inverse 
of the linear isomorphism \eqref{eq:tangent_plus_normal}.
We contemplate the diagram of based $H$-maps:
\[ \xymatrix@C=20mm{
S^U\sm S^V \ar@/^1pc/[dr]^(.6){S^U\sm  S^\varphi}\ar@<-6ex>@/_3pc/[ddd]
\ar[d]_{S^U\sm c}   \\
S^U\sm (G\ltimes_H S^W)  \ar[r]_-{S^U\sm( l_H^G\sm_H S^W)}\ar[d]_{\text{shear}}^\iso 
& S^U\sm S^L\sm S^W \ar[d]^{S^U\sm\tau_{L,W}} \\
 G\ltimes_H ( S^U\sm S^W) \ar[d]^{G\ltimes_H(f\diamond W)}\ar[r]_-{l_{S^U\sm S^W}} &
 S^U\sm S^W\sm S^L \ar[d]^{(f\diamond W)\sm S^L}\\
  G\ltimes_H Y(U\oplus V )\ar[r]_-{l_{Y(U\oplus V)}} & Y(U\oplus V)\sm S^L
} \]
The left vertical composite represents the external transfer $G\ltimes_H\td{f}$,
so the composite around the lower left corner represents
$(l_Y)_*(\res^G_H(G\ltimes_H\td{f}))$.
The upper triangle commutes up to $H$-equivariant homotopy by
Proposition \ref{prop:collaps map interaction} (i).
The middle square commutes up to $H$-homotopy by
Proposition \ref{prop:shear and l} and the fact that $l_A:G\ltimes_H A\to A\sm S^L$
is defined as the composite of $(l_H^G\sm_H A)$ and the twist isomorphism
$S^L\sm A\iso A\sm S^L$.
The lower square commutes by naturality of the maps \eqref{eq:define_l_A}.
Upon expanding the definition \eqref{eq:define_f_diamond_W} 
of $f\diamond W$, the diagram shows that 
the class $(l_Y)_*(\res^G_H(G\ltimes_H\td{f}))$ is also represented by the composite
\begin{align*}
 S^{U\oplus V} \ &\xra{ f\sm S^{\varphi} }\ Y(U)\sm S^L\sm S^L\sm S^W \
\xra{ Y(U)\sm S^L\sm \tau_{L,W}}\ Y(U)\sm S^L\sm S^W\sm S^L \\
& \xra{ Y(U)\sm S^{\varphi^{-1}}\sm S^L}\ Y(U)\sm S^V\sm S^L \
\xra{ \sigma^{\op}_{U,V}\sm S^L}\ Y(U\oplus V) \sm S^L \ . 
  \end{align*}
The isometry 
\[ (\varphi^{-1}\oplus L)\circ(L\oplus\tau_{L,W})\circ(L\oplus \varphi) \ : \ 
L\oplus V\ \to \ V\oplus L \]
is {\em not} the twist isometry $\tau_{L,V}$;
however $(\varphi^{-1}\oplus L)\circ(L\oplus\tau_{L,W})\circ(L\oplus \varphi)$
is equivariantly homotopic to the composite
\[ L\oplus V\ \xra{\tau_{L,V}} \ V\oplus L \ \xra{V\oplus(-\Id_L)}\ V\oplus L\ . \]
Hence the class $(l_Y)_*(\res^G_H(G\ltimes_H\td{f}))$ is also represented by the composite
\[ S^{U\oplus V} \ \xra{ f\sm S^V}\ Y(U)\sm S^L\sm S^V \
 \xra{ Y(U)\sm \tau_{L,V}}\ Y(U)\sm S^V\sm S^L \
\xra{ \sigma^{\op}_{U,V}\sm S^{-\Id}}\ Y(U\oplus V) \sm S^L \ . \]
So altogether this shows the desired relation
\[ \Wirth_H^G(G\ltimes_H\td{f}) \ = \ 
(l_Y)_*(\res^G_H(G\ltimes_H\td{f}))\ = \ \varepsilon_L\td{f}\ . \]
In particular, the Wirthm{\"u}ller map is surjective.

Now we show that the Wirthm{\"u}ller map is injective.
We let $f,f':S^U\to G\ltimes_H Y(U)$ be two $G$-maps that represent
classes with the same image under the Wirthm{\"u}ller map
$\Wirth_H^G :\pi^G_0(G\ltimes_H Y) \to \pi^H_0(Y\sm S^L)$.
By increasing the $G$-representation $U$, if necessary, we can assume
that the composites
\[ S^U\ \xra{\ f, f' \ }\ G\ltimes_H Y(i^* U)\
\xra{l_{Y(i^* U)}} \ Y(i^* U)\sm S^L \]
are $H$-equivariantly homotopic.
Then by Proposition \ref{prop:collaps map interaction}~(ii) the two maps
\[ S^{U\oplus V}\ \xra{\ f\sm S^V, f'\sm S^V \ }\ ( G\ltimes_H Y(i^* U))\sm S^V 
\ \iso \ ( G\ltimes_H Y)(U)\sm S^V \]
are $G$-equivariantly homotopic.
The maps remain $G$-homotopic if we furthermore postcompose
with the opposite structure 
map $\sigma^{\op}_{U,V}:(G\ltimes_H Y)(U)\sm S^V\to(G\ltimes_H Y)(U\oplus V)$.
This shows that the stabilizations from the right of $f$ and $f'$ by $V$
become $G$-homotopic. 
Since such a stabilization represents the same class 
in $\pi_0^G(G\ltimes_H Y)$ as the original map, this shows that $\td{f}=\td{f'}$,
i.e., the Wirthm{\"u}ller map is injective.

Now we know that the Wirthm{\"u}ller map and the map $(G\ltimes_H -)\circ\varepsilon_L$
are inverse to each other, no matter which choices of slice $s:D(L)\to G$, 
representation $V$ and $G$-embedding $G/H\to V$ we made.
Since the choices for the Wirthm{\"u}ller map and the choices for
the external transfer are independent of each other, the two resulting maps
are independent of all choices.

The last thing to show is that the Wirthm{\"u}ller map and the external
transfer are additive maps. If $H$ has finite index in $G$, the Wirthm{\"u}ller
map is the composite of a restriction homomorphism and the effect of
a morphism of orthogonal $H$-spectra, both of which are additive.
In general, however, we need an additional argument, namely naturality.
Indeed, for a fixed choice of slice $s:D(L)\to G$, the Wirthm{\"u}ller
map $\pi^G_0(G\ltimes_H Y) \to  \pi^H_0(Y\sm S^L)$ is natural in the orthogonal
$H$-spectrum $Y$.
Since source and target are reduced additive functors
from orthogonal $H$-spectra to abelian groups, any natural transformation
is automatically additive, by Proposition \ref{prop:additivity prop}.
Since the Wirthm{\"u}ller map is additive, so is its inverse, and
hence also the external transfer.
\end{proof}

Theorem \ref{thm:Wirth iso} states the Wirthm{\"u}ller isomorphism
only for 0-dimensional homotopy groups.
We now extend it to homotopy groups in all integer dimensions.
This extension is a rather formal consequence of the fact that 
the Wirthm{\"u}ller maps
commute with the loop\index{subject}{loop isomorphism} 
and suspension isomorphisms\index{subject}{suspension isomorphism}
\[  \alpha\ : \  \pi^G_k(\Omega X)\ \to\ \pi^G_{k+1} (X)  \text{\quad and\quad}
 -\sm S^1 \ : \ \pi^G_k (X)\ \to \ \pi^G_{k+1}(X\sm S^1) \]
defined in \eqref{eq:loop iso} respectively \eqref{eq:suspension iso}.

\begin{prop}\label{prop:Wirth in dimension k}
Let $H$ be a closed subgroup of a compact Lie group $G$ and $Y$
an orthogonal $H$-spectrum. 
\begin{enumerate}[\em (i)]
\item 
The following diagrams commute for
all integers $k$:
\[ \xymatrix@C=15mm@R=8mm{  
\pi_k^G(G\ltimes_H (\Omega Y))  \ar[d]_{\text{\em assembly}}^\iso 
 \ar[r]^-{\Wirth_H^G} &
\pi_k^H( (\Omega Y)\sm S^L)\ar[d]_\iso^{\text{\em assembly}} 
\\
\pi_k^G(\Omega(G\ltimes_H Y))\ar[d]^\iso_\alpha &
\pi_k^H(\Omega(Y\sm S^L)) \ar[d]^\alpha_\iso  
\\
\pi_{k+1}^G(G\ltimes_H Y) \ar[r]_-{\Wirth_H^G}   &
\pi_{k+1}^H(Y\sm S^L)} \]
\[
 \xymatrix@C=15mm@R=8mm{  
\pi_k^G(G\ltimes_H Y)  \ar[d]_{-\sm S^1}^\iso \ar[r]^-{\Wirth_H^G} &
\pi_k^H(Y\sm S^L) \ar[d]_\iso^{-\sm S^1} 
\\
\pi_{k+1}^G((G\ltimes_H Y)\sm S^1)\ar[d]_{b_*}^\iso&
\pi_{k+1}^H(Y\sm S^L\sm S^1) \ar[d]^{(Y\sm \tau_{S^L,S^1})_*}_-\iso  
\\
\pi_{k+1}^G(G\ltimes_H (Y\sm S^1)) \ar[r]_-{\Wirth_H^G}  &
\pi_{k+1}^H(Y\sm S^1\sm S^L)
 }\]
Here $\tau_{S^L,S^1}:S^L\sm S^1\to S^1\sm S^L$ is the twist homeomorphism,
and the $G$-equivariant isomorphism
\[ b \ : \ (G\ltimes_H Y)\sm S^1\ \iso \ G\ltimes_H (Y\sm S^1) \]
is given by $b([g,y]\sm t)=[g, y\sm t]$.
\item The Wirthm{\"u}ller map
\[ \Wirth_H^G \ : \ \pi_k^G(G\ltimes_H Y) \ \to \ \pi_k^H(Y\sm S^L) \]
is an isomorphism for all $k\in\mZ$.
\end{enumerate}
\end{prop}
\begin{proof}
(i) The loop and suspension isomorphisms commute with restriction from $G$ to $H$, 
by direct inspection. So the proof comes down
to checking that the following diagrams commute:
\[ \xymatrix@C=15mm@R=8mm{  
\pi_k^H(G\ltimes_H (\Omega Y))  \ar[d]_{\text{assembly}}^\iso 
 \ar[r]^-{(l_{\Omega Y})_*}  &
\pi_k^H( (\Omega Y)\sm S^L)\ar[d]_\iso^{\text{assembly}} 
\\
\pi_k^H(\Omega(G\ltimes_H Y))\ar[d]^\iso_\alpha &
\pi_k^H(\Omega(Y\sm S^L)) \ar[d]^\alpha_\iso  
\\
\pi_{k+1}^H(G\ltimes_H Y) \ar[r]_-{(l_Y)_*}   &
\pi_{k+1}^H(Y\sm S^L) 
\\
\pi_k^H(G\ltimes_H Y)  \ar[d]_{-\sm S^1}^\iso \ar[r]^-{(l_Y)_*} &
\pi_k^H(Y\sm S^L) \ar[d]_\iso^{-\sm S^1} 
\\
\pi_{k+1}^H((G\ltimes_H Y)\sm S^1)\ar[d]_{b_*}^\iso&
\pi_{k+1}^H(Y\sm S^L\sm S^1) \ar[d]^{(Y\sm \tau_{S^L,S^1})_*}_-\iso  
\\
\pi_{k+1}^H(G\ltimes_H (Y\sm S^1)) \ar[r]_-{(l_{Y\sm S^1})_*}  &
\pi_{k+1}^H(Y\sm S^1\sm S^L)
 }\]
This in turn follows from the fact -- again verified by direct inspection -- 
that for every based $H$-space $A$ the following two squares commute:
\[ \xymatrix{ 
G\ltimes_H(\Omega A) \ar[r]^-{l_{\Omega A}}\ar[d]_{\text{assembly}} & 
(\Omega A)\sm S^L\ar[d]^{\text{assembly}} 
&
(G\ltimes_H A)\sm S^1 \ar[r]^-{l_A\sm S^1}\ar[d]_b & 
A\sm S^L\sm S^1\ar[d]^{A\sm\tau_{S^L.S^1}}\\
\Omega( G\ltimes_H A) \ar[r]_-{\Omega l_A} &
\Omega(A\sm S^L)
&
G\ltimes_H (A\sm S^1) \ar[r]_-{l_{A\sm S^1}} &
A\sm S^1\sm S^L} \]

(ii) We argue by induction over $|k|$, the absolute value of the integer $k$.
The induction starts with $k=0$, where Theorem \ref{thm:Wirth iso} 
provides the desired conclusion.
If $k$ is positive, the compatibility of the Wirthm{\"u}ller map with the loop isomorphism,
established in part~(i), provides the inductive step.
If $k$ is negative, 
the compatibility of the Wirthm{\"u}ller map with the suspension isomorphism,
also established in part~(i), provides the inductive step.
\end{proof}

\begin{construction}[External transfer in integer degrees]
So far we only defined the external transfer for 0-dimensional homotopy groups.
Theorem \ref{thm:Wirth iso} shows that in dimension~0
the external transfer is inverse to the map
$\varepsilon_L\circ \Wirth_H^G:\pi^G_0(G\ltimes_H Y) \to \pi^H_0(Y\sm S^L)$.
We want the same property in all dimensions, and since the
Wirthm{\"u}ller map is an isomorphism by Proposition \ref{prop:Wirth in dimension k}
we simply define the external transfer isomorphism
\begin{equation}  \label{eq:define_external_dim_k}
 G\ltimes_H - \ : \   \pi^H_k(Y\sm S^L) \ \to \ 
\pi^G_k(G\ltimes_H Y)  
\end{equation}
as the composite
\[  \pi^H_k(Y\sm S^L) 
\ \xra{\ \varepsilon_L\ } \  \pi^H_k(Y\sm S^L) \ 
\ \xra{(\Wirth_H^G)^{-1}} \  
\pi^G_k(G\ltimes_H Y) \ .\]
The compatibility of the  Wirthm{\"u}ller isomorphism with the 
loop\index{subject}{loop isomorphism} and suspension isomorphisms\index{subject}{suspension isomorphism}
then directly implies the analogous compatibility for the
external transfer $G\ltimes_H -$, by reading the diagrams of
Proposition \ref{prop:Wirth in dimension k}~(i) backwards.
For easier reference, we explicitly record this compatibility
in Proposition \ref{prop:internal transfer loop suspension} below.
\end{construction}

In the next proposition we use the Wirthm{\"u}ller isomorphism
to show that smashing with a cofibrant $G$-space is homotopical.
In \cite[III Thm.\,3.11]{mandell-may}, Mandell and May give
a different proof of this fact which does not use 
the Wirthm{\"u}ller isomorphism.

\begin{prop}\label{prop:space smash preserves global}
Let $G$ be a compact Lie group and $A$ a cofibrant based $G$-space. 
\begin{enumerate}[\em (i)]
\item Let $f:X\to Y$ be a morphism of orthogonal $G$-spectra with the
following property: for every closed subgroup $H$ of $G$ that 
fixes some non-basepoint of $A$,
the map $\pi_*^H(f):\pi_*^H(X)\to\pi_*^H(Y)$ is an isomorphism.
Then the morphism $f\sm A:X\sm A\to Y\sm A$ is a
$\upi_*$-isomorphism of orthogonal $G$-spectra.
\item
The functor $-\sm A$ preserves $\upi_*$-isomorphisms of orthogonal $G$-spectra.
\end{enumerate}
\end{prop}
\begin{proof}
(i) Smashing with $A$ commutes with mapping cones, so by the long exact
homotopy group sequence of Proposition \ref{prop:LES for homotopy of cone and fibre}
it suffices to show the following special case. Let $X$ be an orthogonal $G$-spectrum
with the following property: for every closed subgroup $H$ of $G$ that occurs
as the stabilizer group of a non-basepoint of $A$,
the groups $\pi_*^H(X)$ vanish.
Then for all closed subgroups $K$ of $G$ 
the equivariant homotopy groups $\pi_*^K(X\sm A)$ vanish.

In a first step we show this when $A$ is a finite-dimensional $G$-CW-complex.
We argue by contradiction. If the claim were false,
we could find a compact Lie group $G$ of minimal dimension for which it fails.
We let $A$ be a $G$-CW-complex whose dimension $n$ is minimal among all counterexamples.
Then $A$ can be obtained from an $(n-1)$-dimensional subcomplex $B$
by attaching equivariant cells $G/H_i\times D^n$, 
for $i$ in some indexing set $I$, where $H_i$ is a closed subgroup of $G$.
Then $X\sm (A/B)$ is isomorphic to a wedge, over the set $I$,
of orthogonal $G$-spectra $X\sm (G/H_i)_+\sm S^n$.
Since equivariant homotopy groups take wedges to sums, the 
suspension isomorphism and the Wirthm{\"u}ller
isomorphism allow us to rewrite the equivariant homotopy groups 
of $X\sm A/B$ as
\begin{align*}
  \pi_*^G(X\sm A/B)\ &\iso \ {\bigoplus}_{i\in I}\  \pi_*^G(X\sm (G/H_i)_+\sm S^n)\\
&\iso {\bigoplus}_{i\in I} \ \pi_{*-n}^G(X\sm (G/H_i)_+)\
\iso {\bigoplus}_{i\in I} \ \pi_{*-n}^{H_i}(X\sm S^{L_i})\ , 
\end{align*}
where $L_i$ is the tangent representation of $H_i$ in $G$.
Since $H_i$ is the stabilizer of a non-basepoint of $A$,
the groups $\pi_*^K(X)$ vanish for all closed subgroups $K$ of $H_i$, by hypothesis.
If $H_i$ has finite index in $G$, then $L_i=0$ and the 
respective summand thus vanishes.
The representation sphere $S^{L_i}$ admits a finite $H_i$-CW-structure,
so if $H_i$ has strictly smaller dimension than $G$, then 
the respective summand vanishes by the minimality of $G$.
So altogether we conclude that the groups 
$\pi_*^G(X\sm A/B)$ vanish. The groups $\pi_*^G(X\sm B)$
vanish by the minimality of $A$. The inclusion of $B$ into $A$ is
an h-cofibration of based $G$-spaces, so the induced morphism
$X\sm B\to X\sm A$ is an h-cofibration of orthogonal $G$-spectra.
Hence the groups $\pi_*^G(X\sm A)$ vanish by the long exact sequence of
Corollary \ref{cor-long exact sequence h-cofibration}~(i).
Since $A$ was supposed to be a counterexample to the proposition,
we have reached the desired contradiction.
Altogether this proves the claim when $A$ admits the structure of a finite-dimensional
$G$-CW-complex.

If $A$ admits the structure of a $G$-CW-complex, possibly infinite
dimensional, we choose a skeleton filtration by $G$-subspaces $A_n$.
Then the $G$-homotopy groups of $X\sm A_n$ vanish for all $n\geq 0$,
and all the morphisms $X\sm A_n\to X\sm A_{n+1}$ are h-cofibrations
of orthogonal $G$-spectra. Since equivariant
homotopy groups commute with such sequential colimits
(see Proposition \ref{prop:sequential colimit closed embeddings} (i)), 
the groups $\pi_*^G(X\sm A)$ vanish as well.
An arbitrary cofibrant based $G$-space is $G$-homotopy equivalent
to a based $G$-CW-complex, so 
the groups $\pi_*^G(X\sm A)$ vanish for all cofibrant $A$.

Now we let $K$ be an arbitrary closed subgroup of $G$. 
The underlying $K$-space of $A$ is again cofibrant
by Proposition \ref{prop:cofibrancy preservers}~(i),
so we can apply the previous reasoning to $K$ instead of $G$ and
conclude that the groups $\pi_*^K(X\sm A)$ vanish.
\end{proof}

\begin{prop}\label{prop:induction is homotopical}
Let $G$ and $K$ be compact Lie groups and $A$ a cofibrant based $(G\times K)$-space
such that the $G$-action is free away from the basepoint
and the $K$-action is free away from the basepoint.
Then the functor
\[ A\sm_K - \ : \ K\spec \ \to \ G\spec \]
takes $\upi_*$-isomorphisms of orthogonal $K$-spectra to
$\upi_*$-isomorphisms of orthogonal $G$-spectra.
\end{prop}
\begin{proof}
The functor $A\sm_K-$ preserves mapping cones, so by the long exact
homotopy group sequence of Proposition \ref{prop:LES for homotopy of cone and fibre}
it suffices to show the following special case. Let $C$ be an orthogonal $K$-spectrum
all of whose equivariant homotopy groups vanish, for all closed subgroups of $K$.
Then all equivariant homotopy groups of the orthogonal $G$-spectrum $A\sm_K C$
vanish, for all closed subgroups of $G$.

In a first step we show that the $G$-equivariant stable homotopy groups
of $A\sm_K C$ vanish.
A cofibrant based $(G\times K)$-space is equivariantly homotopy
equivalent to a based $(G\times K)$-CW-complex, so it is no loss of generality
to assume an equivariant CW-structure with skeleton filtration
\[ \ast\ = \ A_{-1} \ \subset \ A_0 \ \subset \ \dots \ \subset \ A_n \
\subset \ \dots \ .\]
We show first, by induction on $n$, that the orthogonal $G$-spectrum
$(A_n)\sm_K C$ has trivial $G$-equivariant homotopy groups.
The induction starts with $n=-1$, where there is nothing to show. 
For $n\geq 0$ the quotient $A_n/A_{n-1}$ is $(G\times K)$-equivariantly isomorphic
to a wedge of summands of the form $((G\times K)/\Gamma_i)_+\sm S^n$,
for certain closed subgroups $\Gamma_i$ of $G\times K$.
Since the $K$-action on $A$ is free (away from the basepoint), 
each isotropy group $\Gamma_i$ that occurs 
is the graph of a continuous homomorphism $\alpha_i:H_i\to K$
defined on a closed subgroup $H_i$ of $G$.
Since the $G$-action on $A$ is free (away from the basepoint),
all the homomorphisms $\alpha_i$ that occur are injective.

Since equivariant homotopy groups take wedges to sums, the 
suspension isomorphism and the Wirthm{\"u}ller
isomorphism allow us to rewrite the equivariant homotopy groups 
of $(A_n/A_{n-1})\sm_K C$ as
\begin{align*}
  \pi_*^G((A_n/A_{n-1})\sm_K C)\ &\iso \ \bigoplus\,  \pi_*^G(((G\times K)/\Gamma_i)_+\sm_K C\sm S^n)\\
&\iso \ \bigoplus \, \pi_{*-n}^G( G\ltimes_{H_i}\alpha_i^*(C))\
\iso \bigoplus \, \pi_{*-n}^{H_i}(\alpha_i^*(C)\sm S^{L_i})\ , 
\end{align*}
where $L_i$ is the tangent representation of $H_i$ in $G$.
By hypothesis on $C$ and because $\alpha_i$ is injective, the orthogonal $H$-spectrum
$\alpha_i^*(C)$ is $H_i$-stably contractible. So $\alpha_i^*(C)\sm S^{L_i}$ 
is $H_i$-stably contractible 
by Proposition \ref{prop:space smash preserves global}~(ii). 
Altogether this shows that the orthogonal $G$-spectrum $(A_n/A_{n-1})\sm_K C$ 
has vanishing $G$-equivariant homotopy groups.

The inclusion $A_{n-1}\to A_n$ is an h-cofibration of based $(G\times K)$-spaces,
so the induced morphism $( A_{n-1})\sm_K C\to (A_n)\sm_K C$ is an h-cofibration
of orthogonal $G$-spectra, giving rise to a long exact sequence of equivariant
homotopy groups (Corollary \ref{cor-long exact sequence h-cofibration}).
By the previous paragraph and the inductive hypothesis, the
orthogonal $G$-spectrum $(A_n)\sm_K C$ 
has vanishing $G$-equivariant homotopy groups.
This completes the inductive step.

Since $A$ is the sequential colimit, along h-cofibrations
of based $(G\times K)$-spaces, of the skeleta $A_n$,
the orthogonal $G$-spectrum $A\sm_K C$ is the sequential colimit, 
along h-cofibrations of orthogonal $G$-spectra, 
of the sequence with terms $(A_n)\sm_K C$. 
Equivariant homotopy groups commute with such sequential colimits 
(Proposition \ref{prop:sequential colimit closed embeddings}~(i)), 
so also $A\sm_K C$ has vanishing $G$-equivariant homotopy groups.

Now we let $H$ be any closed subgroup of $G$. 
The underlying $(H\times K)$-space of $A$ is again cofibrant
by Proposition \ref{prop:cofibrancy preservers}~(i),
so we can apply the previous reasoning to $H$ instead of $G$ and
conclude that the groups $\pi_*^H(A\sm_K C)$ vanish.
\end{proof}

\begin{cor}\label{cor:induce H2G}
Let $H$ be a closed subgroup of a compact Lie group $G$. Then the induction functor  
\[ G\ltimes_H - \ : \ H\spec\ \to \ G\spec \]
takes $\upi_*$-isomorphisms of orthogonal $H$-spectra to
$\upi_*$-isomorphisms of orthogonal $G$-spectra.
\end{cor}
\begin{proof}
We let $G\times H$ act on $G$ by left and right translation, i.e., via
\[ (g,h)\cdot\gamma\ = \ g \gamma h^{-1}\ . \]
With this action, $G$ is $(G\times H)$-equivariantly isomorphic to
the homogeneous space $(G\times H)/\Delta$
for $\Delta=\{(h,h)\ : \ h\in H\}$.
In particular, $G$ is $(G\times H)$-cofibrant, and both partial actions
are free. So Proposition \ref{prop:induction is homotopical} 
applies to the functor $(G_+)\sm_H- = G\ltimes_H-$ and yields the desired conclusion.
\end{proof}

Now we discuss the transfer maps of equivariant homotopy groups.

\begin{construction}[Transfers]\label{con:transfer}\index{subject}{transfer!on equivariant homotopy groups}
We let $H$ be a closed subgroup of a compact Lie group $G$.
As before we write $L=T_{e H}(G/H)$ for the tangent space of $G/H$ at the coset $e H$,
which inherits an $H$-action from the $H$-action on $G/H$.
We let $X$ be an orthogonal $G$-spectrum. We form the composite
\[ \pi_k^H(X\sm S^L) \ \xra[\iso]{\ G\ltimes_H -\ } \ \pi_k^G(G\ltimes_H X) 
\ \xra{\ \text{act}_*\ } \ \pi_k^G(X) \]
of the external transfer \eqref{eq:define_external_dim_k} and the effect of
the action map (i.e., the adjunction counit) $G\ltimes_H X\to X$. We call this composite
the {\em dimension shifting transfer}\index{subject}{transfer!dimension shifting} 
and denote it\index{symbol}{$\Tr_H^G$ - {dimension shifting transfer}}
\begin{equation}\label{eq:dimension_shifting_transfer}
 \Tr_H^G\ : \ \pi_k^H (X\sm S^L) \ \to \ \pi_k^G (X) \ .
\end{equation}
The (degree zero) transfer is then defined as the composite\index{subject}{transfer!degree zero},
\begin{equation}\label{eq:transfer}
 \tr_H^G \ : \ \pi_k^H(X)\ \xra{\ (X\sm i)_* \ } \  \pi_k^H(X\sm S^L) \ 
\xra{\ \Tr_H^G\ } \ \pi_k^G(X)\ , 
\end{equation}   
where $i:S^0\to S^L$ is the `inclusion of the origin', the based map sending~0 to~0.
Both kinds of transfer are natural for morphisms of 
orthogonal $G$-spectra.\index{symbol}{$\tr_H^G$ - {degree zero transfer}}
For finite index inclusions, $L=0$ and there is no difference between
the dimension shifting transfer and the degree zero transfer.
The  external transfer is additive;
since the dimension shifting and degree zero transfers are obtained from
there by applying morphisms of equivariant spectra, they are also additive.

The key properties of these transfer maps are:
\begin{itemize}
\item transfers are transitive (Proposition \ref{prop:transfer transitive});
\item transfers commute with inflation maps (Proposition \ref{prop:transfer and epi});
\item the restriction of a degree zero transfer to a closed subgroup satisfies
a double coset formula (Theorem \ref{thm:double coset formula}).
\end{itemize}
\end{construction}

\begin{eg}[Infinite Weyl group transfers]\label{eg:L(H) has fix}  
If $H$ has infinite index in its normalizer, then the 
degree zero transfer $\tr_H^G$ is trivial. Indeed, the inclusion of the normalizer
$N_G H$ of $H$ into $G$ induces a smooth embedding
\[ W_G H\ =\ (N_G H) / H \ \to \ G/H  \]
and thus a monomorphism of tangent $H$-representations
\[ T_{e H}(W_G H) \ \to \ T_{e H}( G/H ) = L \ . \]
If $n\in N_G H$ is an element of the normalizer and $h\in H$, then
\[  h\cdot n H \ = \ n\cdot (n^{-1}h n) H \ = \ n H\ ,\]
so $W_G H$ is $H$-fixed inside $G/H$. Consequently, the tangent space
$T_{e H} (W_G H)$ is contained in the $H$-fixed space $L^H$.
If $H$ has infinite index in its normalizer, then the Weyl group $W_G H$ 
and the tangent space $ T_{e H}(W_G H)$ have positive dimension.
In particular, $L$ has non-zero $H$-fixed points.
The point~0 in $S^L$ can thus be moved through $H$-fixed points 
to the basepoint at infinity.
The first map in the composite \eqref{eq:transfer} is thus the zero map,
hence so is the transfer $\tr_H^G$.
\end{eg}

\begin{rk} As the previous example indicates,
the passage from the dimension shifting transfer $\Tr_H^G$ to the
degree zero transfer $\tr_H^G$ loses information.
An extreme case is when the subgroup $H$ is normal in $G$. 
Then the action of the group $H$ on $G/H$ is trivial;
hence also the $H$-action on the tangent space $L$ is trivial.
Upon choosing an isomorphism $L\iso \mR^d$,
the Wirthm{\"u}ller isomorphism identifies
$\pi_0^G(G\ltimes_H Y)$ with $\pi_0^H(Y\sm S^d)$,
where $d=\dim(G/H)=\dim(G)-\dim(H)$.
The dimension shifting transfer then becomes a natural map
\[ \pi_0^H ( Y ) \ \iso \ \pi_d^H(Y\sm S^L) \ \xra{\Tr_H^G} \ \pi_d^G (Y) \ . \]
This transformation is generically non-trivial.
\end{rk}

Now we recall some important properties of the transfer maps.
We start with the compatibility with the loop and suspension isomorphisms.

\begin{prop}\label{prop:internal transfer loop suspension}
Let $H$ be a closed subgroup of a compact Lie group $G$,
\begin{enumerate}[\em (i)]
\item 
For every orthogonal $H$-spectrum $Y$ and all $k\in\mZ$, 
the following diagrams commute:
\[ \xymatrix@C=15mm@R=7mm{  
\pi_k^H( (\Omega Y)\sm S^L) \ar[r]^-{G\ltimes_H -}_-{\iso} 
\ar[d]^\iso_{\text{\em assembly}} &\pi_k^G(G\ltimes_H (\Omega Y))  \ar[d]^{\text{\em assembly}}_\iso 
\\
\pi_k^H(\Omega(Y\sm S^L)) \ar[d]_\alpha^\iso  &
\pi_k^G(\Omega(G\ltimes_H Y))\ar[d]_\iso^\alpha 
\\
\pi_{k+1}^H(Y\sm S^L) \ar[r]_-{G\ltimes_H-}^-{\iso}  &
\pi_{k+1}^G(G\ltimes_H Y)  
}\]
\[ \xymatrix@C=15mm@R=7mm{  
\pi_k^H(Y\sm S^L) \ar[r]^-{G\ltimes_H -}_-{\iso} 
\ar[d]^\iso_{-\sm S^1} &\pi_k^G(G\ltimes_H Y)  \ar[d]^{-\sm S^1}_\iso \\
\pi_{k+1}^H(Y\sm S^L\sm S^1) \ar[d]_{(Y\sm \tau_{S^L,S^1})_*}^-\iso  &
\pi_{k+1}^G((G\ltimes_H Y)\sm S^1)\ar[d]_\iso^{b_*}\\
\pi_{k+1}^H(Y\sm S^1\sm S^L) \ar[r]_-{G\ltimes_H-}^-{\iso}  &
\pi_{k+1}^G(G\ltimes_H (Y\sm S^1))  
}\]
Here $\tau_{S^L,S^1}:S^L\sm S^1\to S^1\sm S^L$ is the twist isomorphism,
and the $G$-equivariant isomorphism
\[ b \ : \ (G\ltimes_H Y)\sm S^1\ \iso \ G\ltimes_H (Y\sm S^1) \]
is given by $b([g,y]\sm t)=[g, y\sm t]$.
\item
For every orthogonal $G$-spectrum $X$ and all $k\in\mZ$, 
the following diagrams commute:
\[ \xymatrix@C=7mm@R=7mm{  
\pi_k^H( (\Omega X)\sm S^L) \ar[r]^-{\Tr_H^G} \ar[d]^\iso_{\text{\em assembly}}
 &\pi_k^G(\Omega X) \ar[dd]_\iso^\alpha &
\pi_k^H(X\sm S^L) \ar[r]^-{\Tr_H^G}
\ar[d]^\iso_{-\sm S^1} &\pi_k^G(X)  \ar[dd]^{-\sm S^1}_\iso \\
\pi_k^H(\Omega(X\sm S^L))  \ar[d]^\iso_\alpha &
&
\pi_{k+1}^H(X\sm S^L\sm S^1) \ar[d]_{(X\sm \tau_{S^L,S^1})_*}^-\iso  &
\\
\pi_{k+1}^H(X\sm S^L) \ar[r]_-{\Tr_H^G} &
\pi_{k+1}^G(X)  
&
\pi_{k+1}^H(X\sm S^1\sm S^L) \ar[r]_-{\Tr_H^G}  &
\pi_{k+1}^G(X\sm S^1)  }\]
\[ \xymatrix@C=10mm{  
\pi_k^H( \Omega X) \ar[r]^-{\tr_H^G} \ar[d]^\iso_\alpha 
 &\pi_k^G(\Omega X) \ar[d]_\iso^\alpha &
\pi_k^H(X) \ar[r]^-{\tr_H^G}
\ar[d]^\iso_{-\sm S^1} &\pi_k^G(X)  \ar[d]^{-\sm S^1}_\iso \\
\pi_{k+1}^H(X) \ar[r]_-{\tr_H^G} &
\pi_{k+1}^G(X)  
&
\pi_{k+1}^H(X\sm S^1) \ar[r]_-{\tr_H^G}  &
\pi_{k+1}^G(X\sm S^1)  }\]
\end{enumerate}
\end{prop}
\begin{proof}
The first two diagrams commute because we can read the diagrams of
Proposition \ref{prop:Wirth in dimension k}~(i) backwards.
The commutativity of the other diagrams then follows by naturality.
\end{proof}

Now we establish the transitivity property for a 
nested triple of compact Lie groups $K\leq H\leq G$.
We continue to denote by $L=T_{e H}(G/H)$ the tangent $H$-representation in $G$,
and we write $\bar L=T_{e K}(H/K)$ for the tangent $K$-representation in $H$.
We choose a slice 
\[ s\ : \ D( L ) \ \to \ G\]
as in the construction of the map $l_H^G:G\to S^L\sm H_+$ in 
Construction \ref{con:Wirthmuller map}; so $s$ is a wide smooth embedding
of the unit disc of $L$ satisfying
\[ s(0)\ =\ 1 \text{\qquad and\qquad}   s(h\cdot l)\ = \ h\cdot s(l)\cdot h^{-1} \]
for all $(h,l)\in H\times D(L)$, and 
the differential at $0\in D(L)$ of the composite
$\text{proj}_H\circ s:D(L)\to  G/H$ is the identity of $L$.
The differential at $0\in D(L)$ of the composite
\[ D(L)\ \xrightarrow{\ s\ }\ G \ \xrightarrow{\text{proj}_K}\ G/K \]
is then a $K$-equivariant linear monomorphism
\[ d(\text{proj}_K\circ s)_0 \ : \ L\ \to \ T_{e K}(G/K) \ = \ L(K,G) \]
that splits the differential at $e K$ of the projection $q:G/K\to G/H$.
So the combined map
\begin{align}  \label{eq:split_transitive_tangent}
   ( d(\text{proj}_K\circ s)_0,\, (d q )_{e H}) \ : \
L\oplus\bar L  \ &= \ T_{e H}(G/H)\oplus T_{e K}(H/K) \\ 
&\to \ T_{e K}(G/K) \ = \ L(K,G) \nonumber
\end{align}
is an isomorphism of $K$-representations.
Upon one-point compactification this isomorphism induces a homeomorphism
of $K$-spaces
\[
 S^L\sm S^{\bar L}\ \iso \ S^{L(K,G)}\ .   
\]
Any two slices are $H$-equivariantly isotopic 
(compare \cite[VI Thm.\,2.6]{bredon-intro}),
so the $K$-equivariant homotopy class of the latter isomorphism is independent
of the choice of slice.

\begin{prop}[Transitivity of transfers]\label{prop:transfer transitive}
Let $G$ be a compact Lie group, $K\leq H\leq G$ nested closed subgroups
and $X$ an orthogonal $G$-spectrum.
Then the composite
\[ \pi_k^K(X\sm S^{L(K,G)}) \ \iso_{\eqref{eq:split_transitive_tangent}} \ 
\pi_k^K(X\sm S^L\sm S^{\bar L}) \ \xra{\ \Tr_K^H\ } \ 
\pi_k^H(X\sm S^L) \ \xra{\ \Tr_H^G\ } \ \pi_k^G(X) \]
agrees with the transfer $\Tr_K^G$.
Moreover, the degree zero transfers satisfy
\[ \tr_H^G\circ\tr_K^H \ = \ \tr_K^G \ : \ \pi_k^K (X) \ \to \ \pi_k^G (X) \ .\]
\end{prop}
\begin{proof}
We start by establishing transitivity of the Wirthm{\"u}ller maps.
We choose a slice for the inclusion of $K$ into $H$, 
i.e., a wide smooth embedding $\bar s : D( \bar L ) \to  H$ satisfying
\[ \bar s(0)\ =\ 1 \text{\qquad and\qquad}   
\bar s(k\cdot \bar l)\ = \ k\cdot \bar s(\bar l)\cdot k^{-1} \]
for all $(k,\bar l)\in K\times D(\bar L)$, and 
such that the differential at $0\in D(\bar L)$ lifts the identity of $\bar L$.
We combine the two slices into a slice for $K$ inside $G$:
the $K$-equivariant map
\[ D( L \oplus \bar L)\ \to \ G\ ,\quad
(l,\bar l)\ \longmapsto\  s(l)\cdot \bar s(\bar l)\]
sends $(0,0)$ to $1$, and its differential at $(0,0)$ is exactly the
identification \eqref{eq:split_transitive_tangent}.
So we can -- and will -- define the map $l_K^G:G\to S^{L(K,G)}\sm K_+$
from the slice
\[ s' \ : \ D(L(K,G)) \ \xra[\eqref{eq:split_transitive_tangent}^{-1}]{\iso} \
D( L \oplus \bar L)\ \xra{(l,\bar l) \mapsto s(l)\cdot \bar s(\bar l)} \ G\ .\]
The maps $l_H^G:G\to S^L\sm H_+$ respectively
$l_K^H:H\to S^{\bar L}\sm K_+$ are the Thom-Pontryagin collapses 
based on the $H^2$-equivariant smooth embedding
\[  \tilde s \ : \  D(L)\times H \ \to \ G \ , \quad (l,h)\ \longmapsto \ s(l)\cdot h \]
respectively 
the $K^2$-equivariant smooth embedding
\[  \hat s \ : \  D(\bar L)\times K \ \to \ H \ , \quad 
(\bar l,k)\ \longmapsto \ \bar s(\bar l)\cdot k\ . \]
The composite $(S^L\sm l_K^H)\circ l_H^G$ is thus $K^2$-equivariantly homotopic
to the Thom-Pontryagin collapse 
based on the $K^2$-equivariant smooth embedding
\[    D(L)\times D(\bar L)\times K \ \to \ G \ , 
\quad (l,\bar l,k)\ \longmapsto \ s(l)\cdot\bar s(\bar l)\cdot k \ .\]
The following diagram of $K^2$-equivariant smooth embeddings
then commutes by construction:
\[ \xymatrix@C=12mm@R=7mm{
D(L\oplus \bar L)\times  K\ar[d]_{\text{incl}}
\ar[r]_-\iso^-{\eqref{eq:split_transitive_tangent}} 
 & D(L(K,G))\times K \ar[dd]^{(l,k)\mapsto s'(l)\cdot k}  \\
D(L)\times D(\bar L)\times  K
\ar[d]_{(l,\bar l, k)\mapsto (l,\bar s(\bar l)\cdot k)} & \\
D(L)\times  H  \ar[r]_-{(l,h)\mapsto s(l)\cdot h} & G  
} \]
So the associated diagram of $K^2$-equivariant collapse maps also commutes:
\[ \xymatrix@C=12mm@R=7mm{
G_+ \ar[dd]_{l_K^G} \ar[r]^-{l_H^G} &
S^L\sm H_+\ar[d]^{S^L\sm l_K^H} \\
& 
S^L\sm S^{\bar L}\sm K_+\ar[d]^{\Psi} \\
 S^{L(K,G)}\sm K_+ \ar[r]_-\iso^-{\eqref{eq:split_transitive_tangent}^{-1}}  
&  S^{L\oplus \bar L}\sm K_+} \]
Here $\Psi:S^L\sm S^{\bar L}\to S^{L\oplus\bar L}$ is the collapse map
for the inclusion $D(L\oplus \bar L)\to D(L)\times D(\bar L)$.
A rescaling homotopy connects $\Psi$ to the canonical homeomorphism
$S^L\sm S^{\bar L}\iso S^{L\oplus\bar L}$, so the following square
commutes up to $K^2$-equivariant based homotopy:
\[ \xymatrix@C=12mm{
G_+ \ar[d]_{l_K^G} \ar[r]^-{l_H^G} &
S^L\sm H_+\ar[d]^{S^L\sm l_K^H} \\
 S^{L(K,G)}\sm K_+ \ar[r]_-\iso^-{\eqref{eq:split_transitive_tangent}^{-1}}  
&  S^{L\oplus \bar L}\sm K_+} \]
We can thus conclude that the map $l_K^G / K$ 
is $K$-equivariantly homotopic to the composite
\begin{align*}
  G / K_+ \ &\iso \ G\ltimes_H(H /K _+) \ \xra{l_H^G\sm_H (H / K_+)} \
(H /K _+)\sm S^L  \\ 
&\iso \ H\ltimes_K S^L \ \xra{l_K^H\sm_K S^L } \ 
S^L\sm  S^{\bar L} \ \iso_\eqref{eq:split_transitive_tangent} \ S^{L(K,G)} \ .
\end{align*}
Naturality and transitivity of restriction maps then show that the following
diagram commutes:
\[ \xymatrix@C=7mm@R=7mm{ 
\pi_k^G(X\sm G / K_+) \ar[rr]^-{\res^G_K} \ar[d]_\iso&& 
\pi_k^K(X\sm G / K_+) \ar[ddd]^{l_K^G / K} \\
\pi_k^G(X\sm G\ltimes_H (H / K_+)) \ar[d]^-{\res^G_H} \ar@<-8ex>@/_2pc/[dd]_(.3){\Wirth_H^G} && \\
\pi_k^H(X\sm G\ltimes_H(H / K_+)) \ar[d]^{l_H^G\sm_H(H / K_+)} && \\
\pi_k^H(X\sm (H/K_+)\sm S^L) \ar[d]^{\iso} &&
\pi_k^K(X\sm S^{L(K,G)}) \ar[d]^\iso \\
\pi_k^H(X\sm H\ltimes_K S^L)) \ar[r]^-{\res^H_K} 
\ar@/_2pc/[rr]_(.8){\Wirth_K^H} &
\pi_k^K(X\sm H\ltimes_K S^L) \ar[r]^-{l_K^H\sm_K S^L} &
\pi_k^K(X\sm S^L\sm S^{\bar L}) 
} \]
By Theorem \ref{thm:Wirth iso}
the Wirthm{\"u}ller map is inverse to the composite
\[ 
\pi^K_k(X\sm S^{L(K,G)}) \ \xra{\varepsilon_{L(K,G)}} \ \pi^K_k(X\sm S^{L(K,G)}) \ 
\xra{G\ltimes_K -} \ \pi^G_k(G\ltimes_K X)  \ . \]
The map $\varepsilon_{L(K,G)}$ is induced by the antipodal map of $S^{L(K,G)}$.
Under the homeomorphism between $S^{L(K,G)}$ and $S^L\sm S^{\bar L}$,
this becomes the smash product of the antipodal maps of $L$ and $\bar L$.
So reading the diagram backwards gives a commutative diagram of external transfers
\begin{equation}\begin{aligned}\label{eq:transitive external transfer}
 \xymatrix@C=7mm@R=7mm{ 
\pi_k^K(X\sm S^L\sm S^{\bar L}) \ar[r]^-{H\ltimes_K-} \ar[dd]_\iso  &
\pi_k^H(X\sm H\ltimes_K S^L)  \ar[r]^-\iso & 
\pi_k^H(X\sm H/K_+\sm S^L) \ar[d]^{G\ltimes_H-} \\
 && \pi_k^G(X\sm G\ltimes_H H/K_+) \ar[d]^{\iso} \\
\pi_k^K(X\sm S^{L(K,G)}) \ar[rr]_{G\ltimes_K-} && \pi_k^G(X\sm G/K_+) }    
\end{aligned}\end{equation}
Postcomposing with the effect of the projection $G/K\to\ast$
and exploiting naturality gives the claim about the dimension shifting transfer.
The second claim follows by precomposing with the inclusion of the origin of $L(K,G)$.
\end{proof}

\begin{eg}\label{eq:projection versus transfer} 
We let $K\leq H\leq G$ be nested closed subgroups and $X$ an orthogonal $G$-spectrum.
We let $p:G/K\to G/H$ denote the projection.
For later reference we show that under the external transfer
isomorphisms the effect of morphism $X\sm p_+:X\sm G/K_+\to X\sm G/H_+$
corresponds to the transfer from $K$ to $H$. For simplicity
we restrict to the case where $\dim(K)=\dim(H)$, i.e., when
$K$ has finite index in $H$; the general case only differs 
by more complicated notation. If $K$ has finite index in $H$,
then the differential of the projection $G/K\to G/H$ 
is an isomorphism from the $K$-representation $T_{e K}(G/K)$
to the underlying $K$-representation of $L=T_{e H}(G/H)$.
We identify these two representations via this isomorphism.
We claim that then the following square commutes:
\[ \xymatrix@C=15mm{ 
\pi_*^K(X\sm S^L)\ar[r]^-{\tr_K^H}\ar[d]_{G\ltimes_K -}^{\iso} &
\pi_*^H(X\sm S^L)\ar[d]^{G\ltimes_H -}_{\iso}\\
\pi_*^G(X\sm G/K_+)\ar[r]_-{ (X\sm p_+)_*} &
\pi_*^G(X\sm G/H_+)} \]
To see this we compose the commutative 
diagram \eqref{eq:transitive external transfer}
in the proof of Proposition \ref{prop:transfer transitive}
that encodes the transitivity of external transfers
with the map $\pi_*^G(X\sm p_+)$ and arrive at another commutative diagram:
\[ 
 \xymatrix@C=7mm@R=7mm{ 
\pi_k^K(X\sm S^L) \ar[d]^{H\ltimes_K-} \ar@{=}[rr] \ar@<-4ex>@/_3pc/[ddd]_(.3){\tr_K^H}&&
\pi_k^K(X\sm S^L) \ar[dd]^{G\ltimes_K-} \\
\pi_k^H(X\sm H\ltimes_K S^L)  \ar[d]^\iso & \\ 
\pi_k^H(X\sm H/K_+\sm S^L) \ar[r]^-{G\ltimes_H-} \ar[d] &
 \pi_k^G(X\sm G\ltimes_H H/K_+) \ar[r]^-{\iso} & 
 \pi_k^G(X\sm G/K_+) \ar[d]^{(X\sm p_+)_*}\\
\pi_k^H(X\sm  S^L) \ar[rr]_-{G\ltimes_H-} && 
\pi_k^G(X\sm G/H_+)  }    
 \]
The lower left vertical map is induced by the projection $H/K\to\ast$
and the bottom part of the diagram commutes by naturality of the
external transfer.
\end{eg}

Now we prove the compatibility of transfers with inflations, i.e., restriction along
continuous epimorphisms $\alpha:K\to G$. For every closed 
subgroup $H$ of $G$, the map
\[ \bar\alpha\ : \ K/J \ \to \ G/H \ , \quad k J\ \longmapsto \ \alpha(k) H \]
is a diffeomorphism, where $J=\alpha^{-1}(H)$.
The differential at the coset $e J$ is an isomorphism
\[ (d\bar\alpha)_{e J} \ : \ \bar L = T_{e J}(K/J) \ \to \ 
(\alpha|_J)^*(T_{e H}(G/H))\ = \ (\alpha|_J)^*( L ) \]
of $J$-representations.
In the statement and proof of the following proposition 
a couple of unnamed isomorphisms occur. One of them is the natural isomorphism
of $K$-spaces
\[ K\ltimes_J (\alpha|_J)^*(A) \ \iso \ \alpha^*(G\ltimes_H A) \ , \quad
[k, a]\ \longmapsto \ [\alpha(k), a]\ .\]

\begin{prop}\label{prop:transfer and epi}
Let $K$ and $G$ be compact Lie groups and $\alpha:K\to G$ 
a continuous epimorphism. 
Let $H$ be a closed subgroup of $G$,
set $J=\alpha^{-1}(H)$, and let $\alpha|_J:J\to H$ 
denote the restriction of $\alpha$.  
\begin{enumerate}[\em (i)]
\item 
For every orthogonal $H$-spectrum $Y$
the following diagram commutes:
\[ \xymatrix@C=10mm{ 
\pi_k^H(Y\sm S^L) \ar[d]_{G\ltimes_H -} \ar[r]^-{(\alpha|_J)^*} &
\pi_k^J( (\alpha|_J)^*(Y\sm S^L))  \ar[r]^-\iso &
\pi_k^J( (\alpha|_J)^*(Y)\sm S^{\bar L}) \ar[d]^{K\ltimes_J- } \\
\pi_k^G(G\ltimes_H Y) \ar[r]_-{\alpha^*} & 
\pi_k^K(\alpha^*(G\ltimes_H Y)) \ar[r]_-{\iso} & 
\pi_k^K(K\ltimes_J(\alpha|_J)^*( Y) )
} \]
\item 
For every orthogonal $G$-spectrum $X$
the following diagram commutes:
\[ \xymatrix@C=10mm{ 
\pi_k^H(X\sm S^L) \ar[d]_{\Tr_H^G} \ar[r]^-{(\alpha|_J)^*} &
\pi_k^J( (\alpha|_J)^*(X\sm S^L))  \ar[r]^-\iso &
\pi_k^J( (\alpha|_J)^*(X)\sm S^{\bar L}) \ar[d]^{\Tr_J^K} \\
\pi_k^G(X) \ar[rr]_-{\alpha^*} && 
\pi_k^K(\alpha^* (X) ) } \]
Moreover, the degree zero transfers satisfy the relation
\[ \alpha^*\circ \tr_H^G \ = \ \tr_J^K \circ (\alpha|_J)^* \]
as maps $\pi_k^H (X)\to \pi_k^K(\alpha^*(X))$.
\item For every orthogonal $G$-spectrum $X$, every $g\in G$
and all closed subgroups $K\leq H$ of $G$ the following diagram commutes:
\[ \xymatrix{ 
\pi_k^{K^g}(X) \ar[d]_{\tr_{K^g}^{H^g}} \ar[r]^-{g_\star} & \pi_k^K(X)  \ar[d]^{\tr_K^H} \\
\pi_k^{H^g}(X) \ar[r]_-{g_\star} & \pi_k^H(X) } \]
\end{enumerate}
\end{prop}
\begin{proof}
(i) The restriction maps commute with the loop and suspension isomorphisms,
and so do the external transfer maps (by Proposition \ref{prop:internal transfer loop suspension}).
So it suffices to prove the claim in dimension $k=0$.
We choose a wide $G$-equivariant embedding $i:G/H\to V$ into some
$G$-representation; from this input data we form the collapse map
\[ c_{G/H}\ : \ S^V \ \to \ G\ltimes_H S^W  \]
and the external transfer $G\ltimes_H-$. 
Then the composite
\begin{equation}\label{eq:K/J2V}
\alpha^*(i)\circ\bar\alpha\ : \  
K/J \ \to\ \alpha^*(G/H)\ \to \ \alpha^*(V)  
\end{equation}
is a wide $K$-equivariant embedding, and we can (and will) base the 
external transfer $K\ltimes_J-$ on this embedding.
Since $i:G/H\to V$ and $\alpha^*(i)\circ\bar\alpha$ have the same image,
they define the same decomposition of $V$ into tangent and normal subspaces,
i.e.,
\[ \bar W\ = \ \alpha^*(V) - d(\alpha^*(i)\circ\bar\alpha)_{e J}(\bar L)
\ = \ (\alpha|_J)^*(W)\ . \]
Moreover, the composite
\[ S^{\alpha^*(V)}\ = \ 
\alpha^*(S^V)\ \xra{\alpha^*(c_{G/H})}\  \alpha^*(G\ltimes_H S^W)
\ \iso \ K\ltimes_J (\alpha|_J)^*(S^W) \ = \ K\ltimes_J S^{\bar W}\]
is precisely the collapse map based on the wide embedding \eqref{eq:K/J2V}.
From here the commutativity of the square is straightforward from the
definitions.

(ii) The dimension shifting transfer $\Tr_H^G$ is defined as the composite
of the external transfer $G\ltimes_H -$ and the effect of the action map
$G\ltimes_H X \to X$, and similarly for $\Tr_J^K$. 
The action map for the orthogonal $K$-spectrum $K\ltimes_J \alpha^*(X)$
coincides with the composite
\[ K\ltimes_J \alpha^*(X)\ \xra{\ \iso \ }\
 \alpha^*(G\ltimes_ H X)\ \xra{\alpha^*(\text{act})} \ \alpha^*(X)\ . \]
So the following diagram commutes by naturality of the restriction map: 
\[ \xymatrix@C=15mm{ 
\pi_k^G(G\ltimes_H X) \ar[d]_{\pi_k^G(\text{act})} \ar[r]^-{\alpha^*} &
\pi_k^K( \alpha^*(G\ltimes_H X)) \ar[dr]_{\pi_k^K(\alpha^*(\text{act}))} \ar[r]^-\iso &
\pi_k^K(K\ltimes_J \alpha^*(X)) \ar[d]^{\pi_k^K(\text{act})} \\
\pi_k^G(X) \ar[rr]_-{\alpha^*} && 
\pi_k^K(\alpha^* (X) ) } \]
Part~(ii) then follows by stacking this commutative diagram to the
one of part~(i).
The second claim follows by precomposing with the inclusion of the origin of $L$.

Part~(iii) follows from part~(ii) for the epimorphism $c_g:H\to H^g$ and
the closed subgroup $K$ of $H$, and naturality of the transfer:
\begin{align*}
  g_\star\circ\tr_{K^g}^{H^g} \ &= \ (l_g^X)_*\circ (c_g)^*\circ \tr_{K^g}^{H^g}\\
_{\text{(ii)}} \ &= \ (l_g^X)_*\circ \tr_K^H\circ (c_g)^*\
= \ \tr_K^H\circ  (l_g^X)_*\circ (c_g)^*\ =\ \tr_K^H\circ g_\star\qedhere 
\end{align*}
\end{proof}

Now we know how transfers compose and interact with inflations.
The remaining compatibility issue is to rewrite the composite 
of a transfer map followed by a restriction map.
The answer is given by the {\em double coset formula}
that we will prove in Theorem \ref{thm:double coset formula} below.

\index{subject}{geometric fixed points|(}

\section{Geometric fixed points}\label{sec:geometric fixed points}

In this section we study the 
{\em geometric fixed point homotopy groups} $\Phi^G_*(X)$ 
of an orthogonal $G$-spectrum $X$.
We establish the isotropy separation sequence \eqref{eq:isotropy_separation_les} 
that is often useful for inductive arguments, 
and we prove that equivariant equivalences can
also be detected by geometric fixed points, 
see Proposition \ref{prop:fix_and_geometric_fix}.
In Proposition \ref{prop:Phi of transfer} we show that 
geometric fixed points annihilate transfers from proper subgroups.
Theorem \ref{thm:G-equivariant pi_0 of Sigma^infty}  provides
a functorial description of the 0-th equivariant stable homotopy groups 
of a $G$-space $Y$ in terms of the path components of the 
$H$-fixed points spaces $Y^H$ for closed subgroups $H$ of $G$ with finite Weyl group.

\medskip

We define the geometric fixed point homotopy groups of an orthogonal $G$-spectrum $X$.
As before we let $s(\Uc_G)$ denote the set 
of finite-dimensional $G$-subrepresentations 
of the complete $G$-universe $\Uc_G$, considered as a poset under inclusion. 
We obtain a functor from $s(\Uc_G)$ to sets by
\[  V \ \longmapsto\ [S^{V^G}, X(V)^G] \ ,\]
the set of (non-equivariant) homotopy classes of based maps
from the fixed point sphere $S^{V^G}$ to the fixed point space $X(V)^G$.
An inclusion $V\subseteq W$ in $s(\Uc_G)$ is sent to the map
\[ [S^{V^G}, X(V)^G] \ \to \ [S^{W^G},X(W)^G]\]
that takes the homotopy class of $f:S^{V^G}\to X(V)^G$
to the homotopy class of the composite
\begin{align*}
  S^{W^G}\ \iso\ S^{(V^\perp)^G}\sm S^{V^G} \
&\xra{\Id\sm f}\  S^{(V^\perp)^G}\sm X(V)^G  \\ 
&= \  (S^{V^\perp}\sm X(V))^G \ \xra{(\sigma_{V^\perp,V})^G} \
 X(V^\perp\oplus V)^G \ =\ X(W)^G \ .
\end{align*}

 \begin{defn}
Let $G$ be a compact Lie group and $X$ an orthogonal $G$-spectrum.
The  0-th {\em geometric fixed point homotopy group}\index{subject}{geometric fixed points} \index{subject}{fixed points!geometric}  
is defined as \index{symbol}{$\Phi^G_k$ -  geometric fixed point homotopy group}  
\[    \Phi_0^G(X) \ = \ \colim_{V\in s(\Uc_G)}\, [S^{V^G}, X(V)^G]  \ .  \]
If $k$ is an arbitrary integer, we define the $k$-th
geometric fixed point homotopy group $ \Phi_k^G(X)$
by replacing $S^V$ by $S^{V\oplus\mR^k}$ (for $k>0$)
or  replacing $X(V)$ by $X(V\oplus\mR^{-k})$ (for $k<0$),
analogous to the definition of $\pi_k^G(X)$ in \eqref{eq:define_pi_k}.
 \end{defn}

The construction comes with a 
{\em geometric fixed point map}\index{subject}{geometric fixed point homomorphism}
\begin{align}\label{eq:geometric_fix_map}
 \Phi^G\ : \ \pi_0^G (X) \quad &\to \qquad \Phi_0^G (X) \\
[f:S^V\to X(V)]\ &\longmapsto \ [f^G:S^{V^G}\to X(V)^G]   \nonumber
\end{align}
from the $G$-equivariant homotopy group to the 
geometric fixed point homotopy group.
For a trivial group, equivariant and geometric fixed point groups
coincide and the geometric fixed point map $\Phi^e:\pi_0^e(X)\to \Phi^e_0(X)$
is the identity.

\begin{eg}[Geometric fixed points of suspension spectra]
\label{eg:geometric and suspension}\index{subject}{geometric fixed points!of suspension spectra} 
If $A$ is any based $G$-space,
then the geometric fixed points $\Phi_*^G(\Sigma^\infty A)$ 
of the suspension spectrum are given by
\[ \Phi_k^G(\Sigma^\infty A) 
\ = \ \colim_{V\in s(\Uc_G)}\, [S^{V^G\oplus\mR^k}, S^{V^G}\sm A^G]  \ .\]
As $V$ ranges over $s(\Uc_G)$ the dimension of the fixed points
grows to infinity. So the composite
\[ \pi_k^e (\Sigma^\infty A^G)\ \xra{\ p_G^*\ } \
\pi_k^G (\Sigma^\infty A^G)\ \xra{\text{incl}_*} \
\pi_k^G (\Sigma^\infty A)\ \xra{\ \Phi^G\ }\ 
\Phi_k^G (\Sigma^\infty A)\]
is an isomorphism, where $p_G^*$ is inflation along the 
unique group homomorphism $p_G:G\to e$.
 We will sometimes refer to this isomorphism by saying that
`geometric fixed points commute with suspension spectra'.
\end{eg}

\begin{construction}\label{con:restriction in geometric fixed}
We let $X$ be a $G$-orthogonal spectrum and $\alpha:K\to G$ a continuous epimorphism.
We define {\em inflation maps}\index{subject}{inflation map!of geometric fixed point homotopy groups}  
\[ \alpha^* \ : \ \Phi^G_0 (X) \ \to \ \Phi^K_0 (\alpha^* X) \]
on geometric fixed point homotopy groups.
We choose a $K$-equivariant linear isometric embedding $\psi:\alpha^*(\Uc_G)\to \Uc_K$
of the restriction along $\alpha$ of the complete $G$-universe into
the complete $K$-universe. We let $f:S^{V^G}\to X(V)^G$ be a based map representing
an element of $\Phi_0^G(X)$, for some $V\in s(\Uc_G)$.
Since $\alpha$ is surjective, $V^G=(\alpha^* V)^K$ and
$X(V)^G=(\alpha^*(X(V)))^K=((\alpha^* X)(\alpha^* V))^K$. We use $\psi$ to identify
$\alpha^* V$ with $\psi(V)$ as $K$-representations, and 
hence also $(\alpha^* V)^K$ with $\psi(V)^K$. 
This turns $f$ into a based map
\[ S^{\psi(V)^K} \ \iso \  S^{(\alpha^* V)^K} \ = \ S^{V^G} \ \xra{\ f\ }\ 
X(V)^G \ = \ ((\alpha^* X)(\alpha^*V))^K \ \iso \ ((\alpha^* X)(\psi(V)))^K \ .\]
Any two equivariant embeddings of  $\alpha^*(\Uc_G)$ into $\Uc_K$
are homotopic through $K$-equivariant linear isometric embeddings, so
the restriction map is independent of the choice of $\psi$.
This latter map represents the element $\alpha^*[f]$ in $\Phi_0^K(\alpha^* X)$.
The element $\alpha^*[f]$ depends only on the class of $f$ in $\Phi_0^G(X)$,
so the inflation map $\alpha^*$ is a well-defined homomorphism.
\end{construction}

The surjectivity of $\alpha$ is essential to obtain an inflation map $\alpha^*$
on geometric fixed point homotopy groups, 
and geometric fixed points do not have natural restriction maps to subgroups.
These inflation maps between the geometric fixed point homotopy groups
are clearly natural in the orthogonal $G$-spectrum.
The next proposition lists the other naturality properties.

\begin{prop}\label{prop:properties geometric fixed maps} 
Let $G$ be a compact Lie group and $X$ an orthogonal $G$-spectrum.
  \begin{enumerate}[\em (i)]
  \item For every pair of composable continuous epimorphisms $\alpha:K\to G$
    and $\beta:L\to K$ we have
    \[ \beta^*\circ \alpha^* \ = \ (\alpha\beta)^* \ : \ 
    \Phi_0^G (X)\ \to \ \Phi_0^L( (\alpha\beta)^* X)\ .\]
  \item For every element $g\in G$ the composite
    \[ \Phi_0^G(X) \ \xra{\ (c_g)^*\ } \Phi_0^G (c_g^* X) \ \xra{\ (l_g^X)_*\ }
    \ \Phi_0^G (X) \]
    is the identity.
  \item For every continuous epimorphism $\alpha:K\to G$
    of compact Lie groups the following square commutes:
    \[ \xymatrix{ 
      \pi_0^G(X) \ar[r]^-{\Phi^G}\ar[d]_{\alpha^*} & \Phi_0^G(X) \ar[d]^{\alpha^*} \\
      \pi_0^K(\alpha^* X) \ar[r]_-{\Phi^K} & \Phi_0^K (\alpha^* X)} \]
  \end{enumerate}
\end{prop}
\begin{proof}
(i) We choose a $K$-equivariant linear isometric embedding $\psi:\alpha^*(\Uc_G)\to \Uc_K$
and an $L$-equivariant linear isometric embedding $\varphi:\beta^*(\Uc_K)\to \Uc_L$.
If we then use the $L$-equivariant linear isometric embedding
\[\varphi\circ\beta^*(\psi) \ : \ (\alpha\circ\beta)^*(\Uc_G)\ \to \ \Uc_L \]
for the calculation of $(\alpha\circ\beta)^*$, the desired equality even holds
on the level of representatives.

(ii) We let $V$ be a finite-dimensional $G$-subrepresentation of $\Uc_G$
and $f:S^{V^G}\to X(V)^G$ a based map representing an element of $\Phi_0^G(X)$.
We use the $G$-equivariant linear isometry $l_g:c_g^*(\Uc_G)\to\Uc_G$
given by left multiplication by $g$. 
Then $c_g^*(V)$ and $l_g(V)$ have the same underlying sets 
and the restriction of $l_g$ to $(c_g^*(V))^G$ is the identity onto $(l_g(V))^G$.

The class $c_g^*[f]$ is represented by the composite
\begin{align*}
  S^{(l_g(V))^G} =  S^{(c_g^* V)^G} = S^{V^G} \ \xra{\ f\ }\ 
X(V)^G &= ((c_g^* X)(c_g^*V))^G \\ 
&\xra{ ((c_g^*X)(l_g^V))^G} \ ((c_g^* X)(l_g(V)))^G \ .
\end{align*}
Consequently, $(l_g^X)_*(c_g^*[f])$ is represented by the composite
\begin{align*}
 S^{(l_g(V))^G} \ = \  S^{V^G} \ &\xra{\ f\ }\ 
X(V)^G \ = \ ((c_g^* X)(c_g^*V))^G \\ 
&\xra{ ((c_g^*X)(l_g^V))^G} \ ((c_g^* X)(l_g(V)))^G 
\ \xra{ ((l_g^X)(l_g(V)))^G} \ (X(l_g(V)))^G  \ .  
\end{align*}
The $G$-action on $X(V)$ is diagonally,
from the external $G$-action on $X$ and the internal $G$-action on $V$.   
Hence the map $l_g^{X(V)}:c_g^*(X(V))\to X(V)$
is the composite of $(c_g^* X)(l_g^V):(c_g^*X)(c_g^* V)\to (c_g^* X)(V)$
and $(l_g^X)(V):(c_g^* X)(V)\to X(V)$.
Since the restriction of $l_g^{X(V)}$ to the $G$-fixed points is the
identity, the composite of $( (c_g^* X)(l_g^V))^G$
and $((l_g^X)(V))^G$ is the identity.
This shows that $(l_g^X)_*(c_g^*[f])$ is again represented by $f$,
and hence $(l_g^X)_*\circ (c_g)^*=\Id$.

(iii) 
We consider a based continuous $G$-map $f:S^V\to X(V)$; 
then $\alpha^*[f]$ is represented by the $K$-map
\[ \alpha^*(f)\ : \ S^{\alpha^*(V)}\ = \ \alpha^*(S^V) \ 
\to \ \alpha^*(X(V))\ = \ (\alpha^* X)(\alpha^*(V))\ ,\]
and so $\Phi^K(\alpha^*[f])$ is represented by
\[ (\alpha^*(f))^K\ : \ (S^{\alpha^*(V)})^K\ = \ (\alpha^*(S^V))^K \ 
\to \ (\alpha^*(X(V)))^K\ = \ ((\alpha^* X)(\alpha^*(V)))^K\ .\]
Since $\alpha$ is surjective, this is the same map as
\[ f^G\ : \ (S^{V})^G\ \to \ X(V)^G\ .\]
To calculate $\alpha^*(\Phi^G([f]))$ 
we choose a $K$-equivariant linear isometric embedding 
$\psi:\alpha^*(\Uc_G)\to \Uc_K$ and conjugate $f^G$ by the isometry
\[ (\psi|_{\alpha^*(V)})^K\ : \ (\alpha^*(V))^K \ \iso \ \psi(V)^K\ . \]
But conjugation by an isometry does not change the stable homotopy class,
by Proposition \ref{prop:invariant description stable}~(ii)
for the trivial group. So 
\[ \Phi^K(\alpha^*[f])\ =\ [\alpha(f)^K]\ = \ [f^G]\ = \ \alpha^*(\Phi^G[f])\ .\qedhere \]
\end{proof}

\begin{rk}[Weyl group action on geometric fixed points]
\label{rk:Weyl group on geometric fixed}
We let $H$ be a closed subgroup of a compact Lie group $G$, and $X$ an
orthogonal $G$-spectrum. Every $g\in G$ gives rise to a conjugation homomorphism
$c_g:H\to H^g$ by $c_g(h)=g^{-1}h g$. Moreover,
left translation by $g$ is a homomorphism of orthogonal $H$-spectra
$l_g^X:c_g^*(X)\to X$. So combining inflation along $c_g$
with the effect of $l_g^X$ gives a homomorphism\index{subject}{conjugation homomorphism!on geometric fixed point homotopy groups}  
\[ g_\star\ : \ \Phi_0^{H^g}(X)\ \xra{\ (c_g)^*\ }\ \Phi_0^H(c_g^*(X))\ \xra{\ (l^X_g)_*\ }\
\Phi_0^H(X)\ . \]
In the special case when $g$ normalizes $H$, this is a self-map
of the geometric fixed point group $\Phi_0^H(X)$.
If moreover $g$ belongs to $H$, then $g_\star$ is the identity
by Proposition \ref{prop:properties geometric fixed maps}~(ii).
So the maps $g_\star$ define an action of the Weyl group
$W_G H=N_G H/H$ on the geometric fixed point homotopy group $\Phi_0^H(X)$.
By the same arguments as for equivariant homotopy groups
in Remark \ref{rk:Weyl on pi_0^G}, the identity path component of the 
Weyl group acts trivially, so the action factors over an action
of the discrete group $\pi_0(W_G H)=W_G H/(W_G H)^\circ$.
Since the geometric fixed point map $\Phi:\pi_0^H(X)\to\Phi_0^H(X)$
is compatible with inflation and natural in $X$, this map
is $\pi_0(W_G H)$-equivariant.
\end{rk}

Now we recall the interpretation of geometric fixed point homotopy groups 
as the equivariant homotopy groups 
of the smash product of $X$ with a certain universal $G$-space.
We denote by $\Pc_G$ the family of proper closed subgroups of $G$,
and by $E\Pc_G$ a universal space for the family $\Pc_G$.\index{subject}{universal space!for the family of proper subgroups}
So $E\Pc_G$ is a cofibrant $G$-space without $G$-fixed points and such that the
fixed point space $(E\Pc_G)^H$ is contractible for every closed proper subgroup $H$
of $G$.
These properties determine $E\Pc_G$ uniquely up to $G$-homotopy equivalence,
see Proposition \ref{prop:universal spaces}.

We denote by $\tilde E\Pc_G$ the reduced mapping cone of the based $G$-map 
$(E\Pc_G)_+\to S^0$ that sends $E\Pc_G$ to the non-basepoint of $S^0$. 
So $\tilde E\Pc_G$ is the unreduced suspension of the universal space $E\Pc_G$.
The $G$-fixed points of $E\Pc_G$ are empty and 
fixed points commute with mapping cones, so the map
$S^0\to (\tilde E\Pc_G)^G$ is an isomorphism. 
For all proper subgroups $H$ of $G$ the map $(E\Pc_G)^H_+\to (S^0)^H=S^0$ 
is a weak equivalence, so $(\tilde E\Pc_G)^H$ is contractible.

\begin{eg}\label{eg:S of U^perp is E tilde}
We let $\Uc_G^\perp=\Uc_G-(\Uc_G)^G$ be the orthogonal complement of
the $G$-fixed points in the complete $G$-universe $\Uc_G$. 
We claim that the unit sphere $S(\Uc_G^\perp)$
of this complement is a universal space $E\Pc_G$.
Indeed, the unit sphere $S(\Uc_G^\perp)$ is $G$-equivariantly homeomorphic
to the space $\bL(\mR,\Uc_G^\perp)$, so it is cofibrant as a $G$-space
by Proposition \ref{prop:K G cofibration}~(ii). 
Since $S(\Uc_G^\perp)$ has no $G$-fixed points, any stabilizer group
is a proper subgroup of $G$, i.e., in the family $\Pc_G$.
On the other hand, for every proper subgroup $H$ of $G$ 
there is a $G$-representation $V$ with $V^G=0$ but $V^H\ne 0$,
see for example \cite[III Prop.\,4.2]{broecker-tomDieck}.
Since $\Uc^\perp_G$ contains infinitely many isomorphic copies of $V$,
the $H$-fixed points
\[ (S(\Uc_G^\perp))^H \ = \  S((\Uc_G^\perp)^H) \]
form an infinite-dimensional sphere, and hence are contractible.
So $S(\Uc_G^\perp)$ is a universal $G$-space for the family of proper subgroups.

Since $\tilde E\Pc_G$ is an unreduced suspension of $E\Pc_G$,
it is equivariantly homeomorphic to 
\[ S(\mR\oplus \Uc_G^\perp)\  ,\]
the unit sphere in $\mR\oplus \Uc_G^\perp$.
So $S(\mR\oplus\Uc_G^\perp)$ is a model for $\tilde E\Pc_G$.
\end{eg}

The inclusion $i:S^0\to\tilde E\Pc_G$ induces an isomorphism of
$G$-fixed points $S^0\iso(\tilde E\Pc_G)^G$. 
So for every based $G$-space $A$ the map $A\sm i :A\to A\sm \tilde E\Pc_G$
induces an isomorphism of $G$-fixed points. Hence also for every orthogonal $G$-spectrum
the induced map of geometric fixed point homotopy groups
\[ \Phi^G_k(X\sm i)\ : \ \Phi_k^G(X) \ \iso \ \Phi_k^G(X\sm\tilde E\Pc_G)\]
is an isomorphism. If we compose the inverse with the geometric
fixed point homomorphism \eqref{eq:geometric_fix_map}, we arrive at a homomorphism
$\Phi:\pi_k^G(X\sm \tilde E\Pc_G)\to\Phi^G_k(X)$.

\begin{prop}\label{prop:geometric as fixed points} 
For every orthogonal $G$-spectrum $X$ and every integer $k$,
the geometric fixed point map
\[ \Phi\ : \ \pi_k^G (X\sm \tilde E\Pc_G) \ \to \ \Phi_k^G (X) \]
is an isomorphism.
\end{prop}
\begin{proof} 
We claim that for every finite based $G$-CW-complex $A$
and every based $G$-space $Y$ the map
\[ (-)^G \ : \ \map_*^G(A,Y\sm \tilde E\Pc_G) \ \to \ \map_*(A^G,Y^G)\]
that takes a $G$-map $f:A\to Y\sm \tilde E\Pc_G$ to the induced map
on $G$-fixed points
\[ f^G\ :\ A^G\ \to\  (Y\sm \tilde E\Pc_G)^G=Y^G\sm (\tilde E\Pc_G)^G\iso Y^G \]
is a weak equivalence and Serre fibration.

Indeed, since $A$ is a $G$-CW-complex, 
the inclusion of fixed points
$A^G\to A$ is a $G$-cofibration and induces a Serre fibration 
of equivariant mapping spaces
\[   \map_*^G(A,Y\sm \tilde E\Pc_G) \ \to \ 
\map_*^G(A^G,Y\sm \tilde E\Pc_G) \ .\]
Since every $G$-map from $A^G$ lands in the $G$-fixed points of
$Y\sm \tilde E\Pc_G$ and because $(Y\sm\tilde E\Pc_G)^G=Y^G$, 
the target space is the non-equivariant mapping space $\map_*(A^G,Y^G)$. 
The $G$-space $A$ is built from its fixed points 
by attaching $G$-cells $G/H\times D^n$ whose isotropy $H$ is a proper subgroup.
Since the $H$-fixed points of $Y\sm \tilde E\Pc_G$ are contractible 
for all proper subgroups $H$ of $G$, the fibration is also a weak equivalence.

Now we consider a finite-dimensional $G$-representation $V$.
When applied to $A=S^V$ and $Y=X(V)$, the previous claim implies that the fixed point map
\[ [S^V,X(V)\sm \tilde E\Pc_G]^G \ \to \ [S^{V^G}, X(V)^G]\]
is bijective for every $G$-representation $V$.
Passing to colimits over the poset $s(\Uc_G)$ proves
the result for $k=0$. The argument in the other dimensions is similar, and
we leave it to the reader.
\end{proof}

A consequence of the previous proposition is the following
{\em isotropy separation sequence}.\index{subject}{isotropy separation sequence} 
The mapping cone sequence of based $G$-space
\[ (E\Pc_G)_+\ \to \ S^0 \ \to \ \tilde E\Pc_G\]
becomes a mapping cone sequence of $G$-spectra 
\[ X\sm (E\Pc_G)_+ \ \to \ X \ \to \ X\sm \tilde E\Pc_G  \]
after smashing with any given orthogonal $G$-spectrum $X$. So
taking equivariant homotopy groups gives a long exact sequence
\begin{align}\label{eq:isotropy_separation_les}
 \cdots \ \to
\pi_k^G (X\sm (E\Pc_G)_+) \ &\to \
\pi_k^G (X) \ \xra{\ \Phi\ }  \\ 
&\Phi_k^G (X)\ \to \ \pi_{k-1}^G (X\sm (E\Pc_G)_+) \ \to \ \cdots   \nonumber
\end{align}
where we exploited the identification 
of Proposition \ref{prop:geometric as fixed points}.

\begin{prop}\label{prop:fix_and_geometric_fix} 
Let $G$ be a compact Lie group.
For a morphism $f:X\to Y$ of orthogonal $G$-spectra the following are equivalent:
\begin{enumerate}[\em (i)]
\item The morphism $f$ is a $\upi_*$-isomorphism.
\item For every closed subgroup $H$ of $G$ and every integer $k$ the map 
$\Phi_k^H(f)$ of geometric $H$-fixed point homotopy groups is an isomorphism.
\end{enumerate}
\end{prop}
\begin{proof} 
(i)$\Longrightarrow$(ii) If $f$ is
is an equivalence of orthogonal $G$-spectra, then so is $f\sm \tilde E\Pc_G$
by Proposition \ref{prop:space smash preserves global}~(ii).
Proposition \ref{prop:geometric as fixed points}
then implies that $\Phi_k^H(f):\Phi_k^H(X) \to \Phi_k^H(Y)$
is an isomorphism for all $k$. 

(ii)$\Longrightarrow$(i) We argue by induction on the size of the group $G$
(i.e., of the dimension of $G$ and the order of $\pi_0 G$).
If $G$ is the trivial group, then the geometric fixed point map 
$\Phi:\pi_k^e(X)\to\Phi_k^e(X)$ does not do anything, and is an isomorphism.
Since $\Phi_*^e(f)$ is an isomorphism, so is $\pi_*^e(f)$.

If $G$ is a non-trivial group we know by induction hypothesis that $f$ is an equivalence
of orthogonal $H$-spectra for every proper closed subgroup $H$ of $G$. 
Since $E\Pc_G$ is a cofibrant $G$-space without $G$-fixed points, 
Proposition \ref{prop:space smash preserves global}~(i)
lets us conclude that $f\sm (E\Pc_G)_+$ is an
equivalence of orthogonal $G$-spectra.
Since $\Phi_*^G(f):\Phi_*^G(X)\to\Phi_*^G(Y)$ 
is also an isomorphism, the isotropy separation sequence
and the five lemma let us conclude that 
$\pi_*^G(f):\pi_*^G(X)\to\pi_*^G(Y)$ is an isomorphism.
\end{proof}

The next proposition shows that `geometric fixed points vanish on transfers'.
In fact, it is often a helpful slogan to think of geometric fixed points
as `dividing out transfers from proper subgroups' -- despite the fact
that the kernel of the geometric fixed point map
$\Phi:\pi_0^G (X) \to \Phi_0^G (X)$ is in general larger than the subgroup
generated by proper transfers. For finite groups $G$, the slogan is 
in fact true {\em up to torsion}, i.e., 
the geometric fixed point map $\Phi:\pi_0^G (X) \to \Phi_0^G (X)$ is 
rationally surjective and its kernel is rationally generated by transfers
from proper subgroups, compare 
Proposition \ref{prop:rational geometric fixed for GSp} below.
\index{subject}{geometric fixed points!and transfers}
A more general version of part~(iii) below will
appear in Proposition \ref{prop:Phi^K after tr_H^G}~(ii).

\begin{prop}\label{prop:Phi of transfer} 
Let $K$ be a closed subgroup of a compact Lie group $G$
and $X$ an orthogonal $G$-spectrum.
\begin{enumerate}[\em (i)]
\item Let $H$ be a closed subgroup of $G$ such that $K$ is not subconjugate to $H$.
Then the composite
\[ \pi_0^G (X) \ \xra{\res^G_K}\  \pi_0^K (X) \ \xra{\ \Phi^K \ } \Phi_0^K (X)  \]
annihilates the images of the dimension shifting transfer 
$\Tr_H^G:\pi_0^H (X\sm S^L) \to\pi_0^G (X)$
and of the degree zero transfer $\tr_H^G:\pi_0^H (X) \to\pi_0^G (X)$. 
\item The geometric fixed point map $\Phi^G:\pi_0^G (X) \to \Phi_0^G (X)$ 
annihilates the images of the dimension shifting transfer 
and of the degree zero transfer from all proper closed subgroups of $G$.
\item
If the Weyl group $W_G K$ is finite, then the relation
  \[ \Phi^K\circ\res^G_K\circ \tr_K^G \ = \ 
  \sum_{ g K\in W_G K }\, \Phi^K\circ g_\star \]
holds as natural transformations from $\pi_0^K(X)$ to $\Phi_0^K(X)$.
\end{enumerate}
\end{prop}
\begin{proof}
(i) Let $V$ be any $G$-representation. The $G$-space
$(G\ltimes_H X)(V)$ is isomorphic to $G\ltimes_H X(i^* V)$.
If $K$ is not subconjugate to $H$, then both $G$-spaces 
have only one $K$-fixed point, the base point. 
So the geometric fixed point homotopy group $\Phi_0^K(G\ltimes_H X)$ vanishes.
The dimension shifting transfer is defined as the composite
\[ \pi_0^H (X\sm S^L) \ \xra{G\ltimes_H-}\ 
\pi_0^G (G\ltimes_H X) \ \xra{\text{act}_*}\ \pi_0^G (X)\ . \]
The geometric fixed point map is natural for $G$-maps,
so the 
composite $\Phi^K\circ\res^G_K\circ\text{act}_*:\pi_0^G (G\ltimes_H X) \to\Phi_0^K (X)$
factors through the trivial group $\Phi_0^K(G\ltimes_H X)$.
Thus the dimension shifting transfer vanishes.
The degree~0 transfer factors through the dimension shifting transfer,
so it vanishes as well.
Part~(ii) is the special case of~(i) for $K=G$.

(iii)
Since the functor $\pi_0^K$ is
represented by the suspension spectrum of $G/K$ 
(in the sense of Proposition \ref{prop:G/H represents pi_0^H}),
is suffices to check the relation for the orthogonal $G$-spectrum $\Sigma^\infty_+ G/K$
and the tautological class $e_K$ defined in \eqref{eq:tautological_e_H}.

For every $g\in N_G K$ and every class $x\in\pi_0^K(X)$ we have 
\begin{align*}
g_\star( \Phi^K(\res^G_K(\tr_K^G(x)))) \ = \ 
\Phi^K( \res^G_K( g_\star(\tr_K^G(x)))) \ = \ 
\Phi^K(  \res^G_K(\tr_K^G(x))) 
\end{align*}
because $g_\star$ is the identity on $\pi_0^G$.
So all classes of the form 
$\Phi^K(  \res^G_K(\tr_K^G(x)))$ are invariant under the action
of the Weyl group $W_G K$ specified in 
Remark \ref{rk:Weyl group on geometric fixed}.
In the universal case this class lives in the group $\Phi_0^K(\Sigma^\infty_+ G/K)$
which is $W_G K$-equivariantly isomorphic to
$\pi_0^e(\Sigma^\infty_+ (G/K)^K)=\pi_0^e(\Sigma^\infty_+ W_G K)$
and hence a free module of rank~1 over the integral group ring of
the Weyl group $W_G K$, generated by the class $\Phi^K(e_K)$.
So the class $\Phi^K(\res^G_K(\tr_K^G(e_K)))$
is an integer multiple of the norm element, i.e.,
\[   \Phi^K(\res^G_K(\tr_K^G(e_K))) \ = \ 
\lambda \cdot \sum_{g K\in W_G K} g_\star (\Phi^K(e_K )) \]
for some $\lambda\in\mZ$.

It remains to show that $\lambda=1$.
We let $1\in\pi_0^K(\mS)$ be the class represented by the identity of $S^0$.
Inspection of Construction \ref{con:transfer}
reveals that the transfer $\tr_K^G(1)$ in $\pi_0^G(\mS)$
is represented by the $G$-map
\[ S^V \ \xra{\ c \ } \ G\ltimes_K S^W\ \xra{\ a \ } \ S^V  \]
where $c$ is the collapse map based on any wide embedding
of $i:G/K\to V$ into a $G$-representation,
$W$ is the orthogonal complement of the image of $T_{e K}(G/K)$
under the differential of $i$, and $a[g,w]=g w$.
The class $\Phi^K_0(\res^G_K(\tr_K^G(1)))$ is then represented by the restriction
to $K$-fixed points of the above composite, i.e., by the map
\[ S^{V^K} \ \xra{\ c^K \ } \ (G\ltimes_K S^W)^K\ \xra{\ a^K \ }
\ S^{V^K} \ .\]
Every $K$-fixed point of $G\ltimes_K S^W$ is of the form
$[g,w]$ with $g\in N_G K$ and $w\in S^{W^K}$, i.e., the map
\[ (W_G K)_+\sm S^{W^K}\ \to \ (G\ltimes_K S^W)^K\ , \quad 
g K\sm w \ \longmapsto \ [g, w] \]
is a homeomorphism.

Since the Weyl group $W_G K$ is finite we have $(T_{e K}(G/K))^K=0$.
Indeed, the translation action of $K$ on the homogeneous space $G/K$ is smooth,
so by the differentiable slice theorem, 
the $K$-fixed point $e K$ has an open $K$-invariant neighborhood
inside $G/K$ that is $K$-equivariantly diffeomorphic 
to the tangent space $T_{e K} (G/K)$, compare \cite[Thm.\,1.6.5]{palais-classification}
or \cite[VI Cor.\,2.4]{bredon-intro}.
Since $W_G K=(G/K)^K$ is finite, $e K$ is an isolated $K$-fixed point in $G/K$,
and hence 0 is an isolated $K$-fixed point in $T_{e K} (G/K)$,
i.e., $(T_{e K} (G/K))^K=0$.

Since $(T_{e K}(G/K))^K=0$ we have $W^K=V^K$.
Under these identifications, the map $c^K$ becomes a pinch map
\[ S^{V^K}\ \to \ (W_G K)_+\sm S^{V^K}\ \iso \ \bigvee_{g K\in W_G K}\, S^{V^K}\ , \]
i.e., the projection to each wedge summand has degree~1.
On the other hand, the map $a^K$ becomes the fold map
\[  (W_G K)_+\sm S^{V^K}\ \to \ S^{V^K} \ . \]
So the degree of the composite $a^K\circ c^K$ is the order of the Weyl group $W_G K$.
We have thus shown that
\[ \Phi^K(\res^G_K(\tr^G_K(1)))\ = \ |W_G K|\cdot  \Phi^K(1)\]
in the group $\Phi^K_0(\mS)$.
On the other hand, the class~1 is invariant under the action of the Weyl group,
and hence
\[ \Phi^K(\res^G_K(\tr^G_K(1)))\ = \ 
\lambda\cdot   \sum_{ g K\in W_G K }\, 
g_\star (\Phi^K(1)) \ = \ \lambda \cdot |W_G K|\cdot \Phi^K(1)\ .\]
Since the abelian group $\Phi_0^K(\mS)$ is freely generated by $\Phi^K(1)$, we can compare
coefficients in the last two expressions and deduce that $\lambda=1$.
\end{proof}

The 0-th equivariant homotopy groups of equivariant {\em spectra} 
have two extra pieces of structure, compared to equivariant {\em spaces}: 
an abelian group structure and transfers.
Theorem \ref{thm:G-equivariant pi_0 of Sigma^infty} 
and Proposition \ref{prop:pi_0 of Sigma^infty} make precise,
first for suspension spectra of $G$-spaces and then for suspension spectra of
orthogonal spaces, that at the level of 0-th equivariant homotopy sets,
the suspension spectrum `freely builds in' the extra structure
that is available stably.

We introduce specific stabilization maps that relate unstable homotopy sets 
to stable homotopy groups.
We let $H$ be a compact Lie group and $Y$ an $H$-space.
We define a map\index{symbol}{$\sigma^G$ - {stabilization map $\pi_0^G (Y)\to\pi_0^G(\Sigma^\infty_+ Y)$}}
\begin{equation}\label{eq:define sigma^H}
  \sigma^H\ : \ \pi_0(Y^H) \ \to \ \pi_0^H(\Sigma^\infty_+ Y)     
\end{equation}
by sending the path component of an $H$-fixed point $y\in Y^H$
to the equivariant stable homotopy class 
$\sigma^H[y]$ represented by the $H$-map
\[ S^0\ \xra{-\sm y} \ S^0\sm Y_+ \ = \ (\Sigma^\infty_+ Y)(0) \ .\]
By direct inspection, the map $\sigma^H$ can be factored as the composition
\[ \pi_0(Y^H) \ \xra{\ \sigma^e\ } \  \pi_0^e(\Sigma^\infty_+ Y^H)  
\  \xra{\ p_H^*\ }\ \pi_0^H(\Sigma^\infty_+ Y^H)  
\  \xra{\text{incl}_*}\ \pi_0^H(\Sigma^\infty_+ Y)  \ ,
 \]
where $p_H:H\to e$ is the unique group homomorphism.

We recall that for every space $Z$ the non-equivariant stable homotopy group
$\pi_0^e(\Sigma^\infty_+ Z)$ is free abelian generated by the 
classes $\sigma^e(y)$ for all $y\in\pi_0(Z)$, i.e,
\begin{equation}\label{eq:pi_0_non-equivariant}
  \pi_0^e(\Sigma^\infty_+ Z)\ \iso \ \mZ\{\pi_0(Z)\} \ .
\end{equation}
Indeed, for all $n\geq 2$ the group $\pi_n(S^n\sm Z_+,\ast)$ is free abelian,
with basis the classes of the maps $-\sm y:S^n\to S^n\sm Z_+$
as $y$ runs over the path components of $Z$, 
see for example \cite[Prop.\,7.1.7]{tomDieck-algebraic topology}.
Passing to the colimit over $n$ proves the claim.

If $H$ is a closed subgroup of a compact Lie group $G$, 
and $Y$ is the underlying $H$-space of a $G$-space,
then the normalizer $N_G H$ leaves $Y^H$ invariant,
and the action of $N_G H$ factors over an
action of the Weyl group $W_G H=N_G H/H$ on $Y^H$.
This, in turn, induces an action of the component group $\pi_0(W_G H)$
on the set $\pi_0(Y^H)$.
For all $g\in G$, the following square commutes, again by direct inspection:
\[ \xymatrix{ 
\pi_0(Y^{H^g})\ar[r]^-{\sigma^{H^g}}\ar[d]_{\pi_0(l_g)}  &
\pi_0^{H^g}(\Sigma^\infty_+ Y)\ar[d]^{g_\star}\\
\pi_0(Y^H)\ar[r]_-{\sigma^H} & \pi_0^H(\Sigma^\infty_+ Y)} \]
Here $l_g:Y^{H^g}\to Y^H$ is left multiplication by $g$.
In particular, the map $\sigma^H$ is equivariant for the
action of the group $\pi_0(W_G H)$.

After stabilizing along the map
$\sigma^H:\pi_0(Y^H)\to\pi_0^H(\Sigma^\infty_+ Y)$, we can then transfer from $H$ to $G$.
For an element $g\in N_G H$ and a class $x\in\pi_0(Y^H)$ we have
\begin{align} \label{eq:Weyl_annihilated_transfer}
 \tr_H^G(\sigma^H(\pi_0(l_g)(x))) \ &= \ \tr_H^G(g_\star(\sigma^H(x))) \\ 
&= \ g_\star(\tr_H^G(\sigma^H(x))) \ = \ \tr_H^G(\sigma^H(x)) \ ,  \nonumber
\end{align}
because transfer commutes with conjugation,
and inner automorphisms act as the identity.
So the composite $\tr_H^G\circ \sigma^H$ coequalizes 
the $\pi_0(W_G H)$-action on $\pi_0(Y^H)$.

Our proof of the following theorem is based 
on an inductive argument with the isotropy separation sequence.
A different proof, based on the tom Dieck splitting, 
can be found in \cite[V Cor.\,9.3]{lms}.

\begin{theorem}\label{thm:G-equivariant pi_0 of Sigma^infty} 
Let $G$ be a compact Lie group and $Y$ a $G$-space.\index{subject}{suspension spectrum!of a $G$-space} 
\begin{enumerate}[\em (i)]
\item 
The equivariant homotopy group $\pi_0^G(\Sigma_+^\infty Y)$ 
is a free abelian group with a basis given by the elements
\[ \tr_H^G(\sigma^H(x))\ , \]
where $H$ runs through all conjugacy classes of closed subgroups of $G$
with finite Weyl group and $x$ runs through a set of representatives 
of the $W_G H$-orbits of the set $\pi_0(Y^H)$.
\item
Let $z\in \pi_0^G(\Sigma^\infty_+ Y)$ be a class
such that for every closed subgroup $K$ of $G$ with finite Weyl group
the geometric fixed point class
\[ \Phi^K( \res^G_K(z))\ \in \ \Phi_0^K(\Sigma^\infty_+ Y) \]
is trivial. Then $z=0$.
\end{enumerate}
\end{theorem}
\begin{proof}
(i)
In \eqref{eq:pi_0_non-equivariant} we recalled property~(i) when $G$ is a trivial group.
For the trivial group the geometric fixed point map
$\Phi^e:\pi_0^e(\Sigma^\infty_+ Y)\to\Phi_0^e(\Sigma^\infty_+ Y)$ 
is the identity, so part~(ii) is tautologically true.

Now we let $G$ be any compact Lie group.
We let $H$ be a closed subgroup of $G$ with finite Weyl group.
By \eqref{eq:pi_0_non-equivariant} we know that the group
$\pi_0^e(\Sigma^\infty_+ Y^H)$ is free abelian 
on the set of path components of $Y^H$; moreover, the Weyl group $W_G H$
permutes the basis elements, i.e., 
$\pi_0^e(\Sigma^\infty_+ Y^H)$ is an integral permutation representation of
the group $W_G H$. For any representation $M$ of a finite group $W$
the norm map 
\[ N \ : \ M \ \to \ M \ , \quad x\ \longmapsto \ \sum_{w\in W} w\cdot x \]
factors over the group of coinvariants
\[ M_W\ = \ M / \td{ x- w x\, |\, x\in M, w\in W} \ .\]
For the integral permutation representation $M=\mZ[S]$ of a $W$-set $S$,
a special feature is that the induced map $\bar N:M_W \to M$ is injective.

So part~(i) is equivalent to the claim that the map
\[  T \ : \ \bigoplus_{(H)}\, ( \pi_0^e(\Sigma^\infty_+ Y^H))_{W_G H}\ \to \ 
\pi_0^G(\Sigma_+^\infty Y) \]
is an isomorphism, where the sum is indexed by representatives of the
conjugacy classes of closed subgroups with finite Weyl group,
and the restriction of $T$ to the $H$-summand is induced by the composite
\[ \pi_0^e(\Sigma^\infty_+ Y^H)\ \xra{\ p_H^*\ } \ 
\pi_0^H(\Sigma^\infty_+ Y^H)\ \xra{\text{incl}_*} \ 
\pi_0^H(\Sigma^\infty_+ Y)\ \xra{\ \tr_H^G\ } \ 
\pi_0^G(\Sigma^\infty_+ Y)\ . \]
We consider the total geometric fixed point homomorphism
\[ \Phi^{\text{total}} \ : \ 
\pi_0^G(\Sigma_+^\infty Y) \ \to \ \prod_{K:\, W_G K\text{ finite}} \Phi_0^K(\Sigma^\infty_+ Y)
\ , \quad y\ \longmapsto \ (\Phi^K(\res^G_K(y)))_K\ .\]
Property~(ii) is equivalent to the claim that $\Phi^{\text{total}}$ is injective.

We show now that the composite $\Phi^{\text{total}}\circ T$ is injective.
We argue by contradiction and suppose that
\[ z \ = \ (z_H)_H\ \in \ 
\bigoplus_{(H)}\, ( \pi_0^e(\Sigma^\infty_+ Y^H))_{W_G H} \]
is a non-zero element in the source of $T$
such that $\Phi^{\text{total}}( T(z))$ vanishes.
As an element in a direct sum, $z$ has only finitely many non-zero components $z_H$.
We let $K$ be of maximal dimension and with maximal
number of path components among all indexing subgroups 
such that $z_K\ne 0$.
Then
\[ T(z) \ = \ \tr_K^G(\text{incl}_*(p_K^*(y)))\ +\ 
{\sum}_{i=1}^m\,  \tr_{H_i}^G(y_i) \]
with $y$ an element of $\pi_0^e(\Sigma^\infty_+ Y^K)$
with non-zero image in the $W_G K$-coinvariants,
and with certain closed subgroups $H_i$ of $G$
that are not conjugate to $K$ and `no larger' in the sense that
either $\dim(H_i) <\dim(K)$, or $\dim(H_i) =\dim(K)$ and $|\pi_0(H_i)|\leq |\pi_0(K)|$.
This means in particular that $K$ is not subconjugate to any of the
groups $H_1,\dots,H_m$. Thus
\[ \Phi^K(\res^G_K(\tr_{H_i}^G(y_i))) \ = \ 0 \]
for all $i=1,\dots,m$, by Proposition \ref{prop:Phi of transfer}~(i).
Hence
\begin{align*}
0 \ = \ \Phi^K(T(z))\ &= \ 
\Phi^K(\res^G_K(\tr_K^G(\text{incl}_*(p_K^*(y))))) \\
& = \ \sum_{g K\in W_G K} \Phi^K(g_\star(\text{incl}_*(p_K^*(y)))) \\
& = \ \sum_{g K\in W_G K} \Phi^K(\text{incl}_*(g_\star(p_K^*(y)))) \\
&= \ \Phi^K(\text{incl}_*(p_K^*({\sum}_{g K\in W_G K}  (l_g)_*(y)))) \ .
\end{align*}
The third equation is Proposition \ref{prop:Phi of transfer}~(iii).
Since the composite
\[ \pi_0^e(\Sigma^\infty_+ Y^K)\ \xra{\ p_K^*\ }\ 
\pi_0^K(\Sigma^\infty_+ Y^K)\ \xra{\text{incl}_*}\ 
\pi_0^K(\Sigma^\infty_+ Y)\ \xra{\ \Phi^K\ }\ 
\Phi_0^K(\Sigma^\infty_+ Y)  \]
is an isomorphism, we conclude that the 
norm of the element $y\in \pi_0^e(\Sigma^\infty_+ Y^K)$ is zero. 
But since $\pi_0^e(\Sigma^\infty_+ Y^K)$ is an integral permutation
representation of the Weyl group, this only happens if
$y$ maps to~0 in the coinvariants, which contradicts our assumption.

Now we show that the classes $\tr_H^G(\sigma^H(x))$ in the statement of~(i)
generate the group $\pi_0^G(\Sigma^\infty_+ Y)$
(i.e., the homomorphism $T$ is surjective).
We argue by induction on the size of $G$, 
i.e., by a double induction over the dimension and number of path components of $G$. 
The induction starts when $G$ is the trivial group, which we dealt with above.
Now we let $G$ be a non-trivial compact Lie group.
We start with the special case $Y=G/K_+$
for a proper closed subgroup $K$ of $G$. 
The composite
\[ G\ltimes_K\mS \ \xra{G\ltimes_K (e K)_*}\ 
G\ltimes_K (\Sigma^\infty_+ G/K) \ \xra{\ \text{act}\ }\ 
 \Sigma^\infty_+ G/K \]
is an isomorphism of orthogonal $G$-spectra.
Hence the composite
\begin{align*}
\pi_0^K(\Sigma^\infty S^V)\ \xra[\iso]{ G\ltimes_K -} \ 
&\pi_0^G(G\ltimes\mS) \ \xra{\pi_0^G(G\ltimes_K (e K)_*)} \\ 
& \pi_0^G( G\ltimes_K (\Sigma^\infty_+ G/K)) \ \xra{\pi_0^G(\text{act})}\ 
\pi_0^G(\Sigma^\infty_+ G/K)   
\end{align*}
is an isomorphism of abelian groups, 
where $V=T_{e K}(G/K)$ is the tangent representation and
the first map is the external transfer 
(an isomorphism by Theorem \ref{thm:Wirth iso}). 

The inclusion $S^0\to S^V$ is an equivariant h-cofibration
and its quotient $S^V/S^0$ is $G$-homeomorphic to the unreduced suspension
of the unit sphere $S(V)$
(with respect to any $K$-invariant scalar product on $V$). So the group
\[ \pi_0^G(\Sigma^\infty (S^V/S^0)) \ \iso \ \pi_0^G(\Sigma^\infty (S(V)_+\sm S^1))
\ \iso \ \pi_{-1}^G(\Sigma^\infty_+ S(V)) \]
vanishes by the suspension isomorphism and 
Proposition \ref{prop:suspension spectrum homotopical}~(ii).
The long exact sequence of 
Corollary \ref{cor-long exact sequence h-cofibration}~(i)
then shows that the map
\[ \text{incl}_* \ : \
\pi_0^K(\Sigma^\infty_+ \{0\})\ \to \ \pi_0^K(\Sigma^\infty S^V) \]
is surjective.
Since $K$ is a proper closed subgroup of $G$,
it either has smaller dimension or fewer path components, so we
know by the inductive hypothesis that the group
$\pi_0^K(\Sigma^\infty_+ \{0\})$ is generated by the elements $\tr_L^K(\sigma^L[0])$
where $L$ runs through all conjugacy classes of closed subgroups of $K$
with finite Weyl group.

Putting this all together lets us conclude 
that the group $\pi_0^G(\Sigma^\infty_+ G/K)$ is generated by the classes 
\[ \left( \pi_0^G(\text{act}\circ (G\ltimes_K(e K)_*))\circ 
(G\ltimes_K-)\circ\tr_L^K\circ\sigma^L\right)[0] \]
for all $K$-conjugacy classes of closed subgroups $L\leq K$ 
that have finite Weyl group in $K$.
However, this long expression in fact defines a familiar class, as we shall now see.
Indeed, the following diagram commutes by the various naturality properties:
\[ \xymatrix@C=11mm@R=8mm{ 
\pi_0(\{0\})\ar@<-4ex>@/_1pc/[ddd]_(.3){\sigma^L} \ar[rr]^-{\pi_0((e K)_*)}
\ar[d]^{p_L^*\circ\sigma^e}  &&
 \pi_0( (G/K)^L) \ar[d]_-{p_L^*\circ \sigma^e} \ar@<0ex>@/^2pc/[ddr]^(.7){\sigma^L} \\
\pi_0^L(\Sigma^\infty S^0)\ar[dd]^{\text{incl}_*} \ar[rr]^-{\pi_0^L((e K)_*)}
&& 
\pi_0^L( \Sigma^\infty_+ (G/K)^L) \ar[dr]_{\text{incl}_*}  \\
&&& \pi_0^L(\Sigma^\infty_+ G/K)\ar[d]^{\tr_L^K} \ar[dl]_{\text{incl}_*} \\ 
\pi_0^L(\Sigma^\infty_+ S^V)\ar[d]_{\tr_L^K} \ar[rr]^-{\pi_0^L((e K)_*\sm S^V)} 
&& \pi_0^L(\Sigma^\infty_+ G/K\sm S^V)\ar[d]_{\tr_L^K} &
\pi_0^K(\Sigma^\infty_+ G/K)\ar[dl]_{\text{incl}_*} \ar[dd]^{\tr_K^G}\\
\pi_0^K(\Sigma^\infty_+ S^V)\ar[d]^\iso_{G\ltimes_K -} 
\ar[rr]^{\pi_0^K((e K)_*\sm S^V)}
&& 
\pi_0^K(\Sigma^\infty_+ G/K\sm S^V)\ar[d]_{G\ltimes_K -} \\
\pi_0^G(G\ltimes_K \mS )\ar[rr]_-{\pi_0^G(G\ltimes_K (e K)_*))} 
&&\pi_0^G(G\ltimes_K \Sigma^\infty_+ G/K)  \ar[r]_-{\pi_0^G(\text{act})} &
 \pi_0^G(\Sigma^\infty_+ G/K) } \]
So 
\begin{align*}
 \left( \pi_0^G(\text{act}\circ (G\ltimes_K(e K)_*))\circ 
(G\ltimes_K-)\circ\tr_L^K\circ\sigma^L\right)[0]& \\
= \  \tr_K^G(\tr_L^K(\sigma^L( (e K)_*[0] )))& \ =\  \tr_L^G( \sigma^L[e K] ) \ .
\end{align*}
So the group $\pi_0^G(\Sigma^\infty_+ G/K)$ is generated by the classes 
$\tr_L^G( \sigma^L[e K] )$ for all $K$-conjugacy classes of closed subgroups $L\leq K$ 
that have finite Weyl group in $K$.
If the Weyl group of $L$ in the ambient group $G$  happens to be infinite, 
then $\tr_L^G=0$ and the generator is redundant. 
Otherwise $e K$ is an $L$-fixed point of $G/K$, so the generator is one
of the classes mentioned in the statement of~(i).
This shows the generating property for the $G$-space $G/K$.

Next we observe that whenever the claim is true for a family
of $G$-spaces, then it is also true for their disjoint union;
this follows from the fact that both fixed points and $\pi_0$
commute with disjoint unions, and
that equivariant homotopy groups take wedges to direct sums.
In particular, the claim holds when $Y$ is a disjoint union
of homogeneous space $G/H$ for varying proper closed subgroups $H$ of $G$.

Now we prove the claim when $Y$ is any $G$-space
without $G$-fixed points.
For every proper closed subgroup $H$ of $G$ we choose representatives 
for the path components of the $H$-fixed point space $Y^H$. 
These choices determine a continuous $G$-map
\[ f\ : \ Z \ = \ \coprod_{H}\ \coprod_{[x]\in \pi_0(Y^H)} \, G/H \ \to \ Y \]
by sending $g H$ in the summand indexed by $x\in Y^H$ to $g x$.
Because $Y$ has no $G$-fixed points, the induced map $\pi_0(f^H):\pi_0(Z^H)\to \pi_0(Y^H)$
is then surjective for every closed subgroup $H$ of $G$, by construction.
In the commutative square
\[ \xymatrix{
 \bigoplus_{(H)}\, ( \pi_0^e(\Sigma^\infty_+ Z^H))_{W_G H}\ar[r]^-T
\ar[d]_{\bigoplus (\pi_0^e(\Sigma^\infty_+ f^H))_{W_G H}} &
\pi_0^G(\Sigma_+^\infty Z) \ar[d]^{\pi_0^G(\Sigma^\infty_+ f)}\\
 \bigoplus_{(H)}\, ( \pi_0^e(\Sigma^\infty_+ Y^H))_{W_G H}\ar[r]_-T &
\pi_0^G(\Sigma_+^\infty Y)}  \]
the right vertical map is then surjective by
Proposition \ref{prop:suspension spectrum homotopical}~(i) for $m=0$,
and the upper horizontal map is surjective by the above.
So the lower horizontal map is surjective as well.

Now we let $Y$ be any $G$-space, possibly with $G$-fixed points. 
The composite
\[ \pi_*^e(\Sigma^\infty_+ Y^G) \ \xra{\ p_G^*\ } \ 
\pi_*^G(\Sigma^\infty_+ Y^G) \ \xra{\text{incl}_*} \ 
\pi_*^G(\Sigma^\infty_+ Y) \ \xra{\ \Phi^G \ } \ 
\Phi_*^G(\Sigma^\infty_+ Y) 
\]
is an isomorphism, i.e., the map $\text{incl}_*\circ p_G^*$ splits
the geometric fixed point homomorphism.
The long exact isotropy separation sequence \eqref{eq:isotropy_separation_les} 
thus decomposes into short exact sequences and the map
\[ 
(q_*,\,\text{incl}_*\circ p_G^*)\ : \
\pi_*^G(\Sigma^\infty_+ (Y\times E\Pc_G) ) \oplus \pi_*^e(\Sigma^\infty_+ Y^G) 
\ \to\ \pi_*^G(\Sigma^\infty_+ Y)   \]
is an isomorphism, where $q:Y\times E\Pc_G\to Y$ is the projection. 
The $G$-space $Y\times E\Pc_G$ has no $G$-fixed points.
So by the previous part, the group 
$\pi_0^G(\Sigma^\infty_+ (Y\times E\Pc_G))$ is generated by
the classes $\tr_H^G(\sigma^H(x) )$ for $H$ with finite Weyl group as above 
and $x\in \pi_0( (Y\times E\Pc_G)^H)$.
If $H$ is a proper subgroup, then $(E\Pc_G)^H$ is contractible and
so the projection $q:Y\times E\Pc_G\to Y$ induces a bijection on $\pi_0((-)^H)$.
So the generators coming from $Y\times E\Pc_G$
map to the desired basis elements for $Y$ that are
indexed by proper subgroups of $G$.
On the other hand, the group $\pi_*^e(\Sigma^\infty_+ Y^G)$ is free
with basis the classes $\sigma^e(x)$ for all components $x\in \pi_0(Y^G)$,
by \eqref{eq:pi_0_non-equivariant}.
Because
\[ \text{incl}_*(p_G^*(\sigma^e(x)))\ = \ \sigma^G(x) \ , \]
the basis of the second summand precisely maps
to the desired basis elements for $Y$ that are
indexed by the group $G$ itself. This completes the argument 
that the map $T$ is surjective.

Now we can complete the proofs of part~(i) and~(ii).
Since the composite $\Phi^{\text{total}}\circ T$ is injective, 
the homomorphism $T$ is injective; since $T$ is also surjective, 
it is bijective, which shows~(i) for the group $G$.
Since $T$ is bijective and $\Phi^{\text{total}}\circ T$ is injective, 
the homomorphism $\Phi^{\text{total}}$ is injective, which shows~(ii) for the group $G$.
\end{proof}

\begin{eg}[Equivariant 0-stems and Burnside rings]
We recall that for a finite group $G$, the {\em Burnside ring} $A(G)$
is the Grothendieck group of the abelian monoid, under disjoint union,
of isomorphism classes of finite $G$-sets.\index{subject}{Burnside ring}\index{subject}{sphere spectrum}
Every finite $G$-set is the disjoint union of transitive $G$-sets,
so $A(G)$ is a free abelian group generated by the classes
of the $G$-sets $G/H$, as $H$ runs over all conjugacy classes
of subgroups of $G$.

The sphere spectrum $\mS$ is also the unreduced suspension spectrum
of the one-point $G$-space, $\mS\iso\Sigma^\infty_+\{0\}$.
So Theorem \ref{thm:G-equivariant pi_0 of Sigma^infty}~(i) 
says that the map
\[ \psi_G \ : \ A(G)\ \to \ \pi_0^G(\mS)\ , \quad [G/H]\ \longmapsto \tr_H^G(1) \]
is an isomorphism, a result that is originally due to Segal \cite{segal-ICM}.
The Burnside rings for different groups are related by restriction 
homomorphisms $\alpha^*:A(G)\to A(K)$
along homomorphisms $\alpha:K\to G$, induced by restriction of
the action along $\alpha$.
The Burnside rings also enjoy transfer maps
\[ \tr_H^G\ : \ A(H)\ \to \ A(G) \]
induced by sending an $H$-set $S$ to the induced $G$-set $G\times_H S$.

The compatibility of transfers with inflations 
(Proposition \ref{prop:transfer and epi}~(iii)) 
implies that for every surjective homomorphism $\alpha$ the relation
\begin{align*}
   \alpha^*(\psi_G[G/H])\ &= \ \alpha^*(\tr_H^G(1))\ = \ 
\tr_L^K((\alpha|_L)^*(1))\\ 
&= \ \tr_L^K(1)\ = \ \psi_K([K/L]) \ = \ \psi_K(\alpha^*[G/H]) 
\end{align*}
holds. In other words, the isomorphisms $\psi_G$ are compatible
with inflation. 
The double coset formula (see Example \ref{eg:double coset finite index} below),
and the double coset formula in the Burnside rings imply
that the isomorphisms $\psi_G$ are also compatible with transfers. 

If $G$ is a compact Lie group of positive dimension, there is no direct
interpretation of $\pi_0^G(\mS)$ in terms of finite $G$-sets;
in that situation, some authors {\em define} the Burnside ring of $G$
as the group $\pi_0^G(\mS)$.
In \cite[Prop.\,5.5.1]{tomDieck-and representation theory}, 
tom Dieck gives an interpretation of this Burnside ring 
in terms of certain equivalence classes of compact $G$-ENRs.
\index{subject}{euclidean neighborhood retract!equivariant} 
Theorem \ref{thm:G-equivariant pi_0 of Sigma^infty}~(i) identifies $\pi_0^G(\mS)$
as a free abelian group with basis the classes $\tr_H^G(1)$
for all conjugacy classes of closed subgroups $H$ with finite Weyl group.

Theorem \ref{thm:G-equivariant pi_0 of Sigma^infty}~(ii) 
provides a convenient detection criterion for elements in the equivariant
0-stems; as we explain now, this can be rephrased as a degree function.
For every compact Lie group $K$ the 
geometric fixed point group $\Phi^K_0(\mS)$
is the 0-th non-equivariant stable stem, and hence free abelian of rank~1,
generated by the class $\Phi^K(1)$. 
Moreover, if $z\in\pi_0^K(\mS)$ is represented by the $K$-map
$f:S^V\to S^V$ for some $K$-representation $V$, then
\[ \Phi^K(z)\ = \ [f^K\ : \ S^{V^K}\to S^{V^K}]\ = \
\deg(f^K)\cdot \Phi^K(1)\ . \]
So in terms of the basis $\Phi^K(1)$, the geometric fixed point
homomorphism $\Phi^K:\pi_0^K(\mS)\to\Phi_0^K(\mS)$
is extracting the degree of the restriction of any representative
to $K$-fixed points.
We let $C(G)$ denote the set of integer valued, conjugation-invariant functions
from the set of closed subgroup of $G$ with finite Weyl group.
We define a homomorphism
\begin{equation}\label{eq:degree_homomorphism}
 \deg \ : \ \pi_0^G(\mS)\ \to \ C(G)  
\end{equation}
by
\[ \deg[f:S^V\to S^V](K)\ = \ \deg(f^K\ : \ S^{V^K}\to S^{V^K})\ . \]
Theorem \ref{thm:G-equivariant pi_0 of Sigma^infty}~(ii) 
then says that this degree homomorphism is injective.
Whenever the group $G$ is non-trivial, the degree homomorphism
is {\em not} surjective. 
In \cite[Prop.\,5.8.5]{tomDieck-and representation theory},
tom~Dieck exhibits an explicit set of congruences that 
together with a certain continuity condition characterize the 
image of the degree homomorphism \eqref{eq:degree_homomorphism}.
When $G$ is finite, these congruences -- combined with the isomorphism $\psi_G$
to the Burnside ring $A(G)$ -- specialize to the
congruences between the cardinalities of the various fixed point sets
of a finite $G$-set, specified for example 
in \cite[Prop.\,1.3.5]{tomDieck-and representation theory}.

We will revisit the equivariant 0-stems from a global perspective in
Example \ref{eg:sphere spectrum} and with a view towards products
and multiplicative power operations in Example \ref{eg:Burnside global power}.
\end{eg}

Theorem \ref{thm:G-equivariant pi_0 of Sigma^infty}
gives a functorial description of the 0-th equivariant stable homotopy
group of an unbased $G$-space. We will now deduce a similar result
for reduced suspension spectra of {\em based} $G$-spaces.
This version, however, needs a non-degeneracy hypothesis,
i.e., we must restrict to well-pointed $G$-spaces.

\begin{theorem}
Let $G$ be a compact Lie group and $Y$ a well-pointed $G$-space.\index{subject}{well-pointed} 
Then the group $\pi_0^G(\Sigma^\infty Y)$ 
is a free abelian group with a basis given by the elements
\[ \tr_H^G(\sigma^H(x))\ , \]
where $H$ runs through all conjugacy classes of closed subgroups of $G$
with finite Weyl group and $x$ runs through a set of representatives 
of the $W_G H$-orbits of the non-basepoint components of the set $\pi_0(Y^H)$.
\end{theorem}
\begin{proof}
If the inclusion of the basepoint $i:\{y_0\}\to Y$ is an unbased h-cofibration,
then the based map $i_+:\{y_0\}_+\to Y_+$ is a based h-cofibration.
The induced morphism of suspension spectra 
$\Sigma^\infty_+ i:\Sigma^\infty_+\{y_0\}\to \Sigma^\infty_+ Y$ 
is then an h-cofibration of orthogonal $G$-spectra,
so it gives rise to a long exact sequence 
of homotopy groups as in Corollary \ref{cor-long exact sequence h-cofibration}~(i).
The cokernel of $i_+:\{y_0\}_+\to Y_+$ is $G$-homeomorphic to $Y$
(with the original basepoint). Also, the map $i$ has a section,
so the long exact sequence degenerates into a short exact sequence:
\[ 0\ \to \ \pi_0^G (\Sigma^\infty_+ \{y_0\}) \
\ \xra{(\Sigma^\infty_+ i)_*} \ \pi_0^G (\Sigma^\infty_+ Y) \
\xra{\text{proj}_*} \ \pi_0^G (\Sigma^\infty Y) \ \to \ 0 \]
Theorem \ref{thm:G-equivariant pi_0 of Sigma^infty}~(i),
applied to the unbased $G$-spaces $\{y_0\}$ and $Y$,
provides bases of the left and middle group; by naturality
the map  $(\Sigma^\infty_+ i)_*$ hits precisely
the subgroup generated by the basis elements coming from basepoint components
of the space $Y^H$.
So the cokernel of $(\Sigma^\infty_+ i)_*$,
and hence the group $\pi_0^G (\Sigma^\infty Y)$, has a basis of the desired form.
\end{proof}

\index{subject}{geometric fixed points|)}

\section{The double coset formula}
\label{sec:double coset}
  
The main  aim of this section is to establish
the double coset formula for the composite of a transfer 
followed by a restriction to a closed subgroup,
see Theorem \ref{thm:double coset formula} below.
We also discuss various examples and special cases in 
Examples \ref{eg:restrict transfer to e} through \ref{eg:double coset for Brauer}.
For finite groups (or more generally for transfers along finite index inclusions),
the statement and proof of the double coset formula are significantly simpler,
see Example \ref{eg:double coset finite index},
as this special case does not need any equivariant differential topology.
We end the section with a discussion of Mackey functors for
finite groups, and prove that rationally and for finite groups, 
geometric fixed point homotopy groups can be obtained from equivariant homotopy groups 
by dividing out transfers
from proper subgroups (Proposition \ref{prop:rational geometric fixed for GSp}).

\medskip

For use in the proof of the double coset formula,
and as an interesting result on its own, we calculate
the geometric fixed points of the restriction of any transfer.
In fact, when $K=e$ is the trivial subgroup of $G$,
then $\Phi_0^e(X)=\pi_0^e(X)$, the map $\Phi^e$ is the identity
and the following proposition reduces to a special case of the double coset formula.
For the statement we use that the action of a closed subgroup $K$
on a homogeneous space $G/H$ by left translation is smooth,
and hence the fixed points $(G/H)^K$ form a disjoint union of
closed smooth submanifolds, possibly of varying dimensions.

We let $K$ be a compact Lie group and $B$ a closed smooth $K$-manifold.
The Mostow-Palais embedding theorem \cite{mostow-embeddings, palais-imbedding}
provides a smooth $K$-equiva\-riant embedding $i: B\to V$, 
for some $K$-representation $V$. 
We can assume without loss of generality that $V$ is a subrepresentation
of the chosen complete $K$-universe $\Uc_K$.
We use the inner product on $V$ to define the normal bundle $\nu$ 
of the embedding at $b\in B$ by
\[ \nu_b \ = \ V - (d i)(T_b B)\ ,\]
the orthogonal complement of the image of the tangent space $T_b B$ in $V$.
By multiplying with a suitably large scalar, if necessary, 
we can assume that the embedding
is {\em wide}\index{subject}{wide embedding}\index{subject}{embedding!wide|see{wide embedding}}
in the sense that the exponential map
\[ D(\nu) \ \to \ V \ , \quad (b,v)\ \longmapsto \ i(b) \, +\,  v \]
is injective on the unit disc bundle of the normal bundle, 
and hence a closed $K$-equivariant embedding.
The image of this map is a tubular neighborhood of radius~1 around $i(B)$,
and it determines a $K$-equivariant Thom-Pontryagin collapse map
\begin{equation}  \label{eq:collaps B}
 c_B\ : \ S^V \ \to \  S^V\sm B_+ \end{equation}
as follows: every point outside of the tubular neighborhood is sent
to the basepoint, and a point $i(b)+ v$, for $(b,v)\in D(\nu)$,
is sent to 
\[  c_B(i(b) +  v) \ = \ \left( \frac{v}{1-|v|} \right) \sm  b \ . \]
The homotopy class of the map $c_B$ is then an element in
the equivariant homotopy group $\pi_0^K(\Sigma^\infty_+ B)$.
The next result determines the image of the class $[c_B]$
under the geometric fixed point map.

In \eqref{eq:define sigma^H} 
we defined the map $\sigma^K:\pi_0(B^K)\to\pi_0^K(\Sigma^\infty_+ B)$ that
produces equivariant stable homotopy classes from fixed point information.

\begin{prop}\label{prop:Phi^K after tr_H^G}
Let $K$ be a compact Lie group.
\begin{enumerate}[\em (i)]
\item  For every closed smooth $K$-manifold $B$ the relation
  \[ \Phi^K[c_B]\ = \ 
  \sum_{ M\in \pi_0( B^K) }\, \chi(M)\cdot \Phi^K(\sigma^K[M])\]
  holds in the group $\Phi_0^K(\Sigma^\infty_+ B)$,
  where the sum runs over the connected components of the fixed point manifold $B^K$.
\item Let $G$ be a compact Lie group containing $K$ and let $H$ be another
  closed subgroup of $G$. Then
  \[ \Phi^K\circ\res^G_K\circ \tr_H^G \ = \ 
  \sum_{ M\in \pi_0( (G/H)^K) }\, \chi(M)\cdot \Phi^K\circ g_\star\circ\res^H_{K^g} \ ,\]
  where the sum runs over the connected components of the fixed point
  space $(G/H)^K$ and $g\in G$ is such that $g H\in M$.
\end{enumerate}
\end{prop}
\begin{proof}
(i)
Since the composite
\[ \pi_0^e(\Sigma^\infty_+ B^K )\  \xra{\ p_K^*\ }\ 
 \pi_0^K(\Sigma^\infty_+ B^K) \ \xra{\text{incl}_*} \ \pi_0^K(\Sigma^\infty_+ B) 
\ \xra{\ \Phi^K\ }\ \Phi_0^K(\Sigma^\infty_+ B) 
\]
is an isomorphism and the source is free abelian on the path components
of the space $B^K$,
the group $\Phi_0^K(\Sigma^\infty_+ B)$ is free abelian with basis
given by the classes
\[ \Phi^K(\sigma^K[M]) \ = \ \Phi^K(\text{incl}_*(p_K^*(\sigma^e[M]))) \]
for $M\in\pi_0(B^K)$.
We let $i: B\to V$ be a wide smooth $K$-equivariant embedding into a $K$-representation. 
Then the class $\Phi^K[c_B]$ is represented by the non-equivariant map
\[ S^{V^K} \ \xra{\ (c_B)^K \ } \ S^{V^K}\sm B^K_+ \ .\]
Since $K$ acts smoothly on $B$, the fixed points $B^K$ 
are a disjoint union of finitely many closed smooth submanifolds.
For each component $M$ of $B^K$ we let $p_M:B^K_+\to M_+$ denote
the projection, i.e., $p_M$ is the identity on $M$ and sends all
other path components of $B^K$ to the basepoint.
Then the composite
\[ S^{V^K} \ \xra{\ (c_B)^K \ } \ S^{V^K}\sm B^K_+ \ \xra{S^{V^K}\sm p_M}\  S^{V^K}\sm M_+ \]
coincides with 
\[ S^{V^K} \ \xra{\ c_M \ } \ S^{V^K}\sm M_+ \ , \]
the collapse map \eqref{eq:collaps B} 
based on the non-equivariant wide smooth embedding $(i^K)|_M:M\to V^K$.
Since $M$ is path connected, the group of based homotopy classes
of maps $[S^{V^K},S^{V^K}\sm M_+]$ is isomorphic to $\mZ$,
by \eqref{eq:pi_0_non-equivariant},
and an element is determined by the degree of the composite
with the projection $S^{V^K}\sm M_+\to S^{V^K}$.
It is a classical fact that the degree of the composite
\[ S^{V^K} \ \xra{\ c_M \ } \ S^{V^K}\sm M_+ \ \xra{\text{proj}} \ S^{V^K} \]
is the Euler characteristic of the manifold $M$, 
see for example \cite[Thm.\,2.4]{becker-gottlieb}.
So the summand indexed by $M$ precisely contributes the term
$\chi(M)\cdot \Phi^K(\sigma^K[M])$, which proves the desired relation.

(ii)
Both sides of the equation are natural transformations 
on the category of orthogonal $G$-spectra from the functor $\pi_0^H$
to the functor $\Phi^K_0$. Since the functor $\pi_0^H$ is
represented by the suspension spectrum of $G/H$ 
(in the sense of Proposition \ref{prop:G/H represents pi_0^H}),
is suffices to check the relation for the orthogonal $G$-spectrum $\Sigma^\infty_+ G/H$
and the tautological class $e_H$ defined in \eqref{eq:tautological_e_H}.
Inspection of the definition in Construction \ref{con:transfer}
reveals that the transfer 
of the tautological class $\tr_H^G(e_H)$ in $\pi_0^G(\Sigma^\infty_+ G/H)$
is represented
by the $G$-map
\[ S^V \ \xra{\ c \ } \ G\ltimes_H S^W\ \xra{\ a \ }
\ S^V\sm G/H_+ \]
where $c$ is the collapse map based on any wide embedding
of $i:G/H\to V$ into a $G$-representation,
and $a[g,w]=g w\sm g H$.
So the class $\res^G_K(\tr_H^G(e_H))$ is represented by the underlying
$K$-map of the above composite, which is precisely the map $c_{G/H}$ 
for the underlying $K$-manifold of $G/H$. Part~(i) thus provides the relation
  \[ \Phi^K(\res^G_K(\tr_H^G(e_H)))\ = \ \Phi^K[c_{G/H}]\ = \ 
  \sum_{ M\in \pi_0( (G/H)^K) }\, \chi(M)\cdot \Phi^K(\sigma^K[M])\]
  where the sum runs over the connected components of $(G/H)^K$. 
  On the other hand, if $g\in G$ is such that $g H\in M\subset (G/H)^K$, then
  $K^g\leq H$ and
  \begin{align*}
    \sigma^K [M] \ &= \ g_\star (\sigma^{K^g}\td{e H})\
     = \ g_\star( \res^H_{K^g}(e_H))\ ;
  \end{align*}
  this proves the desired relation for the universal class $e_H$.
\end{proof}

Now we can proceed towards the double coset formula
for a transfer followed by a restriction.
The double coset formula was first proved by Feshbach 
for Borel cohomology theories \cite[Thm.\,II.11]{feshbach} 
and later generalized to equivariant cohomology theories by 
Lewis and May \cite[IV \S 6]{lms}. 
I am not aware of a published account of the double coset formula
in the context of orthogonal $G$-spectra.
Our method of proof for the double coset formula is different
from the approach of Feshbach and Lewis-May: 
we verify the double coset formula after passage to geometric fixed points 
for all closed subgroups. 
The detection criterion of Theorem \ref{thm:G-equivariant pi_0 of Sigma^infty}~(ii) 
then lets us deduce the  double coset formula for the universal example,
the tautological class $e_H$ in $\pi_0^H(\Sigma^\infty_+ G/H)$.

Before we can state the double coset formula we have to recall 
some additional concepts, such as orbit type submanifolds and 
the internal Euler characteristic.

\begin{defn}
Let $L$ be a closed subgroup of a compact Lie group $K$ and $B$ a $K$-space.
The {\em orbit type space} associated to the conjugacy class of $L$ is
the $K$-invariant subspace
\[ B_{(L)} \ = \
\{  x \in  B \ | \ \text{the stabilizer group of $x$ is conjugate to $L$}\}\ .\]
\end{defn}

The points in the space $B_{(L)}$ are said to have `orbit type' the conjugacy
class $(L)$.
We mostly care about orbit type subspaces for  equivariant smooth manifolds.
Two sources that collect general information about orbit type manifolds
are Chapters~IV and~VI of Bredon's book \cite{bredon-intro}
and Section~I.5 of tom Dieck's book \cite{tomDieck-transformation}.
If $B$ is a smooth $K$-manifold, then the subspace $B_{(L)}$
is a smooth submanifold of $B$, and locally closed in $B$, 
compare \cite[I Prop.\,5.12]{tomDieck-transformation}
or \cite[VI Cor.\,2.5]{bredon-intro}.\index{subject}{orbit type manifold}
Thus $B_{(L)}$ is called the {\em orbit type manifold}.
If the smooth $K$-manifold is compact, 
then it only has finitely many orbit types --
this was first shown by Yang \cite{yang-montgomery}; 
other references 
are \cite[IV Prop.\,1.2]{bredon-intro}, 
\cite[Thm.\,1.7.25]{palais-classification}, 
or \cite[I Thm.\,5.11]{tomDieck-transformation}.
The orbit type manifolds $B_{(L)}$ need not be closed inside $B$;
but if one orbit type manifold $B_{(L)}$ lies in the closure of another one
$B_{(L')}$, then $L'$ is subconjugate to $L$ in $K$, 
compare \cite[IV Thm.\,3.3]{bredon-intro}.
For every conjugacy class $(L)$, the quotient map $B_{(L)}\to K\bs B_{(L)}$
is a locally trivial smooth fiber bundle with fiber $K/L$ and structure group
the normalizer $N_K L$, compare \cite[VI Cor.\,2.5]{bredon-intro}.
Thus every path component of every orbit space $K\bs B_{(L)}$ of
an orbit type manifold is again a smooth manifold, 
in such a way that the quotient maps are smooth submersions.
In particular, if the action happens to have only one orbit type
(i.e., all stabilizer groups are conjugate), then the quotient space $K\bs B$ 
is a manifold and inherits a smooth structure from $B$.

The terms in the double coset formula will be indexed by path components
of the orbit type orbit manifolds $K\bs B_{(L)}$.
These manifolds need not be connected, and different components may
have varying dimensions, but each $K\bs B_{(L)}$
has only finitely many path components.
Indeed, by \cite[Cor.\,3.7 and 3.8]{verona-triangulation}
the orbit space $K\bs B$ admits a triangulation
such that the orbit type is constant on every open simplex.
So if a path component $M$ of $K\bs B_{(L)}$ intersects a particular open simplex,
then that simplex is entirely contained in $M$, i.e.,
$M$ is a union of open simplices in such a triangulation.
Since $K\bs B$ is compact, any triangulation has only finitely many simplices.

\begin{construction}[Internal Euler characteristic]
We continue to consider a compact Lie group $K$ acting smoothly on
a closed smooth manifold $B$. We let $L$ be a closed subgroup of $K$
and $M\subset K\bs B_{(L)}$ a connected component of
the orbit type manifold for the conjugacy class $(L)$.
We let $\bar M$ denote the closure of $M$ inside $K\bs B$
and $\delta M=\bar M-M$ the complement of $M$ inside its closure. 
Since $K\bs B$ is compact, so is $\bar M$. Since $K\bs B_{(L)}$
is locally closed in $K\bs B$, the set $M$ is open inside its closure $\bar M$.
So the complement $\delta M$ is closed in $\bar M$, hence compact.
One should beware that while $M$ is a topological manifold (without boundary,
and typically not compact),
$\bar M$ need not be a topological manifold with boundary.

By \cite[Cor.\,3.7 and 3.8]{verona-triangulation}
the orbit space $K\bs B$ admits a triangulation (necessarily finite)
such that the orbit type is constant on every open simplex.
So if the orbit type component $M$ intersects a particular open simplex,
then that simplex is entirely contained in $M$, i.e.,
$M$ is a union of open simplices in such a triangulation.
This means that $\bar M$ is obtained from $M$ by adding
some simplices of smaller dimension. 
Hence $M$ is a simplicial subcomplex of any such triangulation,
and $\delta M$ is a subcomplex of $M$.
In particular, $\bar M$ and $\delta M$ admit the structure of finite simplicial complexes
of dimension less than or equal to the dimension of $B$.
Thus the integral singular homology groups 
of $\bar M$ and $\delta M$ are finitely generated in every degree, and they
vanish above the dimension of $B$.
So the Euler characteristics of $\bar M$ and $\delta M$
are well-defined integers. 
The {\em internal Euler characteristic}\index{subject}{Euler characteristic!internal}
of $M$ is then defined as 
\[ \chi^\sharp(M)\ =\ \chi(\bar M)\ -\ \chi(\delta M)\ .\]
The internal Euler characteristic $\chi^\sharp(M)$ also
has an intrinsic interpretation that does not refer to the ambient space $K\bs B$.
Indeed since the integral homology groups of $\bar M$ and $\delta M$
are finitely generated and vanish for almost all degrees, 
the same is true for the relative singular homology groups $H_*(\bar M,\delta M;\mZ)$,
and the internal Euler characteristic satisfies
\[ \chi^\sharp(M)\ =\ \sum_{n\geq 0} (-1)^n\cdot\text{rank}(H_n(\bar M,\delta M;\mZ))
\ =\ \sum_{n\geq 0} (-1)^n\cdot\text{rank}(H^n(\bar M,\delta M;\mZ))\ .\]
The second equality uses that by the universal coefficient theorem,
we can use {\em co-}homology instead of homology to calculate the Euler characteristic.
Since $\delta M$ is a simplicial subcomplex of $\bar M$,
it is a neighborhood deformation retract inside $\bar M$.
So the relative cohomology group $H^n(\bar M,\delta M;\mZ)$
is isomorphic to the reduced cohomology group of the quotient space
$\bar M/\delta M$, by excision. This quotient space is homeomorphic
to $\hat M$, the one-point compactification of $M=\bar M-\delta M$.
Since $M$ is a topological manifold, it is locally compact, and so 
the cohomology groups of $\hat M$ are isomorphic 
to the compactly supported cohomology groups of $M$. 
Altogether, this provides isomorphisms
\[  H^n(\bar M,\delta M;\mZ)\ \iso \ H^n_c(M;\mZ)\ . \]
So we conclude that the internal Euler characteristic $\chi^\sharp(M)$ 
can also be defined as the compactly supported Euler characteristic of $M$, 
i.e., 
\[ \chi^\sharp(M)\ = \ \sum_{n\geq 0} (-1)^n\cdot\text{rank}(H^n_c(M;\mZ))\ . \]
\end{construction}

Now we let $K$ and $H$ be two closed subgroups of a compact Lie group $G$.
Then the homogeneous space $G/H$ is a smooth manifold
and the $K$-action on $G/H$ by left translations is smooth.
The {\em double coset space}\index{subject}{double coset space} 
$K\backslash G/H$\index{symbol}{$K\backslash G/H$ - {double coset space}} 
is the quotient space of $G/H$ by this $K$-action, 
i.e., the quotient of $G$ by the left $K$- and right $H$-action by translation.
In contrast to a homogeneous space $G/H$, the double coset space is in general
not a smooth manifold. However, the orbit type decomposition expresses
the double coset space as the union of certain subspaces that are
manifolds of varying dimensions.
Indeed, $G/H$ decomposes into orbit type manifolds with respect to $K$:
\[ K\bs G/H \ = \ \bigcup_{(L)} \, K\bs (G/H)_{(L)} \ ,\]
where the set-theoretic union is indexed by conjugacy classes of closed
subgroups of $L$; all except finitely many of the subspaces $K\bs (G/H)_{(L)}$ are empty.
The double coset formula expresses the composite $\res^G_K\circ\tr_H^G$
as a sum of terms, indexed by all connected components $M$ 
of orbit type orbit manifolds $K\bs (G/H)_{(L)}$.
The coefficient of the contribution of $M$ is the 
internal Euler characteristic $\chi^\sharp(M)$.

Our next result is a technical proposition that contains the essential
input to the double coset formula from equivariant differential topology.

\begin{prop}\label{prop:smooth Euler characteristic}
Let $K$ be a compact Lie group, $B$ a closed smooth $K$-manifold and
$L\leq J\leq K$ nested closed subgroups. 
Let $N$ be a connected component of the fixed point space $B^L$,
and $M$ a connected component of the orbit space $K\bs B_{(J)}$.
Let $b\in B$ be any point with stabilizer group $J$ and with $K b\in M$.
Let $W[M,N]$ be the preimage of $N\cap B_{(J)}$ under the map
\[   (K/J)^L \ \to \ (B_{(J)})^L \ ,\quad k J \ \longmapsto \ k b\ .\]
Then
\[ \chi^\sharp(N\cap B_{(J)}) \ = \  
\sum_{M\in \pi_0(  K\bs B_{(J)})}\chi^\sharp(M)\cdot \chi(W[M,N])\ .\]
\end{prop}
\begin{proof}
The projection
\[  B_{(J)} \ \xra{\ p\ } \  K\bs B_{(J)} \]
is a smooth fiber bundle with fiber $K/J$ and structure group $N_K J$, compare
\cite[Thm.\,1.7.35]{palais-classification} or \cite[VI Cor.\,2.5]{bredon-intro}.
Since $L$ is contained in $J$ the composite
\[ (B_{(J)})^L \ \xra{\text{incl}}\ B_{(J)} \ \xra{\ p\ } \  K\bs B_{(J)}  \]
is surjective and the restriction of $p$ to
$(B_{(J)})^L$ is another smooth fiber bundle
\[ (B_{(J)})^L \ \to \  K\bs B_{(J)} \ ,\]
now with fiber $(K/J)^L$, and the map
\[ \Psi_M\ : \ (K/J)^L \ \to \ ( B_{(J)})^L \ , \quad  k J \ \longmapsto \ k b \]
is the inclusion of the fiber over $K b \in M$.

For every connected component $M$ of the base $K\bs B_{(J)}$
the inverse image $p^{-1}(M)$ is open and closed in $(B_{(J)})^L$.
For every connected component $N$ of $B^L$
the intersection $N\cap B_{(J)}$ is open and closed in $(B_{(J)})^L$.
So the subset $p^{-1}(M)\cap N$ is open and closed in $(B_{(J)})^L$.
As $M$ varies, the subsets $p^{-1}(M)\cap N$ 
cover all of $N\cap B_{(J)}$.
Moreover, the restriction of the projection to a map
\[ p \ : \ p^{-1}(M)\cap N \ \to \ M  \]
is a smooth fiber bundle with connected base and with fiber $W[M,N]$,
by definition of the latter.
The internal Euler characteristic is multiplicative on smooth fiber bundles
with closed fiber, so
\[ \chi^\sharp(p^{-1}(M)\cap N)\ =\ \chi^\sharp(M)\cdot\chi(W[M,N])\ . \]
The internal Euler characteristic is additive on disjoint unions, so
\[ 
\chi^\sharp ( N\cap  B_{(J)}  )   \ =  \sum_{M\in\pi_0 ( K\bs B_{(J)})}
 \chi^\sharp ( N\cap p^{-1}(M) )   \ =   
\sum_{M\in\pi_0 ( K\bs B_{(J)})} \chi^\sharp (M)\cdot \chi(W[M,N])\ .\qedhere
 \]
\end{proof}

\begin{eg}\label{eg:chi(B) formula}
For later reference we point out a direct consequence of 
Proposition \ref{prop:smooth Euler characteristic}
for $L=e$, the trivial subgroup of $K$.
We let $J$ be any closed subgroup of $K$, and $M$ a connected component
of $K\bs B_{(J)}$.
The space $K/J$ is the disjoint union of its subspaces $W[M,N]$, 
as $N$ varies over the connected components of $B$.
So summing the formula of Proposition \ref{prop:smooth Euler characteristic}
over all components of $B$ yields
\[ \chi^\sharp(B_{(J)}) \ = \  
\sum_{M\in \pi_0(  K\bs B_{(J)})}\chi^\sharp(M)\cdot \chi( K/J )\ ,\]
by additivity of Euler characteristics for disjoint unions.
Because $\chi(B)$ is the sum of the internal Euler characteristics
$\chi^\sharp(B_{(J)})$ over all conjugacy classes of closed subgroups of $K$,
summing up over conjugacy classes gives the formula
\begin{equation}\label{eq:chi(B)}
 \chi(B) \ = \ \sum_{(J)\leq K}\ \sum_{M\in \pi_0(  K\bs B_{(J)})}\chi^\sharp(M)\cdot \chi( K/J )\ .  \end{equation}
\end{eg}

We let $K$ be a compact Lie group and $B$ a closed smooth $K$-manifold.
The Mostow-Palais embedding theorem \cite{mostow-embeddings, palais-imbedding}
provides a wide smooth $K$-equivariant embedding $i: B\to V$, 
for some $K$-representation $V$. 
The associated collapse map 
\[ c_B\ : \ S^V \ \to \  S^V\sm B_+  \]
was defined \eqref{eq:collaps B}.
The homotopy class of the map $c_B$ is then an element in
the equivariant homotopy group $\pi_0^K(\Sigma^\infty_+ B)$.
Theorem \ref{thm:G-equivariant pi_0 of Sigma^infty}~(i)
exhibits a basis of the  group $\pi_0^K(\Sigma^\infty_+ B)$,
and the next theorem expands the class of $c_B$ in that basis. 
As we shall explain in the proof of Theorem \ref{thm:double coset formula} below, 
the double coset formula for $\res^G_K\circ\tr_H^G$
is essentially the special case $B=G/H$ of the following theorem.
A compact smooth $K$-manifold only has finitely many orbit types,
so the sum occurring in the following theorem is finite.

\begin{theorem}\label{thm:pre double coset}
Let $K$ be a compact Lie group and $B$ a closed smooth $K$-manifold.
Then the relation
\[ [c_B] \ = \ \sum_{(J)\leq K}\sum_{M\in\pi_0(K\bs B_{(J)})}\ 
\ \chi^\sharp(M)\cdot \tr_J^K (\sigma^J \td{b_M} ) \]
holds in the group $\pi_0^K(\Sigma^\infty_+ B)$.
Here the sum runs over all connected components $M$ 
of all orbit type orbit manifolds $K\bs B_{(J)}$, 
and the element $b_M\in B$ that occurs is such that $K b_M\in M$
and with stabilizer group $J$, and $\td{b_M}$ is the path component of $b_M$
in the space $B^J$.
\end{theorem}
\begin{proof}
Because the desired relation lies in the 0-th equivariant homotopy group 
of a suspension spectrum, 
Theorem \ref{thm:G-equivariant pi_0 of Sigma^infty}~(ii) 
applies and shows that we only need to verify the formula
after taking geometric fixed points to any closed subgroup $L$ of $K$.
We let $J$ be another closed subgroup of $K$ containing $L$.
We calculate the contribution of the conjugacy class $(J)$
to the effect of $\Phi^L$ on the right hand side of the formula.
For every connected component $M$ of $K\bs B_{(J)}$ we choose an 
element $b_M\in B$ with stabilizer group $J$ and $K b_M\in M$.
Then
\begin{align*}
 \sum_{M \in \pi_0( K\bs B_{(J)})}\
&\chi^\sharp(M)\cdot \Phi^L (\res^K_L( \tr_J^K(\sigma^J\td{b_M})))\\ 
= \  & \sum_{M \in\pi_0( K\bs B_{(J)})}\  \sum_{W\in\pi_0( (K/J)^L)}\ 
\chi^\sharp(M)\cdot 
 \chi(W)\cdot \Phi^L ( k_\star(\res^J_{L^k}(\sigma^J\td{b_M}))) \\
= \  & \sum_{M \in\pi_0( K\bs B_{(J)})}\  \sum_{W\in\pi_0( (K/J)^L)}\ 
\chi^\sharp(M)\cdot 
 \chi(W)\cdot \Phi^L (\sigma^L(k_\star(\res^J_{L^k}\td{b_M}))) \\
= \ & \sum_{M \in\pi_0( K\bs B_{(J)})} \sum_{N\in \pi_0(B^L)}\
 \chi^\sharp(M)\cdot \chi(W[M,N])\cdot \Phi^L (\sigma^L[N]) \\
= \ & 
\sum_{N\in \pi_0(B^L)} \chi^\sharp ( N\cap  B_{(J)})  \cdot \Phi^L (\sigma^L[N]) \ .
\end{align*}
The first equation is Proposition \ref{prop:Phi^K after tr_H^G}~(ii),
applied to the closed subgroups $L$ and $J$ of $K$.
The element $k\in K$ is chosen so that $k J\in W\subset (K/J)^L$;
in particular this forces the relation $L\leq{^k J}\leq {^k\text{stab}(b_M)}$,
and thus $k b_M\in B^L$.
The third equation uses that for every $M$, 
\[ (K/J)^L \ = \ \coprod_{N\in \pi_0(B^L)}\, W[M,N] \ ,\]
that the Euler characteristic is additive on disjoint unions,
and that for all $k J\in W[M,N]$ the element $k b_M$ 
lies in the path component $N$ of $B^L$, by the very definition of $W[M,N]$.
The fourth equation is Proposition \ref{prop:smooth Euler characteristic}.

Now we sum up the contributions from the different conjugacy classes of
subgroups of $K$. 
The smooth $K$-manifold $B$ is stratified by the 
relatively closed smooth submanifolds $B_{(J)}$,
indexed over the poset of conjugacy classes of subgroup of $K$, 
where only finitely many orbit types occur.
Intersecting the strata with a connected component $N$ 
of the fixed point submanifold $B^L$
gives a stratification of $N$ by relatively closed
submanifolds, still indexed over the conjugacy classes of subgroup of $K$. 
The internal Euler characteristic is 
additive for such stratifications, i.e.,
\[ \chi( N ) \ = \ \sum_{(J)\leq K}\, \chi^\sharp( N\cap  B_{(J)} ) \ .\]
An application of~Proposition \ref{prop:Phi^K after tr_H^G}~(i) gives
\begin{align*}
 \Phi^L(\res^K_L [c_B] ) \ &= \
\sum_{N\in\pi_0(B^L)}  \chi(N)\cdot \Phi^L(\sigma^L[N]) \\  
&= \ 
\sum_{(J)\leq K} \sum_{N\in\pi_0(B^L)} 
\chi^\sharp( N\cap B_{(J)})\cdot \Phi^L(\sigma^L[N]) \\   
&= \ \sum_{(J)\leq K} 
 \sum_{M \in \pi_0( K\bs B_{(J)})}\
\chi^\sharp(M)\cdot \Phi^L(\res^K_L(\tr_J^K(\sigma^J\td{b_M}))) \\
&= \ \Phi^L\left(\res^K_L\left( \sum_{(J)\leq K} 
 \sum_{M \in \pi_0( K\bs B_{(J)})}\ \chi^\sharp(M)\cdot \tr_J^K(\sigma^J\td{b_M}) \right)\right)\ .
\end{align*}
This proves the desired formula after composition
with the geometric fixed point map to any closed subgroup of $K$.
\end{proof}

Now we have all ingredients for the statement and proof of the double
coset formula in place, which gives an expression for
the operation $\res^G_K\circ\tr_H^G$ on the equivariant homotopy
groups of an orthogonal $G$-spectrum $X$, as a sum indexed by the path components $M$
of the orbit type orbit manifolds $K\bs (G/H)_{(L)}$.
The operation associated to $M$ 
involves the transfer $\tr_{K\cap{^g H}}^K$, where $g\in G$ 
is an element such that $K g H\in M$;
since $K\cap{^g H}$ is the stabilizer of the coset $g H$, 
the group $K\cap{^g H}$ lies in the conjugacy class $(L)$.
The group $K\cap{^g H}$ may have infinite index in its normalizer in $K$, 
in which case $\tr_{K\cap{^g H}}^K=0$, and then the corresponding
summand in the double coset formula vanishes.

\begin{theorem}[Double coset formula]\label{thm:double coset formula}
Let $H$ and $K$ be closed subgroups\index{subject}{double coset formula}
of a compact Lie group $G$. Then for every orthogonal $G$-spectrum $X$ the relation
\[ \res^G_K\circ \tr_H^G \ = \ \sum_{M}\
\chi^\sharp(M)\cdot \tr_{K\cap{^g H}}^K  \circ g_\star \circ \res^H_{K^g\cap H}\]
holds as homomorphisms $\pi_0^H(X)\to\pi_0^K(X)$.
Here the sum runs over all connected components $M$ 
of all orbit type orbit manifolds $K\bs  (G/H)_{(L)}$, 
and $g\in G$ is an element such that $K g H\in M$.
\end{theorem}
\begin{proof}
Both sides of the double coset formula are natural transformations 
on the category of orthogonal $G$-spectra from the functor $\pi_0^H$
to the functor $\pi^K_0$. Since the functor $\pi_0^H$ is
represented by the suspension spectrum of $G/H$ 
(in the sense of Proposition \ref{prop:G/H represents pi_0^H}),
is suffices to check the relation for the orthogonal $G$-spectrum $\Sigma^\infty_+ G/H$
and the tautological class $e_H$ defined in \eqref{eq:tautological_e_H}.

Inspection of the definition in Construction \ref{con:transfer}
reveals that the transfer 
of the tautological class $\tr_H^G(e_H)$ in $\pi_0^G(\Sigma^\infty_+ G/H)$
is represented
by the $G$-map
\[ S^V \ \xra{\ c \ } \ G\ltimes_H S^W\ \xra{\ a \ }
\ S^V\sm G/H_+ \]
where $c$ is the collapse map based on any wide embedding
$i:G/H\to V$ into a $G$-representation, and $a[g,w]=g w\sm g H$.
So the class $\res^G_K(\tr_H^G(e_H))$ is represented by the underlying
$K$-map of the above composite, which is precisely the map $c_{G/H}$ 
for the underlying $K$-manifold of $G/H$.

Theorem \ref{thm:pre double coset} thus yields the formula
\[ \res^G_K(\tr_H^G(e_H))\ = \ [c_{G/H}] \ = \ \sum_{(J)\leq K}\sum_{M\in\pi_0(K\bs (G/H)_{(J)})}
\ \chi^\sharp(M)\cdot \tr_J^K (\sigma^J \td{g_M H} ) \]
in the group $\pi_0^K(\Sigma^\infty_+ G/H)$.
Here $g_M\in G$ is such that $K g_M H\in M$
and $K\cap {^{g_M} H}=J$. On the other hand,
\begin{align*}
  \sigma^J \td{g_M H} \ &= \ (g_M)_\star (\sigma^{J^{g_M}}\td{e H})\\
&= \ (g_M)_\star( \res^H_{J^{g_M}}(\sigma^H\td{e H}))\
= \ (g_M)_\star( \res^H_{K^{g_M}\cap H}(e_H))\ ,
\end{align*}
so this proves the double coset formula for the universal class $e_H$.
\end{proof}

\begin{eg}\label{eg:restrict transfer to e} 
A special case of the double coset formula is when $K=e$,
i.e., when we restrict a transfer all the way to the trivial subgroup of $G$.
In this case the orbit type component manifolds are the path components
of the coset space $G/H$, and the sum in the double coset
formula is indexed by these.
Since all path components of $G/H$ are homeomorphic, 
they have the same Euler characteristic $\chi(G/H)/|\pi_0(G/H)|$, 
so the double coset formula specializes to
\begin{align*}
   \res^G_e\circ \tr_H^G \ &= \ 
 \sum_{ M\in \pi_0 (G/H) }\, \chi(M)\cdot  g_\star\circ\res^H_e \\
&= \ \chi(G/H)/|\pi_0(G/H)| \cdot \sum_{ [g H]\in \pi_0 (G/H) }\,  g_\star\circ\res^H_e \ .
\end{align*}
In the special case where the conjugation action of $\pi_0(G)$ on $\pi_0^e(X)$
happens to be trivial (for example for all global stable homotopy types),
all the maps $g_\star$ are the identity and the formula simplifies to
\[ \res^G_e\circ \tr_H^G \ = \ \chi( G/H)\cdot \res^H_e \ .  \]
\end{eg}

\begin{eg}[Double coset formula for finite index transfers]\label{eg:double coset finite index}
If $H$ has finite index in $G$,\index{subject}{double coset formula!finite index} 
then the double coset formula simplifies.
In this situation the homogeneous space $G/H$, and
hence also the double coset space $K\bs G/H$, is discrete,
all orbit type manifold components are points, and all internal Euler characteristics
that occur in the double coset formula are~1.
Since the intersection $K\cap {^g H}$ also has finite index in $K$, 
only finite index transfers show up in the double coset formula, 
which specializes to the relation
\[ \res^G_K\circ \tr_H^G \ = \ \sum_{[g]\in K\backslash G/H}\
\tr_{K\cap{^g H}}^K \circ g_\star\circ \res^H_{K^g\cap  H} \ .\]
\end{eg}

\begin{eg}
We calculate the double coset formula for the maximal torus 
\[ H\ =\ 
\left\lbrace \begin{pmatrix} \lambda & 0\\ 0 & \mu\end{pmatrix}\ : \ 
\lambda,\mu\in U(1) \right\rbrace \]
of $G=U(2)$.
We take $K=N_{U(2)} H$, the normalizer of the maximal torus; 
this is Example~VI.2 in \cite{feshbach}.\index{subject}{unitary group!$U(2)$}
Then $K$ is isomorphic to $\Sigma_2\ltimes H$, the semidirect product 
with $\Sigma_2$ permuting the two diagonal entries of matrices in $H$;
the group $K$ is generated by $H$ and the involution 
$(\begin{smallmatrix}0&1\\1&0\end{smallmatrix})$.

We calculate the double coset space $K \bs U(2) / H$
by identifying the space $H\bs U(2) / H$ 
and the residual action of the symmetric group $\Sigma_2$ on this space.
A homeomorphism from $H\bs U(2) / H$ to the unit interval $[0,1]$ is induced by
\[ h \ : \ U(2)  \ \to \ [0,1] \ , \quad
 A \ \longmapsto \ |a_{11}|^2 \ ,\]
the square of the length of the upper left entry $a_{11}$ of $A$.
Lengths are non-negative and every column of a unitary matrix is a unit vector,
so the map really lands in the interval $[0,1]$.
The number $h(A)$ only depends on the
double coset of the matrix $A$, so the map $h$ factors over
the double coset space $H\bs U(2) / H$.
For every $x\in [0,1]$ 
the vector $(\sqrt{x},\sqrt{1-x})\in \mC^2$ 
has length~1, so it can be complemented to an orthonormal basis,
and so it occurs as the first column of a unitary matrix;
the map $h$ is thus surjective.
On the other hand, if $h(A)=h(B)$ for two unitary matrices $A$ and $B$,
then left multiplication by an element in $H$ makes the first row
of $B$ equal to the first row of $A$. Right multiplication by
an element of $1\times U(1)$ then makes the matrices equal.
So $A$ and $B$ represent the same element in the double coset space.
The induced map
\[ \bar h \ : \ H \bs U(2) / H \ \to \ [0,1] \]
is thus bijective, hence a homeomorphism.
The action of $\Sigma_2$ on the orbit space
$H \bs U(2) / H$ permutes the two rows of a matrix; 
under the homeomorphism $\bar h$, this action thus 
corresponds to the involution of $[0,1]$ sending $x$ to $1-x$.
Altogether this specifies a homeomorphism 
from the double coset space to the interval $[0,1/2]$
that sends a double coset $K\cdot A\cdot H$
to the minimum of $|a_{11}|^2$ and $|a_{21}|^2$.
The inverse homeomorphism is 
\[ g \ : \ [0, 1/2]\ \iso \ K \bs U(2) / H\ , \quad
g(t)\ = \  K\cdot 
\begin{pmatrix} \sqrt{1-t} & \sqrt{t}\\ -\sqrt{t}& \sqrt{1-t}\end{pmatrix}\cdot H\ .\]
The orbit type decomposition is as
\[ \{0\}\cup (0,1/2)\cup \{1/2\} \ .\]
So as representatives of the orbit types we can choose
\begin{align*}
 g(0)& = \begin{pmatrix} 1&0\\0&1\end{pmatrix} \ ,\quad
 g(1/5)= \frac{1}{\sqrt{5}} \begin{pmatrix} 2&1\\-1&2\end{pmatrix} 
\text{\quad and\quad}
 g(1/2)= \frac{1}{\sqrt{2}}\begin{pmatrix} 1&1\\-1&1\end{pmatrix} \ .
\end{align*}
For the intersections (which are representatives of the conjugacy classes
of those subgroups of $K$ with non-empty orbit type manifolds) we get
\begin{alignat*}{3}
 & K &&\cap\  {^{g(0)} H}\ &&= \ H \ , \\  
 & K &&\cap {^{g(1/5)} H}\ &&= \ \Delta\text{ , and}\\  
 & K &&\cap {^{g(1/2)} H}\ &&= \ \Sigma_2\times \Delta\ .
\end{alignat*}
Here $\Delta=\left\lbrace 
(\begin{smallmatrix} \lambda & 0\\ 0 & \lambda\end{smallmatrix}) \ : \ 
\lambda\in U(1) \right\rbrace$
is the diagonal copy of $U(1)$.
The group $\Delta$ is normal in $K$, so 
it has infinite index in its normalizer, and the corresponding 
transfer does not contribute to the double coset formula.
The internal Euler characteristic of a point is~1, so
the double coset formula has two non-trivial summands, and it specializes to 
\[ 
 \res^{U(2)}_K\ \circ \ \tr_H^{U(2)}\ = \  
\tr_H^K\ + \ \tr_{\Sigma_2\times\Delta}^K \circ \gamma_\star\circ \res^H_{(\Sigma_2\times\Delta)^\gamma}
 \ ,  \]
where $\gamma=g(1/2)=\frac{1}{\sqrt{2}}(\begin{smallmatrix} \phantom{-}1 & 1\\ -1 & 1\end{smallmatrix})$.
\end{eg}

\begin{eg}\label{eg:double coset for Brauer} 
For later use in the `explicit Brauer induction'
(see Remark \ref{rk:Brauer induction} below), 
we work out another double coset formula for
the unitary group $G=U(n)$.
We write $U(k,n-k)$ for the subgroup of those elements of $U(n)$ 
that leave the subspaces $\mC^k\oplus 0$
and $0\oplus\mC^{n-k}$ invariant. We also use the analogous notation for more than two
factors. As subgroups we take
\[  H\ =\ U(1,n-1) \text{\qquad and\qquad} 
K \ = \ U(k,n-k)\ ,\]
where $1\leq k\leq n-1$.\index{subject}{unitary group}
Again the relevant double coset space `is' an interval.
Indeed, the continuous map
\[ U(n)  \ \to \ [0,1]  \ , \quad
 A \ \longmapsto \ |a_{1 1}|^2 + \dots +  |a_{k 1}|^2 \ ,\]
the partial length of the entries in the first column of $A$,
is invariant under right multiplication by elements of $U(1,n-1)$,
and under left multiplication by elements of $U(k,n-k)$.
Every column of a unitary matrix is a unit vector,
so the map really lands in the interval $[0,1]$.
The map factors over a continuous surjective map
\[ h \ : \ U(k,n-k) \bs U(n) / U(1,n-1)  \ \to \ [0,1] \ . \]
On the other hand, if $h[A]=h[B]$ for two unitary matrices $A$ and $B$,
then left multiplication by an element in $U(k,n-k)$ makes the first
columns of $A$ and $B$ equal.
Right multiplication by an element of $U(1,n-1)$ then makes the matrices equal.
So $A$ and $B$ represent the same element in the double coset space.
This shows that the map $h$ is bijective, hence a homeomorphism.
The inverse homeomorphism is 
\[ [0, 1]\ \iso \ U(k,n-k) \bs U(n) / U(1,n-1) \ ,\quad
t \ \mapsto\ U(k,n-k)\cdot g(t) \cdot U(1,n-1)   \]
with
\[  g(t)\ = \  
\begin{pmatrix} 
\sqrt{t} & 0 & \cdots & 0 & \sqrt{1-t}\\ 
0 & 1 & \cdots & 0 & 0\\ 
\vdots & \vdots  & \ddots & \vdots & \vdots\\ 
0 & 0 & \cdots & 1 & 0\\ 
-\sqrt{1-t} & 0 & \cdots & 0 & \sqrt{t}\end{pmatrix}
 \ . \]
The orbit type decomposition is as $ \{0\}\cup (0,1)\cup \{1\}$ .
As representatives of the orbit types we can choose 
\[ 
g(0)= \begin{pmatrix} 
0 & 0 & \cdots & 0 & 1\\ 
0 & 1 & \cdots & 0 & 0\\ 
\vdots & \vdots  & \ddots & \vdots & \vdots\\ 
0 & 0 & \cdots & 1 & 0\\ 
-1 & 0 & \cdots & 0 & 0\end{pmatrix}
\ , \quad
g(1/2)= \frac{1}{\sqrt{2}}\begin{pmatrix} 
1 & 0 & \cdots & 0 & 1\\ 
0 & 1 & \cdots & 0 & 0\\ 
\vdots & \vdots  & \ddots & \vdots & \vdots\\ 
0 & 0 & \cdots & 1 & 0\\ 
-1 & 0 & \cdots & 0 & 1\end{pmatrix}
 \]
and the identity $g(1)$.
For the intersections (which are representatives of the conjugacy classes
of those subgroups of $U(k,n-k)$ with non-empty orbit type manifolds) we get
\begin{alignat*}{3}
 &U(k,n-k)\, &&\cap\quad  {^{g(0)} U(1,n-1)}\ 
&&= \ U(k, n-k-1, 1)\\  
 &U(k,n-k)\, &&\cap\, {^{g(1/2)} U(1,n-1)}\ &&= \ \Delta \\  
 &U(k,n-k)\, &&\cap\quad  {^{g(1)} U(1,n-1)}\ 
&&= \ U(1, k-1, n-k)\ .
\end{alignat*}
Here $\Delta$ is the subgroup of $U(n)$ consisting of matrices of the form
\[  \left(\begin{matrix} 
\lambda & 0 & 0 & 0\\ 
    0 & A & 0 & 0\\ 
    0 & 0 & B & 0\\ 
    0 & 0 & 0 & \lambda \end{matrix}\right)\]
for $(\lambda,A,B)\in U(1)\times U(k-1)\times U(n-k-1)$.
The subgroup $\Delta$ is normalized by $U(1,\, k-1,\, n-k-1,\, 1)$,
so it has infinite index in its normalizer, and the corresponding 
transfer does not contribute to the double coset formula.
The double coset formula thus
has two non-trivial summands, and comes out as
\begin{align}\label{eq:U(n)_Brauer_double_coset}
 \res^{U(n)}_{U(k,n-k)}\ \circ \ \tr_{U(1,n-1)}^{U(n)}\ &= \\  
\tr_{U(1,k-1,n-k)}^{U(k,n-k)}&\circ\res^{U(1,n-1)}_{U(1,k-1,n-k)}\  +\  
\tr_{U(k,n-k-1,1)}^{U(k,n-k)}\circ g(0)_\star\circ\res^{U(1,n-1)}_{U(1,k,n-k-1)}\ .\nonumber
\end{align}
\end{eg}

\index{subject}{Mackey functor|(}
For the rest of this section we restrict our attention to finite groups.
For easier reference we record that the equivariant homotopy groups
of an orthogonal $G$-spectrum form a Mackey functor.
We first recall one of several equivalent definitions of this concept.

\begin{defn}\label{def:G-Mackey}\index{subject}{Mackey functor}
Let $G$ be a finite group. A {\em $G$-Mackey functor}\index{subject}{G-Mackey functor@{$G$-Mackey functor}|see{Mackey functor}} 
$M$ consists of the following data:
\begin{itemize}
\item an abelian group $M(H)$ for every subgroup $H$ of $G$,
\item conjugation homomorphisms $g_\star:M(H^g)\to M(H)$ for all $H\leq G$ 
and $g\in G$,\index{subject}{conjugation homomorphism!in a Mackey functor}  
\item restriction homomorphisms $\res^H_K:M(H)\to M(K)$ for all $K\leq H\leq G$, and\index{subject}{restriction homomorphism!in a Mackey functor}  
\item transfer homomorphisms $\tr^H_K:M(K)\to M(H)$ for all $K\leq H\leq G$.\index{subject}{transfer!in a Mackey functor}  
\end{itemize}
This data has to satisfy the following conditions:
\begin{enumerate}[(i)]
\item  (Unit conditions)
\[ \res_H^H \ = \ \tr_H^H\ = \ \Id_{M(H)} \]
for all subgroups $H$, and $h_\star =\Id_{M(H)} $ for all $h\in H$.
\item (Transitivity conditions)
\[ 
 \res^K_L\circ \res^H_K \ = \ \res^H_L 
\text{\quad and \quad}
\tr^H_K\circ\tr_L^K \ = \ \tr^H_L 
 \]
for all $L\leq K\leq H\leq G$.
\item(Interaction conditions)
\[ g_\star\circ g'_\star \ = \ (g g')_\star \ , \quad
\tr^H_K\circ g_\star \ = \ g_\star\circ \tr^{H^g}_{K^g}
 \text{\quad and \quad}
\res^H_K\circ g_\star \ = \ g_\star \circ\res^{H^g}_{K^g} \]
for all $g,g'\in G$ and $K\leq H\leq G$.
\item (Double coset formula)
for every pair of subgroups $K, L$ of $H$ the relation
\[  
\res^H_L\circ \tr_K^H  \ = \  \sum_{[h]\in L\backslash H/K} 
\ \tr^L_{L\cap ^h K}\circ h_\star \circ \res^K_{L^h\cap K}
\]
holds as maps $M(K)\to M(L)$;
here $[h]$ runs over a set of representatives for the
$L$-$K$-double cosets.
\end{enumerate}
A {\em morphism} of $G$-Mackey functors\index{subject}{morphism!of $G$-Mackey functors}
 $f:M\to N$ is a collection
of group homomorphisms $f(H):M(H)\to N(H)$, for all subgroups $H$ of $G$,
such that
\[ \res^H_K\circ f(H) \ =\  f(K)\circ \res^H_K\text{\quad and\quad}   
 \tr^H_K\circ f(K) \ =\ f(H)\circ \tr^H_K\]
for all $K\leq H\leq G$, and
\[ g_\star\circ f(H^g) \ =\  f(H)\circ g_\star  \]
for all $g\in G$.
\end{defn}

We denote the category of $G$-Mackey functors by $G\Mack$.\index{symbol}{  $G\Mack$ - {category of $G$-Mackey functors}} 
To my knowledge, the concept of a Mackey functor goes back to Dress \cite{dress-notes}
and Green \cite{green-axiomatic}.
Green in fact defined what is nowadays called a {\em Green functor}\index{subject}{Green functor}
(which he called {\em $G$-functor} in \cite[Def.\,1.3]{green-axiomatic}),
which amounts to $G$-Mackey functor whose values underlie commutative rings,
where restriction and conjugation maps are ring homomorphisms,
and where the transfer maps satisfy Frobenius reciprocity.
Definition \ref{def:G-Mackey} is the `down to earth' definition
of Mackey functor, and the one that is most useful for concrete calculations.
There are two alternative (and equivalent) definitions that are often used:
\begin{itemize}
\item As a pair $(M_*,M^*)$ of additive functors on the category of finite $G$-sets,
where $M_*$ is covariant and $M^*$ is contravariant;
the two functors must agree on objects
and for every pullback diagram of finite $G$-sets
    \[ \xymatrix{ A \ar[r]^-f \ar[d]_g & B \ar[d]^h\\
      C \ar[r]_-k & D } \]
the relation
\[ M_*(f)\circ M^*(g) \ = \ M^*(h)\circ M_*(k)\ : \ M(C)\ \to \ M(B) \]
holds. This is (a special case of) the definition introduced by Dress \cite[\S 6]{dress-notes}.
\item As an additive functor on the category of spans of finite $G$-sets;
the equivalence with Dress' definition is due to Lindner \cite[Thm.\,4]{lindner-Mackey}.
\end{itemize}

Our main source of examples of $G$-Mackey functors comes from 
orthogonal $G$-spectra: as we verified 
in Sections \ref{sec:equivariant homotopy groups}
and \ref{sec:Wirthmuller and transfer}
and Theorem \ref{thm:double coset formula},
the restriction, conjugation and transfer maps 
make the homotopy groups $\pi_k^H(X)$ for varying subgroups $H$
into a $G$-Mackey functor.

\smallskip

Proposition \ref{prop:Phi of transfer} above
shows that the geometric fixed point map \eqref{eq:geometric_fix_map}
factors over the quotient of $\pi_0^G(X)$ by the subgroup
generated by proper transfers. This should serve as motivation for the
following algebraic interlude about Mackey functors for finite groups,
where we study the process of `dividing out transfers' systematically.

\begin{construction}
We let $G$ be a finite group and $M$ a $G$-Mackey functor.
For a subgroup $H$ of $G$ we let  $t_H M$ be the subgroup of $M(H)$ 
generated by transfers from proper subgroups of $H$, and we define
\[ \tau_H M \ = \ M(H) / t_H M\ .\]
For $g\in G$, the conjugation isomorphism
$g_\star: M(H^g)\to M(H)$ descends to a homomorphism 
\begin{equation}  \label{eq:act_on_prod_tau_H}
 g_\star\ : \ \tau_{H^g} M\ \to \ \tau_H  M \ .  
\end{equation}
Moreover, the transitivity 
relation $g_\star\circ g'_\star=(g g')_\star$
still holds. In particular, the action of the Weyl group $W_G H$ on $M(H)$
descends to a $W_G H$-action on the quotient group $\tau_H M$.

Now we recall how a $G$-Mackey functor $M$ can rationally be
recovered from the $W_G H$-modules $\tau_H M$ for all subgroups $H$ of $G$.
We let
\[ \bar\psi^M_G \ : \ M(G)\ \to \ {\prod}_{H\leq G}\, \tau_H M  \]
be the homomorphism whose $H$-component is the composite
\[ M(G)\ \xra{\res^G_H} \ M(H)\ \xra{\text{proj}} \ \tau_H M  \ ;\]
here the product is indexed over all subgroups $H$ of $G$.
The group $G$ acts on the product via the maps \eqref{eq:act_on_prod_tau_H},
permuting the factors within conjugacy classes.
Since inner automorphisms induce the identity we have
$g_\star\circ\res^G_{H^g}=\res^G_H\circ g_\star=\res^G_H$, 
so the map $\bar\psi^M_G$ lands in the subgroup of $G$-invariant tuples.
We let
\[ \psi^M_G \ : \ M(G)\ \to \ \left( {\prod}_{H\leq G}\, \tau_H M\right)^G  \]
denote the same map as $\bar\psi^M_G$, but now with image the subgroups
of $G$-invariants.
While the above description of the target of $\psi^M_G$ is the most natural one,
it is somewhat redundant.
Indeed, for every $H\leq G$, the image of the restriction map 
$\res^G_H:M(G)\to M(H)$ lands in the subgroup
$M(H)^{W_G H}$ of invariants under the Weyl group $W_G H$.
Moreover, if we choose representatives of the conjugacy 
classes of subgroups, then projection from the full product
(over all subgroups of $G$) to the product indexed by the representatives restricts
to an isomorphism
\[ \left( {\prod}_{H\leq G}\, \tau_H M\right)^G \ \xra{\ \iso \ }\ 
 {\prod}_{(H)}\, \left(\tau_H M\right)^{W_G H}  \ .\]
For explicit calculations, the second description of the target of $\psi^M_G$
is often more convenient.
\end{construction}

For a finite group $G$ we set
\[ d_G\ = \  {\prod}_{(H)} \, |W_G H|\ , \]
the product, over all conjugacy classes of subgroups of $G$,
of the orders of the respective Weyl groups.
Since the order of $W_G H$ divides the order of $G$,
inverting $|G|$ also inverts the number $d_G$.
The following result is well known;
the earliest reference that I am aware of is \cite[Cor.\,4.4]{thevenaz-some results}.

\begin{prop}\label{prop:recover G-Mackey rationally}\index{subject}{Mackey functor}
For every finite group $G$ and every $G$-Mackey functor $M$,
the kernel and cokernel of the homomorphism $\psi^M_G$ are annihilated by $d_G$.
In particular, $\psi^M_G$ becomes an isomorphism after inverting the order of $G$.
\end{prop}
\begin{proof}
We reproduce the proof given in \cite{thevenaz-some results}.
We choose representatives for the conjugacy classes of subgroups and number them
by non-decreasing order
\[ e = H_1 \ , H_2, \dots, H_n = G \ .\]
Then
\begin{itemize}
\item for all $i<j$ and $g\in G$ the group $H_j\cap{^g H_i}$
is a proper subgroup of $H_j$, and
\item  every proper subgroup of $H_i$ is conjugate to one of the groups
occurring before $H_i$.
\end{itemize}
We set
\[ K_j \ = \ \ker(\psi^M_G)\, \cap\, {\sum}_{i=1}^j\, \tr_{H_i}^G(M(H_i)) \ \subseteq 
 \ M(G)\ ;\]
this defines a nested sequence of subgroups
\[ 0 \ = \ K_0\  \ \subseteq \
K_1\  \ \subseteq \
K_2\  \ \subseteq \ \dots \ 
K_{n-1}\  \ \subseteq \ K_n \ = \ \ker(\psi^M_G)\ .
 \]
We show that
\begin{equation}  \label{eq:degenerate_step}
 |W_G H_j|\cdot K_j \ \subseteq \ K_{j-1}    
\end{equation}
for all $j=1,\dots,n$.
Altogether this means that
\[ d_G\cdot K_n \ = \  {\prod}_{i=1}^n |W_G H_i|\cdot K_n \ \subseteq\ K_0 \ = \ 0 \ . \]
In the course of proving \eqref{eq:degenerate_step}
we call an element of $M(H)$ {\em degenerate} if it maps to
zero in $\tau_H M$, i.e., if it is a sum of transfers from proper
subgroups of $H$. So the kernel of $\psi^M_G$ is precisely the
subgroup of those elements of $M(G)$ that restrict to degenerate
elements on all subgroups.

We write any given element $x\in K_j$ as
\[ x \ = \ \tr_{H_j}^G(y) \ + \ \bar x \]
for suitable $y\in M(H_j)$ and with $\bar x$ a sum of transfers from the groups
$H_1,\dots,H_{j-1}$. 
For $i=1,\dots, j-1$ the double coset formula for $\res^G_{H_j}\circ \tr^G_{H_i}$
lets us write $\res^G_{H_j}(\bar x)$ as a sum of transfers from proper subgroups
of $H_j$; this uses the hypothesis on the enumeration of the subgroups.
So the class $\res^G_{H_j}(\tr_{H_i}^G(\bar x))$ is degenerate.
Since $x$ is in the kernel of $\psi_G^M$, the class $\res^G_{H_j}(x)$ is also degenerate.
So the class
\[ \res^G_{H_j}(\tr^G_{H_j}(y))\ = \ \res^G_{H_j}(x)\  - \res^G_{H_j}(\bar x)  \]
is degenerate.
The double coset formula expresses 
$\res^G_{H_j}(\tr_{H_j}^G(y))$ as the sum of the $W_G H_j$-conjugates
of $y$, plus transfers from proper subgroups of $H_j$. So the element
\[  \sum_{g H_j\in W_G H_j} \ g_\star(y)\ \in \ M(H_j)\]
is degenerate and hence the element
\[ |W_G H_j|\cdot \tr_{H_j}^G(y)\ = \ 
\tr_{H_j}^G\left( \sum_{g H_j\in W_G H_j}  g_\star(y) \right) \ \in \ M(G)\]
is a sum of transfers from proper subgroups of $H_j$.
Every proper subgroup of $H_j$ is conjugate to one of the group $H_1,\dots,H_{j-1}$,
so
\[ |W_G H_j|\cdot x \ = \ 
|W_G H_j|\cdot \tr_{H_j}^G(y)\ + \ |W_G H_j|\cdot \bar x\]
is both in the kernel of $\psi^M_G$ and
a sum of transfer from the groups $H_1,\dots,H_{j-1}$.
This proves the claim that $|W_G H_j|\cdot x$ belongs to $K_{j-1}$.
Altogether this finishes the proof that the kernel of $\psi^M_G$
is annihilated by the number $d_G$.

Now we show that the cokernel of $\psi^M_G$ is annihilated by $d_G$.
We let $I_j$ denote the subgroup of $\prod_{i=1}^n (\tau_{H_i} M)^{W_G H_i}$
consisting of those tuples
\[ x \ = \ (x_i)_{1\leq i\leq n} \]
such that $x_{j+1}=x_{j+2}=\dots=x_n=0$.
This defines a nested sequence
\[ 0 \ = \ I_0\  \ \subseteq \
I_1\  \ \subseteq \
I_2\  \ \subseteq \ \dots \ 
I_{n-1}\  \ \subseteq \ I_n \ = \ \prod_{i=1,\dots, n} (\tau_{H_i} M)^{W_G H_i}\ . \]
We show that
\begin{equation}  \label{eq:image_step}
 |W_G H_j|\cdot I_j \ \subseteq \ \text{Im}(\psi^M_G) \ +\ I_{j-1} 
\end{equation}
for all $j=1,\dots,n$.
Altogether this means that
\[ d_G\cdot I_n \ \subseteq \  \text{Im}(\psi^M_G) \  , \]
i.e., the cokernel of $\psi^M_G$ is annihilated by $d_G$.

To prove \eqref{eq:image_step}
we consider a tuple $x=(x_n)$ in $I_j$ and choose a representative
$y\in M(H_j)$ for the `last' non-zero component, i.e., such that
$y$ maps to $x_j$ in $\tau_{H_j}M$.
Since $x_j$ is invariant under the action of the Weyl group $W_G H_j$, the element $y$
is at least $W_G H_j$-invariant modulo transfers from proper subgroups of $H_j$.
The double coset formula for $\res^G_{H_j}\circ\tr_{H_j}^G$ thus
gives that
\[ \res^G_{H_j}(\tr_{H_j}^G(y)) \ \equiv \ \sum_{g H_j\in W_G H_j} g_\star(y) \ \equiv \ 
|W_G H_j|\cdot y \ ,\]
with both congruences modulo transfers from proper subgroups of $H_j$.
So the composite
\[ M(G)\ \xra{\res^G_{H_j}} \ M(H_j)\ \xra{\text{proj}}\ \tau_{H_j}M\]
takes $\tr^G_{H_j}(y)$ to $|W_G H_j|\cdot x_j$.
For $i>j$ the double coset formula
shows that $\res^G_{H_i}(\tr_{H_j}^G(y))$ is a sum of transfers
from proper subgroups of $H_i$. So the element
$\res^G_{H_i}(\tr_{H_j}^G(y))$ maps to~0 in $\tau_{H_i} M$ for all $i=j+1,\dots,n$.
In other words, the tuple $\psi^M_G(\tr_{H_j}^G(y))$ belongs to $I_j$.
Since $\psi^M_G(\tr_{H_j}^G(y))$ and $|W_G H_j|\cdot x$ both belong to $I_j$
and agree at the component of $H_j$, we can thus conclude that 
\[ \psi^M_G(\tr_{H_j}^G(y)) - |W_G H_j|\cdot x \ \in \ I_{j-1} \ .\]
This proves \eqref{eq:image_step} and finishes the proof 
that the cokernel of $\psi^M_G$ is annihilated by the number $d_G$.
\end{proof}

The previous Proposition \ref{prop:recover G-Mackey rationally}
shows that for every $G$-Mackey functor $M$
the group $\mZ[1/|G|]\tensor (\tau_G M)$ is a direct summand of the
group $\mZ[1/|G|]\tensor M(G)$, and the splitting is natural for
morphisms of $G$-Mackey functors. So the functor
\[ G\Mack \ \to \ \Ab\  , \quad
M\longmapsto\ \mZ[1/|G|]\tensor (\tau_G M)\] 
is a natural direct summand of an exact functor. So we can conclude:

\begin{cor}\label{cor:tau_G is exact}
Let $G$ be a finite group. The functor that assigns to a $G$-Mackey functor $M$
the abelian group $\mZ[1/|G|]\tensor (\tau_G M)$ is exact.  
\end{cor}

Now we will show that after inverting the order of the finite group $G$,
the category of $G$-Mackey functors splits as a product of categories,
indexed by the conjugacy classes of subgroups of $G$; 
the factor corresponding to $H\leq G$ is the category of $W_G H$-modules,
i.e., abelian groups with an action of Weyl group of $G$
(also known as modules over the group ring $\mZ[W_G H]$).
The argument is not particularly difficult, but somewhat lengthy
because the details involve a substantial amount of book keeping.
The following result is folklore, but I do not know a reference that 
states it in precisely this form. 
The discussion in Appendix~A of \cite{greenlees-may-Tate} is closely related,
but the exposition of Greenlees and May differs from ours 
in that it emphasizes idempotents in the rationalized Burnside ring, 
as opposed to dividing out transfers.
One could deduce part~(ii) of the following theorem
from the results in \cite[App.\,A]{greenlees-may-Tate}
by using the action of the Burnside ring $A(H)$ on the value
$M(H)$ of any $G$-Mackey functor $M$, and showing that after inverting the group order,
the functor $\tau_H:G\Mack\to W_G H \mo$ becomes isomorphic to
$e_H\cdot M(H)$, where $e_H\in A(H)[1/|G|]$ is the idempotent `supported at $H$'.

The special case $Q=\mQ$ of the following theorem
in particular shows that rational $G$-Mackey functors form a semisimple abelian category, 
i.e., every object is both projective and injective. 

\begin{theorem}\label{thm:split Mackey functors}
Let $G$ be a finite group.
\begin{enumerate}[\em (i)]
\item The functor
\[  \tau \ = \ (\tau_H)_{(H)}\ : \  G\Mack \ \to \ 
\prod_{(H)}\, W_G H\mo \]
has a right adjoint $\rho:\prod_{(H)}\, W_G H\mo \to  G\Mack$.
\item Let $Q$ be a subring of the ring of rational numbers
in which the order of $G$ is invertible.
Then the adjoint functors $(\tau,\rho)$ 
restrict to inverse equivalences between the category
$Q$-local $G$-Mackey functors and the product of the categories
of $Q$-local $W_G H$-modules.
\end{enumerate}
\end{theorem}
\begin{proof}
(i) We start by exhibiting a right adjoint $\rho_H:W_G H\mo\to G\Mack$
to the functor $\tau_H$.
We recall that the Weyl group $W_G H$ acts
on the $H$-fixed points of every $G$-set $S$ by
\[ W_G H\times S^H \ \to \ S^H \ , \quad (g H, s) \ \longmapsto \ g s \ .\]
We let $N$ be a $W_G H$-module and $K$ another subgroup of $G$. We set
\[ (\rho_H N)(K)\ = \ \map^{W_G H}((G/K)^H,N)\ , \]
the abelian group of $W_G H$-equivariant functions from $(G/K)^H$ to $N$,
under pointwise addition of functions.
The conjugation map
$\gamma_\star:(\rho_H N)(K^\gamma)\to (\rho_H N)(K)$,
for $\gamma\in G$, is precomposition with the $W_G H$-map
\[ (G/K)^H \ \to \ (G/K^\gamma)^H \ , \quad 
g\cdot K \ \longmapsto \ g \gamma \cdot K^\gamma\ . \]
For $L\leq K$, the restriction map
$\res^K_L:(\rho_H N)(K)\to (\rho_H N)(L)$ is precomposition with the $W_G H$-map
\[ (G/L)^H \ \to \ (G/K)^H \ , \quad g L\ \longmapsto \ g K\ . \]
The transfer map
$\tr^K_L:(\rho_H N)(L)\to (\rho_H N)(K)$ is given by `summation over preimages', i.e.,
\[ \tr_L^K(f)(g K) \ =  \ 
\sum_{\gamma L\in (G/L)^H\ :\ \gamma K= g K} f(\gamma L)\ , \]
where $f:(G/L)^H\to N$ is $W_G H$-equivariant and $g K\in (G/K)^H$.
We omit the verification that $\tr_L^K(f)$ is again $W_G H$-equivariant.

With these definitions, the verification of properties~(i), (ii) and~(iii) 
of Definition \ref{def:G-Mackey} is straightforward, and we also omit it.
Checking the double coset formula is more involved, so we 
spell out the argument. 
We consider a subgroup $K$ of $G$ and subgroups $J$ and $L$ of $K$. 
We let $g\in G$ be such that $g J\in (G/J)^H$, i.e., $H^g\leq J$.
We let $R$ be a set of representatives of the $J$-$L$-double cosets in $K$. 
For every $k\in R$ we let $S_k$ be a set of representatives $j$
of those $(J\cap {^k L})$-cosets $j (J\cap {^k L})\in (G/(J\cap {^k L}))^H$ such that 
$j J=g J$.
Then the sets $S_k\!\cdot k$ are pairwise disjoint for $k\in R$, and their union
\[ {\bigcup}_{k\in R}\  S_k\!\cdot k \]
is a set of representatives of those $H$-fixed $L$-coset 
$\gamma L\in (G/L)^H$ such that $\gamma K=g K$.
Hence for every $W_G H$-map $f:(G/L)^H\to N$ all $g J\in (G/J)^H$ we obtain the relation
\begin{align*}
  (\res^K_J(\tr_L^K (f)))(g J)\ &= \ \sum_{\gamma L\in (G/L)^H\ :\ \gamma K= g K} f(\gamma L)\\
&= \ \sum_{k\in R} \ \sum_{ j\in S_k}\,  f(j k L)\\ 
&= \ \sum_{k\in R} \ \sum_{ j\in S_k}\,  k_\star(\res^L_{J^k\cap L}(f))(j(J\cap {^k L}))\\ 
&= \ \sum_{k\in R} \tr_{J\cap{^k L}}^J(k_\star(\res^L_{J^k\cap L}(f)))(g J)\ .
\end{align*}
For varying $g J$ and $f$, this proves the double coset relation
\[ \res^K_J\circ \tr_L^K\ =\ 
{\sum}_{k\in R}\, \tr_{J\cap{^k L}}^J\circ k_\star\circ\res^L_{J^k\cap L} \ : \ 
(\rho_H N)(L)\to(\rho_H N)(J) \ . \]
Now we have shown that $\rho_H N$ is a $G$-Mackey functor, 
so this completes the definition of the functor $\rho_H$ on objects;
on morphisms, $\rho_H$ is simply given by postcomposition with
a $W_G H$-module homomorphism.

The action of $W_G H$ on $(G/H)^H$ is free and transitive, so the evaluation map
\begin{equation}  \label{eq:rho_H at H}
 (\rho_H N)(H)\ = \ \map^{W_G H}((G/H)^H,N) \ \to \ N \ , \quad 
f\ \longmapsto \ f(e H)  
\end{equation}
is an isomorphism. Moreover, if $H$ is not subconjugate to $K$, then
the set $(G/K)^H$ is empty, and so $(\rho_H N)(K)=0$.
This holds in particular for all proper subgroups of $H$, so in
the $G$-Mackey functor $\rho_H N$, all transfers from proper subgroups 
to $H$ are trivial. Hence the isomorphism \eqref{eq:rho_H at H}
passes to a natural isomorphism of $W_G H$-modules
\[ \epsilon_H^N \ : \ \tau_H(\rho_H N) \ \iso \ N \ , 
\quad  \epsilon_H^N[ f]\ = \ f(e H)\ .\]
Here $[y]$ denotes the class in $\tau_H M$ of an element $y\in M(H)$.
Given a $G$-Mackey functor $M$, we also define a homomorphism
\[ \eta_H^M \ : \ M \ \to \ \rho_H(\tau_H M) \]
of $G$-Mackey functors as follows. For a subgroup $K$ of $G$
we define
\[ \eta_H^M(K)\ : \ M(K)\ \to \ \map^{W_G H}((G/K)^H,\tau_H M) \ = \
\rho_H(\tau_H M)(K) \]
by
\[ \eta_H^M(K)(x)(g K)\ = \  [ g_\star( \res^K_{H^g}(x) ) ] \ .\]
We show that these additive maps indeed define a morphism of $G$-Mackey functors.
For $L\leq K\leq G$, $x\in M(K)$ and $g L\in (G/L)^H$ we have 
\begin{align*}
  \res^K_L(\eta_H^M(K)(x))(g L)\ &= \   
\eta_H^M(K)(x)(g K)\ = \  [ g_\star(\res^K_{H^g}(x)) ] \\
&= \ [ g_\star(\res^L_{H^g}(\res^K_L(x))) ] \ = \ \eta_H^M(L)( \res^K_L(x))(g L)\ .
\end{align*}
So $\res^K_L\circ\eta_H^M(K)=\eta_H^M(L)\circ\res^K_L$, i.e.,
the $\eta_H^M$-maps are compatible with restrictions.
Similarly, for $\gamma\in G$ and $x\in M(K^\gamma)$ we have 
\begin{align*}
  \gamma_\star(\eta_H^M(K^\gamma)(x))(g\cdot K)\ &= \ 
\eta_H^M(K^\gamma)(x)(g\gamma \cdot K^\gamma)\ 
= \  [ (g\gamma)_\star(\res^{K^\gamma}_{H^{g\gamma}}(x)) ] \\
&= \ [ g_\star(\res^K_{H^g}(\gamma_\star(x))) ] \
= \ \eta_H^M(K)( \gamma_\star(x))(g\cdot K)\ .
\end{align*}
So $\gamma_\star\circ\eta_H^M(K^\gamma)=\eta_H^M(K)\circ\gamma_\star$, i.e.,
the $\eta_H^M$-maps are compatible with conjugations.

For compatibility with transfers we consider $g K\in (G/K)^H$ 
(i.e., $H^g\leq K$), and use the double coset formula 
\begin{align}  \label{eq:double_coset_rho_H}
g_\star\circ \res^K_{H^g}\circ\tr_L^K\ &= \ 
\sum_{(H^g)k L\in H^g\bs K/L} 
g_\star\circ \tr_{H^g\cap{^k L}}^{H^g}\circ k_\star\circ\res^L_{H^{g k}\cap L}  \nonumber\\ 
&\equiv\
\sum_{(H^g) k L\ : \ H^{g k}\leq L}  ( g k)_\star\circ\res^L_{H^{g k}} \nonumber \\
&= \ \sum_{\gamma L\in (G/L)^H\ : \ \gamma K = g K}  \gamma_\star\circ\res^L_{H^{g k}}  
\end{align}
as maps $M(L)\to \tau_{H^g}( M)$. The second relation is modulo
transfers from proper subgroups of $H^g$, i.e., we have dropped all summands
where $H^g$ is not contained in ${^k L}$.
The third relation exploits that sending $(H^g) k L$ to $g k L$
is a bijection between the set of those $H^g$-$L$-double cosets with $H^{g k}\leq L$
and the set of those $\gamma L\in (G/L)^H$ that satisfy $\gamma K=g K$.
So for all $x\in M(L)$ we obtain the relation
\begin{align*}
  \tr^K_L(\eta_H^M(L)(x))(g K)\ &= \  \sum\, \eta_H^M(L)(x)(\gamma L)\ 
= \  \sum\, [ \gamma_\star(\res^L_{H^{\gamma}}(x)) ] \\
 _\eqref{eq:double_coset_rho_H} 
&= \ [ g_\star(\res^K_{H^g}(\tr^K_L(x)))] \ = \  \eta_H^M(K)( \tr^K_L(x))(g K)\ .
\end{align*}
The sums run over those $H$-fixed $L$-cosets $\gamma L\in (G/L)^H$
such that $\gamma K= g K$.
So $\tr_L^K\circ\eta_H^M(L)=\eta_H^M(K)\circ\tr_L^K$, i.e.,
the $\eta_H^M$-maps are compatible with transfers.

It is now straightforward to check that the composite maps
\begin{align*}
   \rho_H N\ &\xra{\eta_H^{\rho_H N}}\  \rho_H(\tau_H(\rho_H N)) \
\xra[\iso]{\rho_H(\epsilon_H^N)} \ \rho_H N \text{\qquad and}\\
 \tau_H M\ &\xra{\tau_H(\eta_H^M)}\  \tau_H(\rho_H(\tau_H M)) \
\xra[\iso]{\epsilon_H^{\tau_H M}} \ \tau_H M 
\end{align*}
are the respective identity morphisms.
Indeed,
\begin{align*}
  \rho_H(\epsilon_H^N)(\eta_H^{\rho_H N}(K)(f) )(g K)\ &= \ 
\epsilon_H^N\big(\eta_H^{\rho_H N}(K)(f)(g K)\big)\ = \ 
\epsilon_H^N[ g_\star(\res^K_{H^g}(f))]\\ 
&= \ g_\star(\res^K_{H^g}(f))(e H)\ = \ \res^K_{H^g}(f)(g H)\ = \ f(g K)
\end{align*}
for all $W_G H$-equivariant $f:(G/K)^H\to N$ and all $g K\in(G/K)^H$.
Moreover,
\begin{align*}
  \epsilon_H^{\tau_H M}(\tau_H(\eta^M_H)[x])\ &= \ 
  \epsilon_H^{\tau_H M}[\eta^M_H(H)(x)]\ = \ \eta^M_H(H)(x)(e H)\ = \ [x]
\end{align*}
for all $x\in M(H)$.
This shows that $\eta_H$ and $\epsilon_H$ are the unit and counit 
of an adjunction between $\tau_H$ and $\rho_H$.

Now we can easily prove part~(i). Since $(\tau_H,\rho_H)$
is an adjoint pair for every subgroup $H$ of $G$,
a right adjoint to the product functor $\tau$ is given by
\[ \rho\ : \ {\prod}_{(H)}\, W_G H\mo \to  G\Mack \ , \quad
\rho( (N_H)_{(H)})\ = \ {\prod}_{(H)} \, \rho_H(N_H)\ ,\]
the product of the $G$-Mackey functors $\rho_H(N_H)$.
The universal property of the product of $G$-Mackey functors
and the previous adjunctions combine into a natural isomorphism
\begin{align*}
  G\Mack(M,\rho ( (N_H)_{(H)}) ) \ \xra{\ \iso\ } \ 
  &{\prod}_{(H)}\, G\Mack(M,\rho_H (N_H))\\ 
  \iso\quad  &{\prod}_{(H)} \,  W_G H\mo(\tau_H M, N_H)\ .  
\end{align*}
This completes the construction of the adjunction between $\tau$ and $\rho$.

(ii)
We show first that for every $G$-Mackey functor $M$, the adjunction unit
\[ \eta^M \ = \ (\eta_H^M)\ : \ M \ \to \ \rho(\tau M) \ = \ 
{\prod}_{(H)}\, \rho_H(\tau_H M) \]
becomes an isomorphism after inverting the order of $G$.
When evaluated at $G$, this is the content 
of Proposition \ref{prop:recover G-Mackey rationally}.
At a general subgroup of $G$, we reduce the claim 
to Proposition \ref{prop:recover G-Mackey rationally} as follows.
We let $H$ and $K$ be two subgroups of $G$.
We let $R_{H,K}$ be a set of representatives of the $K$-conjugacy classes 
of subgroups $L\leq K$ that are $G$-conjugate to $H$.
For each $L\in R_{H,K}$ we choose an element $g_L\in G$ such that $H^{g_L}=L$.
Then $g_L K\in (G/K)^H$, and the map
\[ {\coprod}_{L\in R_{H,K}} \ (W_G H) / (W_{^{g_L} K} H)\ \to \ (G/K)^H\ , \quad
[n H]\ \mapsto \ n g_L K \]
is an isomorphism of $W_G H$-sets. So taking $W_G H$-equivariant maps into $\tau_H M$
provides an isomorphism of abelian groups
\begin{align*}
    \rho_H(\tau_H M)(K)\ = \ 
\map^{W_G H}( (G/K)^H, \tau_H M) \ &\to \ {\prod}_{L\in R_{H,K}}\, (\tau_H M)^{W_{^{g_L} K} H} \\
f \qquad &\longmapsto \quad ( f(g_L K) )_{L} \ .
\end{align*}
Because $H^{g_L}=L$, the map $(g_L)_\star:\tau_L M\to\tau_H M$
defined in \eqref{eq:act_on_prod_tau_H} is an isomorphism of abelian groups,
equivariant for the isomorphism~
$W_K L=W_K H^{g_L}\iso W_{^{g_L} K}H$
induced by conjugation by $g_L$.
So the map
\[ \zeta_H\ : \ \rho_H(\tau_H M)(K)\ \to \ 
{\prod}_{L\in R_{H,K}}\, (\tau_L M)^{W_K L} \ , \quad
f \ \longmapsto \ ( (g_L)_\star^{-1} f(g_L K) )_{L} \]
is an isomorphism.

Now the following diagram commutes by direct inspection, where the right vertical map
is the product of the maps $\zeta_H$ over all conjugacy classes
of subgroups of $G$:
\[ \xymatrix@C=8mm{ 
M(K)\ar[rr]^-{\eta_M(K)}\ar@/_1pc/[drr]_{\psi^M_K} &&
\rho(\tau M)(K)\ar@{=}[r] &\prod_{(H)}\, \rho_H( \tau_H M)(K) \ar[d]^{\prod\zeta_H}_\iso \\ 
&& \prod_{(L)\leq K}  (\tau_L M)^{W_K L}\ar@{=}[r]
& \prod_{(H)\leq G} \prod_{L\in R_{H,K}} (\tau_L M)^{W_K L} } \]
The map $\psi^M_K$ becomes an isomorphism after inverting the
order of $K$, by Proposition \ref{prop:recover G-Mackey rationally},
applied to the underlying $K$-Mackey functor of $M$.
Since the right vertical map is an isomorphism and
the order of $K$ divides the order of $G$, we conclude that the map $\eta_M(K)$ 
becomes an isomorphism after inverting the order of $G$.
This completes the proof that $\eta_M:M\to \rho(\tau M)$ becomes an isomorphism 
after inverting $|G|$. 

Now $|G|$ is invertible in $Q$, so when restricted to $Q$-local objects, 
the adjunction unit $\eta:\Id\to\rho\circ\tau$ is a natural isomorphism.
Hence the restriction of the functor $\tau$ to the full subcategory of 
$Q$-local $G$-Mackey functors is fully faithful.
We claim that additionally, every family of $Q$-local $W_G H$-modules
is the image of a $Q$-local $G$-Mackey functor.
Since the functor $\tau$ is additive, it commutes with finite products,
so it suffices to show that for every individual subgroup $H$ of $G$
and every $Q$-local $W_G H$-module $N$, there is an
$Q$-local $G$-Mackey functor $M$ such that $\tau M$ is isomorphic to $N$.
We let $e:\tau(\rho_H N)\to\tau(\rho_H N)$ be the idempotent
endomorphism whose $H$-component is the identity of $\tau_H(\rho_H N)$,
and such that $e_K=0$ for all $K$ that are not conjugate to $H$.
Since $\rho_H N$ is $Q$-local and $\tau$ is fully faithful on 
$Q$-local $G$-Mackey functors, there exists an idempotent endomorphism
$\tilde e:\rho_H N\to\rho_H N$ such that $e=\tau(\tilde e)$.
Then $\rho_H N$ is isomorphic to the direct sum of the image of $\tilde e$
and the image of $\Id-\tilde e$. Since the functor $\tau$ is additive,
$\tau(\tilde e\cdot \rho_H N)$ is isomorphic to 
$e\cdot \tau(\rho_H N)=\tau_H(\rho_H N)$.
As we showed in the proof of part~(i), the $W_G H$-module $\tau_H(\rho_H N)$ 
is isomorphic to $N$, so this proves that $N$ 
is isomorphic to $\tau(\tilde e\cdot \rho_H N)$.
Altogether this shows that the functor $\tau$ is essentially surjective on
$Q$-local objects. On this subcategory, $\tau$ is also fully faithful,
and hence an equivalence of categories.
\end{proof}

Now we return to equivariant spectra.
For every orthogonal $G$-spectrum $X$, Proposition \ref{prop:Phi of transfer} (ii)
shows that the geometric fixed point map \eqref{eq:geometric_fix_map}
factors over the quotient of $\pi_0^G(X)$ by the subgroup
generated by proper transfers. We denote by
\begin{equation}  \label{eq:reduced_geometric_fixed}
 \bar\Phi \ : \  \pi_k^G(X) / t_G( \upi_k(X))
\ = \ \tau_G( \upi_k(X) )\ \to \ \Phi^G_k(X)   
\end{equation}
the induced map on the factor group, and call it the 
{\em reduced geometric fixed point homomorphism}.\index{subject}{geometric fixed point homomorphism!reduced}
Our next result shows that for finite groups $G$,
the `corrected' (i.e., reduced) geometric fixed point map
becomes an isomorphism after inverting the order of $G$.
The following proposition is well known,
and closely related statements appear in Appendix~A of \cite{greenlees-may-Tate};
however, I am not aware of a reference for the following statement in this form.

\begin{prop}\label{prop:rational geometric fixed for GSp} 
For  every finite group $G$, every orthogonal $G$-spectrum $X$
and every integer $k$ the reduced geometric fixed point map
\[ \bar\Phi \ : \ \tau_G( \upi_k(X) )\ \to \ \Phi^G_k(X) \]
becomes an isomorphism after inverting the order of $G$.\index{subject}{geometric fixed points!rationally}\index{subject}{geometric fixed points!and transfers} 
\end{prop} 
\begin{proof}
We start by showing that for every orthogonal $G$-spectrum $X$ and every
subgroup $H$ of $G$ the transfer map
\[ \tr_H^G\ : \ \pi_k^H(X\sm G/H_+) \ \to \ \pi_k^G(X\sm G/H_+) \]
for the spectrum $X\sm G/H_+$ is surjective.
Indeed, the transfer is defined as the composite
\[ \pi_k^H(X\sm G/H_+) \ \xra[\iso]{G\ltimes_H -} \ \pi_k^G(G\ltimes_H( X\sm G/H_+)) 
\ \xra{\text{act}_*} \ \pi_k^G(X\sm G/H_+) \ . \]
The first map -- the external transfer -- is an isomorphism
by Theorem \ref{thm:Wirth iso}.
The second map is surjective because the action morphism has a $G$-equivariant section
\[  X\sm G/H_+ \ \to \ G\ltimes_H( X\sm G/H_+)\ , \quad
x\sm g H \ \longmapsto \ [g,\, g^{-1} x\sm e H] \ .\]

Now we let $A$ be a $G$-CW-complex without $G$-fixed points.
We claim that for every orthogonal $G$-spectrum $X$,
the entire group $\pi_k^G(X\sm A_+)$ is generated by transfers from
proper subgroups after inverting $|G|$.
In a first step we show this when $A$ is a finite-dimensional $G$-CW-complex.
We argue by induction over the dimension of $A$. The induction starts when $A$
is empty, in which case $X\sm A_+$ is the trivial spectrum and there is nothing to show.
Now we consider $n\geq 0$ and assume the claim for all 
$(n-1)$-dimensional $G$-CW-complexes without $G$-fixed points. 
We suppose that $A$ is obtained
from such an $(n-1)$-dimensional $G$-CW-complex $B$ by attaching equivariant $n$-cells. 
After choosing characteristic maps for the $n$-cells we can identify
the quotient $A/B$ as a wedge
\[ A/B \ \iso \ {\bigvee}_{i\in I}\, S^n\sm (G/H_i)_+ \]
for some index set $I$ and certain subgroups $H_i$ of $G$.
Equivariant homotopy groups commute with wedges, so we can identify
the homotopy group Mackey functor of $X\sm (A/B)$ as
\[ \upi_k(X\sm A/B)\ \iso \ {\bigoplus}_{i\in I}\ \upi_k(X\sm S^n\sm (G/H_i)_+)\ . \]
Since $A$ has no $G$-fixed points, the groups $H_i$ are all proper subgroups of $G$,
so the group $\pi_k^G(X\sm S^n\sm (G/H_i)_+)$ is generated by transfers
from the proper subgroup $H_i$, by the previous paragraph.
So altogether the group $\pi_k^G(X\sm A/B)$ is generated by transfers
from proper subgroups (even before any localization).

The inclusion $B\to A$ is a $G$-equivariant h-cofibration, so it induces an
h-cofibration of orthogonal $G$-spectra $X\sm B_+\to X\sm A_+$
that results in a long exact sequence of homotopy group Mackey functors as in 
Corollary \ref{cor-long exact sequence h-cofibration}~(i).
The functor sending a $G$-Mackey functor $M$ to $\mZ[1/|G|]\tensor(\tau_G M)$
is exact (Corollary \ref{cor:tau_G is exact}), 
the groups $\mZ[1/|G|]\tensor \tau_G(\upi_*(X\sm B_+))$
vanish by induction, and the groups $\tau_G(\upi_*(X\sm A/B))$ 
vanish by the previous paragraph.
So the groups $\mZ[1/|G|]\tensor \tau_G(\upi_*(X\sm A_+))$
vanish by exactness, and this finishes the inductive step.

Now we suppose that $A$ is an arbitrary $G$-CW-complex without $G$-fixed points,
possibly infinite dimensional. 
We choose a skeleton filtration by $G$-subspaces $A^n$.
Then the groups $\mZ[1/|G|]\tensor \tau_G(\upi_*(X\sm A^n_+))$
vanish for all $n\geq 0$, by the previous paragraph.
All the morphisms $X\sm A^n_+\to X\sm A^{n+1}_+$ are h-cofibrations
of orthogonal $G$-spectra, hence levelwise closed embeddings. Since equivariant
homotopy groups and the functor $\mZ[1/|G|]\tensor \tau_G(-)$
both commute with such sequential colimits, this shows that  
the groups $\tau_G(\upi_k(X\sm A_+))$ vanish after inverting $|G|$.

Now we can prove the proposition.
The inclusion $i:S^0\to \tilde E\Pc_G$ gives rise to a commutative square:
\begin{equation} \begin{aligned}\label{eq:tau_Phi_square}
 \xymatrix@C=15mm{ 
 \tau_G( \upi_k(X))\ar[r]^-{\tau_G(\upi_k(X\sm i))}\ar[d]_{\bar\Phi} &
 \tau_G( \upi_k(X\sm\tilde E\Pc_G)) \ar[d]^{\bar\Phi} \\
\Phi^G_k(X)\ar[r]_-{\Phi^G_k(X\sm i)} &  \Phi^G_k(X\sm\tilde E\Pc_G) }   
  \end{aligned}\end{equation}
The lower horizontal map is an isomorphism because $i$ 
identifies $S^0$ with the fixed points $(\tilde E\Pc_G)^G$.
Proposition \ref{prop:geometric as fixed points} 
shows that the geometric fixed point map $\Phi: \pi^G_k(X\sm\tilde E\Pc_G)\to
 \Phi_k^G(X\sm\tilde E\Pc_G)$ is an isomorphism {\em before dividing out transfers}.
So all transfers from proper subgroups are in fact trivial  in 
$\pi^G_k(X\sm\tilde E\Pc_G)$, and the right vertical map
in the square \eqref{eq:tau_Phi_square} is an isomorphism 
(even before any localization).
By the previous paragraph, the group $\pi_k^G(X\sm (E\Pc_G)_+)$ 
is generated by transfers from proper subgroups after inverting $|G|$, so
\[ \mZ[1/|G|]\tensor \tau_G( \upi_k(X\sm (E\Pc_G)_+)\ = \ 0 \ . \]
The functor sending a $G$-Mackey functor $M$ to the group $\mZ[1/|G|]\tensor(\tau_G M)$
is exact (Corollary \ref{cor:tau_G is exact}), 
so the isotropy separation sequence \eqref{eq:isotropy_separation_les} 
shows that the upper horizontal map in the square \eqref{eq:tau_Phi_square}
becomes an isomorphism after inverting $|G|$.
So the left vertical map in the square \eqref{eq:tau_Phi_square}
becomes an isomorphism after inverting $|G|$.
\end{proof}

When $X$ is an orthogonal $G$-spectrum, we can apply the earlier algebraic
Proposition \ref{prop:recover G-Mackey rationally}
to the $G$-Mackey functor $\upi_k(X)$; this allows us -- after inverting $|G|$ -- 
to reconstruct $\pi_k^G(X)$ from the groups $\tau_H(\upi_k(X))$ 
with their Weyl group action.
Moreover, after inverting the group order we can use the reduced
geometric fixed point map \eqref{eq:reduced_geometric_fixed},
for the underlying $H$-spectrum of $X$, to identify
the group $\tau_H ( \upi_k(X))$ with the group $\Phi_k^H(X)$.
Under this identification, the map $\bar\psi^M_G$ 
becomes the product of the maps
\[ \pi_k^G(X) \ \xra{\ \res^G_H\ }\ \pi_k^H(X)\ \xra{\ \Phi^H\ } \
\Phi_k^H(X)\ .\]
So Propositions \ref{prop:recover G-Mackey rationally} 
and \ref{prop:rational geometric fixed for GSp}  together prove:

\begin{cor}\label{cor:equivariant from geometric}
For every finite group $G$, every orthogonal $G$-spectrum $X$
and every integer $k$ the map
\[  (\Phi^H\circ\res^G_H)_H\ : \ 
\pi^G_k(X) \ \to \ {\prod}_{(H)}\, \left( \Phi^H_k(X)\right)^{W_G H} \]
becomes an isomorphism after inverting the order of $G$.
\end{cor} 

\index{subject}{Mackey functor|)}

\section{Products}\label{sec:products}

In this section we recall the smash product of orthogonal spectra
and orthogonal $G$-spectra and study its formal and homotopical properties.
Like the box product of orthogonal spaces, the smash product 
of orthogonal spectra is a special case of Day's convolution product on categories
of enriched functors, compare Appendix \ref{app:enriched functors}.
The main homotopical result about the smash product is the `flatness theorem'
(Theorem \ref{thm:G-flat is flat}), showing that smashing with 
$G$-flat orthogonal $G$-spectra (in the sense of Definition \ref{def:G-flat})
preserves $\upi_*$-isomorphisms.

The smash product is intimately related to pairings of equivariant stable homotopy groups
that we recall in Construction \ref{con:pairing equivalence homotopy};
the main properties of these pairings are summarized 
in Theorem \ref{thm:product properties}. 
When specialized to equivariant ring spectra, these pairings turn 
the equivariant stable homotopy groups into a graded ring,
compare Corollary \ref{cor:product properties ring spectrum}. 

\medskip

The smash product of orthogonal spectra is characterized by a
universal property that we recall now.
The indexing category $\bO$ for orthogonal spectra has a symmetric monoidal
product by direct sum as follows. We denote by $\bO\sm\bO$ the
category enriched in based spaces whose objects are
pairs of inner product spaces and whose morphism spaces
are smash products of morphism spaces in $\bO$. A based continuous functor
\[ \oplus \ : \  \bO\sm\bO \ \to \ \bO\]
is defined on objects by orthogonal direct sum, and on morphism spaces by
\begin{align*}
   \bO(V,W)\sm \bO(V',W')\ &\to \ \bO(V\oplus W,V'\oplus W')\\
(w,\varphi)\sm (w',\varphi')\quad &\longmapsto \quad ((w,w'),\varphi\oplus\varphi')\ .
\end{align*}
A {\em bimorphism} $b:(X,Y)\to Z$\index{subject}{bimorphism!of orthogonal spectra}
from a pair of orthogonal spectra $(X,Y)$
to an orthogonal spectrum $Z$ is a natural transformation 
\[ b \ : \ X\bar\sm Y \ \to \ Z\circ\oplus \]
of continuous functors $\bO\sm\bO\to\bT_*$;
here $X\bar\sm Y$ is the `external smash product'  defined by
$(X\bar\sm Y)(V,W)=X(V)\sm Y(W)$. A bimorphism thus consists of
based continuous maps 
\[ b_{V,W} \ : \ X(V) \sm  Y(W) \ \to \ Z(V\oplus W) \]
for all inner product spaces $V$ and $W$ that form morphisms of orthogonal
spectra in each variable separately.
A smash product of two orthogonal spectra is now a universal example
of a bimorphism from $(X,Y)$.

\begin{defn}
A {\em smash product}\index{subject}{smash product}\index{symbol}{$\sm$ - {smash product of orthogonal spectra}}
of two orthogonal spectra $X$ and $Y$ is a pair $(X\sm Y,i)$ 
consisting of an orthogonal spectrum $X\sm Y$
and a universal bimorphism $i:(X,Y)\to X\sm Y$,
i.e., a bimorphism 
such that for every orthogonal spectrum $Z$ the map
\begin{equation}\label{eq:universal property smash}
 \spec(X\sm Y,Z) \ \to \ \text{Bimor}((X,Y),Z) \ , \
f\mapsto f i = \{f(V\oplus W)\circ i_{V,W}\}_{V,W} 
\end{equation}
is bijective. 
\end{defn}

Since the index category $\bO$ is skeletally small and the base category $\bT_*$
is cocomplete, every pair of orthogonal spectra has a smash product
by Proposition \ref{prop:box exists}. 
Often only the object $X\sm Y$ will be referred to 
as the smash product, but one should keep in mind that it comes equipped
with a specific, universal bimorphism. We will often refer to the 
bijection \eqref{eq:universal property smash} as the
{\em universal property} of the smash product of orthogonal spectra.
\index{subject}{universal property!of the smash product}

While the smash products for pairs of orthogonal spectra are choices, 
there is a preferred way to extend any chosen smash products 
to a functor in two variables
(Construction \ref{con:box functoriality});
this functor has a preferred symmetric monoidal structure, i.e., 
distinguished natural associativity and symmetry isomorphisms
\[  \alpha_{X,Y,Z}\ : \ (X\sm Y)\sm Z \ \to \ X\sm(Y\sm Z)
\text{\quad respectively \quad} 
\tau_{X,Y}\ : \  X\sm Y \ \to \ Y\sm X  \]
(see Construction \ref{con:box coherence}).
Together with strict  unit isomorphisms $\mS\sm X= X= X\sm\mS$,
these satisfy the coherence conditions of a symmetric monoidal category,
compare Day's Theorem \ref{thm:symmetric monoidal}.
The smash product of orthogonal spectra is {\em closed} symmetric monoidal
in the sense that the smash product is adjoint to an
internal Hom spectrum.
We won't use the internal function spectrum, so we don't go into any details.

When a compact Lie group $G$ acts on the orthogonal spectra $X$ and $Y$,
then $X\sm Y$ becomes an orthogonal $G$-spectrum via the diagonal action.
So the smash product lifts to a symmetric monoidal closed structure
\[ \sm \ : \ G\spec\times G\spec \ \to \ G\spec \]
on the category of orthogonal $G$-spectra.\index{subject}{smash product!of orthogonal $G$-spectra}

\begin{eg}[Smash product with a free orthogonal spectrum]
We give a `formula' for the smash product with a free orthogonal spectrum,
i.e., the represented functor $\bO(W,-)$ for some inner product space $W$.
A general feature of Day type convolution products is that
the convolution product of represented functors is represented, 
see Remark \ref{rk:box representable}.
In the case at hand this specializes to a preferred isomorphism
\[ \bO(V,-)\sm \bO(W,-)\ \iso \ \bO(V\oplus W,-)    \]
specified, via the universal property \eqref{eq:universal property smash},
by the bimorphism with $(U,U')$-component
\begin{equation}\label{eq:free_smash_free}
 \oplus \ : \ \bO(V,U)\sm \bO(W,U')\ \to \ \bO(V\oplus W,U\oplus U')\ .   
\end{equation}
Now we let $X$ be any orthogonal spectrum.
We give a formula for the value at $V\oplus W$ of the smash product $X\sm \bO(W,-)$. 
The composite 
\[ X(V)\ \xra{-\sm \Id_W} \ X(V)\sm \bO(W,W)\ \xra{\ i_{V,W}\ }\ 
(X\sm \bO(W,-))(V\oplus W) \]
is $O(V)$-equivariant, where $O(V)$ acts on the target by restriction
along the homomorphism $-\oplus W:O(V)\to O(V\oplus W)$.
So the map extends to an $O(V\oplus W)$-equivariant based map
\begin{equation}  \label{eq:smash_with_free}
 O(V\oplus W)\ltimes_{O(V)} X(V)\ \to \ (X\sm \bO(W,-))(V\oplus W) \ .  
\end{equation}
We claim that this map is a homeomorphism.
To see that we first check the special case where $X=\bO(U,-)$
is itself a represented functor. The isomorphism \eqref{eq:free_smash_free}
turns \eqref{eq:smash_with_free} into the map
\begin{align*}
   O(V\oplus W)\ltimes_{O(V)} \bO(U,V)\ &\to \quad \bO(U\oplus W,V\oplus W) \\
[A, (v,\varphi)]\quad &\longmapsto \ (A\cdot (v,0),A\circ(\varphi\oplus W))
\end{align*}
which is a homeomorphism by inspection.
Source and target of the map \eqref{eq:smash_with_free}
commute with smash products with based spaces
and colimits in the variable $X$. Since every orthogonal spectrum
is a coend of represented orthogonal spectra smashed with based spaces,
the map \eqref{eq:smash_with_free} is an isomorphism in general.
\end{eg}

Our next topic is a flatness result, proving that smashing with
certain orthogonal $G$-spectra preserves $\upi_*$-isomorphisms.
To define the relevant class, the `$G$-flat orthogonal $G$-spectra', 
we introduce the {\em skeleton filtration},
a functorial way to write an orthogonal spectrum as a 
sequential colimit of spectra which are made from the information
below a fixed level. The word `filtration' should be used with caution
because the maps from the skeleta to the orthogonal spectrum
need not be injective. Since an orthogonal $G$-spectrum is the same data as an orthogonal
spectrum with continuous $G$-action, and since skeleta are functorial,
the skeleta of an orthogonal $G$-spectrum automatically come as orthogonal $G$-spectra.

\begin{construction}[Skeleton filtration of orthogonal spectra]\label{con:skeleton spec}
As in the unstable situation in Section \ref{sec:global model structures spaces},
the skeleton filtration is a special case of 
the general skeleton filtration on certain enriched functor categories
that we discuss in Appendix \ref{app:enriched functors}.
Indeed, if we specialize the base category to $\Vc=\bT_*$,
the category of based spaces under smash product, and the index category
to $\Dc=\bO$, then the functor category $\Dc^\ast$ becomes the
category $\spec$ of orthogonal spectra. The dimension function needed in
the construction and analysis of skeleta is the vector space dimension.
For every orthogonal spectrum $X$ and every $m\geq 0$,
the general theory provides 
an {\em $m$-skeleton},\index{subject}{skeleton!of an orthogonal spectrum}
\[ \sk^m X\ = \ l_m(X^{\leq m}) \ ,\]
the extension of the restriction of $X$ to $\bO_{\leq m}$,
and a natural morphism $i_m:\sk^m X\to X$, the counit of
the adjunction $(l_m,(-)^{\leq m})$. The value
\[ i_m(V)\ :\ (\sk^m X)(V)\ \to \ X(V) \]
is an isomorphism for all inner product spaces $V$ of dimension at most $m$.
The {\em $m$-th latching space}\index{subject}{latching space!of an orthogonal spectrum}
of $X$ is the based $O(m)$-space
\[ L_m X \ = \ (\sk^{m-1} X)(\mR^m) \ ;\]
it comes with a natural based $O(m)$-equivariant map
\[  \nu_m=i_{m-1}(\mR^m)\ :\ L_m X\ \to \ X(\mR^m) \ , \]
the {\em $m$-th latching map}. 
We also agree to set $\sk^{-1} X=\ast$, the trivial orthogonal spectrum,
and $L_0 X=\ast$, a one-point space. 

The different skeleta are related by natural morphisms
$j_m:\sk^{m-1} X\to \sk^m X$, for all $m\geq 0$, 
such that $i_m\circ j_m=i_{m-1}$.
The sequence of skeleta stabilizes to $X$ in a 
very strong sense. For every inner product space $V$, the maps
$j_m(V)$ and $i_m(V)$ are isomorphisms as soon as $m >\dim(V)$.
In particular, $X(V)$ is a colimit, with respect to the maps $i_m(V)$, 
of the sequence of maps $j_m(V)$. 
Since colimits in the category of orthogonal spectra are created objectwise,
we deduce that the orthogonal spectrum $X$ is a colimit, 
with respect to the morphisms $i_m$, of the sequence of morphisms $j_m$. 
\end{construction}

Now we can define the class of $G$-flat orthogonal $G$-spectra for which smash product
is homotopical. As we remarked earlier,
the skeleta of (the underlying orthogonal spectrum of)
an orthogonal $G$-spectrum inherit a continuous $G$-action by functoriality.
In other words, the skeleta and the various morphisms between them
lift to endofunctors and natural transformations on the category of
orthogonal $G$-spectra.
If $X$ is an orthogonal $G$-spectrum, then the $O(m)$-space $L_m X$ 
comes with a commuting action by $G$, again by functoriality of the latching space.
Moreover, the latching morphism $\nu_m:L_m X\to X(\mR^m)$ is
$(G\times O(m))$-equivariant.

\begin{defn}\label{def:G-flat}
  Let $G$ be a compact Lie group.
  An orthogonal $G$-spectrum $X$ is {\em $G$-flat}\index{subject}{G-flat@$G$-flat!orthogonal $G$-spectrum} if for every $m\geq 0$ the latching map $\nu_m:L_m X \to X(\mR^m)$
  is a $(G\times O(m))$-cofibration.
\end{defn}

When $G$ is the trivial group, we simply speak of {\em flat} 
orthogonal spectra.\index{subject}{flat!orthogonal spectrum}
These flat orthogonal spectra play a special role as the cofibrant objects in the global
model structure on the category of orthogonal spectra, to be established
in Theorem \ref{thm:All global spectra}  below.

\begin{rk}
The $G$-flat orthogonal spectra are the cofibrant objects in a certain model
structure on the category of orthogonal $G$-spectra, the  
$\mS$-model structure of Stolz \cite[Thm.\,2.3.27]{stolz-thesis}.
This model structure has more cofibrant objects 
than the one of Mandell-May \cite[III Thm.\,4.2]{mandell-may},
whose `q-cofibrant' orthogonal $G$-spectra are precisely the ones 
for which the latching  map $\nu_m:L_m X \to X(\mR^m)$
is a $(G\times O(m))$-cofibration and in addition the group $O(m)$
acts freely on the complement of the image of $\nu_m$.

A convenient feature of the class of flat equivariant spectra is that they
are closed under various change-of-group functors. 
For example, if $\alpha:K\to G$ is a continuous 
homomorphism and $X$ an orthogonal $G$-spectrum, then the 
orthogonal $K$-spectrum $\alpha^* X$ satisfies
\[ L_m (\alpha^* X) \ = \ (\alpha\times O(m))^*(L_m X)\text{\quad and\quad}
(\alpha^* X)(\mR^m)\ = \ (\alpha\times O(m))^*(X(\mR^m)) \]
as $(K\times O(m))$-spaces. Since restriction along 
$\alpha\times O(m)$ preserves cofibrations 
(see Proposition \ref{prop:cofibrancy preservers}~(i)),
the restriction functor $\alpha^*:G\spec\to K\spec$
takes $G$-flat orthogonal $G$-spectra to
$K$-flat orthogonal $K$-spectra. Of particular relevance for global considerations
is the special case of the unique homomorphism $G\to e$ to the
trivial group. In this case we obtain that a flat orthogonal spectrum 
is in particular $G$-flat when we give it the trivial $G$-action.

A similar argument applies to induction of equivariant spectra. Indeed, if $H$
is a closed subgroup of $G$ and $Y$ an orthogonal $H$-spectrum, then
\begin{align*}
   L_m (G\ltimes_H Y) \ &\iso \ (G\times O(m))\ltimes_{H\times O(m)} (L_m Y)
\text{\quad and}\\
(G\ltimes_H Y)(\mR^m)\ &\iso \ (G\times O(m))\ltimes_{H\times O(m)} Y(\mR^m)\ . 
\end{align*}
Since induction from $H\times O(m)$ to $G\times O(m)$
preserves cofibrations (see Proposition \ref{prop:cofibrancy preservers}~(ii)),
the induction functor $G\ltimes_H -:H\spec\to G\spec$
takes $H$-flat orthogonal $H$-spectra to
$G$-flat orthogonal $G$-spectra.
\end{rk}

The main purpose of the latching object and latching morphism are to record
how one skeleton of an orthogonal $G$-spectrum is obtained from the previous one.
We denote by
\[ G_m  \ : \ O(m)\bT_* \ \to \ \spec \]
the left adjoint to the evaluation functor $X\mapsto X(\mR^m)$.
Explicitly, this functor sends a based $O(m)$-space $A$ to the orthogonal spectrum
\[ G_m A \ = \ \bO(\mR^m,-)\sm_{O(m)} A\ .\]
Since $G_m$ is a continuous functor, it takes $(G\times O(m))$-spaces
to $G$-orthogonal spectra, by functoriality.
Proposition \ref{prop:basic skeleton pushout square} 
and the fact that pushouts of orthogonal $G$-spectra are created on
underlying orthogonal spectra yield:

\begin{prop}
For every orthogonal $G$-spectrum $X$ and every $m\geq 0$ the commutative square
\[  
\xymatrix@C=12mm{ G_m L_m X \ar[r]^{G_m\nu_m} \ar[d] & G_m X(\mR^m) \ar[d]\\
\sk^{m-1} X \ar[r]_-{j_m} & \sk^m X}  
\]
is a pushout of orthogonal $G$-spectra.
The two vertical morphisms are adjoint to the identity of $L_m X$ 
respectively $X(\mR^m)$.
\end{prop}

The next theorem shows that smashing with a $G$-flat orthogonal $G$-spectrum
preserves $\upi_*$-isomorphisms.
This result is due to Stolz \cite[Prop.\,2.3.29]{stolz-thesis},
because the $G$-flat orthogonal $G$-spectra are precisely the $\mS$-cofibrant
objects in the sense of \cite[Def.\,2.3.4]{stolz-thesis};
a proof can also be found in \cite[Prop.\,2.10.1]{brun-dundas-stolz}. 
The theorem is stronger than the earlier result 
of Mandell and May \cite[III Prop.\,7.3]{mandell-may}
because the class of $G$-flat of orthogonal $G$-spectra
is strictly larger than the cofibrant $G$-spectra 
in the sense of \cite[Ch.\,III]{mandell-may}.

Since Stolz' thesis \cite{stolz-thesis} is not published, the notation and 
level of generality in \cite{brun-dundas-stolz}
is different from ours, 
and the characterization of flat objects in terms of latching maps is not
explicitly mentioned in \cite{brun-dundas-stolz} nor \cite{stolz-thesis},
we spell out the argument.

\begin{theorem}[Flatness theorem]\label{thm:G-flat is flat}\index{subject}{flatness theorem!for orthogonal $G$-spectra} 
Let $G$ be a compact Lie group and $X$ a $G$-flat orthogonal $G$-spectrum.
Then the functor $-\sm X$ preserves $\upi_*$-isomorphisms of orthogonal $G$-spectra.
\end{theorem}
\begin{proof}
We go through a sequence of six steps, 
proving successively more general cases of the theorem.
We call an orthogonal $G$-spectrum {\em $G$-stably contractible}
if all of its equivariant homotopy groups vanish, for all closed subgroups of $G$.
In Step~1 through~5, we let $C$ be a $G$-stably contractible orthogonal $G$-spectrum,
and we show for successively more general $X$ that $C\sm X$ 
is again $G$-stably contractible.

Step~1:
We let $K$ be another compact Lie group, $W$ a $K$-representation
and $A$ a based $(G\times K)$-space.
We define an orthogonal $G$-spectrum
\[  A\triangleright_{K,W} C \ = \ A\sm_K (\sh^W C) \]
at an inner product space $U$ as
\[  (A\triangleright_{K,W} C)(U)  \ = \ A\sm_K C(U\oplus W)  \ ; \]
the $K$-action on $C(U\oplus W)$ is through the $K$-action on $W$;
the $G$-action on the smash product is diagonally, 
through the $G$-actions on $A$ and on $C$.

Now we assume in addition that $A$ is a finite based $(G\times K)$-CW-complex and
the $K$-action is free (away from the basepoint). 
We claim that then $A\triangleright_{K,W} C$ is also $G$-stably contractible.
We argue by the number of equivariant cells, starting with $A=\ast$, in which case
there is nothing to show. Now we assume that the claim holds for $A$;
we let $B$ be obtained from $A$ by attaching an equivariant cell 
$(G\times K)/\Gamma\times D^k$,
where $\Gamma$ is a closed subgroup of $G\times K$.
The inclusion $A\to B$ is an h-cofibration of $(G\times K)$-spaces, 
hence the induced morphism
$A\triangleright_{K,W} C \to B\triangleright_{K,W} C$
is an h-cofibration of orthogonal $G$-spectra. 
By the long exact homotopy group sequence 
(Corollary \ref {cor-long exact sequence h-cofibration}~(i))
it suffices to show the claim for the quotient:
\begin{align*}
  ( B \triangleright_{K,W} C) / ( A \triangleright_{K,W} C) \ &\iso \  
  (A/B) \triangleright_{K,W} C  \ 
 \iso \    S^k \sm  (( G\times K)/\Gamma)_+ \triangleright_{K,W} C  
\end{align*}
By the suspension isomorphism we may show that
$( (G\times K)/\Gamma)_+ \triangleright_{K,W} C $ 
is $G$-stably contractible.
Since the $K$-action on $B$ is free, the stabilizer group $\Gamma$ 
is the graph of a continuous homomorphism
$\alpha:H\to K$ from some closed subgroup $H$ of $G$ (namely the projection
of $\Gamma$ to $G$). So we can rewrite 
\[   (  (G\times K)/\Gamma)_+ \triangleright_{K,W} C \ \iso \
   G\ltimes_H ( \sh^{\alpha^*(W)} C) \ .  \]
Since $C$ is $G$-stably contractible, it is also $H$-stably contractible.
The orthogonal $H$-spectrum  $C\sm S^{\alpha^*(W)}$
is then $H$-stably contractible
by Proposition \ref{prop:space smash preserves global} (ii).
Since $\sh^{\alpha^*(W)}C$ is $\upi_*$-isomorphic to 
 $C\sm S^{\alpha^*(W)}$ (by~Proposition \ref{prop:lambda upi_* isos}~(ii)),
it is $H$-stably contractible.
Since the induction functor $G\ltimes_H-$ preserves
$\upi_*$-isomorphisms (by Corollary \ref{cor:induce H2G}), 
we conclude that $(  (G\times K)/\Gamma)_+ \triangleright_{K,W} C$ 
is $G$-stably contractible.
This finishes the inductive argument, and hence the proof
that the orthogonal $G$-spectrum $A\triangleright_{K,W} C$
is $G$-stably contractible.

Step~2:
We let $G$ be a compact Lie group, $K$ a closed subgroup of $O(m)$
and $\alpha:G\to W_{O(m)} K$ a continuous homomorphism to the Weyl group of $K$
in $O(m)$.
The represented orthogonal spectrum $\bO(\mR^m,-)$
has a right $O(m)$-action through the tautological action on $\mR^m$.
The semifree orthogonal spectrum 
\[ G_m ( O(m)/K_+) \ = \ \bO(\mR^m,-)/K \]
has a residual action of the Weyl group $W_{O(m)} K$, 
and we restrict scalars along $\alpha$ to obtain the orthogonal $G$-spectrum 
\[ \alpha^*(\bO(\mR^m,-)/K)\ . \]
We claim that then the orthogonal $G$-spectrum
$C\sm \alpha^*(\bO(\mR^m,-)/K)$ is $G$-stably contractible.
We show that $\pi_k^H( C\sm \alpha^*(\bO(\mR^m,-)/K) )=0$
for $k\geq 0$ and every closed subgroup $H$ of $G$; for $k<0$ we exploit that
$\pi_k^H( C\sm \alpha^*(\bO(\mR^m,-)/K) )$ is isomorphic to
$\pi_0^H( (C\sm S^{-k})\sm \alpha^*(\bO(\mR^m,-)/K) )$,
so applying the argument to $C\sm S^{-k}$ instead of $C$ shows the
vanishing in negative dimensions.
Since the input data is stable under passage to closed subgroups of $G$
(just restrict $\alpha$ to such a subgroup), it is no loss of generality
to assume $H=G$.

We can represent every class of $\pi_k^G(C\sm \alpha^*(\bO(\mR^m,-)/K ))$
by a based $G$-map
\[  f \ : \  S^{V\oplus\mR^{m+k}}\ \to\  
(C\sm \alpha^*(\bO(\mR^m,-)/K ))(V\oplus\mR^m)\]
for some $G$-representation $V$.
In \eqref{eq:smash_with_free} we exhibited an isomorphism 
\[  O(V\oplus \mR^m)\ltimes_{O(V)} C(V)\ \xra{\ \iso\ } \ 
(C\sm \bO(\mR^m,-))(V\oplus \mR^m) \ .   \]
Under this isomorphism, the $K$-action on $\bO(\mR^m,-)$
becomes the right translation $K$-action on $O(V\oplus\mR^m)$,
so passage to $K$-orbits yields an isomorphism
\[  O(V\oplus \mR^m)/K_+\sm_{O(V)} C(V)\ \iso \ 
(C\sm \alpha^*(\bO(\mR^m,-)/K))(V\oplus \mR^m) \ . \]
The $G$-action on the right hand side is diagonally from three ingredients:
the $G$-action on $C$, the $G$-action on $V$, and the $G$-action
on $\bO(\mR^m,-)/K$ through $\alpha$. 
So $G$ acts on the left hand side diagonally, also through the actions
on $V$ (by left translation on $O(V\oplus\mR^m)$), on $\mR^m$ through $\alpha$,
and on $C$.

The target of $f$ is thus $G$-equivariantly homeomorphic to  
\[  O(V\oplus\mR^m)/K _+\sm_{O(V)}  C(V) \ = \ 
\left( O(V\oplus\mR^m)/K_+ \triangleright_{O(V), V}  C\right)(0) \ . \]
Via this isomorphism, we view $f$ as representing a $G$-equivariant homotopy class in
\begin{equation}\label{eq:complicated homotopy group}
 \pi_{m+k}^G\left( \Omega^V(O(V\oplus\mR^m)/K_+\triangleright_{O(V), V}  C)\right)\ . 
\end{equation}
The homogeneous space $O(V\oplus\mR^m)/K$ has a left $G$-action by translation
through the composite 
\[ G\ \to \ O(V)\ \xra{-\oplus\mR^m}\ O(V\oplus\mR^m) \]
and another $G$-action by restriction of the residual $W_{O(m)}(K)$-action 
along $\alpha$; altogether, the group $G$ acts diagonally.
This space also has a commuting free right
action of $O(V)$ by right translation.
Since $O(V\oplus\mR^m)/K$ is a smooth manifold
and the $(G\times O(V))$-action is smooth, 
Illman's theorem \cite[Cor.\,7.2]{illman} 
provides a finite $(G\times  O(V))$-CW-structure on $O(V\oplus\mR^m)/K$;
so $O(V\oplus\mR^m)/K_+$ is cofibrant as a based $(G\times O(V))$-space,
with free right $O(V)$-action.
Since $C$ is $G$-stably contractible, 
so is $O(V\oplus\mR^m)/K_+ \triangleright_{O(V), V}  C$, by Step~1.
The orthogonal $G$-spectrum
$\Omega^V\left( O(V\oplus\mR^m)/K_+ \triangleright_{O(V), V}  C\right)$
is then $G$-stably contractible by Proposition \ref{prop:map(A,-) preserves global}~(ii).
In particular, the group \eqref{eq:complicated homotopy group} is trivial,
so there is a $G$-representation $U$ such that the composite
\begin{align*}
   S^{U\oplus V\oplus\mR^{m+k}}\ 
&\xra{S^U\sm f}\  S^U\sm (O(V\oplus\mR^m)/K_+) \sm_{O(V)}  C(V) \\ 
&\xra{\sigma_{U,V}}\  O(V\oplus \mR^m)/K_+ \sm_{O(V)}  C(U\oplus V)   
\end{align*}
is $G$-equivariantly null-homotopic.
The structure map
\begin{align*}
   \sigma_{U, V\oplus\mR^m} \ : \ S^U\sm &(C\sm \alpha^*(\bO(\mR^m,-)/K))(V\oplus\mR^m) \\
&\to \ (C\sm \alpha^*(\bO(\mR^m,-)/K ))(U\oplus V\oplus \mR^m)  
\end{align*}
of the $G$-spectrum $C\sm \alpha^*(\bO(\mR^m,-)/K )$ factors as the composite
\begin{align*}
S^U\sm O(V\oplus\mR^m)/K_+ \sm_{O(V)} C(V)
&\xra{ O(V\oplus\mR^m)/K_+ \sm_{O(V)} \sigma_{U,V}} \\ 
O(V\oplus\mR^m)/K_+ \sm_{O(V)} C(U\oplus V) \ 
&\to \ O(U\oplus V\oplus\mR^m)/K_+\sm_{O(U\oplus V)} C(U\oplus V)\ .
\end{align*}
So also the stabilization $\sigma_{U,V}\circ(S^U\sm f)$
in the orthogonal $G$-spectrum $C\sm\alpha^*(\bO(\mR^m,-)/K)$ 
is equivariantly null-homotopic.
Since the stabilization represents the same element as $f$, this shows that 
$\pi_k^G(C\sm \alpha^*(\bO(\mR^m,-)/K))=0$.

Step 3:
We let $\Gamma$ be a closed subgroup of $G\times O(m)$.
We claim that smashing with the orthogonal $G$-spectrum
\[ G_m\left( (G\times O(m))/\Gamma_+\right)  \ = \ 
\bO(\mR^m,-)\sm_{O(m)} (G\times O(m))/\Gamma_+\]
preserves $G$-stably contractible orthogonal $G$-spectra.
To see that we interpret $\Gamma$ as a `generalized graph'.
We let $H\leq G$ be the projection of $\Gamma$ onto the first factor.
We let
\[ K\ = \ \{ A\in O(m) \ : \ (1,A)\in\Gamma\}\]
be the trace of $\Gamma$ in $O(m)$. Then $\Gamma$ is the graph of the continuous 
homomorphism
\[ \alpha\ : \ H \ \to \ W_{O(m)} K\ , \quad 
\alpha(h)\ = \ \{ A\in O(m)\ : \ (h,A)\in \Gamma\}\ , \]
in the sense that
\[ \Gamma \ = \ {\bigcup}_{h\in H} \ \{h\}\times \alpha(h)\ . \]
Moreover, the orthogonal $G$-spectrum 
$ G_m\left( (G\times O(m))/\Gamma_+\right)$ is isomorphic to
\[ G\ltimes_H(\alpha^*(\bO(\mR^m,-)/K))\ . \]
So for every orthogonal $G$-spectrum $C$, the shearing isomorphism
becomes an isomorphism of  orthogonal $G$-spectra
\begin{align*}
 C\sm (G_m ( (G\times O(m))/\Gamma_+) ) \ &\iso \ 
C\sm ( G\ltimes_H(\alpha^*(\bO(\mR^m,-)/K))) \\ 
&\iso\  G\ltimes_H(\res^G_H(C )\sm \alpha^*(\bO(\mR^m,-)/K) )\ .  
\end{align*}
If $C$ is $G$-stably contractible,
then it is also $H$-stably contractible;
so $\res^G_H(C)\sm \alpha^*(\bO(\mR^m,-)/K)$ is  $H$-stably contractible by Step~2.
Since the induction functor $G\ltimes_H-$ takes $\upi_*$-isomorphisms
of orthogonal $H$-spectra to $\upi_*$-isomorphisms
of orthogonal $G$-spectra (Corollary \ref{cor:induce H2G}), this shows that
$C\sm G_m((G\times O(m))/\Gamma_+)$ is $G$-stably contractible.

Step 4:
We let $A$ be a cofibrant based $(G\times O(m))$-space.
We claim that smashing with the orthogonal $G$-spectrum $G_m A$
preserves $G$-stably contractible orthogonal $G$-spectra.
A cofibrant based $(G\times O(m))$-space is equivariantly homotopy
equivalent to a based $(G\times O(m))$-CW-complex, so it is no loss of generality
to assume an equivariant CW-structure with skeleton filtration
\[ \ast\ = \ A_{-1} \ \subset \ A_0 \ \subset \ \dots \ \subset \ A_n \
\subset \ \dots \ .\]
We show first, by induction on $n$, that 
smashing with the orthogonal $G$-spectrum
$G_m A_n$ preserves $G$-stably contractible orthogonal $G$-spectra.
The induction starts with $n=-1$, where there is nothing to show. 
For $n\geq 0$ the inclusion
$A_{n-1}\to A_n$ is an h-cofibration of based $(G\times O(m))$-spaces,
so the induced morphism $C\sm G_m A_{n-1}\to C\sm G_m A_n$ is an h-cofibration
of orthogonal $G$-spectra.
By induction and the long exact sequence of equivariant
homotopy groups (Corollary \ref{cor-long exact sequence h-cofibration})
it suffices to show that the quotient 
\begin{align*}
   ( C\sm G_m A_n ) / ( C\sm G_m A_{n-1})\ &\iso \
C\sm G_m(A_n / A_{n-1}) 
 \end{align*}
is $G$-stably contractible.
Now $A_n/A_{n-1}$ is $(G\times O(m))$-equivariantly isomorphic
to a wedge of summands of the form $S^n\sm (G\times O(m))/\Gamma_+$,
for various closed subgroups $\Gamma$ of $G\times O(m)$.
So the quotient spectrum is isomorphic to the $n$-fold suspension 
of a wedge of orthogonal $G$-spectra of the form $C\sm G_m((G\times O(m))/\Gamma_+)$.
These wedge summands are $G$-stably contractible by Step~3.
Since equivariant homotopy groups take wedges to direct sums
and suspension shifts equivariant homotopy groups,
this proves that $C\sm  G_m( A_n/ A_{n-1})$
is $G$-stably contractible. This completes the inductive step,
and the proof for all finite-dimensional
based $(G\times O(m))$-CW-complexes.

Since $A$ is the sequential colimit, along h-cofibrations
of based $(G\times O(m))$-spaces, of the skeleta $A_n$,
the orthogonal $G$-spectrum
$C\sm G_m A$ is the sequential colimit, along h-cofibrations
of orthogonal $G$-spectra, of the sequence with terms $C\sm G_m A_n$. 
Since h-cofibrations are in particular levelwise closed embeddings, 
equivariant homotopy groups commute with such sequential colimits 
(Proposition \ref{prop:sequential colimit closed embeddings}~(i)), 
so also $C\sm G_m A$ is $G$-stably trivial.

Step 5:
We let $X$ be a $G$-flat orthogonal $G$-spectrum.
We show first, by induction on $m$, 
that $C\sm (\sk^m\! X)$ is $G$-stably contractible,
where $\sk^m \! X$ is the $m$-skeleton 
in the sense of Construction \ref{con:skeleton spec}.
The induction starts with $m=-1$, where there is nothing to show. 
For $m\geq 0$ the morphism
$j_m:\sk^{m-1}\! X\to \sk^m\! X$ is an h-cofibration of orthogonal $G$-spectra, 
hence so is the induced morphism $C\sm j_m:C\sm (\sk^{m-1}\! X)\to C\sm ( \sk^m\! X) $.
By induction and the long exact sequence of equivariant
homotopy groups (Corollary \ref{cor-long exact sequence h-cofibration})
we may show that the quotient 
\begin{align*}
   (C\sm \sk^m X ) / ( C\sm \sk^{m-1} X )\ &\iso \
C \sm (\sk^m X / \sk^{m-1} X)  \ \iso \ C \sm G_m(X(\mR^m)/L_m X) 
 \end{align*}
is $G$-stably contractible.
Since $X$ is $G$-flat, $X(\mR^m)/L_m X$ is a cofibrant based $(G\times O(m))$-space,
and the $G$-equivariant homotopy groups of
$C\sm G_m(X(\mR^m)/L_m X) $ vanish by Step~4.
This finishes the inductive proof that
$C\sm (\sk^m\! X)$ is $G$-stably contractible.
The spectrum $C\sm X$ is the sequential colimit,
along h-cofibrations, of the orthogonal $G$-spectra $C\sm (\sk^m \! X)$;
equivariant homotopy groups commute with such colimits, so 
$C\sm X$ is $G$-stably contractible.

Step~6:
We let $X$ be a $G$-flat orthogonal $G$-spectrum and $f:A\to B$ 
a $\upi_*$-isomorphism of orthogonal $G$-spectra.
Then the mapping cone $C(f)$ is $G$-stably contractible
by the long exact homotopy group sequence
(Proposition \ref{prop:LES for homotopy of cone and fibre}).
So the smash product $C(f)\sm X$ is $G$-stably contractible by Step~5.
Since $C(f)\sm X$ is isomorphic to the mapping cone of $f\sm X:A\sm X\to B\sm X$,
the morphism $f\sm X$ induces isomorphisms on $G$-equivariant stable homotopy groups,
again by the long exact homotopy group sequence.
\end{proof}

Now we turn to products on equivariant homotopy groups.

\begin{construction}\label{con:pairing equivalence homotopy}
Given a compact Lie group $G$ and
two orthogonal $G$-spectra $X$ and $Y$, 
we endow the equivariant homotopy groups with 
a pairing\index{subject}{product!on global homotopy groups}
\begin{equation}\label{eq:def_dot}
\times \ : \ \pi_k^G(X) \ \times \ \pi_l^G (Y) \ \to \ \pi_{k+l}^G(X\sm Y) \ ,
\end{equation}
where $k$ and $l$ are integers. We let 
\[ f\ :\ S^{U\oplus \mR^{m+k}}\ \to\  X(U\oplus\mR^m) \text{\quad and\quad}
g\ :\ S^{V\oplus \mR^{n+l}}\ \to\  Y(V\oplus \mR^n) \]
represent classes in $\pi^G_k(X)$ respectively $\pi^G_l(Y)$,
for suitable $G$-representations $U$ and $V$.
The class $[f]\times[g]$ in $\pi_{k+l}^G(X\sm Y)$
is then represented by the composite
\begin{align*}
 S^{U\oplus V\oplus\mR^{m+n+k+l}} \iso \ S^{U\oplus \mR^{m+k}}\sm S^{V\oplus \mR^{n+l}} \ 
\xra{\ f\sm g\ }\   &X(U\oplus \mR^m)\sm Y(V\oplus\mR^n)\\  
\xra{i_{U\oplus \mR^m, V\oplus\mR^n}} \
&(X\sm Y)(U\oplus\mR^m\oplus V\oplus \mR^n) \\
\xra{(X\sm Y)(U\oplus\tau_{\mR^m,V}\oplus \mR^n)}\  &(X\sm Y)(U\oplus V\oplus \mR^{m+n})\ .
\end{align*}
The first homeomorphism shuffles the sphere coordinates.
We omit the proof that the class of the
composite only depends on the classes of $f$ and $g$.
\end{construction}

The pairing of equivariant homotopy groups has several expected properties
that we summarize in the next theorem.

\begin{theorem}\label{thm:product properties} 
  Let $G$ be a compact Lie group and $X, Y$ and $Z$ orthogonal $G$-spectra.
  \begin{enumerate}[\em (i)]
  \item {\em (Biadditivity)}
    The product $\times:\pi_k^G(X) \times\pi_l^G (Y) \to\pi_{k+l}^G(X\sm Y)$
    is biadditive.
  \item {\em (Unitality)}
    The class $1\in\pi_0^G(\mS)$ represented by the identity of $S^0$
    is a two-sided unit for the pairing $\times$.
  \item {\em (Associativity)} For all $x\in \pi_k^G ( X )$, $y\in \pi_l^G( Y )$
    and $z\in\pi_m^G( Z )$ the relation
    \[ x\times(y\times z) \ = \ (x\times y)\times z\]
    holds in $\pi^G_{k+l+m}(X\sm Y\sm Z)$.
  \item {\em (Commutativity)} For all $x\in \pi_k^G ( X )$ and $y\in \pi_l^G( Y )$
    the relation
    \[ \tau^{X,Y}_*( x\times y ) \ = \ (-1)^{kl}\cdot ( y\times x)\]
    holds in $\pi^G_{k+l}( Y\sm X)$,
    where $\tau^{X,Y}:X\sm Y\to Y\sm X$ is the symmetry isomorphism
    of the smash product.
  \item {\em (Restriction)} For all $x\in \pi_k^G ( X )$ and $y\in \pi_l^G( Y )$
    and all continuous homomorphisms $\alpha:K\to G$ the relation
    \[ \alpha^*(x)\times\alpha^*(y)\ = \ \alpha^*(x\times y) \]
    holds in $\pi_{k+l}^K(\alpha^*(X\sm Y))$.
   \item {\em (Transfer)} Let $H$ be a closed subgroup of $G$.
     For all $x\in\pi_k^G(X)$ and $z\in\pi_l^H(Y\sm S^L)$, the relation
     \[ x\times \Tr_H^G(z)\ = \ \Tr_H^G(\res^G_H(x)\times z) \]
     holds in $\pi_{k+l}^G(X\sm Y)$, where $L=T_{e H}(G/H)$
     is the tangent $H$-representation.
     For all $y\in\pi_l^H(Y)$, the relation
     \[ x\times \tr_H^G(y)\ = \ \tr_H^G(\res^G_H(x)\times y) \]
     holds in $\pi_{k+l}^G(X\sm Y)$.
  \end{enumerate}
\end{theorem}
\begin{proof}
(i)  We deduce the additivity in the first variable from the
general additivity statement in Proposition \ref{prop:additivity prop}.
We consider the two reduced additive functors
\[  X\ \longmapsto\ \pi_k^G(X)\text{\qquad and\qquad} 
X\longmapsto\  \pi_{k+l}^G( X\sm Y) \]
from the category of orthogonal $G$-spectra to the category of abelian groups.
Proposition \ref{prop:additivity prop} shows that for every
$y\in \pi_l^G(Y)$ the natural transformation
\[ -\times y \ : \  \pi_k^G(-)\ \to\  \pi_{k+l}^G( -\sm Y)  \]
is additive.
Additivity in the second variable is proved in the same way.
Properties~(ii), (iii) and~(v) are straightforward consequences
of the definitions.
The sign $(-1)^{k l}$ in the commutativity relation~(iv)
is the degree of the map that interchanges a $k$-sphere with an $l$-sphere.
Indeed, if
\[ f \ : \ S^{U\oplus \mR^{m+k}}\ \to \ X(U\oplus \mR^m)  \text{\quad and\quad}
 g \ : \ S^{V\oplus \mR^{n+l}}\ \to \ Y(V\oplus \mR^n)  \]
represent classes in $\pi_k^G(X)$ respectively $\pi_l^G(Y)$, then 
the left and right vertical composites in the diagram
\[ \xymatrix@C=20mm@R=6mm{
S^{U\oplus V\oplus \mR^{m+n+k+l}}\ar[d]_{f\times g}
\ar[r]^-{\tau_{U,V}\sm\tau_{m,n}\sm S^{k+l}}&
 S^{V\oplus U\oplus \mR^{n+m+k+l}}\ar[d]^{S^{V\oplus U\oplus\mR^{n+m}}\sm\tau_{k,l}}\\
(X\sm Y)(U\oplus V\oplus\mR^{m+n}) \ar[d]_{\tau^{X,Y}(U\oplus V\oplus\mR^{m+n})}&
 S^{V\oplus U\oplus \mR^{n+m+l+k}}\ar[d]^{g\times f}\\
(Y\sm X)(U\oplus V\oplus\mR^{m+n})\ar[r]_-{(Y\sm X)(\tau_{U,V}\oplus\tau_{m,n})} &
 (Y\sm X)(V\oplus U\oplus\mR^{n+m}) 
} \]
represent $\tau^{X,Y}_*([f]\times[g])$, respectively $(-1)^{k l}\cdot([g]\times[f])$.
Since the two composites differ by conjugation by a $G$-equivariant linear isometry,
they represent the same class 
by Proposition \ref{prop:invariant description stable}~(ii).
 
(vi) The following diagram of abelian groups
commutes by naturality of the pairing \eqref{eq:def_dot}, and because 
restriction from $G$ to $H$ is multiplicative:
   \[ \hspace*{-.7cm}\xymatrix@C=6mm@R=6mm{ 
     \pi_k^G(X) \times \pi_l^G(Y\sm G/H_+) 
     \ar@<-10ex>@/_1pc/[dd]_(.2){\Id\times\Wirth_H^G}
     \ar[r]^-{\times} \ar[d]^{\Id\times\res_H^G} & 
  \pi_{k+l}^G(X\sm Y\sm G/H_+) \ar[d]_{\res^G_H}\ar@/^3pc/[dddr]^-{\Wirth_H^G}\\
  \pi_k^G(X) \times \pi_l^H(Y\sm G/H_+) 
  \ar[d]^{ \Id\times (Y\sm l)_*} \ar[dr]^{\res^G_H\times\Id}&  
  \pi_{k+l}^H(X\sm Y\sm G/H_+)
  \ar@/^2pc/[ddr]_(.8){(X\sm Y\sm l)_*}    &\\
     \pi_k^G(X) \times \pi_l^H(Y\sm S^L) \ar[dr]_{\res^G_H\times\Id}  
     &   \pi_k^H(X)\times \pi^H_l(Y\sm G/H_+) \ar[d]^{\Id\times (Y\sm l)_*}\ar[u]_\times
     &  \\  
     & \pi_k^H(X) \times \pi_l^H(Y\sm S^L)       \ar[r]_-{\times} & 
     \pi_{k+l}^H(X\sm Y\sm S^L)   } 
   \]
Here $l:G/H_+\to S^L$ is the $H$-equivariant collapse map \eqref{eq:define_l_A}.
The two vertical composites are the respective Wirthm{\"u}ller isomorphisms
(Theorem \ref{thm:Wirth iso}).
Since the external transfer is inverse to the Wirthm{\"u}ller isomorphism
(up to the effect of the involution $S^{-\Id}:S^L\to S^L$),
we can read the diagram backwards and conclude that the upper part of the
following diagram commutes:
   \[ \xymatrix@C=10mm@R=6mm{ 
     \pi_k^G(X) \times \pi_l^H(Y\sm S^L)    \ar[r]^-{\res^G_H\times\Id}
     \ar@<-8ex>@/_2pc/[ddd]_(.3){\Id\times\Tr_H^G}
     \ar[dd]^{\Id\times (G\ltimes_H-)}_\iso & 
     \pi_k^H(X) \times \pi_l^H(Y\sm S^L)  \ar[d]^-{\times} \\
     & \pi_{k+l}^H(X\sm Y\sm S^L) \ar[d]^\iso_{G\ltimes_H- } \ar@<6ex>@/^2pc/[dd]^(.3){\Tr_H^G}\\
     \pi_k^G(X)\times \pi_l^G(Y\sm G/H_+) \ar[r]_-{\times} \ar[d]^{\Id\times (Y\sm p)_*}&   
     \pi_{k+l}^G(X\sm Y\sm G/H _+) \ar[d]_{(X\sm Y\sm p)_*}\\ 
     \pi_k^G(X)\times \pi_l^G(Y) \ar[r]_-{\times} &  \pi^G_{k+l} ( X\sm Y)
   } \]
   Here $p:G/H_+\to S^0$ is the projection.
   The lower part of the diagram commutes by naturality of the pairing. 
   This proves the first claim about dimension shifting transfers.
   The formula for the degree zero transfers follows by naturality 
   for the inclusion of the origin $S^0\to S^L$.
\end{proof}

\begin{defn}
An {\em orthogonal ring spectrum}
\index{subject}{orthogonal ring spectrum} 
\index{subject}{ring spectrum!orthogonal|see{orthogonal ring spectrum}} 
is a monoid in the category of orthogonal spectra with respect to the smash product.
For a compact Lie group $G$, an
{\em orthogonal $G$-ring spectrum}
\index{subject}{orthogonal $G$-ring spectrum} 
\index{subject}{G-ring spectrum@$G$-ring spectrum!orthogonal|see{orthogonal $G$-ring spectrum}} 
is a monoid in the category of orthogonal $G$-spectra with respect to the smash product.
\end{defn}

An orthogonal ring spectrum is thus an orthogonal spectrum $R$ 
equipped with a multiplication morphism
$\mu:R\sm R\to R$ and a unit morphism $\eta:\mS\to R$ such that the associativity and
unit diagrams commute (compare \cite[VII.3]{maclane-working}).
A {\em morphism}\index{subject}{morphism!of orthogonal ring spectra}
of orthogonal ring spectra is a morphism $f:R\to S$ of orthogonal spectra
that satisfies $f\circ\mu^R=\mu^S\circ(f\sm f)$ and $f\circ\eta^R=\eta^S$.
Via the universal property of the smash product the data contained in the
multiplication morphism can be made more explicit: $\mu:R\sm R\to R$ corresponds to
a collection of based continuous maps $\mu_{V,W}:R(V)\sm R(W)\to R(V\oplus W)$
that together form a bimorphism. The associativity and unit conditions
can also be rephrased in more explicit forms, and then we are requiring
that the multiplication and unit maps make $R:\bO\to\bT_*$ into a lax monoidal functor.
Most of the time we will specify the data of an orthogonal ring spectrum
in the explicit bimorphism form.

An orthogonal ring spectrum $R$ (respective orthogonal $G$-ring spectrum) is 
{\em commutative}\index{subject}{orthogonal ring spectrum!commutative} if the 
multiplication morphism satisfies $\mu\circ\tau_{R,R}=\mu$. In the explicit form
this is equivalent to the commutativity of the square
\[\xymatrix@C=12mm{ R(V)\sm R(W) \ar[r]^-{\mu_{V,W}} \ar[d]_{\tau_{R(V), R(W)}} &
R(V\oplus W)\ar[d]^-{R(\tau_{V,W})} \\
R(W)\sm R(V) \ar[r]_-{\mu_{W,V}} & R(W\oplus V) }\]
for all inner product spaces $V$ and $W$.
Equivalently, the multiplication and unit maps make $R:\bO\to\bT_*$ 
into a lax {\em symmetric} monoidal functor.
Commutative orthogonal ring spectra already appear, 
with an extra pointset topological hypothesis and
under the name ${\mathscr I}_*$-prefunctor, in \cite[IV Def.\,2.1]{may-quinn-ray}.

Since the smash product of orthogonal $G$-spectra is just
the smash product of the underlying orthogonal spectra, endowed with
the diagonal $G$-action, an orthogonal $G$-ring spectrum is nothing
but an orthogonal ring spectrum equipped with a continuous $G$-action
through homomorphisms of orthogonal ring spectra.
Via the universal property of the smash product,
yet another way to package the data in an 
orthogonal $G$-ring spectrum is as a continuous lax monoidal functor
from the category $\bO$ to the category of based $G$-spaces.

\medskip

Given a compact Lie group $G$ and an orthogonal $G$-ring spectrum $R$, 
we define an internal pairing
\begin{equation}\label{eq:multiplication_upi_ringspectrum}
 \cdot \ : \   \pi_k^G (R) \times \pi_l^G (R) \ \to\ \pi_{k+l}^G(R) 
\end{equation}
on the equivariant homotopy groups of $R$ as the composite
\[
  \pi_k^G (R) \times \pi_l^G (R) \ \xra{\ \times \ } \
\pi_{k+l}^G(R\sm R) \ \xra{\ \mu_*\ }\ \pi_{k+l}^G(R)\ .  
\]
Theorem \ref{thm:product properties} then immediately implies:

\begin{cor}\label{cor:product properties ring spectrum} 
  Let $G$ be a compact Lie group, $R$ an orthogonal $G$-ring spectrum
  and $H$ a closed subgroup of $G$.
  \begin{enumerate}[\em (i)]
  \item The products \eqref{eq:multiplication_upi_ringspectrum}
    make the abelian groups $\{\pi_k^G(R)\}_{k\in\mZ}$ into a graded ring.
    The multiplicative unit is the class of the unit map $S^0\to R(0)$.
  \item If the multiplication of $R$ is commutative, then the relation
    \[ x\cdot y\ = \  (-1)^{kl}\cdot y \cdot x \]
    holds for all classes $x\in \pi_k^G(R)$ and $y\in \pi^G_l(R)$. 
 \item {\em (Restriction)} The restriction maps
    $\res^G_H:\pi_*^G(R)\to\pi_*^H(R)$ form a homomorphism of graded rings.
 \item {\em (Conjugation)} For every $g\in G$
   the conjugation maps $g_\star:\pi_*^{H^g}(R)\to\pi_*^H(R)$ 
   form a homomorphism of graded rings.\index{subject}{conjugation homomorphism!on equivariant homotopy groups}  
  \item {\em (Reciprocity)} 
     For all $x\in\pi_k^G(R)$ and all $y\in\pi_l^H(R)$, the relation
     \[ x\cdot \tr_H^G(y)\ = \ \tr_H^G(\res^G_H(x)\cdot y) \]
     holds in $\pi_{k+l}^G(R)$.
  \end{enumerate}
\end{cor}

\begin{rk}[Norm maps]
A lot more is happening for commutative orthogonal $G$-ring spectra:
strict commutativity of the multiplication not only makes the
homotopy pairings graded-commutative, but it also gives rise to new operations,
usually called {\em norm maps}\index{subject}{norm map}
$N_H^G:\pi_0^H(R)\to \pi_0^G(R)$ for all closed subgroups
$H$ of finite index in $G$.
To my knowledge, these norm maps were first defined by
Greenlees and May \cite[Def.\,7.10]{greenlees-may-completion},
and they were later studied for example in \cite{brun-Witt cobordism, HHR-Kervaire}.

We don't discuss norms maps for commutative orthogonal $G$-ring spectra here,
but norms maps and power operations will be a major topic 
in the context of ultra-commutative ring spectra.
In Section \ref{sec:power ops ring spectra} we define
power operations on the 0-th equivariant homotopy groups of
an ultra-commutative ring spectrum, and we recall in 
Remark \ref{rk:power vs Tambara} how to turn the power operations into norm maps.
\end{rk}

\chapter{Global stable homotopy theory}
\label{ch-stable}

In this chapter we embark on the investigation of
global stable homotopy theory.
In Section \ref{sec:orthogonal spectra model global} we specialize the equivariant 
theory of the previous chapter to global stable homotopy types,
which we model by orthogonal spectra (with no additional action of any groups). 
Section \ref{sec:global functors} introduces the category of {\em global functors},
the natural home of the collection of equivariant
homotopy groups of a global stable homotopy type (i.e., an orthogonal spectrum).
Global functors play the same role
for global homotopy theory as the category of abelian groups
in ordinary homotopy theory, or the category of $G$-Mackey functors
for $G$-equivariant homotopy theory.
Global functors are defined as additive functors on the {\em global Burnside category};
an explicit calculation of the global Burnside category 
provides the link to other notions of global Mackey functors.
In the global context, the pairings on equivariant homotopy groups also
give rise to a symmetric monoidal structure on the global Burnside category,
and to a symmetric monoidal `box product' of global functors.

Section \ref{sec:global model structures} establishes the global model structure
on the category of orthogonal spectra; more generally,
we consider a global family $\Fc$ and define the
$\Fc$-global model structure, where weak equivalences are tested
on equivariant homotopy groups for all Lie groups in $\Fc$.
The  $\Fc$-global model structure is monoidal with
respect to the smash product of orthogonal spectra,
provided that $\Fc$ is closed under products.
Section \ref{sec:generators} collects aspects of global stable homotopy theory that
refer to the triangulated structure of the global stable homotopy category.
Specific topics are  compact generators, Brown representability,
a t-structure whose heart is the category of global functors,
global Postnikov sections and Eilenberg-Mac Lane spectra of global functors.

The final Section \ref{sec:change family} is a systematic study of the effects
of changing the global family. We show that 
the `forgetful' functor between the global stable homotopy categories
of two nested global families has fully faithful left and right adjoints, 
which are part of a recollement.
We provide characterizations of the global homotopy types in the image of the 
two adjoints; for example,
the right adjoint all the way from the non-equivariant to the
global stable homotopy category models Borel cohomology theories.
We also relate the global homotopy category to the $G$-equivariant stable homotopy
category for a fixed compact Lie group $G$;
here the forgetful functor also has both adjoints, but these are no longer
fully faithful if the group is non-trivial.

Also in  Section \ref{sec:change family} we
establish an algebraic model for {\em rational $\Fin$-global}
stable homotopy theory, i.e., for rational
global stable homotopy theory based on the global family of finite groups.
Indeed, spectral Morita theory provides a chain of Quillen equivalences 
to the category of chain complexes of rational global functors on finite groups.
On the algebraic side, the abelian category of rational $\Fin$-global functors 
is equivalent to an even simpler category, namely functors from finite groups 
and conjugacy classes of epimorphisms to $\mQ$-vector spaces.
Under the two equivalences, homology groups of chain complexes
correspond to equivariant stable homotopy groups,
respectively geometric fixed point homotopy groups, of spectra.

\section{Orthogonal spectra as global homotopy types}
\label{sec:orthogonal spectra model global}

In this section we specialize the equivariant stable homotopy
theory of Chapter \ref{ch-equivariant} to global stable homotopy types,
which we model by orthogonal spectra (with no additional group action). 
Given an orthogonal spectrum $X$ and a compact Lie group $G$,
we obtain an orthogonal $G$-spectrum by letting $G$ act trivially
on the values of $X$.
We call this the {\em underlying orthogonal $G$-spectrum} of $X$
and denote it $X_G$.
\index{subject}{underlying $G$-spectrum!of an orthogonal spectrum}
\index{symbol}{$X_G$ - {underlying $G$-spectrum of an orthogonal spectrum}}
To simplify notation we omit the subscript~`$_G$' when we refer
to equivariant homotopy groups, i.e., we simply write $\pi_k^G(X)$
instead of $\pi_k^G(X_G)$.

For global homotopy types (i.e., orthogonal spectra),
the notation related to restriction maps simplifies, and some special features happen.
Indeed, if $X$ is an orthogonal spectrum, then for every continuous homomorphism
$\alpha:K\to G$ we have $\alpha^*(X_G)=X_K$, because both $K$ and $G$ act trivially.
So for global homotopy types the restriction maps become homomorphisms
\[ \alpha^* \ : \ \pi_*^G(X) \ \to \ \pi_*^K(X)\ .\]
In particular, we have these restriction maps when $\alpha$ is an epimorphism;
in that case we refer to $\alpha^*$ as an {\em inflation map}.
\index{subject}{inflation map!of equivariant point homotopy groups}  

If $H$ is a closed subgroup of $G$ and $g\in G$, the
conjugation homomorphism $c_g:H\to H^g$ is given by $c_g(h)=g^{-1}h g$.
We recall that for an orthogonal $G$-spectrum $Y$ the conjugation morphism
$g_\star:\pi_0^{H^g}(Y)\to\pi_0^H(Y)$ is defined as
$g_\star=(l_g^Y)_*\circ (c_g)^\ast$, where $l_g^Y:c_g^*(Y)\to Y$ is left multiplication by $g$,
compare \eqref{eq:define_G_star}.
If $Y=X_G$ for an orthogonal spectrum $X$, then $l_g^X$ is the identity;
so for global homotopy types we have
\begin{equation}\label{eq:c_g=g_star}
 g_\star \ = \ (c_g)^* \ : \  \pi_0^{H^g}(X) \ \to \ \pi_0^H(X)\ .
\end{equation}
A different way to express the significance of the relation \eqref{eq:c_g=g_star}
is that whenever $Y=X_G$ underlies a global homotopy type, 
then the action of the group $W_G H=(N_G H)/H$
on $\pi_0^H(X)$ discussed in Construction \ref{rk:Weyl on pi_0^G}
factors through an action of the quotient group
\[ (N_G H) / (H\cdot C_G H)\ .\]
The difference is illustrated most drastically for $H=e$, the trivial subgroup of $G$.
For a general orthogonal $G$-spectrum the action of $G=W_G e$
on $\pi_0^G(Y)$ is typically non-trivial. If $Y=X_G$ arises from
an orthogonal spectrum, then this $G$-action is trivial.

\begin{rk}[Global homotopy types are split $G$-spectra]\label{rk:homotopy properties}
Obviously, only very special orthogonal $G$-spectra $Y$ are part of
a `global family', i.e., arise as $X_G$ for an orthogonal spectrum $X$. 
However, it is not a priori clear what the homotopical
significance of the pointset level condition
that $G$ must act trivially on the values of $Y$ (at trivial representations) is.
Now we formulate obstructions to `being global' in terms 
of the Mackey functor homotopy groups of an orthogonal $G$-spectrum.\index{subject}{Mackey functor}

The equivariant homotopy groups $\pi_0^G(X)=\pi_0^G(X_G)$ of
an orthogonal spectrum $X$ come equipped with restriction maps along arbitrary
continuous group homomorphisms, not necessarily injective. 
This is in contrast to the situation for
a fixed compact Lie group, where one can only restrict to subgroups,
or along conjugation maps by elements of the ambient group,
but where there are no inflation maps.
One obstruction to $Y$ being part of a `global family'
is that the $G$-Mackey functor structure can be extended to
a `global functor' (in the sense of Definition \ref{def:global functor} below).
In particular, the $G$-Mackey functor homotopy groups
can be complemented by restriction maps 
along arbitrary group homomorphisms between the subgroups of $G$. 
So every global functor takes isomorphic
values on a pair of isomorphic subgroups of $G$;
for $G$-Mackey functors this is true when the subgroups are conjugate in $G$, 
but not in general when they are merely abstractly isomorphic.

As an extreme case the global structure includes an inflation map
$p_G^*:\pi_0^e(X)\to \pi_0^G(X)$ associated to the unique homomorphism $p_G:G\to e$,
splitting the restriction map $\pi_0^G(X)\to \pi_0^e(X)$. 
So one obstruction to being global is
that this restriction map from $\pi_0^G(X)$ to $\pi_0^e(X)$
needs to be a split epimorphism.
This is the algebraic shadow of the fact that the $G$-equivariant homotopy types
that underlie global homotopy types (i.e., are represented by orthogonal spectra)
are `split'\index{subject}{split $G$-spectrum} in the sense that there
is a morphism from the underlying non-equivariant spectrum to
the genuine $G$-fixed point spectrum that splits the restriction map.
\end{rk}

\begin{defn}\label{def:global equivalence}
A morphism $f:X\to Y$ of orthogonal spectra 
is a {\em global equivalence}\index{subject}{global equivalence!of orthogonal spectra}
if the induced map $\pi_k^G(f):\pi_k^G(X) \to \pi_k^G(Y)$ is an isomorphism
for all compact Lie groups $G$ and all integers $k$.
\end{defn}

We define the {\em global stable homotopy category} $\GH$\index{symbol}{  $\GH$ - {global stable homotopy category}}\index{subject}{global stable homotopy category}
by localizing the category of orthogonal spectra at the class of
global equivalences.
The global equivalences are the weak equivalences of the
{\em global model structure} on the category of orthogonal spectra,
see Theorem \ref{thm:All global spectra} below. 
So the methods of homotopical algebra are available for studying global
equivalences and the associated global homotopy category.
In later sections we will also consider a relative notion of
global equivalence, the `$\Fc$-equivalences', 
based on a global family $\Fc$\index{subject}{global family}
of compact Lie groups. There we require that 
the induced map $\pi_k^G (f): \pi_k^G (X)\to\pi_k^G(Y)$ is an isomorphism
for all integers $k$ and all compact Lie groups $G$ that belong to
the global family $\Fc$.

When we specialize Proposition \ref{prop:lambda upi_* isos}~(ii),
Corollary \ref{cor-wedges and finite products},
Proposition \ref{prop:map(A,-) preserves global},
Proposition \ref{prop:space smash preserves global}~(ii)
and Proposition \ref{prop:sequential colimit closed embeddings}
to the underlying orthogonal $G$-spectra of orthogonal spectra, we obtain
the following consequences:

\begin{prop}\label{prop:global equiv preservation}
  \begin{enumerate}[\em (i)]
  \item 
    For every orthogonal spectrum $X$ the morphism
    \[  \lambda_X\ : \  X\sm S^1 \ \to\ \sh X  \ ,
    \text{\quad its adjoint\quad}
    \tilde\lambda_X \ : \ X \ \to \ \Omega \sh X \ ,\]
    the adjunction unit $\eta_X:X\to\Omega(X\sm S^1)$ 
    and the adjunction counit $\epsilon_X:(\Omega X)\sm S^1\to X$ 
    are global equivalences.
  \item 
    For every finite family of orthogonal spectra
    the canonical morphism from the wedge to the product 
    is a global equivalence.
  \item
    For every finite based CW-complex $A$, the functor $\map_*(A,-)$
    preserves global equivalences of orthogonal spectra.
  \item For every cofibrant based space $A$, 
    the functor $-\sm A$ preserves global equivalences of orthogonal spectra.
  \item Let $e_n:X_n\to X_{n+1}$ and $f_n:Y_n\to Y_{n+1}$ be morphisms 
    of orthogonal spectra that are levelwise closed embeddings, for $n\geq 0$. 
    Let $\psi_n:X_n\to Y_n$ be global equivalences
    that satisfy $\psi_{n+1}\circ e_n=f_n\circ\psi_n$ for all $n\geq 0$.
    Then the induced morphism $\psi_\infty:X_\infty\to Y_\infty$ 
    between the colimits of the sequences is a global equivalence.
  \item Let $f_n:Y_n\to Y_{n+1}$ be global equivalences of orthogonal spectra
   that are levelwise closed embeddings, for $n\geq 0$. 
   Then the canonical morphism 
   $f_\infty:Y_0\to Y_\infty$ to a colimit of the sequence $\{f_n\}_{n\geq 0}$
   is a global equivalence.
  \end{enumerate}
\end{prop}

For later use we record another closure property of global equivalences.

\begin{cor}\label{cor-closure h-cof and global equi}
The class of h-cofibrations of orthogonal spectra
that are simultaneously global equivalences
is closed under cobase change, coproducts and sequential compositions.
\end{cor}
\begin{proof}
The class of h-cofibrations is closed under coproducts, cobase change
and composition (finite or sequential), 
compare Corollary \ref{cor-h-cofibration closures}~(i).
The class of global equivalences is closed under coproducts
because equivariant homotopy groups take wedges to direct sums 
(Corollary \ref{cor-wedges and finite products}~(i)).
The cobase change of an h-cofibration that is also a global equivalence is 
another global equivalence 
by Corollary \ref{cor-cobase change h-cofibration}~(i).
Every h-cofibration of orthogonal spectra is in particular levelwise
a closed embedding. So the class of h-cofibrations that are also global equivalences
is closed under sequential 
composition by Proposition \ref{prop:global equiv preservation}~(vi).
\end{proof}

For an orthogonal spectrum $X$ the 0-th equivariant homotopy groups $\pi_0^G(X)$
and the restriction maps between them coincide with the homotopy Rep-functor,
in the sense of \eqref{eq:define_pi_0_set},
of a certain orthogonal space $\Omega^\bullet X$ that we now recall.

\begin{construction}\label{con:Omega^bullet}
We introduce the functor
\[ \Omega^\bullet \ : \ \spec \ \to \ \spc \]
from orthogonal spectra to orthogonal spaces.\index{symbol}{$\Omega^\bullet$ - {underlying orthogonal space of an orthogonal spectrum}}
Given an orthogonal spectrum $X$, the value of 
$\Omega^\bullet X$ at an inner product space $V$ is
\[ (\Omega^\bullet X)(V)\ = \ \map_*(S^V,X(V))\ . \]
If $\varphi:V\to W$ is a linear isometric embedding, the induced map
\[ \varphi_* \ : \ (\Omega^\bullet X)(V)\ = \ \map_*(S^V,X(V)) \ \to \ 
 \map_*(S^W,X(W))\ = \ (\Omega^\bullet X)(W) \]
was defined in \eqref{eq:define varphi_* f}.
In particular, the orthogonal group $O(V)$ acts on 
$(\Omega^\bullet X)(V)=\map_*(S^V,X(V))$ by conjugation.
If $\psi:U\to V$ is another isometric embedding, then we have
$ \varphi_*(\psi_* f)  =  (\varphi\psi)_*f$.
The assignment $(\varphi,f)\mapsto \varphi_*f$ is continuous in both variables,
so we have really defined an orthogonal space $\Omega^\bullet X$.
The construction is clearly functorial in the orthogonal spectrum $X$;
moreover, $\Omega^\bullet$ has a left adjoint `unreduced suspension spectrum'
functor $\Sigma^\infty_+$ that we discuss 
in Construction \ref{con:suspension spectrum} below.

If $G$ acts on $V$ by linear isometries, then the 
$G$-fixed subspace of  $(\Omega^\bullet X)(V)$ is the space
of $G$-equivariant based maps from $S^V$ to $X(V)$:
\[ ((\Omega^\bullet X)(V))^G\ = \ \map_*^G(S^V,X(V))\ . \]
The path components of this space are precisely the equivariant homotopy
classes of based $G$-maps, i.e.,
\[ \pi_0\left(((\Omega^\bullet X)(V))^G\right)\ = \ 
\pi_0 \left( \map_*^G(S^V,X(V)) \right) \ = \ [S^V,X(V)]^G\ . \]
Passing to colimits over the poset $s(\Uc_G)$ gives
\[ \pi_0^G(\Omega^\bullet X)\ = \ \pi_0^G(X)\ , \]
i.e., the $G$-equivariant homotopy group of the orthogonal spectrum $X$
equals the $G$-equivariant homotopy set 
of the orthogonal space $\Omega^\bullet X$.
A direct inspection shows that the spectrum level restriction maps
defined in Construction \ref{con:restriction map G-spectra}\index{subject}{restriction map!of equivariant stable homotopy groups}
coincide with the restriction maps for orthogonal spaces
introduced in \eqref{eq:define_alpha^*}.\end{construction}

\begin{construction}[Suspension spectra of orthogonal spaces]\label{con:suspension spectrum} 
To every orthogonal space $Y$ we can associate an unreduced\index{symbol}{$\Sigma^\infty$ - {suspension spectrum of an orthogonal space}}
{\em suspension spectrum}\index{subject}{suspension spectrum!of an orthogonal space}
 $\Sigma^\infty_+Y$ whose value on an inner product space
\[ (\Sigma^\infty_+ Y)(V)\ = \ S^V\sm Y(V)_+\ ; \]
here the orthogonal group $O(V)$ acts diagonally and the
structure map
\[  \sigma_{V,W}\ : \ 
S^V\sm (S^W\sm Y(W)_+) \ \to \  S^{V\oplus W}\sm Y(V\oplus W)_+  \]
is the combination of the canonical homeomorphism $S^V\sm S^W\iso S^{V\oplus W}$
and the map $Y(i_W):Y(W)\to Y(V\oplus W)$.
If $Y$ is the constant orthogonal space with value a topological space $A$,
then $\Sigma^\infty_+Y=\Sigma^\infty_+A$ specializes to the suspension spectrum
of $A$ with a disjoint basepoint added.
The functor
\[ \Sigma^\infty_+ \ : \ \spc \ \to \ \spec \]
is left adjoint to the functor $\Omega^\bullet$ 
of Construction \ref{con:Omega^bullet}.
\end{construction}

Let $G$ be a compact Lie group and $\{V_i\}_{i\geq 1}$ 
an exhaustive sequence\index{subject}{exhaustive sequence!of representations} 
of finite-dimensional $G$-subrepresentations of the complete universe $\Uc_G$. 
Given an orthogonal space $Y$, we denote by
$\tel_i Y(V_i)$ the mapping telescope of the sequence of $G$-spaces
\[ Y(V_1) \ \to \ Y(V_2) \ \to \ \cdots \ \to \ Y(V_i) \ \to \ \cdots \ ;\]
the maps in the sequence are induced by the inclusions,
so they are $G$-equi\-variant, and the telescope inherits
a natural $G$-action. The canonical maps $Y(V_j)\to \tel_i Y(V_i)$
induce maps of equivariant homotopy classes
\begin{align*}
  [S^{V_j\oplus\mR^k}, (\Sigma^\infty_+ Y)(V_j)]^G \ = \  
  &[S^{V_j\oplus\mR^k},S^{V_j}\sm Y(V_j)_+]^G \\ \to \  
  &[S^{V_j\oplus\mR^k},S^{V_j}\sm \tel_i Y(V_i)_+]^G \
  \to \quad \pi_k^G(\Sigma^\infty_+ \tel_i Y(V_i))\ ,
\end{align*}
where the second map is the canonical one to the colimit.
These maps are compatible with stabilization in the source when we increase $j$.
Since the exhaustive sequence is cofinal in the poset $s(\Uc_G)$,
the colimit over $j$ calculates the $k$-th equivariant homotopy group
of $\Sigma^\infty_+ Y$. So altogether, the maps assemble into a natural 
group homomorphism
\[
\pi_k^G(\Sigma^\infty_+ Y) \ \to \  \pi_k^G(\Sigma^\infty_+ \tel_i Y(V_i) ) \ ,
 \]
for $k\geq 0$. For negative values of $k$, we obtain a similar map
by inserting $\mR^{-k}$ in second variable of the above sets of equivariant homotopy
classes.

\begin{prop}\label{prop:G-homotopy of spc}
Let $G$ be a compact Lie group and $\{V_i\}_{i\geq 1}$ an exhaustive sequence 
of $G$-representations. Then for every orthogonal space $Y$ 
and every integer $k$ the map
\[
\pi_k^G(\Sigma^\infty_+ Y) \ \to \  \pi_k^G(\Sigma^\infty_+ \tel_i Y(V_i) ) 
 \]
is an isomorphism.
\end{prop}
\begin{proof}
  The mapping telescope is the colimit, along h-cofibrations,
  of the truncated mapping telescopes $\tel_{[0,n]} Y(V_i)$.
  So the space $S^{V_j}\sm \tel_i Y(V_i)_+$ 
  is the co\-limit, along a sequence of closed embeddings,
  of the spaces $S^{V_j}\sm \tel_{[0,n]} Y(V_i)_+$. 
  For every compact based $G$-space $A$ the canonical map
  \[ \colim_{n\geq 1}\, [A, S^{V_j}\sm \tel_{[0,n]} Y(V_i)_+ ]^G \ \to \ 
  [A, S^{V_j}\sm \tel_i Y(V_i)_+ ]^G  \]
  is thus bijective. Indeed, every continuous map from the compact spaces $A$
  and $A\sm [0,1]_+$ to $S^{V_j}\sm \tel_i Y(V_i)_+$ factors through 
  $S^{V_j}\sm \tel_{[0,n]} Y(V_i)$ some $n$, 
  for example by Proposition \ref{prop:filtered colim preserve weq}.
  The canonical map $Y(V_n)\to \tel_i Y(V_i)$
  factors through $Y(V_n)\to \tel_{[0,n]} Y(V_i)$, and this factorization
  is an equivariant homotopy equivalence. So we can replace 
  $\tel_{[0,n]} Y(V_i)$ by $Y(V_n)$ and conclude that the map
  \[ \colim_{n\geq 1}\, [A, S^{V_j}\sm Y(V_n)_+ ]^G \ \to \ 
  [A, S^{V_j}\sm \tel_i Y(V_i)_+ ]^G  \]
  is bijective.
  We specialize to $A=S^{V_j\oplus\mR^k}$ and pass to the colimit
  over $j$ by stabilization in source and target. The result is a bijection
  \begin{align*}
    \colim_{j\geq 1}  \colim_{n\geq 1}\, &[S^{V_j\oplus\mR^k}, S^{V_j}\sm Y(V_n)_+ ]^G \\
  &\to \ \colim_{j\geq 1}    [S^{V_j\oplus\mR^k}, S^{V_j}\sm \tel_i Y(V_i)_+ ]^G \ . 
  \end{align*}
  The source is isomorphic to the diagonal colimit 
  \[ \colim_{j\geq 1}   [S^{V_j\oplus\mR^k}, S^{V_j}\sm Y(V_j)_+ ]^G  \]
  by cofinality.
  Since the exhaustive sequence is cofinal in the poset $s(\Uc_G)$,
  the two colimits over $j$ also calculate the colimits over $s(\Uc_G)$,
  and hence the $k$-th equivariant homotopy group of 
  $\Sigma^\infty_+ Y$ respectively $\Sigma^\infty_+ \tel_i Y(V_i)$.
  This shows the claim for $k\geq 0$. 
  For negative values of $k$, the argument is similar: 
  we insert $\mR^{-k}$ into the second variable 
  of the sets of equivariant homotopy classes.
\end{proof}

\begin{cor}\label{cor:suspension spectrum globally homotopical}
The unreduced suspension spectrum functor
$\Sigma^\infty_+$ takes global equivalences of orthogonal spaces
to global equivalences of orthogonal spectra.   
\end{cor}
\begin{proof}
Let $f:X\to Y$ be a global equivalence of orthogonal spaces and $G$
a compact Lie group. We choose an exhaustive sequence 
of $G$-representations $\{V_i\}_{i\geq 1}$. Then the $G$-map
\[ \tel_i f(V_i)\ : \ \tel_i X(V_i)\ \to \ \tel_i Y(V_i) \]
is a $G$-weak equivalence
by Proposition \ref{prop:telescope criterion}~(iii). 
So the map of suspension spectra $\Sigma^\infty_+\tel_i f(V_i)$
induces isomorphisms on all $G$-equivariant stable homotopy groups, by
Proposition \ref{prop:suspension spectrum homotopical}.
Since the group $\pi_k^G(\Sigma^\infty_+\tel_i X(V_i))$
is naturally isomorphic to the group $\pi_k^G(\Sigma^\infty_+ X)$
(by Proposition \ref{prop:G-homotopy of spc}),
the morphism of orthogonal suspension 
spectra $\Sigma^\infty_+f:\Sigma^\infty_+X\to\Sigma^\infty_+ Y$
induces an isomorphism on $\pi_*^G$.
\end{proof}

Now we let $Y$ be an orthogonal space and $K$ a compact Lie group.
We define a map
\begin{equation}  \label{eq:sigma_map}
 \sigma^K\ : \ \pi_0^K (Y) \ \to \ \pi_0^K(\Sigma^\infty_+ Y)  
\end{equation}
as the effect of the adjunction unit $Y\to\Omega^\bullet(\Sigma^\infty_+ Y)$
on the  $K$-equivariant homotopy set $\pi_0^K$, using that
$\pi_0^K(\Omega^\bullet(\Sigma^\infty_+ Y))=\pi_0^K(\Sigma^\infty_+ Y)$.
If we unravel the definitions, this comes out as follows: 
if $V$ is a finite-dimensional $K$-subrepresentation
of the complete $K$-universe $\Uc_K$ and $y\in Y(V)^K$ a $K$-fixed point, then
$\sigma^K[y]$ is represented by the $K$-map
\[ S^V\ \xra{-\sm y} \ S^V\sm Y(V)_+ \ = \ (\Sigma^\infty_+ Y)(V) \ .\]
As $K$ varies, the maps $\sigma^K$ are compatible with restriction
along continuous homomorphisms,
since they arise from a morphism of orthogonal spaces.
By the same argument as for orthogonal $K$-spectra 
in \eqref{eq:Weyl_annihilated_transfer},
transferring from a closed subgroup $L$ to $K$ 
annihilates the action of the Weyl group on $\pi_0^L(Y)$.

The next proposition is a global analog
of Theorem \ref{thm:G-equivariant pi_0 of Sigma^infty}~(i). 
We call an orthogonal spectrum {\em globally connective}\index{subject}{globally connective}
if its $G$-equivariant homotopy groups vanish in negative dimensions,
for all compact Lie groups $G$.

\begin{prop}\label{prop:pi_0 of Sigma^infty} 
Let $Y$ be an orthogonal space. 
Then the suspension spectrum $\Sigma^\infty_+ Y$ is globally connective.
Moreover, for every compact Lie group $K$ the equivariant homotopy group
$\pi_0^K(\Sigma^\infty_+ Y)$ is a free abelian group with a basis 
given by the elements
\[ \tr_L^K(\sigma^L(x)) \ ,\]
where $L$ runs through all conjugacy classes of closed subgroups of $K$
with finite Weyl group and $x$ runs through a set of representatives 
of the $W_K L$-orbits of the set $\pi_0^L(Y)$.
\end{prop}
\begin{proof}
We consider the functor on the product poset $s(\Uc_K)^2$ sending
$(U,V)$ to the set $[S^V, S^V\sm Y(U)_+]^K$.
The diagonal is cofinal in $s(\Uc_K)^2$, so the induced map 
\begin{align*}
   \pi_0^K(\Sigma^\infty_+ Y)\ = \ \colim_{V\in s(\Uc_G)} 
 &[S^V, S^V\sm Y(V)_+]^K \\ 
&\to \ \colim_{(V,U)\in s(\Uc_K)^2}  [S^V, S^V\sm Y(U)_+]^K 
\end{align*}
is an isomorphism.
The target can be calculated in two steps, hence the group we are after
is isomorphic to
\[ \colim_{U\in s(\Uc_K)} \left( \colim_{V\in s(\Uc_K)}  [S^V, S^V\sm Y(U)_+]^K \right)
\ = \ \colim_{U\in s(\Uc_K)}  \pi_0^K\left(\Sigma^\infty_+ Y(U) \right)\ .\]
We may thus show that the latter group is free abelian with the specified basis.

Theorem \ref{thm:G-equivariant pi_0 of Sigma^infty}~(i) identifies 
the equivariant homotopy group $\pi_0^K\left(\Sigma^\infty_+ Y(U) \right)$
as the free abelian group with basis the classes $\tr_L^K(\sigma^L(x))$,
where $L$ runs through all conjugacy classes of closed subgroups of $K$
with finite Weyl group and $x$ runs through a set of representatives 
of the $W_K L$-orbits of the set $\pi_0(Y(U)^L)$.
Passage to the colimit over $U\in s(\Uc_K)$ yields the proposition.
\end{proof}

We let $G$ be a compact Lie group and  $V$ a  $G$-representation.
We recall from \eqref{eq:tautological_class} the tautological class
\[ u_{G,V}\ \in \ \pi_0^G (\bL_{G,V}) \]
in the $G$-equivariant homotopy set of the semifree orthogonal space $\bL_{G,V}$.
The {\em stable tautological class}\index{subject}{stable tautological class}\index{subject}{tautological class!stable|see{stable tautological class}}\index{symbol}{$e_{G,V}$ - {stable tautological class in $\pi_0^G(\Sigma^\infty_+ \bL_{G,V})$}} is
\begin{equation}\label{eq:define_stable_tautological}
 e_{G,V}\ = \ \sigma^G(u_{G,V})\ \in \ \pi_0^G(\Sigma^\infty_+ \bL_{G,V}) \ .  
\end{equation}
Explicitly, $e_{G,V}$ is the homotopy class of the $G$-map
\[  S^V \ \to \ S^V\sm (\bL(V,V)/G)_+ \ = \ (\Sigma^\infty_+ \bL_{G,V})(V) \ , \quad 
v\ \longmapsto v\sm (\Id_V\cdot G) \ .\]

In the following corollary we index certain homotopy classes by pairs $(L,\alpha)$
consisting of a closed subgroup $L$ of a compact Lie group $K$
and a continuous homomorphism $\alpha:L\to G$.
The {\em conjugate} of $(L,\alpha)$ by a pair $(k,g)\in K\times G$ is the pair
$(L^k,c_g\circ \alpha\circ c_k^{-1})$ consisting of the
conjugate subgroup $L^k$ and the composite homomorphism
\[  L^k \ \xra{\ c_k^{-1}\ }\ L \ \xra{\ \alpha\  }\ 
G \ \xra{\ c_g\ } \ G \ . \]

\begin{cor}\label{cor-pi_0 of B_gl}\index{subject}{global classifying space}  
Let $G$ and $K$ be compact Lie groups and $V$ a faithful $G$-representation. 
Then the homotopy group $\pi_0^K(\Sigma^\infty_+ \bL_{G,V})$
is a free abelian group with basis given by the classes
\[ \tr_L^K(\alpha^*(e_{G,V})) \]
as  $(L,\alpha)$ runs over a set of representatives of
all $(K\times G)$-conjugacy classes of pairs consisting of a closed subgroup $L$ of $K$ 
with finite Weyl group and a continuous homomorphism $\alpha:L\to G$. 
\end{cor}
\begin{proof} 
The map
\[  \Rep(K,G) \ \to \ \pi_0^K(\bL_{G,V})\ , \quad
[\alpha:K\to G] \ \longmapsto \ \alpha^*(u_{G,V})\]
is bijective according to Proposition \ref{prop:fix of global classifying}~(ii).
Proposition \ref{prop:pi_0 of Sigma^infty} thus says that
$\pi_0^K(\Sigma^\infty_+ \bL_{G,V})$
is a free abelian group with a basis given by the elements
\[ \tr_L^K(\sigma^L(\alpha^*(u_{G,V})))\ = \ 
 \tr_L^K(\alpha^*(\sigma^G(u_{G,V})))\ = \  \tr_L^K(\alpha^*(e_{G,V})) \]
where $L$ runs through all conjugacy classes of closed subgroups of $K$ 
with finite Weyl group and $\alpha$ runs through a set of representatives of the 
$W_KL$-orbits of the set $ \Rep(L,G)$.
The claim follows because $(K\times G)$-conjugacy classes of such pairs
$(L,\alpha)$ biject with pairs 
consisting of a conjugacy class of subgroups $(L)$ and
a $W_KL$-equivalence class in $\Rep(L,G)$.
\end{proof}

\begin{eg}
We discuss a specific example of Corollary \ref{cor-pi_0 of B_gl},
with $G=A_3$ the alternating group on three letters and 
$K=\Sigma_3$ the symmetric group on three letters.
The group $\Sigma_3$ has four conjugacy classes of subgroups, with representatives
$\Sigma_3, A_3, (12)$ and $e$. 
The groups $\Sigma_3, (12)$ and $e$ admit only trivial homomorphisms to $A_3$,
whereas the alternating group $A_3$ also has two automorphisms.
None of the three endomorphisms of $A_3$ are conjugate,
so the set $\Rep(A_3,A_3)$ has three elements.
However, the Weyl group $W_{\Sigma_3}A_3$ has two elements, and its action
realizes the two automorphisms of $A_3$.
So while $\pi_0^{A_3}(B_{\gl}A_3)\iso \Rep(A_3,A_3)$ has three elements,
it only contributes two generators to the stable group
$\pi_0^{\Sigma_3}(\Sigma^\infty_+ B_{\gl}A_3)$.
A basis for the free abelian group $\pi_0^{\Sigma_3}(\Sigma^\infty_+ B_{\gl}A_3)$
is thus given by the classes
\[ p_{\Sigma_3}^*(1) \ , \quad \tr_{A_3}^{\Sigma_3}(e_{A_3})
\ , \quad \tr_{A_3}^{\Sigma_3}(p_{A_3}^*(1))
\ , \quad \tr_{(12)}^{\Sigma_3}(p_{(12)}^*(1))
\text{\quad and\quad} \tr_e^{\Sigma_3}(p_e^*(1))\ . \]
Here the faithful $A_3$-representation $V$ is unspecified and 
we write $e_{A_3}=e_{A_3,V}\in\pi_0^{A_3}(\Sigma^\infty_+ B_{\gl}A_3)$ 
for the stable tautological class.
Moreover $p_H:H\to e$ denotes the unique homomorphism to the trivial group
and $1=\res^{A_3}_e(e_{A_3})$ is the restriction of the class $e_{A_3}$ to the
trivial group.
\end{eg}

Now we discuss how multiplicative features related to the smash 
product and the homotopy group pairings work out for global homotopy types.
In Construction \ref{con:Omega^bullet} we associated
an orthogonal space $\Omega^\bullet X$ to every orthogonal spectrum $X$.
This functor is compatible with the smash product of orthogonal
spectra and the box product of orthogonal spaces, 
in the sense of a lax symmetric monoidal transformation
\begin{equation}\label{eq:monoidal map of Omega^bullet} 
( \Omega^\bullet X)\boxtimes( \Omega^\bullet Y) \ \to \  \Omega^\bullet ( X\sm Y) \ .
\end{equation}
This morphism is associated to a bimorphism from $(\Omega^\bullet X,\Omega^\bullet Y)$ to
$\Omega^\bullet ( X\sm Y)$ with $(V,W)$-component the composite
\begin{align*}
 \map_*(S^V,X(V)) \times \map_*(S^W,Y(W))\ \xra{\ \sm \ }\ 
& \map_*(S^{V\oplus W},X(V)\sm Y(W))  \\
\xra{\map_*(S^{V\oplus W}, i_{V,W})}\ 
& \map_*(S^{V\oplus W},(X\sm Y)(V\oplus W)) \ . 
\end{align*}
The morphism is unital, associative and symmetric. 
Finally, the homotopy group pairing \eqref{eq:def_dot} 
for $k=l=0$ coincides with the composite
\begin{align*}
\pi_0^G(X)  \times  \pi_0^G (Y) \ = \ 
\pi_0^G(\Omega^\bullet X)  \times  \pi_0^G (\Omega^\bullet Y) \ \xra{\ \times\ } \ 
&\pi_0^G(\Omega^\bullet X \boxtimes  \Omega^\bullet Y) \\ 
\xra{\eqref{eq:monoidal map of Omega^bullet} } \ 
& \pi_0^G(\Omega^\bullet(X\sm Y))  \ = \  \pi_0^G(X\sm Y) \ .  
\end{align*}

\begin{eg}[Orthogonal ring spectra and orthogonal monoid spaces]
\label{eg:suspension spectrum of orthogonal monoid space} 
The monoidal morphism \eqref{eq:monoidal map of Omega^bullet} 
for the functor $\Omega^\bullet$ is unital, associative and symmetric. 
In particular the orthogonal space $\Omega^\bullet R$
associated to an orthogonal ring spectrum becomes an orthogonal monoid space
via the composite
\[ ( \Omega^\bullet R)\boxtimes( \Omega^\bullet R) \ \to \ 
 \Omega^\bullet ( R\sm R) \ \xra{\Omega^\bullet\mu}\ \Omega^\bullet R\ ,\]
and this passage preserves commutativity of multiplications.
The bimorphism corresponding to the induced product on $\Omega^\bullet R$ 
thus has as $(V,W)$-component the composite
\begin{align*}
 \map_*(S^V, R(V)) \times \map_*(S^W, R(W)) \ \xra{\ - \sm - \ } \  
&\map_*(S^{V\oplus W}, R(V)\sm R(W)) \\ 
\xra{\map_*(S^{V\oplus W},\, \mu_{V,W})} \ &\map_*(S^{V\oplus W}, R(V\oplus W)) \ . 
\end{align*}

The suspension spectrum functor (see Construction \ref{con:suspension spectrum})
takes the box product of orthogonal spaces to the smash product of orthogonal spectra.
In more detail: for all inner product spaces $V$ and $W$ the maps
\begin{align*}
 (\Sigma^\infty_+ X )(V)\sm(\Sigma^\infty_+ Y )(W) \ = \ 
&( S^V\sm X(V)_+ ) \sm (S^W\sm Y(W)_+)\\ 
\iso \ &S^{V\oplus W}\sm  (X(V)\times Y(W))_+\ \xra{S^{V\oplus W}\sm  i_{X,Y}} \\
&S^{V\oplus W}\sm (X\boxtimes Y)(V\oplus W)_+ \ = \ 
(\Sigma^\infty_+ (X\boxtimes Y))(V\oplus W)  
\end{align*}
form a bimorphism, so they correspond to a morphism of orthogonal spectra
\begin{equation}  \label{eq:box2smash}
 (\Sigma^\infty_+ X)\sm (\Sigma^\infty_+ Y) \ \to \ 
 \Sigma^\infty_+ (X\boxtimes Y) \ . 
\end{equation}
\end{eg}

\begin{prop}\label{prop:box to smash}
For all orthogonal spaces $X$ and $Y$ 
the morphism \eqref{eq:box2smash} is an isomorphism. Together with the 
unique isomorphism $\Sigma^\infty_+\ast\iso\mS$ 
this makes $\Sigma^\infty_+$ into a strong symmetric monoidal functor
from the category of orthogonal spaces to the category of orthogonal spectra.  
\end{prop}
\begin{proof}
We consider the composite
\begin{align*}
X\boxtimes Y\ &\xra{\eta_X\boxtimes\eta_Y} \   
\Omega^\bullet( \Sigma^\infty_+ X)\boxtimes \Omega^\bullet( \Sigma^\infty_+ Y)\
\xra{\eqref{eq:monoidal map of Omega^bullet}}\ \Omega^\bullet( (\Sigma^\infty_+ X)\sm (\Sigma^\infty_+ Y))\ ,\end{align*}
where $\eta:\Id\to\Omega^\bullet\circ\Sigma^\infty_+$ is the adjunction unit.
This composite is adjoint to a morphism of orthogonal spectra
$\Sigma^\infty_+(X\boxtimes Y)\to (\Sigma^\infty_+ X)\sm (\Sigma^\infty_+ Y)$.
We omit the verification that this morphism is indeed inverse to 
the morphism \eqref{eq:box2smash}.
\end{proof}

\begin{construction}[Orthogonal ring spectra from orthogonal monoid spaces]\label{con:monoid ring spectrum}
The suspension spectrum $\Sigma^\infty_+ M$ 
of an orthogonal monoid space $M$
becomes an orthogonal ring spectrum via the multiplication map
\[ (\Sigma^\infty_+ M )\sm(\Sigma^\infty_+ M )\ \xra[\iso]{\eqref{eq:box2smash}} \ 
\Sigma^\infty_+ ( M\boxtimes M )
\ \xra{\Sigma^\infty_+\mu_M}\ \Sigma^\infty_+ M \ .\]
If the multiplication of $M$ is commutative, then so is the resulting
multiplication on $\Sigma^\infty_+ M$.
The functor pair $(\Sigma^\infty_+,\Omega^\bullet)$ is again an adjoint
pair when viewed as functors between the categories of 
orthogonal monoid spaces and orthogonal ring spectra. 

This construction includes spherical monoid ring spectra:\index{subject}{monoid ring spectrum} 
if $M$ is a topological monoid, then the constant orthogonal space
with value $M$ inherits an associative and unital product from $M$.
The suspension spectrum
of such a constant multiplicative functor is the monoid ring spectrum.
\end{construction}

\begin{construction}
In Construction \ref{con:pairing equivalence homotopy}
we introduced pairings on the equivariant homotopy groups
of orthogonal $G$-spectra.
Now we consider two orthogonal spectra $X$ and $Y$.
In the global context, the equivariant homotopy groups also support inflation maps,
which we can use to define another pairing
\begin{equation}\label{eq:def_boxtimes}
\boxtimes \ : \
\pi_k^G(X) \ \times \ \pi_l^K (Y) \ \to \ \pi_{k+l}^{G\times K}(X\sm Y) \ .
\end{equation}
Here $G$ and $K$ are compact Lie groups and $k,l\in\mZ$.
We define this pairing as the composite
\begin{align*}
  \pi_k^G(X) \times \pi_l^K (Y)\ \xra{p_G^*\times p_K^*}\ 
 \pi_k^{G\times K}(X) \times \pi_l^{G\times K} (Y)\ \xra{\ \times \ }\ 
 \pi_{k+l}^{G\times K}(X\sm Y)\ ,
\end{align*}
where $p_G:G\times K\to G$ and $p_K:G\times K\to K$
are the projections to the two factors.
Theorem \ref{thm:product properties}
and the additivity of inflation maps imply 
that this pairing is biadditive, and that it satisfies
certain associativity, commutativity and restriction properties.
We do take the time to spell out the most important
properties of these pairings in the next theorem.

As $G$ and $K$ vary through all compact Lie groups,
the $\boxtimes$-pairings form a bimorphism of global functors 
in the sense of Construction \ref{con:Box product} below.
Moreover, the passage from the `diagonal' pairings
to the `external' pairings \eqref{eq:def_boxtimes} can be reversed
by taking $K=G$ and restricting to the diagonal;
suitably formalized, diagonal and external pairings contain the same
amount of information. We refer 
to Remark \ref{rk:diagonal versus external products} below for more details.
\end{construction}

\begin{theorem}\label{thm:external product properties}
  Let $G,K$ and $L$ be compact Lie groups and $X, Y$ and $Z$ orthogonal spectra.
  \begin{enumerate}[\em (i)]
  \item {\em (Biadditivity)}
    The product $\boxtimes:\pi_k^G(X) \times\pi_l^K (Y) \to\pi_{k+l}^{G\times K}(X\sm Y)$
    is biadditive.
  \item {\em (Unitality)}
    Let $1\in\pi_0^e(\mS)$ denote the class represented by the identity of $S^0$.
    The product is unital in the sense that $1\boxtimes x=x=x\boxtimes 1$ 
    under the identifications $\mS\sm X=X=X\sm \mS$ and $e\times G\iso G\iso G\times e$.
  \item {\em (Associativity)} For all $x\in \pi_k^G ( X )$, $y\in \pi_l^K( Y )$
    and $z\in\pi_m^L( Z )$ the relation
    \[ x\boxtimes(y\boxtimes z) \ = \ (x\boxtimes y)\boxtimes z\]
    holds in $\pi^{G\times K\times L}_{k+l+m}(X\sm Y\sm Z)$.
  \item {\em (Commutativity)} For all $x\in \pi_k^G ( X )$ and $y\in \pi_l^K( Y )$
    the relation
    \[ \tau^{X,Y}_*( x\boxtimes y ) \ = \ (-1)^{kl}\cdot \tau_{G,K}^*( y\boxtimes x)\]
    holds in $\pi^{G\times K}_{l+k}( Y\sm X)$,
    where $\tau^{X,Y}:X\sm Y\to Y\sm X$ is the symmetry isomorphism
    of the smash product and $\tau_{G,K}:G\times K \to K\times G$
    interchanges the factors.
  \item {\em (Restriction)} For all $x\in \pi_k^G ( X )$ and $y\in \pi_l^K( Y )$
    and all continuous homomorphisms $\alpha:\bar G\to G$
    and $\beta:\bar K\to K$ the relation
    \[ \alpha^*(x)\boxtimes\beta^*(y)\ = \ (\alpha\boxtimes\beta)^*(x\boxtimes y) \]
    holds in $\pi_{k+l}^{\bar G\times\bar K}(X\sm Y)$.
   \item {\em (Transfer)} For all closed subgroups $H\leq G$ and $L\leq K$ the square 
    \[ \xymatrix@C=15mm{ \pi_k^H(X) \times \pi^L_l(Y) 
      \ar[r]^-{\boxtimes} \ar[d]_{\tr_H^G\times\tr_L^K} & 
      \pi^{H\times L}_{k+l} (X\sm Y) \ar[d]^{\tr_{H\times L}^{G\times K}}\\
      \pi^G_k (X) \times  \pi^K_l(Y) \ar[r]_-{\boxtimes} & \pi^{G\times K}_{k+l} (X\sm Y)} \]
    commutes.
  \end{enumerate}
\end{theorem}
\begin{proof}
Parts~(i) through~(v) follow from the 
respective parts of Theorem \ref{thm:product properties} 
by naturality.
For part (vi) we start with two special cases, namely $L=K$ respectively $H=G$.
The two proofs are analogous, so we only treat the case $L=K$:
\begin{align*}
  \tr_{H\times K}^{G\times K}(x\boxtimes y)\ 
&= \  \tr_{H\times K}^{G\times K}( p_H^*(x)\times p_K^*(y))\ 
 = \  \tr_{H\times K}^{G\times K}(p_H^*(x)\times \res^{G\times K}_{H\times K}(p_K^*(y)))\\ 
&= \  \tr_{H\times K}^{G\times K}(p_H^*(x))\times p_K^*(y)\ 
= \ p_G^*(\tr_H^G(x))\times p_K^*(y)\ =\tr_H^G(x)\boxtimes y
\end{align*}
We slightly abuse notation by writing $p_K$ for both the projections
of $H\times K$ and of $G\times K$ to $K$.
The third equation is the reciprocity
relation of Theorem \ref{thm:product properties}~(vi). 
The fourth equation is the compatibility of transfers with inflations
(Proposition \ref{prop:transfer and epi}~(ii)).

The general case is now obtained by combining the two special cases:
\begin{align*}
 \tr_H^G(x)\boxtimes\tr_L^K(y) \ &= \ 
\tr_{H\times K}^{G\times K}(x\boxtimes\tr_L^K(y))\\ 
&= \ \tr_{H\times K}^{G\times K}(\tr_{H\times L}^{H\times K}(x\boxtimes y))\ = \ 
\tr_{H\times L}^{G\times K}(x\boxtimes y)\qedhere  
\end{align*}
\end{proof}

For all orthogonal spectra (i.e., all global homotopy types), 
the collection of equivariant homotopy groups $\{ \pi_0^G (X) \}_G$
come with restriction and transfer maps, and this data together forms
a `global functor', compare Definition \ref{def:global functor} below. 
The geometric fixed point homotopy groups have fewer natural operations,
and they do {\em not} allow restriction to subgroups.
However, geometric fixed points still have inflation maps, i.e., 
restriction maps along epimorphisms.
Indeed, in Construction \ref{con:restriction in geometric fixed}
we defined inflation maps on geometric fixed point homotopy groups,
associated to a continuous epimorphism $\alpha:K\to G$ between compact Lie groups.
When $X$ is an orthogonal spectrum, representing a global homotopy type,
then $\alpha^*(X_G)=X_K$, and the inflation maps become homomorphisms\index{subject}{inflation map!of geometric fixed point homotopy groups}  
\[ \alpha^* \ : \ \Phi^G_k (X) \ \to \ \Phi^K_k (X) \ .\]
These inflation maps between the geometric fixed point homotopy groups
are clearly natural in the orthogonal spectrum.

\begin{construction}[Semifree orthogonal spectra]\label{con:free orthogonal}
Given a compact Lie group $G$ and a $G$-representation $V$, the functor
\[ \ev_{G,V}\ : \ \spec \ \to \ G\bT_* \]
that sends an orthogonal spectrum $X$ to the based $G$-space $X(V)$
has a left adjoint\index{subject}{semifree orthogonal spectrum}\index{subject}{orthogonal spectrum!semifree}\index{symbol}{$F_{G,V}$ - {semifree orthogonal spectrum generated by $(G,V)$}}
\begin{equation}\label{eq:define F_{G,V}}
 F_{G,V} \ : \ G\bT_* \ \to \ \spec  \ .  
\end{equation}
The semifree orthogonal spectrum generated by a based $G$-space $A$ in level $V$ is 
\[ F_{G,V} A\ = \ \bO(V,-)\sm_G A\ ; \]
the value at an inner product space $W$ is thus given by
\[ (F_{G,V}A)(W) \ = \ \bO(V,W)\sm_G A\ . \]
We note that $F_{G,V} A$ consists of a single point in all levels below
the dimension of $V$.
The `freeness' property of $F_{G,V}A$
is a consequence of the enriched Yoneda lemma,\index{subject}{Yoneda lemma!enriched}
see Remark \ref{rk:enriched Yoneda} or \cite[Sec.\,1.9]{kelly-enriched category};
it means that for every orthogonal spectrum $X$ and every based $G$-map $f:A\to X(V)$ 
there is a unique morphism 
$f^\flat:F_{G,V}A\to X$ of orthogonal spectra such that the composite
\[  A \xra{\ \Id\sm -\ }\ 
\bO(V,V)\sm_G A = (F_{G,V}A)(V) \ \xra{f^\flat(V)} \ X(V) \]
is $f$.
Indeed, the map $f^\flat(W)$ is the composite
\[ \bO(V,W)\sm_G A \ \xra{\Id\sm_G f} \
\bO(V,W)\sm_G X(V) \ \xra{\ \text{act}\ } \ X(W)\ .\]
\end{construction}

\begin{rk}[Semifree spectra as global Thom spectra]\label{rk:free are Thom}
The underlying non-equivariant stable homotopy type of $F_{G,V}$ 
is the Thom spectrum of the negative of the bundle
\[ \bL(V,\mR^\infty)\times_G V \ \to \ \bL(V,\mR^\infty)/G \ = \ B G \ , \]
the vector bundle over $B G$ associated to the $G$-representation $V$. 
So one should think of the semifree orthogonal spectrum $F_{G,V}=F_{G,V} S^0$
as the `global Thom spectrum' associated to a `virtual global vector bundle',
namely the negative of the vector bundle over $B_{\gl}G$ associated to the
$G$-representation $V$. 

The special case $G=O(m)$ of the orthogonal group
with $V=\nu_m$, i.e., $\mR^m$ with the tautological $O(m)$-action,
will feature prominently in the rank filtration of the global
Thom spectrum $\bmO$ in Section \ref{sec:global Thom}.
Non-equivariantly, $F_{O(m),\nu_m}$ is the Thom spectrum 
of the negative of the tautological $m$-plane bundle
over the Grassmannian $G r_m(\mR^\infty)$; 
the traditional notation for this Thom spectrum is $M T O(m)$
or simply $M T (m)$.
Indeed,
 \[ F_{O(m),\nu_m}(\mR^{m+n}) \ = \ \bO(\mR^m,\mR^{m+n})/O(m) \]
is the Thom space of the orthogonal complement of
the tautological $m$-plane bundle over $\bL(\mR^m,\mR^{m+n})/O(m)=G r_m(\mR^{m+n})$,
and this is precisely the $(m+n)$-th space of $M T O(m)$, see
for example \cite[Sec.\,3.1]{galatius-madsen-tillmann-weiss}.
Similarly, the non-equivariant homotopy type
underlying $F_{S O(m),\nu_m}$ is an oriented version of $M T O(m)$,
which is usually denoted $M T S O(m)$ or sometimes $M T(m)^+$.
\end{rk}

\begin{eg}[Smash products of semifree orthogonal spectra]
The smash product of two semifree orthogonal spectra is again
a semifree orthogonal spectrum.\index{subject}{smash product!of semifree spectra} 
In more detail, we consider
\begin{itemize}
\item  two compact Lie groups $G$ and $K$, 
\item a $G$-representation $V$ and  a $K$-representation $W$, and 
\item a based $G$-space $A$ and a based $K$-space $B$.
\end{itemize}
Then $V\oplus W$ is a $(G\times K)$-representation and
$A\sm B$ is a $(G\times K)$-space via
\[ (g,k)\cdot (v,w)\ = \ (g v,k w) \text{\quad respectively\quad}
(g,k)\cdot (a\sm b)\ = \ g a\sm k b \ . \]
We claim that the smash product
$(F_{G,V}A)\sm (F_{K,W}B)$ is canonically isomorphic
to the semifree orthogonal spectrum generated by the $(G\times K)$-space $A\sm B$
in level $V\oplus W$. 
Indeed, a morphism
\begin{equation}\label{eq:F_V_smash_F_W}
(F_{G,V}A)\sm (F_{K,W}B)\ \to\  F_{G\times K,V\oplus W}(A\sm B) 
\end{equation}
is obtained by the universal property \eqref{eq:universal property smash}
from the bimorphism with $(U,U')$-component
\begin{align*}
 (F_{G,V}A)(U)\sm (F_{K,W}B)(U') \ = \ 
(\bO(V,U)&\sm_G A)\sm(\bO(W,U')\sm_K B) \\ 
\xra{\ \oplus \ } \ &\bO(V\oplus W,U\oplus U')\sm_{G\times K}(A\sm B) \\ 
= \ &((F_{G\times K,V\oplus W})(A\sm B))(U\oplus U') \ . 
\end{align*}
In the other direction, a morphism 
$F_{G\times K,V\oplus W}(A\sm B)\to F_{G,V}A\sm F_{G,W}B$ 
is freely generated by the $(G\times K)$-map
\[ A\sm B \ \to \ (F_{G,V}A)(V) \sm (F_{K,W}B)(W) \ \xra{i_{V,W}}\ 
 (F_{G,V}A\sm F_{K,W}B)(V\oplus W) \ . \]
These two maps are inverse to each other.
\end{eg}

In Proposition \ref{prop:free_orthogonal_space}~(ii) 
we have seen that for every compact Lie group $G$,
every $G$-representation $V$ and every faithful $G$-representation $W$ 
the restriction morphism of orthogonal spaces
$\rho_{V,W}/G: \bL_{G,V\oplus W} \to \bL_{G,W}$ is a global equivalence.
One consequence is that the semifree orthogonal space $\bL_{G.W}$ 
has a well-defined unstable global homotopy type, independent of which faithful
$G$-representation is used. 
Another consequence is that the induced morphism
\[ \Sigma^\infty_+ \rho_{V,W}/G\ :\  
\Sigma^\infty_+ \bL_{G,V\oplus W}\ \to\  \Sigma^\infty_+ \bL_{G,W} \]
of suspension spectra is a global equivalence of orthogonal spectra,
by Corollary \ref{cor:suspension spectrum globally homotopical}.
For an inner product space $W$, the untwisting homeomorphisms \eqref{eq:untwisting_homeo}
descend to homeomorphisms on $G$-orbit spaces
\[     \bO(V,W)\sm_G S^V \ \iso \ S^W\sm \bL(V,W)/G_+\ .\]
As $W$ varies, these form the `untwisting isomorphism'
of orthogonal spectra\index{subject}{untwisting homeomorphism}
\[ F_{G,V}S^V \ \iso \  \Sigma^\infty_+ \bL_{G,V} \ .   \]
So suspension spectra of semifree orthogonal spaces are semifree orthogonal spectra.
We will now prove a generalization of the fact that $\Sigma^\infty_+ \rho_{V,W}/G$
is a global equivalence for these global Thom spectra.
Given $G$-representations $V$ and $W$,
we define a restriction morphism of orthogonal spectra
\begin{equation}\label{eq:define_lambda}\index{symbol}{$\lambda_{G,V,W}$ - {fundamental global equivalence of orthogonal spectra}}
 \lambda_{G,V,W}\ : \ F_{G,V\oplus W}S^V\ \to\  F_{G,W}    
\end{equation}
as the adjoint of the based $G$-map
\[ S^V \ \to \ \bO(W,V\oplus W)/G \ = \  F_{G,W}(V\oplus W) \ , \quad
v \ \longmapsto \  ((v,0) , i)\cdot G \ ,\]
where $i:W\to V\oplus W$ is the embedding of the second summand.
The value of $\lambda_{G,V,W}$ at an inner product space $U$ is then
\begin{align*}
 \lambda_{G,V,W}(U)\ : \ \bO(V\oplus W, U)\sm_G S^V\ &\to\quad  
\bO(W,U)/G \\
[(u,\varphi)\sm v] \qquad &\longmapsto \  (u+\varphi(v),\varphi\circ i)\cdot G\ .  
\end{align*}

\begin{theorem}\label{thm:faithful independence}
Let $G$ be a compact Lie  group, $V$ a $G$-representation
and $W$ a faithful $G$-representation. Then the morphism
\[ \lambda_{G,V,W} \ : \ F_{G,V\oplus W}S^V \ \to \ F_{G,W} \]
is a global equivalence of orthogonal spectra.
\end{theorem}
\begin{proof}
To simplify the notation we abbreviate the restriction morphism of orthogonal spaces
to
\[ \rho \ = \ \rho_{V,W}\ : \ \bL(V\oplus W,-)\ \to \ \bL(W,-)\ . \]
We let $K$ be another compact Lie group and $U\in s(\Uc_K)$ 
a finite-dimensional $K$-subrepresentation of the complete $K$-universe $\Uc_K$.
In a first step we produce a $K$-representation $U'\in s(\Uc_K)$ 
with $U\subseteq U'$ and
a continuous $(K\times G)$-equivariant map 
\[ h\ :\  \bL(W,U) \ \to \ \bL(V\oplus W,U') \] 
such that the lower right triangle in the diagram
\begin{equation}\begin{aligned}\label{eq:base_square} 
\xymatrix@C=15mm{
\bL(V\oplus W,U) \ar[r]^{i_*}\ar[d]_{\rho(U)} & 
\bL(V\oplus W,U') \ar[d]^{\rho(U')} \\
\bL(W,U) \ar[r]_{i_*} \ar@{-->}[ur]^(.4)h & \bL(W,U') } 
\end{aligned}\end{equation}
commutes, and the upper left triangle
commutes up to $(K\times G)$-equivariant fiberwise homotopy over $\bL(W,U')$,
where $i:U\to U'$ is the inclusion.

Since $G$ acts faithfully on $W$ (and hence on $V\oplus W$),
both $\bL(W,\Uc_K)$ and $\bL(V\oplus W,\Uc_K)$ are universal spaces
for the same family of subgroups of $K\times G$,
namely the family $\Fc(K;G)$ of graph subgroups,\index{subject}{graph subgroup}
compare Proposition \ref{prop:free_orthogonal_space}~(i). 
If $\Gamma$ is the graph of
a continuous homomorphism $\alpha:L\to G$ defined on some closed subgroup $L$ of $K$,
then the $\Gamma$-fixed points of $\bL(W,\Uc_K)$ are given by
\[ \bL(W,\Uc_K)^\Gamma \ = \   \bL^L(\alpha^* W,\res^K_L( \Uc_K)) \ , \]
the space of $L$-equivariant linear isometric embeddings from $\alpha^* W$.
The same is true for $V\oplus W$ instead of $W$, and so the
$\Gamma$-fixed point map $\rho(\Uc_K)^\Gamma:\bL(V\oplus W,\Uc_K)^\Gamma \to \bL(W,\Uc_K)^\Gamma$
is the restriction map
\[
\bL^L(\alpha^* V\oplus \alpha^* W,\res^K_L(\Uc_K))\ \to \ 
\bL^L(\alpha^* W,\res^K_L(\Uc_K)) \]
to the summand $\alpha^* W$. This map is a locally trivial fiber bundle,
hence a Serre fibration. We conclude that the restriction map
$\rho(\Uc_K):\bL(V\oplus W,\Uc_K) \to \bL(W,\Uc_K)$ is both
a $(K\times G)$-weak equivalence and a $(K\times G)$-fibration.

Since $\bL(W,U)$ is cofibrant as a $(K\times G)$-space 
(by Proposition \ref{prop:K G cofibration}~(ii)), 
the $(K\times G)$-map $i_*:\bL(W,U)\to \bL(W,\Uc_K)$  
thus admits a $(K\times G)$-equivariant lift 
$h:\bL(W,U)\to \bL(V\oplus W,\Uc_K)$  such that $\rho(\Uc_K)\circ h=i_*$.
Since the space $\bL(W,U)$ is compact and 
$\bL(V\oplus W,\Uc_K)$ is the filtered union of the closed subspaces
$\bL(V\oplus W, U')$ for $U'\in s(\Uc_K)$, the lift $h$
lands in the subspace $\bL(V\oplus W, U')$ for suitably large $U'\in s(\Uc_K)$,
and we may assume that $U\subseteq U'$.

The two maps 
\[ h\circ\rho(U)\ , \ i_* \ : \  \bL(V\oplus W,U) \ \to \ \bL(V\oplus W, U')   \]
become equal after applying $\rho(U'):\bL(V\oplus W, U') \to \bL(V, U')$,
hence the composites with $i_*:\bL(V\oplus W,U') \to \bL(V\oplus W, \Uc_K)$
become equal after applying the $(K\times G)$-equivariant 
acyclic fibration $\rho(\Uc_K):\bL(V\oplus W, \Uc_K) \to \bL(W, \Uc_K)$.
Since $\bL(V\oplus W,U)$ is also $(K\times G)$-cofibrant, there
is a fiberwise $(K\times G)$-equivariant homotopy between 
$h\circ\rho(U)$ and $i_*$ in $\bL(V\oplus W,\Uc_K)$.
Again by compactness, the homotopy has image in 
$\bL(V\oplus W, U'')$ for suitably large $U''\in s(\Uc_K)$.
So after increasing $U'$, if necessary, we have proved the claim
subsumed in the diagram \eqref{eq:base_square}.

Now we lift the data produced in the first step to the Thom spaces
of the orthogonal complement bundles. The diagram \eqref{eq:base_square} 
is covered by morphisms of $(K\times G)$-vector bundles:
\begin{equation}\begin{aligned}\label{eq:bundle_square} 
\xymatrix@C=12mm{
(U'-U) \times \xi(V\oplus W,U)\times V 
 \ar[r]^-{\bar i}\ar[d]_{(U'-U)\times \bar \rho(U)} & 
\xi(V\oplus W,U')\times V \ar[d]^{\bar \rho(U')} \\
(U'-U) \times \xi(W,U) \ar[r]_-{\bar i} \ar@{-->}[ur]^(.4){\bar h} & \xi(W,U') } 
\end{aligned}\end{equation}
The maps on the total spaces of the bundles are defined as follows:
the right vertical map is
\[ \bar \rho(U') \ : \ \xi(V\oplus W,U')\times V\ \to \ \xi(W,U')\ , \quad 
((u',\varphi),\, v)\ \longmapsto \ ( u' +\varphi(v),\,\varphi|_W)\ . \]
The map $\bar\rho(U)$ is defined in the same way.
The lower horizontal map is 
\[ \bar i \ : \ (U'-U)\times \xi(W,U)
\ \to \ \xi(W,U')\ , \quad 
(u',(u,\varphi))\ \longmapsto \ (u'+u,\,\varphi)\ . \]
The upper horizontal map is defined in the same way, but with $V\oplus W$ 
instead of $W$ and multiplied by the identity of $V$.
These four outer morphisms in \eqref{eq:bundle_square} 
are all fiberwise linear isomorphisms; so each of these four bundle
maps expresses the source bundle as a pullback of the target bundle.
In particular, the square
\begin{equation}\begin{aligned}\label{eq:fiber pullback}
 \xymatrix{ 
\xi(V\oplus W, U')\times V \ar[r]^-{\bar\rho(U')} \ar[d] & \xi(W,U') \ar[d]\\
\bL(V\oplus W, U') \ar[r]_-{\rho(U')} & \bL(W,U')}   
\end{aligned}\end{equation}
is a pullback;  so the composite
\[ (U'-U)\times \xi(W,U) \ \to \ \bL(W,U)\ \xra{\ h\ } \ \bL(V\oplus W,U')  \]
and the map of total spaces $\bar i:(U'-U)\times\xi(W,U)\to\xi(W,U')$
assemble into a map 
\[ \bar h \ : \  (U'-U)\times \xi(W,U) \ \to \  \xi(V\oplus W,U')  \]
that covers $h$ and is a fiberwise linear isomorphism.

In \eqref{eq:bundle_square} (as in \eqref{eq:base_square})
the outer square and the lower right triangle commute,
but the upper left triangle does {\em not} commute.
We will now show that the upper left triangle commutes up to homotopy
of $(K\times G)$-equivariant bundle maps.
For this purpose we let
\[ H \ : \  \bL(V\oplus W,U) \times [0,1] \ \to \ \bL(V\oplus W, U') \]
be a  $(K\times G)$-equivariant homotopy from the
map $i_*$ to $h\circ \rho(U)$, such that
$\rho(U')\circ H:\bL(V\oplus W,U) \times [0,1] \to \bL(W, U')$
is the constant homotopy from $\rho(U')\circ i_* = i_*\circ \rho(U)$ to itself.
Again because the square \eqref{eq:fiber pullback} is a pullback, the composite
\[\xi(V\oplus W,U)\times V \times (U'-U) \times [0,1]\ \to \ 
\bL(V\oplus W,U) \times [0,1]\ \xra{\ H\ } \  \bL(V\oplus W, U')  \]
and the map of total spaces 
\begin{align*}
  (U'-U) \times \xi(V\oplus W,U)\times V\times  [0,1]\ &\xra{\text{proj}} \ 
 (U'-U)\times \xi(V\oplus W,U)\times V  \\ 
&\xra{\bar\rho(U')\circ\bar i=\bar i\circ\bar\rho(U)} \  \xi(W,U')  
\end{align*}
assemble into a map 
\[ \bar H \ : \  (U'-U) \times \xi(V\oplus W,U) \times V \times  [0,1] 
\ \to \  \xi(V\oplus W,U')\times V  \]
that covers the homotopy $H$. This lift $\bar H$ is a 
$(K\times G)$-equivariant homotopy of vector bundle morphisms,
and for every $t\in [0,1]$, the relation
\begin{align*}
\bar\rho(U') \circ \bar H(-,t) \ &= \  \bar\rho(U')\circ\bar i \ 
= \ \bar\rho(U')\circ (\bar h\circ\bar\rho(U))
\end{align*}
holds by definition of $\bar H$.
For $t=0$ this shows that $\bar H$ starts with 
$\bar i:(U'-U)\times \xi(V\oplus W,U) \times V  \to 
 \xi(V\oplus W,U')\times V$; for  $t=1$ this shows that $\bar H$ 
ends with $\bar h\circ\bar\rho(U)$, one more time 
because \eqref{eq:fiber pullback} is a pullback.
We conclude that $\bar H$ makes the upper left triangle in \eqref{eq:bundle_square} 
commute up to equivariant homotopy of vector bundle maps.

Passing to Thom spaces in \eqref{eq:bundle_square} 
gives a diagram of $(K\times G)$-equivariant based maps:
\[ \xymatrix@C=20mm{
S^{U'-U}\sm \bO(V\oplus W,U)\sm S^V
 \ar[r]^{\sigma_{U,U'-U}\sm S^V}\ar[d]_{S^{U'-U}\sm \lambda_{V,W}(U)} & 
\bO(V\oplus W,U')\sm S^V \ar[d]^{\lambda_{V,W}(U')} \\
S^{U'-U}  \sm \bO(W,U) \ar[r]_{\sigma_{U,U'-U}} \ar@{-->}[ur]^(.4){\bar h} & \bO(W,U') } 
 \]
Again, the lower right triangle commutes,
and the upper left triangle commutes up to  $(K\times G)$-equivariant 
based homotopy.
We pass to $G$-orbit spaces and obtain a diagram of based $K$-spaces
\[ \xymatrix@C=20mm{
S^{U'-U} \sm (F_{G,V \oplus W} S^V)(U)
\ar[r]^{\sigma_{U,U'-U}}\ar[d]_{S^{U'-U}\sm \lambda_{G,V,W}(U)} & 
(F_{G,V \oplus W} S^V)(U') \ar[d]^{\lambda_{G,V,W}(U')} \\
S^{U'-U}\sm F_{G,W}(U) \ar[r]_{\sigma_{U,U'-U}} \ar[ur]^(.4){\bar h/G} & F_{G,W}(U') } 
 \]
whose lower right triangle commutes, and whose  upper left triangle
commutes up to $K$-equivariant based homotopy.
Since we had started with an arbitrary $K$-subrepresentation $U\in s(\Uc_K)$,
this implies that for every based $K$-space $A$ 
and $G$-representation $\bar U$ the map on colimits
\begin{align*}
 \colim_{U\in s(\Uc_K)}\, &[S^U\sm A, (F_{G,V \oplus W} S^V)(U\oplus\bar U)]^K \\ 
&\to \quad \colim_{U\in s(\Uc_K)} \, [S^U\sm A, F_{G,W}(U\oplus\bar U) ]^K 
\end{align*}
induced by the morphism $\lambda_{G,V,W}$ is bijective. 
For $A=S^k$ and $\bar U=0$ this shows that $\pi_k^G(\lambda_{G,V,W})$
is an isomorphism.
For $A=S^0$ and $\bar U=\mR^k$ this shows that $\pi_{-k}^G(\lambda_{G,V,W})$
is an isomorphism.
So $\lambda_{G,V,W}$ is a global equivalence.
\end{proof}

\section{Global functors}\label{sec:global functors}

This section is devoted to the category of {\em global functors},
the natural home of the collection of equivariant
homotopy groups of a global stable homotopy type.
The category of global functors is a symmetric monoidal abelian category
with enough injectives and projectives that plays the same role
for global homotopy theory that is played by the category of abelian groups
in ordinary homotopy theory, or by the category of $G$-Mackey functors
for $G$-equivariant homotopy theory.

We introduce  global functors in Definition \ref{def:global functor} 
as additive functors on the {\em global Burnside category}, 
the category of natural operations between equivariant stable homotopy groups.
The abstract definition ensures that equivariant homotopy groups
of orthogonal spectra are tautologically global functors.
A key result is Theorem \ref{thm:Burnside category basis} 
that describes explicit bases of the morphism groups of the 
global Burnside category in terms of transfers and restriction operations.
This calculation is the key to comparing our notion of global functors
to other kinds of global Mackey functors, 
as well as for all concrete calculations with global functors.
Example \ref{eg:global functor examples} lists interesting
examples of global functors:
the Burnside ring global functor,
represented global functors, constant global functors,
the representation ring global functor, and Borel type global functors.
Many more examples of global functors are discussed 
in the remaining sections of this book.

An abstract way to motivate the appearance of global functors is as follows.
One can consider the globally connective (respectively globally coconnective)
orthogonal spectra, i.e., those where all equivariant homotopy groups
vanish in negative dimensions (respectively in positive dimensions).
Then the full subcategories of globally connective 
respectively globally coconnective spectra
define a non-degenerate t-structure on the triangulated global stable homotopy category,
and the heart of this t-structure is (equivalent to) 
the abelian category of global functors; we refer the reader to
Theorem \ref{thm:t-structure on GH} below for details.

In the global context, the external pairing \eqref{eq:def_boxtimes}
of equivariant homotopy groups gives rise to a symmetric monoidal structure 
on the global Burnside category,
see Theorem \ref{thm:bA monoidal}. Hence the abelian category of
global functors can also be endowed with a Day type convolution product,
the box product of global functors, see Construction \ref{con:Box product}.

\medskip

\index{subject}{Burnside category|(}

\begin{construction}[Burnside category]
We define the pre-additive {\em Burnside category}\index{subject}{Burnside category} $\bA$.\index{symbol}{$\bA$ - {Burnside category}}
The objects of $\bA$ are all compact Lie groups;
morphisms from a group $G$ to $K$ are defined as
\[ \bA(G,K) \ = \ \Nat(\pi_0^G,\pi_0^K) \ ,\]
the set of natural transformations of functors, from orthogonal spectra to sets,
between the equivariant homotopy group functors $\pi_0^G$ and $\pi_0^K$. 
Composition in the category $\bA$ is composition of natural transformations.
\end{construction}

It is not a priori clear that the natural transformations 
from $\pi_0^G$ to $\pi_0^K$ form a set (as opposed to a proper class),
but this follows from Proposition \ref{prop:B_gl represents} below. 
The Burnside category $\bA$ is skeletally small: isomorphic compact Lie groups
are also isomorphic as objects in the category $\bA$, and every compact Lie
group is isomorphic to a closed subgroup of an orthogonal group $O(n)$.
The functor $\pi_0^K$ is abelian group valued,
so the set $\bA(G,K)$ is an abelian group under objectwise addition of transformations.
Proposition \ref{prop:additivity prop}, 
applied to the category of orthogonal spectra,
shows that set valued natural transformations between 
the two reduced additive functors $\pi_0^G$ and $\pi_0^K$ are automatically
additive. So composition in the Burnside category is additive in each variable, 
and $\bA$ is indeed a pre-additive category.

\begin{defn}\label{def:global functor} 
A {\em global functor}\index{subject}{global functor} is an additive
functor from the Burnside category $\bA$ to the category of abelian groups.
A morphism of global functors is a natural transformation.
We write $\GF$ for the category of global functors.
\end{defn}\index{symbol}{  $\GF$ - {category of global functors}} 
\index{symbol}{  $\Ab$ - {category of abelian groups}}

We discuss various explicit examples of interesting global functors 
in Example \ref{eg:global functor examples}.

\begin{eg}
The definition of the Burnside category $\bA$ is made so that
the collection of equivariant homotopy groups of an orthogonal spectrum 
is tautologically a global functor.
Explicitly, the homotopy group global functor $\upi_0(X)$
of an orthogonal spectrum $X$ is defined on objects by
\[ \upi_0(X)(G)\ = \ \pi_0^G (X) \]
and on morphisms by evaluating natural transformations at $X$.
It is less obvious that conversely every global functor is isomorphic
to the global functor $\upi_0(X)$ of some orthogonal spectrum $X$;
we refer the reader to Remark \ref{rk:general Eilenberg Mac Lane} 
below for the construction of Eilenberg-Mac\,Lane spectra from global functors.
\end{eg}

As a category of additive functors out of a skeletally small 
pre-additive category, the category $\GF$ of global functors has some
immediate properties that we collect in the following proposition.

\begin{prop}
The category $\GF$ of global functors is a Grothendieck abelian
category with enough injectives and projectives.
\end{prop}
\begin{proof}
Any category of additive functors out of a skeletally small additive category
is Grothendieck abelian with objectwise notion of exactness,
see \cite[IV Prop.\,7.2]{stenstrom} and \cite[V Ex.\,1.2]{stenstrom}.
A set of projective generators is given by the represented global functors
$\bA(G,-)$ where $G$ runs through a set of representatives of the
isomorphism classes of compact Lie groups, see \cite[IV Cor.\,7.5]{stenstrom}.
A set of injective cogenerators is given similarly by the global functors
\[ \Hom(\bA(-,K),\mQ/\mZ) \ : \ \bA \ \to \ \Ab \ .\qedhere \] 
\end{proof}

As we shall explain in Construction \ref{con:Box product} below, 
the category $\GF$ has a closed symmetric monoidal product 
$\Box$ that arises as a convolution product for a certain symmetric monoidal structure
on the Burnside category $\bA$.

\medskip

Our definition of the Burnside category is made so that 
every orthogonal spectrum $X$ gives rise to a homotopy group global functor
without further ado, but it is not clear from the definition how to describe
the morphism groups of $\bA$ explicitly. Our next aim 
is to show that each morphism group $\bA(G,K)$ is 
a free abelian group with an explicit basis given by certain composites 
of a restriction and a transfer morphism.
This calculation has two ingredients: 
We identify natural transformations
from $\pi_0^G$ to $\pi_0^K$ with the group $\pi_0^K(\Sigma^\infty_+ B_{\gl}G)$,
and then we exploit the explicit calculation of the latter group 
in Corollary \ref{cor-pi_0 of B_gl}.

\begin{prop}\label{prop:B_gl represents}\index{subject}{global classifying space}   
Let $G$ and $K$ be compact Lie groups and $V$ a faithful $G$-representation.
Then evaluation at the stable tautological class \eqref{eq:define_stable_tautological}
is a bijection
\[  \bA(G,K)\ = \ \Nat^{\spec}(\pi_0^G,\pi_0^K) 
\ \to \ \pi_0^K(\Sigma^\infty_+ \bL_{G,V}) \ , \quad
\tau\ \longmapsto \ \tau(e_{G,V})\]
to the 0-th $K$-equivariant homotopy group 
of the orthogonal spectrum $\Sigma^\infty_+ \bL_{G,V}$.
In other words, the morphism  of global functors
\[ \bA(G,-)\ \to \ \upi_0(\Sigma^\infty_+ \bL_{G,V}) \]
classified by the stable tautological class $e_{G,V}$ is an isomorphism.
\end{prop}
\begin{proof}
  We apply the representability result 
  of Proposition \ref{prop:pi_0^G representability}
  to the category of orthogonal spectra and
  the adjoint functor pair $(\Sigma^\infty_+,\Omega^\bullet)$. 
  If $G$ is a compact Lie group, $V$ a $G$-representation and $W$
  a faithful $G$-representation, then the restriction morphism
  $\rho_{G,V,W}:\bL_{G,V\oplus W}\to \bL_{G,W}$ is a global equivalence
  of orthogonal spaces (Proposition \ref{prop:free_orthogonal_space}). 
  So the induced morphism of unreduced suspension spectra
  $\Sigma^\infty_+\rho_{G,V,W}$ is a global equivalence
  by Corollary \ref{cor:suspension spectrum globally homotopical};
  in particular, the morphism of $\Rep$-functors
  $\upi_0(\Sigma^\infty_+\rho_{G,V,W})$ is an isomorphism. So 
  Proposition \ref{prop:pi_0^G representability} applies and shows that
  evaluation at the tautological class is bijective. This proves the claim.
\end{proof}

Now we name a basis of the group $\bA(G,K)$.
For a pair $(L,\alpha)$ consisting of a closed subgroup $L$ of $K$ 
and a continuous group homomorphism $\alpha:L\to G$ we define
\[ [L,\alpha]\ = \ \tr_L^K\circ\alpha^*\ \in \ \bA(G,K) \ ,\]
the natural transformation whose value at $X$ is the composite
\[  \pi_0^G (X) \ \xra{\ \alpha^*\ }\ \pi_0^L(X) \ \xra{\ \tr_L^K\ } \ \pi_0^K(X)\]
of restriction along $\alpha$ with transfer from $L$ to $K$.

If $L$ has infinite index in its normalizer, then the transfer map
$\tr_L^K$, and hence also the element $[L,\alpha]$, 
is zero by Example \ref{eg:L(H) has fix}.  
The {\em conjugate} of $(L,\alpha)$ by a pair 
$(k,g)\in K\times G$ of group elements is the pair
$(L^k,c_g\circ \alpha\circ c_k^{-1})$ consisting of the
conjugate subgroup $L^k$ and the composite homomorphism
\[  L^k \ \xra{\ c_k^{-1}\ }\ L \ \xra{\ \alpha\  }\ 
G \ \xra{\ c_g\ } \ G \ . \]
Since inner automorphisms induce the identity on equivariant homotopy groups
(compare Proposition \ref{prop:inner is identity}), 
\[ \tr_{L^k}^K \circ (c_g\circ \alpha \circ c_k^{-1})^*\ = \ 
 \tr_{L^k}^K \circ k_\star^{-1} \circ \alpha^* \circ g_\star\ = \ 
k_\star^{-1} \circ  \tr_L^K \circ \alpha^* \circ g_\star\ = \ 
 \tr_L^K \circ \alpha^* \ . \]
So the transformation $[L,\alpha]$ only depends on
the conjugacy class of $(L,\alpha)$, i.e.,
\[ [ L^k, c_g\circ \alpha\circ c_k^{-1} ]\ = \  [L,\alpha]
\text{\quad in $\bA(G,K)$.} \]

\begin{theorem}\label{thm:Burnside category basis} 
Let $K$ and $G$ be compact Lie groups. The morphism group $\bA(G,K)$ 
in the Burnside category is a free abelian group 
with basis the transformations $[L,\alpha]$,
where $(L,\alpha)$ runs over all $(K\times G)$-conjugacy classes of pairs consisting of
\begin{itemize}
\item a closed subgroup $L\leq K$ whose Weyl group $W_K L$ is finite, and
\item a continuous group homomorphism $\alpha:L\to G$.
\end{itemize}  
\end{theorem}
\begin{proof}
We let $V$ be any faithful $G$-representation.
By Proposition \ref{cor-pi_0 of B_gl} the composite
\[ \mZ\{[L,\alpha]\ |\ \ |W_KL| <\infty,\, \alpha:L\to G\}
\ \to \ \Nat(\pi_0^G,\pi_0^K) 
\ \xra{\ \ev\ } \ \pi_0^K(\Sigma^\infty_+ \bL_{G,V})\]
is an isomorphism, where the first  map takes $[L,\alpha]$ to $\tr_L^K\circ\alpha^*$
and the second map is evaluation at the
stable tautological class $e_{G,V}$.
The evaluation map is an isomorphism by 
Proposition \ref{prop:B_gl represents}, so the first map
is an isomorphism, as claimed.  
\end{proof}

Theorem \ref{thm:Burnside category basis} amounts to a complete calculation
of the Burnside category, because we know how to
express the composite of two operations, each given in
the basis of Theorem \ref{thm:Burnside category basis}, as a sum of basis elements.
Indeed, restrictions are contravariantly functorial and transfers are transitive,
and we also know how to expand a transfer followed by a restriction:
every group homomorphism is the composite of an epimorphism and a subgroup
inclusion. Inflations commute with transfers according to
Proposition \ref{prop:transfer and epi}, and the restriction of a transfer
can be rewritten via the double coset formula
as in Theorem \ref{thm:double coset formula}.

Theorem \ref{thm:Burnside category basis} 
tell us what data is necessary to specify a global functor $M:\bA\to\Ab$. 
For this, one needs to give the values $M(G)$ at all compact Lie groups $G$,
restriction maps $\alpha^*:M(G)\to M(L)$ for all continuous group homomorphisms
$\alpha:L\to G$ and transfer maps $\tr_L^K:M(L)\to M(K)$ for all
closed subgroup inclusions $L\leq K$.
This data has to satisfy the same kind of relations that the 
restriction and transfer maps for equivariant homotopy groups satisfy,
namely:
\begin{itemize}
\item the restriction maps are contravariantly functorial;
\item inner automorphisms induce the identity;
\item transfers are transitive and $\tr_K^K$ is the identity;
\item the transfer $\tr_L^K$ is zero if the Weyl group $W_K L$ is infinite;
\item transfer along an inclusion $H\leq G$ 
interacts with inflation along an epimorphism $\alpha:K\to G$ according to
\[ \alpha^*\circ \tr_H^G \ = \ 
\tr_L^K \circ (\alpha|_L)^* \ : \ M(H)\ \to \ M(K) \ ,\]
where $L=\alpha^{-1}(H)$;
\item for all pairs of closed subgroups $H$ and $K$ of $G$,
the double coset formula holds, see Theorem \ref{thm:double coset formula}.
\end{itemize}
The hypothesis that inner automorphisms act as the identity
implies that the restriction map $\alpha^*$ only 
depends on the homotopy class of $\alpha$.
Indeed, suppose that $\alpha,\alpha':K\to G$ 
are homotopic through continuous group homomorphisms. 
Then $\alpha$ and $\alpha'$ belong to the same path component of the space $\hom(K,G)$
of continuous homomorphisms, so they are conjugate by an element of $G$,
compare Proposition \ref{prop:components of hom(K,G)}.

This explicit description allows us to relate our notion of global functor
to other `global' versions of Mackey functors.
For example, our category of global functors is equivalent to the
category of {\em functors with regular Mackey structure}
in the sense of Symonds \cite[\S 3, p.\,177]{symonds-splitting}.
Our global functors are {\em not} equivalent to the global Mackey functors
in the sense of tom Dieck \cite[Ch.\,VI (8.14), Ex.\,5]{tomDieck-transformation};
indeed, in the indexing category $\Omega$ for tom Dieck's global Mackey functors
the group $\Hom_\Omega(G,K)$ has a $\mZ$-basis indexed by $(G\times K)$-conjugacy classes
$[L,\alpha]$ where the Weyl group $W_K L$ is allowed to be infinite.
As we shall explain in Remark \ref{rk:A^fin and A^c} below, 
global functors defined on finite groups
are equivalent to {\em inflation functors} in the sense of Webb \cite{webb}
and to `global $(\emptyset,\infty)$-Mackey functors' in 
the sense of Lewis \cite{lewis-projective not flat}.

\begin{eg}[Sphere spectrum]\label{eg:sphere spectrum}\index{subject}{sphere spectrum}\index{symbol}{$\mS$ - {sphere spectrum}}
The {\em sphere spectrum} $\mS$ is given by $\mS(V)= S^V$,
the one-point compactification of the inner product space $V$.
The orthogonal group acts as the one-point compactification
of its action on $V$. The structure map
$\sigma_{V,W}:S^V\sm S^W\to S^{V\oplus W}$ is the canonical homeomorphism.
The equivariant homotopy groups of the sphere spectrum are
the equivariant stable stems.
The sphere spectrum is the suspension spectrum of the
constant one-point orthogonal space $\bL_{e,0}$,
\[ \mS \ \iso \ \Sigma^\infty_+ \bL_{e,0}\ . \]
The trivial representation is faithful as a representation of the trivial group,
so $\bL_{e,0}=B_{\gl}e$ is a global classifying space for the trivial group.
The class  $1\in\pi_0(\mS)$
represented by the identity of $S^0$ is the stable tautological class $e_{e,0}$
(compare \eqref{eq:define_stable_tautological}).
By Corollary \ref{cor-pi_0 of B_gl} the group $\pi_0^K(\mS)$
is a free abelian group with basis the classes $\tr_L^K(p_L^*(1))$
where $L$ runs over all conjugacy classes of closed subgroups of $K$ 
with finite Weyl group, and where $p_L:L\to e$ is the unique homomorphism.
For finite groups, this is originally due
to Segal \cite{segal-ICM}, and for general compact Lie groups 
to tom\,Dieck, as a corollary to his splitting theorem 
(see~Satz~2 and Satz~3 of \cite{tomDieck-OrbittypenII}).
By Proposition \ref{prop:B_gl represents},
the action on the unit $1\in\pi_0(\mS)$ is an isomorphism of global functors
\[  \mA\ = \  \bA(e,-) \ \to \ \upi_0 ( \mS )\]
from the Burnside ring global functor $\mA$ to the 0-th homotopy
global functor of the sphere spectrum. 
\end{eg}

\begin{eg}\label{eg:global functor examples} 
(i) The {\em Burnside ring global functor}\index{subject}{Burnside ring global functor}\index{symbol}{$\mA$ - {Burnside ring global functor}}
is the represented global functor
$\mA=\bA(e,-)$ of morphisms out of the trivial group $e$.  
By Theorem \ref{thm:Burnside category basis}, the value
$\mA(K)=\bA(e,K)$ at a compact Lie group $K$ is a free abelian group
with basis the set of conjugacy classes of closed subgroups $L\leq K$ 
with finite Weyl group.
When $K$ is finite, then the Weyl group condition is
vacuous and $\mA(K)$ is canonically isomorphic to the Burnside ring of $K$,
by sending the operation $[L,p_L]=\tr_L^K\circ p_L^*\in \mA(K)$
to the class of the $K$-set $K/L$ (where $p_L:L\to e$ is the unique homomorphism). 
As we discussed in Example \ref{eg:sphere spectrum}, the Burnside ring 
global functor $\mA$ is realized by the sphere spectrum $\mS$.
More generally, the represented functors $\bA(G,-)$ are other examples 
of global functors, and we have seen in Proposition \ref{prop:B_gl represents}  
that the represented global functor $\bA(G,-)$ is
realized by the suspension spectrum of the global classifying
space $B_{\gl}G$.\index{subject}{global classifying space}

(ii) In \cite[Sec.\,5.5]{tomDieck-and representation theory},
tom Dieck gives a very different construction of the Burnside ring global functor.
We let $G$ be  a compact Lie group.
A {\em $G$-ENR} is a $G$-space equivariantly homeomorphic
to a $G$-retract of a $G$-invariant open subset of a finite-dimensional
$G$-representation. 
The acronym `ENR' stands for {\em euclidean neighborhood retract}. 
Examples of $G$-ENRs are smooth compact $G$-manifolds and finite $G$-CW-complexes.

Tom Dieck calls two compact $G$-ENRs $X$ and $Y$ {\em equivalent}
if for every closed subgroup $H$ of $G$ the Euler characteristics
of the $H$-fixed point spaces $X^H$ and $Y^H$ coincide;
here Euler characteristics are taken with respect to compactly supported
Alexander-Spanier cohomology, and there is some work involved in showing
that a compact $G$-ENR has a well-defined Euler characteristic.
Then $A(G)$ is defined as the set of equivalence classes of compact $G$-ENRs.
The set $A(G)$ is naturally a commutative ring, with addition induced by
disjoint union and multiplication induced by cartesian product of $G$-ENRs.

Tom Dieck shows in \cite[Prop.\,5.5.1]{tomDieck-and representation theory}, 
that $A(G)$ is a free abelian group with basis the classes of the homogeneous spaces
$G/H$, where $H$ runs over the conjugacy classes of closed subgroups with
finite Weyl group. Moreover, the class of a general compact $G$-ENR $X$
is expressed in terms of this basis by the formula
\[  [X]\ = \ {\sum}_{(H)}\ \chi^{\text{AS}}(G\bs X_{(H)})\cdot [G/H]\ ;\]
the sum is over conjugacy classes of closed subgroups, 
$X_{(H)}$ is the orbit type subspace, and $\chi^{\text{AS}}$
is the Euler characteristic based on Alexander-Spanier cohomology.
A compact $G$-ENR has only finitely many orbit types, so the sum is in fact finite.
Restriction of scalars along a continuous homomorphism $\alpha:K\to G$
induces a ring homomorphism $\alpha^*:A(G)\to A(K)$.
If $H$ is a closed subgroup of $G$, then extension of scalars --
sending an $H$-ENR $Y$ to $G\times_H Y$ --
induces an additive transfer map $\tr_H^G:A(H)\to A(G)$
that satisfies reciprocity with respect to restriction from $G$ to $H$.
By explicit comparison of bases, the homomorphisms
\[ A(G)\ \to \ \mA(G) \ , \quad [G/H]\ \longmapsto \ \tr_H^G\circ p_H^* \]
define an isomorphism of global functors.

(iii)
Given an abelian group $M$, the {\em constant global functor} $\underline{M}$
\index{subject}{global functor!constant} 
is given by $\underline{M}(G)=M$ and all restriction maps
are identity maps. The transfer $\tr_H^G:\underline{M}(H)\to\underline{M}(G)$
is multiplication by the Euler characteristic of the homogeneous space $G/H$. 
In particular, if $H$ has finite index in $G$, then
$\tr_H^G$ is multiplication by the index $[G:H]$.
In this example, the double coset formula
is a special case of the Euler characteristic formula \eqref{eq:chi(B)},
namely for the $K$-manifold $B=G/H$. 

\smallskip

\Danger There is a well-known point set level model 
of an Eilenberg-Mac\,Lane spectrum $\Hc M$\index{subject}{Eilenberg-Mac\,Lane spectrum!of an abelian group}
that we recall in Construction \ref{con:HM} below.
However, contrary to what one may expect,
the 0-th homotopy group global functor $\upi_0(\Hc M)$ is {\em not}
isomorphic to the constant global functor $\underline{M}$.
More precisely, the restriction map $p_G^*: \pi_0^e (\Hc M)\to \pi_0^G(\Hc M)$
is an isomorphism for finite groups $G$, but not for general compact Lie groups,
see Example \ref{eg:HZ for abelian G}.

\smallskip

(iv)
The {\em unitary representation ring global functor}\index{subject}{representation ring!unitary}
$\bRU$ assigns to a compact Lie group $G$ the 
unitary representation ring $\bRU(G)$, i.e., the Grothendieck
group of finite-dimensional complex $G$-representations,
with product induced by tensor product of representations.
The restriction maps $\alpha^*:\bRU(G)\to \bRU(K)$
are induced by restriction of representations along 
a continuous homomorphism $\alpha:K\to G$.
The transfer maps $\tr_H^G:\bRU(H)\to \bRU(G)$ along a closed subgroup inclusion
$H\leq G$ are given by the {\em smooth induction} 
of Segal \cite[\S\,2]{segal-representation}.
If $H$ is a subgroup of finite index of $G$, then this induction sends the class 
of an $H$-representation to the induced $G$-representation $\map^H(G,V)$;
in general, induction may send actual representations to virtual representations.
In the generality of compact Lie groups, 
the double coset formula for $\bRU$ was proved by Snaith \cite[Thm.\,2.4]{Snaith-Brauer}.
We look more closely at the representation ring global functor
in Example \ref{eg:RU_as_global_power},
and we show in Theorem \ref{thm:pi_0 KU is RU} that $\bRU$ 
is realized by the global $K$-theory spectrum $\bKU$
(see Construction \ref{con:global KU}).

(v)
Given any generalized cohomology theory $E$ (in the non-equivariant sense),
we can define a global functor $\underline{E}$  by setting
\[   \underline{E}(G) \ = \ E^0(B G) \ , \]
the 0-th $E$-cohomology of a classifying space of the group $G$.
The contravariant functoriality in group homomorphisms $\alpha:K\to G$
comes from the covariant functoriality of the classifying space construction.
The transfer map for a subgroup inclusion $H\leq G$ 
comes from the stable transfer map (i.e., Becker-Gottlieb transfer)\index{subject}{Becker-Gottlieb transfer}\index{subject}{transfer!Becker-Gottlieb|see{Becker-Gottlieb transfer}} 
\[ \Sigma^\infty_+ B G \ \to \  \Sigma^\infty_+ B H \ .\]
The double coset formula was proved 
in this context by Feshbach \cite[Thm.\,II.11]{feshbach}.
We will show in Proposition \ref{prop:global homotopy of b E}
that the global functor $\un{E}$ is realized by an orthogonal spectrum $b E$,
the `global Borel theory' associated to $E$.
\end{eg}

\begin{rk}\index{subject}{Mackey functor} 
If we fix a finite group $G$ and let $H$ run through all subgroups of $G$, 
then the collection of $H$-equivariant homotopy groups
$\pi_0^H(X)$ of an orthogonal spectrum $X$ forms a
{\em $G$-Mackey functor} (see Definition \ref{def:G-Mackey}),
with respect to restriction to subgroups, conjugation and transfer maps.
As we already discussed in Remark \ref{rk:homotopy properties},
not all $G$-Mackey functors arise this way.
To illustrate this we compare Mackey functors 
for the group $C_3=\{1,\tau,\tau^2\}$ with three elements
to additive functors on the full subcategory of $\bA$ spanned by
the group $e$ and $C_3$. Generating operations can be displayed as
follows:
\[ \xymatrix@R=1mm@!C=4mm{ 
&&&&&&&\\
 F(C_3)  \ar[rr]_{\res} &&
F(e) \ar@<.9ex>@/^1pc/[ll]^{p^*} \ar@/_1pc/[ll]_{\tr}
&& F(C_3)\ar[rr]_{\res} &&
F(e)  \ar@/_1pc/[ll]_{\tr} \\
 \ar@<2ex>@(dl,ul)[uu]^{\alpha^*} &&&&&&  \ar@<-2ex>@(dr,ur)[uu]_{\tau} \\
&&&\\
&\text{global functor on $C_3$ and $e$} &&&&
\text{$C_3$-Mackey functor}
 } \]
Here $\res=\res^{C_3}_e$ and $\tr=\tr_e^{C_3}$ are the restriction and
transfer maps that are present in both cases.
A global functor also comes with inflation maps along the epimorphism 
$p:C_3\to e$ and along the automorphism $\alpha:C_3\to C_3$ with $\alpha(\tau)=\tau^2$, 
and the relations are
\[ \res\circ p^*\ = \ \Id \text{\qquad and\qquad} \res\circ \tr\ = \ 3\cdot \Id \]
as well as $\alpha^*\circ\alpha^*=\Id$,
$\alpha^*\circ p^*=p^*$, $\res\circ\alpha^*=\res$ and
$\alpha^*\circ \tr=\tr$.
In contrast, $C_3$-Mackey functors have an additional action of $C_3$
(the Weyl group of $e$ in $C_3$) on $F(e)$, and this action satisfies the relation
\[ \res\circ \tr\ = \ \Id\ +\ \tau \ +\ \tau^2\ .\]
\end{rk}

We will now use the exterior homotopy group pairings \eqref{eq:def_boxtimes}
\[  \boxtimes \ : \ \pi_0^G(X) \ \times \ \pi_0^K (Y) \ \to \ \pi_0^{G\times K}(X\sm Y)  \]
to define a biadditive functor
\[ \times \ : \ \bA\times\bA\ \to \ \bA \]
that is given on objects by the product of Lie groups,
and that extends to a symmetric monoidal structure on the Burnside category $\bA$. 

\begin{prop}
Let $G, G', K$ and $K'$ be compact Lie groups.
Given operations $\tau\in \bA(G,K)$ and $\psi\in \bA(G',K')$, 
there is a unique operation
\[  \tau\times\psi\ \in\ \bA(G\times G',\, K\times K')\]
with the following property:
for all orthogonal spectra $X$ and $Y$ and all classes $x\in\pi_0^G(X)$
and $y\in\pi_0^{G'}(Y)$ the relation
\begin{equation}\label{eq:characterize tau x tau}
 (\tau\times\psi)(x\boxtimes y)\ = \ 
\tau(x)\boxtimes \psi(y)
\end{equation}
holds in $\pi_0^{K\times K'}(X\sm Y)$.
\end{prop}
\begin{proof}
We choose a faithful $G$-representation $V$ and a faithful $G'$-representation $V'$,
which have associated stable tautological 
classes \eqref{eq:define_stable_tautological}\index{subject}{stable tautological class}
\[  e_{G,V}\ \in \ \pi_0^G (\Sigma^\infty_+ \bL_{G,V}) 
\text{\qquad and\qquad}  e_{G',V'}\ \in \ \pi_0^{G'}(\Sigma^\infty_+ \bL_{G',V'}) \ .\]
Combining the isomorphism \eqref{eq:box2smash}
with the one from Example \ref{eg:box of free orthogonal} shows that
the orthogonal spectrum $\Sigma^\infty_+ \bL_{G,V}\sm \Sigma^\infty_+ \bL_{G',V'}$
is isomorphic to $\Sigma^\infty_+ \bL_{G\times G',V\oplus V'}$ in a way that matches the class
$e_{G,V}\boxtimes e_{G',V'}$ with the class $e_{G\times G',V\oplus V'}$.
Proposition \ref{prop:B_gl represents} then shows that the pair
$(\Sigma^\infty_+ \bL_{G,V}\sm \Sigma^\infty_+ \bL_{G',V'},e_{G,V}\boxtimes e_{G',V'})$ 
represents the functor $\pi_0^{G\times G'}$.
There is thus a unique operation $\tau\times\psi\in\bA(G\times G',K\times K')$
that satisfies
\begin{equation}\label{eq:universal case tau x psi}
   (\tau\times\psi)(e_{G,V}\boxtimes e_{G',V'})\ = \ 
\tau(e_{G,V})\boxtimes \psi(e_{G',V'})   
\end{equation}
in $\pi_0^{K\times K'}(\Sigma^\infty_+ \bL_{G,V}\sm \Sigma^\infty_+ \bL_{G',V'})$.

The relation \eqref{eq:universal case tau x psi} 
is a special case of \eqref{eq:characterize tau x tau},
and it remains to show that the operation $\tau\times\psi$
satisfies the relation \eqref{eq:characterize tau x tau}
in complete generality. As we already argued in the proof 
of Proposition \ref{prop:B_gl represents}, there is a $G$-representation $W$
and a morphism of orthogonal spectra
\[ f \ : \ \Sigma^\infty_+ \bL_{G,V\oplus W}\ \to \ X  \]
that satisfies $f_*(e_{G,V\oplus W})=x$.
Similarly, there is a $G'$-representation $W'$
and a morphism of orthogonal spectra
\[ f' \ : \ \Sigma^\infty_+ \bL_{G',V'\oplus W'}\ \to \ Y  \]
that satisfies $f'_*(e_{G',V'\oplus W'})=y$. The morphism
\[ \Sigma^\infty_+\rho_{G,V,W}\ : \ \Sigma^\infty_+\bL_{G,V\oplus W}\ \to \ 
\Sigma^\infty_+\bL_{G,V}\text{\quad satisfies\quad}
 (\Sigma^\infty_+\rho_{G,V,W})(e_{G,V\oplus W})\ = \ e_{G,V}\ , \]
and similarly for the triple $(G',V',W')$.
Naturality then yields
\begin{align*}
(\Sigma^\infty_+\rho_{G,V,W}\sm\Sigma^\infty_+&\rho_{G',V',W'})_*( (\tau\times\psi)(e_{G,V\oplus W}\boxtimes e_{G',V'\oplus W'}))\\ 
&= \ (\tau\times\psi)( (\Sigma^\infty_+\rho_{G,V,W})_*(e_{G,V\oplus W})\boxtimes 
(\Sigma^\infty_+\rho_{G',V',W'})_*( e_{G',V'\oplus W'}))\\
&= \ (\tau\times\psi)( e_{G,V}\boxtimes  e_{G',V'})\\
_\eqref{eq:universal case tau x psi} 
&= \ \tau( e_{G,V})\boxtimes \psi(e_{G',V'})\\
&= \ \tau((\Sigma^\infty_+\rho_{G,V,W})_*(e_{G,V\oplus W}))\boxtimes 
\psi((\Sigma^\infty_+\rho_{G',V',W'})_*( e_{G',V'\oplus W'}))\\
&= \ 
(\Sigma^\infty_+\rho_{G,V,W}\sm \Sigma^\infty_+\rho_{G',V',W'})_*
(\tau(e_{G,V\oplus W})\boxtimes \psi(e_{G',V'\oplus W'}))\ .
\end{align*}
The morphism 
$\Sigma^\infty_+\rho_{G,V,W}\sm\Sigma^\infty_+\rho_{G',V',W'}$
is isomorphic to $\Sigma^\infty_+\rho_{G\times G',V\oplus V',W\oplus W'}$.
The morphism $\rho_{G\times G',V\oplus V',W\oplus W'}$ is a global equivalence of
orthogonal spaces (by Proposition \ref{prop:free_orthogonal_space}~(ii)),
and so the morphism 
$\Sigma^\infty_+\rho_{G,V,W}\sm\Sigma^\infty_+\rho_{G',V',W'}$
is a global equivalence of orthogonal spectra 
(by Corollary \ref{cor:suspension spectrum globally homotopical}).
So in particular it induces an isomorphism on $\pi_0^{K\times K'}$, 
and we can conclude that
\[   (\tau\times\psi)(e_{G,V\oplus W}\boxtimes e_{G',V'\oplus W'})\ = \ 
\tau(e_{G,V\oplus W})\boxtimes \psi(e_{G',V'\oplus W'}) \ .\]
Now the relation \eqref{eq:characterize tau x tau} follows by simple naturality:
\begin{align*}
 (\tau\times\psi)(x\boxtimes y)\ &= \
 (\tau\times\psi)(f_*(e_{G,V\oplus W})\boxtimes f'_*(e_{G',V'\oplus W'}))\\ 
&= \ 
 (\tau\times\psi)((f\sm f')_*(e_{G,V\oplus W}\boxtimes e_{G',V'\oplus W'}))\\ 
&= \ 
(f\sm f')_*((\tau\times\psi)(e_{G,V\oplus W}\boxtimes e_{G',V'\oplus W'}))\\ 
&= \ 
(f\sm f')_*(\tau(e_{G,V\oplus W})\boxtimes\psi(e_{G',V'\oplus W'}))\\ 
&= \ 
f_*(\tau(e_{G,V\oplus W}))\boxtimes f'_*(\psi(e_{G',V'\oplus W'}))\\ 
&= \ 
\tau(f_*(e_{G,V\oplus W}))\boxtimes \psi(f'_*(e_{G',V'\oplus W'}))\
= \ \tau(x)\boxtimes \psi(y)  \qedhere
\end{align*}
\end{proof}

\begin{eg} We calculate the product of two generating operations in the Burnside category.
We recall that for a pair $(L,\alpha)$ consisting of a closed subgroup $L$ of $K$ 
and a continuous homomorphism $\alpha:L\to G$ we defined
\[ [L,\alpha]\ = \ \tr_L^K\circ\alpha^* \ \in \ \bA(G,K) \ .\]
By Theorem \ref{thm:Burnside category basis}, a certain subset
of these operations forms a basis of the abelian group $\bA(G,K)$. 
Using parts~(v) and~(vi) of Theorem \ref{thm:external product properties} we deduce
\begin{align*}
   [L\times L',\alpha\times\alpha'](x\boxtimes y)\ &= \ 
   \tr_{L\times L'}^{K\times K'}( (\alpha\times\alpha')^*(x\boxtimes y))\\ 
&= \  \tr_{L\times L'}^{K\times K'}( \alpha^*(x)\boxtimes (\alpha')^*(y))\\ 
&= \  \tr_L^K( \alpha^*(x))\boxtimes \tr_{L'}^{K'}((\alpha')^*(y))\
= \  [L,\alpha](x)\boxtimes [L',\alpha'](y)\ .
\end{align*}
So the operation $[L\times L',\alpha\times\alpha']$ has the property
that characterizes the operation $ [L,\alpha]\times [L',\alpha']$.
Hence the monoidal product in $\bA$ satisfies
\begin{equation}\label{eq:times_on_[L,a]}
 [L,\alpha]\times [L'\alpha']\ = \ [L\times L',\alpha\times\alpha']\ .
  \end{equation}
\end{eg}

Now we are ready for the monoidal structure of the Burnside category.

\begin{theorem}\label{thm:bA monoidal} 
Let $G,G',K$ and $K'$ be compact Lie groups.
  \begin{enumerate}[\em (i)]
  \item The map
    \[ \times \ : \ \bA(G,K)\times\bA(G',K')\ \to \ \bA(G\times G',K\times K') \]
    is biadditive.
  \item As the Lie groups vary, the maps of~{\em (i)}
    form a functor $\times:\bA\times\bA\to \bA$.
  \item 
    The restriction operations along the group isomorphisms
    \begin{align*}
 a_{G,G',G''}\ : \ G\times (G'\times G'') \ &\iso \  (G\times G')\times G''\ ,\\ 
    \tau_{G,G'}\ : \ G\times G'\ \iso \  G'\times  G&\text{\quad respectively\quad} 
    G\times e\ \iso \ G\ \iso \ e\times G       
    \end{align*}
    make the functor $\times$ into a symmetric monoidal structure on
    the global Burnside category.
  \end{enumerate}
\end{theorem}
\begin{proof}
(i) We show additivity in the first variable, the other case being analogous. 
The relation
\begin{align*}
((\tau\times \psi)+(\tau'\times\psi))(x\boxtimes y)
 &= \  (\tau\times \psi)(x\boxtimes y)+(\tau'\times\psi)(x\boxtimes y)\\
 &= \  (\tau(x)\boxtimes\psi(y))+(\tau'(x)\boxtimes \psi(y))  \\
 &= \  (\tau(x)+\tau'(x))\boxtimes \psi(y)  \
 = \   (\tau+\tau')(x)\boxtimes \psi(y)  
\end{align*}
shows that the operation $(\tau\times\psi')+(\tau'\times\psi)$
has the property that characterizes the operation 
$(\tau +\tau')\times\psi$. So
\[ (\tau+\tau')\times \psi\ = \ (\tau\times\psi)+(\tau'\times\psi)\ . \]

(ii) The relation
\begin{align*}
 ((\tau'\times\psi')\circ(\tau\times\psi))(x\boxtimes y)\ &= \
 (\tau'\times\psi')( (\tau\times\psi)(x\boxtimes y))  \\ 
&= \ (\tau'\times\psi')( \tau(x)\boxtimes \psi(y) )  \\
&= \
 \tau'(\tau(x))\boxtimes \psi'(\psi(y))  \
= \ (\tau'\circ\tau)(x)\boxtimes (\psi'\circ\psi)(y)
\end{align*}
shows that the operation $(\tau'\times\psi')\circ(\tau\times\psi)$
has the property that characterizes the operation 
$(\tau'\circ\tau)\times(\psi'\circ \psi)$. So
\[ (\tau'\circ\tau)\times(\psi'\circ \psi)\ = \ 
(\tau'\times\psi')\circ(\tau\times\psi)\ . \]
A similar (but even shorter) argument shows that $\Id_G\times\Id_{G'}=\Id_{G\times G'}$.

(iii) We start with naturality of the associativity isomorphism.
We consider operations $\tau\in\bA(G,K)$, $\psi\in\bA(G',K')$ 
and $\kappa\in\bA(G'',K'')$. The relation
\begin{align*}
(\tau\times(\psi\times \kappa))(a_{G,G',G''}^*((x\boxtimes y)\boxtimes z)) \ 
&= \ ((\tau\times(\psi\times \kappa))(x\boxtimes (y\boxtimes z))) \\ 
&= \ \tau(x)\boxtimes(\psi(y)\boxtimes \kappa(z))\\
&= \ a_{K,K',K''}^*( (\tau(x)\boxtimes\psi(y))\boxtimes \kappa(z))\\
&= \ a_{K,K',K''}^*(((\tau\times\psi)\times\kappa)((x\boxtimes y)\boxtimes z))
\end{align*}
shows that the operation 
$(a_{K,K',K''}^*)^{-1}\circ(\tau\times(\psi\times \kappa))\circ a_{G,G',G''}^*$
has the property that characterizes the operation 
$(\tau\times\psi)\times\kappa$. So
\[ (\tau\times(\psi\times \kappa))\circ a_{G,G',G''}^* \ = \ 
a_{K,K',K''}^*\circ((\tau\times\psi)\times\kappa)\ . \]
The arguments for naturality of the unit and symmetry isomorphisms
are similar, and we omit them.

The unit, associativity (pentagon) and symmetry (hexagon)
coherence relations in $\bA$ follow from the corresponding coherence relations
for the product of groups, and the fact that passage from
group homomorphisms to restriction operations is functorial.
\end{proof}

\begin{rk}\label{rk:A^fin and A^c}
The full subcategory $\bA_\Fin$ of the Burnside category $\bA$ spanned 
by {\em finite} groups has a different, more algebraic description, as we shall now
recall. This alternative description is in terms of `bisets',
and is often taken as the definition in algebraic treatments of global functors
(which are then sometimes called `biset functors').\index{subject}{biset functor}
The category of `global functors on finite groups',
i.e., additive functors from $\bA_\Fin$ to abelian groups, is thus equivalent to
the category of `inflation functors' in the sense of Webb \cite[p.271]{webb}
and to the `global $(\emptyset,\infty)$-Mackey functors' in 
the sense of Lewis \cite{lewis-projective not flat}. 

We define the additive {\em biset category} $\mA^\text{c}$. 
The objects of $\mA^\text{c}$ are all finite groups. The abelian group $\mA^\text{c}(G,K)$ 
of morphisms from a group $G$ to $K$ is the Grothendieck
group of finite $K$-$G$-bisets where the right $G$-action is free.
In the special case $G=e$ of the trivial group as source we
obtain $\mA^\text{c}(e,K)$, the Burnside ring of finite $K$-sets.
Composition
\[ \circ \ : \ \mA^\text{c}(K,L) \times \mA^\text{c}(G,K) \ \to \ \mA^\text{c}(G,L)\]
is induced by the balanced product over $K$, i.e., it is the
biadditive extension of 
\[ (S,T) \ \longmapsto \ S\times_K T \ .\] 
Here $S$ has a left $L$-action and a commuting free right $K$-action,
whereas $T$ has a left $K$-action and a commuting free right $G$-action.
The balanced product $S\times_K T$ than inherits a left $L$-action from $S$
and a free right $G$-action from $T$.
Since the balanced product is associative up to isomorphism,
this defines a pre-additive category.
So the category $\mA^\text{c}$ is the `group completion' of the category $\mA^+_\Fin$, 
the restriction of the effective Burnside category $\mA^+$ 
(compare Construction \ref{con:define P and A^+}) to finite groups.

We define additive maps
\[ \Psi_{G,K}\ : \ \bA(G,K)\ \to \ \mA^\text{c}(G,K) \]
that form an additive equivalence of categories (restricted to finite groups).
The map $\Psi_{G,K}$ sends a basis element $[L,\alpha]$ to the
class of the $K$-$G$-biset
\[ K\times_{(L,\alpha)} G \ = \ K\times G/(kl,g)\sim(k,\alpha(l)g) \]
whose right $G$-action is free.
Every transitive $G$-free $K$-$G$-biset is isomorphic to one of this form, 
and $K\times_{(L,\alpha)} G$ is isomorphic, as a $K$-$G$-biset,
to $K\times_{(L',\alpha')} G$ if and only if $(L,\alpha)$ is conjugate to
$(L',\alpha')$. So the map $\Psi_{G,K}$ sends the basis of
$\bA(G,K)$ of Theorem \ref{thm:Burnside category basis} 
to a basis of $\mA^\text{c}(G,K)$, and it is thus an isomorphism.

We claim that the maps $\Psi_{G,K}$ form a functor as $G$ and $K$ vary through
all finite groups; this then shows that $\Psi$ is an additive equivalence
of categories from the full subcategory of $\bA_\Fin$ to $\mA^\text{c}$.
The functoriality boils down to the fact that in both categories
restriction, inflation and transfer interact with each other in exactly the
same way. We omit the details.

The restriction of the monoidal structure on the Burnside category
to finite groups has an interpretation in terms of the cartesian product of
bisets: under the equivalence of categories $\Psi:\bA_\Fin\iso \mA^c$, 
it corresponds to the monoidal structure
\[  \mA^\text{c}(G,K)\times \mA^\text{c}(G',K')\ \to \ 
\mA^\text{c}(G\times G', K\times K')\ , 
\ ([S], [S']) \ \longmapsto [S\times S']\ .\]
Here $S$ is a right free $K$-$G$-biset
and $S'$ is a right free $K'$-$G'$-biset; the cartesian product
$S\times S'$ is then a right free $(K\times K')$-$(G\times G')$-biset.
Indeed, the equivalence $\Psi:\bA_\Fin\iso \mA^c$ 
sends the basis element $[L,\alpha]\in\bA(G,K)$ 
to the class of the biset $K\times_{(L,\alpha)} G$.
So the equivalence is monoidal because 
\[ (K\times_{(L,\alpha)} G)\times (K'\times_{(L',\alpha')} G')
\text{\quad and\quad} (K\times K')\times_{(L\times L',\alpha\times \alpha')} G\times G' \]
are isomorphic as  $(K\times K')$-$(G\times G')$-bisets.
\end{rk}

\index{subject}{Burnside category|)}

\begin{construction}[Box product of global functors]\label{con:Box product}
Since global functors are additive functors on the Burnside category $\bA$,
the symmetric monoidal product on $\bA$ gives rise to a
symmetric monoidal convolution product on the category of global functors.
This is a special case of the general construction of Day \cite{day-closed} that
we review in Appendix \ref{app:enriched functors}.
We now make this convolution product more explicit.
We denote by $\bA\tensor\bA$ the pre-additive category whose objects
are pairs of compact Lie groups, and with morphism groups
\[ (\bA\tensor\bA)((G,G'),(K,K')) \ = \  \bA(G,K) \tensor \bA(G',K') \ . \]
We let $F, F'$ and $F''$ be global functors.
We denote by $F\tensor F':\bA\tensor\bA\to \Ab$ the objectwise
tensor product given on objects by
\[ (F\tensor F')(G,G')\ = \ F(G)\tensor F'(G') \ .\]
A {\em bimorphism} is a natural transformation 
\[ F\tensor F'\ \to \ F''\circ \times  \]
of additive functors on the category $\bA\tensor\bA$.\index{subject}{bimorphism!of global functors}
Since the morphism groups in the Burnside category are
generated by transfers and restrictions, this means more explicitly
that a bimorphism is a collection of group homomorphisms
\[ b_{G,G'} \ : \ F(G) \tensor F'(G') \ \to \ F''(G\times G') \]
for all compact Lie groups $G$ and $G'$, such that for all continuous homomorphisms
$\alpha:K\to G$ and $\alpha':K'\to G'$ 
and for all closed subgroups $H\leq G$ and $H'\leq G'$ the  diagram
\begin{equation*}  \begin{aligned}
\xymatrix@C=15mm@R=10mm{ 
F(H) \tensor F'(H') \ar[r]^-{\tr_H^G\tensor \tr_{H'}^{G'}} \ar[d]_{b_{H,H'}}&
F(G) \tensor F'(G') \ar[r]^-{\alpha^*\tensor (\alpha')^*} \ar[d]_{b_{G,G'}}&
F(K) \tensor F'(K') \ar[d]^{b_{K,K'}}  \\
F''(H\times H') \ar[r]_-{\tr_{H\times H'}^{G\times G'}} &
F''(G\times G') \ar[r]_-{(\alpha\times\alpha')^*} &  F''(K\times K')}
  \end{aligned}\end{equation*}
commutes.
Here we exploited that the generating operations multiply as 
described in \eqref{eq:times_on_[L,a]}.
Equivalently: for every compact Lie group $G$ the maps $\{b_{G,G'}\}_{G'}$
form a morphism of global functors $F(G)\tensor F'(-)\to F''(G\times-)$
and for every compact Lie group $G'$ the maps $\{b_{G,G'}\}_G$
form a morphism of global functors $F(-)\tensor F'(G')\to F''(-\times G')$.

A {\em box product}\index{subject}{box product!of global functors|(}\index{symbol}{$\Box$ - {box product of global functors}} 
of $F$ and $F'$ is a universal example of a global functor
with a bimorphism from $F$ and $F'$.
More precisely, a box product 
is a pair $(F\Box F',i)$ consisting of a global functor $F\Box F'$
and a universal bimorphism $i$,
i.e., such that for every global functor $F''$ the map
\[ \GF(F\Box F',F'') \ \to \ \text{Bimor}((F,F'),F'') \ , \quad f\longmapsto f i  \]
is bijective. Box products exist by the general theory
(see Proposition \ref{prop:box exists}), 
and they are unique up to preferred isomorphism
(see Remark \ref{rk:box unique}). 
The universal property guarantees
that given any collection of choices of box product $F\Box F'$
for all pairs of global functors, $F\Box F'$ extends to
an additive functor in both variables.
Moreover, there are preferred associativity and commutativity isomorphisms
that make the box product into a symmetric monoidal structure on the category
of global functors, with the Burnside ring global functor $\mA$ 
as a strict unit object, compare Theorem \ref{thm:symmetric monoidal}.
The box product of representable global functors is again representable,
by Remark \ref{rk:box representable}.
\end{construction}

\begin{rk}\label{rk:Box is right exact}
We claim that the box product is right exact in both variables. To see this we 
recall from Remark \ref{rk:box as enriched coend} that the values 
of a Day type convolution product can be described as an enriched coend. 
In our present situation this says that
$F\Box M$ is a cokernel of a certain homomorphism of global functors
\begin{align*}
 d\ : \ 
\bigoplus_{H,H',G,G'} 
\bA(H\times H',-)\tensor\bA(G,H)&\tensor\bA(G',H')\tensor F(G)\tensor M(G') \\
\to \  &\bigoplus_{G,G'} \bA(G\times G',-)\tensor F(G)\tensor M(G')  \ .\nonumber
\end{align*}
The left sum is indexed over all quadruples $(H,H',G,G')$ in a set
of representatives of isomorphism classes of compact Lie groups,
the right sum is indexed over pairs $(G,G')$ of such groups.
The map $d$ is the difference of two homomorphisms; one of them
sums the tensor products of 
\[ \bA(H\times H',-)\tensor\bA(G,H)\tensor\bA(G',H')\ \to \ \bA(G\times G',-)\ 
, \  \varphi\tensor\tau\tensor\tau' \ \longmapsto \ \varphi\circ(\tau\times\tau') \]
and the identity on $F(G)\tensor M(G')$.
The other map sums the tensor product of the identity of 
the global functor $\bA(H\times H',-)$
and the action maps $\bA(G,H)\tensor F(G)\to F(H)$
respectively  $\bA(G',H')\tensor M(G')\to M(H')$.
Cokernels of global functors are calculated objectwise, so
the value $(F\Box M)(K)$ is a cokernel of the morphism of abelian groups
that we obtain by plugging $K$ into the free variable above.
Theorem \ref{thm:Burnside category basis} 
describes explicit free generators for the morphism groups
in the Burnside category; using this, the value $(F\Box M)(K)$ 
can be expanded into a cokernel of a morphism between two huge sums
of tensor products of values of $F$ and $M$.

Now we consider a short exact sequence of global functors 
\[ 0 \to M \to M' \to M''\to 0 \ .\]
This gives rise to a commutative diagram
of global functors
\[ \xymatrix@C=5mm@R=5mm{  \bigoplus  \bA(H\times H',-)\tensor \bA(G,H)\tensor\bA(G',H') \tensor F(G)\tensor M(G')\ar[r]\ar[d] &
\bigoplus \bA(G\times G',-)\tensor F(G)\tensor M(G') \ar[d] \\
 \bigoplus \bA(H\times H',-)\tensor \bA(G,H)\tensor\bA(G',H')\tensor F(G) \tensor M'(G')\ar[r]\ar[d] &
\bigoplus \bA(G\times G',-)\tensor F(G)\tensor M'(G') \ar[d] \\
 \bigoplus \bA(H\times H',-)\tensor \bA(G,H)\tensor\bA(G',H')\tensor F(G) \tensor M''(G')\ar[r]\ar[d] &
\bigoplus \bA(G\times G',-)\tensor F(G)\tensor M''(G') \ar[d] \\
0 & 0 } \]
with exact columns. The induced sequence of horizontal cokernels
\[F\Box M \to F\Box M' \to F\Box M''\to 0 \]
is thus also exact.
\end{rk}

\Danger While the box product of global functors shares many properties with
the tensor product of modules over a commutative ring,
the constructions differ fundamentally in one aspect:
{\em projectives are not generally flat} in the category of global functors. 
In other words, for most projective global functors $P$,
the functor $-\Box P$ does {\em not} send monomorphisms to monomorphisms. 
This kind of phenomenon has been analyzed in
great detail by Lewis in \cite{lewis-projective not flat}; 
Lewis' notation for the category $\mA^c$ is $\mathcal{B}_*(\emptyset,\infty)$,
our $\Fin$-global functors are his `global $(\emptyset,\infty)$-Mackey functors'
and the category of $\Fin$-global functors is denoted $\mathfrak{M}_*(\emptyset,\infty)$.
Theorem~6.10 of \cite{lewis-projective not flat} shows that 
the representable functor $\bA_{C_p}$ is not flat, 
where $C_p$ is a cyclic group of prime order $p$.

\index{subject}{box product!of global functors|)}

\medskip

We now remark that bimorphisms of global functors can be identified
with another kind of structure that we call `diagonal products'.

\begin{defn}\label{def:diagonal product}
Let $X,Y$ and $Z$ be global functors. A {\em diagonal product}
is a natural transformation $X\tensor Y\to Z$ of Rep-functors
to abelian groups that satisfies reciprocity,
where $X\tensor Y$ is the objectwise tensor product.
\end{defn}

More explicitly, a diagonal product consists of additive maps
\[ \nu_G \ : \ X(G)\tensor Y(G) \ \to \ Z(G)\]
for every compact Lie group $G$ that are natural for restriction along
continuous homomorphisms and satisfy the reciprocity relation
\[ \tr^G_H (\nu_H(x \tensor \res^G_H(y))) \ = \ \nu_G(\tr^G_H(x) \tensor y) \]
for all closed subgroups $H$ of $G$ and all classes $x\in X(H)$ and $y\in Y(G)$, 
and similarly in the other variable.

\begin{rk}\label{rk:diagonal versus external products}
Any bimorphism $\mu:(X, Y)\to Z$ gives rise to a diagonal product
as follows. For a group $G$ we define $\nu_G$ as the composite
\[ X(G)\otimes Y(G) \ \xra{\mu_{G,G}} \
 Z(G\times G)\ \xra{\Delta^*_G} \ Z(G) \]
where $\Delta_G:G\to G\times G$ is the diagonal. For a group homomorphism
$\alpha:K\to G$ we have $\Delta_G\circ\alpha=(\alpha\times\alpha)\circ\Delta_K$,
so the following diagram commutes:
\[\xymatrix@C=12mm{
X(G)\otimes Y(G) \ar[r]^-{\mu_{G,G}} \ar[d]_{\alpha^*\tensor\alpha^*} &
Z(G\times G)\ar[r]^-{\Delta_G^*} \ar[d]^{(\alpha\times\alpha)^*} & 
Z(G)\ar[d]^{\alpha^*}\\
X(K)\otimes Y(K) \ar[r]_-{\mu_{K,K}} &
Z(K\times K)\ar[r]_-{\Delta_K^*}  & Z(K)}\]
Since there is only one double coset for the left $\Delta_G$-action
and the right $(H\times G)$-action on $G\times G$, the double coset formula becomes 
\[ \Delta_G^*\circ\tr_{H\times G}^{G\times G} \ = \ \tr_H^G \circ 
\Delta_H^*\circ\res^{H\times G}_{H\times H}\ .\]
We conclude that
\begin{align*}
 \tr^G_H(\nu_H(x\tensor\res^G_H(y))) \ &= \ 
 \tr_H^G(\Delta_H^*(\res^{H\times G}_{H\times H}(\mu_{H,G}(x\tensor y)))) \\ 
&= \ \Delta_G^*(\tr_{H\times G}^{G\times G}(\mu_{H,G}(x\tensor y))) \\ 
&= \ \Delta_G^*(\mu_{G,G}(\tr^G_H(x)\tensor y)) \ = \  \nu_G(\tr^G_H(x)\tensor y)\ ,
\end{align*}
the reciprocity relation for the diagonal product $\nu$.
The reciprocity in the other variable is similar.

Conversely, given a diagonal product $\nu$, we define a bimorphism 
with $(G,K)$-component as the composite
\[ X(G)\tensor Y(K) \ \xra{\ p_G^*\tensor p_K^*\ } \
 X(G\times K)\tensor Y(G\times K)\  \xra{\ \nu_{G\times K}\ }\
 Z(G\times K)\ , \]
where $p_G:G\times K\to G$ and $p_K:G\times K\to K$ are the projections.
If the diagonal product $\nu$ was defined from an external product $\mu$ as above, then
\begin{align*}
  \nu_{G\times K}\circ (p_G^*\tensor p_K^*)\ &= \ 
  \Delta_{G\times K}^* \circ\mu_{G\times K,G\times K}\circ (p_G^*\tensor p_K^*)\\ 
&= \   \Delta_{G\times K}^*\circ (p_G\times p_K)^*\circ\mu_{G,K}\ = \ \mu_{G,K}
\end{align*}
because the composite $(p_G\times p_K)\circ\Delta_{G\times K}$
is the identity. So the external product can be recovered from
the diagonal product.

Given homomorphisms $\alpha:G\to G'$ and
$\beta:K\to K'$, we have $p_{G'}\circ(\alpha\times\beta)=\alpha\circ p_G$ and
$p_{K'}\circ(\alpha\times\beta)=\beta\circ p_K$, so the left part of
the diagram
\[\xymatrix@C=12mm{
 X(G)\tensor Y(K) \ar[r]^-{p_G^*\tensor p_K^*} \ar[d]_{\alpha^*\tensor\beta^*}&
 X(G\times K)\tensor Y(G\times K) \ar[r]^-{\nu_{G\times K}} 
\ar[d]^{(\alpha\times\beta)^*\tensor(\alpha\times\beta)^*} &
 Z(G\times K)\ar[d]^{(\alpha\times\beta)^*} \\ 
 X(G')\tensor Y(K') \ar[r]_-{p_{G'}^*\tensor p_{K'}^*} &
 X(G'\times K')\tensor Y(G'\times K') \ar[r]_-{\nu_{G'\times K'}} &
 Z(G'\times K') }\]
commutes. The right part commutes by naturality of the diagonal
product $\nu$.

For naturality with respect to transfers we let $H$ be a closed subgroup of $G$,
and we consider classes $x\in X(H)$ and $y\in Y(K)$. Then
\begin{align*}
 \tr_{H\times K}^{G\times K}(\mu_{H,K}(x\tensor y)) \ &= \ 
\tr_{H\times K}^{G\times K}(\nu_{H\times K}(p_H^*(x)\tensor \res^{G\times K}_{H\times K}(p_K^*(y)))) \\ 
&=\  \nu_{G\times K}(\tr_{H\times K}^{G\times K}(p_H^*(x))\tensor p_K^*(y)) \\ 
&=\
\nu_{G\times K}(p_G^*(\tr_H^G(x))\tensor p_K^*(y))\ = \ \mu_{G,K}(\tr_H^G(x)\tensor y) \ .
\end{align*}
The second equality is reciprocity, the third is compatibility of transfer 
and inflation.
The argument for transfer naturality in the $K$-variable is similar.
\end{rk}

\section{Global model structures for orthogonal spectra}
\label{sec:global model structures}

In this section we establish the strong level and global model structures
on the category of orthogonal spectra. Many arguments are parallel
to the unstable analogs in  Section \ref{sec:global model structures spaces},
so there is a certain amount of repetition. 
The main model structure of interest for us is
the {\em global model structure},
see Theorem \ref{thm:All global spectra}.
The weak equivalences in this model structure are the global equivalences
and the cofibrations are the flat cofibrations.
More generally, we consider a global family $\Fc$ and define the
{\em  $\Fc$-global model structure},
see Theorem \ref{thm:F-global spectra} below,
with weak equivalences the $\Fc$-equivalences, 
i.e., those morphisms inducing isomorphisms
of $G$-equivariant homotopy groups for all $G$ in $\Fc$.
Proposition \ref{prop:ExF ppp} shows that the global model structure
is monoidal with respect to the smash product of orthogonal spectra;
more generally, the $\Fc$-global model structure is monoidal,
provided that $\Fc$ is closed under products.

\medskip

As we explained in more detail in Construction \ref{con:skeleton spec},
there is a functorial way to write an orthogonal spectrum $X$ as 
a sequential colimit of spectra which 
are made from the information below a fixed level, the {\em skeleta} 
\[ \sk^m X\ = \ l_m(X^{\leq m}) \ ,\]
the extension of the restriction of $X$ to $\bO_{\leq m}$.
The skeleton comes with a natural morphism $i_m:\sk^m X\to X$, the counit of
the adjunction $(l_m,(-)^{\leq m})$.
The value $i_m(V):(\sk^m X)(V)\to X(V)$ is an isomorphism 
for all inner product spaces $V$ of dimension at most $m$.
The word `filtration' should be used with caution
because the morphism $i_m$ need not be injective.\index{subject}{skeleton!of an orthogonal spectrum}
The {\em $m$-th latching space}\index{subject}{latching space!of an orthogonal spectrum}
of $X$ is the based $O(m)$-space $ L_m X  =  (\sk^{m-1} X)(\mR^m)$;
it comes with a natural based $O(m)$-equivariant map
\[  \nu_m=i_{m-1}(\mR^m)\ :\ L_m X\ \to \ X(\mR^m) \ , \]
the {\em $m$-th latching map}. 

\begin{eg}
  We let $G$ be a compact Lie group and $V$ a $G$-representation of
  dimension $n$.
  Then the semifree orthogonal spectrum \eqref{eq:define F_{G,V}}
  $F_{G,V} A$ generated by a based $G$-space $A$ in level $V$ 
  is `purely $n$-dimensional' in the following sense.
  The space $(F_{G,V}A)_m$ is trivial for $m < n$, and hence the 
  latching space $L_m (F_{G,V}A)$ is trivial for $m\leq n$.
  For $m>n$ the latching map $\nu_m:L_m(F_{G,V}A) \to (F_{G,V}A)(\mR^m)$ 
  is an isomorphism. So the skeleton $\sk^m (F_{G,V}A)$ is trivial for $m<n$
  and $\sk^m (F_{G,V}A)=F_{G,V}A$ is the entire spectrum for $m \geq n$.
\end{eg}

Proposition \ref{prop:general level model structure} 
is a fairly general construction of level model structures.
We specialize the general recipe to the category of orthogonal spectra.
For a morphism $f:X\to Y$  of orthogonal spectra and $m\geq 0$
we have a commutative square of $O(m)$-spaces:
\[\xymatrix@C=12mm{ L_m X\ar[r]^{L_m f}\ar[d]_{\nu_m^X} &  L_m Y \ar[d]^{\nu_m^Y}\\
X(\mR^m)\ar[r]_{f(\mR^m)} & Y(\mR^m)} \]
We thus get a natural morphism of based $O(m)$-spaces
\[ \nu_m f\ = \ f(\mR^m)\cup \nu_m^Y \ : \ X(\mR^m)\cup_{L_m X}L_m Y \ \to \ Y(\mR^m) \ . \]

\begin{defn}
  A morphism $f:X\to Y$ of orthogonal spectra is
  a {\em flat cofibration}\index{subject}{flat cofibration!of orthogonal spectra} 
  if the latching morphism
  $\nu_m f:X(\mR^m)\cup_{L_m X}L_m Y\to Y(\mR^m)$ is an $O(m)$-cofibration  
  for all $m\geq 0$.
  An orthogonal spectrum $Y$ is {\em flat}\index{subject}{flat!orthogonal spectrum}\index{subject}{orthogonal spectrum!flat} if the morphism from the trivial spectrum to it
  is a flat cofibration, i.e., for every $m\geq 0$ the latching map $\nu_m:L_m Y \to Y(\mR^m)$
  is an $O(m)$-cofibration.
\end{defn}

Flatness as just defined is thus a special case of `$G$-flatness'
in the sense of Definition \ref{def:G-flat}, for $G$ a trivial compact Lie group.

\medskip

We let $\Fc$ be a global family\index{subject}{global family} 
in the sense of Definition \ref{def:global family},
i.e., a non-empty class of compact Lie groups
that is closed under isomorphisms, closed subgroups and quotient groups.
As in the unstable situation 
in Section \ref{sec:global families unstable},
we now develop the {\em  $\Fc$-level model structure}
on the category of orthogonal spectra,
in which the $\Fc$-level equivalences are the weak equivalences.
This model structure has a `global' (or `stable') version, 
see Theorem \ref{thm:F-global spectra} below. 

We recall that $\Fc\cap G$ denotes the family of those closed subgroups
of a compact Lie group $G$ that belong to the global family $\Fc$.
Moreover, $\Fc(m)=\Fc\cap O(m)$ is the family of those closed subgroups
of the orthogonal group $O(m)$ that belong to the global family $\Fc$.

\begin{defn} Let $\Fc$ be a global family.
  A morphism $f:X\to Y$ of orthogonal spectra is
  \begin{itemize}
  \item  an {\em $\Fc$-level equivalence}\index{subject}{F-level equivalence@$\Fc$-level equivalence!of orthogonal spectra} 
    if the map $f(\mR^m):X(\mR^m)\to Y(\mR^m)$ is an $\Fc(m)$-equivalence for all $m\geq 0$;
  \item an {\em $\Fc$-level fibration}\index{subject}{F-level fibration@$\Fc$-level fibration!of orthogonal spectra} 
    if the map $f(\mR^m):X(\mR^m)\to Y(\mR^m)$ is an $\Fc(m)$-fibration for all $m\geq 0$;
  \item an $\Fc$-cofibration if the latching morphism
    $\nu_m f:X(\mR^m)\cup_{L_m X}L_m Y\to Y(\mR^m)$ is an $\Fc(m)$-cofibration 
    or all $m\geq 0$.
\end{itemize}
\end{defn}

In other words, $f:X\to Y$ is an $\Fc$-level equivalence 
(respectively $\Fc$-level fibration) if for every $m\geq 0$ 
and every subgroup $H$ of $O(m)$ that belongs to the family $\Fc$ 
the map $f(\mR^m)^H:X(\mR^m)^H\to Y(\mR^m)^H$ is a weak equivalence
(respectively Serre fibration).

We let $G$ be any group from the family $\Fc$ and $V$ a faithful $G$-representation
of dimension $m$.
We let $\alpha:\mR^m\to V$ be a linear isometry and
define a homomorphism $c_\alpha:G\to O(m)$ by `conjugation by $\alpha$',
i.e., we set 
\[  (c_\alpha(g))(x) \ = \ \alpha^{-1}(g\cdot\alpha(x)) \]
for $g\in G$ and $x\in\mR^m$.
We restrict the $O(m)$-action on $X(\mR^m)$ to a $G$-action
along the homomorphism $c_\alpha$. Then the map
\[  X(\alpha)\ : \ c_\alpha^*( X(\mR^m) )\ \to \ X(V)  \]
is a $G$-equivariant homeomorphism, natural in $X$;
it restricts to a natural homeomorphism from $X(\mR^m)^{\bar G}$ to $X(V)^G$,
where $\bar G=c_\alpha(G)$ is the image of $c_\alpha$. This implies:

\begin{prop}\label{prop:Fc level equivalences} 
Let $\Fc$ be a global family and $f:X\to Y$ a morphism of orthogonal spectra.
\begin{enumerate}[\em (i)]
\item The morphism $f$ is an $\Fc$-level equivalence
  if and only if for every compact Lie group $G$ 
  and every faithful $G$-representation $V$ the map
  $f(V):X(V)\to Y(V)$ is an $(\Fc\cap G)$-equivalence.
  \item The morphism $f$ is an $\Fc$-level fibration
    if and only if for every compact Lie group $G$ 
    and every faithful $G$-representation $V$ the map
    $f(V):X(V)\to Y(V)$ is an $(\Fc\cap G)$-fibration.
\end{enumerate}
\end{prop}

Now we are ready to establish the $\Fc$-level model structure.

\begin{prop}\label{prop:F-level spectra}
  Let $\Fc$ be a global family.
     The $\Fc$-level equivalences, $\Fc$-level fibrations 
    and $\Fc$-cofibrations form a model structure,
    the {\em  $\Fc$-level model structure},\index{subject}{F-level model structure@$\Fc$-level model structure!for orthogonal spectra}
    on the category of orthogonal spectra.
   The  $\Fc$-level model structure is topological and cofibrantly generated.
\end{prop}
\begin{proof}
The first part is a special case of 
Proposition \ref{prop:general level model structure}, in the following way.
We let $\Cc(m)$ be the $\Fc(m)$-projective
model structure on the category of based $O(m)$-spaces,
i.e., the based version of the model structure of 
Proposition~ \ref{prop:proj model structures for G-spaces}.
With respect to these choices of model structures $\Cc(m)$,
the classes of level equivalences, level fibrations and cofibrations
in the sense of Proposition \ref{prop:general level model structure} 
become the $\Fc$-level equivalences, $\Fc$-level fibrations and $\Fc$-cofibrations.

The consistency condition (see Definition \ref{def:consistency condition})
becomes the following condition:
for all $m,n\geq 0$ and every acyclic cofibration $i:A\to B$ in 
the $\Fc(m)$-projective model structure on based $O(m)$-spaces,
every cobase change, in the category of based $O(m+n)$-spaces, of the map 
\[  \bO(\mR^m,\mR^{m+n})\sm_{O(m)} i \ : \
 \bO(\mR^m,\mR^{m+n})\sm_{O(m)} A \ \to\ \bO(\mR^m,\mR^{m+n})\sm_{O(m)}B  \]
is an $\Fc(m+n)$-weak equivalence.
We show a stronger statement, namely that the functor
\[ \bO(\mR^m,\mR^{m+n})\sm_{O(m)} - \ : \ O(m)\bT_* \ \to\ O(m+n)\bT_* \]
takes acyclic cofibrations in 
the $\Fc(m)$-projective model structure 
to acyclic cofibrations in the projective model structure
on the category of based $O(m+n)$-spaces
(i.e., the $\All$-projective model structure, where $\All$ is
the family of all closed subgroups of $O(m+n)$). 
Since the functor under consideration
is a left adjoint, it suffices to prove the claim for the generating
acyclic cofibrations, i.e., the maps
\[   ( O(m)/H \times j_k)_+\]
for all $H\in\Fc(m)$ and all $k\geq 0$,
where $j_k:D^k\times\{0\}\to D^k\times [0,1]$ is the inclusion.
Up to isomorphism, the functor takes this map to
\[  O(m+n)\ltimes_{H\times O(n)} S^n \sm (j_k)_+  \]
where $H$ acts trivially on $S^n$.
Since the projective model structure on the category of based $O(m+n)$-spaces
is topological, it suffices to show that 
$O(m+n)\ltimes_{H\times O(n)} S^n$ is cofibrant in this model structure.
Since $S^n$ is $O(n)$-equivariantly homeomorphic to the reduced mapping cone
of the map $O(n)/O(n-1)_+\to S^0$, it suffices to show that
the two $O(m+n)$-spaces
\[ O(m+n)\times_{H\times O(n)}(O(n)/O(n-1)) \text{\qquad and\qquad}
 O(m+n)/(H\times O(n))  \]
are cofibrant in the unbased sense.
Both of these $O(m+n)$-spaces are homogeneous spaces, hence $O(m+n)$-cofibrant.

We describe explicit sets of generating cofibrations 
and generating acyclic cofibrations for the  $\Fc$-level model structure.
As before we denote by
\[ G_m \ = \ F_{O(m),\mR^m} \ : \ O(m)\bT_* \ \to \ \spec \]
 the left adjoint to the evaluation functor $X\mapsto X(\mR^m)$,
i.e., the semifree functor \eqref{eq:define F_{G,V}}
indexed by the tautological $O(m)$-representation.
We let $I_{\Fc}$ be the set of all morphisms $G_mi$ for $m\geq 0$ and for $i$ in the set
of generating cofibrations for the $\Fc(m)$-projective model
structure on the category of $O(m)$-spaces specified 
in \eqref{eq:I_for_F-proj_on_GT}.
Then the set $I_{\Fc}$ detects the acyclic fibrations 
in the  $\Fc$-level model structure 
by Proposition \ref{prop:general level model structure}~(iii). 
Similarly, we let $J_{\Fc}$ be the set of all morphisms $G_m j$ 
for $m\geq 0$ and for $j$ in the set
of generating acyclic  cofibrations for the $\Fc(m)$-projective model
structure on the category of $O(m)$-spaces specified 
in \eqref{eq:J_for_F-proj_on_GT}.
Again by Proposition \ref{prop:general level model structure}~(iii),
$J_{\Fc}$ detects the fibrations in the  $\Fc$-level model structure. 
The $\Fc$-level model structure is topological 
by Proposition \ref{prop:topological criterion},
where we take $\Gc$ as the set of semifree orthogonal spectra $F_{H,\mR^m}$
for all $m\geq 0$ and all $H\in\Fc(m)$, and $\Zc=\emptyset$
as the empty set.
\end{proof}

When $\Fc=\td{e}$ is the minimal global family consisting of all trivial groups,
then the $\td{e}$-level equivalences (respectively $\td{e}$-level fibrations) 
are the level equivalences (respectively level fibrations) of
orthogonal spectra in the sense of Definition~6.1 of \cite{mmss}.
Hence the $\td{e}$-cofibrations are the `q-cofibrations' 
in the sense of \cite[Def.\,\,6.1]{mmss}.
For the minimal global family, the  $\td{e}$-level model structure 
thus specializes to the level model structure of \cite[Thm.\,6.5]{mmss}.

For easier reference we spell out the case $\Fc=\All$ of the maximal global family 
of all compact Lie groups. In this case $\All(m)$ is the family of all closed
subgroups of $O(m)$, and the $\All$-cofibrations specialize to the flat cofibrations. 
We introduce special names for the $\All$-level equivalences
and the $\All$-level fibrations,
analogous to the unstable situation 
in Section \ref{sec:global model structures spaces}.

\begin{defn}
A morphism $f:X\to Y$ of orthogonal spectra is
a {\em strong level equivalence}\index{subject}{strong level equivalence!of orthogonal spectra}\index{subject}{level equivalence!strong|see{strong level equivalence}} 
(respectively {\em strong level fibration})\index{subject}{strong level fibration!of orthogonal spectra} 
if for every $m\geq 0$ the map $f(\mR^m):X(\mR^m)\to Y(\mR^m)$ 
is an $O(m)$-weak equivalence (respectively an $O(m)$-fibration).
\end{defn}

For the global family $\Fc=\All$, Proposition \ref{prop:F-level spectra} 
specializes to:

\begin{prop}
    The strong level equivalences, strong level fibrations 
    and flat cofibrations form a model structure,
    the {\em strong level model structure},\index{subject}{strong level model structure!for orthogonal spectra}
    on the category of orthogonal spectra.
    The strong level model structure is topological and cofibrantly generated. 
\end{prop}

Now we introduce and discuss an important class of orthogonal spectra.

\begin{defn}\label{def:global Omega} 
An orthogonal spectrum $X$ is a 
{\em global $\Omega$-spectrum}\index{subject}{global Omega-spectrum@global $\Omega$-spectrum}
if for every compact Lie group $G$, every faithful $G$-representation $W$
and an arbitrary $G$-representation $V$ the adjoint structure map 
\[ \tilde\sigma_{V,W}\ :\  X(W)\ \to \ \map_*(S^V,X(V\oplus W))\]
is a $G$-weak equivalence.
\end{defn}

The global $\Omega$-spectra will turn out to be
the fibrant objects in the global model structure on orthogonal spectra, see
Theorem \ref{thm:All global spectra} below.
This means that global $\Omega$-spectra abound, because every
orthogonal spectrum admits a global equivalence to a global $\Omega$-spectrum.

\begin{rk}
Global $\Omega$-spectra are a very rich kind of structure,
because they encode compatible equivariant
infinite loop spaces for all compact Lie groups at once. 
For a global $\Omega$-spectrum $X$ and a compact Lie group $G$ 
the associated orthogonal $G$-spectrum  $X_G$ is 
`eventually an $\Omega$-$G$-spectrum' in the sense that
the $\Omega$-$G$-spectrum condition of \cite[III Def.\,3.1]{mandell-may}
holds for all `sufficiently large' (i.e., faithful) $G$-representations.
However, if $G$ is a non-trivial group,
then the associated orthogonal $G$-spectrum $X_G$ is in general {\em not}
an $\Omega$-$G$-spectrum since there is no control over
the $G$-homotopy type of the values at non-faithful representations.  

For every compact Lie group $G$ and every faithful $G$-representation $W$, 
the $G$-space 
\[ X[G] \ =\ \Omega^W X(W)    \]
is a `genuine' equivariant infinite loop space, i.e., deloopable in
the direction of every representation.
Indeed, for every $G$-representation $V$, the composite
\[ 
X[G]\ = \ \Omega^W X(W)\ \xra{\Omega^W(\tilde\sigma_{V,W})} \ 
\Omega^W (\Omega^V X(V\oplus W)) \ \iso \
\Omega^V (\Omega^W X(V\oplus W)) \]
is a $G$-weak equivalence, so the global $\Omega$-spectrum $X$ provides 
a $V$-deloop $\Omega^W X(V\oplus W)$ of $X[G]$.
The $G$-space $X[G]$ is also independent, up to $G$-weak equivalence,
of the choice of faithful $G$-representation. Indeed, if $\bar W$ is another
faithful $G$-representation, then the $G$-maps
\begin{align*}
   \Omega^W X(W)\ \xra{\Omega^W (\tilde\sigma_{W',W})} &\
\Omega^W(\Omega^{W'} X(W'\oplus W)) \\ 
&\iso \ \Omega^{W'}(\Omega^W,X(W\oplus W')) \ \xla{\Omega^{W'}(\tilde\sigma_{W,W'})}
\Omega^{W'} X(W') 
\end{align*}
are $G$-weak equivalences.

As $G$ varies, the equivariant infinite loop spaces $X[G]$ are 
closely related to each other. For example, if $H$ is a closed subgroup of $G$, 
then any faithful $G$-representation is also faithful as an $H$-representation.
So $X[H]$ is $H$-weakly equivalent to the restriction of the 
$G$-equivariant infinite loop space $X[G]$.
\end{rk}

\begin{rk} 
Let $X$ be a global $\Omega$-spectrum.
Specialized to the trivial group, the condition in Definition \ref{def:global Omega} 
says that $X$ is in particular a non-equivariant $\Omega$-spectrum 
in the sense that the adjoint structure map 
$\tilde\sigma_{\mR,W}:X(W)\to \Omega X(\mR\oplus W)$ is a weak equivalence of 
(non-equivariant) spaces for every inner product space $W$.

If $X$ is a global $\Omega$-spectrum, then so is the shifted
spectrum $\sh X$ and the function spectrum $\map_*(A,X)$
for every cofibrant based space $A$.
Indeed, if $W$ is a faithful $G$-representation, then 
$W\oplus\mR$, the sum with a trivial 1-dimensional representation,
is also faithful. So the adjoint structure map
\[\tilde\sigma_{V,W\oplus\mR}^X \ : \ X(W\oplus\mR) \ \to \ 
\Omega^V(X(V\oplus W\oplus\mR))  \]
is a $G$-weak equivalence for every $G$-representation $V$.
But this map is also the adjoint structure map
\[  \tilde\sigma_{V,W}^{\sh X}\ : \ (\sh X)(W) \ \to \ \Omega^V ((\sh X)(V\oplus W)) \]
of the shifted spectrum.
The argument for mapping spectra is similar.
\end{rk}

\begin{defn}\index{subject}{global fibration!of orthogonal spectra}\label{def:global fibration}
    A morphism $f:X\to Y$ of orthogonal spectra is a {\em global fibration}
    if it is a strong level fibration
    and for every compact Lie group $G$,
    every $G$-representation $V$ and every faithful $G$-representation $W$ 
    the square 
    \begin{equation}  \begin{aligned}
        \label{eq:global fibration spectra}
        \xymatrix@C=13mm{ X(W)^G \ar[d]_{\quad f(W)^G} \ar[r]^-{(\tilde\sigma_{V,W})^G} & 
          \map_*^G(S^V, X(V\oplus W)) \ar[d]^{\map_*^G(S^V,f(V\oplus W))} \\
          Y(W)^G \ar[r]_-{(\tilde\sigma_{V,W})^G} & \map_*^G(S^V, Y(V\oplus W))}
      \end{aligned}\end{equation}
    is homotopy cartesian.  
\end{defn}

We state a useful criterion for checking when a morphism is a global fibration.
We recall from Construction \ref{con:cone and fiber}
the homotopy fiber $F(f)$ of a morphism $f:X\to Y$ of orthogonal spectra 
along with the natural map $i:\Omega Y\to F(f)$.
We apply this for the shifted morphism
$\sh f:\sh X\to \sh Y$\index{subject}{shift!of an orthogonal spectrum}
and precompose with the morphism $\tilde\lambda_Y$,
the adjoint to $\lambda_Y:Y\sm S^1\to\sh Y$ defined in \eqref{eq:defn lambda_n}.
We denote by $\xi(f)$ the resulting composite
\[ Y \ \xra{\ \tilde\lambda_Y\ } \ \Omega ( \sh Y)\  \xra{\ i\ }\  F(\sh f)  \ . \]
We note that the orthogonal spectrum $F(\Id_{\sh X})$ has a preferred
contraction, so part~(ii) of the next proposition is a way to make precise that
the sequence 
\[ X \ \xra{\ f\ }\ Y \ \xra{\ \xi(f)\ }\ F(\sh f) \]
is a `strong level homotopy fiber sequence'.

\begin{prop}\label{prop:global fibration criterion}
Let $f:X\to Y$ be a strong level fibration of orthogonal spectra. 
Then the following two conditions are equivalent.
\begin{enumerate}[\em (i)]
\item The morphism $f$ is a global fibration.
\item
The strict fiber of $f$ is a global $\Omega$-spectrum and
the commutative square
\[ \xymatrix@C=15mm{ 
X \ar[d]_f \ar[r]^-{\xi(\Id_X)} & F(\Id_{\sh X}) \ar[d] \\
Y \ar[r]_-{\xi(f)} &  F(\sh f)  } \]
is homotopy cartesian in the strong level model structure.
\end{enumerate}
\end{prop}
\begin{proof}
(i)$\Longrightarrow$(ii)
Since the square \eqref{eq:global fibration spectra} 
is homotopy cartesian and the vertical maps are Serre fibrations,
the induced map on strict fibers $(f^{-1}(*)(W))^G\to \map_*^G(S^V,f^{-1}(*)(V\oplus W))$
is a weak equivalence. So the strict fiber is a global $\Omega$-spectrum.

The square in (ii) factors as the composite of two commutative squares:
\[
\xymatrix{ X \ar[d]_f \ar[r]^-{\tilde\lambda_X} & 
  \Omega \sh X  \ar[d]^{\Omega \sh f} \ar[r]^-i & F(\Id_{\sh X}) \ar[d] \\
  Y \ar[r]_-{\tilde\lambda_Y} & \Omega\sh Y \ar[r]_i & F(\sh f ) } \]
Specializing the global fibration property \eqref{eq:global fibration spectra}
for $W=\mR$ with trivial $G$-action shows that the left square is 
homotopy cartesian in the strong level model structure.
For every continuous based map $g:A\to B$ the square
\[ \xymatrix{ 
\Omega A \ar[r]^-i \ar[d]_{\Omega g} & F(\Id_A)\ar[d] \\
\Omega B \ar[r]_-i & F(g)} \]
is homotopy cartesian; applying this to  $g=((\sh f)(V))^G$
for all $G$-representa\-tions $V$ shows that 
the right square is homotopy cartesian  in the strong level model structure.
So the composite square is homotopy cartesian in the strong level model structure.

(ii)$\Longrightarrow$(i)
We let $G$ be a compact Lie group, $V$ a $G$-representation,
and $W$ a faithful $G$-representation. We consider the commutative diagram
\begin{equation*}
    \xymatrix@C=10mm{ X(W)^G \ar[d]_{f(W)^G} \ar[r]^-{(\tilde\sigma_{V,W})^G} & 
      \map_*^G(S^V, X(V\oplus W)) \ar[d]^{\map_*^G(S^V,f(V\oplus W))} \\
      Y(W)^G \ar[d]_{(\xi(f)(W))^G}\ar[r]_-{(\tilde\sigma_{V,W})^G} & 
      \map_*^G(S^V, Y(V\oplus W)) \ar[d]^{\map_*^G(S^V,\xi(f)(V\oplus W))}\\
          F(f)(W\oplus\mR)^G \ar[r]^-\sim_-{(\tilde\sigma_{V,W})^G} & \map_*^G(S^V, F(f)(V\oplus W\oplus \mR))}
      \end{equation*}
Since $f$ is a strong level fibration, the natural morphism from
the strict fiber to the homotopy fiber is a strong level equivalence.
So the lower horizontal map is a weak equivalence because the strict fiber,
and hence the homotopy fiber, is a global $\Omega$-spectrum.
Moreover, the two vertical columns are homotopy fiber sequences by hypothesis,
so the upper square is homotopy cartesian, both with respect to the strong
level model structure.
\end{proof}

\begin{defn}
  Let $\Fc$ be a global family. 
  \begin{itemize}
  \item 
    A morphism $f:X\to Y$
    of orthogonal spectra is an {\em $\Fc$-equivalence}\index{subject}{F-equivalence@$\Fc$-equivalence!of orthogonal spectra}
    if the induced map $\pi_k^G (f) : \pi_k^G (X) \to \pi_k^G (Y)$ is an isomorphism
    for all $G$ in $\Fc$ and all integers $k$.
  \item
    A morphism $f:X\to Y$ of orthogonal spectra is an $\Fc$-global fibration
    if it is an $\Fc$-level fibration
    and for every compact Lie group $G$ in $\Fc$,
    every $G$-representation $V$\index{subject}{F-global fibration@$\Fc$-global fibration!of orthogonal spectra} 
    and every faithful $G$-representation $W$ the square \eqref{eq:global fibration spectra} 
        is homotopy cartesian.  
      \item 
        An orthogonal spectrum $X$ is an 
        {\em $\Fc$-$\Omega$-spectrum}\index{subject}{F-Omega-spectrum@$\Fc$-$\Omega$-spectrum}
if for every compact Lie group $G$ in $\Fc$, every $G$-representation $V$
and every faithful $G$-representation $W$ the adjoint structure map 
\[ \tilde\sigma_{V,W}\ :\  X(W)\ \to \ \map_*(S^V,X(V\oplus W))\]
is a $G$-weak equivalence.
\end{itemize}
\end{defn}

When $\Fc=\All$ is the maximal global family of all compact Lie groups,
then an $\All$-equivalence is just a global equivalence
in the sense of Definition \ref{def:global equivalence},
and the $\All$-global fibrations are the global fibrations in the
sense of Definition \ref{def:global fibration}.
Also, an $\All$-$\Omega$-spectrum is the same as a global $\Omega$-spectrum 
in the sense of Definition \ref{def:global Omega}.
When $\Fc=\td{e}$ is the minimal global family of all trivial groups,
then the $\td{e}$-equivalences are just the traditional
non-equivariant stable equivalences of orthogonal spectra,
also known as `$\pi_*$-isomorphisms'.
The $\td{e}$-$\Omega$-spectra are just the traditional non-equivariant $\Omega$-spectra.

The following proposition collects various useful relations
between $\Fc$-equi\-valences, $\Fc$-level equivalences and $\Fc$-$\Omega$-spectra.

\begin{prop}\label{prop:stable zero is level zero}
  Let $\Fc$ be a global family.
  \begin{enumerate}[\em (i)]
  \item 
    Every $\Fc$-level equivalence of orthogonal spectra is an $\Fc$-equivalence.
  \item
    Let $X$ be an $\Fc$-$\Omega$-spectrum. Then for every $G$ in $\Fc$, 
    every faithful $G$-representation $V$ and every $k\geq 0$ the stabilization map
    \[ [S^{V\oplus\mR^k}, X(V)]^G \ \to \ \pi_k^G(X) \ , \quad [f]\ \longmapsto \ \td{f}\]
    is bijective.
  \item 
    Let $X$ be an $\Fc$-$\Omega$-spectrum such that $\pi^G_k (X)=0$
    for every integer $k$ and all $G$ in $\Fc$. Then for every group $G$ 
    in $\Fc$ and every faithful $G$-representation $V$ the space
    $X(V)$ is $G$-weakly contractible.
  \item  
    Every $\Fc$-equivalence between $\Fc$-$\Omega$-spectra
    is an $\Fc$-level equivalence.
  \item  
    Every $\Fc$-equivalence
    that is also an $\Fc$-global fibration is an $\Fc$-level equivalence.
  \end{enumerate}
\end{prop}
\begin{proof}
(i) We let $f:X\to Y$ be an $\Fc$-level equivalence, and we need to show that the map
$\pi_k^G(f): \pi_k^G (X)\to\pi_k^G(Y)$ is an isomorphism for all integers $k$ and
all $G$ in $\Fc$. We start with the case $k=0$.
We let $G$ be a group from $\Fc$ and $V$ a finite-dimensional
faithful $G$-subrepresentation of the complete $G$-universe $\Uc_G$.
By Proposition \ref{prop:Fc level equivalences}~(i) 
the map $f(V):X(V)\to Y(V)$ is a $G$-weak equivalence.
Since the representation sphere $S^V$
can be given a $G$-CW-structure, the induced map on $G$-homotopy classes
\[  [S^V,f(V)]^G \ : \ [S^V,X(V)]^G\ \to\  [S^V,Y(V)]^G \] 
is bijective. The faithful $G$-representations are cofinal in
the poset $s(\Uc_G)$, so taking the colimit over $V\in s(\Uc_G)$ shows that 
$\pi_0^G(f):\pi_0^G(X)\to\pi_0^G (Y)$ is an isomorphism.
For $k>0$ we exploit that $\pi_k^G(X)$ is naturally isomorphic
to $\pi_0^G(\Omega^k X)$ and the $k$-fold loop
of an $\Fc$-level equivalence is again an $\Fc$-level equivalence.
For $k<0$ we exploit that $\pi_k^G(X)$ is naturally
isomorphic to $\pi_0^G(\sh^{-k} X)$ and that every shift
of an $\Fc$-level equivalence is again an $\Fc$-level equivalence.

(ii) We may assume that $V$ is a $G$-invariant subspace
of the complete $G$-universe $\Uc_G$ used to define $\pi_k^G(X)$.
Since $V$ is faithful, so is every $G$-represen\-tation that contains $V$.
So the directed system whose colimit is  $\pi_k^G(X)$
consists of isomorphisms `above $V$'. Hence the canonical map
\[ [S^{V\oplus\mR^k}, X(V)]^G \ \to \ \pi_k^G(X) \]
is bijective.

(iii) We adapt an argument of Lewis-May-Steinberger \cite[I 7.12]{lms} to our context.
Every $\Fc$-$\Omega$-spectrum is in particular a non-equivariant $\Omega$-spectrum;
every non-equivariant $\Omega$-spectrum with trivial homotopy groups
is levelwise weakly contractible, so this takes care of trivial groups.

Now we let $G$ be a non-trivial group in $\Fc$.
We argue by a nested induction over the `size' of $G$:
we induct over the dimension of $G$ and, for fixed dimension, 
over the cardinality of the finite set of path components of $G$.
Every proper closed subgroup $H$ of $G$ either has strictly smaller dimension than $G$,
or the same dimension but fewer path components.
So we know by induction that the fixed point space $X(V)^H$ is weakly contractible 
for every proper closed subgroup $H$ of $G$.
So it remains to analyze the $G$-fixed points of $X(V)$.

We let $W=V-V^G$ be the orthogonal complement of the fixed subspace $V^G$
of $V$. Then $W$ is a faithful $G$-representation with trivial fixed points.
The cofiber sequence of $G$-CW-complexes
\[ S(W)_+ \ \to \ D(W)_+ \ \to \ S^W\]
induces a fiber sequence of equivariant mapping spaces
\[ \map_*^G(S^W,X(V)) \ \to \ \map^G(D(W),X(V)) \ \to \
\map^G(S(W),X(V)) \ . \]
Since $W^G=0$, the $G$-fixed points $S(W)^G$ are empty,
so any $G$-CW-structure on $S(W)$ uses only equivariant cells of the
form $G/H\times D^n$ for proper subgroups $H$ of $G$.
But for proper subgroups $H$, all the fixed point spaces $X(V)^H$ 
are weakly contractible by induction. 
Hence the space $\map^G(S(W),X(V))$ is weakly contractible.
Since the unit disc $D(W)$ is equivariantly contractible,
the space $\map^G(D(W),X(V))$ is homotopy equivalent to $X(V)^G$,
and we conclude that evaluation at the fixed point $0\in W$ is a weak equivalence
\[ \map^G_*(S^W,X(V)) \ \simeq \ X(V)^G\ . \]
We have $X(V)=X(W\oplus V^G)=(\sh^{V^G}X)(W)$.
The stabilization map
\begin{align*}
  \pi_k\left( \map^G_*(S^W,X(V))\right)\ &=\ [S^{W\oplus\mR^k}, (\sh^{V^G}X)(W)]^G\\
 &\to \ \pi_k^G(\sh^{V^G}X) \iso \pi_{k-\dim(V^G)}^G (X)  
\end{align*}
is bijective for all $k\geq 0$, by part~(ii). Since we assumed that the
$G$-equivariant homotopy groups of $X$ vanish, the space $\map_*^G(S^W, X(V))$
is path connected and has vanishing homotopy groups, i.e., 
it is weakly contractible.
So $X(V)^G$ is weakly contractible, and this completes the proof of~(iii).

(iv) Let $f:X\to Y$ be an $\Fc$-equivalence 
between $\Fc$-$\Omega$-spectra.
We let $F$ denote the homotopy fiber of $f$, and we let $G$ be a group from $\Fc$. 
For every $G$-representation $V$ the $G$-space $F(V)$ is the homotopy fiber
of $f(V):X(V)\to Y(V)$. So $F$ is again an $\Fc$-$\Omega$-spectrum.
The long exact sequence of homotopy groups 
(see Proposition \ref{prop:LES for homotopy of cone and fibre})
implies that $\pi^H_*(F)=0$ 
for all $H$ in $\Fc$.

If $G$ acts faithfully on $V$, then by the $\Fc$-$\Omega$-spectrum property, 
the space $X(V)$ is $G$-weakly equivalent
to $\Omega X(V\oplus\mR)$ and similarly for $Y$.
So the map $f(V)$ is $G$-weakly equivalent to 
\[ \Omega f(V\oplus\mR) \ : \ \Omega X(V\oplus\mR)
\ \to \ \Omega X(V\oplus\mR)\ . \]
Hence we have a homotopy fiber sequence of $G$-spaces
\[X(V)\ \xra{\ f(V)\ }\ Y(V) \ \to \ F(V\oplus\mR) \ . \]
Since $F$ is an $\Fc$-$\Omega$-spectrum with vanishing equivariant
homotopy groups for groups in $\Fc$, the space $F(V\oplus\mR)$ 
is $G$-weakly contractible by part (iii). 
So $f(V)$ is a $G$-weak equivalence. 
The morphism $f$ is then an $\Fc$-level equivalence 
by the criterion of Proposition \ref{prop:Fc level equivalences} (i). 

(v) We let $f:X\to Y$ be an $\Fc$-equivalence and an $\Fc$-global fibration.
Then the strict fiber $f^{-1}(\ast)$ of $f$ is an $\Fc$-$\Omega$-spectrum
with trivial $G$-equivariant homotopy groups for all $G\in\Fc$.
So $f^{-1}(\ast)$ is $\Fc$-level equivalent to the trivial spectrum by
part~(iii). Since $f$ is an $\Fc$-level fibration, 
the embedding $f^{-1}(\ast)\to F(f)$ of the strict fiber into the homotopy fiber 
is an $\Fc$-level equivalence, so the  homotopy fiber $F(f)$
is $\Fc$-level equivalent to the trivial spectrum.
The homotopy cartesian square of Proposition \ref{prop:global fibration criterion}~(ii)
then shows that $f$ is an $\Fc$-level equivalence.
\end{proof}

Our next result, in fact the main result of this section,
is the $\Fc$-global model structure.
This model structure is cofibrantly generated, and we spell out
an explicit set of generating cofibrations and generating acyclic cofibrations.
The set $I_\Fc$ was defined in 
the proof of Proposition \ref{prop:F-level spectra} as the set of morphisms $G_m i$ 
for $m\geq 0$ and for $i$ in the set of generating cofibrations 
for the $\Fc(m)$-projective model structure on the category of $O(m)$-spaces specified 
in \eqref{eq:I_for_F-proj_on_GT}.
Similarly, the set $J_\Fc$ is the set of morphisms $G_m j$ 
for $m\geq 0$ and for $j$ in the set of generating acyclic cofibrations 
for the $\Fc(m)$-projective model structure on the category of $O(m)$-spaces specified 
in \eqref{eq:J_for_F-proj_on_GT}.

Given any compact Lie group $G$ and $G$-representations $V$ and $W$
we recall from \eqref{eq:define_lambda} the morphism
\[ \lambda_{G,V,W} \ : \ F_{G,V\oplus W}S^V \ \to \ F_{G,W} \ . \]
If the representation $W$ is faithful, 
then this morphism is a global equivalence 
by Theorem \ref{thm:faithful independence}. 
We set
\begin{equation}\label{eq:define K_F}
 K_\Fc \ = \ \bigcup_{G,V,W} \Zc(\lambda_{G,V,W}) \ ,  
\end{equation}
the set of all pushout products of sphere inclusions $i_m:\partial D^m\to D^m$
with the mapping cylinder inclusions of the morphisms $\lambda_{G,V,W}$
(compare Construction \ref{con:define Z(j)});
here the union is over a set of representatives
of the isomorphism classes of triples $(G,V,W)$ consisting of
a compact Lie group $G$ in the family $\Fc$, a $G$-representation $V$ 
and a faithful $G$-representation $W$.

\begin{theorem}[$\Fc$-global model structure]\label{thm:F-global spectra} 
  Let $\Fc$ be a global family.
  \begin{enumerate}[\em (i)]
  \item 
    The $\Fc$-equivalences, $\Fc$-global fibrations and $\Fc$-cofibrations 
    form a model structure on the category of orthogonal spectra, 
    the {\em  $\Fc$-global model structure}.\index{subject}{F-global model structure@$\Fc$-global model structure!for orthogonal spectra}\index{subject}{global model structure!$\Fc$-}
  \item
    The fibrant objects in the  $\Fc$-global model structure 
    are the $\Fc$-$\Omega$-spectra.\index{subject}{F-Omega-spectrum@$\Fc$-$\Omega$-spectrum}
  \item A morphism of orthogonal spectra is:
    \begin{itemize}
    \item an acyclic fibration in the  $\Fc$-global model structure
      if and only if it has the right lifting property 
      with respect to the set $I_{\Fc}$; 
    \item 
      a fibration in the  $\Fc$-global model structure
      if and only if it has the right lifting property 
      with respect to the set $J_\Fc \cup K_\Fc$.  
    \end{itemize}
  \item
    The $\Fc$-global model structure is cofibrantly generated, 
    proper and topological. 
  \item The adjoint functor pair
    \[ \xymatrix{ \Sigma^\infty_+ \ : \ \spc\ \ar@<.4ex>[r] & 
      \ \spec\ : \ \Omega^\bullet \ar@<.4ex>[l]  } \]
    is a Quillen pair for the two $\Fc$-global model structures
    on orthogonal spaces and orthogonal spectra.
  \end{enumerate}
\end{theorem}
\begin{proof}
The category of orthogonal spectra is complete and cocomplete 
(MC1), the $\Fc$-equivalences satisfy the 2-out-of-3 property (MC2)
and the classes of $\Fc$-equivalences, $\Fc$-global fibrations and 
$\Fc$-cofibrations are closed under retracts (MC3).
The $\Fc$-level model structure (Proposition \ref{prop:F-level spectra})
shows that every morphism of orthogonal spectra
can be factored as $f\circ i$ for an $\Fc$-cofibration $i$
followed by an $\Fc$-level equivalence $f$ that is also an $\Fc$-level fibration.
For every $G$ in $\Fc$, every $G$-representation $V$ 
and every faithful $G$-representation $W$,
both vertical maps in the commutative square \eqref{eq:global fibration spectra} 
are then weak equivalences, so the square is homotopy cartesian.
The morphism $f$ is thus an $\Fc$-global fibration and an $\Fc$-equivalence,
so this provides one of the factorizations as required by MC5.

The morphism $\lambda_{G,V,W}$ represents the map
\[ (\tilde\sigma_{V,W})^G\ :\ X(W)^G\ \to\ \map_*^G(S^V,X(V\oplus W))^G\ ; \]
by Proposition \ref{prop:hocartesian via RLP}, 
the right lifting property with respect to the union $J_\Fc\cup K_\Fc$ thus
characterizes the $\Fc$-global fibrations. We apply the small object argument 
(see for example \cite[7.12]{dwyer-spalinski} or \cite[Thm.\,2.1.14]{hovey-book})
to the  set $J_\Fc\cup K_\Fc$.
All morphisms in $J_\Fc$ are flat cofibrations and $\Fc$-level equivalences;
$\Fc$-level equivalences are $\Fc$-equivalences by Proposition \ref{prop:stable zero is level zero}~(i).
Since $F_{G,V\oplus W}S^V$ and $F_{G,W}$ are flat, the morphisms in $K_\Fc$
are also flat cofibrations, and they are $\Fc$-equivalences
because the morphisms $\lambda_{G,V,W}$ are.
The small object argument provides a functorial factorization
of every morphism $\varphi:X\to Y$ of orthogonal spectra
as a composite
\[ X \ \xra{\ i \ }\ W \ \xra{\ q \ }\ Y \]
where $i$ is a sequential composition of cobase changes of coproducts
of morphisms in $J_\Fc\cup K_\Fc$, 
and $q$ has the right lifting property with respect 
to $J_\Fc\cup K_\Fc$; in particular, the morphism $q$ 
is an $\Fc$-global fibration.
All morphisms in $J_\Fc\cup K_\Fc$ are $\Fc$-equivalences and $\Fc$-cofibrations,
hence also h-cofibrations 
(by Proposition \ref{cor-h-cofibration closures} 
applied to the strong level model structure).
By Corollary \ref{cor-closure h-cof and global equi}
(or rather its modification for $\Fc$-equivalences),
the class of h-co\-fibrations that are simultaneously $\Fc$-equivalences 
is closed under coproducts, cobase change and sequential composition.
So the morphism $i$ is an $\Fc$-cofibra\-tion and an $\Fc$-equivalence.

Now we show the lifting properties MC4. 
By Proposition \ref{prop:stable zero is level zero}~(v)
a morphism that is both an $\Fc$-equivalence and an $\Fc$-global fibration
is an $\Fc$-level equivalence, and hence an acyclic fibration
in the $\Fc$-level model structure. So every morphism that is
simultaneously an $\Fc$-equivalence  and an $\Fc$-global fibration has the
right lifting property with respect to $\Fc$-cofibrations.
Now we let $j:A\to B$ be an $\Fc$-cofibration that is also an $\Fc$-equivalence and 
we show that it has the left lifting property with respect 
to $\Fc$-global fibrations.
We factor $j=q\circ i$, 
via the small object argument for $J_\Fc\cup K_\Fc$,
where $i:A\to W$ is a $(J_\Fc\cup K_\Fc)$-cell complex
and $q:W\to B$ is an $\Fc$-global fibration.
Then $q$ is an $\Fc$-equivalence since $j$ and $i$ are,
so $q$ is an acyclic fibration in the $\Fc$-level model structure, by the above.
Since $j$ is an $\Fc$-cofibration, a lifting in
\[\xymatrix{
A \ar[r]^-i \ar[d]_j & W \ar[d]^q_(.6)\sim \\
B \ar@{=}[r] \ar@{..>}[ur] & B }\]
exists. Thus $j$ is a retract of the morphism $i$ that has the left lifting property
with respect to all $\Fc$-global fibrations.
But then $j$ itself has this lifting property.
This finishes the verification of the model category axioms.
Alongside we have also specified sets of generating flat cofibrations $I_\Fc$
and generating acyclic cofibrations $J_\Fc\cup K_\Fc$.
Sources and targets of all morphisms in these sets are small with
respect to sequential colimits of flat cofibrations. So the 
 $\Fc$-model structure is cofibrantly generated.

Left properness of the $\Fc$-global model structure follows from
Corollary \ref{cor-cobase change h-cofibration}~(i) and the fact that
$\Fc$-cofibrations are h-cofibrations.
Right properness follows from
Corollary \ref{cor-cobase change h-cofibration}~(ii) and the fact that
$\Fc$-global fibrations are in particular $\Fc$-level fibrations.

The $\Fc$-global model structure is topological by 
Proposition \ref{prop:topological criterion};
here we take $\Gc$ to be the set of semifree orthogonal spectra $G_m(O(m)/H)$
for all $m\geq 0$ and all $H\in \Fc(m)$,
and we take $\Zc$ as the set of mapping cylinder inclusions $c(\lambda_{G,V,W})$
of the morphisms $\lambda_{G,V,W}$ indexed as in the set $K_\Fc$.

Part~(v) is straightforward from the definitions.
\end{proof}

In the case $\Fc=\td{e}$ of the minimal global family of trivial groups,
the $\td{e}$-equivalences are the (non-equivariant) stable equivalences
of orthogonal spectra, and the  $\td{e}$-global model structure
coincides with the stable model structure established by Mandell, May, Shipley and the
author in \cite[Thm.\,9.2]{mmss}.
For easier reference we spell out the special case $\Fc=\All$
for the maximal family of all compact Lie groups,
resulting in the {\em global model structure}
on the category of orthogonal spectra.

\begin{theorem}[Global model structure]\label{thm:All global spectra} 
    The global equivalences, global fibrations and flat cofibrations 
    form a model structure,
    the {\em global model structure}\index{subject}{global model structure!for orthogonal spectra}
    on the category of orthogonal spectra. 
    The fibrant objects in the global model structure are the global $\Omega$-spectra.
    The global model structure is proper, topological and cofibrantly generated.
  \end{theorem}

The same proof as in the unstable situation 
in Corollary \ref{cor-characterize F-equivalences spaces}
applies to prove the following characterization of $\Fc$-equivalences.

\begin{cor}
Let $f:A\to B$ be a morphism of orthogonal spectra
and $\Fc$ a global family.
Then the following conditions are equivalent.
\begin{enumerate}[\em (i)]
\item The morphism $f$ is an $\Fc$-equivalence.\index{subject}{F-equivalence@$\Fc$-equivalence!of orthogonal spectra}
\item The morphism can be written as $f=w_2\circ w_1$ for an
  $\Fc$-level equivalence $w_2$ and a global equivalence $w_1$.
\item For some (hence any) $\Fc$-cofibrant approximation
  $f^c:A^c\to B^c$ in the  $\Fc$-level model structure
and every $\Fc$-$\Omega$-spectrum $Y$ the induced map
\[ [f^c,Y] \ : \ [B^c,Y] \ \to \ [A^c,Y] \]
on homotopy classes of morphisms is bijective.
\end{enumerate}
\end{cor}

\begin{rk}[Mixed global model structures]
Cole's `mixing theorem' for model structures \cite{cole-mixed}  
allows to construct many more $\Fc$-model structures 
on the category of orthogonal spectra.
We consider two global families such that $\Fc\subseteq \Ec$. 
Then every $\Ec$-equivalence is an $\Fc$-equivalence and every 
fibration in the  $\Ec$-global model structure
is a fibration in the  $\Fc$-global model structure.
By Cole's theorem \cite[Thm.\,2.1]{cole-mixed} 
the $\Fc$-equivalences and the fibrations of the
 $\Ec$-global model structure
are part of a model structure,
the {\em $\Ec$-mixed $\Fc$-global model structure}\index{subject}{mixed global model structure}\index{subject}{global model structure!mixed}
on the category of orthogonal spectra.
By \cite[Prop.\,3.2]{cole-mixed} 
the cofibrations in the $\Ec$-mixed $\Fc$-global model structure 
are precisely the retracts of all composites
$h\circ g$ in which $g$ is an $\Fc$-cofibration and $h$ is simultaneously
an $\Ec$-equivalence and an $\Ec$-cofibration.
In particular, an orthogonal spectrum is cofibrant in the
$\Ec$-mixed $\Fc$-global model structure if it is $\Ec$-cofibrant and
$\Ec$-equivalent to an $\Fc$-cofibrant 
orthogonal spectrum \cite[Cor.\,3.7]{cole-mixed}.
The $\Ec$-mixed $\Fc$-global model structure is again proper 
(Propositions~4.1 and~4.2 of \cite{cole-mixed}).

When $\Fc=\td{e}$ is the minimal family of trivial groups,
this provides infinitely many $\Ec$-mixed model structures on
the category of orthogonal spectra that are all Quillen equivalent.
In the extreme case, the $\All$-mixed $\td{e}$-model structure
is the $\mS$-model structure of Stolz \cite[Prop.\,1.3.10]{stolz-thesis}.
\end{rk}

\begin{rk}[$\Fin$-global homotopy theory via symmetric spectra]
We denote by $\Fin$\index{symbol}{$\Fin$ - {global family of finite groups}}
the global family of finite groups.\index{subject}{global family!of finite groups} 
The $\Fin$-global stable homotopy theory
has another very natural model, namely the category of 
symmetric spectra\index{subject}{symmetric spectrum}
in the sense of Hovey, Shipley and Smith \cite{HSS}.
In \cite{hausmann-master thesis, hausmann-global} 
Hausmann has established a global model structure 
on the category of symmetric spectra, 
and he showed that the forgetful functor is a right Quillen equivalence
from the category of orthogonal spectra with the $\Fin$-global model structure
to the category of symmetric spectra with the global model structure.
Symmetric spectra cannot model global homotopy types for all compact Lie groups,
basically because compact Lie groups of positive dimensions do not have
any faithful permutation representations.
\end{rk}

One of the points of model category structures is that they facilitate
the analysis of and constructions in the homotopy category (which only depends
on the class of weak equivalences). An example of this is the existence
and constructions of products and coproducts in a homotopy category.
We take the time to make this explicit for the $\Fc$-global stable homotopy category,
for which we write $\GH_\Fc$.

\begin{prop}\label{prop:sums and product in GH} 
Let $\Fc$ be a global family.
\begin{enumerate}[\em (i)]
\item The coproduct (i.e., wedge) of every family  of orthogonal spectra
is also a coproduct in the homotopy category $\GH_\Fc$. 
\item Let $\{X_i\}_{i\in I}$ be a family of orthogonal spectra such that the
canonical map
\[ \pi_k^G\big({\prod}_{i\in I}\, X_i\big)\ \to \ {\prod}_{i\in I}\, \pi_k^G(X_i) \]
is an isomorphism for every compact Lie group $G$ in $\Fc$ and every integer $k$.
Then the product $\prod_{i\in I}  X_i$ of orthogonal spectra
is also a product of $\{X_i\}_{i\in I}$ 
in the homotopy category $\GH_\Fc$. 
\end{enumerate}
In particular, the $\Fc$-global stable homotopy category $\GH_\Fc$ 
has all set indexed coproducts and products.
\end{prop}
\begin{proof}
(i) Since the $\Fc$-equivalences are part of the $\Fc$-global model structure,
general model category theory guarantees that coproducts in $\GH_\Fc$ can be constructed
by taking the pointset level coproduct of $\Fc$-cofibrant approximations,
see for example \cite[Ex.\,1.3.11]{hovey-book}.
Since equivariant homotopy groups take wedges of orthogonal spectra
to direct sums (Corollary \ref{cor-wedges and finite products}~(i)),
any wedge of $\Fc$-equivalences is again an $\Fc$-equivalence.
So the pointset level wedge maps by 
an $\Fc$-equivalence to the wedge of the $\Fc$-cofibrant approximations,
and these two are isomorphic in $\GH_\Fc$.

(ii) This part is essentially dual to~(i), with the caveat that
infinite products of $\Fc$-equivalences need not be $\Fc$-equivalences in general.
To construct a product in $\GH_\Fc$ of the given family,
the abstract recipe is to choose $\Fc$-equivalences $f_i:X_i\to X_i^\Fc$
to $\Fc$-$\Omega$-spectra, and then form the pointset level
product of the replacements.
For all groups $G$ in $\Fc$ we consider the commutative square:
\[ \xymatrix{
 \pi_k^G\left(\prod_{i\in I} X_i\right)\ar[r]\ar[d]_{(\prod f_i)_*} &
\prod_{i\in I} \pi_k^G(X_i) \ar[d]^{\prod (f_i)_*} \\
 \pi_k^G\left(\prod_{i\in I} X^\Fc_i\right)\ar[r] &\prod_{i\in I} \pi_k^G(X^\Fc_i) 
} \]
The upper map is an isomorphism by hypothesis, and the right map
is an isomorphism since each $f_i$ is an $\Fc$-equivalence.
The lower map is also an isomorphism: since all $X_i^\Fc$ 
are $\Fc$-$\Omega$-spectra, the colimit system
\begin{align*}
  \pi_k^G\big( {\prod}_{i\in I}\, X^\Fc_i\big)\ &= \
\colim_{V\in s(\Uc_G)} \, \left[ S^{V\oplus\mR^k}, {\prod}_{i\in I}\, X_i^\Fc(V) \right]^G \\
 &= \
\colim_{V\in s(\Uc_G)} \, \left( {\prod}_{i\in I}\, [ S^{V\oplus\mR^k},  X_i^\Fc(V) ]^G\right) 
\end{align*}
consists of isomorphisms starting at any {\em faithful} $G$-representation $V$.
Since faithful representations are cofinal in the poset $s(\Uc_G)$,
in this particular situation the colimit commutes with the product.
Altogether we can conclude that in our situation the morphism
\[  \prod f_i \ : \  {\prod}_{i\in I}\, X_i\ \to \ {\prod}_{i\in I}\,  X^\Fc_i\]
is an $\Fc$-equivalence.
Since the right hand side is a product in $\GH_\Fc$ of the family $\{X_i\}_{i\in I}$,
so is the left hand side.
\end{proof}

We turn to the interaction of the smash product of orthogonal spectra
with the level and global model structures.
Given two morphisms $f:A\to B$ and $g:X\to Y$ of orthogonal spectra,
the {\em pushout product} morphism\index{subject}{pushout product} is defined as
\[ f\Box g = (f\sm Y)\cup(B\sm g) \ : \ 
A\sm Y\cup_{A\sm X}B\sm X \ \to \ B\sm Y\ .\]
We introduce another piece of notation that is convenient for
the discussion of the monoidal properties: when $I$ and $I'$ are two sets
of morphisms of orthogonal spectra, then we define
\[ I\  \sqsubset \ I'\ \]
to mean that every morphism in $I$ is isomorphic to a morphism in $I'$.
Moreover, we write $I\Box I'$ for the set of pushout product morphisms
$f\Box g$ for all $f\in I$ and all $g\in I'$.
If $\Ec$ and $\Fc$ are global families, then we denote by
$\Ec\times\Fc$ the smallest global family that contains all
groups of the form $G\times K$ for $G\in\Ec$ and $K\in\Fc$.

The sets $I_\Fc$ and $J_\Fc$ of generating cofibrations 
respectively generating acyclic cofibrations 
for the $\Fc$-level model structure
were defined in the proof of Proposition \ref{prop:F-level spectra}.
They consist of the morphisms
\[  F_{H,\mR^m}\sm (i_k)_+
\text{\qquad respectively\qquad}
F_{H,\mR^m} \sm (j_k)_+ \]
for all $k\geq 0$, all $m\geq 0$ and all compact Lie groups $H$ in $\Fc(m)$.
Here we used that the orthogonal spectrum $F_{G,V}( G/H_+)$
is isomorphic to $F_{H,V}$.
The set $K_\Fc$ of acyclic cofibrations for the $\Fc$-global model structure
was defined in \eqref{eq:define K_F}.

\begin{prop}\label{prop:basic box relations} 
Let $\Ec$ and $\Fc$ be two global families.
Then the following relations hold for the sets of generating
cofibrations and acyclic cofibrations:
\[   I_\Ec\, \Box\, I_\Fc \ \sqsubset \ I_{\Ec\times\Fc} \ ,\quad
  I_\Ec\, \Box\, J_\Fc \ \sqsubset \ J_{\Ec\times\Fc} \text{\quad and\quad}
  I_\Ec\, \Box\,  K_\Fc \ \sqsubset \ K_{\Ec\times\Fc} \ . \]
\end{prop}
\begin{proof}
  We start with two key observations concerning 
  the generating cofibrations $i_k:\partial D^k\to D^k$
  and the generating acyclic cofibrations $j_k: D^k\times\{0\}\to D^k\times [0,1]$
  for the model structure of spaces: the pushout product $i_k\Box i_m$
  of two sphere inclusions is homeomorphic to the map $i_{k+m}$;
  similarly, the pushout product $i_k\Box j_m$ is homeomorphic to the map $j_{k+m}$.

  The first relation   $I_\Ec\Box I_\Fc \sqsubset I_{\Ec\times\Fc}$
  is then a consequence of the compatibilities between  
  smash products and the isomorphism \eqref{eq:F_V_smash_F_W} 
  for the smash product of two semifree orthogonal spectra:
  \[ ( F_{G,V}\sm (i_k)_+)\Box ( F_{K,W}\sm (i_m)_+) \ \iso \ 
  ( F_{G,V}\sm F_{K,W})\sm (i_k\Box i_m)_+\ \iso \ 
  F_{G\times K,V\oplus W} \sm (i_{k+m})_+ \]
  The second relation $I_\Ec\Box J_\Fc \sqsubset J_{\Ec\times\Fc}$
  is proved in the same way:
  \[ ( F_{G,V}\sm (i_k)_+)\Box ( F_{K,W}\sm (j_m)_+) \ \iso \ 
  ( F_{G,V}\sm F_{K,W})\sm (i_k\Box j_m)_+ \ \iso \ 
  F_{G\times K,V\oplus W} \sm (j_{k+m})_+ \]
  For the third relation we recall that $c_{K,U,W}$ is the mapping cylinder inclusion 
  of the global equivalence
  \[ \lambda_{K,U,W} \ : \ F_{K,U\oplus W} S^U \ \to \ F_{K,W} \ . \]
  The claim $I_\Ec\Box K_\Fc\sqsubset K_{\Ec\times\Fc}$ then follows from
  \begin{align*}
    ( F_{G,V}\sm (i_k)_+)\Box ( c_{K,U,W}\Box (i_m)_+) \ &\iso \ 
  ( F_{G,V}\sm c_{K,U,W})\Box (i_k\Box i_m)_+ \\ 
  &\iso \   c_{G\times K,V\oplus U, W} \sm (i_{k+m})_+\ .\qedhere
  \end{align*}
\end{proof}

Certain pushout product properties are now formal consequences.
In the special case where $\Ec=\Fc=\td{e}=\Ec\times\Fc$
are the trivial global families, part~(iii) of the previous proposition 
specializes to Proposition~12.6 of \cite{mmss}.

\begin{prop} \label{prop:ExF ppp}
  Let $\Ec$ and $\Fc$ be two global families.
  \begin{enumerate}[\em (i)]
  \item The pushout product of an $\Ec$-cofibration with an $\Fc$-cofibration
    is an $(\Ec\times\Fc)$-cofibration.  
  \item The pushout product of a flat cofibration 
    with an $\Fc$-cofibration that is also an $\Fc$-equivalence
    is a flat cofibration and global equivalence.
  \item Let $\Fc$ be a multiplicative global family,\index{subject}{global family!multiplicative|see{multiplicative global family}}\index{subject}{multiplicative global family}
    i.e., $\Fc\times\Fc =\Fc$. Then the  $\Fc$-global model structure  
    satisfies the pushout product property with respect to the smash product
    of orthogonal spectra.
  \end{enumerate}
\end{prop}
\begin{proof} 
  (i) It suffices to show the claim for a set of generating cofibrations,
  where it follows from the relation  $I_\Ec\Box I_\Fc \sqsubset I_{\Ec\times\Fc}$
  established in Proposition \ref{prop:basic box relations}.

  (ii) Again it suffices to check the pushout product of   
  any generating flat cofibration with a generating acyclic cofibration 
  for the $\Fc$-global model structure.
  For generators, the claim follows from the relation
  \begin{align*}
    I_{\All}\Box( J_\Fc\cup K_\Fc) \ &= \ (I_{\All}\Box J_\Fc)\cup(I_{\All}\Box K_\Fc) \\
  &\sqsubset \ J_{\All\times\Fc} \cup K_{\All\times\Fc} \ = \     J_{\All} \cup K_{\All}  
  \end{align*}
  established in Proposition \ref{prop:basic box relations}.

  (iii) The pushout product of two $\Fc$-cofibrations is an $\Fc$-cofibration by part~(i)
  and the hypothesis that $\Fc$ is multiplicative.
  Since $\Fc$-cofibrations are in particu\-lar flat cofibrations, 
  the pushout product of two $\Fc$-cofibrations 
  one of which is also an $\Fc$-equivalence
  is another $\Fc$-equivalence by part~(ii).
\end{proof}

The sphere spectrum $\mS$ is the unit object for the smash product 
of orthogonal spectra, and it is `free', i.e., $\td{e}$-cofibrant.
Thus $\mS$ is cofibrant in the  
$\Fc$-global model structure for every global family $\Fc$.
So if $\Fc$ is multiplicative, then with respect to the smash product, 
the $\Fc$-global model structure is a 
symmetric monoidal model category 
in the sense of \cite[Def.\,4.2.6]{hovey-book}.
A corollary is that the
homotopy category $\GH_\Fc$, i.e., the localization of the category
of orthogonal spectra at the class of $\Fc$-equivalences,
inherits a closed symmetric monoidal structure, 
compare \cite[Thm.\,4.3.3]{hovey-book}. The 
{\em derived smash product}\index{subject}{derived smash product}\index{subject}{smash product!derived|see{derived smash product}}\index{symbol}{$\sm^\mL_\Fc$ - {derived smash product}} 
 \begin{equation}\label{eq:derived_smash}
\sm^\mL_\Fc \ : \  \GH_\Fc \times \GH_\Fc \ \to \ \GH_\Fc \ ,    
\end{equation}
i.e., the induced symmetric monoidal product on $\GH_\Fc$,
is any total left derived functor of the smash product.

\begin{cor}\label{cor-derived smash} 
For every multiplicative global family $\Fc$, the $\Fc$-global 
homotopy category $\GH_\Fc$\index{symbol}{  $\GH_\Fc$ - {global stable homotopy category with respect to the global family $\Fc$}}
is closed symmetric monoidal under the derived smash product \eqref{eq:derived_smash}. \end{cor}

The value of the derived smash product 
at a pair $(X,Y)$ of orthogonal spectra can be calculated as
\[ X\sm^\mL_\Fc Y \ = \ X^c \sm Y^c\ , \]
where  $X^c\to X$ and $Y^c\to Y$ are
cofibrant replacements in the $\Fc$-global model structure,
i.e., $\Fc$-equivalences with $\Fc$-cofibrant sources.
By the following 'flatness theorem',
it actually suffices to `resolve' only one of the factors,
and it is enough to require the source of the `resolution'
to be flat (as opposed to being $\Fc$-cofibrant).

\begin{theorem}\label{thm:flat is flat}\index{subject}{flatness theorem!for orthogonal spectra} 
Let $\Fc$ be a global family.
\begin{enumerate}[\em (i)]
\item Smashing with a flat orthogonal spectrum preserves $\Fc$-equivalences.
\item Smashing with any orthogonal spectrum preserves $\Fc$-equivalences
between flat orthogonal spectra.
\end{enumerate}
\end{theorem}
\begin{proof}
(i) We let $G$ be a compact Lie group.
Every cofibration of $O(m)$-spaces is also a cofibration 
of $(G\times O(m))$-spaces when we let $G$ act trivially.
So if we let $G$ act trivially on a flat orthogonal spectrum $X$, 
then it is $G$-flat in the sense of Definition \ref{def:G-flat}.
Now we let $f:Y\to Z$ be an $\Fc$-equivalence of orthogonal spectra.
If $G$ belongs to the global family $\Fc$, then $f$ is in particular
a $\upi_*$-isomorphism of orthogonal $G$-spectra with respect
to the trivial $G$-actions. So $X\sm f:X\sm Y\to X\sm Z$ 
is a $\upi_*$-isomorphism of orthogonal $G$-spectra by
Theorem \ref{thm:G-flat is flat}, because $X$ is $G$-flat.
This proves that $X\sm f$ is again an $\Fc$-equivalence.

(ii) 
We let $f:X\to Y$ be an $\Fc$-equivalence between flat orthogonal spectra.
We choose a global equivalence $\varphi:A^\flat \to A$ with flat source
and consider the commutative square:
\[ \xymatrix{ 
A^\flat\sm X\ar[r]^-{\varphi\sm X}\ar[d]_{A^\flat\sm f}&
A\sm X\ar[d]^{A\sm f}\\
A^\flat\sm Y\ar[r]_-{\varphi\sm Y} & A\sm Y
} \]
The two horizontal morphisms are global equivalences by part~(i),
because $X$ and $Y$ are flat. The left vertical morphism is an $\Fc$-equivalence
by part~(i), because $A^\flat$ is flat. 
So $A\sm f$ is also an $\Fc$-equivalence.
\end{proof}

Now we prove another important relationship between
the global model structures and the smash product, 
the {\em monoid axiom} \cite[Def.\,3.3]{schwede-shipley-monoidal}.
As in the unstable situation in Proposition \ref{prop:monoid orthogonal spaces}
we only discuss the weaker form of the monoid axiom with
sequential (as opposed to more general transfinite) compositions.

\begin{prop}[Monoid axiom]\label{prop:monoid axiom}\index{subject}{monoid axiom!for the smash product of orthogonal spectra}
We let $\Fc$ be a global family.
For every flat cofibration $i:A\to B$ that is also an $\Fc$-equivalence
and every orthogonal spectrum $Y$ the morphism
\[ i\sm Y \ : \ A\sm Y \ \to \ B\sm Y \]
is an h-cofibration and an $\Fc$-equivalence.\index{subject}{h-cofibration}
Moreover, the class of h-cofibrations that are also $\Fc$-equivalences
is closed under cobase change, coproducts and sequential compositions.
\end{prop}
\begin{proof}
Given Theorem \ref{thm:flat is flat}, this is a standard argument, 
similar to the proofs of the monoid axiom in the non-equivariant
or the $G$-equivariant context, compare 
\cite[Prop.\,12.5]{mmss}, \cite[Prop.\,1.3.10]{stolz-thesis},
\cite[III Prop.\,7.4]{mandell-may} or \cite[Prop.\,2.3.27]{stolz-thesis}.
Every flat cofibration is an h-cofibration 
(Corollary \ref{cor-h-cofibration closures}~(iii)), 
and h-cofibra\-tions are closed under smashing with any orthogonal spectrum, 
so $i\sm Y$ is an h-cofibration.
Since $i$ is a h-cofibration and $\Fc$-equivalence, its cokernel $B/A$
is $\Fc$-stably contractible by the long exact homotopy group sequence
(Corollary \ref {cor-long exact sequence h-cofibration}).
But $B/A$ is also flat as a cokernel of a flat cofibration,
so the spectrum $(B/A)\sm Y$ is $\Fc$-stably contractible by 
Theorem \ref{thm:flat is flat}. 
Since $i\sm Y$ is an h-cofibration with cokernel isomorphic to $(B/A)\sm Y$, 
the long exact homotopy group sequence then shows that $i\sm Y$ is an $\Fc$-equivalence.

The proof that the class of h-cofibrations that are also $\Fc$-equivalences
is closed under cobase change, coproducts and sequential compositions 
is the same as for for global equivalences
in Corollary \ref{cor-closure h-cof and global equi}.
\end{proof}

Every  $\Fc$-cofibration is in particular a flat cofibration.
So the monoid axiom implies 
the monoid axiom in the  $\Fc$-global model structure.
If the global family $\Fc$ is closed under products,
Theorem \cite[Thm.\,4.1]{schwede-shipley-monoidal} applies
to the  $\Fc$-global model structure and shows the following lifting results.
The additional claims in part~(i) about the behavior of
the forgetful functor on the cofibrations are proved 
as in the unstable analog in Corollary \ref{cor-lift to modules spaces}.

\begin{cor}\label{cor-lift modules} 
Let $R$ be an orthogonal ring spectrum and $\Fc$ 
  a multiplicative global family.\index{subject}{multiplicative global family}
  \begin{enumerate}[\em (i)]
  \item The category of $R$-modules admits the {\em $\Fc$-global model structure}\index{subject}{F-global model structure@$\Fc$-global model structure!for $R$-modules}  
    in which a morphism is an equivalence (respectively fibration)
    if the underlying morphism of orthogonal spectra is an $\Fc$-equivalence
    (respectively $\Fc$-global fibration).
    This model structure is cofibrantly generated.
    Every cofibration in this $\Fc$-global model structure 
    is an h-cofibration of underlying orthogonal spectra.
    If the underlying orthogonal spectrum of $R$ is $\Fc$-cofibrant,
    then every cofibration of $R$-modules is an $\Fc$-cofibration of underlying
    orthogonal spectra.
  \item If $R$ is commutative, then with respect to $\sm_R$ 
    the $\Fc$-global model structure of $R$-modules is a monoidal model category 
    that satisfies the monoid axiom.
  \item If $R$ is commutative, then the category of $R$-algebras 
    admits the {\em  $\Fc$-global model structure}\index{subject}{F-global model structure@$\Fc$-global model structure!for $R$-algebras}  
    in which a morphism is an equivalence (respectively fibration)
    if the underlying morphism of orthogonal spectra is an $\Fc$-equivalence
    (respectively $\Fc$-global fibration).
    Every cofibrant $R$-algebra is also cofibrant as an $R$-module.
  \end{enumerate}
\end{cor}

Strictly speaking, Theorem~4.1 of \cite{schwede-shipley-monoidal}
does not apply verbatim to the $\Fc$-global model structure
because the hypothesis that every object is small (with respect to some regular cardinal)
is not satisfied. However, in our situation the sources of the generating cofibrations
and generating acyclic cofibrations are small with respect to
sequential compositions of h-cofibrations, and this suffices to
run the countable small object argument
(compare also Remark~2.4 of \cite{schwede-shipley-monoidal}).

\begin{prop}\label{prop:sm cofibrant R-mod} 
Let $R$ be an orthogonal ring spectrum and $N$ a right $R$-module
that is cofibrant in the $\All$-global model structure 
of Corollary  {\em \ref{cor-lift modules} (i)}. 
Then for every global family $\Fc$,
the functor $N\sm_R -$ takes $\Fc$-equivalences of left $R$-modules
to $\Fc$-equivalences of orthogonal spectra.
\end{prop}
\begin{proof}
The argument is completely parallel to the unstable
precursor in Proposition \ref{prop:box cofibrant R-mod}. 
We call a right $R$-module $N$ {\em homotopical}
if the functor $N\sm_R -$ takes $\Fc$-equivalences of left $R$-modules
to $\Fc$-equivalences of orthogonal spectra.  
Since the $\All$-global model structure on the category of right $R$-modules
is obtained by lifting the global model structure of orthogonal spectra
along the free and forgetful adjoint functor pair, every cofibrant right $R$-module
is a retract of an $R$-module that arises as the colimit of a sequence
\begin{equation}\label{eq:R-module spectra sequence}
 \ast = M_0 \ \to \ M_1 \ \to \dots \ \to \ M_k \ \to \ \dots 
\end{equation}
in which each $M_k$ is obtained from $M_{k-1}$ as a pushout
\[ \xymatrix{ 
 A_k\sm R \ar[r]^-{f_k\sm R}\ar[d] & 
 B_k \sm R \ar[d] & \\
M_{k-1}\ar[r] & M_k } \]
for some flat cofibration $f_k:A_k\to B_k$ between flat orthogonal spectra.
For example, $f_k$ can be chosen as a wedge of morphisms 
in the set $I_{\All}$ of generating flat cofibrations.
We show by induction on $k$ that each module $M_k$ is homotopical.
The induction starts with the trivial $R$-module, where there is nothing to show.
Now we suppose that $M_{k-1}$ is homotopical, and we claim that then $M_k$
is homotopical as well.
To see this we consider an $\Fc$-equivalence of left $R$-modules
$\varphi:X\to Y$. Then $M_k\sm_R \varphi$
is obtained by passing to horizontal pushouts in the following commutative diagram
of orthogonal spectra:
\[ \xymatrix{ 
M_{k-1}\sm_R X \ar[d]_{M_{k-1}\sm_R \varphi} & 
 A_k \sm X\ar[d]^{A_k\sm \varphi}\ar[l]\ar[r]^-{f_k\sm X}
& B_k \sm X\ar[d]^{B_k\sm \varphi}\\
M_{k-1}\sm_R Y &  A_k \sm Y\ar[l]\ar[r]_-{f_k\sm Y}& B_k \sm Y
} \]
Here we have exploited that $(A_k\sm R)\sm_R X$
is naturally isomorphic to $A_k\sm X$.
In the diagram, the left vertical morphism is an $\Fc$-equivalence by hypothesis. 
The middle and right vertical morphisms are 
$\Fc$-equivalences because smash product with a flat orthogonal spectrum
preserves $\Fc$-equivalences (Theorem \ref{thm:flat is flat}~(i)).
Moreover, since the morphism $f_k$ is a flat cofibration, it is an h-cofibration
(by Corollary \ref{cor-h-cofibration closures}~(iii)),
and so the morphisms $f_k\sm X$ and $f_k\sm Y$ are h-cofibrations.
Corollary \ref{cor-cobase change h-cofibration}~(i)
then implies that the induced morphism on horizontal
pushouts $M_k\sm_R \varphi$ is again an $\Fc$-equivalence.

Now we let $M$ be a colimit of the sequence \eqref{eq:R-module spectra sequence}.
Then $M\sm_R X$ is a colimit of the sequence $M_k \sm_R X$.
Moreover, since $f_k:A_k\to B_k$ is an h-cofibration, so is the 
morphism $f_k\sm R$, and hence also its cobase change $M_{k-1}\to M_k$.
So the sequence whose colimit is $M\sm_R X$ consists of h-cofibrations,
which are in particular levelwise closed embeddings.
The same is true for $M\sm_R Y$.
Since each $M_k$ is homotopical and colimits of orthogonal spectra
along closed embeddings are homotopical
(see Proposition \ref{prop:sequential colimit closed embeddings}~(ii)),
we conclude that the morphism $M\sm_R\varphi:M\sm_R X\to M\sm_R Y$ 
is an $\Fc$-equivalence, so that $M$ is homotopical.
Since $\Fc$-equivalences are closed under retracts, 
the class of homotopical $R$-modules is closed under retracts, and so every
cofibrant right $R$-module is homotopical.
\end{proof}

For later use we record a positive version of 
the global model structure of Theorem \ref{thm:All global spectra}.
This positive variant is the one that can be lifted to the category
of ultra-commutative ring spectra, 
see Theorem \ref{thm:global ultra-commutative} below

\begin{defn} A morphism $f:X\to Y$ of orthogonal spectra is a {\em positive cofibration}\index{subject}{cofibration!positive}\index{subject}{positive cofibration!of orthogonal spectra}
if it is a flat cofibration and the map $f(0):X(0)\to Y(0)$ is a homeomorphism.  
An orthogonal spectrum is a {\em positive global $\Omega$-spectrum}\index{subject}{positive global $\Omega$-spectrum}\index{subject}{global Omega-spectrum@global $\Omega$-spectrum!positive}
if for every compact Lie group $G$, every $G$-representation $V$
and every faithful $G$-representation $W$ with $W\ne 0$ 
the adjoint structure map 
\[ \tilde\sigma_{V,W}\ :\  X(W)\ \to \ \map_*(S^V,X(V\oplus W))\]
is a $G$-weak equivalence.
\end{defn}

If $G$ is a non-trivial compact Lie group, then any faithful $G$-representation
is automatically non-trivial. So a positive global $\Omega$-spectrum
is a global $\Omega$-spectrum (in the absolute sense) 
if in addition the adjoint structure map $\tilde \sigma_{0,\mR}:X(0)\to \Omega X(\mR)$ 
is a non-equivariant weak equivalence.

\begin{prop}[Positive global model structure]\label{prop:positive global spectra} 
  The global equivalences and positive cofibrations are part of 
  a proper, cofibrantly generated, topological model structure,
  the {\em positive global model structure}\index{subject}{positive global model structure!for orthogonal spectra}\index{subject}{global model structure!positive}
  on the category of orthogonal spectra. 
  A morphism $f:X\to Y$ of orthogonal spectra is a fibration 
  in the positive global model structure if and only if
  for every compact Lie group $G$, every $G$-representation $V$
  and every faithful $G$-representation $W$ with $W\ne 0$ 
  the map $f(W)^G:X(W)^G\to Y(W)^G$ is a Serre fibration and the square
    \begin{equation*} 
             \xymatrix@C=13mm{ X(W)^G \ar[d]_{f(W)^G} \ar[r]^-{(\tilde\sigma_{V,W})^G} & 
          \map_*^G(S^V, X(V\oplus W)) \ar[d]^{\map_*^G(S^V,f(V\oplus W))} \\
          Y(W)^G \ar[r]_-{(\tilde\sigma_{V,W})^G} & \map_*^G(S^V, Y(V\oplus W))}
      \end{equation*}
    is homotopy cartesian.  
    The fibrant objects in the positive global model structure 
    are the positive global $\Omega$-spectra.
    The model structure is monoidal with respect 
    to the smash product of orthogonal spectra.
\end{prop}
\begin{proof}
  We start by establishing a {\em positive strong level model structure}.
  A morphism $f:X\to Y$ of orthogonal spectra is a
  {\em positive strong level equivalence}\index{subject}{level equivalence!positive strong!of orthogonal spectra} 
  (respectively {\em positive strong level fibration})\index{subject}{level fibration!positive strong!of orthogonal spectra}
  if for every inner product space $V$ with $V\ne 0$
  the map $f(V):X(V)\to Y(V)$  is an $O(V)$-weak equivalence
  (respectively an $O(V)$-fibration).
  We claim that the positive strong level equivalences, 
  positive strong level fibrations and positive cofibrations form a model structure
  on the category of orthogonal spectra.

  The proof is another application of the general construction of
  level model structures in Proposition \ref{prop:general level model structure}. 
  Indeed, we let $\Cc(0)$ be the degenerate model structure on 
  the category $\bT_*$ of based spaces in which every morphism
  is a weak equivalence and a fibration, but only the isomorphisms
  are cofibrations. 
  For $m\geq 1$ we let $\Cc(m)$ be the projective
  model structure (for the family of all closed subgroups)
  on the category of based $O(m)$-spaces.
  With respect to these choices of model structures $\Cc(m)$,
  the classes of level equivalences, level fibrations and cofibrations
  in the sense of Proposition \ref{prop:general level model structure}
  precisely become the positive strong level equivalences,
  positive strong level fibrations and positive cofibrations.
  The consistency condition (Definition \ref{def:consistency condition})
  is now strictly weaker than for the strong level model structure,
  so it holds. 
  The positive strong level model structure is topological 
  by Proposition \ref{prop:topological criterion},
  where we take $\Gc$ as the set of semifree orthogonal spectra $F_{H,\mR^m}$
  for all $m\geq 1$ and all closed subgroup $H$ of $O(m)$.
  
  We obtain the positive global model structure for orthogonal spectra by `mixing' 
  the positive strong level model structure 
  with the global model structure of Theorem \ref{thm:All global spectra}.
  Every positive strong level equivalence is a global equivalence and every 
  positive cofibration is a flat cofibration.
  The global equivalences and the positive cofibrations
  are part of a model structure by Cole's mixing theorem \cite[Thm.\,2.1]{cole-mixed},
  which is our first claim. 
  By \cite[Cor.\,3.7]{cole-mixed} (or rather its dual formulation),
  an orthogonal spectrum is fibrant in the positive global model structure 
  if it is equivalent in the positive strong level model structure to 
  a global $\Omega$-spectrum; this is equivalent to
  being a positive global $\Omega$-spectrum.
  The proof that the positive global model structure is proper and topological 
  is similar as for the global model structure.
  The proof of the pushout product property is as in the absolute global model
  structure (see Proposition \ref{prop:ExF ppp}); the only new ingredient is
  that the class of generators $F_{G,V}$ with $V\ne 0$ for the positive cofibrations
  is closed under the smash product of orthogonal spectra.

  Finally, the positive global model structure is cofibrantly generated:
  we can simply take the same sets of generating cofibrations 
  and generating acyclic cofibrations as for the global model structure, 
  {\em except} that we  omit all morphisms freely generated in level~0.
  \end{proof}

\section{Triangulated global stable homotopy categories}
\label{sec:generators}

\index{subject}{global stable homotopy category|(}

As the homotopy category of a stable model structure, 
the global stable homotopy category $\GH$
comes with the structure of a triangulated category.
The shift functor is the suspension of orthogonal spectra,
and the distinguished triangles arise from mapping cone sequences.
In this section we collect the aspects of global stable homotopy theory that
are best expressed in terms of the triangulated structure.

More generally, we work in the triangulated $\Fc$-global
stable homotopy category $\GH_\Fc$, where $\Fc$ is any global family. 
Theorem \ref{thm:Bgl G weak generators} identifies
a set of compact generators of $\GH_\Fc$, 
namely the suspension spectra of global classifying spaces of the groups in $\Fc$.
A consequence is Brown representability for cohomological
and homological functors out of $\GH_\Fc$, 
compare Corollary \ref{cor-generators for GH_F}.
Theorem \ref{thm:t-structure on GH} shows that
the classes of $\Fc$-connective and $\Fc$-coconnective spectra 
form a non-degenerate t-structure on $\GH_\Fc$; 
moreover, taking 0-th equi\-variant homotopy groups
is an equivalence from the heart of this t-structure to the
category of $\Fc$-global functors.
Immediate consequences are global Postnikov sections
and the existence of Eilenberg-Mac Lane spectra that realize global functors.
Proposition \ref{prop:Box vs smash} establishes another connection between
the smash product and the algebra of global functors: 
for globally connective orthogonal spectra the box product of global functors 
calculates the 0-th homotopy groups of a derived smash product.

The last topic of this section are certain distinguished triangles 
in the global stable homotopy category that 
arise from special representations,
namely when a compact Lie group $G$ acts faithfully and transitively on the unit sphere.
The stabilizer of a unit vector is then a closed subgroup $H$ 
such that $G/H$ is diffeomorphic to a sphere, 
and Theorem \ref{thm:fundamental triangle} 
exhibits a distinguished triangle that relates the global classifying space
of $G$ to certain semifree orthogonal spectra of $G$ and $H$.
The main examples of this situation are the tautological 
representations of the groups $O(m)$, $S O(m)$, $U(m)$, $S U(m)$ and $S p(m)$;
these will show up again in Section \ref{sec:global Thom} 
in the rank filtrations of global Thom spectra.

\medskip

The $\Fc$-global stable homotopy category $\GH_\Fc$ is the homotopy category
of a stable model structure, so it is naturally a triangulated category,
for example by \cite[Sec.\,7.1]{hovey-book} or \cite[Thm.\,A.12]{schwede-topological}.
The shift functor is modeled by the pointset level suspension of orthogonal
spectra. More precisely, the suspension functor $-\sm S^1:\spec\to\spec$ 
of Construction \ref{con:shift} preserves $\Fc$-equivalences,
so it descends to a functor on the $\Fc$-global stable homotopy category
\[ -\sm S^1 \ : \ \GH_\Fc \ \to \ \GH_\Fc \]
for which we use the same name. 
The distinguished triangles are defined from mapping cone sequences,
i.e., a triangle is distinguished if and only if it is isomorphic, 
in $\GH_\Fc$, to a sequence of the form
\[ X \ \xra{\ f\ }\ Y \ \xra{\ i\ } \ C f \ \xra{\ p \ }\ X\sm S^1 \]
for some morphism of orthogonal spectra $f:X\to Y$; here the morphisms $i$ 
and $p$ were defined in \eqref{eq:define i and p}.

\begin{eg}[Shift preserves distinguished triangles]\label{eg:shift preserves triangles}
The shift functor $\sh:\spec\to\spec$ 
of Construction \ref{con:shift} preserves $\Fc$-equivalences,
so it descends to a functor on the $\Fc$-global stable homotopy category
\[ \sh \ : \ \GH_\Fc \ \to \ \GH_\Fc \]
for which we use the same name. Moreover, shifting commutes with smashing
with a based space {\em on the nose}, i.e., $(\sh X)\sm A=\sh(X\sm A)$;
so we can (and will) leave out parentheses in such expressions.
Since the suspension functor on $\GH_\Fc$ is induced by smashing with $S^1$,
the shift functor commutes with the suspension functor, again on the nose,
both on the pointset level and also on the level of the $\Fc$-global stable homotopy
category. We will now argue that shifting also preserves distinguished 
triangles on the nose; equivalently, the derived shift is an exact
functor of triangulated categories if we equip it with the identity isomorphism
$\sh\circ(-\sm S^1)=(-\sm S^1)\circ\sh$.
  
To prove our claim we consider a distinguished triangle in $\GH_\Fc$:
\[ A \ \xra{\ f\ }\ B \ \xra{\ g\ } \ C  \ \xra{\ h \ }\ A\sm S^1 \]
The morphism $\lambda_A:A\sm S^1\to \sh A$ is a natural global equivalence
by Proposition \ref{prop:global equiv preservation} (i), so all
vertical morphisms in the following diagram are isomorphisms in $\GH_\Fc$:
\[ \xymatrix@C=10mm{ 
 A \sm S^1 \ar[r]^-{f\sm S^1}\ar[d]_{\lambda_A} &
 B \sm S^1 \ar[r]^-{g\sm S^1} \ar[d]_{\lambda_B} & 
C\sm S^1 \ar[r]^-{- h\sm S^1}\ar[d]^{\lambda_C}&
A\sm S^1\sm S^1 \ar[d]^{\lambda_A\sm S^1} \\
\sh A \ar[r]_-{\sh f} & \sh B \ar[r]_-{\sh g} & \sh C \ar[r]_-{\sh h}&  \sh A\sm S^1 } \]
The claim now follows from the following three observations:
\begin{itemize}
\item the suspension functor in any triangulated category 
preserves distinguished triangles {\em up to a sign};
so the upper triangle above is distinguished;
\item the left and middle squares commute by naturality of the $\lambda$-morphisms;
\item the right square commutes because the two morphisms
$\lambda_{A\sm S^1},\lambda_A\sm S^1:A\sm S^1\sm S^1 \to\sh A\sm S^1$
differ by the twist involution of $S^1\sm S^1$; since this involution has degree~-1,
we obtain
\[  (\sh h)\circ\lambda_C \ = \ 
\lambda_{A\sm S^1}\circ (h\sm S^1)\ = \ - (\lambda_A\sm S^1)\circ (h\sm S^1)\ .\]
\end{itemize}
\end{eg}

We recall the notion of compactly generated triangulated categories.
Compact generation has strong formal consequences, see Theorem \ref{thm:Brown rep} below.
Theorem \ref{thm:Bgl G weak generators} shows that the 
$\Fc$-global stable homotopy category enjoys this special property.

\begin{defn}\label{def:compactly generated}
Let $\Tc$ be a triangulated category with infinite sums.
An object $C$ of $\Tc$ is {\em compact}\index{subject}{compact!object in a triangulated category} 
if for every family $\{X_i\}_{i\in I}$ of objects the canonical map
\[ {\bigoplus}_{i\in I}\, \Tc(C,\, X_i) \ \to \ \Tc(C,\, {\oplus}_{i\in I}\, X_i)
\]
is an isomorphism. 
A set $\Gc$ of objects of $\Tc$ is called a set of
{\em weak generators}\index{subject}{generator!weak}
if the following condition holds: if 
the groups $\Tc(G[k],X)$ are trivial for all $k\in\mZ$ and all $G\in\Gc$, 
then $X$ is a zero object.
The triangulated category $\Tc$ is {\em compactly generated}\index{subject}{compactly generated!triangulated category}
if it has sums and a set of compact weak generators.
\end{defn}

If $G$ is from a global family $\Fc$, then the functor $\pi_0^G:\spec\to\Ab$
takes $\Fc$-equivalences to isomorphisms. So the universal property
of a localization provides a unique factorization
\[ \pi_0^G\ :\ \GH_\Fc \ \to\ \Ab \]
through the $\Fc$-global stable homotopy category. We will abuse notation
and use the same symbol for the equivariant homotopy group functor
on the category of orthogonal spectra and for its `derived'
functor defined on $\GH_\Fc$. This abuse of notation is mostly harmless,
but there is one point where it can create confusion, namely in the context of
infinite products; 
we refer the reader to Remark \ref{rk:homotopy of infinite product} for this issue.

We recall from Definition \ref{def:global classifying} 
that the global classifying space $B_{\gl}G$\index{subject}{global classifying space} 
of a compact Lie group $G$ is the semifree orthogonal space
$\bL_{G,V}=\bL(V,-)/G$, for some faithful $G$-representation $V$.
The choice of faithful representation is omitted from the notation 
because the global homotopy type of $B_{\gl}G$ does not depend on it.
The suspension spectrum of $B_{\gl} G$ comes with a stable tautological class
\[  e_G \ = \  e_{G,V} \in \ \pi_0^G(\Sigma^\infty_+ B_{\gl} G) \]
defined in \eqref{eq:define_stable_tautological}.
\index{subject}{stable tautological class}

In the proof of the next theorem
we will start using the shorthand notation\index{symbol}{$\gh{X,Y}_\Fc$ - { morphism group in the $\Fc$-global stable homotopy category}}
\[ \gh{ X,Y }_\Fc \ = \ \GH_\Fc(X,Y)\]
for the abelian group of morphisms
in the triangulated $\Fc$-global stable homotopy category.

\begin{theorem}\label{thm:Bgl G weak generators}
Let $\Fc$ be a global family and $G$ a compact Lie group in $\Fc$.
\begin{enumerate}[\em (i)]
\item The pair $(\Sigma^\infty_+ B_{\gl}G, e_G)$ represents the functor
$\pi_0^G:\GH_\Fc \to \emph{(sets)}$.
\item The orthogonal spectrum $\Sigma^\infty_+ B_{\gl}G$ 
is a compact object of the $\Fc$-global stable homotopy category $\GH_\Fc$.
\item As $G$ varies through a set of representatives
of isomorphism classes of groups in $\Fc$, the spectra
$\Sigma^\infty_+ B_{\gl}G$ form a set of weak generators for 
the $\Fc$-global stable homotopy category $\GH_\Fc$. 
\end{enumerate}
In particular, the $\Fc$-global stable homotopy category $\GH_\Fc$ 
is compactly generated.
\end{theorem}
\begin{proof}
(i) We need to show that for every orthogonal spectrum $X$ the map
\[ \gh{ \Sigma^\infty_+ B_{\gl}G, X }_\Fc \ \to \ \pi_0^G (X)\ , \quad 
f\longmapsto f_*(e_G)\]
is bijective. Since both sides take $\Fc$-equivalences in $X$ to bijections,
we can assume that $X$ is an $\Fc$-$\Omega$-spectrum, and hence fibrant
in the $\Fc$-global model structure.
For $G$ in the family $\Fc$, 
the orthogonal spectrum $\Sigma^\infty_+ B_{\gl}G$ is $\Fc$-cofi\-brant.
So the localization functor induces a bijection
\[ \spec(\Sigma^\infty_+ B_{\gl}G, X) / \text{homotopy} \ \to \ 
\gh{ \Sigma^\infty_+ B_{\gl}G, X }_\Fc \] 
from the set of homotopy classes of morphisms of orthogonal spectra
to the set of morphisms in $\GH_\Fc$.

We let $V$ be the faithful $G$-representation that is implicit in the
definition of the global classifying space $B_{\gl}G$. 
By the freeness property of $B_{\gl} G=\bL_{G,V}$, morphisms from 
$\Sigma^\infty_+ B_{\gl}G$ to $X$ biject with based $G$-maps $S^V\to X(V)$,
and similarly for homotopies.
The composite
\[ [ S^V, X(V) ]^G  \ \xra{\ \iso\ } \
 \gh{ \Sigma^\infty_+ B_{\gl}G, X }_\Fc \ \xra{f\mapsto f_*(e_G)} \ \pi_0^G (X)  \]
is the stabilization map, and hence bijective by 
Proposition \ref{prop:stable zero is level zero}~(ii).
Since the left map and the composite are bijective, so is the evaluation
map at the stable tautological class.

(ii)
By Proposition \ref{prop:sums and product in GH} 
the wedge of any family of orthogonal spectra is a coproduct in $\GH_\Fc$. 
The vertical maps in the commutative square
\[\xymatrix{  {\bigoplus_{i\in I}} \gh{ \Sigma^\infty_+ B_{\gl}G,\, X_i }_\Fc \ar[r]\ar[d] &
\gh{ \Sigma^\infty_+ B_{\gl}G,\, \bigoplus_{i\in I}X_i }_\Fc \ar[d] \\
{\bigoplus_{i\in I}} \pi^G_0(X_i) \ar[r] &
\pi_0^G\left( \bigvee_{i\in I}X_i\right) }\]
are evaluation at the stable tautological class,
and hence isomorphisms by part~(i).
The lower horizontal map is an isomorphism by 
Corollary \ref{cor-wedges and finite products}~(i),
hence so is the upper horizontal map.
This shows that $\Sigma^\infty_+ B_{\gl}G$ is compact 
as an object of the triangulated category $\GH_\Fc$.

(iii) If $X$ is an orthogonal spectrum such that the group
$\gh{ \Sigma^\infty_+ B_{\gl}G[k],\, X}$ is trivial for every $G$ in $\Fc$
and all integers $k$, 
then $X$ is $\Fc$-equivalent to the trivial orthogonal spectrum by part~(i);
so $X$ is a zero object in $\GH_\Fc$. This proves that the
spectra $\Sigma^\infty_+ B_{\gl}G$ form a set of weak generators 
for $\GH_\Fc$ as $G$ varies over $\Fc$.
\end{proof}

A covariant functor $E$ from a triangulated category $\Tc$ 
to the category of abelian groups is called
{\em homological}\index{subject}{homological functor} if 
for every distinguished
triangle $(f,g,h)$ in $\Tc$ the sequence of abelian groups
\[ E(A) \ \xra{E(f)} \ E(B) \ \xra{E(g)} \ E(C)  \ \xra{E(h)} \ E(A[1])\]
is exact.
A contravariant functor $E$ from $\Tc$ to the category of abelian groups 
is called {\em cohomological}\index{subject}{cohomological functor} if 
for every distinguished
triangle $(f,g,h)$ in $\Tc$ the sequence of abelian groups
\[ E(A[1]) \ \xra{E(h)} \ E(C) \ \xra{E(g)} \ E(B)  \ \xra{E(f)} \ E(A)\]
is exact.

\begin{theorem}[Brown representability]\label{thm:Brown rep}
Let $\Tc$ be a compactly generated triangulated category.
\begin{enumerate}[\em (i)]
\item Every cohomological functor $\Tc^{\op}\to\Ab$ that takes sums in $\Tc$ 
to products of abelian groups is representable, 
i.e., isomorphic to the functor $\Tc(-,X)$ for some object $X$ of $\Tc$.
\item Every homological functor $\Tc\to\Ab$ that takes products in $\Tc$ 
to products of abelian groups is representable, 
i.e., isomorphic to the functor $\Tc(Y,-)$ for some object $Y$ of $\Tc$.
\item An exact functor $F:\Tc\to\Sc$ to another triangulated category
  has a right adjoint if and only if it takes sums in $\Tc$ to
  sums in $\Sc$.
\item An exact functor $F:\Tc\to\Sc$ to another triangulated category
  has a left adjoint if and only if it takes products in $\Tc$ to
  products in $\Sc$.
\end{enumerate}
\end{theorem}

A proof of part~(i) of this form of Brown representability can be found in
 \cite[Thm.\,3.1]{neeman-Grothendieck}
or \cite[Thm.\,A]{krause-Brown coherent}.
A proof of part~(ii) of this form of Brown representability can be found 
in \cite[Thm.\,8.6.1]{neeman-triangulated categories}
or \cite[Thm.\,B]{krause-Brown coherent}.
Part~(iii) is a formal consequence of part~(i):
if $F$ preserves sums, then for every object $X$ of $\Tc$ the functor
\[ \Sc(F(-), X) \ : \ \Tc^{\op} \ \to \ \Ab \]
is cohomological and takes sums to products.
Hence the functor is representable by an object $R X$ in $\Tc$
and an isomorphism
\[  \Tc(A, R X) \ \iso \ \Sc(F A, X) \ , \]
natural in $A$.
Once this representing data is chosen, the assignment $X\mapsto R X$ extends 
canonically to a functor $R:\Sc\to\Tc$ that is right adjoint to $F$.
In much the same way, part~(iv) is a formal consequence of part~(ii).

\medskip

We let $\Tc$ be a triangulated category with sums.
A {\em localizing subcategory}\index{subject}{localizing subcategory} 
of $\Tc$ is a full subcategory $\Xc$ which is closed under sums and under 
extensions in the following sense: if
\[ A \ \xra{\ f \ }\ B \ \xra{\ g\ }\ C \  \xra{\ h\ }\ A[1] \]
is a distinguished triangle in $\Tc$
such that two of the objects $A$, $B$ or $C$ 
belong to $\Xc$, then so does the third.
A set of {\em compact} objects is a set of weak generators in the sense
of Definition \ref{def:compactly generated} if and only if the smallest
localizing subcategory containing the set is all of $\Tc$,
see for example \cite[Lemma 2.2.1]{schwede-shipley-modules}.
Hence Theorem \ref{thm:Bgl G weak generators} entitles us to the 
following corollary.

\begin{cor}\label{cor-generators for GH_F}
Let $\Fc$ be a global family and $\Sc$ a triangulated category.
\begin{enumerate}[\em (i)]
\item 
Every localizing subcategory of the $\Fc$-global stable homotopy category 
that contains the spectrum $\Sigma^\infty_+ B_{\gl}G$ for every group $G$ of $\Fc$
is all of $\GH_\Fc$.
\item Every cohomological functor on $\GH_\Fc$ 
that takes sums to products is representable.  
\item Every homological functor on $\GH_\Fc$ that takes products
to products is representable. 
\item An exact functor $F:\GH_\Fc\to\Sc$  has a right adjoint if 
  and only if it preserves sums.
\item An exact functor $F:\GH_\Fc\to\Sc$  has a left adjoint 
  if and only if it preserves products.
\end{enumerate}
\end{cor}

\begin{rk}[Equivariant homotopy groups of infinite products]\label{rk:homotopy of infinite product}\index{subject}{equivariant homotopy group!of a product}
In Corollary \ref{cor-wedges and finite products}~(ii) 
we showed that 
for every compact Lie group $G$
the functor $\pi_0^G:\spec\to\Ab$ preserves {\em finite} products.
However, it is {\em not} true that $\pi_0^G$, as a functor on the category of 
orthogonal spectra, preserves infinite products in general.

On the other hand, the `derived' functor $\pi_0^G:\GH\to\Ab$ is representable
by the spectrum $\Sigma^\infty_+ B_{\gl}G$, so it preserves infinite products.
This is no contradiction because
an infinite product of orthogonal spectra is {\em not} in general
a product in the global homotopy category. To calculate a
product in $\GH$ of a family $\{X_i\}_{i\in I}$ of orthogonal spectra,
one has to choose stable equivalences $f_i:X_i\to X_i^\text{f}$
to global $\Omega$-spectra.
For an infinite indexing set, the morphism
\[ {\prod}_{i\in I}\, f_i\ : \ {\prod}_{i\in I}\, X_i\ \to \ {\prod}_{i\in I}\, X_i^\text{f} \]
may fail to be a global equivalence, and then the target,
but not the source, of this map is a product in $\GH$ 
of the family $\{X_i\}_{i\in I}$.
So when considering infinite products it is important to be aware 
of our abuse of notation and to remember that the symbol $\pi_0^G$
has two different meanings.
\end{rk}

The preferred set of compact generators $\{\Sigma^\infty_+ B_{\gl}G\}$
of the global stable homotopy category has another special property,
it is `positive' in the following sense: for all compact Lie groups $G$ and $K$
and all $n>0$ the group
\[ \gh{\Sigma^\infty_+ B_{\gl} G,\ \Sigma^\infty_+ B_{\gl} K[n]} \ \iso \ 
 \pi_0^G(\Sigma^\infty_+ B_{\gl} K[n]) 
\ \iso \ \pi_{-n}^G(\Sigma^\infty_+ B_{\gl} K)  \]
is trivial because every orthogonal suspension spectrum
is globally connective (see Proposition \ref{prop:pi_0 of Sigma^infty}). 
A set of positive compact generators in this sense has strong implications,
as we shall now explain. 

A {\em t-structure} as introduced by 
Beilinson, Bernstein and Deligne in \cite[Def.\,1.3.1]{BBD}
axiomatizes the situation in the derived category of an abelian category
given by cochain complexes whose cohomology vanishes 
in positive respectively negative dimensions.
We are mainly interested in spectra, where a homological
(as opposed to {\em co\,}homological) grading is more common.
So we adapt the definition of a t-structure to homological notation.

\begin{defn}\label{def:t-structure}
A {\em t-structure}\index{subject}{t-structure} on a triangulated category $\Tc$
is a pair $(\Tc_{\geq 0},\Tc_{\leq 0})$ of full subcategories satisfying the
following three conditions, where $\Tc_{\geq n}= \Tc_{\geq 0}[n]$ 
and $\Tc_{\leq n}=\Tc_{\leq 0}[n]$.
\begin{enumerate}[(a)]
\item For all $X\in\Tc_{\geq 0}$ and all $Y\in\Tc_{\leq -1}$ we have $\Tc(X,Y)=0$.
\item $\Tc_{\geq 0}\subset \Tc_{\geq -1}$ and $\Tc_{\leq 0}\supset \Tc_{\leq -1}$. 
\item For every object $X$ of $\Tc$ there is a distinguished triangle
\[ A \ \to \ X \ \to \ B \ \to \ A[1]  \]
such that $A\in\Tc_{\geq 0}$ and $B\in\Tc_{\leq -1}$.
\end{enumerate}
The t-structure is {\em non-degenerate}
if every object in the intersection $\bigcap_{n\in \mZ} \Tc_{\geq n}$ is a zero object
and every object in the intersection $\bigcap_{n\in \mZ} \Tc_{\leq n}$ is a zero object.
\end{defn}

Some of the basic results of Beilinson, Bernstein and Deligne about t-struc\-tures are
(in our homological notation):
\begin{itemize}
\item  the inclusion $\Tc_{\geq n}\to \Tc$ has a right adjoint
$\tau_{\geq n}:\Tc\to\Tc_{\geq n}$, 
and the inclusion $\Tc_{\leq n}\to \Tc$ 
has a left adjoint $\tau_{\leq n}:\Tc\to\Tc_{\leq n}$ \cite[Prop.\,1.3.3]{BBD};
\item given choices of adjoints as above, then for all
$m\leq n$ there is a preferred natural isomorphism of functors 
between $\tau_{\geq m}\circ\tau_{\leq n}$ and $\tau_{\leq n}\circ\tau_{\geq m}$ \cite[Prop.\,1.3.5]{BBD};
\item 
the {\em heart}\index{subject}{heart! of a t-structure} 
\[ \mathscr H \ = \ \Tc_{\geq 0}\cap \Tc_{\leq 0} \ ,  \]
viewed as a full subcategory of $\Tc$, is an abelian category
and $\tau_{\leq 0}\circ\tau_{\geq 0}:\Tc\to\mathscr H$ is a homological 
functor \cite[Thm.\,1.3.6]{BBD}.
\end{itemize}

We will now show that the global stable homotopy category
has a preferred t-structure whose heart is equivalent
to the category of global functors.

\begin{defn}
Let $\Fc$ be a global family.
An orthogonal spectrum $X$ is {\em $\Fc$-connective}\index{subject}{F-connective@$\Fc$-connective}
if the homotopy group $\pi_n^G(X)$ is trivial for 
every group $G$ in $\Fc$ and every $n<0$.
An orthogonal spectrum $X$ is {\em $\Fc$-coconnective}\index{subject}{F-coconnective@$\Fc$-coconnective}
if the homotopy group $\pi_n^G(X)$ is trivial 
for every group $G$ in $\Fc$ and every $n>0$.
\end{defn}

An {\em $\Fc$-global functor}\index{subject}{F-global functor@$\Fc$-global functor}
is an additive functor to abelian groups
from the full subcategory of the Burnside category
with objects the global family $\Fc$.
A morphism of $\Fc$-global functors is a natural transformation.
We write $\GF_\Fc$ for the category of $\Fc$-global functors.
\index{symbol}{  $\GF_\Fc$ - {category of $\Fc$-global functors}}

\begin{theorem}\label{thm:t-structure on GH}
  For every global family $\Fc$,
  the classes of $\Fc$-connective spectra and $\Fc$-coconnective spectra 
  form a non-degenerate t-structure on $\GH_\Fc$.
  The heart of this t-structure 
  consists of those orthogonal spectra $X$ such that $\pi_n^G(X)=0$
  for all $G\in\Fc$ and all $n\ne 0$.
  The functor
  \[ \upi_0 \ : \ \mathscr H \ \to \ \GF_\Fc \]
  is an equivalence of categories from the heart of the t-structure to the
  category of $\Fc$-global functors.
\end{theorem}
\begin{proof}
We deduce this from the more general arguments 
of Beligiannis and Reiten \cite[Ch.\,III]{beligiannis-reiten} 
who systematically investigate torsion pairs 
and t-structures in triangulated categories that are generated by compact objects. 
By Theorem \ref{thm:Bgl G weak generators} the set 
\[ \Pc\ =\ \{\Sigma^\infty_+ B_{\gl} G\}_{[G]\in\Fc} \]
is a set of compact weak generators for the triangulated category $\GH_\Fc$,
where $[G]$ indicates that we choose representatives from the
isomorphism classes of compact Lie groups in $\Fc$.

We let $\Yc$ be the class of orthogonal spectra $Y$ such that 
\[ \gh{P[n], Y}_\Fc \ = \ 0  \]
for all $P\in\Pc$ and all $n\geq 0$.
The representability result of Theorem \ref{thm:Bgl G weak generators}~(i)
shows that this is precisely the class of those $Y$ such that the
group $\pi_n^G(Y)$ vanishes for all $G\in \Fc$ and all $n\geq 0$.
So $\Yc[1]$ is the class of $\Fc$-coconnective orthogonal spectra.
We let $\Xc$ be the `left orthogonal' to $\Yc$, i.e., 
the class of orthogonal spectra $X$ such that 
$\gh{X,Y}_\Fc=0$ for all $Y\in\Yc$.
Since the objects of $\Pc$ are compact in $\GH_\Fc$
by Theorem \ref{thm:Bgl G weak generators}~(ii), 
\cite[Thm.\,III.2.3]{beligiannis-reiten} shows that 
the pair $(\Xc,\Yc)$ is a `torsion pair'
in the sense of \cite[Def.\,I.2.1]{beligiannis-reiten}.
This simply means that the pair $(\Xc,\Yc[1])$ is a t-structure
in the sense of Definition \ref{def:t-structure}, 
see \cite[Prop.\,I.2.13]{beligiannis-reiten}.

It remains to show that $\Xc$ coincides with the class of
$\Fc$-connective orthogonal spectra.
This needs the `positivity property' of the set $\Pc$ of compact generators.
Representability (Theorem \ref{thm:Bgl G weak generators}~(i))
and the suspension isomorphism provide an isomorphism
\[ \gh{\Sigma^\infty_+ B_{\gl} G, \Sigma^\infty_+ B_{\gl} K[n]}_\Fc\ \iso \  
\pi_0^G(\Sigma^\infty_+ B_{\gl}K[n]) \ \iso \ \pi_{-n}^G(\Sigma^\infty_+ B_{\gl}K) \ .\]
The spectrum $\Sigma^\infty_+ B_{\gl}K$ is globally connective by 
Proposition \ref{prop:pi_0 of Sigma^infty},
so this latter group vanishes for all $n\geq 1$ and all $G,K\in\Fc$.
So \cite[Prop.\,III.2.8]{beligiannis-reiten}
shows that $\Xc$ coincides with the class of 
those orthogonal spectra $X$ such that
$\gh{\Sigma^\infty_+ B_{\gl} G, X[n]}_\Fc=0$
for all $G\in\Fc$ and $n\geq 1$.
Since the group $\gh{\Sigma^\infty_+ B_{\gl} G, X[n]}_\Fc$ is isomorphic
to $\pi_{-n}^G(X)$, this shows that $\Xc$ is precisely the class of
$\Fc$-connective orthogonal spectra.
The t-structure is non-degenerate because spectra with trivial $\Fc$-equivariant
homotopy groups are zero objects in $\GH_\Fc$.

We denote by $\End(\Pc)$ the `endomorphism category' of the set $\Pc$, 
i.e., the full pre-additive subcategory of $\GH_\Fc$ with object set $\Pc$.
By an $\End(\Pc)$-module we mean an additive functor
\[ \End(\Pc)^{\op} \ \to \ \Ab \]
from the opposite category of $\End(\Pc)$.
The tautological functor
\begin{equation}\label{eq:tautological functor}
\GH_\Fc \ \to \ \Mod\End(\Pc)   
\end{equation}
takes an object $X$ to the restriction of the 
contravariant Hom-functor $\gh{-,X}_\Fc$ to the full subcategory $\End(\Pc)$.
Because $\Pc$ is a set of positive, compact generators 
for the triangulated category $\GH_\Fc$, \cite[Thm.\,III.3.4]{beligiannis-reiten}
shows that the restriction of the tautological functor \eqref{eq:tautological functor}
to the heart is an equivalence of categories
\[ \ \mathscr H \ \xra{\ \iso \ } \ \Mod\End(\Pc) \ .\]
So to establish the last claim it suffices to show that $\End(\Pc)$ 
is anti-equivalent to the full subcategory of the Burnside category $\bA$ with
object class $\Fc$, in such a way that the tautological functor 
corresponds to the functor $\upi_0$.
The equivalence $\End(\Pc)^{\op}\to\bA_\Fc$ is given by the inclusion on objects, 
and on morphisms by the isomorphisms
\[ \gh{\Sigma^\infty_+ B_{\gl} G ,\Sigma^\infty_+ B_{\gl} K }_\Fc \ \iso \ 
\pi_0^G(\Sigma^\infty_+ B_{\gl} K) \ \iso \ \bA(K,G)\]
specified in Theorem \ref{thm:Bgl G weak generators}~(i) respectively
Proposition \ref{prop:B_gl represents}.  
\end{proof}

\begin{rk}[Postnikov sections]\index{subject}{Postnikov section!of an orthogonal spectrum}
For the standard t-structure on the global stable homotopy category 
(i.e., Theorem \ref{thm:t-structure on GH} for $\Fc=\All$)
the truncation functor
\[ \tau_{\leq n} \ : \ \GH \ \to \ \GH_{\leq n} \ , \]
left adjoint to the inclusion, provides `global Postnikov sections':
for every orthogonal spectrum $X$ the spectrum $\tau_{\leq n}X$
satisfies $\upi_k(\tau_{\leq n}X)=0$ for $k>n$ and the adjunction unit
$X\to X_{\leq n}$ induces an isomorphism on the global functor $\upi_k$
for every $k\leq n$.
\end{rk}

\begin{rk}[Eilenberg-Mac\,Lane spectra]\label{rk:general Eilenberg Mac Lane}
In the case $\Fc=\All$ of the maximal global family,
Theorem \ref{thm:t-structure on GH}
in particular provides an Eilenberg-Mac\,Lane spectrum\index{subject}{Eilenberg-Mac\,Lane spectrum!of a global functor}\index{symbol}{$H M$ - {Eilenberg-Mac\,Lane spectrum of a global functor}} 
for every global functor $M$, i.e., an orthogonal spectrum $HM$
such that $\upi_k(HM)=0$ for all $k\ne 0$
and such that the global functor $\upi_0(HM)$ is isomorphic to $M$;
moreover, these properties characterize $HM$ up to preferred isomorphism in $\GH$. 
Indeed, a choice of inverse to the equivalence $\upi_0$ of
Theorem \ref{thm:t-structure on GH}~, composed with the inclusion of
the heart, provides an Eilenberg-Mac\,Lane functor
\[ H \ : \ \GF \ \to \ \GH \]
to the global stable homotopy category.
\end{rk}

We let $\Tc$ be a triangulated category with infinite sums and
we let $\Cc$ be a class of objects of $\Tc$. 
We denote by $\td{\Cc}_+$ the smallest class of objects of $\Tc$
that contains $\Cc$, is closed under sums (possibly infinite)
and is {\em closed under cones} in the following sense: if
\[ A \to B \to C  \to  A[1]    \]
is a distinguished triangle such that $A$ and $B$ belong to the class,
then so does $C$. 
Any non-empty class of objects that is closed under cones 
contains all zero objects (because a zero object
is a cone of any identity morphism) and is closed under suspension 
(because $A[1]$ is a cone of the morphism from $A$ to a zero object).

\begin{prop}\label{prop:plus generating GH}
The class $\td{\Sigma^\infty_+ B_{\gl} G\colon\text{ $G$ compact Lie}}_+$ 
in the triangulated category $\GH$ 
coincides with the class of globally connective spectra.  
\end{prop}
\begin{proof}
The spectra $\Sigma^\infty_+ B_{\gl} G$ are globally connective by
Proposition \ref{prop:pi_0 of Sigma^infty}.
The class of globally connective orthogonal spectra is closed under sums
because equivariant homotopy groups commute with sums 
(Corollary \ref{cor-wedges and finite products}~(i)).
The class of globally connective orthogonal spectra is closed under cones
by the long exact homotopy group sequence of Proposition \ref{prop:LES for homotopy of cone and fibre}. So the class $\td{\Sigma^\infty_+ B_{\gl} G}_+$ 
is contained in the class of globally connective spectra.
  
For the converse we choose a set of representatives 
of the isomorphism classes of compact Lie groups and set
\[ \Pc\ = \ \{ (\Sigma^\infty_+ B_{\gl}G)[k]\}_{[G], k\geq 0} \ ,\]
indexed by the chosen representatives.
We let $X$ be any orthogonal spectrum.
By induction on $n$ we construct objects $A_n$ in $\td{\Pc}_+$
and morphisms $i_n:A_n\to A_{n+1}$ and  $u_n:A_n\to X$ such that
$u_{n+1} i_n=u_n$. We start with 
\[ A_0 \ = \ \bigoplus_{P\in \Pc, x\in \gh{P,X}} \, P\ .\]
Then  $A_0$ belongs to $\td{\Pc}_+$ and the canonical morphism
$u_0:A_0\to X$ (i.e., the morphism $x$ on the summand indexed by $x$)
induces a surjection $\gh{P,u_0}:\gh{P,A_0}\to\gh{P,X}$
for all $P\in\Pc$.

In the inductive step we suppose that $A_n$ and $u_n:A_n\to X$
have already been constructed. We define
\[ D_n \ = \ \bigoplus_{P\in \Pc, x\in \ker\gh{P,u_n}} \, P\ ,\]
which comes with a tautological morphism 
$\tau:D_n\to A_n$, again given by $x$ on the summand indexed by $x$.
We choose a distinguished triangle
\[ D_n \xra{\ \tau\ } A_n \xra{\ i_n\ } A_{n+1} \to D_n [1] \ . \]
Since $D_n$ and $A_n$ belong to the class $\td{\Pc}_+$, we also have $A_{n+1}\in\td{\Pc}_+$.
Since $u_n\tau=0 $ (by definition), 
we can choose a morphism $u_{n+1}:A_{n+1}\to X$ such that
$u_{n+1}i_n=u_n$. This completes the inductive construction.

Now we choose a homotopy colimit $(A,\{\varphi_n:A_n\to A\}_n)$
of the sequence of morphisms $i_n:A_n\to A_{n+1}$, i.e.,
a distinguished triangle
\[ {\bigoplus}_{n\geq 0}\, A_n \ \xra{1-\sh}\ 
 {\bigoplus}_{n\geq 0}\, A_n \ \xra{\qquad}\ 
A \ \xra{\qquad}\  {\bigoplus}_{n\geq 0} A_n [1]\ . \]
Since all the objects $A_n$ are in $\td{\Pc}_+$, so is $A$.
Since a homotopy colimit in $\GH$ is a weak colimit,
we can choose a morphism $u:A\to X$ such that $u\varphi_n=u_n$ 
for all $n\geq 0$.
The map
\[\gh{P,u_0}\ = \ \gh{P,u}\circ \gh{P,\varphi_0}\ : \  \gh{P,A_0}\ \to\ \gh{P,X} \]
is surjective for $P\in\Pc$, 
hence so is $\gh{P,u}:\gh{P,A}\to\gh{P,X}$.

We claim that $\gh{P,u}$ is also injective for all $P\in\Pc$.
We let $\alpha:P\to A$ be a morphism such that $u \alpha=0$. 
Since $P$ is compact, there is an $n\geq 0$ and a morphism $\alpha':P\to A_n$ 
such that $\alpha=\varphi_n\alpha'$.
Then $u_n\alpha'=u\varphi_n\alpha'=u\alpha=0$.
So $\alpha'$ indexes one of the summands of $D_n$.
Thus $\alpha'$ factors through the tautological 
morphism $\tau:D_n\to A_n$ as $\alpha'=\tau\alpha''$, and hence 
\[ \alpha \ = \ \varphi_n\alpha' \ = \ \varphi_{n+1}i_n\tau\alpha'' \ = \ 0 \]
since $i_n$ and $\tau$ are consecutive morphisms in a distinguished triangle.
So $\gh{P,u}$ is also injective, hence bijective.
By the choice of the set $\Pc$ and the natural isomorphism
\[  \gh{(\Sigma^\infty_+ B_{\gl}G)[k],A}\ \iso \ \pi_k^G(A)\ , \]
this shows that the morphism $u:A\to X$ induces isomorphisms
on all equivariant homotopy groups in non-negative dimensions.

So far the spectrum $X$ was arbitrary. 
If we now assume that $X$ is globally connective,
then the morphism $u:A\to X$ is a global equivalence
(because $A$ is globally connective by the previous paragraph).
So $X$ is isomorphic in $\GH$ to $A$, and hence $X\in\td{\Pc}_+$.
\end{proof}

Given orthogonal spectra $X$ and $Y$, the external product maps \eqref{eq:def_boxtimes}
\[ \boxtimes \ : \ \pi_0^G (X) \tensor \pi_0^K (Y) \ \to \ \pi_0^{G\times K}(X\sm Y)\]
form a bimorphism of global functors by Theorem \ref{thm:external product properties}. 
The box product of global functors was introduced in
Construction \ref{con:Box product}.
The universal property of the box product produces a morphism of
global functors
\begin{equation}\label{eq:canonical box smash morphism}
 \upi_0(X) \, \Box \, \upi_0 (Y) \ \to  \ \upi_0(X\sm Y) \ .
\end{equation}
We recall from \eqref{eq:derived_smash} that 
the symmetric monoidal derived smash product $\sm^\mL$ 
on the global stable homotopy category is obtained as the total left derived
functor of the smash product of orthogonal spectra.
When applied to flat replacements of $X$ and $Y$, 
the morphism \eqref{eq:canonical box smash morphism}
becomes the morphism of the following proposition.

\begin{prop}\label{prop:Box vs smash} 
For all globally connective spectra $X$ and $Y$ 
the orthogonal spectrum $X\sm^\mL Y$ is globally connective 
and the natural morphism\index{subject}{smash product!derived}\index{subject}{box product!of global functors}
\[ (\upi_0 X) \, \Box \, (\upi_0 Y) \ \to \ \upi_0(X\sm^\mL Y) \]
is an isomorphism of global functors.
\end{prop}
\begin{proof} We fix a compact Lie group $K$ and let $\Xc$ be the class
of globally connective orthogonal spectra $X$ such that 
$X\sm^\mL \Sigma^\infty_+ B_{\gl}K$ is globally connective and the natural morphism
\[ \upi_0 X \, \Box\,  \upi_0(\Sigma^\infty_+ B_{\gl} K)
 \ \to \ \upi_0(X\sm^\mL \Sigma^\infty_+ B_{\gl}K ) \]
is an isomorphism of global functors. 
The class $\Xc$ is closed under sums and we claim that it is also closed under cones.
We let
\[  A \ \to\  B\  \to\  C\   \to\  \Sigma A  \]
be a distinguished triangle in $\GH$ such that $A$ and $B$ belong to $\Xc$.
Since $A$ is globally connective the global functor $\upi_{-1}(A)$ vanishes;
since $A$ belongs to $\Xc$ the global functor 
 $\upi_{-1}(A \sm^\mL \Sigma^\infty_+ B_{\gl} K )$ vanishes.
Since $-\Box\upi_0(\Sigma^\infty_+ B_{\gl} K)$ is right exact
(by Remark \ref{rk:Box is right exact}),
the upper row in the commutative diagram
\[ \xymatrix@C=5mm{ 
(\upi_0 A) \Box\upi_0(\Sigma^\infty_+ B_{\gl} K)  \ar[r]\ar[d] &
 (\upi_0 B) \Box\upi_0(\Sigma^\infty_+ B_{\gl} K) \ar[r]\ar[d] &
(\upi_0 C) \Box \upi_0(\Sigma^\infty_+ B_{\gl} K) \ar[r]\ar[d] &
0  \\
\upi_0(A\sm^\mL \Sigma^\infty_+ B_{\gl} K) \ar[r] &
\upi_0(B\sm^\mL \Sigma^\infty_+ B_{\gl} K)\ar[r] &
\upi_0(C\sm^\mL \Sigma^\infty_+ B_{\gl} K)\ar[r] &
0 } \]
is exact. The lower row is exact because smashing with $\Sigma^\infty_+ B_{\gl}K$ preserves
distinguished triangles. The two left vertical maps are isomorphisms because
$A$ and $B$ belong to the class $\Xc$. So the right vertical map
is an isomorphism and $C$ belongs to $\Xc$ as well.
Moreover, we have 
\begin{align*}
  \upi_0(\Sigma^\infty_+ B_{\gl} G)\, \Box \, \upi_0(\Sigma^\infty_+ B_{\gl} K) \ &\iso \ 
\bA(G,-)\, \Box \, \bA(K,-) \\ 
&\iso \ \bA(G\times K,-) \ \iso \ \upi_0(\Sigma^\infty_+ B_{\gl} (G\times K)) 
\end{align*}
by Proposition \ref{prop:B_gl represents} and Remark \ref{rk:box representable},
and 
\[  (\Sigma^\infty_+ B_{\gl}G)\sm (\Sigma^\infty_+ B_{\gl}K) \ \iso \ 
\Sigma^\infty_+ (B_{\gl}G\boxtimes B_{\gl} K)\ \iso \ \Sigma^\infty_+ B_{\gl}(G\times K) \]
by Proposition \ref{prop:box to smash}, respectively by \eqref{eq:boxtimes_of_B_gl}. 
This shows that the class $\Xc$ is closed under sums and cones
and contains the suspension spectra $\Sigma^\infty_+ B_{\gl}G$ for all
compact Lie groups $G$.
Proposition \ref{prop:plus generating GH}
then shows that $\Xc$ is the class of all globally connective orthogonal spectra.
This proves the proposition in the special case $Y=\Sigma^\infty_+ B_{\gl} K$

Now we perform the same argument in the other variable. 
We fix a globally connective spectrum $X$ and let $\Yc$ denote the class
of globally connective orthogonal spectra $Y$ such that 
$X\sm^\mL Y$ is globally connective and the natural morphism
of the proposition is an isomorphism of global functors. 
The class $\Yc$ is again closed under sums and cones, by the same
arguments as above. 
Moreover, for every compact Lie group $K$
the suspension spectrum $\Sigma^\infty_+ B_{\gl}K$ belongs
to the class $\Yc$ by the previous paragraph. 
Again Proposition \ref{prop:plus generating GH}
shows that $\Yc$ is the class of all globally connective orthogonal spectra.
\end{proof}

The semifree orthogonal spectrum  $F_{G,V}=\bO(V,-)/G$
generated by a $G$-representation $V$
was defined in Construction \ref{con:free orthogonal}.
Now we identify the functor represented by $F_{G,V}$
in the global stable homotopy category. The element
\[ (0,\Id_V)\cdot G \ \in \ \bO(V,V)/G = \ F_{G,V}(V) \]
is a $G$-fixed point, so the $G$-map
\[ S^V \ \to \ \bO(V,V)/G\sm S^V \ = \ (F_{G,V}\sm S^V)(V) \ , \quad
v \ \longmapsto \ (0,\Id_V)\cdot G \sm (-v) \]
represents an equivariant homotopy class
\begin{equation}  \label{eq:define_a_G,G}
a_{G,V}\ \in \ \pi_0^G(F_{G,V}\sm S^V) \ .
\end{equation}
The $\upi_*$-isomorphism of orthogonal $G$-spectra $\lambda^V_E:E\sm S^V\to \sh^V E$ 
was defined in \eqref{eq:defn lambda_n};
its value at an inner product space $U$ is
the opposite structure map $\sigma^{\op}_{U,V}:E(U)\sm S^V \to E(U\oplus V) = (\sh^V E)(U)$.
The following lemma is a direct consequence of
Proposition \ref{prop:lambda upi_* isos}~(i).

\begin{lemma}\label{lemma:normalize a_G,V} 
The morphism $\lambda^V_{F_{G,V}}:F_{G,V}\sm S^V\to \sh^V F_{G,V}$
takes the element $a_{G,V}$ to the class in $\pi_0^G(\sh^V F_{G,V})$ 
that is represented by the $G$-fixed
point $(0,\Id_V)\cdot G$ of $(\sh^V F_{G,V})(0)$.
\end{lemma}

We emphasize that the representative of $a_{G,V}$
involves the involution $S^{-\Id}:S^V\to S^V$,
which represents a certain unit in the ring $\pi_0^G(\mS)$ that squares to~1.
Lemma \ref{lemma:normalize a_G,V} is our reason 
for choosing this particular normalization.
We could remove the involution from the definition of $a_{G,V}$,
thereby changing the class $a_{G,V}$ into $\varepsilon_V(a_{G,V})$,
but then the unit would instead appear in other formulas later on.

\begin{theorem}\label{thm:free spectrum represents} 
Let $G$ be a compact Lie group and $V$ a faithful $G$-represen\-tation.
Then the pair $(F_{G,V}, a_{G,V})$ represents the functor
\[ \GH\to \emph{(sets)}\ , \quad E \ \longmapsto \ \pi_0^G(E\sm S^V)\ . \]
The semifree orthogonal spectrum $F_{G,V}$ is a compact object 
of the global stable homotopy category $\GH$.
\end{theorem}
\begin{proof}
We define
\[ b_{G,V}\ = \ (\lambda^V_{F_{G,V}})_*(a_{G,V})\ \in\ \pi_0^G (\sh^V F_{G,V})\ ; \]
Lemma \ref{lemma:normalize a_G,V} shows that the class $b_{G,V}$
is represented by the $G$-fixed point
\[ (0,\Id_V)\cdot G \ \in \ \bO(V,V)/G\ = \ (\sh^V F_{G,V})(0) \ .\]
Since the morphism $\lambda^V_{F_{G,V}}$ is a $\upi_*$-isomorphism
of orthogonal $G$-spectra,
we may thus show that the pair $(F_{G,V}, b_{G,V})$ represents the functor
\[ \GH\to \text{(sets)}\ , \quad E \ \longmapsto \ \pi_0^G(\sh^V E)\ , \]
i.e., for every orthogonal spectrum $E$ the map
\[ \gh{ F_{G,V}, E } \ \to \ \pi_0^G (\sh^V E)\ , \quad 
f\longmapsto (\sh^V f)_*(b_{G,V})  \]
is bijective. Since both sides take global equivalences in $E$ to bijections,
we can assume that $E$ is a global $\Omega$-spectrum.
The orthogonal spectrum $F_{G,V}$ is flat, and hence cofibrant
in the global model structure.
So the localization functor induces a bijection
\[ \spec(F_{G,V}, E) / \text{homotopy} \ \to \ \gh{ F_{G,V}, E } \] 
from the set of homotopy classes of morphisms of orthogonal spectra
to the set of morphisms in $\GH$.
By the freeness property, morphisms from $F_{G,V}$ to $E$ biject 
with $G$-fixed points of $E(V)$, and homotopies correspond to
paths of $G$-fixed points.
The composite
\[ \pi_0( E(V)^G )  \ \xra{\ \iso\ } \
 \gh{ F_{G,V}, E } \ \xra{f\mapsto (\sh^V f)_*(b_{G,V})} \ \pi_0^G (\sh^V E)  \]
is the stabilization map of the $G$-$\Omega$-spectrum $\sh^V E$, 
and hence bijective.
Since the left map and the composite are bijective, so is the evaluation
map at the class $b_{G,V}$.

Now we prove that $F_{G,V}$ is a compact object in $\GH$.
In the commutative square
\[\xymatrix{  {\bigoplus_{i\in I}} \gh{ F_{G,V},\, X_i } \ar[r]\ar[d] &
\gh{ F_{G,V},\, \bigoplus_{i\in I}X_i } \ar[d] \\
{\bigoplus_{i\in I}} \pi_0^G(X_i\sm S^V) \ar[r] &
\pi_0^G\left( \bigvee_{i\in I}X_i\sm S^V\right) }\]
both vertical maps are evaluation at the class $a_{G,V}$,
which are isomorphisms by part~(i).
The lower horizontal map is an isomorphism by 
Corollary \ref{cor-wedges and finite products}~(i),
hence so is the upper horizontal map.
This shows that $F_{G,V}$ is compact as an object of the triangulated category $\GH$.
\end{proof}

Now we construct certain distinguished triangles in the global stable homotopy
category that arise from special representations of compact Lie groups $G$,
namely when $G$ acts faithfully and transitively on the unit sphere.
Equivalently, we are looking for ways to present spheres as homogeneous spaces.
The distinguished triangles relate two semifree orthogonal spectra for the
compact Lie group $G$ and a semifree orthogonal spectrum for a closed subgroup $H$ 
that occurs as the stabilizer of a unit vector 
in the transitive faithful $G$-representation. 
The main case we care about is $G=O(m)$\index{subject}{orthogonal group}
acting tautologically on $\mR^{m}$ with stabilizer group $H=O(m-1)$.
This special case will become relevant in Section \ref{sec:global Thom}
when we analyze the rank filtration of the global Thom spectrum $\bmO$.
Similarly, we can let $G=S O(m)$\index{subject}{special orthogonal group} act tautologically
on $\mR^{m}$ with stabilizer group $H=S O(m-1)$, and this
shows up in the rank filtration of the global Thom spectrum $\bmSO$; 
here we need $m\geq 2$, since $S O(1)$ does not act transitively on $S(\mR)$.
Other examples are $U(m)$\index{subject}{unitary group} or $S U(m)$\index{subject}{special unitary group}
 (the latter for $m\geq 2$)
acting on the underlying $\mR$-vector
space of $\mC^{m}$, with stabilizer groups $U(m-1)$ respectively $S U(m-1)$. 
Similarly, we can consider the tautological representation of $S p(m)$\index{subject}{symplectic group}
on the underlying $\mR$-vector space of $\mH^{m}$ with stabilizer group $S p(m-1)$.
These examples show up in the rank filtrations of the global Thom spectra $\bmU$,
$\mathbf{mSU}$ and $\mathbf{mSp}$.
There are also more exotic examples, such as the exceptional 
Lie group $G_2$,\index{subject}{G2@$G_2$ exceptional Lie group}
the group of $\mR$-algebra automorphism of the octonions $\mathbb O$.
Here we take $V$ as the tautological 7-dimensional representation
on the imaginary octonions, with stabilizer group isomorphic to $S U(3)$.
For a complete list of examples we refer the reader to \cite[Ch.\,7.B, Ex.\,7.13]{besse}
and the references given therein.

\begin{construction}
We consider a representation $V$ of a compact Lie group $G$
such that $G$ acts faithfully and transitively on the unit sphere $S(V)$.
We choose a unit vector $v\in S(V)$ and let $H$ be the stabilizer group of $v$.
Then the underlying $H$-representation of $V$ decomposes as
\[ V \ = \ L\oplus (\mR\cdot v) \]
where $L$ is the orthogonal complement of $v$.
We use the letter $L$ for this complement because it essentially `is' 
the tangent representation $T_{e H}(G/H)$.
More precisely, the differential at $e H$ of the smooth $G$-equivariant embedding
\[  G/H \ \to \ V \ , \quad g H\ \longmapsto g v \]
is an $H$-equivariant linear isomorphism from $T_{e H}(G/H)$ onto $L$.

We define a based $G$-map
\begin{align}\label{eq:define r}
 r \ : \ S^V \ &\ \to \ \bO(L,V)/H\ = \  F_{H,L}(V)\\
\text{ by\qquad}
 r(g\cdot t\cdot v) & \ = \
g\cdot( (t^2-1)/t\cdot v,\text{incl})\cdot H\ ,\nonumber
\end{align}
where $g\in G$ and $t\in[0,\infty]$.
We let 
\[ T\ :\ \Sigma^\infty_+ B_{\gl} G\ = \ \Sigma^\infty_+ \bL_{G,V}\ \to\  F_{H,L}\]
denote the adjoint of $r$.
We define a morphism of orthogonal spectra
\[ i \ : \ F_{H,L}\ \to \ \sh F_{G,V} \]
as the adjoint of the $H$-fixed point
\[ \psi^{-1}\cdot G\ \in \ \bO(V,L\oplus\mR)/G \ = \ (\sh F_{G,V})(L) \ ,\]
where the $H$-equivariant isometry $\psi:L\oplus\mR\iso V$
is defined as $\psi(x,t)= x+ t\cdot v$. 
The value of the morphism $i$ at an inner product space $W$ is then the map
\begin{align*}
  i(W) \ : \ F_{H,L}(W) = \bO(L,W)/H
\ &\to\ \bO(V,W\oplus \mR)/G \ = \  (\sh F_{G,V})(W) \\
(w,\varphi)\cdot H \qquad &\longmapsto\qquad 
((w,0),(\varphi\oplus\mR)\circ\psi^{-1})\cdot G \ .
\end{align*}
We define a morphism of orthogonal spectra 
$a:F_{G,V}\to \Sigma^\infty_+ \bL_{G,V}=\Sigma^\infty_+ B_{\gl}G$
as the adjoint of the $G$-fixed point
\[ 0 \sm \Id\cdot G\ \in\  S^V\sm \bL(V,V)/G_+ \ = \ 
(\Sigma^\infty_+ \bL_{G,V})(V)\ .\]
\end{construction}

The morphism $\lambda_X:X\sm S^1\to \sh X$ 
was defined in \eqref{eq:defn lambda_n} and is a global equivalence 
by Proposition \ref{prop:global equiv preservation}~(i).
The stable tautological class\index{subject}{stable tautological class}
\[ e_G \ = \ e_{G,V}\  \in \ \pi_0^G(\Sigma^\infty_+ B_{\gl}G )  \]
was defined in \eqref{eq:define_stable_tautological}.
In the next theorem we abuse notation and use the same symbols
for morphisms of orthogonal spectra and their images in the global
stable homotopy category.

\begin{theorem}\label{thm:fundamental triangle}
Let $G$ be a compact Lie group and $V$ a $G$-representation
such that $G$ acts faithfully and transitively on the unit sphere $S(V)$;
let $H$ be the stabilizer group of a unit vector of $V$.
\begin{enumerate}[\em (i)]
\item 
The sequence
\[ F_{G,V} \ \xra{\ a\ }\ \Sigma^\infty_+ B_{\gl} G\ \xra{\ T \ } \
F_{H,L}\ \xra{-\lambda_{F_{G,V}}^{-1}\circ i \ } \  F_{G,V}\sm S^1 \]
is a distinguished triangle in the global stable homotopy category.
\item
The morphism $T$ satisfies the relation
\[ T_*(e_{G}) \ = \ \Tr_H^G(a_{H,L})\]
in the group $\pi_0^G(F_{H,L})$.
\end{enumerate}
\end{theorem}
\begin{proof}
We consider the based map $q:G/H_+\to S^0$ that sends $G/H$
to the non-basepoint.
We identify the mapping cone of $q$
with the sphere $S^{V}$ via the $G$-equivariant homeomorphism
\[ h \ : \ C q\ \iso \ S^{V} \]
that is induced by the map
\[ G/H\times[0,1] \ \to \ S^V \ , \quad
(g  H,x) \ \longmapsto\ g \cdot (1-x)/x \cdot v \ .\]
Under this identification the mapping cone inclusion $i:S^0\to C q$ becomes the
inclusion $\iota:\{0,\infty\}=S^0\to S^{V}$,
and the projection $p:C q\to G/H_+\sm S^1$ becomes the map
\[
 p\circ h^{-1}\ : \ S^V\ \to \ G/H_+\sm S^1 \ ,\quad
g\cdot t\cdot v \ \longmapsto \ g H\sm (1-t^2)/t\ ,  
  \]
where $g\in G$ and $t\in [0,\infty]$.
The semifree functor $F_{G,V}$ takes mapping cone sequences of
based $G$-spaces to mapping cone sequences of orthogonal spectra, so the sequence
\[ F_{G,V}(G/H_+)\ \xra{F_{G,V} q}\ F_{G,V} \ \xra{F_{G,V}\iota} \ 
F_{G,V}S^{V}\ \xra{F_{G,V}(p\circ h^{-1})}\ 
F_{G,V}( G/H_+\sm S^1)  \]
is a distinguished triangle in the global stable homotopy category.

We define a morphism $\Lambda:F_{G,V}( G/H_+\sm S^1)\to F_{H,L}$
as the adjoint of the $G$-map
\begin{align*}
 G/H_+\sm S^1\ &\to \ \bO(L,V)/H \ = \ F_{H,L}(V) \\
g  H\ \sm t \quad &\longmapsto \ g \cdot ( - t\cdot v,\, \text{incl} )\cdot H \ .  
\end{align*}
The morphism $\Lambda$ is the composite of an isomorphism
$F_{G,V}( G/H_+\sm S^1)\iso F_{H,L\oplus \mR}S^1$
and the morphism $\lambda_{H,L,\mR}:F_{H,L\oplus\mR}S^1\to F_{H,L}$
defined in \eqref{eq:define_lambda}. The morphism
$\lambda_{H,L,\mR}$ is a global equivalence
by Theorem \ref{thm:faithful independence}, hence $\Lambda$ is 
a global equivalence as well. The composite
\[ 
\Sigma^\infty_+ B_{\gl} G\ \xra{\text{untwist}^{-1}}\ 
F_{G,V} S^V\ \xra{F_{H,L}(p\circ h^{-1})} \ F_{G,V}(G/H_+\sm S^1)
\ \xra{\ \Lambda\ }\ F_{H,L} \]
coincides with the morphism $T$, by direct inspection
of the effects at the inner product space $V$.

Our next claim is that the following diagram 
of orthogonal spectra commutes up to homotopy:
\begin{equation}\begin{aligned}\label{eq:F_{H,L} square}
\xymatrix@C=13mm{ 
F_{G,V}(G/H_+)\sm S^1\ar[r]^-{(F_{G,V}q)\sm S^1}\ar[d]_\Lambda & 
F_{G,V}\sm S^1\ar[d]^{\lambda_{F_{G,V}}}\\
F_{H,L} \ar[r]_-i & \sh F_{G,V}} 
\end{aligned}\end{equation}
The representing property of $F_{G,V}$ reduces this to evaluating 
the square at $V$ and showing that the two resulting 
$G$-maps from $G/H_+\sm S^1$ to $(\sh F_{G,V})(V)$
are equivariantly homotopic. Since $G$-maps out of
$G/H$ correspond to $H$-fixed points,
this in turn reduces to the claim that the two maps
\[ S^1 \ \to \ \left( \bO(V,V\oplus\mR)/G\right)^{H}
\ = \ ( (\sh F_{G,V})(V))^{H} \]
that send $t\in S^1$ to
\[ ( (-t\cdot v,0), i_0)\cdot G \text{\qquad respectively\qquad} 
( (0,t), i_1)\cdot G  \]
are based homotopic, where $i_0,i_1:V\to V\oplus\mR$ are given by 
\[  i_0(x)\ = \ (x- \td{x,v}\cdot v,\, \td{x,v}) 
\text{\quad respectively\quad}
i_1(x)\ = \  (x,0)\ . \]
A homotopy that witnesses this is
\begin{align*}
  S^1\times [0,1] \ &\to \quad  (\bO(V,V\oplus\mR)/G)^H \\
(t,s)\quad &\longmapsto \ 
( (-t\cdot \cos(\pi s/2)\cdot v, t\cdot \sin(\pi s/2)),\, i_s)\cdot G\ ,
\end{align*}
with
\[ i_s \ : \  V\ \to \ V\oplus\mR \ , \
i_s(x)\ = \ (x + ( \sin(\pi s/2)-1)\cdot\td{x,v}\cdot v, \cos(\pi s/2)\cdot\td{x,v})\ .\]
This shows the claim that the square \eqref{eq:F_{H,L} square} 
is homotopy commutative.

So altogether we have produced a commutative diagram 
in the global stable homotopy category:
\[ \xymatrix@C=13mm{
F_{G,V} \ar@{=}[d] \ar[r]^-{F_{G,V}\iota} & F_{G,V} S^{V} 
\ar[d]_{\text{untwist}}^\iso \ar[r]^-{F_{G,V}(p\circ h^{-1})} &
F_{G,V}(G/H_+\sm S^1) \ar[d]^{\Lambda}_\sim\ar[r]^-{-(F_{G,V}q\sm S^1)} & 
F_{G,V} \ar[d]^{\lambda_{F_{G,V}}}_\sim\sm S^1\\
F_{G,V} \ar[r]_-a & \Sigma^\infty_+ B_{\gl}G\ar[r]_-T &
F_{H,L} \ar[r]_-{-i}& \sh F_{G,V} } \]
The upper sequence is obtained by rotating a distinguished triangle,
so it is distinguished.
Since all vertical morphisms are isomorphisms in $\GH$, we conclude that
the sequence $(a,T,-\lambda^{-1}_{F_{G,V}}\circ i)$ is a distinguished triangle in $\GH$.

It remains to analyze the effect of the morphism $T$ on the
stable tautological class $e_G$.
We consider the wide $G$-equivariant embedding
\[ G/H \ \to \ V \ , \quad g  H\ \longmapsto \ g v\ .\]
This embedding was already used to identify the tangent space $T_{e H}(G/H)$
at the preferred coset with the subspace $L$ inside $V$;
the inclusion $L\to V$ corresponds to the differential at $e H$.
The orthogonal complement of $L$ is precisely the line 
spanned by the vector $v$, and the subgroup $H$ acts trivially on this complement.

As described more generally in Construction \ref{con:external transfer},
the wide embedding gives rise to a
$G$-equivariant Thom-Pontryagin collapse map \eqref{eq:TPcollaps}
\[ c \ : \ S^V\ \to\  G/H_+\sm S^1 \ .\] 
In the more general context of an arbitrary closed subgroup, 
the target of the collapse map is $G\ltimes_H S^{V-L}$; 
but in our situation $H$ acts trivially on $\mR\cdot v=V-L$, 
so we can (and will) identify $V-L$ with $\mR$ via $x\mapsto \td{x,v}$.
The collapse map $c$ sends the origin and all vectors of length at least~2 to
the basepoint; if $x\in V$ satisfies $0<|x|<2$, then 
$x=g \cdot t\cdot v $ for some $g\in G$ and $t\in(0,2)$.
The collapse map $c$ sends such a vector $x$ to
\[ g  H\sm  \frac{t-1}{1-|t-1|}\ \in \ G/H_+\sm S^1\  .\]

Now we calculate the class $\Tr_H^G(a_{H,L})$.
We factor the external transfer $G\ltimes_H-:\pi_0^G(E\sm S^L)\to\pi_0^G(G\ltimes_H E)$
as the composite
\[ \pi_0^H(E\sm S^L)\ \xra{ (\lambda^L_E)_*}\ 
\pi_0^H(\sh^L E) \ \xra{G\boxtimes_H -} \ \pi_0^G(G\ltimes_H E)\ .\]
Here $\lambda^L_E:E\sm S^L\to\sh^L E$ is the morphism of orthogonal $H$-spectra
defined in \eqref{eq:defn lambda_n};
the second homomorphism is a variation of the external transfer, 
and is defined as follows.
We let $U$ be an $H$-representation and $f:S^U \to E(U\oplus L)=(\sh^L E)(U)$ 
an $H$-equivariant 
based map that represents a class in $\pi_0^H(\sh^L E)$.
By enlarging $U$, if necessary, we can assume that it
is the underlying $H$-representation of a $G$-representation.
We denote by $f\diamond\mR$ the composite $H$-map
\[ S^U\sm S^1 \ \xra{f\sm S^1} \ E(U\oplus L)\sm S^1\ \xra{\sigma^{\op}_{U\oplus L,\mR}} \
 E(U\oplus L\oplus\mR)\ \xra[\iso]{E(U\oplus\psi)} \
 E(U\oplus V) \ . \]
The class $G\boxtimes_H\td{f}$ in $\pi_0^G(G\ltimes_H E)$
is then represented by the composite $G$-map
\begin{align*}
  S^{U\oplus V} \ \xra{\ S^U\sm c\ }\  
&S^U\sm (G/H_+\sm S^1) \ \xra[\iso]{\text{shear}}\   
G\ltimes_H( S^U\sm S^1) \\ 
&\xra{G\ltimes_H(f\diamond\mR)}\   
G\ltimes_H  E(U\oplus V)\ \iso \ (G\ltimes_H E)(U\oplus V) \ .
\end{align*}
The relation $G\boxtimes_H((\lambda^L_E)_*\td{f})=G\ltimes_H \td{f}$
is straightforward from the definitions.

By Lemma \ref{lemma:normalize a_G,V}, the class
$(\lambda_{F_{H,L}}^L)_*(a_{H,L})$
is represented by the $H$-fixed point
\[ (0,\Id_L)\cdot H \ \in \  \bO(L,L)/H \ = \ (\sh^L F_{H,L})(0) \ . \]
So the class $\Tr_H^G(a_{H,L})=\text{act}_*(G\boxtimes_H((\lambda^L_{F_{H,L}})_*(a_{H,L})))$
is represented by the composite $G$-map
\begin{align*}
  S^V \ &\xra{\ c\ }\  
G/H_+\sm  S^1  \ 
\xra{G\ltimes_H( (0,\Id_L)H \diamond\mR)}\  G\ltimes_H ( \bO(L,V)/H) \ \xra{\text{act}}\ 
\bO(L,V)/H \ .
\end{align*}
Expanding all the definitions identifies this composite 
as the map
\[ S^V \ \to \ \bO(L,V)/H \ , \quad 
g\cdot t\cdot v \ \longmapsto \ g\cdot(\xi(t)\cdot v,\text{incl})\cdot H\ ,\]
for $g\in G$ and $t\in [0,\infty]$, with
\[ \xi\ : \ [0,\infty]\ \to \ S^1 \ ,
 \quad \xi(t)\ = \ 
\begin{cases}
 \frac{t-1}{1-|t-1|}
& \text{ if $t\in (0,2)$, and}\\   
\qquad \ast & \text{ else.}
\end{cases}\] 
The function $\xi$ is homotopic, relative to $\{0,\infty\}$, 
to the function sending $t$ to $(t^2-1)/t$. 
Any relative homotopy between these two functions
induces a based $G$-equivariant homotopy between the representative
of the class $\Tr_H^G(a_{H,L})$ and the map $r$
defined in \eqref{eq:define r}. So $r$ itself is a representative
of the class $\Tr_H^G(a_{H,L})$.
Since $T:\Sigma^\infty_+ B_{\gl}G\to F_{H,L}$ is adjoint to the $G$-map $r$,
we deduce the desired relation $\Tr_H^G(a_{H,L})=[r]=T_*(e_G)$.
\end{proof}

Part~(ii) of Theorem \ref{thm:fundamental triangle}
can be interpreted as saying that for global stable homotopy types,
the morphism of orthogonal spectra
$T:\Sigma^\infty_+ B_{\gl}G\to F_{H,L}$ represents the dimension shifting transfer
from $H$ to $G$.
Indeed, the pair $(\Sigma^\infty_+ B_{\gl}G, e_G)$ represents
the functor $\pi_0^G$ (by Theorem \ref{thm:Bgl G weak generators}~(i)), 
and the pair $(F_{H,L},a_{H,L})$ 
represents the functor $\pi_0^H(-\sm S^L)$
(by Theorem \ref{thm:free spectrum represents}).
Naturality of transfers then promotes the relation $T_*(e_G)=\Tr_H^G(a_{H,L})$
to a commutative diagram:
\[ \xymatrix@C=15mm{ 
\gh{F_{H,L}, E }\ar[r]^-{\gh{T,E}} \ar[d]_{f\mapsto f_*(a_{H,L})}^\iso & 
\gh{\Sigma^\infty_+ B_{\gl}G, E } \ar[d]^{f\mapsto f_*(e_G)}_\iso \\
\pi_0^H(E\sm S^L)\ar[r]_-{\Tr_H^G } & \pi_0^G( E) } \]

\section{Change of families}\label{sec:change family}

In this section we compare the global stable homotopy categories for two
different global families $\Fc$ and $\Ec$, where we suppose that $\Fc\subseteq \Ec$. 
Then every $\Ec$-equivalence is also an $\Fc$-equivalence, 
so we get a `forgetful' functor on the homotopy categories
\[ U = U^\Ec_\Fc\ : \ \GH_\Ec \ \to \ \GH_\Fc \]
from the universal property of localizations.
The global model structures are stable, so the two global homotopy categories
$\GH_\Ec$ and $\GH_\Fc$ have preferred triangulated structures,
and the forgetful functor is canonically an exact functor of
triangulated categories. 
We show in Theorem \ref{thm:change families} that this forgetful functor
has a left and a right adjoint, both fully faithful,
and that this data is part of a recollement of triangulated categories.
If the global families $\Ec$ and $\Fc$ are multiplicative, 
then the smash product of orthogonal spectra can be derived to
symmetric monoidal products on $\GH_\Ec$ and on $\GH_\Fc$ 
(see Corollary \ref{cor-derived smash}). The forgetful functor
is strongly monoidal with respect to these derived smash products.
Indeed, the derived smash product in $\GH_\Ec$ can be calculated by
flat approximation up to $\Ec$-equivalence;
every $\Ec$-equivalence is also an $\Fc$-equivalence, so these 
flat approximations can also be used to calculate 
the derived smash product in $\GH_\Fc$.
Theorem \ref{thm:change families} below also exhibits symmetric monoidal structures
on the two adjoints of the forgetful functor.

The special case with $\Fc$ the trivial global family
allows a quick calculation of the Picard group of $\GH_{\Ec}$,
with the possibly disappointing answer that it is free abelian of rank~1, 
generated by the suspension of the global sphere spectrum 
(see Theorem \ref{thm:Picard of GH}).
Propositions \ref{prop:reflexive} and \ref{prop:right induced criterion} 
provide characterizations of the global homotopy types in the image of the 
left and right adjoints.
As we explain in Example \ref{eg:global Borel},
the `absolute' right adjoint from the non-equivariant to the
full global stable homotopy category models Borel cohomology theories.
Construction \ref{cor:Borel b} exhibits a particularly nice
lift of the global Borel functor to the pointset level,
i.e., a lax symmetric monoidal endofunctor of the category of orthogonal spectra.

In Theorem \ref{thm:change to groups} we relate 
the global stable homotopy category to the $G$-equivariant stable homotopy
category $G\text{-}\SH$ for a fixed compact Lie group $G$.
There is another forgetful functor $\GH\to G\text{-}\SH$
which is an exact functor of triangulated categories and
has both a left adjoint and a right adjoint.
Here, however, the adjoints are {\em not} fully faithful as soon
as the group $G$ is non-trivial.

The global family $\Fin$ of finite groups is an important example
to which the discussion of this section applies.
We will show that rationally, the associated $\Fin$-global stable homotopy category
admits an algebraic model: Theorem \ref{thm:rational SH}
provides a chain of Quillen equivalences between
the category of orthogonal spectra with the 
rational $\Fin$-global model structure and the category of chain complexes of
rational global functors on finite groups.
Under this equivalence, the homotopy group global functor for spectra
corresponds to the homology group global functor for complexes.

The algebraic model can be simplified further.
Rational $G$-Mackey functors for a fixed finite group
naturally split into contributions indexed by conjugacy classes of subgroups,
see Proposition \ref{prop:recover G-Mackey rationally}.\index{subject}{Mackey functor}
While the analogous global context is {\em not} semisimple, the
category of rational $\Fin$-global functors is equivalent
to a simpler one, namely contravariant functors 
from the category $\Out$ of finite groups and
conjugacy classes of epimorphisms to $\mQ$-vector spaces,
by Theorem \ref{thm:rational global functors} below.
So rational $\Fin$-global stable homotopy theory is also modeled
by the complexes of such functors. Under the composite equivalence,
the geometric fixed point homotopy group functors for spectra
corresponds to the homology group $\Out$-functors for complexes,
see Corollary \ref{cor:rational geometric fixed}.

\medskip

\begin{theorem}\label{thm:change families}
  Let $\Fc$ and $\Ec$ be two global families such that $\Fc\subseteq \Ec$. 
  \begin{enumerate}[\em (i)]
  \item The forgetful functor
    \[ U \ : \ \GH_\Ec \ \to \ \GH_\Fc \]
    has a left adjoint $L$ and a right adjoint $R$, and both adjoints are fully faithful.
  \item  If the global families $\Ec$ and $\Fc$ are multiplicative,\index{subject}{multiplicative global family}  
    then the right adjoint has a preferred lax symmetric monoidal structure
    $( R A ) \sm^\mL_\Ec ( R B ) \to R(  A \sm_\Fc^\mL  B )$.\index{subject}{derived smash product}     \item  If the global families $\Ec$ and $\Fc$ are multiplicative, then
    the left  adjoint has a preferred strong symmetric monoidal structure, i.e., 
    a natural isomorphism between $L(  A \sm_\Fc^\mL  B )$ and 
    $( L A ) \sm^\mL_\Ec ( L B )$.
  \end{enumerate}
\end{theorem}
\begin{proof}
By Proposition \ref{prop:sums and product in GH} 
the categories $\GH_\Ec$ and $\GH_\Fc$ have infinite sums
and infinite products, and the forgetful functor preserves both.

(i) As spelled out in Corollary \ref{cor-generators for GH_F}~(v),
the existence of the left adjoint is a formal consequence
of the fact that $\GH_\Ec$ is compactly generated 
and that the functor $U$ preserves products.
But instead of arguing by hand that $U$ preserves products,
we give an alternative construction 
of the left adjoint by model category theory.
Indeed, it is immediate from the definitions of $\Fc$-equivalences
and $\Fc$-global fibrations that the identity functor is a right Quillen functor 
from the $\Ec$-global to the $\Fc$-global model structure.
So its derived functor has a (derived) left adjoint $L:\GH_\Fc\to \GH_\Ec$, 
see for example
\cite[I.4, Thm.\,3]{Q} or \cite[Lemma 1.3.10]{hovey-book}.
Since the right Quillen functor (i.e., the identity) 
preserves all weak equivalences, the adjunction 
unit $A\to U(L A)$ is an isomorphism in the $\Fc$-global homotopy category. 
So the left adjoint is fully faithful.

By Proposition \ref{prop:sums and product in GH} sums in
both $\GH_\Ec$ and $\GH_\Fc$ are represented by the wedge of orthogonal spectra. 
So the forgetful functor $U$ preserves sums.
As spelled out in Corollary \ref{cor-generators for GH_F}~(iv),
the existence of the right adjoint is then a formal consequence
of the fact that $\GH_\Ec$ is compactly generated. 

The fact that $R$ is fully faithful is also a formal consequence of properties
of the adjoint pair $(L,U)$. As we saw above,
the adjunction unit $\eta_A:A\to U(L A)$ is an isomorphism
for every object $A$ of $\GH_\Fc$.
So the map 
\[ (\eta_A)^* \ : \ \gh{U(L A),X}_\Fc\ \to \ \gh{A,X}_\Fc \]
is bijective for every object $X$ of $\GH_\Fc$. The adjunction $(U,R)$ lets
us rewrite the left hand side as $\gh{L A, R X}_\Ec$,
and the adjunction $(L,U)$ lets
us rewrite this further to $\gh{A, U(R X)}_\Fc$.
Under these substitutions, the map $(\eta_A)^*$ becomes the map
induced by the adjunction counit $\epsilon_X:U(R X)\to X$.
This adjunction counit is thus a natural isomorphism, and so $R$
is also fully faithful.

(ii) The lax monoidal structure of the right adjoint $R$ is a formal
consequence of the strong monoidal structure of the forgetful functor $U$.
Indeed, for every pair of orthogonal spectra $A$ and $B$
this strong monoidal structure and the adjunction counits provide a morphism
\[ U \left( ( R A ) \sm^\mL_\Ec ( R B ) \right) \ \iso \  U (R A ) \sm_\Fc^\mL U (R B) 
\ \xra{\epsilon_A\sm^\mL \epsilon_B }\  A \sm_\Fc^\mL  B \]
whose adjoint $( R A ) \sm^\mL_\Ec ( R B ) \to R(  A \sm_\Fc^\mL  B )$
is associative, commutative and unital.

(iii) 
The strong symmetric monoidal structure on $U$ and the adjunction units
provide a morphism
\[ A \sm^\mL_\Fc B\  \xra{\eta_A \sm^\mL \eta_B }\ U( L A) \sm^\mL_\Fc U ( L B )
\ \iso \ U \left( ( L A) \sm^\mL_\Ec  ( L B ) \right) \]
whose adjoint $\lambda_{A,B}:L(  A \sm_\Fc^\mL  B ) \to ( L A ) \sm^\mL_\Ec ( L B ) $
is associative, commutative and unital.
We claim that the morphism  $\lambda_{A,B}$ is in fact an isomorphism in $\GH_\Ec$
(and hence it can be turned around).
We can assume that $A$ and $B$ are $\Fc$-cofibrant so that $LA=A$ and $LB=B$. 
Since $\Fc$ is multiplicative, the pointset level smash product $A\sm B$
is again $\Fc$-cofibrant by Proposition \ref{prop:ExF ppp}~(i).
So the value of the left adjoint $L$ on $A\sm B$ is also given by $A\sm B$.
\end{proof}

\begin{rk}[Recollements]
Theorem \ref{thm:change families} implies that for all pairs of
nested global families $\Fc\subseteq \Ec$ the diagram
 of triangulated categories and exact functors
\[\xymatrix@C=15mm{
\GH(\Ec;\Fc) \ar@<-.3ex>[r]^-{i_*}  & 
\GH_\Ec\ \ar@<-.3ex>[r]^-U  \ \ar@<.4ex>@/^1pc/[l]^-{i^!} \ar@<-.4ex>@/_1pc/[l]_-{i^*} & 
\GH_\Fc \ \ar@<.4ex>@/^1pc/[l]^-R \ar@<-.4ex>@/_1pc/[l]_-L }
 \]
is a {\em recollement}\index{subject}{recollement} in the sense of \cite[Sec.\,1.4]{BBD}.
Here $\GH(\Ec;\Fc)$ denotes the `$\Ec$-global homotopy category with
support outside $\Fc$', i.e., the full subcategory of $\GH_\Ec$
of spectra all of whose $\Fc$-equivariant homotopy groups vanish.
The functor $i_*:\GH(\Ec;\Fc)\to\GH_\Ec$ is the inclusion, and
$i^*$ (respectively $i^!$) is a left adjoint (respectively right adjoint) of $i_*$.
\end{rk}

\begin{rk}
In Theorem \ref{thm:change families}
the left adjoint $L:\GH_\Fc\to \GH_\Ec$ of the forgetful functor $U:\GH_\Ec \to \GH_\Fc$
is obtained as the total left derived functor of the identity functor
on orthogonal spectra with respect to change of model structure
from the $\Fc$-global to the  $\Ec$-global model structure.
So one can calculate the value of the left adjoint $L$ on an orthogonal spectrum $X$ 
by choosing any cofibrant replacement in the  $\Fc$-global model structure,
i.e., an $\Fc$-equivalence $X_\Fc\to X$ with $\Fc$-cofibrant source.
The global homotopy type (so in particular the $\Ec$-global homotopy type) 
of $X_\Fc$ is then well-defined.
Indeed, since the identity is a left Quillen functor
from the $\Fc$-global to the global model structure
(the special case $\Ec=\All$),
every acyclic cofibration in the $\Fc$-global model structure
is a global equivalence. Ken Brown's lemma 
(see the proof of \cite[I.4 Lemma 1]{brown-abstract},
\cite[Lemma 1.1.12]{hovey-book} or \cite[Cor.\,7.7.2]{hirschhorn}) then implies that
every $\Fc$-equivalence between $\Fc$-cofibrant orthogonal spectra
is a global equivalence. 

It seems worth spelling out the extreme case
when $\Fc=\td{e}$ is the minimal global family of trivial groups and
when $\Ec=\All$ is the maximal global family of all compact Lie groups:
An orthogonal spectrum $X$ is $\td{e}$-cofibrant if for every
$m\geq 0$ the latching morphism $L_m X\to X(\mR^m)$ is an $O(m)$-cofibration
and $O(m)$ acts freely on the complement of the image. These are precisely the
orthogonal spectra called `q-cofibrant' in \cite{mmss}.
So every non-equivariant stable equivalence between 
q-cofibrant orthogonal spectra is a global equivalence.
\end{rk}

The formal properties of the change-of-family functors 
established in Theorem \ref{thm:change families} facilitate an
easy and rather formal argument to identify the Picard group
of the global stable homotopy category. The result is that there
are no `exotic' invertible objects, i.e., the only smash invertible
objects of $\GH$ are the suspensions and desuspensions of the global sphere spectrum.
The same is true more generally for the $\Fc$-global stable homotopy
category relative to any multiplicative global family $\Fc$, 
see Theorem \ref{thm:Picard of GH} below.

We recall that an object $X$ of a monoidal category $(\Cc,\Box,I)$
is {\em invertible} if there is another object $Y$ such that both
$X\Box Y$ and $Y\Box X$ are isomorphic to the unit object $I$.
If the isomorphism classes of invertible objects form a set, 
then the {\em Picard group}\index{subject}{Picard group}
$\Pic(\Cc)$ is this set, with the group structure induced
by the monoidal product.
Every strong monoidal functor between  monoidal categories takes invertible objects
to invertible objects, and thus induces a group homomorphism of Picard groups.

Part of the calculation of $\Pic(\GH)$ involves a very general argument that we spell
out explicitly. 

\begin{prop}\label{prop:Picard idempotent}
Let $(\Cc,\Box,I)$ be a monoidal category,
$P:\Cc\to\Cc$ a strong monoidal functor  
and $\epsilon:P\to\Id_\Cc$ a monoidal transformation
such that $\epsilon_I:P I\to I$ is an isomorphism.
Then the induced endomorphism
\[ \Pic(P)\ : \ \Pic(\Cc)\ \to \ \Pic(\Cc) \]
is the identity.
\end{prop}
\begin{proof}
We let $X$ be any invertible object, and $Y$ an inverse of $X$.  
Since $\epsilon$ is a monoidal transformation, the composite
\[ (P X)\Box (P Y) \ \xra{\ \iso \ }\ 
P ( X\Box Y) \ \xra{\epsilon_{X\Box Y}} \ X\Box Y \]
agrees with the morphism $\epsilon_X\Box\epsilon_Y$,
where the first isomorphism is the strong
monoidal structure on $P$.
Since $X\Box Y$ is isomorphic to the unit object and $\epsilon_I$ is
an isomorphism, the morphism
$\epsilon_{X\Box Y}$ is also an isomorphism.
So the composite $\epsilon_X\Box \epsilon_Y$ is an isomorphism. Since
\[  \epsilon_X\Box \epsilon_Y\ = \ (X\Box\epsilon_Y)\circ (\epsilon_X\Box P Y)
\ = \ (\epsilon_X\Box  Y)\circ (P X\Box\epsilon_Y)\]
we conclude that  $\epsilon_X\Box P Y$ has a left inverse 
and $\epsilon_X\Box  Y$ has a right inverse. 
Since $Y$ and $P Y$ are invertible objects, the functors
$-\Box Y$ and $-\Box P Y$ are equivalences of categories, so 
already $\epsilon_X:P X\to X$ has both a left and a right inverse. 
So $\epsilon_X$ is an isomorphism, and hence $\Pic(P)[X]=[P X]=[X]$.
\end{proof}

Now we have all necessary ingredients to determine the Picard
group of the $\Fc$-global stable homotopy category.

\begin{theorem}\label{thm:Picard of GH}
  For every multiplicative global family $\Fc$, 
  the Picard group of the $\Fc$-global stable 
  homotopy category is free abelian of rank~1, generated by the suspension
  of the global sphere spectrum.\index{subject}{Picard group!of the $\Fc$-global stable homotopy category}
\end{theorem}
\begin{proof}
  The forgetful functor 
  \[ U\ = \ U^\Fc_{\td{e}}\ :\ \GH_\Fc \ \to \ \GH_{\td{e}}\ = \ \SH \]
  to the non-equivariant stable homotopy category and its left adjoint
  \[ L\ :\ \SH \ \to \ \GH_\Fc \]
  both have strong monoidal structures, 
  the latter by Theorem \ref{thm:change families}~(iii).
  Moreover, the adjunction counit $\epsilon:L U\to \Id$ is a monoidal transformation,
  and the morphism $\epsilon_\mS:L (U\mS)\to\mS$ is an isomorphism.
  We apply Proposition \ref{prop:Picard idempotent} to the composite endofunctor $L U$
  of $\GH_\Fc$ and conclude that the composite homomorphism
  \[  \Pic(\GH_\Fc) \ \xra{\Pic(U)} \ \Pic(\SH) \ \xra{\Pic(L)} \ \Pic(\GH_\Fc) \]
  is the identity. In particular, the homomorphism $\Pic(U)$ induced by
  the forgetful functor is injective.
  The Picard group of the non-equivariant stable homotopy category 
  is free abelian of rank~1, generated by the suspension 
  of the non-equivariant sphere spectrum. This generator is the image 
  of the suspension of the global sphere spectrum. So the homomorphism
  $\Pic(U)$ is surjective, hence an isomorphism. So $\Pic(\GH_\Fc)$ is
  also free abelian of rank~1, generated by the suspension 
  of the global sphere spectrum.
  \end{proof}

Now we develop criteria that characterize global homotopy
types in the essential image of one of the adjoints to a
forgetful change-of-family functor.
The following  terminology is convenient here.

\begin{defn}\label{def:left_right_induced}
Let $\Fc$ be a global family. An orthogonal spectrum is
{\em left induced}\index{subject}{left induced} 
from $\Fc$ if it is in the essential image of the left adjoint
$L_\Fc:\GH_\Fc\to\GH$. Similarly, an orthogonal spectrum is 
{\em right induced}\index{subject}{right induced}  
from $\Fc$ if it is in the essential image of the right adjoint $R_\Fc:\GH_\Fc\to\GH$.   
\end{defn}

We start with a criterion, for certain `reflexive' global families,
that characterizes the left induced homotopy types in terms of geometric fixed points.

\begin{defn}
A global family $\Fc$ is {\em reflexive}\index{subject}{reflexive global family}\index{subject}{global family!reflexive|see{reflexive global family}}
if for every compact Lie group $K$ 
there is a compact Lie group $u K$, belonging to $\Fc$,
and a continuous homomorphism $p:K\to u K$ that is initial among
continuous homomorphisms from $K$ to groups in $\Fc$.  
\end{defn}

In other words, $\Fc$ is reflexive if and only if
the inclusion into the category of all compact Lie groups has a left adjoint.
As always with adjoints, the universal pair $(u K,p)$ 
is then unique up to unique isomorphism under $K$. 
Moreover, the universal homomorphism $p:K\to u K$ is necessarily surjective. 
Indeed, the image of $p$ is a closed subgroup of $u K$,
hence also in the global family $\Fc$. So if the image of $p$ were strictly smaller
than $K$, then $p$ would not be initial among morphisms into groups from $\Fc$.
Some examples of reflexive global families are
the minimal global family $\td{e}$ of trivial groups, 
the global family $\Fin$ of finite groups and
the global family of abelian compact Lie groups.
The maximal family of all compact Lie groups is also reflexive, 
but in this case the following proposition has no content.

A reflexive global family $\Fc$ is in particular 
multiplicative.\index{subject}{multiplicative global family}
Indeed, for $G,K\in \Fc$ the projections $p_G:G\times K\to G$
and $p_K:G\times K\to K$ factor through continuous homomorphisms
$q_G:u(G\times K)\to G$ respectively $q_K:u(G\times K)\to K$.
The composite
\[ G\times K \ \xra{\ p\ }\ u(G\times K)\ \xra{(q_G,q_K)}\ G\times K \]
is then the identity, so the universal homomorphism
$p:G\times K\to u(G\times K)$ is injective. Since $u(G\times K)$
belongs to $\Fc$, so does $G\times K$.

\begin{prop}\label{prop:reflexive}
Let $\Fc$ be a reflexive global family.
Then an orthogonal spectrum $X$ is left induced from $\Fc$ if and only if
for every compact Lie group $K$ the inflation map
\[ p^* \ : \ \Phi_*^{u K}( X )\ \to \ \Phi_*^K ( X ) \]
associated to the universal morphism $p:K\to u K$ is an isomorphism 
between the geometric fixed point homotopy groups for $u K$ and $K$.
\end{prop}
\begin{proof}
We let $\Xc$ be the full subcategory of $\GH$ consisting
of the orthogonal spectra $X$ such that for every compact Lie group $K$
the inflation map $p^*: \Phi_*^{u K}( X )\to \Phi_*^K ( X )$ is an isomorphism.
We need to show that $\Xc$ coincides with the class of spectra left induced from $\Fc$.

Geometric fixed point homotopy groups commute with sums and take exact triangles
to long exact sequences. So $\Xc$ is closed under sums and triangles, i.e.,
it is a localizing subcategory of the global homotopy category. 
Now we claim that for every group $G$ in $\Fc$
the suspension spectrum of the global classifying space $B_{\gl} G$
belongs to $\Xc$. 
Since  $p:K\to u K$ is initial among morphisms into groups from $\Fc$,
precomposition with $p$ is a bijection between the sets of
conjugacy classes of homomorphisms into $G$; moreover, the image of a homomorphism
$\alpha:u K\to G$ agrees with the image of $\alpha\circ p:K\to G$, 
because $p$ is surjective. 
Proposition \ref{prop:fix of global classifying}~(i) 
identifies the fixed points of the orthogonal space $B_{\gl} G$
as a disjoint union, over conjugacy classes of homomorphisms,
of centralizers of images.
So the restriction map along $p$ is a weak equivalence of fixed points spaces
\[ p^* \ : \ \left((B_{\gl} G)(\Uc_{u K}) \right)^{u K} \ = \ 
\left((B_{\gl} G)(\Uc_{u K}) \right)^K \ \simeq \ 
\left((B_{\gl} G)(\Uc_K) \right)^K \ .\]
Geometric fixed points commute with suspension spectra
(see Example \ref{eg:geometric and suspension}), 
in the sense of an isomorphism
\begin{align*}
 \Phi^K_* (\Sigma^\infty_+ B_{\gl} G)\ &\iso \ 
\pi_*^e \left( \Sigma^\infty_+ \left( (B_{\gl} G)(\Uc_K)\right)^K\right) \ ,
\end{align*}
natural for inflation maps.
So together this implies the claim for the suspension spectrum of $B_{\gl} G$.

Now we have shown that $\Xc$ is a localizing subcategory of the global
stable homotopy category that contains the suspension spectra 
of global classifying spaces of all groups in $\Fc$.
The left adjoint $L:\GH_\Fc\to \GH$ is fully faithful and 
$\GH_\Fc$ is generated by the suspension spectra of the global classifying spaces
in $\Fc$ (by Theorem \ref{thm:Bgl G weak generators}).
So $L:\GH_\Fc\to \GH$ is an equivalence onto the full triangulated subcategory generated by
the suspension spectra $\Sigma^\infty_+ B_{\gl} G$ for all $G\in \Fc$.
So the image of $L$ is contained in $\Xc$.

Now suppose that conversely $X$ is an orthogonal spectrum in $\Xc$.
The adjunction counit $\epsilon_X:L( U X )\to X$ is an $\Fc$-equivalence, 
so it induces isomorphisms of geometric fixed point groups for all groups in $\Fc$. 
By the hypothesis on $X$ and naturality of the inflation maps $p^*$, 
the morphism $\epsilon_X$ induces isomorphisms 
of geometric fixed point homotopy groups for all compact Lie groups.
So $\epsilon_X$ is a global equivalence, 
and in particular $X$ is left induced from $\Fc$.
\end{proof}

\begin{rk}
The same proof as in Proposition \ref{prop:reflexive}
yields the following relative version of the proposition.
We let $\Fc\subset\Ec$ be global families and assume that $\Fc$ 
is {\em reflexive relative to $\Ec$}, i.e.,
for every compact Lie group $K$ from the family $\Ec$ there is
a compact Lie group $u K$, belonging to $\Fc$,
and a continuous homomorphism $p:K\to u K$ that is initial among
homomorphisms to groups in $\Fc$.
Then an orthogonal spectrum $X$ is in the essential image of the relative left adjoint
$L:\GH_\Fc\to \GH_\Ec$ if and only if
for every compact Lie group $K$ in $\Ec$ the universal morphism $p:K\to u K$
induces isomorphisms
\[ p^* \ : \ \Phi_*^{u K}( X )\ \to \ \Phi_*^K ( X ) \]
between the geometric fixed point homotopy groups of $u K$ and $K$.
\end{rk}

\begin{eg}\label{eg:reflexive example trivial}
The minimal global family $\Fc=\td{e}$ of trivial groups is reflexive, 
and the unique morphism $K\to e$ to any trivial group is universal.
So Proposition \ref{prop:reflexive} characterizes the global homotopy types
in the essential image of the left adjoint $L:\SH=\GH_{\td{e}}\to\GH$ from the
non-equivariant stable homotopy category to the global stable homotopy category:
an orthogonal spectrum $X$ is left induced from the trivial family
if and only if for every  compact Lie group $K$ the unique homomorphism
$p_K:K\to e$ induces an isomorphism
\[ p_K^*\ : \ \Phi^e_*( X ) \ \to \ \Phi^K_*( X ) \ .\]
The geometric fixed point homotopy groups $\Phi^e_*( X )$ with respect
to the trivial group are isomorphic to $\pi_*^e ( X )$, the stable homotopy groups of the
underlying non-equivariant spectrum.
So the global homotopy types in the essential image of the left adjoint $L:\SH\to\GH$
are precisely the orthogonal spectra with `constant geometric fixed points'.
\end{eg}

Here are some specific examples of left induced global homotopy types.

\begin{eg}[Suspension spectra]\label{eg:suspensions are left induced}
The global sphere spectrum $\mS$ and the suspension spectrum
of every based space are left induced from the trivial global family $\td{e}$.
Indeed, geometric fixed points commute with suspension spectra in the
following sense: 
if $A$ has trivial $G$-action, then
\[ \pi_* ( \Sigma^\infty A ) \ \iso \ \Phi^G_*(\Sigma^\infty A)\ , \]
compare Example \ref{eg:geometric and suspension}. 
So the suspension spectrum $\Sigma^\infty A$ has `constant geometric fixed points',
and it is left induced from the trivial family by the criterion
of Example \ref{eg:reflexive example trivial}.
\end{eg}

\begin{eg}[Global classifying spaces and semifree orthogonal spectra]
If $G$ is a compact Lie group from a global family $\Fc$, 
then the suspension spectrum of the global classifying space $B_{\gl} G$ 
is left induced from $\Fc$. To see this, we can refer to the proof
of Proposition \ref{prop:reflexive}; alternatively, we may show that
$\Sigma^\infty_+ B_{\gl}G$ is $\Fc$-cofibrant, i.e., has the left lifting property
with respect to morphisms that are both $\Fc$-level equivalences and
$\Fc$-level fibrations. We recall that $B_{\gl} G = \bL_{G,V}\ = \ \bL(V,-)/G$ 
is a semifree orthogonal space, where $V$ is any faithful $G$-representation.
So morphisms $\Sigma^\infty_+ B_{\gl}G \to X$ of orthogonal spectra
biject with continuous based $G$-maps $S^V \to X(V)$;
since  $S^V$ can be given the structure of a based $G$-CW-complex,
it has the left lifting property with respect to 
$G$-weak equivalences that are also $G$-fibrations, and the claim follows by
adjointness.

The same kind of reasoning shows that the semifree orthogonal spectra $F_{G,V}$
introduced in Construction \ref{con:free orthogonal}
are left induced from $\Fc$ whenever $G$ belongs to $\Fc$ and $V$ is
a faithful $G$-representation.
\end{eg}

\begin{eg}[$\bGamma$-spaces]\index{subject}{Gamma-space@$\bGamma$-space}
  We let $\bGamma$ denote the category whose objects are the based sets
  $n_+=\{0,1,\dots,n\}$, with basepoint~0, and with morphisms all based maps.
  A {\em $\bGamma$-space} is a functor from $\bGamma$ 
  to the category of based spaces which is reduced
  (i.e., the value at $0_+$ is a one-point space).

  A $\bGamma$-space $F:\bGamma\to\bT_*$ can be extended 
  to a continuous functor on the category
  of based spaces by a coend construction
  \begin{equation}  \label{eq:extended_Gamma_space}
    F(K)\ = \ \int^{n_+\in\bGamma} F(n_+)\times  K^n \ ;
  \end{equation}
  here $K$ is a based space and we use that $K^n=\map_*(n_+,K)$ 
  is contravariantly functorial in $n_+$.
  We refer to the extended functor
  as the {\em prolongation}\index{subject}{prolongation!of a $\bGamma$-space} of $F$ 
  and denote it by the same letter.
  This abuse of notation is justified by the fact that the value of the prolongation
  at $n_+$ is canonically homeomorphic to the original value, 
  see Remark \ref{rk:abuse prolongation}.
  The coend can be calculated by a familiar quotient
  space construction in the ambient category of all topological spaces,
  compare Proposition \ref{prop:prolongation in bT}:
  $F(K)$ can be obtained from the disjoint union 
  of the spaces $F(n_+)\times  K^n$, for $n\geq 0$,
  by dividing out the equivalence relation generated by
  \[ (F(\alpha)(x);\,k_1,\dots,k_n)\ \sim\ (x;\, k_{\alpha(1)},\dots,k_{\alpha(m)}) \]
  for all $x\in F(m_+)$, all $(k_1,\dots,k_n)$ in $K^n$,  and all 
  morphisms $\alpha:m_+\to n_+$ in $\bGamma$.
  Here $k_{\alpha(i)}$ is to be interpreted as the basepoint of $K$ whenever $\alpha(i)=0$.
  The non-obvious fact, proved in Proposition \ref{prop:prolongation in bT}~(ii),
  is that this quotient space is automatically compactly generated, and hence
  a coend in the category $\bT_*$.

  We write $[x;\,k_1,\dots,k_n]$ for the equivalence class 
  of a tuple $(x;\,k_1,\dots,k_n)\in F(n_+)\times K^n$ in $F(K)$.
  The prolongation is continuous and comes with a continuous, based {\em assembly map}
  \index{subject}{assembly map!of a $\bGamma$-space}
\[  \alpha\ : \ K\sm F(L) \ \to \ F(K\sm L) \ , \
    \alpha(k\sm [x;\, l_1,\dots,l_n])\ = \ [x;\,k\sm l_1,\dots,k\sm l_n]\ . \]
  The assembly map is natural in all three variables and associative and unital.

   We can now define an orthogonal spectrum $F(\mS)$
   \index{symbol}{$F(\mS)$ - {evaluation of a $\bGamma$-space $F$ on spheres}}
   \index{subject}{Gamma-space@$\bGamma$-space!evaluation on spheres}
  by evaluating the $\bGamma$-space $F$ on spheres. In other words, 
  the value of $F(\mS)$ at an  inner product space $V$ is 
  \[  F(\mS)(V) \ = \ F(S^V)   \]
  and the structure map $\sigma_{V,W}:S^V\sm F(\mS)(W)\to F(\mS)(V\oplus W)$
  is the assembly map for $K=S^V$ and $L=S^W$,
  followed by the effect of $F$ on the canonical homeomorphism 
  $S^V\sm S^W\iso S^{V\oplus W}$.
  The $O(V)$-action on $F(\mS)(V)$ is via the action on $S^V$ and the continuous
  functoriality of $F$.
\end{eg}

\begin{prop}\label{prop:geometric fix of Gamma spaces} 
  Let $F$ be a $\bGamma$-space and $G$ a compact Lie group. 
  \begin{enumerate}[\em (i)]
  \item The projection $p:G\to \pi_0 G=\bar G$ to the group of path components 
    induces an isomorphism    of geometric fixed point homotopy groups 
    \[ p^* \ : \  \Phi^{\bar G}_*( F(\mS) ) \ \to \ \Phi^G_*( F(\mS) ) \]
    of the orthogonal spectrum $F(\mS)$.
  \item The orthogonal spectrum $F(\mS)$ obtained by evaluation of $F$ on spheres
    is left induced from the global family $\Fin$ of finite groups.\index{subject}{global family!of finite groups} 
  \end{enumerate}
\end{prop}
\begin{proof}
  (i) 
  If $G$ is connected, then for every based $G$-space $K$ the map $F(K^G)\to(F(K))^G$
  induced by the fixed point inclusion $K^G\to K$ is a homeomorphism
  by Proposition \ref{prop:Gamma fixed points}, where $G$ acts trivially 
  on the $\bGamma$-space $F$.
  We can calculate $G$-fixed points by first taking fixed points with respect
  to the normal subgroup $G^\circ$ (the path component of the identity) and then
  fixed points with respect to the quotient $\bar G=G/G^\circ=\pi_0 G$. 
  So for $k\geq 0$, the $G$-geometric fixed points 
  of the orthogonal spectrum $F(\mS)$ can be
  rewritten as
  \begin{align*}
    \Phi^G_k ( F(\mS) )\ &= \  \colim_{V\in s(\Uc_G)}\, [S^{V^G\oplus \mR^k},\, F(S^V)^G] \\ 
    &\iso \ \colim_{V\in s(\Uc_G)}\, [ (S^{V^{G^\circ}\oplus\mR^k})^{\bar G},\, F(S^{V^{G^\circ}})^{\bar G}] \\ 
    &\iso \ \colim_{W\in s(\Uc_{\bar G})}\, [ S^{W^{\bar G}\oplus\mR^k},\, F(S^W)^{\bar G}]
    \ =\ \Phi^{\bar G}_k ( F(\mS))\ .
  \end{align*}
  The third step uses that $(\Uc_G)^{G^\circ}$ is a complete universe for
  the finite group $\bar G$ and as $V$ runs through $s(\Uc_G)$,
  the $G^\circ$-fixed points $V^{G^\circ}$ exhaust $(\Uc_G)^{G^\circ}$.
  The composite isomorphism is inverse to the inflation map $p^*$.
  The argument for $k<0$ is similar.
  
  (ii)
  The global family $\Fin$ of finite groups is reflexive, and
  for every compact Lie group $K$ the projection $K\to \pi_0 K$ 
  to the finite group of path components is universal with respect to $\Fin$.
  Part~(i) verifies the geometric fixed point criterion, 
  so by Proposition \ref{prop:reflexive} 
  the orthogonal spectrum $F(\mS)$ is left induced from the global family of
  finite groups.    
\end{proof}

Now we look more closely at right induced global homotopy types.
For a global family $\Fc$ and a compact Lie group $G$ we
denote by $\Fc\cap G$ the family of those closed subgroups of $G$ that belong to $\Fc$, 
and $E(\Fc\cap G)$ is a universal $G$-space for the family $\Fc\cap G$.  
We also need the equivariant cohomology theory represented by an orthogonal spectrum $X$.
If $A$ is a $G$-space, we define the $G$-cohomology $X^k_G(A)$ as
\[ X^k_G(A) \ = \ \gh{\Sigma^\infty_+ \bL_{G,V}A, X[k]}\ ,\]
the group of degree $k$ maps in $\GH$ from the suspension spectrum
of the semifree orthogonal space $\bL_{G,V} A$ to $X$.\index{subject}{equivariant cohomology theory!of an orthogonal spectrum}
Here $V$ is an implicitly chosen faithful $G$-repre\-sen\-tation.
By the adjunction between the global stable homotopy category
and the $G$-equivariant stable homotopy category that we discuss in
Theorem \ref{thm:change to groups} below, the group $X^k_G(A)$
is isomorphic to the value at $A$ of the $G$-cohomology theory
represented by the underlying $G$-spectrum $X_G$.

\begin{prop}\label{prop:right induced criterion}
An orthogonal spectrum $X$ is right induced from a global family $\Fc$
if and only if for every compact Lie group $G$
and every cofibrant $G$-space $A$ the map
\[ X^*_G(A)\ \to \ X^*_G(A\times E(\Fc\cap G)) \]
induced by the projection $A\times E(\Fc\cap G)\to A$ is an isomorphism.
\end{prop}
\begin{proof}
For every $G$-space $A$ the projection from $A\times E(\Fc\cap G)$ to $A$
is an $(\Fc\cap G)$-equivalence; moreover, if $A$ is cofibrant, 
then the source is $(\Fc\cap G)$-projective, so the suspension spectrum
of the semifree orthogonal space
\[ \Sigma^\infty_+ \bL_{G,V} (A\times E(\Fc\cap G)) \]
is left induced from the global family $\Fc$. 
This implies that 
\[  L_\Fc(U_\Fc( \Sigma^\infty_+ \bL_{G,V} A))\ \iso \ 
\Sigma^\infty_+ \bL_{G,V} (A\times E(\Fc\cap G))\]
in the global stable homotopy category. Hence
\begin{align*}
  X^k_G( A\times &E(\Fc\cap G)) \ = \ 
\gh{\Sigma^\infty_+ \bL_{G,V}( A\times E(\Fc\cap G)),X[k]} \\
&\iso \  \gh{L_\Fc(U_\Fc(\Sigma^\infty_+ \bL_{G,V} A)) ,X[k]} \
\iso \  \gh{\Sigma^\infty_+ \bL_{G,V} A ,R_\Fc(U_\Fc (X))[k]} \ .
\end{align*}
for every orthogonal spectrum $X$.
Under this composite isomorphism, the map of the proposition becomes 
the map 
\[ X_G^k(A)\ = \ \gh{\Sigma^\infty_+ \bL_{G,V} A, X[k]} \ \to\
 \gh{\Sigma^\infty_+ \bL_{G,V} A, R_\Fc(U_\Fc(X))[k]}  \]
induced by the adjunction unit $X\to R_\Fc(U_\Fc(X))$.

If $X$ is right induced from $\Fc$,
then this adjunction unit is an isomorphism,
hence so is the map $X^*_G(A)\to X^*_G(A\times E(\Fc\cap G))$.
If, on the other hand, this map is an isomorphism for all $G$-spaces $A$,
then for $A=\ast$ we deduce that the map
\[ \gh{\Sigma^\infty_+ B_{\gl} G, X[k]} \ \to\
 \gh{\Sigma^\infty_+ B_{\gl} G, R_\Fc(U_\Fc(X))[k]}  \]
is an isomorphism. Since the suspension spectrum of $B_{\gl} G$
represents $\pi_0^G$ (by Theorem \ref{thm:Bgl G weak generators} (i)), 
this shows that the adjunction unit  $X\to R_\Fc(U_\Fc(X))$ is a global equivalence.
So $X$ is right induced from $\Fc$.
\end{proof}

\begin{rk}
Essentially the same proof also shows the
following relative version of Proposition \ref{prop:right induced criterion}.
We let $\Fc\subseteq\Ec$ be two nested global families.
Then an orthogonal spectrum $X$ is in the essential image
of the relative right adjoint $R:\GH_\Fc\to\GH_\Ec$ if and only if 
for every group $G$ in $\Ec$ and every cofibrant $G$-space $A$ the map
\[ X^*_G(A)\ \to \ X^*_G(A\times E(\Fc\cap G)) \]
induced by the projection $A\times E(\Fc\cap G)\to A$ is an isomorphism.
\end{rk}

\begin{eg}\label{eg:cofree is cofree}
We let $X$ be a global $\Omega$-spectrum with the property that for every
inner product space $V$, 
the $O(V)$-space $X(V)$ is cofree,\index{subject}{cofree equivariant space} i.e.,
for some (hence any) universal free $O(V)$-space $E$ the map
\[ \text{const}\ : \  X(V)\ \to \ \map(E, X(V) ) \]
that sends a point to the corresponding constant map is an $O(V)$-weak equivalence.
We claim that then the orthogonal spectrum $X$ 
is right induced from the trivial global family $\td{e}$.
We use the criterion of Proposition \ref{prop:right induced criterion}
and show that for every compact Lie group $G$,
every cofibrant $G$-space $A$ and every integer $k$ the map $X^k_G(\Pi)$
induced by the projection $\Pi:A\times E G\to A$ is an isomorphism.

We start with the case $k=0$.
We let $V$ be any faithful $G$-representation.
The projection $\Pi$ is a weak equivalence
of underlying non-equivariant spaces, and source and target are cofibrant
as $G$-spaces. So the $G$-map
\[ S^V\sm \Pi_+\ : \ S^V\sm (A\times E G)_+\ \to \  S^V\sm A_+  \]
is also a weak equivalence
of underlying non-equivariant spaces, and source and target are $G$-cofibrant
in the based sense. 
Since $X(V)$ is cofree as an $O(V)$-space, it is also cofree as a $G$-space,
where $G$ acts via the representation homomorphism $G\to O(V)$.
So the induced map
\[  [ S^V\sm \Pi_+ , X(V) ]^G \ :\
[ S^V\sm A_+ , X(V) ]^G \ \to \ [ S^V\sm(A\times E G)_+, X(V)]^G \]
is bijective.
Since $X$ is a global $\Omega$-spectrum
and the orthogonal suspension spectrum $\Sigma^\infty_+\bL_{G,V}A$ is flat, 
the localization map
\[  \spec(\Sigma^\infty_+ \bL_{G,V}A,\, X) / \text{homotopy}\ \to \ 
\gh{\Sigma^\infty_+ \bL_{G,V}A,\, X} \ = \  X^0_G(A) \]
is bijective.
By the freeness property, the left hand side bijects with the set $[S^V\sm A_+, X(V)]^G$.
So by the previous paragraph the map
$X_G^0(\Pi):X^0_G(A)\to X^0_G(A\times E G)$
is bijective.

For $k>0$ we apply the same argument to the
global $\Omega$-spectrum $\sh^k X$ (which also has cofree levels)
and exploit the natural isomorphism
\[
 (\sh^k X)^0_G(A)\ = \ \gh{\Sigma^\infty_+ \bL_{G,V}A,\, \sh^k X}\ \iso \
\gh{\Sigma^\infty_+ \bL_{G,V}A,\, X[k]}\ = \  X^k_G(A)\ .\]
 For $k<0$ we apply the same argument to the
global $\Omega$-spectrum $\Omega^{-k} X$ (which also has cofree levels)
and exploit the natural isomorphism
\[
 (\Omega^{-k} X)^0_G(A)\ = \ \gh{\Sigma^\infty_+ \bL_{G,V}A,\, \Omega^{-k} X}\ \iso \
\gh{\Sigma^\infty_+ \bL_{G,V}A,\, X[k]}\ = \  X^k_G(A)\ .\]
\end{eg}

As we explain in more detail in Example \ref{eg:Global K is right induced}
below, the global $K$-theory spectrum $\bKU$ 
is right induced from the global family
of finite cyclic groups. This fact is a reformulation of the generalization, 
due to Adams, Haeberly, Jackowski 
and May \cite{adams-haeberly-jackowski-may:atiyah-segal},
of the Atiyah-Segal completion theorem.

\begin{eg}[Global Borel theories]\label{eg:global Borel}
We let $E$ be a non-equivariant generalized cohomology theory.
Then we obtain a global functor $\underline{E}$  by setting
\[ \underline{E}(G) \ = \ E^0(B G) \ , \] 
the 0-th $E$-cohomology of the classifying space of the group $G$.
The contravariant functoriality in group homomorphisms $\alpha:K\to G$
comes from the covariant functoriality of the classifying space construction.
The transfer map for a subgroup inclusion $H\subset G$ 
is the Becker-Gottlieb transfer \cite{becker-gottlieb}
\[ \tr\ : \ \Sigma^\infty_+ B G \ \to \ \Sigma^\infty_+ B H \]
associated to the fiber bundle $E G/H\to E G/G =B G$
with fiber $G/H$, using that $E G/H$ is a classifying space for $H$.
Strictly speaking, in \cite{becker-gottlieb}\index{subject}{Becker-Gottlieb transfer} 
Becker and Gottlieb only define a stable transfer map for 
locally trivial fiber bundles with smooth compact manifold fiber
whenever the base is a finite CW-complex. To get the transfer above
one approximates $E G$ (and hence $B G$) by its finite-dimensional skeleta.
The verification of the double coset formula for this global functor
is due to Feshbach \cite[Thm.\,II.11]{feshbach}.

More generally, we can consider the Borel equivariant cohomology
theory represented by $E$. For a compact Lie group $G$ and a cofibrant $G$-space $A$,
its value is
\[  E^*(E G\times_G A) \ , \]
the $E$-cohomology of the Borel construction 
(also known as homotopy orbit construction).
Here $E G$ is a universal free $G$-space, which is unique up to equivariant
homotopy equivalence.
We claim that these Borel cohomology theories associated to $E$
are represented by a specific global homotopy type, namely the result
of applying the right adjoint 
\[ R \ : \ \SH \ \to \ \GH \]
to the forget functor $U:\GH\to\SH$ 
to any non-equivariant spectrum that represents $E$.
To verify this claim we choose a faithful $G$-representation $V$
and recall from Proposition \ref{prop:free_orthogonal_space}~(i)
that the $G$-space $\bL(V,\mR^\infty)$ is a universal free $G$-space.
So for every $G$-space $A$ the underlying non-equivariant homotopy type 
of the closed orthogonal space $\bL_{G,V} A$ is 
\[ (\bL_{G,V}A)(\mR^\infty) \ = \ \bL(V,\mR^\infty)\times_G A \ = \ E G\times_G A\ .\]
The adjunction between $U$ and $R$ thus provides an isomorphism
\begin{align*}
 (R E)^0_G(A)\ = \ \gh{\Sigma^\infty_+ \bL_{G,V} A, R E} \ 
&\iso \ \SH( U(\Sigma^\infty_+ \bL_{G,V} A), E ) \\
&= \ \SH( \Sigma^\infty_+ (E G\times_G A), E ) \ = \ E^0(E G\times_G A)\ .  
\end{align*}
When $A$ is a one-point $G$-space, then $E G\times_G \ast= B G$,
and this bijection gives rise to a composite
isomorphism
\begin{equation}  \label{eq:pi_0^G(RE)}
  \pi_0^G(R E)\ \iso\ \gh{\Sigma^\infty_+ \bL_{G,V} , R E} \ \iso \  E^0(B G)\ ,     
\end{equation}
where the first one is inverse to evaluation at the stable tautological class
$e_G=e_{G,V}\in \pi_0^G(\Sigma^\infty_+ \bL_{G,V})$.
We claim that the isomorphisms \eqref{eq:pi_0^G(RE)} 
are compatible with restriction maps
arising from continuous group homomorphisms $\alpha:K\to G$.
For this purpose we also choose a faithful $K$-representation $W$.
This data gives rise to a composite morphism of global classifying spaces
\[ B_{\gl}\alpha \ : \ B_{\gl} K \ = \ \bL_{K,\alpha^*(V)\oplus W} \ \xra{\rho_{\alpha^*(V),W}/K}\ 
\bL_{K,\alpha^*(V)} \ \xra{\ \text{proj}\ }\ \bL_{G,V} \ = \ B_{\gl}G\ .\]
The first morphism restricts a linear isometric embedding from $\alpha^*(V)\oplus W$
to $\alpha^*(V)$, and the second morphism is the quotient map from
$K$-orbits to $G$-orbits.
On the underlying non-equivariant homotopy types
(i.e., after evaluating at $\mR^\infty$), the morphism $B_{\gl}\alpha$
classifies the homomorphism $\alpha$. Moreover, the morphism
has the `correct' behavior on the unstable tautological classes. i.e.,
\[ (B_{\gl}\alpha)_*(u_K) \ = \ \alpha^*(u_G) \text{\quad in $\pi_0^K(B_{\gl} G)$} 
\ ,\]
by direct verification from the definitions.
The analogous relation for the stable tautological classes
\begin{align*}
 (\Sigma^\infty_+ B_{\gl}\alpha)_*(e_K) \ &= \ 
 (\Sigma^\infty_+ B_{\gl}\alpha)_*(\sigma^K(u_K)) \ = \ 
 \sigma^K( (B_{\gl}\alpha)_*(u_K)) \\ 
&= \ \sigma^K( \alpha^*(u_G)) \ 
=\  \alpha^*(\sigma^G(u_G)) \ =\ \alpha^*(e_G) 
\end{align*}
in the stable group $\pi_0^K(\Sigma^\infty_+ B_{\gl} G)$
then follows by naturality 
of the suspension maps $\sigma^K:\pi_0^K(Y)\to\pi_0^K(\Sigma^\infty_+ Y)$
and its compatibility with restriction. This proves the compatibility 
of the isomorphisms \eqref{eq:pi_0^G(RE)} with restriction maps.
Compatibility with transfers is essentially built in,
as both the transfer in Construction \ref{con:external transfer} 
and the Becker-Gottlieb transfer in \cite{becker-gottlieb}
are defined as the Thom-Pontryagin construction based on smooth equivariant
embedding of $G/H$ into a $G$-representation;
we omit the formal proof. In any case, the group isomorphisms \eqref{eq:pi_0^G(RE)} 
together form an isomorphism of global functors 
between $\upi_0(R E)$ and $\underline{E}$.
This proves the claim that the `global Borel theories' are 
precisely the ones right induced from non-equivariant stable homotopy theory.
\end{eg}

\begin{construction}\label{cor:Borel b}
We introduce a specific pointset level lift
\[ b \ : \ \spec \ \to \ \spec \]
of the right adjoint $R:\SH\to \GH$ to the category of orthogonal spectra.
Given an orthogonal spectrum $E$  
we define a new orthogonal spectrum $b E$\index{symbol}{$b E$ - {global Borel theory}} 
as follows.
For an inner product space $V$ we set\index{subject}{global Borel theory|(} 
\[ (b E)(V)\ = \ \map(\bL(V,\mR^\infty),E(V)) \ ,   \]
the space of all continuous maps from $\bL(V,\mR^\infty)$ to $E(V)$.
The orthogonal group $O(V)$ acts on this mapping space by conjugation, through
its actions on $\bL(V,\mR^\infty)$ and on $E(V)$.
We define structure maps $\sigma_{V,W}:S^V\sm (b E)(W)\to (b E)(V\oplus W)$
as the composite
\begin{align*}
 S^V\sm \map(\bL(W,\mR^\infty),E(W)) \quad \xra{\text{assembly}} \qquad 
&\map(\bL(W,\mR^\infty),S^V\sm E(W)) \\ 
\xra{\map(\res_W, \sigma_{V,W}^E)}\ &\map(\bL(V\oplus W,\mR^\infty),E(V\oplus W)) 
\end{align*}
where $\res_W:\bL(V\oplus W,\mR^\infty)\to \bL(W,\mR^\infty)$
restricts an isometric embedding from $V\oplus W$ to $W$.
In the functorial description of orthogonal spectra, the structure maps
are given by
\begin{align*}
 \bO(V,W)\sm \map(\bL(V,\mR^\infty),E(V))\ &\to \  
\map(\bL(W,\mR^\infty),E(W)) \\  
(w,\varphi)\sm \quad f \hspace*{2cm} &\longmapsto \ 
\{\ \psi \mapsto X(w,\varphi)(f(\psi\circ\varphi))\ \}\ .
\end{align*}

The endofunctor $b$ on the category of orthogonal spectra
comes with a natural transformation
\[ i_E \ : \ E \ \to \ b E  \]
whose value at an inner product space $V$ sends a point
$x\in E(V)$ to the constant map $\bL(V,\mR^\infty)\to E(V)$ with value $x$.
Since $\bL(V,\mR^\infty)$ is contractible, the morphism
$i_E:E\to b E$ is a non-equivariant level equivalence,
hence a non-equivariant stable equivalence.
\end{construction}

The next result shows that the global Borel construction $b$ 
takes $\Omega$-spectra to global $\Omega$-spectra,
and that the functor $b$ realizes, in a certain precise way, 
the right adjoint to the forgetful functor
from the global to the non-equivariant stable homotopy category.
Since the morphism $i_E:E\to b E$ is a non-equivariant stable
equivalence, it becomes invertible in the non-equivariant stable
homotopy category $\SH$.
Part~(ii) of the following proposition shows 
that the morphism $i_E^{-1}:b E\to E$ is the counit of the adjunction $(U,R)$.

\begin{prop}\label{prop:global homotopy of b E}
Let $E$ be an orthogonal $\Omega$-spectrum.
\begin{enumerate}[\em (i)]
\item The orthogonal spectrum
$b E$ is a global $\Omega$-spectrum and 
right induced from the trivial global family.
\item For every orthogonal spectrum $A$ both of the two group homomorphisms
\[ \gh{A,b E}\ \xra{\ U \ } \ \SH(A,b E)\ \xra{\ (i_E)^{-1}_*}\ \SH(A, E)\ .\]
are bijective. 
\end{enumerate}
\end{prop}
\begin{proof} (i) We let $G$ be a compact Lie group and $V$ and $W$ two
$G$-representations such that $W$ is faithful. Since $E$ is an 
$\Omega$-spectrum, the adjoint structure map
\[ \tilde\sigma^E_{V,W}\ : \ E(W)\ \to \ \Omega^V E(V\oplus W) \]
is a non-equivariant weak equivalence. 
The $G$-space $\bL(W,\mR^\infty)$ 
is cofibrant by Proposition \ref{prop:K G cofibration}~(ii). 
Because $W$ is a faithful $G$-representation, the induced $G$-action on
$\bL(W,\mR^\infty)$ is free. So the induced map
\begin{align*}
  \map(\bL(W,\mR^\infty),\tilde\sigma^E_{V,W})\ : \ 
(b E)(W) \ =\  &\map(\bL(W,\mR^\infty),E(W))\\ 
&\to \ \map(\bL(W,\mR^\infty),\Omega^V E(V\oplus W)) 
\end{align*}
is a $G$-weak equivalence. Moreover, the restriction map
$\res_W:\bL(V\oplus W,\mR^\infty)\to \bL(W,\mR^\infty)$
is a $G$-homotopy equivalence (by Proposition \ref{prop:free_orthogonal_space}~(ii)), 
hence it induces another $G$-homotopy equivalence
\begin{align*}
 \map(\res_W,\Omega^V E(V\oplus W))\ : \ 
\map&(\bL(W,\mR^\infty),\Omega^V E(V\oplus W)) \\ 
&\to \ 
\map(\bL(V\oplus W,\mR^\infty),\Omega^V E(V\oplus W))  \end{align*}
on mapping spaces.
The target of this last map is $G$-homeomorphic to
\[\map_*(S^V,\map(\bL(V\oplus W,\mR^\infty),E(V\oplus W))) \ = \
\Omega^V( (b E)(V\oplus W) )\ ; \]
under this homeomorphism, the composite of the two $G$-weak equivalences
becomes the adjoint structure map
\[ \tilde\sigma^{b E}_{V,W}\ : \ (b E)(W)\ \to \ \Omega^V( (b E)(V\oplus W) )\ . \]
So we have shown that $\tilde\sigma^{b E}_{V,W}$ is a $G$-weak equivalence,
and that means that $b E$ is a global $\Omega$-spectrum.

The $O(V)$-space $\bL(V,\mR^\infty)$ is a universal free $O(V)$-space by
Proposition \ref{prop:free_orthogonal_space}~(i).
So the $O(V)$-space $(b E)(V)=\map(\bL(V,\mR^\infty),E(V))$
is cofree. 
Since $b E$ is also a global $\Omega$-spectrum,
the criterion of Example \ref{eg:cofree is cofree}
shows that it is right induced from the trivial global family.

Part~(ii) is a formal consequence of~(i): since $b E$ is right induced
from the trivial global family, the forgetful functor induces a bijection
$U:\gh{A,b E}\iso \SH(A,b E)$. Since the morphism $i_E:E\to b E$ 
becomes an isomorphism in $\SH$, it induces another bijection on $\SH(A,-)$.
\end{proof}

We endow the functor $b$ with a lax symmetric monoidal transformation
\[ \mu_{E,F}\ : \ b E\sm b F \ \to \ b(E\sm F)\ . \]
To construct $\mu_{E,F}$ we start from the $(O(V)\times O(W))$-equivariant maps
\begin{align*}
 \map(\bL(V,\mR^\infty), E(V)) \sm &\map(\bL(W,\mR^\infty), F(W)) \\ 
\ \xra{\quad \sm \quad } \quad 
&\map(\bL(V,\mR^\infty)\times \bL(W,\mR^\infty), E(V)\sm F(W)) \\
\xra{ \map(\res_{V,W},i_{V,W})}\ 
&\map(\bL(V\oplus W,\mR^\infty), (E\sm F)(V\oplus W))
\end{align*}
that constitute a bimorphism from $(b E,b F)$ to $b(E\sm F)$.
Here
\[ \res_{V,W}\ :\ \bL(V\oplus W,\mR^\infty)\ \to\ \bL(V,\mR^\infty)\times \bL(W,\mR^\infty) \]
maps an embedding of $V\oplus W$ to its restrictions to $V$ and $W$.
The morphism $\mu_{E,F}$ is associated to this bimorphism
via the universal property of the smash product \eqref{eq:universal property smash}.

\begin{rk}
Various `completion' maps\index{subject}{completion map} 
(also called `bundling maps') fit in here as follows. 
For this we suppose that $E$ is a commutative orthogonal 
ring spectrum and a positive $\Omega$-spectrum (in the non-equivariant sense).
Then the morphism $i_E:E\to b E$ is a kind of `global completion map'. 
For every compact Lie group $G$ it induces a ring homomorphism
of $G$-equivariant homotopy groups
\[ \pi_0^G(E) \ \to \ \pi_0^G( b E ) \ \iso \ E^0( B G ) \ .\]
When  $E=\mS$ is the sphere spectrum and $G$ is finite, 
Carlsson's theorem \cite{carlsson-segal conjecture} 
(proving the Segal conjecture) shows that the map
\[ \mA(G)\iso \pi_0^G(\mS) \ \to \ \pi^0(B G)\]
is completion of the Burnside ring 
at the augmentation ideal.\index{subject}{augmentation ideal!of the Burnside ring}
The sphere spectrum is the suspension spectrum of a 
global classifying space of the trivial group; 
more generally, for the global classifying space $B_{\gl} K$ 
of a finite group $K$ the `forgetful' map
\[ \bA(K,G)\iso \pi_0^G(\Sigma^\infty_+ B_{\gl}K) \ \to \ 
\SH( \Sigma^\infty_+ B G,\,\Sigma^\infty_+ B K ) \]
is again completion at the augmentation ideal of the Burnside ring $A(G)$,
see \cite[Thm.\,A]{lewis-may-mcclure-classifying}.

Since the global Borel theory functor $b:\spec\to\spec$
is lax symmetric monoidal,
it takes orthogonal ring spectra to orthogonal ring spectra, 
in a way preserving commutativity and module structures.
Since the transformation $i_E$ is monoidal, it becomes a homomorphism 
of orthogonal ring spectra when $E$ is an orthogonal ring spectrum.
We let $\mS\to\mS^{\text{f}}$ be a `positively fibrant replacement',
i.e., a morphism of commutative orthogonal ring spectra
that is a non-equivariant stable equivalence and whose target
is a positive $\Omega$-spectrum;
such a replacement exists and is homotopically unique by the
positive model structure for commutative orthogonal ring spectra
of \cite[Thm.\,15.1]{mmss}.
The spectrum $\hat \mS =b(\mS^\text{f})$\index{symbol}{$\mS$@$\hat\mS$ - {completed sphere spectrum}} thus comes with a commutative ring spectrum structure,
and we call it the {\em completed sphere spectrum}.\index{subject}{completed sphere spectrum}\index{subject}{sphere spectrum!completed|see{completed sphere spectrum}}
Moreover, for every $\mS^{\text{f}}$-module spectrum $E$ the map
\[ \hat \mS\sm b E\ = \ b(\mS^\text{f}) \sm b E \ \xra{\mu_{\mS^{\text{f}},E}} \ 
b(\mS^\text{f}\sm E) \ \xra{b(\text{act})} \ b E \]
makes the orthogonal spectrum $b E$ into a module spectrum 
over the completed sphere spectrum. 
Since $\hat \mS$ is non-equivariantly stably equivalent to $\mS$,
this shows that for every group $G$ the equivariant homotopy group
\[ \pi_k^G(b E) \ \iso \  E^{-k}(B G)  \]
is naturally a module over the commutative ring $\pi_0^G( \hat\mS )$. 

For the global $K$-theory spectrum 
(compare Construction \ref{con:global KU} below)
and any compact Lie group $G$, the bundling map
\[ \bRU(G)\iso \pi_0^G ( \bKU ) \ \to \ K U^0(B G)\]\index{subject}{augmentation ideal!of the unitary representation ring}
takes a virtual $G$-representation to the class of the associated virtual vector bundle
over $B G$. The Atiyah-Segal completion theorem \cite{atiyah_segal-completion}
shows that this map is completion at the augmentation ideal of the representation ring.
For the Eilenberg-Mac\,Lane spectrum $\Hc\mZ$
(see Construction \ref{con:HM} below), 
the global functor $G\mapsto H^0(B G;\mZ)$ is constant with value $\mZ$; the map 
\[  \pi_0^G(\Hc\mZ) \ \to \ H^0(B G;\mZ)\]
is surjective and an isomorphism modulo torsion 
for all compact Lie groups whose identity path component 
is commutative (compare Example \ref{eg:HZ for abelian G}).
\end{rk}
\index{subject}{global Borel theory|)}

Now we fix a compact Lie group $G$ and relate the global stable homotopy category 
to the $G$-equivariant stable homotopy category (based on a complete universe). 
We denote by $G\text{-}\SH$\index{symbol}{  $G\text{-}\SH$ - {$G$-equivariant stable homotopy category}} 
the {\em $G$-equivariant stable homotopy category},\index{subject}{stable homotopy category!$G$-equivariant}
i.e., any localization of the category $G\spec$ of orthogonal $G$-spectra
at the class of $\upi_*$-isomorphisms.
Various stable model structures have been established with 
 $\upi_*$-isomorphisms as weak equivalences, for example by
Mandell-May \cite[III Thm.\,4.2]{mandell-may}, Stolz \cite[Thm.\,2.3.27]{stolz-thesis}
and Hill-Hopkins-Ravenel \cite[Prop.\,B.63]{HHR-Kervaire}.
In particular, as for every stable model category,
the homotopy category  $G\text{-}\SH$
comes with a preferred structure of a triangulated category.

A functor
\[ (-)_G \ : \ \spec \ \to \ G\spec \ , \quad X \ \longmapsto \ X_G \]
from orthogonal spectra to orthogonal $G$-spectra
is given by endowing an orthogonal spectrum with the trivial $G$-action.
Since the trivial action functor takes global equivalences 
to $\upi_*$-isomorphisms,
the universal property of localizations provides 
a `forgetful' functor on the homotopy categories
\[ U = U_G\ : \ \GH \ \to \ G\text{-}\SH  \ .\]
Moreover, $U$ is canonically an exact functor of triangulated categories. 
We will show now that the forgetful functor has both a left and a right adjoint.

The `equivariant' smash product of orthogonal $G$-spectra
is simply the smash product of the underlying non-equivariant
orthogonal spectra with diagonal $G$-action.
So the trivial action functor $(-)_G:\spec\to G\spec$ is
strong symmetric monoidal.
The smash product of orthogonal spectra and of orthogonal $G$-spectra 
can be derived to symmetric monoidal products on $\GH$ and on $G\text{-}\SH$ 
(see Corollary \ref{cor-derived smash}). The forgetful functor
is strongly monoidal with respect to these derived smash products.
Indeed, the derived smash product in $\GH$ can be calculated by
flat approximation up to global equivalence;
a flat orthogonal spectrum, endowed with trivial $G$-action,
is $G$-flat, and hence can be used to calculate the derived smash product
 $G\text{-}\SH$, by the flatness theorem (Theorem \ref{thm:G-flat is flat}).

When $G=e$ is a trivial group, the next theorem reduces to the change of family
functor of Theorem \ref{thm:change families}, with $\Ec=\All$ and $\Fc=\td{e}$.

\begin{theorem}\label{thm:change to groups}
  For every compact Lie group $G$ the forgetful functor
  \[ U \ : \ \GH \ \to \ G\text{-}\SH \]
  preserves sums and products, and it has a left adjoint and a right adjoint.
  The left adjoint has a preferred lax symmetric comonoidal structure.
  The right adjoint has a preferred lax symmetric monoidal structure.
\end{theorem}
\begin{proof}
  Sums in $\GH$ and $G\text{-}\SH$ are represented 
  in both cases by the pointset level wedge.
  For $\GH$ we state this explicitly in Proposition \ref{prop:sums and product in GH}~(i);
  for $G\text{-}\SH$ we can run the argument based on the stable model
  structure for orthogonal $G$-spectra established in \cite[III Thm\,4.2]{mandell-may}.
  On the pointset level, the forgetful functor preserves wedges,
  so the derived forgetful functor preserves sums.
  As spelled out in Corollary \ref{cor-generators for GH_F}~(iv),
  the existence of the right adjoint is a formal consequence
  of the fact that $\GH$ is compactly generated 
  and that the functor $U$ preserves sums.

  The forgetful functor also preserves infinite products,
  but the argument here is slightly more subtle because products in $\GH$
  are not generally represented by the pointset level product,
  and because equivariant homotopy groups do not in general commute
  with infinite pointset level products,
  compare Remark \ref{rk:homotopy of infinite product}.
  We let $\{X_i\}_{i\in I}$ be a family of orthogonal spectra.
  By replacing each factor by a globally equivalent spectrum, if necessary,
  we can assume without loss of generality that each $X_i$ is a global $\Omega$-spectrum.
  Since global $\Omega$-spectra are the fibrant objects
  in a model structure underlying $\GH$, 
  the pointset level product $\prod_{i\in I} X_i$ then represents
  the product in $\GH$.
  
  Even though $X_i$ is a global $\Omega$-spectrum,
  the underlying orthogonal $G$-spec\-trum $(X_i)_G$ need {\em not}
  be a $G$-$\Omega$-spectrum. However, as we spell out in the proof of
  Proposition \ref{prop:sums and product in GH}~(ii), the natural map
  \[ \pi_k^G({\prod}_{i\in I} \, X_i) \ \to \ {\prod}_{i\in I} \, \pi_k^G(X_i)\]
  is an isomorphism for all integers $k$.
  Again we can run the argument of Proposition \ref{prop:sums and product in GH}~(ii)
  in the stable model structure 
  for orthogonal $G$-spectra \cite[III Thm\,4.2]{mandell-may}, and conclude that
  in this situation, the pointset level product is also a product in $G\text{-}\SH$.
  So the derived forgetful functor preserves products.
  The existence of the left adjoint is then again a formal consequence
  of the fact that $\GH$ is compactly generated,  
  compare Corollary \ref{cor-generators for GH_F}~(v).

  The same formal argument as in part~(iii) of Theorem \ref{thm:change families}
  shows how to turn the strong monoidal structure of the forgetful functor $U$
  into a lax comonoidal structure 
  $L(  A \sm ^\mL  B ) \to ( L A ) \sm^\mL ( L B ) $ of the left adjoint.
  In contrast to Theorem \ref{thm:change families}~(iii), however,
  this morphism is usually {\em not} an isomorphism, so we cannot turn it
  around into a monoidal structure on $L$.
  The same formal argument as in Theorem \ref{thm:change families}~(ii)
  constructs the lax monoidal structure on $R$ from the strong monoidal
  structure of the forget functor $U$.
\end{proof}

\Danger Theorem \ref{thm:change to groups} looks similar to the change-of-family
Theorem \ref{thm:change families}, but there is one important difference:
if the group $G$ is non-trivial, then neither of the two adjoints
to the forgetful functor  $U:\GH \to G\text{-}\SH$ is fully faithful.

\begin{rk}
  We mention an alternative way to construct the two adjoints
  to the forgetful functor $U : \GH \to G\text{-}\SH$, by exhibiting the
  pointset forgetful functor $(-)_G:\spec\to G\spec$
  as a left respectively a right Quillen functor for suitable model structures. 
  We sketch this for the left adjoint, where we can use 
  the stable model structure of orthogonal $G$-spectra
  established by Mandell and May in \cite[III Thm\,4.2]{mandell-may}.
  However, we cannot argue directly with the functor $(-)_G:\spec\to G\spec$,
  since it is {\em not} a right Quillen functor.
  Indeed, if it were a right Quillen functor, then it would  preserve
  fibrant objects. However, a global $\Omega$-spectrum is typically
  {\em not} a $G$-$\Omega$-spectrum when given the trivial $G$-action.
  
  What saves us is that a global $\Omega$-spectrum is `eventually'
  a $G$-$\Omega$-spectrum, i.e., starting at faithful representations.
  This lets us modify $(-)_G$ into a right Quillen functor as follows. 
  We choose a faithful $G$-representation $V$ and let
  \[ \Omega^V \sh^V \ : \ \spec \ \to \ G\spec \]
  denote the functor that takes an orthogonal spectrum $X$ to the
  orthogonal $G$-spectrum with $U$-th level
  \[ ( \Omega^V \sh^V X )(U) \ = \ \map_*(S^V, X(U\oplus V))\ . \]
  We emphasize that the $G$-action on  $\Omega^V \sh^V X$ is non-trivial,
  despite the fact that we started with an orthogonal spectrum without a $G$-action.
  A natural morphism of orthogonal $G$-spectra 
  $\tilde\lambda^V_X:X\to \Omega^V \sh^V X$ is given by the adjoint
  of the morphism $\lambda^V_X:X\sm S^V\to \sh^V X$
  defined in \eqref{eq:defn lambda_n};
  the morphism $\tilde\lambda^V_X$ 
  is a $\upi_*$-isomorphism by Proposition \ref{prop:lambda upi_* isos}~(ii).
  In particular, the functor $ \Omega^V \sh^V $ also takes
  global equivalences of orthogonal spectra to $\upi_*$-isomorphisms
  of orthogonal $G$-spectra, and the derived functor of $\Omega^V \sh^V$ 
  is naturally isomorphic to the forgetful functor $U:\GH\to G\text{-}\SH$.
  The argument can then be completed by showing that  
  the functor $\Omega^V \sh^V$ is a right Quillen functor
  from the global model structure on orthogonal spectra
  to the stable model structure on orthogonal $G$-spectra
  from \cite[III Thm\,4.2]{mandell-may}.
  Then by general model category theory, the derived functor
  of $\Omega^V \sh^V$, and hence also the forgetful functor $U$,
  has a left adjoint.
  
  The existence of the right adjoint to $U$ can also be established
  by model category reasoning. For this one can use the $\mathbb S$-model
  structure on orthogonal $G$-spectra 
  constructed by Stolz \cite[Thm.\,2.3.27]{stolz-thesis}.
  We leave it to the interested reader to verify that the cofibrant objects
  in the $\mathbb S$-model structure are precisely 
  the $G$-flat orthogonal $G$-spectra in the sense of Definition \ref{def:G-flat}.
  This shows that the forgetful functor
  $(-)_G:\spec\to G\spec$ is a left Quillen functor from the
  global model structure on orthogonal spectra to the stable $\mathbb S$-model
  structure on orthogonal $G$-spectra.
\end{rk}

The left adjoint $L:G\text{-}\SH\to\GH$ to the forgetful functor
is an exact functor of
triangulated categories that preserves infinite sums.
The $G$-equivariant stable homotopy category is compactly generated
by the unreduced suspension spectra of all the coset spaces $G/H$,
for all closed subgroups $H$ of $G$. So $L$ is essentially determined by
its values on these generators.
The sequence of natural bijections
\[ \GH(L(\Sigma^\infty_+ G/H), X) \ \iso \
G\text{-}\SH(\Sigma^\infty_+ G/H,\, U X) \ \iso \ \pi_0^H (X) \ \iso \ 
\GH(\Sigma^\infty_+ B_{\gl}H , X) \]
shows that the left adjoint $L$ takes the unreduced suspension spectrum
of the coset space $G/H$ to the suspension spectrum  of the 
global classifying space of $H$.
In the special case $H=G$ the spectrum $\Sigma^\infty_+ G/H$
is the equivariant sphere spectrum $\mS_G$, and we obtain that
\[ L(\mS_G)\ \iso\ \Sigma^\infty_+ B_{\gl}G \ .  \]
Now
\[ G\text{-}\SH(\mS_G,\mS_G) \ \iso \ \pi_0^G (\mS) \ \iso \ \mA(G) \ = \ \bA(e,G) \]
is the Burnside ring, whereas
\[ \GH(L(\mS_G),L(\mS_G)) \ \iso \  \pi_0^G (\Sigma^\infty_+ B_{\gl} G)  
\ \iso \ \ \bA(G,G)  \]
is the double Burnside ring. 
The map $L:G\text{-}\SH(\mS_G,\mS_G) \to \GH(L(\mS_G),L(\mS_G))$ 
corresponds to the ring homomorphism
\[ \mA(G)\ = \ \bA(e,G) \ \to \ \bA(G,G) \ , \quad
\tr_H^G\circ p_H^* \ \longmapsto \ \tr_H^G\circ \res^G_H \]
from the Burnside ring to the double Burnside ring; 
this homomorphism is never surjective unless $G$ is trivial, 
so the left adjoint is not full.

\begin{rk} The discussion in this section could be done relative
to a global family $\Fc$, as long as $\Fc$ contains the compact Lie group
$G$ under consideration (and hence also all its subgroups).
Indeed, if $\Fc$ contains $G$, then every $\Fc$-equivalence
of orthogonal spectra is a $\upi_*$-isomorphism of 
underlying orthogonal $G$-spectra.
Hence the trivial action functor descends to
a `forgetful' functor on the homotopy categories
\[  U_G^\Fc\ : \ \GH_\Fc \ \to \ G\text{-}\SH  \]
by the universal property of localizations. The same arguments as
in Theorem \ref{thm:change to groups} show the existence of both
adjoints to this forgetful functor, with the same kind of monoidal properties. 

Theorem \ref{thm:change to groups} discusses the maximal case
of the global family $\Fc=\All$ of all compact Lie groups.
The minimal case is the global family $\td{G}$ generated by $G$,
i.e., the class of compact Lie groups that are isomorphic to
a quotient of a closed subgroup of $G$. All the forgetful functors $U^\Fc_G$
then factor as composites 
\[ \GH_\Fc \ \xra{U_{\td{G}}^\Fc} \ \GH_{\td{G}} \ \xra{U_G^{\td{G}}}\  G\text{-}\SH\]
of a change-of-family functor and a family-to-group functor.
The various adjoints then compose accordingly.
\end{rk}

\Danger 
Whenever $G$ is non-trivial, then the global homotopy category
$\GH_{\td{G}}$ associated to the global family generated by $G$ is different
from the $G$-equivariant stable homotopy category $G\text{-}\SH$.
In other words, if $G$ is non-trivial, then the forgetful family-to-group functor 
$U_G^{\td{G}}: \GH_{\td{G}}\to G\text{-}\SH$ is {\em not} an equivalence,
and neither of its adjoints is fully faithful.
\index{subject}{global stable homotopy category|)}

\medskip

\index{subject}{global stable homotopy category!for finite groups|(}
\index{subject}{global family!of finite groups|(} 
In the rest of this section we turn to the global family $\Fin$ of finite groups
and describe the associated global stable homotopy category $\GH_\Fin$ rationally. 
By our previous results, $\GH_\Fin$ is a compactly generated triangulated category with
a symmetric monoidal derived smash product.
We call an object $X$ of the category $\GH_\Fin$ {\em rational}\index{subject}{orthogonal spectrum!rational}
if the equivariant homotopy groups $\pi_k^G(X)$ are uniquely divisible 
(i.e., $\mQ$-vector spaces) for all finite groups $G$.
In this section we will give an algebraic model of the 
{\em rational $\Fin$-global stable homotopy category},
i.e., the full subcategory $\GH^\mQ_\Fin$ of rational spectra in $\GH_\Fin$.
Theorem \ref{thm:rational SH} below shows that the homotopy types in $\GH^\mQ_\Fin$
are determined by a chain complex of global functors,
up to quasi-isomorphism.
More precisely, we construct an equivalence of triangulated categories
from  $\GH^\mQ_\Fin$ to the unbounded derived category of rational global functors
on finite groups.

We let $G$ and $K$ be compact Lie groups.
We recall from Proposition \ref{prop:B_gl represents} 
that the evaluation map
\[ \bA(G,K) \ \to \ \pi_0^K(\Sigma^\infty_+ B_{\gl} G) \ , \quad 
\tau \longmapsto \tau(e_G)\]
is an isomorphism, where $e_G\in \pi_0^G(\Sigma^\infty_+ B_{\gl} G)$
is the stable tautological class.
More precisely, the definition of the global classifying space $B_{\gl} G$
involves an implicit choice of faithful $G$-representation $V$ that is omitted 
from the notation, and $e_G$ is the class denote $e_{G,V}$ 
in \eqref{eq:define_stable_tautological}.
The unreduced suspension spectrum of every orthogonal spaces is globally connective
(see Proposition \ref{prop:pi_0 of Sigma^infty}),
so the group $\pi_k^K(\Sigma^\infty_+ B_{\gl} G)$ is trivial for $k<0$.

\begin{prop}\label{prop:pi_0 B G is A(-,G)}\index{subject}{global classifying space!of a finite group}   
Let $G$ and $K$ be finite groups. Then for every $k>0$, 
the equivariant homotopy group $\pi_k^K(\Sigma^\infty_+ B_{\gl} G)$ is torsion.
\end{prop}
\begin{proof} 
We show first that the geometric fixed point group 
$\Phi_k^K(\Sigma^\infty_+ B_{\gl} G)$ is torsion for $k>0$.\index{subject}{geometric fixed points}
This part of the argument needs $G$ to be finite, but $K$ could be any compact Lie group.
Geometric fixed points commute with suspension spectra, i.e., 
the groups $\Phi_*^K(\Sigma^\infty_+ B_{\gl} G)$
are isomorphic to the non-equivariant stable homotopy groups
of the fixed point space $((B_{\gl}G)(\Uc_K))^K$.
Proposition \ref{prop:fix of global classifying}~(i) 
identifies these fixed points as
\[ (( B_{\gl}G)(\Uc_K))^K \ \simeq \ 
\coprod_{[\alpha]\in\text{Rep}(K,G)} \quad B C(\alpha)\ , \]
where the disjoint union is indexed by conjugacy
classes of homomorphisms from $K$ to $G$, 
and $C(\alpha)$ is the centralizer of the image of $\alpha:K\to G$.
Since $G$ is finite, so are all the centralizers $C(\alpha)$,
hence the classifying space $B C(\alpha)$ has no rational homology,
hence no rational stable homotopy, in positive dimensions.
So we conclude that the rationalized stable homotopy groups 
of the space $((B_{\gl}G)(\Uc_K))^K$ vanish in positive dimensions. 

If $K$ is also finite, then the $k$-th rationalized equivariant 
stable homotopy group of any orthogonal $K$-spectrum can be recovered from the 
$k$-th rationalized geometric fixed point homotopy groups
for all subgroups $L$ of $K$, as described in 
Corollary \ref{cor:equivariant from geometric}.
So when both $G$ and $K$ are finite, then 
also the equivariant homotopy group $\pi_k^K(\Sigma^\infty_+ B_{\gl} G)$ is torsion
for all $k>0$.
\end{proof}

\Danger The conclusion of Proposition \ref{prop:pi_0 B G is A(-,G)} 
is no longer true if we drop the finiteness hypothesis on one
of the two groups $G$ or $K$.
For example, for $G=e$ we have $\Sigma^\infty_+ B_{\gl} G=\mS$, 
and the dimension shifting transfer $\Tr_e^{U(1)}(1)$ is an element 
of infinite order in the group $\pi_1^{U(1)}(\mS)$.
On the other hand, for $K=e$ the group $\pi_k^K(\Sigma^\infty_+ B_{\gl} G)$
is the non-equivariant stable homotopy group of the 
ordinary classifying space $B G$. For $G=U(1)$ the group
$\pi_k^e(\Sigma^\infty_+ B U(1))$ contains a free summand of rank~1
whenever $k\geq 0$ is even.

\medskip

Now we can establish the algebraic model for the rational $\Fin$-global
homotopy category.
We let $\Ac$ be a pre-additive category, such as the 
$\Fin$-Burnside category $\bA_\Fin$.
We denote by $\Ac\mo$ the category of additive functors from
$\Ac$ to the category of $\mQ$-vector spaces.
This is an abelian category, and the represented functors
$\Ac(a,-)$, for all objects $a$ of $\Ac$, form a set
of finitely presented projective generators of $\Ac\mo$.
The category of $\mZ$-graded chain complexes in the abelian category $\Ac\mo$
then admits the {\em projective model structure} with the quasi-isomorphisms
as weak equivalences. The fibrations in the projective model structure
are those chain morphisms that are surjective in every chain complex degree
and at every object of $\Ac$.
This projective model structure for complexes of $\Ac$-modules is a
special case of \cite[Thm.~5.1]{Christensen-Hovey-relative}.
Indeed, the projective (in the usual sense) $\Ac$-modules together
with the epimorphisms form a projective class 
(in the sense of \cite[Def.~1.1]{Christensen-Hovey-relative}),
and this class is determined
(in the sense of \cite[Sec.~5.2]{Christensen-Hovey-relative})
by the set of represented functors.

We also need the rational version of the $\Fc$-global model structure,
for a global family $\Fc$.
We call a morphism $f:X\to Y$ of orthogonal spectra a
{\em rational $\Fc$-equivalence}\index{subject}{F-equivalence@$\Fc$-equivalence!rational}
if the map
\[ \mQ\tensor \pi_k(f)\ : \ \mQ\tensor\pi_k^G(X)\ \to\ \mQ\tensor\pi_k^G(Y) \]
is an isomorphism for all integers $k$ and all groups $G$ in the family $\Fc$.

\begin{theorem}[Rational $\Fc$-global model structure]
  \label{thm:rational F-global spectra} 
  Let $\Fc$ be a global family.
  \begin{enumerate}[\em (i)]
  \item 
    The rational $\Fc$-equivalences and the $\Fc$-cofibrations 
    are part of a model structure on the category of orthogonal spectra, 
    the {\em rational $\Fc$-global model structure}.\index{subject}{F-global model structure@$\Fc$-global model structure!rational}
  \item
    The fibrant objects in the rational $\Fc$-global model structure 
    are the $\Fc$-$\Omega$-spectra $X$ such that for all 
    $G\in\Fc$ the equivariant homotopy groups $\pi_*^G(X)$ are uniquely divisible.
  \item
    The rational $\Fc$-global model structure is cofibrantly generated, 
    proper and topological. 
  \end{enumerate}
\end{theorem}

Theorem \ref{thm:rational F-global spectra} is obtained
by Bousfield localization of the $\Fc$-global model structure on
orthogonal spectra, and one can use a similar proof as for the rational stable
model structure on sequential spectra in \cite[Lemma 4.1]{schwede-shipley-uniqueness}.
We omit the details. 

\begin{theorem}\label{thm:rational SH} 
There is a chain of Quillen equivalences between
the category of orthogonal spectra with the 
rational $\Fin$-global model structure and the category of chain complexes of
rational global functors on finite groups. In particular, this induces
an equivalence of triangulated categories
\[  \GH^\mQ_\Fin\ \xra{\ \iso \ } \ \Dc\left( \GF^\mQ_\Fin \right) \ . \]
The equivalence can be chosen so that the homotopy group global functor
on the left hand side corresponds to the homology global functor 
on the right hand side.
\end{theorem}
\begin{proof}
We prove this as a special case of the `generalized tilting theorem'
of Brooke Shipley and the author. 
Indeed, by Theorem \ref{thm:Bgl G weak generators} 
the suspension spectra of the global classifying spaces $B_{\gl}G$ 
are compact generators of the global homotopy category
$\GH_{\Fin}$ as $G$ varies through all finite groups.
So the rationalizations  $(\Sigma^\infty_+ B_{\gl}G)_\mQ$ 
are compact generators of the rational global homotopy category $\GH_{\Fin}^\mQ$. 
If $k$ is any integer, then 
the morphism vector spaces between two such objects are given by
\begin{align*}
\gh{(\Sigma^\infty_+ B_{\gl} K)_\mQ[k],\ (\Sigma^\infty_+ B_{\gl} G)_\mQ} \ &\iso \
\pi_k^K((\Sigma^\infty_+ B_{\gl} G)_\mQ) \ 
\iso \ \mQ\tensor \pi_k^K( \Sigma^\infty_+ B_{\gl} G ) \\ 
&\iso \ 
\begin{cases}
  \mQ\tensor \bA(G,K) &\text{\quad for $k=0$, and}\\
\qquad 0 &\text{\quad for $k\ne0$.}
\end{cases}
\end{align*}
The vanishing for $k>0$ is Proposition \ref{prop:pi_0 B G is A(-,G)};
the vanishing for $k<0$ is Proposition \ref{prop:pi_0 of Sigma^infty}.

The rational $\Fin$-global model structure on orthogonal spectra is 
topological (hence simplicial), cofibrantly generated, proper and stable; 
so we can apply the Tilting Theorem \cite[Thm.\,5.1.1]{schwede-shipley-modules}. 
This theorem yields a chain
of Quillen equivalences between orthogonal spectra 
in the rational $\Fin$-global model structure and the category of chain complexes of
$\mQ\tensor\bA_\Fin$-modules, i.e., 
additive functors from the rationalized Burnside category $\mQ\tensor\bA_\Fin$
to abelian groups. 
This functor category is equivalent to the category of additive functors
from $\bA_\Fin$ to $\mQ$-vector spaces, and this proves the theorem.
\end{proof}

\begin{rk} There is an important homological difference between
global functors on finite groups and Mackey functors for one fixed finite group.
Indeed, for a finite group $G$, 
the category of rational $G$-Mackey functors is equivalent to a product, 
indexed over conjugacy classes $(H)$ of subgroups of $G$, of the
module categories over the rational group rings $\mQ[W_G H]$
of the Weyl groups, see Theorem \ref{thm:split Mackey functors}~(ii). 
In particular, the abelian  category of rational $G$-Mackey functors
is semisimple, every object is projective and injective
and the derived category is equivalent, by taking homology,
to the category of graded rational $G$-Mackey functors.

There is no analog of this for rational $\Fin$-global functors.
For example, the rationalized augmentation
\[ \mQ\tensor\mA= \mQ\tensor \bA(e,-) \ \to \ \underline{\mQ}\]
from the rationalized Burnside ring global functor to the constant global functor
for the group $\mQ$ does {\em not} split on finite groups.
The new phenomenon is that any splitting would have to be natural for inflation maps.
Let us be even more specific.
In the constant global functor $\underline{\mQ}$ we have
\[ \tr_e^{C_2}(1)\ =\ 2\cdot p^*(1)
 \text{\qquad in\qquad $\underline{\mQ}(C_2)=\mQ$\ ,} \]
where $p:C_2\to e$ is the unique group homomorphism.
So for every morphism of global functors $\varphi:\underline{\mQ}\to N$
the image $\varphi(e)(1)$ of the unit element under the map
$\varphi(e):\mQ\to N(e)$ must satisfy
\[ \tr_e^{C_2}(\varphi(e)(1)) \ = \ 2\cdot p^*(\varphi(e)(1))\ . \]
But in the Burnside ring $\mA(e)$, and also in its rationalization,
0 is the only element in the kernel of $\tr_e^{C_2}-2\cdot p^*$; 
so every morphism of global functors from $\underline{\mQ}$ to $\mQ\tensor \mA$~is zero.
\end{rk}

While the abelian category $\GF^\mQ_\Fin$ is not semisimple,
we can still `divide out transfers' and thereby replace $\GF^\mQ_\Fin$ by an equivalent, 
but simpler category.
We let $\Out$\index{symbol}{  $\Out$ - {category of finite groups and conjugacy classes
of epimorphisms}} denote the category of finite groups and conjugacy classes
of surjective group homomorphisms. We write
\[ \Mod\Out \ = \ \Fc(\Out^{\op},\Ab) \]
for the category of `right $\Out$-modules', i.e., contravariant functors
from $\Out$ to the category of abelian groups.
To a $\Fin$-global functor $M:\bA_\Fin\to\Ab$ we can associate a right $\Out$-module
$\tau M:\Out^{\op}\to\Ab$, the {\em reduced functor} as follows.
On objects we set
\[ (\tau M)(G) \ = \ \tau_G M\ = \ M(G)/ t_G M \ ,\]
the quotient of the group $M(G)$ by the subgroup $t_G M$
generated by the images of all transfer maps $\tr_H^G:M(H)\to M(G)$ 
for all proper subgroups $H$ of $G$.
If $\alpha:K\to G$ is a surjective group homomorphism and
$H\leq G$ a proper subgroup, then $L=\alpha^{-1}(H)$ is a proper
subgroup of $K$ and the relation
\[ \alpha^*\circ \tr_H^G \ = \ 
\tr_L^K \circ (\alpha|_L)^*\ : \ M(H)\to M(K) \]
shows that the inflation
map $\alpha^*:M(G)\to M(K)$ passes to a homomorphism
$\alpha^*:(\tau M)(G)\to (\tau M)(K)$ of quotient groups.
We will now argue that the reduction functor
\[ \tau \ : \ \GF_\Fin \ \to \ \Mod\Out \]
is rationally an equivalence of categories, 
compare Theorem \ref{thm:rational global functors}.
By construction, the projection maps $M(G)\to\tau M(G)$ 
form a natural transformation from the restriction of
the global functor $M$ to the category $\Out^{\op}$ to $\tau M$.

\begin{eg} We let $A$ be an abelian group and $\underline{A}$
the constant global functor with value $A$,
compare Example \ref{eg:global functor examples}~(iii).\index{subject}{global functor!constant} 
Then $(\tau\underline A)(G)= A/c A$,
where $c$ is the greatest common divisor of the indices of all proper subgroups
of $G$. If $G$ is not a $p$-group for any prime $p$, then this
greatest common divisor is~1. If $G$ is a non-trivial $p$-group, then
$G$ has a proper subgroup of index $p$. So we have
\[ (\tau\underline A)(G) \ = \ 
\begin{cases}
\  A & \text{if $G=e$,}\\
  A/p A & \text{if $G$ is a non-trivial $p$-group, and}\\
\  0 & \text{else.}
\end{cases}\]
The inflation maps in $\tau\underline{A}$ are quotient maps.
\end{eg}

\begin{eg} The Burnside ring $\mA(K)=\bA(e,K)$\index{subject}{Burnside ring global functor} 
of a finite group $K$ is
freely generated, as an abelian group, by the transfers $\tr_L^K(1)$
where $L$ runs through representatives of the conjugacy classes
of subgroups of $K$. So $(\tau \mA)(K)$ is free abelian of rank~1,
generated by the class of the multiplicative unit~1. All restriction maps
preserve the unit, so the reduced functor $\tau\mA$ of the global Burnside
ring functor is isomorphic to the constant functor $\Out^{\op}\to\Ab$
with value $\mZ$.

We denote by $\bA_G=\bA(G,-)$ the global functor represented by a finite group $G$.
We will now present
$\tau(\bA_G)$ as an explicit quotient of a sum of representable $\Out$-modules.
For every closed subgroup $H$ of $G$ the restriction map $\res_H^G$
is a morphism in $\bA(G,H)$. If $H$ is finite, the 
Yoneda lemma provides a unique morphism
\[\Out_H \ \to \tau(\bA_G)\]
of $\Out$-functors from the representable functor
$\Out_H=\mZ[\Out(-,H)]$ to $\tau(\bA_G)$ that sends the identity of $H$
to the class of $\res_H^G$ in $\tau (\bA_G)(H)$.
For every element $g\in G$ the conjugation isomorphism
$c_g:H\to H^g$ given by $c_g(\gamma)=g^{-1}\gamma g$ induces
an isomorphism
\[g_\star\circ - \ : \ \bA_G(H^g) \ \to \ \bA_G(H) \]
by postcomposition. We have 
\[ g_\star\circ [\res_{H^g}^G] \ = \ [g_\star\circ \res_{H^g}^G] \ = \ 
[\res_H^G\circ g_\star] \ = \ [\res_H^G] \]
in $(\tau \bA_G)(H)$; so the Yoneda lemma shows that the triangle of $\Out$-functors
\[\xymatrix@R=6mm@C=6mm{ \Out_H \ar[dr] \ar[rr]^-{c_g\circ -}&&\Out_{H^g} \ar[dl]\\
& \tau(\bA_G) }\]
commutes. The direct sum of the transformations $\Out_H\to\tau(\bA_G)$
thus factors over a natural transformation
\begin{equation}  \label{eq:presentation_of_tau_A_G}
  \left({\bigoplus}_{H\leq G} \ \Out_H\right)/G \ \to \ \tau (\bA_G) \ . 
\end{equation}
The source of this morphism can be rewritten if we choose representatives
of the conjugacy classes of subgroups in $G$:
\[ \bigoplus_{(H)} \  ( \Out_H / W_G H )\ \to \ \tau (\bA_G) \ .   \]
Now the sum is indexed by conjugacy classes.
\end{eg}

\begin{prop}\label{prop:tau of A_G}
For every finite group $G$, the morphism \eqref{eq:presentation_of_tau_A_G}
is an isomorphism of $\Out$-modules.
\end{prop}
\begin{proof}
By Theorem \ref{thm:Burnside category basis} 
the abelian group $\bA_G(K)=\bA(G,K)$ is
freely generated by the elements $\tr_L^K\circ\alpha^*$
where $(L,\alpha)$ runs through representatives of the $(K\times G)$-conjugacy classes
of pairs consisting of a subgroup $L$ of $K$
and a homomorphism $\alpha:L\to G$. 
So $(\tau \bA_G)(K)$ is a free abelian group with basis the
classes of $\alpha^*$ for all conjugacy classes of
homomorphisms $\alpha:K\to G$. 

On the other hand, the group 
$( \Out_H / W_G H )(K)$ is free abelian with basis
given by $W_G H$-orbits of conjugacy classes of epimorphisms $\alpha:K\to H$.
The map \eqref{eq:presentation_of_tau_A_G} sends the basis
element represented by $\alpha$ to the basis element represented
by the composite of $\alpha$ with the inclusion $H\to G$.
So the homomorphism \eqref{eq:presentation_of_tau_A_G} 
takes a basis of the source to a basis of the target, and is thus an isomorphism.
\end{proof}

We recalled in Proposition \ref{prop:recover G-Mackey rationally}
above how the value of a $G$-Mackey functor $M$,
for a finite group $G$, can rationally be recovered
from the groups $(\tau M)(H)$ for all subgroups $H$ of $G$:
the map
\[ \psi^M_G \ : \ M(G)\ \to \ \left( {\prod}_{H\leq G}\, (\tau M)(H)\right)^G  \]
whose $H$-component is the composite
\[ M(G)\ \xra{\ \res^G_H\ } \ M(H)\ \xra{\text{proj}}  (\tau M)(H)\]
becomes an isomorphism after tensoring with $\mQ$.
When applied to a $\Fin$-global functor $M$, we see that $M$ can rationally be
recovered from the $\Out$-module $\tau M$. This is the key input to the following
equivalence of categories.
I suspect that the following result is well known,
but I have been unable to find an explicit reference.

\begin{theorem}\label{thm:rational global functors} 
The restriction of the functor $\tau$ to the full subcategory
of rational $\Fin$-global functors is  an equivalence of categories
\[ \tau \ : \ \GF_\Fin^\mQ \ \to \ \Mod\Out_\mQ \ = \ \Fc(\Out^\text{\rm op},\mQ)\]
onto the category of rational $\Out$-modules.
\end{theorem}
\begin{proof}
Since global functors are an enriched functor category and the functor
\[ \tau \ : \ \GF_\Fin\ \to \ \Mod\Out \]
commutes with colimits, $\tau$ has a right adjoint
\[ \rho \ : \ \Mod\Out\ \to \ \GF_\Fin \ .\]
If $\rho X$ is an $\Out$-module, the value of the $\Fin$-global functor $\rho X$
at a finite group $G$ is necessarily given by
\[ (\rho X)(G)\ = \ \Mod\Out(\tau(\bA_G), X)\ , \]
the group of $\Out$-module homomorphisms from $\tau(\bA_G)$ to $X$.
The global functoriality in $G$ is via $\tau$, i.e., as the composite
\begin{align*}
  \bA(G,K)\tensor  (\rho X)(G)\
&\to \ \GF(\bA_K,\bA_G)\tensor  (\rho X)(G)\\
&\xra{\tau\tensor\Id} \ 
\Mod\Out(\tau(\bA_K),\tau(\bA_G))\tensor  (\rho X)(G)\
\xra{\ \circ \ } \ (\rho X)(K)\ .
\end{align*}
We rewrite the definition of $(\rho X)(G)$
using the description of the $\Out$-module $\tau(\bA_G)$
given in \eqref{eq:presentation_of_tau_A_G}.
Indeed, by Proposition \ref{prop:tau of A_G},
precomposition with \eqref{eq:presentation_of_tau_A_G} induces an isomorphism
\[ (\rho X)(G)\ \iso \ 
\Mod\Out\left(  \left({\bigoplus}_{H\leq G} \ \Out_H\right)/G, X\right)\ 
\ \iso \ \left( {\prod}_{H\leq G} X(H)\right)^G\ . \]
So for every global functor $M$,
this description shows that 
$\rho(\tau M)(G)$ is isomorphic to the target of the morphism $\psi^M_G$.
A closer analysis reveals that for $X=\tau M$, 
the above isomorphism identifies $\psi^M_G$
with the value of the adjunction unit $\eta:M\to\rho(\tau M)$ at $G$.
So Proposition \ref{prop:recover G-Mackey rationally}
shows that for every rational $\Fin$-global functor $M$
the adjunction unit $\eta_M:M\to\rho(\tau M)$ is an isomorphism.
This implies that the restriction of the left adjoint $\tau$ 
to the category of rational $\Fin$-global functors is fully faithful.

Now we consider a rational $\Out$-module $X$. Then $\rho X$
is a rational $\Fin$-global functor, so $\eta_{\rho X}:\rho X\to \rho(\tau(\rho X))$
is an isomorphism by the previous paragraph.
Since $\eta_{\rho X}$ is right inverse to $\rho(\epsilon_X)$,
the morphism $\rho(\epsilon_X):\rho(\tau(\rho X))\to\rho X$
is an isomorphism of $\Fin$-global functors.
By~Proposition \ref{prop:tau of A_G} the represented $\Out$-module
$\Out_G$ is a direct summand of $\tau(\bA_G)$.
So the group
\[ \Mod\Out(\Out_G,X) \ \iso \ X(G) \]
is a direct summand of the group
\[ \Mod\Out(\tau(\bA_G),X) \ \iso \ \GF(\bA_G,\rho X) \ \iso \ (\rho X)(G) \ ,\]
and this splitting is natural for morphisms of $\Out$-modules in $X$.
In particular, the morphism $\epsilon_X(G):(\tau(\rho X))(G)\to X(G)$
is a direct summand of the morphism
\[ \rho(\epsilon_X)(G)\ :\ (\rho(\tau(\rho X)))(G)\ \to\ (\rho X)(G) \ .  \]
The latter is an isomorphism by the previous paragraph,
so the morphism $\epsilon_X(G)$ is also an isomorphism. 
This shows that for every rational $\Out$-module $X$
the adjunction counit $\epsilon_X:\tau(\rho X)\to X$ is an isomorphism.

Altogether we have now seen that when restricted to rational objects on both sides,
the unit and counit of the adjunction $(\tau,\rho)$ are isomorphisms. 
This proves the theorem.
\end{proof}

The rational equivalence $\tau$ of abelian categories prolongs
to an equivalence of derived categories by applying $\tau$ dimensionwise
to chain complexes. The combination with the equivalence of triangulated
categories of Theorem \ref{thm:rational SH} is then a chain of
two exact equivalences of triangulated categories 
\begin{equation}\label{eq:combined_triangulated_equivalence}
  \GH^\mQ_\Fin\ \xra[\ \iso\ ] \ \Dc\left(\GF^\mQ_\Fin\right)
\ \xra[\iso]{\Dc(\tau)}\ \Dc(\Mod\Out_\mQ)\ . 
\end{equation}
The next proposition shows that this composite equivalence is
an algebraic model for the geometric fixed point homotopy groups.

For every orthogonal spectrum $X$ and every compact Lie group $G$,
the geometric fixed point map $\Phi:\pi^G_0(X)\to\Phi_0^G(X)$ 
annihilates all transfers from proper subgroups 
by Proposition \ref{prop:Phi of transfer}. So the geometric fixed
point map factors over a homomorphism
\[ \bar\Phi\ : \  \tau (\upi_k(X))(G) \ \to \ \Phi_k^G(X)  \]
that we called the reduced geometric fixed point map above.\index{subject}{geometric fixed point homomorphism!reduced}
The geometric fixed point maps are compatible with inflations
(Proposition \ref{prop:properties geometric fixed maps}~(iii)),
so as $G$ varies among finite groups, the reduced geometric fixed
point maps form a morphism of $\Out$-modules.
When we apply Proposition \ref{prop:rational geometric fixed for GSp} 
to the underlying orthogonal $G$-spectrum of an orthogonal spectrum,
it specializes to the following:

\begin{cor}\label{cor:rational geometric fixed} 
For every orthogonal spectrum $X$, every finite group $G$ and every integer $k$
the map\index{subject}{geometric fixed points!rationally} 
\[ \bar\Phi \ : \ \tau( \upi_k(X) ) (G) \ \to \ \Phi^G_k(X) \]
becomes an isomorphism after tensoring with $\mQ$.
So for varying finite groups $G$, these maps
form a rational isomorphism of $\Out$-functors
$\tau(\upi_k(X))\iso \underline{\Phi}_k(X)$.
\end{cor}

As a corollary we obtain that the combined equivalence $\kappa$ of triangulated 
categories \eqref{eq:combined_triangulated_equivalence}
from the rational finite global homotopy category $\GH^\mQ_\Fin$
to the derived category of the abelian category $\Mod\Out_\mQ$
comes with a natural isomorphism
\[ \Phi_*^G( X)\ \iso \   H_*(\kappa(X))   \ ,\]
for every object $X$ of $\GH^\mQ_\Fin$,
between the geometric fixed point homotopy groups and
the homology $\Out$-modules of $\kappa(X)$.

\index{subject}{global stable homotopy category!for finite groups|)}
\index{subject}{global family!of finite groups|)}

\chapter{Ultra-commutative ring spectra}
\label{ch:ultra}

This chapter is devoted to ultra-commutative ring spectra,
our model for extremely highly structured, multiplicative global stable
homotopy types.
On the point set level, these objects are simply 
commutative orthogonal ring spectra;
we use the term {\em `ultra-commutative'} to emphasize 
that we care about their homotopy theory with respect to
multiplicative morphisms that are global equivalences.\index{subject}{ring spectrum!ultra-commutative|see{ultra-commutative ring spectrum}}\index{subject}{ultra-commutative ring spectrum|(}
We refer to the introduction of Chapter \ref{ch-umon} 
for further justification of the adjective `ultra-commutative'.
In short, the slogan ``$E_\infty$=commutative'' is not true globally, 
and a strictly commutative multiplication encodes 
a large amount of additional structure that deserves a special name.

Section \ref{sec:power ops ring spectra}
introduces the formal setup for power operations on ultra-commu\-ta\-tive ring spectra.
We define global power functors as global Green functors 
equipped with additional power operations,
satisfying a list of axioms reminiscent of the properties
of the power maps $x\mapsto x^m$ in a commutative ring.
We show in Theorem \ref{thm:pi_0 is global power functor} 
that the global functor $\upi_0(R)$ of an ultra-commutative ring spectrum $R$
supports such power operations, and is an example of a global power functor.

Section \ref{sec:algebra global power} 
is primarily of algebraic nature, and gives both a monadic and
a comonadic description of the category of global power functors.
We introduce the comonad of `exponential sequences'
on the category of global Green functors, and show that its
coalgebras are equivalent to global power functors.
A formal consequence is that the category of global power functors
has all limits and colimits, and that they are created in the category of
global Green functors.
We discuss localization of global Green functors and global power functors
at a multiplicative subset  of the underlying ring, including
rationalization of global power functors.

In Section \ref{sec:GPF examples} we discuss various examples of 
global power functors, such as the Burnside ring global power functor, 
the global functor represented by an abelian compact Lie group,
free global power functors, constant global power functors,
and the complex representation ring global functor.
In Section \ref{sec:ucom global model} we set up the global model structure 
on the category of ultra-commutative ring spectra, compare 
Theorem \ref{thm:global ultra-commutative}.
In Theorem \ref{thm:generators global power functor}  we calculate the algebra
of natural operations on the 0-th homotopy groups of ultra-commutative ring spectra:
we show that these operations are freely generated by restrictions,
transfers and power operations. 
Theorem \ref{thm:realize global power functor} shows that every global
power functor can be realized by an ultra-commutative ring spectrum.

\section{Power operations}
\label{sec:power ops ring spectra}
\index{subject}{global power functor|(}

In this section we introduce the formal setup for encoding
the power operations on ultra-commutative ring spectra.
In Definition \ref{def:power functor} we define global power functors,
which are global Green functors equipped with additional power operations,
satisfying a list of axioms reminiscent of the properties
of the power maps $x\mapsto x^m$ in a commutative ring.
Theorem \ref{thm:pi_0 is global power functor} shows 
that the global functor $\upi_0(R)$ of an ultra-commutative ring spectrum $R$
supports power operations, and is an example of a global power functor.
As we shall see in Theorem \ref{thm:generators global power functor} below,
all the natural operations on $\upi_0(R)$ are generated by restrictions,
transfers and power operations,
so we are not missing any additional structure.
Moreover, every global power functor is realized by an
ultra-commutative Eilenberg-Mac\,Lane ring spectrum,
see Theorem \ref{thm:realize global power functor} below.
For a different perspective of global power functors
(restricted to finite groups), 
including the relationship to the concepts of $\lambda$-rings, 
$\tau$-rings and $\beta$-rings,
we refer the reader to Ganter's paper \cite{ganter-power}.

\begin{defn}
An {\em ultra-commutative ring spectrum} is a commutative orthogonal ring spectrum.
We write $ucom$ for the category of ultra-commutative ring spectra.\index{symbol}{  $ucom$ - {category of ultra-commutative ring spectra}}
\end{defn}
\index{subject}{ring spectrum!ultra-commutative|see{ultra-commutative ring spectrum}}
\index{subject}{ultra-commutative ring spectrum}

In Section \ref{sec:power op spaces} we introduced 
power operations and transfers on the equivariant homotopy sets of
ultra-commutative monoids.
For every ultra-commutative ring spectrum $R$, 
the orthogonal space $\Omega^\bullet R$ 
inherits a commutative multiplication, making it an ultra-commutative monoid
(compare Example \ref{eg:suspension spectrum of orthogonal monoid space}).
Moreover, $\pi_0^G(\Omega^\bullet R)=\pi_0^G(R)$,
so this endows the 0-th stable equivariant homotopy groups $\pi_0^G(R)$ 
with multiplicative power operations\index{subject}{power operation!in an ultra-commutative ring spectrum}
and transfers, natural for homomorphisms of ultra-commutative ring spectra.
Since these operations come from the multiplicative 
(as opposed to the `additive' structure) of the ring spectrum, 
we now switch to a multiplicative notation and write
\begin{equation}\label{eq:power_op_ring_spectra} 
P^m\  : \  \pi_0^G ( R )\ \to\ \pi^{\Sigma_m\wr G}_0 ( R )
\end{equation}
(instead of $[m]$) for the multiplicative power operations, and we write $N_H^G$
(instead of $\tr_H^G$) for multiplicative transfers.
Multiplicative transfers are often referred to as {\em norm maps},
and that is also the terminology we will use. 
Since we will work with power operations a lot, 
we take the time to expand the definition: the operation $P^m$ takes the class
represented by a based $G$-map $f:S^V\to R(V)$, for some $G$-representation $V$,
to the class of the $(\Sigma_m\wr G)$-map
\[ S^{V^m} = (S^V)^{\sm m} \ \xra{\ f^{\sm m}\ } \
R(V)^{\sm m} \ \xra{\ \mu_{V,\dots,V}\ } R(V^m) \ , \]
where $\mu_{V,\dots,V}$
is the iterated, $(\Sigma_m\wr G)$-equivariant multiplication map of $R$.

\begin{defn}
A {\em global Green functor}\index{subject}{global Green functor}
\index{subject}{Green functor!global|see{global Green functor}}
is a commutative monoid in the category of global functors
under the monoidal structure given by the box product 
of Construction \ref{con:Box product}.\index{subject}{box product!of global functors} 
We write $\GlGre$ for the category of global Green functors.\index{symbol}{  $\GlGre$ - {category of global Green functors}} 
\end{defn}

As we explain after Definition \ref{def:diagonal product},
the commutative multiplication on a global Green functor $R$
can be made more explicit in two equivalent ways:
\begin{itemize}
\item as a commutative ring structure on the group $R(G)$
for every compact Lie group, subject to the requirement
that all restrictions maps are ring homomorphisms and the transfer maps
satisfy reciprocity;
\item as a unit element $1\in R(e)$ and biadditive, commutative, associative and unital
external pairings $\times :R(G)\times R(K)\to R(G\times K)$ that are
morphisms of global functors in each variable separately.
\end{itemize}

We clarify next how the power operations 
of ultra-commutative ring spectra interact with the other structure
on equivariant stable homotopy groups, such as the addition,
restriction and transfer maps.
For $m\geq 2$ the power operation $P^m$ is {\em not} additive, but it satisfies
various properties reminiscent of the map $x\mapsto x^m$ in a commutative ring.
We formalize these properties into the concept of a {\em global power functor}.
Conditions~(i) through~(vi) in the following definition 
express the fact that a global power functor has an underlying `multiplicative'
global power monoid, 
in the sense of Definition \ref{def:power monoid},
if we forget the additive structure.
The definition makes use of certain embeddings between products
and wreath products:
\begin{align}\label{eq:wreath_sum spectra}  
\Phi_{i,j}\ : \ (\Sigma_i\wr G)\ \times\ (\Sigma_j\wr G)\hspace*{1.3cm} &\to\quad \Sigma_{i+j}\wr G \\
 ((\sigma;\, g_1,\dots,g_i),\, (\sigma';\, g_{i+1},\dots,g_{i+j})) \ &\longmapsto \
(\sigma+\sigma';\, g_1,\dots,g_{i+j})\nonumber
\end{align}
and
\begin{align}\label{eq:wreath_iterate spectra} 
\Psi_{k,m}\ : \ \Sigma_k\wr (\Sigma_m\wr G) \qquad
&\to\qquad \Sigma_{k m}\wr G \\
 (\sigma;\, (\tau_1;\, h^1),\dots,(\tau_k;\, h^k)) \ &\longmapsto \
(\sigma\natural(\tau_1,\dots,\tau_k);\, h^1+\dots+h^k)\ .\nonumber
\end{align}
These monomorphisms were defined
in Construction \ref{con:wreath product morphisms}.

\begin{defn}\label{def:power functor} \index{subject}{global power functor}
A {\em global power functor} is a global Green functor $R$ 
equipped with maps
\[ P^m \ : \ R(G) \ \to \ R(\Sigma_m\wr G)\]
for all compact Lie groups $G$ and $m\geq 1$,
called {\em power operations}, that satisfy the following relations.
\begin{enumerate}[(i)]
\item (Unit)  $P^m(1)=1$ for the unit $1\in R(e)$.
\item (Identity) $P^1=\Id$ under the identification $G\iso\Sigma_1\wr G$
sending $g$ to $(1;g)$.
\item (Naturality) 
For every continuous homomorphism $\alpha:K\to G$ between compact Lie groups
and all $m\geq 1$ the relation
\[ P^m\circ \alpha^* \ = \ (\Sigma_m\wr\alpha)^*\circ P^m\]
holds as maps $R(G)\to R(\Sigma_m\wr K)$.
\item (Multiplicativity) 
For all compact Lie groups $G$, 
all $m\geq 1$ and all classes $x,y\in R(G)$ the relation
\[ P^m(x\cdot  y) \ = \ P^m(x)\cdot P^m(y) \]
holds in the group $R(\Sigma_m\wr G)$. 
\item (Restriction)
For all compact Lie groups $G$, all $m >k>0$ and all $x\in R(G)$
the relation
\[  \Phi^*_{k,m-k}( P^m(x) )\ = \ P^k(x)\times P^{m-k}(x) \]
holds in $R((\Sigma_k\wr G)\times(\Sigma_{m-k}\wr G))$.
\item (Transitivity)
For all compact Lie groups $G$, all $k,m\geq 1$ and all $x\in R(G)$ the relation
\[  \Psi^*_{k,m}( P^{k m}(x) )\ = \ P^k( P^m (x) ) \]
holds in $R(\Sigma_k\wr(\Sigma_m \wr G))$.
\item (Additivity) 
For all compact Lie groups $G$, all $m\geq 1$, and all $x,y\in R(G)$
the relation
\[ P^m(x+y) \ = \ \sum_{k=0}^m\ \tr_{k,m-k} (P^k(x)\times P^{m-k}(y)) \]
holds in $R(\Sigma_m\wr G)$, where $\tr_{k,m-k}$ is the transfer associated
to the embedding $\Phi_{k,m-k}:(\Sigma_k\wr G)\ \times\ (\Sigma_{m-k}\wr G)\to\Sigma_m\wr G$ defined in \eqref{eq:wreath_sum spectra}. Here $P^0(x)$ is the multiplicative unit~1.
\item (Transfer) For every closed subgroup $H$ of 
a compact Lie group $G$  and every $m\geq 1$ the relation
\[ P^m\circ \tr_H^G \ = \ \tr_{\Sigma_m\wr H}^{\Sigma_m\wr G}\circ P^m\]
holds as maps $R(H)\to R(\Sigma_m\wr G)$.
\end{enumerate}
A {\em morphism} of global power functors is a morphism of
global Green functors that also commutes with the power operations.\index{subject}{morphism!of global power functors}
We write $\GlPow$ for the category of global power functors.\index{symbol}{  $\GlPow$ - {category of global power functors}} 
\end{defn}

In a global power functor the relation $P^m(0)=0$
also holds for every $m\geq 1$ and all $G$.
Indeed, the additivity and unit relations give
\[   \sum_{k=1}^{m-1}\ \tr_{k,m-k} (P^k(0)\times 1 )\ = \ P^m(0+ 1) - P^m(0) - P^m(1)
 \ = \ -P^m(0)\ . \]
Starting from $P^1(0)=0$ this shows inductively that $P^m(0)=0$.

\begin{rk}[Global power functors versus global Tambara functors]
\label{rk:power vs Tambara}
A global power functor gives rise to two underlying global power monoids, 
the {\em additive} and the {\em multiplicative} one.
As we explained in Construction \ref{con:transfer map}, 
applied to the multiplicative global power monoid,
the power operations $P^m$ lead to `multiplicative transfers',
$N_H^G:R(H)\to R(G)$ that are called {\em  norm maps},
for every subgroup $H$ of finite index in $G$.
For the convenience of the reader, we recall the construction of the norm maps.
We suppose that $[G:H]=m$,
and we choose an $H$-basis of $G$, i.e., an ordered $m$-tuple 
$\bar g=(g_1,\dots,g_m)$ of elements in disjoint $H$-orbits such that
\[ G \ = \ {\bigcup}_{i=1}^m \ g_i H \ .\]
The wreath product $\Sigma_m\wr H$ acts freely and transitively 
from the right on the set of all such $H$-bases of $G$, by the formula
\[ (g_1,\dots,g_m)\cdot (\sigma;\,h_1,\dots,h_m)\ = \ 
(g_{\sigma(1)} h_1,\dots,g_{\sigma(m)}h_m)\ .  \]
We obtain a continuous homomorphism $\Psi_{\bar g}:G\to\Sigma_m\wr H$ by requiring that
\[ \gamma \cdot\bar g \ = \ \bar g\cdot\Psi_{\bar g}(\gamma)  \ .\]
The norm $N_H^G:R(H)\to R(G)$\index{subject}{norm map}\index{symbol}{$N_H^G$ - {norm map}} is then the composite
\[ R(H) \ \xra{\ P^m\ } \ R(\Sigma_m\wr H) \ \xra{\ \Psi_{\bar g}^*\ }\ R(G) \ .\]
Any other $H$-basis is of the form
$\bar g\omega$ for a unique $\omega\in\Sigma_m\wr H$. 
We have $\Psi_{\bar g\omega}=c_\omega\circ \Psi_{\bar g}$ as maps $G\to\Sigma_m\wr H$,
where $c_\omega(\gamma)=\omega^{-1}\gamma\omega$. Since inner automorphisms
induce the identity in any Rep-functor, we have
\[ \Psi_{\bar g}^* = \Psi_{\bar g\omega}^* \ : \ R(\Sigma_m\wr H) \ \to \ R(G)\ . \]
So the norm $N_H^G$ does not depend on the choice of basis $\bar g$. 

The norms maps satisfy a number of important relations,
by Proposition \ref{prop:transfer} applied to the multiplicative monoid
of the global power functor $R$. There relations -- turned into multiplicative notation --
are as follows.
\begin{enumerate}[(i)]
\item {(Transitivity)} 
We have $N_G^G=\Id_{R(G)}$ and for nested subgroups
$H\subseteq G\subseteq F$ of finite index the relation
\[ N_G^F \circ N_H^G \ = \  N_H^F  \]
holds as maps  $R(H)\to R(F)$.
\item {(Multiplicative double coset formula)}
For every subgroup $K$ of $G$ (not necessarily of finite index)
the relation\index{subject}{double coset formula!for norm map}
\[ \res^G_K\circ N_H^G  \ = \ \prod_{[g]\in K\backslash G/H}
N^K_{K\cap{^g H}}\circ g_\star \circ \res^H_{K^g\cap H} \]
holds as maps $R(H)\to R(K)$.
Here $[g]$ runs over a set of representatives of the finite
set of $K$-$H$-double cosets.
\item {(Inflation)}
For every continuous epimorphism  $\alpha:K\to G$ of compact Lie groups
the relation
\[ \alpha^*\circ N_H^G \ =\  N_L^K\circ (\alpha|_L)^* \]
holds as maps from $R(H)\to R(K)$, where $L=\alpha^{-1}(H)$.
\item 
For every $m\geq 1$ the $m$-th power operation can be recovered as
\begin{equation}\label{eq:norm2power}
 P^m \ =\  N_K^{\Sigma_m\wr G}\circ q^* \ ,
\end{equation}
where $K$ is the subgroup of $\Sigma_m\wr G$
consisting of all $(\sigma;\,g_1,\dots,g_m)$ such that $\sigma(m)=m$
and $q:K\to G$ is defined by $q(\sigma;\,g_1,\dots,g_m)=g_m$.
\end{enumerate}
In particular, the power operations define the norm maps,
but they can also be reconstructed from the norm maps.
So the information in a global power functor could be packaged
in an equivalent but different way using norm maps instead of
power operations. The algebraic structure that arises then
is the global analog of a {\em TNR-functor}\index{subject}{TNR-functor} 
in the sense of Tambara \cite{tambara}, nowadays also called
a {\em Tambara functor};\index{subject}{Tambara functor|see{TNR-functor}}
here the acronym stands form `Transfer, Norm and Restriction'.
This observation can be stated as an equivalence of categories 
between global power functors and a certain category of `global TNR functors'; 
we shall not pursue this further.
Our reason for favoring power operations over norm maps
is that they satisfy explicit and intuitive formulas 
with respect to the rest of the structure 
(restriction, transfer, sum, product,\dots).
The norm maps also satisfy universal formulas when applied to sums
and transfers, but I find them harder to describe and to remember.

For a fixed finite group $G$, Brun \cite[Sec.\,7.2]{brun-Witt cobordism} 
has constructed norm maps on the 0-th equivariant homotopy group Mackey functor 
of every commutative orthogonal $G$-ring spectrum, and he showed
that this structure is a TNR-functor. 
So when restricted to finite groups, the global power functor structure
on $\upi_0(R)$ for an ultra-commutative ring spectrum,
obtained in the following Theorem \ref{thm:pi_0 is global power functor}, 
could also be deduced by using Brun's TNR-structure 
for the underlying orthogonal $G$-ring spectrum $R_G$, and then turning the norm maps
into power operations as \eqref{eq:norm2power}.
However, Brun's construction is rather indirect and this would hide
the simple and explicit nature of the power operations.
\end{rk}

In order to show that the 0-th equivariant homotopy group functor
of an ultra-commutative ring spectrum satisfies the transfer axiom~(viii)
of a global power functor, we study the interplay between power operations,
the Wirth\-m{\"u}ller isomorphism and the degree shifting transfer.
To state the results we first have to generalize power operations
from equivariant homotopy groups to equivariant homology theories.

We let $R$ be an orthogonal spectrum, $G$ a compact Lie group and $A$
a based $G$-space. We define the $G$-equivariant $R$-homology group of $A$ 
as the group\index{subject}{equivariant homology group!of an orthogonal spectrum}
\[ R_0^G(A) \ = \ \pi_0^G(R\sm A)\ = \ \colim_{V\in s(\Uc_G)} \, [S^V, R(V)\sm A ]^G\ .\]
Every continuous group homomorphism $\alpha:K\to G$ induces a 
restriction homomorphism
\[ \alpha^* \ : \ R_0^G(A) \ \to \ R_0^K(\alpha^*(A)) \]
that generalizes the restriction homomorphism $\alpha^*:\pi_0^G(R)\to\pi_0^K(R)$.
Again $\alpha^*$ is defined by applying restriction of scalars to any representative
of a given equivariant homology class.

Now we let $R$ be an orthogonal ring spectrum (not necessarily commutative).
Then the equivariant homology theories represented by $R$ inherit 
multiplications in the form of bilinear maps
\[ \times \ : \ R_0^G(A) \times R_0^G(B)\ \to \ R_0^G(A\sm B) \ .\]
We define this pairing simply as the composite
\begin{align*}
  R_0^G(A) \times R_0^G(B)\ = \   \pi_0^G(R\sm A) \times \pi_0^G&(R\sm B)\ 
\xra{\ \times\  } \ \pi_0^G((R\sm A) \sm (R\sm B))\\ 
&\xra{\ \mu_* \ } \ \pi_0^G(R\sm A\sm  B)\ = \ R_0^G(A\sm B) \ ,
\end{align*}
where the first map is the homotopy group pairing of 
Construction \ref{con:pairing equivalence homotopy}
and $\mu:(R\sm A) \sm (R\sm B)\to R\sm A\sm  B$
stems from the multiplication of $R$.
For $A=B=S^0$ this construction reduces to the pairings 
of equivariant homotopy groups \eqref{eq:multiplication_upi_ringspectrum}.

Now we suppose that the ring spectrum $R$ is ultra-commutative.
Given a based $G$-space $A$, we write $A^{(m)}=A^{\sm m}$ for its $m$-fold
smash power, which is naturally a based $(\Sigma_m\wr G)$-space.
Then we define the $m$-th power operation
\[ P^m \ : \ R_0^G(A)\ \to \ R_0^{\Sigma_m\wr G}(A^{(m)}) \]
by the obvious generalization of \eqref{eq:power_op_ring_spectra}: 
the operation $P^m$ takes the class
represented by a based $G$-map $f:S^V\to R(V)\sm A$, for some $G$-representation $V$,
to the class of the $(\Sigma_m\wr G)$-map
\begin{align*}
  S^{V^m} = (S^V)^{(m)} \ \xra{\ f^{(m)}\ } \
&(R(V)\sm A)^{(m)} \\ 
\xra{\ \text{shuffle}\ } \ &R(V)^{(m)}\sm A^{(m)} \ 
\xra{\mu_{V,\dots,V}\sm A^{(m)}} R(V^m)\sm A^{(m)} \ ,
\end{align*}
where $\mu_{V,\dots,V}$ is the iterated multiplication map of $R$.
\index{subject}{power operation!in equivariant homology groups}
We omit the straightforward verification that the power operations 
in equivariant $R$-homology are compatible with restriction maps:
for every continuous homomorphism $\alpha:K\to G$ between compact Lie groups
and every based $G$-space $A$, the relation
\[ P^m\circ \alpha^* \ = \ (\Sigma_m\wr \alpha)^*\circ P^m \]
holds as maps from $R_0^G(A)$ to $R_0^{\Sigma_m\wr K}( \alpha^*(A)^{(m)})$,
exploiting that $(\Sigma_m\wr\alpha)^*(A^{(m)})=\alpha^*(A)^{(m)}$
as $(\Sigma_m\wr K)$-spaces.

The following proposition makes precise how the power operations 
in equivariant $R$-homology interact with 
the Wirthm{\"u}ller isomorphism of Theorem \ref{thm:Wirth iso}.
To give the precise statement we have to introduce additional notation.
We let $H$ be a closed subgroup of a compact Lie group $G$.
As before we let $L = T_{e H}( G/H)$ denote the tangent $H$-representation, 
the tangent space of $G/H$ at the distinguished coset $e H$. 
We write
\[ \gamma \ : \ (G/H)^m \ \iso \ (\Sigma_m\wr G) / (\Sigma_m\wr H) \ , \
(g_1 H,\dots, g_m H)\ \longmapsto \ (1;\,g_1,\dots, g_m)\cdot(\Sigma_m\wr H) \]
for the distinguished $(\Sigma_m\wr G)$-equivariant diffeomorphism.
The differential of $\gamma$ at $(e H,\dots,e H)$ 
is a $(\Sigma_m\wr H)$-equivariant linear isometry
\[  (d \gamma)_{(e H,\dots,e H)}\ : \ L^m \ \iso \ 
T_{e(\Sigma_m\wr H)} \big( (\Sigma_m\wr G) / (\Sigma_m\wr H) \big)\  .\]
In the next proposition and its corollaries, we will use this
equivariant isometry to identify $L^m$ with the tangent representation
of $\Sigma_m\wr H$ inside $\Sigma_m\wr G$.

\begin{prop}\label{prop:power and Wirthmuller} 
Let $R$ be an ultra-commutative ring spectrum and $H$ 
a closed subgroup of a compact Lie group $G$. Then the following diagram commutes
\[ \xymatrix@C=15mm@R=7mm{ 
R_0^G(G/H_+) \ar[r]^-{\Wirth^G_H}_-\iso \ar[d]_{P^m} &
R_0^H(S^L) \ar[d]^{P^m} \\
R_0^{\Sigma_m\wr G}(( G / H)^m_+)  \ar[d]^\iso_{ \gamma_*} &  
R_0^{\Sigma_m\wr H}( ( S^L)^{(m)}) \ar[d]^\iso\\
R_0^{\Sigma_m\wr G}(  (\Sigma_m\wr G) /  (\Sigma_m\wr H)_+ ) 
\ar[r]_-{\Wirth^{\Sigma_m\wr G}_{\Sigma_m\wr H}}^-\iso  & R_0^{\Sigma_m\wr H}(S^{L^m}) } \]
where the horizontal maps are the respective Wirthm{\"u}ller isomorphisms.\index{subject}{Wirthm{\"u}ller isomorphism}
\end{prop}
\begin{proof}
We choose a slice as in the definition of the Wirthm{\"u}ller map in
Construction \ref{con:Wirthmuller map}, 
i.e., a smooth embedding $s:D(L)\to G$ that satisfies
\[ s(0)=1\ , \quad s(h\cdot l)=h\cdot s(l)\cdot h^{-1}
\text{\qquad and\qquad} s(-l)=s(l)^{-1} \]
for all $(l,h)\in D(L)\times H$, and such that the differential at~0 of the composite
\[ D(L) \ \xra{\ s\ }\ G \ \xra{\text{proj}}\ G/H\]
is the identity.
The collapse map
\[ \lambda_H^G\ = \ l_H^G/H \ : \ G/H_+ \ \to \ S^L \]
is then given by the formula
\[ \lambda_H^G(g H) \ = \ 
\begin{cases}
  l / (1-|l|)   & \text{ if $g=s(l)\cdot h$ with $(l,h)\in D(L)\times H$, and }\\
\quad \ast & \text{ if $g$ is not of this form.}
\end{cases} \]
We define a slice for the pair $(\Sigma_m\wr G,\Sigma_m\wr H)$
from the slice $s$ for the pair $(G,H)$, namely as the smooth embedding
\[ \bar s \ : \ D(L^m)\ \to \ \Sigma_m\wr G \ , \quad
\bar s(l_1,\dots, l_m)\ = \ (1;\, s(l_1),\dots, s(l_m))\ .\]
Clearly, $\bar s(0,\dots,0)$ is the multiplicative unit,
\begin{align*}
    \bar s(-l_1,\dots,-l_m)\ &= \ (1;\, s(- l_1),\dots, s(-l_m))\\ 
&= \ (1;\, s(l_1)^{-1},\dots, s(l_m)^{-1})\ = \ \bar s(-l_1,\dots,-l_m)^{-1} \ ,
\end{align*}
and
\begin{align*}
  \bar s((\sigma^{-1};\,h_1,\dots,&h_m)\cdot(l_1,\dots, l_m))\ 
=\ \bar s(h_{\sigma(1)}l_{\sigma(1)},\dots,h_{\sigma(m)} l_{\sigma(m)})\\ 
&=\ (1;\,s(h_{\sigma(1)}l_{\sigma(1)}),\dots,s(h_{\sigma(m)} l_{\sigma(m)}))\\ 
&=\ (1;\,h_{\sigma(1)}s(l_{\sigma(1)})h_{\sigma(1)}^{-1},\dots,
h_{\sigma(m)}s(l_{\sigma(m)})h_{\sigma(m)}^{-1})\\ 
&=\  (\sigma^{-1};\,h_1,\dots,h_m)\cdot (1;\,s(l_1),\dots, s(l_m))\cdot 
(\sigma;\,h_{\sigma(1)}^{-1},\dots,h_{\sigma(m)}^{-1})\\
&=\  (\sigma^{-1};\,h_1,\dots,h_m)\cdot \bar s(l_1,\dots, l_m)\cdot 
(\sigma^{-1};\,h_1,\dots,h_m)^{-1}\ ,
\end{align*}
for all $(l_1,\dots,l_m)\in D(L^m)$ and all
$(\sigma^{-1};\,h_1,\dots,h_m)\in\Sigma_m\wr H$.
Finally, the differential of the composite
\[ D(L^m)\ \xra{\ \bar s \ }\ \Sigma_m\wr G \ \xra{\text{proj}}
(\Sigma_m\wr G) / (\Sigma_m\wr H) \ \xra{\gamma^{-1}} \ (G/H)^m \]
is the identity, so we have indeed defined a slice.
We let
\[ \lambda_{\Sigma_m\wr H}^{\Sigma_m\wr G}\ : \
(\Sigma_m\wr G) /  (\Sigma_m\wr H)_+\ \to \ S^{L^m} \]
denote the collapse map based on the slice $\bar s$.
The composite
\[ (G / H )^m_+ \ \xra{\ \gamma_+\ } \
(\Sigma_m\wr G) /  (\Sigma_m\wr H)_+ \
\xra{\lambda_{\Sigma_m\wr H}^{\Sigma_m\wr G}}  \ S^{L^m}   \]
sends a point $( s(l_1)H,\dots, s(l_m)H)$ with $(l_1,\dots,l_m)\in D(L^m)$ to the point
\[ 
\left( 1-\sqrt{|l_1|^2+\dots+|l_m|^2} \right)^{-1}\cdot(l_1,\dots,l_m)  \ ;
\]
all other points of $(G/H)^m$ are sent to the basepoint at infinity.
A scaling homotopy thus witnesses that the following
diagram commutes up to $(\Sigma_m\wr H)$-equivariant based homotopy:
\[ \xymatrix@C=15mm{ 
(G / H )^m_+ \ar[r]^-{(\lambda_H^G)^{(m)}}\ar[d]^\iso_{\gamma_+} &
(S^L )^{(m)} \ar[d]^\iso \\
(\Sigma_m\wr G) /  (\Sigma_m\wr H)_+
\ar[r]_-{\lambda_{\Sigma_m\wr H}^{\Sigma_m\wr G} } & S^{L^m} }  \]
Now we contemplate the diagram:
\[ \xymatrix@C=8mm@R=7mm{ 
R_0^G( G / H_+) \ar[r]^-{\res^G_H} \ar[d]_{P^m} &
R_0^H( G / H_+) \ar[r]^-{(\lambda_H^G)_*} \ar[d]_{P^m} &
R_0^H(S^L) \ar[d]^{P^m} \\
R_0^{\Sigma_m\wr G}(( G / H ) ^m_+) 
\ar[r]^-{\res^{\Sigma_m\wr G}_{\Sigma_m\wr H}}  \ar[d]^\iso_{ (\gamma_+)_*} &  
R_0^{\Sigma_m\wr H}(( G / H)^m_+) 
\ar[r]^-{ ((\lambda_H^G)^m)_*}   \ar[d]_{(\gamma_+)_*}^\iso & 
R_0^{\Sigma_m\wr H}(( S^L)^{(m)}) \ar[d]^\iso\\
R_0^{\Sigma_m\wr G}( (\Sigma_m\wr G) / (\Sigma_m\wr H) _+ ) 
\ar[r]_-{\res^{\Sigma_m\wr G}_{\Sigma_m\wr H}}  &
R_0^{\Sigma_m\wr H}((\Sigma_m\wr G) / (\Sigma_m\wr H) _+) 
\ar[r]_-{(\lambda_{\Sigma_m\wr H}^{\Sigma_m\wr G})_*} & R_0^{\Sigma_m\wr H}( S^{L^m}) } \]
The two squares on the left commute by compatibility of power operations
with restriction and by naturality of restriction.
The upper right square commutes by naturality of power operations.
The lower right square commutes by the previous paragraph.
Since the upper and lower horizontal composites are 
the respective Wirthm{\"u}ller maps, this proves the proposition.
\end{proof}

A direct consequence of the previous proposition is that power
operations are compatible with dimension shifting 
and degree zero transfer maps.

\begin{cor}\label{cor-power and transfers}
Let $R$ be an ultra-commutative ring spectrum and $H$ 
a closed subgroup of a compact Lie group $G$. Then the following two diagrams commute:
\[ \xymatrix@C=10mm@R=6mm{ 
R_0^H(S^L) \ar[d]_{P^m} \ar[r]^-{\Tr_H^G} &\pi_0^G(R)  \ar[dd]^{P^m}&
\pi_0^H(R)   \ar[dd]_{P^m} \ar[r]^-{\tr_H^G} &\pi_0^G(R) \ar[dd]^{P^m} \\
R_0^{\Sigma_m\wr H}( (S^L)^{(m)})   \ar[d]_\iso &  \\
R_0^{\Sigma_m\wr H}(S^{L^m})  \ar[r]_-{\Tr_{\Sigma_m\wr H}^{\Sigma_m\wr G}} &
\pi_0^{\Sigma_m\wr G}(R) &
\pi_0^{\Sigma_m\wr H}(R)  \ar[r]_-{\tr_{\Sigma_m\wr H}^{\Sigma_m\wr G}} & \pi_0^{\Sigma_m\wr G}(R) 
} \]
\end{cor}
\begin{proof}
By Theorem \ref{thm:Wirth iso}, the external transfer
\[ G\ltimes_H- \ : \ R_0^H(S^L)\ = \ \pi_0^H(R\sm S^L) 
\ \to \ \pi_0^G(R\sm G/H_+) \ = \ R_0^G(G/H_+)\]
is inverse to the Wirthm{\"u}ller isomorphism, 
up to the effect of the antipodal map of $S^L$.
The identification $(S^L)^{(m)}\iso S^{L^m}$ takes $(S^{-\Id_L})^{(m)}$
to $S^{-\Id_{L^m}}$, so the following diagram commutes by naturality of power operations:
\[ \xymatrix@C=12mm@R=7mm{ 
R_0^H(S^L) \ar[r]^-{\varepsilon_L}_-\iso \ar[d]_{P^m} &
R_0^H( S^L) \ar[d]^{P^m} \\
R_0^{\Sigma_m\wr H}( ( S^L)^{(m)})  \ar[d]_\iso &  
R_0^{\Sigma_m\wr H}(( S^L)^{(m)}) \ar[d]^\iso\\
R_0^{\Sigma_m\wr H}( S^{L^m}) 
\ar[r]_-{\varepsilon_{L^m}}^-\iso  & R_0^{\Sigma_m\wr H}(S^{L^m}) } \]
Stacking this diagram next to the commutative diagram of Wirthm{\"u}ller
isomorphisms given in Proposition \ref{prop:power and Wirthmuller},
and reading the composite diagram backwards yields the commutativity 
of the left part of the following diagram:
\[ \xymatrix@C=12mm@R=7mm{ 
R_0^H(S^L) \ar[d]_{P^m} \ar[r]_-{G\ltimes_H-} 
\ar@/^1pc/[rr]^(.3){\Tr_H^G}& R_0^G( (G/H)_+)  \ar[d]^{P^m}
 \ar[r]_-{(p_+)_*} & \pi_0^G(R)   \ar[d]^{P^m}  \\
R_0^{\Sigma_m\wr H}( (S^L)^{(m)})   \ar[d]_\iso &  
R_0^{\Sigma_m\wr G}((G/H)^m_+)\ar[d]^{\gamma_*}_\iso \ar[r]_-{(p^m_+)_*} &
\pi_0^{\Sigma_m\wr G}(R) \ar@{=}[d] \\
R_0^{\Sigma_m\wr H}(S^{L^m})  \ar[r]^-{\Sigma_m\wr G\ltimes_{\Sigma_m\wr H}-} 
\ar@/_1pc/[rr]_(.3){\Tr_{\Sigma_m\wr H}^{\Sigma_m\wr G}}&
R_0^{\Sigma_m\wr G}(  (\Sigma_m\wr G) / (\Sigma_m\wr H)_+ ) \ar[r]^-{(p_+)_*} &
\pi_0^{\Sigma_m\wr G}(R) } \]
The dimension shifting transfer $\Tr_H^G : R_0^H(S^L) \to \pi_0^G(R)$
is the composite of the external transfer 
and the effect of the unique $G$-map $p:G/H\to \ast$.
So the first claim follows from the commutativity of the right part of the
above diagram.
The degree zero transfer is obtained from the dimension shifting transfer
by precomposing with the effect of the map $S^0\to S^L$,
the inclusion of the origin into the tangent representation.
If we raise the inclusion of the origin of $S^L$ to the $m$-th power,
the canonical homeomorphism $(S^L)^{(m)}\to S^{L^m}$ identifies it with
the inclusion of the origin of $S^{L^m}$.
So the power operations are also compatible with degree zero transfers.
\end{proof}

Much of the next result is contained, at least implicitly,
in Greenlees' and May's construction 
of norm maps \cite[Sec.\,7-9]{greenlees-may-completion},
simply because an ultra-commutative ring spectrum is an example of
a `$\mathcal G\mathcal I_*$-FSP' 
in the sense of \cite[Def.\,5.5]{greenlees-may-completion}.

\begin{theorem}\label{thm:pi_0 is global power functor} 
Let $R$ be an ultra-commutative ring spectrum.
The power operations \eqref{eq:power_op_ring_spectra} make the
global functor $\upi_0(R)$ into a global power functor.
\end{theorem}
\begin{proof} The properties~(i) through~(vi) only involve the
multiplication, power operations and restriction maps,
so they are special cases of Proposition \ref{prop:power op unstable} 
for the multiplicative ultra-commutative monoid $\Omega^\bullet R$.
The transfer relation (viii) is proved in Corollary \ref{cor-power and transfers}.
The most involved argument remaining is required for 
the additivity formula (vii), identifying the behavior of power operations on sums.

We first show an external version of the additivity relation.
We consider two equivariant homology classes $x, y\in \pi_0^G(R)$.
Since equivariant homotopy groups take wedges to direct sums there is
a unique class
\[ x\oplus y\ \in \ \pi_0^G(R\sm \{1,2\}_+)  \]
 such that
\[ p^1_*(x\oplus y) \ = \ x \text{\qquad and\qquad}
p^2_*(x\oplus y) \ = \ y \ ,\]
where $p^1,p^2:\{1,2\}_+\to S^0$ are the projections determined by $(p^i)^{-1}(0)=\{i\}$.
We write $\Sigma_{k,m-k}\wr G$ for the image of the homomorphism
$\Phi_{k,m-k}:(\Sigma_k\wr G)\times(\Sigma_{m-k}\wr G)\to\Sigma_m\wr G$
and let
\[ \psi^k\ : \ (\Sigma_m\wr G) / (\Sigma_{k,m-k}\wr G)\ \to \ \{1,2\}^m \]
be the embedding of $(\Sigma_m\wr G)$-sets that sends the preferred coset
$e(\Sigma_{k,m-k}\wr G)$ to the point
\[ (\underbrace{1,\dots,1}_k,\underbrace{2,\dots,2}_{m-k})\ \in \ \{1,2\}^m\ .\]
Here $\Sigma_m\wr G$ acts on $\{1,2\}^m$ through the projection to $\Sigma_m$,
by permutation of coordinates.
The orbits of the $(\Sigma_m\wr G)$-set $\{1,2\}^m$ are precisely 
the images of the maps $\psi^0,\dots,\psi^m$.
We will show the relation
\begin{equation}  \label{eq:power_of_external_sum}
 P^m(x\oplus y)\ = \ \sum_{k=0}^m \, (\psi^k_+)_*( 
(\Sigma_m\wr G)\ltimes_{\Sigma_{k,m-k}\wr G} 
 (P^k(x)\times P^{m-k}(y)))   
\end{equation}
in the group $\pi_0^{\Sigma_m\wr G}(R\sm \{1,2\}^m_+)$.
Since $\pi_0^{\Sigma_m\wr G}(R\sm -)$ is additive on wedges, it suffices to show the relation
after projection to each  $(\Sigma_m\wr G)$-orbit of $\{1,2\}^m$.
So we let 
\[ \bar\psi^k\ : \ \{1,2\}^m_+\ \to \ 
(\Sigma_m\wr G) / (\Sigma_{k,m-k}\wr G)_+\]
be the right inverse to $\psi^k_+$ that sends the other orbits to the
basepoint, i.e., $\bar\psi^k\circ\psi^j_+$ is constant for $j\ne k$.
The relation \eqref{eq:power_of_external_sum} thus follows if we can show
\begin{equation}  \label{eq:component_of_external_sum}
\bar\psi^k_*( P^m(x\oplus y)) \ = \ 
(\Sigma_m\wr G)\ltimes_{\Sigma_{k,m-k}\wr G} (P^k(x)\times P^{m-k}(y))   
\end{equation}
in the group 
\[ \pi_0^{\Sigma_m\wr G}( R\sm (\Sigma_m\wr G)/(\Sigma_{k,m-k}\wr G)_+) \]
for all $0\leq k\leq m$. 
We apply the Wirthm{\"u}ller isomorphism, i.e., the composite
\begin{align*}
\pi_0^{\Sigma_m\wr G}(R\sm (\Sigma_m\wr G)/(\Sigma_{k,m-k}\wr G)_+)  
\ &\xra{\res^{\Sigma_m\wr G}_{\Sigma_{k,m-k}\wr G}}\\
\pi_0^{\Sigma_{k,m-k}\wr G}(R\sm &(\Sigma_m\wr G)/(\Sigma_{k,m-k}\wr G)_+) \
\xra{ (l_k)_*}\ \pi_0^{\Sigma_{k,m-k}\wr G}(R) \ .
\end{align*}
Here $l_k:(\Sigma_m\wr G)/(\Sigma_{k,m-k}\wr G)_+ \to S^0$
is the projection to the preferred coset.
We obtain
\begin{align*}
   (l_k)_*(\res^{\Sigma_m\wr G}_{\Sigma_{k,m-k}\wr G} &(\bar\psi^k_*( P^m(x\oplus y)))) \ = \ 
(l_k\circ\bar\psi^k)_*(\res^{\Sigma_m\wr G}_{\Sigma_{k,m-k}\wr G} ( P^m(x\oplus y))) \\ 
&= \  ( (p^1)^{(k)}\sm (p^2)^{(m-k)})_*( P^k(x\oplus y)\times P^{m-k}(x\oplus y))\\
&= \  P^k( p^1_*(x\oplus y) )\times P^{m-k}(p^2_*(x\oplus y) )\
= \  P^k(x)\times P^{m-k}(y)\ .
\end{align*}
Since the Wirthm{\"u}ller isomorphism is inverse to the external transfer
(compare Theorem \ref{thm:Wirth iso}),
this proves \eqref{eq:component_of_external_sum},
and hence \eqref{eq:power_of_external_sum}. 

Now we obtain the additivity relation by naturality for the fold map
$\nabla:\{1,2\}_+\to S^0$ with $\nabla(1)=\nabla(2)=0$.
Then
\begin{align*}
P^m(x+y)\ &= \ P^m(\nabla_*(x\oplus y)) \ = \ 
(\nabla^{(m)})_*(P^m(x\oplus y)) \\
_\eqref{eq:power_of_external_sum} &= \  
\sum_{k=0}^m \, (\nabla^{(m)}\circ \psi^k_+)_*( (\Sigma_m\wr G)\ltimes_{\Sigma_{k,m-k}\wr G}
 (P^k(x)\times P^{m-k}(y)))   \\
&= \  
\sum_{k=0}^m \, \tr_{\Sigma_{k,m-k}\wr G}^{\Sigma_m\wr G} (P^k(x)\times P^{m-k}(y))\ .  
\end{align*}
The third equation uses that the composite
\[ (\Sigma_m\wr G) / (\Sigma_{k,m-k}\wr G)_+\ \xra{\ \psi^k_+\ } \ \{1,2\}^m_+\
\xra{\nabla^{(m)}} (S^0)^{(m)}\ = \ S^0 \]
sends all of $(\Sigma_m\wr G) / (\Sigma_{k,m-k}\wr G)$ to the non-basepoint.
 \end{proof}

\index{subject}{H-infinity@$H_\infty$-ring spectrum|(}
\begin{rk}[Relation to classical power operations]\label{rk:classical power} 
We owe an explanation how power operations for ultra-commutative ring spectra
refine the classical power operations
defined in the (non-equivariant) cohomology theory represented by an
$H_\infty$-ring spectrum. 
We recall from \cite[I.\S 4]{bruner-may-mcclure-steinberger}
that an $H_\infty$-structure is an algebra structure over the monad
\[ L\mP\ : \ \SH \ \to\ \SH  \]
on the stable homotopy category that can be obtained
by deriving the `symmetric algebra' monad
\[ \mP \ : \ \spec \ \to \ \spec \]
on the category of orthogonal spectra (whose algebras are
commutative orthogonal ring spectra).
This is not the full truth, because 
Bruner, May, McClure and Steinberger use a different model for
the stable homotopy category, so strictly speaking one would have 
to translate the relevant parts of \cite{bruner-may-mcclure-steinberger}
to the context of orthogonal spectra.
If we did that, the derived functor of the $m$-symmetric power of orthogonal spectra
would be modeled by the 
$m$-th {\em extended power}\index{subject}{extended power!of an orthogonal spectrum}
\[ D_m X \ = \ (E\Sigma_m)_+\sm_{\Sigma_m} X^{\sm m} \ .\]
Specifying an $H_\infty$-structure on an orthogonal spectrum $E$ 
thus amounts to specifying morphisms,
in the non-equivariant stable homotopy category,
\[ \mu_m \ : \ D_m E \ \to \ E \]
from the $m$-th extended power to $E$;
the algebra structure over the monad $L\mP$ then translates
into a specific collection of relations among the morphisms $\mu_m$
that are spelled out in \cite[Ch.\,I Def.\,3.1]{bruner-may-mcclure-steinberger}.

For every space $A$, the $H_\infty$-structure gives rise to power operations
\begin{equation}  \label{eq:H_infty_P_m}
 \mathcal P_m \ : \ E^0(A)\ \to \ E^0(B\Sigma_m\times A)   
\end{equation}
in $E$-cohomology 
defined in \cite[Ch.\,I Def.\,4.1]{bruner-may-mcclure-steinberger} as 
the following composite:
\begin{align*}
 E^0(A)\ = \ [\Sigma^\infty_+ A&, E] \ \xra{\ D_m\ } \ [D_m( \Sigma^\infty_+ A), D_m E] \ 
\xra{[D_m(\Sigma^\infty_+ A),\mu_n]} \
 [D_m(\Sigma^\infty_+ A), E] \\ 
&\iso\ [\Sigma^\infty_+ (D_m A), E] \ = \ E^0(D_m A) \ 
\xra{E^0(E\Sigma_m\times_{\Sigma_m}\Delta)} \  E^0(B\Sigma_m\times A)\ .  
\end{align*}
Here $[-,-]$ denotes the morphism group in the stable homotopy category $\SH$,
$D_m A=E\Sigma_m\times_{\Sigma_m}A^m$ is the space level extended power,  
and $\Delta:A\to A^m$ is the diagonal.
Depending on the context, the power operations \eqref{eq:H_infty_P_m}
are often processed further; in favorable cases,
$E^0(B\Sigma_m\times A)$ can be explicitly described as a functor of $E^0(A)$,
and the power operations can be translated into a specific kind of algebraic structure.

Now we let $R$ be an ultra-commutative ring spectrum. 
Then the underlying $H_\infty$-structure
is given by the composite morphism
 \[ D_m R \ = \ (E\Sigma_m)_+\sm_{\Sigma_m} R^{\sm m} \ \to 
\Sigma_m\bs R^{\sm m} \ \xra{\text{mult}} \ R\]
where the first morphism collapses $E\Sigma_m$ to a point and the second
map is induced by the iterated multiplication $R^{\sm m}\to R$.
The definition \eqref{eq:power_op_ring_spectra} of the power operations
on equivariant homotopy groups
directly extends to power operations
\[ P^m \ : \ R^0_G(A) \ \to \ R^0_{\Sigma_m\wr G}(A^m) \]
in the equivariant cohomology of a $G$-space $A$.
The operation $P^m$ takes the class
represented by a based $G$-map $f:S^V\sm A_+\to R(V)$, for some $G$-representation $V$,
to the class of the $(\Sigma_m\wr G)$-map
\[ S^{V^m}\sm A^m_+\ \iso\  (S^V\sm A_+)^{\sm m} \ \xra{\ f^{\sm m}\ } \
R(V)^{\sm m} \ \xra{\ \mu_{V,\dots,V}\ } R(V^m) \ , \]
where $\mu_{V,\dots,V}$
is the iterated, $(\Sigma_m\wr G)$-equivariant multiplication map of $R$.
A forgetful homomorphism
\[ R^0_G(A)\ \to \ R^0(E G\times_G A) \]
is defined as the composite
\begin{align*}
    R^0_G(A)\ = \ \gh{\Sigma^\infty_+ \bL_{G,V}A, R}\ \xra{\ U \ } \
&[U(\Sigma^\infty_+ \bL_{G,V}A), R] \\ 
= \quad &[\Sigma^\infty_+ (E G\times_G A), R] \ = \ R^0(E G\times_G A)  \ ,
\end{align*}
where $U:\GH\to\SH$ is the forgetful functor from the 
global stable homotopy category to the non-equivariant stable homotopy category.
Then the following diagram commutes:
\[ \xymatrix@C=12mm@R=8mm{ 
R^0(A) \ar[d]_{P^m} \ar[r]^-{D_m} &
[D_m( \Sigma^\infty_+ A), D_m R]  \ar[d]^{(\mu_m)_*}\\
R^0_{\Sigma_m}(A^m) \ar[r]^-U \ar[d]_{R^0_{\Sigma_m}(\Delta)} & 
R^0(E\Sigma_m\times_{\Sigma_m} A) \ar[d]^{R^0(E\Sigma_m\times_{\Sigma_m}\Delta)}\\
R^0_{\Sigma_m}(A) \ar[r]_-U  & R^0(B\Sigma_m\times A) 
} \]
The composite through the upper right corner is the $H_\infty$-power operation $\Pc_m$;
so this diagram makes precise in which way power operations 
for ultra-commu\-tative ring spectra refine power operations for $H_\infty$-ring spectra.
\end{rk}

\begin{rk}[$G_\infty$-ring spectra]
As we recalled in the previous remark, non-equivariant power operations on a ring valued
cohomology theory already arise from a weaker structure 
than a commutative multiplication (or equivalently $E_\infty$-multiplication):
all that is needed is an $H_\infty$-structure,
compare \cite[I.\S 4]{bruner-may-mcclure-steinberger}.
This suggests a global analog of an $H_\infty$-structure that, for lack of better name,
we call a {\em $G_\infty$-structure}.\index{subject}{G-infinity@$G_\infty$-ring spectrum}
For the formal definition we exploit the 
global model structure for ultra-commutative ring spectra,
to be established in Theorem \ref{thm:global ultra-commutative} below.
This model structure is obtained by lifting the positive global model structure 
on the category of orthogonal spectra 
(see Proposition \ref{prop:positive global spectra})
along the forgetful functor. In particular, the
free and forgetful functor form a Quillen adjoint functor pair:
\[ \xymatrix@C=10mm{ \mP\ : \ \spec\ \ar@<.4ex>[r] & 
  \ ucom \ : \ U \ar@<.4ex>[l]} \]
Every Quillen adjoint pair between model categories derives to
an adjoint functor pair between the homotopy categories, 
see \cite[I.4, Thm.\,3]{Q} or \cite[Lemma 1.3.10]{hovey-book}. 
In our situation this provides the derived adjunction:
\[ \xymatrix@C=10mm{ L_{\gl}\mP \ : \ \GH \ \ar@<.4ex>[r]   & 
  \ \Ho(ucom) \ : \ \Ho(U) \ar@<.4ex>[l]  } \]
Since the forgetful functor is fully homotopical, it does not even have to be derived.
The composite
\[ \mG\ = \ \Ho(U)\circ L_{\gl}\mP\ : \ \GH \ \to\ \GH  \]
is then canonically a monad on the global stable homotopy category, 
whose algebras we call $G_\infty$-ring spectra.

The underlying non-equivariant homotopy type of a $G_\infty$-ring spectrum
comes with an $H_\infty$-structure, so we arrive at a square of 
forgetful functors between categories of structured ring spectra
with different degrees of commutativity:
\[ \xymatrix{ 
\Ho(ucom) \ar[r]\ar[d] & 
(\text{$G_\infty$-ring spectra}) \ar[d] \\
\Ho^{\pi_*\text{-iso}}(\text{comm.\,ring spectra}) \ar[r] & 
(\text{$H_\infty$-ring spectra}) } \]
We emphasize that, like $H_\infty$-ring spectra,
the category of $G_\infty$-ring spectra is not the homotopy category
of any natural model category.
\end{rk}
\index{subject}{H-infinity@$H_\infty$-ring spectrum|)}

\begin{eg}[Units of an ultra-commutative ring spectrum]\index{subject}{units!of an ultra-commutative ring spectrum}
In Example \ref{eg:R^naive times}\index{subject}{units!of an orthogonal monoid space!naive}
we defined the naive units of an orthogonal monoid space.
When $R$ is an orthogonal ring spectrum, then the naive units
of the multiplicative orthogonal monoid space $\Omega^\bullet R$ satisfy
\[ \pi_0^G( (\Omega^\bullet R)^{n\times}) \ = \ 
\{ x \in \pi_0^G(R)\ | \ \text{$\res^G_e(x)$ is a unit in $\pi_0^e(R)$}\}\ ,\]
the multiplicative submonoid of $\pi_0^G(R)$ of elements that become
invertible when restricted to the trivial group.
One should beware that these naive units may contain non-invertible elements,
i.e., the orthogonal monoid space $(\Omega^\bullet R)^{n\times}$ need not be group-like.

When the ring spectrum $R$ is ultra-commutative, then there is a more
refined construction 
\[ G L_1(R) \ = \ (\Omega^\bullet R)^\times\ , \]
the {\em global units}\index{symbol}{$G L_1(R)$ - {global units of the ultra-commutative ring spectrum $R$}} 
of $R$. Indeed, if $R$ is ultra-commutative, then $\Omega^\bullet R$ is an
ultra-commutative monoid, so we can form the `true' global units,
the homotopy fiber of the multiplication morphism, 
see Construction \ref{con:R^times}.
Then $G L_1(R)$ is a group-like ultra-commutative monoid 
and for every compact Lie group $G$,
\[ \pi_0^G( G L_1(R) ) \ = \ 
(\pi_0^G(R))^\times \ , \]
the multiplicative submonoid of units of the commutative ring $\pi_0^G(R)$.
Moreover, the power operations in $\upi_0(R)$ correspond to the power operations
in $\upi_0(G L_1(R))$.

In the non-equivariant context,
$G L_1(R)$ is an infinite loop space, i.e., weakly equivalent
to the 0-th space of an $\Omega$-spectrum of units.
This fact has a global generalization as follows.
As we hope to explain elsewhere, every ultra-commu\-tative monoid $M$ has a
{\em  global delooping}\index{subject}{global delooping} $\bB M$,
an orthogonal spectrum that is a $\Fin$-global $\Omega$-spectrum. 
It also comes with a natural morphism of orthogonal spaces
$M\to \Omega^\bullet (\bB M)$ that is a $\Fin$-global equivalence 
whenever $M$ is group-like.
Since $G L_1(R)$ is an ultra-commutative monoid, it
has a global delooping\index{subject}{global delooping} 
\[ g l_1(R)\ = \ \bB( G L_1(R) ) \ .\]
Since the ultra-commutative monoid $G L_1(R)$ is group-like, the morphism 
\[ \xi \ : \  G L_1(R) \ \to \   \Omega^\bullet ( \bB (G L_1(R) ) )\ = \ 
 \Omega^\bullet ( g l_1(R) ) \]
is a $\Fin$-global equivalence of orthogonal spaces. 
For every compact Lie group $G$, this induces a map
\begin{align*}
    (\pi_0^G( R))^\times \ &\iso  \  \pi_0^G(G L_1(R)) \\ 
&\iso  \ \pi_0^G( ( G L_1(R))^c) \ \xra{\ \upi_0(\xi)\ } \  
\pi_0^G(\Omega^\bullet ( g l_1(R) )) \ = \  \pi_0^G(g l_1(R)) \ .
\end{align*}
These maps are compatible with restriction
along continuous homomorphisms and they are bijective whenever $G$ is finite.
Moreover, the maps take the multiplication respectively norm operations 
in $(\upi_0( R))^\times$ 
to addition respectively finite index transfers in $\upi_0(g l_1(R))$.
\end{eg}

\begin{rk}[Picard groups] The global units $G L_1(R)$ of an 
ultra-commu\-tative ring spectrum $R$ ought to have an interesting delooping $\pic(R)$\index{symbol}{$\pic(R)$ - {Picard monoid of an ultra-commutative ring spectrum}}
\index{subject}{Picard group! of an ultra-commutative ring spectrum}
that records the information about invertible modules over the
equivariant ring spectra underlying $R$. At present I have no construction
of this delooping as an ultra-commutative monoid, but I 
describe the evidence for expecting its existence.

For every compact Lie group $G$ the underlying orthogonal $G$-ring spectrum $R_G$
of $R$ has a symmetric monoidal model category of modules,
i.e., orthogonal $G$-spectra with an action by $R$ (where $G$ acts trivially on $R$).
The equivalences we consider here are $R$-linear morphisms 
that are $\upi_*$-isomorphisms of underlying orthogonal $G$-spectra;
the construction of such a symmetric mo\-noidal model 
category can be found in \cite[III Thm.\,7.6]{mandell-may}.
We let
\[ \Pic(R)(G)\ = \ \Pic(\Ho(R_G\text{-mod})) \]
be the resulting Picard group, i.e., the set of isomorphism
classes, in the homotopy category of $R_G$-modules,
of objects that are invertible under the derived smash product.
For a continuous group homomorphism $\alpha:K\to G$ the restriction functor
$\alpha^*:R_G\text{-mod}\to R_K\text{-mod}$
derives to a strong symmetric monoidal functor
\[ R\alpha^*\ : \ \Ho(R_G\text{-mod})\ \to \ \Ho(R_K\text{-mod})\ . \]
So $R\alpha^*$ preserves invertibility and induces a group homomorphism
\[ \alpha^*\ : \ \Pic(R)(G)  \to \ \Pic(R)(K)\ . \]
For a second homomorphism $\beta:L\to K$ the functors
$(R\beta^*)\circ(R\alpha^*)$ and $R(\alpha\circ\beta)^*$
are naturally isomorphic.
Moreover, for every element $g\in G$ the restriction functor $R (c_g)^*$
is naturally isomorphic to the identity functor of
$\Ho(R_G\text{-mod})$, via left multiplication by $g$.
So the assignment $G\mapsto \Pic(R)(G)$ becomes a functor
\[ \Pic(R)\ : \ \Rep^{\op} \ \to \ \Ab \ . \]
But the ultra-commutativity gives more. For every finite index subgroup $H\leq G$,
the norm construction of Hill, Hopkins and Ravenel
derives to a strong symmetric monoidal functor
\[ N_H^G \ : \ \Ho(R_H\text{-mod})\ \to \ \Ho(R_G\text{-mod})\ , \]
compare \cite[Prop.\,B.105]{HHR-Kervaire}.
So also $N_H^G$ preserves invertibility and induces a group homomorphism
\[ N_H^G\ : \ \Pic(R)(H)  \to \ \Pic(R)(G)\ . \]
These norm maps are transitive and they extend the abelian $\Rep$-monoid
to a global power monoid $\Pic(R)$. 

We expect that there is a `natural' ultra-commutative monoid $\pic(R)$
such that $\upi_0(\pic(R))\iso \Pic(R)$ as global power monoids
and such that $\Omega (\pic(R))$ is globally equivalent, as an ultra-commutative
monoid, to $G L_1(R)$. The $G$-fixed points $(\pic(R)(\Uc_G))^G$
ought to have the homotopy type, as an $E_\infty$-space,
of the nerve of the category of invertible $R_G$-modules
and $\upi_*$-isomorphisms. Despite the strong evidence for its existence,
I cannot presently construct $\pic(R)$ as an ultra-commutative monoid
in our formalism.
\end{rk}

\begin{eg}[Free global power functors]\index{subject}{global power functor!free}\label{eg:free global power functor}  
For a compact Lie group $K$ we construct a 
{\em free global power functor} $C_K$ generated by $K$. 
The underlying global functor is
\[ C_K \ = \ {\bigoplus}_{m\geq 0}\ \bA(\Sigma_m\wr K,-) \ ,\]
the direct sum of the global functors represented by the wreath products 
$\Sigma_m\wr K$, including the trivial group $\Sigma_0\wr K=e$.
The multiplication $\mu: C_K\Box C_K\to C_K$
that makes this into a global Green functor restricted to the $(m,n)$-summand
is the morphism
\[  \bA(\Sigma_m\wr K,-)\Box \bA(\Sigma_n\wr K,-) \ \to \ C_K \]
that corresponds, via the universal property of the box product, 
to the bimorphism with $(G,G')$-component
\begin{align*}
   \bA(\Sigma_m\wr K,G)\Box \bA(\Sigma_n\wr K,G') \ \xra{\ \times \ } \ 
&\bA((\Sigma_m\wr K)\times(\Sigma_n\wr K),G\times G') \\ 
\xra{\bA(\Phi^*_{m,n},G\times G')} \ 
&\bA(\Sigma_{m+n}\wr K,G\times G') \ \xra{\text{incl}} \ C_K(G\times G')\ ;
\end{align*}
here $\Phi^*_{m,n}$
is the restriction map associated to the embedding \eqref{eq:wreath_sum spectra}   
\[  \Phi_{m,n}\ : \ (\Sigma_m\wr K)\times(\Sigma_n\wr K) \ \to \ \Sigma_{m+n}\wr K \ .
 \]
The multiplication is associative because
\begin{align*}
  \Phi_{k+m,n}\circ(\Phi_{k,m}\times(\Sigma_n\wr K))\ &=\ 
\Phi_{k,m+n}\circ((\Sigma_k\wr K)\times\Phi_{m,n}) \ : \\
&(\Sigma_k\wr K)\times(\Sigma_m\wr K)\times (\Sigma_n\wr K)\ 
\to\ \Sigma_{k+m+n}\wr K\ .
\end{align*}
The multiplication is commutative because the group homomorphisms
\[ \Phi_{m,n}\ , \
\Phi_{n,m}\circ\tau_{\Sigma_m\wr K, \Sigma_n\wr K} \ : \
(\Sigma_m\wr K)\times(\Sigma_n\wr K) \ \to\ \Sigma_{m+n}\wr K\]
are conjugate, so they represent the same morphism in
$\bA((\Sigma_m\wr K)\times(\Sigma_n\wr K),\ \Sigma_{m+n}\wr K)$.
The unit is the inclusion $\bA(e,-)\to C_K$ of the summand indexed by $m=0$.

The global Green functor $C_K$ can be made into a
global power functor in a unique way such that the relation
\[ P^m(1_K) \ = \ 1_{\Sigma_m\wr K}\]
holds in the $m$-th summand of $C_K(\Sigma_m\wr K)$, 
where $1_K\in \bA(K,K)$ and $1_{\Sigma_m\wr K}\in\bA(\Sigma_m\wr K,\Sigma_m\wr K)$
are the identity operations.
Indeed, $C_K$  is generated as a global functor by the classes 
$1_{\Sigma_m\wr K}$ for all $k\geq 0$, so there is at most one such global power structure,
and every morphism of global power functors out of $C_K$ is determined
by its value on the class $1_K$.
The existence of a global power structure on $C_K$ with this property
could be justified purely algebraically, but we show it
by realizing $C_K$ by an ultra-commutative ring spectrum.
As we shall make precise in Proposition \ref{prop:what C_K represents}~(ii)
below, $C_K$ is indeed {\em freely} generated, as a global power functor,
by the class $1_K$.

The unreduced suspension spectrum\index{subject}{global classifying space}  
\[ \Sigma^\infty_+ \mP(B_{\gl}K ) \ \iso \ 
{\bigvee}_{m\geq 0} \,\Sigma^\infty_+ B_{\gl}(\Sigma_m\wr K)\]
of the free ultra-commutative monoid
(compare~Example \ref{eg:free umon})
generated by a global classifying space of $K$
is an ultra-commutative ring spectrum. 
According to Proposition \ref{prop:B_gl represents}, 
its 0-th homotopy group global functor is given additively by
\[ \upi_0 ( \Sigma^\infty_+ \mP(B_{\gl}K ) ) \ \iso \ 
{\bigoplus}_{m\geq 0} \, \upi_0( \Sigma^\infty_+ B_{\gl}(\Sigma_m\wr K) ) \ \iso \ 
{\bigoplus}_{m\geq 0} \, \bA(\Sigma_m\wr K,-)  \ . \]
Under this isomorphism, the stable tautological class 
$e_K\in \pi_0^K ( \Sigma^\infty_+ B_{\gl}K)$ maps to the 
generator $1_K\in \bA(K,K)$.
The class $e_K$ is the stabilization of the unstable tautological class
$u_K\in \pi_0^K( B_{\gl} K)$, whose $m$-th power is $[m](u_K)= u_{\Sigma_m\wr K}$
in $\pi_0^{\Sigma_m\wr K}(\mP (B_{\gl} K))$, see \eqref{eq:P^m_of_u_G}.
The stabilization map $\sigma:\upi_0(B_{\gl}K)\to\upi_0(\Sigma^\infty_+ B_{\gl}K)$
commutes with power operations, so this shows that
\[ P^m(e_K) \ = \  P^m(\sigma^K(u_K))\ = \ 
\sigma^{\Sigma_m\wr K}([m](u_K)) \ = \ 
\sigma^{\Sigma_m\wr K}(u_{\Sigma_m\wr K}) \ = \ e_{\Sigma_m\wr K}\]
in the group $\pi_0^{\Sigma_m\wr K}(\Sigma^\infty_+ \mP(B_{\gl}K))$.
\end{eg}

\Danger We warn the reader that the previous example
of the free global power functor $C_K$ is {\em not}
the symmetric algebra, with respect to the box product of global functors,
of the represented global functor $\bA(K,-)$. The issue is that the global functors
\[ \bA(K,-)^{\Box m}/\Sigma_m\ \iso\ \bA(K^m,-)/\Sigma_m
\text{\qquad and\qquad}  \bA(\Sigma_m\wr K,- )\]
are typically not isomorphic.
The restriction map $\res^{\Sigma_m\wr K}_{K^m}\in\bA(\Sigma_m\wr K,K^m)$
induces a morphism of represented global functors
\[ -\circ \res^{\Sigma_m\wr K}_{K^m} \ : \  \bA(K^m,-) \ \to  \bA(\Sigma_m\wr K,- )\]
that equalizes the $\Sigma_m$-action on the source because
every permutation of the factors of $K^m$ becomes an inner automorphism
in $\Sigma_m\wr K$.
So the morphism factors over a morphism of global functors
\[ \bA(K^m,-)/\Sigma_m \ \to \ \bA(\Sigma_m\wr K,- )\]  
which, however, is generally {\em not} an isomorphism
(already for $K=e$ and $m=2$).
The box product symmetric algebra
\[ {\bigoplus}_{m\geq 0} \, \bA(K,-)^{\Box m}/\Sigma_m \ \iso \ 
 {\bigoplus}_{m\geq 0} \, \bA(K^m,-)/\Sigma_m  \]
also has a universal property: it is the free global Green functor
generated by $K$. However, this box product symmetric algebra
does not seem to have natural power operations.

\section{Comonadic description of global power functors}
\label{sec:algebra global power}

In this section we show that the category of global power functors 
is both monadic and comonadic over the category of global Green functors.
We introduce the functor of exponential sequences 
and make it into a comonad on the category of global Green functors.
For a global Green functor $R$ and a compact Lie group $G$,
Construction \ref{con:exp for GPF} introduces the commutative ring $\exp(R;G)$ 
of exponential sequences.
Construction \ref{con:Green on exp} connects these rings for
varying Lie groups $G$, making the entire data into a new global Green functor $\exp(R)$,
compare Proposition \ref{prop:exp(R) is global Green}.
Theorem \ref{thm:exp is comonad} extends the functor of
exponential sequences to a comonad on the category of global Green functors,
and Theorem \ref{thm:power comonadic} shows that the category of its coalgebras 
is isomorphic to the category of global power functors.
Proposition \ref{prop:GPF are monadic} shows that the category of
global power functors is also monadic over the category of global Green functors.
A formal consequence is that the category of global power functors
has all limits and colimits, and that they are created in the category of
global Green functors.
So the relationship of global power functors to global Green functors
is formally similar to the situation for $\lambda$-rings,
which are both monadic and comonadic over the category of commutative rings.
When restricted to finite groups, most of the results about the comonad
of exponential sequences are contained in the PhD thesis
of J.\,Singer \cite{singer}, a former student of the author.
Also for finite groups (as opposed to compact Lie groups),
this comonadic description has independently been obtained by Ganter \cite{ganter-power}.

Proposition \ref{prop:universal prop localization} discusses localization
of global Green functors at a multiplicative subset of the underlying ring;
while the structure of global Green functor always `survives localization',
this is not generally true for power operations.
Theorem \ref{thm:localize global power functor}
exhibits a necessary and sufficient condition so that a localization
of a global power functor inherits power operations.
For localization at a set of integer primes this condition is always satisfied
(Example \ref{eg:rationalize GPF}), so global power functors can in particular 
be rationalized.

\begin{construction}\label{con:exp for GPF}
We let $R$ be a global Green functor and $G$ a compact Lie group. We let
\[ \exp(R;G) \ \subset \ {\prod}_{m\geq 0}\ R(\Sigma_m\wr G)\]
be the set of {\em exponential sequences},\index{subject}{exponential sequence!in a global Green functor}\index{symbol}{$\exp(R)$ - {global Green functor of exponential sequences}}
i.e., of those families $(x_m)_m$ that satisfy $x_0=1$ in $R(\Sigma_0\wr G)=R(e)$ and
\[ \Phi_{k,m-k}^*(x_m) \ = \ x_k\times x_{m-k}\]
in $R((\Sigma_k\wr G)\times(\Sigma_{m-k}\wr G))$ for all $0 < k < m$,
where 
\[  \Phi_{k,m-k}\ : \ (\Sigma_k\wr G)\ \times\ (\Sigma_{m-k}\wr G)\ \to\ \Sigma_m\wr G  \]
is the monomorphism \eqref{eq:wreath_sum spectra}.   
We define a multiplication on the set $\exp(R;G)$ by
coordinatewise multiplication in the rings $R(\Sigma_m\wr G)$, i.e.,
\[ (x \cdot y )_m  \ = \ x_m\cdot y_m\ . \]
We introduce another binary operation $\oplus$ on $\exp(R;G)$ by
\[ ( x \oplus y )_m  \ = \ \sum_{k=0}^m \ \tr_{k,m-k} (x_k\times y_{m-k}) \ ,    \]
where $x=(x_m)$, $y=(y_m)$ and
$\tr_{k,m-k}:R( (\Sigma_k\wr G)\times(\Sigma_{m-k}\wr G))\to R(\Sigma_m\wr G)$ 
is the transfer associated to the monomorphism $\Phi_{k,m-k}$. 
\end{construction}

\begin{prop}\label{prop:exp(R,G) is ring} 
  Let $R$ be a global Green functor and $G$ a compact Lie group.
  Then the addition $\oplus$ and the componentwise multiplication
  make the set $\exp(R;G)$ of exponential sequences into a commutative ring.
\end{prop}
\begin{proof}
If $x$ and $y$ are exponential sequences, then the relation
\begin{align*}
  \Phi_{k,m-k}^*( x_m \cdot y_m )\ &= \ 
  \Phi_{k,m-k}^*( x_m ) \cdot  \Phi_{k,m-k}^*(  y_m )\\ 
&= \   (x_k\times  x_{m-k})\cdot  (y_k\times y_{m-k})\ 
= \ (x_k\cdot y_k)\times (x_{m-k}\cdot y_{m-k})
\end{align*}
holds in $R( (\Sigma_k\wr G)\times(\Sigma_{m-k}\wr G))$;
so the product $x\cdot y$ is again exponential.
The product is associative and commutative since 
all the multiplications in the rings $R(\Sigma_m\wr G)$ 
have this property. The exponential sequence $(1)_{m \geq 0}$ is a multiplicative unit.

Now we show that the sum, with respect to $\oplus$,
of two exponential sequences is again exponential.
The key step in this verification is an application of the double coset formula,
for which we need to understand
the $(\Sigma_i\times\Sigma_{m-i})$-$(\Sigma_k\times\Sigma_{m-k})$-double cosets
inside $\Sigma_m$.
We parametrize these double cosets by pairs $(a,b)$ of natural numbers satisfying
\begin{equation}  \label{eq:(a,b)_constraints}
 0\leq a \leq i \ , \quad  0 \leq b\leq m-i  \text{\qquad and\qquad} a + b = k \ .  
\end{equation}
For each such pair we define a permutation $\chi(a,b)\in\Sigma_m$ by
\[ \chi(a,b)(j) \ = \ 
\begin{cases}
\  j      & \text{ for $1\leq j\leq a$, }\\
\  j-a+i  & \text{ for $a+1\leq j\leq a+b$, }\\
\  j-b  & \text{ for $a+b+1\leq j\leq i+b$, }\\
\  j & \text{ for $i+b+1\leq j\leq m$. }
\end{cases} \]
In other words, $\chi(a,b)$ is the unique $(k,m-k)$-shuffle such that
\[ \chi(a,b)(\{1,\dots,k\}) \ = \  \{1,\dots, a\}   \ \cup \ \{i+1,\dots, i+b\} \ .\]
The permutations $\chi(a,b)$ form a set of double coset representatives
for the subgroups $\Sigma_i\times\Sigma_{m-i}$ and $\Sigma_k\times\Sigma_{m-k}$
inside $\Sigma_m$, for all pairs $(a,b)$ subject to \eqref{eq:(a,b)_constraints}.

When applying the double coset formula we will need the relations
\[ (\Sigma_i\times\Sigma_{m-i})^{\chi(a,b)}\ \cap\ (\Sigma_k\times\Sigma_{m-k}) \ = \ 
\Sigma_a\times\Sigma_b\times\Sigma_{i-a}\times\Sigma_{m-i-b} \]
and
\[ (\Sigma_i\times\Sigma_{m-i})\ \cap\ {^{\chi(a,b)} (\Sigma_k\times\Sigma_{m-k})} \ = \ 
\Sigma_a\times\Sigma_{i-a}\times\Sigma_b\times\Sigma_{m-i-b} \ .\]
Thus the double coset formula becomes
\begin{align}\label{eq:double coset double symmetric}
 \Phi_{i,m-i}^*&(\tr_{k,m-k} (x\times y ))  \\ 
&= \ 
\sum_{a,b}
\tr_{(\Sigma_i\times\Sigma_{m-i})\cap{^{\chi(a,b)}(\Sigma_k\times\Sigma_{m-k})}}^{\Sigma_i\times\Sigma_{m-i}} \left( 
\chi(a,b)_\star \left(
\res_{(\Sigma_i\times\Sigma_{m-i})^{\chi(a,b)}\cap(\Sigma_k\times\Sigma_{m-k})}^{\Sigma_k\times\Sigma_{m-k}}(x\times y) \right)\right)\nonumber  \\ 
 &= \ 
\sum_{a,b}
\tr_{\Sigma_a\times\Sigma_{i-a}\times\Sigma_b\times\Sigma_{m-i-b}}^{\Sigma_i\times\Sigma_{m-i}} \left( 
\chi(a,b)_\star \left( \res^{\Sigma_k}_{\Sigma_a\times\Sigma_b}(x)\times\res^{\Sigma_{m-k}}_{\Sigma_{i-a,m-i-b}}( y)\right)\right) \nonumber
\end{align}
The two sums run over all pairs $(a,b)$ of natural numbers 
satisfying \eqref{eq:(a,b)_constraints}.
Now we consider exponential sequences $x,y\in\exp(R;G)$ and calculate
\begin{align*}
\Phi_{i,m-i}^*( (x\oplus y)_m) \ &= \ 
\sum_{k=0}^m \, \Phi_{i,m-i}^*(\tr_{k,m-k} (x_k\times y_{m-k})) 
\\ 
_\eqref{eq:double coset double symmetric} \ &= \ \sum_{a,b}
\tr_{\Sigma_i\times\Sigma_{i-a}\times\Sigma_b\times\Sigma_{m-i-b}}^{\Sigma_i\times\Sigma_{m-i}} \left( 
\chi(a,b)_\star ( x_a\times x_b\times y_{i-a}\times y_{m-i-b}) \right) \\ 
 &= \ \sum_{a,b}
\tr_{a,i-a} (x_a\times y_{i-a})\times  \tr_{b,m-i-b} (x_b\times y_{m-i-b}) \\
&= \  (x\oplus y)_i\times (x\oplus y)_{m-i} \ .
\end{align*}
Here the last two sums run over all pairs $(a,b)$ of natural numbers satisfying
$0\leq a\leq i$ and $0\leq b\leq m-i$.
This shows that the sequence $x\oplus y$ is again exponential.

The following square of group monomorphisms commutes:
\[ \xymatrix@C=15mm{ 
(\Sigma_j\wr G)\times (\Sigma_k\wr G)\times(\Sigma_l\wr G) 
\ar[r]^-{(\Sigma_j\wr G)\times \Phi_{k,l}} 
\ar[d]_{\Phi_{j,k}\times (\Sigma_l\wr G)} &
(\Sigma_j\wr G) \times (\Sigma_{k+l}\wr G) \ar[d]^{\Phi_{j,k+l}}\\
(\Sigma_{j+k}\wr G)\times(\Sigma_l\wr G) 
\ar[r]_-{\Phi_{j+k,l}} & \Sigma_{j+k+l}\wr G} \]
So for all $x\in R(\Sigma_j\wr G)$, $y\in R(\Sigma_k\wr G)$,
and $z\in R(\Sigma_l\wr G)$, the relation
\[ \tr_{j,k+l}(x\times \tr_{k,l}(y\times z)) \ = \ 
\tr_{j+k,l}(\tr_{j,k}(x\times y)\times z) \]
holds in the group $R(\Sigma_{j+k+l}\wr G)$. 
By unraveling the definitions, this becomes the associativity 
of the operation $\oplus$.

Also, the following square of group monomorphisms commutes:
\[ \xymatrix@C=15mm{ 
(\Sigma_k\wr G)\times(\Sigma_l\wr G) 
\ar[r]^-{\Phi_{k,l}} 
\ar[d]_{\text{twist}} & 
\Sigma_{k+l}\wr G \ar[d]^{c_\chi} \\
(\Sigma_l\wr G)\times(\Sigma_k\wr G) 
\ar[r]_-{\Phi_{l,k}} & \Sigma_{l+k}\wr G} \]
Here $\chi=(\chi_{k,l};1,\dots,1)$, 
for the shuffle permutation $\chi_{k,l}\in\Sigma_{k+l}$.
So for all $x\in R(\Sigma_k\wr G)$ and $y\in R(\Sigma_l\wr G)$, the relation
\[ \tr_{l,k}(y\times x) \ = \ 
\chi_\star( \tr_{k,l}(x\times y)) \ = \ \tr_{k,l}(x\times y) \]
holds in the group $R(\Sigma_{k+l}\wr G)$. 
By unraveling the definitions, this implies the commutativity 
of the operation $\oplus$.

The sequence $\underline 0$ with ${\underline 0}_0=1$ and
${\underline 0}_m = 0$ for $m\geq 1$ is a neutral element for $\oplus$. 
Given an exponential sequence $x$ we define a sequence $y$ inductively by
$y_0=1$ and by
\[  y_m \ = \  - {\sum}_{k=1}^m \ \tr_{k,m-k} (x_k\times y_{m-k})  \]
for $m\geq 1$.
To see that the sequence $y$ is again exponential,
we show the relation
\[ \Phi_{i,m-i}^*(  y_m) \ = \ y_i\times y_{m-i} \]
by induction on $m$. The induction starts with $m=1$, where there is nothing to show.
Now we assume that $m\geq 2$; then 
\begin{align*}
- \Phi_{i,m-i}^*(  y_m) \ &= \ 
 {\sum}_{k=1}^m \, \Phi_{i,m-i}^*(\tr_{k,m-k} (x_k\times y_{m-k})) \\ 
_\eqref{eq:double coset double symmetric} \ &= \ 
\sum_{ a+b\geq 1}
\tr_{\Sigma_a\times\Sigma_{i-a}\times\Sigma_b\times\Sigma_{m-i-b}}^{\Sigma_i\times\Sigma_{m-i}} \left( 
\chi(a,b)_\star ( x_a\times x_b\times y_{i-a}\times y_{m-i-b} )\right) \\ 
 &= \ \sum_{a,b\geq 1}
\tr_{a,i-a} (x_a\times y_{i-a})\times  \tr_{b,m-i-b} (x_b\times y_{m-i-b}) \\
&+ \ \sum_{a=1}^i \tr_{a,i-a}(x_a\times y_{i-a}) \times y_{m-i}
\ +\  \sum_{b=1}^{m-i} y_i\times \tr_{b,m-i-b}(x_b\times y_{m-i-b})\\
 &= \ 
 y_i \times y_{m-i}\ -\ y_i\times y_{m-i}\ -\  y_i \times  y_{m-i} \ = \  
- y_i\times y_{m-i} \ .
\end{align*}
In the sums, $(a,b)$  runs over all pairs of natural numbers satisfying
$a\leq i$ and $b\leq m-i$ plus the conditions attached to the summation symbols.
The second equation uses the inductive hypothesis.
This shows that the sequence $y$ is again exponential.
The relation $x\oplus y=\underline 0$ holds by construction,
so $y$ is an inverse of $x$ with respect to $\oplus$.
Altogether this shows that $\exp(R;G)$ is an abelian group under $\oplus$. 

It remains to show distributivity of multiplication over $\oplus$.
We let $x,y, z\in\exp(R;G)$ be exponential sequences. Then
\begin{align*}
( ( x\oplus y)\cdot z)_m  \ &= \ 
\left({\sum}_{k=0}^m\, \tr_{k,m-k}(x_k\times y_{m-k}) \right)\cdot z_m \\
&= \ 
{\sum}_{k=0}^m\, \tr_{k,m-k} ( (x_k\times y_{m-k}) \cdot \Phi_{k,m-k}^*(z_m) ) \\
&= \ 
{\sum}_{k=0}^m\, \tr_{k,m-k} ( (x_k\times y_{m-k}) \cdot (z_k\times  z_{m-k} )) \\
&= \ 
{\sum}_{k=0}^m\, \tr_{k,m-k} ( (x_k\cdot z_k ) \times (y_{m-k}\cdot z_{m-k}) ) \ = \
( ( x\cdot z)\oplus (y\cdot z))_m  \ .
\end{align*}
The second equation is reciprocity in the global Green functor $R$.
\end{proof}

The ring $\exp(R;G)$ is covariantly functorial in $R$:
given a morphism of global Green functors $\varphi:R\to S$, we define a map
\[ \exp(\varphi;G) \ : \ \exp(R;G) \ \to \  \exp(S;G) \text{\quad by\quad}
 \left( \exp(\varphi;G)(x) \right)_m \ = \ \varphi(x_m)\ . \]
We omit the straightforward verifications that
for $x\in\exp(R;G)$ the sequence $\exp(\varphi;G)(x)$ is again exponential,
that $\exp(\varphi;G)$ is a ring homomorphism,
and that $\exp(-;G)$ is compatible with identity and composition.

Now we turn to the functoriality of $\exp(R;G)$ in the
Lie group variable, which is slightly more subtle. 
The ultimate aim is to make $\exp(R;-)$ into another global Green functor. 
Proposition \ref{prop:exp(R,G) is ring} already provides the ring structures
on the values, so the missing data is to give $\exp(R;-)$
the structure of a global functor. Our approach is indirect and not based on
explicit formulas; rather we exploit that the functor $R\mapsto \exp(R;G)$
is representable by the underlying global Green functor
of the free global power functor $C_G$ discussed in 
Example \ref{eg:free global power functor}.
The next proposition will be used for establishing this representability property,
see Proposition \ref{prop:what C_K represents} below.

\begin{prop}\label{prop:good elements}
Let $R$ and $S$ be global power functors and $f:R\to S$
a morphism of global Green functors.
Then the collection of elements $x\in R(G)$, for varying compact Lie groups $G$,
that satisfy
\[ f(P^m(x))\ = \ P^m(f(x)) \]
in $S(\Sigma_m\wr G)$  for all $m\geq 0$ form a global power subfunctor of $R$.
\end{prop}
\begin{proof}
  This is a straightforward consequence of the various relations enjoyed by
  power operations. For the course of this proof we call a pair $(G,x)$
  consisting of a compact Lie group $G$ and an element $x\in R(G)$ {\em good}
  if the relation $f(P^m(x))= P^m(f(x))$ holds for all $m\geq 0$.
  Then the naturality relation $P^m\circ \alpha^* = (\Sigma_m\wr\alpha)^*\circ P^m$
  for a continuous homomorphism $\alpha:K\to G$ implies
  \begin{align*}
    f(P^m(\alpha^*(x))) \ &= \ f((\Sigma_m\wr\alpha)^*(P^m(x))) \ = \ 
    (\Sigma_m\wr\alpha)^*(f^*(P^m(x))) \\ 
    &= \ (\Sigma_m\wr\alpha)^*(P^m(f^*(x))) \ = \   P^m(\alpha^*(f^*(x))) \ 
    = \   P^m(f^*(\alpha^*(x))) \ . 
  \end{align*}
  So if $(G,x)$ is good, then so is $(K,\alpha^*(x))$.
  By the analogous calculation, the transfer relation
  $P^m\circ \tr_H^G = \tr_{\Sigma_m\wr H}^{\Sigma_m\wr G}\circ P^m$ for a closed subgroup $H$
  of $G$ implies that if $(H,x)$ is good, then so is $(G,\tr_H^G(x))$.
  
  The relation $P^m(1)=1$ and the fact $f(1)=1$ show that $(e,1)$ is good.
  Since good pairs are closed under restriction, the pair $(G,1)$ 
  is good for every compact Lie group $G$.
  The multiplicativity $P^m(x\cdot  y)= P^m(x)\cdot P^m(y)$
  and the hypothesis that $f$ is multiplicative show that
  for all $G$, the good elements are closed under multiplication in the ring $R(G)$.
  Closure under restriction then implies that good elements
  are also closed under the external pairing
  \[ \times \ : \ R(G)\times R(K)\ \to \ R(G\times K)\ . \]
  The additivity relation 
  \[ P^m(x+y) \ = \ {\sum}_{k=0}^m\ \tr_{k,m-k} (P^k(x)\times P^{m-k}(y)) \ ,\]
  the hypothesis that $f$ is additive and the already established 
  closure under transfers and external multiplication then shows that for all $G$, 
  the good elements are closed under addition in the ring $R(G)$.
  The additivity relation and the fact that $P^m(1+(-1))=P^m(0)=0$
  imply that
  \[ P^m(-1)\ = \  - \sum_{k=1}^m\ \tr_{k,m-k} (P^k(1)\times P^{m-k}(-1)) 
  \ = \  - \sum_{k=1}^m\ \tr_{k,m-k} ( 1\times P^{m-k}(-1)) \ .\]
  So induction on $m$ and the previously established closure properties
  show that $f(P^m(-1))=P^m(f(-1))$, so $(e,-1)$ is good. Closure under restriction
  implies that also $(G,-1)$ is good for every compact Lie group $G$.
  Since the good elements are closed under multiplication and contain $-1$,
  they are closed under additive inverses. Since $(G,1)$ is good, so is 
  $(G,1-1)=(G,0)$. This completes the proof that the good elements 
  form a global Green subfunctor of $R$.
  
  The closure property under power operations is then a consequence of
  the transitivity relation and closure under restriction: if $(G,x)$ is good, then
  \begin{align*}
    f(P^k(P^m(x)))\ &= \  f(\Psi_{k,m}^*(P^{k m}(x)))\ = \ \Psi_{k,m}^*(f(P^{k m}(x)))\\ 
    &= \ \Psi_{k,m}^*( P^{k m}(f(x)))\ = \ P^k( P^m(f(x)))\ = \ P^k(f(P^m(x)))
  \end{align*}
  holds in $S(\Sigma_k\wr(\Sigma_m \wr G))$, for all $k\geq 0$.
  This proves that the pair $(\Sigma_m\wr G, P^m(x))$ is also good.
  Altogether this shows that the good elements form a  
  global power subfunctor of $R$.
\end{proof}

The free global power functor $C_K$ associated to a compact Lie group $K$ 
was introduced in Example \ref{eg:free global power functor}.
The next proposition justifies the adjective `free',
and also shows that the underlying global Green functor of $C_K$
represents the functor $\exp(-;K)$ of exponential sequences.
We recall that
\[ 1_{\Sigma_m\wr K}\ \in \ \bA(\Sigma_m\wr K,\Sigma_m\wr K) \ \subset \ 
C_K(\Sigma_m\wr K) \]
denotes the identity operation of the functor $\pi_0^{\Sigma_m\wr K}$.
Two key relations among these elements in the global power functor $C_K$ are
\[ P^m(1_K)\ = \ 1_{\Sigma_m\wr K} \text{\qquad and\qquad}
 \Phi^*_{k,m-k}(1_{\Sigma_m\wr K})\ = \ 1_{\Sigma_k\wr K}\times 1_{\Sigma_{m-k}\wr K} \ .\]
The second set of relations shows that the tuple $(1_{\Sigma_m\wr K})_{m\geq 0}$
is exponential, i.e., an element of $\exp(C_K;K)$.
In fact, part~(i) of the following proposition shows that it is a universal
exponential element.

\begin{prop}\label{prop:what C_K represents}
Let $K$ be a compact Lie group.
\begin{enumerate}[\em (i)]
\item For every global Green functor $R$
the map
\begin{align*}
  \epsilon_K \ : \  \GlGre(C_K,R)\ &\to \ \exp(R;K) \\
f \qquad &\longmapsto \ \exp(f;(1_{\Sigma_m\wr K})_{m\geq 0})\ = \ ( f(1_{\Sigma_m\wr K}))_{m\geq 0}
\end{align*}
is bijective.  
\item For every global power functor $R$ the map
\[  \GlPow(C_K,R)\ \to \ R(K) \ , \quad
f \ \longmapsto \ f(1_K) \]
is bijective.  
\end{enumerate}
\end{prop}
\begin{proof}
(i)
The underlying global functor of $C_K$ is the direct sum of the 
represented global functors $\bA(\Sigma_m\wr K,-)$, for $m\geq 0$.  
So the enriched Yoneda lemma (see Remark \ref{rk:enriched Yoneda})
shows that evaluation at the universal elements is a bijection
\[ \GF(C_K,R)\ \to \ {\prod}_{m\geq 0}\, R(\Sigma_m\wr K)\ , \quad
f \ \longmapsto \ ( f(1_{\Sigma_m\wr K}))_{m\geq 0}\ .\]
A morphism of global functors $f:C_K\to R$
is a morphism of global Green functors if and
only if it is also multiplicative and unital. Unitality corresponds
to the condition $f(1_{\Sigma_0\wr K})=f(\Id_e)=1$ in $R(e)$.
Multiplicativity means that the following square of global functors commutes:
\[ \xymatrix@C=15mm@R=7mm{
C_K\Box C_K\ar[r]^{f\Box f}\ar[d]_{\text{mult}} & R\Box R\ar[d]^{\text{mult}} \\
C_K\ar[r]_f & R} \]
Since the box product of global functors is biadditive and
the box product of two represented functors is represented 
(see Remark \ref{rk:box representable}), the global functor $C_K\Box C_K$
is isomorphic to
\[ {\bigoplus}_{m,n\geq 0}\, \bA(\Sigma_m\wr K,-)\Box\bA(\Sigma_n\wr K,-)\ \iso\
{\bigoplus}_{m,n\geq 0}\, \bA( (\Sigma_m\wr K)\times (\Sigma_n\wr K),-)\ .\]
So commutativity of the above square can be tested by evaluation at the
universal classes, by another application of the enriched Yoneda lemma.
The multiplication of $C_K$ takes the exterior product of
the classes $1_{\Sigma_m\wr K}$ and $1_{\Sigma_n\wr K}$ to the class
$\Phi^*_{m,n}( 1_{\Sigma_{m+n}\wr K})$, so the multiplicativity condition becomes 
the exponential condition
\[ \Phi^*_{m,n}( f(1_{\Sigma_{m+n}\wr K}))\ = \  
f(1_{\Sigma_m\wr K})\times f(1_{\Sigma_n\wr K}) \ .\]

(ii)  Every morphism $f:C_K\to R$ of global power functors is in particular
a morphism of global Green functors. So by part~(i), $f$ is
determined by the exponential sequence $( f(1_{\Sigma_m\wr K}))_{m\geq 0}$. 
The relation $P^m(1_K)=1_{\Sigma_m\wr K}$ holds in the global power functor $C_K$;
so since $f$ also commutes with power operations, it is already determined
by the element $f(1_K)$. This shows that the map in question
is injective.

Now we show that evaluation at $1_K$ is also surjective. We let $y\in R(K)$
be any element. Then the sequence $(P^m(y))_{m\geq 0}$ is exponential,
so part~(i) provides a unique morphism of global Green functors $f:C_K\to R$
satisfying
\[ f(1_{\Sigma_m\wr K}) \ = \ P^m(y)\text{\quad in } R(\Sigma_m\wr K)\]
for all $m\geq 0$. We must show that $f$ is also compatible with power operations. 
By Theorem \ref{thm:Burnside category basis}  the abelian group
\[ C_K(G)\ = \ {\bigoplus}_{m\geq 0}\bA(\Sigma_m\wr K,G) \]
is generated by the classes
\[  \tr_H^G(\alpha^*(P^m(1_K)))  \]
for $m\geq 0$, and where $(H,\alpha)$
runs over pairs consisting of a closed subgroup $H$ of $G$
and a continuous homomorphism $\alpha:H\to\Sigma_m\wr K$.
In particular, $C_K$ is generated, as a global power functor,
by the single element $1_K$ in $C_K(K)$.

Proposition \ref{prop:good elements} shows that the collection of elements $x$
of $C_K$ that satisfy $f(P^m(x))= P^m(f(x))$ for all $m\geq 0$
form a global power subfunctor of $C_K$.
The element $1_K$ is among these because 
$f(P^m(1_K))=f(1_{\Sigma_m\wr K})= P^m(y)=P^m(f(1_K))$.
Since $1_K$ generates $C_K$ as a global power functor, all elements of $C_K$ 
have this property.
So $f$ is a morphism of global power functors.
\end{proof}

Now we can make the rings $\exp(R;G)$ of exponential sequences into
a global functor for varying $G$.
As an auxiliary tool we introduce a functor
\[ \Gamma \ : \ \bA^{\op}\ \to \  \GlPow \  . \]
The functor is given  by $\Gamma(K)=C_K$ on objects;
on morphisms, the freeness property of Proposition \ref{prop:what C_K represents}~(ii)
allows us to define
\[ \Gamma \ : \ \bA(K,G)\ \to \  \GlPow(C_G,C_K) \]
by the requirement
\[ \Gamma(\tau)(1_G)\ = \ \tau \ \in \ \bA(K,G)\ \subset \ C_K(G) \ .\]
The contravariant functoriality is rather formal: for $\psi\in\bA(L,K)$ we have
\[ \Gamma(\tau\circ\psi)(1_G) \ = \ \tau\circ\psi\ = \ 
\Gamma(\psi)(\tau) \ = \ 
 (\Gamma(\psi)\circ\Gamma(\tau))(1_G)\ . \]
So $\Gamma(\tau\circ\psi)=\Gamma(\psi)\circ\Gamma(\tau)$ by freeness,
i.e., Proposition \ref{prop:what C_K represents}~(ii).
We emphasize that while $\bA$ is a preadditive category,
$\Gamma$ is just a plain functor, i.e., not additive in any sense.

\begin{construction}\label{con:Green on exp}
As before we consider a global Green functor $R$.
For compact Lie groups $G$ and $K$ we define a map
\[ \exp(R;-)\ : \ \bA(G,K)\times \exp(R;G)\ \to \ \exp(R;K) \ , \
(\tau,x)\ \longmapsto \ \exp(R;\tau)(x)\ .\]
The definition exploits the representability of $\exp(R;G)$, i.e.,
that the map
\[   \epsilon_G \ : \  \GlGre(C_G,R)\ \to \ \exp(R;G)  \]
sending $f$ to the exponential sequence $( f(1_{\Sigma_m\wr G}))_{m\geq 0}$
is bijective, by Proposition \ref{prop:what C_K represents} (i).
We define $\exp(R;\tau)$
by requiring that for every morphism of global Green functors
$f:C_G\to R$ the following relation holds:
\[ \exp(R;\tau)(\epsilon_G(f))\ = \ \epsilon_K(f\circ\Gamma(\tau)) \ .\]
Since $\epsilon_G$ is bijective, this is a legitimate definition.
It is rather straightforward to see that the assignments
\[ G\ \longmapsto \ \exp(R;G)\text{\qquad and\qquad}
\tau\ \longmapsto \ \exp(R;\tau) \]
define a functor from the global Burnside category to the category of sets.
Indeed, for $\tau\in\bA(K,G)$, $\psi\in\bA(L,K)$
and a morphism of global Green functors $f:C_G\to R$, we have
\begin{align*}
  \exp(R;\tau\circ\psi)&(\epsilon_G(f))\ = \ 
\epsilon_L(f\circ\Gamma(\tau\circ\psi)) \ = \ 
\epsilon_L(f\circ\Gamma(\psi)\circ\Gamma(\tau)) \\ 
&= \ \exp(R;\tau)(\epsilon_K(f\circ\Gamma(\psi))) \ 
= \ \exp(R;\tau)(\exp(R;\psi)(\epsilon_G(f))) \ .
\end{align*}
Since every exponential sequence in $\exp(R;G)$ is of the form
$\epsilon_G(f)$, this proves that $\exp(R;\tau\circ\psi)=\exp(R;\tau)\circ\exp(R;\psi)$.
It is not completely obvious, though, that this construction is additive
in both variables, but we will show that in the next proposition.

The construction of exponential sequences is in fact a functor in two variables:
for every operation $\tau\in\bA(G,K)$ 
and all morphisms of global Green functors $f:C_G\to R$ we have
\begin{align*}
 \left( ( \exp(\varphi;K)\circ \exp(R;\tau))(\epsilon_G(f)) \right)_m\ &=\ 
\left( ( \exp(\varphi;K) (\epsilon_K(f\circ\Gamma(\tau)))) \right)_m\\ 
&=\ \varphi\left( (\epsilon_K(f\circ\Gamma(\tau)))_m \right)\\
&=\ (\varphi\circ f\circ \Gamma(\tau))(1_{\Sigma_m\wr K})\\ 
&= \ \left( \exp(S;\tau)(\epsilon_G(\varphi \circ f)) \right)_m\\
&= \ \left( (\exp(S;\tau)\circ\exp(\varphi;G))(\epsilon_G(f)) \right)_m\ .
\end{align*}
Every exponential sequence in $\exp(R;G)$ is of the form
$\epsilon_G(f)$, so
\[ \exp(\varphi;K)\circ \exp(R;\tau)\ =\ \exp(S;\tau)\circ\exp(\varphi;G) \ : \
\exp(R;G)\ \to \ \exp(S;K)\ .\]
\end{construction}

The next proposition makes the abstract definition of the map
$\exp(R;\tau)$ more explicit by giving a concrete formula when
$\tau$ is a restriction or a transfer.

\begin{prop}\label{prop:exp(R) is global Green}
  Let $R$ be a global Green functor.
  \begin{enumerate}[\em (i)]
  \item For every continuous group homomorphism $\alpha:K\to G$ between
    compact Lie groups the map 
    \[ \exp(R;\alpha^*)\ :\ \exp(R;G)\ \to\ \exp(R;K) \]
    is a ring homomorphism and satisfies
    \[ \left( \exp(R;\alpha^*)(x)\right)_m \ = \  (\Sigma_m\wr\alpha)^*(x_m)\]
    in $R(\Sigma_m\wr K)$ for all $x\in \exp(R;G)$.
  \item For every closed subgroup $H$ of a compact Lie group $G$ 
    the map
    \[  \exp(R;\tr_H^G)\ :\ \exp(R;H)\ \to\ \exp(R;G)\]
    is additive, satisfies reciprocity with respect to restriction from $G$ to $H$
    and is given by
    \[ \left( \exp(R;\tr_H^G)(x)\right)_m \ = \  \tr_{\Sigma_m\wr H}^{\Sigma_m\wr G}(x_m)\]
    for all $x\in\exp(R;H)$.
  \item
    For all compact Lie groups $G$ and $K$, the map 
    \[ \exp(R;-)\ : \ \bA(G,K)\times \exp(R;G)\ \to \ \exp(R;K) \]
    is biadditive, i.e., $\exp(R;-)$ becomes a global functor in the Lie group variable.
  \item
    As the Lie group $G$ varies, the ring structures and the functoriality
    in the global Burnside category make $\exp(R)$ into a global Green functor.
    \end{enumerate}
\end{prop}
\begin{proof}
(i) We let $f:C_G\to R$ be any morphism of global Green functors. Then
\begin{align*}
 \left( \exp(R;\alpha^*)(\epsilon_G(f))\right)_m \ &= \  
\left( \epsilon_K(f\circ \Gamma(\alpha^*))\right)_m\
= \ (f\circ \Gamma(\alpha^*))(1_{\Sigma_m\wr K}) \\
&= \ f ( \Gamma(\alpha^*)(P^m(1_K))) \
= \ f ( P^m( \Gamma(\alpha^*)(1_K))) \\
&= \ f ( P^m( \alpha^*(1_G))) \
= \ f ( (\Sigma_m\wr \alpha)^*(P^m(1_G))) \\
&= \ f ( (\Sigma_m\wr \alpha)^*(1_{\Sigma_m\wr G})) \
= \  (\Sigma_m\wr \alpha)^*(f(1_{\Sigma_m\wr G})) \\
&= \ (\Sigma_m\wr\alpha)^*( \epsilon_G(f)_m)\ .
\end{align*}
Since every exponential sequence in $\exp(R;G)$ is of the form
$\epsilon_G(f)$, this proves the formula for $\exp(R;\alpha^*)$.
Since the multiplication in $\exp(R;G)$ is coordinatewise
and the original restriction maps for $R$ are ring homomorphisms,
the formula for $\exp(R;\alpha^*)$ shows that $\exp(R;\alpha^*)$
is multiplicative and preserves the multiplicative unit.
The additivity of $\exp(R;\alpha^*)$ uses the relation
\begin{equation}  \label{eq:alpha_after_tr}
 (\Sigma_m\wr \alpha)^*\circ \tr_{k,m-k} \ = \ \tr_{k,m-k}\circ ((\Sigma_k\wr\alpha)\times(\Sigma_{m-k}\wr\alpha))^*   
\end{equation}
as maps from $R((\Sigma_k\wr G)\times(\Sigma_{m-k}\wr G))$
to $R(\Sigma_m\wr K)$.
To prove \eqref{eq:alpha_after_tr} we distinguish two cases.
If $\alpha$ is surjective, then so is $\Sigma_m\wr\alpha$, and 
\[ (\Sigma_m\wr\alpha)^{-1}( (\Sigma_k\wr G)\times(\Sigma_{m-k}\wr G)) \ = \ 
 (\Sigma_k\wr K)\times(\Sigma_{m-k}\wr K)\ . \]
So for epimorphisms, the relation \eqref{eq:alpha_after_tr} 
is a special case of compatibility of transfer with inflation.
If $H$ is a closed subgroup of $G$, then $\Sigma_m\wr G$ consists of a single
double coset for the left  $(\Sigma_m\wr H)$-action
and right $( (\Sigma_k\wr G)\times(\Sigma_{m-k}\wr G))$-action, and
\[(\Sigma_m\wr H)\ \cap \  ( (\Sigma_k\wr G)\times(\Sigma_{m-k}\wr G)) 
  \ = \ 
 (\Sigma_k\wr H)\times(\Sigma_{m-k}\wr H)  \ .\]
So the double coset formula specializes to 
\[  \res^{\Sigma_m\wr G}_{\Sigma_m\wr H} \circ \tr_{k,m-k} \ = \ 
\tr_{k,m-k} \circ  \res^{ (\Sigma_k\wr G)\times(\Sigma_{m-k}\wr G)}_{(\Sigma_k\wr H)\times(\Sigma_{m-k}\wr H)} 
 \ ,\]
which is precisely the relation \eqref{eq:alpha_after_tr} 
for the inclusion $H\to G$. 
Every homomorphism factors as an epimorphism followed by a subgroup
inclusion, so relation \eqref{eq:alpha_after_tr} follows in general.
We can then conclude that $\exp(R;\alpha^*)$ is additive:
\begin{align*}
 (\exp(R;\alpha^*)(x\oplus y) )_m \ &= \ (\Sigma_m\wr \alpha)^*( (x\oplus y)_m) \
 = \ \sum_{k=0}^m \ (\Sigma_m\wr\alpha)^*(\tr_{k,m-k} (x_k\times y_{m-k})     )\\
_\eqref{eq:alpha_after_tr}\ &= \ \sum_{k=0}^m \ 
\tr_{k,m-k} \left( (  (\Sigma_k\wr \alpha)\times(\Sigma_{m-k}\wr\alpha))^*(x_k\times y_{m-k})\right)   \\
&= \ \sum_{k=0}^m \ 
\tr_{k,m-k} \left( (\Sigma_k\wr \alpha)^*(x_k)\times(\Sigma_{m-k}\wr\alpha)^*(y_{m-k})\right)\\
&= \ \sum_{k=0}^m \ 
\tr_{k,m-k} \left( (\exp(R;\alpha^*)(x))_k \times (\exp(R;\alpha^*)(y))_{m-k}\right)\\
&= \  \left( \exp(R;\alpha^*)(x)\oplus \exp(R;\alpha^*)(y) \right)_m \ .
\end{align*}

Part~(ii) is similar to part~(i).
We let $f:C_H\to R$ be any morphism of global Green functors. Then
\begin{align*}
 \left( \exp(R;\tr_H^G)(\epsilon_H(f))\right)_m \ &= \  
\left( \epsilon_G(f\circ \Gamma(\tr_H^G))\right)_m\
= \ (f\circ \Gamma(\tr_H^G))(1_{\Sigma_m\wr G}) \\
&= \ f ( \Gamma(\tr_H^G)(P^m(1_G))) \
= \ f ( P^m( \Gamma(\tr_H^G)(1_G))) \\
&= \ f ( P^m( \tr_H^G(1_H))) \
= \ f ( \tr_{\Sigma_m\wr H}^{\Sigma_m\wr G}(P^m(1_H))) \\
&= \ f ( \tr_{\Sigma_m\wr H}^{\Sigma_m\wr G}(1_{\Sigma_m\wr H})) \
= \  \tr_{\Sigma_m\wr H}^{\Sigma_m\wr G}(f(1_{\Sigma_m\wr H})) \
= \  \tr_{\Sigma_m\wr H}^{\Sigma_m\wr G}( \epsilon_H(f)_m)\ .
\end{align*}
Since every exponential sequence is of the form $\epsilon_H(f)$, 
this proves the formula for $\exp(R;\tr_H^G)$.
To see that the map $\exp(R;\tr_H^G)$ is additive we observe that
\[  \tr_{\Sigma_m\wr H}^{\Sigma_m\wr G}\circ \tr_{k,m-k} \ = \
\tr_{k,m-k} \circ \tr_{(\Sigma_k\wr H)\times(\Sigma_{m-k}\wr H)}^{(\Sigma_k\wr G)\times(\Sigma_{m-k}\wr G)}\ ,
 \]
by transitivity of transfers. Thus 
\begin{align*}
 ( \exp(R;\tr_H^G)(x \oplus y) )_m \ &= \ 
\sum_{k=0}^m \ \tr_{\Sigma_m\wr H}^{\Sigma_m\wr G}(\tr_{k,m-k} (x_k\times y_{m-k})   )\\
 &= \ \sum_{k=0}^m \ 
\tr_{k,m-k} \left( \tr_{(\Sigma_k\wr H)\times(\Sigma_{m-k}\wr H)}^{(\Sigma_k\wr G)\times(\Sigma_{m-k}\wr G)}(x_k\times y_{m-k})\right)   \\
&= \ \sum_{k=0}^m \ 
\tr_{k,m-k} \left( \tr_{\Sigma_k\wr H}^{\Sigma_k\wr G}(x_k)\times \tr_{\Sigma_{m-k}\wr H}^{\Sigma_{m-k}\wr G}(y_{m-k})\right)\\
&= \ ( \exp(R;\tr_H^G)(x) \oplus \exp(R;\tr_H^G)(y) )_m \ .
\end{align*}
The reciprocity for restriction and transfer is now a direct
consequence of the reciprocity for the global Green functor $R$:
\begin{align*}
 ( \exp(R;\tr_H^G)(x)\cdot y )_m \ &= \ 
\tr_{\Sigma_m\wr H}^{\Sigma_m\wr G}(x_m)\cdot y_m \\ 
&= \ 
\tr_{\Sigma_m\wr H}^{\Sigma_m\wr G}(x_m\cdot \res^{\Sigma_m\wr G}_{\Sigma_m\wr H}( y_m)) \
= \ ( \exp(R;\tr_H^G)(x \cdot \res_H^G(y) ))_m \ . 
\end{align*}

(iii)
We start with additivity of the map $\exp(R;-)$ in the variable $\tau$.
We let $\tau'\in\bA(G,K)$ be another natural transformation.
Then the morphism of global power functors $\Gamma(\tau+\tau'):C_K \to C_G$
satisfies
\begin{align*}
 \Gamma(\tau+\tau')(1_{\Sigma_m\wr K}) \ &= \   
 \Gamma(\tau+\tau')(P^m(1_K)) \ = \   P^m( \Gamma(\tau+\tau')(1_K)) \\ 
&= \   P^m( \tau+\tau') \ =\
\sum_{k=0}^m \tr_{k,m-k}(P^k(\tau)\times P^{m-k}(\tau')) \\ 
&= \  
\sum_{k=0}^m \tr_{k,m-k}( P^k( \Gamma(\tau)(1_K))\times P^{m-k}(\Gamma(\tau')(1_K))) \\ 
&= \  
\sum_{k=0}^m \tr_{k,m-k}(  \Gamma(\tau)(1_{\Sigma_k\wr K})\times \Gamma(\tau')(1_{\Sigma_{m-k}\wr K})) \ .
\end{align*}
For a morphism of global Green functors $f:C_G\to R$ this implies
\begin{align*}
(f\circ\Gamma(\tau+\tau'))(1_{\Sigma_m\wr K}) \ &= \ 
\sum_{k=0}^m \tr_{k,m-k}(  f(\Gamma(\tau)(1_{\Sigma_k\wr K}))\times f(\Gamma(\tau')(1_{\Sigma_{m-k}\wr K})))\ .
\end{align*}
This establishes the $m$-th component of the relation
\[ \epsilon_K(f\circ \Gamma(\tau+\tau'))\ = \ 
\epsilon_K(f\circ\Gamma(\tau)) \oplus \epsilon_K(f\circ\Gamma(\tau'))\ ,\]
and this in turn proves
\[ \exp(R;\tau+\tau')(\epsilon_G(f))\ = \ 
\exp(R;\tau)(\epsilon_G(f))\ \oplus \ \exp(R;\tau')(\epsilon_G(f))
\ .\]
Now we come to additivity of the map in the exponential sequence,
i.e., that for fixed $\tau\in\bA(G,K)$ the map
\[ \exp(R;\tau)\ : \ \exp(R;G)\ \to \exp(R;K) \]
is additive with respect to $\oplus$.
By Theorem \ref{thm:Burnside category basis} the abelian group $\bA(G,K)$ 
is generated by operations of the form $\tr_L^K\circ\alpha^*$
where $L$ is a closed subgroup of $K$ and $\alpha:L\to G$
a continuous homomorphism.
Since we have already established additivity in $\tau$, we can assume
that $\tau=\tr_L^K\circ\alpha^*$ is one of the generating operations.
Since
\[ \exp(R;\tr_L^K\circ\alpha^*)\ = \ \exp(R;\tr_L^K)\circ \exp(R;\alpha^*) \ ,\]
this in turn follows from additivity for restrictions and transfers,
which we showed in part~(i) respectively~(ii).

Part~(iv) is just the summary of the properties in parts~(i), (ii) and~(iii).
\end{proof}

The previous proposition establishes the functor $\exp$
of exponential sequences as an endofunctor of the category of global Green functors.
When iterating the construction, 
we encounter iterated wreath products, and we will save a 
substantial number of parentheses by using the short hand notation
\[ \Sigma_k\wr\Sigma_m\wr G \ = \ \Sigma_k\wr(\Sigma_m\wr G)\ . \]
Now we make this endofunctor into a comonad.
A natural transformation of global Green functors
\[ \eta_R\ : \ \exp(R) \ \to \ R  \]
is given by $\eta(x)=x_1$, using the identification $G\iso \Sigma_1\wr G$
via $g\mapsto (1;g)$.
A natural transformation 
\[  \kappa_R\ : \ \exp(R) \ \to \ \exp(\exp(R))   \]
is given at a compact Lie group $G$ by 
\[ ( \kappa(x)_m )_k \ = \ \Psi_{k,m}^*(x_{k m}) \ \in\  R(\Sigma_k\wr\Sigma_m\wr G)\ ;\]
here the restriction is along the monomorphism \eqref{eq:wreath_iterate spectra} 
\begin{align*} 
\Psi_{k,m}\ : \ \Sigma_k\wr \Sigma_m\wr G \qquad
&\to\qquad \Sigma_{k m}\wr G \\
 (\sigma;\, (\tau_1;\, h^1),\dots,(\tau_k;\, h^k)) \ &\longmapsto \
(\sigma\natural(\tau_1,\dots,\tau_k);\, h^1+\dots+h^k)\ .
\end{align*}

The analog of the following for finite groups is Satz~2.17 in \cite{singer}.

\begin{theorem}\label{thm:exp is comonad}
  \begin{enumerate}[\em (i)]
    Let $R$ be a global Green functor.
  \item 
     For every compact Lie group $G$ 
    and for every exponential sequence $x\in \exp(R;G)$,
    the sequence $\kappa(x)$ is an element of $\exp(\exp(R);G)$.
  \item
    As the group varies, the maps $\kappa$ form a morphism of global Green functors
    $\kappa_R:\exp(R)\to\exp(\exp(R))$, natural in $R$.
  \item
    The natural transformations 
    \[ \eta\ : \ \exp \ \to\ \Id \text{\qquad and\qquad}
    \kappa\ :\ \exp\ \to\ \exp\circ\exp \]
    make the functor $\exp$ into a comonad on the category of global Green functors.
  \end{enumerate}
\end{theorem}

\begin{proof}
(i) Because the square of group homomorphisms
\[ \xymatrix@C=13mm@R=7mm{ 
(\Sigma_j\wr\Sigma_m\wr G) \times (\Sigma_{k-j}\wr\Sigma_m\wr G) 
\ar[r]^-{\Phi_{j,k-j}}\ar[d]_{\Psi_{j,m}\times\Psi_{k-j,m}} &
\Sigma_k\wr\Sigma_m\wr G \ar[d]^{\Psi_{k,m}} \\
(\Sigma_{j m}\wr G) \times (\Sigma_{(k-j)m}\wr G) \ar[r]_-{\Phi_{j m,(k-j)m}} & 
\Sigma_{k m}\wr G } \]
commutes, the elements $\Psi_{k,m}^*(x_{k m})$ satisfy
\begin{align*}
  \Phi_{j,k-j}^*( \Psi_{k,m}^*(x_{k m}) )\ 
&= \  (\Psi_{j,m}\times\Psi_{k-j,m})^*(\Phi_{j m,(k-j)m}^*(x_{k m}))\\
&= \  (\Psi_{j,m}\times\Psi_{k-j,m})^*(x_{j m}\times x_{(k-j)m})\\ 
&= \  \Psi_{j,m}^*(x_{j m})\times\Psi_{k-j,m}^*(x_{(k-j)m})\ .
\end{align*}
This shows that for fixed $m\geq 0$, the sequence
$\kappa(x)_m=(\Psi_{k,m}^*(x_{k m}))_{k\geq 0}$ is exponential, i.e., an element
of the ring $\exp(R;\Sigma_m\wr G)$.

The square of group homomorphisms
\begin{align}\begin{aligned}\label{eq:iterated_wreat_diagram}
 \xymatrix@C=13mm@R=7mm{ 
\Sigma_k\wr( (\Sigma_i\wr G) \times (\Sigma_{m-i}\wr G) ) 
\ar[r]^-{\Sigma_k\wr \Phi_{i,m-i}} \ar[d]_{\Delta}  &
\Sigma_k\wr\Sigma_m\wr G \ar[dd]^{\Psi_{k,m}} \\
(\Sigma_k\wr \Sigma_i\wr G) \times (\Sigma_k\wr \Sigma_{m-i}\wr G ) 
\ar[d]_{\Psi_{k,i}\times \Psi_{k,m-i}}  &\\
(\Sigma_{k i}\wr G) \times (\Sigma_{k(m-i)}\wr G) \ar[r]_-{\Phi_{k i,k(m-i)}} & 
\Sigma_{k m}\wr G }     
\end{aligned}\end{align}
does {\em not} commute; we invite the reader to check
the case $k=m=2$, $i=1$ and $G=e$, where the phenomenon is already visible
in the fact that the square
\[ \xymatrix@C=13mm@R=7mm{ 
\Sigma_2 \ar[r]^-{\Sigma_2\wr \Phi_{1,1}} \ar[d]_{\Delta}  &
\Sigma_2\wr \Sigma_2 \ar[d]^{\natural} \\
\Sigma_2 \times \Sigma_2\ar[r]_-{+} &  \Sigma_4} \]
does {\em not} commute.
However, the square \eqref{eq:iterated_wreat_diagram} 
does commute up to conjugation by an element of $\Sigma_{k m}\wr G$.
Since inner automorphisms are invisible through the eyes of a Rep-functor,
we conclude that the relation
\begin{align*}
\left(  \Phi_{i,m-i}^*( \kappa(x)_m ) \right)_k \ &= \ 
(\Sigma_k\wr\Phi_{i,m-i})^*( \Psi_{k,m}^*(x_{k m}) ) \\ 
&= \  \Delta^*( ( \Psi_{k,i}\times \Psi_{k,m-i})^*(\Phi_{k i, k(m-i)}^*( x_{k m})))\\
&= \  \Delta^*( \Psi_{k,i}^*(x_{k i}) \times \Psi_{k,m-i}^*(x_{k(m-i)}))\\
&= \ \Delta^*( (\kappa(x)_i )_k \times (\kappa(x)_{m-i})_k)\
= \ \left(\kappa(x)_i \times \kappa(x)_{m-i}\right)_k
\end{align*}
holds in $R(\Sigma_k\wr ((\Sigma_i\wr G)\times (\Sigma_{m-i}\wr G)))$.
For varying $k\geq 0$, this shows that
\[   \Phi_{i,m-i}^*( \kappa(x)_m ) \ = \ \kappa(x)_i\times \kappa(x)_{m-i}
\text{\qquad in\quad $\exp(R;(\Sigma_i\wr G)\times (\Sigma_{m-i}\wr G))$.} \]
In other words, the sequence $\kappa(x)=(\kappa(x)_m)_{m\geq 0}$ is 
itself exponential.

(ii) The relations $\kappa(1)=1$ and
$\kappa(x\cdot y)=\kappa(x)\cdot \kappa(y)$
are straightforward from the definitions, using that multiplication
is defined coordinatewise and that the restriction map $\Psi_{k,m}^*$
is multiplicative and unital.
The verification of the additivity of $\kappa_R(G):\exp(R;G)\to\exp(\exp(R);G)$
is more involved,
and forces us to confront the double coset formula for the subgroups
\[ \Sigma_k\wr'\Sigma_m \ = \  \Psi_{k,m}(\Sigma_k\wr\Sigma_m\wr G) 
\text{\quad and\quad} 
\Sigma_i\times'\Sigma_{k m-i}\ = \ \Phi_{i,k m-i}( (\Sigma_i\wr G)\times(\Sigma_{k m-i}\wr G))
\]
of the wreath product $\Sigma_{k m}\wr G$. To a large extent, the group $G$
acts like a dummy, which is why we omit it from the notation for
$\Sigma_k\wr'\Sigma_m$ and $\Sigma_i\times'\Sigma_{k m-i}$.
We specify a bijection between the set of double cosets 
\[  \Sigma_k\wr'\Sigma_m \bs \Sigma_{k m}\wr G / \Sigma_i\times'\Sigma_{k m-i} \]
and the set of tuples $(a_0,\dots,a_m)$ of natural numbers that satisfy 
\begin{equation}  \label{eq:two_conditions}
 a_0+\dots+a_m \ = \ k \text{\qquad and\qquad}
{\sum}_{j=0}^m\, j\cdot a_j \ = \ i \ .   
\end{equation}
For $1\leq b\leq k$ we let $j_b\in\{1,\dots, m\}$ be the unique number such that
\[ a_0+\dots+a_{j_b-1} \ < \ b\  \leq\  a_0+\dots+a_{j_b-1}+a_{j_b} \ .\]
We define
\[ A(a_0,\dots,a_m)\ = \ {\bigcup}_{b=1}^k\, \{ (b-1)m+1,\dots, (b-1)m+ j_b \} \ 
\subseteq \ \{1,\dots, k m\}\ ; \]
because
\[ {\sum}_{b=1}^k \, j_b \ = \ 
 {\sum}_{j=1}^m \, a_j\cdot j\ = \  i\ , \]
the set $A(a_0,\dots,a_m)$ has exactly $i$ elements.
We let $\bar\sigma=\bar\sigma(a_0,\dots,a_m)\in\Sigma_{k m}$
be any permutation such that
\[ \{\bar\sigma(1),\dots,\bar\sigma(i)\} \ = \ A(a_0,\dots,a_m)\ . \]
Then the elements
\[ \sigma(a_0,\dots,a_m)\ = \ (\bar\sigma(a_0,\dots,a_m);\,1,\dots,1) \ \in \ \Sigma_{k m}\wr G  \]
are a complete set of double coset representatives as $(a_0,\dots,a_m)$
ranges over those tuples that satisfy \eqref{eq:two_conditions}.
Moreover,
\begin{align*}
  (\Sigma_k\wr'\Sigma_m)^{\sigma(a_0,\dots,a_m)}&\cap (\Sigma_i\times'\Sigma_{k m-i}) 
\ = \ 
 \left( {{\prod}'}_{j=0}^m \Sigma_{a_j}\wr'\Sigma_j \right)
\times'
\left( {{\prod}'}_{j=0}^m\Sigma_{a_j}\wr'\Sigma_{m-j}\right)   
\end{align*}
and
\begin{align*}
 (\Sigma_k\wr'\Sigma_m)&\cap {^{\sigma(a_0,\dots,a_m)}}(\Sigma_i\times'\Sigma_{k m-i}) 
\ = \ 
{{\prod}'}_{j=0}^m (\Sigma_{a_j}\wr'\Sigma_j)\times' (\Sigma_{a_j}\wr'\Sigma_{m-j})\ .
\end{align*}
For $x_i\in R(\Sigma_i\wr G)$ and $y_{k m-i}\in R(\Sigma_{k m-i}\wr G)$,
the double coset formula thus becomes the relation
\begin{align*}
 \Psi_{k,m}^*&( \tr_{i,k m-i}( x_i \times y_{k m-i} ))\\ =   
&\sum  \tr^{\Sigma_k\wr\Sigma_m\wr G}_{\prod' (\Sigma_{a_j}\wr'\Sigma_j)\times' (\Sigma_{a_j}\wr'\Sigma_{m-j})} \left( \sigma(a_0,\dots,a_m)_\star\left( \res^{(\Sigma_i\wr G)\times(\Sigma_{k m-i}\wr G)}_{ \left( \prod' \Sigma_{a_j}\wr'\Sigma_j \right)
\times'
\left( \prod'\Sigma_{a_j}\wr'\Sigma_{m-j}\right)}( x_i \times y_{k m-i} ) \right)\right)
\end{align*}
where the sum ranges over all tuples $(a_0,\dots,a_m)$ 
that satisfy \eqref{eq:two_conditions}.
If $x_i$ and $y_{k m-i}$ are the respective components of two exponential sequences 
$x, y\in \exp(R;G)$, then
\begin{align*}
\sigma(a_0,\dots,a_m)_\star&\left(\res^{(\Sigma_i\wr G)\times(\Sigma_{k m-i}\wr G)}_{ \left( {\prod'} \Sigma_{a_j}\wr'\Sigma_j \right)
\times'
\left( \prod'\Sigma_{a_j}\wr'\Sigma_{m-j}\right)}( x_i \times y_{k m-i} )\right) \\
&= \  
\sigma(a_0,\dots,a_m)_\star\left(\res^{\Sigma_i\wr G}_{ \prod' \Sigma_{a_j}\wr'\Sigma_j }( x_i) \times
\res^{\Sigma_{k m-i}\wr G}_{ \prod'\Sigma_{a_j}\wr'\Sigma_{m-j}}( y_{k m-i} )\right)\\
&= \ \sigma(a_0,\dots,a_m)_\star\left( {\prod}_{j=0}^m \Psi_{a_j,j}^*(x_{ a_j j}) \times
  {\prod}_{j=0}^m \Psi_{a_j,m-j}^*(y_{ a_j(m-j)})\right)\\
&= \  {\prod}_{j=0}^m \left( \Psi_{a_j,j}^*(x_{a_j j})\times
\Psi_{a_j,m-j}^*(y_{a_j(m-j)})\right)\ .
\end{align*}
Given $x,y\in\exp(R;G)$ we have
\[  ( \kappa(x) \oplus \kappa(y) )_m  \ = \ 
{\oplus}_{j=0}^m\ \tr_{j,m-j}( \kappa(x)_j \times \kappa(y)_{m-j} )\ ,\]
where the sum on the right is taken in the group $\exp(R;\Sigma_m\wr G)$
under $\oplus$. Expanding this further we arrive at the expression
\begin{align*}
( ( \kappa(x) &\oplus \kappa(y) )_m)_k  \ = \ 
\sum_{a_0+\dots+a_m=k} \tr_{a_0,\dots,a_m}\left(
{\prod}_{j=0}^m \tr_{j,m-j}( \kappa(x)_j \times \kappa(y)_{m-j} )_{a_j} \right)  \\
&= \ 
\sum_{a_0+\dots+a_m=k} \tr_{a_0,\dots,a_m}\left(
{\prod}_{j=0}^m \tr_{\Sigma_{a_j}\wr\Sigma_{j,m-j}}^{\Sigma_{a_j}\wr\Sigma_m}( \Psi_{a_j,j}^*( x_{a_j j}) \times \Psi^*_{a_j,m-j}(y_{a_j(m-j)} ))\right) \\
&= \  
\sum_{a_0+\dots+a_m=k}
  \tr^{\Sigma_k\wr\Sigma_m\wr G}_{\prod' (\Sigma_{a_j}\wr'\Sigma_j)\times' (\Sigma_{a_j}\wr'\Sigma_{m-j})} 
\left(  {\prod}_{j=0}^m \Psi_{a_j,j}^*(x_{a_j j})\times
\Psi_{a_j,m-j}^*(y_{a_j(m-j)})\right)\\
&= \ \sum_{a_0+\dots+a_m=k}
  \tr^{\Sigma_k\wr\Sigma_m\wr G}_{\prod' (\Sigma_{a_j}\wr'\Sigma_j)\times' (\Sigma_{a_j}\wr'\Sigma_{m-j})} 
\left(  {\prod}_{j=0}^m \Psi_{a_j,j}^*(x_{a_j j})\times
\Psi_{a_j,m-j}^*(y_{a_j (m-j)})\right)\\
&= \ 
\sum_{i=0}^{k m} \Psi_{k,m}^*( \tr_{i,k m-i}( x_i \times y_{k m-i} ))\ =\
  ( \kappa(x \oplus y)_m )_k  
\end{align*}
in the group $R(\Sigma_k\wr\Sigma_m\wr G)$.
Here $\tr_{a_0,\dots,a_m}$ is shorthand notation for the transfer along the monomorphism
\[  \Phi_{a_0,\dots,a_m}\ : \ (\Sigma_{a_0}\wr\Sigma_m\wr G)\times\dots\times
 (\Sigma_{a_j}\wr\Sigma_m\wr G)\times\dots\times (\Sigma_{a_m}\wr\Sigma_m\wr G)
\ \to \  \Sigma_k\wr\Sigma_m\wr G  \ ,\]
the analog of the monomorphism \eqref{eq:wreath_sum spectra}
with multiple inputs (and for the group $\Sigma_m\wr G$ instead of $G$).   
So the map $\kappa_R(G):\exp(R;G)\to \exp(\exp(R);G)$ is additive,
and hence a ring homomorphism.

It remains to check that the maps $\kappa_R$ commute with restrictions
and transfers.
For every continuous homomorphism $\alpha:K\to G$ the relations
\begin{align*}
( \kappa(\alpha^*(x))_m )_k \ &= \ \Psi_{k,m}^*((\Sigma_{k m}\wr \alpha)^* (x_{k m})) \
  = \   (\Sigma_k\wr\Sigma_m\wr \alpha)^* (\Psi_{k,m}^*(x_{k m})) \\ 
 &= \   (\Sigma_k\wr\Sigma_m\wr \alpha)^* ((\kappa(x)_m)_k) \ = \ 
 \left( ( \alpha^*(\kappa(x)))_m \right)_k 
\end{align*}
hold in $R(\Sigma_k\wr\Sigma_m\wr K)$. 
So $\kappa_R(K)\circ\exp(R;\alpha^*)=\exp(\exp(R);\alpha^*)\circ\kappa_R(G)$.

The compatibility of $\kappa$ with transfers needs another application 
of a double coset formula.
Indeed, for every closed subgroup $H$ of $G$, the group 
$\Sigma_{k m}\wr G$ consists of a single
double coset for the left $ (\Sigma_k\wr \Sigma_m\wr G)$-action
and right $(\Sigma_{k m}\wr H)$-action, and
\[  (\Sigma_k\wr \Sigma_m\wr G) \ \cap \ (\Sigma_{k  m}\wr H)  \ = \ 
 \Sigma_k\wr \Sigma_m\wr H  \ .\]
So for every $x\in\exp(R;H)$, the relations
\begin{align*}
( \kappa(\tr_H^G(x))_m )_k\ = \ 
\Psi_{k,m}^* \left(\tr_{\Sigma_{k m}\wr H}^{\Sigma_{k m}\wr G} (x_{k m}) \right) \ 
&= \  \tr_{\Sigma_k\wr\Sigma_m\wr H}^{\Sigma_k\wr\Sigma_m\wr G} (\Psi_{k,m}^*(x_{k m})) \\ 
&= \  \tr_{\Sigma_k\wr\Sigma_m\wr H}^{\Sigma_k\wr\Sigma_m\wr G} ((\kappa(x)_m)_k) \ = \
 \left( ( \tr_H^G(\kappa(x)))_m \right)_k 
\end{align*}
holds in $R(\Sigma_k\wr\Sigma_m\wr G)$. Hence $\kappa\circ\tr_H^G=\tr_H^G\circ\kappa$.
Altogether this shows that the maps $\kappa_R(G)$ are
ring homomorphisms and compatible with restrictions and transfers,
so they form a morphism of global Green functors.

(iii)
We have to show that the transformation $\kappa$ is coassociative, 
and counital with respect to $\eta$, and 
these are all straightforward from the definitions.
The counitality relations
\[ \eta_{\exp(M)}\circ \kappa_M \ = \ \Id_{\exp(M)}\ = \ \exp(\eta_M)\circ \kappa_M\]
come down to the facts that the homomorphism $\Psi_{k,1}$
is the result of applying $\Sigma_k\wr -$
to the preferred isomorphism $\Sigma_1\wr G\iso G$,
and that the homomorphism $\Psi_{1,n}$ is 
the preferred isomorphism $\Sigma_1\wr \Sigma_n\wr G\iso\Sigma_n\wr G$.
The coassociativity relation
\[ \exp(\kappa_M)\circ\kappa_M \ = \ \kappa_{\exp(M)}\circ\kappa_M \]
ultimately boils down to the observation that the following square of monomorphisms
commutes:
\[
\begin{gathered}[b]
\xymatrix@C=15mm@R=7mm{ 
\Sigma_k\wr\Sigma_m\wr\Sigma_n\wr G\ar[r]^-{\Sigma_k\wr \Psi_{m,n}}
\ar[d]_{\Psi_{k,m}} &
\Sigma_k\wr\Sigma_{m n}\wr G\ar[d]^{\Psi_{k, m n}} \\
\Sigma_{k m}\wr\Sigma_n\wr G\ar[r]_-{\Psi_{k m, n}} &
\Sigma_{k m n}\wr G} \\[-\dp\strutbox]
\end{gathered}
\qedhere
\]
\end{proof}

Now we can finally get to the main result of this section,
identifying global power functors with coalgebras over the comonad
of exponential sequences.
We suppose that $R$ is a global Green functor and $P:R\to\exp(R)$
a morphism of global Green functors. 
For every compact Lie group $G$, a sequence of operations $P^m:R(G)\to R(\Sigma_m\wr G)$
is then defined by
\[ P^m(x) \ = \ ( P(x) )_m \ , \]
i.e., $P^m(x)$ is the $m$-th component of the exponential sequence $P(x)$.

\begin{theorem}[Comonadic description of global power functors]\label{thm:power comonadic} \ 
  \begin{enumerate}[\em (i)]
  \item 
    Let $R$ be a global Green functor and
    $P:R\to\exp(R)$ a morphism of global Green functors that makes $R$ into
    a coalgebra over the comonad $(\exp,\eta,\kappa)$.
    Then the operations  $P^m:R(G)\to R(\Sigma_m\wr G)$
    make $R$ into a global power functor.
  \item
    The functor
    \begin{align*}
    (\exp\text{\em-coalgebras}) \ &\to \ \GlPow \ ,\quad
    (R, P) \ \longmapsto \ (R, \{P^m\}_{m\geq 0})
    \end{align*}
    is an isomorphism of categories.
  \end{enumerate}
\end{theorem}
\begin{proof}
(i) The fact that $P:R\to\exp(R)$ takes values in exponential sequences
is equivalent to the restriction condition of the power operations.
The fact that $P:R\to\exp(R)$ is a morphism of global Green functors
encodes simultaneously the unit, contravariant naturality, transfer,
multiplicativity and additivity relations of a global power functor.
The identity relation $P^1=\Id$ is equivalent to the counit condition of a coalgebra,
i.e., that the composite
\[ R \ \xra{\ P\ } \ \exp(R) \ \xra{\ \eta_R\ } \ R\]
is the identity. The transitivity relation is equivalent to 
\[ \exp(P)\circ P  \ = \ \kappa_R \circ P \ ,  \]
the coassociativity condition of a coalgebra.

Part~(ii) is essentially reading part~(i) backwards, and we omit the details.
\end{proof}

The interpretation of global power functors as coalgebras over
a comonad has some useful consequences.
In general, the forgetful functor from any category of coalgebras 
to the underlying category has a right adjoint `cofree' functor. In particular,
colimits in a category of coalgebras are created in the underlying category. 
In our situation that means:

\begin{cor}\label{cor-colimits in power functors}
  \begin{enumerate}[\em (i)]
  \item 
    Colimits in the category of global power functors
    exist and are created in the underlying category of global Green functors.  
  \item
    For every global Green functor $R$, the maps
    \[ P^m \ : \ \exp(R;G) \ \to  \exp(R;\Sigma_m\wr G)\ ,\
    P^m(x)\ = \ \kappa(x)_m=(\Psi_{k,m}^*(x_{k m}))_{k\geq 0}\]
    make the global Green functor $\exp(R)$
    into a global power functor.
  \item When viewed as a functor to the category of
    global power functors as in~{\em (ii)}, the functor $\exp$ is right adjoint to the
    forgetful functor.
  \end{enumerate}
\end{cor}

\begin{eg}[Coproducts]\label{eg:coproduct of power functors}
We let $R$ and $S$ be two global Green functors. Global Green functors
are the commutative monoids, with respect to the box product,
in the category of global functors. 
So the box product $R\Box S$ is the coproduct in the category of global
Green functors, with multiplication defined as the composite
\[ R\Box S\Box R\Box S \ \xra{R\Box\tau_{S,R} \Box S} \
R\Box R\Box S\Box S \ \xra{\mu_R\Box\mu_S} \ R\Box S \ .  \]
If $P:R\to \exp(R)$ and $P':S\to\exp(S)$
are global power structures on $R$ and $S$, 
then $R\Box S$ has preferred power operations specified by
the morphism of global Green functors
\[ R\Box S \ \xra{\ P\Box P'\ } \
 \exp(R)\Box \exp(S) \ \to \  \exp(R\Box S) \ ,\]
where the second morphism is the canonical one from the coproduct
of $\exp$ to the values of $\exp$ at a coproduct.
With these power operations, $R\Box S$ becomes a coproduct 
of $R$ and $S$ in the category of global power functors,
by Corollary \ref{cor-colimits in power functors}~(i).

This abstract definition of the power operations on $R\Box S$ 
can be made more explicit. Indeed, the power operations on $R\Box S$ 
are determined by  the formula
\[ P^m(x\times y) \ = \ \Delta^*(P^m(x)\times P^m(y)) \] 
for all compact Lie groups $G$ and $K$ and classes $x\in R(G)$ and $y\in S(K)$, 
and by the relations of the power operations.
Here $\Delta:\Sigma_m\wr(G\times K)\to(\Sigma_m\wr G)\times(\Sigma_m\wr K)$
is the diagonal monomorphism \eqref{eq:wreath_diagonal}.

The coproduct of global power functors is realized by the coproduct
of ultra-commutative ring spectra in the following sense.
If $E$ and $F$ are ultra-commutative ring spectra, then the ring spectra
morphisms $E\to E\sm F$ and $F\to E\sm F$ induce morphisms
of global power functors $\upi_0(E) \to \upi_0(E\sm F)$
and $\upi_0 (F) \to \upi_0(E\sm F)$, and together they define a morphism
from the coproduct of global power functors
\[ \upi_0 (E)\, \Box\, \upi_0 (F) \ \to \ \upi_0(E\sm F) \ . \]
If $E$ and $F$ are globally connective and at least one of them is flat
as an orthogonal spectrum, then this is an isomorphism of
global functors by Proposition \ref{prop:Box vs smash}, 
hence an isomorphism of global power functors.
\end{eg}

\begin{eg}[Localization of global power functors]
Now we discuss localizations of global power functors.
\index{subject}{localization!of global Green functors}
We first consider a global Green functor $R$ and a multiplicative
subset $S\subseteq R(e)$ in the `underlying ring', i.e., the value
at the trivial group. We define a global Green functor $R[S^{-1}]$
and a morphism of global Green functors $i:R\to R[S^{-1}]$.
The value at a compact Lie group $G$ is the ring
\[ R[S^{-1}](G)\ = \ R(G)[p_G^*(S)^{-1}]\ , \]
the localization of the ring $R(G)$ at the multiplicative subset
obtained as the image of $S$ under the inflation homomorphism $p_G^*:R(e)\to R(G)$.
The value of the morphism $i$ at $G$
is the localization map $R(G)\to R(G)[p_G^*(S)^{-1}]$.
If $\alpha:K\to G$ is a continuous homomorphism, the relation
$p_G\circ \alpha=p_K$ implies that the ring homomorphism
\[ \alpha^* \ : \ R(G)\ \to \ R(K) \]
takes the set $p_G^*(S)$ to the set $p_K^*(S)$.
So the universal property of localization provides a unique ring homomorphism
\[ \alpha[S^{-1}]^*\ : \   R[S^{-1}](G)\ = \ R(G)[p_G^*(S)^{-1}] \ \to \ R(K)[p_K^*(S)^{-1}]
\ = \ R[S^{-1}](K) \]
such that $\alpha[S^{-1}]^*\circ i(G)=i(K)\circ\alpha^*$.
Again by the universal property of localizations, this data produces 
a contravariant functor from the category $\Rep$ to the category of commutative rings.

Now we let $H$ be a closed subgroup of $G$. We consider $R(H)$
as a module over $R(G)$ via the restriction homomorphism $\res^G_H:R(G)\to R(H)$.
Because $\res^G_H(p_G^*(S))=p_H^*(S)$,
the localization $R(H)[p_H^*(S)^{-1}]=R[S^{-1}](H)$
is also a localization of $R(H)$ at $p_G^*(S)$ as an $R(G)$-module.
The reciprocity formula means that the transfer map $\tr_H^G:R(H)\to R(G)$
is a homomorphism of $R(G)$-modules. The composite $R(G)$-linear map
\[ R(H)\ \xra{\ \tr_H^G\ }\ R(G)\ \xra{i(G)} \ R(G)[p_G^*(S)^{-1}]\ = \ R[S^{-1}](G)\]
thus extends over a unique $R[S^{-1}](G)$-linear map
\[ \tr[S^{-1}]_H^G \ : \
R[S^{-1}](H)\ = \ R(H)[p_G^*(S)^{-1}]\ \to \ R[S^{-1}](G)\]
such that $\tr[S^{-1}]_H^G\circ i(H)=i(G)\circ\tr_H^G$.
Reciprocity for $\tr[S^{-1}]_H^G$ is equivalent to $R[S^{-1}](G)$-linearity.
The other necessary properties of transfers, such as transitivity, 
compatibility with inflations, vanishing for infinite Weyl groups, 
and the double coset formula all follow from corresponding properties
in the global Green functor $R$ and the universal property of localization.
So altogether this shows that the objectwise localizations assemble into
a new global Green functor $R[S^{-1}]$; the homomorphisms $i(G)$
altogether form a morphism of global Green functors $i:R\to R[S^{-1}]$
by construction. The following universal property is also straightforward
from the universal property of localizations of commutative rings and modules.
\end{eg}

\begin{prop}\label{prop:universal prop localization}\index{subject}{localization!of global power functors}
Let $R$ be a global Green functor and $S$ a multiplicative subset of
the underlying ring $R(e)$. Let $f:R\to R'$ be a morphism of global Green functors
such that all elements of the set $f(e)(S)$ are invertible in the ring $R'(e)$.
Then there is a unique homomorphism of global Green functors $\bar f:R[S^{-1}]\to R'$
such that $\bar f i=f$.
\end{prop}

Now we let $R$ be a global power functor. If we want the localization
$R[S^{-1}]$ to inherit power operations, then we need an extra hypothesis
on the multiplicative subset $S$.

\begin{theorem}\label{thm:localize global power functor}
Let $R$ be a global power functor and $S$ a multiplicative subset of
the underlying ring $R(e)$. Suppose that the multiplicative subset
\[ P^m(S) \ \subset \ R(\Sigma_m)\]
becomes invertible in the ring $R[S^{-1}](\Sigma_m)$, for every $m\geq 1$.  
\begin{enumerate}[\em (i)]
\item 
There is a unique extension of the global Green functor $R[S^{-1}]$
to a global power functor such that the morphism $i:R\to R[S^{-1}]$
is a morphism of global power functors.
\item
Let $f:R\to R'$ be a morphism of global power functors
such that all elements of the set $f(e)(S)$ are invertible in the ring $R'(e)$.
Then there is a unique homomorphism of global power functors $\bar f:R[S^{-1}]\to R'$
such that $\bar f i=f$.
\end{enumerate}
\end{theorem}
\begin{proof}
(i) We use the comonadic description of global power functors given in  
Theorem \ref{thm:power comonadic}. This exhibits the power operations of $R$
as a morphism of global Green functors $P:R\to\exp(R)$.

The multiplication in the
ring $\exp(R;e)$ is componentwise, and for every $s\in S$
the element $P^m(s)$ has a multiplicative inverse $t_m\in R[S^{-1}](\Sigma_m)$ 
by hypothesis.
We claim that the sequence $t=(t_m)_{m\geq 0}$ is again exponential,
and hence an element of the ring $\exp(R[S^{-1}],e)$. Indeed,
for $0<i<m$ we have
\begin{align*}
 \Phi_{i,m-i}^*(t_m)\cdot (P^i(s)\times P^{m-i}(s))\ &= \   
 \Phi_{i,m-i}^*(t_m)\cdot \Phi_{i,m-i}^*(P^m(s))\\ 
&= \ \Phi_{i,m-i}^*(t_m\cdot P^m(s))\ = \    \Phi_{i,m-i}^*(1)\ = \   1\ .
\end{align*}
Since $t_i\times t_{m-i}$ is also inverse to $P^i(s)\times P^{m-i}(s)$
and inverse are unique, we conclude that $\Phi^*_{i,m-i}(t_m)=t_i\times t_{m-i}$.
This shows that the inverses form another exponential sequence.

The relation
\[ P(s)\cdot t \ = \ (P^m(s)\cdot t_m)_m\ = \ 1 \]
holds in the ring $\exp(R[S^{-1}];e)$, by construction. So the composite morphism
\[ R \ \xra{\ P \ }\ \exp(R)\ \xra{\exp(i)}\ \exp(R[S^{-1}]) \]
takes $S$ to invertible elements. The universal property of 
Proposition \ref{prop:universal prop localization} thus provides 
a unique morphism of global Green functors
\[ P[S^{-1}] \ : \ R[S^{-1}]\ \to \ \exp(R[S^{-1}]) \]
such that $P[S^{-1}]\circ i=\exp(i)\circ P$.
Now we observe that
\[  \eta_{R[S^{-1}]}\circ P[S^{-1}]\circ i \ = \ 
 \eta_{R[S^{-1}]}\circ \exp(i)\circ P \ = \ 
 i\circ \eta_R\circ P \ = \ i \ ; \]
the uniqueness clause in the universal property of 
Proposition \ref{prop:universal prop localization} then implies that
$\eta_{R[S^{-1}]} \circ P[S^{-1}]$ is the identity of $R[S^{-1}]$.
Similarly,
\begin{align*}
  \exp(P[S^{-1}])\circ P[S^{-1}]\circ i \ &= \ 
  \exp(P[S^{-1}])\circ \exp(i)\circ P\
= \    \exp(P[S^{-1}]\circ i)\circ P\\
&= \    \exp(\exp(i)\circ P)\circ P\
= \    \exp(\exp(i))\circ \exp(P)\circ P\\
&= \    \exp(\exp(i))\circ \kappa_R\circ P\
= \    \kappa_{R[S^{-1}]}\circ \exp(i)\circ  P\\
&= \ \kappa_{R[S^{-1}]}\circ P[S^{-1}]\circ i\ .
\end{align*}
The uniqueness clause then implies that
\[   \exp(P[S^{-1}])\circ P[S^{-1}]\ = \  \kappa_{R[S^{-1}]}\circ P[S^{-1}] \ .\]
So the morphism $P[S^{-1}]$ is a coalgebra structure over the $\exp$ comonad.
The relation $P[S^{-1}]\circ i=\exp(i)\circ P$ then says that $i$ is a morphism
of $\exp$-coalgebras, so this completes the proof.

(ii) Since morphisms of global power functors are in particular morphisms
of global Green functors, the uniqueness clause follows from
Proposition \ref{prop:universal prop localization}.
If $f:R\to R'$ is a morphism of global power functors
such that all elements of $f(e)(S)$ are invertible in the ring $R'(e)$,
then Proposition \ref{prop:universal prop localization} provides a
homomorphism of global Green functors $\bar f:R[S^{-1}]\to R'$
such that $\bar f i=f$. We need to show that $\bar f$ is also compatible
with the power operations.
We let
\[ P[S^{-1}] \ : \ R[S^{-1}]\ \to \ \exp(R[S^{-1}]) \text{\qquad and\qquad}
 P' \ : \ R'\ \to \ \exp(R') \]
be the morphisms of global Green functors that encode the $\exp$-coalgebra structures.
We observe that
\begin{align*}
  \exp(\bar f)\circ P[S^{-1}]\circ i \ &= \ 
  \exp(\bar f)\circ \exp(i)\circ P\\
&= \    \exp(f)\circ P\ = \    P'\circ f\ = \ P'\circ \bar f\circ i\ ;
\end{align*}
the third equation is the hypothesis that $f$ is a morphism of global power functors.
This is an equality between morphisms of global Green functors, so the 
uniqueness of Proposition \ref{prop:universal prop localization}
shows that $\exp(\bar f)\circ P[S^{-1}]=P'\circ \bar f$, i.e.,
$\bar f$ is a morphism of $\exp$-coalgebras.
\end{proof}

\begin{eg}[Localization at a subring of $\mQ$]\label{eg:rationalize GPF}
We use Theorem \ref{thm:localize global power functor}
to show that global power functors can always be rationalized;
more generally, power operations `survive' localization 
at any subring of the ring $\mQ$ of rational numbers.
We consider a global power functor $R$ and a natural number $n\geq 2$.
We claim that for every $m\geq 1$ the element $P^m(n)$ of $R(\Sigma_m)$ becomes
invertible in the ring $R(\Sigma_m)[1/n]$. We argue by induction over $m$;
we start with $m=1$, where the relation $P^1(n)=n$ shows the claim.
Now we consider $m\geq 2$ and assume that for all $0\leq i < m$
the element $P^i(n)$ has an inverse $t_i$ in the ring $R(\Sigma_i)[1/n]$. 
Then for every $n$-tuple $(i_1,\dots,i_n)$ that
satisfies $0\leq i_j <m$ and $i_1+\dots+i_n=m$ we get the relation
\begin{align*}
\res^{\Sigma_m}_{\Sigma_{i_1}\times\dots\times\Sigma_{i_n}} (P^m(n))\cdot(t_{i_1}&\times\dots\times t_{i_n} )\ 
= \ (P^{i_1}(n)\times\dots\times P^{i_n}(n)) \cdot(t_{i_1}\times\dots\times t_{i_n} )\\  
&= \ (P^{i_1}(n)\cdot t_{i_1})\times\dots\times (P^{i_n}(n)\cdot t_{i_n} )\  
= \ 1\times\dots\times 1
\end{align*}
in the ring $R(\Sigma_{i_1}\times\dots\times\Sigma_{i_n})[1/n]$. Thus 
\begin{align*}
  P^m(n)\ &= \ P^m(\underbrace{1+\dots+1}_n)\ = \ 
\sum_{i_1+\dots+i_n=m}  \tr_{\Sigma_{i_1}\times\dots\times\Sigma_{i_n}}^{\Sigma_m}( P^{i_1}(1)\times\dots\times P^{i_n}(1)) \\
&= \ 
n\cdot P^m(1)\ + \sum_{i_1+\dots+i_n=m,\ i_j<m}  \tr_{\Sigma_{i_1}\times\dots\times\Sigma_{i_n}}^{\Sigma_m}( 1\times\dots\times 1) \\
&= \ 
n\ + \sum_{i_1+\dots+i_n=m,\ i_j<m}  \tr_{\Sigma_{i_1}\times\dots\times\Sigma_{i_n}}^{\Sigma_m}(\res^{\Sigma_m}_{\Sigma_{i_1}\times\dots\times\Sigma_{i_n}} (P^m(n))\cdot(t_{i_1}\times\dots\times t_{i_n} )) \\
&= \ 
n\ + \sum_{i_1+\dots+i_n=m,\ i_j<m} P^m(n)\cdot \tr_{\Sigma_{i_1}\times\dots\times\Sigma_{i_n}}^{\Sigma_m}(t_{i_1}\times\dots\times t_{i_n} ) 
\end{align*}
Rearranging the terms gives
\begin{align*}
  P^m(n)\cdot & \left(1- \sum_{i_1+\dots+i_n=m,\ i_j< m} 
\tr_{\Sigma_{i_1}\times\dots\times\Sigma_{i_n}}^{\Sigma_m}( t_{i_1}\times\dots\times t_{i_n})\right) \
 = \ n \ .
\end{align*}
So $P^m(n)$ has an inverse in the ring $R(\Sigma_m)[1/n]$,
and this completes the inductive step.

Now we let $S$ be any multiplicative subset of the ring of integers.
By the previous paragraph, Theorem \ref{thm:localize global power functor}
applies and provides a unique structure of global power functor on the
global Green functor $R[S^{-1}]$ such that the morphism $i:R\to R[S^{-1}]$
is a morphism of global power functors.\index{subject}{rationalization!of global power functors}
In particular, if we let $S$ be the set of all positive integers, 
we can conclude that the rationalization $\mQ\tensor R$ of $R$
has a unique structure of global power functor such that the localization
map $R\to \mQ\tensor R$ is a  morphism of global power functors.  
\end{eg}

\begin{rk}[Monadic description of global power functors] 
Now we explain that the category of global power functors in not only comonadic,
but also monadic over the category of global Green functors.
In fact, both categories are examples of algebras over 
{\em multisorted algebraic theories} (also called {\em colored theories}).
The `sorts' (or `colors') are the compact Lie groups and the 
content of this claim is that the structure of
global Green functors respectively global power functors
can be specified by giving the values $R(G)$ at every compact Lie group,
together with $n$-ary operations for different $n\geq 0$ and
varying inputs and output, and relations between composites of
those operations.

In the case of global Green functors, the operations to be specified are
\begin{itemize}
\item the constants given by the additive and multiplicative units in the rings $R(G)$,
\item the unary operations given by the additive inverse map in $R(G)$,
the restriction maps and transfers,
\item and the binary operations specifying the addition and multiplication in
the rings $R(G)$.
\end{itemize}
The relations include, among others, the neutrality, associativity and
commutativity of addition and multiplication; the distributivity
in the rings $R(G)$; the additivity of restriction and transfers;
the functoriality of restrictions and transitivity of transfers;
and the double coset and reciprocity formulas.

Global power functors have additional unary operations, the power operations,
and additional relations as listed in Definition \ref{def:power functor}.
\end{rk}

\begin{prop}\label{prop:GPF are monadic}
The forgetful functor from the category of global power functors
to the category of global Green functors has a left adjoint.  
The category of global power functors is isomorphic to the
category of algebras over the monad of this adjunction.
\end{prop}
\begin{proof}
For the existence of the left adjoint we have to show that for every
global Green functor $R$ the functor
\[ \GlPow\ \to \ \text{(sets)} \ , \quad
S \ \longmapsto \ \GlGre(R, S)\]
is representable. This is a formal consequence of the existence of free
global power functors, colimits of global power functors and the fact
that global power functors are a multi-sorted theory. 
We explain this in more detail, without completely formalizing the argument. 

We choose a set of compact Lie groups $\{K_i\}_{i\in I}$
that contains one compact Lie group from every isomorphism class.
We form the global power functor
\[ L \ = \ \Box_{ i\in I, x\in R(K_i)}  C_{i,x}\ ,\]
a coproduct, indexed by all pairs $(i,x)$
consisting of an index $i\in I$ and an element $x\in R(K_i)$,
of free global power functors
\[ C_{i,x} \ = \ C_{K_i}\]
generated by the compact Lie groups $K_i$.
On this free global power functor we impose the minimal amount of relations
so that the maps $R(K_i)\to L(K_i)$ that send $x\in R(K_i)$
to the generator indexed by $(i,x)$ becomes a morphism of global Green functors.
Here `imposing relations' means that we form 
another box product $L'$  of free global power functors,
with one box factor for each relation between elements in the various sets $R(K_i)$.
For example, we include one factor for the sum of each pair 
of elements in the same set $R(K_i)$,
another factor for the product of each pair 
of elements in the same set $R(K_i)$,
and more factors for zero elements, multiplicative units, 
all restriction relations  and all transfers relations. 
Then we form a coequalizer,
in the category of global power functors
\[ \xymatrix{ L' \ar@<.4ex>[r] \ar@<-.4ex>[r] & L \ar[r] & F } \]
where the two morphisms from $L'$ to $L$ restrict, on each box factor,
to the morphism that represents the respective relation.
The resulting global power functor then represents 
the functor $\GlGre(R,-)$,
so we can take $F$ as the value of the left adjoint on $R$.

Since global power functors are equivalent to
the coalgebras over the $\exp$-comonad, the forgetful functor
creates all colimits, in particular coequalizers.
So by Beck's monadicity theorem (see for example \cite[VI.7 Thm.\,1]{maclane-working}),
the tautological functor from global power functors to algebras over
the adjunction monad is an isomorphism of categories.
\end{proof}

\begin{eg}[Limits]
The category of global power functors has limits, and they are 
defined `groupwise'. 
A special case of a limit is the product of global power functors,
which is realized by the product of ultra-commutative ring spectra.
Indeed, if $E$ and $F$ are ultra-commutative ring spectra, then so is
the product $E\times F$ of the underlying orthogonal spectra,
and the canonical map
\[ \upi_0(E\times F) \ \to \ \upi_0( E) \times \upi_0( F) \]
is an isomorphism of global power functors
(by Corollary \ref{cor-wedges and finite products}~(ii)).
\end{eg}

\section{Examples}\label{sec:GPF examples}

In this section we discuss various examples of 
and constructions with global power functors, and how these are
realized topologically by ultra-commutative ring spectra.
These examples include the Burnside ring global power functor
(Example \ref{eg:Burnside global power}), the global functor represented
by an abelian compact Lie group (Proposition \ref{prop:unique power bA(A,-)}),
and constant global power functors (Example \ref{eg:constant global Green}).
The orthogonal spectrum $\Hc A$ consisting of the $A$-linearizations
of spheres (Construction \ref{con:HM}) tries to be an Eilenberg-Mac Lane spectrum
for the constant global functor, and it is so on finite groups.
Closely related to $\Hc\mZ$ is the infinite symmetric product spectrum $\Spinf$,
see Example \ref{eg:infinite symmetric product}; this ultra-commutative
ring spectrum comes with a filtration by finite symmetric products,
realizing an interesting filtration of the Burnside ring global
functor on $\upi_0$, see \eqref{eq:pi_0_Sp^n}. 
We close this section with a discussion of
the complex representation ring global functor (Example \ref{eg:RU_as_global_power}),
and a global view on `explicit Brauer induction' (Remark \ref{rk:Brauer induction}).

\begin{eg}[Burnside ring global functor]\index{subject}{Burnside ring global functor}\label{eg:Burnside global power}
The Burnside ring global functor $\mA=\bA(e,-)$ 
is the unit object for the box product of global functors,
and hence an initial object in the category of global Green functors.
Initial objects are examples of colimits,
so Corollary \ref{cor-colimits in power functors}~(i) 
implies that $\mA$ has a unique structure of global power functor. 
Indeed, there is a unique morphism $P:\mA\to\exp(\mA)$ of global Green functors
(since $\mA$ is initial), and the coalgebra diagrams commute
(again since $\mA$ is initial).
With these power operations, $\mA$ is also an initial global power functor.

We can make the power operations in the Burnside ring global functor more explicit. 
Indeed, the group $\mA(G)$ is free abelian with a basis
given by the elements $t_H=\tr_H^G(p_H^*(1))$ for every conjugacy class
of closed subgroups $H\leq G$ with finite Weyl group, where $p_H:H\to e$
is the unique homomorphism. On these generators,
the naturality properties of a global power functor force the power operations to be 
\begin{align}\label{eq:P^m on t_H}
 P^m(t_H) \ &= \ P^m(\tr_H^G(p_H^*(1)))\ = \ 
\tr_{\Sigma_m\wr H}^{\Sigma_m\wr G}((\Sigma_m\wr p_H)^*(P^m(1)))\\ 
&= \  \tr_{\Sigma_m\wr H}^{\Sigma_m\wr G}((\Sigma_m\wr p_H)^*(p_{\Sigma_m}^*(1)))\
= \  \tr_{\Sigma_m\wr H}^{\Sigma_m\wr G}(p_{\Sigma_m\wr H}^*(1))\ = \   
t_{\Sigma_m\wr H}\ .\nonumber
\end{align}
This determines the power operations in general by the additivity property,
and also shows the uniqueness.

When restricted to {\em finite} groups, the ring $\mA(G)$ is isomorphic to 
the Grothen\-dieck group of finite $G$-sets, and in this description 
the power operations are given by raising a finite $G$-set to a power, i.e.,
the power map
\[ P^m \ : \ \mA(G) \ \to \ \mA(\Sigma_m\wr G)\]
takes the class of a finite $G$-set $S$ to the
class of the $(\Sigma_m\wr G)$-set $S^m$.
Indeed, for the additive generator $[G/H]=t_H$
of $\mA(G)$ this is the relation \eqref{eq:P^m on t_H},
and for general finite $G$-sets it follows from the additivity formula
for power operations and the fact that for two finite $G$-sets $S$ and $T$
the power $(S\amalg T)^m$ is $(\Sigma_m\wr G)$-equivariantly isomorphic to
the coproduct 
\[ {\coprod}_{k=0}^m \ 
(\Sigma_m\wr G)\times_{(\Sigma_k\wr G)\times (\Sigma_{m-k}\wr G)} (S^k\times T^{m-k})\ .\]
The canonical power operations in the Burnside ring global
functor correspond to the  homotopy theoretic power operations 
for the global sphere spectrum.\index{subject}{sphere spectrum}
Indeed, since $\mA$ is initial in both the category
of global Green functors and in the category of global power functors,
any isomorphism of global Green functors is automatically
compatible with power operations. In other words, we can conclude that the square
\[\xymatrix{ \mA(G) \ar[r]^-{P^m} \ar[d]_\iso & 
\mA(\Sigma_m\wr G) \ar[d]^\iso\\
\pi_0^G(\mS) \ar[r]_-{P^m} &  \pi_0^{\Sigma_m\wr G} (\mS) } \]
commutes for all $G$ and $m$ without having to go back to the
definition of the operations in $\upi_0(\mS)$;
the vertical maps are the action on the multiplicative units.
\end{eg}

The previous example generalizes to the represented global functor
$\bA_A=\bA(A,-)$ for every {\em abelian} compact Lie group,
which has a preferred structure of a global power functor.
The following proposition is a stable analog of 
Proposition \ref{prop:unique power Rep(-,A)}, which describes
a structure  of global power monoid on the $\Rep$-functor represented by $A$.

\begin{prop}\label{prop:unique power bA(A,-)}\index{subject}{compact Lie group!abelian}
Let $A$ be an abelian compact Lie group.
The represented global functor $\bA_A$
has a unique structure of global power functor subject to the following two
conditions:
\begin{itemize}
\item The multiplication is the composite
\[ \bA_A\Box\bA_A \ \iso \ \bA_{A\times A}\ \xra{\bA(\mu^*,-)}\ 
\bA_A\ ,\]
where $\mu:A\times A\to A$ is the multiplication
and $\mu^*\in \bA(A,A\times A)$ is the associated restriction morphism.
\item The power operations satisfy
\[ P^m(1_A)\ = \ p_m^* \text{\quad in\quad} \bA(A,\Sigma_m\wr A) \ ,   \]
the inflation operation of the continuous homomorphism 
\[  p_m\ :\ \Sigma_m\wr A\ \to\  A \ , \quad
(\sigma;\,a_1,\dots,a_m)\ \longmapsto \ a_1\cdot\ldots\cdot a_m \ . \]
\end{itemize}
Moreover, for every global power functor $R$ the evaluation map  
\[ \GlPow(\bA_A,R) \ \to \ R(A)\ , \quad f \ \longmapsto \ f_A(1_A)\]
is injective with image the set of those $x\in R(A)$
such that $P^m(x)=p_m^*(x)$ for all $m\geq 1$.
\end{prop}
\begin{proof}
Since the convolution product of two represented functors is again represented,
morphisms of global functors $\bA_A\Box\bA_A\to\bA_A$
are determined by their effect on 
the element $1_A\Box 1_A\in (\bA_A\Box\bA_A)(A\times A)$.
The first condition fixes this to be the operation $\mu^*$.
So there is at most one structure of global Green functor on $\bA_A$
that satisfies the first condition.
Since this representable global functor is freely generated by the identity
$1_A$ in $\bA(A,A)$, the power operations are all determined by naturality
from the effect on this generator. So together we have shown that there is at most
one structure of global power functor that satisfies the two conditions.

Now we construct the desired structure.
The associativity, commutativity and unitality 
of the multiplication $\mu:A\times A\to A$ readily imply that
the associated multiplication of $\bA_A$ is 
associative, commutative and unital. So this multiplication makes $\bA_A$
into a global Green functor. The relations
\[  p_m\circ \Phi_{k,m-k}\ = \ p_k\cdot p_{m-k}\ : \ 
(\Sigma_k\wr A)\times(\Sigma_{m-k}\wr A) \ \to \ A\]
imply that the sequence $(p_m^*)_{m\geq 0}$ is exponential.
So the Yoneda lemma provides a unique morphism of global functors
\[ P \ =\ (P^m)_{m\geq 0}\ : \ \bA_A\ \to \ \exp(\bA_A) \text{\qquad such that\qquad}
P^m(1_A)\ = \  p_m^* \ .\]
We claim that $P$ is a morphism of global Green functors.
The unitality follows from the fact that the restriction of $p_m$
to $\Sigma_m\wr e$ is the trivial homomorphism.
Since the global functor $\bA_A\Box\bA_A$
is representable by $A\times A$, the Yoneda lemma reduces 
the multiplicativity property to the relation
\[ \exp(\bA_A,\mu^*)((p^*_m)_m) \ = \ 
\exp(\bA_A,q_1^*)((p_m^*)_m)\ \cdot\ \exp(\bA_A,q_2^*)((p_m^*)_m)  \]
in the ring $\exp(\bA_A,A\times A)$, where $q_1,q_2:A\times A\to A$
are the two projections. This relation, in turn,
is a consequence of the commutativity of
the following diagram of group homomorphisms:
\[ \xymatrix@C=12mm{ 
\Sigma_m\wr(A\times A)\ar[d]_\Delta \ar[r]^-{\Sigma_m\wr \mu} &
\Sigma_m\wr A\ar[d]^{p_m}\\
(\Sigma_m\wr A)\times (\Sigma_m\wr A)\ar[r]_-{p_m\cdot p_m} & A
} \]
Here the left vertical map is a diagonal, given by
\[\Delta(\sigma;\,(a_1,\bar a_1),\dots,(a_m,\bar a_m))\ =\
((\sigma;\,a_1,\dots,a_m),(\sigma;\,\bar a_1,\dots,\bar a_m)) \ .  \]
We conclude the construction of the global power structure by showing 
that the morphism $P:\bA_A\to\exp(\bA_A)$ provides a coalgebra structure
over the $\exp$-comonad.
It suffices to show
the relations $\eta_{\bA_A}\circ P=\Id$ and $\exp(P)\circ P=\kappa_{\bA_A}\circ P$ 
on the universal element $1_A$ in $\bA_A(A)$,
one more time by the Yoneda lemma.
For the first one, this boils down to the fact 
that $p_1:\Sigma_1\wr A\to A$ is the preferred isomorphism.
The second one follows from the commutativity of the following square:
\[
\begin{gathered}[b]
 \xymatrix@C=12mm{
\Sigma_k\wr\Sigma_m\wr A\ar[r]^-{\Psi_{k,m}}\ar[d]_{\Sigma_k\wr p_m} &
\Sigma_{k m}\wr A\ar[d]^{p_{k m}} \\
\Sigma_k\wr A\ar[r]_-{p_k}& A}\\[-\dp\strutbox]
\end{gathered}
\qedhere
 \]
\end{proof}

The global power functor $\bA_A$ described
in Proposition \ref{prop:unique power bA(A,-)} is realized by 
an ultra-commutative ring spectrum.
Indeed, Construction \ref{con:multiplicative B G}
provides a global classifying space $R(B A)$ for the
abelian compact Lie group $A$ with an ultra-commutative multiplication.\index{subject}{global classifying space!of an abelian compact Lie group}  
The associated unreduced suspension spectrum
$\Sigma^\infty_+ R(B A)$ is then an ultra-commutative ring spectrum.
Theorem \ref{prop:B_gl represents}
shows that $\upi_0(\Sigma^\infty_+ R(B A))$ is 
freely generated, as a global functor,
by the stable tautological class $e_A\in \pi_0^A(\Sigma^\infty_+ R(B A))$,
the stabilization of the unstable tautological class $u_A$.
The characterization of the multiplication and power operations on $\bA_A$
only involve restriction maps (but no transfers); so the fact that
$\upi_0(\Sigma^\infty_+ R(B A))$ 
realizes the global power structure of
Proposition \ref{prop:unique power bA(A,-)} follows from the
corresponding unstable relations in the global power monoid
$\upi_0( R(B A))$, as explained in Example \ref{con:multiplicative B G}. 

\medskip

For a compact Lie group $G$ we define
\[ u_G^{ucom}\ \in \ \pi_0^G(\Sigma^\infty_+\mP(B_{\gl}G))\]
as the image of the unstable tautological class
$u_G\ \in \ \pi_0^G( B_{\gl}B)$, defined in \eqref{eq:tautological_class},
under the two composites in the commutative square
\[ \xymatrix@C=15mm{ 
\pi_0^G( B_{\gl}B)\ar[r]^-{\pi_0^G(\eta)}\ar[d]_{\sigma^G} & 
\pi_0^G( \mP(B_{\gl}B))\ar[d]^{\sigma^G}\\
\pi_0^G(\Sigma^\infty_+ B_{\gl}B) \ar[r]_-{\pi_0^G(\Sigma^\infty_+ \eta)}&
\pi_0^G(\Sigma^\infty_+ \mP(B_{\gl}B))} \]
where $\eta:B_{\gl} G\to \mP(B_{\gl} G)$ is the adjunction unit, i.e., the
inclusion of the homogeneous summand for $m=1$. 
The stabilization map $\sigma^G$ was defined in \eqref{eq:sigma_map}.
Here, as usual, we have made an implicit choice of non-zero faithful 
$G$-representation $V$, and the global classifying space is $B_{\gl}G=\bL_{G,V}$.
We already know that
\begin{itemize}
\item the class $u_G$ freely generates $\upi_0(B_{\gl}G)$ as a Rep-functor
  (Proposition \ref{prop:fix of global classifying}~(ii)),
\item the class $u_G^{umon}=\pi_0^G(\eta)(u_G)$ freely generates
  $\upi_0(\mP(B_{\gl}G))$ as a global power monoid (Theorem \ref{thm:umon operation}~(ii)), and
\item the class $e_G=\sigma^G(u_G)$ freely generates
  $\upi_0(\Sigma^\infty_+ B_{\gl}G)$ as a global functor
  (Proposition \ref{prop:B_gl represents}).\index{subject}{global classifying space}
\end{itemize}
\begin{itemize}
\item $\upi_0(B_{\gl}G)$ is the freely generated by $u_G$ as a  Rep-functor
  (Proposition \ref{prop:fix of global classifying}~(ii)),
\item $\upi_0(\mP(B_{\gl}G))$ is the freely generated by $u_G^{umon}=\pi_0^G(\eta)(u_G)$ 
  as a global power monoid (Theorem \ref{thm:umon operation}~(ii)), and
\item $\upi_0(\Sigma^\infty_+ B_{\gl}G)$ is the freely generated by $e_G=\sigma^G(u_G)$ 
  as a global functor
  (Proposition \ref{prop:B_gl represents}).\index{subject}{global classifying space}
\end{itemize}
The next theorem complements these results and shows that
the global power functor $\upi_0(\Sigma^\infty_+\mP(B_{\gl}B))$
is freely generated by the class $u_G^{ucom}$.
The free global power functor\index{subject}{global power functor!free}
$C_G$ generated by the compact Lie group $G$ 
was introduced in Example \ref{eg:free global power functor}.
The underlying global functor is
\[ C_G \ = \ {\bigoplus}_{m\geq 0}\ \bA(\Sigma_m\wr G,-) \ , \]
and there is a preferred element $1_G\in C_G(G)$, the identity operation
of $\pi_0^G$ in the summand indexed by $m=1$.
The adjective `free' is justified by Proposition \ref{prop:what C_K represents}~(ii),
which says that for every global power functor $R$ and every element $x\in R(G)$
there is a unique morphism of global power functors $f:C_G\to R$ such that $f(1_G)=x$.

\begin{prop}\label{cor:B_gl free ucom}
  Let $G$ be a compact Lie group.
  \begin{enumerate}[\em (i)]
  \item The unique morphism of global power functors
\[ \varphi\ : \ C_G \ \to \ \upi_0(\Sigma^\infty_+ \mP(B_{\gl} G)) \]
satisfying $\varphi(1_G)=u_G^{ucom}$ is an isomorphism.
  \item
 For every global power functor $R$ and every element $x\in R(G)$
  there is a unique morphism of global power functors
  $f:\upi_0(\Sigma^\infty_+ \mP(B_{\gl} G))\to R$ such that $f(u_G^{ucom}) =  x$.
  \end{enumerate}
\end{prop}
\begin{proof}
  (i)
  The ultra-commutative monoid $\mP(B_{\gl}G)$ is the disjoint union
  of the orthogonal spaces $\mP^m(B_{\gl}G)$ for $m\geq 0$.
  Moreover, $\mP^m(B_{\gl}G)$ is a global classifying space for the group
  $\Sigma_m\wr G$, in such a way that the class $[m](u_G)\in\pi_0^G(\mP^m(B_{\gl}G))$
  matches the unstable tautological class $u_{\Sigma_m\wr G}$, 
  by Example \ref{eg:P B_gl G}.
  So by Proposition \ref{prop:B_gl represents} the morphism of global functors
  \[ \psi_m \ : \ \bA(\Sigma_m\wr G,-)\ \to \ \upi_0(\Sigma^\infty_+\mP^m(B_{\gl}G))\]
  classified by $\sigma^{\Sigma_m\wr G}([m](u_G))=P^m(\sigma^G(u_G))$ in 
  $\pi_0^{\Sigma_m\wr G}(\Sigma^\infty_+\mP^m(B_{\gl}G))$
  is an isomorphism.
  The functor $\Sigma^\infty_+$ takes a disjoint union of orthogonal
  spaces to a wedge of orthogonal spectra, and
  equivariant stable homotopy groups take wedges to sums.
  Hence the morphism  
  \begin{align*}
    \psi\ = \ {\bigoplus}_{m\geq 0}\psi_m \ : \ 
  C_G \ = \ &{\bigoplus}_{m\geq 0}\,\bA(\Sigma_m\wr G,-)\\ 
  \to \ &{\bigoplus}_{m\geq 0 }\,\upi_0(\Sigma^\infty_+\mP^m(B_{\gl}G))\ = \ 
  \upi_0(\Sigma^\infty_+\mP(B_{\gl}G))
  \end{align*}
  is an isomorphism of global functors.

  The power operations in $C_G$ satisfy $[m](1_G)=1_{\Sigma_m\wr G}$;
  since $\varphi$ is a morphism of global power functors, it also satisfies
  \[ \varphi(1_{\Sigma_m\wr G})\ = \  \varphi([m](1_G))\  = \ [m](\varphi(1_G))
  \  = \  [m](u^{umon}_G)\ = \ \sigma^G([m](u_G))\ .\]
  So the morphisms $\varphi$ and $\psi$ coincide on the classes $1_{\Sigma_m\wr G}$
  for all $m\geq 0$. Since these classes generated $C_G$ as a global functor,
  we conclude that $\varphi=\psi$. Since $\psi$ is an isomorphism, this proves the claim.
  
  Part~(ii) is the combination of~(i) and the freeness property of $C_G$
  from Proposition \ref{prop:what C_K represents}~(ii).
\end{proof}

\begin{eg}[Monoid rings] 
We let $R$ be a global Green functor and $M$ a commutative monoid.
We denote by $R[M]$ the {\em monoid ring functor}; its value at
a compact Lie group $G$ is given by
\[ (R[M])(G) \ = \ R(G)[M] \ ,\]
the monoid ring of $M$ over $R(G)$. The structure as global functor
is induced from the structure of $R$ and constant in $M$.
The multiplication and unit are
induced from the multiplication and units of $R$ and $M$.
The global Green functor $R[M]$ can be characterized as follows
by the functor that it represents:
for every global Green functor $S$, morphisms $R[M]\to S$ biject with pairs consisting 
of a morphism of global Green functors $R\to S$ and a monoid homomorphism
$M\to (S(e),\cdot)$ to the multiplicative monoid of the underlying ring $S(e)$.

Now suppose that $R$ is even a global power functor.
Then $R[M]$ inherits a natural structure as global power functor:
we define the power operation
\[ P^m \ : \  R(G)[M] \ \to \ R(\Sigma_m\wr G)[M] \]
by $P^m(r\cdot x)=P^m(r)\cdot x^m$ for $r\in R(G)$ and $x\in M$, 
and then we extend this
by additivity to general elements in $R(G)[M]$.
In this upgraded setting, the global power functor $R[M]$ 
has a similar characterization as in the previous paragraph:
for a global power functor $S$, we let 
\[ \Mon(S) \ = \ \{ x\in S(e)\ | \ \text{$P^m(x)=(p_{\Sigma_m})^*(x)$ for all $m\geq 1$}\}\]
be the set of {\em monoid-like elements} of $S$;
here $p_{\Sigma_m}:\Sigma_m\to e$ is the unique homomorphism.
The set $\Mon(S)$ is a submonoid of the multiplicative monoid of $S(e)$.
Then morphisms $R[M]\to S$ biject with pairs consisting 
of a morphism of global power functors $R\to S$ and a monoid homomorphism
$M\to \Mon(S)$ to the monoid-like elements of $S$.

Now we let $E$ be an ultra-commutative ring spectrum.
Then the monoid ring spectrum $E[M]=E\sm M_+$\index{subject}{monoid ring spectrum}
is another ultra-commutative ring spectrum, 
see Construction \ref{con:monoid ring spectrum}.
The unit of $M$ induces a morphism of ultra-commutative ring spectra $E\to E[M]$
that induces a morphism of global power functors 
$\upi_0(E)\to \upi_0(E[M])$. Similarly, the unit of $E$ induces a monoid homomorphism
$M\to E(0)\sm M_+$, and this induces a monoid homomorphism $M\to \Mon(\upi_0(E[M]))$
to the  monoid-like elements in the multiplicative monoid of the ring $\pi_0^e(E[M])$.
The universal property of the algebraic monoid ring combines these
two pieces of data into a morphism of global power functors
\begin{equation}\label{eq:monoid ring isomorphism}
 (\upi_0 E)[M]  \ \to \ \upi_0(E[M]) \ .  
\end{equation}
Additively, the left hand side is a direct sum of copies
of $\upi_0(E)$, indexed by the elements of $M$.
Similarly, the orthogonal spectrum $E[M]$ is a wedge of copies
of $E$, indexed by the elements of $M$.
So the right hand side $\upi_0(E[M])$ is also a direct sum of copies of $\upi_0(E)$
(by Corollary \ref{cor-wedges and finite products}~(i)).
We conclude that the morphism of global power functors \eqref{eq:monoid ring isomorphism} is an isomorphism. 
\end{eg}

\begin{eg}[Constant global power functors and Eilenberg-Mac Lane spectra]\label{eg:constant global Green}\index{subject}{Eilenberg-Mac\,Lane spectrum!of a commutative ring}
We let $B$ be a commutative ring. Then the constant global\index{subject}{global power functor!constant}  
functor $\underline B$ (see Example \ref{eg:global functor examples}~(iii)) 
is naturally a global Green functor, via the ring structure of $B$.
For every compact Lie group $G$, every exponential sequence $x\in\exp(\underline{B};G)$ 
satisfies $x_n=(x_1)^n$; so an exponential sequence
is completely determined by the element $x_1$, which can be chosen arbitrarily. 
Hence the morphism
\[ \eta_{\underline B}\ : \ \exp(\underline{B}) \ \to \ \underline{B} \]
is an isomorphism of global Green functors. So when restricted
to constant global Green functors, the exp comonad is isomorphic to the
identity. Thus $\underline{B}$ has a unique structure of coalgebra
over the comonad $\exp$, with structure morphism the inverse of $\eta_{\underline B}$.

The preferred global power structure on the constant global functor $\underline B$
can also be derived more directly.
In fact, since the power operation $P^m$ has to be an equivariant refinement
of the $m$-th power map in $R(G)$ (compare Remark \ref{rk:power operation remarks}), 
the only possibility to define $P^m$ is as
\[ P^m\ : \ \underline B(G)\ = B \ \to \ B = \underline B (\Sigma_m\wr G)\ , \quad
b \ \longmapsto \ b^m \  \ ,\]
the $m$-th power in the ring $B$.

Since $\underline{B}$ is constant, 
every morphism $R\to \underline{B}$ of global functors is determined
by the map $R(e)\to \underline{B}(e)=B$. Moreover, every ring homomorphism
$\psi:R(e)\to B$ extends uniquely to a morphism of global power functors
$\hat\psi:R\to \underline{B}$ by defining its value at $G$ as the composite
\[ R(G)\ \xra{\ \res^G_e\ } \ R(e)\ \xra{\ \psi\ } \ B\ .  \]
In other words, the functor 
\[ \text{(commutative rings)}  \ \to \  \GlPow\ , \quad
B \ \longmapsto \ \underline{B} \]
is right adjoint to the functor that takes a global power functor $R$ 
to the underlying ring $R(e)$.
\end{eg}

The author does not know an explicit pointset level model
for an ultra-commutative ring spectrum that realizes the  constant global power functor
$\underline B$ of a commutative ring $B$. 
The Eilenberg-Mac\,Lane spectrum $\Hc B$, 
discussed in the next construction,
is an ultra-commutative ring spectrum and tries to realize $\underline B$:
the global power functor $\upi_0(\Hc B)$
is indeed constant {\em on finite groups}, but the restriction maps
are not generally isomorphisms, compare Example \ref{eg:HZ for abelian G}.
The morphism of global power functors
$\upi_0(\Hc B)\to\underline{B}$ adjoint to the identification
$\pi^e_0(\Hc B)\iso B$ is thus an isomorphism at finite groups
(and some other compact Lie groups), but {\em not} an isomorphism in general.

\index{subject}{Eilenberg-Mac\,Lane spectrum!of an abelian group|(} 
\index{subject}{Eilenberg-Mac\,Lane spectrum!of an abelian group} 

\begin{construction}\label{con:HM}
Let $A$ be an abelian group.
The {\em Eilenberg-Mac\,Lane spectrum}
$\Hc A$\index{symbol}{$\Hc A$ - {Eilenberg-Mac\,Lane spectrum of an abelian group}} 
is defined at an inner product space $V$ by
\[ (\Hc A)(V) \ = \ A[S^V] \ ,  \]
the reduced $A$-linearization\index{subject}{linearization!of a space} 
of the $V$-sphere. 
The orthogonal group $O(V)$ acts through the action on $S^V$ and
the structure map
$\sigma_{V,W}:S^V\sm (\Hc A)(W)\to (\Hc A)(V\oplus W)$ is given by
\[  S^V\sm A[S^W] \ \to\  A[S^{V\oplus W}] \ , \quad
v\, \sm\, \big( {\sum}_i\ m_i\cdot w_i \big) \ \longmapsto \
{\sum}_i\ m_i\cdot \, (v\sm w_i) \  . \]
Non-equivariantly, $A[S^V]$ is an Eilenberg-Mac\,Lane space of type $(A,n)$, 
where $n=\dim(V)$. 
Hence the underlying non-equivariant homotopy type of $\Hc A$ is that of
an Eilenberg-Mac\,Lane spectrum for $A$.
\end{construction}

As we shall now discuss, the equivariant homotopy groups of $\Hc A$ are 
in general {\em not} concentrated in dimension~0, and hence $\Hc A$ is 
not the Eilenberg-Mac\,Lane spectrum of a global functor.
However, on {\em finite} groups, the equivariant behavior of $\Hc A$ is as expected.
We recall from Definition \ref{def:left_right_induced}
that an orthogonal spectrum is
{\em left induced}\index{subject}{left induced} 
from the global family $\Fin$\index{subject}{global family!of finite groups} 
of finite groups if it is in the essential image of the left adjoint
$L_\Fc:\GH_\Fin\to\GH$ from the $\Fin$-global homotopy category.
The inclusion of generators is a homeomorphism $A\iso A[S^0]=(\Hc A)(0)$
and this induces a bijection
\[ A \ \iso \ \pi_0(A) \ = \ [S^0, (\Hc A)(0)] \ .\]
Since $\Hc A$ is in particular a non-equivariant $\Omega$-spectrum
(for example by Proposition \ref{prop:H M is Fin-Omega} below),
the composite of this bijection with the stabilization map is an 
isomorphism of abelian groups
\[ A \ \iso \ \pi_0^e(\Hc A) \ .\]
The restriction maps
\[ \res^G_e \ : \  \pi_0^G(\Hc A) \ \to \ 
\pi_0^e(\Hc A)\ \iso\  A \ = \ \underline{A}(G) \]
form a morphism of global functors $\upi_0(\Hc A)\to \underline{A}$
to the constant global functor with value $A$.
Since $\Hc A$ is globally connective (by the next proposition), 
there is a unique morphism 
\[ \rho \ : \ \Hc A \ \to \ H\underline{A} \]
in the global stable homotopy category that realizes the morphism on $\pi_0^e$.

\begin{prop}\label{prop:H M is Fin-Omega}
For every abelian group $A$ the Eilenberg-Mac\,Lane spectrum~ $\Hc A$ is 
globally connective, 
left induced from the global family $\Fin$ of finite groups 
and a $\Fin$-$\Omega$-spectrum. 
The morphism $\rho: \Hc A \to H\underline{A}$ is a $\Fin$-global equivalence.
\end{prop}
\begin{proof}
The orthogonal spectrum $\Hc A$ is obtained by evaluation of a 
$\bGamma$-space\index{subject}{Gamma-space@$\bGamma$-space} $\underline{A}$ on spheres,
where $\underline{A}(T)=A[T]$ is the reduced $A$-linearization of a finite based set $T$.
So $\Hc A$ is left induced from the global family $\Fin$
by Proposition \ref{prop:geometric fix of Gamma spaces}~(ii). 

The $\bGamma$-space $\underline{A}$ is actually a $\bGamma$-{\em set},
i.e., all its values have the discrete topology.
The latching map $l_n:\colim_{U\subsetneq \{1,\dots,n\}} \, A[U]\ \to \  A[\{1,\dots,n\}]$
defined in \eqref{eq:Gamma latching map} is thus an injective map 
(by Proposition \ref{prop:Gamma latching closed embedding}~(ii)) 
between discrete $\Sigma_n$-spaces, hence a $\Sigma_n$-cofibration.
It is then a $(G\times \Sigma_n)$-cofibration
if we let a compact Lie group $G$ act trivially.
Proposition \ref{prop:Gamma on S^V}~(ii) then shows that
the orthogonal $G$-spectrum $\underline{A}(\mS)=\Hc A$ (with trivial $G$-action)
is equivariantly connective. Since $G$ was arbitrary, this proves that $\Hc A$
is globally connective.

Dos Santos shows in \cite{dosSantos} that for every representation $V$ 
of a finite group $G$, the $G$-space $\Hc A(V)=A[S^V]$ 
is an equivariant Eilenberg-Mac\,Lane space of type $(A,V)$, i.e., the $G$-space
$\map_*(S^V, A[S^V])$ has homotopically discrete fixed points for 
all subgroups $H$ of $G$ and the natural map
\[A \ \to \ [S^V,A[S^V]]^H\ =\ \pi_0\left( \map_*^H(S^V,A[S^V]) \right) \]
sending $m\in A$ to the homotopy class of $m\cdot -:S^V\to A[S^V]$
is an isomorphism.
This shows that $\Hc A$ is a $\Fin$-$\Omega$-spectrum
for the constant global functor $\underline A$.

We offer an independent proof of the $\Fin$-$\Omega$-property 
via the $\bGamma$-$G$-space techniques of Segal and Shimakawa \cite{segal-some_equivariant, shimakawa},\index{subject}{Gamma-space@$\bGamma$-space!equivariant} 
in our adaptation of Theorem \ref{thm:prolonged delooping}.
We let $G$ act trivially on the $\bGamma$-space $\underline{A}$.
The $G$-maps 
\[ P_S \ : \ \underline{A}(S_+) \ \to \ \map(S,\underline{A}(1_+)) \ = \ A^S \ , \quad
P_S(x)(s)\ = \ \underline{A}(p_s)(x) \ ,\]
are homeomorphisms, so in particular $G$-weak equivalences.
So $\underline{A}$ is a special $\bGamma$-$G$-space in the sense 
of Shimakawa \cite[Def.\,1.3]{shimakawa}, 
see also Definition \ref{def:special G-Gamma} below. 
Since $\pi_0(\underline{A}(1_+))$ is isomorphic to $A$, 
and hence a group (as opposed to a monoid only),
the $\bGamma$-$G$-space $\underline{A}$ is very special.
Since $\underline{A}$ is also $G$-cofibrant,
 Theorem \ref{thm:prolonged delooping} shows that
for every finite group $G$ and every pair of
$G$-representations $V$ and $W$, the adjoint structure maps 
$\tilde\sigma_{V,W}: A[S^W]\to\map_*(S^V, A[S^{V\oplus W}])$
is a $G$-weak equivalence. 
\end{proof}

\Danger The properties mentioned in the previous proposition
do {\em not} generalize to compact Lie groups of positive dimension,
i.e., contrary to what one may expect at first, 
$(\Hc A)_G$ is not generally a $G$-$\Omega$-spectrum,
not all restriction maps in dimension~0 are isomorphisms 
(see Example \ref{eg:HZ for abelian G}),
and the groups $\pi_*^G ( \Hc A )$ need not be concentrated in dimension~0
(see Theorem \ref{thm:pi_1^T HA}).

\medskip

We now consider the important special case $A=\mZ$.
The equivariant homotopy group $\pi_0^G( \Hc\mZ )$ 
may be larger than\index{symbol}{$\Hc\mZ$ - {Eilenberg-Mac\,Lane spectrum of the integers}} 
a single copy of the integers, and we are now going to give a presentation
of $\pi_0^G( \Hc\mZ )$. Before we do so, we compare $\Hc\mZ$ to
the `infinite symmetric product of the sphere spectrum'.

\begin{eg}[Infinite symmetric product]\label{eg:infinite symmetric product}\index{subject}{infinite symmetric product spectrum}
There is no essential difference if we consider
the infinite symmetric product $S\! p^\infty$ (i.e., the reduced free abelian monoid)  
instead of the reduced free abelian group $\mZ[-]$ generated by representation spheres:
the levelwise inclusions of the free abelian monoids into the free abelian groups
provide a morphism of ultra-commutative ring spectra\index{symbol}{$\Spinf$ - {infinite symmetric product}}
\begin{equation}\label{eq:Sp^infty2HZ} 
 \Spinf \ = \ \{\Spinf(S^V) \}_V\ \to \
\{\mZ[S^V] \}_V\ = \  \Hc\mZ\ , 
\end{equation}
and this morphism is a global equivalence by the following proposition.

In his unpublished preprint \cite{segal-some_equivariant},
Segal argues that for every finite group $G$ 
and every $G$-representation $V$ with $V^G\ne 0$ the map
\[  \Spinf(S^V) \ \to \ \mZ[S^V]   \]
is a $G$-weak equivalence. A published proof of this fact
appears as Proposition\,A.6 in Dugger's paper \cite{dugger-real K}. 
This generalizes to compact Lie groups:
\end{eg}\index{subject}{Eilenberg-Mac\,Lane spectrum!of the integers}

\begin{prop}\label{prop:Sp2HZ}
Let $G$ be a compact Lie group.
\begin{enumerate}[\em (i)]
  \item  For every $G$-representation $V$ such that $V^G\ne 0$, 
    the natural map $\Spinf(S^V)\to \mZ[S^V]$ is a $G$-weak equivalence.
  \item The morphism \eqref{eq:Sp^infty2HZ} is a global equivalence 
    of orthogonal spectra from $S\! p^\infty$ to $\Hc\mZ$.
\end{enumerate}
\end{prop}
\begin{proof}
(i) By application to all closed subgroups of $G$, it suffices to show 
that the fixed point map $(S\! p^{\infty}(S^V))^G\to (\mZ[S^V])^G$ is a 
non-equivariant weak equivalence. 
We let $G^\circ\leq G$ denote the connected component of the identity element and
$\bar G=G/G^\circ$ the finite group of components of $G$.
To calculate $G$-fixed points we can first take $G^\circ$-fixed points and
then $\bar G$-fixed points. 
Proposition \ref{prop:Gamma fixed points} provides homeomorphisms
\[  (S\! p^\infty(S^V))^G \ \iso \  (S\! p^\infty(S^{V^{G^\circ}}))^{\bar G}
\text{\qquad and\qquad}
(\mZ[S^V])^G \iso \   (\mZ[S^{V^{G^\circ}}])^{\bar G} \ .\]
Since $V^{G^\circ}$ is an orthogonal representation of
the finite group $\bar G$ the map
\[  (S\! p^\infty(S^{V^{G^\circ}}))^{\bar G}\ \to \  (\mZ[S^{V^{G^\circ}}])^{\bar G}\ ,\]
is a weak equivalence by \cite[Prop.\,A.6]{dugger-real K}. 
Strictly speaking, Dugger's proposition is stated only for
geometric realizations of $\bar G$-simplicial sets; 
since the representation sphere $S^{V^{G^\circ}}$ 
admits the structure of a finite $\bar G$-CW-complex, 
it is $\bar G$-homotopy equivalent to the realization of a $\bar G$-simplicial set,
compare Proposition \ref{prop:equivariant simplicial realization}~(ii).
So we can also apply Dugger's result in our situation. 
Part~(ii) is then immediate from~(i).
\end{proof}

The infinite symmetric product of any based space $X$ has a natural filtration
\[ X = S\! p^1(X)\ \subseteq\ S\! p^2(X)\ \subseteq\ \ldots\ 
\subseteq\ S\! p^n(X)\ \subseteq\ \ldots \]
by the finite symmetric products. Evaluation on spheres
provides orthogonal spectra $S\! p^n  =  S\! p^n(\mS)$ and a filtration
\[ \mS=S\! p^1\ \subseteq\ S\! p^2\ \subseteq\ \ldots\ \subseteq \
S\! p^n\ \subseteq\ \ldots \]
of $S\! p^\infty$ by orthogonal subspectra. We study the equivariant
and global properties of this filtration in \cite{schwede-infinite symmetric}.
The 0-th homotopy group global functor of the orthogonal spectrum $S\! p^n$
has a very compact description as follows.
We recall that the Burnside ring $\mA(G)$ of a compact Lie group $G$ is
freely generated by the classes $t_H^G=\tr_H^G(p_H^*(1))$
where $H$ runs through a set of representatives of the conjugacy classes
of subgroups of $G$ with finite Weyl group.
Theorem~3.12 of \cite{schwede-infinite symmetric} shows that
the global functor $\upi_0(S\! p^n)$ is the quotient
of the Burnside ring global functor by the global subfunctor generated by
the element $n\cdot 1- t_{\Sigma_{n-1}}^{\Sigma_n}$ in $\mA(\Sigma_n)$,
\begin{equation}  \label{eq:pi_0_Sp^n}
 \upi_0(S\! p^n) \ \iso \ \mA/\td{n\cdot 1 \ - \ t_{\Sigma_{n-1}}^{\Sigma_n} } \ .  
\end{equation}
Letting $n$ go to infinity and combining this with part~(ii)
of Proposition \ref{prop:Sp2HZ} gives the following calculation of
the global functor $\upi_0(\Hc\mZ)$.
We let $I_\infty$ denote the global subfunctor of the Burnside ring
global functor $\mA$ generated by the classes 
$n\cdot 1- t_{\Sigma_{n-1}}^{\Sigma_n}$ for all $n\geq 1$.
We show in \cite[Thm.\,3.12]{schwede-infinite symmetric}
that the unit map $\mS\to S\! p^\infty$ induces an isomorphism of global functors
\[ \mA/I_\infty\ \iso \ \upi_0(S\! p^\infty)\ \iso \ \upi_0(\Hc\mZ)\ . \]
This calculation can be made even more explicit.
Elementary algebra (see \cite[Prop.\,4.1]{schwede-infinite symmetric})
identifies the value $I_\infty(G)$ at a compact Lie group $G$
as the subgroup of $\mA(G)$ generated by the classes
\[  [H:K]\cdot t_H^G \ - \ t_K^G \]
for all nested sequences of 
closed subgroup $K\leq H\leq G$ such that $W_G H$ is finite
and $K$ has finite index in $H$. The possibility that $K$ has infinite Weyl group
in $G$ is allowed, in which case $t_K^G=0$.
This presents the equivariant homotopy group
as an explicit quotient of the Burnside ring of $G$:
\[ \pi_0^G(\Hc\mZ)\ \iso \ \pi_0^G(S\! p^\infty)\ \iso \ \mA(G)/I_\infty(G) \]

\begin{eg}\label{eg:HZ for abelian G}
For every finite group $G$, the group $\pi_0^G(\Hc\mZ)$ is free of rank~1,
generated by the multiplicative unit, 
and the inflation map $p_G^*:\pi_0^e(\Hc\mZ)\to\pi_0^G(\Hc\mZ)$ is an isomorphism.
This does {\em not} persist to general compact Lie groups of positive dimension.
An explicit example for which $\pi_0^G( \Hc\mZ )$ has rank bigger than one 
is $G=S U(2)$.\index{subject}{special unitary group!$S U(2)$}
Then the classes $1$ and $\tr_N^{S U(2)}(1)$ are a $\mZ$-basis for
$\pi_0^{S U(2)}( \Hc\mZ )$ modulo torsion,
see \cite[Ex.\,4.16]{schwede-infinite symmetric},
where $N$ is a maximal torus normalizer. 
\end{eg}

To illustrate that the orthogonal spectrum $\Hc A$
is typically {\em not} an Eilenberg-Mac\,Lane spectrum of any global functor,
we calculate the $U(1)$-equivariant homotopy group $\pi_1^{U(1)}(\Hc A)$.

\begin{construction}
As we shall now see, all classes in the group $\pi_1^{U(1)}(\Hc A)$ 
arise in a systematic way via dimension shifting transfers 
from finite subgroups of $U(1)$.\index{subject}{unitary group!$U(1)$}
We organize these transfers into a natural isomorphism from $A\tensor\mQ$
to $\pi_1^{U(1)}(\Hc A)$.
 
We let $C$ be any finite subgroup of $U(1)$.
We identify $\mR$ with the tangent space 
of the distinguished coset $e C$ in $U(1)/C$ via the differential of the smooth curve
\[ \mR \ \to \ U(1)/C \ , \quad t\ \longmapsto \ e^{2\pi i t}\cdot C\ .\]
The upshot of this identification $\mR\iso T_{e C}( U(1)/C)$
is that we can view the dimension shifting transfer \eqref{eq:dimension_shifting_transfer}
as a homomorphism
\[ \Tr_C^{U(1)}\ : \ \pi_1^C(E\sm S^1) \ \to \ \pi_1^{U(1)}(E)\ .\]
We can then define an additive map $\psi_C:A\to\pi_1^{U(1)}(\Hc A)$ as the composite
\[ A \ \iso \ \pi_0^e(\Hc A)\ \xra[\iso]{\ p_C^*\ } \ 
\pi_0^C(\Hc A)\ \xra[\iso]{-\sm S^1} \
\pi_1^C(\Hc A\sm S^1)\ \xra{\ \Tr_C^{U(1)}\ } \ \pi_1^{U(1)}(\Hc A) \]
of inflation, the suspension isomorphism and the dimension shifting transfer.

Now we consider a subgroup $C'$ of $C$ of index $m$. 
Since the restriction of $\upi_0(\Hc A)$ to finite groups is a 
constant global functor, the relation
\begin{align*}
  \Tr_{C'}^{U(1)}\circ (-\sm S^1)&\circ p_{C'}^* \ = \ 
  \Tr_C^{U(1)}\circ\tr_{C'}^C\circ (-\sm S^1)\circ p_{C'}^* \\ 
&= \   \Tr_C^{U(1)}\circ (-\sm S^1)\circ\tr_{C'}^C\circ p_{C'}^* \ 
= \   m\cdot \Tr_C^{U(1)}\circ (-\sm S^1)\circ p_C^* 
\end{align*}
holds as homomorphisms $\pi_0^e(\Hc A)\to\pi_1^{U(1)}(\Hc A)$,
by transitivity of transfers (Proposition \ref{prop:transfer transitive}).
Hence $\psi_{C'}=m\cdot\psi_C$ and we obtain a well-defined homomorphism
\[ \psi \ : \ A\tensor\mQ \ \to \ \pi_1^{U(1)}(\Hc A)\text{\qquad by\qquad}
\psi(a\tensor r/s ) \ = \ r\cdot \psi_{C_s}(a)\ ,\]
where $C_s\subset U(1)$ is group of $s$-th roots of unity.
\end{construction}

\begin{theorem}\label{thm:pi_1^T HA} 
The homomorphism $\psi : A\tensor\mQ \to\pi_1^{U(1)}(\Hc A)$ is an isomorphism
for every abelian group $A$.
\end{theorem}
\begin{proof}
Since the orthogonal spectrum $\Hc A$ is obtained from a $\bGamma$-space
by evaluation on spheres, the inflation homomorphism
\[ p^*\ : \ \pi_*^e(\Hc A)\ = \ \Phi_*^e(\Hc A)\ \to \ \Phi_*^{U(1)}(\Hc A) \]
is an isomorphism of geometric fixed point homotopy groups,
by Proposition \ref{prop:geometric fix of Gamma spaces}~(ii).
The non-equivariant homotopy groups of $\Hc A$ are concentrated in dimension~0,
so the $U(1)$-equivariant geometric fixed point homotopy groups $\Phi_k^{U(1)}(\Hc A)$
vanish for all $k\ne 0$.
The isotropy separation sequence \eqref{eq:isotropy_separation_les}
thus shows that the map $E\Pc\to \ast$ induces an isomorphism 
\[ \pi_1^{U(1)}(\Hc A\sm E\Pc_+)\ \iso \  \pi_1^{U(1)}(\Hc A) \ , \]
where $E\Pc$ is a universal $U(1)$-space for the family of proper closed (i.e., finite)
subgroups of $U(1)$.

We let $X$ be mapping telescope of the sequence of projections
\[ U(1)\ \to \ U(1)/C_{2} \ \to \ U(1)/C_{6} \ \to \ \dots
\ \to \ U(1)/C_{n !} \ \to \ \dots \ .\]
This is a $U(1)$-CW-complex without $U(1)$-fixed points and
such that the fixed point space $X^C$ is path connected for every finite subgroup $C$
of $U(1)$. So we can build a universal space $E\Pc$ by attaching
equivariant cells with finite isotropy groups of dimension~2 and larger to $X$.
The Wirthm{\"u}ller isomorphisms
\[ \pi_k^{U(1)}( \Hc A\sm (U(1)/C)_+\sm S^n) \ \iso \ 
\pi_k^C(\Hc A\sm S^1\sm S^n) \ \iso \ \pi_{k-1-n}^C(\Hc A)  \]
and the fact that the $C$-equivariant homotopy groups of $\Hc A$
are concentrated in dimension~0 mean that the group
$\pi_k^{U(1)}( \Hc A\sm E\Pc/X)$ is trivial for all $k\leq 2$. 
Hence the inclusion $X\to E\Pc$ induces an isomorphism
\[ \pi_1^{U(1)}(\Hc A\sm  X_+) \ \xra{\ \iso \ }\ 
 \pi_1^{U(1)}(\Hc A\sm  E\Pc_+) \ . \]

We let $f^n:U(1)/C_{n !}\to X$ be the inclusion.
The dimension shifting transfer map $\Tr_C^{U(1)}$ factors as the composite
\begin{align*}
  \pi_1^C(\Hc A\sm S^1) \ \xra{U(1)\ltimes_C -} \ &\pi_1^{U(1)}(\Hc A\sm (U(1)/C_{n !})_+)\\ 
  \xra{\quad f^n_*\quad } \ \ &\pi_1^{U(1)}(\Hc A\sm X_+) \ \to \ \pi_1^{U(1)}(\Hc A) \ .
\end{align*}
So the map $\psi$ factors through $\pi_1^{U(1)}(\Hc A\sm X_+)$.
Since $\Hc A\sm X_+$ is a mapping telescope, its equivariant homotopy
groups can be calculated as the colimit of the sequence:
\[ \pi_1^{U(1)}(\Hc A\sm U(1)_+)
\ \to \ \dots \ \to \ \pi_1^{U(1)}(\Hc A\sm (U(1)/C_{n !})_+) \ \to \ \dots \]
We rewrite this using the Wirthm{\"u}ller and suspension isomorphisms as
\[ \pi_1^{U(1)}(\Hc A\sm (U(1)/C_k)_+) \ \iso\ \pi_1^{C_k}(\Hc A\sm S^1)
\ \iso \ \pi_0^{C_k}(\Hc A)\ \iso \ A \ . \]
Example \ref{eq:projection versus transfer} 
shows that for all $k,m\geq 1$ the following diagram commutes:
\[ \xymatrix@C=20mm{ 
\pi_1^{C_k}(\Hc A\sm S^1)\ar[r]^-{\tr_{C_k}^{C_{k m}}}\ar[d]_{U(1)\ltimes_{C_k} -}^{\iso} &
\pi_1^{C_{k m}}(\Hc A\sm S^1)\ar[d]^{U(1)\ltimes_{C_{k m}} -}_{\iso}\\
\pi_1^{U(1)}(\Hc A\sm (U(1)/C_k)_+)\ar[r]_-{ (\Hc A\sm \text{proj})_*} &
\pi_1^{U(1)}(\Hc A\sm (U(1)/C_{k m})_+)} \]
Transfers commute with the suspension isomorphism and 
in the global functor $\upi_0(\Hc A)$, 
the transfer $\tr_{C_k}^{C_{k m}}$ is multiplication by
the index $[C_{k m}:C_k]=m$;
so the above homotopy group sequence becomes the sequence
\[
 A \ \xra{\ \cdot 2\ }\ A\ \xra{\ \cdot 3\ }\ \cdots \ \to \ A 
\xra{\ \cdot n\ }\ A \ \to \ \cdots \ \ . \]
The colimit of this sequence is isomorphic to $A\tensor\mQ$,
compatible with the map $\psi:A\tensor\mQ\to \pi_1^{U(1)}(\Hc A\sm X_+)$. 
This proves the claim.
\end{proof}
\index{subject}{Eilenberg-Mac\,Lane spectrum!of an abelian group|)}

\begin{eg}\index{subject}{global Borel theory}  
The right adjoint $R:\SH\to \GH$ to the forgetful functor
from the global to the non-equivariant stable homotopy category
is modeled on the pointset level by the functor
$b:\spec\to\spec$ discussed in Construction \ref{cor:Borel b}.
The global homotopy type of $b E$ is
that of a Borel cohomology theory, and in particular,
\[ \pi_0^G( b E) \ \iso \ E^0(B G)\ ,\]
natural in $G$ for transfers and restriction maps.
The functor $b$ is lax symmetric monoidal, so it takes
an ultra-commutative ring spectrum $R$ to an ultra-commutative
ring spectrum $b R$; the power operations
\[ P^m\ : \ \pi_0^e(b E) \ \to \ \pi_0^{\Sigma_m}(b E) \]
correspond to the classical power operations
\[ \Pc_m\ : \ E^0(\ast) \ \to \ E^0(B\Sigma_m) \ , \]
compare the more general Remark \ref{rk:classical power}.
\end{eg}

\begin{eg}[Representation ring global power functor]\label{eg:RU_as_global_power}
As $G$ varies over all compact Lie groups, the unitary representation rings $\bRU(G)$
form the {\em unitary representation ring global functor}\index{subject}{representation ring!unitary} $\bRU$, compare Example \ref{eg:global functor examples}~(iv).
This is classical in the restricted realm of finite groups, but somewhat less familiar
for compact Lie groups in general.
The restriction maps $\alpha^*:\bRU(G)\to \bRU(K)$
are induced by restriction of representations along a homomorphism $\alpha:K\to G$.
The transfer maps $\tr_H^G:\bRU(H)\to \bRU(G)$ along a closed subgroup inclusion
$H\leq G$ are given by the {\em smooth induction} 
of Segal \cite[\S\,2]{segal-representation}.
If $H$ is a subgroup of finite index of $G$, then this induction 
sends the class of an $H$-representation $V$ to the induced $G$-representation
$\map^H(G,V)$;
in general, induction sends actual representations to virtual representations.
In the generality of compact Lie groups,
the double coset formula for $\bRU$ was proven by Snaith \cite[Thm.\,2.4]{Snaith-Brauer}.

The representation rings also have well-known power operations
\[ P^m \ : \ \bRU(G)\ \to \ \bRU(\Sigma_m\wr G) \ ;\]
on the class of a $G$-representation $V$, the power operation is
represented by the tensor power
\[ P^m[V] \ = \ [V^{\tensor m}]\ ,\]
using the canonical action of $\Sigma_m\wr G$ on $V^{\tensor m}$.
Since power operations are not additive, one has to argue why
this assignment extends to virtual representations.
The standard way is to assemble all power operations 
on representations into a map
\[ P \ : \ \bRU^+(G)\ \to \ \exp(\bRU;G) \ , \quad
P( [V] )\ = \ \left( [V^{\tensor m}]\right)_{m\geq 0} \]
from the monoid of isomorphism classes of $G$-representations.
If $W$ is another $G$-representation, then
\[ (V\oplus W)^{\tensor m} \text{\qquad and\qquad}
{\bigoplus}_{k=0}^m \, \tr_{ (\Sigma_k\wr G)\times(\Sigma_{m-k}\wr G)}^{\Sigma_m\wr G} 
(V^{\tensor k}\tensor W^{\tensor (m-k)})\]
are isomorphic as $(\Sigma_m\wr G)$-representations,
because tensor product distributes over direct sum.
This means that the total power map $P$ is a monoid homomorphism
from $\bRU^+(G)$ to the group $\exp(\bRU;G)$ under $\oplus$.
So the total power operation extends uniquely to a group homomorphism
\[ P \ : \ \bRU(G)\ \to \ \exp(\bRU;G) \]
on the representation ring. The operation $P^m$ is then the $m$-th factor
of this extension.
The representation ring global power functor $\bRU$ is
realized by the ultra-commutative ring spectrum $\bKU$,
the unitary global $K$-theory spectrum,
see Theorem \ref{thm:pi_0 KU is RU} below.\index{subject}{K-theory@$K$-theory!periodic global}\end{eg}

\begin{rk}[Brauer induction]\label{rk:Brauer induction}
By Brauer's theorem \cite[Thm.\,I] {brauer} 
the complex representation ring of a finite group
is generated, as an abelian group, by representations that are induced from
1-dimensional representations of subgroups.
Segal generalized this result 
to compact Lie groups in \cite[Prop.\ 3.11 (ii)]{segal-representation}, 
where `induction' refers to the smooth induction. 
In fact, in the world of compact Lie group, Segal's smooth induction
for not necessarily finite index subgroups
makes the proof quite transparent, as we shall now recall.
In our language the statement can be expressed by saying that
the representation ring global functor $\bRU$ is `cyclic' in the sense that it is
generated by a single element, the class
$x\in \bRU(U(1))$ of the tautological 1-dimensional representation 
of the circle group $U(1)$.\index{subject}{unitary group!$U(1)$}
Equivalently, the morphism of global functors
\[ \ev_x \ : \ \bA(U(1),-)\ \to \ \bRU \]
classified by the element $x$ is an epimorphism. 
We recall the argument: we let $i:U(1) \times U(n-1)\to U(n)$ be the block sum embedding
and $q:U(1)\times U(n-1)\to U(1)$ the projection to the first factor.
The character formula \cite[p.\,119]{segal-representation}
for induced representations shows that the smooth transfer
\begin{equation}  \label{eq:tautological_as_induction}
 i_!(q^*(x))\ \in \ \bRU(U(n))   
\end{equation}
has the same character as the tautological $n$-dimensional representation
of $U(n)$. Since characters determine unitary representations of compact Lie groups,
$i_!(q^*(x))$ equals the class of the tautological representation $\tau_n$ of $U(n)$. 
Any unitary representation of a compact Lie group $G$ of dimension $n$
is isomorphic to $\alpha^*(\tau_n)$ for a continuous homomorphism
$\alpha:G\to U(n)$; so the class of such a representation equals
\[ \alpha^*(i_!(q^*(x)))\ \in \ \bRU(G) \ .\]
So the global functor $\bRU$ is generated by the single class $x=\tau_1$.

An interesting line of investigation, dubbed {\em explicit Brauer induction},\index{subject}{Brauer induction!explicit}
started with Snaith's paper \cite{Snaith-Brauer}.
Informally speaking, an `explicit Brauer induction' is a section to the map
\[ \bA(U(1),G) \ \to \ \bRU(G) \]
that is specified by a direct recipe, for example an explicit formula, 
and has naturality properties as the group $G$ varies.
So such a map is an `explicit and natural' way to write a (virtual) representation
as a sum of induced representations of 1-dimensional representations.
The first explicit Brauer induction 
was Snaith's formula \cite[Thm.\,(2.16)]{Snaith-Brauer};
however, Snaith's maps are not additive and not compatible with restriction to subgroups.
Later Boltje \cite{boltje-canonical brauer}
specified a different explicit Brauer induction
formula by purely algebraic means;
Symonds \cite{symonds-splitting} gave a topological interpretation of Boltje's
construction. The Boltje-Symonds maps are additive and natural for
restriction along group homomorphisms; the maps are not (and in fact cannot be)
in general compatible with transfers. 
In our present language, the Boltje-Symonds maps form a natural transformation of
$\Rep$-abelian groups from $\bRU$ to $\bA(U(1),-)$.
We recall the construction of these maps; 
we follow Symonds' approach, as his reasoning is very much 
in the spirit of global homotopy theory.

The starting point is the formula \eqref{eq:tautological_as_induction}
that expresses the class of the tautological $U(n)$-representation
in $\bRU(U(n))$ as a smooth induction of a specific 1-dimensional
representation of the subgroup $U(1)\times U(n-1)$, namely the one whose character
is the projection $q:U(1)\times U(n-1)\to U(1)$
to the first factor. 
Formula \eqref{eq:tautological_as_induction} suggest a class in $\bA(U(1),U(n))$ 
as the image of the tautological $U(n)$-representation; 
if we also want naturality and additivity, then this fixes things completely:
\end{rk}

\begin{theorem} [Boltje \cite{boltje-canonical brauer}, Symonds \cite{symonds-splitting}]
There is a unique natural transformation of
$\Rep$-abelian groups
\[ b \ : \ \bRU \ \to \ \bA(U(1),-) \]
that satisfies
\[ b_{U(n)}(\tau_n)\ = \ \tr_{U(1)\times U(n-1)}^{U(n)}\circ q^*\text{\qquad in\quad} \bA(U(1),U(n))  \ . \]
Moreover:
\begin{enumerate}[\em (i)]
\item 
The transformation $b$ is a section to the evaluation morphism
of global power functors $\ev_x:\bA(U(1),-)\to\bRU$ at the class $x\in\bRU(U(1))$.
\item 
The value of $b_G$ at the 1-dimensional representation with
character $\chi:G\to U(1)$ is given by
\[ b_G[\chi^*(\tau_1)]\ = \ \chi^*  \ \in \ \bA(U(1),G)\ .\]
\end{enumerate}
\end{theorem}
\begin{proof}
Every class in $\bRU(G)$ is a formal difference of classes of
actual representations, and every $n$-dimensional representation
is the restriction of $\tau_n$ along some continuous homomorphism
$G\to U(n)$. So uniqueness is a consequence of naturality and additivity.  

Conversely, this also suggests how to define the transformation.
If $V$ is any $n$-dimensional unitary $G$-representation,
then $V$ is isomorphic to $\alpha^*(\mC^n)$ 
for some continuous homomorphism $\alpha:G\to U(n)$,
unique up to conjugacy. So we set
\[ b_G[V] \ = \  b_G(\alpha^*(\tau_n)) \ = \ 
\alpha^*\circ\tr_{U(1)\times U(n-1)}^{U(n)}\circ q^*\  \in \ \bA(U(1),G)  \ . \]
This defines a set theoretic map
\[ b_G \ : \ \bRU^+(G) \ \to \ \bA(U(1),G) \]
from the monoid of isomorphism classes of unitary $G$-representations,
and these maps are automatically compatible with restriction along group homomorphisms.
The double coset formula \eqref{eq:U(n)_Brauer_double_coset}
for $\res^{U(n+m)}_{U(n)\times U(m)}\circ\tr_{U(1)\times U(n+m-1)}^{U(n+m)}$
implies the relation
\[ \res^{U(n+m)}_{U(n)\times U(m)} ( b_{U(n+m)}(\tau_{n+m}) ) \ = \  
b_{U(n)}(\tau_n) \ +\ b_{U(m)}(\tau_m)\]
in the group $\bA(U(1),U(n)\times U(m))$.
Hence the maps $b_G$ are additive, and so they extend uniquely to
a group homomorphism
\[ b_G \ : \ \bRU(G) \ \to \ \bA(U(1),G) \]
on the Grothendieck group, for which we use the same name. 
These homomorphisms are still compatible
with restriction along continuous group homomorphisms.

It remains to show the additional properties.
The transformation $\ev_x\circ b:\bRU\to\bRU$ is additive
and natural for restriction along continuous homomorphisms,
so for property~(i) it suffices to show the relation $\ev_x\circ b=\Id$
in the universal examples, i.e., for the tautological representations
of the unitary groups $U(n)$. This universal example is taken care of
by the formula \eqref{eq:tautological_as_induction}.

Property~(ii) holds because for $n=1$ the map $q$ is the identity of $U(1)$.
\end{proof}

While the above construction of the explicit Brauer map $b:\bRU\to\bA(U(1),-)$ 
is slick, it is not yet particularly explicit. To write the class $b_G[V]$
as a $\mZ$-linear combination of transfers of 1-dimensional representations
of subgroups of $G$, one would now have to write 
the classifying homomorphism $\alpha:G\to U(n)$ for $V$
as the composite of an epimorphism and a subgroup inclusion and then 
expand the term $\alpha^*\circ\tr^{U(n)}_{U(1)\times U(n-1)}$
using the compatibility of transfers with inflation, and the
double coset formula for the restriction of a transfer.
Readers desperate for a truly explicit formula can find one in
\cite[Thm.\,(2.1)]{boltje-canonical brauer}
or \cite[Thm.\,2.24 (e)]{boltje-snaith-symonds}.

\index{subject}{global power functor|)}

\section{Global model structure}
\label{sec:ucom global model}

In this section we construct a model structure on the category of
ultra-commu\-ta\-tive ring spectra with global equivalences as the weak equivalences,
see Theorem \ref{thm:global ultra-commutative}.
The strategy is the same is in the unstable situation 
in Section \ref{sec:global model monoid spaces}:
we lift the positive version of the global model structure 
to commutative monoid objects with the help of the
general lifting theorem \cite[Thm.\,3.2]{white-commutative monoids}.

We also calculate the algebra
of natural operations on the homotopy groups of ultra-commutative ring spectra:
we show that these operations are freely generated by restrictions,
transfers and power operations,
see Theorem \ref{thm:generators global power functor} for the precise statement.
Then we show in Theorem \ref{thm:realize global power functor} that every global
power functor can be realized by an ultra-commutative ring spectrum;
the realization can moreover be chosen as a global Eilenberg-Mac Lane spectrum,
i.e., such that all equivariant homotopy groups in non-zero dimensions vanish.

\medskip

The construction of the global model structure
on the category of ultra-commutative ring spectra is largely parallel
to the unstable precursor, the global model structure for ultra-commutative monoids 
in Theorem \ref{thm:global umon}.
The analogous arguments as in the unstable case 
in Corollary \ref{cor:co-limits in umon} show 
that the category of ultra-commutative ring spectra is complete and cocomplete
and that the forgetful functor 
to the category of orthogonal spectra creates all limits, 
all filtered colimits and those coequalizers that are reflexive in
the category of orthogonal spectra.
The category of ultra-commutative ring spectra is
tensored and cotensored over the category $\bT$ of spaces,
by the same arguments in Construction \ref{con:umon enrich}
for ultra-commutative monoids.

We recall from \cite[Def.\,3]{gor-gul-symmetric} 
the notion of a {\em symmetrizable cofibration}\index{subject}{symmetrizable cofibration}\index{subject}{symmetrizable acyclic cofibration}
respectively {\em symmetrizable acyclic cofibration},
see also Definition \ref{def:symmetrizable}.
For a morphism $i:A\to B$ of orthogonal spectra we let
$K^n(i)$ denote the $n$-cube of orthogonal spectra whose
value at a subset $S\subseteq \{1,2,\dots,n\}$ is
\[ K^n(i)(S) \ = C_1 \sm C_2 \sm \dots \sm C_n \text{\qquad with\qquad}
 C_j \ = \left\{ \begin{array}{l@{\quad}l} A & \mbox{if } j \not\in  
S \\
B & \mbox{if } j \in S. \end{array} \right. \]
All morphisms in the cube $K^n(i)$ are smash products of identities and
copies of the morphism $i:A\to B$.
We let $Q^n(i)$ denote the colimit of the punctured $n$-cube,
i.e., the cube $K^n(i)$ with the terminal vertex removed, and
$i^{\Box n}:Q^n(i)\to K^n(i)(\{1,\dots,n\})=B^{\sm n}$ the canonical map.
So the morphism $i^{\Box n}$ is an iterated  pushout product morphism.
The symmetric group $\Sigma_n$ acts on $Q^n(i)$ and $B^{\sm n}$ 
by permuting the smash factors,
and the iterated pushout product morphism
$i^{\Box n}:Q^n(i)\to B^{\sm n}$ is $\Sigma_n$-equivariant.
The morphism $i:A\to B$ of orthogonal spectra 
is a {\em symmetrizable cofibration} \index{subject}{symmetrizable cofibration}\index{subject}{symmetrizable acyclic cofibration}
(respectively a {\em symmetrizable acyclic cofibration}) if the morphism
\[ i^{\Box n}/\Sigma_n \ :\ Q^n(i)/\Sigma_n \ \to \ B^{\sm n}/\Sigma_n = \mP^n(B) \]
is a cofibration (respectively an acyclic cofibration) for every $n\geq 1$.  
Since the morphism $i^{\Box 1}/\Sigma_1$ is the original morphism $i$, 
every symmetrizable cofibration
is in particular a cofibration, and similarly for acyclic cofibrations.

The next theorem says that in the category of orthogonal spectra,
all cofibrations and acyclic cofibrations in the positive global model structure
are symmetrizable with respect to the monoidal structure given by the smash product.

\begin{theorem}\label{thm:symmetrizable in orthogonal spectra}
  \begin{enumerate}[\em (i)]
  \item 
    Let $i:A\to B$ be a flat cofibration of orthogonal spectra.
    Then for every $n\geq 1$ the morphism 
    \[ i^{\Box n}/\Sigma_n \ : \ Q^n(i)/\Sigma_n \ \to \ B^{\sm n}/\Sigma_n \]
    is a flat cofibration. 
    In other words, all cofibrations in the global model structure
    of orthogonal spectra are symmetrizable.
  \item
    Let $i:A\to B$ be a positive flat cofibration of orthogonal spectra 
    that is also a global equivalence.
    Then for every $n\geq 1$ the morphism 
    \[ i^{\Box n}/\Sigma_n \ : \ Q^n(i)/\Sigma_n \ \to \ B^{\sm n}/\Sigma_n \]
    is a global equivalence.
    In other words, all acyclic cofibrations in the positive global model structure
    of orthogonal spectra are symmetrizable.
  \end{enumerate}
\end{theorem}
\begin{proof}
(i) We recall from Theorem \ref{thm:F-global spectra}~(iii) 
the set 
\[ I_{\All} \ = \ \{\ G_m( ( O(m)/H\times i_k )_+) \ | \ m,k \geq 0, H\leq O(m)\} \]
of generating flat cofibrations of orthogonal spectra, where $i_k:\partial D^k\to D^k$
is the inclusion.
The set $I_{\All}$ detects the acyclic fibrations 
in the strong level model structure of orthogonal spectra.
In particular, every flat cofibration is a retract of an $I_{\All}$-cell complex.
By \cite[Cor.\,9]{gor-gul-symmetric} or \cite[Lemma A.1]{white-commutative monoids} 
it suffices to show that the
generating flat cofibrations in $I_{\All}$ are symmetrizable.

For a space $A$, the orthogonal spectrum $G_m( (O(m)/H\times A)_+ )$ 
is isomorphic to $F_{H,\mR^m}\sm A_+$;
so we show more generally that every morphism of the form
\[ F_{G,V} \sm (i_k)_+ \ : \  F_{G,V}\sm \partial D^k_+ \ \to \ F_{G,V}\sm D^k_+ \]
is a symmetrizable cofibration, where $V$ is any representation 
of a compact Lie group $G$.
The symmetrized iterated pushout product
\begin{equation}  \label{eq:sym_iterated_ppp spectra}
 ( F_{G,V}\sm (i_k)_+ )^{\Box n} / \Sigma_n \ : \ 
Q^n(F_{G,V}\sm (i_k)_+ ) /\Sigma_n\ \to \  (F_{G,V} \sm D^k_+)^{\sm n}/\Sigma_n 
\end{equation}
is isomorphic to 
\[ F_{\Sigma_n\wr G,V^n}((i_k^{\Box n})_+) \ : \ 
 F_{\Sigma_n\wr G,V^n} ( Q^n (i_k)_+ ) \ \to \  F_{\Sigma_n\wr G,V^n}( (D^k)^n_+) \ ,\]
where
\[ i_k^{\Box n}\ : \ Q^n (i_k) \ \to \ (D^k)^n \]
is the $n$-fold pushout product of the inclusion $i_k:\partial D^k\to D^k$, 
with respect to the cartesian product of spaces.
The map $i_k^{\Box n}$ is $\Sigma_n$-equivariant, 
and we showed in the proof of the analogous unstable result in
Theorem \ref{thm:symmetrizable in orthogonal spaces}~(i)
that $i_k^{\Box n}$ is a cofibration of $\Sigma_n$-spaces.
Proposition \ref{prop:cofibrancy preservers}~(i)
then shows that $i_k^{\Box n}$ is also a cofibration of $(\Sigma_n\wr G)$-spaces,
with respect to the action by restriction along the projection
$(\Sigma_n\wr G)\to\Sigma_n$. 
So the morphism \eqref{eq:sym_iterated_ppp spectra} is a flat cofibration.

(ii) Theorem \ref{thm:F-global spectra}~(iii) describes a 
set $J_{\All}\cup K_{\All}$ of generating acyclic cofibrations for the global model structure
on the category of orthogonal spectra.
From this we obtain a set $J^+\cup K^+$ 
of generating acyclic cofibration for the {\em positive} 
global model structure of Proposition \ref{prop:positive global spectra} 
by restricting to those morphisms in $J_{\All}\cup K_{\All}$ 
that are positive cofibrations, 
i.e., homeomorphisms in level~0; so explicitly, we set
\[ J^+ \ = \ \{\ G_m( ( O(m)/H\times j_k)_+ )\ | \ m\geq 1, k \geq 0, H\leq O(m)\} \ ,\]
where $j_k:D^k\times\{0\}\to D^k\times [0,1]$ is the inclusion, and
\[ K^+\ = \ \bigcup_{G,V,W\ : \  V\ne 0} \Zc(\lambda_{G,V,W}) \ ,\]
the set of all pushout products of sphere inclusions $i_k$
with the mapping cylinder inclusions of the global equivalences
$\lambda_{G,V,W}: F_{G,V\oplus W}S^V \to F_{G,W}$.
Here $(G,V,W)$ runs through a set of representatives
of the isomorphism classes of triples consisting of a compact Lie group $G$,
a $G$-representation $V$ and non-zero faithful $G$-representation $W$.
By \cite[Cor.\,9]{gor-gul-symmetric} or \cite[Lemma A.1]{white-commutative monoids} 
it suffices to show that all morphisms in $J^+\cup K^+$
are symmetrizable acyclic cofibrations. 

We start with a morphism $G_m( ( O(m)/H\times j_k)_+ )$ in $J^+$.
For every $n\geq 1$, the morphism
\[ (G_m( ( O(m)/H\times j_k)_+ ))^{\Box n}/\Sigma_n \]
is a flat cofibration by part~(i), and a homeomorphism in level~0
because $m\geq 1$. Moreover, the morphism $j_k$ is a homotopy
equivalence of spaces, so the morphism $G_m( ( O(m)/H\times j_k)_+ )$ 
is a homotopy equivalence of orthogonal spectra; the morphism
$\mP^n(G_m( ( O(m)/H\times j_k)_+ ))$ is then again a homotopy equivalence 
for every $n\geq 1$,
by Proposition \ref{prop: Zc is symmetrizable}~(i).
Then \cite[Cor.\,23]{gor-gul-symmetric} shows that $G_m( ( O(m)/H\times j_k)_+ )$
is a symmetrizable acyclic cofibration. 
This takes care of the set $J^+$.

Now we consider the morphisms in the set $K^+$.
Since $G$ acts faithfully on the non-zero inner product space $W$, the action
of the wreath product $\Sigma_n\wr G$ on $W^n$ is again faithful.
So the morphism
\[ \lambda_{\Sigma_n\wr G,V^n, W^n} \ : \ F_{\Sigma_n\wr G, V^n\oplus W^n}\, S^{V^n} \ \to \
F_{\Sigma_n\wr G, W^n} \]
is a global equivalence by Theorem \ref{thm:faithful independence}.
The vertical morphisms in the commutative square
\[ \xymatrix@C=17mm{ 
F_{\Sigma_n\wr G, V^n\oplus W^n} \, S^{V^n}
\ar[r]^-{ \lambda_{\Sigma_n\wr G,V^n, W^n} } \ar[d]_\iso &
F_{\Sigma_n\wr G, W^n} \ar[d]^\iso\\
\mP^n(F_{G,V\oplus W} S^V) \ar[r]_-{\mP^n(\lambda_{G,V,W}) } & \mP^n(F_{G,W})
} \]
are isomorphisms; so the morphism  $\mP^n(\lambda_{G,V,W})$ is a global equivalence. 
Proposition \ref{prop: Zc is symmetrizable}~(iii)
then shows that all morphisms in $\Zc(\lambda_{G,V,W})$ 
are symmetrizable acyclic cofibrations. 
\end{proof}

Now we put all the pieces together.
We call a morphism of ultra-commutative ring spectra a
{\em global equivalence} (respectively {\em positive global fibration})
if the underlying morphism of orthogonal spectra is a
global equivalence (respectively fibration in the positive global model structure
of Proposition \ref{prop:positive global spectra}).

\begin{theorem}[Global model structure for ultra-commutative ring spectra]\label{thm:global ultra-commutative} \quad
  \begin{enumerate}[\em (i)]
  \item 
    The global equivalences and positive global fibrations are part of a 
    cofibrantly generated, proper, topological model structure
    on the category of ultra-commutative ring spectra,
    the {\em global model structure}.\index{subject}{global model structure!for ultra-commutative ring spectra}
  \item Let $j:R\to S$ be a cofibration in the global model structure 
  of ultra-commutative ring spectra.
  \begin{enumerate}[\em (a)]
  \item The morphism of $R$-modules underlying $j$ is a cofibration
    in the global model structure of $R$-modules 
    of Corollary {\em \ref{cor-lift modules}~(i)}. 
  \item The morphism of orthogonal spectra underlying $j$ is an h-cofibration.
  \item If the underlying orthogonal spectrum of $R$ is flat, then 
    $j$ is a flat cofibration of orthogonal spectra.
\end{enumerate} 
  \end{enumerate}
\end{theorem}
\begin{proof}
The proof is completely parallel to the proof of the analogous unstable
theorem, Theorem \ref{thm:global umon}.

(i) The positive global model structure of orthogonal spectra
(Proposition \ref{prop:positive global spectra}) is monoidal 
and cofibrantly generated.
The `unit axiom' also holds: we let $f:\mS^{\text{c}}\to\mS$ 
be any positive flat replacement
of the orthogonal sphere spectrum.
Then for every orthogonal spectrum $X$ the induced morphism
$f\sm X:\mS^{\text{c}}\sm X\to \mS\sm X$ is a global equivalence by
Theorem \ref{thm:flat is flat}~(ii).
The monoid axiom holds by Proposition \ref{prop:monoid axiom}.
Cofibrations and acyclic cofibrations are symmetrizable by
Theorem \ref{thm:symmetrizable in orthogonal spectra},
so the model structure satisfies the `commutative monoid axiom'
of \cite[Def.\,3.1]{white-commutative monoids}.
The symmetric algebra functor $\mP$ commutes with filtered colimits
by the analog of Corollary \ref{cor:co-limits in umon} 
for ultra-commutative ring spectra.
Theorem~3.2 of \cite{white-commutative monoids} thus shows that the positive global
model structure of orthogonal spectra lifts to 
a cofibrantly generated model structure on the category of 
ultra-commutative ring spectra.
The global model structure is topological 
by Proposition \ref{prop:topological criterion};
here we take $\Gc$ as the set of free ultra-commutative 
ring spectra $\mP(F_{H,\mR^m})$ for all $m\geq 1$ 
and all closed subgroups $H$ of $O(m)$,
and we take $\Zc$ as the set of acyclic cofibrations
$\mP(c(\lambda_{G,V,W}))$ for the mapping cone inclusions $c(\lambda_{G,V,W})$ of
the global equivalences $\lambda_{G,V,W}:F_{G,V\oplus W}S^W\to F_{G,W}$,
indexed by representatives as in the definition of the set $K^+$.
Since weak equivalences, fibrations and pullbacks of ultra-commutative ring spectra are
are created on underlying orthogonal spectra,
right properness is inherited from the positive global model structure
of orthogonal spectra (Proposition \ref{prop:positive global spectra}).
We defer the proof of left properness until after the proof of part~(ii).

Part~(ii) is proved in exactly the same way as in the unstable case
in Theorem \ref{thm:global umon};
whenever the unstable proof refers to the model category of $R$-modules
in Corollary \ref{cor-lift to modules spaces}, the stable 
proof instead uses Corollary \ref{cor-lift modules}.
For the symmetrizability of the cofibrations,
the stable proof uses Theorem \ref{thm:symmetrizable in orthogonal spectra}~(i)
instead of Theorem \ref{thm:symmetrizable in orthogonal spaces}~(i).
We refrain from repeating the remaining details.

It remains to show that the model structure is left proper.
A pushout square of ultra-commutative ring spectra has the form
\[ \xymatrix@C=12mm{
R \ar[d]_j \ar[r]^-f_-\simeq & T\ar[d]^{j\sm_R T}\\ 
S\ar[r]_-{S\sm_R f} & S\sm_R T} \]
where $S$ and $T$ are
considered as $R$-modules by restriction along $j$ respectively $f$.
For left properness we now suppose that $j$ is a cofibration and $f$ 
is a global equivalence.
By part~(a) of~(ii), the morphism $j$ is then a cofibration of $R$-modules 
in the global model structure of Corollary \ref{cor-lift modules}~(i). 
Since $R$ is cofibrant in that model structure, also $S$ is
cofibrant as an $R$-module.
Proposition \ref{prop:sm cofibrant R-mod} then shows that the functor
$S\sm_R -$ preserves global equivalences.
So the cobase change $S\sm_R f$ of $f$ is a global equivalence.
This shows that the global model structure of ultra-commutative ring spectra 
is left proper. 
\end{proof}

Now we can show that the restriction maps, (additive) transfer maps and
(multiplicative) power operations generate all natural operations between the 0-th
equivariant homotopy groups of ultra-commutative ring spectra. 
The strategy is the one that we have employed several times before:
the functor $\pi_0^G$ from  ultra-commutative ring spectra to sets
is representable, namely by $\Sigma^\infty_+ \mP(B_{\gl}G)$,
the unreduced suspension spectrum of the free ultra-commutative monoid
generated by $B_{\gl} G$.
So we have to determine the equivariant homotopy groups 
$\pi_0^K(\Sigma^\infty_+ \mP(B_{\gl}G))$, which just means assembling various
results already proved.

\begin{theorem}\label{thm:generators global power functor} 
Let $G$ and $K$ be compact Lie groups.
The group of natural transformations  $\pi_0^G \to \pi_0^K$
of set valued functors on the category of ultra-commutative ring spectra 
is a free abelian group with basis the operations
\[   \tr_L^K\circ\, \alpha^*\circ P^m\ : \ \pi_0^G \ \to \pi_0^K \]
for all $m\geq 0$ and all $(K\times(\Sigma_m\wr G))$-conjugacy classes of pairs 
$(L,\alpha)$ consisting of a closed subgroup $L$ of $K$ with finite Weyl group
and a continuous homomorphism $\alpha:L\to\Sigma_m\wr G$.
\end{theorem}
\begin{proof}
  We let $W$ be any non-zero faithful $G$-representation and write
  $B_{\gl}G=\bL_{G,W}$ for the global classifying space of $G$ based on $W$.
  We denote by $u^{ucom}_G$ the class in $\pi_0^G( \Sigma^\infty_+ \mP( B_{\gl} G))$
  obtained by pushing forward the tautological class $u_{G,W}$ 
  along the adjunction unit $\eta:B_{\gl}G\to \mP( B_{\gl} G)$ 
  and then applying the stabilization map $\sigma^G$ of \eqref{eq:sigma_map}.

  We apply the representability result 
  of Proposition \ref{prop:pi_0^G representability}
  to the category of ultra-commutative ring spectra and the adjoint functor pair: 
  \[ \xymatrix@C=10mm{ 
\Sigma^\infty_+\circ \mP \ : \ \spc \ \ar@<.4ex>[r] & 
    \ ucom \ : \ U\circ\Omega^\bullet \ar@<.4ex>[l]  } \]
  If $V$ is any $G$-representation, then the restriction morphism
  $\rho_{G,V,W}:\bL_{G,V\oplus W}\to \bL_{G,W}$ is a global equivalence
  between positive flat orthogonal spaces. We showed in the proof of
  Theorem \ref{thm:symmetrizable in orthogonal spaces}~(ii)
  that the induced morphism of free ultra-commutative monoids
  $\mP(\rho_{G,V,W}):\mP(\bL_{G,V\oplus W})\to \mP(\bL_{G,W})$ is a global equivalence.
  So the induced morphism of unreduced suspension spectra
  $\Sigma^\infty_+\mP(\rho_{G,V,W})$ is a global equivalence by 
  Corollary \ref{cor:suspension spectrum globally homotopical}.
  In particular, the morphism of $\Rep$-functors
  $\upi_0(\Sigma^\infty_+\mP(\rho_{G,V,W}))$ is an isomorphism. So 
  Proposition \ref{prop:pi_0^G representability} applies and shows that
  evaluation at the tautological class is a bijection\index{subject}{ultra-commutative ring spectrum!free|see{free ultra-commutative ring spectrum}}\index{subject}{free ultra-commutative ring spectrum|(}
  \[   \Nat^{ucom}(\pi_0^G,\pi_0^K) 
  \ \to \ \pi_0^K(\Sigma^\infty_+ \mP(B_{\gl} G)) \ , \quad
  \tau\ \longmapsto \ \tau(u^{ucom}_G) \]
  to the 0-th $K$-equivariant homotopy group 
  of the ultra-commutative ring spectrum $\Sigma^\infty_+ \mP(B_{\gl} G)$.

  The identification \eqref{eq:P( B_gl G)}
  and Proposition \ref{prop:fix of global classifying}~(ii) 
  show that the homotopy set $\pi_0^L(\mP(B_{\gl} G))$ bijects with the set
  \[ {\coprod}_{m\geq 0}\, \Rep(L,\Sigma_m\wr G) \]
  by sending the conjugacy class of $\alpha:L\to\Sigma_m\wr G$ 
  to $\alpha^*([m](u_G))=\alpha^*(u_{\Sigma_m\wr G})$.
  Proposition \ref{prop:pi_0 of Sigma^infty} 
  then implies that the group $\pi_0^K(\Sigma^\infty_+ \mP(B_{\gl} G))$
  is a free abelian group, and it specifies a basis consisting of the elements
  \[ \tr_L^K(\sigma^L(x)) \]
  where $L$ runs through all conjugacy classes of closed subgroups of $K$
  with finite Weyl group and $x$ runs through a set of representatives 
  of the $W_K L$-orbits of the set $\pi_0^L (\mP(B_{\gl} G))$.
  So together this shows that 
  $\pi_0^K(\Sigma^\infty_+ \mP(B_{\gl} G))$ is a free abelian group with
  basis the classes 
  \begin{align*}
    \tr_L^K(\sigma^L(\alpha^*([m](u_G)))) \ &= \ 
    \tr_L^K(\alpha^*(\sigma^{\Sigma_m\wr G}([m](u_G)))) \\ 
    &= \  \tr_L^K(\alpha^*(P^m(\sigma^G(u_G)))) 
    \ = \  \tr_L^K(\alpha^*(P^m(u^{ucom}_G)))   
  \end{align*}
  for all $(m,L,\alpha)$ as in the statement of the theorem.
\end{proof}

Our last major goal in this chapter is to show that every global power functor
is realized by an ultra-commutative ring spectrum that can moreover
be chosen to have all its equivariant homotopy groups concentrated in
dimension~0, see Theorem \ref{thm:realize global power functor} below.
Towards this aim we will show various auxiliary results, which may be of
independent interest.
We start with a proposition that describes the effect on the homotopy group
global power functors of `coning off' a free ultra-commutative ring
spectrum. We can only analyze this process algebraically if all spectra involved
are globally connective, for otherwise the symmetric powers $\mP^m(A)$
generate non-trivial equivariant homotopy in arbitrarily low negative dimensions.

We recall from Example \ref{eg:coproduct of power functors}
that the coproduct in the category of global power functors
is given by the box product of underlying global functors.
So if $f:R\to F$ and $g:S\to F$ are morphisms of global power functors,
we let $f\diamond g:R\Box S\to F$ denote their categorical sum, i.e., the composite
\[ R\Box S \ \xra{\ f\Box g\ }\ F\Box F \ \xra{\text{multiplication}}\ F\ . \]
If $A$ is an orthogonal spectrum and $R$ an ultra-commutative ring spectrum,
then the zero morphism of orthogonal spectra from $A$ to $R$
freely generates a morphism of ultra-commutative ring spectra $\tau_A:\mP A\to R$
that we call the `trivial morphism'.
In terms of the wedge decomposition $\mP A=\bigvee_{m\geq 0}\mP^m(A)$,
this trivial morphism is the unit morphism of $R$ on $\mP^0(A)=\mS$, 
and trivial on all other wedge summands.
We recall that $C A=A\sm [0,1]$ is the cone of an orthogonal spectrum $A$,
and $i=(-\sm 1):A\to C A$ is the cone inclusion.

\begin{prop}\label{prop:umon cone}
Let $A$ be a positively flat orthogonal spectrum such that
the free ultra-commutative ring spectrum $\mP A$ is globally connective.
Let $R$ be a globally connective ultra-commutative ring spectrum and $\rho:\mP A\to R$
a morphism of ultra-commutative ring spectra.
We define $T$ as the following pushout, 
in the category of ultra-commutative ring spectra:
\[ \xymatrix{ \mP A\ar[r]^-{\mP i}\ar[d]_{\rho} & \mP (C A) \ar[d]\\ 
R \ar[r]_\psi & T
} \]
Then $T$ is globally connective, 
the morphism of global power functors $\upi_0(\psi):\upi_0(R)\to\upi_0(T)$
is surjective and 
a coequalizer, in the category of global power functors, of the two morphisms
\[ \upi_0(R)\diamond\upi_0(\rho)\ , \ \upi_0(R)\diamond\upi_0(\tau_A) \ : \ 
\upi_0(R)\Box\upi_0(\mP A)\ \to \ \upi_0(R)\ . \]
\end{prop}
\begin{proof}
The cone $C A=A\sm [0,1]$ is isomorphic to the realization 
of the simplicial orthogonal spectrum $B_\bullet(A,A,\ast)$,
the bar construction with respect to wedge.
Explicitly, 
\[ B_m(A,A,\ast)\ = \ A\sm \Delta[1]_m\ , \]
where $\Delta[1]$ is the simplicial 1-simplex, pointed by the vertex~0;
the simplicial structure is induced by that of $\Delta[1]$.
Since $\Delta[1]$ has $(m+1)$ non-basepoint simplices of dimension $m$, 
the spectrum $B_m(A,A,\ast)$ is a wedge of $m+1$ copies of $A$.

Since the free functor $\mP:\spec\to ucom$ is an enriched left adjoint,
it preserves coends and cotensors with spaces. 
So $\mP(C A)$ is isomorphic to the realization,
internal to ultra-commutative ring spectra,
of the simplicial ultra-commutative ring spectrum $\mP(B_\bullet(A,A,\ast))$.
Under this isomorphism, the morphism $\mP i:\mP A\to \mP(C A)$ identifies
$\mP A$ with the degeneracies of the object of 0-simplices.
As a left adjoint, the free functor preserves coproducts, 
i.e., it  takes wedges of orthogonal spectra
to smash products of ultra-commutative ring spectra.
So $\mP(B_\bullet(A,A,\ast))$ is isomorphic, 
as a simplicial ultra-commutative ring spectrum,
to $B_\bullet^{\sm}(\mP A,\mP A,\mS)$,
the bar construction of $\mP A$ with respect to smash product.
\index{subject}{realization!of simplicial ultra-commutative ring spectra}

Since colimits commute among themselves and with cotensors,
we can describe the pushout $T=R\sm _{\mP A}\mP(C A)$ as the realization 
of the simplicial ultra-commutative ring 
spectrum $R\sm_{\mP A} B_\bullet^{\sm}(\mP A,\mP A,\mS)$
where we take pushouts in every simplicial dimension.
The smash product over $\mP A$ `cancels' with the 0-th smash factor
in the bar construction, so we can take $T$ as the internal realization of 
the bar construction $B_\bullet^{\sm}(R,\mP A,\mS)$.

By the stable analog of Proposition \ref{prop:U preserves realization}, 
the realization internal to the category of ultra-commutative ring spectra 
can equivalently be taken in the underlying category of orthogonal spectra. 
The underlying realization
is the sequential colimit of partial realizations $B^{[n]}$,
i.e., `skeleta' in the simplicial direction, defined as the coend
\[ B^{[n]}\ = \ \int^{[m]\in\bDelta_\leq n}  B_m^{\sm}(R,\mP A,\mS) \sm \Delta^m_+ \]
of the restriction to the full subcategory $\bDelta_{\leq n}$
of $\bDelta$ with objects all $[m]$ with $m\leq n$.
We refer to \cite[VII Sec.\ 1-3]{goerss-jardine}
for background on realization of simplicial objects and the skeleton filtration.
The realization $|B_\bullet^{\sm}(R, \mP A, \mS)|$ 
is the colimit of the sequence of orthogonal spectra
\begin{equation}\label{eq:filtration of ucom po}
 R\ =\  B^{[0]} \ \to \ B^{[1]}\ \to \ \dots \ 
\to \ B^{[n]}\ \to \ \dots \ .   
\end{equation}
Moreover, the $n$-skeleton $B^{[n]}$ is obtained from the previous stage
as a pushout of orthogonal spectra:
\begin{equation}  \begin{aligned}\label{eq:simplicial skeleton pushout}
\xymatrix@C=13mm{ 
L_n^\bDelta\sm \Delta^n_+\cup_{ L_n^\bDelta\sm \partial\Delta^n_+}
B_n^{\sm}(R,\mP A,\mS)\sm\partial\Delta^n_+\ar[r]^-{\text{incl}}\ar[d] &
B_n^{\sm}(R,\mP A,\mS)\sm\Delta^n_+\ar[d] \\
B^{[n-1]}\ar[r] &B^{[n]} }   
  \end{aligned}\end{equation}
see~(3.8) of \cite[VII Sec.\,3]{goerss-jardine}.
Here the orthogonal spectrum $L_n^\bDelta$ 
is the $n$-th latching object in the simplicial direction, 
compare~(1.5) of \cite[VII Sec.\,1]{goerss-jardine},
or \eqref{eq:simplicial_latching}.
In our situation, the latching morphism
$ L_n^\bDelta\to B_n^{\sm}(R,\mP A,\mS)$ is the morphism
  \[ R\sm i^{\Box n} \ : \ R\sm Q^n(i)\ \to \ 
R\sm (\mP A)^{\sm n} \ = \ B_n^{\sm}(R,\mP A,\mS)\ ,\]
where $i^{\Box n}$ is the iterated pushout product, with respect to
smash product of orthogonal spectra,
of the unit morphism $\mS\to \mP A$. Since 
flat cofibrations are symmetrizable 
(Theorem \ref{thm:symmetrizable in orthogonal spectra}~(i)), this unit morphism
is a flat cofibration.
The pushout product property of the flat cofibrations 
(Proposition \ref{prop:ExF ppp})
thus  shows that $i^{\Box n}$, and hence the latching morphism, is a flat cofibration.

So both horizontal morphisms 
in the pushout square \eqref{eq:simplicial skeleton pushout}
are h-cofibrations, and the cokernel of the inclusion $B^{[n-1]}\to B^{[n]}$
is isomorphic to 
\[ (B_n^{\sm}(R,\mP A,\mS) / L_n^{\bDelta})\sm S^n \ \iso \ 
R\sm (\bar\mP A)^{\sm n}\sm S^n\]
where $\bar\mP A=(\mP A)/\mS=\bigvee_{m\geq 1}\mP^m(A)$.
Since $A$ is flat and flat cofibrations are symmetrizable,
the underlying orthogonal spectrum of $\bar\mP A$ is flat.
Since $\mP A$ is globally connective, so is $\bar\mP A$.
Since $R$ is also globally connective, so is
$R\sm (\bar\mP A)^{\sm n}$, by Proposition \ref{prop:Box vs smash}. 
This shows that the cokernel of the inclusion $B^{[n-1]}\to B^{[n]}$
is globally $(n-1)$-connected. Since the inclusion is
also an h-cofibration,
Corollary \ref{cor-long exact sequence h-cofibration}~(i)
shows that all the morphisms in the sequence \eqref{eq:filtration of ucom po}
induce isomorphisms on the global functors $\upi_k$ for all $k< 0$
and an epimorphism on $\upi_0$. Moreover, starting from $B^{[1]}$ onward,
the morphisms even induce an isomorphism on $\upi_0$.
Each morphism in the sequence is an h-cofibration,
and so levelwise a closed embedding.
Since equivariant homotopy groups commute with sequential colimits over closed embeddings
(Proposition \ref{prop:sequential colimit closed embeddings}~(i)),
we conclude that also the morphism $\upi_k(\psi):\upi_k(R)\to \upi_k(T)$ 
is an isomorphism of global functors for all $k < 0$
and an epimorphism for $k=0$.
Since $R$ is globally connective, 
this shows in particular that the colimit $T$ is globally connective.
Moreover, the inclusion $B^{[1]}\to |B_\bullet^{\sm}(R,\mP A,\mS)|=T$ 
induces an isomorphism on $\upi_0$.

It remains to identify $\upi_0(\psi):\upi_0(R)\to\upi_0(T)$ as a coequalizer.
The latching object $L_1^\bDelta$ is simply the spectrum $B_0^{\sm}(R,\mP A,\mS)=R$, 
and for $n=1$ the pushout \eqref{eq:simplicial skeleton pushout}
specializes to a pushout of orthogonal spectra
\begin{equation}
  \begin{aligned}\label{eq:bar 1 in ucom}
\xymatrix@C=12mm{ 
 R\sm \mP A\sm\{0,1\}_+\ar[r]^-{\text{incl}}\ar[d]_{d_0+ d_1} &
 R\sm \mP A\sm [0,1]_+\ar[d] \\
R\ar[r] &B^{[1]} }   
  \end{aligned}\end{equation}
The orthogonal spectrum $R\sm \mP A\sm\{0,1\}_+$
is the wedge of two copies of $R\sm \mP A$.
On one copy the left vertical morphism is the simplicial face morphism
$d_0:R\sm \mP A\to R$, which is the product, with respect to
the multiplication of $R$, of the identity and the morphism $\rho:\mP A\to R$.
On the other copy the left vertical morphism is 
the simplicial face morphism $d_1:R\sm \mP A\to R$, which is the product
of the identity of $R$ and the trivial morphism $\tau_A:\mP A\to R$,
adjoint to the zero morphism from $A$ to $R$.

The orthogonal spectra underlying $R$ and $\mP A$ are globally connective
by hypothesis. Moreover, $\mP A$ is flat because $A$ is and because
flat cofibrations are symmetrizable.
So by Proposition \ref{prop:Box vs smash},
the orthogonal spectrum $R\sm \mP A$ is again globally connective. 
Since the two horizontal morphisms in the pushout square \eqref{eq:bar 1 in ucom}
are h-cofibrations, the equivariant homotopy groups of the pushout
$B^{[1]}$ participate in an exact Mayer-Vietoris sequence
\[ \upi_0(R\sm \mP A)\oplus \upi_0(R\sm \mP A)
\ \to \ 
\upi_0(R)\oplus\upi_0(R\sm \mP A)\ \to \ \upi_0(B^{[1]}) \ \to \ 0 \ , \]
where the first map is given by
\[ (x,y)\ \longmapsto \ (\upi_0(d_0)(x)+\upi_0(d_1)(y),\, x+y)\ . \]
Using the identification $\upi_0(B^{[1]})\iso\upi_0(T)$
induced by the inclusion $B^{[1]}\to T$, we can rewrite this as an
exact sequence of global functors
\begin{equation}\label{eq:exact_for_ucom_pushout}
 \upi_0(R\sm \mP A)\ \xra{\upi_0(d_0) - \upi_0(d_1)} \ 
\upi_0(R)\ \xra{\ \upi_0(\psi)\ } \ \upi_0(T) \ \to \ 0 \ . \end{equation}
Since $R$ and $\mP A$ are globally connective and $\mP A$ is flat, the canonical morphism
\[ \upi_0(R)\Box \upi_0(\mP A)\ \to \ \upi_0(R\sm \mP A) \]
is an isomorphism of global power functors
by Example \ref{eg:coproduct of power functors}. 
Under this identification, the effects of the two face morphisms $d_0,d_1:R\sm \mP A\to R$
on $\upi_0$ correspond to the morphisms of global power functors
\[ \upi_0(R)\diamond\upi_0(\rho)\text{\qquad respectively\qquad}
\upi_0(R)\diamond\upi_0(\tau_A)  \ . \]
So we may show that $\upi_0(\psi)$ is a coequalizer
of the two morphisms $\upi_0(d_0)$ and $\upi_0(d_1)$.
To this end we consider a morphism of global power functors
$\varphi:\upi_0(R)\to F$
satisfying $\varphi\circ \upi_0(d_0)=\varphi\circ \upi_0(d_1)$.
By exactness of the sequence \eqref{eq:exact_for_ucom_pushout}
there is a unique morphism of global functors $\bar\varphi:\upi_0(T)\to F$
such that $\bar\varphi\circ \upi_0(\psi)=\varphi$.
Since $\upi_0(\psi)$ is a surjective homomorphism of global power functors
and the composite $\bar\varphi\circ \upi_0(\psi)$ is a homomorphism
of global power functors, the morphism $\bar\varphi$ is also compatible
with multiplication and power operations, hence a morphism of 
global power functors. This establishes the universal property of a coequalizer,
and hence completes the proof.
\end{proof}

Now we have all the necessary tools to show that the 
topological construction of free ultra-commutative ring spectra
realizes the algebraic construction of free global power functors,
at least for sufficiently cofibrant and globally connective spectra.

\begin{defn}
  A {\em free global power functor}\index{subject}{global power functor!freely generated by a global functor|(}  
  generated by a global functor $F$ is a pair $(R,i)$ consisting of 
  a global power functor $R$ and a morphism of global functors $i:F\to R$ 
  with the following property:
  for every global power functor $S$ and every morphism of 
  global functors $j:F\to S$
  there is a unique morphism of global power functors $\varphi:R\to S$
  such that $\varphi\circ i=j$.
\end{defn}
  
In other words, $(R,i)$ is a free global power functor for $F$ 
if and only if it represents the functor
\[   \GlPow\ \to \ \text{(sets)} \ , \quad S \ \longmapsto \ \GF(F,S) \ . \]

For every orthogonal spectrum $X$, the free ultra-commutative ring spectrum $\mP X$
gives rise to a global power functor $\upi_0(\mP X)$.
The adjunction unit $\eta:X\to\mP X$ (i.e., the inclusion of the `linear' wedge summand)
is a morphism of orthogonal spectra, and thus induces a morphism of global functors
\[ \upi_0(\eta)\ : \ \upi_0(X)\ \to \ \upi_0(\mP X)\ . \]

\begin{theorem}\label{thm:free global power}
Let $X$ be a positively flat globally connective orthogonal spectrum.
Then the free ultra-commutative ring spectrum $\mP X$ is globally connective and
the pair $(\upi_0(\mP X),\upi_0(\eta))$ is a free global power functor generated
by the global functor $\upi_0(X)$.
\end{theorem}
\begin{proof} We let $\Yc$ be the class
of globally connective orthogonal spectra $Y$ such that 
for some (hence any) global equivalence $Y^c\to Y$ from a 
positively flat orthogonal spectrum the claim of the theorem holds for $Y^c$.
We let $G$ be a compact Lie group and $V$ a non-zero faithful $G$-representation.
Then $B_{\gl}G=\bL_{G,V}$ is a global classifying space for $G$.
The free ultra-commutative ring spectrum
$\mP(\Sigma^\infty_+ B_{\gl}G)$ is isomorphic to the unreduced suspension spectrum
of the free ultra-commutative monoid $\mP(B_{\gl}G)$; 
so $\mP(\Sigma^\infty_+ B_{\gl}G)$ is globally connective by 
Proposition \ref{prop:pi_0 of Sigma^infty}.
By Corollary \ref{cor:B_gl free ucom} the global power functor
\[ \upi_0(\mP(\Sigma^\infty_+ B_{\gl}G)) \]
is freely generated by the class $u_G^{ucom}$ in $\pi_0^G(\mP(\Sigma^\infty_+ B_{\gl}G))$.
Moreover,  the global functor $\upi_0(\Sigma^\infty_+ B_{\gl}G)$ is
freely generated by the stable tautological class $e_G$,
by Proposition \ref{prop:B_gl represents}.
So altogether this shows that the orthogonal spectrum $\Sigma^\infty_+ B_{\gl}G$
belongs to the class $\Yc$.

Now we establish two closure properties for the class $\Yc$.
We let $\{X_i\}_{i\in I}$ be a family of orthogonal spectra in $\Yc$;
we claim that the wedge $\bigvee_{i\in I}X_i$ also belongs to $\Yc$.
Since wedges of orthogonal spectra preserve global equivalences 
(by Corollary \ref{cor-wedges and finite products}~(i)),
we can assume that each $X_i$ is positively flat. 
Then the wedge $\bigvee_{i\in I}X_i$ is also positively flat.
The underlying orthogonal spectrum of $\mP(\bigvee_{i\in I}X_i)$
is the wedge, indexed by functions
$f: I\to \mN$ with finite support (i.e., $f(i)=0$ for almost all $i\in I$)
of smash products of the finite symmetric powers:
\[ \mP\left(  {\bigvee}_{i\in I}\, X_i \right) \ \iso \ 
\bigvee_{f:I\to\mN, \text{finite supp}} \, \bigwedge_{f(i)\ne 0} \mP^{f(i)}(X_i)\ . \]
Since flat cofibrations are symmetrizable and $X_i$ is flat, the
symmetric power $\mP^{f(i)}(X_i)$ is flat. 
Since $\mP(X_i)$ is globally connective, so is its
retract $\mP^{f(i)}(X_i)$. So for every finitely supported function
$f:I\to\mN$, the morphism
\[ {\bigbox}_{f(i)\ne 0}\, \upi_0\left(\mP^{f(i)}(X_i)\right) \ \to \ 
\upi_0\left({\bigwedge}_{f(i)\ne 0}\, \mP^{f(i)}(X_i) \right)  \]
from the iterated box product is an isomorphism of global functors by 
Proposition \ref{prop:Box vs smash}. 
Since equivariant homotopy groups take wedges to sums by 
Corollary \ref{cor-wedges and finite products}~(i), 
this establishes a sequence of isomorphisms of global functors 
\begin{align*}
\colim_{J\subset I,\ J\text{ finite}}  \bigbox_{j\in J}  \upi_0(\mP (X_j))
\ &\iso \ 
\bigoplus_{f:I\to\mN, \text{finite supp}} \, \bigbox_{f(i)\ne 0} \, \upi_0(\mP^{f(i)}(X_i))\\
&\iso \ \bigoplus_{f:I\to\mN, \text{finite supp}} \, 
\upi_0\left({\bigwedge}_{f(i)\ne 0}\, \mP^{f(i)}(X_i) \right)\\
&\xra{\ \iso \ } \ \upi_0\left( \mP\left(  {\bigvee}_{i\in I}\, X_i \right)\right)   
\end{align*}
We claim that the source of this isomorphism is the coproduct,
in the category of global power functors, of the global power functors
$\upi_0(\mP(X_i))$.
Indeed, the coproduct of two global power functors is given by the
box product of underlying global functors, 
see Example \ref{eg:coproduct of power functors}.
As in any category, a general (possibly infinite) $I$-indexed coproduct 
can be obtained as the filtered colimit,
formed over the poset of finite subsets of $I$, of the finite coproducts.
Moreover, filtered colimits in the category of global Green functors,
and hence also of global power functors,
are created in the underlying category of global functors,
see Proposition \ref{prop:filtered colim P-alg}.
This proves the claim.

Since $\upi_0(\mP(X_i))$ is freely generated by the global functor $\upi_0(X_i)$
and coproducts of free objects are free,
the right hand side of the isomorphism is a free global power functor
generated by the global functor $\bigoplus_{i\in I}\upi_0(X_i )$.
Equivariant homotopy groups take wedges to direct sums, so the canonical morphism
\[ {\bigoplus}_{i\in I} \, \upi_0(X_i)\ \to \ \upi_0\left({\bigvee}_{i\in I} X_i \right) \]
is an isomorphism of global functors.
So altogether we obtain that $\upi_0(\mP(  \bigvee_{i\in I} X_i))$
is a free global power functor generated by the global functor
 $\upi_0(\bigvee_{i\in I} X_i)$.
This concludes the proof that the class $\Yc$ is closed under wedges.

Now we claim that the class $\Yc$ is also closed under cones.
In other words, we let
\[  A \ \to \ B \ \to \ C \ \to \ A \sm S^1 \]
be a distinguished triangle in the global stable homotopy category 
such that $A$ and $B$ belong to $\Yc$; we must show that $C$ also belongs to $\Yc$.
After replacing the distinguished triangle by an isomorphic one,
we can assume that $A$ and $B$ are positively flat and $C=C f$ is
the mapping cone of a morphism $f:A\to B$ of orthogonal spectra, i.e., the
pushout of the diagram of orthogonal spectra
\[ B \ \xla{\ f\ }\ A \ \xra{\ i \ }\ C A \ , \]
where the right morphism is the cone inclusion.
As a left adjoint, the free functor $\mP$ preserves pushouts,
so $\mP(C f)$ is a pushout of the diagram of ultra-commutative ring spectra
\[ \mP B \ \xla{\ \mP f\ }\ \mP A \ \xra{\ \mP i \ }\ \mP( C A )\ . \]
Since $A$ and $B$ belong to $\Yc$, both $\mP B$ and $\mP A$
are globally connective. So Proposition \ref{prop:umon cone} applies
and shows that the morphism
\[ \mP j\ : \ \mP B \ \to \ \mP( C f) \]
is a coequalizer, in the category of global power functors, of the two morphisms
\[ \upi_0(\mP B)\diamond\upi_0(\mP f)\ , \ \upi_0(\mP B)\diamond\upi_0(\tau_A) \ : \ 
\upi_0(\mP B)\Box\upi_0(\mP A)\ \to \ \upi_0(\mP B)\ . \]
Since $A$ and $B$ belong to the class $\Yc$, 
the global power functors $\upi_0(\mP A)$ and $\upi_0(\mP B)$ are freely 
generated by the global functors $\upi_0(A)$ respectively $\upi_0(B)$.
So the coequalizer property witnesses that the global power functor $\upi_0(\mP( C f))$
is freely generated by the cokernel of the morphism of global functors 
$\upi_0(f):\upi_0(A)\to\upi_0(B)$.
Since $A$ is globally connective,
\[  \upi_0(A)\ \xra{\upi_0(f)} \ \upi_0(B)\ \to \ \upi_0(C f)\ 
\ \to\  0  \]
is an exact sequence of global functors
by Proposition \ref{prop:LES for homotopy of cone and fibre}.
This concludes the proof that the pair $(\upi_0(\mP(C f)), \upi_0(\eta))$
is a free global power functor generated by the global functor $\upi_0(C f)$.
So we have also shown that the class $\Yc$ is closed under cones.

Altogether this shows that the class of orthogonal spectra $\Yc$ 
is closed under wedges and cones 
and contains the suspension spectra $\Sigma^\infty_+ B_{\gl}G$ for all
compact Lie groups $G$.
Proposition \ref{prop:plus generating GH} 
then shows that $\Yc$ is the class of all globally connective orthogonal spectra.
This proves the theorem.
\end{proof}
\index{subject}{global power functor!freely generated by a global functor|)}

\begin{theorem}\label{thm:power connectivity}
Let $X$ be a positively flat globally $(n-1)$-connected orthogonal spectrum,
for some $n\geq 0$. Then for every $m\geq 1$ the orthogonal spectrum
$\mP^m(X)$ is globally $(n-1)$-connected.
\end{theorem}
\begin{proof}
We choose a cofibrant replacement $q:Y\to \Omega^n X$ of $\Omega^n X$
in the positive global model structure of 
Proposition \ref{prop:positive global spectra}.
Then $Y$ is globally connective (i.e., globally $(-1)$-connected).
The adjoint $q^\flat$ of $q$ factors as the composite
\[ Y\sm S^n \ \xra{q\sm S^n} \ (\Omega^n X)\sm S^n \ \xra{\ \epsilon_X\ }\ X \ . \]
The first morphism is a global equivalence since suspension is fully homotopical,
and the second morphism is a global equivalence
by Proposition \ref{prop:lambda upi_* isos}~(ii).
So the morphism $q^\flat$ is a global equivalence.

We show that $\mP^m(Y\sm S^n)$ is globally $(n-1)$-connected.
We let
\[ \nu_m\ = \ \{(x_1,\dots,x_m)\in\mR^m\ | \ x_1+\dots+x_m=0\} \]
denote the reduced natural $\Sigma_m$-representation,
where $\Sigma_m$ acts by permutation of coordinates.
The smash power $(S^n)^{\sm m}$ is $\Sigma_m$-equivariantly homeomorphic to
the smash product of $S^n$ (with trivial $\Sigma_m$-action)
and the representation sphere $S^{\mR^n\tensor\nu_m}$.
This decomposition induces an isomorphism of orthogonal spectra
\[ \mP^m(Y\sm S^n)\ \iso \ 
 (Y^{\sm m}\sm (S^n)^{\sm m})/\Sigma_m \ \iso \  
(Y^{\sm m}\sm S^{\mR^n\tensor\nu_m})/\Sigma_m \sm S^n \ .  \]
The sphere $S^{\mR^n\tensor\nu_m}$ admits a finite $\Sigma_m$-CW-structure.
Moreover, all isotropy groups are conjugate to subgroups of the form
$\Sigma_{i_1}\times\dots\times\Sigma_{i_k}$ with $i_1+\dots+i_k=m$.
So the orthogonal spectrum $\mP^m(Y\sm S^n)$ admits a finite filtration,
along h-cofibrations, with subquotients isomorphic to 
\begin{align*}
  (Y^{\sm m}\sm \Sigma_m/(\Sigma_{i_1}\times\dots\times\Sigma_{i_k})_+\sm S^k)/\Sigma_m\sm S^n \ &\iso\  
 Y^{\sm m}/(\Sigma_{i_1}\times\dots\times\Sigma_{i_k})\sm S^{k+n} \\ 
&\iso\    \mP^{i_1}(Y)\sm \dots\sm \mP^{i_k}(Y)\sm S^{k+n} \ . 
\end{align*}
Since $Y$ is positively flat and globally connective,
Theorem \ref{thm:free global power}
says that the free ultra-commutative ring spectrum $\mP Y$ is 
globally connective. Since $\mP^{i_k}(Y)$ is a retract of $\mP Y$,
it is also globally connective. 
Since flat cofibrations are symmetrizable
(Theorem \ref{thm:symmetrizable in orthogonal spectra}~(i)), 
each of the orthogonal spectra $\mP^{i_j}(Y)$ is again flat. 
So the smash product
$\mP^{i_1}(Y)\sm\dots\sm \mP^{i_k}(Y)$ is globally connective by
Proposition \ref{prop:Box vs smash}.
Altogether this shows that each of the subquotients of $\mP^m(Y\sm S^n)$
is globally $(n-1)$-connected. So the spectrum $\mP^m(Y\sm S^n)$ itself
is globally $(n-1)$-connected.

The free functor $\mP:\spec\to ucom$ is a left Quillen functor
from the positive global model structure on orthogonal spectra
to the global model structure of ultra-commutative ring spectra.
So $\mP$ preserves global equivalences 
between positively flat orthogonal spectra, hence so does its retract $\mP^m$.
So the morphism $\mP^m(q^\flat):\mP^m(Y\sm S^n)\to \mP^m(X)$ is a global equivalence,
and $\mP^m(X)$ is globally $(n-1)$-connected.
\end{proof}

The next theorem shows that we can `kill higher homotopy groups' in
the world of ultra-commutative ring spectra. Again, the construction
needs the hypothesis of global connectivity.

\begin{theorem}\label{thm:kill homotopy ucom} 
Let $R$ be a globally connective ultra-commutative ring spectrum and $n\geq 1$.
Then there is a cofibration of ultra-commutative ring spectra
$\psi:R\to T$ with the following properties:
\begin{itemize}
\item the induced morphism of global functors $\upi_k(\psi):\upi_k(R)\to\upi_k(T)$
is bijective for $k<n$, and
\item the global functor $\upi_n(T)$ is trivial.
\end{itemize}
\end{theorem}
\begin{proof}
We choose an index set $J$, compact Lie groups $G_j$ and
elements $y_j \in \pi_n^{G_j}(R)$, for $j\in J$, 
that altogether generate the global functor $\upi_n(R)$.
We represent each class $y_j$ as a morphism of
orthogonal spectra
\[  f_j \ : \ \Sigma^\infty_+ B_{\gl} G_j \sm S^n\ \to \ R  \]
that sends the $n$-fold suspension of the
stable tautological class $e_{G_j}$ to $y_j$; this involves an
implicit choice of non-zero faithful $G_j$-representation.
We form the wedge of all these morphisms 
\[ F \ : \ X \ = \ {\bigvee}_{j\in J}\, \Sigma^\infty_+ B_{\gl} G_j \sm S^n  \ \to \ R\ .\]
All we need to remember about $F$ is that its source $X$ is globally $(n-1)$-connected
and positively flat, and that the morphism of global functors
\[ \upi_n(F)\ : \ \upi_n(X)\ \to \ \upi_n(R) \]
is surjective.
We extend $F$ freely to a morphism of ultra-commutative ring spectra
$\tilde F:\mP X \to R$.
We let $T$ be the pushout, in the category
of ultra-commutative ring spectra, of the diagram
\[ R \ \xla{\ \tilde F \ }
\  \mP  X  \ \xra{\ \mP i\ } \ \mP(C X) \ , \]
where $C X$ is the cone of the orthogonal spectrum $X$,
and $i:X\to C X$ is the cone inclusion.
Since $X$ is positively flat, the cone inclusion is a positive cofibration
of orthogonal spectra, and so $\mP i$ is a cofibration of ultra-commutative ring
spectra. The cobase change $\psi:R\to T$ of $\mP i$ is then also
a cofibration of ultra-commutative ring spectra.

We recall that the cobase change
of a free morphism $\mP i:\mP X\to\mP(C X)$
comes with a filtration by `number of factors from $\mP(C X)$',
see for example \cite[Prop.\,B.2]{white-commutative monoids}.
More precisely, the morphism $\psi:R\to T$ is the sequential composite
of  morphisms of $R$-module spectra
\[ R \ = \ P_0 \  \xra{\ \psi_1\ } \ P_1 \ \xra{\ \psi_2\ } \ P_2 \
\xra{\ \psi_3\ } \ \cdots  \ ,   \]
such that each step in the filtration can be obtained from the previous one
as a pushout of $R$-modules:
\[ \xymatrix@C=12mm{ 
R\sm Q^m(i)/\Sigma_m\ar[r]^-{i^{\Box m}/\Sigma_m}\ar[d] &
R\sm (C X)^{\sm m}/\Sigma_m\ar[d]\\
P_{m-1}\ar[r]_-{\psi_m}& P_m} \]
Since $i:X\to C X$ is a flat cofibration, and flat cofibrations
are symmetrizable (Theorem \ref{thm:symmetrizable in orthogonal spectra}~(i)),
the upper horizontal morphism is a cofibration of $R$-modules,
hence so is its cobase change $\psi_m:P_{m-1}\to P_m$.
Moreover, the pushout square witnesses that the cokernel of $\psi_m$ is isomorphic to
\[ 
R\sm  \text{coker}(i^{\Box m})/\Sigma_m\ \iso \ 
R\sm (C X/X )^{\sm m}/\Sigma_m \ \iso \ R\sm \mP^m(X\sm S^1) \ .\]
Since $X$ is globally $(n-1)$-connected, its suspension $X\sm S^1$
is globally $n$-connected, and so the $m$-th symmetric power
$\mP^m(X\sm S^1)$ is also globally $n$-connected for all $m\geq 1$,
by Theorem \ref{thm:power connectivity}.
Since $R$ is connective and $\mP^m(X\sm S^1)$ is flat,
the smash product $R\sm \mP^m(X\sm S^1)$ is again
globally $n$-connected by Proposition \ref{prop:Box vs smash}.
In particular, the morphism~ $\upi_k(\psi_m):\upi_k(P_{m-1})\to \upi_k(P_m)$
is bijective for $k<n$ and surjective for $k=n$.
Since all morphisms $\psi_m:P_{m-1}\to P_m$ are h-cofibrations of orthogonal spectra,
they are levelwise closed embeddings, and so that canonical map
\[ \colim_m\, \upi_k(P_m)\ \to \ \upi_k(\colim_m P_m) \ = \ \upi_k(T)    \]
is an isomorphism by Proposition \ref{prop:sequential colimit closed embeddings}~(i).
Altogether this shows that the morphism
$\upi_k(\psi):\upi_k(R)\to\upi_k(T)$
is bijective for $k<n$ and surjective for $k=n$.

The diagram
\[ \xymatrix@C=15mm{ 
\upi_n(X) \ar[r]^-{\upi_n(\eta_X)} \ar[dr]_-{\upi_n(F)} &
\upi_n(\mP X) \ar[r]^-{\upi_n(\mP(i))} 
 \ar[d]^{\upi_n(\tilde F)}& \upi_n(\mP(C X))\ar[d]\\
& \upi_n(R)\ar[r]_-{\upi_n(\psi)}&\upi_n(T)} \]
commutes and the composite through the upper right corner is the zero map
because $C X$ is contractible.
Since $\upi_n(F)$ and $\upi_n(\psi)$ are surjective, we conclude that
the global functor $\upi_n(T)$ is trivial.
\end{proof}

Now we can finally show that every global power
functor is realized by an ultra-commutative ring spectrum.
The analogous result in equivariant stable homotopy theory for a fixed
finite group has been obtained by Ullman \cite{ullman-tambara}.
More is true: the next theorem effectively constructs a right adjoint functor
\[ H \ : \  \GlPow \ \to \ \Ho^{\text{gl.\,connective}}(ucom) \]
to the functor $\upi_0$ such that the adjunction counit is an isomorphism
$\upi_0(H F)\iso F$ of global power functors.
In other words, for every globally connective ultra-commutative ring spectrum $T$,
the functor $\upi_0$ restricts to a bijection
\[ \upi_0\ : \  \Ho(ucom)(T, H F)\ \iso \
 \GlPow(\upi_0(T), F)\ . \]
However, we won't prove these more general facts.

\begin{theorem}\label{thm:realize global power functor} 
Let $F$ be a global power functor.\index{subject}{Eilenberg-Mac\,Lane spectrum!of a global power functor} 
There is an ultra-commu\-ta\-tive ring spectrum $H F$
such that $\upi_k(H F)=0$ for all $k\ne 0$  and an isomorphism of global power functors
\[ \upi_0(H F)\ \iso \ F\ . \]
\end{theorem}
\begin{proof}
The t-structure on the global stable homotopy category
of Theorem \ref{thm:t-structure on GH} 
(for the global family of all compact Lie groups) 
lets us choose an Eilen\-berg-Mac\,Lane spectrum for the underlying global functor
of $F$, see Remark \ref{rk:general Eilenberg Mac Lane}.
Cofibrant replacement in the positive global model structure 
of Proposition \ref{prop:positive global spectra} 
then provides a globally connective positive flat orthogonal spectrum $X$ 
and an isomorphism of global functors $\upi_0(X)\iso F$.

We form the free ultra-commutative ring spectrum $\mP X$.
Theorem \ref{thm:free global power} then shows that $\mP X$ is globally connective 
and $\upi_0(\mP X)$ is a free global power functor generated
by the global functor $\upi_0(X)$, with respect to
the morphism  $\upi_0(\eta):\upi_0(X)\to \upi_0(\mP X)$.
In particular, there is a unique morphism of global power functors
\[ \epsilon\ : \ \upi_0(\mP X)\ \to \ F \]
such that $\epsilon\circ\upi_0(\eta)$ is the previous isomorphism.
The morphism $\epsilon$ is then surjective.

We choose a positively cofibrant globally connective orthogonal spectrum $A$
together with a morphism of orthogonal spectra $\rho':A\to \mP X$
such that the image of the morphism $\upi_0(\rho'):\upi_0(A)\to\upi_0(\mP X)$
coincides with the kernel of $\epsilon$.
For example, we can choose compact Lie groups $G_j$ and
elements $y_j \in \pi_0^{G_j}(\mP X)$
that altogether generate this kernel as a global functor.
Then we represent each class $y_j$ as a morphism of
orthogonal spectra
\[  f_j \ : \ \Sigma^\infty_+ B_{\gl} G_j \ \to \ \mP X  \]
that sends the stable tautological class $e_{G_j}$ to $y_j$,
and let $\rho':A\to\mP X$ be the wedge of all these morphisms.

We let $\rho:\mP A\to\mP X$ be the free extension of $\rho'$ to
a morphism of ultra-commutative ring spectra. Then we
define $T$ as the following pushout, 
in the category of ultra-commutative ring spectra:
\[ \xymatrix{ \mP A\ar[r]^-{\mP i}\ar[d]_{\rho} & \mP (C A) \ar[d]\\ 
\mP X \ar[r]_\psi & T
} \]
Since $\mP A$ and $\mP X$ are both globally connective by
Theorem \ref{thm:free global power}, we can apply Proposition \ref{prop:umon cone}
and conclude that $T$ is globally connective, 
the morphism of global power functors $\upi_0(\psi):\upi_0(\mP X)\to\upi_0(T)$
is surjective, and a coequalizer, 
in the category of global power functors, of the two morphisms
\[ \upi_0(\mP X)\diamond\upi_0(\rho)\ , \ \upi_0(\mP X)\diamond\upi_0(\tau_A) \ : \ 
\upi_0(\mP X)\Box\upi_0(\mP A)\ \to \ \upi_0(\mP X)\ . \]
Now $\upi_0(\mP A)$ is the free global power functor generated by
the global functor $\upi_0(A)$, by Theorem \ref{thm:free global power}.
So the effect of the coequalizer is to annihilate the global power ideal 
generated by the image of $\upi_0(\rho'):\upi_0(A)\to\upi_0(\mP X)$. 
This global power ideal is, by construction, the kernel of $\epsilon:\upi_0(\mP X)\to F$.
So $\epsilon$ descends to an isomorphism of global power functors
between $\upi_0(T)$ and $F$.

Now we use Theorem \ref{thm:kill homotopy ucom} 
to kill the higher homotopy groups of $T$. More precisely,
we construct a sequence of cofibrations of ultra-commutative ring spectra
\begin{equation}  \label{eq:HR_colimit_sequence}
 T \ = \ T_0 \ \to \ T_1 \ \to \ \dots \ \to \ T_n \ \to \ \dots  
\end{equation}
by induction on $n$. We obtain $T_n$ by applying Theorem \ref{thm:kill homotopy ucom}
in dimension $n$ to $R=T_{n-1}$. 
Then $T_n$ is globally connective, $\upi_0(T_n)$ is isomorphic to
$\upi_0(T)$, and hence to $F$, 
and the global functors $\upi_k(T_n)$ are trivial for all $1\leq k\leq n$.

Finally, we define $H F$ as the colimit of the sequence \eqref{eq:HR_colimit_sequence}
of ultra-commu\-ta\-tive ring spectra.
Each morphism in the sequence is a cofibration of ultra-commutative ring spectra,
hence an h-cofibration of underlying orthogonal spectra
(see Theorem \ref{thm:global ultra-commutative}~(ii)),
and so levelwise a closed embedding.
Since equivariant homotopy groups commute with sequential colimits over closed embeddings
(see Proposition \ref{prop:sequential colimit closed embeddings}~(i)),
the colimit has the desired properties of an
ultra-commutative Eilenberg-Mac Lane spectrum for $F$.
\end{proof}
\index{subject}{free ultra-commutative ring spectrum|)}

\chapter{Global Thom and \texorpdfstring{$K$}{K}-theory spectra}
\label{ch:examples}

The final chapter of this book is devoted to an in-depth study 
of interesting examples of ultra-commutative ring spectra, in particular
global Thom spectra and ultra-commutative global models for various flavors
of equivariant $K$-theory spectra.
In Section \ref{sec:global Thom} we discuss global refinements of
the classical Thom spectrum $M O$, as well as Thom spectra for the other families
of classical Lie groups, such as $M S O$ and $M U$. 
We carefully analyze two different global forms 
of the Thom spectrum $M O$ that represents unoriented bordism,
namely the ultra-commutative Thom ring spectrum $\bMO$, 
and an  $E_\infty$-orthogonal ring spectrum $\bmO$.
Both Thom spectra are the homogeneous degree~0 summands in certain
$\mZ$-graded periodic extensions $\bMOP$ respectively $\bmOP$.
The Thom spectrum $\bmO$ is the natural target for the Thom-Pontryagin map
from geometric equivariant bordism.
The ultra-commutative ring spectrum $\bMO$ is 
the $\mR$-analog of tom\,Dieck's equivariant homotopical bordism.
The equivariant homology theory represented by $\bMO$ is a localization
of the one represented by $\bmO$, formed by inverting certain `inverse Thom classes',
see Corollary \ref{cor:MO is localized mO}.
The Thom spectrum $\bmO$ comes with a rank filtration
that we use to show that $\bmO$ is globally connective,
and to calculate the global functor $\upi_0(\bmO)$, see Theorem \ref{thm:pi_0 of mO}.

Section \ref{sec:equivariant bordism}
is devoted to equivariant bordism and its relation to
the equivariant homology theories represented by the global Thom spectra
of Section \ref{sec:global Thom}.
We recall various facts about equivariant bordism groups of smooth
compact $G$-manifolds in some detail, highlighting the global perspective.
The main result of this section is Theorem \ref{thm:TP is iso} 
which says that when $G$ is isomorphic to a product of
a finite group and a torus, then the Thom-Pontryagin map
is an isomorphism from $G$-equivariant bordism to 
$G$-equivariant $\bmO$-homology.
This result is usually credited to Wasserman,
because it can be derived from his equivariant transversality
theorem \cite[Thm.\,3.11]{wasserman}.
Theorem \ref{thm:stable TP} gives a localized version
of this result: the Thom-Pontryagin map
is an isomorphism from {\em stable} equivariant bordism to 
$\bmO[1/\tau]$-theory, without any restriction on the compact Lie group.
Given that $\bmO[1/\tau]$-theory is isomorphic to $\bMO$-theory
(by Corollary \ref{cor:MO is localized mO}),
this is equivalent to a result of Br{\"o}cker and Hook \cite[Thm.\,4.1]{broecker-hook}
that identifies stable equivariant bordism with equivariant $\bMO$-homology.

Section \ref{sec:connective global K} 
discusses the ultra-commutative ring spectrum $\bku$,
the {\em connective global $K$-theory spectrum}, 
see Construction \ref{con:connective global K-theory}.
The degree zero equivariant cohomology theory represented by $\bku$
tries hard to be equivariant $K$-theory (see Theorem \ref{thm:K_G_is_ku} for the
precise statement),
and there is a ring homomorphism $\bRU(G)\to\pi_0^G(\bku)$
from the complex representation ring that is 
an isomorphism whenever $G$ is finite, and almost a morphism
of global power functors, compare Theorem \ref{thm:R(G) to pi ku}.
We also introduce the Bott class in the group $\pi_2^e(\bku)$ 
(see Construction \ref{con:Bott class}) and the more general
equivariant Bott classes associated to  $G$-$\Spin^c$-representations
(see Construction \ref{con:equivariant Bott class}).

The final Section \ref{sec:global K} reviews $\bKU$, the
{\em periodic global $K$-theory} spectrum, see Construction \ref{con:global KU}.
This is an ultra-commutative ring spectrum whose $G$-homotopy type realizes 
$G$-equivariant periodic $K$-theory, in the sense that
the equivariant cohomology represented by $\bKU$ on finite $G$-CW-complexes is
isomorphic to equivariant $K$-theory, see Corollary \ref{cor:global_K_represents}.
Moreover, $\upi_0(\bKU)$ is isomorphic, as a global power functor, 
to the complex representation ring functor $\bRU$,
by Theorem \ref{thm:pi_0 KU is RU}.
Periodic global $K$-theory receives a morphism of ultra-commutative ring spectra
$j:\bku\to\bKU$ from connective global $K$-theory; 
the image of the Bott class under $j$ becomes a unit
(Theorem \ref{thm:Bott and inverse}),
justifying the adjective `periodic' for $\bKU$.
Our final example of an ultra-commutative ring spectrum is $\bkuc$,
the {\em global connective $K$-theory}
(see Construction \ref{con:global connective ku}),
a certain global refinement 
of Greenlees' equivariant connective $K$-theory \cite{greenlees-connective}.

\section{Global Thom spectra}
\label{sec:global Thom}

\index{subject}{global Thom spectrum|(}

In this section we discuss two different global forms 
of the Thom spectrum $M O$ that represents unoriented bordism,
namely the ultra-commutative Thom ring spectrum $\bMO$, 
and a variation $\bmO$ that is only $E_\infty$-commutative. 
Both Thom spectra are the homogeneous degree~0 summands in certain
$\mZ$-graded periodic extensions $\bMOP$ respectively $\bmOP$.
All four orthogonal spectra are Thom spectra over
certain orthogonal spaces defined in Section \ref{sec:forms of BO}.
The partners are easy to identify from the notation: 
the relevant orthogonal spaces either have
a $\mathbf B$ or $\mathbf b$ in their name,
and in the corresponding Thom spectrum this letter is replaced by 
an $\mathbf M$ respectively $\mathbf m$.

While the spectrum $\bmO$ has less structure than $\bMO$
(it is `only' $E_\infty$, not ultra-commutative),
it is closely related to geometry. Indeed, the Thom spectrum $\bmO$
is the natural target for the Thom-Pontryagin map
from geometric equivariant bordism.
The ultra-commutative ring spectrum $\bMO$ is 
the $\mR$-analog of equivariant homotopical bordism, 
due to tom\,Dieck \cite{tomDieck-bordism};
the unoriented version $\bMO$ is studied in detail in \cite{broecker-hook}.
The following commutative diagram gives a schematic overview of the relevant
equivariant homology theories:
\[ \xymatrix@C=10mm{ 
\mathcal N_*^G \ar[r]^-{\Theta^G} \ar[d] & 
\bmO_*^G \ar[d] \ar[r] &
 \bmOP_*^G \ar[d] \\
\mathfrak N_*^{G:S}\ar[r]_-{\Theta^G}^\iso &
\bMO_*^G  \ar[r] & \bMOP_*^G  } \]
The two theories in the middle column are the equivariant homology
theories represented by the orthogonal Thom spectra $\bmO$ and $\bMO$.
The vertical transformations are isomorphisms for $G=e$, but not in general.
In fact, when $G$ has more than one element, then $\pi_*^G(\bMO)$
has non-trivial elements in negative degrees, while $\bmO$ is globally connective.
The vertical transformations are localizations, i.e., they become
isomorphisms after inverting all inverse Thom classes,
compare Corollary \ref{cor:MO is localized mO}.
The two theories $\bmOP$ and $\bMOP$ in the last column are the periodic versions
of $\bmO$ and $\bMO$; each is a wedge of all suspensions and desuspensions
of the non-periodic version.

In the left column, $\mathcal N_*^G$ is geometrically defined equivariant bordism,
and $\mathfrak N_*^{G:S}$ is {\em stable} equivariant bordism,
a certain localization of $\mathcal N_*^G$.
So the two theories in the left column are not represented by
orthogonal spectra, but they are defined from bordism classes of $G$-manifolds; 
we will recall these geometric theories in Section \ref{sec:equivariant bordism}.
The transformations labeled $\Theta^G$ are the equivariant Thom-Pontryagin construction
and its `stabilization'. The upper Thom-Pontryagin
map $\Theta^G:\mathcal N_*^G\to\pi_*^G(\bmO)$ is an isomorphism whenever $G$ is isomorphic
to a product of a finite group and a torus; 
this result seems to have been folklore
at some point, and we give a proof in Theorem \ref{thm:TP is iso}. 
The upper Thom-Pontryagin map is not an isomorphism for more general
compact Lie groups; in fact the geometric bordism theory $\mathcal N_*^G$
cannot in general be represented by an orthogonal $G$-spectrum because
the Wirthm{\"u}ller map fails to be an isomorphism for those closed
subgroups $H$ of $G$ that act non-trivially on the tangent space $T_{e H} (G/H)$,
compare Remark \ref{rk:Wirthmuller fails}.
The stabilized Thom-Pontryagin
map $\Theta^G:\mathfrak N_*^{G:S}\to\pi_*^G(\bMO)$ is an isomorphism
in complete generality, by a theorem of Br{\"o}cker 
and Hook \cite[Thm.\,4.1]{broecker-hook}; we derive this fact in a different way,
see Remark \ref{rk:proof of Brocker Hook}.

The Thom spectrum $\bmO$ comes with a rank filtration, the subquotients
of which we identify in Theorem \ref{thm:rank mO triangles}.
We use the rank filtration to show that $\bmO$ is globally connective,
and we calculate the global functor $\upi_0(\bmO)$ in Theorem \ref{thm:pi_0 of mO}.

\begin{eg}\label{eg:defined MGr}
We start with the ultra-commutative ring spectrum $\bMGr$,\index{symbol}{$\bMGr$ - {Thom spectrum over the additive Grassmannian}}
the Thom spectrum over the additive Grassmannian $\bGr$,
discussed in Example \ref{eg:Gr additive}.\index{subject}{additive Grassmannian}
We recall that the value of $\bGr$ at an inner product space $V$ is
\[ \bGr(V)\ = \ {\coprod}_{n\geq 0}\, Gr_n(V) \ , \]
the disjoint union of all Grassmannians in $V$.
Over the space $\bGr(V)$ sits a tautological euclidean vector bundle 
(of non-constant rank) whose total space consists of pairs 
$(x,U)\in V\times\bGr(V)$ such that $x\in U$.
We define $\bMGr(V)$ as the Thom space of this tautological vector bundle, i.e., 
the one-point compactification of the total space.
The structure maps are given by
\[ \bO(V,W) \sm \bMGr(V) \ \to \ \bMGr(W)\ ,\
(w,\varphi) \sm  (x,U) 
 \ \longmapsto \  ( w +\varphi(x),\, \varphi^\perp \oplus \varphi(U) )\ , \]
where $\varphi^\perp=W-\varphi(V)$ is the orthogonal complement 
of the image of $\varphi:V\to W$.
Multiplication maps
are defined by direct sum, i.e., 
\begin{align*}
 \mu_{V,W}\ : \ \bMGr(V) \sm \bMGr(W) \ &\to \quad \bMGr(V\oplus W)   \\
(x,U)\sm (x',U')\quad &\longmapsto\ ( (x,x'),\, U\oplus U')\ .   
\end{align*}
Unit maps are defined by
\[ \eta(V)\ : \ S^V\ \to\  \bMGr(V) \ , \quad  v\ \longmapsto \ (v,V)\ .  \]
The multiplication maps are binatural, associative, commutative and unital,
and all this structure makes $\bMGr$ into an ultra-commutative ring spectrum.
The orthogonal spectrum $\bMGr$ is graded, with $k$-th homogeneous summand given by
\[  \bMGr^{[k]}(V)\ = \  T h(G r_{|V|+k}(V)) \ . \]
So $\bMGr$ is concentrated in non-positive gradings, i.e.,
$\bMGr^{[k]}$ is trivial for $k >0$.
The unit morphism $\eta:\mS\to\bMGr$ is an isomorphism onto the summand $\bMGr^{[0]}$.

Now we let $V$ be a representation of a compact Lie group $G$.
We define the {\em inverse Thom class}\index{subject}{inverse Thom class!in $\bMGr$} 
\begin{equation}  \label{eq:inverse_Thom_MGr}
 \tau_{G,V}\ \in \ \bMGr_0^G(S^V)   
\end{equation}
as the class represented by the $G$-map 
  \begin{align*}
t_{G,V}\ : \  S^V \ &\to \ T h(G r(V))\sm S^V \ = \ 
\bMGr(V)\sm S^V \ ,\quad
v\ \longmapsto \  (0, \{0\})\sm (-v)   \ .
  \end{align*}
If $V$ has dimension $m$, then the class $\tau_{G,V}$ 
has internal degree $-m$, i.e., it lies in the homogeneous summand $\bMGr^{[-m]}$.
The justification for the name `inverse Thom class' is that in
the theory $\bMOP$, the image of the class $\tau_{G,V}$
becomes invertible, and its inverse is the Thom class of $V$
(considered as $G$-vector bundle over a point). 
We explain this in more detail in Theorem \ref{thm:Thom iso for MOP} below.
So while the theory $\bMGr$ is not globally orientable,
and does not have Thom isomorphisms for equivariant bundles,
informally speaking the inverses of the prospective Thom classes 
are already present in $\bMGr$. 
\end{eg}

\begin{rk}[$\bMGr$ as a wedge of semifree spectra]
We recall from Construction \ref{con:free orthogonal}
that the semifree orthogonal spectrum generated
by the tautological $O(m)$-representation $\nu_m$ on $\mR^m$ is 
given by\index{subject}{semifree orthogonal spectrum}
\[ F_{O(m),\nu_m} \ = \ \bO(\nu_m,-)/O(m)\ .\]
We claim that the spectrum $\bMGr$ is isomorphic
to the wedge of these semifree orthogonal spectra. Indeed, the maps
\[ \bO(\nu_m,V)/O(m)\ \to \ \bMGr^{[-m]}(V) \ ,\quad
(v,\varphi)\cdot O(m)\ \longmapsto \ (v, \varphi^\perp ) \]
define an isomorphism of orthogonal spectra from $F_{O(m),\nu_m}$
to the homogeneous summand $\bMGr^{[-m]}$. 
This isomorphism takes the equivariant homotopy class
$a_{O(m),\nu_m}$ in $\pi_0^G(F_{O(m),\nu_m}\sm S^{\nu_m})$ defined in \eqref{eq:define_a_G,G}
to the inverse Thom class $\tau_{O(m),\nu_m}$ of the tautological $O(m)$-representation,
by direct inspection of the definitions.
\end{rk}

Now we suppose that $\varphi:V\to W$ is an isomorphism of
orthogonal $G$-representations.
Then $\varphi$ compactifies to a $G$-equivariant homeomorphism 
$S^\varphi:S^V\to S^W$ and hence induces an isomorphism 
\[ (S^\varphi)_* \ : \  \bMGr_0^G(S^V) \ \to \ \bMGr_0^G(S^W) \ .\]
The following properties of the inverse Thom classes $\tau_{G,V}$
are all straightforward from the definition.
The power operations that show up in part~(v) 
are defined in \eqref{eq:power_op_ring_spectra}; 
the norm maps of part~(vi) are introduced in Remark \ref{rk:power vs Tambara}.

\begin{prop}\label{prop:inverse thom properties}
The inverse Thom classes $\tau_{G,V}$ have the following properties.  
\begin{enumerate}[\em (i)]
\item For every isomorphism $\varphi:V\to W$ of orthogonal $G$-representations,
the induced isomorphism $(S^\varphi)_*$ 
takes the class $\tau_{G,V}$ to the class $\tau_{G,W}$.
\item The class $\tau_{G,0}$ of the trivial 0-dimensional $G$-representation
is the multiplicative unit $1$ in $\bMGr_0^G(S^0)=\pi_0^G(\bMGr)$.
\item For every continuous homomorphism $\alpha:K\to G$ 
of compact Lie groups the relation
\[ \alpha^*(\tau_{G,V}) \ = \ \tau_{K,\alpha^*V}\]
holds in  $\bMGr_0^K(S^{\alpha^*(V)})$. 
\item For all orthogonal $G$-representations $V$ and $W$
the relation
\[ \tau_{G,V}\cdot \tau_{G,W} \ = \ \tau_{G,V\oplus W}\]
holds in  $\bMGr_0^G(S^{V\oplus W})$.
\item For all orthogonal $G$-representations $V$ and all $k\geq 1$
the relation
\[ P^k(\tau_{G,V})\ = \ \tau_{\Sigma_k\wr G,V^k}\]
holds in  $\bMGr_0^{\Sigma_k\wr G}(S^{V^k})$.
\item For every closed subgroup $H$ of $G$ of finite index 
and all orthogonal $H$-representations $W$ the relation
\[ N_H^G(\tau_{H,W})\ = \ \tau_{G,\Ind_H^G W}\]
holds in  $\bMGr_0^G(S^{\Ind_H^G W})$.
\end{enumerate}
\end{prop}

The next proposition shows that multiplication by the inverse Thom class
$\tau_{G,V}$ is realized by a certain morphism
of orthogonal $G$-spectra $j^V_\bMGr:\bMGr\to\sh^V\bMGr$. 
The value at an inner product space $U$ is the map
\begin{align}\label{eq:define j^V for MGr}
 j_\bMGr^V(U)\ : \ \bMGr(U)\ &\to \ \bMGr(U\oplus V)\ = \ (\sh^V\bMGr)(U)\\
(x,L)\quad &\longmapsto \ ((x,0), L\oplus 0)  \ . \nonumber
\end{align}
If $V$ has dimension $m$, then $j^V_\bMGr$ is homogeneous of degree $-m$
in terms of the internal grading of $\bMGr$, i.e.,
$j^V_\bMGr$ takes the wedge summand $\bMGr^{[k]}$ to the summand $\sh^V\bMGr^{[k-m]}$.  
The $\upi_*$-isomorphism $\lambda^V_X:X\sm S^V\to \sh^V X$ 
was defined in \eqref{eq:defn lambda_n}.

\begin{prop}\label{prop:j^V is inverse Thom MGr} 
Let $V$ be a representation of a compact Lie group $G$.
The composite
\[ \pi_k^G(\bMGr\sm A)\ \xra{-\cdot \tau_{G,V}} \ \pi_k^G(\bMGr\sm A\sm S^V) 
\ \xra{(\lambda^V_{\bMGr\sm A})_*} \ \pi_k^G(\sh^V \bMGr\sm A) \]
coincides with the effect of the morphism $j_\bMGr^V\sm A:\bMGr\sm A\to \sh^V\bMGr\sm A$.
\end{prop}
\begin{proof}
In \eqref{eq:shift2suspension} 
we defined an isomorphism $\psi^V_X:\pi_k^G(\sh^V X)\to\pi_k^G(X\sm S^V)$,
natural for morphisms of orthogonal $G$-spectra $X$, essentially given by
smashing with the identity of $S^V$.
We also defined $\varepsilon_V:\pi_k^G(X\sm S^V)\to\pi_k^G(X\sm S^V)$
as the effect of the antipodal map of $S^V$.
Now we observe that  multiplication by the class $\tau_{G,V}$
factors as the composite
\begin{align*}
   \pi_k^G(\bMGr\sm A)\ &\xra{(j^V_\bMGr\sm A)_*} 
\ \pi_k^G(\sh^V \bMGr\sm A) \\ 
&\xra{\ \psi^V_{\bMGr\sm A}\ } \ \pi_k^G(\bMGr\sm A\sm S^V) \
\xra{\ \varepsilon_V\ }\ \pi_k^G(\bMGr\sm A\sm S^V) \ .
\end{align*}
Besides the definitions,
this uses that the map $j^V_\bMGr(U)\sm S^V$ equals the composite
\begin{align*}
   \bMGr(U)\sm S^V\ \xra{\bMGr(U)\sm t_{G,V}} \  \bMGr(U)&\sm\bMGr(V)\sm S^V\\ 
&\xra{\mu_{U,V}\sm S^{-\Id}} \  \bMGr(U\oplus V)\sm S^V\ ,
\end{align*}
where $t_{G,V}$ is the defining representative of the class $\tau_{G,V}$.
The claim then follows from the fact,
established in Proposition \ref{prop:lambda upi_* isos}~(i),
that $(\lambda^V_{\bMGr\sm A})_*$ is inverse to $\varepsilon_V\circ\psi^V_{\bMGr\sm A}$.
\end{proof}

\begin{eg}\label{eg:MO}
We define two ultra-commutative ring spectra $\bMO$ 
and $\bMOP$. The latter is a periodic version of the former, 
and the former is the homogeneous degree~0 summand with respect
to a natural $\mZ$-grading of the latter.
Non-equivariantly, $\bMO$ is a version of the unoriented Thom spectrum $M O$,
and it is a global refinement of
equivariant homotopical bordism, due to tom\,Dieck \cite{tomDieck-bordism};
tom\,Dieck first considered the unitary version in \cite{tomDieck-bordism},
and the paper \cite{broecker-hook} by Br{\"o}cker and Hook
studies the orthogonal version.

The spectrum $\bMOP$ is a Thom spectrum over the orthogonal space $\bBOP$
discussed in Example \ref{eg:BOP}.
We recall that the value of $\bBOP$ at an inner product space $V$ is
\[ \bBOP(V)\ = \ {\coprod}_{n\geq 0}\, Gr_n(V^2) \ , \]
the disjoint union of all Grassmannians in $V^2$.
Over the space $\bBOP(V)$ sits a tautological euclidean vector bundle 
(again of non-constant rank) with total space consisting of pairs 
$(x,U)\in V^2\times\bBOP(V)$ such that $x\in U$.
We define $\bMOP(V)$\index{symbol}{$\bMOP$ - {periodic global Thom spectrum}}
as the Thom space of this tautological vector bundle, i.e., the
one-point compactification of the total space.
The structure maps are given by
\begin{align*}
  \bO(V,W) \sm \bMOP(V) \ &\to \qquad \bMOP(W)\\
(w,\varphi) \sm  (x,U)  \quad &\longmapsto \  ((w,0)+\bBOP(\varphi)(x),\,\bBOP(\varphi)(U))\ . 
\end{align*}
Multiplication maps
\[  \mu_{V,W}\ : \ \bMOP(V) \sm \bMOP(W) \ \to \ \bMOP(V\oplus W)   \]
are defined by sending $(x,U)\sm (x',U')$ to 
$(\kappa_{V,W}(x,x'),\kappa_{V,W}(U\oplus U'))$ 
where $\kappa_{V,W}:V^2\oplus W^2\iso (V\oplus W)^2$ is the preferred isometry 
defined by
\[ \kappa_{V,W}((v,v'),(w,w'))\ = \  ((v,w),(v',w')) \ . \]
Unit maps are defined by
\[ S^V\ \to\  \bMOP(V) \ , \quad  v\ \longmapsto \ ((v,0),\,V\oplus 0)\ .  \]
The multiplication maps are binatural, associative, commutative and unital,
and all this structure makes $\bMOP$ into an ultra-commutative ring spectrum.

The orthogonal space $\bBOP$ is $\mZ$-graded, with $k$-th homogeneous summand
\[  \bBOP^{[k]}(V)\ = \  Gr_{|V|+k}(V^2) \ . \]
The spectrum $\bMOP$ inherits a $\mZ$-grading,
where the summand $\bMOP^{[k]}(V)$ of degree $k$ 
is the Thom space of the tautological $(|V|+k)$-plane bundle over $\bBOP^{[k]}(V)$;
then $\bMOP(V)$ is the one-point union of the
Thom spaces  $\bMOP^{[k]}(V)$ for $-|V|\leq k\leq |V|$, and thus
\begin{equation}\label{eq:MOP decomposition}
 \bMOP \ = \ {\bigvee}_{k\in\mZ}\ \bMOP^{[k]}
\end{equation}
as orthogonal spectra.

We define $\bMO=\bMOP^{[0]}$ as the homogeneous degree zero wedge summand
of $\bMOP$; this is then an ultra-commutative ring spectrum in its own right.
Explicitly, $\bMO(V)$ is the Thom space of the tautological $|V|$-plane
bundle over $G r_{|V|}(V^2)$.\index{symbol}{$\bMO$ - {global Thom spectrum}}
\end{eg}

\begin{rk}
Certain variations of the construction of $\bMO$ and $\bMOP$ are possible, and
have been used at other places in the literature.
Indeed, if $U$ is any euclidean vector space, finite or infinite dimensional, 
and $u\in S(U)$ a unit vector,
we obtain an ultra-commutative ring spectrum $\bMO_{U,u}$ in exactly the same way
as above, with value at $V$ given by 
the Thom space of the tautological vector bundle over $Gr_{|V|}(V\tensor U)$.
The chosen vector $u$ enters in the definition of the unit and structure maps.
For $U=\mR^2$ and $u=(1,0)$, the construction specializes to $\bMO$ as above.

If the dimension of $U$ is at least~2, then we always get the same global homotopy type.
Indeed: any linear isometric embedding $\psi:U\to U'$ such that $\psi(u)=u'$
induces a morphism of ultra-commutative ring spectra $\psi_*:\bMO_{U,u}\to\bMO_{U',u'}$.
If the dimension of $U$ is at least~2, this morphism is a global equivalence. 
\end{rk}

Since $\bMOP$ is an ultra-commutative ring spectrum, 
the equivariant homotopy groups $\upi_0(\bMOP)$ form a global power functor.
The global power functor $\upi_0(\bMOP)$ is an interesting algebraic structure,
but a complete algebraic description does not seem to be known.
Since $2=0$ in $\pi_0^e(\bMOP)$,
the global power functor $\upi_0(\bMOP)$ takes values in $\mF_2$-vector spaces.
In $\upi_0(\bMOP)$, the stronger relation $\tr_e^{C_2}(1)=0$ holds,
compare Theorem \ref{thm:pi_0 of mO} below.

The orthogonal spectrum underlying $\bMOP$ comes with a $\mZ$-grading, 
i.e., a wedge decomposition \eqref{eq:MOP decomposition}
into summands $\bMOP^{[k]}$. The geometric splitting induces a direct sum decomposition
of $\pi_0^G(\bMOP)$ for every compact Lie group $G$ 
and makes it into a commutative $\mZ$-graded ring.
The $m$-th power operation takes the summand $\bMOP^{[k]}$ to the summand $\bMOP^{[m k]}$.

\medskip

We move on to explain the periodicity of the ultra-commutative ring spectrum $\bMOP$. 
We let $t\in\pi_{-1}(\bMOP^{[-1]})$
be the class represented by the point
\begin{equation} \label{eq:define_t_MOP}
(0,\{0\}) \ \in \ T h(G r_0(\mR^2)) \ = \  \bMOP^{[-1]}(\mR) \ .\end{equation}
We let $\sigma\in\pi_1(\bMOP^{[1]})$
be the class represented by the map
\begin{equation} \label{eq:define_sigma_MOP}
S^2 \ \to \ T h(G r_2(\mR^2)) \ = \  \bMOP^{[1]}(\mR) \ ,\quad
x\ \longmapsto \ (x,\mR^2)\ .\end{equation}
As the next proposition shows, $\bMOP$ is periodic in the sense
that $t$ is a unit in the graded ring $\pi_*^e(\bMOP)$, with inverse $\sigma$.

The orthogonal spectrum $\bMOP$ has an even stronger kind of
`$R O(G)$-graded' periodicity.
We define the {\em inverse Thom class}\index{subject}{inverse Thom class!in $\bMOP$} 
\begin{equation}  \label{eq:inverse_Thom_MOP}
 \tau_{G,V}\ \in \ \bMOP_0^G(S^V)   
\end{equation}
as the class represented by the $G$-map 
  \begin{align*}
t_{G,V}\ : \  S^V \ &\to \ T h(G r(V^2))\sm S^V \ = \ \bMOP(V)\sm S^V \\
v\quad &\longmapsto \  ( (0,0), 0\oplus 0)\sm (-v)   \ .
  \end{align*}
Here we abuse notation by denoting the inverse Thom class in $\bMOP$-theory
by the same symbol as the  inverse Thom class in $\bMGr$-theory
defined in \eqref{eq:inverse_Thom_MGr}. The justification
for this abuse is that the homomorphism $c:\bMGr\to\bMOP$
introduced in \eqref{eq:MGr_to_bMOP} below takes one inverse Thom class
to the other.
If $V$ has dimension $m$, then the class $\tau_{G,V}$ 
has internal degree $-m$, i.e., it lies in the homogeneous summand $\bMOP^{[-m]}$.
The justification for the name `inverse Thom class' is that
it is inverse to the Thom class $\sigma_{G,V}$ in $\bMOP^0_G(S^V)$, 
defined in \eqref{eq:define_sigma_G,V for MOP} below.

The periodicity class $t$ of \eqref{eq:define_t_MOP} is essentially
the inverse Thom class of the 1-dimensional representation of the trivial group.
More precisely, $t\sm S^1 = \tau_{e,\mR}$, i.e., the suspension isomorphism
\[ -\sm S^1 \ : \  \pi_{-1}^e(\bMOP)\ \xra{\ \iso\ } \ \pi_0^e(\bMOP\sm S^1)\ = \ 
\bMOP_0^e(S^1)   \]
takes the periodicity class $t$ to the inverse Thom class $\tau_{e,\mR}$.
Indeed, the suspension of the defining representative \eqref{eq:define_t_MOP}
for $t$ differs from the  defining representative
for $\tau_{e,\mR}$ by the inversion map $-\Id:S^1\to S^1$.
So $t\sm S^1=-\tau_{e,\mR}$; however, since $2=0$ in $\pi_0^e(\bMOP)$,
this yields the claim.

The next proposition shows that multiplication by any inverse Thom class
is invertible in $\bMOP$-theory. Moreover, multiplication by $\tau_{G,V}$
is realized, in a certain precise way,
by a periodicity morphism of orthogonal $G$-spectra $j_\bMOP^V:\bMOP\to\sh^V\bMOP$:
the value at an inner product space $U$ is the map
\begin{align}\label{eq:define j^V for MOP}
 j_\bMOP^V(U)\ : \ \bMOP(U)\ &\to \ \bMOP(U\oplus V)\ = \ (\sh^V\bMOP)(U)\nonumber\\
(x,L)\quad &\longmapsto \quad (i(x),i(L))   
\end{align}
induced by the linear isometric embedding
$i:U\oplus U\to U\oplus V\oplus U\oplus V$ with $i(u,u')=(u,0,u',0)$.
The morphism $j^V_\bMOP$ is even a homomorphism of left $\bMOP$-module spectra.
If $V$ has dimension $m$, then $j^V_\bMOP$ is homogeneous of degree $-m$ in terms
of the $\mZ$-grading of $\bMOP$, i.e., it restricts to a morphism
of orthogonal $G$-spectra
\[ j_\bMOP^V\ : \ \bMOP^{[k+m]}\ \to\ \sh^V\bMOP^{[k]}  \]
where $k$ is any integer.
In the special case $V=\mR$ with trivial $G$-action, the map
\[  j \ = \ j_\bMOP^\mR \ : \ \bMOP\ \to \ \sh \bMOP \]
is a morphism of orthogonal spectra (with trivial $G$-action).
As usual we denote by $p_G:G\to e$ the unique group homomorphism
to the trivial group, and $p_G^*$ is the associated inflation homomorphism.

\begin{theorem}\label{thm:j^V is upi_*-iso}
  \begin{enumerate}[\em (i)]
  \item
    The relation $t\cdot \sigma =1$ holds in $\pi_0^e(\bMOP)$.
    For every compact Lie group $G$, every based $G$-space $A$
    and all $k\in\mZ$, the maps 
    \[  \bMOP_{k+1}^G(A)  \ \xra{\ \cdot p_G^*(t)\  }    \bMOP_k^G(A) 
    \text{\quad and\quad}
    \bMOP_k^G(A)   \ \xra{\ \cdot p_G^*(\sigma)\  }  \bMOP_{k+1}^G(A) \]
    are mutually inverse isomorphisms.
  \item 
    For every representation $V$ of a compact Lie group $G$, the morphism
    \[ j^V_\bMOP\ : \ \bMOP \ \to \ \sh^V \bMOP \]
    is a $\upi_*$-isomorphism of orthogonal $G$-spectra.
    In particular, the morphism $j=j_\bMOP^\mR :\bMOP \to \sh \bMOP$
    is a global equivalence.
  \item
    For every based $G$-space $A$, the composite
    \[ \pi_0^G(\bMOP\sm A)\ \xra{-\cdot \tau_{G,V}} \ \pi_0^G(\bMOP\sm A\sm S^V) 
    \ \xra{(\lambda^V_{\bMOP\sm A})_*} \ \pi_0^G(\sh^V \bMOP\sm A) \]
    coincides with the effect of the morphism 
    $j^V_\bMOP\sm A:\bMOP\sm A\to \sh^V\bMOP\sm A$.
    In particular, exterior multiplication by the inverse Thom class $\tau_{G,V}$ 
    is invertible in equivariant $\bMOP$-homology.
  \end{enumerate}
\end{theorem}
\begin{proof}
(i) The class $t\cdot \sigma$ is represented by the composite
\[
  S^2 \ \xra{x\mapsto (0,\{0\})\sm (x,\mR^2)}\ \bMOP^{[-1]}(\mR)\sm\bMOP^{[1]}(\mR)
\ \xra{\mu_{\mR,\mR}}\ \bMOP^{[0]}(\mR\oplus\mR)
\]
where the first map is the smash product of the defining representatives
for $t$ and $\sigma$. Expanding the definition of $\mu_{\mR.\mR}$
identifies this composite as the map
\[   S^2 \ \to\  \bMOP^{[0]}(\mR\oplus\mR) \ , \quad
x\ \longmapsto \ ( (\mR\oplus\tau_{\mR,\mR}\oplus\mR)(0,0,x), 0\oplus\mR\oplus 0\oplus\mR )\ .\]
This differs from the representative of the unit $1\in\pi_0^G(\bMOP)$
by the action of the linear isometry
\[ \mR^4 \ \to \ \mR^4 \ , \quad (a,b,c,d)\ \longmapsto \ 
(b,d,c,a) \ . \]
This isometry has determinant~1, so we conclude that $t\cdot\sigma=1$
in $\pi_0^e(\bMOP)$.

(ii) We show first that $j_\bMOP^V$ induces an isomorphism on 0-th 
$G$-equivariant homotopy groups. To this end we define a map
\[ \Phi\ : \ \pi_0^G(\sh^V \bMOP)\ \to \ \pi_0^G(\bMOP) \]
in the opposite direction as follows.
We define a $G$-map 
\[   s_{G,V}\ : \ S^{V\oplus V} \ \to \ T h(G r(V\oplus V)) \ = \ \bMOP(V) 
\text{\, by\, } s_{G,V}(v,w)\ =\ ( (v,w), V\oplus V)  \ .   \]
If $f:S^U\to \bMOP(U\oplus V)=(\sh^V\bMOP)(U)$ 
represents a class in $\pi_0^G(\sh^V \bMOP)$,
then we define $\Phi[f]$ as the class of the composite
\begin{align*}
S^{U\oplus V\oplus V}  \ \xra{f\sm s_{G,V}} \
\bMOP(U\oplus V)\sm \bMOP(V)  \ \xra{\mu_{U\oplus V,V}} \
\bMOP(U\oplus V\oplus V)  \ .
\end{align*}
This recipe is compatible with stabilization, so $\Phi$ is indeed well-defined.

The composite
\[ S^{V\oplus V} \ \xra{\ s_{G,V}\ }\ \bMOP(V) \ \xra{j_\bMOP^V(V)}\ \bMOP(V\oplus V)\]
is $G$-equivariantly homotopic to the unit 
map $\eta(V\oplus V):S^{V\oplus V}\to \bMOP(V\oplus V)$.
So 
\begin{align*}
  (j_\bMOP^V)_*(\Phi[f])\ &= \ [j_\bMOP^V(U\oplus V\oplus V)\circ \mu_{U\oplus V,V}\circ(f\sm s_{G,V})] \\
&= \ [\mu_{U\oplus V,V\oplus V}\circ (\bMOP(U\oplus V)\sm j_\bMOP^V(V))\circ (f\sm s_{G,V})] \\
&= \ [\mu_{U\oplus V,V\oplus V}\circ (f\sm ( j_\bMOP^V(V)\circ s_{G,V}))] \\
&= \ [\mu_{U\oplus V,V\oplus V}\circ (f\sm  \eta(V\oplus V))] \ = \ [f]\ .
\end{align*}
The second equation is the fact that $j_\bMOP^V$ is a homomorphism of left $\bMOP$-modules.
This proves that $(j_\bMOP^V)_*\circ \Phi$ is the identity.

On the other hand, the composite
\[  \bMOP(U)\sm S^{V\oplus V} \ \xra{j_\bMOP^V(U)\sm s_{G,V}}\ 
 \bMOP(U\oplus V)\sm \bMOP(V) \ \xra{\mu_{U\oplus V,V}}\ \bMOP(U\oplus V\oplus V) \]
agrees with the opposite structure map $\sigma^{\op}_{U,V\oplus V}$
of the spectrum $\bMOP$.
So if $\varphi:S^U\to\bMOP(U)$ represents a class in $\pi_0^G(\bMOP)$, then
\begin{align*}
   \Phi( (j_\bMOP^V)_*[\varphi])\ &= \ [\mu_{U\oplus V,V}\circ(j_\bMOP^V(U)\circ \varphi)\sm s_{G,V}] \\
&= \ [\mu_{U\oplus V,V}\circ(j_\bMOP^V(U)\sm s_{G,V})\circ (\varphi\sm S^{V\oplus V})] \\
&= \ [\sigma^{\op}_{U,V\oplus V}\circ(\varphi\sm S^{V\oplus V})]\ = \ [\varphi]\ .
\end{align*}
This proves that $\Phi\circ (j_\bMOP^V)_*$ is the identity,
and thus completes the proof that $\pi_0^G(j_\bMOP^V)$ is an isomorphism.
Since $j_\bMOP^V$ is a homomorphism of left $\bMOP$-modules,
its effect on homotopy groups is $\pi_*^G(\bMOP)$-linear.
In particular, it commutes with the action of the element $p_G^*(t)$,
where $t\in\pi_{-1}^e(\bMOP)$ is the periodicity element defined 
in \eqref{eq:define_t_MOP}. Since $p_G^*(t)$ is invertible by part~(i)
and $\pi_0^G(j_\bMOP^V)$ is an isomorphism, the map $\pi_k^G(j_\bMOP^V)$ is 
then an isomorphism for every integer $k$.
Applying the above to a closed subgroup $H$ of $G$ and the underlying $H$-representation
of $V$ shows that $j_\bMOP^V$ induces isomorphisms of equivariant stable homotopy groups
for all closed subgroups of $G$. So $j_\bMOP^V$ is a $\upi_*$-isomorphism.

(iii) The proof of the first claim proceeds in the same way as its analog for
$\bMGr$ in Proposition \ref{prop:j^V is inverse Thom MGr}.
Since the morphisms $j^V_{\bMOP}$ and $\lambda^V_{\bMOP\sm A}$
are both $\upi_*$-isomorphisms of orthogonal $G$-spectra,
they induce isomorphisms on $\pi_0^G$.
So exterior multiplication by the inverse Thom class $\tau_{G,V}$ 
is invertible.
\end{proof}

\Danger The orthogonal spectra $\bMGr$ and $\bMOP$ admit `shift morphisms'
$j_\bMGr^V: \bMGr \to \sh^V \bMGr$ and $j_\bMOP^V : \bMOP\to \sh^V\bMOP$
defined in \eqref{eq:define j^V for MGr} respectively \eqref{eq:define j^V for MOP}.
However, the morphism $j^V_\bMGr$  is {\em not} a $\upi_*$-isomorphism,
whereas the morphism $j^V_\bMOP$ {\em is} a $\upi_*$-isomorphism,
by Theorem \ref{thm:j^V is upi_*-iso}~(ii).
This is a reflection of the fact that the inverse Thom classes $\tau_{G,V}$
are {\em not} invertible in equivariant $\bMGr$-homology,
whereas their $\bMOP$-counterparts are.

\begin{construction}[Thom classes for representations]
The Thom spectrum $\bMOP$ comes with distinguished Thom classes for representations.
We let $V$ be a representation of a compact Lie group $G$.
We consider the $G$-map
\[  s_{G,V}\ : \ S^{V\oplus V} \ \to \ T h(G r(V\oplus V)) \ = \ \bMOP(V) \ ,
\quad (v,w)\ \longmapsto \ ( (v,w), V\oplus V)  \ .  \]
If $V$ has dimension $m$, then $s_{G,V}$ is a homeomorphism onto the 
homogeneous summand $\bMOP^{[m]}(V)$.
The adjoint $S^V\to\map_*(S^V,\bMOP(V))$ of $s_{G,V}$ represents the {\em Thom class}
\index{symbol}{$\sigma_{G,V}$ - {Thom class in $\bMOP^0_G(S^V)$}}\index{subject}{Thom class!in $\bMOP$}
\begin{equation}\label{eq:define_sigma_G,V for MOP}
 \sigma_{G,V}\ \in \ \bMOP^0_G(S^V) \ = \ \pi_0^G(\map_*(S^V,\bMOP ))
\end{equation}
in the  $G$-equivariant $\bMOP$-cohomology of $S^V$.
\end{construction}

The following theorem is a special case of a Thom isomorphism 
in the equivariant cohomology theory represented by $\bMOP$.
It also makes precise in which way the inverse Thom class $\tau_{G,V}$
is inverse to the Thom class $\sigma_{G,V}$.
This relation between $\tau_{G,V}$ and $\sigma_{G,V}$
is the ultimate justification for naming $\tau_{G,V}$ the `inverse Thom class'.

\begin{theorem}\label{thm:Thom iso for MOP}\index{subject}{Thom isomorphism!for $\bMOP$}
Let $V$ be a representation of a compact Lie group $G$.
Then the composite
\[ \bMOP^0_G(S^V)\ = \ \pi_0^G(\Omega^V\bMOP )\ \xra{\ \cdot\tau_{G,V}}\
\pi_0^G((\Omega^V\bMOP)\sm S^V)\ \xra{(\epsilon^V_{\bMOP})_*}\
\pi_0^G(\bMOP)  \]
is inverse to multiplication by $\sigma_{G,V}$,
where $\epsilon^V_\bMOP:(\Omega^V\bMOP)\sm S^V\to\bMOP$ is the evaluation morphism.
In particular, $\bMOP^0_G(S^V)$ is a free module 
of rank~1 over the ring $\pi_0^G(\bMOP)$,
and the Thom class $\sigma_{G,V}$ is a generator.
\end{theorem}
\begin{proof}
We consider a $G$-map $f:S^U\to\Omega^V \bMOP(U)$
that represents a class in $\bMOP^0_G(S^V)$.
Then $\epsilon^V_\bMOP([f]\cdot\tau_{G,V})$ is represented by the following composite:
\begin{align*}
S^{U\oplus V} \ \xra{f\sm t_{G,V}}\ &(\Omega^V \bMOP(U))\sm \bMOP(V)\sm S^V \\ 
\xra{\text{evaluate}}\  &\bMOP(U)\sm \bMOP(V) \ 
\xra{\ \mu^{\bMOP}_{U,V}\ }\  \bMOP(U \oplus V) \ ,
\end{align*}
where $t_{G,V}:S^V\to\bMOP(V)\sm S^V$ 
is the defining representative for the class $\tau_{G,V}$
from \eqref{eq:inverse_Thom_MOP}. 
If we let $f=s_{G,V}^\sharp:S^V\to\Omega^V \bMOP(V)$ be adjoint to
the defining representative for $\sigma_{G,V}$, then the composite comes out as the map
\[ S^{V\oplus V} \ \to\  \bMOP(V \oplus V) \ , \quad
(v,w)\ \longmapsto \ ( (v,0,-w,0),V\oplus 0\oplus V\oplus 0)\ .\]
This composite is equivariantly homotopic to the map
\[ (v,w)\ \longmapsto \ ( (v,w,0,0) ,\  V\oplus V\oplus 0\oplus 0) \]
which represents the multiplicative unit. So we have shown that
\[ \epsilon^V_{\bMOP}(\sigma_{G,V}\cdot\tau_{G,V})\ =\ 1 \ .\]
All maps in sight are left $\pi_0^G(\bMOP)$-linear, so we deduce that
\[ \epsilon^V_{\bMOP}(x\cdot \sigma_{G,V}\cdot\tau_{G,V})\ =\ 
x\cdot \epsilon^V_{\bMOP}(\sigma_{G,V}\cdot\tau_{G,V})\ =\ x \]
for every class $x\in\pi_0^G(\bMOP)$.
On the other hand, the composite
\[  \pi_0^G(\Omega^V\bMOP )\ \xra{\ \cdot\tau_{G,V}}\
\pi_0^G((\Omega^V\bMOP)\sm S^V)\ \xra{(\lambda^V_{\Omega^V\bMOP})_*}\
\pi_0^G(\sh^V \Omega^V\bMOP) \]
is the effect of the morphism $\Omega^V j^V_\bMOP:\Omega^V\bMOP \to\sh^V\Omega^V\bMOP$, 
by the same reasoning as in Theorem \ref{thm:j^V is upi_*-iso}~(iii). 
Since $j^V_\bMOP$ is a $\upi_*$-isomorphism, so is $\Omega^V j^V_\bMOP$ by 
Proposition \ref{prop:map(A,-) preserves global}~(ii).
On the other hand, $\lambda^V_{\Omega^V\bMOP}$ is a  $\upi_*$-isomorphism 
by Proposition \ref{prop:lambda upi_* isos}~(ii).
So multiplication by the class $\tau_{G,V}$ is an isomorphism.
The morphism $\epsilon^V_\bMOP$ is a $\upi_*$-isomorphism,
also by Proposition \ref{prop:lambda upi_* isos}~(ii).
So the composite $(\epsilon^V_\bMOP)_*\circ(-\cdot\tau_{G,V})$
is bijective. Since it is also left inverse to multiplication by $\sigma_{G,V}$,
this proves the first claim.
\end{proof}

\begin{construction}[Thom classes for equivariant vector bundles]
The Thom spectrum $\bMOP$ comes with a distinguished orientation, 
given by Thom classes for equivariant vector bundles.
These Thom classes generalize the classes $\sigma_{G,V}$ 
defined in \eqref{eq:define_sigma_G,V for MOP},
when we view a $G$-representation as a $G$-vector bundle over a one-point
$G$-space.

We recall the definition of the Thom classes.
Given a compact Lie group $G$, an orthogonal $G$-spectrum $X$ and a
compact based $G$-space $B$, we define the {\em $G$-equivariant $X$-cohomology group}
of $B$ as
\[ X^0_G(B) \ = \ \pi_0^G(\map_*(B, X) ) \ .\]
We let $\xi:E\to B$ be a $G$-equivariant vector bundle.
The bundle has a classifying $G$-map $\psi:B\to G r(V)$
for some $G$-representation $V$, i.e., such that $\xi$ is isomorphic
to the pullback of the tautological $G$-vector bundle over the Grassmannian.
We let $\bar \psi:E\to V$ be a map that covers $\psi$, i.e., $\bar\psi$ is
fiberwise linear and satisfies $\bar \psi(e)\in \psi(\xi(e))$. 
We define a based $G$-map
\[ S^V \sm T h(\xi)\ \to \ \bMOP(V) = T h(G r(V^2))\text{\, by\,}
v\sm e \ \longmapsto \ ((v,\bar \psi(e)), V\oplus \psi(\xi(e)))\ .\]
We denote the equivariant cohomology class represented by the adjoint 
$S^V\to \map_*(T h(\xi), \bMOP(V))$ of this map by
\[ \sigma_G(\xi)\ \in \ \bMOP_G^0(T h(\xi)) \]
and refer to it as the 
{\em Thom class}\index{subject}{Thom class!in $\bMOP$!of a $G$-vector bundle} 
of the $G$-vector bundle $\xi$.
If the bundle $\xi$ has constant rank $m$,
then the image of the map lies in the wedge summand $\bMOP^{[m]}$.
It is straightforward to see that the Thom classes just defined
are natural for pullback of bundles, compatible with
restriction along continuous homomorphisms, and 
the Thom class of an exterior product of bundles is
the exterior product of the Thom classes.

The Thom diagonal of the $G$-vector bundle $\xi:E\to B$ is the map
\[ \Delta\ : \ 
T h(\xi)\ \to \ B_+\sm T h(\xi) \ , \quad e\ \longmapsto \ \xi(e)\sm e \ . \]
This diagonal induces an action map of equivariant cohomology groups
\[ \bMOP^0_G(B_+) \times \bMOP^0_G(T h(\xi)) 
\ \xra{\ \times \ } \  \bMOP^0_G(B_+\sm T h(\xi)) \ \xra{\ \Delta^*\ } \
 \bMOP^0_G(T h(\xi)) \ .\]
Since the diagonal is coassociative and counital, the action map
makes the group $\bMOP^0_G(T h(\xi))$ into a left module over the
commutative ring $\bMOP^0_G(B_+)$.
The Thom isomorphism then says that whenever $B$ admits the structure
of a finite $G$-CW-complex, then $\bMOP^0_G(T h(\xi))$ is a free module 
of rank~1 over $\bMOP^0_G(B_+)$, with the Thom class $\sigma_G(\xi)$
as a generator.
Theorem \ref{thm:Thom iso for MOP} is the special case when the base
consists of a single point.
We refrain from giving the proof of the general Thom isomorphism 
in equivariant $\bMOP$-theory.
\end{construction}

\begin{construction}[Euler classes]
We let $G$ be a compact Lie group and $V$ an $m$-dimensional $G$-representation.
As usual, the Thom classes $\sigma_{G,V}$ give rise to
Euler classes\index{symbol}{$e(V)$ - {Euler class in $\pi_0^G(\bMOP)$}}\index{subject}{Euler class!in $\bMO$}
by `restriction to the zero section', i.e.,
\[ e(V)\ = \ i^*(\sigma_{G,V})\ \in \ \bMOP^0_G(S^0)\ = \ \pi_0^G(\bMOP) \ . \]
Here $i:S^0\to S^V$ is the inclusion
of the fixed points~0 and $\infty$, and
$i^*:\bMOP^0_G(S^V)\to\bMOP^0_G(S^0)$ the induced map on equivariant cohomology groups.
The Euler class is thus represented by the based $G$-map
\[ S^V \ \to \  \bMOP(V) \ ,\quad v\ \longmapsto \ ( (v,0),\, V\oplus V) \ .  \]
Since the Thom class lives in the homogeneous summand $\bMOP^{[m]}$,
so does the Euler class.
If $V$ has nonzero $G$-fixed points, then the inclusion $i:S^0\to S^V$
is $G$-equivariantly null-homotopic, so $e(V)=0$ whenever $V^G\ne 0$.
\end{construction}

\begin{rk}[Shifted Thom and Euler classes in $\bMO$]\index{subject}{shifted inverse Thom class!in $\bMO$} 
The author thinks that the periodic theory $\bMOP$ is the most
natural home for the Thom classes, the Euler classes and the inverse Thom classes, 
but the more traditional place
to host them is the degree~0 wedge summand $\bMO=\bMOP^{(0]}$.
Indeed, for an $m$-dimensional $G$-representation $V$,
the Thom class $\sigma_{G,V}$ and the Euler class $e(V)$ lie 
in the homogeneous summand $\bMOP^{[m]}$,
and we can use the periodicity of $\bMOP$ to move the classes into $\bMO$,
at the expense of shifting their degrees by $m$.
In other words, by multiplying by a suitable power of the
periodicity class $t\in \pi_{-1}^e(\bMOP^{[-1]})$, we define
\[ \bar \sigma_{G,V}\ = \ p_G^*(t^m)\cdot \sigma_{G,V}  \in \ \bMO^m_G(S^V) 
\text{\quad and\quad}
\bar e(V)\ = \ p_G^*(t^m)\cdot e(V)  \in \ \pi_{-m}^G(\bMO) \ , \]
where $p_G:G\to e$ is the unique group homomorphism.
The Thom isomorphism theorem for $\bMO$ then says that $\bMO^*_G(S^V)$ 
is a free graded module of rank~1 over the graded ring $\pi_*^G(\bMO)$,
and the shifted Thom class $\bar\sigma_{G,V}$ is a generator.
This version of the Thom isomorphism follows directly 
from Theorem \ref{thm:Thom iso for MOP} because 
$\bMOP$ is globally equivalent to the wedge of all suspensions and desuspensions
of $\bMO$. 
More precisely, the maps 
\[ \bigoplus_{n\in\mZ} \, \pi_n^G(\bMO)\ \to \ \pi_0^G(\bMOP)  
\text{\quad and\quad}
 \bigoplus_{n\in\mZ} \, \bMO^{-n}_G(S^V)\ \to \ \bMOP^0_G(S^V) \ , \]
given on the $n$-th summand by multiplication by $p_G^*(t^n)$,
are isomorphisms; moreover the latter isomorphism 
takes the shifted Thom class $\bar\sigma_{G,V}$
to the original Thom class $\sigma_{G,V}$.

Similarly, we can use the periodicity of $\bMOP$ to move the inverse Thom class 
$\tau_{G,V}\in \bMOP_0^G(S^V)$ into $\bMO$,
at the expense of shifting it from degree~0 to homological degree $m$.
The class $\tau_{G,V}$ lies in the homogeneous summand $\bMOP^{[-m]}$,
so by multiplying by a suitable power of the
periodicity class $\sigma\in \pi_1^e(\bMOP^{[1]})$ we define
\begin{equation}\label{eq:MO inverse Thom class}
 \bar \tau_{G,V}\ = \ p_G^*(\sigma^m)\cdot \tau_{G,V}  \in \ \bMO_m^G(S^V) \ .
\end{equation}
Theorem \ref{thm:j^V is upi_*-iso}~(iii) then implies that
for every based $G$-space $A$, the map
\[-\cdot \bar\tau_{G,V}\ : \ \pi_*^G(\bMO\sm A)\ \to \ \pi_{*+m}^G(\bMO\sm A\sm S^V) \]
is an isomorphism.
\end{rk}

Our next task is to show that the Thom spectrum $\bMOP$ is
a localization of the Thom spectrum $\bMGr$, 
obtained by formally inverting all inverse Thom classes.
This result can be viewed as a `thomification' of the fact
that the morphism of ultra-commutative monoids $i:\bGr\to\bBOP$
induces a group completion of abelian monoids
$[A,i]^G : [A,\bGr]^G\to[A,\bBOP]^G$
for every compact Lie group $G$ and every $G$-space $A$, 
compare Proposition \ref{prop:i is group completion}.

The ultra-commutative ring spectra $\bMGr$ and $\bMOP$ are connected by a homomorphism
\begin{equation}  \label{eq:MGr_to_bMOP}
a \ : \  \bMGr\ \to\ \bMOP
\end{equation}
whose value at an inner product space $V$ is
\[ a(V)\ : \ \bMGr(V)\ \to \ \bMOP(V) \ ,\quad  
(x,L)\ \longmapsto \ ((x,0), L\oplus 0)\ . \]
For varying $V$, these maps form a morphism of 
$\mZ$-graded ultra-commutative ring spectra.
The morphism \eqref{eq:MGr_to_bMOP}
induces natural transformations of equivariant homology theories
\[ (a\sm A)_* \ : \ \bMGr_*^G(A)\ \to\ \bMOP_*^G(A)  \]
for all compact Lie groups $G$ and all based $G$-spaces $A$.
We observe that 
\[ (a\sm S^V)_*(\tau_{G,V})\ = \ \tau_{G,V} \ ,\]
i.e., the morphism $a$ takes the $\bMGr$-inverse Thom class
to the $\bMOP$-inverse Thom class with the same name.
This relation is immediate from the explicit
representatives of the inverse Thom classes in \eqref{eq:inverse_Thom_MGr}
respectively \eqref{eq:inverse_Thom_MOP}.

We define a localized version of equivariant $\bMGr$-homology by
\[  \bMGr^G_k(A)[1/\tau]\ = \ 
\colim_{V\in s(\Uc_G)}\,  \bMGr_k^G(A\sm S^V) \ ;\]
for $V\subset W$, the structure map in the colimit system is the multiplication
\[ \bMGr_k^G(A\sm S^V) \ \xra{-\cdot\tau_{G,W-V}}\ 
\bMGr_k^G(A\sm S^V\sm S^{W-V}) \ \iso\ \bMGr_k^G(A\sm S^W)  \ .\]
In equivariant $\bMOP$-theory, the inverse Thom classes become
invertible by Theorem \ref{thm:j^V is upi_*-iso}~(iii).
So for a $G$-representation $V$ we can consider the map
\[ \bMGr_k^G(A\sm S^V) \ \xra{\ (a\sm A\sm S^V)_*\ }\ 
\bMOP^G_k(A\sm S^V)\ \xra{ (-\cdot \tau_{G,V})^{-1}} \ \bMOP^G_k(A)\ . \]
By the multiplicativity of $a$,
these maps are compatible as $V$ varies over the poset $s(\Uc_G)$,
so they assemble into a natural transformation
\[ a^\sharp \ : \ \bMGr_k^G( A)[1/\tau] \ =\ 
\colim_{V\in s(\Uc_G)}\,  \bMGr_k^G(A\sm S^V) \ \to\ \bMOP^G_k(A)\ . \]

\begin{theorem}\label{thm:MOP is localized MGr}
For every compact Lie group $G$, every based $G$-space $A$ 
and every integer $k$ the map
\[ a^\sharp \ : \ \bMGr_k^G(A)[1/\tau] \ \to \ \bMOP_k^G(A) \]
is an isomorphism.
\end{theorem}
\begin{proof}
To simplify the exposition we prove the claim for $k=0$ only, 
the argument for a general integer being essentially the same.
Alternatively, we can observe that source and target of $a^\sharp$
are periodic, so it suffices to establish bijectivity in a single dimension.
We show two separate statements that amount to the injectivity
respectively surjectivity of the map $a^\sharp$.

(a) 
We show that for every class $x$ in the kernel 
of the map $(a\sm A)_*:\bMGr_0^G(A) \to \bMOP^G_0(A)$,
there is a $G$-representation $V$ such that $x\cdot\tau_{G,V}=0$.
Indeed, we can represent any such class $x$ by a based $G$-map
$f:S^V\to \bMGr(V)\sm A$, for some $G$-representation $V$,
such that the composite
\[ S^V \ \xra{\ f\ }\ \bMGr(V)\sm A\ \xra{a(V)\sm A}\ \bMOP(V)\sm A\]
is equivariantly null-homotopic. 
In \eqref{eq:define j^V for mOP} we defined a morphism of orthogonal $G$-spectra
$j^V_\bMGr:\bMGr\to\sh^V\bMGr$. 
We observe that $\bMOP(V)=\bMGr(V\oplus V)$ and $a(V)=j^V_\bMGr(V)$.
So we can apply Proposition \ref{prop:j^V is inverse Thom MGr} and conclude that
\begin{align*}
(\lambda^V_{\bMGr\sm A})_*(x\cdot\tau_{G,V})\ &= \ 
 (j_\bMGr^V\sm A)_*(x) \\ 
&= \ [ (j_\bMGr^V(V)\sm A)\circ f]\ = \ [ (a(V)\sm A)\circ f]\ = \ 0 \ .
\end{align*}
Since $\lambda^V_{\bMGr\sm A}:\bMGr\sm A\sm S^V\to\sh^V\bMGr\sm A$ is
a $\upi_*$-isomorphism (by Proposition \ref{prop:lambda upi_* isos}~(ii)),
this implies the desired relation $x\cdot\tau_{G,V}=0$.

(b) 
We show that for every class $y$ in $\bMOP_0^G(A)$,
there is a $G$-representation $V$ and a class $x$ in $\bMGr_0^G(A\sm S^V)$
such that $y\cdot \tau_{G,V}= (a\sm A\sm S^V)_*(x)$.
To this end we represent $y$ by a based $G$-map
$f:S^V\to\bMOP(V)\sm A$. Because $\bMOP(V) = \bMGr(V^2)=(\sh^V \bMGr)(V)$,
the map $f$ also defines a class in $\pi_0^G(\sh^V\bMGr\sm A)$.
Since $\lambda^V_{\bMGr\sm A}:\bMGr\sm A\sm S^V\to\sh^V\bMGr\sm A$ is
a $\upi_*$-isomorphism by Proposition \ref{prop:lambda upi_* isos}~(ii),
there is a unique class $x\in\bMGr_0^G(A\sm S^V)$ such that
\[ (\lambda^V_{\bMGr\sm A})_*(x)\ = \ [f]\ . \]
On the other hand, the map $a(V^2):\bMGr(V^2)\to \bMOP(V^2)$
is {\em equal} to the map $j^V_{\bMOP}(V):\bMOP(V)\to\bMOP(V^2)$.
So
\begin{align*}
(\lambda^V_{\bMOP\sm A})_*(  (a\sm A\sm S^V)_*(x))\ &= \  
  (\sh^V a\sm A)_*((\lambda^V_{\bMGr\sm A})_*(x))\\ 
&= \  (\sh^V a\sm A)_*[f]\ = \  [(a(V^2)\sm A)\circ  f]\\
&= \  [(j^V_{\bMOP}(V)\sm A)\circ  f]\ = \ (j^V_\bMOP\sm A)_* (y) \\ 
&= \ (\lambda^V_{\bMOP\sm A})_* (y\cdot\tau_{G,V}) \ .
\end{align*}
The sixth equation is Theorem \ref{thm:j^V is upi_*-iso}~(iii),
the others are either definitions or naturality properties.
Since $\lambda^V_{\bMOP\sm A}$ is a $\upi_*$-isomorphism,
we can conclude that 
\[ (a\sm A\sm S^V)_*(x)\ = \  y\cdot\tau_{G,V}\ .\qedhere \]
\end{proof}

\begin{eg}[The global Thom spectra $\bmO$ and $\bmOP$]\label{eg:geometric global bordism}
We define two $E_\infty$-orthogonal ring spectra $\bmO$\index{symbol}{$\bmO$ - {global Thom spectrum}} 
and $\bmOP$,\index{symbol}{$\bmOP$ - {periodic global Thom spectrum}} 
the Thom spectra over the orthogonal spaces $\bbO$ and $\bbOP$
defined in Examples \ref{eg:define bO}
respectively \ref{eg:define bOP}.
The spectrum $\bmOP$ is a periodic version of $\bmO$,
and conversely $\bmO$ is the homogeneous degree~0 summand with respect
to a certain $\mZ$-grading of $\bmOP$.
Non-equivariantly, $\bmO$ is another version of the unoriented Thom spectrum $M O$.
The equivariant homology theory represented by $\bmO$ is the natural target
of the equivariant Thom-Pontryagin map from equivariant bordism, 
and that map is trying hard to be an isomorphism, see Theorem \ref{thm:TP is iso}  below.

We recall that the value of $\bbOP$ at an inner product space $V$ is
\[ \bbOP(V)\ = \  {\coprod}_{n\geq 0}\, G r_n(V\oplus\mR^\infty) \ , \]
the disjoint union of all the Grassmannians in $V\oplus\mR^\infty$.
The map $\bbOP(\varphi):\bbOP(V)\to\bbOP(W)$ induced by a 
linear isometric embedding $\varphi:V\to W$ is defined as
\[ \bbOP(\varphi)(L) \ = \ 
(\varphi\oplus\mR^\infty)(L) \ + \ ((W-\varphi(V))\oplus 0) \ .\]
In other words: we apply the linear isometric 
embedding $\varphi\oplus\mR^\infty:V\oplus\mR^\infty\to W\oplus\mR^\infty$ 
to the subspace $L$ and add the orthogonal complement of the image of $\varphi$
(sitting in the first summand of $W\oplus\mR^\infty$).

Over the space $\bbOP(V)$ sits a tautological euclidean vector bundle 
(again of non-constant rank);
the total space of this bundle consist of pairs 
$(x,U)\in (V\oplus\mR^\infty)\times \bbO(V)$  such that $x\in U$.
We define $\bmOP(V)$ as the Thom space of this tautological vector bundle.
The structure maps are given by
\begin{align*}
  \bO(V,W)\sm  \bmOP(V) \ &\to \quad \bmOP(W)\\
(w,\varphi)\sm (x,U) \quad &\longmapsto \ ((w,0)+\bbOP(\varphi)(x), \bbOP(\varphi)(U))\ . 
\end{align*}
As we explained in Remark \ref{rk:bO versus BO},
the orthogonal spaces $\bbO$ and $\bbOP$ have natural $E_\infty$-structures.
Correspondingly, the orthogonal spectra $\bmO$ and $\bmOP$
have natural $E_\infty$-structures,
by which we mean an action of the linear isometries operad. 
This multiplication is, however, {\em not} ultra-commutative.
Multiplication maps
\[ \mu_{V,W}\ : \ \bL((\mR^\infty)^2,\mR^\infty)_+ \sm \bmOP(V) \sm \bmOP(W) \ \to \ 
\bmOP(V\oplus W)   
 \]
are defined by sending $\psi\sm (x,U)\sm (x',U')$ 
to $(\psi_\sharp(x,x'),\psi_\sharp(U\oplus U'))$,
where $\psi_\sharp$ is the linear isometric embedding
\begin{equation}  \label{eq:define_psi_sharp}
 V\oplus\mR^\infty\oplus W\oplus\mR^\infty\ \to \ 
V\oplus W\oplus \mR^\infty\  ,\quad \psi_\sharp(v,y,w,z)\ = \ (v,w,\psi(y,z))\ .  
\end{equation}
Unit maps are defined by
\[ S^V\ \to\  \bmOP(V) \ , \quad  v\ \longmapsto \ ((v,0),\, V\oplus 0)\ .  \]
All this structure makes $\bmOP$ into an $E_\infty$-orthogonal ring spectrum.

The orthogonal spectrum $\bmOP$ is $\mZ$-graded, where
the summand $\bmOP^{[k]}(V)$ of degree $k$ is defined as 
the Thom space of the tautological
$(|V|+k)$-plane bundle over $\bbOP^{[k]}(V)=G r_{|V|+k}(V\oplus\mR^\infty)$;
then $\bmOP(V)$ is the one-point union of the
Thom spaces  $\bmOP^{[k]}(V)$ for $|V|+k\geq 0$.
So we have a wedge decomposition
\[ \bmOP \ = \ {\bigvee}_{k\in\mZ}\ \bmOP^{[k]} \]
as orthogonal spectra.
We define $\bmO=\bmOP^{[0]}$ as the homogeneous wedge summand of degree~0.
\end{eg}

In the rest of this section we will also use products
on the equivariant homology theories represented by the Thom spectra $\bmO$
and $\bmOP$.
As we just explained,
the spectrum $\bmOP$ comes with an $E_\infty$-multiplication
which is, however, neither strictly associative, nor strictly commutative.
So we briefly explain how to define these products.
We choose a linear isometric embedding $\psi:\mR^\infty\oplus\mR^\infty\to\mR^\infty$.
We define continuous maps
\begin{align}\label{eq:psi_for_mu_bmO}
  \psi_{V,W}\ : \ \bmOP(V) \sm \bmOP(W) \ &\to \quad  \bmOP(V\oplus W)  \text{\qquad by} \\
\psi_{V,W}((x,U),(x',U'))\ &\ = \ (\psi_\sharp(x,x'),\psi_\sharp(U\oplus U')) \ ,  
\nonumber
\end{align}
where $\psi_\sharp$ was defined in \eqref{eq:define_psi_sharp}.
These maps form a bimorphism, which corresponds to a morphism of orthogonal spectra
\[ \psi_\star\ : \  \bmOP\sm\bmOP\ \to \ \bmOP \]
by the universal property of the smash product.
Since the space $\bL((\mR^\infty)^2,\mR^\infty)$
of linear isometric embeddings is contractible,
the morphism $\psi_\star$ is independent up to homotopy of the choice of $\psi$.
Even though the multiplication map $\psi_\star$ is neither associative
nor commutative, the contractibility of the space $\bL((\mR^\infty)^3,\mR^\infty)$ 
implies that the square
\[ \xymatrix@C=18mm{ 
\bmOP\sm\bmOP\sm\bmOP \ar[r]^-{\bmOP\sm \psi_\star}\ar[d]_{\psi_\star\sm\bmOP}&
\bmOP\sm\bmOP \ar[d]^-{\psi_\star}\\
\bmOP\sm\bmOP \ar[r]_-{\psi_\star} & \bmOP } \]
commutes up to homotopy,
and the contractibility of the space $\bL((\mR^\infty)^2,\mR^\infty)$
implies that the composite
\[ \bmOP\sm\bmOP \ \xra{\tau_{\bmOP,\bmOP}} \ \bmOP\sm\bmOP \ \xra{\ \psi_\star\ } \ \bmOP \]
is homotopic to $\psi_\star$.
So whenever we pass to induced maps on equivariant homotopy groups,
an $E_\infty$-multiplication is as good as a strictly
associative and commutative multiplication. 
However, an $E_\infty$-multiplication does {\em not} entitle us to power operations.

Given a compact Lie group $G$ and based $G$-spaces $A$ and $B$,
we define a multiplication
\[ \cdot \ : \ \bmOP^G_k(A)\times \bmOP_l^G(B) \ \to \ \bmOP_{k+l}^G(A\sm B) \]
as the composite
\begin{align*}
\pi_k^G(\bmOP\sm A)\times \pi_l^G(\bmOP\sm B) \ &\xra[\eqref{eq:def_dot}]{\ \times \ } \ 
\pi_{k+l}^G(\bmOP\sm A\sm \bmOP\sm B) \\ 
&\xra{(\text{twist})_*} \  \pi_{k+l}^G(\bmOP\sm\bmOP\sm A\sm B) \\
&\xra{(\psi_\star\sm A\sm B)_*} \  \pi_{k+l}^G(\bmOP\sm A\sm B) \ .
\end{align*}

We move on to explain the periodicity property of $\bmOP$. 
As the theories $\bMGr$ and $\bMOP$, the theory $\bmOP$
also has its own inverse Thom classes and shift morphisms.
We define the {\em inverse Thom class}\index{subject}{inverse Thom class!in $\bmOP$} 
\begin{equation}  \label{eq:inverse_Thom_mOP}
 \tau_{G,V}\ \in \ \bmOP_0^G(S^V)   
\end{equation}
as the class represented by the $G$-map 
\[ S^V \ \to \ T h(G r(V\oplus\mR^\infty))\sm S^V \ = \ 
\bmOP(V)\sm S^V \ ,\quad
v\ \longmapsto \  ( (0,0), 0\oplus 0)\sm (-v)   \ . \]
Here we abuse notation one more time and also denote the inverse Thom class 
in $\bmOP$-theory
by the same symbol as its counterparts in $\bMGr$-theory
and $\bMOP$-theory defined in \eqref{eq:inverse_Thom_MGr}
respectively \eqref{eq:inverse_Thom_MOP}.
The justification for this abuse is that the inverse Thom classes
match up under certain homomorphisms relating $\bMGr$, $\bMOP$ and $\bmOP$.
As usual, if $V$ has dimension $m$, then the class $\tau_{G,V}$ 
lies in the homogeneous summand $\bMOP^{[-m]}$.
A shift morphism of orthogonal $G$-spectra $j_\bmOP^V:\bmOP\to\sh^V\bmOP$ 
is defined as for $\bMGr$ and $\bMOP$: the value at an inner product space $U$ is the map
\begin{align}\label{eq:define j^V for mOP}
 j_\bmOP^V(U)\ : \ \bmOP(U)\ &\to \ \bmOP(U\oplus V)\ = \ (\sh^V\bmOP)(U)\nonumber\\
(x,L)\quad &\longmapsto \quad (i(x),i(L))   
\end{align}
induced by the linear isometric embedding
$i:U\oplus \mR^\infty\to U\oplus V\oplus \mR^\infty$ with $i(u,x)=(u,0,x)$.
If $V$ has dimension $m$, then $j^V_\bmOP$ is homogeneous of degree $-m$.
In the special case $V=\mR$ with trivial $G$-action, the map
\[ j \ = \ j_\bmOP^\mR \ : \ \bmOP \ \to \ \sh \bmOP  \]
is a morphism of orthogonal spectra (with trivial $G$-action).

On the level of homotopy groups, the periodicity is realized by
multiplication with a periodicity element $t\in\pi_{-1}^e(\bmOP^{[-1]})$
represented by the point
\begin{equation}\label{eq:t in mOP}
 (0,\{0\})\ \in \ T h(G r_0(\mR\oplus\mR^\infty))\ = \ \bmOP^{[-1]}(\mR) \ .  
\end{equation}

\begin{prop}\label{thm:j^V for mOP}
  \begin{enumerate}[\em (i)]
  \item 
    The morphism $j=j_\bmOP^\mR :\bmOP \to \sh \bmOP$
    is a homotopy equivalence of orthogonal spectra, hence a global equivalence.
    For every compact Lie group $G$ and every based $G$-space $A$, the induced map
    \[ (j\sm A)_*\ : \ \pi_*^G(\bmOP\sm A) \ \to\ \pi_*^G(\sh \bmOP\sm A) \]
    is an isomorphism.
  \item    
    For every compact Lie group $G$, every $G$-representation $V$ and every
    based $G$-space $A$, the composite
    \[ \pi_0^G(\bmOP\sm A)\ \xra{-\cdot \tau_{G,V}} \ \pi_0^G(\bmOP\sm A\sm S^V) 
    \ \xra{(\lambda^V_{\bmOP\sm A})_*} \ \pi_0^G(\sh^V \bmOP\sm A) \]
    coincides with the effect of the morphism 
    $j^V_\bmOP\sm A:\bmOP\sm A\to \sh^V\bmOP\sm A$.
    In particular, exterior multiplication by the inverse Thom class $\tau_{G,\mR}$ 
    of the trivial 1-dimensional $G$-representation
    is invertible in equivariant $\bmOP$-homology.
  \item
    For every compact Lie group $G$, every based $G$-space $A$
    and every integer $k$ the multiplication map
    \[  \bmOP_{k+1}^G(A)  \ \xra{\ \cdot p_G^*(t)\  } \   \bmOP_k^G(A) \]
    is an isomorphism. In particular, the class $t\in\pi_{-1}^e(\bmOP)$
    is a unit in the graded homotopy ring $\pi_*^e(\bmOP)$.
  \end{enumerate}
\end{prop}
\begin{proof}
(i) The morphism $j$ is based on the 
linear isometric embedding $i:\mR^\infty\to \mR\oplus \mR^\infty$ 
defined by $i(x)=(0,x)$. This linear isometric embedding is homotopic, 
through linear isometric embeddings,
to the linear isometry $\mR^\infty\iso \mR\oplus \mR^\infty$ 
sending $(x_1,x_2,x_3,\dots)$ to $(x_1,(x_2,x_3,\dots))$.
Such a homotopy induces a homotopy from the morphism $j$ to an
isomorphism between $\bmOP$ and $\sh\bmOP$. So $j$ is homotopic to an isomorphism,
hence a homotopy equivalence.

The first statement in~(ii) is proved by the same argument as the analogous statement 
for $\bMGr$ in Proposition \ref{prop:j^V is inverse Thom MGr}.
For $V=\mR$, the trivial 1-dimensional $G$-representation,
the morphisms $j^\mR_{\bmOP}\sm A$ and $\lambda_{\bmOP\sm A}$
are both global equivalences by part~(i) 
respectively Proposition \ref{prop:global equiv preservation}~(i).
So they induce isomorphisms on $\pi_0^G$.
Hence exterior multiplication by the class $\tau_{G,\mR}$ is invertible.

(iii) The composite
\[ \bmOP_{k+1}^G(A)\ \xra{-\cdot p_G^*(t)} \
 \bmOP_k^G(A)\ \xra[\iso]{-\sm S^1} \ \bmOP_{k+1}^G(A\sm S^1) \]
differs from multiplication by the inverse Thom class $\tau_{G,\mR}$
by the effect of the involution $\bmOP\sm A\sm S^{-\Id}$ of $\bmOP\sm A\sm S^1$.
Multiplication by $\tau_{G,\mR}$ is an isomorphism by part~(ii),
and the suspension isomorphism is bijective by
Proposition \ref{prop:loop and suspension isomorphisms}.
So multiplication by $p_G^*(t)$ is bijective as well.
\end{proof}

In the earlier Theorem \ref{thm:MOP is localized MGr}
we showed that equivariant $\bMOP$-homology is obtained from 
equivariant $\bMGr$-homology by inverting all inverse Thom classes.
Now we add $\bmOP$ to this picture, which turns out to be
an intermediate localization. 
As we will now explain, $\bmOP$-theory is obtained from
$\bMGr$-theory by inverting the inverse Thom classes of all {\em trivial} representation.
Then $\bMOP$-theory is obtained from
$\bmOP$-theory by inverting the remaining inverse Thom classes,
i.e., the ones of non-trivial representations.
Informally speaking, the first localization turns $\bMGr$-theory 
into a theory that is periodic in the $\mZ$-graded sense;
the second localization then turns $\bmOP$-theory 
into a theory that is periodic in the $R O(G)$-graded sense.
Schematically:
\[\xymatrix@C=17mm{ 
 \bMGr_*^G(A) \quad  \ar@{=>}[r]^-{\text{invert $\tau_{G,\mR}$}} & \quad \bmOP_*^G(A)
\quad \ar@{=>}[r]^-{\text{invert all $\tau_{G,V}$}} & \quad \bMOP_*^G(A) } \]

The ring spectra $\bMGr$ and $\bmOP$ are connected by a homomorphism
\begin{equation}  \label{eq:MGr_to_bmOP}
b \ : \  \bMGr\ \to\ \bmOP
\end{equation}
whose value at an inner product space $V$ is
\[ b(V)\ : \ \bMGr(V)\ \to \ \bmOP(V) \ ,\quad  
(x,L)\ \longmapsto \ ((x,0), L\oplus 0)\ . \]
For varying $V$, these maps form a morphism of 
$\mZ$-graded orthogonal $E_\infty$-ring spectra.
The morphism $b:\bMGr\to\bmOP$
induces natural transformations of equivariant homology theories
\[ (b\sm A)_* \ : \ \bMGr_*^G(A)\ \to\ \bmOP_*^G(A)  \]
for all compact Lie groups $G$ and all based $G$-spaces $A$.
We observe that 
\[ (b\sm S^V)_*(\tau_{G,V})\ = \ \tau_{G,V} \ ,\]
i.e., the morphism $b$ takes the $\bMGr$-inverse Thom class
to the $\bmOP$-inverse Thom class with the same name.
This relation is again immediate from the explicit
representatives of the inverse Thom classes in \eqref{eq:inverse_Thom_MGr}
respectively \eqref{eq:inverse_Thom_mOP}.

We define a localized version of equivariant $\bMGr$-homology by
\[  \bMGr^G_k(A)[\tau_{G,\mR}^{-1}]\ = \ 
\colim_{n\geq 0}\,  \bMGr_k^G(A\sm S^n) \ ,\]
the colimit of the sequence
\[  \bMGr_k^G(A) \ \xra{-\cdot \tau_{G,\mR}}\ 
 \bMGr_k^G(A\sm S^1) \ \xra{-\cdot \tau_{G,\mR}}\
 \bMGr_k^G(A\sm S^2) \ \xra{-\cdot \tau_{G,\mR}}\ \dots
\]
along multiplication by the inverse Thom class $\tau_{G,\mR}\in \bMGr_0^G(S^1)$.
In equivariant $\bmOP$-theory, the class $\tau_{G,\mR}$
becomes invertible by Theorem \ref{thm:j^V for mOP}~(ii).
So we can consider the maps
\[ \bMGr_k^G(A\sm S^n) \ \xra{\ (b\sm A\sm S^n)_*\ }\ 
\bmOP^G_k(A\sm S^n)\ \xra{ (-\cdot \tau_{G,\mR}^n)^{-1}} \ \bmOP^G_k(A)\ . \]
By the multiplicativity of $b$, these maps are compatible,
so they assemble into a natural transformation
$b^\sharp : \bMGr_k^G(A)[\tau_{G,\mR}^{-1}]\to\bmOP^G_k(A)$.

\begin{theorem}\label{thm:mOP is localized MGr}
For every compact Lie group $G$, every based $G$-space $A$ and every integer $k$
the map 
\[ b^\sharp \ : \ \bMGr_k^G(A)[\tau_{G,\mR}^{-1}] \ \to \ \bmOP_k^G(A) \]
is an isomorphism.
\end{theorem}
\begin{proof}
The `standard' linear isometric embedding
\[ \mR^n \ \to \ \mR^\infty\ , \quad 
(x_1,\dots,x_n)\ \longmapsto \ (x_1,\dots,x_n,0,0,\dots) \]
induces a continuous map
\[ \psi^n(V)\ : \ (\sh^n \bMGr)(V)\ = \ T h(G r(V\oplus\mR^n))
\ \to \ T h(G r(V\oplus\mR^\infty))\ = \ \bmOP(V)\ . \]
As $V$ varies, these maps form a morphism of orthogonal spectra
$\psi^n:\sh^n\bMGr\to\bmOP$. Since $\psi^n= \psi^{n+1}\circ(\sh^n j_{\bMGr})$, 
the morphisms $\psi^n$ are compatible with the sequence of morphisms
\[ \bMGr \ \xra{\ j_\bMGr\ }\ 
\sh \bMGr \ \xra{ \sh j_\bMGr }\ 
\sh^2 \bMGr \ \xra{ \sh^2 j_\bMGr }\ \cdots\  \ .\]
Moreover, the morphisms $\psi^n$ express $\bmOP$ as the colimit of this sequence.
So $\bmOP\sm A$ is the colimit 
of the sequence of orthogonal $G$-spectra $\sh^n\bMGr\sm A$.
The map $j_\bMGr(V)$ is an h-cofibration of based $O(V)$-spaces.
So if $G$ acts on $V$ by linear isometries,
then $j_\bMGr(V)\sm A$ is an h-cofibration of based $G$-spaces, hence a closed embedding
by Proposition \ref{prop:h-cof is closed embedding}.
The morphism $j_\bMGr\sm A$ and all its shifts are thus levelwise closed embeddings.

Proposition \ref{prop:sequential colimit closed embeddings}~(i)
shows that equivariant homotopy groups commute with sequential colimits
over closed embeddings; so the canonical map
\[ \colim_{n\geq 0}\, \pi_k^G(\sh^n\bMGr\sm A)\ \to \ \pi_k^G(\bmOP\sm A) \]
is an isomorphism.
The diagram
\[ \xymatrix@C=12mm{ 
\pi_k^G(\bMGr\sm A\sm S^{n-1})\ar[r]^-{\cdot\tau_{G,\mR}}\ar[d]_{\lambda^{n-1}_{\bMGr\sm A}}^\iso &
\pi_k^G(\bMGr\sm A\sm S^n)\ar[d]^{\lambda^n_{\bMGr\sm A}}_\iso 
\ar[r]^-{\pi_k^G(b\sm A\sm S^n)}&  \pi_k^G(\bmOP\sm A\sm S^n)\ar[d]_\iso^{(-\cdot\tau_{G,\mR}^n)^{-1}}\\
\pi_k^G(\sh^{n-1} \bMGr\sm A)\ar[r]_-{\pi_k^G(\sh^{n-1} j_\bMGr\sm A)} &
\pi_k^G(\sh^n\bMGr\sm A)\ar[r]_-{\pi_k^G(\psi^n\sm A)} &
\pi_k^G(\bmOP\sm A)
} \]
commutes by Proposition \ref{prop:j^V is inverse Thom MGr} and because
$(\sh^{n-1}j_\bMGr)\circ\lambda^{n-1}_{\bMGr} = \lambda^{n-1}_{\sh \bMGr}\circ (j_\bMGr\sm S^{n-1})$.
The three vertical maps are isomorphisms.
So $\pi_k^G(\bMOP)$ is also a colimit of the sequence of
multiplication maps by $\tau_{G.\mR}$, with respect to the maps
that define $b^\sharp$.
\end{proof}

There is one localization result left: it remains to exhibit 
$\bMOP$-theory as the localization of $\bmOP$-theory, 
by inverting the inverse Thom classes of arbitrary representations.
This is in fact a formal consequence of 
Theorem \ref{thm:MOP is localized MGr}
and Theorem \ref{thm:mOP is localized MGr} which exhibit both
$\bMOP_*^G(A)$ and $\bmOP_*^G(A)$ as localizations of $\bMGr_*^G(A)$,
the former being a more drastic localization than the latter.

We spell this out in more detail.
We define a localized version of equivariant $\bmOP$-homology by
\[  \bmOP^G_0(A)[1/\tau]\ = \ 
\colim_{V\in s(\Uc_G)}\,  \bmOP_0^G(A\sm S^V) \ ;\]
for $V\subset W$, the structure map in the colimit system
is the multiplication
\[ \bmOP_0^G(A\sm S^V) \ \xra{-\cdot\tau_{G,W-V}}\ 
\bmOP_0^G(A\sm S^V\sm S^{W-V}) \ \iso\ \bmOP_0^G(A\sm S^W)  \ .\]
In \eqref{eq:MGr_to_bMOP} we introduced a morphism 
of ultra-commutative ring spectra $a:\bMGr\to\bMOP$.
In \eqref{eq:MGr_to_bmOP} we introduced the morphism
of $E_\infty$-ring spectra $b:\bMGr\to\bmOP$.
These morphisms induce multiplicative natural transformations
\begin{align*}
   (a\sm A)_* \ &: \ \bMGr_*^G(A)\ \to\ \bMOP_*^G(A)  \text{\quad respectively}\\
 (b\sm A)_* \ &: \ \bMGr_*^G(A)\ \to\ \bmOP_*^G(A)  
\end{align*}
for all compact Lie groups $G$ and all based $G$-spaces $A$.
Moreover, the morphisms match up the inverse Thom classes
in the sense that
\[ (a\sm S^V)_*(\tau_{G,V})\ = \ \tau_{G,V} \text{\qquad and\qquad}
 (b\sm S^V)_*(\tau_{G,V})\ = \ \tau_{G,V} \ ;\]
indeed, this is the justification for our abuse of notation of
using the same name for the inverse Thom classes in $\bMGr$, $\bMOP$ and $\bmOP$.

Theorem \ref{thm:mOP is localized MGr} says that the map $(b\sm A)_*$
becomes an isomorphism after inverting the inverse Thom class $\tau_{G,\mR}$
of the trivial 1-dimensional $G$-represen\-tation.
So it also becomes an becomes after inverting the inverse Thom classes
of all representations, i.e., the induced transformation
\[ (b\sm A)_*[1/\tau] \ : \ \bMGr_*^G(A)[1/\tau]\ \to\ \bmOP_*^G(A)[1/\tau]  \]
is an isomorphism.
On the other hand, the transformation $(a\sm A)_*$ induces an
isomorphism $a^\sharp$ from $\bMGr_*^G(A)[1/\tau]$ to $\bMOP_*^G(A)$, by
Theorem \ref{thm:MOP is localized MGr}.
So combining these two theorems yields:

\begin{cor}\label{cor:MOP is localized mOP}
For every compact Lie group $G$, every based $G$-space $A$ 
and every integer $k$ the map
\[ a^\sharp\circ ((b\sm A)_*[1/\tau])^{-1}\ : \ \bmOP_k^G(A)[1/\tau] \ \to \ \bMOP_k^G(A) \]
is an isomorphism.
\end{cor}

While the authors thinks that the periodic theories $\bmOP$ and $\bMOP$ give the
most convenient formulation of the localization result, 
the more traditional formulation is in terms of the degree~0 summands
$\bmO=\bmOP^{[0]}$ and $\bMO=\bMOP^{[0]}$. So we also take the
time to reformulate Corollary \ref{cor:MOP is localized mOP}
in terms of $\bmO$ and $\bMO$.
Since the inverse Thom classes do {\em not} lie in the degree zero summands,
we instead use the shifted inverse Thom classes\index{subject}{shifted inverse Thom class!in $\bmO$}\index{subject}{inverse Thom class!shifted|see{shifted inverse Thom class}}  
\[  \bar \tau_{G,V}\ = \ p_G^*(\sigma^m)\cdot \tau_{G,V}  \in \ \bMO_m^G(S^V)  \]
introduced in \eqref{eq:MO inverse Thom class}, 
and its $\bmO$-analog
\begin{equation}\label{eq:shifted inverse_Thom_mO}
  \bar \tau_{G,V}\ = \ p_G^*(\sigma^m)\cdot \tau_{G,V}  \in \ \bmO_m^G(S^V)  \ .
\end{equation}
In the $\bMOP$-case, the periodicity class $\sigma\in\pi_1^e(\bMOP^{[1]})$
was defined in \eqref{eq:define_sigma_MOP}, and it is inverse to the class $t$. 
In the $\bmOP$-case, the periodicity class $\sigma\in\pi_1^e(\bmOP^{[1]})$
was not yet defined, and we take it to be the inverse
of the $\bmOP$ periodicity class $t$ from \eqref{eq:t in mOP}.
Both theories $\bmOP$ and $\bMOP$ are periodic (in the $\mZ$-graded sense),
i.e., the maps
\[ \bigoplus_{m\in\mZ} \, \bMO_m^G(A)\ \to \ \bMOP_0^G(A)  \text{\quad and\quad}
 \bigoplus_{m\in\mZ} \, \bmO_m^G(A)\ \to \ \bmOP_0^G(A)  \]
that multiply by the appropriate powers of the periodicity classes are isomorphisms.

We define a localized version of equivariant $\bmO$-homology by
\[  \bmO^G_k(A)[1/\bar\tau]\ = \ 
\colim_{V\in s(\Uc_G)}\,  \bmO_{k+|V|}^G(A\sm S^V) \ ;\]
for $V\subset W$, the structure map in the colimit system is the multiplication
\[ \bmO_{k+|V|}^G(A\sm S^V) \ \xra{-\cdot\bar\tau_{G,W-V}}\ 
\bmO_{k+|W|}^G(A\sm S^V\sm S^{W-V}) \ \iso\
\bmO_{k+|W|}^G(A\sm S^W)  \ .\]
The periodicity of $\bmOP$ is inherited by the localized theory, 
i.e., the upper horizontal map in the following commutative square
\[ \xymatrix{
\bigoplus_{m\in\mZ}  \bmO_m^G(A)[1/\bar\tau]\ar[r] 
\ar[d]_{a^\sharp\circ ((b\sm A)_*[1/\bar\tau])^{-1}}  &
\bmOP_0^G(A)[1/\tau] \ar[d]^{a^\sharp\circ ((b\sm A)_*[1/\tau])^{-1}} \\
\bigoplus_{m\in\mZ}  \bMO_m^G(A)\ar[r] &
\bMOP_0^G(A) } \]
is an isomorphism.
The lower horizontal map is an isomorphism by the periodicity of $\bMOP$, and
the right vertical map is an isomorphism by Corollary \ref{cor:MOP is localized mOP}.
So the left vertical map is an isomorphism. Since the left map
is homogeneous with respect to the $\mZ$-grading, it is an
isomorphism in every degree. So we conclude:

\begin{cor}\label{cor:MO is localized mO}
For every compact Lie group $G$, every based $G$-space $A$ 
and every integer $m$ the map
\[ a^\sharp\circ ( (b\sm A)_*[1/\bar\tau])^{-1}\ : \ 
\bmO_m^G(A)[1/\bar\tau] \ \to \ \bMO_m^G(A) \]
is an isomorphism.
\end{cor}

Now we investigate the global homotopy type of the Thom spectrum $\bmO$
in more detail; 
one reason for wanting to understand $\bmO$ better is the close connection
to equivariant bordism, compare Theorem \ref{thm:TP is iso}  below. 
Our main tool is the following `rank filtration';
there is an analog of the rank filtration for $\bmOP$, 
but we refrain from making it explicit.

\begin{construction}[Rank filtration of $\bmO$]\index{subject}{rank filtration! of $\bmO$}\label{con:mO_m}
The orthogonal spectrum $\bmO$ is a global Thom spectrum over
an orthogonal space $\bbO$;   
in Proposition \ref{prop:bbO as hocolim}
we identified $\bbO$ as a certain global homotopy colimit
of the global classifying spaces $B_{\gl} O(m)$. More precisely,
the filtration of $\mR^\infty$ by the subspaces $\mR^m$ induces a
filtration of $\bbO$ by orthogonal subspaces $\bbO_{(m)}$,
and $\bbO_{(m)}$ receives a global equivalence from $B_{\gl} O(m)$.
We now define and study the corresponding orthogonal Thom spectrum
$\bmO_{(m)}$ over $\bbO_{(m)}$, which turns out to be an $m$-fold suspension
of the orthogonal spectrum $M_{\gl} T(m)$, the global refinement
of the spectrum traditionally denoted $M T (m)$.

Construction \ref{con:free orthogonal} defines 
the semifree orthogonal spectrum $F_{G,V}$
generated by a $G$-representation $V$. We are interested in the
tautological $O(m)$-representation $\nu_m$ with underlying inner product space $\mR^m$,
and to simplify the notation we abbreviate the corresponding semifree spectrum to 
\[ F_m\ =\ F_{O(m),\nu_m}\ . \]
The shift functor $\sh^m=\sh^{\mR^m}$ 
by the inner product space $\mR^m$ was defined in Construction \ref{con:shift}.
We set
\[ \bmO_{(m)}\ = \ \sh^m F_m\ , \]
the $m$-th shift of $F_m$.\index{symbol}{$\bmO_{(m)}$ - {truncated global Thom spectrum}}  
Unpacking this definition reveals the value of $\bmO_{(m)}$
at an inner product space $V$ as the space
\[ \bmO_{(m)}(V)\ = \ \bO(\nu_m,V\oplus \mR^m)/O(m)\ . \]
To justify the notation $\bmO_{(m)}$ we clarify the connection
to the orthogonal space $\bbO_{(m)}$ defined in \eqref{eq:define_bmO_m}.
The value of $\bbO_{(m)}$ at $V$ is
\[ \bbO_{(m)}(V)\ = \ G r_{|V|}(V\oplus\mR^m) \ ,\]
the Grassmannian of $|V|$-planes in $V\oplus\mR^m$.
Over the space $\bbO_{(m)}(V)$ sits a tautological euclidean $|V|$-plane bundle, 
with total space consisting of pairs 
$(x,U)\in (V\oplus\mR^m)\times \bbO_{(m)}(V)$  such that $x\in U$.  
Passage to orthogonal complements provides a homeomorphism:
\begin{align*}
\bmO_{(m)}(V)\ = \  \bO(\nu_m,V\oplus\mR^m)/O(m)\ &\iso \  T h(G r_{|V|}(V\oplus \mR^m))  \\
[x, \varphi]\quad &\longmapsto \ (x,\varphi(\nu_m)^\perp )
\end{align*}
In this sense, $\bmO_{(m)}(V)$ `is' the Thom space over $\bbO_{(m)}(V)$.
Just as the orthogonal spaces $\bbO_{(m)}$ form an exhaustive filtration
of $\bbO$, the orthogonal spectra $\bmO_{(m)}$ form an exhaustive filtration
of the Thom spectrum $\bmO$, compare Proposition \ref{prop:mO global hocolim} below.
\end{construction}

As we explained in Remark \ref{rk:free are Thom},
the semifree orthogonal spectrum $F_m=F_{O(m),\nu_m}$
is a global refinement of the Thom spectrum $M T(m)$ of the
negative of the tautological $m$-plane bundle over $G r_m(\mR^\infty)$.
Since shift and suspension are globally equivalent
(by Proposition \ref{prop:global equiv preservation}~(i)),
$\bmO_{(m)}$ is globally equivalent to the
$m$-fold suspension of the orthogonal spectrum $M_{\gl} T (m)=F_{O(m),\nu_m}$ 
\[ \bmO_{(m)}\ = \ \sh^m F_m\ \simeq_{\gl} \ F_m \sm S^m\ = \  M_{\gl} T(m) \sm S^m\ .\]

We define a morphism of orthogonal spectra
\[ i \ : \ F_m\ \to \ \sh F_{m+1} \ ;\]
the value at an inner product space $V$ is the closed embedding
\begin{align*}
  i(V) \ : \ F_m(V) = \bO(\nu_m,V)/O(m)
 &\xra{ -\oplus\mR}  \bO(\nu_{m+1},V\oplus \mR)/O(m+1) = (\sh F_{m+1})(V) \\
[x,\varphi] \qquad &\longmapsto\qquad [(x,0),\varphi\oplus\mR] \,
\end{align*}
where we identify
$\nu_m\oplus\mR$ with $\nu_{m+1}$ by 
sending $((x_1,\dots,x_m),y)$ to $(x_1,\dots,x_m,y)$.
In fact, there are very few morphisms of orthogonal spectra from $F_m$ to $\sh F_{m+1}$:
by the representing property of the semifree spectrum $F_m$, such
morphisms biject with $O(m)$-fixed points
of $(\sh F_{m+1})(\nu_m)=\bO(\nu_{m+1},\nu_m\oplus\mR)/O(m+1)$;
this space only has two elements, and the morphism $i$
corresponds to the non-basepoint element.
We define 
\[  j^m\ = \ \sh^m i \ :\ \bmO_{(m)} = \sh^m F_m \ \to \ 
\sh^m(\sh F_{m+1})\ = \ \sh^{m+1}F_{m+1}\ = \ \bmO_{(m+1)}\ .\]
We define a morphism
\[ \psi^m \ : \ \bmO_{(m)}\ = \ \sh^m F_m \ \to \ \bmO    \]
at an inner product space $V$ as the map
\begin{align*}
 \psi^m(V) \ : \ \bO(\nu_m,V\oplus\mR^m)/O(m)\ &\to \ T h(G r_{|V|}(V\oplus\mR^\infty))\\
[x,\varphi] \quad &\longmapsto \quad ( i(x), i(\varphi^\perp)  )\ ,  
\end{align*}
where $\varphi^\perp=(V\oplus \mR^m)-\varphi(\nu_m)$
is the orthogonal complement of the image of $\varphi$ and
 $i:V\oplus \mR^m \to V\oplus\mR^\infty$ is the `standard' embedding given by
\[ i(v,x_1,\dots,x_m)\ = \ (v,x_1,\dots,x_m,0,0,\dots) \ .\]
The rank filtration expresses the orthogonal spectrum $\bmO$ 
as the colimit of the sequence of closed embeddings
\begin{equation}  \label{eq:bmO_filtration}
 \mS \ \iso \ \bmO_{(0)}\ \xra{j^0} \ \bmO_{(1)}\ \xra{j^1} \ \dots\ \to \ \bmO_{(m)}\ \xra{j^m} \ \dots   \ .
\end{equation}
The following proposition is straightforward from the definitions.
The periodic analog was already used in the proof of 
Theorem \ref{thm:mOP is localized MGr},
since it expresses $\bmOP$ as a sequential colimit of the spectra $\sh^m\bMGr$.
We omit the proof.

\begin{prop}\label{prop:mO global hocolim} 
For every $m\geq 0$ the morphism 
$j^m:\ \bmO_{(m)} \to \bmO_{(m+1)}$ is levelwise a closed embedding. 
These morphisms satisfy $\psi^{m+1}\circ j^m= \psi^m$.  
With respect to the morphisms $\psi^m:\bmO_{(m)}\to\bmO$,
the orthogonal spectrum $\bmO$ is a colimit of the sequence \eqref{eq:bmO_filtration}.
\end{prop}

Since colimits along sequences of closed embeddings
are invariant under global equivalences 
(Proposition \ref{prop:global equiv preservation}~(v)),
Proposition \ref{prop:mO global hocolim} says that $\bmO$ is a homotopy colimit
of the sequence \eqref{eq:bmO_filtration}.
The underlying non-equivariant statement, i.e., that $M O$ is
a homotopy colimit of the spectra $\Sigma^m M T(m)$,
can for example be found in \cite[Sec.\,3]{galatius-madsen-tillmann-weiss}. 
The identification of $\bmO$ as a homotopy colimit of semifree orthogonal spectra
now allows an algebraic description of $\gh{\bmO,E}$, 
the group of morphisms in the global stable homotopy category $\GH$,
into any orthogonal spectrum $E$.

We define a distinguished class
\[ \tau_m\ \in \ \pi_m^{O(m)}(\bmO_{(m)}\sm S^{\nu_m}) \]
in the $m$-th $O(m)$-equivariant $\bmO_{(m)}$-homology group of $S^{\nu_m}$
as the class represented by the based $O(m)$-map
\begin{align*}
  S^{\nu_m\oplus\mR^m}\ &\to \ \bO(\nu_m,\nu_m\oplus\mR^m)/O(m)\sm S^{\nu_m} 
\ = \ \bmO_{(m)}(\nu_m)\sm S^{\nu_m}\\
(v,x)\ &\longmapsto \qquad [ (0,x),i] \sm (-v) \ ,\nonumber
\end{align*}
where $i:\nu_m\to\nu_m\oplus\mR^m$ is the embedding as the first summand.

\begin{prop}\label{prop:i on universal classes}
  Let $m\geq 0$ be a natural number.
  \begin{enumerate}[\em (i)]
  \item  
    The pair $(\bmO_{(m)},\tau_m)$ represents the functor
    \[ \GH\to \text{\em (sets)}\ , \quad E \ \longmapsto \ E_m^{O(m)}(S^{\nu_m})
    \ = \ \pi_m^{O(m)}(E\sm S^{\nu_m}) \ . \]
  \item  
    The morphism $j^m:\bmO_{(m)}\to\bmO_{(m+1)}$ satisfies the relation
    \[ (j^m\sm S^{\nu_m})_*(\tau_m) \sm S^1 \ = \  \res^{O(m+1)}_{O(m)}(\tau_{m+1}) \]
    in the group  $\pi_{m+1}^{O(m)}(\bmO_{(m+1)}\sm S^{\nu_m}\sm S^1)$.
  \item The morphism $\psi^m:\bmO_{(m)}\to\bmO$ sends $\tau_m$
    to the shifted inverse Thom class $\bar\tau_{O(m),\nu_m}$
    defined in \eqref{eq:shifted inverse_Thom_mO}.\index{subject}{shifted inverse Thom class!in $\bmO$} 
 \end{enumerate}
\end{prop}
\begin{proof}
(i) 
In \eqref{eq:define_a_G,G} we defined a distinguished equivariant homotopy class
\[  a_m \ = \ a_{O(m),\nu_m} \ \in \ \pi_0^{O(m)}(F_m\sm S^{\nu_m})  \ .\] 
Expanding the definitions of $a_m$ and of the morphism
\[  \lambda^m_{F_m\sm S^{\nu_m}}\ :\ F_m\sm S^{\nu_m}\sm S^m\ \to\ \sh^m F_m\sm S^{\nu_m}
\ = \ \bmO_{(m)}\sm S^{\nu_m}\]
shows that
\[ \tau_m\ = \ (\lambda^m_{F_m\sm S^{\nu_m}})_*(a_m\sm S^m) \ .\]
By Theorem \ref{thm:free spectrum represents},
the pair $(F_m,a_m)$ represents the functor
\[ \GH\to \text{(sets)}\ , \quad E \ \longmapsto \ E_0^{O(m)}(S^{\nu_m})\ = \
\pi_0^{O(m)}(E\sm S^{\nu_m})\ . \]
We claim that the following composite
\begin{align*}
\gh{\bmO_{(m)}, E} \ &= \ 
\gh{\sh^m F_m, E} \ \xra{\gh{\lambda^m_{F_m},E}} \ \gh{F_m\sm S^m, E} \\ 
&\xra{\text{adjunction}} \ \gh{F_m,\Omega^m E} \ 
\xra{\text{eval at $a_m$}} \ \pi_0^{O(m)}((\Omega^m E)\sm S^{\nu_m}) \\
&\xra{\text{assembly}_*} \ \pi_0^{O(m)}(\Omega^m (E\sm S^{\nu_m})) \
\xra{\ \alpha^m\ } \ \pi_m^{O(m)}(E\sm S^{\nu_m}) 
\end{align*}
coincides with evaluation at the class $\tau_m$. 
Here the fourth map is induced by the assembly morphism
\[ (\Omega^m E)\sm S^{\nu_m}\ \to \ \Omega^m( E\sm S^{\nu_m}) \ ,\]
and
\[ \alpha^m\ :\ \pi_0^{O(m)}(\Omega^m Y)\ \iso \ \pi_m^{O(m)}(Y) \]
is the analog of the loop isomorphism \eqref{eq:loop iso}. 
Since all the maps are natural in $E$, it suffices to check this claim
for the identity of the universal example $E=\bmO_{(m)}=\sh^m F_m$.
Since the square
\[ \xymatrix@C=15mm{ 
E\sm S^{\nu_m} \ar[r]^-{\eta} \ar[d]_{\tilde\lambda^m_E\sm S^{\nu_m}} &
\Omega^m(E\sm S^{\nu_m}\sm S^m) \ar[d]^{\Omega^m(\lambda^m_{E\sm S^{\nu_m}})}\\
(\Omega^m \sh^m E)\sm S^{\nu_m} \ar[r]_{\text{assembly}} & 
\Omega^m (\sh^m E\sm S^{\nu_m}) } \]
commutes, we obtain
\begin{align*}
  \alpha^m(\text{assembly}_*((\tilde\lambda^m_{F_m}\sm S^{\nu_m})_*(a_m))) \ &= \ 
  \alpha^m( \Omega^m(\lambda^m_{F_m\sm S^{\nu_m}})_*(\eta_*(a_m))) \\ 
&= \  (\lambda^m_{F_m\sm S^{\nu_m}})_*(\alpha^m(\eta_*(a_m))) \\ 
&= \  (\lambda^m_{F_m\sm S^{\nu_m}})_*( a_m\sm S^m)\ = \ \tau_m\ .
\end{align*}
This verifies the relation in the universal example.
Since all the individual maps in the above composite are bijective,
so is the composite, which proves the representability property
of the pair $(\bmO_{(m)},\tau_m)$.

(ii)
The class $ (j^m\sm S^{\nu_m})_*(\tau_m) \sm S^1$ is represented by the composite:
\begin{align*}
 S^{\nu_m\oplus\mR^{m+1}}\ &\xra{t_m\sm S^1} \ 
\bO(\nu_m,\nu_m\oplus\mR^m)/O(m)\sm S^{\nu_m} \sm S^1 \\
 &\xra{j^m(\nu_m\oplus\mR^m)\sm S^{\nu_m}\sm S^1} \ 
\bO(\nu_{m+1},\nu_m\oplus\mR^{m+1})/O(m+1)\sm S^{\nu_m} \sm S^1\\
 (v,x,s)\ &\longmapsto \ [ (0,x),i]\sm (-v)\sm s
\ \longmapsto \ [ (0,x,0),i\oplus\mR ]\sm (-v)\sm s  \ .
\end{align*}
If we stabilize this representative along the linear isometric embedding
\[ j\ : \ \nu_m\oplus\mR^{m+1}\ \to \ \res^{O(m+1)}_{O(m)}(\nu_{m+1})\oplus\mR^{m+1} \ , \quad
 (v,x,s)\ \longmapsto \ (v,0,x,s) \]
we obtain another representative, namely
\begin{align*}
 S^{\nu_m\oplus\mR\oplus\mR^{m+1}}\ &\to \ 
\bO(\nu_{m+1},\res^{O(m+1)}_{O(m)}(\nu_{m+1})\oplus\mR^{m+1})/O(m+1)\sm S^{\nu_m} \sm S^1\\
 (v,u,x,s)\ &\longmapsto \quad [ (0,u,x,0),i\oplus\mR ] \sm (-v)\sm s \ .
\end{align*}
Here $i\oplus\mR:\nu_{m+1}\to\res^{O(m+1)}_{O(m)}(\nu_{m+1})\oplus\mR^{m+1}$
sends $(v,u)$ to $(v,0,0,u)$.

On the other hand, 
 $\res^{O(m+1)}_{O(m)} (\tau_{m+1})$ is represented by the underlying
$O(m)$-map of the $O(m+1)$-map
\begin{align*}
 S^{\nu_{m+1}\oplus\mR^{m+1}}\ &\xra{\ t_{m+1}\ } \ 
\bO(\nu_{m+1},\nu_{m+1}\oplus\mR^{m+1})/O(m+1)\sm S^{\nu_{m+1}}\\
 (v,u,x,s)\ &\longmapsto \quad [ (0,0,x,s),i' ]\sm (-v, -u)\ ,
\end{align*}
where $i':\nu_{m+1}\to\res^{O(m+1)}_{O(m)}(\nu_{m+1})\oplus\mR^{m+1}$
sends $(v,u)$ to $(v,u,0,0)$.
The two representatives differ by conjugation with the $O(m)$-equivariant
linear isometry
\[ \nu_m\oplus\mR \oplus\mR^m\oplus\mR \ \to \ 
\nu_m\oplus\mR \oplus\mR^m\oplus\mR  \ , \quad
(v,u,x,s)\ \longmapsto \ (v,-s,x,u)\ , \]
so they represent the same class in 
the group  $\pi_{m+1}^{O(m)}(\bmO_{(m+1)}\sm S^{\nu_m}\sm S^1)$.

(iii) Substituting the definitions of $\tau_m$ and of the morphism $\psi^m$
shows that $(\psi^m\sm S^{\nu_m})_*(\tau_m)$ is represented by the map
\[    S^{\nu_m\oplus\mR^m}\ \to \  \bmO(\nu_m)\sm S^{\nu_m}\ ,\quad
(v,x)\ \longmapsto \quad [ (0,i(x)), 0\oplus i(\mR^m)] \sm (-v) \ , \]
where $i:\mR^m\to\mR^\infty$ is the `standard' embedding as the leading $m$ coordinates.
The periodicity class $\sigma$ is represented by the map
\[ s\ : \ S^1 \ \to \ T h(G r(\mR^\infty)) \ = \ \bmOP(0)\ , \quad
x\ \longmapsto \ ((x,0,0,\dots),\mR\oplus 0)\ .\]
The multiplication of $\bmOP$ is an $E_\infty$-multiplication,
so as explained in \eqref{eq:psi_for_mu_bmO},
multiplication of classes in $\bmOP$-theory
involves choices of linear isometric embeddings.
We choose a linear isometric embedding
\[ \psi \ : \ (\mR^\infty)^m \ \to \ \mR^\infty \]
that satisfies 
\[ \psi( (x_1,0,0,\dots),\dots,(x_m,0,0,\dots))\ = \ 
(x_1,\dots,x_m,0,0,\dots) \ .\]
If we base the multiplication of $\bmOP$ on such a choice $\psi$,
then the product of the defining representative for $\tau_{O(m),\nu_m}$
and $m$ factors of the map $s$ above is precisely the 
previous representative for the class $(\psi^m\sm S^{\nu_m})_*(\tau_m)$.
This proves the relation $(\psi^m\sm S^{\nu_m})_*(\tau_m)=\tau_{O(m),\nu_m} \cdot p_{O(m)}^*(\sigma^m)$ in$~\bmO_m^{O(m)}(S^{\nu_m})$.
\end{proof}

The fact that $\bmO$ is the global homotopy colimit of the sequence of
orthogonal spectra $\bmO_{(m)}$ (see Proposition \ref{prop:mO global hocolim})
has the following consequence. 

\begin{cor}\label{cor:Milnor for maps out of bmO} 
For every orthogonal spectrum $E$ the following sequence is short exact:
\[ 0 \ \to \ {\lim_m}^1\,  E_{m+1}^{O(m)}(S^{\nu_m})
\ \to \ \gh{\bmO, E}\ \to \ 
\lim_m  \, E_m^{O(m)}(S^{\nu_m})\ \to \ 0 \]
Here  the right map is evaluation at the shifted inverse Thom
classes $\bar\tau_{O(m),\nu_m}$ and 
the inverse limit and derived limit are formed along the maps
\[ 
E_{m+1}^{O(m+1)}(S^{\nu_{m+1}})\ \xra{\res^{O(m+1)}_{O(m)}} \ E_{m+1}^{O(m)}(S^{\nu_m}\sm S^1)\ 
\xra{(-\sm S^1)^{-1}} \ E_m^{O(m)}(S^{\nu_m}) \ .\]
\end{cor}
\begin{proof}
Since $\bmO$ is the sequential homotopy colimit, 
in the triangulated global stable homotopy category, 
of the sequence of orthogonal spectra \eqref{eq:bmO_filtration}, 
the Milnor exact sequence takes the form:
\[ 0 \ \to \ {\lim_m}^1\,  \gh{\bmO_{(m)}\sm S^1,E}
\ \to \ \gh{\bmO, E}\ \to \ \lim_m  \, \gh{\bmO_{(m)},E}\ \to \ 0 \]
By  Proposition \ref{prop:i on universal classes}
the pair $(\bmO_{(m)},\tau_m)$ then represents the functor
\[ \GH\to \text{(sets)}\ , \quad E \ \longmapsto \ E_m^{O(m)}(S^{\nu_m})\ , \]
and the morphism $j^m:\bmO_{(m)}\to\bmO_{(m+1)}$ has the correct behavior
on the universal classes.
\end{proof}

Now that we recognized $\bmO$ as the global homotopy colimit of
the rank filtration \eqref{eq:bmO_filtration}, we study how
one filtration term is obtained from the previous one.
The answer given by Theorem \ref{thm:rank mO triangles} below
takes the form of a distinguished triangle
in the global stable homotopy category, witnessing that
the mapping cone of $j^m:\bmO_{(m-1)}\to \bmO_{(m)}$
`is' the $m$-fold suspension of the suspension spectrum
of the global classifying space~ $B_{\gl}O(m)$.

We define a morphism of orthogonal spectra
\[ T_m\ =\ T_{O(m)}^{O(m+1)} \ : \ 
 \Sigma^\infty_+ B_{\gl}O(m+1)\ \to \ F_m \]
as the adjoint of the $O(m+1)$-equivariant map
\begin{align*}
 r \ : \  S^{\nu_{m+1}} \ &\ \to \ \bO(\nu_m,\nu_{m+1})/O(m) 
\ = \ F_m(\nu_{m+1}) \\
 r(A\cdot(0,\dots,0,t))\ &= \
A\cdot( (0,\dots,0, (t^2-1)/t),\text{incl})\cdot O(m)\ ,
\end{align*}
where $A\in O(m+1)$ and $t\in[0,\infty]$.
Theorem \ref{thm:fundamental triangle}~(ii) shows that 
the morphism $T_m$ represents the dimension shifting transfer
from $O(m)$ to $O(m+1)$, in the sense of the relation
\begin{equation}  \label{eq:T_M_of_e_m+1}
 (T_m)_*(e_{O(m+1)})\ = \ \Tr_{O(m)}^{O(m+1)}(a_m)   
\end{equation}
between the tautological classes.
From $T_m$ we define another morphism of orthogonal spectra
\begin{align*}
  \partial \ = \ \lambda_{F_m}^m\circ (T_m\sm S^m)
 \ = \ &(\sh^m T_m)\circ \lambda_{\Sigma^\infty_+ B_{\gl} O(m+1)}^m
\ : \\ 
 &\Sigma^\infty_+ B_{\gl}O(m+1)\ \to \ \sh^m F_m \ = \ \bmO_{(m)}\ .
\end{align*}
The space
\[ ((\Sigma^\infty_+ \bL_{O(m+1),\nu_{m+1}})(\nu_{m+1}))^{O(m+1)} 
\ = \ (S^{\nu_{m+1}}\sm \bL(\nu_{m+1},\nu_{m+1})/O(m+1) )^{O(m+1)} \]
has two points, the basepoint and $0\sm \Id_{\nu_{m+1}}\cdot O(m+1)$.
So there is a unique non-trivial morphism of orthogonal spectra 
\[ a\ :\ F_{m+1}\ \to\ \Sigma^\infty_+ \bL_{O(m+1),\nu_{m+1}}=\Sigma^\infty_+ B_{\gl}O(m+1) \ .  \]
For every orthogonal spectrum $E$
the morphism $\lambda_E^{m+1}:E\sm S^{m+1}\to\sh^{m+1} E$ 
is a global equivalence
(by an iteration of Proposition \ref{prop:global equiv preservation}~(i)),
so it becomes invertible in the global stable homotopy category.
We can thus define a morphism in $\GH$ as
\[ q \ = \ (\lambda_{\Sigma^\infty_+ B_{\gl} O(m+1)}^{m+1})^{-1}\circ (\sh^{m+1 }a)
\ : \ \bmO_{(m+1)}\ \to \ \Sigma^\infty_+ B_{\gl}O({m+1})\sm S^{m+1}\ .  \]

\begin{theorem}\label{thm:rank mO triangles}\index{subject}{global classifying space!of $O(n)$}  
The sequence
\[ 
 \Sigma^\infty_+ B_{\gl}O({m+1})\sm S^{m}\ \xra{\ \partial \ } \
\bmO_{(m)}\ \xra{j^{m}} \ \bmO_{({m+1})} \ \xra{\ q\ }\ 
\Sigma^\infty_+ B_{\gl}O({m+1})\sm S^{m+1}
 \]
is a distinguished triangle in the global stable homotopy category.
The behavior of the first morphism on the stable tautological class is
given by
\[ \partial_*( e_{O({m+1})}\sm S^m) \ = \  \Tr_{O(m)}^{O(m+1)}(\tau_m)\ .  \]
\end{theorem}
\begin{proof}
The tautological $O(m+1)$-representation $\nu_{m+1}$ is faithful,
and the action of $O(m+1)$ on the unit sphere $S(\nu_{m+1})$ is transitive.
The stabilizer group of the unit vector $(0,\dots,0,1)$
identifies with the group $O(m)$, and the orthogonal complement of
this vector becomes the tautological representation $\nu_m$. 
We can thus apply Theorem \ref{thm:fundamental triangle}
and obtain a distinguished triangle:
\[ F_{m+1} \ \xra{\ a\ }\ \Sigma^\infty_+ B_{\gl}O({m+1})\ \xra{\ T_m \ } \
F_m\ \xra{-\lambda_{F_{m+1}}^{-1}\circ i \ } \  F_{m+1} \sm S^1\]
By Example \ref{eg:shift preserves triangles}
shifting preserves distinguished triangles; 
so the following sequence is also distinguished:
\[ \sh^m F_{m+1} \ \xra{\sh^m a}\ \sh^m \Sigma^\infty_+ B_{\gl}O({m+1})\ \xra{\sh^m T_m } \
\sh^m F_m\ \xra{-\sh^m(\lambda_{F_{m+1}}^{-1}\circ i) } \  \sh^m F_{m+1} \sm S^1\]
The rotation of this triangle is the lower sequence in the following diagram:
\[ \xymatrix@C=10mm{ 
\Sigma^\infty_+ B_{\gl}O(m+1)\sm S^{m}\ar[r]^-{\partial} 
\ar[d]_{\lambda^m_{\Sigma^\infty_+ B_{\gl}O(m+1)}} &
\bmO_{(m)} \ar[r]^-{j^m} \ar@{=}[d]&  \bmO_{(m+1)} \ar[d]^{\sh^m(\lambda_{F_{m+1}}^{-1})} \ar[r]^-q&
\Sigma^\infty_+ B_{\gl}O(m+1) \sm S^{m+1} 
\ar[d]^{\lambda^m_{\Sigma^\infty_+ B_{\gl}O(m+1)}\sm S^1}\\ 
\sh^m \Sigma^\infty_+ B_{\gl}O(m+1)\ar[r]_-{\sh^m T_m} &
\sh^m F_m \ar[r]_-{\sh^m(\lambda_{F_{m+1}}^{-1}\circ i)} & 
\sh^m F_{m+1}\sm S^1   \ar[r]_-{\sh^m a\sm S^1} &
\sh^m \Sigma^\infty_+ B_{\gl}O(m+1) \sm S^1 } \]
The left and middle squares commute by definition of the morphisms $\partial$
respectively $j^m$; the right square commutes by the relation
\begin{align*}
(\lambda^m_{\Sigma^\infty_+ B_{\gl}O(m+1)}\sm S^1)\circ q \ &= \ 
(\lambda^m_{\Sigma^\infty_+ B_{\gl}O(m+1)}\sm S^1)\circ (\lambda_{\Sigma^\infty_+ B_{\gl} O(m+1)}^{m+1})^{-1}\circ (\sh^{m+1 }a)\\
&= \ 
(\sh^m \lambda_{\Sigma^\infty_+ B_{\gl} O(m+1)})^{-1}\circ (\sh^{m+1}a) \\
 &= \ \sh^m(\lambda_{\Sigma^\infty_+ B_{\gl} O(m+1)}^{-1}
\circ(\sh a)) \
= \ \sh^m((a\sm S^1)\circ\lambda_{F_{m+1}}^{-1}) \ .
\end{align*}
So the upper sequence is a distinguished triangle.

The final relation is a consequence of various other previously established relations:
\begin{align*}
\partial_*( e_{O({m+1})}\sm S^m ) \ &= \ 
(\lambda_{F_m}^m)_*( (T_m\sm S^m )_*(e_{O({m+1})}\sm S^m ) ) \\ 
 &= \ (\lambda_{F_m}^m)_*( (T_m)_*(e_{O({m+1})})\sm S^m )  \\ 
_\eqref{eq:T_M_of_e_m+1} &= \ (\lambda_{F_m}^m)_*(\Tr_{O(m)}^{O({m+1})}(a_m)\sm S^m )  \\ 
 &= \ (\lambda_{F_m}^m)_*(\Tr_{O(m)}^{O(m+1)}((F_m\sm\tau_{\nu_m,\mR^m})_*(a_m\sm S^m )))  \\ 
 &= \ \Tr_{O(m)}^{O(m+1)}((\lambda_{F_m}^m\sm S^{\nu_m})_*((F_m\sm\tau_{\nu_m,\mR^{m}})_*(a_m\sm S^m )))  \\ 
&= \  \Tr_{O(m)}^{O({m+1})}((\lambda^m_{F_m\sm S^{\nu_m}})_*(a_m\sm S^m))
\ = \  \Tr_{O(m)}^{O(m+1)}(\tau_m)\ . 
\end{align*}
The second and fifth equations are naturality.
The fourth equation is the compatibility of transfer and suspension isomorphism,
see Proposition \ref{prop:internal transfer loop suspension}.
\end{proof}

For calculations of equivariant homotopy groups of $\bmO$
we also need to understand the composite:
\[   \Sigma^\infty_+ B_{\gl} O(m+1) \sm S^m
\ \xra{\ \partial\ } \ \bmO_{(m)} \ \xra{\ q\ } \ 
 \Sigma^\infty_+ B_{\gl} O(m) \sm S^m\]
We start from the calculation
\begin{align*}
  \left(a\circ T_m\right)_*(e_{O(m+1)})\ &= \ 
  a_*\left( \Tr_{O(m)}^{O(m+1)}(a_m ) \right)\ = \
 \Tr_{O(m)}^{O(m+1)} \left(  a_*(a_m ) \right)\\ 
&= \ \Tr_{O(m)}^{O(m+1)} \left(  a\cdot e_{O(m)}) \right)\ = \
   \tr_{O(m)}^{O(m+1)}(e_{O(m)})\ .
\end{align*}
On the other hand,
\begin{align*}
  q\circ \partial \ &= \ 
(\lambda_{\Sigma^\infty_+ B_{\gl} O(m)}^m)^{-1}\circ (\sh^m a)\circ
\lambda_{F_m}^m\circ (T_m\sm S^m)  \\
&= \ (\lambda_{\Sigma^\infty_+ B_{\gl} O(m)}^m)^{-1}\circ 
\lambda_{\Sigma^\infty_+ B_{\gl}O(m)}^m\circ (a\sm S^m)\circ (T_m\sm S^m)  \
= \  (a\circ T_m)\sm S^m\ ,
\end{align*}
by definition.
Combining these two facts gives
\begin{align}\label{eq:partial q}
  ( q\circ \partial)_*(e_{O(m+1)}\sm S^m) \ &= \ 
((a\circ T_m)\sm S^m)_*(e_{O(m+1)}\sm S^m) \\ 
&= \ 
(a\circ T_m)_*(e_{O(m+1)})\sm S^m\ = \    \tr_{O(m)}^{O(m+1)}(e_{O(m)})\sm S^m\ .\nonumber
\end{align}
In other words, the composite $q\circ\partial$
represents the degree zero transfer $\tr_{O(m)}^{O(m+1)}$.

Now we can easily show that $\bmO$ is globally connective
and describe the global functor $\upi_0(\bmO)$.
We denote by $\td{\tr_e^{O(1)}}$ the global subfunctor of the Burnside ring functor
$\mA$ generated by $\tr_e^{O(1)}\in\mA(O(1))$.

\begin{theorem}\label{thm:pi_0 of mO} 
The orthogonal spectrum $\bmO$ is globally connective
and the action of the Burnside ring global functor on the unit element
$1\in \pi_0^e(\bmO)$ induces an isomorphism of global functors
\[ \mA /\td{\tr_e^{O(1)}} \ \iso\ \upi_0(\bmO) \ .\]
\end{theorem}
\begin{proof}
The suspension spectrum $\Sigma^\infty_+ B_{\gl}O(m+1)$ is globally connective
(Proposition \ref{prop:pi_0 of Sigma^infty});
so the distinguished triangle of Theorem \ref{thm:rank mO triangles}
implies that the morphism $j^m:\bmO_{(m)}\to\bmO_{(m+1)}$
induces an isomorphism of global functors 
\[ \upi_k(\bmO_{(m)})\ \iso \ \upi_k(\bmO_{(m+1)})\]
for $k\leq m-1$ and an exact sequence of global functors
\[ \bA(O(m+1),-)\ \to \ \ \upi_m(\bmO_{(m)})\ \to \ \upi_m(\bmO_{(m+1)})\ \to \ 0 \ .\]
Here we used the isomorphisms
\begin{align*}
 \bA(O(m+1),-) \ &\xra[\iso]{\tau\mapsto\tau(e_{O(m+1)})} \ 
\upi_0(\Sigma^\infty_+ B_{\gl}O(m+1)) \\
&\xra[\iso]{\ -\sm S^m\ } \ \upi_m(\Sigma^\infty_+ B_{\gl}O(m+1)\sm S^m) \  \ .
\end{align*}
Since $\bmO_{(0)}$ is isomorphic to the global sphere spectrum,
which is globally connective, we conclude inductively that $\bmO_{(m)}$
is globally connective for all $m\geq 0$, and that the inclusion
$\bmO_{(1)}\to\bmO_{(m)}$ induces an isomorphism on $\upi_0$ for all $m\geq 1$.
Since $\bmO$ is a colimit of the sequence of closed embeddings
$j^m:\bmO_{(m)}\to \bmO_{(m+1)}$, the map
\[ \colim_{m}\, \upi_k(\bmO_{(m)})\ \to \ \upi_k(\bmO) \]
induced by the morphisms $\psi^m:\bmO_{(m)}\to\bmO$
is an isomorphism of global functors for every integer $k$.
So $\bmO$ is globally connective and the morphism $\psi^1:\bmO_{(1)}\to\bmO$
induces an isomorphism
\[ \upi_0(\bmO_{(1)})\ \iso \ \upi_0(\bmO)\ . \]
The unit morphism identifies $\bmO_{(0)}$ with the global sphere spectrum,
so  the action on the class $1\in\bmO_{(0)}$
is an isomorphism of global functors $\mA\iso \upi_0(\bmO_{(0)})$.
For $m=1$, the exact sequence thus becomes an exact sequence of global functors
\[ \bA(O(1),-)\ \xra{\bA(\tr_e^{O(1)},-)} \ \bA(e,-)\ \to \ \upi_0(\bmO_{(1)})\ \to \ 0 \ ;\]
here we used \eqref{eq:partial q}
to identify the first morphism as the one induced by the transfer $\tr_e^{O(1)}$.
This proves the claim about $\upi_0(\bmO)$.
\end{proof}

Theorem \ref{thm:pi_0 of mO} gives a nice compact description of
the 0-th equivariant homotopy groups of the Thom spectrum $\bmO$,
but we may still ask for a more explicit calculation of the
group $\pi_0^G(\bmO)$ for an individual compact Lie group $G$.
Given the presentation of $\upi_0(\bmO)$ as the quotient of
$\mA$ by the global subfunctor generated by $\tr_e^{O(1)}$,
this is a purely algebraic exercise.
Since the element $\tr_e^{O(1)}(1)$ is trivial in the group $\pi_0^{O(1)}(\bmO)$, also
\[ 2 \ = \ \res^{O(1)}_e(\tr^{O(1)}_e(1))\ = \ 0  \]
in $\pi_0^e(\bmO)$;
thus all equivariant homotopy groups $\pi_*^G(\bmO)$ are $\mF_2$-vector spaces.
The next proposition pins down an $\mF_2$-basis 
of $\pi_0^G(\bmO)$ in terms of the subgroup structure of $G$.
For a closed subgroup $H$ of $G$ we use the familiar notation
\[ t_H^G \ = \ \tr_H^G(p_H^*(1)) \ \in \ \pi_0^G(\bmO) \ , \]
where $p_H:H\to e$ is the unique homomorphism.

\begin{prop}\label{prop:explicit mO kernel} 
  For every compact Lie group $G$, an $\mF_2$-basis of $\pi_0^G(\bmO)$
  is given by the classes $t_H^G$, indexed by conjugacy classes
  of those closed subgroups $H$ of $G$ whose Weyl group is finite and of odd order.
\end{prop}
\begin{proof}
  We abbreviate $C=O(1)$.
  The group $\mA(G)$ is free abelian with basis the classes $t_H^G$ 
  for all conjugacy classes of subgroups $H$ with finite Weyl group.
  So the claim follows if we can show that the value $\td{\tr_e^C}(G)$ 
  of the global functor $\td{\tr_e^C}$ at $G$ is the subgroup of $\mA(G)$ 
  generated by $2\cdot\mA(G)$ and the classes $t_H^G$ 
  for those closed subgroups $H$ whose Weyl group is finite of even order.

  By Theorem \ref{thm:Burnside category basis}, 
  the group $\bA(C,G)$ is freely generated by the elements
  $\tr_L^G\circ\alpha^*$ indexed by $(G\times C)$-conjugacy classes
  of pairs $(L,\alpha)$ where $L$ is a closed subgroup of $G$ with finite Weyl group
  and $\alpha:L\to C$ is a continuous group homomorphism. 
  So $\td{\tr_e^C}(G)$ is generated as an abelian group by the elements
  $\tr_L^G\circ\alpha^*\circ \tr_e^C$.

  A homomorphism to $C$ is either trivial or surjective,
  and the generating elements come in two flavors.
  If $\alpha$ is the trivial homomorphism, then
  \[ \tr_L^G\circ\alpha^*\circ \tr_e^C \ = \ 
  \tr_L^G\circ p_L^*\circ \res^C_e\circ \tr_e^C \ = \ 
  2\cdot \tr_L^G\circ p_L^*\ = \ 2\cdot t_L^G\ .  \]
  These elements generate the subgroup $2\cdot \mA(G)$.
  If $\alpha$ is surjective with kernel $H$, then
  \[ \tr_L^G\circ\alpha^*\circ \tr_e^C \ = \ 
  \tr_L^G\circ \tr_H^L\circ p_H^* \ = \ \tr_H^G\circ p_H^* \ = \  t_H^G\ .  \]
  The group $H$ that arises in this way is normal of index~2 in $L$.
  So $L$ is contained in the normalizer $N_G H$, and~2 divides $[N_G H:H]$. 
  So the order of the Weyl group $W_G H$ is even.
  Conversely, if $H$  is a subgroup of $G$ with finite Weyl group of even order,
  then we can choose a subgroup $C\leq W_G H$ of order~2.
  The preimage $L$ of $C$ under the projection $N_G H\to W_G H$
  then contains $H$ as an index~2 subgroup. By the above, the class $t_H^G$ 
  is then one of the generating elements of $\td{\tr_e^C}(G)$.
\end{proof}

\begin{eg}[Geometric fixed points of $\bmO$ and $\bmOP$]\label{eg:geometric fixed points of bmO_m}
\index{subject}{geometric fixed points!of $\bmOP$}
\index{subject}{geometric fixed points!of $\bmO$}
We give a description of the geometric fixed points of $\bmO$ and $\bmOP$.
The geometric fixed points of other global Thom spectra such as
$\bmO_{(m)}$, $\bMO$ or $\bMOP$ can be worked out in a similar fashion, 
but we leave that to the interested reader.
We let $V$ be a representation of a compact Lie group $G$.
We write $V^\perp=V-V^G$ for the orthogonal complement of the $G$-fixed subspace.
A point $((v,x),U)\in \bmOP(V)=T h(G r(V\oplus\mR^\infty))$ 
is $G$-fixed if and only if the subspace $U$
of $V\oplus\mR^\infty$ is $G$-invariant and the vector $v\in V$ is $G$-fixed.
The first condition guarantees that $U=U^G\oplus (U\cap V^\perp)$,
and $U^G$ is a subspace of $V^G\oplus\mR^\infty$.
So we obtain a homeomorphism
\begin{align}  \label{eq:split_fix}
 \bmOP(V^G)\sm \bGr(V^\perp)^G\ &\xra{\ \iso\ } \quad  \bmOP(V)^G \\
(x, U)\sm W\quad  &\longmapsto\ ( x ,\, U\oplus (V^\perp- W)) \ .\nonumber   
\end{align}
Under this identification, the structure map
\[  (\sigma_{V,W})^G \ : \ S^{V^G}\sm \bmOP(W)^G \ \to \ \bmOP(V\oplus W)^G \]
becomes the  smash product of the structure map 
\[ \sigma_{V^G,W^G} \ : \ S^{V^G}\sm \bmOP(W^G) \ \to \ \bmOP(V^G\oplus W^G) \]
with the map
\[ \bGr(i)^G \ : \ \bGr(W^\perp)^G \ \to \  \bGr(V^\perp\oplus W^\perp)^G  \]
induced by the embedding $W^\perp\to V^\perp\oplus W^\perp$
as the second summand.
So in the colimit over $V\in s(\Uc_G)$ this gives an isomorphism
\begin{equation}\label{eq:Phi^G of mOP}
 \Phi_*^G(\bmOP) \ \iso \ \bmOP_*\left(\bGr(\Uc_G^\perp)^G_+\right)
\end{equation}
to the non-equivariant $\bmOP$-homology groups 
of the $G$-fixed point space of $\bGr(\Uc_G^\perp)$,
the disjoint union of all Grassmannians in $\Uc_G^\perp=\Uc_G-(\Uc_G)^G$.

The formula \eqref{eq:Phi^G of mOP} is a compact way 
to express the geometric fixed points of $\bmOP$, 
but it can be decomposed and rewritten further, thereby making it more explicit.
The orthogonal spectrum $\bmOP$ is a $\mZ$-indexed wedge of homogeneous
summands $\bmOP^{[k]}$. The space $\bGr(\Uc_G^\perp)$, and hence also
its $G$-fixed points, is the disjoint union indexed by the dimension of
the subspaces, i.e.,
\[ \bGr(\Uc_G^\perp)^G \ = \ {\coprod}_{j\geq 0}\, \left(G r_j(\Uc_G^\perp)\right)^G 
\ = \ {\coprod}_{j\geq 0} \, G r_j^{G,\perp} \ . \]
The two decompositions induce a direct sum decomposition of the right
hand side of the isomorphism \eqref{eq:Phi^G of mOP} as
\begin{align}\label{eq:split Phi^G(mOP)}
 \bmOP_*\left(\bGr(\Uc_G^\perp)^G_+\right)\ &= \ 
{\bigoplus}_{k\in\mZ}{\bigoplus}_{j\geq 0}\ 
\bmOP^{[k]}_*\left( (G r_j^{G,\perp})_+\right)\\ 
&\iso\
{\bigoplus}_{k\in\mZ,\ j\geq 0}\ 
\bmO_{*-k}\left( (G r_j^{G,\perp})_+\right)\ . \nonumber  
\end{align}
The second isomorphism uses the periodicity of $\bmOP$
to identify $\bmOP^{[k]}$ with $\bmO\sm S^k$.

The condition $\dim(U\oplus (V^\perp-W))=\dim(V)$
is equivalent to $\dim(U)=\dim(V^G)+\dim(W)$.
So the homeomorphism \eqref{eq:split_fix} identifies $\bmO(V)^G$
with the wedge of the spaces $\bmOP^{[j]}(V^G)\sm G r_j(V^\perp)$ for $j\geq 0$. 
The composite of the isomorphism \eqref{eq:Phi^G of mOP}
and the isomorphism \eqref{eq:split Phi^G(mOP)}
thus takes the wedge summand of $\Phi^G_*(\bmOP)$
corresponding to $\bmO=\bmOP^{[0]}$ to the sum of the terms with $k=j$.
So the isomorphisms restrict to an isomorphism
\begin{equation}\label{eq:Phi^G of mO}
 \Phi^G_*(\bmO)\ \iso\ {\bigoplus}_{j\geq 0}\ 
\bmO_{*-j}\left( (G r_j^{G,\perp})_+\right)\ .   
\end{equation}
The space $G r_j^{G,\perp}$ can be decomposed further:
every $G$-invariant subspace of $\Uc_G^\perp$ is the direct sum of
its isotypical components, indexed by the non-trivial irreducible $G$-representations. 
The irreducibles come in three flavors (real, complex or quaternionic),
and so the space $G r_j^{G,\perp}$ is a disjoint union
of products of classifying spaces of the groups $O(n)$, $U(n)$, and $S p(n)$
for various $n$.

The reader may want to compare the previous description of the geometric
$G$-fixed points of $\bmO$ with Proposition \ref{prop:pi_0 of bbO}~(i),
which identifies the $G$-fixed points of the orthogonal space $\bbO$ as
\[ \bbO(\Uc_G)^G\ \simeq \ {\coprod}_{j\geq 0} \,   B O\times G r_j^{G,\perp} 
 \ .\]
This illustrates the general phenomenon that `geometric $G$-fixed points of
a global Thom spectrum are the Thom spectrum over the $G$-fixed points'.
\end{eg}

\begin{eg}
The Thom spectrum $\bmO$ has an oriented analog.
For $m\geq 0$ we define an orthogonal spectrum $\bmSO_{(m)}$ by
\[ \bmSO_{(m)}\ = \ \sh^m F_{S O(m),\nu_m} \ ,\]
the $m$-th shift of the semifree orthogonal spectrum generated
by the $S O(m)$-representation $\nu_m$.
In much the same way in Construction \ref{con:mO_m}, 
$\bmSO_{(m)}(V)$ `is' (by passage to orthogonal complements)
the Thom space over the tautological
bundle over the oriented Grassmannian $G r^+_{|V|}(V\oplus\mR^m)$
of oriented $|V|$-planes in $V\oplus\mR^m$. 

We define a morphism
\[ i \ : \ F_{S O(m),\nu_m}\ \to \ \sh F_{S O(m+1),\nu_{m+1}} \]
at an inner product space $V$ as the closed embedding
\begin{align*}
  i(V) \ : \bO(\nu_m,V)/ S O(m)
\ &\xra{- \oplus\mR} \ \bO(\nu_{m+1},V\oplus\mR)/S O(m+1)  \\
[x,\varphi] \qquad &\longmapsto\qquad [(x,0),\varphi\oplus\mR]\ .
\end{align*}
Then we set
\[  j^m\ = \ \sh^m i \ :\ \bmSO_{(m)} = \sh^m F_{S O(m),\nu_m} \ \to \ 
 \sh^{m+1} F_{S O(m+1),\nu_{m+1}}\ = \ \bmSO_{(m+1)}\]
and we define $\bmSO$ as the colimit of the sequence of closed embeddings
\[ \bmSO_{(0)} \ \xra{\ j^0 \ } \ \bmSO_{(1)} \ \xra{\ j^1\ } \ \dots \ \to \ \bmSO_{(m)} 
\ \xra{\ j^m\ } \ \dots \ . \]
The orthogonal spectrum $\bmSO$ supports 
inverse Thom classes\index{subject}{inverse Thom class}\index{subject}{Thom class!inverse|see{inverse Thom class}}
for {\em oriented} representations,
i.e., oriented inner product spaces $V$ 
equipped with an orientation preserving isometric action of a compact Lie group $G$.
\end{eg}

Now we mention the unitary analogs $\bmU$ and $\bMU$ 
of the Thom spectra $\bmO$ and $\bMO$.
Beside the complexification, there is an extra twist to the
unitary definitions, because we need to `loop by imaginary spheres'
to really get orthogonal spectra.
The unitary Thom spectra have periodic versions $\bmUP$ respectively $\bMUP$, 
but we will not go into any details about those.\index{symbol}{$\bMUP$ - {periodic unitary global Thom spectrum}}

\begin{eg}
We define an orthogonal spectrum $\bmU$ in analogy
with $\bmO$.\index{symbol}{$\bmU$ - {unitary global Thom spectrum}}
Non-equivariantly, $\bmU$ is the unitary Thom spectrum $M U$. 
The spectrum $\bmU$ is essentially a Thom spectrum over the orthogonal space $\bbU$,
the complex analog of the $E_\infty$-orthogonal monoid space $\bbO$
discussed in Example \ref{eg:define bO}.
As before we denote by $V_\mC=\mC\tensor_\mR V$ the complexification 
of a real inner product space $V$. 
The value of $\bbU$ at $V$ is
\[ \bbU(V)\ = \  Gr_{|V|}^\mC(V_\mC\oplus\mC^\infty) \ , \]
the Grassmannian of $\mC$-linear subspaces of $V_\mC\oplus\mC^\infty$
of the same dimension as $V$.
The map $\bbU(\varphi):\bbU(V)\to\bbU(W)$ induced by a 
linear isometric embedding $\varphi:V\to W$ is defined as
\[ \bbU(\varphi)(L) \ = \ 
(\varphi_\mC\oplus\mC^\infty)(L) \ + \ ((W_\mC-\varphi_\mC(V_\mC))\oplus 0) \ .\]
In other words: we apply the linear isometric 
embedding $\varphi_\mC\oplus\mC^\infty:V_\mC\oplus\mC^\infty\to W_\mC\oplus\mC^\infty$ 
to the subspace $L$ and add the orthogonal complement of the image of $\varphi_\mC$
(sitting in the first summand of $W_\mC\oplus\mC^\infty$).

We let $i V$ denote the imaginary subspace of $V_\mC$, i.e., the
$\mR$-subspace consisting of the elements $i\tensor v$ for $v\in V$.
Over the space $\bbU(V)$ sits a tautological hermitian vector bundle 
with total space consisting of the pairs 
$(x,U)\in (V_\mC\oplus\mC^\infty)\times \bbU(V)$  such that $x\in U$.
We define $\bmU(V)$ as 
\[ \bmU(V) \ = \ \map_*(S^{i V}, T h( \bbU(V)))\ ,\]
the space of based maps from the imaginary sphere $S^{i V}$
to the Thom space of this tautological vector bundle.
The structure map of the spectrum $\bmU$ starts from the 
$(O(V)\times O(W))$-equivariant map
\begin{align*}
 S^{V_\mC}\sm T h(\bbU( W))\ &\to \quad T h(\bbU (V\oplus W)) \\
 v \sm  (x, U) \qquad &\longmapsto \  
( (v,x),\ U\oplus (0\oplus W_\mC\oplus 0))\ .   
\end{align*}
The structure map 
\[ \sigma_{V,W}\ : \ S^V\sm \bmU(W) \ \to \ \bmU(V\oplus W)\]
is then adjoint to the composite
\begin{align*}
 S^V\sm \map_*(S^{i W}, T h(\bbU(W)))&\sm S^{i(V\oplus W)} \\ 
\iso \quad  &S^{V_\mC}\sm \map_*(S^{i W}, T h( \bbU(W)))\sm S^{i W}\\
\xra{S^{V_\mC}\sm \ev} \ &S^{V_\mC}\sm T h(\bbU( W))\ \to \ T h(\bbU (V\oplus W)) \ .&
\end{align*}
Looping by $S^{i V}$ is essential for obtaining an orthogonal spectrum;
without it, we would end up with a structure one may call a `unitary spectrum'.
\end{eg}

Most of our results about $\bmO$ have analogues for $\bmU$.
The orthogonal spectrum $\bmU$ comes with an $E_\infty$-structure,
which is, however,  {\em not} ultra-commuta\-tive.
There are unitary versions of the shifted inverse Thom classes 
\index{subject}{shifted inverse Thom class!in $\bmU$}
\[  \bar\tau^U_{G,V} \ \in \ \bmU_{2 n}^G(S^V) \ , \]
in the $G$-equivariant $\bmU$-homology groups of $S^V$, defined
for {\em unitary} representations $V$ of a compact Lie group $G$.
Here $n$ is the complex dimension of $V$, so that $2 n$ is its real dimension.
The orthogonal spectrum $\bmU$ is the union of a
sequence of orthogonal subspectra
\[ 
\bmU_{(0)}\ \subset \ \bmU_{(1)}\ \subset \ \dots\ \subset \ \bmU_{(m)}\ \subset \ \dots \
 \]
and $\bmU$ is also the global homotopy colimit of this sequence.
The unit morphism is a global equivalence $\mS\simeq\bmU_{(0)}$
and $\bmU_{(m)}$ is globally equivalent to the $2 m$-th 
suspension of the semifree orthogonal spectrum generated by
the tautological unitary representation $\nu_m^U$ of $U(m)$ on $\mC^m$:
\[ \bmU_{(m)} \ \simeq \ F_{U(m),\nu_m^U}\sm S^{2 m} \]
There are distinguished triangles
in the global stable homotopy category:\index{subject}{global classifying space!of $U(n)$}  
\[ \Sigma^\infty_+ B_{\gl} U(m+1)\sm  S^{2 m + 1}
 \ \to \ \bmU_{(m)} \ \to \  \bmU_{(m+1)} \ \to \ 
\Sigma^\infty_+ B_{\gl} U(m+1) \sm S^{2 m +2 } \]
and the first map is classified by the $U(m+1)$-equivariant homotopy class
\[ \Tr_{U(m)}^{U(m+1)}(\tau_{U(m),\nu_m^U} \sm S^1)
\text{\qquad in\quad $\pi_{2 m +1}^{U(m+1)}(\bmU_{(m)})$.} \]
This uses that the tangent $U(m)$-representation 
of $U(m+1)/U(m)$ is isomorphic to $\nu_m^U\oplus\mR$.
So loosely speaking, $\bmU_{(m+1)}$ is obtained from $\bmU_{(m)}$
by coning off this transfer class.
A consequence is then that all $\bmU_{(m)}$ and $\bmU$ are globally connective.
Since the first unitary relation $\Tr_e^{U(1)}(1\sm S^1)=0$
in $\pi_1^{U(1)}(\bmU)$ lives in a positive dimension, 
the description of the global functor $\upi_0(\bmU)$ 
is easier than the corresponding
calculation of $\upi_0(\bmO)$ in Theorem \ref{thm:pi_0 of mO}.
Indeed, the action of the Burnside ring global functor on the element
$1\in \pi_0^e(\bmU)$ induces an isomorphism of global functors
\[ \mA  \ \iso\ \upi_0(\bmU) \ .\]
Moreover, there is an exact sequence of global functors
\[ \bA(U(1),-)\ \to \ \upi_1(\mS)\ \to \ \upi_1(\bmU)\ \to \ 0 \ ,\]
where the first morphism is the action on the class
$\Tr_e^{U(1)}(1\sm S^1)$ in $\pi_1^{U(1)}(\mS)$.

Corollary \ref{cor:Milnor for maps out of bmO} has a unitary analog that
describes morphisms in the global stable homotopy category from $\bmU$:
for every orthogonal spectrum $E$ the sequence 
\[ 0 \ \to \ {\lim_m}^1\,  E_{2 m+1}^{U(m)}(S^{\nu_m^U})
\ \to \ \gh{\bmU, E}\ \to \ 
\lim_m  \, E_{2 m}^{U(m)}(S^{\nu_m^U})\ \to \ 0 \]
is short exact.
The right map is evaluation at the unitary shifted inverse Thom classes,
and the inverse and derived limits are formed along the maps
\[ 
E_{2 m+2}^{U(m+1)}(S^{\nu_{m+1}^U})\ \xra{\res^{U(m+1)}_{U(m)}} \ E_{2 m+2}^{U(m)}(S^{\nu_m^U}\sm S^2)\ 
\xra{(-\sm S^2)^{-1}} \ E_{2 m}^{U(m)}(S^{\nu_m^U}) \ .\]
We leave it to the interested reader to formulate the analogous properties
of the rank filtrations for the special unitary and symplectic
global Thom spectra $\mathbf{mSU}$ and $\mathbf{mSp}$.

\begin{eg}\index{subject}{unitary global bordism}
We define another unitary global Thom spectrum $\bMU$,\index{symbol}{$\bMU$ - {unitary global Thom spectrum}} 
an ultra-commutative ring spectrum and the unitary analog of $\bMO$.
Non-equi\-variantly, $\bMU$ is another version of the complex bordism spectrum $MU$;
for a compact Lie group $G$, the orthogonal $G$-spectrum $\bMU_G$ 
is a model for tom Dieck's homotopical equivariant bordism \cite{tomDieck-bordism}.
Closely related, strictly commutative ring spectrum models
for these homotopy types have been discussed in various places, 
for example \cite{may-quinn-ray}, \cite[Ex.\,5.8]{greenlees-may-completion},
\cite[App.\,A]{strickland-realizing formal groups} or \cite[Sec.\,8]{brun-Witt cobordism}.

For an inner product space $V$ we consider the complex Grassmannian 
\[ \bBU(V)\ = \  Gr_{|V|}^\mC(V_\mC^2) \ . \]
Over the space $\bBU(V)$ sits a tautological hermitian vector bundle and we set
\[ \bMU(V)\ = \ \map_*(S^{i V}, T h(\bBU(V))) \ ,\]
the $i V$-loop space of the Thom space of this tautological vector bundle.
The structure maps are defined in a similar way as for $\bmU$,
and a commutative multiplication is given by
\[ \mu_{V,W}\ : \ \bMU(V) \sm \bMU(W) \ \to \ \bMU(V\oplus W) \ , \quad
f\sm g \ \longmapsto \ f\cdot g \ .  \]
Here $f:S^{i V}\to T h(\bBU(V))$, $g:S^{i W}\to T h(\bBU(W))$,
and $f\cdot g$ denotes the composite
\begin{align*}
 S^{i V}\sm S^{i W} \ \xra{f\sm g} \  &T h(\bBU(V))\sm T h(\bBU(W))\\ 
&\xra{(x,U)\sm (x',U')\mapsto (\kappa_{V,W}(x,x'),\kappa_{V,W}(U\oplus U'))} \ 
T h(\bBU(V\oplus W)) \ ,    
\end{align*}
where $\kappa_{V,W}$ is the preferred isometry from
$V_\mC^2\oplus W_\mC^2\iso (V\oplus W)_\mC^2$ sending $(v,v',w,w')$ to $(v,w,v',w')$.
Unit maps are defined by
\[ S^V\ \to\  \bMU(V) \ , \quad  v\ \longmapsto \ 
[ v'\mapsto \ ( (v+v',0),\,  V_\mC\oplus 0 ) ]\ .  \]
These multiplication maps unital, associative and commutative,
and make $\bMU$ into an ultra-commutative ring spectrum.
\end{eg}

The global Thom spectrum $\bMU$ comes with distinguished Thom classes
\index{symbol}{$\sigma^U_{G,V}$ - {unitary Thom class in $\bMU_V^G(S^{2 n})$}} 
\[ \sigma^U_{G,V} \ \in \ \bMU^{2 n}_G(S^V) \ , \]
for {\em unitary} representations $V$ of compact Lie groups $G$,
where $n$ is the complex dimension of $V$.
As in the orthogonal situation in Theorem \ref{thm:Thom iso for MOP}
the Thom class $\sigma_{G,V}^U$
is inverse to the image of the inverse Thom class $\tau_{G,V}^U$,
and $\sigma_{G,V}^U$ restricts to a unitary Euler class
in $\pi_{-2 n}^G(\bMU)$.\index{subject}{Euler class!in $\bMU$}
The proof of Corollary \ref{cor:MO is localized mO} generalizes to the
unitary situation and proves that the $G$-equivariant homology
represented by $\bMU$ is the localization 
at the unitary inverse Thom classes of the theory represented by $\bmU$.

Some general facts about the equivariant homotopy groups of $\bMU$
are known for {\em abelian} compact Lie groups $A$.
In that case, $\pi_*^A(\bMU)$ is a free module on even dimensional generators
over the non-equivariant homotopy ring $\pi_*^e(\bMU)$;
this calculation was announced by L{\"o}ffler in \cite{loeffler-budva},
and a proof by Comeza{\~n}a can be found in \cite[XXVIII Thm.\,5.3]{may-alaska}.
Since the graded ring $\pi_*^e(\bMU)$ is concentrated in even degrees
and the periodic version $\bMUP$ is a wedge of even suspensions of $\bMU$,
the groups $\pi_*^A(\bMUP)$ are concentrated in even degrees.
Since these groups are also 2-periodic, for abelian compact Lie groups
all the information is concentrated in the ring $\pi_0^A(\bMUP)$. 
The non-equivariant homotopy groups $\pi_0^e(\bMUP)$
are a polynomial ring in countably many generators.
For cyclic groups of prime order, Kriz \cite{kriz-Z/p cobordism}
has described $\pi_0^{C_p}(\bMUP)$ as a pullback of two explicit
ring homomorphisms. 
For the cyclic group of order~2, Strickland \cite{strickland-MU involution} 
has turned this into an explicit
presentation of $\pi_0^{C_2}(\bMUP)$ as an algebra over $\pi_0^e(\bMUP)$. 

\index{subject}{global Thom spectrum|)}

\section{Equivariant bordism}
\label{sec:equivariant bordism}
\index{subject}{equivariant bordism|(}

In this section we recall equivariant bordism groups and their
relationship to the equivariant homology groups defined by the global Thom spectrum $\bmO$
introduced in Example \ref{eg:geometric global bordism}. 
The main result is Theorem \ref{thm:TP is iso} 
which says that when $G$ is isomorphic to a product of
a finite group and a torus, then the Thom-Pontryagin map
is an isomorphism from $G$-equivariant bordism to 
$G$-equivariant $\bmO$-homology.
Theorem \ref{thm:TP is iso} is usually credited to Wasserman
because it can be derived from his equivariant transversality
theorem \cite[Thm.\,3.11]{wasserman}; as far as I know, the only place
where the translation is spelled out
in detail is the unpublished part of Costenoble's thesis,
see \cite[Thm.\,11.1]{costenoble-thesis}. 
Wasserman's theorem is based on equivariant differential topology;
for finite groups, tom Dieck \cite[Satz 5]{tomDieck-OrbittypenII} gives
a different proof of Theorem \ref{thm:TP is iso}
based on the geometric and homotopy theoretic isotropy separation sequences.
We generalize tom Dieck's proof, translated into our present language,
with an emphasis on global aspects.

In Theorem \ref{thm:stable TP} we also present a localized version
of this result: the Thom-Pontryagin map
is an isomorphism from {\em stable} equivariant bordism to 
$\bmO[1/\tau]$-theory, without any restriction on the compact Lie group.
Given that $\bmO[1/\tau]$-theory is isomorphic to $\bMO$-theory
(by Corollary \ref{cor:MO is localized mO}),
this is equivalent to a result of Br{\"o}cker and Hook \cite[Thm.\,4.1]{broecker-hook}
that identifies stable equivariant bordism with equivariant $\bMO$-homology.

\medskip

The serious study of equivariant bordism groups was initiated
by the work \cite{conner-floyd} of Conner and Floyd.
We recall equivariant bordism as a homology theory for $G$-spaces $X$,
where $G$ is a compact Lie group. 
A {\em singular $G$-manifold}\index{subject}{singular $G$-manifold}
over $X$ is a pair $(M, h)$ consisting of a closed smooth $G$-manifold $M$
and a continuous $G$-map $h: M\to X$.
Two singular $G$-manifolds $(M, h)$ and $(M', h')$ are {\em bordant}
if there is a triple $(B, H,\psi)$ consisting of a compact smooth $G$-manifold $B$,
a continuous $G$-map $H:B\to X$ and an equivariant diffeomorphism
\[ \psi\ : \ M \cup M' \ \iso \ \partial B \   \]
such that $(H\circ\psi)|_M=h$  and to $(H\circ\psi)|_{M'}=h'$.

Bordism of singular $G$-manifolds over $X$ is an equivalence relation.
Reflexivity and symmetry are straightforward; transitivity is established by 
gluing two bordisms along a common piece of the boundary.
To get a smooth structure on the glued bordism that is compatible 
with the $G$-action one needs smooth equivariant collars; 
the existence of such collars is guaranteed 
by \cite[Thm.\,21.2]{conner-floyd}.

We denote by $\mathcal N_n^G(X)$\index{symbol}{$\mathcal N_n^G(X)$ - {equivariant geometric bordism}}
the set of bordism classes of $n$-dimensional singular $G$-manifolds over $X$.
This set becomes an abelian group under disjoint union. Every element $x$ of
$\mathcal N_n^G(X)$ satisfies $2 x=0$.
The groups $\mathcal N_n^G(X)$ are covariantly functorial
in continuous $G$-maps, by postcomposition.

\medskip

\begin{prop}\label{prop:bordism homotopy invariance}
  Let $G$ be a compact Lie group. 
  \begin{enumerate}[\em (i)]
  \item
    Let $\varphi,\varphi':X\to Y$ be equivariantly homotopic continuous $G$-maps.
    Then $\varphi_*=\varphi'_*$ as homomorphisms
    from $\mathcal N_n^G(X)$ to $\mathcal N_n^G(Y)$.
  \item For every $G$-weak equivalence
    $\varphi:X\to Y$ the induced homomorphism $\varphi_*:
    \mathcal N_n^G(X)\to\mathcal N_n^G(Y)$ is an isomorphism.
  \item  Let $\{X_i\}_{i\in I}$ be a family of $G$-spaces. 
    Then the canonical map
    \[ {\bigoplus}_{i\in I} \, \mathcal N_n^G(X_i) \ \to \ 
    \mathcal N_n^G\left({\coprod}_{i\in I}\, X_i\right)\]
    is an isomorphism.
\end{enumerate}
\end{prop}
\begin{proof}
(i) We let $H:X\times [0,1]\to Y$ be an equivariant homotopy from $\varphi$
to $\varphi'$ and $(M,h)$ a singular $G$-manifold over $X$.
Then $(M\times[0,1], H\circ (h\times[0,1]),\psi)$ is a bordism 
from $(M,\varphi h)$ to $(M,\varphi' h)$, 
where $\psi:M\cup M\to \partial(M\times [0,1])$
identifies one copy of $M$ with $M\times \{0\}$
and the other copy with $M\times \{1\}$.
So $\varphi_*[M,h]=[M,\varphi\circ h]=[M,\varphi'\circ h]=\varphi'_*[M,h]$.

(ii) For surjectivity of $\varphi_*$ we consider any singular $G$-manifold
$(M,g)$ over $Y$. Illman's theorem \cite[Cor.\,7.2]{illman}
shows that $M$ admits the structure of a $G$-CW-complex.
Since $\varphi$ is a $G$-weak equivalence there exists
a continuous $G$-map $h:M\to X$ such that $\varphi h$ 
is equivariantly homotopic to $g$. 
Part (i) then shows that $\varphi_*[M,h]=(\varphi h)_*[M,\Id_M]=g_*[M,\Id_M]=[M,g]$.

The argument for injectivity is similar. 
We consider a singular $G$-manifold $(M,h)$ over $X$
that represents an element in the kernel of $\varphi_*$. 
There is then a null-bordism $(B,H,\psi)$ of $(M,\varphi h)$.
By Illman's theorem \cite[Cor.\,7.2]{illman} and the discussion immediately
preceding it there is a $G$-CW-structure on $B$ for which
the boundary is an equivariant subcomplex. So the map 
$\psi:M\to B$ that identifies $M$ with the boundary of $B$ is a $G$-cofibration.
Since $\varphi$ is a $G$-weak equivalence, there exists
a continuous $G$-map $H':B\to X$ such that $H'\circ\psi =h$.
The triple $(B,H',\psi)$ thus witnesses that $[M,h]=0$.
Since $\varphi_*$ is a group homomorphism, it is injective. 

Property~(iii) holds because compact manifolds only have finitely 
many connected components, so all
continuous reference maps from singular manifolds or bordisms have image
in a finite union. 
\end{proof}

Now we state the key exactness property of equivariant bordism 
in the form of a Mayer-Vietoris sequence. 
The definition of the boundary map needs the existence of $G$-invariant
separating functions as provided by the following lemma.

\begin{lemma}\label{lemma:separate} 
Let $M$ be a compact smooth $G$-manifold, $C$ and $C'$
two disjoint, closed, $G$-invariant subsets of $M$, and
\[ s \ : \ \partial M \ \to \ \mR \]
a smooth $G$-invariant map such that
\[  C\cap\partial M \subseteq s^{-1}(0) \text{\qquad and\qquad}
 C'\cap\partial M \subseteq s^{-1}(1) \ . \]
Then there exists a smooth $G$-invariant extension 
$r:M\to\mR$ of $s$ such that $C\subseteq r^{-1}(0)$ and
$C'\subseteq r^{-1}(1)$.
\end{lemma}
\begin{proof}
Since $M$ is compact, hence normal, the Tietze extension theorem provides
a continuous map $r_0:M\to \mR$ that extends $s$ and satisfies
$C\subseteq r_0^{-1}(0)$ and $C'\subseteq r_0^{-1}(1)$.
We average $r_0$ to make it $G$-invariant, i.e., we define $r_1:M\to\mR$ by 
\[ r_1(x)\ = \ \int_G r_0(g\cdot x)\,  d g\ . \]
The integral is taken with respect to the normalized invariant measure
(Haar measure) on $G$. The new map $r_1$ is again continuous,
compare \cite[Ch.\ 0 Prop.\ 3.2]{bredon-intro}.
Since $r_0$ is already $G$-invariant on $C\cup C'\cup \partial M$,
the new map $r_1$ coincides with $r_0$ on this subset.
In particular, $r_1$ is smooth on $C\cup C'\cup\partial M$.
By \cite[VI Thm.\,4.2]{bredon-intro} we can then find 
a smooth $G$-invariant map $r:M\to \mR$ that coincides with $r_1$ 
on $C\cup C'\cup\partial M$; this is the desired separating function.
\end{proof}

\begin{construction}[Boundary map in equivariant bordism]\label{con:bordism boundary}
We define a boundary homomorphism for a Mayer-Vietoris sequence.
We let $X$ be a $G$-space and $A,B\subset X$ open $G$-invariant subsets
with $X=A\cup B$. Then a homomorphism
\[ \partial \ : \ \mathcal N_n^G(X)\ \to \ \mathcal N_{n-1}^G(A\cap B) \]
is defined as follows.

We let $(M,h)$ be a singular $G$-manifold that represents
a class in $\mathcal N_n^G(X)$. The sets $h^{-1}(X-A)$ and $h^{-1}(X-B)$ are 
$G$-invariant, disjoint closed subsets of $M$;
we let $r:M\to\mR$ be a $G$-invariant smooth separating function
as provided by Lemma \ref{lemma:separate}, i.e., 
such that $h^{-1}(X-A)\subseteq r^{-1}(0)$ and
$h^{-1}(X-B)\subseteq r^{-1}(1)$.
We let $t\in (0,1)$ be any regular value of $r$. Then
\[ M_t \ =  \ r^{-1}(t) \]
is a smooth closed $G$-submanifold of $M$ of dimension $n-1$ (possibly empty),
and $h_t=h|_{M_t}$ lands in $A\cap B$;
so $(M_t, h_t)$ is a singular $G$-manifold over $A\cap B$.
\end{construction}

\begin{prop}
In the situation above, the bordism class $[M_t,h_t]$ 
is independent of the choice of regular value $t$,
of the choice of separating function and of the representative
for the given class in $\mathcal N_n^G(X)$.
The resulting map
\[ \partial\ : \  \mathcal N_n^G(X)\ \to \ \mathcal N_{n-1}^G(A\cap B)  \ , \quad
[M,h]\ \longmapsto \ [M_t, h_t]\]
is a group homomorphism.  
\end{prop}
\begin{proof}
We let $t<t'$ be two regular values in $(0,1)$ of the separating function $r$.
Then
\[ (r^{-1}[t,t'],\, h|_{r^{-1}[t,t']},\,\text{incl} ) \]
is a bordism from $(r^{-1}(t),h|_{r^{-1}(t)})$ to $(r^{-1}(t'),h|_{r^{-1}(t')})$,
so the bordism class does not depend on the regular value.

Now we let $(M,h)$ and $(N,g)$ be two singular $G$-manifolds over $X$
in the same bordism class, and we let $(B,H,\psi)$ be a $G$-bordism
from $(M,h)$ to $(N,g)$.
We let $r:M\to\mR$ and $\bar r:N\to\mR$ be two $G$-invariant separating
functions. Lemma \ref{lemma:separate} lets us extend this data 
to a smooth $G$-invariant separating function
\[ \Psi\ : \ B \ \to \ \mR \]
such that $\Psi\circ\psi|_M=r$, $\Psi\circ\psi|_N=\bar r$,
\[ H^{-1}(X-A)\subseteq \Psi^{-1}(0) \text{\quad and\quad}
 H^{-1}(X-B)\subseteq \Psi^{-1}(1) \ .\]
We choose a simultaneous regular value $t\in(0,1)$ for $\Psi$, $r$ and $\bar r$. Then 
\[ (\Psi^{-1}(t),\, H|_{\Psi^{-1}(t)},\, \psi|_{r^{-1}(t)\cup \bar r^{-1}(t)}) \]
is a bordism from 
$(r^{-1}(t),h|_{r^{-1}(t)})$ to $(\bar r^{-1}(t),g|_{\bar r^{-1}(t)})$.
This shows at the same time that the bordism class is independent of
the choice of separating function and of the choice of representing
singular $G$-manifold.
Additivity of the resulting boundary map is then clear: a separating function
for a disjoint union can be taken as the union of separating functions
for each summand.
\end{proof}

Now we formulate the property that makes equivariant bordism 
a homology theory for $G$-spaces.

\begin{prop}\label{prop:excision in bordism} 
  Let $G$ be a compact Lie group, $X$ a $G$-space 
  and $A,B\subset X$ open $G$-invariant subsets
  with $X=A\cup B$. Let $i^A:A\cap B\to A$, $i^B:A\cap B\to B$,
  $j^A:A\to X$ and $j^B:B\to X$ denote the inclusions.
  Then the following sequence of abelian groups is exact:
  \[ \dots \to \mathcal N_n^G(A\cap B) \ \xra{(i^A_*,i^B_*)} \ 
  \mathcal N_n^G(A) \oplus \mathcal N_n^G(B) \ 
  \xra{\genfrac(){0pt}{}{j^A_*}{-j^B_*}} \  
  \mathcal N_n^G(X) \ \xra{\ \partial\ } \ 
  \mathcal N_{n-1}^G(A\cap B) \to \dots \]
\end{prop}

While Proposition \ref{prop:excision in bordism} has been frequently used,
there does not seem to be any published account with a complete proof
in the generality of compact Lie groups;
the argument is somewhat involved, but it proceeds along the lines of
the non-equivariant argument as given for example in
\cite[Prop.\,21.1.7]{tomDieck-algebraic topology}.
As input one needs that certain basic tools from differential topology 
generalize to the $G$-equivariant context, such as for example the
existence of equivariant collars and bicollars, and that 
non-equivariant smoothing of corners is compatible with $G$-actions.
We honor the tradition in this area and refrain from giving further details.

\medskip

We define the {\em reduced bordism group}\index{subject}{equivariant bordism!reduced}
of a based $G$-space $X$ as
\[ \widetilde{\mathcal N}_n^G(X) \ = \ 
\coker\left( \mathcal N_n^G(\ast)  \to \mathcal N_n^G(X) \right)\ ,\]
the cokernel of the homomorphism induced by the basepoint inclusion.
If $(M, h)$ is a singular $G$-manifold over $X$, 
then we use the notation
$\gh{M,h}$ for the class it represents in the reduced bordism group
$\widetilde{\mathcal N}^G_n(X)$.
The unique map $u:X\to\ast$ is retraction to the basepoint inclusion, so the map
\[ (\text{proj}, u_*)\ : \ 
\mathcal N_n^G(X) \ \to \ 
\widetilde{\mathcal N}_n^G(X)\oplus \mathcal N_n^G(\ast) \]
is an isomorphism. 
On the other hand, if we add a disjoint $G$-fixed
basepoint to an unbased $G$-space $Y$, then the composite 
\[ \mathcal N_n^G(Y) \ \xra{\text{incl}_*} \  \mathcal N_n^G(Y_+) \ \xra{\text{proj}} \ 
\widetilde{\mathcal N}_n^G(Y_+) \]
is an isomorphism.

\begin{construction}
We consider a continuous $G$-map $f:X\to Y$ and let
\[ C f\ = \ C X \cup_f Y \ = \ ( X\times[0,1]\cup_f Y ) / X\times\{0\} \]
denote its unreduced mapping cone. The two open sets
\[ A \ = \ X\times [0,1) / X\times \{0\}  \text{\qquad and\qquad} 
B \ = \ X\times (0,1]\cup_f Y \]
are $G$-invariant and together cover the mapping cone.
The intersection $A\cap B$ is homeomorphic to $X\times(0,1)$,
so the open covering has an associated boundary homomorphism
\[ \partial \ : \ \mathcal N_n^G(C f) \ \to \  \mathcal N_{n-1}^G(X\times (0,1)) \]
as in Construction \ref{con:bordism boundary}.
We take the cone point as the basepoint of $C f$; this is contained
in the subset $A$, so the map $\iota:\mathcal N_*(\ast)\to\mathcal N_*(C f)$
induced by the basepoint inclusion factors through 
$j^A_*:\mathcal N_n^G(A)\to \mathcal N_n^G(C f)$,
and the composite $\partial\circ \iota$ is trivial by exactness
of the excision sequence. The boundary map thus factors over 
the reduced bordism group. We define a `reduced boundary map' $\bar\partial$
as the composite\index{subject}{connecting homomorphism!in equivariant bordism}
\[ \widetilde{\mathcal N}_n^G(C f)\ \xra{\ \partial\ } \ 
\mathcal N_{n-1}^G(X\times(0,1)) \ \xra{\text{proj}_*} \ \mathcal N_{n-1}^G(X) \ .\]
\end{construction}

\begin{prop}\label{prop:bordism cone sequence} 
  Let $G$ be a compact Lie group and $f:X\to Y$ a continuous $G$-map. 
  Then the following sequence of abelian groups is exact:
  \[ \dots\ \to\ \mathcal N_n^G(X) \ \xra{\ f_*\ } \ 
  \mathcal N_n^G(Y)  \ \xra{\ i_*\ } \ 
  \widetilde{\mathcal N}_n^G(C f) \ \xra{\ \bar\partial\ } \ 
  \mathcal N_{n-1}^G(X)\  \to\ \dots \]
\end{prop}
\begin{proof}
  We use the open covering of the mapping cone $C f$ 
  as in the definition of the boundary map.
  In the diagram
  \[ \xymatrix@C=15mm{ 
    X  \ar[r]^-f \ar[d]_{x\mapsto (x,1/2)} & Y \ar[r]^-i  \ar[d]_-i &  C f \ar@{=}[d]  \\
    A\cap B  \ar[r]_-{\text{incl}} & B \ar[r]_-{\text{incl}} & C f} \]
  the right square commutes and the left square commutes up to equivariant homotopy.
  Moreover, all vertical maps are equivariant homotopy equivalences, so
  they induce isomorphisms in equivariant bordism,
  by Proposition \ref{prop:bordism homotopy invariance}.
  So the resulting diagram of bordism groups commutes:
  \[\xymatrix{
    \mathcal N_n^G(X) \ar[r]^-{f_*} \ar[d]_\iso &
    \mathcal N_n^G(Y) \ar[r]^-{i_*} \ar[d]_\iso^-{i_*}&
    \mathcal N_n^G(C f) \ar@{=}[d]\\
    \mathcal N_n^G(A\cap B) \ar[r] &
    \mathcal N_n^G(B) \ar[r] &
    \mathcal N_n^G(C f)  } \]
  Moreover, all vertical maps in this diagram are isomorphisms,
  so we can substitute $\mathcal N_n^G(X)$ and $\mathcal N_n^G(Y)$
  into the long exact excision sequence 
  of Proposition \ref{prop:excision in bordism}.
  Since $A$ is equivariantly contractible to the cone point,
  we can also replace the corresponding summand by the coefficient group,
  and the result is an exact sequence
  \[ \dots \to \mathcal N_n^G(X) \ \xra{\ (f_*,u_*)\ } \ 
  \mathcal N_n^G(Y) \oplus \mathcal N_n^G(\ast) \ \to \  
  \mathcal N_n^G(C f) \ \xra{\ \partial\ } \ 
  \mathcal N_{n-1}^G(X) \to \dots \]
  The sequence then remains exact if we divide out the summand
  $\mathcal N_n^G(\ast)$ and replace the absolute bordism group of $C f$
  by the reduced group $\widetilde{\mathcal N}_n^G(C f)$.
\end{proof}

If $f:A\to B$ is an h-cofibration of $G$-spaces, then 
the projection $q:C f\to B/A$
from the mapping cone to the quotient is a based equivariant homotopy equivalence.
So we can substitute $\widetilde{\mathcal N}_*^G(B/A)$
into the long exact mapping cone sequence of Proposition \ref{prop:bordism cone sequence} 
and obtain a long exact sequence of abelian groups:
    \[ \dots\ \to\ \mathcal N_n^G(A) \ \xra{\ f_*\ } \ 
    \mathcal N_n^G(B)  \ \xra{\ q_*\ } \ 
    \widetilde{\mathcal N}_n^G(B/A) \ \to \ 
    \mathcal N_{n-1}^G(A)\ \to \ \dots \]

The next proposition says, loosely speaking, that 
in a reduced bordism group the part of
a singular $G$-manifold that sits over the basepoint can be ignored. 

\begin{prop}\label{prop:relative ignores basepoint} 
Let $G$ be a compact Lie group
and $h:N\to X$ a singular $G$-manifold over a based $G$-space $X$. 
Let $V$ be a $G$-representation and $M$ a closed smooth $G$-manifold
such that $\dim(M)+\dim(V)=\dim(N)$.
Let $j:M\times D(V)\to N$ be a smooth $G$-equivariant embedding.
Suppose that $h$ sends $N-j(M\times\mathring{D}(V))$ to the basepoint.
Then 
\[ \gh{N,h}\ = \ \gh{M\times S(\mR\oplus V), f} 
\text{\quad in\quad $\widetilde{\mathcal N}_n^G(X)$,} \]
where $f:M\times S(\mR\oplus V)\to X$ is defined by
\[ f(m, (x,v))\ = \
\begin{cases}
h(j(m,v)) & \text{ if $x\leq 0$, and}\\  
\quad \ast & \text{ if $x\geq 0$.}  
\end{cases}\]
\end{prop}
\begin{proof}
We define a $G$-space by
\[ B \ = \ ( N\times[-1,0]  \ \cup\  N\times [0,1] ) / \sim\ , \]
where the equivalence relation identifies the two copies
of $(n,0)$ for every $n\in N- j(M\times\mathring{D}(V))$.
The group $G$ acts by $g\cdot[n,t]=[g n,t]$.
The space $B$ is a topological $(n+1)$-manifold whose boundary
consists of three disjoint parts that we now parametrize.
There are two obvious embeddings
\[ \psi\ , \ \psi' \ : \ N \ \to \ B \text{\qquad by\qquad}
\psi(n)=[n,-1] \text{\qquad and\qquad} \psi'(n)=[n,1]  \]
as the two endpoints in the direction of the internal $[-1,1]$.
We define another embedding 
\begin{align*}
 i \ : \ M\times S(\mR\oplus V)\ &\to \ B \text{\qquad by}\\
 i(m,(x,v))\ &= \
\begin{cases}
\quad  [(j(m,v), 0]^\text{left} & \text{ for $x\leq 0$, and}\\
\quad  [(j(m,v), 0]^\text{right}  & \text{ for $x\geq 0$.}
\end{cases}  
\end{align*}
Here the superscripts `left' and `right' indicate whether we refer to the point
$[n,0]$ as the image of $(n,0)$ in $N\times [-1,0]$ or in $N\times [0,1]$.
The manifold boundary of $B$ is then the disjoint union of the images
of $\psi, \psi'$ and $i$.

The topological manifold $B$ admits a smooth structure for which the given $G$-action
is smooth and such that the embeddings $\psi, \psi'$ and $i$ are smooth;
the construction involves `smoothing of corners' 
(also called `straightening of angles') near the image of $j(M\times S(V))\times 0$
and is explained for example
in Construction 15.10.3 of tom Dieck's textbook \cite{tomDieck-algebraic topology}.
Tom Dieck has no group actions around, 
but we also have to ensure that the given $G$-action on $B$
is smooth. This can be arranged by insisting that the collars used 
in \cite[15.10.3]{tomDieck-algebraic topology}
are $G$-equivariant collars, which is possible for example
by \cite[Thm.\,21.2]{conner-floyd}.

Now we define a continuous $G$-map
$H:B \to X$ by $H(n,t)=h(n)$ on $N\times[-1,0]$
and as the constant map to the base point of $X$ on $N\times[0,1]$.
Then
\[ H\circ\psi \ = \ h \text{\qquad and\qquad } H\circ i \ = \ f \ ,\]
so the bordism $(B,H,\psi+ \psi'+ i)$ witnesses the relation
\[ [N,h] \ = \ [N,H\circ\psi'] \ + \  [M\times S(\mR\oplus V),f]  \]
in the unreduced bordism group $\mathcal N_n^G(X)$. 
Since $H\circ\psi'$ is constant to the basepoint of $X$,
the class $[N,H\circ\psi']$ vanishes in the reduced bordism group 
$\widetilde{\mathcal N}_n^G(X)$; this proves the claim. 
\end{proof}

The equivariant bordism groups come with natural products, given by the biadditive maps
\[ \times \ : \ \mathcal N_m^G(X)\ \times\ \mathcal N_n^G(Y) \ \to \ 
\mathcal N_{m+n}^G(X\times Y)\ , \quad
[M,h]\times [N,g]\ = \ [M\times N, h\times g]\ .\]
These products are suitably associative, commutative and unital.
The product pairing descends to a pairing on reduced bordism groups
if the $G$-spaces $X$ and $Y$ are based.
Indeed, the composite
\begin{align*}
\mathcal N_m^G(X) \tensor \mathcal N_n^G(Y) \ &\xra{\ \times\ } \   
\mathcal N_{m+n}^G(X\times Y) \ \xra{\ q_*\ } \ 
\mathcal N_{m+n}^G(X\sm Y) \ \xra{\text{proj}} \
\widetilde{\mathcal N}_{m+n}^G(X\sm Y) 
\end{align*}
annihilates the image of $\mathcal N_m^G(\ast)\tensor \mathcal N_n^G(Y)$
and the image of $\mathcal N_m^G(X)\tensor\mathcal N_n^G(\ast)$,
where $q:X\times Y\to X\sm Y$ is the quotient map;
so the composite factors uniquely over a homomorphism
\[ \sm \ : \ \widetilde{\mathcal N}_m^G(X)\ \otimes\ \widetilde{\mathcal N}_n^G(Y) 
\ \to \ \widetilde{\mathcal N}_{m+n}^G(X\sm Y)\ .\]

We will frequently use certain distinguished bordism classes
associated to representations.
We let $V$ be an $m$-dimensional representation of a compact Lie group $G$.
Stereographic projection is a $G$-equivariant 
homeomorphism\index{subject}{stereographic projection}
\[ \Pi_V\ : \ S(\mR\oplus V)\ \xra{\ \iso\ } \  S^V 
\ , \quad (x,v) \ \longmapsto \ \frac{v}{1-x} \]
from the unit sphere of $\mR\oplus V$ to the one-point compactification $S^V$.
We define a reduced $G$-bordism class over $S^V$ by\index{symbol}{$d_{G,V}$ - {equivariant bordism class of a representation}}
\begin{equation}  \label{eq:define_d_G}
 d_{G,V} \ = \ \gh{ S(\mR\oplus V),\Pi_V } \quad \in\ \widetilde{\mathcal N}_m^G(S^V)\ .  
\end{equation}
The following multiplicativity property is expected, but not entirely obvious.

\begin{prop}\label{prop:d-classes associative}
Let $V$ and $W$ be two representations of a compact Lie group $G$.
Then the relation
\[  d_{G,V}\sm\ d_{G,W}\ =\ d_{G,V\oplus W}  \]
holds in $\widetilde{\mathcal N}_{m+n}^G(S^{V\oplus W})$.
\end{prop}
\begin{proof} 
We define a `distorted' version 
$\tau_V : S(\mR\oplus V) \to S^V$
of the stereographic projection as the composite
\[  S(\mR\oplus V) \ \xra{\ \Pi_V\ } \ S^V \ \xra{\ J \ }\ S^V\ , \]
where the second map $J$ is given by
\[ J(v)\ = \  
\begin{cases}
\frac{v}{1-|v|} & \text{for $|v| < 1$, and}\\
\ \infty  & \text{for $|v|\geq 1$.}
\end{cases} \]
The map $J$ is equivariantly based homotopic to the identity of $S^V$,
so the pair $(S(\mR\oplus V),\tau_V)$ 
is another representative for the bordism class $d_{G,V}$.
The map
\begin{align*}
 j\ : \ V\oplus W \ &\to \ S(\mR\oplus V)\times S(\mR\oplus W)\\
j(v,w)\ &= \ \left( \frac{|v|^2-1}{|v|^2+1}, \frac{2 v}{|v|^2+1},
  \frac{|w|^2-1}{|w|^2+1},\frac{2 w}{|w|^2+1}
 \right)
\end{align*}
is a smooth $G$-equivariant embedding.
The map $j$ is the product of the inverse stereographic projections 
$\Pi_V^{-1}\times\Pi_W^{-1}:S^V\times S^W\to S(\mR\oplus V)\times S(\mR\oplus W)$,
restricted to $V\oplus W$.

If a quadruple $(x,v,y,w)$ is in the image $j(D(V\oplus W))$ of the unit disc, 
then in particular $x\leq 0$ and $y\leq 0$.
Equivalently, the points $(x,v,y,w)$ of $S(\mR\oplus V)\times S(\mR\oplus W)$
that are in the complement of $j(D(V\oplus W))$ have $x > 0$ or $y > 0$.
So the map
\[  q\circ(\tau_V\times\tau_W)\ :\ S(\mR\oplus V)\times S(\mR\oplus W)\ \to\ S^{V\oplus W}\]
sends the complement of $j(\mathring{D}(V\oplus W))$ to the basepoint at infinity,
where $q:S^V\times S^W\to S^{V\oplus W}$ is the projection.
Proposition \ref{prop:relative ignores basepoint} thus shows that
\begin{align*}
 q_*(d_{G,V}&\sm d_{G,W}) \ = \ 
q_*(\gh{S(\mR\oplus V),\tau_V}\sm \gh{S(\mR\oplus W),\tau_W}) \\ 
&= \ \gh{ S(\mR\oplus V)\times S(\mR\oplus W),\ q\circ(\tau_V\times\tau_W) }  \
= \ \gh{ S(\mR\oplus V\oplus W),\ f }  \ ,
\end{align*}
where
$f:S(\mR\oplus V\oplus W)\to S^{V\oplus W}$ is defined by
\[ f(x,v,w)\ = \
\begin{cases}
q(\tau_V\times\tau_W(j(v,w))) & \text{ if $x\leq 0$, and}\\  
\qquad \infty & \text{ if $x\geq 0$.}  
\end{cases}\]
We claim that for all $(x,v,w)\in S(\mR\oplus V\oplus W)$ the relation
\[ f(x,v,w)\ \ne \  \Pi_{V\oplus W}(-x, - v, -w) \]
holds in $S^{V\oplus W}$.
Assuming this claim, we can finish the proof as follows:
since the $G$-map $\Pi_{V\oplus W}^{-1}\circ f$
never takes a point to its antipode,
the linear homotopy between  $\Pi_{V\oplus W}^{-1}\circ f$ and the
identity in the ambient vector space $\mR\oplus V\oplus W$ can be
normalized to land in the unit sphere. 
This yields an equivariant based homotopy between $\Pi_{V\oplus W}^{-1}\circ f$ 
and the identity of $S(\mR\oplus V\oplus W)$.
Hence $f$ is equivariantly based homotopic 
to the stereographic projection $\Pi_{V\oplus W}$,
and so $(S(\mR\oplus V\oplus W),f)$ represents the bordism class $d_{G,V\oplus W}$.

It remains to prove the claim.
The only point of $S(\mR\oplus V\oplus W)$ 
that $\Pi_{V\oplus W}$ sends to the point at infinity is $(1,0,0)$.
Since $f(-1,0,0)=(0,0)$, the claim is true for all
$(x,v,w)$ such that $f(x,v,w)=\infty$.
It remains to consider those tuples for which $f(x,v,w)\ne \infty$,
which means that  $x < 0$ and 
$|v|< 1$ and $|w| < 1$.
On such points the map $f$ is given by
\begin{align*}
 f(x,v,w) \ &= \
 q \left( \tau_V\left(\frac{|v|^2-1}{|v|^2+1},\frac{2 v}{|v|^2+1}\right),\ 
\tau_W\left(\frac{|w|^2-1}{|w|^2+1},\frac{2 w}{|w|^2+1}\right)\right)\\
&= \ q(J(v),J(w))\ , 
\end{align*}
whereas
\[ \Pi_{V\oplus W}(-x,-v,-w) \ = \  
  \left( \frac{-v}{1+x},\quad \frac{-w}{1+x}    \right) \ .
\]
If these two expressions were the same, then
\[ \frac{v}{1-|v|}\ = \ \frac{-v}{1+x}\text{\qquad and\qquad}
 \frac{w}{1-|w|}\ = \ \frac{-w}{1+x}\ . \]
For $v\ne 0$ this implies
\[ |v|-1\ = \ \ 1+x \ ,\]
which is impossible since $|v|<1$ and $x\geq -1$.
If $w\ne 0$ we obtain the same kind of contradiction.
The final case is $(x,v,w)=(-1,0,0)$, in which case
\[ f(1-,0,0) \ = \ (0,0) \ \ne \ \infty \ = \ \Pi_{V\oplus W}(1,0,0)\ . \qedhere \]
\end{proof}

We shall now recall that if $G$ acts trivially on $V$, 
then exterior multiplication by $d_{G,V}$ is an isomorphism,
the suspension isomorphism in reduced equivariant bordism.
In general, however, the class $d_{G,V}$ is not invertible
and exterior multiplication by $d_{G,V}$ need not be an isomorphism.
The theory obtained by formally inverting all
the classes $d_{G,V}$ is {\em stable equivariant bordism},
to which we return in Remark \ref{rk:local mO and stable bordism} below.
The projection $p:C f\to X\sm S^1$ from the reduced mapping to the suspension
was defined in \eqref{eq:define_t}.

\begin{prop}\label{prop:suspension iso in bordism}
If $G$ acts trivially on $V$ and $X$ is a cofibrant based $G$-space,
then the  exterior product map
\[ -\sm  d_{G,V} \ : \   \widetilde{\mathcal N}_n^G(X) \ \to \  
\widetilde{\mathcal N}_{n+|V|}^G(X\sm S^V) \]
is an isomorphism.  
For every continuous $G$-map $f:X\to Y$ between based $G$-spaces
the connecting homomorphism in the mapping cone sequence equals the composite
\[  \widetilde{\mathcal N}_n^G(C f)\ \xra{\ p_* \ } \
 \widetilde{\mathcal N}_n^G(X\sm S^1)\ \xra{(-\sm d_{G,\mR})^{-1}} \
 \widetilde{\mathcal N}_{n-1}^G(X)\ . \]
\end{prop}
\begin{proof}
We start with the special case $V=\mR$.
We apply Proposition \ref{prop:bordism cone sequence}
to the map $f:X\to \ast$ to a one-point $G$-space.
The cone of this map is
\[ X^\diamond \ = \ X\times [0,1] / \sim \ , \]
the unreduced suspension of $X$, where $X\times\{0\}$ and $X\times\{1\}$ 
are collapsed to one point each.
Since $X$ has a $G$-fixed point, 
the map $f_*:\mathcal N_*^G(X)\to\mathcal N_*^G(\ast)$
is a split epimorphism. So the long exact sequence provided by 
Proposition \ref{prop:bordism cone sequence} reduces to a short exact sequence:
\[ 0  \ \to  \  \widetilde{\mathcal N}_{n+1}^G(X^\diamond) \ \xra{\ \bar\partial\ } \ 
\mathcal N_n^G(X)\ \xra{\ f_*\ } \ 
\mathcal N_n^G(\ast)  \ \to \ 0  \]
Since $X$ is cofibrant in the based sense, the projection 
\[ \psi \ : \ X^\diamond \ \to \ X\sm S^1 \ , \quad
q[x,s]\ = \ x\sm \frac{2s-1}{s(1-s)}   \]
that collapses $\{x_0\}\times [0,1]$ is an equivariant homotopy equivalence.

Then  the composite
\[ \widetilde{\mathcal N}_{n+1}^G(X\sm S^1) \ \xra[\iso]{\ \psi_*^{-1}\ }\ 
 \widetilde{\mathcal N}_{n+1}^G(X^\diamond) \ \xra{\ \bar\partial\ }\ 
\mathcal N_n^G(X) \ \xra{\text{proj}}\ 
\widetilde{\mathcal N}_n^G(X)  \]
is an isomorphism.
We claim that the relation
\[ \text{proj}(\bar\partial(\psi_*^{-1}(x\sm d_{G,\mR}))) \ = \ x\]
holds for all classes $x\in\widetilde{\mathcal N}_n^G(X)$.
Since $\text{proj}\circ\bar\partial\circ \psi_*^{-1}$ 
is an isomorphism, so is smash product with the class $d_{G,\mR}$.

This relation, in turn, is a consequence of the geometric origin
of the class $d_{G,\mR}$, the product in bordism  and the boundary map.
In more detail, we suppose that $x=\gh{M,h}$ for a singular
$G$-manifold $(M,h)$ over $X$.
We define a continuous map $H:M\times S(\mR\oplus\mR)\to X^\diamond$ by
\[ H(m,(x,y)) \ = \
\begin{cases}
 \ [h(m), (y+1)/2] & \text{ for $x \leq 0$, and }\\
 \quad [x_0, (y+1)/2] & \text{ for $x \geq 0$.}
\end{cases} \]
Then the following square commutes up to $G$-equivariant homotopy:
\[ \xymatrix{ 
M\times S(\mR\oplus\mR) \ar[r]^-H \ar[d]_{h\times\Pi_\mR} &
X^\diamond\ar[d]^\psi\\
X\times S^1 \ar[r]_-q & X\sm S^1
} \]
Hence
\[ 
\psi_*\gh{ M\times S(\mR\oplus\mR), H }\ = \ 
\gh{ M\times S(\mR\oplus\mR), q\circ(h\times\Pi_\mR) } \ = \ \gh{M,h}\sm d_{G,\mR}\ ,
\]
and thus
\[ \bar\partial(\psi_*^{-1}(\gh{M,h}\sm d_{G,\mR})) \ = \
\bar\partial\gh{ M\times S(\mR\oplus\mR), H } \ .\]
To calculate this geometric boundary we use the smooth separating function
\[ r\ : \ M\times S(\mR\oplus\mR) \ \to \ [0,1] \ , \quad
r(m,(x,y))\ =(y+1)/2\ .\]
Then $1/2$ is a regular value of this separating function, and 
the preimage over this regular value is 
$r^{-1}(1/2)=M\times\{(1,0), (-1,0)\}$, two disjoint copies of $M$.
The function $H$ takes the copy $M\times\{(1,0)\}$ to the basepoint
of $X$, so this copy does not contribute to the reduced bordism group.
The restriction of $H$ to the other copy  $M\times\{(-1,0)\}$ 
is the original map $h$, so we obtain
\[  \bar\partial[M\times S(\mR\oplus\mR), H]  \ = \ 
 \gh{M\times \{(1,0),(-1,0)\}, H|_{M\times \{(1,0),(-1,0)\}}}\ \equiv \ \gh{M,h}\]
in the reduced bordism group of $X$.

The general case now follows easily. Since the claim is true for $V=\mR$,
it also holds for $V=\mR^n$ by the associativity of the smash product pairing
and the classes $d_{G,\mR}$, compare Proposition \ref{prop:d-classes associative}.
If $G$ acts trivially on $V$, then it is equivariantly isomorphic to $\mR^n$
for some $n$.
\end{proof}

The bordism theories $\mathcal N_*^G$ for different compact Lie groups
are related by geometrically defined restriction and induction maps.
Every continuous group homomorphism $\alpha:K\to G$ is automatically smooth
(see for example \cite[Prop.\,I.3.12]{broecker-tomDieck}),
and thus induces a restriction homomorphism
\[ \alpha^* \ : \ \mathcal N_n^G(X)\ \to \ \mathcal N_n^K(\alpha^*(X)) \ , \quad
 [M,h]\ \longmapsto \ [\alpha^*(M),\alpha^*(h)] \]
by restricting all actions along $\alpha$.
Restriction maps preserve the distinguished classes \eqref{eq:define_d_G}
in the sense that $\alpha^*(d_{G,V})=d_{K,\alpha^*(V)}$.
The product in equivariant bordism is compatible with restriction maps in the sense that
\[ \alpha^*(x\times y) \ = \ \alpha^*(x)\times \alpha^*(y)\ .\]

For every closed subgroup $H$ of $G$ and every closed smooth $H$-manifold $M$
of dimension $n$, the induced space $G\times_H M$
has a preferred smooth structure making it 
a smooth closed $G$-manifold of dimension $d+n$
with $d=\dim(G/H)=\dim(G)-\dim(H)$.
Indeed, since the diagonal $H$-action on $G\times M$ is smooth and free,
there is a unique smooth structure on $G\times_H M$
such that the projection $G\times M\to G\times_H M$
is a submersion, see for example \cite[Thm.\,15.3.4]{tomDieck-algebraic topology},
and the $G$-action is smooth with respect to this smooth structure.
We can also apply $G\times_H-$ to bordisms, so this
gives a well-defined induction homomorphism
\[ G\times_H- \ : \ \mathcal N_n^H(Y) \ \to  \ \mathcal N_{d+n}^G(G\times_H Y)\ , \quad 
[M,h]\ \longmapsto \ [G\times_H M, G\times_H h]\ .    \]
For $Y=\ast$ we can compose the induction map with the effect of
the projection $G\times_H\ast \to \ast$ and arrive at an induction homomorphism
on coefficient groups
\[ \ind_H^G \ : \  \mathcal N_n^H(\ast) \ \to  \ \mathcal N_{d+n}^G(\ast)\ , \quad
[M]\ \longmapsto \ [G\times_H M]\ .\]
The induction map $\ind_H^G$ is compatible with inflations,
in the sense of the formula
\[ \alpha^*\circ \ind_H^G \ = \ \ind_L^K \circ (\alpha|_L)^*\ : \ 
\mathcal N_n^H(\ast) \ \to \ \mathcal N_{d+n}^K(\ast) \]
for a continuous epimorphism $\alpha:K\to G$, where $L=\alpha^{-1}(H)$.
If $H$ has finite index in $G$, then the induction $\ind_H^G$ preserves
the dimension, and then it satisfies the double coset formula.
So for fixed $n\geq 0$ the coefficient groups $\mathcal N_n^G(\ast)$
{\em almost} form a global functor; the only missing structure are the
transfer maps for closed inclusions that are not of finite index.

Multiplication, restriction and induction satisfy a reciprocity relation.
We let $H$ be a closed subgroup of $G$, 
$X$ an $H$-space and $Y$ a $G$-space.
Then the following diagram commutes:
\[ \xymatrix@C=13mm@R=6mm{ 
\mathcal N_m^H(X)\times \mathcal N_n^G(Y)
\ar[d]_{\Id\times\res^G_H}\ar[r]^-{(G\times_H-)\times\Id} &
\mathcal N_{d+m}^G(G\times_H X)\times \mathcal N_n^G(Y)
\ar[d]^{\times} \\
\mathcal N_m^H(X)\times \mathcal N_n^H(\res^G_H(Y)) \ar[d]_{\times}&
\mathcal N_{d+m+n}^G((G\times_H X) \times Y) \ar[d]^{\chi_*} \\
\mathcal N_{m+n}^H(X\times \res^G_H(Y)) \ar[r]_-{G\times_H-}&
\mathcal N_{d+m+n}^G(G\ltimes_H( X\times \res^G_H(Y))) } \]
Here 
\[ \chi \ : \ 
(G\times_H X) \times Y \ \to \ 
G\ltimes_H( X\times \res^G_H(Y) ) \ , \quad 
([g,x],y)\ \longmapsto \ [g, (x, g^{-1} y)] \]
is the $G$-equivariant shearing isomorphism.
The proof of the commutativity is straightforward from the definitions,
using the shearing diffeomorphism for equivariant manifolds.
If we specialize to $X=Y=\ast$ and postcompose with the projection
to the one-point $G$-space, we obtain the reciprocity formula 
\[ \ind_H^G(x\times \res^G_H(y)) \ = \ \ind_H^G(x)\times y \]
for classes $x\in\mathcal N_m^H(\ast)$ and $y\in\mathcal N_n^G(\ast)$
in the coefficient groups.

\medskip

The next proposition shows that the distinguished bordism classes $d_{G,V}$ 
measure the failure of the Wirthm{\"u}ller isomorphism in equivariant bordism,
see Remark \ref{rk:Wirthmuller fails} below.
We consider a closed subgroup $H$ of a compact Lie group $G$ 
and continue to write
\[ L\ = \ T_{e H}(G/H) \]
for the tangent $H$-representation.
For an $H$-space $Y$ the $H$-equivariant continuous  map
\[ l_Y \ : \ G\times_H Y \ \to \ Y_+\sm S^L \]
was defined in Construction \ref{con:Wirthmuller map}.

\begin{prop}\label{prop:geometric Wirthmuller} 
For every closed subgroup $H$ of a compact Lie group $G$ and every
$H$-space $Y$, the composite
\[ \mathcal N_n^H(Y) \ \xra{\ G\times_H-\ } \
 \mathcal N_{n+d}^G(G\times_H Y) \ \xra{\ \res^G_H\ } \
\mathcal N_{n+d}^H(G\times_H Y) \ \xra{\ (l_Y)_*\ } \
\widetilde{\mathcal N}_{n+d}^H(Y_+\sm S^L) \]
is exterior multiplication by the class $d_{H,L}\in \widetilde{\mathcal N}_d^H(S^L)$,
where $d=\dim(G/H)$.
\end{prop}
\begin{proof}
We let $(M,h)$ be a singular $H$-manifold that
represents a class in $\mathcal N_n^H(Y)$.
The $G\times_H h:G\times_H M\to G\times_H Y$ represents the class $G\times_H[M,h]$.
As in the construction of the collapse map $l_Y$ we choose 
a slice around $1\in G$ orthogonal to $H$, i.e., a smooth embedding
$s: D(L) \to  G$ satisfying $s(0)=1$ and
\[ s(h\cdot l) \ = \ h \cdot s(l) \cdot h^{-1}  \]
for all $(h,l)\in H\times D(L)$ and such that the differential at 0 of the composite
\[ D(L)\ \xrightarrow{\ s\ }\ G \ \xrightarrow{\text{proj}}\ G/H \]
is the identity of $L$.
The slice property implies that the map 
\[  \bar s \ : \ D(L)\times H \ \to \ G \ , \quad (l,h)\ \longmapsto \ s(l)\cdot h \]
is a tubular neighborhood of $H$ inside $G$. Moreover, this 
embedding is equivariant for the action of $H^2$, acting on source
and target by
\[ (h_1,h_2)\cdot (l,h)\ = \ (h_1 l,\, h_1 h h_2^{-1}) 
\text{\quad respectively \quad}
(h_1,h_2)\cdot g\ = \  h_1 g h_2^{-1} \ .\]
The map 
\[ l_H^G \ : \ G_+ \ \to \ S^L\sm H_+ \]
was defined as the $H^2$-equivariant collapse map 
with respect to the tubular neighborhood $\bar s$.
So explicitly,
\[ l_H^G(g) \ = \ 
\begin{cases}
 ( l / (1-|l|) )\sm h & \text{ if $g=s(l)\cdot h$ with $(l,h)\in D(L)\times H$, and }\\
\quad \infty & \text{ if $g$ is not in the image of $\bar s$.}
\end{cases} \]
We obtain a smooth $H$-equivariant embedding
\[ j\ :\ M\times D(L) \ \to\ G\times_H M \text{\qquad by\qquad}
j(m,l) \ = \ [ s(l), m] \ ,\]
where $H$ acts diagonally on the left.
The composite 
\[ 
G\times_H M \ \xra{\ G\times_H h\ } \ G\times_H Y\ \xra{\ l_H^G\sm_H Y_+\ }\  
(S^L\sm  H_+)\sm_H Y_+  \ \iso \ Y_+\sm S^L  \]
sends the complement of $j(M\times D(L))$ to the basepoint at infinity.
Proposition \ref{prop:relative ignores basepoint} (for $(H,L)$ instead of $(G,V)$)
shows that then
\begin{align*}
((l_Y)_*\circ\res^G_H\circ (G\times_H-))[M,h] \ &= \ 
\gh{G\times_H M,\ l_Y\circ (G\times_H h)} \\
&= \ \gh{M\times S(\mR\oplus L), f} 
\end{align*}
in the reduced bordism group of $Y_+\sm S^L$,
where $f:M\times S(\mR\oplus L)\to Y_+\sm S^L$ is defined by
\[ f(m,(x,l))\ = \
\begin{cases}
( l_Y\circ(G\times_H h))(j(m,l)) & \text{ if $x\leq 0$, and}\\  
\qquad \infty  & \text{ if $x\geq 0$.}  
\end{cases}\]
So $f$ equals the composite
\[ M\times S(\mR\oplus L)\ \xra{h\times \Psi}\  Y\times S^L\ \xrightarrow{\ q\ }\ 
Y_+\sm S^L\]
where 
\[  \Psi(x,l)\ = \ l/(1-|l|)  \]
for $x\leq 0$, and $\Psi(x,l)=\infty$ for $x\geq 0$.
The map $\Psi:S(\mR\oplus L)\to S^L$ is homotopic, in the equivariant based sense,
to the stereographic projection $\Pi_L$.
We can thus conclude that
\begin{align*}
((l_Y)_*\circ\res^G_H\circ (G\times_H-))&[M,h] \ = \ 
\gh{M\times S(\mR\oplus L), f} \\  
&= \  [M,h]\sm \gh{S(\mR\oplus L),\Psi}  \ = \ [M,h]\sm d_{H,L}\ .\qedhere
\end{align*}
\end{proof}

\begin{rk}[Failure of the Wirthm{\"u}ller isomorphism
in equivariant bordism]\label{rk:Wirthmuller fails}
\index{subject}{Wirthm{\"u}ller isomorphism!in equivariant bordism}
\index{subject}{Wirthm{\"u}ller map!in equivariant bordism}
As we now explain, there is no general Wirthm{\"u}ller isomorphism
in equivariant bordism, i.e., the Wirthm{\"u}ller map 
\[ \Wirth_H^G\ = \ (l_Y)_*\circ\res^G_H \ : \ \mathcal N_*^G(G\times_H Y)\ \to \ 
\widetilde{\mathcal N}_*^H(Y_+\sm S^L) \]
is {\em not} in general bijective for closed subgroups $H$ of $G$. 
Here $Y$ is an $H$-space and $l_Y$ is the $H$-equivariant collapse map
defined in Construction \ref{con:Wirthmuller map}.
This immediately shows that there {\em cannot} be a natural isomorphism,
compatible with restriction to subgroups,
between equivariant bordism and the equivariant homology theory represented by
any orthogonal $G$-spectrum. Indeed, Proposition \ref{prop:geometric Wirthmuller}
says that the composite
\[ \Wirth_H^G\circ(G\times_H -)\ : \ \mathcal N_n^H(Y) \ \to
\ \widetilde{\mathcal N}_{n+d}^H(Y_+\sm S^L) \]
is exterior multiplication by the class $d_{H,L}\in \widetilde{\mathcal N}_d^H(S^L)$.
Since the induction map $G\times_H -$ is an isomorphism,
the Wirthm{\"u}ller map is bijective if any only if the class $d_{H,L}$
acts invertibly on $\mathcal N_*^H(Y)$.
In this sense, the distinguished bordism class $d_{H,L}$ 
precisely measures the failure of the Wirthm{\"u}ller isomorphism in equivariant bordism.
One can show that the class $d_{H,L}$ is {\em not} 
invertible as soon as $H$ acts non-trivially on $L$,
i.e., in this situation the class $d_{H,L}$ 
does {\em not} generate $\widetilde{\mathcal N}_*^H(S^L)$
as an $\mathcal N_*^H$-module.

If $G$ is a product of a finite group and a torus, then $H$ acts
trivially on $L$ for every closed subgroup $H$ of $G$, and hence
the Wirthm{\"u}ller isomorphism {\em does hold} for such compact Lie groups $G$.
And in fact, for this class of groups, the Thom-Pontryagin construction
provides an natural isomorphism between equivariant bordism and
the equivariant homology theory represented by the Thom spectrum $\bmO$,
compare Theorem \ref{thm:TP is iso}  below.
\end{rk}

Our argument to compare the geometric bordism theory with the equivariant homotopy
groups of the Thom spectrum $\bmO$ is based on the isotropy separation sequence.
We will now identify the `geometric fixed point term'
in the isotropy separation sequence for equivariant bordism.

\begin{construction}
As before we let $E\Pc$ be a universal space for the family of proper closed subgroups
of $G$. So $E\Pc$ is a cofibrant $G$-space without $G$-fixed points, 
and $(E\Pc)^H$ is contractible for every proper closed subgroup of $G$.\index{subject}{universal space!for the family of proper subgroups}
We let $\tilde E\Pc$ denote the unreduced suspension of $E\Pc$.
Then $(\tilde E\Pc)^G=S^0$ consists of the two cone points,
and $(\tilde E\Pc)^H$ is contractible for every proper subgroup of $G$.

We recall the `geometric fixed point' homomorphism,\index{subject}{geometric fixed point homomorphism!in equivariant bordism}
see \eqref{eq:define_geometric_fix} below, that identifies the
reduced bordism group $\widetilde{\mathcal N}^G_n(\tilde E\Pc)$
in terms of non-equivariant bordism groups.
Loosely speaking, the class $\Phi_\text{geom}\gh{M,h:M\to\tilde E\Pc}$ remembers the
bordism class of the part of the fixed point manifold
that lies over the fixed point~0,
together with the normal data of the embedding into $M$.
As before, we use the abbreviation
\[ G r_j^{G,\perp}\ = \ \left( G r_j(\Uc_G^\perp) \right)^G \]
for the space of $j$-dimensional $G$-invariant subspaces of $\Uc_G^\perp=\Uc_G-(\Uc_G)^G$.
We first define a fixed point map 
\[ \bar\Phi\ : \ \mathcal N_n^G(\tilde E\Pc)\ \to \ 
{\bigoplus}_{j\geq 0}\, \mathcal N_{n-j} (G r_j^{G,\perp} ) \]
into the unreduced, non-equivariant bordism group as follows. 
We let $(M,h:M\to\tilde E\Pc)$ be an
$n$-dimensional singular $G$-manifold.
Since $G$ acts smoothly, $M^G$ is a disjoint union of smooth submanifolds
of varying dimensions, compare \cite[VI Cor.\,2.5]{bredon-intro}.
Since $\tilde E\Pc$ has exactly two fixed points, and both are isolated,
$h$ must map every path component of the fixed point manifold $M^G$
to either~0 or $\infty$. We denote by
\[ M_0^G \ = \ M^G \cap h^{-1}(0) \]
the union of those components of $M^G$ that lie over the point~0.
We let $M^{(j)}$ be the union of all $(n-j)$-dimensional components of $M^G_0$.

The Mostow-Palais embedding theorem \cite{mostow-embeddings, palais-imbedding}
provides a smooth $G$-equivariant embedding $i: M\to V$, 
for some $G$-representation $V$; we can assume that $V$ is in fact
a subrepresentation of the complete universe $\Uc_G$. 
Then for every fixed point $x\in M^{(j)}$, 
\[ (d i)(T_x (M^{(j)})) \ = \  ( (d i)(T_x M) )^G \ ,\]
i.e., the tangent space inside $M^{(j)}$ `is' the $G$-fixed part
of the tangent space in $M$. 
So we can define a continuous map
\[ \nu_j \ : \ M^{(j)}\ \to \ ( G r_j(V^\perp))^G\  \xra{\text{incl}}\ 
 G r_j^{G,\perp} \]
by sending a fixed point $x\in M^{(j)}$ to the orthogonal complement
of $(d i)(T_x (M^{(j)}))$ inside $(d i)(T_x M)$.
By its very construction, the map $\nu_j$ classifies 
the normal bundle of the inclusion $M^{(j)}\to M$.
The geometric fixed point map is then given by
\[ \bar\Phi[M,h] \ = \ {\sum}_{j=0}^n \ [M^{(j)}, \nu_j]\ .\]
Since the image of the map $i_*:\mathcal N_n^G(\ast)\to \mathcal N_n^G(\tilde E\Pc)$
is concentrated over the basepoint $\infty$, it is annihilated by
the map $\bar\Phi$. So $\bar\Phi$ factors
over the reduced bordism group of $\tilde E\Pc$ as a homomorphism
\begin{equation}  \label{eq:define_geometric_fix}
 \Phi_\text{geom}\ :  \ \widetilde{\mathcal N}_n^G(\tilde E\Pc)\ \to \ 
{\bigoplus}_{j\geq 0}\, \mathcal N_{n-j} (G r_j^{G,\perp} )\ .  
\end{equation}
The fact that the geometric fixed point map \eqref{eq:define_geometric_fix}
is an isomorphism is ubiquitous in calculations of equivariant bordism groups,
and it goes back to Conner and Floyd \cite{conner-floyd}.
In the classical literature, the reduced equivariant bordism groups 
of $\tilde E\Pc$ usually appear in a different guise, namely as
the groups $\mathcal N_*^G[\All,\Pc]$ of bordism classes of
smooth compact $G$-manifolds with boundary, but where there
are no fixed points on the boundary, 
compare Remark \ref{rk:reinterpret geometric fixed} below.
\end{construction}

\begin{prop}\label{prop:geometric isotropy}
For every compact Lie group $G$ and every $n\geq 0$
the geometric fixed point map \eqref{eq:define_geometric_fix}
is an isomorphism.
\end{prop}
\begin{proof}
The proof is by explicit geometric constructions.
We start with surjectivity.
We consider a non-equivariant closed smooth $(n-j)$-manifold $N$
and a continuous map $f:N\to G r_j^{G,\perp} $.
The space $G r_j^{G,\perp}$ is the filtered colimit of its
closed subspaces $(G r_j(V))^G$ for $V\in s(\Uc_G^\perp)$.
Since $N$ is compact, the image of $f$ lands in $(G r_j(V))^G$ 
for some finite-dimensional $G$-representation $V$ with $V^G=0$.
The space $(G r_j(V))^G$ is a smooth manifold, 
and by smooth approximation 
we can assume without loss of generality that $f$ is a smooth map.
We define a closed smooth $n$-dimensional manifold by
\[ M\ = \ \{ (n,x,v)\in N \times S(\mR\oplus V) \ |\  
v\in f(n) \} \ . \]
Another way to say this is that $M$ is a double of
the unit disc bundle of the pullback 
of the tautological $j$-plane bundle along $f:N\to (G r_j(V))^G$. 
The group $G$ acts smoothly on $M$ by
$g\cdot (n,x,v) = (n, x, g v)$.

At this point it will be convenient to use a specific model 
for the space $\tilde E\Pc$, namely
the unit sphere in $\mR\oplus\Uc_G^\perp$, 
compare Example \ref{eg:S of U^perp is E tilde}.
We can then define a continuous $G$-map
\[ h \ : \ M  \ \to \ \tilde E\Pc\ =\ S(\mR\oplus\Uc_G^\perp) \text{\qquad by\qquad} 
 h(n,x,v)\ = \ (x,v) \ . \]
So the pair $(M,h)$ represents a bordism class 
\[ [M,h] \quad \in\  \mathcal N_n^G(\tilde E\Pc)\ . \]
Since $V^G=0$, the $G$-fixed points of $M$
are a disjoint union of two copies of $N$ embedded as
$ N\times (-1,0)$ respectively $N\times (1,0)$, 
and $h$ maps one copy to each of the two fixed points of $\tilde E\Pc$.
The normal bundle of the copy of $N$ over the non-base point
is the bundle classified by the original map $f$, so we obtain
\[ \Phi_\text{geom}[ M,h]\ = \ [N,f]\ . \]
This shows that every class in $\mathcal N_{n-j}(G r_j^{G,\perp})$
is in the image of the geometric fixed point map, and so
the map \eqref{eq:define_geometric_fix} is surjective.

Now we consider a reduced equivariant bordism class $\gh{M,h}$ over $\tilde E\Pc$
in the kernel of the fixed point map \eqref{eq:define_geometric_fix}.
Being in the kernel of $\Phi_\text{geom}$ means that 
for every $0\leq j\leq n$ there is a non-equivariant null-bordism $B_j$ with
$\partial B_j= M^{(j)}$ and a continuous map $F_j:B\to G r_j^{G,\perp}$
whose restriction to $M^{(j)}$ classifies the normal bundle
of $M^{(j)}$ inside $M$.
As in the first part of the proof we can compress $F_j$
to a $G$-map $B\to (G r_j(V))^G$ for some finite-dimensional $G$-subrepresentation
$V$ of $\Uc_G^\perp$ and replace this factorization by a homotopic smooth
map $F_j:B\to (G r_j(V))^G$.
We use this data to `cut out' the fixed points $M_0^G$ from $M$ by replacing a
tubular neighborhood by the sphere bundles of the maps $F_j$; this produces
a new singular $G$-manifold over $\tilde E\Pc$, bordant to $(M,h)$,
that has no more fixed points over~0. 

The construction is done separately
and disjointly over each of the components $M^{(j)}$ of the fixed points $M^G_0$.
To simplify the exposition we restrict to the special case where 
the fixed points $M^G_0$ are of constant dimension $n-j$
(i.e., all other components $M^{(i)}$ for $i\ne j$ are empty).
Then we proceed as follows.
We let $\nu$ denote the normal bundle of $M^{(j)}$ inside $M$.
The equivariant tubular neighborhood theorem provides
a smooth $G$-equivariant embedding
\[ \psi\ : \ D(\nu) \ \to \ M \]
of the unit disc bundle of $\nu$ whose composite with 
the zero section $s:M^{(j)}\to D(\nu)$ is the inclusion,
see for example \cite[VI Thm.\,2.2]{bredon-intro}.
In particular, the image of $\psi$ is a closed $G$-invariant 
tubular neighborhood of $M^{(j)}$.
By shrinking the neighborhood, if necessary, we can make its image
disjoint from all other components of $M^G$ except $M^{(j)}$.
The bundle arising by pulling back the tautological bundle 
along $F:B\to (G r_j(V))^G$ has its own disc bundle with total space
\[ D(F) \ = \ \{ (b,v)\in B\times V \ | \ v\in F(b),\, |v|\leq 1 \}\ .\]
The sphere bundle $S(F)$ is then a smooth compact $G$-manifold with boundary
\[ \partial(S(F))\ = \ S(F|_{\partial B}) \ = \ S(\nu)\ .\]
Now we form the $G$-manifold
\[ \bar M \ = \ ( M- \psi(\mathring{D}(\nu))\cup_{S(\nu)} S(F)\ ,\]
where the gluing uses the restriction $\psi: S(\nu) \to M$
of the tubular neighborhood to the sphere bundle.
A smooth structure on $\bar M$ is provided by choices of $G$-equivariant
collars of $\psi(S(\nu))$ inside $M- \psi(\mathring{D}(\nu))$
and of $S(\nu)$ inside $S(F)$.

The boundary of the disc bundle of $F$ decomposes as
\[\partial(D(F)) \ = \ S(F) \cup_{S(F|_{\partial B})} D(F|_{\partial B}) 
\ = \ S(F) \cup_{S(\nu)} D(\nu) \ . \]
By equivariant smoothing of corners the disc bundle $D(F)$
can be given a smooth structure such that the given $G$-action is smooth
and that the embeddings of $D(\nu)$ and $S(F)$ into $D(F)$ are smooth.
We define a bordism as the $G$-space
\[  W\ = \ \left( M\times [0,1]\right)\cup_{D(\nu)} D(F) \]
where the gluing is along $(\psi,1):D(\nu)\to M\times [0,1]$.
The space $W$ is a topological $(n+1)$-manifold whose boundary
is the union of two disjoint parts that we now parametrize.
An obvious embedding is given by
\[ \psi\ : \ M \ \to \ W \ ,\quad  \psi(m)=[m,0] \ .\]
A second embedding 
\[ i \ : \ \bar M\ = \ ( M-\psi(\mathring{D}(\nu)))\cup_{S(\nu)} S(F) \ \to \ W \]
identifies $M-\psi(\mathring{D}(\nu))$
with $M-\psi(\mathring{D}(\nu))\times 1$
and includes the sphere bundle $S(F)$ into the disc bundle $D(F)$.
The boundary $\partial W$ is then the disjoint union of the images of $\psi$ and $i$.

The topological manifold $W$ admits a smooth structure for which the given $G$-action
is smooth and such that the embeddings $\psi$ and $i$ are smooth;
in the non-equivariant version, this is explained
in Construction 15.10.3 of \cite{tomDieck-algebraic topology}.
To ensure that the given $G$-action on $W$
is smooth, we must insist that the collars involved in the construction
are $G$-equivariant collars, which is possible for example
by \cite[Thm.\,21.2]{conner-floyd}.

Now we have to arrange the equivariant reference maps to $\tilde E\Pc$.
In fact, this extra data goes along for the ride, as the only homotopical information
it encodes is a decomposition of the $G$-fixed points into
two disjoint subspaces, the preimages of the two fixed points of $\tilde E\Pc$. 
In more detail, we let $W$ be a compact smooth $G$-manifold, possibly with boundary.
Then $W$ admits the structure of a finite $G$-CW-complex
by Illman's triangulation theorem \cite[Thm.\,7.1]{illman}. 
By the proof of Proposition \ref{prop:geometric as fixed points}, 
the fixed point map
\[ (-)^G \ : \ \map^G(W,\tilde E\Pc)\ \to \ \map(W^G,\{0,\infty\})\]
is a weak equivalence and Serre fibration.
A continuous map to the discrete space $\{0,\infty\}$
is equivalent to a decomposition into two disjoint open subsets.
So altogether we conclude that for every disjoint union decomposition $W^G=A\cup B$
there is a $G$-map $h:W\to \tilde E\Pc$ with $h(A)=\{0\}$
and $h(B)=\{\infty\}$, and any two such maps are equivariantly homotopic.
We use this property three times, namely for the $G$-manifolds $W$,
$M$ and $\bar M$.

The $G$-fixed points of $W$ are a disjoint union of
$ M_\infty^G\times [0,1]$ and $(M^{(j)}\times [0,1])\cup_{M^{(j)}\times 1} B$.
So there is a $G$-map $H: W\to \tilde E\Pc $
such that $H( M_\infty^G\times [0,1])=\{\infty\}$
and $H(( M^{(j)}\times [0,1])\cup_{M^{(j)}\times 1} B )=\{0\}$.
The triple $(W,H,\psi+i)$ is then a bordism that witnesses the relation
$\gh{M,H \psi}=\gh{\bar M,H i}$ in the group
$\widetilde{\mathcal N}_n^G(\tilde E\Pc)$.
The map $H i:\bar M\to \tilde E\Pc$ sends all of $\bar M^G$ to the fixed point $\infty$,
so $H i$ is equivariantly homotopic to the constant map 
with value $\infty$. By homotopy invariance, 
the class $\gh{\bar M,H i}$ thus vanishes in the reduced bordism group.
On the other hand, the $G$-map $H\psi$ agrees with the original map $h$
on the fixed points $M^G$. So $H\psi$ is equivariantly homotopic to $h$,
and homotopy invariance yields
\[\gh{M,h}\ = \ \gh{M,H \psi} \ = \ \gh{\bar M,H|_{\bar M}} \ = \ 0  \]
in the group $\widetilde{\mathcal N}_n^G(\tilde E\Pc)$.
This shows that the map $\Phi_\text{geom}$ is injective.
\end{proof}

\begin{rk}\label{rk:reinterpret geometric fixed}
In the classical literature, the reduced equivariant bordism groups 
of $\tilde E\Pc$ usually appear in an isomorphic form, namely as
the groups $\mathcal N_*^G[\All,\Pc]$ of bordism classes of
smooth compact $G$-manifolds with boundary, but where there
are no fixed points on the boundary. A homomorphism
\begin{equation}  \label{eq:P,All_to_EP}
 \mathcal N_*^G[\All,\Pc]\ \to \ \widetilde{\mathcal N}_n^G(\tilde E\Pc)  
\end{equation}
is defined as follows. We let $M$ be a smooth compact $G$-manifold without
$G$-fixed points on the boundary.
The double of $M$ is the smooth closed $G$-manifold
\[ D M \ = \ M\cup_{\partial M} M \]
obtained by gluing two copies of $M$ along their boundary.
A smooth structure in the neighborhood of the gluing locus
is provided by a choice of $G$-equivariant collar
(see \cite[Thm.\,21.2]{conner-floyd}).
By Illman's theorem \cite[Cor.\,7.2]{illman},
$D M$ admits the structure of a finite $G$-CW-complex.
Since the original manifold $M$ had no $G$-fixed points on the boundary,
$(D M)^G$ is the disjoint union of two copies of $M^G$, one from each of 
the two copies of $M$ in the double.
There is thus a continuous $G$-map
$h:D M\to \tilde E\Pc$, unique up to equivariant homotopy,
that takes the `left' copy of $M^G$ to the fixed point~0
and the `right' copy of $M^G$ to the fixed point $\infty$.
The pair $(D M,h)$ is then a singular $G$-manifold over $\tilde E\Pc$,
and it represents a bordism class 
\[ [D M,h] \quad \in\  \mathcal N_n^G(\tilde E\Pc)\ . \]
The fact that the map \eqref{eq:P,All_to_EP} is an isomorphism
is Satz~3 in \cite{tomDieck-OrbittypenI}.

For finite groups, Stong shows in \cite[Cor.\,5.1]{stong-unoriented finite}
that the geometric fixed point map
\[  \mathcal N_*^G[\All,\Pc]\ \to \ 
{\bigoplus}_{j\geq 0}\, \mathcal N_{n-j} (G r_j^{G,\perp} ) \]
is an isomorphism.
Combined with tom Dieck's isomorphism \eqref{eq:P,All_to_EP}
this gives a different proof of Proposition \ref{prop:geometric isotropy}
in the special case of finite groups.
I was unable to find a convenient reference for 
Proposition \ref{prop:geometric isotropy} in the 
generality of compact Lie groups, which is the main reason for including a proof.
\end{rk}

\begin{eg}[Bordism of manifolds with involution]
We look at the geometric isotropy separation sequence 
in the simplest non-trivial case of the two-element group $G=C_2=C$,
the case originally considered by  Conner and Floyd \cite[Thm.\,28.1]{conner-floyd}.
In this case $\Pc$ consists only of the trivial subgroup,
and so $\mathcal N_*^C(E\Pc)=\mathcal N_*^C(E C)$
is the bordism ring of manifolds with free $C$-action.
If $C$ acts freely and smoothly on $M$, then we can form the smooth manifold
\[ [-1,1] \times_C M \ = \ ([-1,1] \times M ) / (x,m)\sim (-x,\tau m)  \]
with $C$-action by $\tau \cdot [x, m] = [-x,\tau m]$;
the boundary of this manifold is equivariantly diffeomorphic to the original
manifold $M$. So every $C$-manifold with free action is null-bordant,
and thus the forgetful map $\mathcal N_*^C(E\Pc)\to\mathcal N_*^C$
is zero. 
The long exact mapping cone sequence of the $C$-map $E C\to\ast$
(see Proposition \ref{prop:bordism cone sequence}) thus 
decomposes into short exact sequences.
For every compact Lie group $G$, the map
\[ \mathcal N_n^G(E G) \ \to \ \mathcal N_n( B G) \ ,\quad
[M,h]\ \longmapsto \ [G\bs M, G\bs h ] \]
is an isomorphism from the bordism group of $G$-manifolds with free
action to the non-equivariant bordism group of the classifying space $B G=G\bs E G$.
So in the case at hand, $\mathcal N_*^C(E\Pc)=\mathcal N_*^C(E C)$
is isomorphic to $\mathcal N_*(\mR P^\infty)$.

We use Proposition \ref{prop:geometric isotropy}
to replace the group $\widetilde{\mathcal N}_*^C(\tilde E C)$
by the direct sum of non-equivariant bordism groups.
Since there is only one non-trivial irreducible
$C$-representation, the 1-dimensional sign representation,
every linear subspace of $\Uc_C^\perp$ is $C$-invariant.
Hence $G r_j^{C,\perp}$ is just a Grassmannian of $j$-planes
in an infinite dimensional $\mR$-vector space,
hence a classifying space of the orthogonal group $O(j)$.
Altogether, this yields a short exact sequence
\[ 0\ \to \ \mathcal N_*^C \ \xra{\ \Phi_\text{geom}\ }\ 
{\bigoplus}_{j\geq 0}\, \mathcal N_{*-j} ( G r_j(\mR^\infty) ) \ \xra{\ J\ } \ 
\mathcal N_{*-1}(\mR P^\infty) \ \to \ 0 \ .\]
The map $J$ is defined as follows. Given an $(n-j)$-dimensional 
singular manifold $(F,\eta:F\to G r_j(\mR^\infty))$,
we factor the continuous map $\eta$ through a finite-dimensional Grassmannian
$F\to G r_j(\mR^m)$, choose a homotopic smooth approximation
$\eta':F\to G r_j(\mR^m)$ and let $P(\eta')$ denote the 
projectivized bundle of the pullback of the tautological $j$-plane bundle along $\eta'$. 
The total space of $P(\eta')$ is then a smooth closed $(n-1)$-manifold
equipped with a map to $\mR P^\infty$ that classifies 
the tautological line bundle over $P(\eta')$; this data represents the class $J[F,\eta]$.

One can deduce from the above short exact sequence that 
$\mathcal N_*^C$ is free as a module over the non-equivariant 
bordism ring $\mathcal N_*$.
This is in fact true much more generally: Stong shows 
in \cite[Prop.\,9.4]{stong-unoriented finite}
that $\mathcal N_*^G(X)$ is free as an $\mathcal N_*$-module
for every finite group $G$ and every $G$-space $X$.
The paper \cite{alexander} by Alexander exhibits an explicit 
$\mathcal N_*$-basis of $\mathcal N_*^C$.
Some basic $C$-bordism classes are represented by the projective spaces
$\mR P^n$ equipped with the involution
\[ \tau\cdot [x_0:x_1:\ldots:x_n] \ = \ [-x_0:x_1:\ldots:x_n] \ .\]
We denote by $ y_n = [\mR P^n,\tau]$
the bordism class of this $C$-manifold in $\mathcal N_n^C$.
An $\mathcal N_*$-linear map
\[ \Gamma \ : \ \mathcal N_*^C \ \to \ \mathcal N_{*+1}^C \]
of degree~1 is given by sending the class of
a manifold $M$ with involution $\tau:M\to M$ to the manifold
\[ S(\mC)\times_C M \ = \ (S(\mC)\times M ) / (z,m)\sim (-z,\tau m) \ .\]
So $S(\mC)\times_C M$ is diffeomorphic to the
mapping torus of the involution $\tau$. The involution on this manifold is by
\[ \tau \ : \ S(\mC)\times_C M \ \to \ S(\mC)\times_C M \ ,\quad
\tau\cdot [z,m] \ = \ [\bar z, \tau m]\ . \]
In our present notation, the operator can also be expressed as
\[ \Gamma \ = \ \res^{O(2)}_C \circ \ind_{O(1)\times O(1)}^{O(2)}\circ p^*\ , \]
where we embed $C$ into $O(2)$ as a reflection, and
where $p:O(1)\times O(1)\to C$ is the epimorphism with kernel $e\times O(1)$.

Alexander shows in \cite[Thm.\,1.1]{alexander} that the multiplicative unit~1
together with the classes
\[ \Gamma^n(y_{i_1}\cdot\ldots\cdot y_{i_r}) \]
for all $n\geq 0$, $r\geq 1$ and $i_j\geq 2$ form a basis of $\mathcal N_*^C$
as a module over $\mathcal N_*$.
Besides the trivial group and $C$, 
the equivariant bordism groups have been calculated for various
finite abelian groups, see \cite{beem-2mod4, beem-C_4, beem-rowlett, firsching}. 
\end{eg}

\index{subject}{global Thom spectrum|(}

Now we work our way towards the equivariant Thom-Pontryagin construction
that assigns to every equivariant bordism class over a $G$-space $X$
an equivariant homology class in $\bmO_*^G(X_+)$.
We break the construction up into two steps, 
and we first discuss the {\em normal class}, 
a basic invariant associated to a closed  smooth $G$-manifold.

\index{subject}{normal class!of a $G$-manifold|(}
\begin{construction}[Normal class of a $G$-manifold]
In Example \ref{eg:defined MGr} we defined the ultra-commutative 
ring spectrum $\bMGr$ that consists of the Thom spaces of the tautological
vector bundles over Grassmannians.
To every smooth closed $G$-manifold $M$ we associate
a {\em normal class}\index{subject}{normal class!of a $G$-manifold}
\[ \td{M}\ \in \  \bMGr_0^G(M_+)\ .\]
This class records the equivariant homotopical information in the
stable normal bundle of $M$, and
it is the geometric input for the Thom-Pontryagin map
to equivariant $\bmO$-homology.
If $M$ has dimension $m$, then the class lives in the homogeneous summand
$\bMGr^{[-m]}$ of $\bMGr$.

The construction starts from the 
Mostow-Palais embedding theorem \cite{mostow-embeddings, palais-imbedding}
that provides a smooth $G$-equivariant embedding $i: M\to V$, 
for some $G$-representation $V$. 
We can assume without loss of generality that $V$ is a subrepresentation
of the chosen complete $G$-universe $\Uc_G$.
We use the inner product on $V$ to define the normal bundle $\nu$ 
of the embedding at $x\in M$ by
\[ \nu_x \ = \ V - (d i)(T_x M)\ ,\]
the orthogonal complement of the image of the tangent space $T_x M$ in $V$.
By multiplying with a suitably large scalar, if necessary, 
we can assume that the embedding is {\em wide}\index{subject}{wide embedding}
in the sense that the exponential map
\[ D(\nu) \ \to \ V \ , \quad (x,v)\ \longmapsto \ i(x) \, +\,  v \]
is injective on the unit disc bundle of the normal bundle, 
and hence a closed $G$-equivariant embedding.
The image of this map is a tubular neighborhood of radius~1 around $i(M)$,
and it determines a $G$-equivariant Thom-Pontryagin collapse map
\[ c_M\ : \ S^V \ \to \  T h(G r(V) )\sm M_+ \ = \ \bMGr(V)\sm M_+ \]
as follows: every point outside of the tubular neighborhood is sent
to the basepoint, and a point $i(x)+ v$, for $(x,v)\in D(\nu)$,
is sent to 
\[ 
c_M(i(x) +  v) \ = \ \left( \frac{v}{1-|v|},\,  \nu_x  \right) \sm  x \ .
\]
The normal class $\td{M}$ is the homotopy class of the collapse map $c_M$.
\end{construction}

\begin{prop} The normal class of a smooth closed $G$-manifold
is independent of the choice of wide embedding into a $G$-representation.
\end{prop}
\begin{proof}
If we enlarge the embedding $i:M\to V$ by postcomposition
with a direct summand embedding $(0,-):V\to U\oplus V{}$,
then the collapse map associated to the composite embedding $(0,-)\circ i$
is equivariantly homotopic to the composite
\[  S^{U\oplus V} \ \xra{S^U\sm c_M} \  S^U\sm \bMGr(V) \sm M_+
\ \xra{\sigma_{U,V} } \ \bMGr(U\oplus V)\sm M_+ \ .\]
So the resulting class in $\bMGr_0^G(M_+)$ does not change.
Two classes based on two different wide embeddings $i:M\to V$ and $j:M\to W$
can be compared by passing to $V\oplus W$; in this larger representation,
the map
\[ M\times [0,1]\ \to \ V\oplus W \ , \quad 
(m,t) \ \longmapsto \ ( t\cdot i(m),\ (1-t)\cdot j(m)) \]
is a smooth isotopy through wide embeddings. This isotopy induces a
homotopy between the two collapse maps and shows that altogether
the normal class $\td{M}$ is independent of the wide embedding.
\end{proof}

Part~(ii) of the following proposition refers to an external
multiplication morphism
\begin{align*}
   \mu_{A,B}\ :\ (\bMGr\sm A_+)\sm (\bMGr\sm B_+)\ \iso\  
&(\bMGr\sm \bMGr)\sm (A\times B)_+\\ 
&\xra{\mu\sm (A\times B)_+}\ \bMGr\sm(A\times B)_+ \ , 
\end{align*}
where $A$ and $B$ are two $G$-spaces
and $\mu:\bMGr\sm\bMGr\to\bMGr$ is the multiplication morphism.
In part~(iv) we consider the $k$-th power $M^k$ of a $G$-manifold $M$
as a $(\Sigma_k\wr G)$-manifold via the action
\[ (\sigma;\,g_1,\dots,g_k) \cdot (x_1,\dots, x_k) \ = \ 
(g_{\sigma^{-1}(1)}x_{\sigma^{-1}(1)},\dots,g_{\sigma^{-1}(k)} x_{\sigma^{-1}(k)}) \ . \]
The morphism $j=j^\mR_{\bMGr}:\bMGr\to\sh\bMGr$ that shows up in part~(v) was defined 
in \eqref{eq:define j^V for MGr}.

\begin{prop}\label{prop:normal class properties}
Let $G$ be a compact Lie group, and let $M$ and $N$ be smooth closed $G$-manifolds. 
  \begin{enumerate}[\em (i)]
  \item Let $i^1:M\to M\cup N$ and $i^2:N\to M\cup N$ denote the inclusions
    into a disjoint union. Then the relation
    \[ \td{M \cup N} \ = \ i^1_*\td{M} \ +\ i^2_*\td{N} \]
    holds in the group $\bMGr^G_0((M\cup N)_+)$.
  \item 
    The relation
    \[ \td{M\times N}\ = \  (\mu_{M,N})_*(\td{M}\times \td{N}) \]
    holds in the group $\bMGr_0^G((M\times N)_+)$.
  \item
    For every continuous homomorphism $\alpha:K\to G$ of compact Lie groups
    and every smooth closed $G$-manifold $M$ the relation
    \[ \td{ \alpha^* M} \ = \ \alpha^*\td{M}      \]
    holds in $\bMGr^K_0(\alpha^*(M)_+)$.
  \item For every $k\geq 0$, the relation\index{subject}{power operation!for $G$-manifolds}
    \[ \td{ M^k} \ = \ P^k\td{M} \]
    holds in the group $\bMGr_0^{\Sigma_k\wr G}(M^k)$.
  \item Let $B$ be a compact smooth $G$-manifold with boundary $\partial B$.
    Then the class $\td{\partial B}$ is in the kernel of the homomorphism
    \[ (j\sm \iota_+)_* \ : \ \bMGr_0^G(\partial B_+) \ \to \ (\sh \bMGr)_0^G(B_+)  \ ,\]
    where $\iota:\partial B\to B$ is the inclusion.
\end{enumerate}
\end{prop}
\begin{proof}
(i) 
We let $p^1:(M\cup N)_+\to M_+$ and $p^2:(M\cup N)_+\to N_+$ 
denote the two projections.
We choose any wide smooth equivariant embedding $i:M\cup N\to V$
and observe that the composite
\[  S^V \ \xra{c_{M\cup N}} \  \bMGr(V)\sm (M\cup N)_+ \ \xra{\bMGr(V)\sm p^1} \ 
 \bMGr(V)\sm M  \]
is on the nose the collapse map for $M$ based on the restriction of
the embedding $i$ to $M$.
We then obtain
\begin{align*}
  p^1_* \td{M \cup N}  \ = \ \td{M}
\ = \ p^1_*( i^1_*\td{M}  +i^2_*\td{N}) \ ;
\end{align*}
the second relation uses that $p^1\circ i^1$ is the identity,
$p^2\circ i^1$ is the trivial map and $p^1_*$ is additive.
The analogous argument shows that 
$p^2_* \td{M \cup N}\ = \ p^2_*( i^1_*\td{M}  +i^2_*\td{N})$.
Since equivariant homotopy groups are additive on wedges, the map 
\[ (p^1_*,p^2_*)\ : \ \bMGr^G_0((M\cup N)_+)\ \to \ 
\bMGr^G_0(M_+)\times \bMGr^G_0( N_+) \]
is bijective, and this proves the claim.

(ii) We choose wide smooth equivariant embeddings
\[ i \ : \ M\ \to \ V \text{\qquad and\qquad} j \ : \ N \ \to \ W \]
into $G$-representations.
The product map
\[ i\times j \ : \ M\times N \ \to \ V\oplus W \]
is then another wide smooth equivariant embedding that we use for the
Thom-Pontryagin construction of $M\times N$.
The normal bundle of $i\times j$ is the
exterior direct sum of the normal bundles of $i$ and $j$.
The unit disc $D(V\oplus W)$ of the direct sum is contained in the
product $D(V)\times D(W)$ of the unit discs, so the exponential tubular neighborhood
for $i\times j$ is contained in the product of the 
exponential tubular neighborhoods for $i$ and $j$.
The collapse map
\[ S^{V\oplus W}\ \xra{\ c_{M\times N}\ } \ 
\bMGr(V\oplus W)\sm (M\times N)_+ \]
is equivariantly homotopic to the composite
\begin{align*}
   S^V\sm S^W\ \xra{c_M\sm c_N} \ 
&( \bMGr(V) \sm M_+)\sm  (\bMGr(W)\sm  N_+) \\
&\xra{\ \mu_{M,N}^{V,W}\ } \  \bMGr(V\oplus W)\sm (M\times N)_+\ .
\end{align*}
This shows the desired relation.
Part~(iii) is straightforward from the definitions.

(iv) We choose a wide smooth equivariant embedding $i:M\to V$
into a $G$-representation. Then
\[ i^k \ : \ M^k \ \to \ V^k \]
is a $(\Sigma_k\wr G)$-equivariant wide smooth embedding
that we use to calculate the class $\td{M^k}$.
The collapse map 
\[  c_{M^k}\ : \ S^{V^k} \ \to \  T h(G r(V^k) )\sm M^k_+  \]
based on $i^k$ is $(\Sigma_k\wr G)$-equivariantly homotopic to the composite
\begin{align*}
  S^{V^k} \ &\xra{(c_M)^{\sm k}} \  (\bMGr(V)\sm M_+)^{\sm k} \\ 
&\xra{\text{shuffle}} \  \bMGr(V)^{\sm k}\sm M^k_+ \ 
\xra{\mu_{V,\dots,V}^{\bMGr}\sm M^k_+} \  \bMGr(V^k)\sm M^k_+ \ .\end{align*}
This latter composite represents the power operation $P^k\td{M}$,
so altogether this shows the desired relation.

(v) We let $C( j\sm B_+)$ denote the mapping cone of the morphism 
$j\sm B_+:\bMGr\sm B_+\to\sh \bMGr\sm B_+$.
We define a relative normal class
\[  \td{B}^{\rel} \ \in \ \pi_1^G( C (j\sm B_+)) \]
such that the relation
\begin{equation}\label{eq:bordism_boundary_relation}
\partial( \td{B}^{\rel} ) \ = \ (\bMGr\sm\iota_+)_* \td{\partial B}
\end{equation}
holds in the group $\pi_0^G(\bMGr\sm B_+)$, where
$\partial$ is the connecting homomorphism \eqref{def:boundary} 
in the long exact homotopy group sequence of the mapping cone.
Because two consecutive maps in the long exact homotopy
group sequence compose to zero, we can then conclude that
\begin{align*}
 (j\sm \iota_+)_* \td{\partial B}\ 
= \   (j\sm B_+)_*( \partial (\td{B}^{\rel}))\  = \ 0\ .
\end{align*}
It remains to construct the class $\td{B}^{\rel}$ and establish 
the relation \eqref{eq:bordism_boundary_relation}.
We choose an equivariant collar, i.e., a smooth $G$-equivariant embedding
\[ c\ : \ \partial B\times [0,1) \ \to \ B \]
such that $c(-,0):\partial B\to B$ is the inclusion and the image of $c$ is
an open neighborhood of the boundary inside $B$.
Then we choose a smooth function 
\[ \kappa\ : \ [0,1]\ \to \ [0,1] \]
that is the identity on $[0,1/3]$, identically~1 on $[2/3,1]$
and whose restriction to $[0,2/3)$ is injective.
We define the smooth function
\[  \psi\ : \ [0,1)\ \to \ [0,1) \text{\qquad by\qquad}
\psi(t)\ = \ \frac{\kappa(t)-1}{t-1} \ ; \]
then $\psi(t)=1$ for $t\in [0,1/3]$ and $\kappa(t)=0$ for $t\in[2/3,1)$.

The Mostow-Palais embedding theorem
provides a wide smooth $G$-equivariant embedding $j: B\to V$, 
for some $G$-representation $V$. Then the smooth $G$-map
\[ i \ : \ B \ \to \  V\oplus \mR \]
defined by
\[  i(b) = 
\begin{cases}
(\psi(t)\cdot  j(x)+ (1-\psi(t))\cdot j(b),\,\kappa(t)) & \text{for $b=c(x,t)$ with $(x,t)\in \partial B\times [0,1)$,}\\
\qquad (j(b),1) & \text{for $b\not\in c(B\times [0,1))$,}
\end{cases}
\]
is a new wide smooth equivariant embedding which satisfies
\[ i(\partial B)\ \subset \  V\times \{0\} \]
and which is `orthogonal to $V$ near the boundary', i.e.,
the set $U=c(\partial B\times [0,1/3))$ is
an open neighborhood of $\partial B$ in $B$, and
\[  i(U) \ = \  i(\partial B) \times [0,1)\ .    \]
Since the embedding $i$ is wide, the exponential map
\[ D(\nu) \ \to \ V\oplus \mR \ , \quad (x,v,t)\ \longmapsto \ i(x) \, +\,  (v,t) \]
is injective on the unit disc bundle of the normal bundle $\nu$ of $j$, 
and hence a closed $G$-equivariant embedding.
We define a continuous map
\[  \kappa \ : \  
D(\nu) \ \to \ C\left( j(V)\colon\bMGr(V)  \to  \bMGr(V\oplus \mR)\right) \sm B_+
\]
to the reduced mapping cone of the embedding $j(V)$ of
$\bMGr(V)\sm B_+$ into  $\bMGr(V\oplus \mR)\sm B_+$ as follows.
We consider $(b,v,t)\in D(\nu)$ where $b\in B$ and $(v,t)\in V\oplus \mR$ is
normal to $i(B)$ at $i(b)$. 
If $b\in U$, then the normal vector must lie in $V\oplus 0$, i.e., $t=0$.
The map $\kappa$ then takes $(b,v,0)$ to
\[ \left[ \left(\frac{v}{1-|v|}, \nu_b \right) ,  i_2(b)  \right] \sm b\]
in the cone of $\bMGr(V)\sm B_+$,
where $i_2(b)\in[0,1)$ is the second component of $i(b)$.
For $b\not\in U$, the map $\kappa$ sends $(b,v,t)$ to
$\left(\frac{(v,t)}{1-|v,t|}, \nu_b\right)\sm b$ 
in $\bMGr(V\oplus\mR)\sm B_+$.

The total space of the disc bundle $D(\nu)$ is a topological manifold with boundary,
and its boundary is the union of the sphere bundle $S(\nu)$ and 
the subspace $D(\nu)|_{\partial B}$, the part sitting over the boundary of $B$.
The map $\kappa$ sends the subspace $D(\nu)|_{\partial B}$
to the cone point in the mapping cone, and it sends
the sphere bundle $S(\nu)$ to the basepoint.
So $\kappa$ sends the entire boundary of $D(\nu)$ to the basepoint
of the reduced mapping cone.
So we can extend $\kappa$ continuously to $S^{V\oplus \mR}$
by sending the complement of $D(\nu)$ to the basepoint.
The result is a continuous based $G$-map
\[ \tilde c_B \ : \ S^{V\oplus \mR}\ \to \ 
C\big( j(V)\sm B_+ : \bMGr(V)\sm B_+ \to  (\sh \bMGr)(V) \sm B_+\big) \ .\]
The map $\tilde c_B$ represents the relative normal class $\td{B}^{\rel}$  
in $\pi_1^G( C (j\sm B_+))$.

It remains to establish the relation \eqref{eq:bordism_boundary_relation}.
The composite
\[  S^{V\oplus\mR}\ \xra{\tilde c_B} \ C( j(V)  \sm B_+) \
\xra{\ p \ } \ \bMGr(V)\sm B_+  \sm S^1 \]
is equivariantly homotopic to the map $((\bMGr(V)\sm\iota_+)\circ c_{\partial B})\sm S^1$, 
where
\[ c_{\partial B} \ : \ S^V\ \to \ \bMGr(V) \sm (\partial B)_+ \]
is the collapse map for $\partial B$
based on the restriction of $i$ to an embedding $\partial B\to V$.
Thus
\[ p_*(\td{B}^{\rel}) \ = \ \iota_* \td{\partial B}\sm S^1\]
in the group $\pi_1^G(\bMGr\sm B_+\sm S^1)$.
The relation \eqref{eq:bordism_boundary_relation} thus follows from
the definition of the connecting homomorphism \eqref{def:boundary}
as the composite of $p_*$ and the inverse suspension isomorphism.
\end{proof}

The inverse Thom class\index{subject}{inverse Thom class} 
\[ \tau_{H,W}\ \in \ \bMGr_0^H(S^W) \]
of an $H$-representation $W$ was defined in \eqref{eq:inverse_Thom_MGr}.
The next theorem shows how the normal class of an
induced equivariant manifold $G\times_H M$
is determined by the normal class of the $H$-manifold $M$ 
and the inverse Thom class of the tangent $H$-representation $L=T_{e H}(G/H)$.
The Wirthm{\"u}ller isomorphism
\[ \Wirth_H^G\ : \ \bMGr_0^G( (G\times_H M)_+) \ \xra{\ \iso \ }\ 
\bMGr_0^H(M_+\sm S^L) \]
was established in Theorem \ref{thm:Wirth iso}.

\begin{theorem}\label{thm:normal of induced}
Let $H$ be a closed subgroup of a compact Lie group $G$ and $M$ a closed smooth
$H$-manifold. 
Then the relation
\[ \Wirth_H^G\td{G\times_H M} \ = \ \td{M}\sm\tau_{H,L} \]
holds in the group $\bMGr_0^H( M_+\sm S^L)$, where $L=T_{e H} (G/H)$
is the tangent $H$-representation.
\end{theorem}
\begin{proof}
The Wirthm{\"u}ller map is the composite
\begin{align*}
 \bMGr_0^G( (G\times_H M)_+) \ &\xra{\res^G_H}\
\bMGr_0^H((G\times_H M)_+) \ \xra{ (\bMGr\sm l_M)_*}\ \bMGr_0^H(M_+\sm S^L) \ .
\end{align*}
Here $l_H^G:G\to S^L\sm H_+$ is the $H^2$-equivariant collapse maps
for the embedding of $H$ into $G$, and $l_M:(G\times_H M)_+\to M_+\sm S^L$
is the composite
\[ (G\times_H M)_+\ \xra{l_H^G\sm_H M}\ S^L\sm M_+ \ \iso \  M_+\sm S^L \ ,\]
compare Construction \ref{con:Wirthmuller map}.
So the claim is equivalent to the relation
\[  (l_H^G\sm_H M)_*(\res^G_H \td{G\times_H M} ) 
\ = \ \tau_{H,L}\sm\td{M} \]
in the group $\bMGr_0^H(S^L\sm  M_+)$.

We choose a $G$-equivariant wide smooth embedding $i:G/H\to V$
into a $G$-representation; 
we also choose an $H$-equivariant wide smooth embedding $j:M\to W$ 
into an $H$-representation underlying some other $G$-representation.
Then the map
\[ \psi\ : \ G\times_H M\ \to \ V\oplus W \ , \quad 
[g,m]\ \longmapsto \ (i(g H),\ g\cdot j(m)) \]
is a $G$-equivariant wide smooth embedding. We base the collapse map for
the $G$-manifold $G\times_H M$ on the embedding $\psi$.

The differential at the coset $H$ of the embedding $i$ is a linear embedding
\[ L \ = \ T_{e H} (G/H) \ \xra{\ d i\ }\ V \ ; \]
we define a scalar product on $L$ so that this embedding becomes isometric.
Now we choose a slice as in the construction of the map $l_H^G$,
i.e., a smooth embedding
\[ s\ : \ D( L ) \ \to \ G\]
of the unit disc of $L$ with $s(0)=1$, such that
$s(h \cdot l)= h\cdot s(l)\cdot h^{-1}$
for all $(h,l)\in H\times D(L)$, and such that 
the differential at $0$ of the composite
\[ D(L)\ \xrightarrow{\ s\ }\ G \ \xrightarrow{\text{proj}}\ G/H \]
is the identity of $L$.
After scaling the slice, if necessary, the map 
\[  \bar s \ : \ D(L)\times H \ \to \ G \ , \quad (l,h)\ \longmapsto \ s(l)\cdot h \]
is a smooth embedding whose image is an $H^2$-equivariant tubular neighborhood
of $H$ inside $G$. Here the group $H^2$ acts on the source  by
\[ (h_1,h_2)\cdot (l, h)\ = \ (h_1 l,\, h_1 h h_2^{-1})  \ .\]
The map $l_H^G: G \to S^L\sm H_+$ is then the collapse map 
for this tubular neighborhood.

The following $H$-equivariant composite is thus a representative for the class
$(l_H^G\sm_H M)_*(\res^G_H \td{G\times_H M})$
\[ S^{V\oplus W}\ \xra{\ c_{G\times_H M}\ }\ \bMGr(V\oplus W)\sm (G\times_H M_+)
\xra{\Id\sm (l_H^G\sm_H M)}\ 
\bMGr(V\oplus W)\sm S^L\sm  M_+\ . \]
The map $l_H^G$ takes the complement of the tubular neighborhood $\bar s(D(L)\times H)$
to the basepoint. So the map
\[ l_H^G\sm_H M_+\ : \ (G\times_H M)_+ \ \to \  S^L\sm  M_+ \]
takes the complement of the subset $\bar s(D(L)\times H)\times_H M$
to the basepoint. The composite
\[ D(L)\times M \ \xra{\ \bar s\ } \ G \ \xra{\ \psi \ } \ V\oplus W \]
is given by the formula
\[ (l,m)\ \longmapsto \ ( i(s(l)\cdot H),\, s(l)\cdot j(m))\ . \]
We define a homotopy of smooth wide $H$-equivariant embeddings
\[ [0,1]\times D(L)\times M \ \to \ V\oplus W
\text{\quad by\quad}  (t,l,m)\ \longmapsto \ ( i(s(l)\cdot H),\, s(t l)\cdot j(m)) \ .\]
For every time $t\in[0,1]$, the corresponding embedding has an associated
collapse map 
\[ S^{V\oplus W}\ \to \ \bMGr(V\oplus W)\sm S^L\sm M_+ \ . \]
So the homotopy induces a homotopy of based $H$-equivariant maps between
the above representative for $(l_H^G\sm_H M)_*(\res^G_H \td{G\times_H M})$
and the collapse map for the product embedding
\[  D(L)\times M \ \to \ V\oplus W
\text{\quad by\quad}  (l,m)\ \longmapsto \ ( i(s(l)\cdot H),\, j(m)) \ .\qedhere\]
\end{proof}

\begin{eg}
We identify the normal classes of some equivariant manifolds
in terms of classes that were previously defined.
In the case where $M=\ast$ is a single point, Theorem \ref{thm:normal of induced}
specializes to the relation
\[ \Wirth_H^G\td{G/H}\ = \ \tau_{H,L}  \]
in the group $\bMGr_0^H(S^L)$, where $L=T_{e H} (G/H)$ is the tangent $H$-representation.
In other words, the normal class of the homogeneous $G$-manifold $G/H$
is the inverse, under the Wirthm{\"u}ller isomorphism,
of the inverse Thom class.

We let $V$ be a representation of a compact Lie group $G$.
Then the unit sphere $S(\mR\oplus V)$ in $\mR\oplus V$ is a smooth closed
$G$-manifold. We recall that
\[ \Pi_V\ :\ S(\mR\oplus V)\ \to\  S^V
\ , \quad (x,v) \ \longmapsto \ \frac{v}{1-x} \]
is the $G$-equivariant stereographic projection.\index{subject}{stereographic projection}
We claim that 
\begin{equation}  \label{eq:normal class of S(V+R)}
\td{S(\mR\oplus V)} \ = \  (\bMGr\sm\Pi_V^{-1})_*(\tau_{G,V}) \ .   
\end{equation}
Since both sides of equation \eqref{eq:normal class of S(V+R)} 
commute with restriction along continuous
homomorphisms, it suffices to show the relation
for the tautological $m$-dimensional $O(m)$-representation $\nu_m$.
The composite
\[ O(m+1)/O(m)\ \xra{\ \psi\ }\ S(\mR\oplus \nu_m) \ \xra{\Pi_{\nu_m}}\ S^{\nu_m}\]
is $O(m)$-equivariantly homotopic to the map $l_{O(m)}^{O(m+1)}:O(m+1)/O(m)\to S^{\nu_m}$
that appears in the Wirthm{\"u}ller isomorphism,
where $\psi(A\cdot O(m))=A\cdot (0,\dots,0,1)$.
So we can argue:
\begin{align*}
(\bMGr\sm \Pi_{\nu_m})_* &\td{S(\mR\oplus\nu_m)} \\ 
&= \ (\bMGr\sm \Pi_{\nu_m})_*((\bMGr\sm \psi_*)\td{\res^{O(m+1)}_{O(m)} O(m+1)/O(m)}) \\ 
&= \ (\bMGr\sm l_{O(m)}^{O(m+1)})_*(\res^{O(m+1)}_{O(m)} \td{O(m+1)/O(m)}) \\ 
&= \ \Wirth_{O(m)}^{O(m+1)}\td{O(m+1)/O(m)} \ = \ \tau_{O(m),\nu_m} \  .
\end{align*}
Inverting $(\bMGr\sm \Pi_{\nu_m})_*$ gives the desired relation \eqref{eq:normal class of S(V+R)} 
in the universal example.
\end{eg}
\index{subject}{normal class!of a $G$-manifold|)}

\begin{construction}[Equivariant Thom-Pontryagin construction]
\index{subject}{Thom-Pontryagin construction}
The equivariant Thom-Pontryagin construction defines a natural transformation
of $G$-homology theories
\[ \Theta^G \ = \ \Theta^G(X)\ : \  \widetilde{\mathcal N}_*^G(X)  \ \to \ \bmO_*^G(X) \ , \]
as we now recall.
We let $(M,h)$ be an $m$-dimensional singular $G$-manifold over a based $G$-space $X$.
The way we have set things up, all the geometry is already encoded
in the normal class $\td{M}\in \bMGr_0^G(M_+)$;
the rest is a formal procedure: we push the normal class forward 
along the morphism $b:\bMGr\to\bmOP$ defined in \eqref{eq:MGr_to_bmOP},
and use the periodicity of $\bmOP$ to move into the homogeneous
summand $\bmO$ of degree~0.
While the normal class is {\em not} yet a bordism invariant, 
pushing it forward to $\bmOP$ makes it one, 
see Proposition \ref{prop:Theta bordism invariant} below.
The periodicity class $t\in\pi_{-1}^e(\bmOP^{[-1]})$ was defined in \eqref{eq:t in mOP},
and we let $\sigma\in\pi_1^e(\bmOP^{[1]})$ be its inverse.
We define
\[ \Theta^G[M,h]\ = \ (b\sm h)_* \td{M}\cdot p_G^*(\sigma^m) \ \in \ \bmO_m^G(X)\ , \]
i.e., we take the image of the normal class of $M$ under the homomorphism
\[  (b\sm h)_*\ : \ \bMGr_0^G(M_+)\ \to \ \bmOP_0^G(X)\]
and multiply by the unit $p_G^*(\sigma^m)$ in $\pi_m^G(\bmOP^{[m]})$.
Since $M$ has dimension $m$, the normal class lies in the 
homogeneous summand $\bMGr^{[-m]}$, whereas $\sigma^m$ lies in the
summand $\bmOP^{[m]}$; so the product indeed lies in the
homogeneous degree~0 summand $\bmO=\bmOP^{[0]}$.
\end{construction}

\begin{prop}\label{prop:Theta bordism invariant} 
The class $\Theta^G[M,h]$ in $\bmO_m^G(X)$
  only depends on the bordism class of the singular $G$-manifold $(M,h)$.
\end{prop}
\begin{proof}
The morphisms $j_{\bMGr}=j^\mR_\bMGr:\bMGr\to\sh\bMGr$
and $j_{\bmOP}=j_\bmOP^\mR:\bmOP\to\sh\bmOP$ 
were defined in \eqref{eq:define j^V for MGr} respectively \eqref{eq:define j^V for mOP}.
The square
\[ \xymatrix@C=15mm{ 
\bMGr \ar[r]^-{j_{\bMGr}}\ar[d]_b &\sh \bMGr \ar[d]^{\sh b} \\
\bmOP \ar[r]_-{j_{\bmOP}}^-\simeq & \sh \bmOP} \]
commutes and the lower horizontal morphism is a homotopy equivalence of
orthogonal spectra by Proposition \ref{thm:j^V for mOP}~(i).
Now we consider a smooth compact $G$-manifold $B$ with boundary $\partial B$. 
 Proposition \ref{prop:normal class properties}~(v) shows that
\begin{align*}
   (j_\bmOP\sm B_+)_*( (b\sm \iota_+)_*\td{\partial B}) \ 
&= \    (\sh b\sm B_+)_* ( (j_{\bMGr}\sm \iota_+)_*\td{\partial B}) \ =\ 0\ , 
\end{align*}
where $\iota:\partial B\to B$ is the inclusion.
Since $j_{\bmOP}$ is a homotopy equivalence, this implies 
$(b\sm \iota_+)_*\td{\partial B}=0$.
We let $(M,h:M\to X)$ be a singular $G$-manifold that is null-bordant.
We choose a null-bordism $(B,H:B\to X,\psi:M\iso\partial B)$, 
so that $H\circ\iota\circ\psi=h$. Then
\begin{align*}
(b\sm h)_*\td{M} \ &= \   (b\sm (H\circ\iota\circ\psi))_*\td{M} \ 
= \ (\bmOP\sm H)_*( (b\sm \iota_+)_*\td{\partial B} )\ = \ 0\ .
\end{align*}
Multiplying by $p_G^*(\sigma^m)$ gives that $\Theta^G[M,h]=0$.
Since the normal class is additive on disjoint unions 
(Proposition \ref{prop:normal class properties}~(i)), naturality
then implies that $\Theta^G[M,h]$ only depends on the bordism class of $(M,h)$.
\end{proof}

\begin{eg}\label{eg:geometric goes to inverse Thom}
We let $G$ be a compact Lie group and $V$ a $G$-representation of dimension $m$.
We claim that then 
\begin{equation}  \label{eq:match_classes}
 \Theta^G(d_{G,V}) \ = \ \bar\tau_{G,V}    
\end{equation}
in the group $\bmO_m^G(S^V)$.
In other words, the Thom-Pontryagin construction matches the distinguished 
geometric bordism class $d_{G,V}$ in $\widetilde{\mathcal N}_m^G(S^V)$
with the shifted inverse Thom class in the 
Thom spectrum $\bmO$.\index{subject}{shifted inverse Thom class!in $\bmO$} 

Indeed, the class $d_{G,V}$ is represented by the singular
$G$-manifold $(S(\mR\oplus V),\Pi_V)$, where $\Pi_V:S(\mR\oplus V)\to S^V$
is the stereographic projection.\index{subject}{stereographic projection}
So
\begin{align*}
  \Theta^G(d_{G,V}) \ &= \ (b\sm \Pi_V)_*\td{S(\mR\oplus V)}\cdot p_G^*(\sigma^m)\\ 
_\eqref{eq:normal class of S(V+R)}
&= \ (b\sm \Pi_V)_*(  (\bMGr\sm\Pi_V^{-1})_*(\tau_{G,V}) )\cdot p_G^*(\sigma^m)\\
&= \ (b\sm S^V)_*( \tau_{G,V}) \cdot p_G^*(\sigma^m) \ = \  \bar\tau_{G,V} \ .
\end{align*}
\end{eg}

The next theorem says, roughly speaking, that the Thom-Pontryagin construction
is compatible `with all global structure'.
There is one caveat, though, namely in how the geometric induction 
in equivariant bordism compares with the homotopy theoretic transfer.
Indeed, the geometric induction
increases the dimension by the dimension of $G/H$,
whereas the Wirthm{\"u}ller isomorphism
increases the dimension in a twisted way,
namely by the sphere of the tangent representation $L=T_{e H}(G/H)$ of $H$ in $G$.
Multiplication by the inverse Thom class is needed to compensate
this `twist' on the homotopy theory side.
In the special case where $H$ has finite index in $G$, then the
tangent representation is zero, so in this special case $\bar\tau_{H,L}=1$ and part~(v)
of the following theorem specializes to the simpler relation
    \[ \Theta^G( G\times_H  y) \ = \ G\ltimes_H \Theta^H(y)\ .\]

\begin{theorem} \label{thm:compatibility of Theta}
  \begin{enumerate}[\em (i)]
  \item 
    The Thom-Pontryagin map $\Theta^G$ is additive.  
  \item 
    The Thom-Pontryagin map is multiplicative, i.e., for
    all classes $x\in\widetilde{\mathcal N}_m^G(X)$ and
    $y\in\widetilde{\mathcal N}_n^G(Y)$, the relation
    \[ \Theta^G(x\sm y)\ = \ \Theta^G(x) \sm\Theta^G(y) \]
    holds in $\bmO_{m+n}^G(X\sm Y)$.
  \item 
    The Thom-Pontryagin map is compatible 
    with the boundary maps in the mapping cone sequences in equivariant bordism
    and equivariant $\bmO$-homology, i.e., $\Theta^G$ is a transformation of
    equivariant homology theories.
  \item
    For every continuous homomorphism $\alpha:K\to G$ of compact Lie groups,
    every based $G$-space $X$ and all $x\in\widetilde{\mathcal N}_m^G(X)$
    the relation
    \[ \Theta^K(\alpha^*(x))  \ = \ \alpha^*(\Theta^G(x))       \]
    holds in $\bmO_m^K(\alpha^*(X))$.
  \item
    If $H$ is a closed subgroup of $G$, then
    for every $H$-space $Y$ and all $y\in \mathcal N_m^G(Y)$, 
    the relation
    \[ \Theta^G(G\times_H y) \ = \ G\ltimes_H (\Theta^H(y)\sm \bar\tau_{H,L} ) \]
    holds in $\bmO_{m+d}^G( (G\times_H Y)_+)$, where $L=T_{e H}(G/H)$ is the 
    tangent $H$-representation and $d=\dim(G/H)$.
  \end{enumerate}
\end{theorem}
\begin{proof}
(i) The two functors
\[  X\ \longmapsto\ \widetilde{\mathcal N}_n^G(X) \text{\qquad and\qquad} 
X\longmapsto\  \bmO_n^G(X) \]
from the category of based $G$-spaces to the category of abelian groups
are reduced and additive.
Proposition \ref{prop:additivity prop} thus shows that the
Thom-Pontryagin map is additive.

Part~(ii) is a formal consequence of the multiplicativity of the normal classes.
We consider singular $G$-manifolds
$(M,h:M\to X)$ and $(N,g:N\to Y)$. The class $\gh{M,h}\sm\gh{N,g}$
is then represented by the singular $G$-manifold
$(M\times N,q\circ (h\times g))$,
where $q:X\times Y\to X\sm Y$ is the quotient map. 
Then
\begin{align*}
(b\sm (q\circ(h\times g)))_* \td{M\times N} \ &= \ 
(b\sm (q\circ(h\times g)))_* ( (\mu_{M,N})_*(\td{M}\times \td{N}) )\\ 
&= \ (\mu_{M,N})_* ((b\sm h)_* \td{M}\times (b\sm g)_* \td{N})  \\ 
&= \ (b\sm h)_* \td{M} \sm  (b\sm g)_* \td{N}  
\end{align*}
in the group $\bmOP_0^G(X\sm Y)$.
The first equation is Proposition \ref{prop:normal class properties}~(ii);
the second equation exploits that $b:\bMGr\to\bmOP$ 
is a homomorphism of $E_\infty$-orthogonal ring spectra.
Now we multiply with the class $p_G^*(\sigma^{m+n})$ and obtain
the desired relation 
\begin{align*}
\Theta^G\gh{ M\times N,q\circ (h\times g) } \
&= \ (b\sm (q\circ(h\times g)))_* \td{M\times N}\cdot p_G^*(\sigma^{m+n}) \\ 
&= \  ( (b\sm h)_* \td{M} \sm  (b\sm g)_* \td{N} ) \cdot p_G^*(\sigma^{m+n}) \\ 
&= \ ((b\sm h)_* \td{M}\cdot p_G^*(\sigma^m)) \sm ( (b\sm g)_* \td{N}\cdot p_G^*(\sigma^n) ) \\ 
&= \ \Theta^G\gh{M,h} \sm\Theta^G\gh{N,g} \ .
\end{align*}
    
(iii) 
We let $f:X\to Y$ be a continuous $G$-map.
Compatibility of the Thom-Pontryagin construction with the boundary homomorphism 
amounts to the commutativity of the following square:
\[ \xymatrix@C=18mm{
 \widetilde{\mathcal N}_{m+1}^G(C f)\ar[r]^-{p_*} \ar[d]_{\Theta^G} &
 \widetilde{\mathcal N}_{m+1}^G(X_+\sm S^1)\ar[r]^-{(-\sm d_{G,\mR})^{-1}}_-{\iso}
\ar[d]_{\Theta^G} &
 \mathcal N_m^G(X)\ar[d]^{\Theta^G}\\
 \bmO_{m+1}^G(C f)\ar[r]_-{p_*} &
 \bmO_{m+1}^G(X_+\sm S^1)\ar[r]_-{(-\sm\bar\tau_{G,\mR})^{-1}}^-{\iso} & \bmO_m^G(X_+)
} \]
Here $p:C f\to X_+\sm S^1$ is the projection defined in \eqref{eq:define_t}.
Indeed, the upper composite agrees with the boundary map in equivariant bordism by 
Proposition \ref{prop:suspension iso in bordism};
the lower composite is the homotopy theoretic boundary map by the definition in
Construction \ref{con:cone and fiber} and the fact that the 
suspension isomorphism in $G$-equivariant $\bmO$-homology is
exterior multiplication with the class $\bar\tau_{G,\mR}\in\bmO_1^G(S^1)$.

The Thom-Pontryagin construction is natural for continuous $G$-maps,
so it remains to show the commutativity of the right square above.
However, multiplicativity and \eqref{eq:match_classes} give
\[
\Theta^G(x\sm d_{G,\mR})  \ = \ \Theta^G(x)\sm \Theta^G(d_{G,\mR})
\ = \ \Theta^G(x)\sm  \bar\tau_{G,\mR}
\]
for all $x\in \widetilde{\mathcal N}_m^G(X)$.

Part~(iv) is straightforward from the definitions.  

(v)
This is a formal consequence of the formula for the normal class
of an induced manifold in Theorem \ref{thm:normal of induced}.
We consider a singular $H$-manifold
$(M,h:M\to Y)$. The class $G\times_H[M,h]$
is then represented by the singular $G$-manifold
$(G\times_H M,G\times_H h)$. 
Then
    \begin{align*}
\Wirth_H^G ( (b\sm (G\ltimes_H h))_*\td{G\times_H M}) \
 &= \ (b\sm h\sm S^L)_*( \Wirth_H^G \td{G\times_H M}) \\
 &= \ (b\sm h\sm S^L)_*( \td{M}\sm\tau_{H,L}) \\
 &= \  (b\sm  h)_*\td{M} \sm  (b\sm S^L)_*(\tau_{H,L}) \\
 &= \ (b\sm h)_*\td{M}\sm \tau_{H,L} \ .
  \end{align*}
Now we multiply with the class $p_G^*(\sigma^{m+n})$ and obtain
    \begin{align*}
         \Wirth_H^G ( \Theta^G(G\times_H y)) \ &= \ 
     \Wirth_H^G \left( (b\sm (G\ltimes_H h))_*\td{G\times_H M}\cdot p_G^*(\sigma^{m+d})\right) \\ 
 &= \ \Wirth_H^G ( (b\sm (G\ltimes_H h))_*\td{G\times_H M}) \cdot p_H^*(\sigma^{m+d})\\ 
&= \ ( (b\sm h)_*\td{M}\sm \tau_{H,L} ) \cdot p_H^*(\sigma^{m+d})\\ 
&= \ ( (b\sm h)_*\td{M}\cdot p_H^*(\sigma^m))\sm 
( \tau_{H,L} \cdot p_H^*(\sigma^d))\\ 
&= \ \Theta^H(y)\sm \bar\tau_{H,L} \ .
  \end{align*}
By Theorem \ref{thm:Wirth iso} the Wirthm{\"u}ller isomorphism is inverse 
to the composite of $\varepsilon_L:\bmO_{m+d}^H(Y_+\sm S^L) \to \bmO_{m+d}^H(Y_+\sm S^L)$
(the effect of $-\Id_L$ on $S^L$) and the exterior transfer. 
Proposition \ref{prop:inverse thom properties}~(i) implies that
the inverse Thom class $\tau_{H,L}$ is fixed by the involution $\varepsilon_L$.
So this last relation is equivalent to
the desired relation
$\Theta^G(G\times_H y) = G\ltimes_H( \Theta^H(y)\sm \bar\tau_{H,L})$.
\end{proof}

As before we let $E\Pc$ be a universal $G$-space 
for the family of proper closed subgroups of $G$,
and $\tilde E\Pc$ denotes its unreduced suspension.
Then $(\tilde E\Pc)^G=S^0$, consisting of the two cone points,
and $(\tilde E\Pc)^H$ is contractible for every proper subgroup of $G$.

\begin{prop}\label{prop:Theta for tilde EP}
The Thom-Pontryagin map
\[ \Theta^G \ : \ \widetilde{\mathcal N}_*^G(\tilde E\Pc)\ \to \ 
 \bmO_*^G(\tilde E\Pc)  \]
is an isomorphism for every compact Lie group $G$. 
\end{prop}
\begin{proof}
We claim that the following diagram commutes:
\[ \xymatrix@C=12mm{ 
 \widetilde{\mathcal N}_n^G(\tilde E\Pc)  \ar[d]_{\Theta^G} 
  \ar[rr]^-{\Phi_\text{geom}}_-\iso  &&
\bigoplus_{j\geq 0}\, \mathcal N_{n-j}( G r_j^{G,\perp} ) \ar[d]^{\bigoplus\Theta^e} \\
 \bmO_n^G(\tilde E\Pc)\ar[r]_-{\Phi}^-\iso &
\Phi_n^G(\bmO) \ar[r]^-\iso_-{\eqref{eq:Phi^G of mO}}  &
\bigoplus_{j\geq 0}\, \bmO_{n-j}( (G r_j^{G,\perp} )_+) } \]
Here $G r_j^{G,\perp}=( G r_j(\Uc_G^\perp))^G$,
and the geometric fixed point homomorphism $\Phi_{\text{geom}}$
was defined in \eqref{eq:define_geometric_fix}.
To see this we let $(M,h)$ be an $n$-dimensional singular
$G$-manifold over $\tilde E\Pc$.
We choose a smooth $G$-equivariant wide embedding $i:M\to V$ into
a $G$-representation and let
\[ c_M \ : \  S^V \ \to \ 
T h(G r_{|V|-n}(V))  \sm  \tilde E\Pc_+ \]
be the associated collapse map.
As before we let $M^{(j)}$ denote the $(n-j)$-dimensional component of the
fixed point manifold $M^G_0$ over the point~0, and
\[ \nu_j \ : \ M^{(j)}\ \to \ ( G r_j(V^\perp))^G \]
the classifying map for the normal bundle of $M^{(j)}$ inside $M$,
sending a fixed point $x\in M^{(j)}$ to the orthogonal complement
of $(d i)(T_x (M^{(j)}))$ inside $(d i)(T_x M)$.

The trick is now to base the Thom-Pontryagin construction for $M^{(j)}$ 
on the non-equivariant embedding
\[ M^{(j)} \ \xra{\text{incl}} \ M^G_0 \ \xra{\ i^G \ } \ V^G \]
into the $G$-fixed points of $V$.
As explained in Example \ref{eg:geometric fixed points of bmO_m},
the $G$-fixed points of the Thom space that occurs in the target of $c_M$ 
decompose as a wedge. The composite 
\[ S^{V^G}\ \xra{(c_M)^G }  ( T h(G r_{|V|-n}(V)) )^G \]
with the projection to the $j$-th summand is then on the nose the map 
\[ S^{V^G}\ \xra{c_{M^{(j)}}\sm (\nu_j)_+} \  
 T h(G r_{\dim(V^G)+j}(V^G))\sm ( G r_j(V^\perp))^G_+ \ ,\]
the smash product of the collapse map for the non-equivariant manifold $M^{(j)}$,
based on the embedding $i^G$, and $\nu_j$.
This shows that $\Phi(\Theta^G\gh{M,h})$ is the 
non-equivariant Thom-Pontryagin construction applied to $\Phi_\text{geom}\gh{M,h}$.

The upper map in the original diagram is an isomorphism
by Proposition \ref{prop:geometric isotropy}.
The lower left homotopical geometric fixed point map $\Phi$ identifies
the equivariant homotopy groups of $\bmO\sm \tilde E\Pc$
with the geometric fixed point groups $\Phi_*^G(\bmO)$,
as shown in Proposition \ref{prop:geometric as fixed points}. 
The  lower right isomorphism is the calculation of these geometric fixed point groups 
in \eqref{eq:Phi^G of mO}.
The right vertical map is a direct sum of non-equivariant Thom-Pontryagin maps,
hence an isomorphism by Thom's theorem \cite[Thm.\,IV.8]{thom-quelques}.
Since the diagram commutes, the left vertical map is also an isomorphism.
\end{proof}

Now we come to the main result of this section, showing that the
equivariant homology theory represented by the global Thom spectrum $\bmO$
is equivariant bordism, at least for products of finite groups and tori.
The result is usually credited to Wasserman,
because it can be derived from his equivariant transversality
theorem \cite[Thm.\,3.11]{wasserman}.

\begin{theorem}[Wasserman]\label{thm:TP is iso} 
Let $G$ be a compact Lie group that is isomorphic to a product of
a finite group and a torus. Then for every cofibrant $G$-space $X$,
the Thom-Pontryagin map
\[ \Theta^G(X)\ : \  \mathcal N_*^G(X)    \ \to \ \bmO_*^G(X_+) \] 
is an isomorphism.
\end{theorem}
\begin{proof}
We adapt tom Dieck's proof given in \cite[Satz 5]{tomDieck-OrbittypenII} to our setting.
Tom Dieck only considers finite groups, where homotopy theoretic transfer  
and geometric induction match up under the Thom-Pontryagin construction.
The general case has a new ingredient, namely the observation that for compact Lie groups
of positive dimensions the difference between homotopy theoretic transfer  
and geometric induction is controlled by the inverse Thom class
of the tangent representation,
compare Theorem \ref{thm:normal of induced} and 
Theorem \ref{thm:compatibility of Theta}~(v).

We prove the statement by double induction over the dimension and the
number of path components of $G$. The induction starts with the trivial
group, i.e., the non-equivariant statement, which is Thom's 
celebrated theorem \cite[Thm.\,IV.8]{thom-quelques}.
Now we let $G$ be a non-trivial compact Lie group that is a product of a finite 
group and a torus, and we assume that the theorem has been established 
for all such groups of smaller dimension, 
and for all groups of the same dimension but with fewer path components.

Every cofibrant $G$-space is equivariantly homotopy equivalent to
a $G$-CW-complex.
We can thus assume without loss of generality that $X$ is a $G$-CW-complex.
To show that $\Theta^G$ is an isomorphism,
we exploit that $\mathcal N_*^G$ and $\bmO_*^G$ 
are both equivariant homology theories and $\Theta^G$ 
is a morphism of homology theories. This reduces the claim to
the special case $X=G/H$ of an orbit for a closed subgroup $H$ of $G$.
The argument for an orbit falls into two cases, depending on whether $H$ 
is a proper subgroup or $H=G$.

When $H$ is a proper closed subgroup of $G$, the following diagram
commutes by Proposition \ref{thm:compatibility of Theta}~(v)
for the one-point $H$-space:
\[ \xymatrix@C=12mm{
\mathcal N_m^H(\ast) \ar[rr]^-{G\times_H -}\ar[d]_{\Theta^H(\ast)} && 
\mathcal N_{m+d}^G( G/H) \ar[d]^{\Theta^G(G/H)} \\
\bmO_m^H(S^0)  \ar[r]_-{-\cdot\bar\tau_{H,L}} &
\bmO_{m+d}^H(S^L) \ar[r]_-{ G\ltimes_H -} & \bmO_{m+d}^G(G/H_+)  } \]
In that diagram, the left vertical map is an isomorphism by the inductive hypothesis.
The induction map $G\times_H -:\mathcal N_m^H\to\mathcal N_{m+d}^G(G/H)$
is an isomorphism, with inverse given by sending 
a singular $G$-manifold $f:M\to G/H$ to the fiber over
the coset $H$ of an equivariant smooth approximation of $f$.
The lower right horizontal map is an isomorphism
by Theorem \ref{thm:Wirth iso}.\index{subject}{Wirthm{\"u}ller isomorphism}
Now we use the hypothesis that the group $G$ is a product
of a finite group and a torus.
In this situation the group $H$ acts trivially on the tangent representation $L$, 
so multiplication by the class $\bar\tau_{H,L}$ in $\bmO_d^H(S^L)$ 
is the suspension isomorphism, hence bijective, 
by Proposition \ref{prop:suspension iso in bordism}.
Hence the right vertical Thom-Pontryagin map $\Theta^G(G/H)$ is an isomorphism.

Now we treat the case $H=G$. 
We let $E\Pc$ be a universal $G$-space for the family of proper subgroups of $G$,
and $\tilde E\Pc$ its unreduced suspension.
We then get compatible long exact isotropy separation sequences:
\[ \xymatrix@C=5mm{
\cdots \ar[r] & \mathcal N_*^G(E\Pc) \ar[r]^-{p_*} \ar[d]_{\Theta^G} &
 \mathcal N_*^G(\ast) \ar[r]^-{i_*}\ar[d]_{\Theta^G} &
 \widetilde{\mathcal N}_*^G(\tilde E\Pc) \ar[r]^-{\partial} \ar[d]^{\Theta^G} &
 \mathcal N_{*-1}^G(E\Pc) \ar[r] \ar[d]^{\Theta^G} & \cdots \\
\cdots \ar[r] &  \bmO_*^G(E\Pc_+) \ar[r]_{p_*} &
 \bmO_*^G(S^0) \ar[r]_-{i_*} &
 \bmO_*^G(\tilde E\Pc) \ar[r]_-\partial &
 \bmO_{*-1}^G(E\Pc_+) \ar[r] & \cdots } \]
The map $\Theta^G$ is an isomorphism 
for the cofibrant and fixed point free $G$-space $E\Pc$ by the previous paragraph,
and $\Theta^G$ is an isomorphism for $\tilde E\Pc$ by Proposition \ref{prop:Theta for tilde EP}.
So the five lemma shows that the Thom-Pontryagin map for the one-point $G$-space 
is an isomorphism.
This completes the inductive step, and hence the proof of the theorem.
\end{proof}

\begin{rk}
In dimension~0, Theorem \ref{thm:pi_0 of mO} gives an explicit description
of the global functor $\upi_0(\bmO)$
as the quotient of the Burnside ring global functor $\mA$ 
by the global subfunctor generated by $\tr_e^{O(1)}\in\mA(O(1))$.
Equivariant manifolds of dimension~0 are easy to understand,
and this allows us to present the groups $\mathcal N_0^G$ in a global fashion 
very similar to (but different from) this presentation of $\upi_0(\bmO)$.
We use this description to give a direct verification that the map
\[ \Theta^G\ : \  \mathcal N_0^G    \ \to \ \pi_0^G(\bmO) \] 
is an isomorphism when $G$ is a product of a finite group and a torus.
We will also use the calculation of $\upi_0(\bmO)$
to show that the map $\Theta^G$ is {\em not} an isomorphism in general.

As we mentioned above, the groups $\mathcal N_0^G$
enjoy the structure of a restricted global functor,
i.e., a group-like global power monoid.
Moreover, the interval $[-1,1]$ with $C_2$-action by reflection at the origin
is a $C_2$-equivariant null-bordism of the free transitive $C_2$-set.
This shows that
\[ \ind_e^{C_2}(1) \ = \ 0 \text{\quad in \quad} \mathcal N_0^{C_2}\ ,\]
where $1\in \mathcal N_0^e$ is the bordism class of a point.
The action on the class~1 thus factors over a morphism of restricted global functors
\begin{equation}  \label{eq:identify_N_0}
   \mA^{\res} /\td{\ind_e^{C_2}}^{\res} \ \to\ \underline{\mathcal N}_0 \ ,
\end{equation}
from the quotient of the represented restricted global functor $\mA^{\res}=\bA^{\res}(e,-)$
by the restricted global subfunctor generated by $\ind_e^{C_2}\in\mA^{\res}(C_2)$.

We claim that the morphism \eqref{eq:identify_N_0} is an isomorphism.
Indeed, a smooth $G$-manifold of dimension~0 is just a finite $G$-set,
so  \eqref{eq:identify_N_0} is surjective. 
By the same algebraic argument as for unrestricted global functors in 
Proposition \ref{prop:explicit mO kernel},
the value of the restricted global functor $\td{\ind_e^{C_2}}^{\res}$ at $G$
is the subgroup of $\mA^{\res}(G)$ generated by $2\cdot\mA^{\res}(G)$ and
the classes $\ind_H^G\circ p_H^*$ for those finite index subgroups $H$
with Weyl group of even order.
So the source of the map \eqref{eq:identify_N_0} at a group $G$
is an $\mF_2$-vector space with basis the classes $\ind_H^G\circ p_H^*$ 
for those conjugacy classes of finite index subgroups $H$ with Weyl group of odd order.
The same classes form a basis for the bordism group $\mathcal N_0^G$;
this is shown for finite groups in \cite[Prop.\,13.1]{stong-unoriented finite},
and the general case follows because restriction along the projection $p:G\to\pi_0(G)$ 
induces an isomorphism $p^*:\mathcal N_0^{\pi_0(G)}\to\mathcal N_0^G$
for bordism of 0-manifolds.

Summing up, we have calculated both sides of the Thom-Pontryagin map $\Theta$
in dimension~0; under these isomorphisms $\Theta$ becomes the map
\[   \mA^{\res} /\td{\ind_e^{C_2}}^{\res} \ \to\  \mA /\td{\tr_e^{C_2}} \ ,\]
i.e., the same kind of quotient in the category of restricted versus unrestricted 
global functors. So the only difference between $\underline{\mathcal N}_0$
and $\upi_0(\bmO)$ is that the left side only has finite index transfers,
whereas the right hand side also has transfers for infinite index inclusions
with finite Weyl group. In finite and abelian compact Lie groups,
every subgroup inclusion with finite Weyl group is necessarily of finite index,
so for finite and abelian compact Lie groups, there is no difference in
the two kinds of quotients. 
This is an independent verification of Theorem \ref{thm:TP is iso} 
in dimension~0. Moreover, we conclude that the Thom-Pontryagin map
$\Theta^G:\mathcal N_0^G\to\pi_0^G(\bmO)$ in dimension~0
is always injective.

On the other hand, the map $\Theta^G$ is {\em not} generally surjective
in dimension~0.
A specific example is the group $G=S U(2)$:\index{subject}{special unitary group!$S U(2)$}
the normalizer $N=N_{S U(2)} T$ of a maximal torus $T$ of $S U(2)$
is self-normalizing, so $N$ has trivial Weyl group in $S U(2)$.
So the classes $1$ and $\tr_N^{S U(2)}(1)$ are linearly independent in $\pi_0^{S U(2)}(\bmO)$.
On the other hand $\mathcal N_0^{S U(2)}=\mZ/2$ because $S U(2)$
is connected.
\end{rk}

\begin{rk}[Stable equivariant bordism and localized $\bmO$]
\label{rk:local mO and stable bordism}
We showed in Corollary \ref{cor:MO is localized mO} that the localized 
equivariant $\bmO$-theory $\bmO^G_*[1/\tau]$
is isomorphic to the theory $\bMO^G_*$.
Our next task is to show that $\bmO[1/\tau]$,
and hence also $\bMO$, 
has a geometric interpretation as {\em stable equivariant bordism}.
\index{subject}{equivariant bordism!stable}
In contrast to Theorem \ref{thm:TP is iso}, 
this stable interpretation works for all compact Lie groups,
not only for products of finite groups and tori.

In \cite{broecker-hook},
Br{\"o}cker and Hook define the {\em stable equivariant bordism groups}
$\tilde{\mathfrak  N}^{G:S}_*(X)$ of a based $G$-space $X$
as the localization of
the geometric bordism group $\tilde{\mathcal N}_*^G(X)$
by formally inverting all the classes $d_{G,V}$.
More precisely, their definition comes down to
\[ \tilde{\mathfrak  N}^{G:S}_m(X)\ = \ 
\colim_{V\in s(\Uc_G)}\,  \widetilde{\mathcal N}_{m+|V|}^G(X\sm S^V) \ ;   \]
for $V\subset W$, the structure map in the colimit system is the multiplication
\[ \widetilde{\mathcal N}_{m+|V|}^G(X\sm S^V) \ \xra{-\sm d_{G,W-V}}\ 
\widetilde{\mathcal N}_{m+|W|}^G(X\sm S^V\sm S^{W-V}) \ \iso\
\widetilde{\mathcal N}_{m+|W|}^G(X\sm S^W)  \ .\]
As we explained in Example \ref{eg:geometric goes to inverse Thom},
the Thom-Pontryagin construction takes the distinguished 
geometric bordism class $d_{G,V}$ to the shifted inverse Thom class, i.e., 
$ \Theta^G(d_{G,V}) = \bar\tau_{G,V}$ in $\bmO_{|V|}^G(S^V)$.
Since the Thom-Pontryagin maps take the geometric product in equivariant
bordism to the homotopy theoretic product in $\bmO$, we conclude that
for every compact Lie group $G$, 
every $G$-representation $V$ and every based $G$-space $X$, the following
square commutes:
\[ \xymatrix@C=12mm{
 \widetilde{\mathcal N}_m^G(X) \ar[d]_{-\sm d_{G,V}} \ar[r]^{\Theta^G} & 
\bmO_m^G(X) \ar[d]^{-\cdot\bar\tau_{G,V}}\\
\widetilde{\mathcal N}_{m+|V|}^G(X\sm S^V) 
\ar[r]_-{\Theta^G} & \bmO_{m+|V|}^G(X\sm S^V) } \]  
The Thom-Pontryagin maps thus assemble into a natural transformation
\[ \Theta^G \ : \  \tilde{\mathfrak  N}^{G:S}_m(X)\  \to \ 
 \bmO^G_m(X)[1/\tau]\]
between the localized theories, for which we use the same letter.
\end{rk}

If $G$ is a product of a finite group and a torus, then the next theorem
is a direct consequence of Theorem \ref{thm:TP is iso}.
The point, however, is that the following localized version holds
without any restriction on the compact Lie group $G$.
Morally, the reason for this is that formally inverting the classes $d_{G,V}$ 
forces the Wirthm{\"u}ller isomorphism to hold,
so in stable equivariant bordism this potential obstruction
to representability by a global homotopy type vanishes. 

\begin{theorem}\label{thm:stable TP} 
For every compact Lie group $G$ and every cofibrant based $G$-space $X$,
the map
\[ \Theta^G \ : \  \tilde{\mathfrak  N}^{G:S}_*(X)\  \to \ 
 \bmO^G_*(X)[1/\tau]\]
is an isomorphism of graded abelian groups.  
\end{theorem}
\begin{proof}
Since filtered colimits of abelian groups are exact, the localized theories 
$\tilde{\mathfrak  N}^{G:S}_*(-)$ and $\bmO^G_*(-)[1/\tau]$
are both equivariant homology theories for every fixed group $G$.
We can run the same inductive argument
as in the proof of Proposition \ref{thm:TP is iso}, i.e., show the
claim by double induction over the dimension and the
number of path components of $G$. 
But this time the case of an orbit $X=G/H$ for a proper closed subgroup $H$ of $G$
works without any restriction on $G$ and $H$. Indeed, in the commutative diagram
\[ \xymatrix@C=12mm{
\tilde{\mathfrak N}_m^{H:S}(S^0) \ar[rr]^-{G\times_H -}\ar[d]_{\Theta^H} && 
\tilde{\mathfrak N}_{m+d}^{G:S}( G/H_+) \ar[d]^{\Theta^G} \\
\bmO_m^H(S^0)[1/\tau]  \ar[r]_-{-\cdot\bar\tau_{H,L}} &
\bmO_{m+d}^H(S^L)[1/\tau] \ar[r]_-{ G\ltimes_H -} & \bmO_{m+d}^G(G/H_+)[1/\tau]  } \]
multiplication by the class $\bar\tau_{H,L}$ is now invertible
in the localized theory $\bmO_*^G[1/\tau]$. So all horizontal maps in the diagram
are isomorphisms. Since the left vertical map is an isomorphism by induction, 
the right vertical map is an isomorphism as well.

The case $X=G/G$ is essentially taken care of
by Proposition \ref{prop:Theta for tilde EP}.
We let $V$ be any $G$-representation, 
and observe that the fixed point inclusion $i:V^G\to V$
induces a $G$-homotopy equivalence
\[\tilde E\Pc\sm i \ : \  \tilde E\Pc\sm S^{V^G} \ \to \ \tilde E\Pc\sm S^V \ .\]
In the commutative diagram
\[ \xymatrix@C=13mm@R=7mm{ 
\widetilde{\mathcal N}_{m+|V^\perp|}^G(\tilde E\Pc)\ar[r]^-{\Theta^G} 
\ar[d]_{-\sm d_{G,V^G}} &
 \bmO_{m+|V^\perp|}^G(\tilde E\Pc) \ar[d]^{-\cdot\tau_{G,V^G}}\\
\widetilde{\mathcal N}_{m+|V|}^G( \tilde E\Pc\sm S^{V^G})\ar[r]^-{\Theta^G} 
\ar[d]_{(\tilde E\Pc\sm i)_*} &
 \bmO_{m+|V|}^G(\tilde E\Pc\sm S^{V^G}) \ar[d]^{(\tilde E\Pc\sm i)_*}\\
\widetilde{\mathcal N}_{m+|V|}^G(\tilde E\Pc\sm S^V)\ar[r]_-{\Theta^G} &
 \bmO_{m+|V|}^G(\tilde E\Pc\sm S^V) }
 \]
the lower vertical maps are thus isomorphisms.
The upper vertical maps are suspension isomorphisms,
hence both vertical composites are isomorphisms.
The upper horizontal map is an isomorphism by Proposition \ref{prop:Theta for tilde EP},
hence so is the lower horizontal map.
As a colimit of isomorphisms, the localized Thom-Pontryagin map
\[ \Theta^G\ :\ \mathfrak N_*^{G:S}(\tilde E\Pc)\ \to\
 \bmO_*^G(\tilde E\Pc)[1/\tau]\]
is also an isomorphism.
Now we finish the argument as in Proposition \ref{thm:TP is iso}:
we compare the isotropy separation sequence for
$\tilde{\mathfrak  N}^{G:S}_*(-)$ to that for $\bmO^G_*(-)[1/\tau]$,
and the five lemma concludes the inductive step.
\end{proof}

\begin{rk}[Stable equivariant bordism and $\bMO$]\label{rk:proof of Brocker Hook}
Corollary \ref{cor:MO is localized mO} and Theorem \ref{thm:stable TP} 
provide natural isomorphisms
\[  \tilde{\mathfrak  N}^{G:S}_*(X)\  \xra{\ \Theta^G\ }
 \bmO^G_*(X)[1/\tau]\ \xra{\ \iso \ }\  \bMO^G_*(X) \]
for cofibrant based $G$-spaces $X$.
Hence the composite is an isomorphism, which provides an alternative proof 
that stable equivariant bordism
agrees with equivariant $\bMO$-homology, which is the main result 
of the paper \cite{broecker-hook} by Br{\"o}cker and Hook. 
Strictly speaking there is a bit more work involved
in the translation, because our group $\bMO_*^G(X)$ is not literally the same
as the homotopy theoretic equivariant bordism group $\tilde N_n^G(X)$
in \cite{broecker-hook}; we invite the reader to spell out
an isomorphism, which boils down to a certain rewriting of colimits,
such that our composite isomorphism 
becomes the map $\Phi^S$ considered by Br{\"o}cker and Hook.
\end{rk}

\index{subject}{global Thom spectrum|)}
\index{subject}{equivariant bordism|)}

\section{Connective global \texorpdfstring{$K$}{K}-theory}\label{sec:connective global K}

In this section we define and discuss the ultra-commutative ring spectrum $\bku$,
the {\em connective global $K$-theory spectrum}, 
see Construction \ref{con:connective global K-theory}.
Our construction is an elaboration of a model of non-equivariant connective $K$-theory 
by Segal \cite{segal-K-homology},
constructed from certain $\bGamma$-spaces of 
`orthogonal subspaces in the symmetric algebra'.
The degree zero equivariant cohomology theory represented by $\bku$
tries hard to be equivariant $K$-theory.
Indeed, Theorem \ref{thm:K_G_is_ku} exhibits a natural transformation
from the $G$-equivariant $K$-group that is an isomorphism
for finite $G$-CW-complexes with finite isotropy groups.
A special case is a ring homomorphism $\bRU(G)\to\pi_0^G(\bku)$
from the complex representation ring that is an isomorphism whenever $G$ is finite.
For varying $G$, these homomorphisms almost form a morphism of global power functors,
in the sense that they are compatible
with finite index transfers, but {\em not} with general degree zero transfers
across an infinite index inclusion, see Remark \ref{rk:homotopy vs smooth transfer}.   

The {\em Bott class} is a preferred element $\beta\in\pi_2^e(\bku)$ 
that we define in Construction \ref{con:Bott class}.
This class becomes invertible in the periodic global $K$-theory spectrum $\bKU$,
by Theorem \ref{thm:Bott and inverse} below.
More generally, every $G$-$\Spin^c$-re\-presen\-tation has an associated
{\em equivariant Bott class}
\[ \beta_{G,W} \ \in \ \bku^0_G(S^W) \ ,\]
see Construction \ref{con:equivariant Bott class},
that becomes an $R O(G)$-graded unit in the periodic theory $\bKU$.

\begin{construction} 
We let $\Uc$ be a complex vector space of countable dimension (finite or infinite)
equipped with a hermitian inner product.
We recall a certain $\bGamma$-space $\bk(\Uc)$ of `orthogonal subspaces in $\Uc$',
due to Segal \cite[Sec.\,1]{segal-K-homology}.\index{subject}{Gamma-space@$\bGamma$-space}
For a finite based set $A$ we let $\bk(\Uc,A)$ be the space of tuples
$(E_a)$, indexed by the non-basepoint elements of $A$,
of finite-dimensional, pairwise orthogonal $\mC$-subspaces of $\Uc$. 
The topology on $\bk(\Uc,A)$ 
is that of a disjoint union of subspaces of a product
of Grassmannians.
The basepoint of  $\bk(\Uc,A)$ is the tuple where each $E_a$
is the zero subspace.
For a based map $\alpha:A\to  B$ the induced map
$\bk(\Uc,\alpha):\bk(\Uc,A)\to \bk(\Uc,B)$ sends $(E_a)$ to $(F_b)$ where
\[ F_b \ = \ {\bigoplus}_{\alpha(a)=b}\, E_a \ .  \]
Then $\bk(\Uc)$ is a $\bGamma$-space whose underlying space is 
\[ \bk(\Uc,1_+) \ = \ {\coprod}_{n\geq 0}\, G r_n^\mC(\Uc)\ ,\]
the disjoint union of the different Grassmannians of $\Uc$.

Every $\bGamma$-space can be evaluated on a based space 
by the coend construction \eqref{eq:extended_Gamma_space}.
We write $\bk(\Uc,K)=\bk(\Uc)(K)$ for the value of the $\bGamma$-space
$\bk(\Uc)$ on a based space $K$.
When $A$ is a finite based set with discrete topology, 
then the prolongation $\bk(\Uc)(A)$ is canonically
homeomorphic to the original space $\bk(\Uc,A)$, 
compare Remark \ref{rk:abuse prolongation}.
Elements of $\bk(\Uc,K)$ can be interpreted as `labeled configurations': 
a point is represented by an unordered tuple 
\[ [ E_1,\dots,E_n;\,k_1,\dots,k_n ] \]
where $(E_1,\dots, E_n)$ is an $n$-tuple of finite-dimensional,
pairwise orthogonal subspaces of $\Uc$, 
and $k_1,\dots,k_n$ are points of $K$, for some $n$. 
The topology is such that,
informally speaking, the labels sum up whenever two points collide,
and a label disappears whenever a point approaches the basepoint of $K$.
\end{construction}

\begin{rk}\label{rk:conf versus C^ast}
When $K$ is compact, the space $\bk(\Uc,K)$ can be described differently,
compare again \cite[Sec.\,1]{segal-K-homology}, namely as the space
\[  \colim_{V\in s(\Uc)}\,  C^\ast(C_0(K),\End_\mC(V) ) \ ,\]
where the colimit runs over all finite-dimensional $\mC$-subspaces of $\Uc$.
Here  $C_0(K)$ is the $C^\ast$-algebra of continuous $\mC$-valued functions
on $K$ that vanish at the basepoint,
$\End_\mC(V)$ is the endomorphism $C^\ast$-algebra of $V$
(which is isomorphic to a matrix algebra over $\mC$),
and $C^\ast(-,-)$ is the space of $C^\ast$-algebra homomorphisms,
with the zero homomorphism as the basepoint.
For $U\subset V$ the map in the colimit system is induced by
the $\ast$-homomorphism $\End_\mC(U)\to \End_\mC(V)$
that extends an endomorphism by~0 on the orthogonal complement.
A homeomorphism
\[ \bk(\Uc,K)\ \to \ \colim_{V\in s(\Uc)} \,  C^\ast(C_0(K),\End_\mC(V)) \]
is given by sending a configuration $[E_1,\dots,E_n;\,k_1,\dots,k_n]$
in $\bk(\Uc,K)$ to the homomorphism that takes a function $\varphi\in C_0(K)$ to 
\[ {\sum}_{i=1}^n\ \varphi(k_i)\cdot p_{E_i}\ ;\]
here $V$ is chosen large enough to contain all the spaces $E_i$, and
$p_{E_i}:V\to V$ is the orthogonal projection onto the subspace $E_i$.
\end{rk}

\begin{rk}[Eigenspace decomposition]\label{rk:eigenspace decomposition}
If $\Uc$ is infinite dimensional, then the $\bGamma$-space $\bk(\Uc)$ is special,
compare Theorem \ref{thm:bk special} (i) below.
The orthogonal spectrum $\bk(\Uc)(\mS)$
is then a positive $\Omega$-spectrum by the general theory. In particular,
the space $\bk(\Uc)(\mS)(\mR)=\bk(\Uc,S^1)$ is an infinite loop space. 
In fact, $\bk(\Uc,S^1)$ is a familiar space, 
namely the infinite unitary group $U(\Uc)$, i.e., the group of 
linear self-isometries of $\Uc$ that are the identity on the orthogonal complement
of some finite-dimensional subspace.
This {\em eigenspace decomposition}\index{subject}{eigenspace decomposition}
works for every hermitian inner product space $\Uc$, of finite or countable
dimension, as follows.
As before we identify $S^1$ with the unit circle $U(1)$ in the complex numbers via
the Cayley transform\index{subject}{Cayley transform}
\[ c \ : \ S^1\ \iso\   U(1)\ , \quad x \ \longmapsto \ (x+i)(x-i)^{-1}\ .\]
This homeomorphism sends the basepoint at infinity to~1. 
Given a tuple of pairwise orthogonal subspaces $(E_1,\dots,E_n)$ 
of $\Uc$ and a point $(x_1,\dots,x_n)\in (S^1)^n$,
we let $\psi(E_1,\dots,E_n;\,x_1,\dots,x_n)$
be the isometry of $\Uc$ that is multiplication by $c(x_i)$ on $E_i$ 
and the identity on the orthogonal complement of $\bigoplus_{i=1}^n E_i$. 
In other words: $E_i$ is the eigenspace 
of $\psi(E_1,\dots,E_n;\,x_1,\dots,x_n)$ for the eigenvalue $c(x_i)$.
As $n$ varies, these maps are compatible
with the equivalence relation and so they assemble into a continuous map
\begin{align}\label{eq:C(S^1)_is_U}\index{subject}{eigenspace decomposition}
 \bk(\Uc,S^1) \ = \ 
\int^{n_+\in\bGamma} \bk(\Uc,n_+) \times (S^1)^n\ &\to \quad U(\Uc)\\
[E_1,\dots,E_n;\,x_1,\dots,x_n]\ &\longmapsto \
\psi(E_1,\dots,E_n;\,x_1,\dots,x_n)\ .  \nonumber
\end{align}
This map is a homeomorphism because every unitary transformation
is diagonalizable with eigenvalues in $U(1)$ 
and pairwise orthogonal eigenspaces. 
\end{rk}

Now we let $\Uc$ and $\Vc$ be two hermitian vector spaces
(of countable dimension, but possibly infinite dimensional).
Here and below we write $\tensor=\tensor_\mC$ for the tensor product over $\mC$.
We endow  $\Uc\tensor\Vc$ with a hermitian scalar product by declaring
\begin{equation}\label{eq:scalar_product_on_tensors}
 (u\tensor v,\, u'\tensor v')\ = \ (u,u')\cdot (v,v')   
\end{equation}
on elementary tensors and extending biadditively.
If $E,F$ are orthogonal subspaces of $\Uc$
and  $E'$ is a subspace of $\Vc$, then
$E\tensor E'$ and $F\tensor E'$ are orthogonal subspaces of $\Uc\tensor \Vc$,
and similarly in the second variable.
For based spaces $K$ and $L$ we can thus define 
a continuous multiplication map
\begin{align}  \label{eq:tensor_product}
 \bk(\Uc,K)\sm \bk(\Vc,L)\ &\to \ \bk(\Uc\tensor\Vc,K\sm L)\nonumber\\
[E_i;k_i]\sm [F_j;l_j]\quad &\longmapsto \quad [E_i\tensor F_j;k_i\sm l_j]\ .
  \end{align}
These multiplication maps are associative, and commutative
in the sense that the following square commutes:
\[ \xymatrix@C=12mm{
 \bk(\Uc,K)\sm \bk(\Vc,L)\ar[r]\ar[d]_{\tau_{\bk(\Uc,K),\bk(\Vc,L)}} & 
\bk(\Uc\tensor\Vc,K\sm L)\ar[d]^{\bk(\tau_{\Uc,\Vc},\tau_{K,L})}\\
 \bk(\Vc,L)\sm \bk(\Uc,K)\ar[r] & \bk(\Vc\tensor\Uc,L\sm K)
} \]

Our construction of connective global $K$-theory
needs induced inner products on symmetric powers.
We explain the complex version; the real version 
works in much the same way.
For a $\mC$-vector space $V$ we denote by
\[ \Sym^n(V) \ = \ V^{\tensor n}/\Sigma_n \]
the $n$-th symmetric power of $V$ and by $\Sym(V)=\bigoplus_{n\geq 0}\Sym^n(V)$
the symmetric algebra of $V$.
If $W$ is another $\mC$-vector space, then
the two direct summand embeddings of $V$ and $W$ into $V\oplus W$
induce homomorphisms of symmetric algebras 
that combine (by multiplying in the target) into a natural
$\mC$-algebra isomorphism
\begin{equation}\label{eq:tensor and Sym} 
  \Sym(V)\tensor\Sym(W)\ \iso \ \Sym(V\oplus W) \ .  
\end{equation}
If $V$ is equipped with a hermitian inner product, 
then the symmetric powers inherit a preferred inner product:

\begin{prop}\label{prop:Sym induced inner product}
  For every hermitian inner product space $V$ there is a unique 
  inner product on $\Sym^n(V)$ that satisfies 
  \[ 
  ( v_1\cdot\ldots\cdot v_n,\,\bar v_1\cdot\ldots\cdot \bar v_n) \ = \ 
  \sum_{\sigma\in\Sigma_n} 
  (v_1,\bar v_{\sigma(1)})\cdot\ldots\cdot (v_n,\bar v_{\sigma(n)})
  \]
  for all $v_i,\bar v_i\in V$.
  This inner product on $\Sym^n(V)$ is natural for $\mC$-linear isometric embeddings
  and it makes the algebra isomorphism \eqref{eq:tensor and Sym} into an isometry.
\end{prop}
\begin{proof}
  Uniqueness of the scalar product follows from the fact that the
  symmetric products $v_1\cdot\ldots\cdot v_n$ generate $\Sym^n(V)$ 
  as a $\mC$-vector space.
   The tensor product $V\tensor W$ of two hermitian inner product spaces $V$ and $W$ 
   has a preferred hermitian inner product as in \eqref{eq:scalar_product_on_tensors}.
   By iteration, the $n$-fold tensor product $V^{\tensor n}$ inherits an inner product. 
   There is thus a unique inner product
   on $\Sym^n(V)$ that makes
   the normalized $\mC$-linear `symmetrization' embedding
  \begin{align*}
    \Sym^n(V) \ &\to \qquad    V^{\tensor n} \\
    v_1\cdot\ldots\cdot v_n \ &\longmapsto \ \frac{1}{\sqrt{n!}}\cdot
    \sum_{\sigma\in\Sigma_n} v_{\sigma(1)}\tensor\ldots\tensor v_{\sigma(n)}\nonumber
   \end{align*}
   a linear isometric embedding. 
   We omit the straightforward verification that this inner product is indeed
   given by the formula in the statement of the proposition,
   and that it is natural for linear isometric embeddings.

  The algebra isomorphism \eqref{eq:tensor and Sym}  is the sum of the embeddings
  \begin{align*}
    \Sym^m(V)\tensor\Sym^n(W)\quad &\to \quad \Sym^{m+n}(V\oplus W) \\
    v_1\cdot\ldots\cdot v_m\tensor  w_1\cdot\ldots\cdot w_n\ &\longmapsto\
    (v_1,0)\cdot\ldots\cdot (v_m,0)\cdot(0,w_1)\cdot\ldots\cdot(0, w_n)\  .  
  \end{align*}
  The relation
  \begin{align*}
    ( v_1\cdot\ldots\cdot v_m&\tensor  w_1\cdot\ldots\cdot w_n,\,\bar v_1\cdot\ldots\cdot \bar v_m\tensor \bar w_1\cdot\ldots\cdot\bar w_n) \\ 
    &= \
    (v_1\cdot\ldots\cdot v_m,\,\bar v_1\cdot\ldots\cdot \bar v_m)
    \cdot (w_1\cdot\ldots\cdot w_n,\,\bar w_1\cdot\ldots\cdot\bar w_n)\\ 
    &= \  \sum_{(\sigma,\tau)\in\Sigma_m\times\Sigma_n} 
  ( v_1,\bar v_{\sigma(1)})\cdot\ldots\cdot (v_m,\bar v_{\sigma(m)})
  \cdot  (w_1,\bar w_{\tau(1)})\cdot\ldots\cdot (w_n,\bar w_{\tau(m)})\\ 
  &= \ \big( (v_1,0)\cdot\ldots\cdot (v_m,0)\cdot(0,w_1)\cdot\ldots\cdot(0, w_n),\\
& \hspace*{3cm} ( \bar v_1,0)\cdot\ldots\cdot (\bar v_m,0)\cdot(0,\bar w_1)\cdot\ldots\cdot(0,\bar w_n) \big)
\end{align*}
then proves that \eqref{eq:tensor and Sym} preserves the inner product.
The last equation uses that $(v,0)$ and $(0,w)$ are orthogonal
in $V\oplus W$ for all $v\in V$ and $w\in W$.
\end{proof}

  The case $n=2$ gives an idea of the induced inner product on $\Sym^n(V)$:
  if $\{e_1,\dots,e_k\}$ is an orthonormal basis of $V$, then the vectors
   \[ 1/\sqrt{2}\cdot e_i^2 \quad (1\leq i\leq k)\text{\qquad and\qquad} 
   e_i\cdot e_j \quad (1\leq i< j\leq k)\]
   form an orthonormal basis of $\Sym^2(V)$.

\index{subject}{K-theory@$K$-theory!connective global|(}

 \begin{construction}[Connective global $K$-theory]\label{con:connective global K-theory}
We can now define an ultra-commutative ring spectrum $\bku$, 
the {\em connective global $K$-theory spectrum}.\index{subject}{connective global $K$-theory}\index{subject}{K-theory@$K$-theory!connective global}\index{symbol}{$\bku$ - {connective global $K$-theory}}  
The value of $\bku$ on a euclidean inner product space $V$ is
\[  \bku(V)\ = \ \bk(\Sym(V_\mC),S^V)\ ,   \]
the value of the $\bGamma$-space $\bk(\Sym(V_\mC))$
on the one-point compactification of $V$.
Here $V_\mC$ is the complexification of $V$ with the
induced hermitian inner product, and the inner product on the symmetric algebra described 
in Proposition \ref{prop:Sym induced inner product}.
The action of $O(V)$ on $V$ then extends to a unitary action on $\Sym(V_\mC)$. 
We let the orthogonal group $O(V)$ act diagonally, 
via the action on the sphere $S^V$ and the action 
on the $\bGamma$-space $\bk(\Sym(V_\mC))$. 
Using the tensor product pairing \eqref{eq:tensor_product}
we define an $(O(V)\times O(W))$-equivariant multiplication map
\begin{align*}
 \mu_{V,W}\ : \  \bku(V)\sm  \bku(W)\ = \quad 
&\bk(\Sym(V_\mC),S^V)\sm \bk(\Sym(W_\mC),S^W)\\
_\eqref{eq:tensor_product}  \to\  &\bk(\Sym(V_\mC)\tensor\Sym(W_\mC),S^V\sm S^W)\\ 
_\eqref{eq:tensor and Sym} \,
\iso \quad  &\bk(\Sym((V\oplus W)_\mC),S^{V\oplus W}) = \bku(V\oplus W)\ . 
\end{align*}
The maps $\mu_{V,W}$  are associative and commutative.
An $O(V)$-equivariant unit map is given by
\[ \iota_V \ : \ S^V\ \to \ \bk(\Sym(V_\mC),S^V) = \bku(V) \ , \quad 
v \ \longmapsto \ [\mC\!\cdot\! 1;\, v]\ ,\]
where $\mC\!\cdot\! 1$ is the homogeneous summand of degree~0
in the symmetric algebra, i.e., the line spanned by the multiplicative unit.
This structure makes $\bku$ into an ultra-commutative ring spectrum.

The space $\bku_0=\bk(\mC\!\cdot\! 1,S^0)$ consists of all subspaces of
$\Sym(0)=\mC\!\cdot\! 1$, so it has two points, the basepoint~0
and the point $\mC\!\cdot\! 1$. The unit map $\iota_0:S^0\to \bku_0$ is 
thus a homeomorphism.
\end{construction}

\begin{construction}[Complex conjugation on $\bku$]\label{con:conjugation on bku}\index{subject}{complex conjugation!on $\bku$}
The ultra-commutative ring spectrum $\bku$ comes with an involution 
by `complex conjugation'
that preserves all the structure available.
Indeed, for every euclidean inner product space $V$ the 
complex symmetric algebra $\Sym(V_\mC)$ of the complexification
is canonically isomorphic to $\mC\tensor_\mR \Sym_\mR(V)$, 
the complexification of the real symmetric algebra of $V$.
So $\Sym(V_\mC)$ comes with an involution $\psi_{\Sym(V)}$
that is $\mC$-semilinear and preserves the orthogonality relation. 
Applying this involution elementwise to tuples of orthogonal subspaces
gives an involution $\bk(\psi_{\Sym(V)}):\bk(\Sym(V_\mC))\to\bk(\Sym(V_\mC))$
of the $\bGamma$-space and hence a homeomorphism
\begin{equation}\label{eq:psi_on_bku}
  \psi(V)\ = \ \bk(\psi_{\Sym(V)},S^V)\ : \ \bku(V) \ \to \ \bku(V)  
\end{equation}
of order~2.
As $V$ varies, the maps $\psi(V)$ form an automorphism 
$\psi:\bku\to\bku$ of the ultra-commutative ring spectrum $\bku$.
\end{construction}

\begin{rk}[Connective real global $K$-theory]
There is a straightforward real analog $\bko$ of the 
complex connective global $K$-theory spectrum $\bku$. 
\index{subject}{connective global $K$-theory!real}\index{subject}{K-theory@$K$-theory!real connective global}\index{symbol}{$\bko$ - {real connective global $K$-theory}}  
The value of $\bko$ on an inner product space $V$ is
\[  \bko(V)\ = \ \bk_\mR(\Sym(V),S^V)\ ,    \]
where now $\Sym(V)$ is the symmetric algebra of $V$ over the real numbers, 
and $\bk_\mR(\Sym(V))$ is the $\bGamma$-space of tuples of
pairwise orthogonal, finite-dimensional $\mR$-subspaces of $\Sym(V)$.
As $V$ varies, the spaces $\bko(V)$ again come with the structure 
of an ultra-commutative ring spectrum. 

Complexification defines a morphism of ultra-commutative ring spectra
\index{subject}{complexification morphism!from $\bko$ to $\bku$}
\begin{equation}\label{eq:c_bko2bku}
 c \ : \ \bko \ \to \ \bku\ .   
\end{equation}
In more detail: if $V$ is a real inner product space, then the map
\[ c(V) \ : \ \bko(V) = \bk_\mR(\Sym(V),S^V) \ \to \ 
 \bk(\Sym(V_\mC),S^V) = \bku(V)  \]
sends a configuration $[E_1,\dots,E_n;\,v_1,\dots,v_n]$ 
to the complexified configuration  
\[  [(E_1)_\mC,\dots,(E_n)_\mC;\,v_1,\dots,v_n]  \ .\]
As $V$ varies, these maps form the morphism \eqref{eq:c_bko2bku} 
of ultra-commutative ring spectra.
The complexification of a real subspace of $\Sym(V)$ is invariant under 
the complex conjugation involution $\psi_{\Sym(V)}$ of $\Sym(V_\mC)$,
so the image of the complexification morphism 
is invariant under the complex conjugation involution $\psi$ of $\bku$
defined in \eqref{eq:psi_on_bku}.
Even more is true: a complex subspace of $\Sym(V_\mC)$
is $\psi$-invariant if and only if it is the complexification of
its real part (the +1 eigenspace of $\psi_{\Sym(V)}$). 
This means that  $\bko$ `is'
the $\psi$-fixed orthogonal ring subspectrum of $\bku$; more formally,
the complexification morphism \eqref{eq:c_bko2bku} is an isomorphism
from $\bko$ onto the $\psi$-fixed orthogonal ring subspectrum of $\bku$.
\end{rk}

The next proposition justifies the adjective `connective'
that we attached to the names of the orthogonal spectra $\bku$ and $\bko$.
The proof uses a certain cofibrancy property of the $\bGamma$-space $\bk(\Sym(V_\mC))$
that we now introduce.

\begin{construction}[Latching map]
We let $\Pc(n)$ denote the power set of the set $\{1,\dots,n\}$, 
i.e., the set of subsets.
We also write $\Pc(n)$ for the associated poset category, i.e., with object set $\Pc(n)$
and exactly one morphism $U\to T$ whenever $U\subseteq T$.
Given a $\bGamma$-space $F$, we obtain a functor from $\Pc(n)$ 
to based spaces by sending a subset $U$
to $F(U_+)$, with the maps $F(U_+)\to F(T_+)$ induced by the inclusions.
We obtain a {\em latching map}\index{subject}{latching map!of a $\bGamma$-space}
\[ l_n\ : \ \colim_{U\subsetneq \{1,\dots,n\}} \, F(U_+)\ \to \  F(n_+) \ .  \]
The map $l_n$ is equivariant for the action of $\Sigma_n$
given on the target by functoriality of $F$. The action of a permutation
$\sigma\in\Sigma_n$ on the source is induced by
sending $F(U_+)$ to $F(\sigma(U)_+)$ via the map
\[ F( (\sigma\cdot -)_+) \ : \ F(U_+)\ \to \ F(\sigma(U)_+)\ .\]
The latching map $l_n$ is always a closed embedding, 
by Proposition \ref{prop:Gamma latching closed embedding}~(ii).
\end{construction}

We will now consider a $\bGamma$-$G$-space for a compact Lie group $G$,
i.e., a reduced functor $F:\bGamma\to G\bT_*$ to the category of based $G$-spaces.
Then $G$ acts on source and target of the latching map $l_n$, which
is thus $(\Sigma_n\times G)$-equivariant.
The following cofibrancy condition on a $\bGamma$-$G$-space ensures that
the prolongation is homotopically meaningful.

\begin{defn}
Let $G$ be a compact Lie group. A $\bGamma$-$G$-space $F$ is {\em $G$-cofibrant}
\index{subject}{Gamma-space@$\bGamma$-space!cofibrant}
if for every $n\geq 1$ the latching map $l_n$ 
is a $(\Sigma_n\times G)$-cofibration.
\end{defn}

\begin{eg}\label{eg:ku cofibrant}  
We let $G$ be a compact Lie group and $\Uc$ a unitary $G$-representation,
of finite or countably infinite dimension.
We argue that the $\bGamma$-$G$-space $\bk(\Uc,-)$ is $G$-cofibrant.
The actions of $G$ respectively $\Sigma_n$ on an $n$-tuple of
orthogonal subspaces are componentwise respectively by permuting the entries, i.e., 
\[ (\sigma,g)\cdot (E_1,\dots,E_n) \ = \ 
(g\cdot E_{\sigma^{-1}(1)},\dots,g\cdot E_{\sigma^{-1}(n)}) \]
for $(\sigma,g)\in\Sigma_n\times G$.
The topology on $\bk(\Uc,n_+)$ is that of a disjoint union of subspaces of a product
of Grassmannians:
\[ \bk(\Uc,n_+) \ = \ \coprod_{(i_1,\dots,i_n)\in\mN^n}\, \bk(\Uc;i_1,\dots,i_n) \ ,\]
where
$\bk(\Uc;i_1,\dots,i_n)$ is the subspace of those tuples
such that $\dim(E_j)=i_j$ for all $1\leq j\leq n$.

The summands are invariant under the $G$-action, 
and $\bk(\Uc;i_1,\dots,i_n)$ is $G$-equivariantly homeomorphic to
\[  \bL^\mC(\mC^{i_1+\dots+i_n},\Uc) /  U(i_1)\times\dots\times U(i_n)\  , \]
where $\bL^\mC(-,-)$ is the space of $\mC$-linear isometric embeddings.
This $G$-space is $G$-cofibrant by the unitary analog of
Proposition \ref{prop:K G cofibration}~(ii).
We let $\Gamma\subset\Sigma_n$ 
be the stabilizer group of the dimension vector $(i_1,\dots,i_n)$,
i.e., the group of those $\sigma\in \Sigma_n$ such that $i_j=i_{\sigma(j)}$ 
for all $1\leq j\leq n$. Then $\Gamma$ acts trivially on $\bk(\Uc;i_1,\dots,i_n)$,
so this space is $(\Gamma\times G)$-cofibrant.
Hence the induced $(\Sigma_n\times G)$-space
\[ \Sigma_n \times_{\Gamma}  \bk(\Uc; i_1,\dots,i_n)  \]
is $(\Sigma_n\times G)$-cofibrant.
The full configuration space $\bk(\Uc,n_+)$
decomposes as the disjoint union of $(\Sigma_n\times G)$-invariant subspaces
indexed by the $\Sigma_n$-orbits on $\mN^n$,
\[ \bk(\Uc,n_+) \ = \ 
\coprod_{i \Sigma_n \in \mN^n/\Sigma_n}\, \Sigma_n\times_{\Gamma_i}  \bk(\Uc; i) \ , \]
where $\Gamma_i\subset\Sigma_n$ is the stabilizer of $i=(i_1,\dots,i_n)$.
Altogether this shows that $\bk(\Uc,n_+)$ is a disjoint union of 
$(\Sigma_n\times G)$-cofibrant spaces. 
The latching map
\[ l_n \ : \ \colim_{U\subsetneq\{1,\dots,n\}} \, \bk(\Uc,U_+)\ \to \  \bk(\Uc,n_+)  \]
is a closed embedding 
by Proposition \ref{prop:Gamma latching closed embedding}~(ii),
and its image is the disjoint union of the
summands $\bk(\Uc; i_1,\dots,i_n)$ such that $i_j=0$ for some $j\in\{1,\dots,n\}$.
Since all additional summands not in the image are $(\Sigma_n\times G)$-cofibrant, 
the map $l_n$ is a $(\Sigma_n\times G)$-cofibration.
\end{eg}

\begin{prop}
The orthogonal spectra $\bku$ and $\bko$ are globally connective.\index{subject}{globally connective}  
\end{prop}
\begin{proof}
  The arguments for $\bku$ and $\bko$ are completely parallel, and we concentrate
  on the complex case.  We let $G$ be any compact Lie group, and we show 
  that the group $\pi_{-k}^G(\bku)$ is trivial for all $k\geq 1$.
  We let $V$ be a $G$-representation and set $\Uc=\Sym((\mR^k\oplus V)_\mC)$,
  a unitary $G$-representation of countably infinite dimension.
  We let $H$ be a closed subgroup of $G$.
  The $\bGamma$-$H$-space $\bk(\Uc,-)$ is $H$-cofibrant by
  Example \ref{eg:ku cofibrant}.  
  So the space $\bk(\Uc,S^{\mR^k\oplus V})^H$
  is $(k+\dim(V^H)-1)$-connected by Proposition \ref{prop:Gamma on S^V}~(i). 

  The cellular dimension of $S^V$ at $H$,
  in the sense of \cite[II.2, p.\,106]{tomDieck-transformation},
  is the topological dimension of the space $(S^{V^H})/N_G H$;
  this cellular dimension is at most $\dim(V^H)$. 
  Because $k$ is positive, the cellular dimension of $S^V$ at $H$ does not exceed
  the connectivity of $\bku(\mR^k\oplus V)^H$.
  So every based continuous $G$-map $S^V\to \bku(\mR^k\oplus V)$ is equivariantly
  null-homotopic by \cite[II Prop.\,2.7]{tomDieck-transformation},
  and the set $[S^V,\bku(\mR^k\oplus V)]^G$ has only one element.
  Passage to the colimit over $V\in s(\Uc_G)$ proves the claim.  
\end{proof}

Every $\mC$-linear isometric embedding $u:\Uc\to \Vc$ 
of complex inner product spaces induces a morphism of $\bGamma$-spaces
$\bk(u):\bk(\Uc)\to\bk(\Vc)$ by applying $u$ elementwise
to a tuple of orthogonal subspaces.
So if a compact Lie group $G$ acts on $\Uc$ by linear isometries
(for example if $\Uc$ is a $G$-universe), 
then the $\bGamma$-space $\bk(\Uc)$ inherits a $G$-action, 
so it becomes a $\bGamma$-$G$-space.

\begin{prop}\label{prop:configuration_invariance}
Let $G$ be a compact Lie group and $\Uc$ and $\Vc$ two isomorphic
complex $G$-universes. Then for every $G$-equivariant linear isometric embedding
$u:\Uc\to \Vc$ and every based $G$-space $K$ the map
\[ \bk(u,K)\ :\ \bk(\Uc,K)\ \to \ \bk(\Vc,K)\]
is a $G$-homotopy equivalence.
\end{prop}
\begin{proof}
We start with the special case where $\Vc=\Uc$.
The space of $G$-equi\-variant linear isometric embeddings from $\Uc$
to itself is contractible. A homotopy from  $u$
to the identity then induces a $G$-homotopy from $\bk(u,K)$
to the identity of the $G$-space $\bk(\Uc,K)$.
In the general case we choose a 
$G$-equivariant linear isometry $v:\Vc\iso \Uc$.
Then the two $G$-maps
\[ \bk(v u,K) \colon \bk(\Uc,K)\ \to \ \bk(\Uc,K)\text{\ \ and\ \ }
 \bk(u v,K) \colon \bk(\Vc,K)\ \to \ \bk(\Vc,K) \]
are $G$-homotopic to the respective identity maps.
\end{proof}

If $F$ is any $\bGamma$-space and $S$ a finite set, then we define the map
\[ P_S \ : \  F(S_+)\ \to \ \map(S, F(1_+))    \]
by $P_S(x)(s)= F(p_s)(x)$, where $p_s:S_+\to 1_+$ sends $s$ to $1$ and all other elements
of $S_+$ to the basepoint. 

Now we suppose that $G$ is a compact Lie group and $F$ is a $\bGamma$-$G$-space.
When we evaluate $F$ (or rather its prolongation)
at a $G$-space $K$, it comes with a $(G\times G)$-action;
one action comes from the `external' action on $F$, the other one from the
$G$-action on $K$ via the continuous functoriality of $F$. In such a situation,
we always consider $F(K)$ as a $G$-space 
via the diagonal action of this $(G\times G)$-action.
The map $P_S$ is natural for bijections in $S$;
so whenever the group $G$ acts on $S$, then $P_S:F(S_+)\to \map(S,F(1_+))$ 
is $G$-equivariant for the diagonal action on the source 
and the conjugation action on the target.

We recall from Definition \ref{def:special G-Gamma}
that a $\bGamma$-$G$-space $F$ is {\em special}\index{subject}{Gamma-space@$\bGamma$-space!special}\index{subject}{Gamma-space@$\bGamma$-space!equivariant}  
if for every closed subgroup $H$ of $G$ and every finite $H$-set $S$ the map
\[ (P_S)^H \ : \ F(S_+)^H\ \to \ \map^H(S,F(1_+)) \]
is a weak equivalence. It goes without saying that $H$ must act continuously
on $S$ with the discrete topology; so $S$ is a disjoint union of $H$-sets
of the form $H/K$ for finite index subgroups $K$ of $H$.
Equivalently, $F$ is special if and only if the map
\[ P_n \ : \ F(n_+)\ \to \ \map(\{1,\dots,n\},F(1_+)) \ = \ F(1_+)^n\]
is an $\Fc(G;\Sigma_n)$-equivalence for every $n\geq 1$, where
$\Fc(G;\Sigma_n)$ is the family of graph subgroups, 
see Proposition \ref{prop:special characterizations}.\index{subject}{graph subgroup}

\begin{theorem}\label{thm:bk special} 
Let $G$ be a compact Lie group and $\Uc$ a complete complex $G$-universe.
\begin{enumerate}[\em (i)]
\item 
The $\bGamma$-$G$-space $\bk(\Uc,-)$ is special.
\item
If $G$ is finite and $V$ and $W$ are two $G$-representations with $W^G\ne \{0\}$, 
then the adjoint assembly map
\[ \tilde\alpha \ : \ \bk(\Uc,S^W)\ \to \  \map_*(S^V, \bk(\Uc,S^{V\oplus W}))  \]
is a $G$-weak equivalence.
\end{enumerate}
\end{theorem}
\begin{proof}
(i)
We show that for every $n\geq 1$ the $(G\times\Sigma_n)$-map
\[  P_n \ : \ \bk(\Uc,n_+)\ \to \ \bk(\Uc,1_+)^n    \]
is an $\Fc(G;\Sigma_n)$-weak equivalence, where $\Fc(G;\Sigma_n)$
is the family of graph subgroups of $G\times\Sigma_n$.
We define a morphism of $(G\times\Sigma_n)$-spaces
\[ \lambda_n\ : \ \bk(\Uc,1_+)^n \ \to \ \bk(\mC^n\tensor\Uc,n_+) \ ;\]
here the $\Sigma_n$-action on the target is diagonally, from
the permutation action on $n_+$ and on the tensor factor $\mC^n$.
The map $\lambda_n$ sends an $n$-tuple $( E_1,\dots, E_n)$
of finite-dimensional $\Uc$-subspaces to the configuration 
\[ ( e_1\tensor E_1,\dots, e_n\tensor E_n ) \]
of pairwise orthogonal subspaces of $\mC^n\tensor\Uc$, 
where $e_j=(0,\dots,0,1,0,\dots,0)$ is the $j$-th vector of the
standard basis of $\mC^n$. 

Now we consider the two composites $\lambda_n\circ P_n$ 
and $P_n'\circ\lambda_n$, where $P_n'$ is the map $P_n$
for the universe $\mC^n\tensor\Uc$ (as opposed to $\Uc$):
\[
\bk(\Uc,n_+)\ \xra{\ P_n\ } \
\bk(\Uc,1_+)^n\ \xra{\ \lambda_n\ } \
 \bk(\mC^n\tensor \Uc,n_+) \ \xra{\ P_n'\ }\ \bk(\mC^n\tensor\Uc,1_+)^n \]
We claim that both composites are $\Fc(G;\Sigma_n)$-weak equivalences.
We start by investigating the composite $\lambda_n\circ P_n$.
For $1\leq j\leq n$, we define a 1-parameter family of unit vectors
\[   u_j \ : \ [0,1] \ \to \ \mC^n
\text{\qquad by\qquad}
u_j(t)\ = \ 
t\cdot e_j + \sqrt{\frac{1-t^2}{n-1}}\cdot\sum_{k\ne j} e_k\ . \]
This provides a homotopy
\[ H \ : \  \bk(\Uc,n_+)\times [0,1]\ \to\ \bk(\mC^n\tensor\Uc,n_+) \]
by
\[ H( (E_1,\dots,E_n),t) \ = \ ( u_1(t)\tensor E_1,\dots,u_n(t)\tensor E_n ) \ .\]
We have
\[ u_1(1/\sqrt{n})\ = \ \dots \ = \ u_n(1/\sqrt{n})\ = \ 
 1/\sqrt{n}\cdot (1,\dots,1) \]
and 
\[  H( (E_1,\dots,E_n),1) \ = \ ( e_1\tensor E_1,\dots,e_n\tensor E_n) \ = \ 
(\lambda_n\circ P_n)(E_1,\dots,E_n)\ . \]
Moreover, each of the maps $H(-,t)$ is $(G\times\Sigma_n)$-equivariant,
so the composite $\lambda_n\circ P_n$ is $(G\times \Sigma_n)$-equivariantly homotopic 
to the map $\bk(u,n_+)$, where
$u:\Uc\to \mC^n\tensor\Uc$ is the $(G\times\Sigma_n)$-equivariant
linear isometric embedding 
\[ u(v) \ = \ 1/\sqrt{n}\cdot(1,\dots,1) \tensor v\ . \]
Now we let $\alpha:H\to \Sigma_n$ be a continuous homomorphism
defined on a closed subgroup $H$ of $G$.
Since $\Uc$ is a complete complex $G$-universe,
both $\Uc$ and $\alpha^*(\mC^n)\tensor\Uc$
are complete complex $H$-universes, 
and $u$ is an $H$-equivariant linear isometric embedding
\[ u \ : \ \res^G_H(\Uc)\ \to \ \alpha^*(\mC^n)\tensor \res^G_H(\Uc)\ . \]
Interpreted in this way, the map $\bk(u,n_+)^H$ is a homotopy equivalence by 
Proposition \ref{prop:configuration_invariance}.
This shows that the morphism $\bk(u,n_+)$ is an $\Fc(G;\Sigma_n)$-weak equivalence,
hence so is the morphism $\lambda_n\circ P_n$.

Now we show that the composite $P_n'\circ\lambda_n$ is an
$\Fc(G;\Sigma_n)$-weak equivalence.
We let $H$ be a closed subgroup of $G$,
$\alpha:H\to \Sigma_n$ a continuous homomorphism,
and $\Gamma=\{ (h,\alpha(h))  \ | \ h\in H\}\leq G\times \Sigma_n$
the graph of $\alpha$.
We let $a_1,\dots, a_k\in \{1,\dots,n\}$ be a set of representatives
of the orbits of the $H$-action on $\{1,\dots,n\}$ through $\alpha$.
We let $H_i\leq H$ be the stabilizer group of $a_i$.
Then projection to the factors indexed by 
$a_1,\dots, a_k$ provides homeomorphisms
\begin{align*}
    (\bk(\Uc,1_+)^n)^\Gamma\quad &\iso \ {\prod}_{i=1}^k \quad  \bk(\Uc,1_+)^{H_i} \text{\quad and}\\
  (\bk(\mC^n\tensor \Uc,1_+)^n)^\Gamma\ &\iso \ 
{\prod}_{i=1}^k\, \bk(\alpha^*(\mC^n)\tensor\Uc,1_+)^{H_i}\ .
\end{align*}
Under these identifications, the map $(P_n'\circ\lambda_n)^\Gamma$ 
becomes the product of the maps
\begin{equation}  \label{eq:H_i_fixed}
 (\bk(\Uc,1_+))^{H_i}\ \to \ \bk(\alpha^*(\mC^n)\tensor\Uc,1_+)^{H_i}   
\end{equation}
induced by the $H_i$-equivariant linear isometric embeddings
\[ \Uc \ \to \ \alpha^*(\mC^n)\tensor\Uc \ ,\quad v \ \longmapsto \ e_{a_i}\tensor v\ .\]
Since $\Uc$ is a complete complex $G$-universe,
both $\Uc$ and $\alpha^*(\mC^n)\tensor\Uc$
are complete complex $H_i$-universes.
So the map \eqref{eq:H_i_fixed} is a homotopy equivalence by 
Proposition \ref{prop:configuration_invariance}.
Since $\lambda_n\circ P_n$ and $P_n'\circ\lambda_n$ are
$\Fc(G;\Sigma_n)$-weak equivalences, 
so is the map $P_n:\bk(\Uc,n_+)\to\bk(\Uc,1_+)^n$.

Part~(ii) is essentially a special case of Segal and Shimakawa's equivariant
delooping machine based on  $\bGamma$-$G$-spaces \cite{segal-some_equivariant,shimakawa}.
There is a bit of extra work necessary, though, because our claim is about
the prolongation $\bk(\Uc,S^{V\oplus W})$ (which is a categorical coend),
whereas Segal and Shimakawa formulate their results in terms of a bar construction
(which amounts to a homotopy coend).  
We provide the additional arguments 
in Theorem \ref{thm:special cofibrant Gamma positive Omega};
to apply this, we need that the $\bGamma$-$G$-space $\bk(\Uc,-)$
is special by part~(i) and $G$-cofibrant by Example \ref{eg:ku cofibrant}.
\end{proof}

The orthogonal spectrum $\bku$ is trying to be a $\Fin$-global $\Omega$-spectrum.
However, the global $\Omega$-spectrum condition on the adjoint structure maps
only holds for `sufficiently large' representations.

\begin{defn}\label{def:ample} 
Let $G$ be a compact Lie group. An orthogonal $G$-re\-presen\-tation $W$
is {\em ample}\index{subject}{ample $G$-representation}\index{subject}{G-representation@$G$-representation!ample|see{ample $G$-representation}}
if the complex symmetric algebra $\Sym(W_\mC)$
is a complete complex $G$-universe.
\end{defn}

\begin{rk} 
Every ample $G$-representation is non-zero and the action of $G$ is faithful.
Examples of ample representations are 
non-zero faithful permutation representations
(which in particular forces the group to be finite). 
Indeed, we let $\mR A$ denote the permutation representation
of a non-empty faithful finite $G$-set $A$.
Then the complexified symmetric algebra $\Sym(\mC A)$
is a complex permutation representation, namely of the infinite
$G$-set $\mN^A=\map(A,\mN)$ of functions from $A$ to $\mN$.
Since $G$ acts faithfully on $A$, every injective map $A\to\mN$ 
generates a free $G$-orbit in the $G$-set $\mN^A$. 
Since $A$ is non-empty, there are infinitely many injections from $A$ to $\mN$ with
pairwise disjoint images, and these generate infinitely many
distinct free $G$-orbits in $\mN^A$. So $\Sym(\mC A)$ contains
infinitely many copies of the complex regular $G$-representation,
and is thus a complete complex $G$-universe.  

A more general class of ample $G$-representations are non-zero representations
containing a vector with trivial isotropy group.
Indeed, given a non-zero $G$-free vector $w_0\in W$, 
we let $\nu$ be the normal space at $w_0$ to the orbit $G w_0$.
A choice of $G$-equivariant tubular neighborhood of the orbit $G w_0$
allows us to $G$-equivariantly embed the Hilbert space $L^2(G;\mC)\tensor \Sym(\nu)$
into $L^2(W;\mC)$.
Since $L^2(G;\mC)\tensor\Sym(\nu)$ contains every finite-dimensional $G$-representation
infinitely often, so does $L^2(W;\mC)$. 
But $L^2(W;\mC)$ is equivariantly isomorphic to 
the Hilbert space completion of $\Sym(W^\ast_\mC)$:
elements of $\Sym(W^\ast_\mC)$ are polynomial functions on $W$,
which we can map to $L^2$-functions by scaling with the function
$w\mapsto \exp(-|w|^2)$. This gives a $G$-equivariant isometric embedding
$\Sym(W^\ast_\mC)\to L^2(W;\mC)$ with dense image.
So if there was any complex $G$-representation that did not embed
equivariantly into $\Sym(W^\ast_\mC)$, then it would also
not embed into $L^2(W;\mC)$, a contradiction.
\end{rk}

\begin{theorem}\label{thm:weak global Omega for bku}
Let $G$ be a finite group and $W$ an ample orthogonal $G$-representation.
Then for every $G$-representation $V$ the adjoint structure map
\[ \tilde\sigma_{V,W}\ : \ \bku(W)\ \to \ \map_*(S^V,\bku(V\oplus W)) \]
is a $G$-weak equivalence. 
\end{theorem}
\begin{proof}
The adjoint structure map $\tilde\sigma_{V,W}$ factors as the composite
\begin{align*}
\bku(W)\ &= \ \bk(\Sym(W_\mC),S^W)\ \xra{\bk(\Sym(i_\mC),S^W) }\
 \bk(\Sym((V\oplus W)_\mC),S^W)\\ 
 &\xra{\ \tilde\alpha\ }\
 \map_*(S^V, \bk(\Sym((V\oplus W)_\mC),S^{V\oplus W}))\ = \ \map_*(S^V,\bku(V\oplus W)); 
\end{align*}
here $i:W\to V\oplus W$ is the inclusion of the second summand and
$\tilde\alpha$ is the adjoint assembly map for the $\bGamma$-$G$-space
$\bk(\Sym((V\oplus W)_\mC))$.
Since $W$ is ample, the map $\Sym(i_{\mC}):\Sym(W_\mC)\to\Sym((V\oplus W)_\mC)$
is an equivariant linear isometric embedding between complete complex
$G$-universes. So the first map is a $G$-homotopy equivalence by 
Proposition \ref{prop:configuration_invariance}.
The second map is a $G$-weak equivalence by Theorem \ref{thm:bk special}~(ii).
\end{proof}

Now we will justify that for every finite group $G$ 
the underlying orthogonal $G$-spectrum of $\bku$ represents 
connective $G$-equivariant topological $K$-theory.

\begin{construction}
We define a morphism of orthogonal spaces
\begin{equation}\label{eq:define_c}
 c \ : \ \bGr^\mC \ \to \ \Omega^\bullet \bku \ ,
\end{equation}
where $\bGr^\mC$ is the complex additive Grassmannian
introduced in Example \ref{eg:Gr additive complex}.\index{subject}{additive Grassmannian}
The value at an inner product space $V$ is the continuous map
\[ c(V) \ : \ G r^\mC(V_\mC) \ \to \ \map_*(S^V,\bku(V))\ = \ (\Omega^\bullet \bku)(V)\]
that sends a complex subspace $L\subset V_\mC$ to the continuous based map
\[ [L;-]\ : \ S^V \ \to \ \bk(\Sym(V_\mC),S^V)\ = \ \bku (V) \ ,\quad
v \ \longmapsto \ [L;\,v]\ .\]
Here we consider $L\subset V_\mC$ as sitting in the linear summand of
the symmetric algebra of $V_\mC$.
\end{construction}

It will be useful to also define a delooping of the morphism $c$, namely
a morphism of orthogonal spaces
\begin{equation}\label{eq:define_eig}
\eig \ : \ \bU \ \to \ \Omega^\bullet(\sh \bku)\ , 
\end{equation}
where $\bU$ is the ultra-commutative monoid of unitary groups
(compare Example \ref{eg:unitary group monoid space}),
and $\sh\bku=\sh^\mR\bku$ is the shift of $\bku$ as 
defined in \eqref{eq:define_shift}.
The name `eig' refers to the fact that the morphism records the eigenspace
decomposition of a unitary isomorphism, compare Remark \ref{rk:eigenspace decomposition}.
The definition uses the homeomorphism
\[ h \ : \ U(1)\ \iso \ S^1 \ , 
\quad h(\lambda)\ = \ i\cdot\frac{\lambda+1}{\lambda-1}\ ,\]
the inverse of the Cayley transform.\index{subject}{Cayley transform}
Every unitary automorphism of a finite-dimensional hermitian inner
product space is diagonalizable with eigenvalues in $U(1)$
and pairwise orthogonal eigenspaces.
So given an inner product space $V$ we define
\[ \eig(V) \ : \ \bU(V)\ = \ U(V_\mC) \ \to \ \map_*(S^V, \bku(V\oplus\mR) )
= \left( \Omega^\bullet(\sh\bku) \right)(V) \]
by
\[ \eig(V)(A)(v) \ = \ 
[E(\lambda_1),\dots, E(\lambda_n);\,(v,h(\lambda_1)),\dots,(v,h(\lambda_n)) ] \ .\]
Here $\lambda_1,\dots\lambda_n\in U(1)$ are the eigenvalues of 
$A$ and $E(\lambda_i)$ is the eigenspace of $\lambda_i$.
Strictly speaking, $E(\lambda_i)$ is a subspace of $V_\mC$, which we embed into
the linear summand of $\Sym((V\oplus\mR)_\mC)$.

\begin{theorem}\label{thm:U_vs_shift_ku}
The morphism
\[ \eig \ : \  \bU \ \to \ \Omega^\bullet(\sh\bku) \]
is a $\Fin$-global equivalence of orthogonal spaces.
\end{theorem}
\begin{proof}
  We let $\bar\bU$ denote the orthogonal space with
\[ \bar\bU(V)\ = \ U(\Sym((V\oplus\mR)_\mC))\ . \]
The structure map induced by $\varphi:V\to W$ 
is given by extending by the identity on the orthogonal complement
of $\Sym((\varphi\oplus\mR)_\mC):\Sym((V\oplus\mR)_\mC)\to \Sym((W\oplus\mR)_\mC)$.
The eigenspace decomposition map then factors as the composite
\[ \bU \ \to \ \bar\bU \ 
\xra{\ \widebar\eig \ } \ \Omega^\bullet(\sh \bku) \ , \]
where $\widebar\eig$ is defined in the same way as $\eig$, recording the 
set of eigenvalues and eigenspaces.
Since $\bar\bU(V)$ is the colimit over $n\geq 0$, along closed embeddings,
of the spaces $U(\Sym^{\leq n}((V\oplus\mR)_\mC))$,
the morphism $\bU\to\bar\bU$ is a global equivalence of orthogonal spaces
by Theorem \ref{thm:general shift of osp}
and Proposition \ref{prop:global equiv basics}~(ix). 
We may thus show that the morphism
\[\widebar\eig  \ : \ \widebar\bU \ \to\ \Omega^\bullet(\sh \bku) \]
is a $\Fin$-global equivalence.
We show the stronger statement that for every finite group $G$ and
every ample $G$-representation $V$ the map
$\widebar\eig(V)$ is a $G$-weak equivalence. 
Indeed, the map $\widebar\eig(V)$ factors as the composite
\begin{align*}
   U(\Sym((V\oplus\mR)_\mC)) \ &\xra{\ \iso\ } \  \bk(\Sym( (V_\mC\oplus\mR)_\mC), S^1) \\ 
&\xra{\ \sim \ } \  \map_*(S^V, \bk(\Sym( (V\oplus\mR)_\mC), S^{V\oplus\mR})) \ .
\end{align*}
The first map is the eigenspace decomposition, hence a homeomorphism
by Remark \ref{rk:eigenspace decomposition}.
The second map is adjoint to the assembly map
 \[ S^V\sm \bk(\Sym((V\oplus\mR)_\mC), S^1) \ \to \  
 \bk(\Sym( (V\oplus\mR)_\mC), S^{V\oplus\mR})   \]
of the $\bGamma$-$G$-space $\bk(\Sym((V\oplus\mR)_\mC))$.
Since $V$ is ample, $\Sym((V\oplus\mR)_\mC)$ is a complete complex $G$-universe.
So the adjoint assembly map is a $G$-weak equivalence 
by Theorem \ref{thm:bk special}~(ii).
\end{proof}

The set $[A,\bGr^\mC]^G$ has an abelian monoid structure arising from
the ultra-commutative multiplication of $\bGr^\mC$ 
as explained in \eqref{eq:addition_on_[A,R]^G}.
The set $[A, \Omega^\bullet\bku]^G$ has an abelian group structure 
as an equivariant stable homotopy group, i.e., through the
adjunction bijection
\[ [A, \Omega^\bullet\bku]^G\ \iso \ \pi_0^G(\map(A,\bku)) \ ;\]
so this group structure arises from concatenation of loops.
The morphism of orthogonal spaces $c: \bGr^\mC\to\Omega^\bullet \bku$
defined in \eqref{eq:define_c} is {\em not} a homomorphism of ultra-commutative monoids,
nor is it a loop map; so it is not a priori clear whether
the induced map on equivariant homotopy sets is a monoid homomorphism.

\begin{theorem}\label{thm:[A,Gr] to [A,bku]}
Let $G$ be a compact Lie group and $A$ a finite $G$-CW-complex.
Then the map
\[ [A,c]^G\ : \ [A,\bGr^\mC]^G \ \to\   [A, \Omega^\bullet\bku]^G\]
is a monoid homomorphism. If all isotropy groups of $A$ are finite, then
$[A,c]^G$ is a group completion of abelian monoids.
\end{theorem}
\begin{proof}
To show that $[A,c]^G$ is a monoid homomorphism we exhibit a `delooping' of $c$,
namely the eigenspace morphism \eqref{eq:define_eig}.
In \eqref{eq:define_beta}
we defined a morphism of ultra-commutative monoids
$\beta:\bGr^\mC\to\Omega \bU$ that is a global group completion 
by Theorem \ref{thm:group completion Gr to Omega U}.
Now we link the monoid homomorphism
\[ [A,\beta]^G\ : \ [A,\bGr^\mC]^G\ \to \ [A,\Omega\bU]^G \]
to the set map $[A,c]^G:[A,\bGr^\mC]^G \to[A, \Omega^\bullet\bku]^G$.
The global equivalence
\[ \tilde\lambda_{\bku}\ : \ \bku \ \to \ \Omega(\sh\bku ) \]
was defined in \eqref{eq:defn lambda_n}.
The adjunction isomorphisms
\begin{align*}
  (\Omega^\bullet (\Omega(\sh \bku)))(V)\ &= \ \map_*(S^V,\Omega \bku(V\oplus\mR))\\ 
&\iso \ 
\map_*(S^{V\oplus \mR}, \bku(V\oplus\mR))\ = \ (\sh_\oplus(\Omega^\bullet \bku))(V) 
\end{align*}
provide an isomorphism of orthogonal spaces between
$\Omega^\bullet (\Omega(\sh \bku))$ and
$\sh_\oplus(\Omega^\bullet \bku)$, where $\sh_\oplus=\sh_\oplus^\mR$
is the additive shift defined in Example \ref{eg:Additive and multiplicative shift}.
Under this isomorphism, the morphism
\[ \Omega^\bullet\tilde\lambda_{\bku} \ : \ 
\Omega^\bullet \bku \ \to \Omega^\bullet(\Omega(\sh \bku)) \]
becomes the morphism
\[ (\Omega^\bullet\bku)\circ i \ : \  \Omega^\bullet \bku \ 
\to \ \sh_\oplus(\Omega^\bullet \bku) \]
given by precomposition with the embeddings $i_V:V\to V\oplus\mR$.
This morphism is a global equivalence by Theorem \ref{thm:general shift of osp}.
So the morphism $\Omega^\bullet\tilde\lambda_{\bku}$ is a global equivalence
of orthogonal spaces as well.

The rigorous statement of the delooping property of the eigenspace morphism 
is the following commutative square of orthogonal spaces:
\[  \xymatrix@C=12mm{
\bGr^\mC \ar[d]_{\beta} \ar[r]^-c & 
\Omega^\bullet \bku\ar[d]^{\Omega^\bullet \tilde\lambda_{\bku}}_\simeq\\
\Omega \bU \ar[r]_-{\Omega\eig} &  \Omega (\Omega^\bullet (\sh\bku)) }     
 \]
This square induces a commutative square of set maps
\begin{equation}\begin{aligned}\label{eq:[A,-]^G square}
  \xymatrix@C=15mm{
[A, \bGr^\mC]^G \ar[d]_{[A,\beta]^G} \ar[r]^-{[A,c]^G} & 
[A,\Omega^\bullet \bku]^G\ar[d]^{[A,\Omega^\bullet \tilde\lambda_{\bku}]^G}_\iso\\
[A,\Omega \bU]^G \ar[r]_-{[A,\Omega \eig]^G} & 
[A,\Omega\left( \Omega^\bullet (\sh \bku)\right)]^G
} 
 \end{aligned}\end{equation}
The set $[A,\Omega\bU]^G$ can be endowed with a monoid structure 
in two ways, via the ultra-commutative multiplication
as in \eqref{eq:addition_on_[A,R]^G}, and by concatenation of loops. 
Since the ultra-commutative monoid structure of $\Omega\bU$ 
is `pointwise' (i.e., induced by the ultra-commutative monoid structure of $\bU$),
these two monoid structures satisfy the interchange law,
so they coincide, and both are abelian group structures. 
The morphism $\beta$ is a homomorphism of ultra-commutative monoids,
so it induces an additive map on $[A,-]^G$.
The morphisms $\Omega\eig$ and $\Omega^\bullet\tilde\lambda_{\bku}$
are loop maps, so they induce homomorphisms with respect to the
group structure by concatenation of loops. This shows that three of the four
maps in \eqref{eq:[A,-]^G square} are homomorphisms of abelian monoids.
Since the right vertical map is an isomorphism, the map 
$[A,c]^G$ is a homomorphism as well.

For the rest of the proof we suppose that the isotropy groups of $A$ are finite.
The morphism $\Omega\eig$ is a $\Fin$-global equivalence
by Theorem \ref{thm:U_vs_shift_ku}, so the map $[A,\Omega\eig]^G$
is bijective, and hence an isomorphism of abelian groups,
by Proposition \ref{prop:[A,Y]^G of closed}~(ii).
The morphism $\Omega^\bullet\tilde\lambda_{\bku}:\Omega^\bullet\bku\to\Omega (\Omega^\bullet(\sh\bku))$
is a global equivalence of orthogonal spaces,
so the map $[A,\Omega^\bullet\tilde\lambda_{\bku}]^G$ is bijective, 
and hence an isomorphism of abelian groups,
again by Proposition \ref{prop:[A,Y]^G of closed}~(ii).
Since $[A,\beta]^G$ is a group completion of abelian monoids
(by Corollary \ref{cor:group completion Gr to Omega U}),
the commutative square \eqref{eq:[A,-]^G square} 
shows that $[A,c]^G$ is a group completion.
\end{proof}

We draw an important consequence of Theorem \ref{thm:[A,Gr] to [A,bku]},
namely that the equivariant cohomology theory represented by the 
connective global $K$-theory spectrum $\bku$ is essentially equivariant $K$-theory.
There is a caveat, however, as this is {\em not} true on arbitrary
finite $G$-CW-complexes, but only under the hypothesis of finite stabilizer groups.

We let $G$ be a compact Lie group and $A$ a compact $G$-space. 
We define the 0-th $G$-equivariant $\bku$-cohomology group of $A$ as
\[ \bku^0_G(A_+)\ = \ [A,\Omega^\bullet \bku]^G \ , \]
the equivariant homotopy set into the orthogonal space $\Omega^\bullet\bku$.
This set is an abelian group by concatenation of loops,
i.e., via the adjunction bijection
\[  [A,\Omega^\bullet \bku]^G \ \iso \ \pi_0^G(\map(A,\bku)) \]
to an equivariant stable homotopy group. 
The set also has a multiplication, i.e., another
commutative binary operation arising from the ring spectrum structure
of $\bku$, which turns $\Omega^\bullet \bku$ into a `multiplicative'
ultra-commutative monoid 
as in Example \ref{eg:suspension spectrum of orthogonal monoid space}. 
A conjugation involution of the ultra-commutative ring spectrum $\bku$
was defined 
in Construction \ref{con:conjugation on bku}.\index{subject}{complex conjugation!on $\bku$} 

We denote by $\bK_G(A)$\index{symbol}{$\bK_G(A)$ - {equivariant $K$-group of $A$}}\index{subject}{equivariant $K$-theory}
the $G$-equivariant $K$-group of $A$, 
i.e., the group completion (Grothendieck group) of the abelian monoid
of isomorphism classes of complex $G$-vector bundles over $A$.
The map 
\[ \td{-}\ : \ [A,\bGr^\mC]^G \ \to \ \bK_G(A) \ , 
\quad [f:A\to G r^\mC(V)]\ \longmapsto \ [f^\star(\gamma_V^\mC)]\]
that takes the pullback of the tautological vector bundles
over Grassmannians is a group completion of abelian monoids by
Theorem \ref{thm:BOP_to_KO_G} (or rather its complex analog).
Since $[A,c]^G:[A,\bGr^\mC]^G\to [A,\Omega^\bullet\bku]^G$
is a monoid homomorphism  (by Theorem \ref{thm:[A,Gr] to [A,bku]})
to an abelian group, the universal property of group completion
provides a unique additive extension 
\begin{equation}\label{eq:K_G_to_ku} 
  [-]\ : \ \bK_G(A)\ \to \ \bku_G^0(A_+) 
\end{equation}
such that the composite
\[ [A,\bGr^\mC]^G \ \xra{\ \td{-}\ }\ \bK_G(A)\ \xra{\ [-]\ } \ \bku_G^0(A_+) 
\ = \ [A,\Omega^\bullet \bku]^G  \]
is the map induced by the morphism of orthogonal spaces $c:\bGr^\mC\to\Omega^\bullet\bku$.

\index{subject}{equivariant $K$-theory|(}
\begin{theorem}\label{thm:K_G_is_ku} 
Let $G$ be a compact Lie group and $A$ a finite $G$-CW-complex.
\begin{enumerate}[\em (i)]
\item For every $G$-representation $V$ and every 
continuous $G$-map $f:A\to G r^\mC(\Sym(V_\mC))$, the homomorphism $[-]$
sends the class of the $G$-vector bundle $f^\star(\gamma^\mC_{\Sym(V_\mC)})$
to the homotopy class of the $G$-map
\[ A \ \xra{[f(-);\, -]} \ \map_*(S^V,\bku(V)) \ ,\quad
a \ \longmapsto \ \{ v\ \longmapsto \ [f(a);v] \} \ .\]
\item The additive map $[-]$ is a ring homomorphism, natural for $G$-maps in $A$,
natural for restriction homomorphisms in $G$,
and compatible with complex conjugation.
\item If $A$ has finite isotropy groups, 
then the homomorphism $[-] :\bK_G(A)\to \bku_G^0(A_+)$ is an isomorphism.
\end{enumerate}
\end{theorem}
\begin{proof}
(i)
We recall from Example \ref{eg:Gr multiplicative real}
the multiplicative Grassmannian $\bGr_\tensor^\mC$\index{subject}{multiplicative Grassmannian}
with values
\[ \bGr_\tensor^\mC(V)\ = \ {\coprod}_{n\geq 0} \, Gr_n^\mC(\Sym(V_\mC)) \ ,\]
the disjoint union of all Grassmannians in the symmetric algebra of $V_\mC$.
We let $i:V_\mC\to \Sym(V_\mC)$ be the embedding as the linear summand of the
symmetric algebra. Then as $V$ varies, the maps
\[ i(V)\ : \ \bGr^\mC(V)\ \to\  \bGr_\tensor^\mC(V)\ , \quad L\ \longmapsto \ i(L) \]
form a global equivalence $i:\bGr^\mC\to\bGr_\tensor^\mC$ of orthogonal spaces,
see Example \ref{eg:Gr multiplicative real}.
The morphism of orthogonal spaces
$c :\bGr^\mC\to \Omega^\bullet \bku$ defined in \eqref{eq:define_c}
has an extension
\[ c_\tensor  \ : \ \bGr_\tensor^\mC \ \to \ \Omega^\bullet \bku \ ,  \]
defined by the same formula as for $c$, namely
\begin{align*}
 c_\tensor(V) \ : \ \bGr_\tensor^\mC(V) \ \to \ 
\map_*(S^V,\bku(V))\ = \ (\Omega^\bullet \bku)(V)\ ,\quad  
L \ \longmapsto \ [L;-]
\end{align*}
with
\[ [L;-]\ : \ S^V \ \to \ \bk(\Sym(V_\mC),S^V)\ = \ \bku (V) \text{\quad given by\quad}
v \ \longmapsto \ [L;\,v]\ .\]
The `pullback bundle' map $\td{-}:[A,\bGr^\mC]^G\to \bK_G(A)$
also has a straightforward extension
\[ \td{-}_\tensor\ :\ [A,\bGr^\mC_\tensor]^G\ \to\ \bK_G(A) \]
again given by the same recipe: we send the class represented by
a $G$-map $f:A\to G r^\mC(\Sym(V_\mC))$ to the class of the pullback
$f^\star(\gamma_{\Sym(V_\mC)}^\mC)$ of the tautological bundle over 
$G r^\mC(\Sym(V_\mC))$.
The outer square in the diagram
\[ \xymatrix@C=15mm{ 
[A,\bGr^\mC]^G\ar[r]^-{\td{-}}\ar[d]_{[A,i]^G}^\iso &
\bK_G(A)\ar[d]^{[-]} \\
[A,\bGr^\mC_\tensor]^G\ar[r]_-{[A,c_\tensor]^G} 
\ar[ur]^(.4){\td{-}_\tensor}&[A,\Omega^\bullet\bku]^G } \]
commutes because $c_\tensor\circ i=c$ and $[A,c]^G=[-]\circ\td{-}$.
The upper left triangle commutes because 
the tautological bundle on $\bGr^\mC_\tensor(V)$
restricts to the tautological bundle on $\bGr^\mC(V)$
along the map $i(V):\bGr^\mC(V)\to \bGr^\mC_\tensor(V)$.
The left vertical map is bijective
by Proposition \ref{prop:[A,Y]^G of closed}~(ii), 
because $i:\bGr^\mC\to\bGr_\tensor^\mC$ is a global equivalence.
So the lower right triangle commutes as well.
The map
\[ [A,c_\tensor]^G\ : \ [A,\bGr^\mC_\tensor]^G \ \to \ 
[A, \Omega^\bullet\bku]^G \ = \ \bku^0_G(A_+)\]
sends the homotopy class of a continuous $G$-map $f:A\to G r^\mC(\Sym(V_\mC))$ 
to the class of the map $[f(-);-]:A\to\map_*(S^V,\bku(V))$, 
by the very definition of $c_\tensor$.
So this proves the claim.

(ii) 
Naturality in $G$-maps is straightforward.
For a $G$-map $h:A'\to A$ and every $G$-map $f:A\to\bGr^\mC(V)$,
the two $G$-vector bundles $h^\star(f^\star(\gamma_V^\mC))$
and $(f h)^\star(\gamma_V^\mC)$ over $A'$ are isomorphic, 
so the following square commutes:
\[ \xymatrix@C=10mm{ 
[A,\bGr^\mC]^G \ar[r]^-{\td{-}}\ar[d]_{[h,\bGr^\mC]^G} &\bK_G(A)\ar[d]^{\bK_G(h)} \\
[A',\bGr^\mC]^G \ar[r]_-{\td{-}}&\bK_G(A')
} \]
Together with the defining property of the homomorphism $[-]$, 
this implies that the two group homomorphisms
\[ [-]\circ\bK_G(h)\ , \ \bku_G^0(h_+)\circ [-]\ : \ 
\bK_G(A)\ \to \ \bku_G^0(A'_+) \]
coincide after precomposition with $\td{-}:[A,\bGr^\mC]^G\to\bK_G(A)$.
Since this map is a group completion of abelian monoids, already
the homomorphisms $[-]\circ\bK_G(h)$ and $\bku_G^0(h_+)\circ [-]$ agree, 
by the universal property of group completions.

The compatibility with complex conjugation and restriction along
group homomorphisms follow the same pattern.
The conjugation morphism of orthogonal spectra $\psi:\bku\to\bku$
deloops the conjugation morphism of ultra-commutative monoids $\psi:\bGr^\mC\to\bGr^\mC$,
in the sense that the square of orthogonal spaces commutes:
\[ \xymatrix{
\bGr^\mC \ar[d]_{\psi} \ar[r]^-c  & \Omega^\bullet \bku \ar[d]^{\Omega^\bullet\psi}
\\
\bGr^\mC \ar[r]_-c  & \Omega^\bullet \bku 
 } \]
In particular, the homomorphism $[A,c]^G:[A,\bGr^\mC]^G\to \bku_G^0(A_+)$
commutes with complex conjugation.
On the other hand, for every $G$-map $f:A\to\bGr^\mC(V)$,
the bundle $(\psi(V)\circ f)^\star(\gamma_V^\mC)$ is isomorphic
to the complex conjugate of the bundle $f^\star(\gamma_V^\mC)$.
So the following square commutes:
\[ \xymatrix@C=10mm{ 
[A,\bGr^\mC]^G \ar[r]^-{\td{-}}\ar[d]_{[A,\psi]^G} &\bK_G(A)\ar[d]^{\psi} \\
[A,\bGr^\mC]^G \ar[r]_-{\td{-}}&\bK_G(A) } \]
Thus the two group homomorphisms
\[ [-] \circ\psi\ , \ \psi\circ [-]\ : \ \bK_G(A)\ \to \ \bku_G^0(A_+) \]
agree after precomposition with a group completion, hence they coincide.
The analogous argument works for restriction
along a continuous homomorphism $\alpha:K\to G$,
using that the underlying $K$-vector bundle of $f^\star(\gamma_V^\mC)$
equals the bundle $\alpha^*(f)(\gamma_{\alpha^*(V)}^\mC)$, and 
the effect of $c$ (being a morphism of orthogonal spaces) commutes with restriction.

Now we show that the additive map $[-]:\bK_G(A)\to\bku_G^0(A_+)$
is a ring homomorphism. The multiplicative unit of $\bK_G(A)$ is the class
of the trivial line bundle $A\times\mC$. 
This bundle is isomorphic to $f^\star(\gamma_{\Sym(0)}^\mC)$ for
the constant map $f:A\to G r^\mC(\Sym(0))$ with value $\mC$,
the constant (and only non-trivial) summand in the symmetric algebra associated 
to the 0-dimensional $G$-representation.
The associated $G$-map $[f(-);-]:A\to \map_*(S^0,\bku(0))$
is constant with image the unit map $\iota_0:S^0\to\bku(0)$;
so by part~(i) the class $\td{A\times\mC}$ of the trivial line bundle
is the multiplicative unit in the ring $\bku_G^0(A_+)$.
So the map $[-]$ preserves multiplicative units.

It remains to show that the map $[-]$ preserves products.
Because $[-]$ is additive and the group $\bK_G(A)$ is generated by
classes of actual vector bundles (as opposed to virtual  bundles),
it suffices to show multiplicativity for two classes in $\bK_G(A)$ represented 
by $G$-vector bundles. We may assume that the two bundles are classified by
continuous $G$-maps
\[ f \ : \ A \ \to \ G r^\mC(V_\mC) \text{\quad respectively\quad}
 g \ : \ A \ \to \ G r^\mC(W_\mC) \ ,\]
where $V$ and $W$ are $G$-representations.
The tensor product bundle $f^\star(\gamma^\mC_V) \tensor g^\star(\gamma^\mC_W)$ 
is then classified by the composite
\begin{align*}
 A \ \xra{\ (f,g)\ } \ G r^\mC(V_\mC)\times G r^\mC(W_\mC)\ 
&\xra{-  \tensor -}\ G r^\mC(V_\mC\tensor W_\mC) \\ 
&\xra{G r^\mC(j)} \ G r^\mC(\Sym( (V\oplus W)_\mC)) \ ,
\end{align*}
where 
\[ j \ : \ V_\mC\tensor W_\mC\ \to \ \Sym((V\oplus W)_\mC) \text{\quad is defined by\quad}
j(v\tensor w)=(v,0)\cdot(0,w)\ .\]
By part~(i) the associated homotopy class 
$[f^\star(\gamma^\mC_V) \tensor g^\star(\gamma^\mC_W)]$
is represented by the $G$-map 
\[ [ (j\circ (f\tensor g))(-);-] \ : \ 
A\ \to \ \map_*( S^{V\oplus W},\bku(V\oplus W))\ . \]
This map is the adjoint of the composite 
\[ 
A_+\sm S^{V\oplus W}\ \xra{a\sm v\sm w\mapsto [f(a);v]\sm [g(a);w]}\ 
\bku(V)\sm \bku(W)\ \xra{\mu_{V,W}}\ \bku(V\oplus W)\ .\]
This composite represents the product
of the classes $[f^\star(\gamma_V^\mC)]$ and $[g^\star(\gamma_W^\mC)]$, 
so this establishes the relation
\[ [f^\star(\gamma_V^\mC)]\ \cdot\ [g^\star(\gamma_W^\mC)]\ = \  
[f^\star(\gamma^\mC_V) \tensor g^\star(\gamma^\mC_W)] \]
in the group $\bku_G^0(A)$.

(iii) If $A$ is a finite $G$-CW-complex with finite isotropy groups,
then the map $[A,c]^G:[A,\bGr^\mC]^G \to \bku^0_G(A_+)$ is a group completion
of abelian monoids by Theorem \ref{thm:[A,Gr] to [A,bku]}.
The map $\td{-}: [A,\bGr^\mC]^G\to \bK_G(A)$ is also a group completion
(by Theorem \ref{thm:BOP_to_KO_G}, or rather its complex analog),
so the unique extension $[-]:\bK_G(A)\to\bku^0_G(A_+)$ is an isomorphism.
\end{proof}
\index{subject}{equivariant $K$-theory|)}

We specialize Theorem \ref{thm:K_G_is_ku} 
to the case $A=\ast$, i.e., when the base is a single point.
In this case the bundle projection is no information,
$G$-vector bundles specialize to $G$-representations,
and the ring $\bK_G(\ast)$ becomes the unitary representation ring $\bRU(G)$.\index{subject}{representation ring!unitary}
On the other hand, $\bku_G^0(S^0)$ specializes to $\pi_0^G(\bku)$.
The ring homomorphism
\[ [-]\ : \ \bRU(G)\ \to \ \pi_0^G(\bku) \]
is then determined by its effect on the classes
of actual representations. If $W$ is a unitary $G$-representation,
we get an explicit representative for the class $[W]$ by choosing
a $G$-equivariant $\mC$-linear isometric embedding $j:W\to V_\mC$
into the complexification of an orthogonal $G$-representations.
For example, the map
\[ 
 j_W \ : \ W \ \to \ ( u W)_\mC\ , \quad 
j_W(w) \ \longmapsto \ 1/\sqrt{2}\cdot ( 1\tensor w - i\tensor(i w)  )
 \]
into the complexification of the underlying orthogonal $G$-representation of $W$
does the job.
Then $[W]$ is the homotopy class of the $G$-map 
\begin{equation}\label{eq:representative [W]}
  S^V\ \to \ 
\bk(\Sym(V_\mC),S^V) \ = \ \bku (V)
\ , \quad v \ \longmapsto \quad [j(W);\, v ] \ .  
\end{equation}
Both $\bRU(G)$ and $\pi_0^G(\bku)$ have restriction maps,
transfers and multiplicative power operations in $G$, 
i.e., they are global power functors
in the sense of Definition \ref{def:power functor}.
The maps $[-]:\bRU(G)\to\pi_0^G(\bku)$ preserve
most of this additional structure, but that is not apparent from
what we discussed so far.

\begin{theorem}\label{thm:R(G) to pi ku}
For every compact Lie group $G$, the map
\begin{equation}\label{eq:R2pi(ku)}
 [-]\ : \ \bRU(G) \ \to \ \pi_0^G(\bku)   
\end{equation}
is a ring homomorphism.
As $G$ varies over all compact Lie groups, the
homomorphisms \eqref{eq:R2pi(ku)}
are compatible with restriction maps, with complex conjugation, 
with finite index transfers, and with multiplicative power operations.
Moreover, the map $[-]$ is an isomorphism whenever the group $G$ is finite.
\end{theorem}
\begin{proof}
The fact that the map is a ring homomorphism, compatible with restrictions,
compatible with complex conjugation
and an isomorphism for finite groups is a special case of 
Theorem \ref{thm:K_G_is_ku} for a one-point $G$-space.

Now we show that the maps $[-]$ are compatible with finite index transfers.
We let $H$ be a finite index subgroup of $G$ and introduce 
an orthogonal $G$-spectrum $\bku[G/H]$ by
\[ \bku[G/H](V) \ = \ \bk(\Sym(V_\mC), S^V\sm (G/H)_+)  \ . \]
The structure maps of $\bku[G/H]$ are defined in much the same way as for $\bku$,
with the extra smash factor $(G/H)_+$ acting as a dummy;
the $G$-action comes entirely from the translation action on $G/H$.
As $V$ varies over all inner product spaces, the assembly maps
\begin{align*}
    \bku(V)\sm (G/H)_+\ = \ &\bk(\Sym(V_\mC),S^V)\sm (G/H)_+ \\ 
&\to \ \bk(\Sym(V_\mC),S^V\sm (G/H)_+) \ = \ \bku[G/H](V) 
\end{align*}
form a morphism of orthogonal $G$-spectra $\alpha:\bku\sm (G/H)_+\to \bku[G/H]$.

We let $\nabla:(G/H)_+\to 1_+$ be the $G$-equivariant fold map
that sends all of $G/H$ to~1. We let $l:G/H_+\to 1_+$ be the
$H$-equivariant `projection' onto the preferred coset, i.e., $l(e H)=1$ 
and $l(g H)=0$ for $g\not\in H$.
Then the following diagram commutes:
\begin{equation}\label{eq:ku[T] diagram} 
 \xymatrix@C=11mm{ 
  \pi_0^H(\bku) &
\pi_0^H(\bku\sm G/H_+)\ar[d]^{\pi_0^H(\alpha)} \ar[l]^-{\pi_0^H(\bku\sm l)}&
\pi_0^G(\bku\sm G/H_+)\ar[d]^{\pi_0^G(\alpha)}
\ar[l]^-{\res_H^G}\ar[r]^-{\pi_0^G(\bku\sm \nabla)}
\ar@<-.4ex>@/_1pc/[ll]_(.7){\Wirth_H^G}&
 \pi_0^G(\bku) \\
&\pi_0^H(\bku[G/H])\ar@/^1pc/[ul]^{\pi_0^H(\bku[l])} &
\pi_0^G(\bku[G/H]) \ar[l]^{\res^G_H}  \ar@/_1pc/[ur]_{\pi_0^G(\bku[\nabla])}& 
}  
\end{equation}
The upper left horizontal composite
is the Wirthm{\"u}ller map \eqref{eq:define_Wirthmuller}.
Since the Wirthm{\"u}ller map is an isomorphism 
(Theorem \ref{thm:Wirth iso}),
the map $\pi_0^H(\bku[l])\circ\res^G_H$ is surjective.
Now we let $f:S^V\to \bku[G/H](V)$ be a $G$-map that represents an element
in the kernel of the map $\pi_0^H(\bku[l])\circ\res^G_H$,
where $V$ is a $G$-representation.
After increasing $V$, if necessary, 
we can assume that $V$ is ample and that the composite
\[  S^V\ \xra{\ f\ }\  \bku[G/H](V)\ \xra{\bku[l](V)}\ \bku(V) \]
is $H$-equivariantly null-homotopic. By adjointness, this means that the composite
\begin{align*}
  S^V\ \xra{\ f\ }\  \bku[G/H](V)\ &= \ \bk(\Sym(V_\mC);S^V\sm (G/H)_+)\\
&\xra{\ l^\flat\ }\ \map^H(G,\bk(\Sym(V_\mC);S^V))\ = \ \map^H(G,\bku(V))  
\end{align*}
is $G$-equivariantly null-homotopic,
where $l^\flat$ is adjoint to the $H$-equivariant map
\[ \bk(\Sym(V_\mC);S^V\sm l)\ : \   \bk(\Sym(V_\mC);S^V\sm (G/H)_+)\ \to \
\bk(\Sym(V_\mC);S^V) \ . \]
The map $l^\flat$ coincides with the Wirthm{\"u}ller map 
\[ \omega_{S^V}\ :\ \bk(\Sym(V_\mC);G\ltimes_H S^V)\ \to\ 
\map^H(G,\bk(\Sym(V_\mC);S^V))\]
defined in \eqref{eq:define_omega_A},
up to the effect of the shearing isomorphism $S^V\sm (G/H)_+\iso G\ltimes_H S^V$.
Since $V$ is ample, the $\bGamma$-$G$-space
$\bk(\Sym(V_\mC))$ is special by Theorem \ref{thm:bk special}~(i), 
and it is $G$-cofibrant by Example \ref{eg:ku cofibrant}.  
So the Wirthm{\"u}ller map $\omega_{S^V}$, and hence also the map $l^\flat$,
is a $G$-weak equivalence by Theorem \ref{prop:prolonged Wirthmuller}~(ii).
Hence already $f$ is $G$-equivariantly null-homotopic. 
Altogether this proves that the map
$\pi_0^H(\bku[l])\circ\res^G_H$ is also injective, hence bijective.
Since $\pi_0^H(\bku[l])\circ\res^G_H$ and the Wirthm{\"u}ller map
are both bijective, the commutativity of \eqref{eq:ku[T] diagram} 
shows that the map $\pi_0^G(\alpha):\pi_0^G(\bku\sm(G/H)_+)\to\pi_0^G(\bku[G/H])$ 
is also bijective.

Now we let $W$ be any unitary $H$-representation
and $j:\map^H(G,W)\to V_\mC$ a $G$-equivariant $\mC$-linear embedding into
the complexification of some $G$-representation.
We define a class
\[ [W]_H^G \ \in \ \pi_0^G(\bku[G/H]) \]
by specifying a representative $G$-map
\[ S^V \ \to \  \bk(\Sym(V_\mC), S^V\sm (G/H)_+)
\ = \ \bku[G/H](V) \]
by
\[ v\ \longmapsto \ [j(\map^H(H g_1^{-1},W)),\dots,j(\map^H(H g_m^{-1},W));\, v\sm g_1 H,\dots,v\sm g_m H] \ .\]
Here $g_1,\dots,g_m$ is a set of coset representatives for $H$ in $G$.
The relation
\begin{align*}
\bk(\Sym(V_\mC), S^V\sm l_+)
[j(\map^H(H g_k,W));\,v\sm g_k H&]_{k=1,\dots,m} \\ 
&=\ [j(\map^H(H,W));\, v]
\end{align*}
shows that
\[ \pi_0^H(\bku[l])(\res^G_H([W]_H^G))\ = \ [W] \ . \]
The commutativity of the left part of diagram \eqref{eq:ku[T] diagram} then shows that
\[ [W]\ = \ \pi_0^H(\bku[l])(\res^G_H([W]_H^G)) \ = \ 
\Wirth_H^G(\pi_0^G(\alpha)^{-1}( [W]_H^G)) \ .\]
Since the Wirthm{\"u}ller map is inverse to the external transfer
(Theorem \ref{thm:Wirth iso}),
this is equivalent to the relation
\[ [W]_H^G \ = \ \pi_0^G(\alpha)(G\ltimes_H [W] ) \]
in the group $\pi_0^G(\bku[G/H])$.

On the other hand, if we fold $S^V\sm (G/H)_+$ onto $S^V$, then
\begin{align*}
\bk(\Sym(V_\mC), S^V\sm \nabla_+)
&[\map^H(H g_1,W),\dots,\map^H(H g_m,W);\, v\sm g_1 H,\dots,v\sm g_m H] \\ 
&=\ [\map^H(H g_1,W),\dots,\map^H(H g_m,W);\, v,\dots,v] \\ 
&=\ [\map^H(H g_1,W)\oplus\dots\oplus\map^H(H g_m,W);\, v] \\ 
&=\ [\map^H(G,W);\, v] \ .
\end{align*}
This shows that
\[ \pi_0^G(\bku[\nabla])([W]_H^G)\ = \ [\map^H(G,W)] \ . \]
The commutativity of the right part of diagram \eqref{eq:ku[T] diagram} 
then yields the desired relation for classes of actual representations:
\begin{align*}
 [\map^H(G,W)]\ = \ \pi_0^G(\bku[\nabla])([W]_H^G) \ 
&= \ \pi_0^G(\bku[\nabla])(\pi_0^G(\alpha)(G\ltimes_H [W] )) \\ 
&= \   \pi_0^G(\bku\sm\nabla)(G\ltimes_H [W]) \ = \ \tr_H^G[W] \ .
\end{align*}
Since transfer maps are additive,
the relation persists to classes of virtual representations.

Now we treat the compatibility with power operations,
i.e., that the following square commutes for every compact Lie group $G$
and all $m\geq 1$:
\[ \xymatrix@C=20mm{ 
\bRU(G)\ar[r]^-{[-]} \ar[d]_{P^m} &\pi_0^G(\bku)\ar[d]^{P^m}\\
\bRU(\Sigma_m\wr G) \ar[r]_-{[-]} &\pi_0^{\Sigma_m\wr G}(\bku)} \]
We consider a unitary $G$-representation $W$ and a $G$-equivariant
$\mC$-linear isometric embedding $j:W\to V_\mC$ into the complexification
of an orthogonal $G$-representation.
The class $[W]$ is then represented by the $G$-map
\[ [j(W);-]\ : \ S^V\ \to\ \bku(V)\ .\]
The map
\begin{align*}
J \ : \ W^{\tensor m} \quad &\to \qquad \Sym(V^m_\mC) \\
w_1\tensor \cdots\tensor w_m &\longmapsto  
(j(w_1),0,\dots,0)\cdot (0,j(w_2),0,\dots,0)\cdot
\ldots\ \cdot (0,\dots,0,j(w_m))
\end{align*}
is a $(\Sigma_m\wr G)$-equivariant $\mC$-linear isometric embedding.
So the class $[W^{\tensor m}]=[P^m\td{W}]$
is represented by the $G$-map
\[ [J(W^{\tensor m});-] \ : \ S^{V^m}\ \to \ \bku( V^m) \ .\]
This map coincides with the composite
\[ S^{V^m}\  \xra{[j(W);-]^{\sm m}} \
\bku(V)^{\sm m} \ \xra{\mu_{V,\dots,V}}\ \bku(V^m) \] 
which represents the power operation $P^m[W]$, so we have shown the relation
\[ P^m[W]\ = \ [W^{\tensor m}]\ = \ [P^m\td{W}] \]
in the group $\pi_0^{\Sigma_m\wr G}( \bku )$.
Since the map $[-]$ is additive, compatible with finite index transfers
and the classes of actual representations generate $\bRU(G)$ as an abelian group, 
the additivity formula for power operations 
\[ P^m(x+y) \ = \ {\sum}_{k=0}^m\ \tr_{k,m-k} (P^k(x)\times P^{m-k}(y)) \]
implies that the map $[-]$ is compatible with
power operations of virtual representations.
\end{proof}

\begin{eg}[Dimension homomorphism]\label{eg:dimension hom}
The {\em dimension homomorphism} is\index{subject}{dimension homomorphism} 
a homomorphism of ultra-commutative ring spectra
\[ \dim\ : \ \bku\ \to \ S\! p^\infty  \]
from the connective global $K$-theory spectrum to the infinite symmetric product spectrum
defined in Example \ref{eg:infinite symmetric product}.\index{subject}{infinite symmetric product spectrum}
The value of this homomorphism at an inner product space $V$ is the map
\begin{align*}
   \dim(V)\ : \ \bku(V)=\bk(\Sym(V_\mC),S^V) \ &\to \quad S\! p^\infty(S^V)  \\
[E_1,\dots, E_n;\, v_1,\dots,v_n]\ &\longmapsto \ {\sum}_{i=1}^n\, \dim(E_i)\cdot v_i\ ,
\end{align*}
i.e., a configuration of vector spaces is mapped to the configurations
of its dimensions (i.e., the sum in the abelian monoid structure of $S\! p^\infty(S^V)$). 
The map is multiplicative because dimension
is multiplicative on tensor products.  

The infinite symmetric product spectrum includes 
by a global equivalence of ultra-commutative ring spectra $S\! p^\infty\to\Hc\mZ$
into the Eilenberg-Mac Lane spectrum
of the integers,\index{subject}{Eilenberg-Mac\,Lane spectrum!of the integers}
see Proposition \ref{prop:Sp2HZ}~(ii).
So we sometimes take $\Hc\mZ$ (instead of $S\! p^\infty$) as the target of
the dimension homomorphism.
\end{eg}

We record that the effect of the dimension homomorphism 
on equivariant homotopy groups `is' the augmentation
\[ \dim \ : \ \bRU(G)\ \to\ \mZ \]
that sends a virtual representation to its virtual dimension.
We warn the reader that the two vertical maps in
the following commutative square are isomorphisms for finite groups,
but {\em not} for general compact Lie groups.

\begin{prop}\label{prop:dim is dim}
For every compact Lie group $G$, 
the following square of commutative rings commutes:
\[ \xymatrix@C=15mm{ 
\bRU(G)\ar[r]^-{\dim} \ar[d]_{[-]} & \mZ \ar[d]^{n\mapsto n\cdot 1} \\
\pi_0^G(\bku)\ar[r]_-{\pi_0^G(\dim)}  & \pi_0^G(S\! p^\infty)} \]
\end{prop}
\begin{proof}
For $n\geq 1$ we consider the `geometric' degree $n$ morphism
\[ \delta^n \ : \ \mS \ \to \ S\! p^\infty \ , \quad x \ \longmapsto \ n\cdot x\ .\]
In level~0, the map
\[ \delta^n(0)\ : \ S^0 \ \to \ S\! p^{\infty}(S^0)\ \iso \ \mN\{0\} \]
sends the non-basepoint $0\in S^0$ to $n\cdot 0$, the $n$-fold
multiple of the corresponding generator of the free abelian monoid $\mN\{0\}$.
Since the ring homomorphism $\mZ\to \pi_0^G(S\! p^\infty)$
sends $n$ to the class of $n\cdot 0$, this shows that
$\delta^n_*(1)= n\cdot 1$ in the group $\pi_0^G(S\! p^\infty)$.

Now we let $W$ be a unitary $G$-representation of dimension $n$.
A representative $[j_W(W);-]$ for the class $[W]$ in $\bRU(G)$
was specified in \eqref{eq:representative [W]}. The composite 
\[ S^{u W}\ \xra{[j_W(W);-]}\  \bku(u W)\ \xra{\dim(u W)} \  S\! p^\infty(S^{u W})\]
is the map 
\[ \delta^n(u W)\ :\ S^{u W}\to S\! p^{\infty}(u W)\ .\]
Hence we conclude that
\[ \pi_0^G(\dim)[W] \ = \ [\delta^n(u W)]\ = \ \delta^n_*(1)\ = \ n\cdot 1\ 
= \ \dim(W)\cdot 1\ .\]
This proves the proposition for the classes of actual representations.
Since all four maps in the square are homomorphisms of abelian groups,
the relation also holds for the classes of virtual representations.
\end{proof}

\begin{rk}\label{rk:homotopy vs smooth transfer}
Segal \cite[\S 2]{segal-representation}
discusses a `smooth transfer' $i_!:\bRU(H)\to \bRU(G)$
where $i:H\to G$ is the inclusion of a closed subgroup of a compact Lie group $G$.
This smooth transfer is {\em not} to be confused with the holomorphic transfer, 
which is defined when $G/H$ is a complex algebraic variety. 
One could hope that the maps $[-]:\bRU(G)\to \pi_0^G ( \bku )$ take the smooth transfer
to the homotopy theoretic transfer 
\[ \tr_H^G \ : \ \pi_0^H ( \bku ) \ \to \ \pi_0^G ( \bku ) \ .\]
This (false!) expectation is suggested by the fact that the smooth
transfer has all the right formal properties, and the complex representations
rings do form a global functor.
As we now illustrate by a specific example,
the collection of ring homomorphisms $[-]:\bRU(G)\to \pi_0^G ( \bku )$ 
comes close to being a morphism of global functors, but it fails to commute with
infinite index transfers.
For {\em periodic} global $K$-theory $\bKU$, 
discussed below in Construction \ref{con:global KU}, 
the two transfers do correspond also for infinite index inclusions;
in fact, $\upi_0 ( \bKU )$ is isomorphic, as a global power functor,
to the representation ring global functor,
see Theorem \ref{thm:pi_0 KU is RU} below.

To illustrate the issue we recall a specific smooth transfer
to $G=S U(2)$ from the normalizer $N=N_{S U(2)}T$ of the diagonal maximal torus
\[ T \ = \ \left\lbrace \left( \begin{smallmatrix} \lambda & 0 \\ 0 & \lambda^{-1} \end{smallmatrix} \right)
\ : \  \lambda\in U(1) \right\rbrace \ .\]\index{subject}{special unitary group!$S U(2)$}
The maximal torus normalizer $N$ is generated by $T$ and the matrix  
$\left(\begin{smallmatrix} \phantom{-}0 & 1\\ -1 & 0\end{smallmatrix}\right)$.
We use Segal's character formula \cite[p.\,119]{segal-representation}
to calculate the character of the smooth transfer $i_!(1)$ of the trivial
1-dimensional $N$-representation, where $i:N\to S U(2)$ is the inclusion. 
An element $g\in S U(2)$ is regular 
in the sense of \cite[Def.\,1.9]{segal-representation}
if and only if it has infinite order, and then the closed subgroup
generated by $g$ is a conjugate $^\gamma T$ of $T$.
The fixed set of such a regular element $g$ on $G/N$ then
consists of a single coset, namely $\gamma N$. So
\[ \chi_{i_!(1)}(g)\ = \  \chi_1(\gamma^{-1}g\gamma)\ = \ 1 \ . \]
Since the character is continuous and regular elements are dense,
the character of $i_!(1)$ is constant with value~1.
Since the character determines the representation up to isomorphism, we have
$i_!(1) = 1$ in $\bRU(S U(2))$.
We claim that in contrast to this relation in $\bRU(S U(2))$,
the elements $\tr_N^{S U(2)}(1)$ and~1 are linearly independent 
(so in particular different) in $\pi_0^{S U(2)} ( \bku )$. 
We can detect this through the 
dimension homomorphism\index{subject}{dimension homomorphism}
$\dim: \bku\to \Hc\mZ$ introduced in Example \ref{eg:dimension hom}.
As we discussed in Example \ref{eg:HZ for abelian G},
the elements $\tr_N^{S U(2)}(1)$ and~1 are linearly independent in $\pi_0^{S U(2)}(H\mZ)$,
so they must also be linearly independent in $\pi_0^{S U(2)} ( \bku )$.

The dimension homomorphism can also be used to detect some
odd-dimen\-sional classes in the coefficient ring $\pi_*^{U(1)}(\bku)$:
we showed in Theorem \ref{thm:pi_1^T HA} that the group $\pi_1^{U(1)}(\Hc\mZ)$
is isomorphic to $\mQ$, and that the dimension shifting transfer
from the trivial group to $U(1)$, applied to the suspension of the multiplicative unit
$1\in\pi_0^e(\Hc\mZ)$, is non-zero in $\pi_1^{U(1)}(\Hc\mZ)$.
So the class
\[ \Tr_e^{U(1)}(1\sm S^1) \ \in \ \pi_1^{U(1)}(\bku)\]
has infinite order, and is in particular non-zero.
\end{rk}

\begin{rk}
There is a  morphism of orthogonal spaces
$c : \bGr \to \Omega^\bullet \bko$
defined in much the same way as its complex analog \eqref{eq:define_c}.
Derived from this are ring homomorphisms
\[  [-]\ : \ \bKO_G(A)\ \to \ \bko_G^0(A_+)  \]
for every compact $G$-space $A$, analogous to the homomorphism \eqref{eq:K_G_to_ku} 
in the complex case.\index{subject}{equivariant $K$-theory}
For $A=\ast$ this specializes to a ring homomorphism
\[ [-]\ : \ \bRO(G)\ \to\ \pi_0^G ( \bko ) \ ,\]
where $\bRO(G)$ is the orthogonal representation ring of $G$.\index{subject}{representation ring!orthogonal}
The analog of Theorem \ref{thm:R(G) to pi ku} then says that these homomorphisms are
compatible with restriction, finite index transfer and multiplicative power operations, 
and an isomorphism for finite groups $G$.
However, in the real situation there is no 
eigenspace decomposition, hence no direct analog of
the delooping $\eig :\bU \to \Omega^\bullet(\sh\bku)$ of the morphism $c$.
So to establish the real analog of Theorem \ref{thm:[A,Gr] to [A,bku]}
(which is used in the proof of Theorem \ref{thm:R(G) to pi ku}),
one has to use a different proof.
\end{rk}

\begin{construction}[Rank filtrations]\index{subject}{rank filtration!of $\bku$}\label{con:rank_filtration}\index{subject}{rank filtration!of $\bko$}
The connective global $K$-theory spectra come with exhaustive filtrations 
\[\bko^{[1]}\to \bko^{[2]}\to \dots \to \bko^{[m]}\to \dots\  \to \ \bko\]
respectively
\[\bku^{[1]}\to \bku^{[2]}\to \dots \to \bku^{[m]}\to \dots\  \to \ \bku\]
by orthogonal subspectra.
We define the rank filtration for $\bku$, the case of $\bko$ being similar.
For $m\geq 1$  we let
\[ \bku^{[m]}(V)\ = \ \bk^{[m]}(\Sym(V_\mC),S^V) \ \subset \ 
\bk(\Sym(V_\mC),S^V)\ = \ \bku(V) \]
be the subspace 
of those configurations $[E_1,\dots, E_n;\,v_1,\dots,v_n]$
such that
\[ {\sum}_{i=1}^n \, \dim(E_i)\ \leq \ m \ .  \]
As $V$ varies, the spaces $\bku^{[m]}(V)$ form an 
orthogonal subspectrum $\bku^{[m]}$ of $\bku$.

The dimension function is multiplicative on tensor products, so
the multiplication of $\bku$ restricts to associative and unital pairings
\[ \bku^{[m]}\sm \bku^{[n]}\ \to \ \bku^{[m n]} \ . \]
For $m=1$ this gives $\bku^{[1]}$ the structure 
of an ultra-commutative ring spectrum 
and it gives $\bku^{[n]}$ a module structure over $\bku^{[1]}$.
The inclusion $\bku^{[1]}\to\bku$ is multiplicative, i.e., a morphism
of ultra-commutative ring spectra.

The first pieces $\bko^{[1]}$ and
$\bku^{[1]}$ are multiplicative models for the suspension spectra
of global classifying space 
of the cyclic group $C_2=O(1)$ respectively the circle group $U(1)$.
Indeed, the complex global projective space $\bP^\mC=\bGr_\tensor^{\mC,[1]}$,
introduced in \eqref{eq:complex_P}, is a multiplicative model of the 
global classifying space of $U(1)$.
A configuration of points labeled by vector spaces 
of total dimension~1 has to be concentrated on at most one point.
So the map
\begin{align*}
  S^V\sm \bP^\mC(V)_+\ = \ S^V\sm P(\Sym(V_\mC))_+\ &\to \ \bk^{[1]}(\Sym(V_\mC),S^V)
\ = \ \bku^{[1]}(V)\\
 v\sm L \qquad &\longmapsto \quad [L; v]  
\end{align*}
is a homeomorphism. 
As $V$ varies through real inner product spaces,
these maps form an isomorphism of ultra-commutative ring spectra
\[ \Sigma^\infty_+ \bP^\mC \ \iso \ \bku^{[1]} \ .\]
\index{subject}{global projective space}\index{symbol}{$\bP$ - {global projective space}}
Similarly, the ultra-commutative ring spectrum $\bko^{[1]}$ 
is globally a suspension spectrum of a global classifying space 
of the group $C_2$.

The underlying non-equivariant spectrum of $\Sigma^\infty_+ B_{\gl} U(1)$,
hence of $\bku^{[1]}$, 
has the stable homotopy type of the suspension spectrum of $\mC P^\infty$; 
moreover, on underlying non-equivariant spectra,  the morphism 
\begin{equation}\label{eq:B_gl T to ku}  
\bku^{[1]}\ \to\ \bku    
\end{equation}
is a rational stable equivalence.
However, the morphism \eqref{eq:B_gl T to ku} 
is {\em not} a rational global equivalence, which can be seen already
on the level of $\upi_0$ for the group $C_2$.
Indeed, by Proposition \ref{prop:B_gl represents}  
the homotopy group global functor $\upi_0( \Sigma^\infty_+ B_{\gl} U(1) )$,
and hence also the global functor $\upi_0( \bku^{[1]})$,
is the represented global functor $\bA( U(1),-)$. The inclusion 
$\bku^{[1]}\to\bku$ induces a morphism of global power 
functors $\bA(U(1),-) \iso \upi_0( \bku^{[1]} )\to \upi_0( \bku )$.
An element of infinite order in the kernel of the map 
$\bA(U(1),C_2)\to \pi_0^{C_2}(\bku)$ is 
\[ \tr_e^{C_2}\circ\res^{U(1)}_e - z^* - \res^{U(1)}_{C_2} \ \in \ \bA(U(1),C_2) \ ,\]
where $z:U(1)\to C_2$ is the trivial homomorphism.
This element maps trivially to $\upi_0(\bku)$ because
the regular representation of $C_2$ splits as the sum of the 1-dimensional
trivial and sign representations.
\end{construction}

For finite groups $G$, 
the ring homomorphism $[-]:\bRU(G)\to\pi_0^G(\bku)$ is an isomorphism.
It will follow from Theorem \ref{thm:pi_0 KU is RU} 
below that the map $[-]:\bRU(G)\to\pi_0^G(\bku)$ is always a split monomorphism,
also when $G$ has positive dimension.
On the other hand, Remark \ref{rk:homotopy vs smooth transfer}
shows that the map is not always surjective, as the class
$\tr_N^{S U(2)}(1)$ in $\pi_0^{S U(2)}(\bku)$ is not in the image.
In \cite{hausmann-ostermayr}, Hausmann and Ostermayr give a complete
calculation of $\upi_0(\bku)$ as a global functor.\index{subject}{special unitary group!$S U(2)$}
The strategy of \cite{hausmann-ostermayr} is to identify the 
global homotopy types of the subquotients of the rank filtration 
and deduce from it a presentation $\upi_0(\bku)$ by generators and relations.
The final answer is that $\upi_0(\bku)$ is generated as a global
functor by the classes
\[ x_m \ = \ [\tau_m]\ \in \ \pi_0^{U(m)}(\bku) \ ,\]
where $\tau_m$ denotes the tautological $U(m)$-representation on $\mC^m$.
There are two kinds of relations; on the one hand, the relations
\[ \res^{U(m+n)}_{U(m)\times U(n)}(x_{m+n}) \ = \ p^*(x_m) + q^*(x_n)
\quad \in \pi_0^{U(m)\times U(n)}(\bku) \]
for all $m,n\geq 1$, where $p:U(m)\times U(n)\to U(m)$
and $q:U(m)\times U(n)\to U(n)$ are the two projections.
These relations follow from the fact that the maps $[-]$
are additive and compatible with restrictions.
The other set of relations equates finite index transfers
of representations with the corresponding finite index transfer in
homotopy theory.

\begin{construction}[Bott class]\label{con:Bott class}\index{subject}{Bott class!in connective $K$-theory}
The {\em Bott class}
\[ \beta\ \in \  \pi_2^e(\bku) \]
is an important non-equivariant homotopy class of the spectrum $\bku$
that we recall now. As we shall show in Theorem \ref{thm:Bott and inverse}
below, the Bott class becomes invertible in the homotopy ring
of the periodic $K$-theory spectrum $\bKU$.\index{symbol}{$\beta$ - {Bott class in $\pi_2^e(\bku)$}}

We define the Bott class by specifying an explicit representative.
We define a continuous map $m:S^{\mR\oplus\mC}\to S U(2)$ by the formula
\begin{equation}  \label{eq:define_m}
 m(v)\ = \ m(x,z)\ = \  
\frac{1}{|v|^2+1} \begin{pmatrix} |v|^2-1 - i 2 x & 
2 i \bar z \\ 2 i z &  |v|^2-1 + i 2 x 
 \end{pmatrix}\ .
\end{equation}
The eigenspace decomposition 
\[ \eig \ : \  U(2) \ \iso \ \bk(\mC^2;S^1) \]
is the inverse of the homeomorphism \eqref{eq:C(S^1)_is_U};\index{subject}{eigenspace decomposition}
it sends a unitary matrix to the configuration of its eigenvalues,
with the inverse Cayley transform applied, 
labeled by the corresponding eigenspaces.
The Bott class is then represented by the composite
\begin{align}  \label{eq:Bott_representative}
 S^3 \ = \ S^{\mR\oplus\mC} \ \xra{\ m\ } \  U(2)\  &\xra[\iso]{\ \eig\ } \ 
\bk(\mC^2;S^1)\\ 
&\xra{\text{incl}}\ \bk(\Sym(\mC);S^1)\ =\ \bku(\mR) \ .  \nonumber
\end{align}
In the last step we identify $\mC^2$ with the subspace of $\Sym(\mC)$ spanned
by the constant and linear summand in the symmetric algebra,
in terms of the preferred basis $1\in \Sym^0(\mC)$ and $1\in\Sym^1(\mC)$.
\end{construction}

\begin{prop}
  The  group $\pi_2^e(\bku)$ is infinite cyclic, 
  and the Bott class $\beta$ is a generator.
\end{prop}
\begin{proof}
The map $m:S^{\mR\oplus\mC}\to S U(2)$ is bijective;
indeed, an explicit formula for the inverse is 
\[ S U(2)\ \to \ S^{\mR\oplus \mC} , \quad
 \begin{pmatrix} 
a & \bar b \\ - b & \bar a
 \end{pmatrix}  \ \longmapsto \ 
\frac{\text{Im}(a)}{\text{Re}(a)-1}\sm \frac{i b}{1-\text{Re}(a)}\ ,\]
where $a,b\in \mC^2$ are complex numbers satisfying $|a|^2+|b|^2=1$.
As a continuous bijection between compact spaces, $m$ is thus a
homeomorphism from $S^3$ to $S U(2)$.
Since the inclusion $S U(2)\to U(2)$ induces an isomorphism on $\pi_3$,
the map $m$ represents a generator of the 
infinite cyclic group $\pi_3(U(2),1)$.
The standard embedding $U(2)\to U$ into the infinite unitary group is 4-connected, 
so the composite \eqref{eq:Bott_representative}
represents a generator of the infinite cyclic group $\pi_3(\bku(\mR),\ast)$.
Theorem \ref{thm:weak global Omega for bku} says that $\bku$ 
is a positive $\Omega$-spectrum (in the non-equivariant sense),
so in particular the stabilization map
\[ \pi_3(\bku(\mR),\ast) \ = \ [S^3,\bku(\mR)]\ \to \ 
\colim_{n\geq 0} \, [S^{n+2},\bku(\mR^n)]\ = \ \pi_2^e(\bku)\]
is an isomorphism. So the group $\pi_2^e(\bku)$ is infinite cyclic, 
and the Bott class is a generator.
\end{proof}

\begin{construction}[Equivariant Bott classes]\label{con:equivariant Bott class}
There are more general {\em equivariant Bott classes}\index{subject}{Bott class!of a $\Spin^c$-representation}\index{symbol}{$\beta_{G,W}$ - {equivariant Bott class of a $\Spin^c$-representation}}
\[ \beta_{G,W} \ \in \ \bku^0_G(S^W) \]
defined for $G$-$\Spin^c$-representations, i.e., 
hermitian inner product spaces $W$
equipped with a lift $G\to \Spin^c(u W)$ of the representation homomorphism $G\to O(u W)$
through the adjoint representation
$\ad(u W) :\Spin^c(u W)\to SO(u W)$, see \eqref{eq:adjoint_rep_spin^c}.
These equivariant Bott classes become
invertible in the $R O(G)$-graded homotopy ring of $\bKU$,
compare Remark \ref{rk:equivariant inverse Bott} below;
this equivariant generalization of Bott periodicity goes back to
Atiyah \cite[Thm.\,4.3]{atiyah-Bott and elliptic}.\index{subject}{Bott periodicity!equivariant} 

We recall the definition as given, 
for example, in \cite[III p.\,44]{karoubi-sur la K-theorie}.
The construction depends heavily on a canonical isomorphism,
for every hermitian inner product space $W$, 
between the complexified Clifford algebra 
$\mC\tensor_\mR\Cl(u W)$ of the underlying euclidean vector space of $W$, 
and $\End_\mC(\Lambda^*(W))$, the endomorphism algebra of the exterior algebra of $W$. 
The underlying $\mR$-vector space $u W$ of the given hermitian inner product space
has a euclidean inner product defined by
\[ \td{w,w'} \ = \ \text{Re} (w,w')\ ,\]
the real part of the complex inner product.
In what follows, we use the hermitian inner product 
on the exterior algebra $\Lambda^*(W)$ characterized by the formula
\[ ( v_1\sm\dots\sm v_n,\,w_1\sm\dots\sm w_n )\ = \   \det( ( v_i,w_j )_{i,j})  \]
for all $v_i,w_j\in W$. Another way to say this is that 
if $(e_i)_{i=1,\dots, k}$ is an orthonormal basis of $W$, then the vectors
\[ e_{i_1}\sm\dots\sm e_{i_n} \]
form an orthonormal basis of $\Lambda^*( W)$ 
as the indices run through all tuples with
$1\leq i_1<\dots< i_n\leq k$.

For $w\in W$ we let
\[ d_w\ :\ \Lambda^*(W)\ \to\ \Lambda^*(W) \ , \quad d_w(x)\ = \ w\sm x  \]
denote left exterior multiplication by $w$. 
We let $d_w^\ast:\Lambda^*(W)\to\Lambda^*(W)$ denote the adjoint of $d_w$,
i.e., the $\mC$-linear map characterized by
\[ (x, d^\ast_w(y))\ = \ ( d_w(x) ,y ) \ = \ ( w\sm x ,y ) \]
for all $x,y\in\Lambda^*(W)$.

The exterior algebra $\Lambda^*(W)$ is $\mZ/2$-graded by
even respectively odd exterior powers. The endomorphism
algebra $\End_\mC(\Lambda^*(W))$ then inherits a $\mZ/2$-grading by
even respectively odd operators.
Exterior multiplication by $w$ takes $\Lambda^n(W)$ to $\Lambda^{n+1}(W)$,
so the operator $d_w$ is odd for every $w\in W$. 
Hence the adjoint $d^*_w$ is odd, and so is the endomorphism
\begin{equation}\label{eq:define_delta_w}
 \delta_w \ = \ d_w - d_w^\ast \ .  
\end{equation}
\end{construction}

\begin{lemma}\label{lemma:delta_w squared}
Let $W$ be a hermitian inner product space.
For every $w\in W$, the endomorphisms $d_w$ and $\delta_w$ 
of $\Lambda^*(W)$ satisfy the relations
\[ d_w\circ d_w^*\ +\ d_w^*\circ d_w \ = \ |w|^2\cdot\Id\text{\qquad and\qquad} 
  \delta_w\circ\delta_w\ = \ -|w|^2\cdot\Id\ .  \]
\end{lemma}
\begin{proof}
If we multiply $w$ by a real scalar $\lambda\in\mR$, then both sides
of both equations scale by $\lambda^2$. So it suffices to show the claims
for unit vectors, i.e., we may assume that $|w|=1$.

If $w$ is a unit vector, then $\Lambda^*(W)$ decomposes as 
an orthogonal direct sum $\Lambda^*(W^\perp)\oplus (w\sm \Lambda^*(W^\perp))$,
where $W^\perp$ is the orthogonal complement of $w$ in $W$.
Moreover, $d_w$ is an isometry from the summand $\Lambda^*(W^\perp)$ 
onto the summand $w\sm \Lambda^*(W^\perp)$, and it vanishes on $w\sm \Lambda^*(W^\perp)$.
So the adjoint $d_w^*$ is inverse to $d_w$ on the summand
 $w\sm \Lambda^*(W^\perp)$, and it vanishes on $\Lambda^*(W^\perp)$.
Hence $d_w^*\circ d_w$ is the orthogonal projection onto the summand $\Lambda^*(W^\perp)$, 
and $d_w\circ d_w^*$ is the orthogonal projection onto the other summand 
$w\sm \Lambda^*(W^\perp)$. This shows that 
$d_w\circ d_w^* + d_w^*\circ d_w$ is the identity.

We have $d_w\circ d_w=0$, hence also $d_w^*\circ d_w^*=0$.
This gives 
\[ \delta_w\circ\delta_w \ = \ (d_w -d_w^*)\circ(d_w -d_w^*)\ = \ 
- ( d_w\circ  d_w^* + d_w^* \circ d_w )\ = \ - \Id\ . \qedhere \]
\end{proof}

Because of the previous lemma, the $\mR$-linear map
\[ u W \ \to \ \End_\mC(\Lambda^*(W)) \ , \quad w \ \longmapsto \ \delta_w \]
satisfies $ \delta_w^2 =  -|w|^2\cdot \Id$;
so the universal property of the Clifford algebra provides a unique
homomorphism of $\mZ/2$-graded $\mC$-algebras 
\begin{equation}\label{eq:Clifford_iso}
\delta_W\ : \  \mC\tensor_\mR\Cl(u W) \ \to\ \End_\mC(\Lambda^*(W)) 
\end{equation}
sending $w\in W$ to $\delta_w$. 
The homomorphism is also compatible with passage to adjoints;
indeed, for $w\in W$ we have
\[ \delta_w^* \ = \ (d_w-d_w^\ast)^\ast\ = \ d^\ast_w-d_w\ = \ 
d_{-w} -d_{-w}^\ast\ = \ \delta_{-w}\ = \ \delta_{w^\ast}\ . \]
Since $\mC\tensor_\mR\Cl(u W)$ is generated by the elements of $u W$, the compatibility
with adjoints holds in general.
The homomorphism $\delta_W$ thus makes the exterior algebra $\Lambda^*(W)$
into a $\mZ/2$-graded $\mC\tensor_\mR\Cl(u W)$-module.

\begin{eg}[The complexified Clifford algebra of $u\mC$]\label{eg:delta C}
We make the homomorphism $\delta_W$ explicit in the simplest non-trivial case,
namely for $W=\mC$. In this case the elements $1, i\in\mC$ form
an orthonormal $\mR$-basis of $u\mC$; we let
\[ e \ = \ 1\tensor [1]\text{\qquad respectively\qquad} f \ = \ 1\tensor [i] \]
be their images in the complexified Clifford algebra $\mC\tensor_\mR\Cl(u\mC)$.
This Clifford algebra then has a $\mC$-basis $(1, e, f, e f)$,
with multiplicative relations
\[ e^2 \ = \ f^2 \ = \ -1 \text{\qquad and\qquad} f e \ = \ - e f\ . \]
On the other hand, the exterior algebra $\Lambda^\ast(\mC)$ 
is two dimensional, and an orthonormal $\mC$-basis is given by
the classes $1\in \Lambda^0(\mC)$ (the multiplicative unit)
and $x\in\Lambda^1(\mC)$, the class of $1\in\mC$.
A direct calculation shows that in terms of the basis $(1,x)$
of $\Lambda^*(\mC)$, the map $\delta_\mC$ is given by
\begin{align}\label{eq:delta of e,f}
 \delta_e \ =\ \begin{pmatrix} 0 & -1 \\ 1 & 0\end{pmatrix} 
\ , \quad\
 \delta_f \ = \ \begin{pmatrix} 0 & i \\ i &  0 \end{pmatrix}
\text{\quad and\quad}
 \delta_e\circ \delta_f \ = \ \begin{pmatrix} - i & 0 \\ 0 &  i \end{pmatrix}\ .  
\end{align}
In particular, the homomorphism $\delta_\mC$ sends the basis $(1,e, f, e f)$
of $\mC\tensor_\mR\Cl(u\mC)$ to a basis of $\End_\mC(\Lambda^*(\mC))$,
so $\delta_\mC$ is an isomorphism.
\end{eg}

We recall that the homomorphism $\delta_W$ is in fact an isomorphism
in general. 
Indeed, given two hermitian inner product spaces $V$ and $W$, the following
square commutes: 
\[ \xymatrix@C=12mm@R=7mm{ 
 \mC\tensor_\mR\Cl(u V\oplus u W) \ar[r]^-{\delta_{V\oplus W}}
\ar[dd]_\iso &
\End_\mC(\Lambda^*(V\oplus W)) \ar[d]_\iso\\
& \End_\mC(\Lambda^*(V)\tensor \Lambda^*(W)) \\
 (\mC\tensor_\mR\Cl(u V))\tensor_\mC (\mC\tensor_\mR\Cl(u W)) 
\ar[r]_-{\delta_V\tensor\delta_W} & 
 \End_\mC(\Lambda^*(V))\tensor\End_\mC(\Lambda^*(W)) \ar[u]_\tensor^\iso
} \]
The left vertical isomorphism sends $1\tensor [v,w]$
to $(1\tensor v)\tensor (1\tensor 1)+(1\tensor 1)\tensor (1\tensor w)$.
The upper right vertical map is induced by the isomorphism 
of $\mZ/2$-graded $\mC$-algebras
\[ \Lambda^*(V\oplus W)\ \xra{\ \iso\ } \  \Lambda^*(V)\tensor \Lambda^*(W)  \]
that extends the linear map
\[  V\oplus W\ \to \  \Lambda^*(V)\tensor \Lambda^*(W)  \ , \quad
(v,w)\ \longmapsto \ v\tensor 1\ + \ 1\tensor w\ .\]
The vertical maps in the above diagram are isomorphisms, and every 
hermitian inner product space is isometrically isomorphic
to $\mC^n$ with standard inner product; so this reduces the
claim to the special case $W=\mC$, in which case the morphism $\delta_\mC$
is an isomorphism by Example \ref{eg:delta C}.

\medskip

To construct the equivariant Bott class
we now start with a hermitian inner product space $W$ and a continuous homomorphism
$G\to \Spin^c(u W)$.
The group $G$ then acts on the underlying euclidean inner product space of $W$
via the adjoint representation $\ad(u W):\Spin^c(u W)\to S O(u W)$.
The equivariant Bott class most naturally lives in the
relative $K$-group $\bK_G(D(W), S(W))$;
elements in this group are represented by triples
$(\xi,\eta,\alpha)$, where $\xi$ and $\eta$ are $G$-vector bundles
over $D(W)$, and $\alpha:\xi|_{S(W)}\iso\eta|_{S(W)}$ is an 
equivariant isomorphism between the restrictions of the two bundles
to the unit sphere $S(W)$.
In this description, the Bott class $\beta_{G,W}$ is represented by the triple
\[ ( D(W)\times \Lambda^{\ev}(W), D(W)\times\Lambda^{\odd}(W),\alpha) \]
consisting of the trivial vector bundles over $D(W)$
with fibers $\Lambda^{\ev}(W)$ respectively $\Lambda^{\odd}(W)$,
and the  equivariant bundle isomorphism $\alpha$ is given 
by the Clifford action \eqref{eq:Clifford_iso}:
\[ \alpha \ : \ S(W)\times \Lambda^{\ev}(W) \ \to \  S(W)\times\Lambda^{\odd}(W) \ , 
\quad (w, x)\ \longmapsto \ (w, \delta_w(x))\ . \]
Here the group $G$ acts on $\Lambda^{\ev}(W)$ 
and $\Lambda^{\odd}(W)$ via the given $\Spin^c$-structure, i.e., 
through the composite
\begin{equation}  \label{eq:G_acts_via_spin^c}  
 G \ \to\ \Spin^c(u W)\ \subset \ \mC\tensor_\mR\Cl(u W)
\ \xra[\iso]{\ \delta_W\ } \ \End_\mC(\Lambda^*(W))\ .
\end{equation}
Equivariant Bott periodicity \cite[Thm.\,4.3]{atiyah-Bott and elliptic}\index{subject}{Bott periodicity!equivariant} 
says that for every compact $G$-space $A$, exterior product
with the Bott class is an isomorphism
\[ -\times \beta_{G,W} \  : \ \bK_G(A)\ \to \ \bK_G(A\times D(W),A\times S(W))\ . \]
To represent the Bott class in $\bku_G^0(S^W)$
we interpret the Clifford action as a clutching function for a vector bundle.
We define a $G$-vector bundle $\xi(W)$ over $S^W$
by gluing the trivial bundle with fiber $\Lambda^{\ev}(W)$ over $D(W)$
with the trivial bundle with fiber $\Lambda^{\odd}(W)$ over $S^W- \mathring D(W)$
using the map
\[ S(W) \ \to \ \bL^\mC(\Lambda^{\ev}(W),\Lambda^{\odd}(W))\ , \quad 
w\ \longmapsto \ \delta_w\ .\]
The homomorphism
\[ [-] \ : \ \bK_G(S^W)\ \to \ \bku_G^0(S^W_+)\]
defined in \eqref{eq:K_G_to_ku} turns this $G$-vector bundle
into an unreduced equivariant $\bku$-cohomology class.
By construction, the fiber of $\xi(W)$ over the basepoint at infinity
is the $G$-representation $\Lambda^{\odd}(W)$,
with $G$ acting via the composite \eqref{eq:G_acts_via_spin^c}.  
So by subtracting the class of the trivial vector bundle with fiber 
$\Lambda^{\odd}(W)$ we obtain the equivariant Bott class as a reduced 
equivariant $\bku$-cohomology class, 
\[  \beta_{G,W}\ = \ [ \xi(W) ] - [ S^W\times \Lambda^{\odd}(W)]  \ \in \ 
\bku_G^0(S^W) \ . \]
We invite the reader to perform a reality check and establish the relation
\[  \beta_{e,\mC}\ = \ - \beta \]
in $\pi_2^e(\bku)$. In other words,
the equivariant Bott class of the trivial group acting on $\mC$
specializes to the non-equivariant Bott class 
as defined in Construction \ref{con:Bott class}, up to a sign.

\index{subject}{K-theory@$K$-theory!connective global|)}

\section{Periodic global \texorpdfstring{$K$}{K}-theory}\label{sec:global K}

Our main object of study in this section is {\em periodic global $K$-theory} $\bKU$, 
an ultra-commutative ring spectrum whose $G$-homotopy type realizes 
$G$-equi\-variant periodic $K$-theory, see Construction \ref{con:global KU}.
The model we use is due to M.\,Joachim \cite{joachim-coherences},
and made of spaces of homomorphisms of $\mZ/2$-graded $C^\ast$-algebras.
Corollary \ref{cor:global_K_represents} shows that the 
equivariant cohomology theory represented by $\bKU$ on finite $G$-CW-complexes `is'
equivariant $K$-theory;
Theorem \ref{thm:pi_0 KU is RU} shows that $\upi_0(\bKU)$ is isomorphic,
as a global power functor, to the complex representation ring functor $\bRU$.
Periodic global $K$-theory receives a morphism of ultra-commutative ring spectra
$j:\bku\to\bKU$ from connective global $K$-theory defined in the previous section.
Theorem \ref{thm:Bott and inverse} shows that $\bKU$ is Bott periodic,
i.e., the Bott class in $\pi_2^e(\bku)$ becomes invertible in $\pi_2^e(\bKU)$.

{\em Global connective $K$-theory} $\bkuc$
is a certain homotopy pullback of the periodic theory $\bKU$, 
its associated global Borel theory,
and the global Borel theory of connective $K$-theory, 
see Construction \ref{con:global connective ku}.
The global homotopy type $\bku^c$ is a refinement 
of Greenlees `equivariant connective $K$-theory' \cite{greenlees-connective}.
One should  note the different order of the adjectives `global' and `connective',
indicating that $\bku$ and $\bku^c$ are quite different global homotopy types
(with the same underlying non-equivariant homotopy type).

\medskip

The spectrum $\bKU$ consists of spaces of homomorphisms 
of $\mZ/2$-graded $C^\ast$-algebras.
Consequently, in this section we will use basic results 
from the theory of $C^\ast$-algebras; the textbooks \cite{murphy-Cstar, wegge-olsen} 
can serve as general references.
Since we are not assuming the reader to be fluent with $C^\ast$-algebras, 
we recall a certain amount of material in some detail.

\begin{construction}[The $C^\ast$-algebra $s$]
  The construction of the orthogonal spectrum $\bKU$ is based
  on a certain $C^\ast$-algebra and on spaces of $\ast$-homo\-morphisms out of it. 
  We consider {\em graded} $C^\ast$-algebras, i.e., $C^\ast$-algebras $A$
  equip\-ped with a $\ast$-automorphism $\alpha:A\to A$ such that $\alpha^2=\Id$.
  We can then decompose $A$ into the $\pm 1$ eigenspaces of $\alpha$
  and obtain a $\mZ/2$-grading of the underlying $\mC$-algebra by setting
  \[ A_{\ev} \ = \ \{ a\in A \ | \ \alpha(a)=a\}\text{\quad and \quad}
  A_{\odd} \ = \ \{ a\in A \ | \ \alpha(a)=-a\}\ .\]
  The conjugation of $A$ preserves the grading into even and odd parts.

  We let $s$ denote the $C^\ast$-algebra of complex valued continuous functions on 
  $\mR$ vanishing at infinity; this is a $\mZ/2$-graded $C^\ast$-algebra with involution 
  $\alpha:s\to s$ defined by
  \[ \alpha(f)(t)\ = \ f(-t)\ . \]
  With respect to this grading, `even' and `odd' have their usual meaning:
  a function $f\in s$ is even (respectively odd)
  if and only if $f(-t)=f(t)$ (respectively $f(-t)=-f(t)$) for all $t\in\mR$.
  
  Graded $\ast$-morphisms out of $s$ correspond to special elements
  in the target $C^\ast$-algebra. The continuous function
  \[ r(x) \ = \ \frac{2i}{x-i} \ : \ \mR \ \to \ \mC\]
  is an element of $s$ that satisfies
  \[    r r^\ast \ + \ r \ + \ r^\ast \ = \ 0 \text{\qquad and\qquad}
    \alpha(r)\ = \ r^\ast \ .    
  \]
  The element $r$ is not homogeneous; its even respectively odd components
  are given by
  \begin{equation} \label{eq:even-odd_of_r}
 r_+(x) \ = \ \frac{-2}{x^2+1}\text{\quad respectively\quad}
  r_-(x) \ = \  \frac{2 x i}{x^2+1}\ .    
  \end{equation}
  Moreover, $s$ is the universal $\mZ/2$-graded $C^\ast$-algebra generated by
  such an element, i.e., evaluation at $r$ is a bijection
  \begin{equation} \label{eq:C(s)_versus_U'}
     C^\ast_{\gr}(s, A) \ \iso \ 
  \{ x\in  A \ : \ x x^\ast = x^\ast x = - x - x^\ast\ ,\ \alpha(x)=x^\ast \}   
  \end{equation}
  for every graded $C^\ast$-algebra $A$.
  Indeed, after adjoining a unit to the algebra $s$ and identifying $S^1$
  with $U(1)$ via the Cayley transform\index{subject}{Cayley transform}
  $c:S^1\iso U(1)$ given by $c(x)=(x+i)(x-i)^{-1}$,
  the algebra $\mC\oplus s$ becomes isomorphic to the unital $C^\ast$-algebra $C(U(1))$
  of continuous $\mC$-valued functions on $U(1)$;
  this isomorphism takes $1+r$ to the inclusion $z:U(1)\to \mC$,
  which is a unitary element in $C(U(1))$.
  So the universal property of $s$
  follows from the fact that $C(U(1))$ is freely generated,
  as an ungraded, unital $C^\ast$-algebra, by the unitary element $z$;
  indeed, unitary elements are in particular normal and have their spectrum
  contained in $U(1)$, so the bijectivity of \eqref{eq:C(s)_versus_U'} 
  becomes a special case of functional calculus for normal elements in unital
  $C^\ast$-algebras, see for example \cite[Thm.\,2.1.13]{murphy-Cstar}.

  In what follows, $\hat\tensor$ denotes the {\em spatial tensor product}\index{subject}{spatial tensor product!of $C^\ast$-algebras}
  of graded $C^\ast$-algebras, see \cite[Sec.\,6.3]{murphy-Cstar}, 
  \cite[Def.\,T.5.16]{wegge-olsen} or \cite[Def.\,1.10]{higson-guentner}.
  The spatial tensor product is the completion of the algebraic tensor product
  with respect to the `spatial norm'; it is also called the {\em minimal tensor product}
  because the spatial norm is minimal among all $C^\ast$-norms on
  the algebraic tensor product, compare \cite[Thm.\,6.4.18]{murphy-Cstar}. 
  A comprehensive discussion of the spatial and other tensor products of
  $C^\ast$-algebras can be found in \cite[App.\,T]{wegge-olsen}.
   We need the {\em graded} tensor product to $\mZ/2$-graded algebras, which involves
  a sign in the formula for multiplication.
  If~ $a,a'\in A$ and $b,b'\in B$  are homogeneous elements of two
  $\mZ/2$-graded algebras, then the multiplication in $A\tensor B$ is defined by
  \[ (a\tensor b)\cdot(a'\tensor b')\ = \ (-1)^{|b||a'|}\cdot (a a')\tensor (b b')\ .  \]
  In other words, if $a\in A$ and $b\in B$ are odd, then
  $a\tensor 1$ and $1\tensor b$ {\em anti-commute} in the graded tensor product.

  \medskip

  \Danger The underlying algebra of the graded tensor of two 
  $\mZ/2$-graded algebras is {\em not} the tensor product of the underlying
  ungraded algebras, as soon as both factors have non-zero odd elements.
  For example, while the ungraded $C^\ast$-algebra underlying $s$ is commutative,
  the underlying $C^\ast$-algebra of the  graded tensor product $s\hat\tensor s$ 
  is {\em not} commutative anymore.  In particular, $s\hat\tensor s$ 
  does not embed into the algebra of continuous functions on any space.

  \medskip

  The algebra $s$ has another important piece of extra structure,
  namely a comultiplication $\Delta:s\to s\hat\tensor s$
  in the category of $\mZ/2$-graded $C^*$-algebras.
  The algebra $s$ is also generated by the self-adjoint functions 
  \begin{equation}\label{eq:s generators u and v}
    u_+(t)\ =\ e^{-t^2} \text{\qquad and\qquad}
    u_-(t)\ =\ t\cdot e^{-t^2} \ ,    
  \end{equation}
  and the comultiplication is completely determined by the values on these;
  indeed, there is a unique graded $\ast$-homomorphism 
  $\Delta:s\to s\hat\tensor s$ that satisfies
  \[ \Delta(u_+)\ = \ u_+\tensor u_+ \text{\quad respectively\quad}
  \Delta(u_-)\ =\ u_-\tensor u_+\ +\ u_+\tensor u_- \ ,\]
  see for example Lemma~1.3.1 and~Remark~1.3.3 of \cite{higson-guentner}.
  This characterization readily implies
  the cocommutativity and coassociativity of $\Delta$.
  The cocommutativity of $\Delta$ is with respect to the
  {\em graded} symmetry automorphism of $s\hat\tensor s$,
  which involves a sign whenever two odd elements are interchanged.
  An explicit definition of the diagonal morphism $\Delta$
  can be found in \cite[Sec.\,1, (13)]{haag}.
\end{construction}

\begin{construction}[Complex Clifford algebras]
  We let $V$ be a euclidean inner product space. We define the complex
  Clifford algebra $\mCl(V)$ by\index{subject}{Clifford algebra|(}\index{symbol}{$\Cl(V)$@$\mCl(V)$ - {complexified Clifford algebra}}
  \[ \mCl(V) \ = \ (T V)_\mC / ( v\tensor v - |v|^2\cdot 1) \ , \]
  the quotient of the complexified tensor algebra of $V$
  by the ideal generated by the elements $v\tensor v -  |v|^2\cdot 1$
  for all $v\in V$.  We write $[-]:V\to \mCl(V)$ for the
  $\mR$-linear and injective composite
  \[ V \ \xra{\text{linear summand}} \  T V\ \xra{1\tensor- }\ (T V)_\mC
  \ \to \  \mCl(V)  \ .\]
  With this notation the relation $[v]^2=|v|^2\cdot 1$ holds in $\mCl(V)$
  for all $v\in V$.
  The Clifford algebra construction is functorial for $\mR$-linear isometric embeddings,
  so in particular $\mCl(V)$ inherits an action of the orthogonal group $O(V)$. 
  The Clifford algebra is $\mZ/2$-graded, coming from the grading of
  the tensor algebra by even and odd tensor powers.

  The complex Clifford algebra is in fact a $\mZ/2$-graded $O(V)$-$C^\ast$-algebra.
  The $\ast$-involution on $\mCl(V)$ is defined by declaring $[v]^\ast=[v]$
  for all $v\in V$ and extending this to a $\mC$-semilinear anti-automorphism. 
  This makes the elements $[v]$ for $v\in S(V)$ into unitary elements of $\mCl(V)$.
  The norm on $\mCl(V)$ arises from an embedding
  into the endomorphism algebra of the exterior algebra $\Lambda^*( V_\mC)$.
  Indeed, we consider the $\mR$-linear map
  \[ V \ \to \ \End_\mC(\Lambda^*(V_\mC))\ , \quad
  v \ \longmapsto \ i\cdot\delta_{1\tensor v} \ , \]
  where $\delta_{1\tensor v}$ was  defined in \eqref{eq:define_delta_w}.
  Lemma \ref{lemma:delta_w squared} provides the relation
  \begin{align*}
    (i\cdot\delta_{1\tensor v})\circ (i\cdot\delta_{1\tensor v}) \ = \ 
    - \delta_{1\tensor v}^2\ = \ |v|^2\cdot \Id\ .
  \end{align*}
  The universal property of $\mCl(V)$ provides a morphism of $\mZ/2$-graded $\mC$-algebras
  \[ \mCl(V)\ \to\  \End_\mC(\Lambda^*(V_\mC))\]
  that sends $[v]$ to $i\cdot\delta_{1\tensor v}$.
  This homomorphism is an embedding and compatible with passage to adjoints.
  The target is a $C^\ast$-algebra via
  the operator norm, so we can endow the Clifford algebra with a $C^\ast$-norm
  via this embedding. For this norm, the embedding 
  $[-]:V\to\mCl(V)$ is isometric. Indeed, 
  \[ \delta_{1\tensor v}^*\circ\delta_{1\tensor v} \ = \ 
  (d^*_{1\tensor v} -d_{1\tensor v})\circ(d_{1\tensor v} - d_{1\tensor v}^*)\ = \ 
   d_{1\tensor v}\circ  d_{1\tensor v}^* + d_{1\tensor v}^* \circ d_{1\tensor v} 
  \ = \ |v|^2\cdot \Id\ , \]
  by Lemma \ref{lemma:delta_w squared}, and so
  \[ \|i\cdot \delta_{1\tensor v}\|^2 
  \ = \ \|\delta_{1\tensor v}^*\circ \delta_{1\tensor v}\| \ = \ |v|^2 \ .\]
\end{construction}

  The Clifford algebra functor sends orthogonal direct sum to 
  graded tensor product, in the following sense.
  The $\mR$-linear map
  \[ V\oplus W \ \to \ \mCl(V)\tensor\mCl(W) \ , \quad 
  (v,w) \ \longmapsto \ [v] \tensor 1 + 1\tensor [w]  \]
  satisfies
  \[ ([v] \tensor 1 + 1\tensor [w])^2 \ = \ [v]^2\tensor 1 + 1\tensor [w]^2 \ = \ 
  (|v|^2+|w|^2)\cdot 1\tensor 1  \]
  because the elements $[v]\tensor 1$ and $1\tensor [w]$ anti-commute
  in $\mCl(V)\tensor\mCl(W)$. The universal property then provides
  a unique extension to a morphism of graded $\mC$-algebras
  \begin{equation}\label{eq:Clifford tensor iso}
    \mu_{V,W}\ : \ \mCl(V\oplus W)\ \iso\ \    \mCl(V)\tensor \mCl(W) 
  \end{equation}
  characterized by $\mu_{V,W}[v,w] =  [v] \tensor 1 + 1\tensor [w]$
  for all $(v,w)\in V\oplus W$. Moreover, this morphism is an isomorphism,
  and the square 
  \[ \xymatrix@C=15mm@R=6mm{ 
    \mCl(V\oplus W)\ar[d]_{\mCl(\tau_{V,W})}\ar[r]^-{\mu_{V,W}} & 
    \mCl(V)\tensor \mCl(W) \ar[d]^\iso  \\
   \mCl(W\oplus V)\ar[r]_-{\mu_{W,V}}  &  \mCl(W)\tensor \mCl(V) } \]
  commutes, where the right vertical map is the symmetry isomorphism
  for $\mZ/2$-graded $\mC$-algebras (which involves a sign whenever two
  odd degree elements interchange places).

\medskip

  \Danger  The complex Clifford algebra $\mCl(V)$
  is isomorphic to the complexification of the real Clifford algebra $\Cl(V)$ considered 
  in the definition of the ultra-commutative monoids $\bPin$ and $\bSpin$ in 
  Example \ref{eg:SPin monoid space}. However, the isomorphism involves
  a slight twist that we want to make explicit.
  Indeed, the real Clifford algebra $\Cl(V)$ was defined from the tensor algebra
  by imposing the relations $v\tensor v=- |v|^2\cdot 1$, so that unit vectors of $V$
  square to $-1$ in $\Cl(V)$. 
  Over the field $\mR$ it makes a difference whether  we make
  $v\tensor v$ equal to $|v|^2\cdot 1$ or $- |v|^2\cdot 1$, 
  but the presence of the imaginary unit makes the difference disappear over $\mC$.
  Indeed, the $\mR$-linear map
  \[ \psi\ : \ V \ \ \to \ \mCl(V) \ , \quad v\ \longmapsto \ i\cdot [v] \]
  satisfies $\psi(v)^2=-|v|^2\cdot 1$, so it 
  extends to a morphism of $\mC$-algebras $\mC\tensor_\mR\Cl(V)\to\mCl(V)$
  by the universal property of the former, and this morphism is an isomorphism.\index{subject}{Clifford algebra|)}

  \medskip

  Another key piece of structure is a continuous based map
  \begin{equation}\label{eq:functional_calculus}
     \fc \ : \ S^V \ \to \  C^\ast_{\gr}(s,\mCl(V)) \ , \quad v\ \longmapsto \ (-)[v] \ ,
  \end{equation}
  often referred to as `functional calculus'.
  For $v\in V$ the $\ast$-homomorphism $\fc(v)$ is given on homogeneous elements
  of $s$ by
  \[ f[v]\ = \ \fc(v)(f)\ = \ 
  \begin{cases}
    \quad f(|v|)\cdot 1 & \text{ when $f$ is even, and}\\
    \frac{ f(|v|)}{|v|}\cdot [v] & \text{ when $f$ is odd.}
  \end{cases}
  \]
  For $v=0$ the formula for odd functions is to be interpreted as $f[0]=0$; 
  this is continuous because for $v\ne 0$, 
  the norm of $f(|v|)/|v|\cdot [v]$ is $f(|v|)$, which tends to $f(0)=0$ 
  if $v$ tends to~0. If the norm of $v$ tends to infinity, then $f(|v|)$ tends to~0, so 
  $\fc(v)$ tends to the constant $\ast$-homomorphism 
  with value~0, the basepoint of $C^\ast_{\gr}(s,\mCl(V))$.
  Hence $\fc$ extends to a continuous
  map on the one-point compactification $S^V$. 
  The functional calculus map is $O(V)$-equivariant for 
  the $O(V)$-action on the mapping space $C^\ast_{\gr}(s,\mCl(V))$
  through the action on the target.
  The functional calculus maps are multiplicative in the sense that
  the following diagram commutes:
  \begin{equation}\begin{aligned}\label{eq:fc_multiplicative}
      \xymatrix@C=5mm@R=8mm{
        S^V \sm S^W \ar[dd]_\iso \ar[rr]^-{\fc\sm\fc} &&
        C^\ast_{\gr}(s,\mCl(V))\sm C^\ast_{\gr}(s,\mCl(W))\ar[d]^{\hat\tensor} \\
        && C^\ast_{\gr}(s\hat\tensor s,\mCl(V)\tensor\mCl(W)) 
        \ar[d]^{\Delta^*} \\
        S^{V\oplus W}\ar[r]_-{\fc} &
        C^\ast_{\gr}(s,\mCl(V\oplus W))\ar[r]^-{\eqref{eq:Clifford tensor iso}}_-{(\mu_{V,W})_*}& 
        C^\ast_{\gr}(s,\mCl(V)\tensor\mCl(W))
} 
    \end{aligned}
   \end{equation}
   To see this we consider $v\in V$ and $w\in W$, chase the element $v\sm w$
   both ways around the square and compare the results 
   on the two functions \eqref{eq:s generators u and v} that generate the algebra $s$.
   Indeed, for the even generating function we have
   \begin{align*}
     \mu_{V,W}(\fc(v,w)(u_+))\ &= \  \mu_{V,W}( e^{-|v|^2-|w|^2}\cdot 1)\ = \ 
     ( e^{-|v|^2}\cdot 1)\tensor (e^{- |w|^2}\cdot 1)\\
     &=\    \fc(v)(u_+)\tensor\fc(w)(u_+)\ 
     = \    (\fc(v)\hat\tensor\fc(w))(u_+\tensor u_+)\\ 
     &= \ (\fc(v)\hat\tensor\fc(w))(\Delta(u_+))\ \ ;
   \end{align*}
   and similarly for the odd generating function:
   \begin{align*}
     \mu_{V,W}(&\fc(v,w)(u_-))\ = \  \mu_{V,W}(  e^{-|v|^2-|w|^2}\cdot [v,w])\\
     &= \    e^{- |v|^2-|w|^2}\cdot ([v]\tensor 1+1\tensor [w])\\ 
     &= \ 
      ( e^{-|v|^2}\cdot [v])\tensor (e^{- |w|^2}\cdot 1)\ +\
    ( e^{-|v|^2}\cdot 1)\tensor ( e^{- |w|^2}\cdot [w])\\
     &=\  \fc(v)(u_-)\tensor\fc(w)(u_+) \ + \   \fc(v)(u_+)\tensor\fc(w)(u_-)\\ 
     &= \    (\fc(v)\hat\tensor\fc(w))(u_-\tensor u_+ + u_+\tensor u_-)\ 
     = \ (\fc(v)\hat\tensor\fc(w))(\Delta(u_-))\  .
   \end{align*}

\index{subject}{K-theory@$K$-theory!periodic global|(}
Now we have all ingredients for the periodic global $K$-theory spectrum $\bKU$.

\begin{construction}[Periodic global $K$-theory]\label{con:global KU}
We let $V$ be a euclidean inner product space.
As we explained in Proposition \ref{prop:Sym induced inner product},
the symmetric algebra $\Sym(V_\mC)$ of the complexification
inherits a hermitian inner product and an $O(V)$-action by $\mC$-linear isometries.
The inner product space $\Sym(V_\mC)$ 
is usually infinite dimensional (unless $V=0$),
but it is not complete. We denote by $\Hc_V$ 
the Hilbert space completion of $\Sym(V_\mC)$.
Since the action of $O(V)$ on $\Sym(V_\mC)$
is by linear isometries, it extends to an analogous
action on the completion $\Hc_V$.
So $\Hc_V$ becomes a complex Hilbert space representation
of the orthogonal group $O(V)$. 
We denote by $\Kc_V$ the $C^\ast$-algebra
of compact operators on the Hilbert space $\Hc_V$,
see for example \cite[Sec.\,2.4]{murphy-Cstar}.\index{subject}{compact operators!on a Hilbert space}

The orthogonal spectrum $\bKU$ assigns to a euclidean inner product space $V$
the space\index{symbol}{$\bKU$ - {periodic global $K$-theory}}\index{subject}{K-theory@$K$-theory!periodic global}
\[ \bKU(V) \ = \ C^\ast_{\gr}(s,\mCl(V)\tensor \Kc_V) \]
of $\mZ/2$-graded $\ast$-homomorphisms from $s$ to
the tensor product, over $\mC$, of $\mCl(V)$ and $\Kc_V$.
Here we consider $\Kc_V$ as evenly graded, so the grading comes
entirely from the grading of the Clifford algebra.
The topology is the topology of pointwise convergence in the 
operator norm of $\Kc_V$; the basepoint is the zero $\ast$-homomorphism.

The continuous action of the orthogonal group $O(V)$ by linear isometries of $\Hc_V$ 
induces an action on the algebra $\Kc_V$ by conjugation,
and together with the action on $\mCl(V)$ it gives 
an $O(V)$-action on $\mCl(V)\tensor\Kc_V$ by graded $\ast$-auto\-morphisms.
This induces an $O(V)$-action on the mapping space $\bKU(V)$.

The multiplication of the spectrum $\bKU$ starts from the 
$(O(V)\times O(W))$-equivariant isometric isomorphism
\begin{align*}
\Sym(V_\mC)\tensor \Sym(W_\mC) \ &\iso\
\Sym(V_\mC\oplus W_\mC)  \ \iso \ \Sym((V\oplus W)_\mC)  \ ,
\end{align*}
compare \eqref{eq:tensor and Sym}.
This extends to an isometry of the Hilbert space completions
\begin{align*}
\Hc_V\hat\tensor \Hc_W \ = \ \widehat{\Sym}(V_\mC)\hat\tensor \widehat{\Sym}(W_\mC) \ 
&\iso\ (\Sym(V_\mC)\tensor \Sym(W_\mC))^\wedge \\ 
&\iso \  \widehat{\Sym}((V\oplus W)_\mC)  \ = \ \Hc_{V\oplus W}\ .
\end{align*}
Tensor product of compact operators and conjugation with this isometry 
is an injective homomorphism of $\mC$-algebras
\begin{equation}\label{eq:multiply Ks}
   \Kc_V\tensor \Kc_W\ \to\ \Kc( \Hc_V\hat\tensor \Hc_W) \ 
\iso\ \Kc(\Hc_{V\oplus W}) = \Kc_{V\oplus W} \ .
\end{equation}
The source of this map is the algebraic tensor product, which is not complete; 
the homomorphism \eqref{eq:multiply Ks} extends to an isomorphism of $C^\ast$-algebras
\begin{equation}\label{eq:tensor compacts}
   \Kc_V\hat\tensor \Kc_W\ \iso\  \Kc_{V\oplus W} 
\end{equation}
from the spatial (minimal) tensor product; this is in fact tautologically
true, as $\Kc_V$ acts on the Hilbert space $\Hc_V$ and
so the spatial norm on the algebraic tensor product can be defined via
the embedding \eqref{eq:multiply Ks}.
Combining \eqref{eq:tensor compacts} with
the isomorphism between $\mCl(V)\tensor \mCl(W)$
and $\mCl(V\oplus W)$ specified in \eqref{eq:Clifford tensor iso},
we obtain an isomorphism of graded $C^\ast$-algebras
\begin{align}\label{eq:compact operator iso}
 (\mCl(V)\tensor \Kc_V)\hat\tensor(\mCl(W)\tensor \Kc_W)\ &\iso \ 
(\mCl(V)\tensor\mCl(W))\tensor(\Kc_V\hat\tensor\Kc_W) \nonumber\\  
&\iso \ \mCl(V\oplus W)\tensor \Kc_{V\oplus W}\ .  
\end{align}
The first step is the symmetry isomorphism,
which in this special case does not involve any signs because $\Kc_V$ 
is concentrated in even grading.
The multiplication map
\[ \mu_{V,W}\ : \ \bKU(V)\sm \bKU(W)\ \to \ \bKU(V\oplus W) \]
is now defined as the composite
\begin{align*}
C^\ast_{\gr}(s,\mCl(V)\tensor \Kc_V)&\sm C^\ast_{\gr}(s,\mCl(W)\tensor \Kc_W)\\ 
&\xra{\ \hat\tensor\ }\ 
C^\ast_{\gr}(s\hat\tensor s,\ (\mCl(V)\tensor \Kc_V)\hat\tensor(\mCl(W)\tensor \Kc_W)) \\
 &\xra{C^\ast_{\gr}(\Delta, \eqref{eq:compact operator iso})}\   
C^\ast_{\gr}(s,\mCl(V\oplus W)\tensor \Kc_{V\oplus W})\ .
\end{align*}
These multiplication maps are associative and commutative.
The unit map 
\[ \eta_V \ : \ S^V \ \to \  C^\ast_{\gr}(s,\mCl(V)\tensor \Kc_V) \ = \ \bKU(V) \]
is defined as the composite
\[ \ S^V \ \xra{\ \fc\ } \ C^\ast_{\gr}(s,\mCl(V))\ 
\xra{(-\tensor p_0)_*} \ C^\ast_{\gr}(s,\mCl(V)\tensor\Kc_V)\ = \ \bKU(V) \ ,  \]
where $p_0\in\Kc_V$ is the orthogonal projection onto the
constant summand in the symmetric algebra. In other words, 
\[  \eta_V(v)(f)\ = \ f[v]\tensor p_0\ .  \]
The multiplicativity of the unit maps follows from 
the multiplicativity \eqref{eq:fc_multiplicative}
of the functional calculus maps and the fact that 
the isomorphism \eqref{eq:tensor compacts} sends $p_0\tensor p_0$ to $p_0$.
\end{construction}

\begin{construction}\label{con:j_ku2KU}
The connective and periodic global $K$-theory spectra are related by a morphism 
of ultra-commutative ring spectra $j:\bku\to\bKU$.\index{subject}{K-theory@$K$-theory!connective global}
To define the value 
\[  j(V)\ : \  \bku(V)\ = \ \bk(\Sym(V_\mC),S^V) \ \to \ 
C^\ast_{\gr}(s,\mCl(V)\tensor \Kc_V)  \ = \ \bKU(V)  \]
at an inner product space $V$ we consider a configuration 
\[ [E_1,\dots,E_n;\,v_1,\dots,v_n] \ \in \  \bk(\Sym(V_\mC),S^V)  \]
of pairwise orthogonal, finite-dimensional subspaces on $S^V$.
The associated $\ast$-homomorphism 
\[ j(V)[E_1,\dots,E_n;\,v_1,\dots,v_n]\ :\ s\ \to\ \mCl(V)\tensor\Kc_V \]
is then defined on a function $f\in s$ by
\[ j(V)[E_1,\dots,E_n;\,v_1,\dots,v_n](f) \ = \ 
{\sum}_{i=1}^n\ f[v_i]\tensor p_{E_i}\ ,\]
where $f[v_i]=\fc(v_i)(f)$ is the functional calculus \eqref{eq:functional_calculus}
and $p_E$ denotes the orthogonal projection onto a subspace $E$.
The fact that the morphism $j(V)[E_k;\,v_k]$ is $\mZ/2$-graded
is a direct consequence of the same property for the functional calculus map.
The verification that the map $j(V)$ is $O(V)$-equivariant
is straightforward, and we omit it.
The verification that the maps $j(V)$ are multiplicative
amounts to the multiplicativity \eqref{eq:fc_multiplicative}
of the functional calculus maps and the fact that 
the isomorphism \eqref{eq:tensor compacts}
sends $p_E\tensor p_F$ to $p_{E\tensor F}$.
Moreover, the composite of $j(V)$ with the unit map $S^V\to\bku(V)$
is the unit map of $\bKU$, so the maps $j$ indeed form a morphism
of ultra-commutative ring spectra. 
\end{construction}

\begin{rk}
We have defined $\bKU(V)$ in a slightly
different way, compared to the presentation of Joachim \cite{joachim-coherences}. 
We use the Hilbert space completion of $\Sym(V_\mC)$,
instead of the Hilbert space $L^2(V)$ of $\mC$-valued square integrable
functions on $V$. These two Hilbert spaces are naturally isomorphic,
as follows. We use the inner product on $V$ to identify it
with its dual space $V^\ast$, and hence the symmetric algebra
$\Sym(V_\mC)$ with $\Sym(V^\ast_\mC)$. Elements of $\Sym(V^\ast_\mC)$
are complex valued polynomial functions on $V$.
We make them square integrable by multiplying with the rapidly
decaying function $\epsilon:V\to\mR$, $\epsilon(v)=e^{-|v|^2}$. 
This provides a linear isometric embedding
\[ \Sym(V^\ast_\mC) \ \to \ L^2(V)\ , \quad
f \ \longmapsto \ f\cdot \epsilon  \]
with dense image. 
Altogether this exhibits $L^2(V)$ as a Hilbert space
completion of $\Sym(V_\mC)$.

\Danger The reader should beware that 
for us the $C^\ast$-algebra $\Kc_V$ 
is evenly graded, whereas Joachim
uses a non-trivial grading on $L^2(V)$ and $\Kc(L^2(V))$ 
by even and odd functions.
Our grading convention is necessary to ensure that the 
map $j(V): \bku(V)\to C^\ast_{\gr}(s,\mCl(V)\tensor\Kc_V)=\bKU(V)$
defined in Construction \ref{con:j_ku2KU} 
takes values in {\em graded} $\ast$-homomorphisms.
\end{rk}

\begin{rk}[Homotopy type of $\bKU(V)$]
Up to isomorphism, there are only
three different graded $C^\ast$-algebras of the form $\mCl(V)\tensor\Kc_V$,
and only three different homeomorphism types of spaces $\bKU(V)$.
The case $V=0$ is degenerate in that $\Sym(0)$ is just the copy of $\mC$ generated by~1.
This is already complete, so $\Hc_0=\mC$ and $\mCl(0)=\Kc_0=\mC$. 
The space $\bKU(0)=C^\ast_{\gr}(s,\mCl(0)\tensor \Kc_0)$ is thus
discrete with two points: the zero homomorphism as basepoint
and the augmentation 
\[ s\ \to \ \mCl(0)\tensor \Kc_0\ , \quad f \ \longmapsto \ f(0)\tensor 1 \ .\]
This non-basepoint is a unit for the multiplication maps $\mu_{V,W}$.
As we explain in Remark \ref{rk:equivariant inverse Bott}
below, the space $\bKU(\mR^n)$ is non-canonically homeomorphic to $\bKU(\mR^{n+2})$
for all $n\geq 1$. So all odd-dimensional terms of $\bKU$ are homeomorphic,
and all even-dimensional terms in positive degrees are homeomorphic.

We can also identify the space $\bKU(\mR)$ as follows. We claim that the map
\[ j(\mR)\ : \ \bku(\mR) \ \to \ \bKU(\mR) \]
is a homotopy equivalence.
Indeed, if $A$ is an (ungraded) $C^\ast$-algebra, then 
$\ast$-homomorphisms from $s$ to $A$ (in the ungraded sense) 
biject with graded $\ast$-homomorphisms 
from $s$ to $\mCl(\mR)\tensor A$, via 
\[  C^\ast_{\gr}(s,\mCl(\mR)\tensor A)  \ \iso \ C^\ast(s, A)
\ ,\quad \varphi \ \longmapsto \ \pi\circ \varphi\ ,\]
where $\pi : \mCl(\mR)\tensor A \to A$ is the $\ast$-homomorphism 
that sends both $1\tensor a$ and $[1]\tensor a$ to $a$.
The composite
\[ S^1 \ \xra{\ \fc \ } \ C^\ast_{\gr}(s,\mCl(\mR)) \ \xra[\iso]{\pi\circ -} \ 
 C^\ast(s,\mC)   \]
sends $x\in S^1$ to evaluation at $x$.
This implies that the following square commutes:
\[ \xymatrix@C=10mm{ 
\bk(\Sym(\mC),S^1)\ = \ 
\bku(\mR)\ar[r]^-{j(\mR)} \ar[d]_\iso & 
\bKU(\mR) = C^\ast_{\gr}(s, \mCl(\mR)\tensor \Kc_\mR) \ar[d]_\iso^{\pi\circ-}\\
\colim_n C^\ast(s, \End_\mC(\Sym^{[n]}(\mC))) \ar[r] & C^\ast(s, \Kc_\mR) } \]
The lower horizontal map is induced by the $\ast$-homomorphisms
$\End_\mC(\Sym^{[n]})\to\Kc_\mR$ that extend an endomorphism by zero
on the orthogonal complement, and 
the left vertical map
sends a configuration $[E_1,\dots,E_n;x_1,\dots,x_n]$ to 
the $\ast$-homomorphism
\[ f\ \longmapsto \ {\sum}_{i=1}^n \, f(x_i)\cdot p_{E_i} \ .\]
The two vertical maps are homeomorphisms
(compare Remark \ref{rk:conf versus C^ast}) and 
the lower horizontal map is a homotopy equivalence 
by \cite[Prop.\,1.2]{segal-K-homology} 
or \cite[Prop.\,4.6]{hohnhold-stolz-teichner}
(or rather its complex analog).
So the map $j(\mR)$ is a homotopy equivalence, and $\bKU(\mR)$
is homotopy equivalent to the infinite unitary group $U$.

The spectrum $\bKU$ is positive $\Omega$-spectrum (in the non-equivariant sense)
by Theorem \ref{thm:KU Omega property} below.
So the adjoint structure map $\bKU(\mR^2)\to\Omega\bKU(\mR^3)$
is a weak equivalence, and the target is homeomorphic to
$\Omega\bKU(\mR)\simeq \Omega U$. So
$\bKU(\mR^2)$ has the weak homotopy type of $\mZ\times B U$,
by Bott periodicity.
\end{rk}

The analysis of the global homotopy type of $\bKU$ depends 
on the operator theoretic formulation of equivariant Bott periodicity
that we now recall.
We let $G$ be a compact Lie group and $V$ an orthogonal $G$-representation.
We denote by $C_0(V, \mCl(V))$
the $G$-$C^\ast$-algebra of continuous $\mCl(V)$-valued functions on $V$
that vanish at infinity.
Since the Clifford algebra is finite-dimensional, the  map
\begin{equation}  \label{eq:C_0and_Cl}
    C_0(V)\tensor \mCl(V)\ \to \ C_0(V, \mCl(V)) \ , \quad 
    f\tensor x \ \longmapsto \ \ f(-)\cdot x
\end{equation}
is an isomorphism of $C^\ast$-algebras.
Functional calculus \eqref{eq:functional_calculus} 
provides a distinguished graded $\ast$-homomorphism
\[ \beta_V \ : \ s \ \to \  C_0(V,\mCl(V)) \ , \quad
\beta_V(f)(v)\ = \ f[v] \ . \]
Since the functional calculus map is $G$-equivariant, 
$\beta_V$ takes values in the $G$-fixed points of $C_0(V,\mCl(V))$ 
for the conjugation action of $G$.

We let $\Hc_G$ be any complete $G$-Hilbert space universe, i.e.,
a Hilbert $G$-representation that is isometrically isomorphic
to the completion of a complete unitary $G$-universe.
We let $\Kc_G$ be the $G$-$C^\ast$-algebra of 
(not necessarily equivariant) compact operators on $\Hc_G$, with
$G$ acting by conjugation.

\begin{theorem}[Equivariant Bott periodicity]\label{thm:equivariant Bott}
Let $G$ be a compact Lie group and $V$ an orthogonal $G$-representation.
Then for every graded $G$-$C^\ast$-algebra $A$, the map
\[ \beta_V\cdot - \ : \ 
C^\ast_{\gr}(s,  A\tensor \Kc_G )\ \to \ 
C^\ast_{\gr}(s, C_0(V,\mCl(V)) \tensor  A\tensor\Kc_G ) \]
is a $G$-weak equivalence.
\end{theorem}

When $G=e$ is a trivial group, $\Kc_G$ reduces to the $C^\ast$-algebra
of compact operators on a separable Hilbert space,
and then this formulation of Bott periodicity can 
be found in \cite[Thm.\,1.14]{higson-guentner}.
Unfortunately, we are lacking a reference in the generality of
compact Lie groups. If one specializes \cite[Sec.\,3, Thm.\,5]{higson-kasparov-trout} 
to finite-dimensional representations of finite groups, 
one obtains a formulation very close to 
(but not exactly the same as) Theorem \ref{thm:equivariant Bott}.

Since the unit map $\eta_V:S^V\to\bKU(V)$ is essentially the functional
calculus map, a direct consequence of the formulation of equivariant Bott periodicity
in Theorem \ref{thm:equivariant Bott} is that $\bKU$ is `eventually'
a global $\Omega$-spectrum, see Theorem \ref{thm:KU Omega property} below.

\medskip

We showed in Proposition \ref{prop:configuration_invariance}
that an equivariant linear isometric embedding between
two complete complex $G$-universes induces a $G$-homotopy equivalence
between configuration spaces with labels in the two universes.
We will now establish an analogous property for the $C^\ast$-algebras
of compact operators on the Hilbert space completions.

\begin{construction}
We recall that a $\mC$-linear isometric embedding $\varphi:\Hc\to\Hc'$
between complex separable Hilbert spaces gives rise to a preferred $\ast$-homomorphism
$\Kc(\varphi)$ between the $C^\ast$-algebras of compact operators.
We refer to $\Kc(\varphi)$ as the {\em conjugation homomorphism}\index{subject}{conjugation homomorphism!between compact operators}
induced by $\varphi$.
We emphasize that there is no further hypothesis on $\varphi$ besides
$\mC$-linearity and the requirement
\[ (x,y)_\Hc \ = \ (\varphi(x),\varphi(y))_{\Hc'} \]
for all $x,y\in\Hc$;
in particular, $\varphi$ is {\em not} assumed to be bounded, complemented, or adjointable.
To characterize $\Kc(\varphi)$ we need the following notation:
for a finite-dimensional subspace $L$ of $\Hc$ we write
$p_L\in\Kc(\Hc)$ for the orthogonal projection onto $L$.
There is then a unique $\ast$-homomorphism
\[ \Kc(\varphi) \ : \ \Kc(\Hc)\ \to \ \Kc(\Hc') 
\text{\qquad such that\qquad}
 \Kc(\varphi)(p_L) \  = \ p_{\varphi(L)}  \]
for all finite-dimensional $L$ in $\Hc$.
The uniqueness is a consequence of the fact that 
the linear span of the finite rank projections is the space of finite rank
operators \cite[Thm.\,2.4.6]{murphy-Cstar}
and the finite rank operators are dense in $\Kc(\Hc)$ \cite[Thm.\,2.4.5]{murphy-Cstar}.
The construction of $\Kc(\varphi)$ can
be found in \cite[Lemma 4.1]{meyer-equivariant generalized}.
 
Informally, one can think of $\Kc(\varphi)$ as `conjugation by $\varphi$',
combined with extension by~0 on the orthogonal complement, which justifies the name.
Indeed, if $\varphi$ happens to be complemented, i.e.,
$\Hc'$ is the direct sum of $\varphi(\Hc)$ and the orthogonal complement
$\varphi(\Hc)^\perp$, then this is literally true.
In particular, if $\varphi$ is a linear isometric isomorphism, 
then $\Kc(\varphi)$ is conjugation by $\varphi$.

The construction is covariantly functorial in the linear isometric embedding.
Indeed, if $\psi:\Hc'\to\Hc''$ is another linear isometric embedding,
then 
\[ \Kc(\psi)\circ\Kc(\varphi)\ = \  \Kc(\psi\circ\varphi) \]
by the uniqueness clause.

The construction generalizes to the 
$G$-equivariant context, where $G$ is any compact Lie group.
More precisely, we assume that $\Hc$ and $\Hc'$ are 
separable complex $G$-Hilbert space representations. 
The $C^\ast$-algebra $\Kc(\Hc)$ then inherits a continuous
$G$-action by conjugation with the $G$-action on $\Hc$. 
If $\varphi:\Hc\to\Hc'$ is $G$-equivariant,
then the functoriality of the homomorphism $\Kc(\varphi)$
as a function of $\varphi$ directly implies that
the $\ast$-homomorphism $\Kc(\varphi):\Kc(\Hc)\to\Kc(\Hc')$
is also $G$-equivariant. 
\end{construction}

\begin{prop}\label{prop:compact operator invariance}
Let $G$ be a compact Lie group and
$u:\Uc\to \Vc$ a $G$-equivariant $\mC$-linear isometric embedding
between two complete complex $G$-universes. Then the morphism of $G$-$C^\ast$-algebras
\[ \Kc(\hat u)\ :\ \Kc(\widehat \Uc)\ \to \ \Kc(\widehat \Vc)\]
is a $G$-equivariant homotopy equivalence.
\end{prop}
\begin{proof}
We start with the special case where $\Vc=\Uc$.
The space of $G$-equivariant linear isometric embeddings from $\Uc$
to itself is contractible.
By passage to completions, a $G$-equivariant homotopy from $u$
to the identity induces a $G$-homotopy
\[ \Phi \ : \ \widehat\Uc\times [0,1]\ \to \ \widehat\Uc \ , \]
i.e., the map $\Phi(x,-):[0,1]\to\widehat\Uc$ is continuous for every $x\in\widehat\Uc$,
and the map $\Phi(-,t)$ is a linear isometric embedding for every $t\in[0,1]$.
The map 
\[ [0,1]\ \to \ C^*(\Kc(\widehat\Uc),\Kc(\widehat\Uc))\ , \quad t \ \longmapsto \ \Kc(\Phi(-,t)) \]
is then continuous, see for example \cite[p.\,211]{meyer-equivariant generalized}.
So $\Kc(\hat u):\Kc(\widehat \Uc)\to\Kc(\widehat \Uc)$
is homotopic, through $G$-equivariant $\ast$-homomorphisms,
to the identity.

In the general case we choose an equivariant linear isometry $v:\Vc\iso \Uc$.
By the previous paragraph the $G$-$C^\ast$-homomorphisms
$\Kc(\hat v)\circ\Kc(\hat u)=\Kc(\widehat{v u})$
and $\Kc(\hat u)\circ\Kc(\hat v)=\Kc(\widehat{u v})$
are $G$-homotopic to the respective identity maps.
\end{proof}

We recall from Definition \ref{def:ample} that an orthogonal $G$-representation
is {\em ample}\index{subject}{ample $G$-representation}
if its complexified symmetric algebra is a complete complex $G$-universe.

\begin{theorem}\label{thm:KU Omega property} 
For every orthogonal $G$-representation $V$ 
and every ample $G$-representation $W$, the adjoint structure map
\[ \tilde\sigma_{V,W}\ : \ \bKU(W) \ \to \  \map_*(S^V, \bKU(V\oplus W) ) \]
is a $G$-weak equivalence.
\end{theorem}
\begin{proof}
We let $u:\Sym(W_\mC)\to \Sym( (V\oplus W)_\mC)$ be 
the linear isometric embedding
induced by the direct summand embedding $W\to V\oplus W$.
Since $W$ is ample, $u$ is an equivariant 
$\mC$-linear isometric embedding between complete 
complex $G$-universes.
So the induced map $\Kc(\hat u):\Kc_W\to \Kc_{V\oplus W}$
of compact operators is a $G$-equivariant homotopy equivalence of $C^\ast$-algebras
by Proposition \ref{prop:compact operator invariance}.

The adjoint structure map $\tilde\sigma_{V,W}$ factors as the composite:
\begin{align*}
  \bKU(W) \ = \  &C^\ast_{\gr}(s, \mCl(W)\tensor\Kc_W)\\ 
\xra{\beta_V\cdot -} \ 
&C^\ast_{\gr}(s, C_0(V,\mCl(V))\tensor \mCl(W)\tensor\Kc_W)\\
\xra{ \Kc(\hat u)_*} \ 
&C^\ast_{\gr}(s, C_0(V,\mCl(V)) \tensor \mCl(W)\tensor\Kc_{V\oplus W})\\
\xra[\iso]{\eqref{eq:C_0and_Cl}}\
 &C^\ast_{\gr}(s, C_0(V)\tensor \mCl(V\oplus W)\tensor\Kc_{V\oplus W})\\
\iso \quad &\map_*(S^V, C^\ast_{\gr}(s, \mCl(V\oplus W)\tensor\Kc_{V\oplus W}))\ =\
 \map_*(S^V, \bKU(V\oplus W) )
\end{align*}
Since $W$ is ample, $\Hc_W$ is a complete Hilbert $G$-universe 
and so $\Kc_W$ is a $\Kc_G$. 
Bott periodicity (Theorem \ref{thm:equivariant Bott})
for the $G$-$C^\ast$-algebra $\mCl(W)$
shows that the first map is a $G$-weak equivalence.
The second map is a $G$-homotopy equivalence by the first paragraph,
so altogether $\tilde\sigma_{V,W}$ is a $G$-weak equivalence.
\end{proof}

The shift $\sh X=\sh^\mR X$ of an orthogonal spectrum $X$ was
defined in \eqref{eq:define_shift}.
Shifting an orthogonal spectrum provides a delooping by 
Proposition \ref{prop:global equiv preservation}~(i).
The eigenspace morphism was defined in \eqref{eq:define_eig}.

\begin{theorem}\label{thm:U_vs_shift_KU}
The composite 
\[ \bU \ \xra{\ \eig \ } \ \Omega^\bullet(\sh \bku) \ \xra{\Omega^\bullet(\sh j)} \
\Omega^\bullet(\sh \bKU)  \]
is a global equivalence of orthogonal spaces.
\end{theorem}
\begin{proof}
As in the proof of Theorem \ref{thm:U_vs_shift_ku}
we let $\bar\bU$ denote the orthogonal space with
$\bar\bU(V)= U(\Sym((V\oplus\mR)_\mC))$,
and we factor the eigenspace decomposition morphism as the composite
\[ \bU \ \to \ \bar\bU \ \xra{\ \widebar\eig \ } \ \Omega^\bullet(\sh \bku) \ , \]
where the first morphism is the global equivalence of orthogonal spaces
induced by the natural linear isometric embedding of $V_\mC$
into the linear summand in $\Sym((V\oplus\mR)_\mC)$.

We claim that composite
\[ \bar\bU \ \xra{\ \widebar\eig \ } \ \Omega^\bullet(\sh \bku) 
\ \xra{\Omega^\bullet(\sh j)} \ \Omega^\bullet(\sh \bKU) \]
is an `eventually strong level equivalence'. More precisely,
we show that for every ample $G$-representation $V$ the map
\begin{align*}
   U(\Sym((V\oplus\mR)_\mC))\ = \ \bar\bU(V) \ &\to 
\ (\Omega^\bullet(\sh \bKU))(V)  \\ 
&= \ \map_*(S^V, C^\ast_{\gr}(s, \mCl(V\oplus\mR)\tensor \Kc_{V\oplus\mR}))  
\end{align*}
is a $G$-weak equivalence.
This map factors through the $G$-map
\begin{align*}
 U(\Sym((V\oplus\mR)_\mC)) \iso \bk&(\Sym((V\oplus\mR)_\mC),S^1) \ \to \ 
\hat \bk(\Hc_{V\oplus\mR},S^1)\\ 
&\iso \  C^\ast(s, \Kc_{V\oplus\mR})\ \iso \ 
C^\ast_{\gr}(s, \mCl(\mR)\tensor\Kc_{V\oplus\mR}) \ .  
\end{align*}
Here $\hat \bk(\Hc_{V\oplus\mR},S^1)$ is the space of configurations,
not necessarily finite, of pairwise orthogonal finite-dimensional 
subspaces of $\Hc_{V\oplus\mR}$,
that are allowed to accumulate around $\infty$.
This map is a $G$-equivariant homotopy equivalence;
the non-equivariant argument (or rather its real analog)
can be found in \cite[Prop.\,1.2]{segal-K-homology} 
or \cite[Prop.\,4.6]{hohnhold-stolz-teichner},
and the explicit homotopies in \cite[Prop.\,4.6]{hohnhold-stolz-teichner}
work just the same way in the presence of an isometric action of a compact Lie group.

By the operator theoretic equivariant 
Bott periodicity theorem (Theorem \ref{thm:equivariant Bott}), 
multiplication by the graded $\ast$-homomorphism $\beta_V:s\to C_0(V,\mCl(V))$ 
is a $G$-weak equivalence 
\[ \beta_V \cdot -\ : \ 
C^\ast_{\gr}(s, \mCl(\mR)\tensor \Kc_{V\oplus\mR} )\ \to \ 
C^\ast_{\gr}(s, C_0(V,\mCl(V))\tensor \mCl(\mR)\tensor\Kc_{V\oplus\mR} ) \ .\]
We use here that $\mCl(\mR)\tensor\Kc_{V\oplus\mR}$ is a $\Kc_G$ because
the $G$-representation $V$ is ample.
This proves the claim because the target is $G$-equivariantly homeomorphic to 
\begin{align*}
 C^\ast_{\gr}(s, C_0(V)\tensor \mCl(V)\tensor\mCl(\mR)&\tensor\Kc_{V\oplus\mR}) \\ 
_\eqref{eq:Clifford tensor iso} &\iso \ 
 C^\ast_{\gr}(s, C_0(V)\tensor \mCl(V\oplus\mR)\tensor\Kc_{V\oplus\mR}) \\ 
&\iso \ \map_*(S^V, C^\ast_{\gr}(s, \mCl(V\oplus\mR)\tensor\Kc_{V\oplus\mR}))  \\ 
&= \ \map_*(S^V, \sh\bKU(V))\ .\qedhere
\end{align*}
\end{proof}

We emphasize that in contrast to Theorem \ref{thm:U_vs_shift_ku},
the composite in Theorem \ref{thm:U_vs_shift_KU}
is a global equivalence for all compact Lie groups (as opposed to 
only a $\Fin$-global equivalence). Another key difference is that the
spectrum $\bku$ is globally connective; in contrast, we
will see in Theorem \ref{thm:Bott and inverse} below that $\bKU$ is Bott periodic,
so applying $\Omega^\bullet$ is losing information.

\begin{rk}[$\bKU$ globally deloops $\bBUP$]
A corollary of the previous theorem is that $\bKU$ globally deloops
the orthogonal space $\bBUP$ defined in Example \ref{eg:BUP and BSpP}.
Indeed,  the global formulation of Bott periodicity in 
Theorem \ref{thm:global Bott periodicity} provides a global equivalence
of ultra-commutative monoids $\bar\beta:\bBUP\to\Omega(\sh_\tensor \bU)$.
Combined with Theorem \ref{thm:U_vs_shift_KU} this yields a chain of global equivalences
of orthogonal spaces
\[ \bBUP \ \xra[\simeq]{\ \bar\beta\ } \ \Omega(\sh_\tensor \bU) \ \xla[\ \simeq\ ] 
\ \Omega \bU \ \xra[\simeq]{\Omega( (\Omega^\bullet(\sh j))\circ\eig)} 
\Omega (\Omega^\bullet(\sh \bKU)) 
\ \xla[\simeq]{\Omega^\bullet(\tilde\lambda_{\bKU})}\  \Omega^\bullet \bKU \ . \]
\end{rk}

We spell out another corollary of Theorem \ref{thm:U_vs_shift_KU}.
We recall that $\bK_G(A)$ denotes the equivariant $K$-group of a $G$-space $A$,
i.e., the Grothendieck group of isomorphism classes of $G$-vector bundles over $A$.
A ring homomorphism $[-]$ from this Grothendieck group
to the equivariant cohomology group $\bku^0_G(A_+)$ 
was defined in \eqref{eq:K_G_to_ku}.

\begin{cor}\label{cor:global_K_represents}
For every compact Lie group $G$ and every finite $G$-CW-complex $A$ the composite
\[
\bK_G(A) \ \xra{\ [-]\ }  \bku^0_G(A_+) \ \xra{\ j_*\ } \ \bKU^0_G(A_+) 
  \]
is an isomorphism.\index{subject}{equivariant $K$-theory}
\end{cor}
\begin{proof}
We contemplate the commutative diagram of orthogonal spaces:
\[  \xymatrix@C=13mm{
\bGr^\mC  \ar[d]_{\beta} \ar[r]^-c & 
\Omega^\bullet \bku\ar[d]^{\Omega^\bullet \tilde\lambda_{\bku}}_\simeq
\ar[r]^-{\Omega^\bullet j}
& \Omega^\bullet \bKU\ar[d]^{\Omega^\bullet \tilde\lambda_{\bKU}}_\simeq\\
\Omega \bU \ar[r]_-{\Omega \eig} & 
\Omega\left( \Omega^\bullet (\sh \bku)\right) 
\ar[r]_-{\Omega (\Omega^\bullet(\sh j))} & 
\Omega( \Omega^\bullet(\sh \bKU)) } \]
The morphism $(\Omega^\bullet(\sh j)) \circ \eig$
is a global equivalence by Theorem \ref{thm:U_vs_shift_KU},
so the lower horizontal composite is also a global equivalence. 
The same argument as for $\Omega^\bullet\tilde\lambda_{\bku}$ in
the proof of Theorem \ref{thm:[A,Gr] to [A,bku]}
shows that the right vertical morphism $\Omega^\bullet\tilde\lambda_{\bKU}$
is a global equivalence of orthogonal spaces.
In the commutative diagram of abelian monoids
\[ 
 \xymatrix@C=12mm{
[A,\bGr^\mC]^G  \ar[d]_{[A,\beta]^G} \ar[r]^-{ [A,c]^G} & 
[A,\Omega^\bullet\bku]^G   \ar[r]^-{ [A,\Omega^\bullet j]^G} & 
[A, \Omega^\bullet \bKU]^G\ =\ \bKU_G^0(A_+) 
\ar[d]^{ [A, \Omega^\bullet \tilde\lambda_{\bKU}]^G }_\iso \\
[A, \Omega \bU]^G 
\ar[rr]_-{  [A, \Omega( (\Omega^\bullet(\sh j))\circ\eig)]^G}^-\iso && 
[A, \Omega (\Omega^\bullet(\sh \bKU)) ]^G }
  \]
the lower horizontal and right vertical maps are thus isomorphisms
by Proposition \ref{prop:[A,Y]^G of closed} (ii). 
The left vertical map $[A,\beta]^G$ is a group completion of abelian monoids
by Corollary \ref{cor:group completion Gr to Omega U}.
So the upper horizontal composite 
$[A,(\Omega^\bullet j)\circ c]^G:[A,\bGr^\mC]^G\to \bKU_G^0(A_+)$
is also a group completion of abelian monoids. 
On the other hand,
the homomorphism $\td{-}:[A,\bGr^\mC]^G\to \bK_G(A)$
is yet another group completion of abelian monoids, by 
the complex analog of Theorem \ref{thm:BOP_to_KO_G}.
The universal property of group completions yields
an isomorphism of abelian groups
\[ \Psi \ : \ \bK_G(A)\ \to \ \bKU_G^0(A_+) \]
such that 
$\Psi\circ\td{-}=[A,(\Omega^\bullet j)\circ c]^G:[A,\bGr^\mC]^G\to \bKU_G^0(A_+)$.
The defining relation $[A, c]^G=  [-]\circ\td{-}$ 
of the homomorphism \eqref{eq:K_G_to_ku} then yields
\[  \Psi\circ\td{-} \ = \ [A, (\Omega^\bullet j)\circ c]^G\ 
= \ [A,\Omega^\bullet j]^G\circ [-]\circ\td{-} \ . \]
Since $\td{-}:[A,\bGr^\mC]^G\to\bK_G(A)$ is a group completion, this
forces the relation $\Psi=[A,\Omega^\bullet j]^G\circ [-]$.
\end{proof}

The special case $A=\ast$ of the previous corollary is worth spelling out explicitly.
In this case the group $\bK_G(\ast)$ becomes
the unitary representation ring $\bRU(G)$, and $\bKU_G^0(\ast)$ becomes
the 0-th equivariant homotopy group $\pi_0^G(\bKU)$.

\begin{theorem}\label{thm:pi_0 KU is RU}\index{subject}{representation ring!unitary} 
As $G$ ranges of all compact Lie groups, the composite maps
\[   \bRU(G)\ \xra{\ [-]\ } \ \pi_0^G(\bku)\ \xra{\ \pi_0^G(j)\ } \ \pi_0^G(\bKU)  \]
form an isomorphism of global power functors between $\bRU$ and $\upi_0(\bKU)$.
\end{theorem}
\begin{proof}
By Corollary \ref{cor:global_K_represents} for $A=\ast$ 
the composite is a ring isomorphism for every compact Lie group $G$.
The second map is induced by a homomorphism $j:\bku\to\bKU$ of ultra-commutative
ring spectra, so the maps $\pi_0^G(j)$ form a morphism of global power functors.

The maps $[-]:\bRU(G)\to\pi_0^G(\bku)$
are ring homomorphisms, compatible with restriction homomorphisms,
with finite index transfers  and multiplicative power
operations by Theorem \ref{thm:R(G) to pi ku}. 
There is still something to show, though, because the maps $[-]$
are definitely {\em not} compatible with general transfers (i.e., of infinite index). 
So an additional argument is needed to see that the composite {\em does}
commute with all transfers.

We consider a closed subgroup $H$ of a compact Lie group $G$, not necessarily of
finite index, and we want to show that
\begin{equation}  \label{eq:general_transfer_in_KU}
 \tr_H^G(j_* [x] )  \ = \ j_*[\tr_H^G(x)]  
\end{equation}
in $\pi_0^G(\bKU)$ for all classes $x\in\bRU(H)$.
Since representations are detected by characters, two classes 
in $\bRU(G)$ are equal already if their restrictions to all
finite abelian subgroups of $G$ coincide. 
Since the composite maps  $\bRU(G)\to\pi_0^G(\bKU)$ 
are all isomorphisms and compatible with restriction,
the analogous property holds for the global functor $\upi_0(\bKU)$.
So it suffices to show that the relation \eqref{eq:general_transfer_in_KU}
holds after restriction to every finite abelian subgroup $A$ of $G$.
The double coset formula
\[  \res^G_A\circ \tr_H^G \ = \ \sum_{[M]}\
\chi^\sharp(M)\cdot \tr_{A\cap{^g H}}^A  \circ g_\star \circ \res^H_{A^g\cap H}\]
holds for $\bRU$ by \cite[Thm.\,2.4]{Snaith-Brauer},
and it holds in the homotopy global functor of every orthogonal spectrum.
Since $A$ is finite, 
the right hand side of the double coset formula only involves finite index transfers. 
Since the maps from $\bRU$ to $\upi_0(\bKU)$ under consideration
do commute with finite index transfers, they commute with the right hand
side of the double coset formula, hence also with the left hand side.
This shows that \eqref{eq:general_transfer_in_KU} holds after restriction
to every finite abelian subgroup, so it holds altogether.
\end{proof}

\begin{rk}
In Construction \ref{con:rank_filtration} we discussed 
the rank filtration\index{subject}{rank filtration!of $\bku$} 
of the connective global $K$-theory spectrum $\bku$,
and we identified the first stage $\bku^{[1]}$ with the suspension spectrum
of the ultra-commutative monoid $\bP^\mC$.
On $\upi_0$, the morphisms of ultra-commutative ring spectra
\[ \Sigma^\infty_+\bP^\mC \ \iso \ \bku^{[1]} \ \xra{\text{incl}} 
\ \bku \ \xra{\ j \ } \ \bKU \]
induce morphisms of global power functors. 
Since $\bP^\mC$ is a global classifying space of the circle group $U(1)$,
the global functor $\upi_0(\Sigma^\infty_+ \bP^\mC)$
is representable by $U(1)$, by Proposition \ref{prop:B_gl represents}.  
On the other hand, $\upi_0(\bKU)$ is isomorphic to the representation
ring global functor $\bRU$, by Theorem \ref{thm:pi_0 KU is RU}.
Under these identifications, the composite morphism 
$\Sigma^\infty_+\bP^\mC \to\bKU$
of ultra-commutative ring spectra becomes the morphism
\[ \ev_x\ : \ \bA(U(1),-)\ \to \ \bRU \]
that sends the generator $1_{U(1)} \in \bA(U(1), U(1))$ to the class of the tautological
$U(1)$-representation on $\mC$. 
As we recalled in Remark \ref{rk:Brauer induction}, this morphism is surjective,
and `explicit Brauer induction' provides a specific section.
\end{rk}

\begin{eg}\label{eg:Global K is right induced}
Our language allows a reformulation of the generalization, 
due to Adams, Haeberly, Jackowski 
and May \cite{adams-haeberly-jackowski-may:atiyah-segal},
of the Atiyah-Segal completion theorem:
the global $K$-theory spectrum $\bKU$ 
is right induced from the global family $\cyc$
of finite cyclic groups.\index{symbol}{$\cyc$ - {global family of finite cyclic groups}}\index{subject}{global family!of finite cyclic groups} 

If $G$ is a compact Lie group, then a virtual $G$-representation
that restricts to zero on every finite cyclic subgroup is already zero.
In other words, the intersection of the kernels 
of all restriction maps $\res^G_C:R(G)\to R(C)$ for all
finite cyclic subgroups $C$ of $G$, is trivial. 
Then by \cite[Cor.\,2.1]{adams-haeberly-jackowski-may:atiyah-segal}, 
the projection $A\times E(\cyc\cap G)\to A$ induces an isomorphism 
\[ \bK^*_G( A )\ \iso \ \bK^*_G( A \times E(\cyc\cap G) )  \]
on equivariant $K$-groups for every finite $G$-CW-complex $A$,
where $E(\cyc\cap G)$ is a universal $G$-space for the family of finite cyclic subgroups
of $G$.
The Milnor short exact sequence lets us extend this to infinite $G$-CW-complexes,
so the criterion provided by Proposition \ref{prop:right induced criterion} 
for being right induced from the global family $\cyc$ is satisfied.
\end{eg}

Our next aim is to establish Bott periodicity for the global $K$-theory spectrum $\bKU$.

\begin{construction}[Inverse Bott class]
The Bott class $\beta\in\pi_2^e(\bku)$
was defined in Construction \ref{con:Bott class}.\index{subject}{Bott class!in connective $K$-theory}
We will now define the {\em inverse Bott class}\index{subject}{inverse Bott class}
$\lambda\in\pi_{-2}^e(\bKU)$
and show that it is multiplicatively inverse to the image of $\beta$ under
the homomorphism of ultra-commutative ring spectra $j:\bku\to\bKU$.

We recall that $u\mC$ denotes the underlying euclidean inner product space of $\mC$.
We let $e=[1]$ and $f=[i]$ be the images in $\mCl(u\mC)$
of the standard orthonormal $\mR$-basis $\{1,i\}$ of $u\mC$.
We denote by $q\in \mCl(u\mC)$ the projection
(i.e., self-adjoint idempotent)
\[ q \ = \ (1 - i\cdot e f)/2\ . \]
The isomorphism $\delta_\mC:\mCl(u\mC)\to\End_\mC(\Lambda^*\mC)$
defined in \eqref{eq:Clifford_iso}
takes the projection $q$ to the orthogonal projection onto $\Lambda^0\mC$, 
the constant summand of the exterior algebra, see Example \ref{eg:delta C}.

Moreover, we denote by $p_0\in \Kc_{u\mC}$ the orthogonal projection onto
the subspace $\mC\cdot 1$ in $\Hc_{u\mC}=\widehat{\Sym}( (u\mC)_\mC)$.
Then the element $q\tensor p_0$ of $\mCl(u\mC)\tensor \Kc_{u\mC}$
is another projection. So the map
\[ s \ \to \ \mCl(u\mC)\tensor \Kc_{u\mC} \ , \quad 
f \ \longmapsto \ f(0)\cdot q\tensor p_0 \]
is a $\mZ/2$-graded $\ast$-homomorphism, i.e., an element in the space
\[ C^\ast_{\gr}(s,\mCl(u\mC)\tensor \Kc_{u\mC}) \ = \ \bKU(u\mC)
\ = \ \bKU(\mR^2)\ . \] 
We denote by
\[ \lambda \ \in \ \pi_{-2}^e(\bKU)  \]
the homotopy class represented by this point.
\end{construction}

\begin{theorem}\label{thm:Bott and inverse}\index{subject}{Bott periodicity}
The relation
$  j_*(\beta)\cdot \lambda =  1 $
holds in $\pi_0^e(\bKU)$.
\end{theorem}
\begin{proof}
Step 1: We start with an identification of
$\mCl(\mR)\tensor M_2$ with $\mCl(\mR\oplus u\mC)$,
where $M_2$ is the $C^\ast$-algebra of $2\times 2$ complex matrices, 
concentrated in even grading.
We let $d=[1]\in \mCl(\mR)$ denote the odd unit corresponding to $1\in\mR$.
The map
\[  \Psi\left(1\tensor \left(\begin{smallmatrix} w & x \\ y &  z \end{smallmatrix}\right)\right)
\ = \ (w+z)/2\cdot 1 \ + \ (z-w)/2\cdot i\cdot e f 
\ + \ (y-x)/2\cdot d e \ - \ (x+y)/2\cdot i\cdot d f \]
identifies the even summand $1\tensor M_2$ of 
$\mCl(\mR)\tensor M_2$ 
with the even summand of $\mCl(\mR\oplus u\mC)$.
The odd element $i\cdot d e f$ of $\mCl(\mR\oplus u\mC)$ 
squares to~1 and commutes with all even elements of $\mCl(\mR\oplus u\mC)$.
So we can extend the isomorphism to the odd summands by setting
\[  \Psi\left(d\tensor \left(\begin{smallmatrix} w & x \\ y &  z \end{smallmatrix}\right)\right)
\ = \ i\cdot d e f \cdot \Psi\left(1\tensor \left(\begin{smallmatrix} w & x \\ y &  z \end{smallmatrix}\right)\right)\ .\]
The result is an isomorphism of $\mZ/2$-graded $C^\ast$-algebras
\[ \Psi \ : \ \mCl(\mR)\tensor M_2\ \xra{\iso}\ \mCl(\mR\oplus u\mC) \ . \]
While $C^\ast$-algebras are not required to have multiplicative units, 
and $\ast$-homo\-morphisms are ignorant of units, these two $C^\ast$-algebras are unital
and $\Psi$ happens to preserve multiplicative units.

Step 2: We let $j':\bk(\mC^2,S^1)\to C^\ast_{\gr}(s,\mCl(\mR)\tensor M_2)$
be the restriction of the map
\[ j(\mR)\ :\ \bku(\mR)=\bk(\Sym(\mC),S^1)\ \to\ 
C^\ast_{\gr}(s,\mCl(\mR)\tensor \Kc_\mR)=\bKU(\mR) \] 
to those configurations contained in $\mC^2\subset \Sym(\mC)$,
defined by the same formula
\[ j'[E_1,E_2;\,v_1,v_2](u) \ = \  u[v_1]\tensor p_{E_1}\ + \ u[v_2]\tensor p_{E_2}\ .\]
Our next claim is that the composite
\begin{align*}
 S^{\mR\oplus\mC} \ &\xra[\iso]{\ m\ } \ S U(2) \ \xra{\ \eig\ } \ \bk(\mC^2,S^1)\\ 
&\xra{\ j'\ } \ 
C^\ast_{\gr}(s,\mCl(\mR)\tensor M_2) \ \xra{\ \Psi_*\ }\ 
C^\ast_{\gr}(s,\mCl(\mR\oplus u\mC))
\end{align*}
is the functional calculus map 
$\fc: S^{\mR\oplus\mC} \to C^\ast_{\gr}(s,\mCl(\mR\oplus u\mC))$,
where $m$ is the homeomorphism \eqref{eq:define_m} defined by the formula
\[ m(v)\ = \ m(x,z)\ = \  
\frac{1}{|v|^2+1} \begin{pmatrix} |v|^2-1 - i 2 x & 
2 i \bar z \\ 2 i z &  |v|^2-1 + i 2 x 
 \end{pmatrix}\ .\]
The map $\eig\circ m$ sends $(x,z)$ to the configuration of Cayley transforms
of the eigenvalues of $m(x,z)$,
labeled by their eigenspaces.
So to identify the composite $j'\circ\eig\circ m$ we must
calculate these eigenvalues and eigenspaces, and their orthogonal projections.
This is a straightforward exercise in linear algebra:
a direct calculation shows that the matrices
\[ p_+ \ = \ 
\frac{1}{2|v|} \begin{pmatrix} |v|- x & \bar z \\ z &  |v| + x  \end{pmatrix}
\text{\quad and\quad}
 p_- \ = \ 
\frac{1}{2|v|} \begin{pmatrix} |v|+ x & -\bar z \\ - z &  |v|- x  \end{pmatrix}\]
are projections and orthogonal to each other, i.e., they satisfy the relations
\[ p_+^2\ = \ p_+^\ast \ = \ p_+ \ , \quad 
 p_-^2\ = \ p_-^\ast \ = \ p_- \text{\quad and\quad}
p_+\cdot p_-\ = \ p_-\cdot p_+ \ = \ 0 \ . \]
Moreover,
\[ c(|v|)\cdot p_+ \ + \ c(-|v|)\cdot p_- \ = \ m(x,z) \ ,\]
so $p_+$ and $p_-$ are the orthogonal projections onto the eigenspaces
of $m(x,z)$, with eigenvalues 
\[ c(|v|)\ = \ (|v|+i)(|v|-i)^{-1} \text{\quad respectively\quad}
 c(-|v|)\ = \ c(|v|)^{-1}\ = \ (|v|-i)(|v|+i)^{-1} \ .\]
The map $j'\circ \eig\circ m$ thus sends $v=(x,z)$ to the $\ast$-homomorphism
\[  (j'\circ \eig\circ m)(x,z)(u)\ =  u[|v|]\tensor p_+ \ + \ u[-|v|]\tensor p_-\ . \]
The $C^\ast$-algebra $s$ is generated by the function $r(x)=2 i(x-i)^{-1}$.
So to verify that two graded $\ast$-homomorphisms with source $s$ agree,
it suffices to show that they coincide on the even and odd 
components of the generator $r$,
which we spelled out explicitly in \eqref{eq:even-odd_of_r}. 
For the even component $r_+$ the calculation is
\begin{align*}
  (\Psi_*\circ j'\circ\eig\circ m)(x,z)(r_+)\ &= \ 
  \Psi\big( r_+[|v|]\tensor p_+ + r_+[-|v|]\tensor p_-\big) \\ 
&= \   \Psi\big( r_+[|v|]\tensor (p_+ +  p_-)\big) \ 
= \    r_+[|v|]\cdot 1\ = \    r_+[x,z]\ .
\end{align*}
For the odd component $r_-$ we first observe that
\begin{align*}
 \Psi\left( d\tensor  
    \left(\begin{smallmatrix}  - x & \bar z \\  z &    x  \end{smallmatrix}\right)\right) \ &= \   
i\cdot d e f \cdot (  x\cdot i e f  + (z-\bar z)/2 \cdot d e
-(z+\bar z )/2\cdot i\cdot d f  )\\
&= \ 
- x\cdot  d e f e f\ -\ \text{Im}(z)\cdot  d e f d e\ + \ \text{Re}(z) \cdot  d e f d f\\
&= \ 
 x\cdot  d + \text{Im}(z)\cdot  f\ + \text{Re}(z)\cdot e\ = \ [x,z]
\end{align*}
in $\mCl(\mR\oplus u\mC)$. So
\begin{align*}
  (\Psi_*\circ j'\circ \eig\circ m)(x,z)(r_-)\ &= \ 
  \Psi\big( r_-[|v|]\tensor p_+ + r_-[-|v|]\tensor p_-\big) \\ 
&= \   \Psi\big( r_-[|v|]\tensor (p_+ -  p_-)\big) \\ 
&= \   \Psi\left( \frac{2 i [|v|]}{|v|^2+1}\tensor \frac{1}{|v|} \begin{pmatrix}-x &\bar z\\ z &  x  \end{pmatrix}\right) \\ 
&= \   \frac{2 i}{|v|^2+1}\cdot \Psi\left( d\tensor \begin{pmatrix}-x &\bar z\\ z & x  \end{pmatrix}\right) \\ 
&= \   \frac{2 i\cdot [x,z]}{|v|^2+1}\ = \  r_-[x,z]\ .
\end{align*}
This completes the verification that $\Psi_*\circ j'\circ\eig\circ m$
is the functional calculus map $\fc: S^{\mR\oplus \mC} \to
C^\ast_{\gr}(s,\mCl(\mR\oplus u\mC))$.

Step 3:
We recall that $q=(1 - i e f)/2$ in the Clifford algebra $\mCl(\mR\oplus u\mC)$, 
which satisfies the relations
\begin{equation}  \label{eq:relations_for_q}
  d q\ = \ q d  \ , \quad  e q   = (1-q)e
\text{\qquad and\qquad} q e q  = 0 \ .  
\end{equation}
We consider the even element $u$ of $\mCl(\mR\oplus u\mC)\tensor M_2$
defined by
\begin{align*}
 u \ = \ q\tensor\left(\begin{smallmatrix}1&0\\0&0\end{smallmatrix}\right)
\ +\ q e d\tensor\left(\begin{smallmatrix}0&0\\1&0\end{smallmatrix}\right)
\ + \    e d q\tensor\left(\begin{smallmatrix}0&1\\0&0\end{smallmatrix}\right)
\ + \  (1-q)\tensor\left(\begin{smallmatrix}0&0\\0&1\end{smallmatrix}\right) \ .
\end{align*}
A direct calculation, making repeated use of the identities \eqref{eq:relations_for_q},
shows that this is a unitary element, i.e.,
\[ u^\ast\cdot u\ = \ u\cdot  u^\ast \ = \  1\ . \]
So conjugation by $u$ is a graded $\ast$-automorphism
\[ \zeta \ : \ \mCl(\mR\oplus u\mC)\tensor M_2 \ \to \ 
\mCl(\mR\oplus u\mC)\tensor M_2 \ , \quad x \ \longmapsto \ u^\ast\cdot x\cdot u \ . \]
Now we observe that the following square of graded $C^\ast$-algebras commutes:
\[ \xymatrix@C=13mm{ 
 \mCl(\mR)\tensor M_2\ar[d]_-{\Psi} \ar[r]^-{-\tensor q\tensor -} &
 \mCl(\mR)\tensor\mCl(u\mC)\tensor M_2
\ar[r]^-{\eqref{eq:Clifford tensor iso}\tensor M_2}&
\mCl(\mR\oplus u\mC)\tensor M_2\ar[d]^\zeta\\
\mCl(\mR\oplus u\mC)\ar[rr]_-{ - \tensor \left(\begin{smallmatrix}1&0\\0&0\end{smallmatrix}\right)} &&
\mCl(\mR\oplus u\mC)\tensor M_2} \]
This, again, is a direct calculation, using the identities \eqref{eq:relations_for_q}.
We claim that the conjugation map $\zeta$
is homotopic, through graded $\ast$-homomorphisms, to the identity
of $\mCl(\mR\oplus u\mC)\tensor M_2$.
To show this we define a continuous path
\[ u\ : \ [0,\pi/2] \ \to \ \mCl(\mR\oplus u\mC)\tensor M_2\]
by the formula
\begin{align*}
 u(t) \ = \ (q+\sin(t)&(1-q))\tensor\left(\begin{smallmatrix}1&0\\0&0\end{smallmatrix}\right)
\ +\ \cos(t)\cdot q e d \tensor\left(\begin{smallmatrix}0&0\\1&0\end{smallmatrix}\right)\\
& + \  \cos(t)\cdot e d q\tensor\left(\begin{smallmatrix}0&1\\0&0\end{smallmatrix}\right)
\ + \ 
(\sin(t)\cdot q + (1-q))\tensor\left(\begin{smallmatrix}0&0\\0&1\end{smallmatrix}\right) \ .
\end{align*}
One more direct calculation shows that all elements in this path are even unitaries.
So conjugation by the elements $u(t)$ is a path of graded $\ast$-automorphism.
Since $u(0)=u$, this path starts with the conjugation map $\zeta$;
since $u(\pi/2)=1$, the path ends with the identity.
Since $\zeta$ is homotopic to the identity, the induced self-map
of $C^\ast_{\gr}(s,\mCl(\mR\oplus u\mC)\tensor M_2)$
is based homotopic to the identity.

Step 4: The class $\lambda$ is represented by the point of $\bKU(u\mC)$
given by the graded $\ast$-homomorphism
\[ s \ \xra{\ \epsilon\ } \ \mC \ \xra{\cdot q\tensor p_0} \ 
\mCl(u\mC)\tensor\Kc_{u\mC}\]
associated to the projection $q\tensor p_0$; here $\epsilon(\varphi)=\varphi(0)$
is the augmentation. So the multiplication map
$ -\cdot \lambda :  \pi_{k+2}^e(\bKU)\to \pi_k^e(\bKU)$
is the effect on homotopy groups of the based continuous map
\begin{align*}
\bKU(V)\ = \  &C^\ast_{\gr}(s,\mCl(V)\tensor\Kc_V)\\
&\xra{-\sm (-\cdot q\tensor p_0)}  \ 
 C^\ast_{\gr}(s,\mCl(V)\tensor\Kc_V)\sm C^\ast_{\gr}(\mC,\mCl(u\mC)\tensor \Kc_{u\mC})  \\ 
&\xra{\Id\sm \epsilon^*} \
 C^\ast_{\gr}(s,\mCl(V)\tensor\Kc_V)\sm C^\ast_{\gr}(s,\mCl(u\mC)\tensor\Kc_{u\mC})  \\ 
&\xra{\mu^\bKU_{V,u\mC}} \
 C^\ast_{\gr}(s,\mCl(V\oplus u\mC)\tensor\Kc_{V\oplus u\mC}) \ = \ \bKU(V\oplus u\mC)\ .
\end{align*}
The composite $s\hat\tensor\epsilon\circ\Delta$ is the identity of $s$, 
so multiplication by $\lambda$ is the effect of postcomposition
with the graded $\ast$-homomorphism
\begin{align*}
 \mCl(V)\tensor\Kc_V\ &\xra{-\tensor q\tensor -\tensor p_0}\ 
\mCl(V)\tensor\mCl(u\mC)\tensor\Kc_V\hat\tensor\Kc_{u\mC}\\
&\xra[\quad \iso\quad]{\eqref{eq:Clifford tensor iso}\tensor\eqref{eq:tensor compacts}}\ 
\mCl(V\oplus u\mC)\tensor\Kc_{V\oplus u\mC}\ .  
\end{align*}

Now we put the pieces together and prove the theorem.
We contemplate the following diagram of continuous based maps:
\[ \xymatrix@C=5mm@R=8mm{ 
 S^{\mR\oplus\mC} \ar[d]_{\eig\circ m} \ar[rr]^-{\fc} && 
C^\ast_{\gr}(s,\mCl(\mR\tensor u\mC)) 
\ar[d]^-{\left(-\tensor \left(\begin{smallmatrix}1&0\\0&0\end{smallmatrix}\right)\right)_*}
\ar@<8ex>@/^2pc/[dd]^-(.2){(-\tensor p_0)_*} \\
 \bk(\mC^2,S^1)\ar[d]_-{\bk(i,S^1)}\ar[r]_-{j'} &
C^\ast_{\gr}(s,\mCl(\mR)\tensor M_2)\ar[r]_-{(-\tensor q\tensor -)_*}
\ar[ur]^(.4){\Psi_*}\ar[d]^-{(-\tensor i_!)_*} &
 C^\ast_{\gr}(s,\mCl(\mR\tensor u\mC)\tensor M_2) \ar[d]^-{(-\tensor i_!)_*} \\
\bk(\Sym(\mC),S^1)\ar@{=}[dd] &
C^\ast_{\gr}(s,\mCl(\mR)\tensor\Kc_\mR) \ar@{=}[dd] \ar[r]_-{(-\tensor q\tensor -)_*} & 
 C^\ast_{\gr}(s,\mCl(\mR\tensor u\mC)\tensor \Kc_{\mR}) 
\ar[d]^{(-\tensor p_0)_*} \\
&&  C^\ast_{\gr}(s,\mCl(\mR\tensor u\mC)\tensor \Kc_{\mR\oplus u\mC}) \ar@{=}[d] \\
\bku(\mR)\ar[r]_-{j(\mR)} &\bKU(\mR)\ar[r]_-{(-\tensor q\tensor -\tensor p_0)_*} &
 \bKU(\mR\oplus u\mC)
} \]
Here $i:\mC^2\to\Sym(\mC)$ is the identification with the
constant and linear summands of the symmetric algebra used throughout,
and $i_!:M_2\to\Kc_\mR$ is the $\ast$-homomorphism given by extension 
by~0 on the higher symmetric powers.
By Step~2, the upper left part commutes on the nose,
and by Step~3 the upper right triangle commutes up to based homotopy.
The lower left part commutes because $j'$ is the restriction of $j(\mR)$. 
The lower right part commutes as well.
So the whole diagram commutes up to based homotopy.
By Step~4, the composite
through the lower left corner represents the class $j_*(\beta)\cdot\lambda$,
whereas the composite along the right side of the diagram is
the unit map $\eta_{\mR\oplus\mC}:S^{\mR\oplus\mC}\to\bKU(\mR\oplus u\mC)$
of the ring spectrum structure. So the diagram witnesses the desired relation
\[ j_*(\beta)\cdot\lambda \ = \ 
[(-\tensor q\tensor - \tensor p_0)_*\circ j(\mR)\circ \bk(\Sym(\mC),S^1)\circ\eig\circ m]
\ = \ [\eta_{\mR\oplus u\mC}] \ = \ 1 \ .\qedhere \]
\end{proof}

\begin{rk}[Thom isomorphism in equivariant $K$-theory]\label{rk:equivariant inverse Bott}
The spectrum $\bKU$ enjoys a much stronger form of periodicity that generalizes
the $\mZ$-graded periodicity manifested by
Theorem \ref{thm:Bott and inverse}. 
The equivariant Bott class $\beta_{G,W}$ associated
to a $G$-$\Spin^c$-representation $W$ was introduced 
in Construction \ref{con:equivariant Bott class},
and $\beta_{G,W}$ is in fact an `$R O(G)$-graded unit' in $\bKU$. 
This kind of equivariant periodicity in $K$-theory goes back to 
Atiyah \cite[Thm.\,4.3]{atiyah-Bott and elliptic};\index{subject}{Bott periodicity!equivariant} 
it is formally similar to the periodicity of $\bMO$ and $\bMOP$
for orthogonal representations (compare Theorem \ref{thm:j^V is upi_*-iso}~(iii)),
or the periodicity of $\bMU$ and $\bMUP$ for unitary representations. 

Even more is true. The equivariant Bott class is a special case of
the Thom class of an equivariant $\Spin^c$-vector bundle,\index{subject}{Thom class!in $\bKU$}
and equivariant Bott periodicity is a special case of a Thom isomorphism 
for such bundles.\index{subject}{Thom isomorphism!for $\bKU$}
This Thom isomorphism for equivariant $\bKU$-theory goes under the name 
`Atiyah-Bott-Shapiro orientation'\index{subject}{Atiyah-Bott-Shapiro orientation} 
and was established in \cite{atiyah-bott-shapiro}.\index{subject}{G-spin c-bundle@$(G,\Spin^c)$-bundle}\index{subject}{spin$^c$ group}\index{subject}{spin$^c$-bundle!equivariant}
In \cite[Thm.\,6.9]{joachim-coherences} 
Joachim shows that the Atiyah-Bott-Shapiro orientation
is in fact extremely highly structured, name\-ly `globally and ultra-commutative'.
Joachim defines a morphism of ultra-commutative ring spectra
$\alpha:\mathbb M \Spin^c\to\bKU$ from a global equivariant version of the
$\Spin^c$-Thom spectrum; his morphism refines the Atiyah-Bott-Shapiro orientation,
in the sense that it takes certain tautological Thom classes 
for $\Spin^c$-vector bundles in $\mathbb M \Spin^c$
to the $\bKU$-theoretic Thom classes of Atiyah-Bott-Shapiro. 
\end{rk}

\index{subject}{K-theory@$K$-theory!periodic global|)}

\index{subject}{K-theory@$K$-theory!global connective|(}

\begin{construction}[Global connective $K$-theory]\label{con:global connective ku}
Now we define {\em global connective $K$-theory} $\bkuc$,
an ultra-commutative ring spectrum whose associated $G$-homotopy type,
for compact Lie groups $G$,
is that of {\em $G$-equivariant connective $K$-theory}
in the sense of Greenlees \cite{greenlees-connective}.
This is {\em not} a connective equivariant theory, i.e.,
the equivariant homotopy groups $\pi_*^G ( \bkuc )$ do {\em not}
vanish in negative dimensions, as soon as the group $G$ is non-trivial.
Hence the order of the adjectives `global' and `connective' matters, i.e.,
`global connective' $K$-theory is different from `connective global' $K$-theory.
Two of the advantages of $\bkuc$ are that it is equivariantly 
(and in fact globally) orientable,
and that $\bkuc$ satisfies a completion theorem: for every
compact Lie group $G$, the completion of the graded ring $\pi_*^G( \bkuc )$ 
at the augmentation ideal of the unitary representation ring\index{subject}{augmentation ideal!of the unitary representation ring}
is the connective $\bku$-cohomology 
of the classifying space $B G$, see \cite[Prop.\,2.2 (i)]{greenlees-connective}.

Our construction of $\bkuc$ is a direct `globalization' of Greenlees' 
definition in \cite[Def.\,3.1]{greenlees-connective}.
We define $\bkuc$\index{symbol}{$\bkuc$ - {global connective $K$-theory}}\index{subject}{K-theory@$K$-theory!global connective} 
as the homotopy pullback in the square of ultra-commutative ring spectra
\begin{equation}\begin{aligned}\label{eq:define ku^c}
\xymatrix{\bkuc\ar[r]^-{}\ar[d] &
b(\bku)\ar[d]^{b j} \\
\bKU\ar[r]_-{i_{\bKU}} & b(\bKU)  } 
\end{aligned}\end{equation}
The morphism $j:\bku\to\bKU$ from connective to periodic global $K$-theory
was defined in Construction \ref{con:global KU};
the Borel theory functor $b$ and the natural transformation
$i:\Id\to b$ were defined in Construction \ref{cor:Borel b}.
So more explicitly, we set 
\[ \bkuc \ = \ \bKU\times_{b(\bKU)} b(\bKU)^{[0,1]}\times_{b(\bKU)} b(\bku) \ . \]
Since the spectra $\bKU$, $b(\bKU)$ and $b(\bku)$
are ultra-commutative ring spectra and the two morphisms
$i_{\bKU}:\bKU\to b(\bKU)$ and $b j:b(\bku)\to b(\bKU)$
are homomorphisms, the homotopy pullback
is canonically an ultra-commutative ring spectrum
and the two morphisms from $\bkuc$ to $\bKU$ and $b(\bku)$
are morphisms of ultra-commutative ring spectra.
As a homotopy pullback, the square \eqref{eq:define ku^c}
does {\em not} commute, but the construction comes with
a preferred homotopy between the two composites around the square.

The morphisms $j:\bku\to \bKU$, $i_{\bku}:\bku\to b(\bku)$
and the constant homotopy provide a morphism of ultra-commutative ring spectra
\[ \iota\ : \ \bku\ \to \ \bkuc \]
from connective global $K$-theory to global connective $K$-theory.
For finite groups, this morphism induces an isomorphism on homotopy global functors
in non-negative dimensions.
The construction of $\bkuc$ endows it with a morphism of ultra-commutative ring spectra
$\bkuc\to\bKU$. As explained by Greenlees 
in \cite[Thm.\,2.1 (iv)]{greenlees-connective},
this morphism becomes a global equivalence after inverting
the Bott class $v\in \pi_2^e( \bkuc )$, 
the image of the Bott class $\beta\in\pi_2^e(\bku)$
defined in Construction \ref{con:Bott class}.\index{subject}{Bott class!in connective $K$-theory}
Indeed, for every compact Lie group $G$, the morphism
$j:\bku\to\bKU$ induces a morphism of graded rings
\[ j^*(B G)\ : \ \bku^*(B G) \ \to \ \bKU^*( B G)\]
that becomes an isomorphism after inverting the Bott class
by \cite[Lemma~1.1.1]{bruner-greenlees-connective_K_finite}.
The natural isomorphisms
\[ \pi_{-*}^G(b(\bku)) \ \xra{\ \iso \ } \  \bku^*(B G) \text{\quad and\quad}
 \pi_{-*}^G(b (\bKU)) \ \xra{\ \iso \ }\  \bKU^*(B G) \]
of Proposition \ref{prop:global homotopy of b E}
then show that the map
\[ \upi_*(b j) \ : \ \upi_*(b(\bku)) \ \to \  \upi_*(b (\bKU)) \]
becomes an isomorphism of graded global functors after inverting the Bott class.
Since the right vertical morphism in the 
defining global homotopy pullback \eqref{eq:define ku^c}
becomes a global equivalence after inverting the Bott class, the same is true for
the left vertical morphism.
\end{construction}

The composite
\[ \dim^c \ : \ \bkuc \ \to \ b(\bku) \ \xra{ b(\dim)} \  b(\Hc\mZ) \]
is a morphism of ultra-commutative ring spectra,
where $\Hc\mZ$ is the Eilenberg-Mac\,Lane spectrum of the integers 
(see Construction \ref{con:HM})\index{subject}{Eilenberg-Mac\,Lane spectrum!of the integers} 
and the dimension homomorphism $\dim:\bku\to \Hc \mZ$
was defined in Example \ref{eg:dimension hom}.\index{subject}{dimension homomorphism} 
The dimension morphism annihilates the Bott class $\beta\in\pi_2^e(\bku)$.
We claim that the sequence
\begin{equation}\label{eq:global_Bott_fiber_sequence}
 \bkuc\sm S^2\ \xra{\ \tilde v \ } \ \bkuc\ \xra{\dim^c} \  b(\Hc\mZ)
\end{equation}
is a `global homotopy cofiber sequence', i.e., part of a distinguished
triangle in the global stable homotopy category, where $\tilde v$
is the extension of the Bott class $v\in\pi_2^e(\bkuc)$ 
to a morphism of $\bkuc$-module spectra.
Indeed, the sequence 
\[ \bku\sm S^2 \ \xra{\ \tilde\beta \ } \ \bku\ \xra{\dim} \ \Hc\mZ \]
is a non-equivariant homotopy fiber sequence;
the Borel theory functor $b$ takes this to the global homotopy fiber sequence
\[ b(\bku)\sm S^2 \ \xra{b(\tilde\beta)} \  b(\bku)
\ \xra{ b(\dim)} \  b(\Hc\mZ) \ .\]
The spectrum $\bKU$ is Bott periodic, i.e., the image of the Bott class
in $\pi_2^e(\bKU)$ is invertible by Theorem \ref{thm:Bott and inverse}.
Since the morphism $\bKU\to b(\bKU)$ is multiplicative, the same goes
for the Borel theory $b(\bKU)$. So the morphisms
\[
\tilde\beta\ : \  \bKU\sm S^2 \ \to\ \bKU 
\text{\quad and\quad}
 b(\tilde\beta)\ : \ b(\bKU)\sm S^2 \ \to\ b(\bKU)
 \]
are global equivalences, and their respective homotopy fibers 
are globally stably contractible.
Passing to homotopy pullbacks gives the desired 
global homotopy fiber sequence \eqref{eq:global_Bott_fiber_sequence}.

The global homotopy fiber sequence \eqref{eq:global_Bott_fiber_sequence}
and the isomorphism
\[ \pi_k^G( b(\Hc\mZ))\ \iso \ H^{-k}(B G, \mZ) \]
of Proposition \ref{prop:global homotopy of b E}~(ii)
give rise to a long exact sequence of global functors
\[ \cdots  \to  \upi_{k+1} ( b(\Hc\mZ) )  \xra{\ \partial\ }  
\upi_{k-2}(\bkuc)  \xra{\ \cdot v \ }  
\upi_k(\bkuc)  \xra{\ \dim^c_*\ }  \upi_k ( b(\Hc\mZ) )
  \to  \cdots\ .  \]
This way Greenlees calculates some of the homotopy group global functors
of equivariant connective $K$-theory in \cite[Prop.\,2.6]{greenlees-connective}.
We review Greenlees' calculations in our language.
The group $H^{-k}(B G,\mZ)$ vanishes for $k>0$.
So multiplication by the Bott class
\[ -\cdot v \ :\ \upi_{k-2}(\bkuc) \ \to\ \upi_k(\bkuc) \]
is an isomorphism for $k>0$ and a monomorphism for $k=0$. 
In particular, we conclude that
\[ \upi_k(\bkuc) \ \iso \
\begin{cases}
  \bRU & \text{\ for $k\geq 0$ and $k$ even, and}\\
  \, 0  & \text{\ for $k\geq -1$ and $k$ odd.}
\end{cases} \]
More precisely, for every $m\geq 0$, the composite
\[ \bRU \ \xra{[-]}\ \upi_0(\bku)\ \xra{\upi_0(\iota)}
\ \upi_0(\bku^c)\ \xra{\ \cdot v^m}\ \ \upi_{2m}(\bku^c)\]
is an isomorphism of global functors.

The global functor $\upi_0(b(\Hc\mZ))$ sends $G$ to $\pi_0^G(b(\Hc\mZ))=H^0(B G,\mZ)$; 
so $\upi_0(b(\Hc\mZ))$ is constant with value $\mZ$ and the morphism 
$\dim^c_*:\upi_0(\bkuc) \to \upi_0 ( b(\Hc\mZ) )$
is isomorphic to the augmentation morphism
$\dim:\bRU\to\underline{\mZ}$ of global functors, by Proposition \ref{prop:dim is dim}.
This is an epimorphism, so the sequence of global functors
\[ 0 \ \to \  \upi_{-2}(\bkuc) \ \xra{\ \cdot v \ } \ 
\upi_0(\bkuc) \ \xra{\ \dim^c_*\ } \ \upi_0 ( b(\Hc\mZ) )\ \to \ 0  \]
is short exact and
\[ \upi_{-2}(\bkuc) \ \iso \ \mathbf{I U}\  =\ \ker(\dim:\bRU \to \underline{\mZ})\]
is the augmentation ideal global functor.\index{subject}{augmentation ideal!of the unitary representation ring}
Again since the map  
$\dim^c_*:\upi_0(\bkuc) \to \upi_0 ( b(\Hc\mZ) )$
is surjective, the global functor $\upi_{-3}(\bkuc)$
injects into $\upi_{-1}(\bkuc)$, which is trivial by the above. 
So we conclude that $\upi_{-3}(\bkuc)  =  0$.

This method can be pushed a little further to also determine 
the global functors $\pi_{-4}(\bkuc)$ and $\pi_{-5}(\bkuc)$;
we refer to \cite[Prop.\,2.6]{greenlees-connective} for the argument.
The result is that
\[ \upi_{-4}(\bkuc) \ \iso \ 
\mathbf{I S U}(G) \ = \ \{x\in \mathbf{I U}(G)\ |\ \det(x)=0\} \]
and that $\upi_{-5}(\bkuc)=0$.
After this point things become less explicit.

\index{subject}{K-theory@$K$-theory!global connective|)}
\index{subject}{ultra-commutative ring spectrum|)}

\backmatter
\appendix

\chapter{Compactly generated spaces}
\label{app:CGWH}

In this appendix we recall some background material about compactly
generated spaces, our basic category to work in.
Compactly generated spaces are in particular `$k$-spaces',
a notion that seems to go back to Kelley's book \cite[p.\ 230]{kelley-general topology}.
Compactly generated spaces were
popularized by Steenrod in his paper \cite{steenrod-convenient}
as a `convenient category of spaces'; however, in contrast to our usage of
the term, Steenrod includes the Hausdorff property in his definition
of `compactly generated'.
But Steenrod already writes that `(...) {\em The Hausdorff property is
imposed to ensure that compact subsets are closed.} (...)';
the {\em weak} Hausdorff condition it thus the next logical step,
as it isolates this relevant property.
The weak Hausdorff condition first appears in print in 
McCord's paper \cite{mccord}, who credits the idea to J.\,C.\,Moore. 
McCord also was the first to use the terminology
`compactly generated spaces' for weak Hausdorff $k$-spaces.
The fact that much of the recent literature in equivariant and stable
homotopy theory uses compactly generated spaces
can be taken as evidence that the weak Hausdorff condition 
is even more convenient than the actual Hausdorff separation property.
Section~7.9 of tom Dieck's textbook \cite{tomDieck-algebraic topology}
is a nice summary of compactly generated spaces (where these are called
$whk$-spaces), and contains most of the material that we discuss here.
Two influential -- but unpublished -- sources about compactly generated
spaces are the Appendix~A of Gaunce Lewis's thesis \cite{lewis-thesis}
and Neil Strickland's preprint \cite{strickland-CGWH}.

I want to emphasize that this appendix does not contain new mathematics
and makes no claim to originality.
I have decided to include it because I found it cumbersome to
collect proofs of all the relevant properties of compactly generated spaces 
from the scattered literature. An additional complication stems
from the fact that the basic 
references \cite{kelley-general topology, mccord, steenrod-convenient}
all work in slightly different categories; so in the interest of a self-contained
treatment, I felt obliged to fill in arguments 
where a reference confines itself to the statement that 
`(...) the proof is analogous to that of (...)'.

\medskip 

We fix some terminology. A topological space is {\em compact}\index{subject}{compact!topological space}
if it is quasi-compact (i.e., every open cover has a finite subcover)
and satisfies the Hausdorff separation property (i.e., every pair of distinct
points can be separated by disjoint open subsets).

\begin{defn}\label{def:cgwh} Let $X$ be a topological space.
  \begin{itemize}
  \item A subset $A$ of $X$ is {\em compactly closed}\index{subject}{compactly closed} 
    if for every compact space $K$ and every continuous map $f:K\to X$, 
    the inverse image $f^{-1}(A)$ is closed in $K$.
  \item $X$ is a {\em $k$-space}\index{subject}{k-space@$k$-space} 
    if every compactly closed subset is closed.
  \item $X$ is {\em weak Hausdorff}\index{subject}{weak Hausdorff space}
    \index{subject}{Hausdorff space!weak|see{weak Hausdorff space}} 
    if for every compact space $K$ and
    every continuous map $f:K\to X$ the image $f(K)$ is closed in $X$.
  \item $X$ is a {\em compactly generated}\index{subject}{compactly generated!topological space} space if it is a $k$-space and weak Hausdorff.
\end{itemize}
\end{defn}

Every closed subset is also compactly closed. One can similarly define
{\em compactly open}\index{subject}{compactly open}  subsets of $X$ by demanding 
that for every compact space $K$ and every continuous map $f:K\to X$, 
the inverse image is open in $K$. A subset is then compactly open 
if and only if its complement is compactly closed.
Thus $k$-spaces can equivalently be defined by the property 
that all compactly open subsets are open.

We denote by $\bSpc$\index{symbol}{$\bSpc$ - {category of topological spaces}}
the category of topological spaces and
continuous maps. We denote by $\bK$\index{symbol}{$\bK$ - {category of $k$-spaces}} 
respectively $\bT$\index{symbol}{$\bT$ - {category of compactly generated spaces}} 
the full subcategories of $\bSpc$  
consisting of $k$-spaces respectively compactly generated spaces.
So we have full embeddings
\[ \bT \ \subset \ \bK \ \subset \ \bSpc \ .  \]
We will see below that the inclusion $\bK\subset \bSpc$
has a right adjoint `Kelleyfication' $k:\bK\to\bSpc$
and the inclusion $\bT\subset \bK$ has a left adjoint  $w:\bT\to\bK$.

\medskip

\Danger 
We follow the terminology introduced by McCord in \cite{mccord}.
We warn the reader that the usage of the terms `$k$-space' and
 `compactly generated' is not consistent throughout the literature. 
For example, some authors define
$k$-spaces by the property that a subset is closed if and only its intersection
with every compact subspace is closed. For Hausdorff spaces,
and more generally weak Hausdorff spaces, that definition agrees with the
one in Definition \ref{def:cgwh} because for such spaces 
the image of a compact space under a continuous map is automatically compact,
by Proposition \ref{prop:properties wH spaces}~(v) below.
Moreover, there are sources (for example \cite{steenrod-convenient})
that require a compactly generated space
to be Hausdorff (as opposed to only weak Hausdorff);
on the other hand, several recent references do not include the weak Hausdorff condition
in `compactly generated'.

\medskip

We recall various useful properties of $k$-spaces, 
weak Hausdorff spaces and compactly generated spaces.
We start by discussing $k$-spaces.
A Hausdorff topological space $X$ 
is {\em locally compact}\index{subject}{locally compact!topological space}
if the following condition holds: for every open subset $U$ of $X$
and every point $x\in U$, there is an open neighborhood $V$ of $x$
whose closure $\bar V$ is compact and contained in $U$.
For example, every compact space is in particular locally compact.
A topological space is {\em first countable}
if every point has a countable basis of neighborhoods.

Because the product of two $k$-spaces need not be a $k$-space in
the usual product topology, we need notation for two different kinds 
of product topologies.
We shall denote by $X\times_0 Y$ the cartesian
product of two spaces $X$ and $Y$, endowed with the product topology
(which makes it a categorical product in the category $\bSpc$ of all topological spaces).
So a basis of the topology $X\times_0 Y$ is given by products of open subsets
in the two factors. We reserve the symbol $X\times Y$ 
for $k(X\times_0 Y)$, the Kelleyfication of the product topology, 
discussed in more detail below. This product $X\times Y$
is a categorical product in the categories $\mathbf K$ and $\bT$.
Part~(vi) of the following proposition says that these two product topologies
coincide if one factor is a $k$-space and the other factor is locally compact Hausdorff.

\begin{prop}\label{prop:properties k-spaces}
  \begin{enumerate}[\em (i)]
  \item Every quotient space of a $k$-space is a $k$-space.
  \item Every closed subset of a $k$-space is a $k$-space in the subspace topology.
  \item Every locally compact Hausdorff space, and hence every compact space,
    is a $k$-space.
  \item Every first countable space is a $k$-space.
  \item Every metric space is first countable, and hence a $k$-space.
  \item If $X$ is a $k$-space and $Y$ a locally compact Hausdorff space, then
    the product $X\times_0 Y$ is a $k$-space in the product topology.
  \end{enumerate}
\end{prop}
\begin{proof}
  (i) Let $p:X\to Y$ be a quotient projection and $B$ a compactly closed subset of $Y$.
  We claim that $p^{-1}(B)$ is compactly closed in $X$. Indeed, if $f:K\to X$
  is a continuous map from a compact space, then $p f:K\to Y$ is continuous, so
  $f^{-1}(p^{-1}(B))=(p f)^{-1}(B)$ is closed because $B$ is compactly closed.
  Since $p^{-1}(B)$ is compactly closed and $X$ is a $k$-space, the set $p^{-1}(B)$ 
  is in fact closed in $X$. So $B$ is closed in $Y$ by definition of
  the quotient topology.

  (ii) We let $Y$ be a closed subset of a $k$-space $X$,
  and $A$ a compactly closed subset of $Y$ with respect to the subspace topology.
  We claim that then $A$ is also compactly closed as a subset of $X$. 
  Indeed, if $f:K\to X$ is a continuous map from a compact space, then
  $L=f^{-1}(Y)$ is closed in $K$, and hence compact. The restriction
  $f|_L:L\to Y$ is then continuous and so $(f|_L)^{-1}(A)=f^{-1}(A)$ is
  closed in $L$ because $A$ was assumed to be compactly closed in $Y$.
  Since $L$ is closed in $K$, the set $f^{-1}(A)$ is closed in $K$.
  This shows that $A$ is compactly closed as a subset of $X$.
  Since $X$ is a $k$-space, $A$ is closed in $X$. But then $A$ is also closed in $Y$
  in the subspace topology, so this concludes the proof that $Y$ is a $k$-space.

  (iii) 
  The argument goes all the way back to 
  Kelley \cite[Ch.\,7, Thm.\,13]{kelley-general topology}.
  We let $A$ be a compactly closed subset of a locally compact Hausdorff space $X$. 
  We let $\bar A$ be the closure of $A$ in $X$ and $x\in\bar A$.
  Since $X$ is locally compact, the point $x$ has a compact neighborhood $K$.
  Since $A$ is compactly closed and the inclusion $K\to X$ is continuous,
  the set $K\cap A$ is closed inside $K$. Since $K$ is compact and $X$ is Hausdorff,
  $K$ is closed in $X$. So $K\cap A$ is closed in $X$.

  Now we claim that $x\in K\cap A$.
  We argue by contraction and suppose that $x\not\in K\cap A$.
  Then $X-(K\cap A)$ is an open neighborhood of $x$, 
  and hence $K\cap(X-(K\cap A))=K\cap(X-A)$
  is another neighborhood of $x$.
  Let $U$ be an open subset of $X$ with $x\in U\subset K\cap(X-A)$. 
  Then $A\subset X-U$, and hence $x\in\bar A\subset X-U$ because $X-U$ is closed.
  But his contradicts the hypothesis $x\in U$.
  Altogether this proves the claim that $x\in K\cap A\subset A$.
  Hence $A=\bar A$, and so $A$ is closed.

  (iv) 
  This, too, goes back to Kelley \cite[Ch.\,7, Thm.\,13]{kelley-general topology}.
  We let $X$ be a first countable space and $A$ a compactly closed subset of $X$.
  We let $z\in\bar A$ be a point in the closure of $A$. 
  The point $z$ has a countable basis of open neighborhoods $\{U_n\}_{n\geq 1}$,
  which we can moreover take to be nested, i.e.,
  \[ U_1 \ \supset \ U_2 \ \supset \ \cdots \ \supset \ U_n \
  \supset \ \cdots \ .\]
  If the intersection of $U_n$ and $A$ were empty, then
  $\bar A\subset X-U_n$ which contradicts the fact that $z\in \bar A\cap U_n$. 
  So for every $n\geq 1$ there is a point $x_n\in U_n\cap A$.
  We define a map
  \[ g\ : \ K\ = \ \{0\}\cup \{1/n \ | \ n\geq 1\}\ \to \ X \]
  by $g(0)=z$ and $g(1/n)=x_n$. The hypotheses
  imply that the map $g$ is continuous if we give the source $K$
  the subspace topology of the interval $[0,1]$.
  In this topology the space $K$ is compact, so $g^{-1}(A)$ is
  closed since $A$ was assumed to be compactly closed.
  On the other hand, all the points $1/n$ are contained in $g^{-1}(A)$,
  and the closure of the set of these point is all of $K$. 
  So $0\in g^{-1}(A)$, which means that $z=g(0)\in A$.
  So $A$ coincides with its closure, i.e., $A$ is closed in $X$.

  (v) 
  In a metric space the balls of radius $1/n$ for all $n\geq 1$
  form a countable neighborhood basis of a given point. 
  So metric spaces are first countable, hence $k$-spaces by~(iv).

  (vi) This argument goes back to D.\,E.\,Cohen \cite[3.2]{cohen-weak topology}
  who thanks J.\,H.\,C. White\-head for `suggesting the subject (...) and for 
  valuable help'.
  We let $A$ be a compactly closed subset of $X\times_0 Y$.
  We let $(x_0,y_0)\in (X\times Y)-A$ be a point in the complement.
  We let
  \[ A_0 \ = \ \{ y\in Y\ | \ (x_0,y)\in A\} \]
  be the `slice' of $A$ through $x_0\in X$.
  We let $N$ be a compact neighborhood of $y_0$ in $Y$.
  Then $A_0\cap N$ is closed in $N$ because $A$ is compactly closed.
  Since $Y$ is Hausdorff, $N$ is closed in $Y$, and hence
  $A_0\cap N$ is closed in $Y$. 
  Since $y_0\not\in A_0$, the set $Y-A_0$ is an open neighborhood of $y_0$.
  Since $Y$ is locally compact Hausdorff, there is a compact neighborhood $K$ 
  of $y_0$ with $K\subset Y-A_0$.
  We let
  \[ B \ = \ \{x \in X \ | \ (\{x\}\times K)\cap A \ne \emptyset \} \]
  be the projection of $(X\times K)\cap A$ to $X$.
  The condition $K\subset Y-A_0$ is then equivalent to $x_0\not\in B$.
  
  We let $f:C\to X$ be a continuous map from a compact space. Then $C\times K$ is compact,
  and so $(f\times K)^{-1}((X\times K)\cap A)$ is closed in $C\times K$ 
  since $A$ is compactly closed.
  Hence $(f\times K)^{-1}((X\times K)\cap A)$ is compact in the subspace topology
  inherited from $C\times K$.
  Since 
  \[ f^{-1}(B) \ = \ \{c \in C \ | \ (\{f(c)\}\times K)\cap A \ne \emptyset \} \]
  is the projection of $(f\times K)^{-1}((X\times K)\cap A)$ onto $C$, 
  the set $f^{-1}(B)$ is closed in $C$.
  Altogether this shows that the set $B$ is compactly closed.
  Since $X$ is a $k$-space, $B$ must be closed in $X$.
  Since $x_0\not\in B$, the set $(X-B)\times K$ is a neighborhood of $(x_0,y_0)$.
  Moreover, $(X-B)\times K$ is disjoint from $A$ by definition of the set $B$.
  This shows that the complement of the original set $A$ is open,
  hence $A$ is closed. This completes the proof that 
  $X\times_0 Y$ is a $k$-space in the product topology.
\end{proof}

The analog of Proposition \ref{prop:properties k-spaces}~(ii) 
is also true for {\em open} subsets, i.e., every open subset 
of a $k$-space is a $k$-space with respect to the subspace topology.
We will not use this, 
so we refer the reader to \cite[Prop.\,7.9.10]{tomDieck-algebraic topology}.

\medskip

If $X$ is any topological space we let $k X$ be the space which has the
same underlying set as $X$, but such that the closed subsets of $k X$
are the compactly closed subsets of $X$. This indeed defines a topology 
which makes $k X$ into a $k$-space
and such that the identity $\Id:k X\to X$ is continuous.
Moreover, any continuous map $Y\to X$ whose source $Y$ is a $k$-space 
is also continuous when viewed as a map to $k X$. 
In more fancy language, the assignment $X\mapsto k X$ extends to a functor
$k:\bSpc \to \mathbf K$ that is right adjoint
to the inclusion of the full subcategory of $k$-spaces.
Since the inclusion $\mathbf K \to\bSpc$ has a right adjoint,
the category $\mathbf K$ of $k$-spaces has small limits and colimits. 
Colimits can be calculated in the ambient category of all topological spaces;
equivalently, any colimit of $k$-spaces is again a $k$-space.
To construct limits, we can first take a limit 
in the ambient category of all topological spaces; this ambient limit
need not be a $k$-space, but applying the Kelleyfication functor  
$k:\bSpc \to \mathbf K$ yields a limit in $\mathbf K$. Since $k$
does not change the underlying set, the categories $\mathbf K$ and
$\bSpc$ share the property that the forgetful functor to sets
preserves all limits and colimits. More loosely speaking, the underlying
set of a limit or colimit in $\mathbf K$ is what one first thinks of.

The discussion about limits above applies in particular to products,
and the product of two $k$-spaces need {\em not} be a $k$-space in
the usual product topology.
As already mentioned, we denote by $X\times_0 Y$ the cartesian
product of two spaces $X$ and $Y$, endowed with the product topology
(which makes it a categorical product in the category $\bSpc$ of all topological spaces).
We denote by $X\times Y=k(X\times_0 Y)$ the Kelleyfication 
of the product topology; if $X$ and $Y$ are $k$-spaces, then $X\times Y$
is a categorical product in the category $\mathbf K$.
Proposition \ref{prop:properties k-spaces}~(vi)
shows that Kelleyfication is unnecessary if one
of the factors is locally compact Hausdorff.

An important example where products in $\bK$ and $\bSpc$ 
can differ is the product of two `sufficiently large' CW-complexes $X$ and $Y$. 
Every CW-complex is a $k$-space,
for example by \cite[Prop.\,1.2.1]{fritsch-piccinini}.
The product $X\times_0 Y$ with the usual product topology 
comes with a filtration $(X\times_0 Y)_{(n)}=\cup_{p+q=n}X_{(p)}\times_0 Y_{(q)}$,
where $X_{(p)}$ is the $p$-skeleton of the CW-structure on $X$.
If $X$ or $Y$ is locally compact, then the product topology
is a $k$-space by Proposition \ref{prop:properties k-spaces}~(vi),
and then the above filtration makes $X\times_0 Y$ into a CW-complex,
as was already noted 
by J.\,H.\,C.\,Whitehead \cite[(H), p.\,227]{whitehead-combinatorialI}.
In general, however, 
$X\times_0 Y$ may not be a $k$-space, and hence cannot have a CW-structure.
The first example of this phenomenon 
was given by Dowker \cite[III.5, p.\,563]{dowker-metric complexes},
namely where $X$ and $Y$ are countably respectively uncountably
infinite wedges of circles.
This is also an example where the topology on $X\times Y=k(X\times_0 Y)$
is strictly finer than the product topology.
The product in the category $\mathbf K$, i.e., the space $X\times Y=k(X\times_0 Y)$,
is always compactly generated and a CW-complex via the above filtration,
see for example \cite[Thm.\,A.6]{hatcher}.

\medskip

A surjective continuous 
map $p:X\to Y$ is a {\em proclusion}\index{subject}{proclusion}
if whenever $O\subset Y$ is such that $p^{-1}(O)$
is open in $X$, then $O$ is already open in $Y$.
So a proclusion is a continuous map that is homeomorphic, under $X$,
to the projection onto a quotient space.
The following proposition is one of the key properties of the category of $k$-spaces,
and is stated without proof as Proposition~2.2 of \cite{mccord}.
The argument goes back, at least, to Steenrod,
who states it in \cite[Thm.\,4.4]{steenrod-convenient}
for Hausdorff $k$-spaces. The analog of the following proposition
does {\em not} hold for general topological spaces in the usual product topology,
and the examples in \cite[Ex.\,7.9.23]{tomDieck-algebraic topology}
illustrate what can go wrong.

\begin{prop}\label{prop:proclusion times Z}
Let $X$ and $Z$ be $k$-spaces and $p:X\to Y$ a proclusion.
Then the map $p\times Z:X\times Z\to Y\times Z$ is a proclusion.
\end{prop}
\begin{proof}
 We adapt Steenrod's argument from \cite[Thm.\,4.4]{steenrod-convenient}
to the more general context, i.e., for $k$-spaces
without any additional separation hypothesis.
We start with the special case when $Z$ is compact. 
In that case $X\times Z=X\times_0 Z$ and
$Y\times Z=Y\times_0 Z$ by Proposition \ref{prop:properties k-spaces}~(vi),
i.e., the usual product topologies coincide with their Kelleyfications.
In this formulation, the special case goes back 
to J.\,H.\,C.\,Whitehead \cite[Lemma~4]{whitehead-note borsuk},
and the proof can also be found in many topology textbooks,
for example \cite[Ch.\ 8, Lemma 8.9]{rotman-intro at}.

Now we reduce the general case to the special case.
We consider a subset $A\subset Y\times Z$
such that $(p\times Z)^{-1}(A)$ is closed in $X\times Z$.  
We let
\[ f \ = \ (f_1,f_2) \ : \ K\ \to \ Y\times Z \]
be a continuous map from a compact space.
Then 
\[ (p\times K)^{-1}((Y\times f_2)^{-1}(A))\ = \ 
 (X\times f_2)^{-1}((p\times Z)^{-1}(A)) \]
is closed in $X\times K$.
The map $p\times K:X\times K\to Y\times K$
is a proclusion by the special case, so 
the set $(Y\times f_2)^{-1}(A)$ is closed in $Y\times K$.
So
\[ f^{-1}(A)\ = \ (f_1,\Id_K)^{-1}( (Y\times f_2)^{-1}(A)) \]
is closed in $K$.
This shows that $A$ is compactly closed, and hence closed in
the $k$-space topology of $Y\times Z$.
\end{proof}

Now we turn to weak Hausdorff spaces, and collect various
useful properties in the following proposition.

\begin{prop}\label{prop:properties wH spaces}
  \begin{enumerate}[\em (i)]
    \item If $i:A\to X$ is a continuous injection and $X$ is weak Hausdorff, then
    $A$ is also weak Hausdorff. In particular, every subspace 
of a weak Hausdorff space is again weak Hausdorff.
  \item Every Hausdorff space is also weak Hausdorff.
  \item Every finite subset of a weak Hausdorff space is closed.
  \item Every continuous bijection from a compact space to weak Hausdorff space 
    is a homeomorphism.
  \item Let $f:K\to X$ be a continuous map from a compact space to a weak Hausdorff
    space. Then the image $f(K)$ is compact in the subspace topology.
  \item Every disjoint union of weak Hausdorff spaces is weak Hausdorff.
  \item Every limit in $\bSpc$ of a functor with values in 
    weak Hausdorff spaces is weak Hausdorff.
  \item If $X$ is a weak Hausdorff space, then its Kelleyfication
    $k X$ is again weak Hausdorff.
  \end{enumerate}
\end{prop}
\begin{proof}
  (i) Let $f:K\to A$ be a continuous map from a compact space.
  Then $i f:K\to X$ is also continuous, hence $(i f)(K)$ is closed in $X$.
  Since $i$ is injective we have $f(K)=i^{-1}( (i f)(K))$, which is thus
  closed in $A$. This shows that $A$ is a weak Hausdorff space.

  (ii) If $K$ is compact and $f:K\to X$ continuous, then
  the image $f(K)$ is always quasi-compact; if $X$ is Hausdorff, then
  any quasi-compact subset such as $f(K)$ is closed. 
  So $X$ is weak Hausdorff.

  (iii) Every one point space is compact, so  every point of any space is the continuous 
 image of a compact space. So in weak Hausdorff spaces, all points and 
 thus all finite subsets are closed.

 (iv) Let $f:X\to Y$ be a continuous bijection from a compact space 
 to a weak Hausdorff space. Every closed subset $A$ of $X$ is compact
 in the subspace topology, so $f(A)$ is closed in $Y$ by the weak Hausdorff
 property. This shows that $f$ is also a closed map, hence a homeomorphism.
 
 (v) The property of being quasi-compact is automatically inherited under
 continuous surjections, so the main issue is the Hausdorff property of $f(K)$,
 for which we reproduce the argument from \cite[Lemma 2.1]{mccord}.
 We let $x,y\in f(K)$ be two distinct points. The sets $\{x\}$
 and $\{y\}$ are closed in $X$ by part~(iii), so
 $f^{-1}(x)$ and $f^{-1}(y)$ are disjoint closed subsets of $K$.
 Since compact spaces are normal, there are disjoint open subsets $U$ and $V$
 of $K$ with $f^{-1}(x)\subset U$ and $f^{-1}(y)\subset V$.
 Since $X$ is weak Hausdorff, the sets $f(K-U)$ and $f(K-V)$
 are closed in $X$, and hence also in $f(K)$.
 But then $f(K)-f(K-U)$ and $f(K)-f(K-V)$ are disjoint open subsets of $f(K)$ 
 that separate $x$ and $y$.

 (vi) We let $\{Y_i\}_{i\in I}$ be a family of weak Hausdorff spaces and
 $f:K\to \coprod_{i\in I} Y_i$ a continuous map from a compact space to
 the disjoint union. Since $K$ is the disjoint union of the closed subspaces
 $K_i=f^{-1}(Y_i)$, each $K_i$ is compact in the subspace topology. 
 The image of the restriction $f|_{K_i}:K_i\to Y_i$ is closed
 since $Y_i$ is weak Hausdorff. So $f(K)=\coprod_{i\in I}f(K_i)$
 is closed in the disjoint union.

 (vii) In a first step we show that a product of any family 
 $\{Y_i\}_{i\in I}$ of weak Hausdorff spaces is weak Hausdorff.
 We let $f=(f_i)_{i\in I}:K\to \prod^0_{i\in I}Y_i$ be a continuous map from a compact space,
 with respect to the product topology on the target. 
 Because $Y_i$ is weak Hausdorff, the subset $f_i(K)$ is closed in $Y_i$, 
 and a Hausdorff space in the subspace topology by part~(v).
 So the subset $\prod_{i\in I} f_i(K)$ is closed in $\prod_{i\in I}^0 Y_i$
 and Hausdorff in the subspace topology. 
 Since $K$ is compact  and the image of $f$ is contained in
 the Hausdorff space $\prod_{i\in I} f_i(K)$, 
 the image $f(K)$ is closed in the subspace topology of $\prod_{i\in I} f_i(K)$.
 Since $\prod_{i\in I} f_i(K)$ is closed in $\prod_{i\in I}^0 Y_i$,
 the image $f(K)$ is also closed in $\prod^0_{i\in I}Y_i$.
 This proves that the product $\prod^0_{i\in I}Y_i$
 is weak Hausdorff.

 A limit in $\bSpc$ of a functor $F:I\to\bSpc$ is a subspace
 of the product of the values $F(i)$ for all object $i\in I$.
 Since a subspace of a weak Hausdorff space is weak Hausdorff by part~(i),
 this proves the claim.

 (viii) The identity is a continuous injection $k X\to X$, so part~(i) proves the claim.
\end{proof}

Now we turn to compactly generated spaces. 
By combining Propositions \ref{prop:properties k-spaces}
and \ref{prop:properties wH spaces}, we get some immediate corollaries:

\begin{prop}\label{prop:properties cgwh spaces}
  \begin{enumerate}[\em (i)]
  \item Every closed subset of a compactly generated space
    is compactly generated in the subspace topology.
  \item Every locally compact Hausdorff space, and hence every compact space,
    is compactly generated.
  \item Every metric space is compactly generated.
  \item Every disjoint union of compactly generated spaces is compactly generated.
    \item If $X$ is compactly generated and $Y$ locally compact Hausdorff, then
      $X\times_0 Y$ is compactly generated in the product topology.
  \end{enumerate}
\end{prop}

There is a useful criterion, due to McCord \cite[Prop.\,2.3]{mccord},
for when a $k$-space is weak Hausdorff (and hence compactly generated).

\begin{prop}\label{prop:wH criterion} 
A $k$-space $X$ is weak Hausdorff if and only if 
the diagonal is closed in $X\times X=k(X\times_0 X)$. 
\end{prop}
\begin{proof}
We suppose first that the  diagonal is closed in $X\times X$. 
We let $f:K\to X$ be a continuous map from a compact space.
Then the map
\[ f\times X \ : \ K\times X \ \to\ X\times X\]
is continuous, and so $(f\times X)^{-1}(\Delta_X)$ is closed in
$K\times X$. Since $K$ is compact, the topology on $K\times X$
is the usual product topology by Proposition \ref{prop:properties k-spaces}~(vi).
Moreover,  the projection $p:K\times X\to X$ away from $K$ is a closed map,
also by compactness.
Hence
\[ f(K)\ = \ p((f\times X)^{-1})(\Delta_X) \]
is closed in $X$. This shows that $X$ is weak Hausdorff.

For the converse we assume that $X$ is weak Hausdorff.
We show that the diagonal $\Delta_X$ is compactly closed in $X\times_0 X$,
and hence closed in $X\times X$. 
We let $f=(f_1,f_2):K\to X\times_0 X$ be a continuous map from a compact space.
Then the map
\[ f_1+f_2 \ : \ K\amalg K \ \to \ X \]
is continuous.
Since $X$ is weak Hausdorff, the subset $A=(f_1+f_2)(K\amalg K)=f_1(K)\cup f_2(K)$
of $X$ is closed and a Hausdorff space in the subspace topology, by
Proposition \ref{prop:properties wH spaces}~(v). 
So the diagonal $\Delta_A$ is closed in $A\times_0 A$.
Since moreover, $f(K)$ is contained in $A\times_0 A$,
the set
$ f^{-1}(\Delta_X)  =  f^{-1}(\Delta_A)$
is closed in $K$. This completes the proof that $\Delta_X$ is compactly closed.
\end{proof}

Proposition \ref{prop:wH criterion} now leads to the following useful criterion
for when a quotient space of a $k$-space is weak Hausdorff (and hence compactly generated).
\begin{prop}\label{prop:closed implies WH}
  Let $X$ be a $k$-space and $E\subset X\times X$ an equivalence relation.
  Then the quotient space $X/E$ is compactly generated if and only if $E$ is
  closed in the $k$-topology of $X\times X$.
\end{prop}
\begin{proof}
We let $p:X\to X/E$ denote the projection.
Any quotient space of a $k$-space is automatically a $k$-space
by Proposition \ref{prop:properties k-spaces}~(i).
If $X/E$ is also weak Hausdorff, then $\Delta_{X/E}$ is closed in 
$X/E\times X/E$ by Proposition \ref{prop:wH criterion}. 
Hence 
\begin{equation}  \label{eq:E_inverse_pxp}
 E \ = \ (p\times p)^{-1}(\Delta_{X/E}) 
\end{equation}
is closed in $X\times X$.
Conversely, suppose that $E$ is closed in $X\times X$.
Since $X$ and $X/E$ are $k$-spaces, the map
\[ p\times p \ : \ X\times X \ \to \ (X/E)\times(X/E) \]
is a proclusion by two applications of Proposition \ref{prop:proclusion times Z}.
Since $E$ is closed by hypothesis, the relation \eqref{eq:E_inverse_pxp} 
shows that $\Delta_{X/E}$ is closed in $(X/E)\times(X/E)$.
So $X/E$ is weak Hausdorff by Proposition \ref{prop:wH criterion}.
\end{proof}

\begin{cor}\label{cor:X/A WH}
Let $X$ be a compactly generated space and $A$ a closed subset of $X$.
Then the quotient topology on $X/A$ is again compactly generated.
\end{cor}
\begin{proof}
Since $A$ is closed in $X$, $A\times A$ is closed in $X\times X$.
Since $X$ is weak Hausdorff, the diagonal $\Delta_X$ is closed $X\times X$.
So the equivalence relation $E=(A\times A)\cup\Delta_X$ is closed in $X\times X$,
and the quotient space $X/A$ is compactly generated 
by Proposition \ref{prop:closed implies WH}.
\end{proof}

Proposition \ref{prop:closed implies WH} 
also suggests how to construct a `maximal weak Hausdorff quotient'
of a $k$-space:

\begin{prop}
Let $X$ be a $k$-space. Let $E_{\min}\subset X\times X$ be the intersection of
all equivalence relations on $X$ that are closed
in the $k$-topology of $X\times X$.
Then $X/E_{\min}$ with the quotient topology is a compactly generated space
and the quotient map $X\to X/E_{\min}$ is initial among
continuous maps from $X$ to a compactly generated space.
\end{prop}
\begin{proof}
The intersection $E_{\min}$ is again an equivalence relation,
and $E_{\min}$ is closed in  $X\times X$ as an intersection of closed subsets; 
so the quotient topology is weak Hausdorff by Proposition \ref{prop:closed implies WH}.

Now we let $f:X\to Y$ be a continuous map to a compactly generated space.
The equivalence relation
\[ E \ = \ \{ (x,x')\in X\times X \ | \ f(x)=f(x') \} \]
is the preimage of the diagonal under the continuous map
\[ f\times f \ : \ X\times X\ \to \ Y\times Y \ . \]
Since the diagonal is closed in $Y\times Y$ by Proposition \ref{prop:wH criterion}, 
the set $E$ is closed in $X\times X$, 
and so $E_{\min}\subseteq E$. So $f$ factors uniquely over
a continuous map from $X/E_{\min}$, by the universal property of the
quotient topology.
\end{proof}

The previous proposition implies that the assignment
\[ X\ \longmapsto \ X/E_{\min} \ = \ w(X) \]
extends canonically to a functor
\[ w \ : \ \bK \ \to \ \bT \]
that is left adjoint to the inclusion of compactly generated spaces into $k$-spaces. 
Moreover, if $X$ is already compactly generated, then the 
diagonal is closed in $X\times X$ by Proposition \ref{prop:wH criterion}; 
every equivalence relation contains the diagonal, 
so $E_{\min}=\Delta_X$ whenever $X$ is 
compactly generated. In this situation the quotient map $X\to X/E_{\min}=w(X)$
is a homeomorphism.

It follows formally from the existence of a left adjoint 
to the inclusion $\mathbf K\subset \bT$ that the category of compactly generated spaces 
has small limits and colimits; 
limits can be calculated in the category $\mathbf K$ of $k$-spaces.\index{subject}{limit!in $\bT$}
To construct a colimit of a diagram in $\bT$, we can first take a colimit 
in the category $\mathbf K$ of $k$-spaces (or equivalently in $\bSpc$);
while a $k$-space, this colimit need not be weak Hausdorff,
but applying the left adjoint $w:\mathbf K\to\bT$ yields a colimit in $\bT$.
\index{subject}{colimit!in $\bT$} 

\medskip

\Danger If the colimit, taken in the category $\bK$ of $k$-spaces,
of a diagram of compactly generated spaces is not already weak Hausdorff,
then the minimal closed equivalence relation on it is strictly larger than the
diagonal, so the `maximal weak Hausdorff quotient' 
identifies at least two distinct points and thus changes the underlying set.
So one has to be especially careful with general colimits
in $\bT$: unlike for $\bSpc$ of $\mathbf K$,
the forgetful functor from $\bT$ to sets need not preserve colimits.
More loosely speaking, the underlying set  of colimit in $\bT$ 
may be smaller than one first thinks.

\medskip

The fact that colimits in $\bT$ may be hard to identify 
may seem like a problem at first. However, the issue is largely irrelevant
for purposes of homotopy theory because we don't expect to have
homotopical control over arbitrary colimits anyhow.
The colimits that we do care about turn out to be `as expected';
in particular,
for pushouts along closed embeddings (see Proposition \ref{prop:pushout in cg}),  
filtered colimits along closed embeddings 
(see Proposition \ref{prop:filtered colim in cg}), 
wedges (see Proposition \ref{prop:compact to wedge}), 
and orbits by actions of compact topological groups 
(see Proposition \ref{prop:G quotient properties}),
it makes no difference if we calculate 
the colimit in the category $\bT$ or in $\bK$ respectively $\bSpc$.

We call a continuous map $i:A\to X$ between topological spaces
a {\em closed embedding}\index{subject}{closed embedding} if $i$ is injective
and a closed map; equivalently, the image $i(A)$ is closed in $X$ 
and $i$ is a homeomorphism onto its image.
By Proposition \ref{prop:pushout in cg} below,
the cobase change, in $\bSpc$, $\bK$ or $\bT$, of a closed embedding is again 
a closed embedding.

\medskip

\Danger There is an ambiguity with the meaning of `embedding' in general, 
due to the fact that a general subset of a $k$-space, endowed with the
subspace topology, need not be a $k$-space,
and so one may or may not want to apply `Kelleyfication' $k:\bSpc\to\bK$
to the subspace topology. However, closed subsets of $k$-spaces
are again $k$-spaces in the usual subspace topology,
so there is no such ambiguity with the notion of `closed embedding'.

\medskip

\begin{prop}\label{prop:product closed embeddings}
Let $i:A\to X$ and $j:B\to Y$ be closed embeddings between topological spaces.
Then the product maps $i\times_0 j:A\times_0 B\to X\times_0 Y$
and $i\times j:A\times B\to X\times Y$ are closed embeddings.
\end{prop}
\begin{proof}
The map $i\times j$ is clearly injective, so must show that it is a closed
map with respect to the product topologies, and with respect to their Kelleyfications.
We may assume that $i$ is the inclusion of a closed subset $A\subset X$
and $j$ is the inclusion of a closed subset $B\subset Y$.
The subspace topology on $A\times B$ induced from $X\times_0 Y$
agrees with the product topology $A\times_0 B$,
which is the claim for $i\times_0 j$, i.e., the usual product topologies.

For the second claim we must show that every compactly closed subset $C$
of $A\times_0 B$ is also compactly closed in $X\times_0 Y$.
So we let $f:K\to X\times_0 Y$ be a continuous map from a compact space.
Since $A\times_0 B$ is closed inside $X\times_0 Y$,
the subset $L=f^{-1}(A\times_0 B)$ is closed in $K$, and hence compact
in the subspace topology. 
The relation
\[ f^{-1}(C)\ = \ (f|_L)^{-1}(C) \]
and the hypothesis that $C$ is compactly closed inside $A\times_0 B$
then show that $f^{-1}(C)$ is closed in $K$. 
So $C$ is compactly closed in $X\times_0 Y$.
\end{proof}

This following proposition is \cite[Lemma 8.1]{lewis-thesis}.

\begin{prop}\label{prop:section is closed embedding}
Let $i:X\to Y$ be a continuous map between compactly generated spaces
that admits a continuous retraction.
Then $i$ is a closed embedding.
\end{prop}
\begin{proof}
Let $r:Y\to X$ be a continuous retraction, i.e., such that $r i=\Id_X$.
Then the map
\[ (\Id_Y,i r) \ : \ Y\ \to \ Y\times Y\]
is continuous and for every subset $A\subset X$ we have
\[ i(A) \ =  \ \{y\in Y\ | \ y=i (r(y)) \text{\ and\ } r(y)\in A\}
\ = \ (\Id_Y,r i)^{-1}(\Delta_Y)\, \cap\, r^{-1}(A)\ .\]
Since $Y$ is compactly generated, the diagonal $\Delta_Y$ is closed
in $Y\times Y$ by Proposition \ref{prop:wH criterion};
hence the set $(\Id_Y,r i)^{-1}(\Delta_Y)$ is closed in $Y$.
If $A$ is closed, then so is $r^{-1}(A)$, and hence also $i(A)$.
This proves that $i$ is a closed embedding.
\end{proof}

\begin{prop}\label{prop:pushout in cg}
    Given a pushout square in the category $\bSpc$ of topological spaces
    \[ \xymatrix{ X \ar[r]^-f \ar[d]_i & Z \ar[d]^j \\
      Y \ar[r]_-g & P } \]
       such that $i$ is a closed embedding, then $j$ is also a closed
    embedding. If moreover $Y$ and $Z$ are compactly generated, then
    so is $P$, and hence the square is a pushout in $\bT$.
\end{prop}
\begin{proof}
We adapt the argument given in \cite[Prop.\,2.5]{mccord}.
The map $j$ is injective because $i$ is. Indeed, we can choose
a set-theoretic retraction $r:Y\to X$ to $i$ (not necessarily continuous),
and then $(f r)\cup\Id: P\to Z$ is a set-theoretic retraction to $j$.

For the other claims we first treat the special case where the map $g$ is a proclusion.
We let $A\subset Z$ be a closed subset; then $f^{-1}(A)$ is closed in $X$. 
Since the map $i$ is injective, the relation
\[ g^{-1}(j(A)) \ = \ i(f^{-1}(A)) \]
holds as subsets of $Y$. So $i(f^{-1}(A))$, and hence $g^{-1}(j(A))$,
is closed in $Y$ because $i$ is a closed embedding.
Since $g$ is a proclusion, this shows that $j(A)$ is closed, 
and hence $j$ is a closed map.
Altogether this proves that $j$ is a closed embedding.

Now we suppose that $Y$ and $Z$ are compactly generated.
Then $P$ is a $k$-space by Proposition \ref{prop:properties k-spaces}~(i)
because $g:Y\to P$ is a proclusion.
Moreover, the diagonals of $Y$ and $Z$
are closed in $Y\times Y$ respectively $Z\times Z$ by 
Proposition \ref{prop:wH criterion}. 
Hence $(f\times f)^{-1}(\Delta_Z)$ is closed in $X\times X$.
Since $i$ is a closed embedding, so is $i\times i:X\times X\to Y\times Y$,
by Proposition \ref{prop:product closed embeddings}.
Because
\[ (g\times g)^{-1}(\Delta_P)\ = \ 
\Delta_Y\cup (i\times i)\left( (f\times f)^{-1}(\Delta_Z)\right) \]
we conclude that $(g\times g)^{-1}(\Delta_P)$ is closed in $Y\times Y$.
Since $g$ is a proclusion, so is $g\times g$, by two applications of
Proposition \ref{prop:proclusion times Z}.
So $\Delta_P$ is closed in $P\times P$, 
and $P$ is weak Hausdorff by the criterion of Proposition \ref{prop:wH criterion}. 

Now we treat the general case.
The pushout $P$ can be constructed as a quotient space of the disjoint union
$Y\amalg Z$, which means that the lower horizontal map 
in the commutative square
\[ \xymatrix@C=12mm{ X\amalg Z \ar[r]^-{f+\Id} \ar[d]_{i\amalg\Id} & Z \ar[d]^j \\
  Y \amalg Z\ar[r]_-{g+j} & P } \]
is a proclusion.
Since the original square is a pushout, so is this new square.
The spaces $Y\amalg Z$ and $Z$ are compactly generated by hypothesis
and Proposition \ref{prop:properties cgwh spaces}~(iv),
and the left vertical map is a closed embedding. 
So we can apply the special case above to this new square, 
and conclude that $P$ is compactly generated.
\end{proof}

The following example of Lewis \cite[App.\,A, p.168]{lewis-thesis} 
shows that a pushout of compactly generated spaces calculated 
in the ambient category $\bSpc$ need not be compactly generated, and that 
a cobase change in $\bT$ of a continuous injection need 
not be injective. We consider the diagram
\[ \{-1,1\} \ \xla{\qquad}\ [-1,0) \cup (0,1] \ \xra{\text{inclusion}} \ [-1,1] \]
where all three spaces have the subspace topology of $\mR$, and the left map 
takes $[-1,0)$ to $-1$ and it takes $(0,1]$ to $1$.
All three spaces are compactly generated, and the pushout $P$ 
in the categories $\bSpc$ or $\bK$ has three points, only one of which is closed.
Any continuous map from $P$ to a weak Hausdorff space must be constant, 
so the space $w(P)$ (which is a pushout in $\bT$) has only one point.

\medskip

We recall that a poset $(P,\leq)$ 
is {\em filtered}\index{subject}{filtered  poset} if 
for all elements $i,j\in P$ there is a $k\in P$ such that $i\leq k$ and $j\leq k$.
The associated category has object set $P$ 
and a unique morphism $(i,j):i\to j$ for every pair of elements such that $i\leq j$.
The filtered poset we mostly care about is $(\mN,\leq)$, the set of natural numbers
under it usual ordering; a functor $F:\mN\to \Cc$ from the associated poset category
is determined by the sequence of morphisms
\[ F(0)\ \xra{F(0,1)}\ F(1)\ \xra{F(1,2)}\ F(2)\ \xra{F(2,3)}\ \dots\ 
\xra{F(n-1,n)}\ F(n)\ \xra{\qquad} \ \dots \]
So colimits of such functors are sequential colimits.
Other filtered posets that come up are the set of finite subsets of an infinite set,
and the set of finite-dimensional vector subspaces of an infinite dimensional
vector spaces, both ordered by inclusion.

\begin{prop}\label{prop:filtered colim in cg}  
  Let $P$ be a filtered poset and $F:P\to \bSpc$ a
  functor from the associated poset category to the category of topological spaces. 
  Let $F_\infty$ be a colimit of $F$ in the category $\bSpc$
  with respect to continuous maps $\kappa_i:F(i)\to F_\infty$.
  \begin{enumerate}[\em (i)]
  \item If the map $F(i,j):F(i)\to F(j)$ is a closed embedding
    for every comparable pair of elements of $P$, then 
    all the maps $\kappa_i:F(i)\to F_\infty$ are closed embeddings.
  \item If the space $F(i)$ is compactly generated for every element $i\in P$
    and the map $F(i,j):F(i)\to F(j)$ in injective 
    for every comparable pair of elements,
    then the colimit $F_\infty$ is compactly generated.
    Hence $F_\infty$ is also a colimit of the functor $F$ in the category $\bT$.
  \end{enumerate}
\end{prop}
\begin{proof}
(i) Colimits in $\bSpc$ are created on underlying sets, so the map
$\kappa_i:F(i)\to F_\infty$  is injective because all the maps $F(i,j)$ are injective.
Now we let $A\subset F(i)$ be closed subset. If $j\in P$ is
another element, we can choose $k\in P$ such that $i\leq k$ and $j\leq k$.
Then
\[ \kappa_k^{-1}(\kappa_k(F(i,k)(A))) \ = \ F(i,k)(A) \]
because the map $\kappa_k$ is injective.
Hence
\begin{align*}
  \kappa_j^{-1}(\kappa_i(A))\ =\ F(j,k)^{-1}(\kappa_k^{-1}(\kappa_k(F(i,k)(A))))\ =\ 
F(j,k)^{-1}(F(i,k)(A)) \ .
\end{align*}
Since $F(i,k)$ is a closed map, this shows that
$\kappa_j^{-1}(\kappa_i(A))$ is closed in $F(j)$ for all $j\in P$.
The map $\coprod_{j\in P}F(j)\to F_\infty$
given by $\kappa_j$ on $F(j)$ is a proclusion, 
so we conclude that $\kappa_i(A)$ is closed in $F_\infty$.
So $\kappa_i$ is a closed embedding.

(ii) The map $\kappa:\coprod_{j\in P}F(j)\to F_\infty$
given by $\kappa_j$ on $F(j)$ is a proclusion, and the source is compactly generated by 
Proposition \ref{prop:properties cgwh spaces}~(iv).
We show that the equivalence relation 
\[ E \ = \ (\kappa\times\kappa)^{-1}(\Delta_{F_\infty}) \ \subset \ 
\left({\amalg}_{j\in P}\, F(j)\right)\times \left({\amalg}_{j\in P}\, F(j)\right) \]
that gives rise to $F_\infty$ is closed; the criterion
of Proposition \ref{prop:closed implies WH} then shows that $F_\infty$
is compactly generated.
The product of ${\coprod}_{j\in P}F(j)$ with itself
is the disjoint union of the subspaces $F(i)\times F(j)$ for all $i,j\in P$,
so it suffices to show that the intersection of $E$ with all
$F(i)\times F(j)$ is closed. We choose $k\in P$ such that $i\leq k$ and $j\leq k$.
Then
\[ E\cap (F(i)\times F(j))\ = \ (F(i,k)\times F(j,k))^{-1}(\Delta_{F(k)}) \]
because the map $\kappa_k:F(k)\to F_\infty$ is injective.
The diagonal is closed in $F(k)\times F(k)$ since $F(k)$ is compactly generated
(by Proposition \ref{prop:closed implies WH}).
So $E\cap (F(i)\times F(j))$ is closed in $F(i)\times F(j)$.
This verifies the criterion of Proposition \ref{prop:closed implies WH} 
and concludes the proof that the colimit $F_\infty$ is compactly generated.
\end{proof}

Many functors $F:P\to\bT$ out of filtered posets that arise in practice 
have the property that a continuous map $f:K\to F_\infty$
from a compact space factors through $F(i)$ for some $i\in P$.
Such filtered colimits then tend to preserve weak equivalences.
We now formulate a precise version of this property for certain lattices.
The most common special case is the poset $(\mN,\leq)$ 
of natural numbers under the usual linear ordering.

We recall that a poset $(P,\leq)$ is a {\em lattice}\index{subject}{lattice}
if every pair of elements $p,q\in P$ has a join $p\vee q$ (least upper bound) 
and a meet $p\wedge q$ (greatest lower bound).
Every lattice is in particular a filtered poset.

\begin{prop} \label{prop:filtered colim preserve weq}
Let $(P,\leq)$ be a lattice with the following property:
for every element $q\in P$ and every countable chain
\[ p_1\ \leq\ p_2\ \leq\ \dots\ \leq\ p_n\ \leq\ \dots \]
in $P$ with $p_n\leq q$ for all $n\geq 1$, 
the sequence is eventually constant.
\begin{enumerate}[\em (i)]
\item 
Let $F:P\to\bT$ be a functor to the category of compactly generated spaces
with the following properties:
\begin{enumerate}[\em (a)]
\item the map $F(i,j):F(i)\to F(j)$ is a closed embedding for all comparable 
pairs $i\leq j$ of elements of $P$, and
\item for all elements $p,q\in P$ the following square is a pullback:
\[ 
\xymatrix{ 
F( p\wedge q) \ar[r]\ar[d] & F(q)\ar[d]\\
F(p)\ar[r] & F(p\vee q)}    
 \]\end{enumerate}
Let $F_\infty$ be a colimit of $F$ with respect to the continuous maps
$\kappa_i:F(i)\to F_\infty$.
Then for every  continuous map $\alpha:K\to F_\infty$ from a compact space $K$
there exists an element $i\in P$ and a  continuous map $\alpha':K\to F(i)$
such that $\alpha=\kappa_i\circ\alpha'$.
\item Let $F':P\to\bT$ be another functor satisfying conditions {\em (a)} 
and {\em (b)} and $\psi:F\to F'$ a natural transformation.
Suppose that the map $\psi(i):F(i)\to F'(i)$ is $m$-connected for all $i\in P$,
for some $m\geq 0$.
Then the map $\psi_\infty:F_\infty\to F'_\infty$ induced on colimits is $m$-connected.
\item Let $F:P\to\bT$ satisfy conditions {\em (a)} 
and {\em (b)} and suppose that the map $F(i,j):F(i)\to F(j)$ is a $m$-connected
for all pairs of comparable elements.
Then the map  $\kappa_i:F(i)\to F_\infty$ is $m$-connected for every $i\in P$.
\end{enumerate}
\end{prop}
\begin{proof}
(i) The canonical maps $\kappa_i:F(i)\to F_\infty$ are
closed embeddings by Proposition \ref{prop:filtered colim in cg}~(i). 
Hence we can, and will, pretend that $F(i)$ is a closed subspace
of $F_\infty$ and the maps $F(i,j):F(i)\to F(j)$ are inclusions.

The image of the continuous map $f:K\to F_\infty$ is compact in the subspace topology
by Proposition \ref{prop:properties wH spaces}~(v). So we may replace
$K$ by its image and pretend that $f$ is the inclusion of a compact subspace
of $F_\infty$. We then have to show that $K\subset F(i)$ for some $i\in P$.

We argue by contradiction and assume that $K$ is not contained in $F(i)$
for any $i\in P$.
We can then choose comparable elements
\[ p_1\ \leq \ p_2 \  \leq \ \dots\ \leq \ p_n \  \leq \ \dots \]
of $P$ and elements $x_i\in K$ with $x_i\in F(p_i)\bs F(p_{i-1})$.
Indeed, we start with any choice $x_1\in F(p_1)$ and continue inductively:
since $K$ is not contained in $F(p_n)$, there is an element $x_{n+1}\in K\bs F(p_n)$.
We must have $x_{n+1}\in F(q)$ for some $q\in P$, and then $p_{n+1}=p_n\vee q$
can serve for the inductive step.

We set $C =  \{x_1,x_2,x_3,\dots\}$, a countably infinite subset of $K$.
Now we claim that for every $q\in P$ the 
intersection $C\cap F(q)$ is finite.
If this were not the case, then there would be infinitely many indices
$i_1<i_2<i_3<\dots$ such that $x_{i_j}\in F(q)$ for all $j\geq 1$.
The commutative square
\[ 
\xymatrix{ F( p_{i_j}\wedge q) \ar[r]\ar[d] & F(q)\ar[d]\\ F(p_{i_j}) \ar[r]& 
F(p_{i_j}\vee q)}    
 \]
is a pullback by hypothesis,
and so $x_{i_j}\in F(p_{i_j}\wedge q)$ for all $j\geq 1$.
Also by hypothesis there is an $N\geq 1$ such that
\[  p_N\wedge q\ = \  p_{N+1}\wedge q\ = \  p_{N+2}\wedge q\ = \ \cdots \ .  \]
So $F( p_N\wedge q)$, and hence also $F(p_N)$, contains infinitely
many of the elements of the set $C$. This contradicts the construction of $C$,
and we can conclude that the intersection $C\cap F(q)$ is finite.

Now we let $T$ be any subset of $C$.
Since $T\cap F(q)$ is finite, it is a closed subset of $F(q)$
by Proposition \ref{prop:properties wH spaces}~(iii).
The map
\[ \coprod \kappa_q \ : \ {\coprod}_{q\in P} F(q)\ \to \ F_\infty \]
is a proclusion, i.e., the space $F_\infty$ carries 
the quotient topology. Since $T\cap F(q)$ is closed for all $q\in P$,
the set $T$ is closed in $F_\infty$. Since $T$ was any subset of $C$,
the subset $C$ is discrete in the subspace topology from $F_\infty$.
So altogether $C$ is an infinite discrete subset of the compact space $K$,
which is a contradiction. This proves that $K$ is contained in $F(i)$
for some $i\in P$.

(ii) We consider a commutative square of continuous maps for $k\leq m$:
\[ \xymatrix{
\partial D^k\ar[r]^\alpha \ar[d]_{\text{incl}} & F_\infty  \ar[d]^{\psi_\infty} \\
D^k\ar[r]_\beta & F'_\infty }\]
Part~(i) provides $i,j\in P$ and continuous maps $\alpha':\partial D^k\to F(i)$
and $\beta':D^k\to F(j)$ such that 
$\kappa_i\circ\alpha'=\alpha$ and $\kappa'_j\circ\beta'=\beta$.
By replacing $i$ and $j$ by their join we may suppose that $i=j$.
Then 
\[ \kappa'_i\circ \psi(i)\circ\alpha'\ =\ 
\psi_\infty\circ \kappa_i\circ \alpha'\ =\ 
\psi_\infty\circ \alpha\ =\ \beta|_{\partial D^k}\ = \ 
\kappa'_i\circ \beta'|_{\partial D^k}\ . \]
Since $\kappa'_i$ is injective this shows that $\psi(i)\circ\alpha'=\beta'|_{\partial D^k}$.
Since $\psi(i)$ is $m$-connected, there is 
a continuous map $\lambda:D^k\to F(i)$ that restricts to $\alpha'$
on $\partial D^k$ and such that $\psi(i)\circ\lambda$ is homotopic, 
relative $\partial D^k$,
to $\beta'$. Then $\kappa_i\circ\lambda:D^k\to F_\infty$
solves the original lifting problem. So the map $F_\infty$ is $m$-connected.

(iii) We consider a commutative square of continuous maps for $k\leq m$:
\[ \xymatrix{
\partial D^k\ar[r]^\alpha \ar[d]_{\text{incl}} & F(i)  \ar[d]^{\kappa_i} \\
D^k\ar[r]_\beta & F_\infty }\]
Part~(i) provides a $j\in P$ and a continuous map $\beta':D^k\to F(j)$ such that 
$\kappa_j\circ\beta'=\beta$.
By replacing $j$ by the join $i\vee j$ 
and $\beta'$ by the map $F(i,i\vee j)\circ\beta':D^k\to F(i\vee j)$,
we can assume without loss of generality that $i\leq j$.
Then 
\[ \kappa_j\circ F(i,j)\circ\alpha\ =\ 
\kappa_i\circ \alpha\ =\ \beta|_{\partial D^k}\ = \ 
\kappa_j\circ \beta'|_{\partial D^k}\ . \]
Since $\kappa_j$ is injective this shows that $F(i,j)\circ\alpha=\beta'|_{\partial D^k}$.
Since $F(i,j)$ is $m$-connected, there is 
a continuous map $\lambda:D^k\to F(i)$ that restricts to $\alpha$
on $\partial D^k$ and such that $F(i,j)\circ\lambda$ is homotopic, 
relative $\partial D^k$, to $\beta'$. Then $\lambda$ also
solves the original lifting problem. So the map $\kappa_i$ is a $m$-connected.
\end{proof}

\begin{eg}\label{eg:linear order}
We let $(P,\leq)$ be a {\em linear order}\index{subject}{linear order},
i.e., a poset in which every two elements are comparable.
Such a poset is in particular a lattice. 
We observe that for every linear order, 
condition~(b) of Proposition \ref{prop:filtered colim preserve weq}~(i)
is trivially satisfied for every functor.
Indeed, any two elements $p,q\in P$ are comparable,
and we suppose that $p\leq q$, the other case being analogous. 
Then $p\wedge q=p$ and $p\vee q=q$, and
the two vertical maps in the commutative square of condition~(b) are
identity maps. Hence the square is a pullback.  
\end{eg}

The most familiar special case of Proposition \ref{prop:filtered colim preserve weq}
is the poset $(\mN,\leq)$ of natural numbers under the usual ordering.
This poset is a linear order, so as we explained in Example \ref{eg:linear order},
condition~(b) of Proposition \ref{prop:filtered colim preserve weq}~(i)
is automatically satisfied for every functor.
For $(\mN,\leq)$, part~(i) of Proposition \ref{prop:filtered colim preserve weq}
can be found in many topology textbooks, for example \cite[Prop.\,2.4.2]{hovey-book}.
For easier reference we spell out the content of part~(ii) and~(iii) 
of Proposition \ref{prop:filtered colim preserve weq} for the poset $(\mN,\leq)$
and for $m=\infty$, i.e., for weak equivalences.

\begin{prop}\label{prop:closed seq colim weak equivalences}\quad
  \begin{enumerate}[\em (i)]
   \item Let $e_n:X_n\to X_{n+1}$ and $f_n:Y_n\to Y_{n+1}$ be closed embeddings
    between compactly generated spaces, for $n\geq 0$. 
    Let $\psi_n:X_n\to Y_n$ be weak equivalences
    that satisfy $\psi_{n+1}\circ e_n=f_n\circ\psi_n$ for all $n\geq 0$.
    Then the induced map $\psi_\infty:X_\infty\to Y_\infty$ 
    between the colimits of the sequences is a weak equivalence.
  \item Let $f_n:Y_n\to Y_{n+1}$ be a closed embedding of compactly generated spaces
   that is also a weak equivalence, for $n\geq 0$. 
   Then the canonical map $Y_0\to Y_\infty$ to the colimit of the sequence
   is a weak equivalence.
  \end{enumerate}
\end{prop}

Besides the linear order $(\mN,\leq)$, 
further examples of relevant posets that satisfy the hypotheses of 
Proposition \ref{prop:filtered colim preserve weq} are 
the set of finite subsets of a given set 
(which features in the following proposition),
and the set of finite-dimensional vector subspaces of a given vector space.

\begin{prop}\label{prop:compact to wedge} 
Let $\{X_i\}_{i\in I}$ be a family of based compactly generated
spaces. Then the wedge (one-point union) $\bigvee_{i\in I}X_i$
is compactly generated, and thus a coproduct in the category $\bT_*$ 
of based compactly generated spaces. 
Moreover, for every compact space $K$ and every continuous map
$f:K\to\bigvee_{i\in I}X_i$ there is a finite subset $J$ of $I$ such that
$f$ factors through the subwedge $\bigvee_{j\in J} X_j$.
\end{prop}
\begin{proof}
The disjoint union $\coprod_{i\in I}X_i$ is compactly generated
by Proposition \ref{prop:properties cgwh spaces} (iv).
Moreover, the equivalence relation that identifies all the basepoints is closed
in $(\coprod_{i\in I}X_i)\times(\coprod_{i\in I}X_i)$.
So the quotient space $\bigvee_{i\in I}X_i$ is compactly generated by 
Proposition \ref{prop:closed implies WH}.

We let $P$ be the poset of finite subsets of $I$, ordered by inclusion;
this is a lattice satisfying the hypotheses of 
Proposition \ref{prop:filtered colim preserve weq}.
We consider the functor $F:P\to\bT$ sending a finite subset $J\subset I$
to the finite wedge $\bigvee_{j\in J}X_j$.
For $J'\subset J\subset I$, the map $\bigvee_{j\in J'}X_j\to \bigvee_{j\in J}X_j$
is a closed embedding 
(by direct inspection, or by Proposition \ref{prop:section is closed embedding}),
and property~(b) of Proposition \ref{prop:filtered colim preserve weq}~(i)
is satisfied.
Proposition \ref{prop:filtered colim preserve weq}~(i) 
thus provides the desired factorization through $\bigvee_{j\in J}X_j$ 
for some finite subset $J$ of $I$.
\end{proof}

Now we let $X$ and $Y$ be compactly generated topological spaces 
equipped with basepoints $x_0\in X$ respectively $y_0\in Y$. 
We define the {\em smash product} as\index{subject}{smash product!of based spaces}
\[ X\sm Y \ =  \ X\times Y / ( X\times\{y_0\}\cup \{x_0\}\times Y) \ ,\]
the quotient of the product (with the $k$-topology)
by the subspace $X\times\{y_0\}\cup \{x_0\}\times Y$.
We write $x\sm y$ for the class of~ $(x,y)\in X\times Y$ in the quotient space.
Since points in compactly generated spaces are closed
(Proposition \ref{prop:properties wH spaces}~(iii)), the map 
\[ X\vee Y\ \to \ X\times Y\]
from the coproduct to the product in $\bT_*$
is a closed embedding with image the subspace 
$X\times\{y_0\}\cup \{x_0\}\times Y$;
so $X\sm Y$ is also a cokernel of this map.

The next proposition is another important reason for working in 
the category $\bT$ of compactly generated spaces. 
Indeed, in the larger category $\bSpc$ of all topological spaces,
the smash product is not generally associative:
if smash products are defined as quotients of the
usual product topology, then the canonical set-theoretic bijection from 
$(X\sm Y)\sm Z$ to $X\sm (Y\sm Z)$ need not be continuous.
An explicit example is mentioned without proof in \cite[5.8]{puppe-homotopiemengen},
namely that $(\mQ\sm \mQ)\sm \mN$ and $\mQ\sm (\mQ\sm \mN)$ 
are not homeomorphic in $\bSpc$, where $\mQ$ has the subspace topology of $\mR$,
and $\mN$ has the discrete topology.
A proof can be found in \cite[Thm.\,1.7.1]{may-sigurdsson}.

\begin{prop}\label{prop:smash associativity} 
Let $X$, $Y$, and $Z$ be based compactly generated spaces.
  \begin{enumerate}[\em (i)]
  \item The space $X\sm Y$ is compactly generated in the quotient topology
    of $X\times Y$.
  \item The maps
    \[ X\sm Y \ \to \ Y\sm X\ , \quad  x\sm y\ \longmapsto \ y\sm x\]
    and
    \[ (X\sm Y)\sm Z \ \to \ X\sm (Y\sm Z)\ , \quad  
    (x\sm y)\sm z\ \longmapsto \ x\sm (y\sm z) \]
    are homeomorphisms.
  \end{enumerate}
\end{prop}
\begin{proof}
(i) In weak Hausdorff spaces all points are closed 
(Proposition \ref{prop:properties wH spaces}~(iii)), hence
$X\times\{y_0\}\cup \{x_0\}\times Y$ is a closed subspace of $X\times Y$.
So the quotient topology is compactly generated by Corollary \ref{cor:X/A WH}. 

The first claim of part~(ii) is straightforward from the fact that the twist homeomorphism
$X\times Y\to Y\times X$ sending $(x,y)$ to $(y,x)$ takes the subspace
$X\times\{y_0\}\cup \{x_0\}\times Y$ to the subspace
$Y\times\{x_0\}\cup \{y_0\}\times X$ and hence descends to a homeomorphism
on quotient spaces.

Since the projection $p_{X,Y}:X\times Y\to X\sm Y$ is a proclusion,
so is $p_{X,Y}\times Z: X\times Y\times Z\to (X\sm Y)\times Z$,
by Proposition \ref{prop:proclusion times Z}. Hence the map
\[ p_{X\sm Y,Z}\circ(p_{X,Y}\times Z) \ : \
X\times Y\times Z\ \to \ (X\sm Y)\sm Z \ , \quad 
(x,y,z)\ \longmapsto \ (x\sm y)\sm z\] 
is a proclusion as the composite of two proclusions.
The composite
\[ X\times Y\times Z\ \xra{p_{X\sm Y,Z}\circ(p_{X,Y}\times Z)} \ 
(X\sm Y)\sm Z \ \xra{(x\sm y)\sm z \mapsto x\sm(y\sm z)}\ X\sm (Y\sm Z) \]
is the map $p_{X,Y\sm Z}\circ(X\times p_{Y,Z})$, and hence continuous.
Since  $p_{X\sm Y,Z}\circ(p_{X,Y}\times Z)$ is a proclusion, the
second map $(X\sm Y)\sm Z \to X\sm (Y\sm Z)$ is continuous.
By the same reasoning the inverse map 
$X\sm (Y\sm Z)\to (X\sm Y)\sm Z$ is continuous. Hence the two maps
are mutually inverse homeomorphisms.
\end{proof}

Another key advantage of the categories $\bK$ and $\bT$
over the category of all topological spaces is that they are cartesian closed,
i.e., categorical product with a fixed object is a left adjoint.
The discovery of these facts was predated 
by similar result in closely related categories;
for example, Brown \cite[\S 3]{brown-function spaces}
discusses internal function spaces in the category of Hausdorff $k$-spaces,
and this was popularized by Steenrod in \cite[Sec.\,5]{steenrod-convenient}.
The first reference to cartesian closedness for $\bK$ and $\bT$ that I am aware of is
the appendix of Lewis' thesis \cite[App.\,A, Thm.\,5.5]{lewis-thesis};
we will closely follow Lewis' exposition.
Internal function spaces in $\bK$ and $\bT$ are given by the set of all continuous maps 
endowed with the Kelleyfication 
of a slight modification of the compact-open topology.
For weak Hausdorff spaces, the `modified' compact-open topology actually coincides with
the compact-open topology.

\begin{construction}\label{con:define C(X,Y)}
For topological spaces $X$ and $Y$, we let $C(X,Y)$
denote the set of continuous maps from $X$ to $Y$. 
For a continuous map $h:K\to X$ from a compact space $K$ and an open subset $U$ of $Y$
we define a subset of $C(X,Y)$ by
\[ N(h,U)\ = \ \{\varphi\in C(X,Y)\ | \ \varphi(h(K))\subset U\}\ . \]
We endow the set $C(X,Y)$ of continuous maps with the topology
generated by the subbasis consisting of all these sets $N(h,U)$.
This topology is covariantly functorial in $Y$ and contravariantly functorial in $X$.
Indeed, for every continuous map $f:Y'\to Y$ we have
\[ C(X,f)^{-1}(N(h,U))\ = \ N(h,f^{-1}(U))\ , \]
which is open in $C(X,Y')$, and so $C(X,f):C(X,Y')\to C(X,Y)$ is continuous.
Similarly, for a continuous maps $g:X\to X'$ and $h:K\to X$,
with $K$ a compact space, we have
\[ C(g,Y)^{-1}(N(h,U))\ = \ N(g h,U)\ , \]
which is open in $C(X',Y)$, and so $C(g,Y):C(X',Y)\to C(X,Y)$ is continuous.

Now we let $X$ and $Y$ be $k$-spaces.
Then $C(X,Y)$ need not be a $k$-space. We fix this by Kelleyfication,
i.e., we define the internal mapping space in $\bK$ as\index{symbol}{$\map(Y,Z)$ - {unbased mapping space}} 
\[ \map(X,Y)\ = \ k C(X,Y) \ .\]
\end{construction}

\begin{rk}[Relation to the compact-open topology]
We recall that a subbasis for the compact-open topology\index{subject}{compact-open topology} 
on $C(X,Y)$ is given by all sets
\[ W(K,U)\ = \ \{\varphi\in C(X,Y)\ | \ \varphi(K)\subset U\}\ , \]
where $K$ is compact subset of $X$ and $U$ is an open subset of $Y$.
Because $W(K,U)=N(\text{incl}:K\to X,U)$,
the topology generated by the sets $N(h,U)$ is finer than the compact-open topology.
In general, the compact-open topology may be strictly coarser.

For weak Hausdorff spaces $X$, the topology on $C(X,Y)$ 
coincides with the compact-open topology.
Indeed, if $X$ is weak Hausdorff, then the image of every continuous map $h:K\to X$
is compact in the subspace topology (Proposition \ref{prop:properties wH spaces}~(v)). 
So in this case $N(h,U)=W(h(K),U)$, and the two subbases coincide.
\end{rk}

We recall that a category 
is {\em cartesian closed}\index{subject}{cartesian closed category}
if it has finite products and product with a fix object is a left adjoint.
The proof of the following theorem can be found 
in \cite[App.\,A, Thm.\,5.5]{lewis-thesis}.

\begin{theorem}\label{thm:internal hom in K} 
For all $k$-spaces $X,Y$ and $Z$, the natural map
\[ \Psi\ : \ \map(X\times Y,Z) \ \to \ \map(X,\map(Y,Z))  \ , \quad
\Psi(f)(x)(y) \ =\  f(x,y) \]
is a homeomorphism. 
In particular, the category of $k$-spaces is cartesian closed.
\end{theorem}
\begin{proof}
We reproduce Lewis' arguments.
We start by checking that the map
\[ \eta_X \ : \ X \ \to \ C(Y,X\times Y)\ , \quad \eta_X(x)(y)\ =\ (x,y)\]
is continuous. 
Since $X$ is a $k$-space, it suffices to check that the composite
$\eta_X\circ f$ is continuous for every continuous map $f:L\to X$ 
from a compact space. Since $\eta_X\circ f=C(Y,f\times Y)\circ\eta_L$
and $C(Y,f\times Y)$ is continuous,
we can assume without loss of generality that $X$ is compact.
If $X$ is compact, then $X\times Y=X\times_0 Y$ 
by Proposition \ref{prop:properties k-spaces}~(vi), i.e.,
$X\times Y$ is already a $k$-space in the usual product topology.

We let $h:K\to Y$ be a continuous map from a compact space and $U$ an open subset 
of $X\times Y=X\times_0 Y$.
Then
\[ \eta_X^{-1}(N(h,U))
 \ = \  \{x\in X\ |\ \{x\}\times K \subset (X\times h)^{-1}(U)\}\ .\]
For every $x\in \eta_X^{-1}(N(h,U))$, the `tube lemma' 
(see for example \cite[Ch.\,8, Lemma 8.9']{rotman-intro at})
provides an open neighborhood $O$ of $x$ in $X$ 
such that $O\times K\subset (X\times h)^{-1}(U)$, 
i.e., such that $O\subset \eta_X^{-1}(N(h,U))$.
So the set $\eta_X^{-1}(N(h,U))$ contains a neighborhood of each of its points,
and is thus open. Altogether this shows that inverse images of
the subbasis of the topology on $C(Y,X\times Y)$ are open; 
hence $\eta_X$ is continuous.
Since $X$ is a $k$-space, the map $\eta_X$ is also continuous
with respect to the Kelleyfied topology on the target,
i.e., when considered as a map to $\map(Y,X\times Y)$.

We denote the evaluation map by 
\[ \epsilon_Z \ : \ \map(Y,Z)\times Y\ = \ 
k( C(Y,Z)\times_0 Y) \ \to \ Z\ , \quad \epsilon(f,y)\ =\ f(y)\ .\]
We show next that for every continuous map
$h:K\to Y$ from a compact space the composite
$\epsilon_Z\circ(C(Y,Z)\times_0 h):C(Y,Z)\times_0 K\to Z$ is continuous. 
We let $U$ be an open subset of $Z$ and consider 
a point $(\varphi,k)\in C(Y,Z)\times_0 K$ in $(\epsilon_Z\circ(C(Y,Z)\times_0 h))^{-1}(U)$;
this simply means that $k\in (\varphi\circ h)^{-1}(U)$.
Since $K$ is compact and $(\varphi\circ h)^{-1}(U)$ is open in $K$,
there is an open neighborhood $V$ of $k$ in $K$ whose closure $\bar V$
is contained in $(\varphi\circ h)^{-1}(U)$.
This closure $\bar V$ is compact in the subspace topology, 
so the set $N(h|_{\bar V}:\bar V\to Y, U)$ is open in $C(Y,Z)$,
and it contains the map $\varphi$.
Hence $N(h|_{\bar V}:\bar V\to Y, U)\times V$ is open in $C(Y,Z)\times_0 K$,
and contains the point $(\varphi,k)$.
Since this open set is also contained in $(\epsilon_Z\circ(C(Y,Z)\times_0 h))^{-1}(U)$,
this latter set is a neighborhood of each of its points, and hence open 
in $C(Y,Z)\times_0 K$. This completes the proof that the map 
$\epsilon_Z\circ (C(Y,Z)\times_0 h)$ is continuous. 

Now we show that the evaluation map $\epsilon_Z$ itself is continuous. 
At this point the Kelleyfication of the source is crucial, 
as $\epsilon_Z$ need not be continuous 
with respect to the topology $C(Y,Z)\times_0 Y$.
So it suffices to show that the composite $\epsilon_Z\circ(g,h):K\to Z$
is continuous for every continuous map $(g,h):K\to C(Y,Z)\times_0 Y$ 
from a compact space.
The map $\epsilon_Z\circ (g,h)$ coincides with the composite
\[ K \ \xra{(g,\Id_K)} \ C(Y,Z)\times_0 K \ \xra{C(Y,Z)\times_0 h} \ 
C(Y,Z)\times_0 Y \ \xra{\ \epsilon_Z \ }\ Z \ . \]
The composite of the last two maps is continuous by the previous paragraph,
so the whole composite $\epsilon_Z\circ (g,h)$ is continuous.

At this point we know that the two maps
\[ \eta_X \ : \ X \ \to \ \map(Y,X\times Y)\text{\quad and\quad}
 \epsilon_Z \ : \ \map(Y,Z)\times Y \ \to \ Z\]
are continuous.
The rest of the argument is formal.
Since $\epsilon_{X\times Y}\circ(\eta_X\times Y)$ is the identity of $X\times Y$
and $\map(Y,\epsilon_Z)\circ\eta_{\map(Y,Z)}$ is the identity of $\map(Y,Z)$,
the natural transformations $\eta$ and $\epsilon$ are the unit respectively counit
of an adjunction $(-\times Y,\map(Y,-))$.
In particular, the map
\[ \Psi\ : \ C(X\times Y,Z) \ \to \ C(X,\map(Y,Z))  \ , \quad
\Psi(f)(x)(y) \ =\  f(x,y) \]
is bijective.
We let $T$ be another $k$-space. Replacing $X$ by $T\times X$ shows that 
\[ \Psi\ : \ C(T\times X\times Y, Z) \ \to \ C(T\times  X,\map(Y,Z)) 
\]
is bijective. The adjunction just established translates this into the fact that 
\[ C(T,\Psi)\ : \ C(T,\map(X\times Y,Z)) \ \to \ C(T,\map(X,\map(Y,Z))) \]
is bijective. Since $T$ is an arbitrary $k$-space, $\Psi$ is a homeomorphism.
\end{proof}

Now we turn to compactly generated spaces.

\begin{theorem}\label{thm:internal hom in T} 
  \begin{enumerate}[\em (i)]
  \item 
    If $Y$ is a weak Hausdorff space and $X$ is any topological space,
    then the spaces $C(X,Y)$ and $\map(X,Y)$ are weak Hausdorff.   
    In particular, if $X$ and $Y$ are compactly generated, then so is $\map(X,Y)$.
  \item For all compactly generated spaces $X,Y$ and $Z$, the natural map
    \[ \Psi\ : \ \map(X\times Y,Z) \ \to \ \map(X,\map(Y,Z))  \ , \quad
    \Psi(f)(x)(y) \ =\  f(x,y) \]
    is a homeomorphism. 
    In particular, the category of compactly generated spaces is cartesian closed.
  \end{enumerate}
\end{theorem}
\begin{proof}
(i) For every point $x\in X$ the evaluation map
\[ \ev_x \ : \ C(X,Y)\ \to \ Y \ , \quad f\ \longmapsto\ f(x) \]
is continuous. Indeed, for an open set $U$ of $Y$ we have
$\ev_x^{-1}(U)=N(\text{incl}:\{x\}\to X,U)$ 
which is open because a one-point space is compact.

The evaluation map
\[ C(X,Y)\ \to \ {\prod}^0_{x\in X}\ Y \ , \quad f \ \longmapsto \ \{f(x)\}_{x\in X}\]
is injective, and it is continuous by the above.
Any product, with the usual product topology, of copies of $Y$ is again weak Hausdorff
by Proposition \ref{prop:properties wH spaces}~(vii).
So $C(X,Y)$ is weak Hausdorff by Proposition \ref{prop:properties wH spaces}~(i).
The space $\map(X,Y)=k C(X,Y)$ is then still weak Hausdorff
by Proposition \ref{prop:properties wH spaces}~(viii),
and hence compactly generated.

(ii)  Part~(i) shows that the category of compactly generated
spaces is closed under the internal function objects in the category $\bK$
of $k$-spaces. Since products in $\bT$ are also formed in the ambient category $\bK$,
part~(ii) becomes a special case of 
Theorem \ref{thm:internal hom in K}.
\end{proof}

\begin{construction}[Based mapping spaces]
Since we work a lot with {\em based} spaces, we record that the based
version of the mapping space is right adjoint to the smash product in $\bT_*$.
In combination with the commutativity and associativity property of
the smash product (Proposition \ref{prop:smash associativity}),
and the natural homeomorphisms $X\sm S^0\iso X\iso S^0\sm X$,
this shows that the smash product provides a closed symmetric monoidal structure
on the category $\bT_*$ of based compactly generated spaces.

The based adjunction is a straightforward consequence of the unbased
version in Theorem \ref{thm:internal hom in T}.
We let $(X,x_0), (Y,y_0)$ and $(Z,z_0)$ be based compactly generated spaces.
We denote by $\map_*(Y,Z)$ the subspace of $\map(Y,Z)$\index{symbol}{$\map_*(Y,Z)$ - {based mapping space}} 
consisting of all based continuous maps $f$, i.e., such that $f(y_0)=z_0$;
the basepoint is the constant map with value $z_0$. 
Since points are closed in compactly generated spaces and
$\map_*(Y,Z)$ is the preimage of $\{z_0\}$ under the continuous evaluation map
\[ \map(Y,Z)\ \to \ Z \ , \quad f\ \longmapsto \ f(y_0)\ , \]
we see that $\map_*(Y,Z)$ is a closed subset of $\map(Y,Z)$.
So $\map_*(Y,Z)$ is compactly generated in the subspace topology 
(Proposition \ref{prop:properties cgwh spaces}~(i)).
The composite
\[ X \ \xra{\ \eta_X\ }\ \map(Y,X\times Y) \ \xra{\map(Y,\text{proj})}
\map(Y,X\sm Y) \]
takes values in the subspace $\map_*(Y,X\sm Y)$, so it restricts to a continuous map
\[ \eta'_X \ : \ X\ \to \ \map_*(Y,X\sm Y)\ ,\]
which is moreover based.
The composite
\[ \map_*(Y,Z)\times Y \ \xra{\text{incl}\times Y} \ 
\map(Y,Z)\times Y \ \xra{\ \epsilon_Z\ } \ Z \]
takes the subspace $\{\text{const}\}\times Y\cup \map_*(Y,Z)\times \{y_0\}$
to the basepoint of $Z$, so it factors over the quotient space
though a continuous based map
\[ \epsilon_Z' \ : \ \map_*(Y,Z)\sm Y \ \to \ Z \ .\]
By direct verification, the composites
\[ \epsilon'_{X\sm Y}\circ(\eta'_X\sm Y)
\text{\quad and\quad}
\map(Y,\epsilon'_Z)\circ \eta'_{\map_*(Y,Z)}\]
are the identity of $X\sm Y$ respectively $\map_*(Y,Z)$.
So $\eta'$ and $\epsilon'$ are the unit respectively counit of an 
adjunction for the pair of endofunctors $(-\sm Y,\map_*(Y,-))$
on the category $\bT_*$.
\end{construction}

Now we specialize to the space of continuous homomorphisms
between Lie groups, with compact source. 
If $(Y,d)$ is a metric space and $K$ a compact topological space,
then the supremum metric
\[ d(f,g)\ = \ \sup_{x\in K}\, \{ d(f(x),g(x))\} \]
is a metric on the set $C(K,Y)$ of continuous maps.
The topology induced by the supremum metric 
agrees with the compact-open topology, see for example \cite[Prop.\,A.13]{hatcher},
and hence also with  the topology of Construction \ref{con:define C(X,Y)}.
Since the topology on $C(K,Y)$ is metrizable, the space $C(X,Y)$ 
is in particular compactly generated,
by Proposition \ref{prop:properties cgwh spaces}~(iii).
Hence $\map(X,Y)=C(X,Y)$, with topology induced by the supremum metric.

The topology of every second countable smooth manifold is metrizable,
so the previous discussion applies in particular when $K$ is a compact Lie group and 
$G$ is a Lie group.
So in that case the compact-open topology
on the space $C(K,G)$ of continuous maps is metrizable,
and hence compactly generated. Hence $C(K,G)=\map(K,G)$,
i.e., the compact-open topology is the topology of the internal mapping
space in the category $\bT$ of compactly generated spaces. 

We let $\hom(K,G)$ denote the set of continuous homomorphisms
with the subspace topology of $\map(K,G)$.
Since $\map(K,G)$ is metrizable, so is its subspace $\hom(K,G)$,
which is then compactly generated by Proposition \ref{prop:properties cgwh spaces}~(iii).
The next proposition in particular shows that the space $\hom(K,G)$
is the topological disjoint union of the orbits, under conjugation,
of the identity component group $G^\circ$.
The key input is a theorem of 
Montgomery and Zippin \cite[Thm.\,1 and Corollary]{montgomery-zippin}
from~1942 which roughly says that in a Lie group `nearby compact subgroups are conjugate'.
The following consequence of Montgomery and Zippin's result
appears in \cite[III, Lemma\,38.1]{conner-floyd}.

\begin{prop}\label{prop:components of hom(K,G)}
Let $K$ and $G$ be Lie groups, and suppose that $K$ is compact.
Then every orbit of the conjugation action by $G^\circ$ on $\hom(K,G)$
is an open subset of the space $\hom(K,G)$.
In particular, the connected components of the space $\hom(K,G)$ coincide
with its path components, and with the $G^\circ$-orbits under the conjugation action.
\end{prop}
\begin{proof}
We reproduce the argument from \cite[III, Lemma\,38.1]{conner-floyd}, in somewhat
expanded form. 
We start by showing that the $G^\circ$-orbits are open subsets of $\hom(K,G)$.
Let $\alpha:K\to G$ be a continuous homomorphism. 
Its graph $\Gamma_\alpha=\{(k,\alpha(k))\ |\ k\in K\}$ is then a compact subgroup
of the Lie group $K\times G$. 
By \cite[Thm.\,1 and Corollary]{montgomery-zippin},
there is an open subset $O$ of $K\times G$ containing $\Gamma_\alpha$
with the following property:
for every closed subgroup $\Delta$ of $K\times G$ with $\Delta\subseteq O$
there is an element $(k,g)\in K^\circ\times G^\circ$ such that
$(k,g)^{-1}\cdot \Delta\cdot(k,g)\subseteq \Gamma_\alpha$.
We set
\[ U \ = \ \{\beta\in\hom(K,G)\ | \ \Gamma_\beta\subseteq O \} \ ;\]
then $\alpha\in U$, and $U$ is contained in the $G^\circ$-orbit of $\alpha$ because
\[ (k,g)^{-1}\cdot \Gamma_{\beta}\cdot(k,g)\ = \ \Gamma_{c_g\circ c_{\beta(k)^{-1}}\circ\beta}\ . \]
We claim that $U$ is an open subset of $\hom(K,G)$.
To see this, we let $\beta:K\to G$ be a continuous homomorphism
such that $\Gamma_\beta\subset O$,
and $k\in K$. Since $O$ is open and $(k,\beta(k))\in O$,
there are open subsets $U_k\subset K$ and $V_k\subset G$ with
\[ (k,\beta(k))\ \subseteq \ U_k\times V_k \ \subseteq \ O \ . \]
Since $\beta$ is continuous, the set $\beta^{-1}(V_k)$ is open in $K$,
hence so is $U_k\cap \beta^{-1}(V_k)$; moreover, $k\in U_k\cap \beta^{-1}(V_k)$.
Since $K$ is compact, there is a compact neighborhood $L_k$ of $k$ 
that is contained in $U_k\cap \beta^{-1} (V_k)$.
Since $L_k$ is a neighborhood of $k$ and $K$ is compact, there
are finitely many points $k_1,\dots,k_n$ such that
$K=L_{k_1}\cup\dots\cup L_{k_n}$.
Since $\beta(L_{k_i})\subset V_{k_i}$ for all $1\leq i\leq n$, the 
homomorphism $\beta$ is contained in the set
\[ W\ = \ \bigcap_{i=1,\dots,n}  N(\text{incl}:L_{k_i}\to K,V_{k_i}) \]
which is an open subset of $\map(K,G)$.
So we are done if we can show that every $\beta\in W\cap \hom(K,G)$
satisfies $\Gamma_\beta\subset O$.
Given $k\in K$, there is an $i\in\{1,\dots,n\}$ with $k\in L_{k_i}$.
Thus $\beta(k)\in V_{k_i}$ by hypothesis on $\beta$, and so
\[ (k,\beta(k))\ \in \ L_{k_i}\times V_{k_i}\ \subseteq \ O \ . \]
Since $k$ was arbitrary, this proves that $\Gamma_\beta\subseteq O$.
Since the set $U$ is open in $\hom(K,G)$ and contained in the $G^\circ$-orbit 
of $\alpha$, the $G^\circ$-orbit of $\alpha$ is open in $\map(K,G)$.

Now we show that every $G^\circ$-orbit is path connected.
For every continuous homomorphism $\alpha:K\to G$, the map 
\[  G \times K \ \to \ G \ , \quad
(g,k)\ \longmapsto \ g^{-1}\cdot\alpha(k)\cdot g \]
is continuous since it can be written as a composite of a diagonal map,
an evaluation map and data from the group structure.
The adjoint 
\[  G \ \to \ \hom(K,G) \ , \quad  g\ \longmapsto \ c_g\circ \alpha \]
is thus continuous as well.
Since the component group $G^\circ$ is path connected and conjugation is continuous,
every $G^\circ$-orbit is path connected, so in particular connected. 
The $G^\circ$-orbits are open by the previous paragraph. Since the complement
of an orbit is a union of other orbits, the $G^\circ$-orbits are also closed.
So the $G^\circ$-orbits are open, closed and path connected; 
hence they coincide with the path components and the connected components.
\end{proof}

For a closed subgroup $H$ of a compact Lie group $G$,
we let $C_G H$ and $N_G H$ denote the centralizer respectively
normalizer of $H$ in $G$, and we write $W_G H=(N_G H)/H$ for the Weyl group.

\begin{prop}\label{prop:1-component of W_G H}
For every closed subgroup $H$ of a compact Lie group $G$,
the composite
\[ (C_G H)^\circ \ \xra{\text{\em incl}} \ 
(N_G H)^\circ\ \xra{\text{\em proj}} \ ( W_G H)^\circ \]
is surjective.
\end{prop}
\begin{proof}
We consider a continuous path $\omega:[0,1]\to W_G H$
starting at $\omega(0)= e H$.
The projection $N_G H\to (N_G H)/H=W_G H$ is a locally trivial fiber bundle,
so the path can be lifted to a continuous path
$\bar \omega:[0,1]\to N_G H$ with $\bar\omega(0)=e$.
Proposition \ref{prop:components of hom(K,G)} shows that the
image of the map
\begin{equation} \label{eq:H^circ_acts}
 H^\circ \ \to \ \hom(H,H)\ , \quad h \ \longmapsto \ c_h   
\end{equation}
is the path component of the identity.
Since $H^\circ$ is compact and $\hom(H,H)$ is a Hausdorff space, 
this map factors over a homeomorphism from
$H^\circ /(H^\circ\cap Z(H))$ onto the path component of the identity,
where $Z(H)$ is the center of $H$.
In particular, the map \eqref{eq:H^circ_acts}
is a locally trivial fiber bundle over the path component of the identity.
The path
\[ [0,1]\ \to \ \hom(H,H) \ , \quad t \ \longmapsto \ c_{\bar\omega(t)}\]
can thus be lifted to a continuous path $\nu:[0,1]\to H^\circ$ satisfying
\[  \nu(0)=e \text{\qquad and \qquad}c_{\bar\omega(t)}\ = \ c_{\nu(t)}\]
for all $t\in [0,1]$.
The second relation means that $\bar\omega(t)\cdot\nu(t)^{-1}$ centralizes $H$,
so 
\[ \bar\omega\cdot\nu^{-1}\ : \ [0,1]\ \to \ C_G H \]
is a path in the centralizer of $H$ that starts at the identity element.
Moreover,
\[ \bar\omega(1)\cdot\nu^{-1}(1)\cdot H \ =\ \bar\omega(1)\cdot H \ = \ \omega(1)\ .\qedhere\]
\end{proof}

The {\em h-cofibrations} are the morphisms with
the homotopy extension property. We will use this concept in various categories,
for example in the category of $G$-spaces, orthogonal spaces 
and orthogonal spectra. So we recall some basic properties of
h-cofibrations in the context of categories enriched over the category of spaces. 
The arguments are all standard and well known, and we include them
for completeness and convenience.

For the discussion of h-cofibrations we work in a cocomplete category $\Cc$ 
that is tensored and cotensored over the category $\bT$ of compactly generated spaces.
We write `$\times$' for the pairing
and $X^K$ for the cotensor of an object $X$ with a compactly generated space $K$.
A {\em homotopy} is then a morphism $H:A\times [0,1]\to X$
defined on the pairing of a $\Cc$-object with the unit interval.
For a homotopy and any $t\in [0,1]$ we denote by $H_t:A\to X$ the
composite morphism
\[ A \ \iso \ A\times\{t\}\ \xra{A\times\text{incl}} \ A\times [0,1]\
\ \xra{\ H \ } \ X \ . \]

\begin{defn}
Let $\Cc$ be a category tensored over the category $\bT$ of spaces. 
A $\Cc$-morphism $f:A\to B$ is an 
{\em h-cofibration}\index{subject}{h-cofibration} 
if it has the homotopy extension property, 
i.e., given a morphism $\varphi:B\to X$ and a homotopy
$H:A\times [0,1]\to X$ such that $H_0=\varphi f$,
there is a homotopy $\bar H:B\times [0,1]\to X$ such that
$\bar H\circ(f\times [0,1])=H$ and $\bar H_0=\varphi$.
\end{defn}

There is a universal test case for the homotopy extension problem, 
namely when $X$ is the pushout:
\[  \xymatrix{ 
A \ar[d]_f \ar[r]^-{(-,0)} & A\times [0,1] \ar[d]^H \\
B \ar[r]_-\varphi & B\cup_f  ( A\times [0,1] )  }\]
So a morphism $f:A\to B$ is an h-cofibration if and only if the canonical morphism
\begin{equation}\label{eq:canonical_h_cof_morphism} 
 B \cup_f ( A\times  [0,1] )\ \to\ B\times [0,1]
\end{equation}
has a retraction. Also, the adjunction between $-\times [0,1]$ and $(-)^{[0,1]}$
lets us rewrite any homotopy extension data $(\varphi,H)$ in adjoint form
as a commutative square:
\[ \xymatrix{ 
A\ar[d]_f \ar[r]^-{\hat H} & X^{[0,1]} \ar[d]^{\ev_0} \\
B \ar[r]_-{\varphi} & X} \]
A solution to the homotopy extension problem is adjoint to a lifting, 
i.e., a morphism $\lambda:B\to X^{[0,1]}$ such that 
$\lambda f=\hat H$ and $\ev_0\circ \lambda=\varphi$.
So a morphism $f:A\to B$ is an h-cofibration if and only if it 
has the left lifting property\index{subject}{lifting property!left}
with respect to the morphisms $\ev_0:X^{[0,1]}\to X$ for all
objects in $\Cc$.

The three equivalent characterizations of h-cofibrations quickly imply
various closure properties.

\begin{cor}\label{cor-h-cofibration closures} 
Let $\Cc$ be a cocomplete category tensored and cotensored over the category 
$\bT$ of compactly generated spaces.
\begin{enumerate}[\em (i)]
\item The class of h-cofibrations in $\Cc$ is closed under
retracts, cobase change, coproducts and sequential compositions.
\item Let $\Cc'$ be another category tensored over the category $\bT$,
and $F:\Cc\to\Cc'$ a continuous functor that commutes with colimits and
tensors with $[0,1]$. Then $F$ takes h-cofibrations
in $\Cc$ to  h-cofibrations in $\Cc'$.
\item If $\Cc$ is a topological model category in which every object is fibrant,
  then every cofibration is an h-cofibration.
\end{enumerate}
\end{cor}
\begin{proof}
(i) Every class of morphisms that can be characterized by the left lifting property 
with respect to some other class has the closure properties listed.

(ii) Let $f:A\to B$ be a cofibration in $\Cc$ and 
$r:B\times [0,1]\to B \cup_f (A\times [0,1])$ a retraction to the canonical morphism.
The composite
\begin{align*}
   F B\times [0,1]  \ \iso\  &F(B\times [0,1]) \ \xra{\ Fr\ }\\ 
&F(B\cup_f ( A\times [0,1]) ) \ \iso \   F B \cup_{Ff} ( F A\times [0,1] )
\end{align*}
is then a retraction to the canonical morphism for $Ff:F A\to F B$.
So $Ff$ is an h-cofibration.

(iii)
Since the model structure is topological, for every cofibration
$f:A\to B$ the canonical morphism \eqref{eq:canonical_h_cof_morphism}
is an acyclic cofibration.
Since every object is fibrant, this morphism has a retraction, 
and so $f$ is an h-cofibration.
\end{proof}

We turn to a key technical property of h-cofibrations in the category $\bT$
of compactly generated spaces, namely that these are always closed embeddings.
The same conclusion also holds for h-cofibrations in the category $\bT_*$
of based compactly generated spaces. \medskip

\Danger 
The forgetful functor $\bT_*\to \bT$ does {\em not} preserve the tensors with
unbased spaces, and it does {\em not} preserve h-cofibrations.
Indeed, the `based tensor' of a based space $(A,a_0)$ with an unbased space $K$ is 
$A\sm K_+$, whereas the unbased tensor of the underlying space is simply $A\times K$.
For every based space $(A,a_0)$ the inclusion 
$\{a_0\}\to A$ is an h-cofibration in the based sense, but not generally 
in the unbased sense. 
By definition, $(A,a_0)$ is well-pointed\index{subject}{well-pointed}
if this inclusion is an h-cofibration of unbased spaces.

\medskip

The following proposition is Lemma~8.2 in \cite[App.\,A]{lewis-thesis}.

\begin{prop}\label{prop:h-cof is closed embedding}
  Every h-cofibration between compactly generated spaces is a closed embedding.
  Every based h-cofibration between compactly generated based spaces 
  is a closed embedding.
\end{prop}
\begin{proof}
Let $f:A\to B$ be an h-cofibration between
compactly generated spaces.
Since the map $(-,0):A\to A\times[0,1]$ is a closed embedding, 
Proposition \ref{prop:pushout in cg} says that the
pushout $B\cup_f ( A\times[0,1])$ in the category $\bT$ can 
in fact be calculated in the ambient category $\bSpc$.
So the map
\[  (-,1)\ : \ A\ \to \ B\cup_f (A\times[0,1]) \ , \quad a \ \longmapsto \ (a,1)\]
is a closed embedding by direct inspection of the topology on the pushout.

The universal example of the homotopy extension property provides 
a continuous retraction to the canonical map 
\[ i\ = \  (-,0)\cup (f\times [0,1])\ :\ B\cup_f (A\times[0,1])\ \to \ B\times [0,1]\ . \]
So $i$ is a closed embedding by Proposition \ref{prop:section is closed embedding}.
In the commutative square
\[ \xymatrix@C=12mm{ 
 A\ar[r]^-{(-,1)}\ar[d]_f & B\cup_f (A\times[0,1]) \ar[d]^i \\
 B\ar[r]_-{(-,1)} &  B\times[0,1]} \]
the two horizontal and the right vertical maps are thus closed embeddings. 
So the map $f$ is a closed embedding as well.

The based case, albeit logically independent, proceeds along the same lines.
We let $f:(A,a_0)\to (B.b_0)$ 
be a based h-cofibration between compactly generated based spaces.
Since points in weak Hausdorff spaces are closed
(Proposition \ref{prop:properties wH spaces}~(iii)), 
the set $\{a_0\}\times [0,1]$ is closed in $A\times [0,1]$, and so the map
\[ (-,t)\ :\ A\to A\sm [0,1]_+\ , \quad a\ \longmapsto \ a\sm t \]
is a closed embedding for every $t\in [0,1]$.
Any pushout of based spaces can be calculated in the ambient category
of unbased spaces; so Proposition \ref{prop:pushout in cg} says that the
pushout $B\cup_f ( A\sm [0,1]_+)$ in the category $\bT_*$ can 
in fact be calculated in the ambient category $\bSpc$.
So the map
\[  (-,1)\ : \ A\ \to \ B\cup_f (A\sm [0,1]_+) \ , \quad a \ \longmapsto \ a\sm 1 \]
is a closed embedding by direct inspection.

The universal example of the homotopy extension property provides 
a continuous retraction to the canonical based map 
\[ i\ = \  (-,0)\cup (f\sm [0,1]_+)\ :\ B\cup_f (A\sm [0,1]_+)\ \to \ B\sm [0,1]_+\ . \]
So $i$ is a closed embedding by Proposition \ref{prop:section is closed embedding}.
In the commutative square
\[ \xymatrix@C=15mm{ 
 A\ar[r]^-{(-,1)}\ar[d]_f & B\cup_f (A\sm[0,1]_+) \ar[d]^i \\
 B\ar[r]_-{(-,1)} &  B\sm[0,1]_+} \]
the two horizontal and the right vertical maps are thus closed embeddings. 
So the map $f$ is a closed embedding as well.
\end{proof}

Now we turn to geometric realization of simplicial spaces,
which is a frequent tool to construct interesting homotopy types.

\begin{construction}\label{con:realize simplicial space}
We recall the geometric realization of a simplicial space.
Geometric realization was originally introduced 
by Milnor \cite{milnor-realization} for simplicial sets 
(which were called `semi-simplicial complexes' at that time);
the version for simplicial spaces is a straightforward generalization.
We let $\bDelta$ denote
the simplicial indexing category, with objects the finite
totally ordered sets $[n]=\{0\leq 1\leq \dots\leq n\}$
for $n\geq 0$. Morphisms in $\bDelta$ are all weakly monotone maps.
We let
\[ \Delta^n \ = \ \{ (t_1,\dots,t_n)\in [0,1]^n \ |\ t_1\leq t_2\leq \dots\leq t_n \} \]
be the topological $n$-simplex.
As $n$ varies, these topological simplices assemble into a covariant functor
\[ \Delta^{\bullet}\ : \ \bDelta \ \to \bSpc \ , \quad [n]\ \longmapsto \ \Delta^n\ ; \]
the coface maps are given by
\[ (d_i)_*(t_1,\dots,t_n) \ = \left\lbrace \begin{array}{ll}
\quad (0,t_1,\dots,t_n) & \mbox{for $i=0$,} \\
(t_1,\dots,t_i,t_i,\dots,t_n)  & \mbox{for $0<i< n$,} \\
\quad (t_1,\dots,t_n,1) & \mbox{for $i=n$.}
\end{array} \right.  \]
The codegeneracy map $(s_i)_*:\Delta^n\to\Delta^{n-1}$ drops the entry $t_{i+1}$.

A {\em simplicial space}\index{subject}{simplicial space}
is functor $X:\bDelta^{\op}\to\bSpc$, i.e., a contravariant functor from $\bDelta$
to the category of topological spaces.
We use the customary notation $X_n=X([n])$ for the value of a simplicial space
at $[n]$. 
The {\em geometric realization}\index{subject}{geometric realization!of a simplicial space}
of $X$ is
\[ |X|\ = \ \left( {\coprod}_{n\geq 0}\, X_n\times_0 \Delta^n \right) / \ \sim\ ,\]
the quotient space of the disjoint union by the equivalence relation generated by
\begin{equation}  \label{eq:generate_realization}
 (x,\alpha_*(t))\ \sim \ (\alpha^*(x),t)   
\end{equation}
for all morphisms $\alpha:[n]\to [m]$ in $\bDelta$ 
and all $(x,t)\in X_m\times_0 \Delta^n$.
A more categorical way to say this is that $|X|$ is a coend of the functor
\begin{equation}\label{eq:Spc coend}
 \bDelta^{\op}\times\bDelta\ \to \ \bSpc\ , \quad 
([m],[n])\ \longmapsto \ X_m\times_0\Delta^n \ .  
\end{equation}
\end{construction}

We recall that the equivalence relation generated by \eqref{eq:generate_realization}
is well understood in terms of `minimal representatives':
\begin{itemize}
\item  every equivalence class has a unique 
representative $(x,t)\in X_l\times\Delta^l$ of minimal dimension $l$;
\item an element $(x,t)\in X_l\times\Delta^l$ is the minimal representative in
its equivalence class if and only if $x$ is non-degenerate 
(i.e., not in the image of any degeneracy map $s_i^*:X_{l-1}\to X_l$)
and $t$ belongs to the interior of $\Delta^l$;
\item if $(x,t)\in X_l\times\Delta^l$ is a minimal representative 
and $(y,s)\in X_m\times \Delta^m$ is equivalent to $(x,t)$, 
then there is a unique surjective morphism $\sigma:[k]\to[l]$, 
a unique injective morphism $\delta:[k]\to[m]$
and a unique $u\in\Delta^k$ such that
\[ \delta^*(y)\ = \ \sigma^*(x)\ , \quad 
s\ = \  \delta_*(u) \text{\quad and\quad} t =\sigma_*(u) \ .\]
\end{itemize}
These facts go back all the way to Milnor \cite[Lemma 3]{milnor-realization}.

The next proposition shows that compactly generated spaces are 
`closed under geometric realization';
more precisely, we consider a simplicial space $X_\bullet$ 
such that $X_n$ is compactly generated for every $n\geq 0$.
The geometric realization, formed in the ambient category $\bSpc$, 
is an enriched coend, hence automatically a $k$-space.
Since colimits of compactly generated spaces are not, in general,
formed in the ambient category, it is not completely obvious, though,
that the ambient realization is again a weak Hausdorff space.
I learned about this fact as a parenthetical remark 
in the introduction of the paper \cite{seguinspazzis-realize}.

\begin{prop}\label{prop:geometric realization in bT}
  Let $X:\bDelta^{\op}\to \bK$ be a simplicial $k$-space.
  \begin{enumerate}[\em (i)]
  \item The the geometric realization $|X|$ is a $k$-space,
    and hence also a coend internal to the category $\bK$,
    of the functor \eqref{eq:Spc coend}.
  \item If $X_n$ is compactly generated for every $n\geq 0$, 
    then the geometric realization $|X|$ is compactly generated,
    and hence also a coend internal to the category $\bT$,
    of the functor \eqref{eq:Spc coend}.
  \item Let $X_n$ be compactly generated for all $n\geq 0$.
    Let $Y$ be a simplicial subspace of $X$ such that
    $Y_n$ is closed in $X_n$ for all $n\geq 0$. 
    Then the inclusion induces a closed embedding $|Y|\to |X|$.
  \end{enumerate}
\end{prop}
\begin{proof}
(i) Since $\Delta^n$ is compact and $X_n$ is a $k$-space,
$X_n\times_0\Delta^n=X\times\Delta^n$ is a $k$-space in the product topology 
by Proposition \ref{prop:properties k-spaces}~(vi),
and so the disjoint union $\coprod_{n\geq 0}\, X_n\times \Delta^n$
is a $k$-space.
As a quotient space, the geometric realization $|X|$ is then 
a $k$-space by Proposition \ref{prop:properties k-spaces}~(i).

(ii) As we argued in part~(i),
the space ${\coprod}_{n\geq 0} X_n\times_0 \Delta^n={\coprod}_{n\geq 0} X_n\times\Delta^n$
is a $k$-space.
We can thus use the criterion given by Proposition \ref{prop:closed implies WH}
to show that the quotient space $|X|$ is compactly generated.
We let $E\subset ({\coprod}_{n\geq 0} X_n\times \Delta^n)^2$
be the equivalence relation generated by \eqref{eq:generate_realization},
which we had simply denoted `$\sim$' above. We will show that
$E$ is closed in the $k$-topology of $({\coprod}_{n\geq 0} X_n\times \Delta^n)^2$.
Since products distribute over disjoint unions, we may show that
\[ E_{m,n}\ = \ E\cap (X_m\times\Delta^m\times X_n\times\Delta^n) \]
is closed in $X_m\times\Delta^m\times X_n\times\Delta^n$ for all $m,n\geq 0$.
If $(y,s,\bar y,\bar s)\in E_{m,n}$, then $(y,s)$ and $(\bar y,\bar s)$
have the same minimal representative $(x,t)\in X_l\times\Delta^l$.
So there are surjective morphisms
$\sigma:[k]\to[l]$ and $\bar\sigma:[\bar k]\to[l]$, 
injective morphisms $\delta:[k]\to[m]$ and $\bar\delta:[\bar k]\to[n]$,
and $u\in\Delta^k$, $\bar u\in\Delta^{\bar k}$ such that
\[ \delta^*(y)\ = \ \sigma^*(x)\ , \quad 
s\ = \  \delta_*(u) \text{\quad and\quad} t =\sigma_*(u) \]
and
\[ \bar\delta^*(\bar y)\ = \ \bar \sigma^*(x)\ , \quad 
\bar s\ = \  \bar \delta_*(\bar u) \text{\quad and\quad} t =\bar \sigma_*(\bar u) \ .\]
Hence $E_{m,n}$ is the union,
indexed over $l,\sigma,\bar\sigma,\delta,\bar\delta$ as above,
of the finite number of sets
\[
 (\delta^*\times\Delta^m\times\bar\delta^*\times\Delta^n)^{-1}\left(
 (\sigma^*\times\delta_*\times\bar\sigma^*\times\bar\delta_*)
\left((X_l\times\sigma_*\times X_l\times\bar\sigma_*)^{-1}
( \Delta_{X_l\times\Delta^l})\right)\right) \ .\]
The diagonal $\Delta_{X_l\times\Delta^l}$ is closed in $(X_l\times\Delta^l)^2$
because $X_l$, and hence $X_l\times\Delta^l$, is compactly generated.
So its inverse image under the continuous map
$X_l\times\sigma_*\times X_l\times\bar\sigma_*$
is closed in $X_l\times\Delta^k\times X_l\times\Delta^{\bar k}$.
Every surjective morphism in $\bDelta$ has a section,
so $\sigma^*:X_l\to X_k$ and $\bar\sigma^*:X_l\to X_{\bar k}$ 
have continuous retractions.
The maps $\delta_*:\Delta^k\to\Delta^m$ and $\bar\delta_*:\Delta^{\bar k}\to\Delta^n$ 
also have continuous retractions, hence so does
$\sigma^*\times\delta_*\times\bar\sigma^*\times\bar\delta_*$.
This map is then a closed embedding by Proposition \ref{prop:section is closed embedding}.
So the set 
\[  (\sigma^*\times\delta_*\times\bar\sigma^*\times\bar\delta_*)
\left((X_l\times\sigma_*\times X_l\times\bar\sigma_*)^{-1}
( \Delta_{X_l\times\Delta^l})\right)  \]
is closed in  $X_k\times\Delta^m\times X_{\bar k}\times\Delta^n$.
Since $E_{m,n}$ is the inverse image of this latter closed set
under a continuous map, this show the claim that $E_{m,n}$
is a closed subset of $X_m\times\Delta^m\times X_n\times\Delta^n$.

  (iii) We write $\iota:Y\to X$ for the inclusion. 
  Our first claim is that the induced map $|\iota|:|Y|\to |X|$ is injective.
  We consider a point $(y,s)\in Y_m\times \Delta^m$ 
  and we let $(x,t)\in X_l\times \Delta^l$ be the minimal representative 
  of the equivalence class of $(y,s)$ 
  in the ambient space $\coprod_{n\geq 0}X_n\times\Delta^n$.
  Then there is a surjective morphism
  $\sigma:[k]\to[l]$, an  injective morphism $\delta:[k]\to[m]$ 
  and $u\in\Delta^k$ such that
  \[ \delta^*(y)\ = \ \sigma^*(x)\ , \quad 
  s\ = \  \delta_*(u) \text{\quad and\quad} t =\sigma_*(u) \ .\]
  We let $\bar\delta:[l]\to [k]$ be a morphism in $\bDelta$ 
  such that $\sigma\bar\delta=\Id_{[l]}$. 
  Then $\bar\delta^*(\delta^*(y))=\bar\delta^*(\sigma^*(x))=x$.
  Since $Y$ is a simplicial subspace of $X$, this shows that $x$ belongs to $Y_l$.
  If $(\bar y,\bar s)\in Y_n\times \Delta^n$ is another pair that
  represents the same point in $|X|$ as $(y,s)$, then
  by the above, the minimal representative of the equivalence class belongs
  to $\coprod_{n\geq 0} Y_n\times\Delta^n$, and so $(y,s)$ and $(\bar y,\bar s)$
  already represent the same point in $|Y|$. So the map $|\iota|$ is injective. 

  It remains to show that the continuous injection $|\iota|$ is a closed map.
  We consider the commutative square
  \[ \xymatrix@C=15mm{
    \coprod_{n\geq 0} Y_n\times \Delta^n\ar[r]^-{\coprod \iota_n\times\Delta^n} \ar[d]_p&
    \coprod_{n\geq 0} X_n\times  \Delta^n\ar[d]^q\\
    |Y|\ar[r]_-{|\iota|} & |X| } \]
  where the vertical maps are the quotient maps.
  We let $(y,s)\in X_m\times\Delta^m$ be a point whose equivalence class
  lies in the image of $|\iota|:|Y|\to|X|$.
  As we showed in the previous paragraph, the representative 
  $(x,t)\in X_l\times \Delta^l$ 
  of minimal dimension in the equivalence class of $(y,s)$ must then
  lie in the simplicial subspace $Y$, i.e., we must have $x\in Y_l$.
  Moreover, there is a surjective morphism
  $\sigma:[k]\to[l]$, an  injective morphism $\delta:[k]\to[m]$ 
  and $u\in\Delta^k$ such that
  \[ \delta^*(y)\ = \ \sigma^*(x)\ , \quad 
  s\ = \  \delta_*(u) \text{\qquad and\qquad} t =\sigma_*(u) \ .\]
  So for every subset $A\subset |Y|$, we have
   \begin{align*}
  ( X_m\times \Delta^m) \cap &q^{-1}( |\iota|(A) )\\ 
  = \ 
  {\bigcup}_{\sigma,\delta}\,  &  (\delta^*\times \Delta^m)^{-1}
\left(( \sigma^*\times\delta_*)
   (  ( X_l\times\sigma_*)^{-1}(( Y_l\times \Delta^l)\cap p^{-1}(A) ) ) \right)\ .
  \end{align*}
  The union is over the finite set of pairs consisting of a surjective 
  morphism $\sigma:[k]\to [l]$ and an injective morphism $\delta:[k]\to [m]$.

  Now we assume that $A$ is closed inside $|Y|$. Because $p$ is continuous
  and $Y_l$ is closed in $X_l$, the set
  $(Y_l\times \Delta^l)\cap p^{-1}(A)$ is closed inside $X_l\times \Delta^l$.
  So
  $(X_l\times\sigma_*)^{-1}(( Y_l\times \Delta^l)\cap p^{-1}(A))$ 
  is a closed subset of $X_l\times \Delta^k$.
  Since the map $\sigma^*\times\delta_*:X_l\times\Delta^k\to X_k\times\Delta^m$ 
  has a continuous retraction, it is a 
  closed embedding by Proposition \ref{prop:section is closed embedding}.
  So the set $ (\sigma^*\times \delta_*) 
  (  ( X_l\times\sigma_*)^{-1}((Y_l\times \Delta^l)\cap p^{-1}(A) ))$
  is closed in $X_k\times \Delta^m$.
  As the inverse image under a continuous map, the set
  \[ (\delta^*\times \Delta^m)^{-1} ( (\sigma^*\times\delta_*) 
  (  ( X_l\times\sigma_*)^{-1}(( Y_l\times \Delta^l)\cap p^{-1}(A) ) ) ) \]
  is then closed in $X_m\times \Delta^m$.

  So each set in the finite union above is closed inside $X_m\times\Delta^m$.
  We conclude that $( X_m\times \Delta^m) \cap q^{-1}( |\iota|(A) )$
  is closed in $X_m\times \Delta^m$ for every $m\geq 0$, 
  hence the set $q^{-1}(|\iota|(A))$ is closed.
  Since $q$ is a quotient map, this shows that $|\iota|(A)$ is closed in $|X|$.
\end{proof}

The next proposition is about the interaction of geometric realization
and products for simplicial $k$-spaces and simplicial compactly generated spaces.
The essential input here is that the categories $\bK$ and $\bT$
are cartesian closed, hence product with a fixed object
is a left adjoint and commutes with colimits and coends.
A direct consequence is that geometric realization 
commutes with colimits and products with a fixed space,
see part~(i) and~(ii) of Proposition \ref{prop:iterated geometric realization} below.

We also recall that different ways to realize a multi-simplicial space are homeomorphic.
To this end we consider a 
{\em bisimplicial $k$-space},\index{subject}{bisimplicial space} i.e., a functor
$Z:\bDelta^{\op}\times\bDelta^{\op}\to \bK$.
The diagonal simplicial $k$-space $\diag Z$ is the composite 
with the diagonal functor $\bDelta^{\op}\to\bDelta^{\op}\times\bDelta^{\op}$.
On the other hand, for every fixed `horizontal' dimension $m$ 
we can define the simplicial space $\partial_m Z$ as the composite 
\[ \bDelta^{\op}\ \xra{\ (\Id_{[m]},-) \ }\ \bDelta^{\op}\times\bDelta^{\op}
\ \xra{\ Z \ } \bK\ , \]
so that $(\partial_m Z)_n=Z_{m,n}$. For varying $m$, these define a simplicial object
in the category of simplicial $k$-spaces; so the assignment
\[ [m]\ \longmapsto \ | \partial_m Z| \]
becomes a simplicial $k$-space. We can geometrically realize this simplicial space
and arrive at the {\em iterated realization}\index{subject}{iterated realization!of a bisimplicial space} that we denote
\begin{equation}  \label{eq:def_iterated_real}
 |Z|^{\text{it}}\ = \ |[m]\mapsto |\partial_m Z| | \ .  
\end{equation}
The diagonal realization and the iterated realization are related
by a natural continuous map, see \eqref{eq:diagonal2iterated} below.
We show in Proposition \ref{prop:iterated geometric realization}~(iii)
that this preferred natural map is a homeomorphism.

I am not sure who deserves credit for parts~(ii) and~(iii)
of the following proposition.
For preservation of products, Milnor \cite[Thm.\,2]{milnor-realization} 
already treats the special case of
simplicial sets (which we can view as discrete simplicial spaces),
where he imposes size conditions to ensure that the ordinary product 
of the two realizations is a $k$-space.
The earliest reference in the present generality I know of 
is \cite[\S 1, p.\,106]{segal-classifying spaces}, 
where Segal states the result without proof.
May \cite[Thm.\,11.5]{may-geometry iterated loop} gives a proof
in the category of Hausdorff $k$-spaces (which he calls 
`compactly generated Hausdorff spaces'), but the proof does not use 
the Hausdorff property. 

\begin{prop}\label{prop:iterated geometric realization}
  \begin{enumerate}[\em (i)]
  \item 
    For every simplicial $k$-space $Y:\bDelta^{\op}\to \bK$ and every $k$-space $K$,
    the canonical map
    \[  |K\times Y|\ \to \ K\times |Y| \ , \quad [k,y,t ]\ \longmapsto \ (k,[y,t])\]
    is a homeomorphism.
  \item
    The realization functors  
    \[ |-|\ :\ \bDelta^{\op}\bK\ \to\ \bK \text{\qquad and\qquad}
    |-|\ :\ \bDelta^{\op}\bT \ \to\ \bT \]
    for simplicial $k$-spaces respectively 
    simplicial compactly generated spaces 
    preserves colimits and finite products.
  \item
    For every bisimplicial $k$-space $Z:\bDelta^{\op}\times\bDelta^{\op}\to \bK$, 
    the map 
    \begin{equation} \label{eq:diagonal2iterated}
      \delta\ : \ |\diag Z| \ \to \ |Z|^{\text{\em it}}\text{\quad given by\quad}
      \delta[z,t] \ = \ [ [z,t],t]   
    \end{equation}
    is a homeomorphism.
  \end{enumerate}
\end{prop}
\begin{proof}
(i)
The category of $k$-spaces is cartesian closed
(Theorem \ref{thm:internal hom in K}), so product with $K$ 
is a left adjoint and commutes with coends. Hence the canonical map
\[ 
|K\times Y |\ = \  \int^{[n]\in\bDelta} \left( K\times Y_n\times\Delta^n\right)
\ \to\
 K\times \left( \int^{[n]\in\bDelta} Y_n\times \Delta^n\right) \ = \ K\times |Y|
\]
is a homeomorphism. 

(ii) We start by showing that the realization functors preserve colimits.
The argument is entirely formal, and literally the same for the
categories $\bK$ and $\bT$, so we only treat the former case.
We let $I$ be a small category and $F:I\to \bDelta^{\op}\bK$ a functor.
The category of $k$-spaces is cartesian closed
(Theorem \ref{thm:internal hom in K}), so product with $\Delta^n$ 
is a left adjoint and commutes with coends. 
Moreover, colimits commute with coends, so the canonical maps
\begin{align*}
\colim_I &\left( |-|\circ F\right)\ = \ 
\colim_I  \left( \int^{[n]\in\bDelta} F_n\times \Delta^n\right)\
\xra{\ \iso \ }\  \int^{[n]\in\bDelta}  \colim_I (F_n \times\Delta^n) \\
&\xra{\ \iso \ }\  \int^{[n]\in\bDelta} ( \colim_I F_n )\times\Delta^n\
=\  \int^{[n]\in\bDelta} ( \colim_I F)_n \times\Delta^n\ = \ |\colim_I F|  
\end{align*}
are homeomorphisms. 

We postpone the proof that the realization functors preserve finite products
until later.

(iii) 
In any presheaf category with values in $\bK$, every functor is canonically
a colimit of representable functors times a fixed space.
For the category $\bDelta\times\bDelta$, 
this means that every bisimplicial $k$-space can be expressed as a colimit of
bisimplicial spaces of the form $(\bDelta\times\bDelta)(-, ([m],[n]))\times K$
for $m,n\geq 0$ and $k$-spaces $K$.
Passage to the diagonal clearly commutes with colimits and products,
so source and target of the map $\delta$ commute with colimits in $\bK$ 
and products with a fixed $k$-space, by parts~(i) and~(ii).
This reduces the claim to the special case 
$Z=\Delta[m,n]=(\bDelta\times\bDelta)(-, ([m],[n]))$,
the bisimplicial discrete space represented by the object $([m],[n])$.
The diagonal of this bisimplicial space is $\Delta[m]\times\Delta[n]$,
the product of the represented simplicial sets $\Delta[m]$ and $\Delta[n]$.

On the other hand, for fixed $k\geq 0$, the simplicial discrete space 
$\partial_k \Delta[m,n]$ is isomorphic to $\bDelta([k],[m])\times \Delta[n]$,
and so its realization is homeomorphic to $\bDelta([k],[m])\times |\Delta[n]|$.
The simplicial space $|\partial_\bullet \Delta[m,n] |$ is thus isomorphic
to $\Delta[m]\times |\Delta[n]|$. Since product with $|\Delta[n]|$
commutes with realization by part~(i), we conclude that
the iterated realization $|\Delta[m,n]|^{\text{it}}$ is homeomorphic to the 
product $|\Delta[m]|\times |\Delta[n]|$.
Moreover, under these identifications, the map $\delta$ for $Z=\Delta[m,n]$ specializes
to the canonical map $|\Delta[m]\times \Delta[n]|\to|\Delta[m]|\times|\Delta[n]|$.
It is a classical fact,
already observed by Milnor \cite[Thm.\,2]{milnor-realization},
that this canonical map is a homeomorphism;
other references are \cite[Ch.\,III 3.4]{gabriel-zisman}
and \cite[Lemma 3.1.8]{hovey-book}.

We still need to show that the realization functors preserve finite products.
The category $\bT$ of compactly generated space is closed under products
in the ambient category $\bK$ of $k$-spaces (because the inclusion $\bK\to\bT$
is a left adjoint), and closed under realization 
by Proposition \ref{prop:geometric realization in bT}~(ii).
So it suffices to treat realization of simplicial $k$-spaces.

We consider the bisimplicial $k$-space $X\bar\times Y$ defined as the composite
\[ \bDelta^{\op}\times\bDelta^{\op}\ \xra{X\times Y}\ \bK\times\bK \ \xra{\ \times\ }\
\bK\ ; \]
then $X\times Y=\diag(X\bar\times Y)$ is its diagonal.
The canonical map
\[ 
|\partial_m (X\bar\times Y) | \ = \ 
| X_m\times Y |
\ \to\ X_m\times \left( \int^{[n]\in\bDelta} Y_n\times \Delta^n\right)
\ = \  X_m\times |Y|
\]
is a homeomorphism by part~(i). 
Taking realizations over varying $m$ and using part~(i) again 
gives two homeomorphisms
\[ 
|X\bar\times Y|^{\text{it}} \ = \ 
|[m]\mapsto | X_m\times Y | |
\ \xra{\ \iso \ }  \ |X \times |Y| | \ \xra{\ \iso\ }\ 
|X |\times |Y| \ .\]
Under these identifications, the canonical map
$(|p_X|,|p_Y|): |X\times Y|\ \to \ |X|\times  |Y|$
becomes the map $\delta:|\diag(X\bar\times Y)|\to |X\bar\times Y|^{\text{it}}$
of part~(iii). Since $\delta$ is a homeomorphism, this shows the claim.
\end{proof}

\Danger For general simplicial topological spaces, i.e., 
functors $X,Y:\bDelta^{\op}\to \bSpc$ geometric realization need not
commute with products in $\bSpc$. Examples already arise when $X$ and $Y$ 
are simplicial sets (considered as simplicial discrete spaces), 
in which case $|X\times Y|$ is compactly generated. 
However, if $X$ and $Y$ are sufficiently large 
(i.e., neither countable nor locally finite), then the product topology
need not make $|X|\times_0 |Y|$ into a $k$-space. A specific example
is given by taking both $X$ and $Y$ as wedges of the simplicial 1-sphere, 
where $X$ has countably infinitely many copies,
and $Y$ has uncountably many copies, 
compare \cite[III.5, p.\,563]{dowker-metric complexes}.

\medskip

Now we discuss the latching spaces\index{subject}{latching space!of a simplicial space} 
of a simplicial space $X$, of which there are competing definitions in the literature.
The definition of the $n$-th latching object that we adopt is 
as a colimit of a certain functor whose values are the subspaces $X_i$ for $i<n$.
A possible point of confusion is that different references work
in slightly different categories (such as general topological spaces, 
$k$-spaces, Hausdorff $k$-spaces or weak Hausdorff $k$-spaces), 
and the presence of a Hausdorff or weak Hausdorff condition
affects the meaning of `colimit'.
Some authors define the $n$-th latching space of $X$ as a subspace of $X_n$,
namely the union of the subspaces $s_i^*(X_{n-1})$ for $i=0,\dots,n-1$,
where $s_i^*:X_{n-1}\to X_n$ is the $i$-th degeneracy map.
I suspect that the different definitions are not generally homeomorphic 
in the context of simplicial topological spaces, 
i.e., functors $X:\bDelta^{\op}\to\bSpc$.

As we will discuss now, all these definitions in fact coincide for
simplicial compactly generated spaces, i.e., functors $\bDelta^{\op}\to\bT$.
We already explained in Proposition \ref{prop:geometric realization in bT}~(ii)
that then the geometric realization (formed in the ambient
category $\bSpc$ of all topological spaces) is automatically compactly generated,
and hence also a coend internal to the category $\bT$.
We will now show that also the latching objects, 
defined in the ambient category $\bSpc$,
are automatically compactly generated, and hence also a colimit in the category $\bT$.
Moreover, the latching map $l_n:L_n X\to X_n$ 
is a closed embedding, with the expected image.

\begin{construction}
For $n\geq 0$ we let $\bDelta(n)$ denote category with objects the weakly 
monotone surjections $\sigma:[n]\to [k]$;
a morphism from $\sigma:[n]\to [k]$ to $\sigma':[n]\to [k']$
is a morphism $\alpha:[k]\to [k']$ in $\bDelta$ with $\alpha\circ\sigma=\sigma'$.
Such a morphism $\alpha$, if it exists, is necessarily unique and also surjective.
We let $\bDelta(n)_\circ$ denote the full subcategory with all objects
{\em except} the identity of $[n]$.

A simplicial topological space $X:\bDelta^{\op}\to \bSpc$ 
can be restricted along the forgetful functor
\[ \bDelta(n)_\circ^{\op}\ \xra{\ u \ } \ \bDelta^{\op} \ ,\quad 
(\sigma:[n]\to [k])\ \longmapsto \ [k]\ , \quad \alpha\ \longmapsto\ \alpha \ .\]
The $n$-th {\em latching space} of $X$ is the colimit of the composite
functor $X\circ u:\bDelta(n)_\circ^{\op}\to \bSpc$.
As $\sigma$ ranges over the objects of $\bDelta(n)^{\op}_\circ$, the maps
\[ \sigma^* \ : \ X_k\ \to \ X_n \]
assemble into a continuous map $l_n:L_n X\to X$,
the $n$-th {\em latching map}.
\end{construction}

 \begin{rk}\label{rk:Delta(n) vs P(n)}
The category $\bDelta(n)^{\op}$ is isomorphic to the poset category $\Pc(n)$ 
of subsets of the set $\{1,\dots,n\}$.
Indeed, an isomorphism $\kappa:\bDelta(n)^{\op}\to\Pc(n)$ is given on objects by
\[ \kappa(\sigma:[n]\to [k]) \ = \ 
\{ i\in\{1,\dots, n\}\ : \ \sigma(i) > \sigma(i-1)\} \ .\]
In the other direction, a subset $U\subset\{1,\dots,n\}$ 
is taken to  the monotone surjection $\kappa^{-1}(U):[n]\to[|U|]$
defined by
\[ \kappa^{-1}(U)(i)\ = \ |U\cap \{1,\dots,i\}|\ . \]
Since $\kappa(\Id_{[n]})=\{1,\dots,n\}$, 
the subcategory $\bDelta(n)^{\op}_\circ$ is taken
to the poset $\Pc(n)^\circ$ of proper subsets of $\{1,\dots,n\}$.

The maximal elements of the poset $\Pc(n)^\circ$ correspond to the morphisms
\[ s_i\ : \ [n]\ \to \ [n-1] \ , \quad 
  s_i(j)\ = \ 
\begin{cases}
\ j & \text{for $0\leq j\leq i$, and}\\
j-1 & \text{for $i+1\leq j\leq n$.}
\end{cases}\]
So the latching space $L_n X$ can also be presented 
as a coequalizer of two continuous maps
\begin{equation}  \label{eq:L_n_as_coequalizer}
\xymatrix{
 \coprod_{0\leq i < j\leq n-1} \, X_{n-2} 
\quad \ar@<.4ex>[r]^-{\alpha}  \ar@<-.4ex>[r]_-{\beta} &
\quad \coprod_{0\leq i\leq  n-1} \, X_{n-1}\ . }  
\end{equation}
The map $\alpha$ sends the $(i,j)$-th summand to the $i$-summand by the degeneracy
map $s_{j-1}^*:X_{n-2}\to X_{n-1}$.
The map $\beta$ sends the $(i,j)$-th summand to the $j$-summand by the degeneracy
map $s_i^*:X_{n-2}\to X_{n-1}$.
For example, the category $\bDelta(0)_\circ$ is empty, so $L_0(X)$ is empty.
The category $\bDelta(1)_\circ$ has a unique object $s_0:[1]\to [0]$,
so $L_1(X) = X_0$ and the latching map is given by $s_0^*:X_0\to X_1$.
The category $\bDelta(2)_\circ$ has three objects and two non-identity morphisms,
and $L_2(X)$ is a pushout of the diagram
\[ X_1 \ \xla{\ s_0^* \ }\ X_0 \ \xra{\ s_0^*\ } X_1 \ .\]
 \end{rk}

\begin{prop}\label{prop:sspaces latching properties}
Let $X:\bDelta^{\op}\to\bSpc$ be a simplicial topological space and $n\geq 0$.  
\begin{enumerate}[\em (i)]
\item The continuous latching map $l_n:L_n X\to X_n$ is injective
and its image is the union of the sets $s_i^*(X_{n-1})$ for $i=0,\dots,n-1$.
\item Suppose that the degeneracy map $s_i^*:X_{n-1}\to X_n$ is
a closed embedding for all $0\leq i\leq n-1$.
Then the latching map $l_n:L_n X\to X_n$ is a closed embedding.
\item Suppose that the space $X_n$ is compactly generated for every $n\geq 0$.
Then the $n$-latching space $L_n X$ formed in the ambient category $\bSpc$
of topological spaces is compactly generated, and hence a latching object
of $X$ internal to the category $\bT$.
Moreover, the latching map $l_n:L_n X\to X_n$ is a closed embedding.
\end{enumerate}
\end{prop}
\begin{proof}
(i) This argument is purely combinatorial and does not use the topology in any way.
For easier book keeping we label the different summands in the source
and target of the coequalizer diagram \eqref{eq:L_n_as_coequalizer}
as $X_{n-2}^{[i,j]}$ respectively $X_{n-1}^{[i]}$.
We let $q:\coprod_{0\leq i\leq n-1}X_{n-1}^{[i]}\to L_n X$ denote the quotient map.
We consider $x\in X_{n-1}^{[i]}$ and $y\in X_{n-1}^{[j]}$ such that
$l_n(q(x))=l_n(q(y))$ in $X_n$. This means that $s_i^*(x)=s_j^*(y)$.
We assume without loss of generality that $i\leq j$.
We write the elements as degeneracies of non-degenerate elements, 
i.e., $x=\sigma^*(z)$ and $y=\bar\sigma^*(\bar z)$
for surjective morphisms $\sigma:[n-1]\to[k],\bar\sigma:[n-1]\to[\bar k]$ 
and non-degenerate simplices $z$ and $\bar z$. Then
\begin{align*}
(\sigma s_i)^*(z)\ = \   s_i^*(\sigma^*(z))\ = \ s_i^*(x)\ =\ s_j^*(y) \ = \  
s_j^*(\bar\sigma^*(\bar z))\ = \ (\bar\sigma s_j)^*(\bar z)\ .  
\end{align*}
By the `Eilenberg-Zilber lemma' (\cite[(8.3)]{eilenberg-zilber-ssc}, 
see also \cite[Sec.\,II.3]{gabriel-zisman}),
the representation of a simplex as a degeneracy of a non-degenerate
element is unique, so $k=\bar k$, $z=\bar z$ 
and $\sigma s_i=\bar\sigma s_j$.
Because $i\leq j$, the second relation means that
\[ \bar\sigma(a)\ = \
\begin{cases}
\ \sigma(a) & \text{for $0\leq a\leq i$,}\\  
\sigma(a-1) & \text{for $i+1\leq a\leq j$, and}\\  
\ \sigma(a) & \text{for $j+1\leq a\leq n-1$.} 
\end{cases}
\]
We define $\tau:[n-2]\to[k]$ by
\[ \tau(a)\ = \
\begin{cases}
\ \bar\sigma(a) & \text{for $0\leq a\leq i$,}\\  
\bar\sigma(a+1) & \text{for $i+1\leq a\leq n-2$.}
\end{cases}
\]
Then
\[ \sigma \ = \ \tau s_{j-1}\text{\qquad and\qquad}\bar\sigma\ = \ \tau s_i \ . \]
Setting $w=\tau^*(z)$ then yields
\[ x \ = \ \sigma^*(z)\ = \ s_{j-1}^*(w)\text{\qquad and\qquad}
y\ = \ \bar\sigma^*(z)\ = \ s_i^*(w)\ .  \]
This means that $x$ and $y$ are equivalent, and hence $q(x)=q(y)$.
So the latching map $l_n:L_n X\to X_n$ is injective.

(ii) The latching map is continuous and injective by part~(i),
so it remains to show that $l_n$ is also a closed map.
The composite
\[ {\coprod}_{0\leq i\leq n-1}\ X_{n-1}^{[i]}\ \xra{\ q\ }\ L_n X \ \xra{\ l_n\ }\ X_n \]
of the quotient map and the latching map is the disjoint union
of the degeneracy maps $s_i^*:X_{n-1}\to X_n$.
So for every subset $A$ of $L_n X$, we have
\[ l_n(A) \ = \ \bigcup_{0\leq i\leq n-1}\, s_i^*\left(q^{-1}(A)\cap X_{n-1}^{[i]}\right)\ .\]
Now we suppose that $A$ is closed in $L_n X$.
Since each $s_i^*$ is a closed embedding,
$l_n(A)$ is a finite union of closed subsets of $X_n$, and hence itself closed.

(iii)
Since the morphism $s_i:[n]\to[n-1]$ has a retraction, 
the degeneracy map $s_i^*:X_{n-1}\to X_n$ has a continuous retraction.
Since $X_{n-1}$ and $X_n$ are compactly generated, 
$s_i^*:X_{n-1}\to X_n$ is a closed embedding 
by Proposition \ref{prop:section is closed embedding}.
So part~(ii) applies and shows that $l_n:L_n X\to X_n$ is a closed embedding.
Since $X_n$ is compactly generated, so is  $L_n X$,
by Proposition \ref{prop:properties cgwh spaces}~(i).
\end{proof}

\Danger
The proof of part~(iii) of the previous proposition
makes critical use of the fact that all degeneracy maps
in a simplicial compactly generated space are closed embeddings.
This is not the case more generally for simplicial $k$-spaces 
or simplicial topological spaces; so while a latching map
is always a continuous injection, it need not be a homeomorphism onto its image.
So if latching spaces are ever considered in this generality,
one has to beware the possible difference between $L_n X$ defined as the colimit
or as the subspace $s_0^*(X_{n-1})\cup\dots\cup s_{n-1}^*(X_{n-1})$ of $X_n$.

\medskip

Now we recall Reedy cofibrancy, which will be our main condition to ensure
good homotopical behavior of geometric realization.

\begin{defn}
A simplicial topological space $X:\bDelta^{\op}\to\bSpc$ 
is {\em Reedy cofibrant}\index{subject}{Reedy cofibrant!simplicial space}
is the latching morphism $l_n:L_n X\to X_n$ 
is a cofibration in the Quillen model structure on the category of topological spaces.
\end{defn}

The next proposition is the key reason why we care about Reedy cofibrancy.

\begin{prop}\label{prop:realization invariance in Top}
  Let $f:X\to Y$ be a morphism between Reedy cofibrant simplicial topological spaces.
  If $f_n:X_n\to Y_n$ is a weak equivalence for every $n\geq 0$, then
  the geometric realization $|f|:|X|\to |Y|$ is a weak equivalence.
\end{prop}
\begin{proof}
We use the Reedy model structure on
the category of simplicial topological spaces.
We endow the category $\bSpc$ of topological space
with the Quillen model structure, see \cite[II.3, Thm.\,1]{Q},
which makes it into a simplicial model category.
Geometric realization thus takes levelwise weak equivalences between
Reedy cofibrant objects to weak equivalences by \cite[VII Prop.\,3.6]{goerss-jardine}.
\end{proof}

\begin{rk}
By the precious proposition, the condition that a simplicial space is Reedy cofibrant
ensures that its geometric realization is `homotopically invariant'.
There are two other widely used conditions for ensuring that 
a geometric realization is well-behaved, the notions of a 
{\em good} and {\em proper} simplicial spaces.
A simplicial space $X$ is {\em proper}
in the sense of May \cite[\S 11]{may-geometry iterated loop} if
the latching map $\nu_n:L_n X\to X_n$ is a closed h-cofibration
for every $n\geq 0$. 
A simplicial space $X$ is {\em good} 
in the sense of Segal \cite[Def.\,A.4]{segal-cat coho}
if the degeneracy map $s_i^*:X_{n-1}\to X_n$ is a closed h-cofibration
for all $0\leq i\leq n-1$. 
One should beware that in the generality of arbitrary topological spaces,
the homotopy extension property does not imply that a map is a closed
embedding, hence `closed h-cofibration' is a stronger condition.
For compactly generated spaces, however, a continuous map with 
the homotopy extension property is automatically a closed embedding,
by Proposition \ref{prop:h-cof is closed embedding}.

For proper simplicial spaces, homotopical invariance of geometric realization 
is proved in \cite[Thm.\,A.4]{may-group completions}.
The case of good simplicial spaces then follows from the fact that
`goodness' implies `properness'; this is implicit in the proof 
of \cite[Lemma\,A.5]{segal-cat coho}, and stated explicitly in
the proof of Lewis' \cite[Cor.\,2.4 (b)]{lewis-cofibration}. 
The argument makes essential use of Lillig's `union theorem' \cite{lillig},
so it is {\em not} of a formal, model category theoretical nature.

Since cofibrations of compactly generated spaces are in particular 
closed h-cofibrations, 
every Reedy cofibrant simplicial compactly generated space 
is in particular proper and good.
Whenever we want to appeal to homotopy invariance of geometric realization,
our simplicial spaces are Reedy cofibrant; so we do not use the more
general conditions `good' and `proper', and we don't elaborate on these any further.
\end{rk}

We recall that the {\em singular complex} $\Sing(K)$ of a topological space $K$
is the simplicial set defined as the composite
\[ \bDelta^{\op}\ \xra{(\Delta^{\bullet})^{\op}}\ 
\bSpc^{\op} \ \xra{\bSpc(-,K)} \ \text{(sets)}\ .\]
In particular, the set of $n$-simplices $\Sing(K)_n$ is the set
of continuous maps from $\Delta^n$ to $K$.
Now we consider a simplicial topological space $X$.
Applying the singular complex functor levelwise yields a simplicial simplicial set
(i.e., a bisimplicial set) $\Sing\circ X$, 
sending $[n]$ to the simplicial set $\Sing(X_n)$.
We apply geometric realization to the simplicial sets $\Sing(X_n)$
for all $n\geq 0$, and end up with a new simplicial space  $|-|\circ\Sing\circ X$.

The singular complex functor is right adjoint the geometric realization
(restricted from simplicial spaces to simplicial sets),
and the adjunction counit is the continuous map
\[ \epsilon_K\ : \ |\Sing(K)| \ \to \ K \ , \quad [f,t]\ \longmapsto \ f(t)\ .\]
As was already noticed Milnor when he introduced the geometric
realization of simplicial sets \cite[Thm.\,4]{milnor-realization}, 
the map $\epsilon_X$ is a weak equivalence.
Since $\epsilon_K$ is natural, these maps assemble into a
morphism of simplicial spaces
\[  \epsilon_X\ : \  |-|\circ \Sing\circ X \ \to \ X \ . \]
We observe that $||-|\circ\Sing\circ X|=|\Sing\circ X|^{\text{it}}$
is the iterated realization \eqref{eq:def_iterated_real} 
of the bisimplicial set $\Sing\circ X$.

\begin{prop}\label{prop:real Sing X}
Let $X:\bDelta^{\op}\to\bSpc$ be a Reedy cofibrant simplicial topological space.
\begin{enumerate}[\em (i)]
\item The map $| \epsilon_X |:|\Sing\circ X|^{\text{\em it}} \to |X|$
is a weak equivalence.
\item If the spaces $X_0,\dots, X_n$ consist only of a single point each,
then the realization $|X|$ is $n$-connected.
\end{enumerate}
\end{prop}
\begin{proof}
(i) We claim that this simplicial space $|-|\circ \Sing\circ X$
sending $[m]$ to the realization of the simplicial set $\Sing(X_m)$
is automatically Reedy cofibrant.
Indeed, for every simplicial set $A$ the latching maps
$L_n A\to A_n$ are automatically injective
(a special case of Proposition \ref{prop:sspaces latching properties}~(i)). 
Hence for fixed $n\geq 0$ and varying $[m]\in\bDelta^{\op}$, the maps 
\[ L_n^{\Sing}(  \Sing(X_m) )\ \to \ \Sing(X_m)_n \]
(with latching object taken in the direction of the singular complex)
are a monomorphism of simplicial sets.
Geometric realization takes monomorphisms of simplicial sets to
relative CW-inclusions, and hence to cofibrations of spaces, 
and it commutes with colimits.
So the simplicial space $|-|\circ\Sing\circ X$ is Reedy cofibrant.
Since we assumed that $X$ is also Reedy cofibrant, 
the map $|\epsilon_X|$ is a weak equivalence 
by Proposition \ref{prop:realization invariance in Top}.

(ii) By part~(i) we may show that the space
$| \Sing\circ X|^{\text{it}}$ is $n$-connected.
The iterated realization is homeomorphic to the diagonal realization, 
by Proposition \ref{prop:iterated geometric realization}~(iii);
so we may show that the space $|\diag( \Sing\circ X)|$ is $n$-connected.
But the simplicial set $\diag( \Sing\circ X)$ has only a single $k$-simplex
for $0\leq k\leq n$, by hypothesis. So its geometric realization 
admits a CW-structure with one 0-cell and no cells in dimension 1 through $n$.
Hence the space $|\diag( \Sing\circ X)|$ is $n$-connected.
\end{proof}

In \cite[Thm.\,B.4]{bousfield-friedlander}, Bousfield and 
Friedlander establish a fairly general criterion to ensure that a dimensionwise
homotopy cartesian square remains homotopy cartesian after realization.
The criterion involves the `$\pi_*$-Kan condition' introduced 
in \cite[B.3]{bousfield-friedlander}. 
We now translate the condition into the context of simplicial spaces;
it roughly says that the simplicial sets obtained by taking $t$-th homotopy groups
dimensionwise satisfy the Kan extension condition.
The precise definition is more elaborate because one needs to properly 
take care of varying basepoints for homotopy groups.

\begin{defn}[Bousfield-Friedlander]\label{def:pi_*-Kan}
Let $X:\bDelta^{\op}\to\bSpc$ be a simplicial topological space, 
$m,t\geq 1$ and $a\in X_m$. 
Then $X$ satisfies the {\em $\pi_t$-Kan condition at $a$} 
if the following condition holds:
for all tuples of elements $x_i\in\pi_t(X_{m-1},d_i^*(a))$
for $i\in\{0,1,\dots,k-1,k+1,\dots,m\}$ satisfying 
\[ d_i^*(x_j)\ = \ d_{j-1}^*(x_i) \text{\quad in\quad $\pi_t(X_{m-2},(d_j d_i)^*(a))$} \]
for all $0\leq i<j\leq m$ with $i\ne k, j\ne k$, there is an element
$y\in \pi_t(X_m,a)$ such that
\[ d_i^*(y) \ = \ x_i \]
for all $0\leq i\leq m$ with $i\ne k$.

The simplicial space $X$ satisfies the 
{\em $\pi_*$-Kan condition}\index{subject}{pi star Kan condition@$\pi_*$-Kan condition}
if it satisfies the $\pi_t$-Kan condition for all $m,t\geq 1$ and all $a\in X_m$.
\end{defn}

\begin{eg}
As Bousfield and Friedlander already remark,
the $\pi_t$-Kan condition is automatically satisfied at $a\in X_m$
for all $t\geq 1$ whenever $a=\sigma_m^*(b)$ 
is the degeneracy of some point $b\in X_0$,
where $\sigma_m:[m]\to[0]$ is the unique morphism.
Indeed, in that special case, the collection of 
homotopy groups $\pi_t(X_n,\sigma_n^*(b))$ forms a simplicial group
as $n$ varies over the objects of $\bDelta$; 
the simplicial sets underlying simplicial groups always 
satisfy the Kan extension condition, an observation due to Moore,
see for example \cite[I Lemma\,3.4]{goerss-jardine}.
\end{eg}

Bousfield and Friedlander work with simplicial sets,
whereas we need their result for topological spaces.
The next theorem does the straightforward translation 
using the singular complex functor.
The reader may observe that \cite[Thm.\,B.4]{bousfield-friedlander}
has no Reedy cofibrancy condition; this should not be surprising 
because latching morphisms of simplicial objects in simplicial sets 
are automatically injective, and hence cofibrations. 
Put differently, bisimplicial spaces are automatically
Reedy cofibrant when considered as simplicial objects in simplicial sets.

A {\em homotopy fiber sequence}\index{subject}{homotopy fiber sequence}
is a pair of composable continuous maps
\[ A \ \xra{\ i \ }\ X \ \xra{\ f\ }\ Y  \]
whose composite $f i:A\to Y$ is constant with value $y_0$,
and such that the induced map
\begin{align*}
  A \ &\to \  F(f)\ = \ 
\{(\lambda,x)\in Y^{[0,1]}\times X\ | \ \lambda(0)=y_0,\, \lambda(1)=f(x)\} \\
a \ &\longmapsto \qquad (\text{const}_{y_0}, i(a))
\end{align*}
from $A$ to the homotopy fiber of $f$ is a weak equivalence.

\begin{theorem}\label{thm:realization preserves cartesian}
Let
\[ A \ \xra{\ i \ }\ X \ \xra{\ f\ }\ Y  \]
be morphisms of Reedy cofibrant simplicial topological spaces
such that the composite $f\circ i$ is constant.
Suppose that the following conditions hold:
\begin{enumerate}[\em (i)]
\item 
For every $n\geq 0$ the sequence $(i_n,f_n)$ is
a homotopy fiber sequence.
\item The simplicial spaces $X$ and $Y$ satisfy the $\pi_*$-Kan condition.
\item The morphism of simplicial sets $\pi_0(f):\pi_0(X)\to\pi_0(Y)$
is a Kan fibration.
\end{enumerate}
Then the sequence $(|i|,|f|)$ is a homotopy fiber sequence.
\end{theorem}
\begin{proof}
We transfer the question into the context of simplicial sets by use of the
singular complex functor, and then quote the theorem of Bousfield and Friedlander.
In the commutative diagram
\[ 
\xymatrix@C=13mm{ 
|\Sing\circ A|^{\text{it}} \ar[r]^-{|\Sing\circ i|^{\text{it}}} \ar[d]_{\epsilon_A} & 
|\Sing\circ X|^{\text{it}} \ar[d]^{\epsilon_X} \ar[r]^-{|\Sing\circ f|^{\text{it}}} &
|\Sing\circ Y|^{\text{it}} \ar[d]^{\epsilon_Y}\\
|A| \ar[r]_{|i|} & |X| \ar[r]_-{|f|} &  |Y| } \]
all vertical maps are weak equivalences by Proposition \ref{prop:real Sing X},
because we assumed that the simplicial spaces $A, X$ and $Y$ are
Reedy cofibrant.
The property of being a homotopy fiber sequence is invariant under 
pointwise weak equivalences, so we may show that the upper row
is a homotopy fiber sequence.
This sequence arises from a commutative square of bisimplicial sets
\[ \xymatrix@C=12mm{ \Sing\circ A \ar[r]^-{\Sing\circ i}\ar[d] & 
\Sing\circ X\ar[d]^{\Sing\circ f} \\
\ast \ar[r] & \Sing\circ Y} \]
by taking iterated realizations. 
In every fixed simplicial dimension, the square of simplicial sets
is a homotopy fiber square in the sense of \cite[Def.\,A.2]{bousfield-friedlander},
by hypothesis~(i).
Moreover, the bisimplicial sets $\Sing\circ X$ and $\Sing\circ Y$
satisfy the $\pi_*$-Kan condition by assumption (ii),
and the morphism of simplicial sets 
$\pi_0^v(\Sing\circ f):\pi_0^v(\Sing\circ X)\to\pi_0^v(\Sing\circ Y)$
is a Kan fibration by assumption~(iii).
So \cite[Thm.\,B.4]{bousfield-friedlander} shows that the
square of diagonal simplicial sets
\[ \xymatrix@C=15mm{ \diag(\Sing \circ A) \ar[r]^-{\diag(\Sing\circ i)}\ar[d] & 
\diag(\Sing\circ X)\ar[d]^{\diag(\Sing\circ f)} \\
\ast \ar[r] & \diag(\Sing\circ Y)} \]
is homotopy cartesian.
After geometric realization we can identify the diagonal realization 
with the iterated realization 
as in Proposition \ref{prop:iterated geometric realization}~(iii),
so this proves the claim.
\end{proof}

At some later stage we will need to know that the product of two
Reedy cofibrant simplicial spaces is again Reedy cofibrant.
The proof of this fact is rather formal and works in a much broader context,
as the following proposition shows.

\begin{prop}\label{prop:box preserves Reedy cofibrant}
Let $\Cc$ be a model category that is equipped with a monoidal product $\boxtimes$
satisfying the pushout product property for cofibrations.
Then for all Reedy cofibrant simplicial objects $X,Y:\bDelta^{\op}\to\Cc$,
the objectwise monoidal product $X\boxtimes Y:\bDelta^{\op}\to\Cc$
is Reedy cofibrant.
\end{prop}
\begin{proof}
I learned this proof from Cary Malkiewich.
We let $\Pc(m)$ denote the power set of the set $\{1,\dots,m\}$, 
i.e., the set of subsets.
We also write $\Pc(m)$ for the associated poset category, under inclusion.
An {\em $m$-cube} in $\Cc$ is simply a functor $F:\Pc(m)\to\Cc$.
We call such an $m$-cube {\em cofibrant} if for every subset $B$ of $\{1,\dots,m\}$
the canonical morphism
\[ l_B \ : \ \colim_{A\subsetneq B} \, F(A)\ \to \  F(B)\]
is a cofibration.

We claim that for every cofibrant $m$-cube $F$ and all
sets $\Yc\subset\Zc\subset\Pc(m)$ that are both closed under passage to subsets,
the canonical morphism
\[ \colim_{A\in \Yc} \, F(A)\ \to \ \colim_{A\in \Zc}\, F(A)\]
is a cofibration. 
We start with the special case where $\Zc=\Yc\cup\{B\}$ 
for some subset $B$ of $\{1,\dots,m\}$ that does not belong to $\Yc$.
Since $\Zc$ is closed under taking subsets, 
every proper subset of $B$ belongs to $\Zc$, and hence to $\Yc$.
Then the square
\[ \xymatrix@C=10mm{ 
 \colim_{A\subsetneq B} \, F(A) \ar[d]\ar[r]^-{l_B} &  F(B)\ar[d]\\
\colim_{A\in \Yc} \, F(A)\ar[r] & \colim_{A\in \Zc}\, F(A)
} \]
is a  pushout in $\Cc$.
The upper horizontal latching morphism is a cofibration by hypothesis;
so the lower horizontal morphism is a cofibration.  

In the general case we choose a chain of intermediate subsets
\[ \Yc \ = \ \Yc_0 \ \subset\ \Yc_1 \ \subset\ \dots\ \subset \Yc_k \ = \ \Zc \]
such that each $\Yc_i$ is closed under taking subsets 
and $\Yc_i$ has exactly one element more than $\Yc_{i-1}$. 
The claim then holds for each pair $(\Yc_i,\Yc_{i-1})$.
Since the composite of two cofibrations 
is a cofibration, this proves the general case.

Now we turn to the proof of the proposition.
In Remark \ref{rk:Delta(n) vs P(n)} we specified an isomorphism of categories
$\kappa:\bDelta(m)^{\op}\iso \Pc(m)$.
Since $X$ is Reedy cofibrant, the $m$-cube
\[ \Pc(m)\ \xra[\iso]{\ \kappa^{-1}} \ \bDelta(m)^{\op}\ \xra{\ u \ } 
\ \bDelta^{\op}\ \ \xra{\ X \ }\ \Cc  \]
is cofibrant for every $m\geq 0$,
and the same is true for $Y$. We observe that the $(m+n)$-cube
\[ X\tensor Y\ : \ \Pc(m+n)\ \to \ \Cc \ , 
\quad U+V\ \longmapsto \ X_{|U|}\boxtimes Y_{|V|} \]
is then again cofibrant, where `+' denotes disjoint union of sets.
Indeed, for $B\in \Pc(m)$ and $B'\in\Pc(n)$,
the latching object
\[ \colim_{A+A'\subsetneq B+B'} X_{|A|}\boxtimes Y_{|A'|} \]
is isomorphic to the pushout
\[ (L_B X)\boxtimes Y_{|B'|} \cup_{(L_B X)\boxtimes (L_{B'}Y)} X_{|B|}\boxtimes (L_{B'}Y) \]
because $\boxtimes$ preserves colimits in both variables.
Under this isomorphism, the latching morphism
\[ l_{B+B'}\ : \ \colim_{A+A'\subsetneq B+B'} X_{|A|}\boxtimes Y_{|A'|} 
\ \to \ X_{|B|}\boxtimes Y_{|B'|} \]
becomes the pushout product $l_B^X\Box l_{B'}^Y$
of the latching morphisms for $X$ respectively $Y$.
This latter morphism is a cofibration by the pushout product property.

Now we take $m=n$ and let $\Xc$ be the subposet of $\Pc(n+n)$ 
of proper diagonal elements, i.e., the sets $U+U$ 
for a proper subset $U$ of $\{1,\dots,n\}$.
The latching object $L_n(X\boxtimes Y)$ is then a colimit
of the functor $X\tensor Y$ over the poset $\Xc$.
We let $\Yc$ be the subposet of $\Pc(n+n)$
consisting of those sets $U+V$ such that $U\cup V\ne \{1,\dots,n\}$.
The latching map for $X\boxtimes Y$ factors as the composite
\begin{align}  \label{eq:latching_as_composite}
 L_n(X\boxtimes Y)\ = \ &\colim_{U+U\in\Xc} \, X_{|U|}\boxtimes Y_{|U|}\nonumber\\ 
\to\  &\colim_{U+V\in\Yc} \, X_{|U|}\boxtimes Y_{|V|}\ \to \ X_n\boxtimes Y_n \ .  
\end{align}
We observe that the inclusion $\Xc\to\Yc$ is final, i.e., 
for every $U\in\Yc$ the comma category $U\downarrow \Xc$
is non-empty and connected. 
So the first morphism in \eqref{eq:latching_as_composite} is an isomorphism.
We set $\Zc=\Pc(n+n)$;  since $\Zc$ has a terminal object,
the object $X_n\boxtimes Y_n$ is a colimit 
of the functor $X\tensor Y$ over $\Zc$.
Both $\Yc$ and $\Zc$ are closed under passage to subsets, 
so the first paragraph shows that 
the second morphism in \eqref{eq:latching_as_composite} is a cofibration.
This proves the claim.  
\end{proof}

\chapter{Equivariant spaces}
\label{app:equivariant spaces}

In this appendix we collect basic results about the equivariant homotopy
theory of $G$-spaces. Initially $G$ can be any topological group, but we will
eventually specialize to compact Lie groups.
We start out by checking that taking fixed points commutes with certain kinds of
colimits, namely pushouts and sequential colimits along closed embeddings,
smash products, geometric realization and latching objects,
see Proposition \ref{prop:G-fix preserves pushouts}.
Proposition \ref{prop:proj model structures for G-spaces}
provides a self-contained proof of the standard model structure 
on the category of $G$-spaces, relative to a set of closed subgroups.
In the realm of compact Lie groups, 
various change of group functors preserve cofibrations,
namely fixed points (see Proposition \ref{prop:cofamily pushout property}),
restriction along a continuous homomorphism, induction,
and orbits (see Proposition \ref{prop:cofibrancy preservers}). 
Proposition \ref{prop:fix of free cofibration}
is a useful decomposition result for $(G\bs X)^K$, 
the $K$-fixed points of the $G$-orbits of a $G$-free 
$(K\times G)$-space.

Then we turn to equivariant $\bGamma$-spaces.
We start with the observation that compactly generated spaces are 
`closed under prolongation of $\bGamma$-spaces'.\index{subject}{Gamma-space@$\bGamma$-space}
More precisely, if $F:\bGamma\to\bK$ is a $\bGamma$-$k$-space
such that $F(n_+)$ is compactly generated for every $n\geq 0$,
then it is not completely obvious whether the prolongation $F(K)$
to a general compactly generated space $K$, defined as a certain quotient space of
$\coprod_{n\geq 0} F(n_+)\times K^n$, is weak Hausdorff.
We show in Proposition \ref{prop:prolongation in bT}
that this is automatically the case.
Because no `weak Hausdorffication' is necessary, 
there are no unexpected identifications, and the underlying set of $F(K)$ 
is what one first thinks of.

For us, the main purpose of $\bGamma$-spaces is to provide spectra
by evaluation on spheres. We want to do this equivariantly, 
in the presence of an action of a compact Lie group.
Proposition \ref{prop:Gamma fixed points} gives a way to calculate 
the fixed points of a prolonged $\bGamma$-$G$-space
under the action of a {\em connected} compact Lie group.
Since $G$-fixed points can be obtained by first taking fixed points
of the identity component $G^\circ$, and then fixed points 
of the finite group $\pi_0 G=G/G^\circ$, this formula effectively reduces
questions about fixed points to the case of finite groups.

In order to analyze prolonged $\bGamma$-$G$-spaces homotopically,
we need a cofibrancy condition that we introduce 
in Definition \ref{def:G-cofibrant Gamma-G-space}.
This kind of cofibrancy is stable under passage to closed subgroups
and taking fixed points (see Proposition \ref{prop:fix preserves Gamma-cofibrant})
and guarantees that evaluation on spheres produces
a $G$-spectrum that is equivariantly connective (see Proposition \ref{prop:Gamma on S^V}).
Moreover, for $G$-co\-fibrant $\bGamma$-$G$-spaces,
prolongation is homotopical, i.e., strict equivalences prolong
to $G$-weak equivalences on finite $G$-CW-complexes
(see Proposition \ref{prop:cofibrant Gamma-G on G-CW}).

We conclude this appendix with a reformulation and generalization 
of the Segal-Shima\-kawa delooping formalism
for equivariant $\bGamma$-spaces, see Theorems \ref{thm:prolonged delooping} 
and \ref{thm:special cofibrant Gamma positive Omega}. 
When restricted to finite groups, these theorems show that
evaluation of a $G$-cofibrant $\bGamma$-$G$-space $F$  on spheres provides
a positive $G$-$\Omega$-spectrum if $F$ is `special' 
(compare Definition \ref{def:special G-Gamma}),
and a full fledged $G$-$\Omega$-spectrum if $F$ is `very special'
(compare Definition \ref{def:very special G-Gamma}).
As we explain in Remark \ref{rk:Gamma limitations}, there is no hope
to obtain a $G$-$\Omega$-spectrum
for compact Lie groups of positive dimension.
However, we do have partial delooping results for non-finite compact Lie groups:
Theorem \ref{thm:special cofibrant Gamma positive Omega} 
effectively says that evaluating a $G$-cofibrant special $\bGamma$-$G$-space
on spheres yields a `$G^\circ$-trivial positive $G$-$\Omega$-spectrum',
where $G^\circ$ is the identity component of $G$.
In other words, evaluating on spheres provides equivariant deloopings with respect to 
all those $G$-representations on which $G^\circ$ acts trivially.
If $G$ is not finite, this is of course a very restricted class of
representations; but for trivial representations we can at least conclude
that evaluation on spheres gives a `naive' $G$-$\Omega$-spectrum.

\medskip

We let $G$ be a topological group, which we take to mean a group
object in the category $\bT$ of compactly generated spaces.
So a topological group is a compactly generated
space equipped with an associative and unital multiplication
\[ \mu \ : \  G \times G \ \to \ G \]
that is continuous with respect to the 
compactly generated product topology, and such that the shearing map
\[   G \times G \ \to \  G \times G \ ,\quad (g,h)\ \longmapsto \ (g, g h)\]
is a homeomorphism (again for the compactly generated product topology). 
This implies in particular that inverses exist in $G$,
and that the inverse map $g\mapsto g^{-1}$ is continuous.
A $G$-space is then a compactly generated space $X$ 
equipped with an associative and unital action 
\[ \alpha\ : \ G  \times X \ \to \ X \]
that is continuous with respect to the compactly generated product topology.
We write $G\bT$\index{symbol}{$G\bT$ - {category of $G$-spaces}}
for the category of $G$-spaces and continuous $G$-maps.

The forgetful functor from $G$-spaces to compactly generated spaces 
has both a left and a right adjoint,
and hence limits and colimits of $G$-spaces are created in the underlying
category $\bT$. This has nothing to do with topology, and is entirely formal,
using only that the underlying category $\bT$ is cartesian closed, 
complete and cocomplete.

\medskip

\Danger Since colimits of $G$-spaces are created in the underlying category $\bT$
of compactly generated spaces, the earlier caveat applies as well.
A colimit in $G$-spaces is calculated by first forming a colimit 
in the category $\bSpc$ of all topological spaces (or equivalently,
in the full subcategory $\bK$ of $k$-spaces); the result is a $k$-space,
but it need not be weak Hausdorff. In that case an application of
the functor $w:\bK\to\bT$ left adjoint to the inclusion produces 
a colimit in $\bT$. This colimit comes with a preferred $G$-action making
it a colimit in the category of $G$-spaces. 
One has to beware that whenever the colimit in $\bSpc$ is not weak Hausdorff,
the functor $w$ changes the underlying set; in particular, the forgetful functor to sets
does {\em not} preserve such colimits.

\medskip

Now we consider a closed subgroup $H$ of a topological group $G$.
Then $H$ is compactly generated in the subspace topology by
Proposition \ref{prop:properties cgwh spaces}~(i), and
hence a topological group (internal to the category $\bT$) in its own right.
For a $G$-space $X$ we denote by
\[   X^H \ = \ \{ x\in X\ | \ h x = x \text{ for all $h\in H$}\}  \]
the subspace of $H$-fixed points. For an individual element $h\in H$
the $h$-fixed subspace $\{x\in X\ | \ h x = x \}$
is the preimage of the diagonal under the continuous map
$(\Id,h\cdot-):X\to X\times X$, so it is a closed subspace of $X$
by Proposition \ref{prop:wH criterion}.
As an intersection of closed subsets, $X^H$ is then closed in $X$, and
hence compactly generated in the subspace topology,
by Proposition \ref{prop:properties cgwh spaces}~(i).

The following proposition records that fixed points preserve certain kinds of colimits.
In part~(iv) we consider a simplicial $G$-space $X:\bDelta^{\op}\to G\bT$,
and we write $X^G$ for the simplicial space consisting of the fixed points of $X$,
i.e., the composite functor
\[ \bDelta^{\op}\ \xra{\ X \ }\ G\bT \ \xra{(-)^G} \ \bT\ . \]

\begin{prop}\label{prop:G-fix preserves pushouts}
Let $G$ be a topological group.
\begin{enumerate}[\em (i)]
\item 
For every pushout square of $G$-spaces on the left
\[ \xymatrix{ A \ar[d]_-i \ar[r]^-f & C \ar[d]^-j && 
A^G \ar[d]_-{i^G} \ar[r]^{f^G} & C^G \ar[d]^-{j^G}\\
B\ar[r]_-h & D && B^G\ar[r]_-{h^G} & D^G } \]
in which the map $i$ is a closed embedding,
the square of fixed point spaces on the right is a pushout.
\item Taking $G$-fixed points commutes with filtered colimits
along continuous $G$-maps that are closed embeddings.
\item  For all based $G$-spaces $X$ and $Y$ 
 the canonical map $X^G\sm Y^G\to (X\sm Y)^G$  is a homeomorphism.  
\item
For every simplicial $G$-space $X:\bDelta^{\op}\to G\bT$ and every $n\geq 0$,
the canonical maps 
\[ |X^G|\to |X|^G\text{\qquad and\qquad} L_n(X^G)\to (L_n X)^G \]
are homeomorphisms.
\end{enumerate}
\end{prop}
\begin{proof}
(i) Pushouts in $G$-spaces are formed on underlying compactly generated spaces.
Since $i$ is a closed embedding, so is its restriction to
fixed points $i^G$. So by Proposition \ref{prop:pushout in cg} 
both pushouts are `as expected', i.e., formed in the ambient category
of all topological spaces. In particular, $D$ is the set-theoretic
disjoint union of the images of $B-i(A)$ and $C$, which are both $G$-invariant.
So $D^G$ is the set-theoretic disjoint union of the images of
$(B-i(A))^G=B^G-(i^G)(A^G)$ and $C^G$. The canonical map
$B^G\cup_{A^G}C^G\to D^G$ is thus a continuous bijection.

Now we show that the canonical map is closed.
We consider the commutative square:
\[ \xymatrix@C=12mm{ 
B^G\amalg C^G \ar[r]^-{\text{incl}}\ar[d]_p & B\amalg C\ar[d]^{h+j} \\
B^G\cup_{A^G}C^G\ar[r]_-{h^G\cup j^G} & D} \]
If $O$ is a closed subset of $B^G\cup_{A^G}C^G$,
then $p^{-1}(O)$ is closed in $B^G\amalg C^G$, and hence in $B\amalg C$
(since fixed points are closed in the ambient space).
The relation
\[ p^{-1}(O)\ = \ (h+j)^{-1}( (h^G\cup j^G)(O)) \]
holds in $B\amalg C$; since the right vertical map $h+j$ 
is a proclusion, $(h^G\cup j^G)(O)$ is closed in $D$, and hence also in $D^G$.

(ii) We let $P$ be a filtered poset and $F:P\to G\bT$ a functor
to the category of $G$-spaces. 
We must show that the canonical continuous map
\[ \kappa \ : \ \colim_P \, F^G \ \to \ (\colim_P F)^G \]
is a homeomorphism, where $F^G=(-)^G\circ F$.
Proposition \ref{prop:filtered colim in cg} ensures that both colimits
are formed in the ambient category of sets, so the map is bijective.
To see that the canonical map is also closed, we consider the commutative square
\[ \xymatrix{
\coprod_{j\in P} F(j)^G \ar[r]^-{\text{incl}}\ar[d]_p & \coprod_{j\in P} F(j) \ar[d]^q\\
\colim_P \, F^G \ar[r]_-\kappa & \colim_P F } \]
in which the vertical maps are proclusions.
The canonical maps $\kappa_j:F(j)\to \colim_P F$ are injective,
so every point of $F(j)$ that becomes $G$-fixed in the colimit
is already $G$-fixed in $F(j)$, i.e., the square is a pullback.
If $A$ is a closed subset of $\colim_P F^G$, then $p^{-1}(A)$
is closed in the coproduct of the $F(j)^G$'s. Since $F(j)^G$
is closed in $F(j)$, the set $q^{-1}(\kappa(A)) = p^{-1}(A)$ 
is closed in the coproduct of the $F(j)$'s. Since $q$ is a proclusion,
the set $\kappa(A)$ is closed in $\colim_P F$, and hence also in $(\colim_P F)^G$.
This shows that $\kappa$ is a closed map.

(iii)
Points of compactly generated spaces are closed 
(Proposition \ref{prop:properties wH spaces}~(iii)),
so the subspace $X\times \{y_0\}\cup \{x_0\}\times Y$ is closed in $X\times Y$.
The claim then follows by applying part~(i) to the pushout of the diagram
\[ \ast \ \longleftarrow\ X\times \{y_0\}\cup \{x_0\}\times Y \ \to \  X\times Y \ . \]

(iv) Since $X_n^G$ is closed inside $X_n$, the inclusion $X^G\subset X$ 
induces a closed embedding
$|X^G|\to |X|$ by Proposition \ref{prop:geometric realization in bT}~(iii).
The image of this map is clearly $G$-fixed. So it remains to show 
that every $G$-fixed point of $|X|$ is in the image of $|X^G|$.
We let $(x,t)\in X_l\times \Delta^l$ be the  minimal representative
of a given $G$-fixed point in $|X|^G$. Then for every $g\in G$ 
the point $(g x,t)$ is equivalent to $(x,t)$.
Since the minimal representative is unique, this forces $g x=x$.
Hence $x\in X_l^G$, which proves the claim.

The latching map $l_n:L_n X\to X_n$
is a closed embedding by Proposition \ref{prop:sspaces latching properties}~(iii).
Hence the restriction to fixed points $(l_n)^G:(L_n X)^G\to (X_n)^G$ 
is also a closed embedding.
The composite
\[ L_n(X^G)\ \to\  (L_n X)^G\ \xra{(l_n)^G}\ (X_n)^G \]
is the latching map for the simplicial compactly generated space $X^G$,
hence also a closed embedding, 
again by Proposition \ref{prop:sspaces latching properties}~(iii).
So the map in question is a closed embedding.

It remains to show that every $G$-fixed point of $L_n X$ 
is the image of a point in $L_n(X^G)$.
By Proposition \ref{prop:sspaces latching properties},
the latching map identifies $L_n X$ with the union
of the subspaces $s_i^*(X_{n-1})$  for $i=0,\dots,n-1$.
So we may show that the composite
\[ L_n(X^G)\ \to\  (L_n X)^G\ \xra[\iso]{\ (l_n)^G\ }\ 
\left({\bigcup}_{i=0,\dots,n-1} s_i^*(X_{n-1}) \right)^G \]
is surjective.
If  $y\in X_{n-1}$ is such that $x=s_i^*(y)$ is $G$-fixed, then $y$ is
$G$-fixed because every degeneracy map is injective and $G$-equivariant.
So $y\in (X_{n-1})^G$ represents a point in $L_n(X^G)$ that maps to $x$.
\end{proof}

\begin{prop}\label{prop:G/H is finite} 
  Let $H$ be a closed subgroup of a topological group $G$.
  \begin{enumerate}[\em (i)]
  \item The orbit space $G/H$ is compactly generated in the quotient topology,
    and hence a $G$-space.  
  \item Let
    \[ X_0 \ \to \ X_1 \ \to \ \dots \ \to \ X_n \ \to \ \dots \]
    be a sequence of closed embeddings of $G$-spaces
    and $X_\infty$ a colimit of the sequence in the category of $G$-spaces.
    Then for every compact space $K$,
    every continuous $G$-map from $G/H\times K$ to $X_\infty$ 
    factors through $X_n$ for some $n$, and that factorization is
    a continuous $G$-map.
  \end{enumerate}
\end{prop}
\begin{proof}
(i) The equivalence relation $E\subset G\times G$
that gives rise to $G/H$ is the inverse image of $H$ under the continuous map
\[ G\times G \ \to \ G \ , \quad (g,\bar g)\ \longmapsto \ g^{-1}\cdot\bar g \ . \]
Because $H$ is closed, so is the equivalence relation, and hence
the quotient topology is compactly generated by 
the criterion given by Proposition \ref{prop:closed implies WH}.

(ii)
We let $f : G/H \times K \to X_\infty$ be a continuous $G$-map. 
Colimits of $G$-spaces are created in
the underlying category $\bT$ of compactly generated spaces;
since the maps in question are closed embeddings, 
the sequential colimit is created in the ambient category $\bSpc$
of all topological spaces, by Proposition \ref{prop:filtered colim in cg}.  
The composite 
\begin{equation}  \label{eq:K2X_infty}
 K \ \xra{(e H,-)} \ G/H \times K \ \xra{\ f\ } \ X_\infty   
\end{equation}
factors through a continuous map $\bar f:K\to X_n$ for some $n\geq 0$,
by Proposition \ref{prop:filtered colim preserve weq}~(i).
The canonical map $X_n\to  X_\infty$ is injective
by Proposition \ref{prop:filtered colim in cg}.
Since the image of the composite \eqref{eq:K2X_infty}
is contained in the $H$-fixed points of $X_\infty$,
the image of the factorization $\bar f:K\to X_n$ is contained
in the $H$-fixed points of $X_n$.
So the composite
\begin{equation}  \label{eq:f^flat}
 G\times K \ \xra{G\times \bar f}\ G\times X_n \ \xra{\text{act}} \ X_n
\end{equation}
factors through a well-defined set theoretic $G$-equivariant map
\[  f' \ : \  G/H\times K \ \to \ X_n \ .\]
The map $\text{proj}\times K:G\times K\to G/H\times K$ is a
proclusion by Proposition \ref{prop:proclusion times Z},
so the continuity of \eqref{eq:f^flat} implies that $f'$ is also continuous,
and hence the desired factorization of the original morphism $f$.
\end{proof}

An argument that we need several times in this book proves
that a certain model structure is topological.
To avoid repeating the same kind of argument, we axiomatize it.
The `classical' model structure on the category of all topological spaces
was established by Quillen in \cite[II.3 Thm.\,1]{Q}.
We use the straightforward adaptation of this model structure to
the category of compactly generated spaces,
which is described for example in \cite[Thm.~2.4.25]{hovey-book}.
In this model structure on the category $\bT$, the weak equivalences 
are the weak homotopy equivalences and fibrations
are the Serre fibrations. The cofibrations are the retracts of 
generalized CW-complexes, i.e., relative cell complexes in which cells
can be attached in any order and not necessarily to cells of lower dimensions.

We consider a model category $\M$ that is also enriched, tensored and cotensored
over the category $\bT$ of compactly generated spaces.
We denote the tensor by $\times$.
Given a morphism $f:X\to Y$ in $\M$  and a continuous map of spaces $g:A\to B$,
we denote by $f\Box g$ the 
{\em pushout product} morphism\index{subject}{pushout product} defined as
\[ f\Box g = (f\times B)\cup(Y\times g) \ : \ 
X\times B\cup_{X\times A}Y\times A \ \to \ Y\times B\ .\]
We recall that the model structure is called {\em topological}\index{subject}{topological model structure}\index{subject}{model structure!topological|see{topological model structure}}
if the following two conditions hold:
\begin{itemize}
\item if $f$ is a cofibration in $\M$ and $g$ is a cofibration of spaces, 
  then the pushout product morphism $f\Box g$ is a cofibration;
\item if in addition $f$ or $g$ is a weak equivalence, then so is the
  pushout product morphism  $f\Box g$.
\end{itemize}
The pushout product condition can also be stated in two different, but equivalent,
adjoint forms, compare \cite[Lemma 4.2.2]{hovey-book}.
In the next proposition, we denote by
\[ i_k \ : \ \partial D^k \ \to \ D^k \text{\qquad and\qquad}
 j_k \ : \  D^k\times\{0\} \ \to \ D^k\times[0,1]  \]
the inclusions. Then $\{i_k\}_{k\geq 0}$ is the standard set
of generating cofibrations for the Quillen model structure on the category
of spaces, and  $\{j_k\}_{k\geq 0}$ is the standard set
of generating acyclic cofibrations, compare \cite[Thm.\,2.4.25]{hovey-book}.

\begin{prop}\label{prop:topological criterion}
Let $\M$ be a model category that is also enriched, tensored and cotensored
over the category $\bT$ of spaces.
Suppose that there is a set of cofibrant objects $\Gc$ 
and a set of acyclic cofibrations $\Zc$ of $\M$ with the following properties:
\begin{enumerate}[\em (a)]
\item The acyclic fibrations are characterized by the right lifting property
with respect to the morphisms $K\times i_k$ for all $K\in\Gc$ and $k\geq 0$.
\item The fibrations are characterized by the right lifting property
with respect to the union of the morphisms $K\times j_k$ for all $K\in\Gc$ and $k\geq 0$
and the pushout products $c\Box i_k$ for all $c\in\Zc$ and $k\geq 0$.
\end{enumerate}
Then the model structure is topological.
\end{prop}
\begin{proof}
Since the tensor bifunctor $\times$ has an adjoint in each variable, 
it preserves colimits in each variable.
So it suffices to check the pushout product properties
when the maps $f$ and $g$ are from the sets of generating (acyclic) cofibrations,
compare \cite[Cor.\,4.2.5]{hovey-book}.
The set of inclusions of spheres into discs is closed under pushout product,
in the sense that $i_k\Box i_l$ is homeomorphic to $i_{k+l}$.
So pushout product with $i_l$ preserves the set of
generating cofibrations $\{K\times i_k\}_{K\in \Gc,k\geq 0}$ (up to isomorphism). 
This takes care of the part
of the pushout product property that involves the cofibrations only.

Similarly, the pushout product of $j_k$ with $i_l$ is isomorphic to $j_{k+l}$.
So pushout product with $i_l$ preserves the set  of
generating acyclic cofibrations $\{K\times j_k\}_{K\in \Gc,k\geq 0}$, and
it preserves the set of generating acyclic cofibrations  $\{c\times i_k\}_{c\in \Zc,k\geq 0}$.
This shows that the pushout product of an acyclic cofibration in $\M$
with a cofibration of spaces is again an acyclic cofibration.

Finally, pushout product with $j_l$ takes the generating cofibrations 
$K\times i_k$ to morphisms of the form
$K\times j_m$ with  $K\in \Gc$ and $m\geq 0$, which are acyclic cofibrations.
This shows that the pushout product of a cofibration in $\M$
with an acyclic cofibration of spaces is again an acyclic cofibration.
\end{proof}

Now we let $\Cc$ be a set of closed subgroups of a topological group $G$.
We call a morphism $f:X\to Y$ of $G$-spaces a {\em $\Cc$-equivalence}
(respectively {\em $\Cc$-fibration})\index{subject}{C-fibration@$\Cc$-fibration} 
if the restriction $f^H:X^H\to Y^H$ to $H$-fixed points is a
weak equivalence (respectively Serre fibration)
of spaces for all subgroups $H$ of $G$ in $\Cc$.
A {\em $\Cc$-cofibration}\index{subject}{C-cofibration@$\Cc$-cofibration} 
is a morphism with the left lifting property 
with respect to all morphisms that are simultaneously
$\Cc$-equivalences and $\Cc$-fibrations.

\begin{prop}[Gluing lemma]\label{prop:gluing lemma G-spaces}
  Let $G$ be a topological group and $\Cc$ a set 
  of closed subgroups of $G$. Consider  a commutative diagram of $G$-spaces 
  \[ \xymatrix{
    C  \ar[d]_\gamma & A \ar[l]_-g \ar[d]^\alpha \ar[r]^-f & B \ar[d]^\beta\\
   \bar C  & \bar A \ar[l]^-{\bar g}\ar[r]_-{\bar f} & \bar B } \]
  such that $f$ and $\bar f$ are h-cofibrations of $G$-spaces.
  Suppose that the maps $\alpha,\beta$ and $\gamma$ are $\Cc$-equivalences.
  Then the induced map of pushouts
    \[ \gamma\cup \beta\ : \ C\cup_A B \ \to \ \bar C\cup_{\bar A} \bar B \]
    is a $\Cc$-equivalence.
\end{prop}
\begin{proof}
  We let $H$ be a closed subgroup from the set $\Cc$,
  and we contemplate the commutative diagram of fixed points:
  \[ \xymatrix@C=10mm{
    C^H  \ar[d]_{\gamma^H} & A^H \ar[l]_-{g^H} \ar[d]^{\alpha^H} \ar[r]^-{f^H} & 
    B^H \ar[d]^{\beta^H}\\
    \bar C^H  & \bar A^H \ar[l]^-{\bar g^H}\ar[r]_-{\bar f^H} & \bar B^H } \]
  Since $f$ and $\bar f$ are h-cofibrations of $G$-spaces,
  $f^H$ and $\bar f^H$ are h-cofibrations of spaces. 
  The three vertical maps are weak equivalences by hypothesis.
  The gluing lemma for weak equivalences and pushouts along 
  h-cofibrations then shows that the induced map on horizontal pushouts 
    \[ \gamma^H\cup \beta^H\ : \ C^H \cup_{A^H} B^H \ \to \ \bar C^H\cup_{\bar A^H} \bar B^H \]
  is a weak equivalence,
  see for example \cite[Appendix, Prop.\,4.8 (b)]{boardman-vogt-homotopy invariant}.
  Since $f$ and $\bar f$ are h-cofibrations of $G$-spaces, they are in particular 
  h-cofibrations of underlying spaces, and hence closed embeddings
  (Proposition \ref{prop:h-cof is closed embedding}~(ii)).
  So taking $H$-fixed points commutes with the horizontal pushout
  (by Proposition \ref{prop:G-fix preserves pushouts}~(i)), and we conclude that
  also the map
    \[ (\gamma\cup \beta)^H\ : \ (C\cup_A B)^H \ \to \ (\bar C\cup_{\bar A} \bar B)^H \]
    is a weak equivalence. This proves the claim.
\end{proof}

The following $\Cc$-projective model structure
is well known and fairly standard, and mentioned in various places in the literature,
for example in \cite[Prop.\,2.11]{fausk-pro} and \cite[III Thm.\,1.8]{mandell-may}.
However, I do not know a reference that is both self-contained and complete, 
so I provide the proof.

\begin{prop}\label{prop:proj model structures for G-spaces}
  Let $G$ be a topological group and $\Cc$ a set 
  of closed subgroups of $G$. 
  Then the $\Cc$-equivalences, $\Cc$-cofibrations and
  $\Cc$-fibrations form a model structure,
  the {\em $\Cc$-projective model structure}
  on the category of $G$-spaces.
  This model structure is proper, cofibrantly generated and topological.
\end{prop}
\begin{proof}
We number the model category axioms as in \cite[3.3]{dwyer-spalinski}.
The category of $G$-spaces is complete and cocomplete and all
limits and colimits are created in the underlying category of compactly generated spaces.
Model category axioms MC2 (2-out-of-3)
and MC3 (closure under retracts) are clear. One half of MC4 
(lifting properties) holds by the definition of $\Cc$-cofibrations.
The proof of the remaining axioms uses Quillen's\index{subject}{small object argument}
small object argument, originally given in \cite[II p.\,3.4]{Q}, 
and later axiomatized in various places,
for example in \cite[7.12]{dwyer-spalinski} or \cite[Thm.\,2.1.14]{hovey-book}.
In the category of (non-equivariant) spaces, 
the set  $\{i_k: \partial D^k \to D^k\}_{k\geq 0}$ of inclusions of spheres into discs 
detects Serre fibrations that are simultaneously weak equivalences. 
By adjointness, the set
\begin{equation}\label{eq:I_for_F-proj_on_GT}
I_\Cc\ = \  \{ G/H\times i_k \ : \ G/H \times \partial D^k\ \to \ G/H \times  D^k \}_{k\geq 0, H\in \Cc}
\end{equation}
then detects acyclic fibrations in the $\Cc$-projective 
model structure on $G$-spaces.
Similarly, the set of inclusions $\{j_k: D^k\times \{0\} \to D^k\times [0,1]\}_{k\geq 0}$ 
detects Serre fibrations; so by adjointness, the set
\begin{equation}\label{eq:J_for_F-proj_on_GT}
J_\Cc \ = \  \{ G/H \times j_k \}_{k\geq 0, H\in \Cc}
\end{equation}
detects fibrations in the $\Cc$-projective  model structure on $G$-spaces.

All morphisms in $I_\Cc$ and $J_\Cc$ are closed embeddings,
and this property is preserved by coproducts, cobase change and
sequential colimits in the category of $G$-spaces.
Proposition \ref{prop:G/H is finite} guarantees that sources and targets
of all  morphisms in $I_\Cc$ and $J_\Cc$ 
are finite (sometimes called `finitely presented') 
with respect to sequences of closed embeddings of $G$-spaces.
In particular, the sources of all these morphisms
are finite with respect to sequences of $I_\Cc$-cell complexes and
$J_\Cc$-cell complexes.

Now we can prove the factorization axiom MC5.
Every morphism in $I_\Cc$ and $J_\Cc$ is a $\Cc$-cofibration by adjointness.  
Hence every $I_\Cc$-cofibration or $J_\Cc$-cofibration is a $\Cc$-cofibration of
$G$-spaces.
The small object argument applied to the set  
$I_\Cc$ gives a factorization of any morphism of $G$-spaces as
a $\Cc$-co\-fibra\-tion followed by a morphism with the right lifting property
with respect to $I_\Cc$. Since $I_\Cc$ detects the $\Cc$-acyclic $\Cc$-fibrations,
this provides the factorizations as cofibrations followed by acyclic fibrations.

For the other half of the factorization axiom MC5
we apply the small object argument to the set $J_\Cc$; 
we obtain a factorization of any morphism of $G$-spaces as
a $J_\Cc$-cell complex followed by a morphism with the right lifting property
with respect to $J_\Cc$. Since $J_\Cc$ detects the $\Cc$-fibrations,
it remains to show that every  $J_\Cc$-cell complex 
is a $\Cc$-equivalence.
To this end we observe that the morphisms in $J_\Cc$ are inclusions of deformation
retracts internal to the category of $G$-spaces.
This property is inherited by coproducts and cobase changes, 
so every morphism obtained by cobase changes of coproducts
of morphisms in $J_\Cc$ is a homotopy equivalence of $G$-spaces,
hence also a $\Cc$-equivalence. We also need to pass to sequential colimits,
which is fine because $J_\Cc$-cell complexes are closed embeddings,
and taking $H$-fixed points commutes with sequential colimits over
closed embeddings (Proposition \ref{prop:G-fix preserves pushouts}~(ii)).

It remains to prove the other half of MC4, i.e., that every $\Cc$-acyclic
$\Cc$-cofibration $f:A\to B$ has the left lifting property
with respect to $\Cc$-fibrations.
The small object argument provides a factorization 
\[  A \ \xra{\ j\ } \ W\ \xra{\ q\ } \ B\]
as a $J_\Cc$-cell complex followed by a $\Cc$-fibration.
In addition, $q$ is a $\Cc$-equivalence since $f$ and $j$ are.
Since $f$ is a $\Cc$-cofibration, a lifting in 
\[\xymatrix{
A \ar[r]^-j \ar[d]_f & W \ar[d]^q_(.6)\sim \\
B \ar@{=}[r] \ar@{..>}[ur] & B }\]
exists. Thus $f$ is a retract of the morphism $j$ that has the left lifting
property for $\Cc$-fibrations. 
So $f$ itself has the left lifting property for $\Cc$-fibrations.

The model structure is topological by Proposition \ref{prop:topological criterion}.
Right properness of the model structure is a straightforward consequence
of right properness of the model structure on spaces,
since the $H$-fixed point functor
preserves pullbacks and takes $\Cc$-fibrations to Serre fibrations.
Since the projective $\Cc$-model structure is topological
and all objects are fibrant, every cofibration is an h-cofibration
by Corollary \ref{cor-h-cofibration closures}~(iii). 
So left properness follows from the gluing lemma for $\Cc$-equivalences
(Proposition \ref{prop:gluing lemma G-spaces}).
\end{proof}

\begin{defn}
Let $G$ be a topological group.
A {\em universal space}\index{subject}{universal space!for a set of subgroups}
for a set $\Cc$ of closed subgroups of $G$
is a $\Cc$-cofibrant $G$-space $E$ such that
for every subgroup $H\in\Cc$ the fixed point space $E^H$ 
is weakly contractible.
\end{defn}

Any cofibrant replacement of the one-point $G$-space 
in the $\Cc$-projective model structure is a universal space for the set $\Cc$,
so universal spaces exist for any set of closed subgroups.
Moreover, any two universal spaces for the same subgroups
are $G$-equivariantly homotopy equivalent:

\begin{prop}\label{prop:universal spaces}
Let $G$ be a topological group, $\Cc$ a set 
of closed subgroups of $G$, and $E$ a universal $G$-space for the set $\Cc$.
\begin{enumerate}[\em (i)]
\item Every $\Cc$-cofibrant $G$-space admits a continuous $G$-map to $E$,
and any two such maps are homotopic as $G$-maps.
\item If $E'$ is another universal $G$-space for $\Cc$, then every continuous $G$-map
from $E'$ to $E$ is a $G$-equivariant homotopy equivalence.
\end{enumerate}
\end{prop}
\begin{proof}
(i) This all follows from the existence of the $\Cc$-projective model
structure described in Proposition \ref{prop:proj model structures for G-spaces}.
We let $A$ be a $\Cc$-cofibrant $G$-space.
The unique map $E\to \ast$ to a one-point $G$-space is a $\Cc$-equivalence and
a $\Cc$-fibration. So the unique morphism from $A$
to the one-point $G$-space lifts to a continuous $G$-map $A\to E$.
The $\Cc$-projective model structure is topological, so 
the inclusion $A\times\{0,1\}\to A\times [0,1]$ is a $\Cc$-cofibration,
and thus has the left lifting property with respect to $E\to\ast$.
Given two $G$-maps $f,f':A\to E$, any solution to the lifting problem
\[ \xymatrix{ 
A\times\{0,1\}\ar[r]^-{f+f'}\ar[d] & E\ar[d]\\
A\times[0,1]\ar[r] &\ast} \]
is a $G$-homotopy from $f$ to $f'$.
Since universal $G$-spaces are $\Cc$-cofibrant, part~(ii) is a consequence of~(i). 
\end{proof}

A morphism of $G$-spaces is a {\em $G$-cofibration}\index{subject}{G-cofibration@$G$-cofibration} 
if it has the left lifting property with respect 
to all morphisms that are simultaneously weak equivalences and Serre
fibrations on the fixed points for all closed subgroups of $G$.
Equivalently, $G$-cofibrations are the $\All$-cofibrations
in the sense of Proposition \ref{prop:proj model structures for G-spaces}
for the maximal set of all closed subgroups of $G$.
Relative $G$-CW-complexes are special kinds of $I_{\All}$-cell complexes
(with $I_{\All}$ defined in \eqref{eq:I_for_F-proj_on_GT}),
namely those where all equivariant cells of the same
dimension are attachable at once, and in order of increasing dimensions.
Thus we have the following implications between the various kinds of
`nice equivariant embeddings':
\[
\text{relative $G$-CW-complex} \ \Longrightarrow \
\text{$G$-cofibration} \ \Longrightarrow \
\text{h-cofibration of $G$-spaces}
\]
The second implication is Corollary \ref{cor-h-cofibration closures},
applied to the $\All$-projective model structure.
Both implications are strict. 

\begin{prop}\label{prop:cofamily pushout property}
  Let $N$ be a closed normal subgroup of a topological group $G$. 
  Then for every $G$-cofibration $i:A\to B$ the map
  \[ i^N\ : \ A^N\ \to \  B^N \]
  is a $G/N$-cofibration, and the map
  \[   i\cup\text{\em incl}\ : \ A\cup_{A^N} B^N \ \to \ B  \]
  is a $G$-cofibration.
\end{prop}
\begin{proof}
  The class of $G$-cofibrations for which the claim holds is clearly closed under
  coproducts and retracts. We consider a pushout square of $G$-spaces on the left
  \[ \xymatrix{A \ar[r]^-i\ar[d]&B\ar[d] &&
    A^N\ar[r]^-{i^N}\ar[d]&B^N\ar[d]\\  
    C \ar[r]_j& D  &&   C^N \ar[r]_{j^N}& D^N } \]
  such that $i$ is a $G$-cofibration for which the claim holds. 
  Then $i$ is a closed embedding,
  and the square on the right is also a pushout of $G$-spaces
  (by Proposition \ref{prop:G-fix preserves pushouts}~(i)).
  In particular, $j^N$ is a $G$-cofibration whenever $i^N$ is.
  This means that $C\cup_{C^N} D^N$ is also a pushout
  of $C$ and $B^N$ over $A^N$, and hence the square
  \[ \xymatrix@C=12mm{ 
    A\cup_{A^N} B^N\ar[r]^-{i\cup\text{incl}} \ar[d] & B \ar[d]\\
    C\cup_{C^N} D^N\ar[r]_-{j\cup\text{incl}} &  D}\]
  is another pushout of $G$-spaces. The upper horizontal map is a $G$-cofibration
  by hypothesis, hence so is the lower horizontal map.
  So the class of $G$-cofibrations satisfying the claim is closed under cobase change.
  Because $N$-fixed points preserve sequential colimits along closed embeddings
  (Proposition \ref{prop:G-fix preserves pushouts}~(ii)),
  the class of $G$-cofibrations satisfying the claim is also closed under
  sequential composites.
  
  Given these closure properties, it suffices to
  verify the claim for the generating $G$-cofibrations of the form
  \[ G/H\times i_k \ : \ G/H\times \partial D^k \ \to \ G/H\times D^k\ . \]
  Because $(G/H)^N$ is either empty (whenever $N$ is not contained in $H$)
  or all of $G/H$ (whenever $N\leq H$),
  the map in question is either the map $G/H\times i_k$
  or the identity of $G/H\times D^k$; in either case it is a $G$-cofibration.
  This proves the claim for the generating cofibrations, and thus concludes the proof. 
\end{proof}

Now we specialize from general topological groups to 
the class of compact topological groups. The first useful
special feature is that in this restricted context,
orbit spaces are always compactly generated without
any need to force the weak Hausdorff condition.
The arguments involved in the following proposition are well-known,
see for example \cite[1.1.3]{palais-classification} or \cite[I Cor.\,1.3]{bredon-intro}
for closely related statements.
Since the references I know of work in slightly different categories,
I spell out the proof.

\begin{prop}\label{prop:G quotient properties}
  Let $G$ be a compact topological group and $X$ a $G$-space.
  \begin{enumerate}[\em (i)]
  \item The orbit space $G\bs X$ is compactly generated in the quotient topology.
  \item The projection $\Pi_X:X\to G\bs X$ is both open and closed.
  \item For every $G$-invariant closed subset $Y$ of $X$, the tautological map
    $ G\bs Y \to G\bs X$ is a closed embedding.
  \end{enumerate}
\end{prop}
\begin{proof}
(i)
Proposition \ref{prop:closed implies WH} reduces the claim to showing that
the orbit equivalence relation
\[ E \ = \ \{  (g x, x)\ | \ g\in G, \ x\in X \} \]
is closed in $X\times X$.
Since $G$ is compact, projection away from $G$ is a closed map;
so the composite
\[ G\times X\times X \ \xra{(g,x,y)\mapsto (g,g x,y)} \
 G\times X\times X \ \xra{\ \text{proj}\ } \ X\times X  \]
is a closed map since the first map is a homeomorphism.
Since $X$ is weak Hausdorff, the diagonal $\Delta_X$ 
is closed in $X\times X$. So $G\times\Delta_X$ is closed in
$G\times X\times X$.
The orbit relation $E$ is the image of $G\times\Delta_X$ under the above composite,
so $E$ is closed in $X\times X$.

(ii)
If $O$ is an open subset of $X$, then $g O$ is open for every $g\in G$, 
since left translation by $g$ is a homeomorphism. So 
\[ \Pi_X^{-1}(\Pi_X(O))\ = \ \bigcup_{g\in G}\, g O \]
is open as a union of open subsets. Hence $\Pi_X(O)$ is open in the quotient topology.
The shearing map
\[ \chi\ : \ G\times X \ \to \ G\times X \ , \quad (g,x)\ \longmapsto \ (g, g x) \]
is a homeomorphism, and the composite
\[ G\times X \ \xra{\ \chi\ }\ G\times X \ \xra{\text{proj}}\ X \]
is the action map $\alpha:G\times X\to X$. 
Since $G$ is compact, projection away from $G$ is a closed map; 
so the action map is also a closed map.
If $A$ is closed in $X$, then $G\times A$ is closed in $G\times X$.
So the set
\[ \Pi_X^{-1}(\Pi_X(A))\ = \  \alpha(G\times A)\]
is closed in $X$. Hence $\Pi_X(A)$ is closed in the quotient topology.

(iii) The tautological map $G\bs Y\to G\bs X$ is continuous 
and injective. We show that the map is also closed.
We denote by $\Pi_Y:Y\to G\bs Y$ the quotient map.
We let $A\subset G\bs Y$ be any closed subset. Then $\Pi_Y^{-1}(A)$ is closed in $Y$,
hence also in $X$.
Since $\Pi_X$ is a closed map by part~(ii), $\Pi_X(A)$ is closed in $G\bs X$.
\end{proof}

Now we specialize our discussion to compact Lie groups.
If $K$ is a closed subgroup of a compact Lie group $G$ of smaller dimension,
then the underlying $K$-space of a $G$-CW-complex 
need not admit a $K$-CW-structure -- 
an example is given by Illman in \cite[Sec.\,2]{illman-restricting equivariance}.
Nevertheless, the underlying $K$-space of a $G$-CW-complex 
is always $K$-homotopy equivalent to a $K$-CW-complex, which can be chosen
to be compact if the original $G$-space is,
see \cite[Thm.\,A]{illman-restricting equivariance}.

The next proposition shows that with respect to restriction of group actions,
the class of cofibrant equivariant spaces
is better behaved than equivariant CW-complexes.
An additional advantage of cofibrant spaces over CW-complexes
is that `cofibrant' is a property, whereas CW-structures are additional data.

\begin{prop}\label{prop:cofibrancy preservers} 
Let $G$ be a compact Lie group.
\begin{enumerate}[\em (i)]
\item For every compact Lie group $K$ and every continuous homomorphism $\alpha:K\to G$
  the restriction functor $\alpha^*:G\bT\to K\bT$ 
  takes $G$-cofibrations to $K$-cofibrations.
\item For every closed subgroup $H$ of $G$ 
  the induction functor $G\times_H -:H\bT\to G\bT$
  takes $H$-cofibrations to $G$-cofibrations.
\item For every closed normal subgroup $N$ of $G$ 
  the orbit functor $N\bs -:G\bT\to (G/N)\bT$
  takes $G$-cofibrations to $G/N$-cofibrations.
\end{enumerate}
\end{prop}
\begin{proof}
(i) The restriction functor $\alpha^*$ preserves colimits,
so we may show that it takes the generating $G$-cofibrations
$G/H\times i_k:G/H\times \partial D^k\to G/H\times D^k$ to $K$-cofibrations, 
for any closed subgroup $H$ of $G$. 
The $K$-action on the smooth compact manifold $G/H$ is by
left translation through $\alpha$; 
continuous homomorphisms between Lie groups are automatically smooth,
compare \cite[Prop.\,I.3.12]{broecker-tomDieck}, so the $K$-action on $G/H$ is smooth.
Illman's theorem \cite[Cor.\,7.2]{illman} provides a finite $K$-CW-structure on $G/H$,
so in particular $G/H$ is cofibrant as a $K$-space. Since the projective
model structure on $K$-spaces (for the set of all closed subgroups) is topological,
the map $G/H\times i_k$ is a $K$-cofibration.

(ii) Since $G\times_H -$ preserves colimits,
it suffices to show that it takes the generating $H$-cofibrations
$H/J\times i_k:H/J \times \partial D^k\to H/J \times D^k$ to $G$-cofibrations, 
for any closed subgroup $J$ of $H$.
This in turn is clear since $G\times_H (H/J)$ is $G$-homeomorphic to $G/J$. 

(iii)  Since $N\bs -$ preserves colimits,
it suffices to show that it takes the generating $G$-cofibrations
$G/H\times i_k:G/H \times \partial D^k\to G/H \times D^k$ to $G/N$-cofibrations, 
for any closed subgroup $H$ of $G$.
This in turn is clear since $N\bs (G/H)$ is $(G/N)$-homeomorphic to $(G/N)/(H/H\cap N)$. 
\end{proof}

Now we show that taking cartesian product preserves equivariant cofibrations.
There are two related questions, 
namely `external products' of a $G$-space and a $K$-space, 
and `internal products' of two $G$-spaces with diagonal action.

\begin{prop}
  Let $G$ and $K$ be topological groups.
  \begin{enumerate}[\em (i)]
  \item The pushout product of a $G$-cofibration with a $K$-cofibration is
    a $(G\times K)$-cofibration.
  \item If $G$ is a compact Lie group, then the pushout product of two $G$-cofibrations
    is a $G$-cofibration with respect to the diagonal $G$-action.
  \end{enumerate}
\end{prop}
\begin{proof}
(i) The product functor 
\[ \times \ : \ G\bT \times K\bT \ \to \ (G\times K)\bT \]
preserves colimits in each variable, so it suffices to check the 
pushout product of a generating $G$-cofibration $G/H\times i_k$ with
a generating $K$-cofibration $K/L\times i_m$, where $i_k:\partial D^k\to D^k$ is
the inclusion. The pushout product of these is isomorphic to
\[ (G\times K)/(H\times L)\times i_{k+m} \]
and hence a cofibration of $(G\times K)$-spaces.
  
(ii) By~(i) the pushout product of two $G$-cofibrations
is a $(G\times G)$-cofibration. Since $G$ is compact Lie,
restriction along the diagonal embedding $G\to G\times G$
preserves cofibrations by Proposition \ref{prop:cofibrancy preservers}~(i). 
\end{proof}

Now we prove a decomposition result for certain kinds of fixed points.
We let $G$ and $K$ be topological groups and $X$ a $(K\times G)$-space. 
We want to describe
the $K$-fixed points $(G\bs X)^K$ of the $G$-orbit space $G\bs X$.
For  a continuous homomorphism $\alpha:K\to G$ we set
\[ X^{\alpha}\ = \  \{x\in X\ |\ (k,\alpha(k)) x = x\text{ for all $k\in K$}\}\ .\]
Equivalently, $X^{\alpha}$ is the fixed point space of the graph of $\alpha$,
which is a closed subgroup of $K\times G$.
The subspace $X^{\alpha}$ of $X$ is invariant under the action of the centralizer 
of the image of $\alpha$, i.e., the group
\[ C(\alpha) \ = \ 
\{g\in G\ |\ g\alpha(k)=\alpha(k)g \text{ for all $k\in K$}\} \ .\]
The inclusion $X^\alpha\to X$ then passes to a continuous map
\[ \alpha^\flat\ : \  C(\alpha) \ \bs X^{\alpha} \ \to \  G\bs X \]
on orbit spaces. For $x\in X^\alpha$ and $k\in K$ the relation
\[ k(G x)\ = \ G (k x) \ = \ G (\alpha(k)^{-1} x) \ = \  G x \]
shows that $\alpha^\flat$ takes values in the $K$-fixed points $(G\bs X)^K$.
Moreover,  for $g\in G$ we have
\begin{equation}\label{eq:g_times_X^alpha}
 g\cdot X^\alpha \ = \ X^{c_g\circ\alpha}   
\end{equation}
as subspaces of $X$. So the maps $\alpha^\flat$ and $(c_g\circ\alpha)^\flat$
arising from conjugate homomorphisms 
have the same image in the orbit space $G\bs X$.
It is relatively straightforward to see that the coproduct of
all the maps $\alpha^\flat$ is bijective if the $G$-action is free; 
some pointset topology is involved in
showing that this continuous bijection is in fact a homeomorphism
if in addition $G$ and $K$ are compact Lie groups.

\begin{prop}\label{prop:fix of free cofibration}
Let $G$ and $K$ be compact Lie groups.
Let $X$ be a $(K\times G)$-space such that $G$-action is free.
Then the map
\[ \coprod \alpha^\flat\ : \  
 {\coprod}_{\td{\alpha}}  \, C(\alpha) \bs X^{\alpha}\ \to \   (G\bs X)^K  \]
is a homeomorphism.
Here the disjoint union runs over conjugacy classes of continuous homomorphisms 
$\alpha:K\to G$ and $C(\alpha)$ is the centralizer, in $G$, of the image of $\alpha$.
\end{prop}
\begin{proof}
We let $\Pi:X\to G\bs X$ denote the quotient map.  
We set
\[ \bar X\ = \ \Pi^{-1}((G\bs X)^K) \ . \]
Since $X$ is compactly generated, so is $G\bs X$ 
by Proposition \ref{prop:G quotient properties}~(i).
Since $(G\bs X)^K$ a closed subset of the orbit space, 
$\bar X$ is a $(K\times G)$-invariant closed subspace of $X$. 
In particular, $\bar X$ is compactly generated in the subspace topology.

For a given continuous homomorphism $\alpha:K\to G$, we set
\[  X^{(\alpha)}\ = \ G\cdot X^{\alpha}\ , \]
the smallest $G$-subspace of $\bar X$ containing $X^\alpha$. 
The relation \eqref{eq:g_times_X^alpha} shows that $X^{(\alpha)}$
is the union of the subsets $X^{\alpha'}$ as $\alpha'$ runs 
over all conjugates of $\alpha$.
We factor the map in question as a composite
\[  {\coprod}_{\td{\alpha}}  \, C(\alpha) \bs X^{\alpha}\ \to \ 
 {\coprod}_{\td{\alpha}}  \, G \bs X^{(\alpha)}\ \to \ G\bs \bar X\ \to \  (G\bs X)^K \ , \]
induced by the inclusions $X^{\alpha}\to X^{(\alpha)}\to \bar X \to X$;
we show that each of the three maps is a homeomorphism.
The third map $G\bs \bar X\to (G\bs X)^K$ is a homeomorphism
by Proposition \ref{prop:G quotient properties}~(iii).

For every $x\in \bar X$ and every $k\in K$ we have $k\cdot (G x)=G(k x)= G x$,
so there exists a $g\in G$ such that $k x= g^{-1} x$.
Since the $G$-action is free, the element $g$ is uniquely determined
by this property. So we can define a map
$\beta_x:K\to G$ by the property $k x=\beta_x(k)^{-1} x$;
equivalently, the characterizing condition for $\beta_x$ is that
$(k,\beta_x(k))\cdot x = x$. 
It is straightforward to see that $\beta_x$ is a group homomorphism. 

By definition, the stabilizer group of $x\in \bar X$ inside $K\times G$
is precisely the graph of the homomorphism $\beta_x:K\to G$.
Since $X$ is compactly generated, this stabilizer group is a closed
subset of $K\times G$, which means that the homomorphism $\beta_x$
is continuous.  
Our next claim is that the assignment
\[ \beta \ : \ \bar X\ \to \ \hom(K,G) \ , \quad x \ \longmapsto \ \beta_x\]
is continuous.
Since $\hom(K,G)$ has the subspace topology of $\map(K,G)$, 
it suffices to show that the adjoint map
\[ \bar\beta \ : \ \bar X \times K \ \to \ G \ , \quad (x,k)\ \longmapsto \ \beta_x(k)\]
is continuous. This, in turn, is equivalent to the claim that the graph of $\bar\beta$
is closed as a subset of $\bar X\times K\times G$.
The graph of $\bar\beta$ is the inverse image of the diagonal under the 
continuous map
\[ \bar X\times K\times G \ \to \ \bar X\times \bar X\ , \quad
(x,k,g)\ \longmapsto \ (k x, g x)\ . \]
So the graph of $\bar\beta$ is closed, hence $\bar\beta$ and $\beta$ are
continuous.

By Proposition \ref{prop:components of hom(K,G)} the space $\hom(K,G)$ 
is the topological disjoint union of the $G^\circ$-orbits under the conjugation action.
In particular, every $G^\circ$-conjugacy class is open and closed.
A $G$-conjugacy class is a finite union of $G^\circ$-conjugacy classes,
so every $G$-conjugacy class in $\hom(K,G)$ is also open and closed.
Since $\beta:\bar X\to \hom(K,G)$ is continuous, 
$\bar X$ is the topological disjoint union of the subsets $X^{(\alpha)}$
as $\alpha$ ranges over all conjugacy classes in $\hom(K,G)$. 
Taking orbits commutes with disjoint unions, so the canonical map
from ${\coprod}_{(\alpha)}  G\bs X^{(\alpha)}$ to $G\bs \bar X$ is a homeomorphism.

The final step is to show that for every continuous homomorphism $\alpha:K\to G$
the canonical continuous map
\begin{equation}  \label{eq:C_bs_X2G_bs_X}
 C(\alpha)\bs X^\alpha\ \to \ G\bs X^{(\alpha)}   
\end{equation}
is a homeomorphism. The map is surjective because $X^{(\alpha)}=G\cdot X^{\alpha}$.
The map is also injective; if $x,y\in X^{\alpha}$ satisfy $G x=G y$,
then $x=g y$ for some $g\in G$, and so $x\in X^{\alpha}\cap X^{c_g\circ\alpha}$.
This implies that $\alpha=c_g\circ\alpha$, and so $g\in C(\alpha)$.
Hence $C(\alpha) x= C(\alpha) g y =C(\alpha) y$.
Finally, the map \eqref{eq:C_bs_X2G_bs_X} is closed:
if $O\subset C(\alpha)\bs X^{\alpha}$ is a closed subset,
then $\Pi_{X^{\alpha}}^{-1}(O)$ is closed in $X^{\alpha}$.
Since $X^{\alpha}$ is closed in $X$, hence also in $X^{(\alpha)}$,
the set $\Pi_{X^{\alpha}}^{-1}(O)$ is also closed in $X^{(\alpha)}$.
Since the projection $X^{(\alpha)}\to G\bs X^{(\alpha)}$ is a closed map,
the image of $O$ in $G\bs X^{(\alpha)}$ is closed. This completes the proof
that the map \eqref{eq:C_bs_X2G_bs_X} is a homeomorphism.
\end{proof}

\begin{prop}
Let $G$ be a compact Lie group and $A$ a free cofibrant $G$-space.
Then the functor $A\times_G -$ takes $G$-maps that are
non-equivariant weak equivalences to weak equivalences.
\end{prop}
\begin{proof}
We let $f:X\to Y$ be a continuous $G$-map
that is a non-equivariant weak equivalence.
We start with the special case where $A$ is a $G$-CW-complex with skeleton filtration
\[ \emptyset = A^{-1} \ \subset \ A^0 \ \subset \ A^1 \ \subset \ \dots \ \subset \ A^n 
\ \subset \ \dots\ .\]
We show by induction over $n$ that the map $A^n\times_G f:A^n\times_G X\to A^n\times_G Y$
is a weak equivalence. The induction start with $n=-1$, where there is nothing to show.
Now we consider $n\geq 0$ and assume the claim for smaller values of $n$.
Since $G$ acts freely on $A$, only free $G$-cells occur in the equivariant CW-structure.
So there is an index set $I$ and a pushout square of $G$-spaces:
\[ \xymatrix{
G\times I \times \partial D^n \ar[r]^-{\text{incl}}\ar[d]_\alpha & 
G\times I \times D^n \ar[d]\\
A^{n-1}\ar[r]& A^n} \]
The map $A^n\times_G f$ can thus be obtained by passing to pushouts 
in the horizontal direction of the commutative diagram
\[ \xymatrix@C=12mm{ 
A^{n-1}\times_G X\ar[d]_{A^{n-1}\times_G f} & 
I\times\partial D^n\times X\ar[d]^{I\times\partial D^n\times f} 
\ar[l]_-{\alpha\times_G X}\ar[r]^-{\text{incl}} &
I\times D^n\times X\ar[d]^{I\times D^n\times f}\\
A^{n-1}\times_G Y& I\times\partial D^n\times Y \ar[l]^{\alpha\times_G Y}\ar[r]_-{\text{incl}} &
I\times D^n\times Y} \]
The left vertical map is a weak equivalence by the inductive hypothesis,
and the middle and right vertical maps are weak equivalences because $f$ is.
The upper and lower right horizontal maps are h-cofibrations,
so the gluing lemma for weak equivalences and pushouts along h-cofibrations 
(see for example \cite[Appendix, Prop.\,4.8 (b)]{boardman-vogt-homotopy invariant})
shows that then the map on pushouts $A^n\times_G f$ is also a weak equivalence.
Instead of the gluing lemma, we can alternatively quote \cite[Lemma A.1]{dugger-isaksen}.

The space $A\times_G X$ is the sequential colimit, along h-cofibrations,
of the spaces $A^n\times_G X$, and similarly for $A\times_G Y$.
In the category $\bT$, all h-cofibrations are closed embeddings 
(Proposition \ref{prop:h-cof is closed embedding}),
and sequential colimits over sequences of closed embeddings in $\bT$
preserve weak equivalences 
(Proposition \ref{prop:closed seq colim weak equivalences}~(i)).
So the map $A\times_G f$ is a weak equivalence.
This completes the proof in the special case when $A$ admits a $G$-CW-structure.
A general cofibrant $G$-space is $G$-homotopy equivalent to a 
$G$-CW-complex, and the functor $-\times_G X$ takes
$G$-homotopy equivalences to non-equivariant homotopy equivalences.
This reduces the general case to the case of $G$-CW-complexes.
\end{proof}

Now we turn to $\bGamma$-spaces and study the categorical and homotopical
properties of prolongation of a $\bGamma$-space to a continuous
functor defined on all based spaces.
In the body of this book, we always work in the category $\bT$ of compactly
generated spaces. However, this full subcategory is not closed under
quotient spaces nor coends inside the ambient category
of all topological spaces; since the construction of the prolongation
involves a coend (quotient space), some care needs to be taken. 
So while we are mostly interested  $\bGamma$-spaces with values 
in the full subcategory $\bT$ of compactly generated spaces, 
we define and discuss prolongation for $\bGamma$-$k$-spaces.

\begin{construction}
We let $\bGamma$ denote the category whose objects are the based sets
$n_+=\{0,1,\dots,n\}$, with basepoint~0, and with morphisms all based maps.
A {\em $\bGamma$-$k$-space} is a functor from $\bGamma$ 
to the category of $k$-spaces that is reduced,
i.e., the value at $0_+$ is a one-point space.

A $\bGamma$-$k$-space $F:\bGamma\to\bK$ 
can be `evaluated' at a based $k$-space $K$ as follows.
We set
\[  F(K) \ = \ \left({\coprod}_{n\geq 0}\,   F(n_+)\times  K^n\right) /\sim \ .\] 
In more detail, $F(K)$ is the quotient space of the disjoint union of the spaces
$F(n_+)\times  K^n$ by the equivalence relation generated by
\begin{equation}  \label{eq:Gamma relation generator}
 (F(\alpha)(x);\,k_1,\dots,k_n)\ \sim\ (x;\, k_{\alpha(1)},\dots,k_{\alpha(m)})   
\end{equation}
for all $x\in F(m_+)$, all $(k_1,\dots,k_n)$ in $K^n$,  and all 
morphisms $\alpha:m_+\to n_+$ in $\bGamma$.
Here $k_{\alpha(i)}$ is to be interpreted as the basepoint of $K$
whenever $\alpha(i)=0$.  
We emphasize that we use the Kelleyfied product topology 
on $F(n^+)\times K^n$, in order to remain inside the category of $k$-spaces.
We mostly care about the special case when $K$ is compact, and then
this coincides with the ordinary product topology 
$F(m^+)\times_0 K^{\times_0 n}$ by Proposition \ref{prop:properties k-spaces}~(vi).
The class of $k$-spaces is closed under disjoint unions and
passage quotient spaces (by Proposition \ref{prop:properties k-spaces}~(i));
so $F(K)$ is indeed a $k$-space.
We refer to the extension as the 
{\em prolongation}\index{subject}{prolongation!of a $\bGamma$-space}
of the $\bGamma$-space $F$.
A more categorical way to describe $F(K)$ is as a coend of the functor
\begin{equation}\label{eq:Gamma coend}
  \bGamma\times\bGamma^{\op}\ \to \ \bK\ , \quad 
  (m^+,n^+)\ \longmapsto \ F(m^+)\times K^n \ .  
\end{equation}
\end{construction}

 \begin{rk}\label{rk:abuse prolongation}
  We want to justify the abuse of notation of not distinguishing 
  between the original $\bGamma$-space and its prolongation:
  the value of the prolongation on $m_+$ is canonically homeomorphic 
  to the original value.   
  Slightly more is true, namely that the coend description 
  exhibits the prolongation as the left Kan extension of $F$ 
  along the inclusion $\bGamma\to \bK_*$; we won't use this, however.

  Only for this argument we denote the prolongation of $F$ by $\hat F$. 
  The continuous maps
  \[ F(n_+)\times (m_+)^n \ \to \ F(m_+)\ , \quad
  (x;\, i_1,\dots,i_n)\ \longmapsto \ F(i)(x)\]
  are compatible with the equivalence relation defining $\hat F(m_+)$;
  here $i:n_+\to m_+$ is the based map with $i(j)=i_j$.
  So these maps assemble into a continuous map $\hat F(m_+)\to F(m_+)$
  that is inverse to the continuous map
  \[ F(m_+)\ \to \  \hat F(m_+)\ , \quad x \ \longmapsto \ [x;\,1,2,\dots,m]\ .\]
  \end{rk}

We have to analyze the equivalence relation
generated by \eqref{eq:Gamma relation generator} more closely.
This part of the argument has nothing to do with topology, is purely combinatorial,
and of a very similar flavor as the analysis of the 
equivalence relation defining the geometric realization of a simplicial set.
We call an element $x\in F(m_+)$ {\em degenerate}
if it is in the image of the map $F(\alpha):F((m-1)_+)\to F(m_+)$
for some morphism $\alpha:(m-1)_+\to m_+$ in $\bGamma$.
We call $x$ {\em non-degenerate} if it is not degenerate.
We let $C_m(K)\subset K^m$ be the set of those tuples
$(k_1,\dots,k_m)$ whose coordinates $k_i$ are pairwise distinct and all different 
from the basepoint of $K$.
Part~(i) of the following proposition 
is a special case of \cite[Prop.\,6.9]{berger-moerdijk-Reedy};
part~(ii) also ought to be well known, but I am not aware of a reference.
The next proposition in particular implies that
the reduction map \eqref{eq:reduction_map} defined in the proof of part~(ii)
is a bijection from $F(K)$ to the set-theoretic disjoint union, 
for $m\geq 0$, of the sets $F(m_+)^{\text{nd}}\times_{\Sigma_m}C_m(K)$,
where $F(m_+)^{\text{nd}}$ is the set of non-generate elements of $F(m_+)$.

\begin{prop}\label{prop:Gamma minimal representative}
Let $F:\bGamma\to\text{\em (sets)}$ be a functor such that $F(0_+)$ has one element.
\begin{enumerate}[\em (i)]
\item For every element $y\in F(m_+)$ there exists an injective morphism
$\delta:l_+\to m_+$ and a non-degenerate element $x\in F(l_+)^{\text{\em nd}}$
such that $y=F(\delta)(x)$.
If moreover $y=F(\bar\delta)(\bar x)$ for 
another injective morphism $\bar\delta:\bar l_+\to m_+$ 
and a non-degenerate element $\bar x\in F(\bar l_+)^{\text{\em nd}}$,
then $l=\bar l$ and there exists a bijective morphism $\lambda:l_+\to l_+$
such that $\delta=\bar\delta\lambda$ and $\bar x=F(\lambda)(x)$.
\item
Let $K$ be a based set. Let $(x,t)\in F(l_+)\times K^l$ be an element
of minimal dimension $l$ within its equivalence class.
Let $(y,s)\in F(m_+)\times K^m$ be equivalent to $(x,t)$.
Then there exists a surjective morphism
$\sigma:m_+\to k_+$, an injective morphism $\delta:l_+\to k_+$
and $u\in C_k(K)$ such that
\[ F(\sigma)(y)\ = \ F(\delta)(x)\ , \quad 
s\ = \  \sigma^*(u) \text{\qquad and\qquad} t =\delta^*(u) \ .\]
\item 
If $(x,t)$ and $(y,s)$ are equivalent elements, both of minimal dimension $l$
in their equivalence class, then
there is a unique isomorphism $\alpha:l_+\to l_+$ such that
$(F(\alpha)(y),s)=(x,\alpha^*(t))$. 
\end{enumerate}
\end{prop}
\begin{proof}
Part~(i) is an analog of the `Eilenberg-Zilber lemma' \cite[(8.3)]{eilenberg-zilber-ssc}
for simplicial sets, and a special case of \cite[Prop.\,6.9]{berger-moerdijk-Reedy}.
For the convenience of the reader, we give a direct proof. 
The existence part is proved by induction on $m$, starting with $m=0$,
where there is nothing to show. If $m$ is positive and $y$ is non-degenerate,
then $x=y$ and $\delta=\Id$ do the job. Otherwise 
$y=F(\alpha)(z)$ for some morphism $\alpha:(m-1)_+\to m_+$ and $z\in F((m-1)_+)$.
We factor $\alpha=\delta'\circ \beta$ such that $\beta:(m-1)_+\to k_+$ is surjective
and $\delta':k_+\to m_+$ is injective. We must have $k<m$, so the inductive
hypothesis provides an injective morphism
$\delta:l_+\to k_+$ and a non-degenerate element $x\in F(l_+)^{\text{nd}}$
such that $F(\beta)(z)=F(\delta)(x)$. But then $\delta'\delta$ is also injective and
\[ y \ = \ F(\alpha)(z)\ = \ F(\delta')(F(\beta)(z))\ = \
F(\delta')(F(\delta)(x))\ = \ F(\delta'\delta)(x)\ . \]
This proves the first statement.

For the second statement we consider injective morphisms $\delta:l_+\to m_+$
and $\bar\delta:\bar l_+\to m_+$, and non-degenerate elements $x\in F(l_+)^{\text{nd}}$
and $\bar x\in F(\bar l_+)^{\text{nd}}$ such that $F(\delta)(x)=F(\bar\delta)(\bar x)$.
We let $\bar\sigma:m_+\to \bar l_+$ be the unique morphism
that sends all elements not in the image of $\bar\delta$ to the basepoint $0$
and satisfies $\bar\sigma\bar\delta=\Id$.
Then
\[ F(\bar\sigma \delta)(x)\ = \ F(\bar\sigma)(F(\delta)(x))\ = \ 
F(\bar\sigma)(F(\bar\delta)(\bar x))\ = \ 
F(\bar\sigma\bar \delta)(\bar x)\ = \ \bar x \ . \]
If $\bar\sigma\delta:l_+\to \bar l_+$ were not surjective,
then it could be factored through $(\bar l-1)_+$, and $\bar x$ would be degenerate,
contradicting the assumptions. So $\bar\sigma\delta$ is surjective,
and hence $l\geq \bar l$. Reversing the roles of $(x,\delta)$ and $(\bar x,\bar\delta)$
gives $l\leq \bar l$, and hence $l=\bar l$. 
Since $l=\bar l$ and $\bar\sigma\delta$ is surjective, it must in fact
be bijective.

Now we claim that the morphisms $\delta,\bar\delta:l_+\to m_+$ have the same image. 
If this were not the case,
then $\bar\sigma\delta:l_+\to l_+$ would not be surjective,
which again would contradict the assumption that $\bar x$ is non-degenerate.
Since $\delta$ and $\bar\delta$ have the same image, there is 
a bijection $\lambda:l_+\to l_+$ such that $\delta=\bar\delta\lambda$.
Hence 
\[  F(\lambda)(x)\ = \ F(\bar\sigma\bar\delta\lambda)(x)\ = \ 
F(\bar\sigma)(F(\delta)(x))\ = \ 
F(\bar\sigma)(F(\bar\delta)(\bar x))\ = \ \bar x\ .\]
So $\lambda$ is the bijection with the desired properties.

(ii)
We call a quadruple $(\sigma,\delta,u,x)$ consisting of
a surjective morphism $\sigma:m_+\to k_+$, 
an injective morphism $\delta:l_+\to k_+$, 
a tuple $u\in C_k(K)$ and a non-degenerate element $x\in F(l_+)$
a {\em reduction datum} for $(y,s)\in F(m_+)\times K^m$ if 
\[ F(\sigma)(y)=F(\delta)(x)\text{\qquad and\qquad}  s=\sigma^*(u)
 \ .\]
Since $\delta$ is injective, the map $\delta^*:K^k\to K^l$ omits
some coordinates and reorders the remaining ones. In particular,
$\delta^*$ does not introduce duplicates or basepoints, 
so it sends the subset $C_k(K)$ to the subset $C_l(K)$, and hence
$\delta^*(u)\in C_l(K)$.

Claim 1:
Every pair $(y,s)$ has a reduction datum.
We let $u=(u_1,\dots,u_k)$ be the non-basepoint elements in the set $\{s_1,\dots,s_m\}$, 
without repetition and in some chosen order.
We define a surjective morphism $\sigma:m_+\to k_+$ by
\[ \sigma(i) \ = \
\begin{cases}
j  & \text{\ if $s_i=u_j$, and}\\
0 &\text{ if $s_i$ is the basepoint of $K$.}
\end{cases}
\]
Then $u\in C_k(K)$ and $s=\sigma^*(u)$. 
Part~(i) now provides an injective morphism $\delta:l_+\to k_+$
and $x\in F(l_+)^{\text{nd}}$ such that $F(\sigma)(y)=F(\delta)(x)$.
So $(\sigma,\delta,u,x)$ is a reduction datum for $(y,s)$.

Claim 2:
We show that the reduction datum is `essentially unique' in the following sense.
If
\[ (\sigma:m_+\to k_+,\delta:l_+\to k_+,u,x)\text{\quad and\quad} 
(\bar\sigma:m_+\to \bar k_+,\bar\delta:\bar l_+\to \bar k_+,\bar u,\bar x) \]
 are reduction data for the same element $(y,s)$,
then $k=\bar k$, $l=\bar l$ and there are bijective morphisms
$\lambda:l_+\to l_+$ and $\beta:k_+\to k_+$ such that
\[ \bar\sigma=\beta \sigma\ , \quad \bar\delta\lambda = \beta\delta\ ,\quad
u = \beta^*(\bar u)\text{\qquad and\qquad}  \bar x =F(\lambda)(x)\ . \]
Indeed, because $s=\sigma^*(u)=\bar\sigma^*(\bar u)$
and $u$ and $\bar u$ don't contain duplicates or basepoints, we must have $k=\bar k$
and $u$ and $\bar u$ can only differ by the ordering.
So there is a bijective morphism $\beta:k_+\to k_+$ such that
$u=\beta^*(\bar u)$. For all $i\in \{1,\dots,m\}$ we then have
\[ \bar u_{\bar\sigma(i)}\ = \ s_i\ = \ u_{\sigma(i)} \ = \ \bar u_{\beta(\sigma(i))}\ .\]
Since the coordinates of $\bar u$ are all distinct, this implies
that $\bar\sigma=\beta\sigma$.
Using that both quadruples are reduction data, we know that 
\[ F(\bar\delta)(\bar x)\ = \
F(\bar\sigma)(y)\ = \  F(\beta)(F(\sigma)(y)) \ = \ 
F(\beta)(F(\delta)(x))\ = \ F(\beta\delta)(x)\ . \]
Since $\bar\delta$ and $\beta\delta$ are both injective and
$x$ and $\bar x$ are both non-degenerate, then essential uniqueness
statement in part~(i)
shows that $l=\bar l$ and provides a bijective morphism $\lambda:l_+\to l_+$
such that $\beta\delta=\bar\delta\lambda$ and $F(\lambda)(x)=\bar x$.

\medskip

We can now define a {\em reduction map}
\begin{equation}  \label{eq:reduction_map}
 \rho \ : \ {\coprod}_{n\geq 0}\, F(n_+)\times K^n\ \to \ 
{\coprod}_{l\geq 0}\, F(l_+)^{\text{nd}}\times_{\Sigma_l} C_l(K)   
\end{equation}
by choosing a reduction datum $(\sigma,\delta,u,x)$ 
for a given element $(y,s)$ and setting
\[ \rho(y,s)\ = \ [ x,\delta^*(u) ] \ ,\]
where $[-,-]$ denotes the $\Sigma_l$-orbit.
If $(\sigma,\delta,u,x)$ is another reduction datum for $(y,s)$, then
\[ [x,\delta^*(u)] \ = \ 
[x,\delta^*(\beta^*(\bar u))] \ = \ 
[x,\lambda^*(\bar \delta^*(\bar u))] \ = \ 
[F(\lambda)(x),\bar\delta^*(\bar u)] \ = \ [\bar x,\bar\delta^*(\bar u)]\ ,\]
by Claim~2. So $\rho(y,s)$ is well-defined.

\medskip

Claim 3:
If $(y,s)\in F(m_+)\times K^m$ and $(\bar y,\bar s)\in F(\bar m_+)\times K^{\bar m}$
are equivalent, then $\rho(y,s)=\rho(\bar y,\bar s)$.
It suffices to show the claim whenever $(y,s)$ and $(\bar y,\bar s)$
are related by a generating relation \eqref{eq:Gamma relation generator},
i.e., we can assume that $y=F(\alpha)(\bar y)$ and $\bar s=\alpha^*(s)$
for some morphism $\alpha:\bar m_+\to m_+$.
We let $(\sigma,\delta,u,x)$ be a reduction datum for $(y,s)$.
We choose a factorization
\[ \sigma\circ\alpha \ = \ \bar\delta\circ \bar\sigma \]
as a surjective morphism $\bar\sigma:\bar m_+\to \bar k_+$
followed by an injective morphism $\bar\delta:\bar k_+\to k_+$.
Using part~(i) we write
\[ F(\bar\sigma)(\bar y) \ = \ F(\delta')(\bar x)\]
for an injective morphism $\delta':l'_+\to \bar k_+$
and a non-degenerate element $\bar x\in F(l'_+)^{\text{nd}}$.
Then
\begin{align*}
  F(\delta)(x)\ &= \ F(\sigma)(y)\ = \ F(\sigma)(F(\alpha)(\bar y)) \\ 
&= \ F(\bar\delta)(F(\bar\sigma)(\bar y)) \ 
= \ F(\bar\delta)(F(\delta')(\bar x)) \ = \ F(\bar\delta\delta')(\bar x) \ . 
\end{align*}
Since $x$ and $\bar x$ are non-degenerate and $\delta$ and $\bar\delta\delta'$
are injective, the essential uniqueness of part~(i)
shows that $l=l'$ and provides a bijection $\lambda:l_+\to l_+$
such that $\delta=\bar\delta\delta'\lambda$ and $\bar x=F(\lambda)(x)$.
We conclude that
\[
F(\bar\sigma)(\bar y)\ = \ F(\delta')(\bar x)\ = \ 
F(\delta')(F(\lambda)(x))\ = \ 
F(\delta'\lambda)(x)\]
 and
\[
\bar s\ = \ \alpha^*(s)\ = \ \alpha^*(\sigma^*(u))\ = \ 
\bar\sigma^*( \bar\delta^*(u)) \ .\]
Since $u\in C_k(K)$ and $\bar\delta$ is injective, $\bar\delta^*(u)$
belongs to $C_{\bar k}(K)$.
This shows that $(\bar\sigma, \delta'\lambda, \bar\delta^*(u),x)$
is a reduction datum for $(\bar y,\bar s)$.
So
\[ \rho(\bar y,\bar s)\ = \ [x, (\delta'\lambda)^*(\bar\delta^*(u))] \ = 
[x, \delta^*(u)] \ = \ \rho(y,s)\ .\]
This finishes the proof of Claim~3.

Now we can prove part~(ii) of the proposition.
If $(x,t)\in F(l_+)\times K^l$ is of minimal dimension in its equivalence class, then $x$
is non-degenerate and $t\in C_l(K)$, for otherwise $(x,t)$ would be equivalent
to an element of smaller dimension. So $(\Id,\Id,t,x)$ is a reduction datum
for $(x,t)$, and hence $\rho(x,t)=[x,t]$.
Now we let $(y,s)$ be equivalent to $(x,t)$, 
and we let $(\sigma,\delta,u,x')$ be a reduction datum for $(y,s)$.
Then
\[ [x',\delta^*(u)] \ = \ \rho(y,s)\ = \ \rho(x,t)\ = \ [x,t] \ ,\]
where the second equality is Claim~3. So there is a bijective morphism
$\lambda:l_+\to l_+$ such that
\[ (x',\lambda^*(\delta^*(u))) \ = \ (F(\lambda)(x),t)\ .\]
Hence
\[ F(\sigma)(y)\ = \ F(\delta)(x')\ = \ 
F(\delta\lambda)(x)\ , \quad 
s\ = \  \sigma^*(u) \text{\quad and\quad} t =(\delta\lambda)^*(u) \ .\]

(iii)
By the minimality assumptions, the data provided by part~(ii)
must satisfy $m=k=l$, and the morphisms $\sigma$ and $\delta$ 
must both be isomorphisms. So $(F(\alpha)(y),s)=(x,\alpha^*(t))$ 
with $\alpha=\delta^{-1}\sigma:l_+\to l_+$.
The coordinates of $t$ are distinct, by minimality of dimension,
so there is only one permutation with $s=\alpha^*(t)$.
\end{proof}

Now that we understand the equivalence relation generated by \eqref{eq:Gamma relation generator}, we can analyze the topology of the space $F(K)$
when $F$ is a $\bGamma$-space and $K$ is a based topological space.
In particular, we show now that the prolongation of a $\bGamma$-space
with values in the category $\bT$ of compactly generated spaces
to a compactly generated space is automatically compactly generated,
with the weak Hausdorff property being the issue.
The statement and proof of the following proposition
are analogous to the ones in Proposition \ref{prop:geometric realization in bT}.

\begin{prop}\label{prop:prolongation in bT}
  Let $F:\bGamma\to \bT$ be a $\bGamma$-space with values in the
  category of compactly generated spaces, and $K$ a compactly generated based space.
  \begin{enumerate}[\em (i)]
  \item 
    The space $F(K)$ is compactly generated, and hence 
    a coend internal to the category $\bT$,
    of the functor \eqref{eq:Gamma coend}.
  \item  Let $E$ be a $\bGamma$-subspace of $F$ such that
    $E(n_+)$ is closed in $F(n_+)$ for every $n\geq 0$.
    Let $L$ be a closed based subset of $K$. 
    Then the inclusions $E\to F$ and $L\to K$
    induce a closed embedding $\iota:E(L)\to F(K)$.
  \end{enumerate}
\end{prop}
\begin{proof}
  (i) 
  The category $\bT$ is closed under products inside the category $\bK$,
  so the functor \eqref{eq:Gamma coend} takes values in $\bT$;
  the issue is that a priori, the quotient topology need not be weak Hausdorff.
  We use Proposition \ref{prop:closed implies WH}
  to show that the quotient space $F(K)$ is compactly generated.
  We let $E\subset ({\coprod}_{n\geq 0} F(n_+)\times K^n)^2$
  be the equivalence relation generated by \eqref{eq:Gamma relation generator},
  which we had simply denoted `$\sim$' above. We will show that $E$ 
  is closed in the $k$-topology of $({\coprod}_{n\geq 0} F(n_+)\times K^n)^2$.
  Since products distribute over disjoint unions, we may show that
  \[ E_{m,n}\ = \ E\cap (F(m_+)\times K^m\times F(n_+)\times K^n) \]
  is closed in $F(m_+)\times K^m\times F(n_+)\times K^n$ for all $m,n\geq 0$.
  We consider $(y,s,\bar y,\bar s)\in E_{m,n}$, i.e., the pairs $(y,s)$
  and $(\bar y,\bar s)$ are equivalent.
  Proposition \ref{prop:Gamma minimal representative}~(ii) provides
  surjective morphisms
  $\sigma:m_+\to k_+$ and $\bar\sigma:n_+\to \bar k_+$, 
  injective morphisms $\delta:l_+\to k_+$ and $\bar\delta:l_+\to \bar k_+$
  and $u\in K^k$, $\bar u\in K^{\bar k}$  and $(x,t)\in F(l_+)\times K^l$
  such that
  \[ F(\sigma)(y)\ = \ F(\delta)(x)\ , \quad 
  s\ = \  \sigma^*(u) \text{\qquad and\qquad} t =\delta^*(u) \]
  and
  \[ F(\bar\sigma)(\bar y)\ = \ F(\bar \delta)(x)\ , \quad 
  \bar s\ = \  \bar \sigma^*(\bar u) \text{\qquad and\qquad} 
  t =\bar \delta^*(\bar u) \ .\]
  Equivalently, $E_{m,n}$ is the union,
  indexed over $l,\sigma,\bar\sigma,\delta,\bar\delta$ as above,
  of the finite number of sets
  \begin{align*}
  (F(\sigma)\times &K^m\times F(\bar\sigma)\times K^n)^{-1}\\
  &\left(
    (F(\delta)\times\sigma^*\times F(\bar\delta)\times\bar\sigma^*)
    \left((F(l_+)\times\delta^*\times F(l_+)\times\bar\delta^*)^{-1}
      ( \Delta_{F(l_+)\times K^l})\right)\right) \ .
  \end{align*}
  The diagonal $\Delta_{F(l_+)\times K^l}$ is closed in $(F(l_+)\times K^l)^2$
  because $F(l_+)\times K^l$ is compactly generated.
  So its inverse image under the continuous map
  $F(l_+)\times\delta^*\times F(l_+)\times\bar\delta^*$
  is closed in $F(l_+)\times K^k\times F(l_+)\times K^{\bar k}$.
  The map $F(\delta)\times\sigma^*\times F(\bar\delta)\times\bar\sigma^*$
  has a continuous retraction, and is thus a closed embedding by
  Proposition \ref{prop:section is closed embedding}.
  So the set 
  \[
  (F(\delta)\times\sigma^*\times F(\bar\delta)\times\bar\sigma^*)
  \left((F(l_+)\times\delta^*\times F(l_+)\times\bar\delta^*)^{-1}
    ( \Delta_{F(l_+)\times K^l})\right) \]
  is closed in  $F(k_+)\times K^m\times F(\bar k_+)\times K^n$.
  Since $E_{m,n}$ is the inverse image of this latter closed set
  under a continuous map, this show the claim that $E_{m,n}$
  is a closed subset of $F(m_+)\times K^m\times F(n_+)\times K^n$.

  (ii)
  Our first claim is that the map $\iota:E(L)\to F(K)$ is injective.
  We let $(x;\,l_1,\dots,l_m)\in E(m_+)\times L^m$ 
  be a minimal representative of an element of $E(L)$.
  So the element $x$ is non-degenerate, and the $l_i$'s are pairwise distinct
  and different from the basepoint of $L$.
  We claim that $x$ is also non-degenerate when viewed as an element of the
  ambient $\bGamma$-space $F$. To see this, we argue by contradiction and assume
  that there was an injective based map $\delta:k_+\to n_+$ with $k<n$ 
  and an element $y\in F(k_+)$ such that $F(\delta)(y)=x$. 
  We let $\sigma:n_+\to k_+$ be a retraction to $\delta$.
  Then $y=F(\sigma)(F(\delta)(y))=F(\sigma)(x)$. Since $E$ is a $\bGamma$-subspace
  of $F$, we conclude that $y\in E(k_+)$. Then $x=E(\delta)(y)$,
  contradicting the non-degeneracy of $x$.
  Altogether this shows that a representative
  of minimal dimension for a class in $E(L)$ remains a minimal representative
  for its class in $F(K)$. Because minimal representatives are unique
  up to a permutation (Proposition \ref{prop:Gamma minimal representative}~(iii)), 
  this implies that the map $\iota:E(L)\to F(K)$ is injective.

  It remains to show that the continuous injection $\iota$ is a closed map.
  We consider the commutative square
  \[ \xymatrix@C=12mm{
    \coprod_{n\geq 0} E(n_+)\times L^n\ar[r]^-{\text{incl}} \ar[d]_p&
    \coprod_{n\geq 0} F(n_+)\times  K^n\ar[d]^q\\
    E(L)\ar[r]_-{\iota} & F(K) } \]
  where the vertical maps are the quotient maps.
  We consider a point $(y,s)\in F(m_+)\times K^m$ whose class in $F(K)$
  is in the image of $\iota$. 
  Then $(y,s)$ is equivalent to an element $(x,t)\in E(l_+)\times L^l$, 
  which we can take as a minimal representative in its $(E,L)$-equivalence class.
  As we argued in the injectivity statement, $(x,t)$ is then also 
  a minimal representative in its $(F,K)$-equivalence class.

  Proposition \ref{prop:Gamma minimal representative}~(ii) 
  provides a surjective morphism $\sigma:m_+\to k_+$, 
  an injective morphism $\delta:l_+\to k_+$ and $u\in K^k$ such that
  \[ F(\sigma)(y)\ = \ F(\delta)(x)\ , \quad 
  s\ = \  \sigma^*(u) \text{\qquad and\qquad}  t \ =\ \delta^*(u) \ .\]
  So for every subset $A\subset E(L)$, we have
  \begin{align*}
  ( F(m_+&)\times K^m) \cap q^{-1}( \iota(A) ) \ =   \\  {\bigcup}_{\sigma,\delta}\,  
&  (F(\sigma)\times K^m)^{-1}
 \left( (F(\delta)\times \sigma^*) 
 (  (F(l_+)\times\delta^*)^{-1}( (E(l_+)\times L^l)\cap p^{-1}(A) ) ) \right)\ .
  \end{align*}
  The union is over the finite set of pairs $(\sigma:m_+\to k_+,\delta:l_+\to k_+)$
  as above.

  Now we assume that $A$ is closed inside $E(L)$. Because $p$ is continuous,
  $E(l_+)$ is closed in $F(l_+)$ and $L$ is closed in $K$,
  the set $(E(l_+)\times L^l)\cap p^{-1}(A)$ is then closed inside $F(l_+)\times K^l$.
  So
  $(F(l_+)\times\delta^*)^{-1}((E(l_+)\times L^l)\cap p^{-1}(A))$ 
  is a closed subset of $F(l_+)\times K^k$.
  Since $F(\delta)\times\sigma^*$ has a continuous retraction,
  it is a closed embedding by Proposition \ref{prop:section is closed embedding}.
  So the set $ (F(\delta)\times \sigma^*) 
  (  (F(l_+)\times\delta^*)^{-1}((E(l_+)\times L^l)\cap p^{-1}(A) ))$
  is closed in $F(k_+)\times K^m$.
  As the inverse image of a closed set under a continuous map,
  each set in the finite union above is closed inside $F(m_+)\times K^m$.
  We conclude that $( F(m_+)\times K^m) \cap q^{-1}( \iota(A) )$
  is closed in $F(m_+)\times K^m$ for every $m\geq 0$, 
  hence the set $q^{-1}(\iota(A))$ is closed.
  Since $q$ is a quotient map, this shows that $\iota(A)$ is closed in $F(K)$.
\end{proof}

We let $F:\bGamma\to\bK$ be a $\bGamma$-$k$-space and $K$ and $L$ two based $k$-spaces.
The prolongation comes with a continuous, 
based {\em assembly map}\index{subject}{assembly map!of a $\bGamma$-space}
\[ 
  \alpha\ : \ K\sm F(L) \ \to \ F(K\sm L) \ , \
  \alpha(k\sm [x;\, l_1,\dots,l_n])\ = \ [x;\,k\sm l_1,\dots,k\sm l_n]\ . 
 \]
The assembly map is natural in all three variables and associative and unital.
We define a `shifted' 
$\bGamma$-space $F_K=F\circ(K\sm -)$ as the composite
\[ \bGamma \ \xra{K\sm-}\ \bK_* \ \xra{\ F\ }\ \bK_*\ . \]
Then for every based $k$-space $L$, we consider the maps
\[ \phi_n \ : \ F_K(n_+)\times L^n  \ =\ 
F(K\sm n_+)\times L^n  \ \to \ F(K\sm L)\]
defined by
\[ \phi_n\big( [x;k_1\sm i_1,\dots,k_m\sm i_m];\,l_1,\dots,l_n \big) \ = \ 
[x;\, k_1\sm l_{i_1},\dots, k_m\sm l_{i_m}] \ ,\]
where $x\in F(m_+)$, $k_1,\dots,k_m\in K$ and $i_1,\dots,i_m\in\{1,\dots,n\}$.
Another way to say this is that $\phi_n$ is the composite 
\begin{align*}
    F(K\sm n_+)\times L^n  \, \xra{\text{proj}} \, F(K\sm n_+)\sm  L^n  \, 
\xra{\text{ass}}\,  F(K\sm n_+\sm L^n)  \, \xra{F(K\sm\epsilon)}\, F(K\sm L)\ ,
\end{align*}
where
\[ \epsilon\ : \ n_+\sm L^n\ \to \ L \ ,\quad i\sm (l_1,\dots,l_n)\ \longmapsto \ l_i\ .\]
is the evaluation map.
These maps $\phi_n$ are compatible with the equivalence relation defining $F_K(L)$,
so they combine into a continuous map
\[ \phi\ : \ F_K(L)\ \to \ F(K\sm L) \ .  \]
The map $\phi$ is natural in $F$, $K$ and $L$.

\begin{prop}\label{prop:Fubini}
Let $F:\bGamma\to \bK$ be a $\bGamma$-$k$-space and $K$ and $L$ based $k$-spaces. 
Then the natural map $\phi: F_K(L) \to  F(K\sm L)$ is a homeomorphism.
\end{prop}
\begin{proof}
We let $q:K\times L\to K\sm L$ denote the projection.
In the category $\bK$, product with any space preserves proclusions 
(Proposition \ref{prop:proclusion times Z}), so the map
\[ F(m_+)\times q^m\ : \ F(m_+)\times (K\times L)^m  \ \to \ F(m_+)\times (K\sm L)^m  \]
is a proclusion for every $m\geq 0$.
The continuous map
\[ \psi_m \ : \ F(m_+)\times (K\times L)^m  \ \to \ F_K(L)\]
defined by
\[ \psi_m(x;(k_1,l_1),\dots,(k_m,l_m)) \ = \ 
[[x;\, k_1\sm 1,\dots, k_m\sm m];\,l_1,\dots,l_m] \]
is constant on the fibers of $F(m_+)\times q^m$,
so it factors over a continuous map
\[ \bar\psi_m \ : \ F(m_+)\times (K\sm L)^m  \ \to \ F_K(L)\ .\]
These maps are compatible with the equivalence relation defining $F(K\sm L)$,
so they combine into a continuous map
\[ \psi\ : \ F(K\sm L)\ \to \ F_K(L) \ .  \]
Since $\phi$ and $\psi$ are inverse to each other, this completes the proof.
\end{proof}

\begin{construction}
Now we prove an interchange relation between geometric realization and 
prolongation of a $\bGamma$-space. 
We let $A:\bDelta^{\op}\to\bK_*$ be a simplicial based $k$-space and
$F:\bGamma\to\bK$ a $\bGamma$-$k$-space. 
The composite
\[  \bDelta^{\op}\ \xra{\ A\ }\ \bK_* \ \xra{\ F \ } \ \bK_*  \]
is a simplicial space, where the second functor is the prolongation of $F$,
denoted by the same symbol. The composite has a geometric realization 
\[ |F\circ A|\ = \ | [n]\longmapsto \ F(A_n)| \ .\]
We exhibit a natural homeomorphism from $|F\circ A|$ to $F(|A|)$,
the value of $F$ at the realization of $A$.
We let $\kappa_n:A_n\sm\Delta^n_+\to |A|$
be the canonical based map sending a point $(a,t)$ to its equivalence class in $A$.
Then the maps
\[ F(A_n)\times \Delta^n\ \xra{\text{assembly}} \
 F(A_n\sm \Delta^n_+)\ \xra{F(\kappa_n)} \ F(|A|)  \]
are compatible with the equivalence relation defining $|F\circ A|$,
as $[n]$ varies over the objects of $\bDelta$.
So the maps define a continuous map
\[ \kappa \ : \ |F\circ A|\ \to \ F(|A|)\ .\]
In terms of elements, $\kappa$ is thus given by
\[ \kappa[[x;\,a_1,\dots,a_m],t] \ = \ [x;\,[a_1,t],\dots,[a_m,t]]\ ,\]
for $x\in F(m_+)$, $a_1,\dots,a_m\in A_n$ and $t\in\Delta^n$.
\end{construction}

The following proposition is \cite[Lemma 1.9]{woolfson}.

\begin{prop}\label{prop:prolong versus realize}
For every simplicial based $k$-space $A:\bDelta^{\op}\to\bK_*$ and every 
$\bGamma$-space $F:\bGamma\to\bK$,
the map $\kappa:|F\circ A|\to F(|A|)$ is a homeomorphism.
\end{prop}
\begin{proof}
This is a special case of the `Fubini theorem' for coends.
In more detail we consider the functor
\[ \bGamma\times\bGamma^{\op}\times\bDelta^{\op}\times\bDelta\ \to \ \bK \ , \quad
(k_+,l_+,[m],[n])\ \longmapsto \ F(k_+)\times A_m^l\times\Delta^n\ . \]
If we fix objects $[m]$ and $[n]$ of $\bDelta$, 
and exploit that coends commute with product with $\Delta^n$,
we obtain a homeomorphism
\begin{align*}
  \int^{k_+\in\bGamma} F(k_+)\times A_m^k\times\Delta^n\ &\iso \ 
\left(  \int^{k_+\in\bGamma} F(k_+)\times A_m^k\right) \times\Delta^n\ \iso \ 
F(A_m)\times \Delta^n\ . 
\end{align*}
These homeomorphisms are natural as $([m],[n])$ vary in $\bDelta^{\op}\times\bDelta$.
So taking coends over $\bDelta$ gives a homeomorphism
\begin{align*}
  \int^{[m]\in\bDelta} \int^{k_+\in\bGamma} F(k_+)\times A_m^k\times\Delta^m\ &\iso \ 
\int^{[m]\in\bDelta} F(A_m)\times\Delta^m\ = \ |F\circ A|\ . 
\end{align*}
On the other hand, if we fix $k$ and $l$ and exploit that geometric
realization commutes with product with the space $F(k_+)$ and with products 
(see parts~(i) and~(ii) of Proposition \ref{prop:iterated geometric realization}),
we obtain a homeomorphism
\begin{align*}
  \int^{[m]\in\bDelta} F(k_+)\times A_m^l\times\Delta^m\ &\iso \ 
F(k_+)\times \int^{[m]\in\bDelta} A_m^l\times\Delta^m\\ 
&= \ F(k_+)\times |A^l| \ \iso \ F(k_+)\times |A|^l\ . 
\end{align*}
These homeomorphisms are natural as $(k_+,l_+)$ vary in $\bGamma\times\bGamma^{\op}$.
So taking coends over $\bGamma$ gives a homeomorphism
\begin{align*}
\int^{k_+\in\bGamma}  \int^{[m]\in\bDelta} F(k_+)\times A_m^k\times\Delta^m\ &\iso \ 
\int^{k_+\in\bGamma} F(k_+)\times |A|^k\ = \ F(|A|)\ . 
\end{align*}
The `Fubini theorem' for iterated coends 
(the dual to \cite[IX.8, Corollary]{maclane-working})
says that these two ways of iteratively taking coends are canonically isomorphic;
the isomorphism is in fact the map $\kappa$ of the proposition.
\end{proof}

\index{subject}{Gamma-space@$\bGamma$-space!equivariant|(}  
Now we add group actions to the discussion of $\bGamma$-spaces.\index{subject}{Gamma-G-space@$\bGamma$-$G$-space|see{$\bGamma$-space, equivariant}}
By a {\em $\bGamma$-$G$-space}, for a topological group $G$,
we simply mean a reduced functor from $\bGamma$ to the category $G\bT$ of $G$-spaces.
The previous proposition reduces the study of prolonged $\bGamma$-$G$-spaces
on the realizations of $G$-simplicial sets to the study of simplicial $G$-spaces
of the form $F\circ A$.
A sufficient condition for the realization of a simplicial $G$-space 
to be homotopically meaningful is Reedy cofibrancy. This explains why we look
for a practical condition to ensure that simplicial $G$-spaces of the form $F\circ A$
are Reedy cofibrant; the concept of `$G$-cofibrancy' introduced in
Definition \ref{def:G-cofibrant Gamma-G-space} below does the job.

\begin{construction}
We let $\Pc(n)$ denote the power set of $\{1,\dots,n\}$, 
i.e., the set of subsets.
We also write $\Pc(n)$ for the associated poset category, i.e., with object set $\Pc(n)$
and exactly one morphism $U\to T$ whenever $U\subseteq T$.
Given a $\bGamma$-space $F$, we obtain a functor from the
poset category $\Pc(n)$ to based spaces by sending a subset $U$
to $F(U_+)$, with the maps $F(U_+)\to F(T_+)$ induced by the inclusions.
We obtain a {\em latching map}\index{subject}{latching map!of a $\bGamma$-space}
\begin{equation}\label{eq:Gamma latching map} 
 l_n\ : \ \colim_{U\subsetneq \{1,\dots,n\}} \, F(U_+)\ \to \  F(n_+) \ . 
\end{equation}
The symmetric group $\Sigma_n$ acts on the target of the latching map 
by functoriality of $F$. The symmetric group also acts on the source of $l_n$,
by letting $\sigma\in\Sigma_n$ send $F(U_+)$ to $F(\sigma(U)_+)$ via the map
\[ F( (\sigma\cdot -)_+) \ : \ F(U_+)\ \to \ F(\sigma(U)_+)\ .\]
The latching map is equivariant for these two $\Sigma_n$-actions.

We will often consider a $\bGamma$-$G$-space for a topological group $G$.
Then $G$ also acts on source and target of the latching map $l_n$, which
is thus $(G\times\Sigma_n)$-equivariant.
\end{construction}

Now we observe that the latching map
for a $\bGamma$-space $F$ is also a latching map 
for a certain simplicial space, whence the terminology and notation.
We recall that the `simplicial circle' $\mathbf S^1:\bDelta^{\op}\to \bGamma$\index{symbol}{$\mathbf S^1$ - {simplicial circle}}\index{subject}{simplicial circle}  
is given on objects by $\mathbf S^1_n  =  n_+$,
with face maps $d_i:n_+\to(n-1)_+$ given by
\[ d_i(j) \ = \
\begin{cases}
j-1 & \text{ for $i <j$, and}\\  
\quad j & \text{ for $i \geq j$ and $j\ne n$,}\\
\quad 0 & \text{ for $i=j=n$,}
\end{cases}
\]
and degeneracy maps $s_i:n_+\to(n+1)_+$ given by
\[  
 s_i(j) \ = \
\begin{cases}
j+1 & \text{ for $i  <j$, and}\\  
\quad j & \text{ for $i\geq j$.}
\end{cases}\]
The simplicial set $\mathbf S^1$ is isomorphic to $\Delta[1]/\partial\Delta[1]$,
and its realization is homeomorphic to a circle, whence the name.

\begin{prop}\label{prop:Gamma latching closed embedding} 
Let $F$ be a $\bGamma$-space and $n\geq 0$.
\begin{enumerate}[\em (i)]
\item There is a homeomorphism
\[ 
L_n(F\circ \mathbf S^1)\ \xra{\ \iso \ } \ \colim_{U\subsetneq \{1,\dots,n\}} \,  F(U_+)
 \]
whose composite with the latching map \eqref{eq:Gamma latching map} 
is the latching map $l_n:L_n(F\circ \mathbf S^1)\to F(\mathbf S^1_n)=F(n_+)$
of the simplicial space $F\circ\mathbf S^1$.
\item The latching map $l_n$ \eqref{eq:Gamma latching map} is a closed embedding.  
\item
Suppose that $F$ is a $\bGamma$-$G$-space for a topological group $G$.
Then the canonical map
\[ \colim_{U\subsetneq \{1,\dots,n\}} \, ( F(U_+))^G\ \to \ 
\left( \colim_{U\subsetneq \{1,\dots,n\}} \, F(U_+)\right)^G   \]
is a homeomorphism.
\end{enumerate}
\end{prop}
\begin{proof}
(i) We recall that $\bDelta(n)$ is the category with objects the weakly 
monotone surjections $\sigma:[n]\to [k]$;
a morphism from $\sigma$ to $\sigma'$
is a morphism $\alpha:[k]\to [k']$ in $\bDelta$ with $\alpha\circ\sigma=\sigma'$.
Moreover $\bDelta(n)_\circ$ is the full subcategory with all objects except 
the identity of $[n]$.

As we recalled in Remark \ref{rk:Delta(n) vs P(n)},
the category $\bDelta(n)^{\op}$ is isomorphic to the
poset category $\Pc(n)$.
Indeed, an isomorphism is given on objects by
\[ \kappa(\sigma:[n]\to [k]) \ = \ 
\{ i\in\{1,\dots, n\}\ : \ \sigma(i) > \sigma(i-1)\} \ .\]
In the other direction, a subset $U\subset\{1,\dots,n\}$ 
is taken to  the monotone surjection $\sigma_U:[n]\to[|U|]$
defined by
\[ \sigma_U(i)\ = \ |U\cap \{1,\dots,i\}|\ . \]
The structure map $\sigma^*:\mathbf S^1_k\to\mathbf S^1_n=n_+$
is injective with image
\[ \sigma^*(\mathbf S^1_k) \ = \ \kappa(\sigma)\cup\{0\}\ .\]
This means that the maps $\sigma^*:\mathbf S^1_k=k_+\to \sigma^*(k_+)=\kappa(\sigma)_+$ 
define a natural isomorphism
between the two composites in the square of functors:
\[ \xymatrix@C=2mm@R=2mm{ 
\bDelta(n)^{\op}\ar[ddd]_u\ar[rrr]^-\kappa_-\iso 
&&&\Pc(n)\ar[ddd]^{U\mapsto U_+}\\
&&&\\ 
&\ar@{=>}[ur]^{\sigma^*} &&\\
 \bDelta^{\op}\ar[rrr]_-{\mathbf S^1} 
 &&& \bGamma } \]
So as $\sigma$ varies over the objects of $\bDelta(n)_\circ$,
the homeomorphisms
\[ (F\circ\mathbf S^1\circ u)(\sigma) \ = \ F(k_+)\ \xra[\iso]{F(\sigma^*)}\ 
F(\kappa(\sigma)_+)\ \to\
 \colim_{U\subsetneq \{1,\dots,n\}} \,  F(U_+) \]
assemble into the desired isomorphism from $L_n(F\circ\mathbf S^1)$.

(ii) Part~(i) shows that the latching map \eqref{eq:Gamma latching map} 
is an instance of the latching map of a simplicial compactly generated space.
Claim~(ii) is then a special case of 
Proposition \ref{prop:sspaces latching properties}~(iii).

(iii) We contemplate the commutative square:
\[\xymatrix{ 
L_n\left( (F\circ\mathbf S^1)^G\right)\ar[d]_\iso \ar[r]
& (L_n(F\circ\mathbf S^1))^G\ar[d]^\iso \\
 \colim_{U\subsetneq \{1,\dots,n\}} \, ( F(U_+))^G\ar[r] &
\left( \colim_{U\subsetneq \{1,\dots,n\}} \, F(U_+)\right)^G }  \]
The left vertical map is the homeomorphism of part~(i)
for the non-equivariant $\bGamma$-space $(F\circ \mathbf S^1)^G$.
The right vertical map is the effect on $G$-fixed points of the homeomorphism of part~(i).
The upper horizontal map is a homeomorphism by 
Proposition \ref{prop:G-fix preserves pushouts}~(iv),
for the simplicial $G$-space $F\circ\mathbf S^1$.
So the lower horizontal map is a homeomorphism.
\end{proof}

\begin{defn}\label{def:G-cofibrant Gamma-G-space}
Let $G$ be a compact Lie group. A $\bGamma$-$G$-space $F$ is {\em $G$-cofibrant}
if for every $n\geq 1$ the latching map
\[ l_n\ : \ \colim_{U\subsetneq \{1,\dots,n\}} \, F(U_+)\ \to \  F(n_+)  \]
is a $(G\times\Sigma_n)$-cofibration.\index{subject}{Gamma-space@$\bGamma$-space!cofibrant}
\end{defn}

We discuss a relevant example of a cofibrant $\bGamma$-$G$-space 
in Example \ref{eg:ku cofibrant}, namely the $\bGamma$-$G$-space 
of pairwise orthogonal, finite-dimensional subspaces of 
a complete complex $G$-universe. 

\begin{eg}[Equivariant $\bGamma$-simplicial sets]
We let $G$ be a finite group.
In the literature about equivariant $\bGamma$-spaces with values in simplicial sets,
no cofibrancy condition is needed to ensure good homotopical behavior.
One way to explain this is to observe that the geometric realization
of a $\bGamma$-$G$-simplicial set $E:\bGamma\to G\textbf{sset}$ 
is automatically $G$-cofibrant.
Indeed, applying Proposition \ref{prop:Gamma latching closed embedding}~(ii) 
to the $\bGamma$-set in any given simplicial dimension
shows that the latching morphism
 \[ l_n\ : \ \colim_{U\subsetneq \{1,\dots,n\}} \, E(U_+)\ \to \  E(n_+)\  , \]
taken internal to simplicial sets, is dimensionwise injective.
Geometric realization commutes with colimits and takes injective morphisms
of $(G\times\Sigma_n)$-simplicial sets to relative $(G\times\Sigma_n)$-CW-inclusions, 
which are in particular
$(G\times\Sigma_n)$-cofibrations. So the map $|l_n|$ 
is a $(G\times\Sigma_n)$-cofibration without any hypotheses on $E$.  

Our notion of `$G$-cofibrant' should not be confused with 
cofibrancy in the strict model structure that Bousfield and Friedlander
introduce for non-equi\-variant $\bGamma$-simplicial sets 
in \cite[Thm.\,3.5]{bousfield-friedlander},
and generalized to $\bGamma$-$G$-simplicial sets 
by Ostermayr \cite[Thm.\,4.12]{ostermayr}. 
Being `strictly cofibrant' includes the condition that the symmetric group
$\Sigma_n$ acts freely on the complement of the image of the latching map $l_n$;
for our purposes, no such freeness is necessary.
\end{eg}

As we shall now see, the notion of cofibrancy for equivariant $\bGamma$-spaces
is stable under passage to closed subgroups and fixed points by normal subgroups.
If $F$ is a $\bGamma$-$G$-space and $H$ a closed normal subgroup of $G$,
we obtain a $\bGamma$-$G/H$-space $F^H$ by taking $H$-fixed points objectwise.
In other words, $F^H$ is the composite functor
\[ \bGamma  \ \xra{\ F \ } \ G\bT \ \xra{(-)^H}\ (G/H)\bT \ .\]

\begin{prop}\label{prop:fix preserves Gamma-cofibrant}
Let $H$ be a closed subgroup of a compact Lie group $G$ 
and $F$ a $G$-cofibrant $\bGamma$-$G$-space.
\begin{enumerate}[\em (i)]
\item The underlying $\bGamma$-$H$-space of $F$ is $H$-cofibrant.
\item If $H$ is normal, then
the $\bGamma$-$G/H$-space $F^H$ is $(G/H)$-cofibrant.
\end{enumerate}
\end{prop}
\begin{proof}
Part~(i) is clear because restriction from $G\times\Sigma_n$ to $H\times\Sigma_n$
preserves colimits, and it preserves equivariant cofibrations 
by Proposition \ref{prop:cofibrancy preservers}~(i).

(ii) 
The $n$-th latching map for $F^H$ factors as the composite
\[ \colim_{U\subsetneq \{1,\dots,n\}} \, F^H(U_+)\ \to \ 
\left( \colim_{U\subsetneq \{1,\dots,n\}} \, F(U_+)\right)^H \ \xra{\ (l_n)^H} \  F(n_+)^H  \]
of the canonical map and the restriction of the latching map for $F$ 
to $H$-fixed points.
The first map is an isomorphism of $(G/H\times\Sigma_n)$-spaces 
by Proposition \ref{prop:Gamma latching closed embedding}~(iii). 
Since $l_n$ is a $(G\times\Sigma_n)$-cofibration, the
second map $(l_n)^H$ is a $(G/H\times\Sigma_n)$-cofibration by 
Proposition \ref{prop:cofamily pushout property}.  
\end{proof}

\begin{prop}\label{prop:strong stronger cofibrant} 
Let $G$ be a compact Lie group and $F$ a $G$-cofibrant $\bGamma$-$G$-space.
Let $K$ be a finite group and $T$ a finite $K$-set.
Let $\Yc\subset\Pc(T)$ be a $K$-invariant set
of subsets of $T$ that is closed under passage to subsets.
Then the canonical map
\[ \colim_{A\in \Yc} \, F(A_+)\ \to \  F(T_+)\]
is a $(G\times K)$-cofibration.
\end{prop}
\begin{proof}
We prove the following more general statement.
We let $\Yc\subset\Zc\subset\Pc(T)$ be two $K$-invariant sets 
that are both closed under passage to subsets.
We show that then the canonical morphism
\[ \colim_{A\in \Yc} \, F(A_+)\ \to \ \colim_{A\in \Zc}\, F(A_+)\]
is a $(G\times K)$-cofibration. The claim is the special case
$\Zc=\Pc(T)$, which has a terminal object $T$.

We start with the special case where 
$\Zc=\Yc\cup\{ k\cdot B\ : \ k\in K\}$ 
for some subset $B$ of $T$ that does not belong to $\Yc$.
Since $\Zc$ is closed under taking subsets, 
every proper subset of $B$ then belongs to $\Zc$, and hence to $\Yc$.
We let $L\leq K$ be the stabilizer group of $B$, i.e., the subgroup
of those $k\in K$ that map $B$ to itself.
Then the square
\[ \xymatrix@C=15mm{ 
K\times_L \colim_{U\subsetneq B} \, F(U_+) \ar[d]\ar[r]^-{K\times_L l_B} & 
K\times_L F(B_+)\ar[d]\\
\colim_{A\in \Yc} \, F(A_+)\ar[r] & \colim_{A\in \Zc}\, F(A_+)
} \]
is a  pushout of $(G\times K)$-spaces. 
The latching map
\[ l_B \ : \ \colim_{U\subsetneq B} \, F(U_+)\ \to \ F(B_+)\]
is a $(G\times\Sigma_B)$-cofibration by the hypothesis that $F$ is $G$-cofibrant. 
The $L$-action on $B$ specifies a homomorphism $L\to \Sigma_B$, 
so the latching map $l_B$ is a $(G\times L)$-cofibration by
Proposition \ref{prop:cofibrancy preservers}~(i). 
Hence the upper horizontal map in the square is a 
$(G\times K)$-cofibration by Proposition \ref{prop:cofibrancy preservers}~(ii).
Since equivariant cofibrations are stable under cobase change,
the lower horizontal morphism is a $(G\times K)$-cofibration.  

In the general case we choose a chain of intermediate $K$-invariant subsets
\[ \Yc \ = \ \Yc_0 \ \subset\ \Yc_1 \ \subset\ \dots\ \subset \Yc_n \ = \ \Zc \]
such that each $\Yc_i$ is closed under taking subsets and $\Yc_i$ has exactly one
$K$-orbit more than $\Yc_{i-1}$. 
The claim then holds for each pair $(\Yc_i,\Yc_{i-1})$.
Since $(G\times K)$-cofibrations are stable under composition, 
this proves the general case.
\end{proof}

\begin{prop}\label{prop:G-cof Reedy cof}
  Let $G$ be a finite group, $F$ a $G$-cofibrant $\bGamma$-$G$-space,
  and $A$ a simplicial finite based $G$-set.
  \begin{enumerate}[\em (i)]
  \item For all $m,n\geq 0$ the `double latching map'
    \begin{align}  \label{eq:double_latch}
      L_m \left( F\circ(A\sm n_+)\right) \cup_{\colim_{V\subsetneq\{1,\dots,n\}} L_m( F\circ(A\sm V_+))}
      &\left(\colim_{V\subsetneq\{1,\dots,n\}} F(A_m\sm V_+)\right) \nonumber\\ 
      &\to \  F(A_m\sm n_+)  
    \end{align}
    is a $(G\times\Sigma_n)$-cofibration.
  \item The shifted $\bGamma$-$G$-space $F_{|A|}$ is $G$-cofibrant.
  \item  The simplicial $G$-space $F\circ A$ is Reedy $G$-cofibrant.
  \item  For every subgroup $H$ of $G$ the simplicial space $(F\circ A)^H$ 
is Reedy cofibrant.
  \end{enumerate}
\end{prop}
\begin{proof}
(i) 
The double latching space is the colimit of the functor
\[ ( \bDelta(m)^{\op}\times\Pc(n) )^\circ\  
\xra{(\sigma:[m]\to[k],V)\mapsto A_k\sm V_+\ }\ \bGamma\ 
\xra{\ F \ }\ G\bT_* \ ; \]
here $( \bDelta(m)^{\op}\times\Pc(n) )^\circ$ is the `punctured' poset category with 
objects those pairs $(\sigma:[m]\to[k],V)$ 
that are not both maximal, i.e., such that $k<m$ or $V$ is a proper subset
of $\{1,\dots,n\}$.
As we explained in Remark \ref{rk:Delta(n) vs P(n)},
the category $\bDelta(m)^{\op}$ is isomorphic to the poset category $\Pc(m)$,
with $U\subset\{1,\dots,m\}$ corresponding to $\sigma_U:[m]\to[|U|]$
defined as 
\[ \sigma_U(i)\ = \ |U\cap \{1,\dots,i\}|\ . \]
Under this isomorphism of categories, the first functor above
becomes isomorphic to the functor
\begin{equation}  \label{eq:index_for_double_latch}
   ( \Pc(m)\times\Pc(n) )^\circ\  
\xra{(U,V)\mapsto \sigma_U^*(A_k)\sm V_+\ }\ \bGamma\ .
\end{equation}
We rewrite the functor \eqref{eq:index_for_double_latch}.
We set $S=A_m\bs\{\ast\}$, identify $S_+$ with $A_m$, and set
\[ I_U \ = \ \text{Im}(\sigma_U^*: A_{|U|}\to A_m)\bs\{\ast\}\ , \]
which is a $G$-invariant subset of $S$.
We let $\Yc$ denote the subposet of
$\Pc(S\times\{1,\dots,n\})$ consisting of all subsets that are
contained in $I_U\times \{1,\dots,n\}$ 
for some $(U,V)\in (\Pc(m)\times\Pc(n))^\circ$.
Then the functor \eqref{eq:index_for_double_latch} factors as the composite
\[ ( \Pc(m)\times\Pc(n) )^\circ\ \xra{\ \varphi\ } \ 
\Yc \ \xra{\text{incl}}\ \Pc(S\times\{1,\dots,n\})\ \xra{\ (-)_+\ }\ \bGamma\ , \]
where $\varphi(U,V)=I_U\times V$.
Altogether this exhibits the double latching space
as the colimit of the composite 
\[ ( \Pc(m)\times\Pc(n) )^\circ\ \xra{\ \varphi\ } \ 
\Yc \ \xra{\ (-)_+\ }\ \bGamma\ \xra{\ F \ }\ G\bT_* \ . \]

We claim that the poset map $\varphi$ is final, 
i.e., for every subset $B\in\Yc$ the comma category $B\downarrow \varphi$ 
is non-empty and connected.
The comma category $B\downarrow \varphi$ is non-empty by the very definition of $\Yc$.
Now we let $(U,V)$ and $(U',V')$ be elements of the poset 
$( \Pc(m)\times\Pc(n) )^\circ$ such that
\[ B\ \subseteq \ I_U\times V \text{\qquad and\qquad}
 B\ \subseteq \ I_{U'}\times V' \ .\]
We let $\sigma=\sigma_U:[m]\to[k]$ and $\sigma'=\sigma_{U'}:[m]\to[k']$ 
be the corresponding monotone surjections.
There is then a unique pushout in the category $\bDelta$:
\begin{equation} \begin{aligned}\label{eq:Delta po square}
 \xymatrix{ 
[m]\ar[r]^-{\sigma'}\ar[d]_{\sigma} & [k']\ar[d]^{\alpha'} \\
[k] \ar[r]_-{\alpha} & [l]} 
  \end{aligned}\end{equation}
The morphisms $\alpha$ and $\alpha'$ are again surjective.
We set $\tau=\alpha\circ\sigma=\alpha'\circ\sigma'$ and claim that
\begin{equation}  \label{eq:intersect_degeneracies}
  \tau^*(A_l)\ = \ \sigma^*(A_k) \cap (\sigma')^*(A_{k'})\  .
\end{equation}
This is essentially the content of \cite[3.2]{gabriel-zisman},
slightly reformulated; 
we reproduce the argument for the convenience of the reader.
Since $\tau^*(A_l)=\sigma^*(\alpha^*(A_l))$, the set $\tau^*(A_l)$ 
is contained in $\sigma^*(A_k)$,
and similarly for $(\sigma')^*(A_{k'})$.
Conversely, we let $a\in A_m$ be a simplex such that
$a=\sigma^*(x)=(\sigma')^*(y)$ for some $x\in A_k$ and $y\in A_{k'}$.
We write $x=\beta^*(z)$ and $y=(\beta')^*(\bar z)$
for surjective homomorphisms $\beta:[k]\to[p],\beta':[k']\to[p']$ 
and non-degenerate simplices $z$ and $\bar z$. Then
\begin{align*}
(\beta\sigma)^*(z)\ = \  a \ = \  (\beta'\sigma')^*(\bar z)\ .
\end{align*}
By the `Eilenberg-Zilber lemma' (\cite[(8.3)]{eilenberg-zilber-ssc}, 
see also \cite[Sec.\,II.3]{gabriel-zisman}),
the representation of a simplex as a degeneracy of a non-degenerate
element is unique, so $p=p'$, $z=\bar z$ 
and $\beta\sigma=\beta'\sigma'$.
Since the square \eqref{eq:Delta po square} is a pushout, 
there is a unique morphism $\lambda:[l]\to[p]$ 
such that $\lambda\alpha=\beta$ and $\lambda\alpha'=\beta'$.
So
\[ a \ = \ (\beta\sigma)^*(z) 
\ = \ (\lambda\alpha\sigma)^*(z) 
\ = \ (\lambda\tau)^*(z) 
\ = \ \tau^*(\lambda^*(z))\ \in \ \tau^*(A_l)\ . \]
This completes the proof of the relation \eqref{eq:intersect_degeneracies}.

We define
\[  T\ = \ \{ i\in\{1,\dots, m\}\ : \ \tau(i) > \tau(i-1)\} \ ,\]
so that $\tau^*(A_l\bs\{\ast\})=I_T$. We now observe that
\[ B\ \subseteq \ (I_U\times V) \cap (I_{U'}\times V') \ = \ 
 \tau^*(A_l\bs\{\ast\})\times (V\cap V') \ = \
 I_T\times (V\cap V') \ .\]
We have thus found an element $(T,V\cap V')$ 
in the poset $( \Pc(m)\times\Pc(n) )^\circ$ that is less than or equal
to both $(U,V)$ and $(U,V')$, 
and such that $B\subseteq \varphi(T,V\cap V')$.
So the comma category  $B\downarrow \varphi$ is connected.
Since the functor $\varphi$ is final, we can conclude that
the canonical $(G\times\Sigma_n)$-equivariant map
\[ \colim_{ (U,V)\in ( \Pc(m)\times\Pc(n) )^\circ} \, 
F( (I_U\times V)_+)\ \to \ \colim_{B\in\Yc} \,  F(B_+)\]
is a homeomorphism.
On the other hand, the set $\Yc$ is closed under passage to subsets
and $(G\times\Sigma_n)$-invariant inside $\Pc(S\times\{1,\dots,n\})$.
So we can apply Proposition \ref{prop:strong stronger cofibrant} 
with $K=G\times\Sigma_n$ and $T=S\times\{1,\dots,n\}$.
We conclude that the canonical map
\[ \colim_{B\in \Yc} \, F(B_+)\ \to \  F( (S\times\{1,\dots,n\})_+)\ \iso\
 F(A_m\sm n_+)\]
is a $(G\times G\times\Sigma_n)$-cofibration.
When we restrict to the diagonal $G$-action,
the same map is a $(G\times\Sigma_n)$-cofibration by 
Proposition \ref{prop:cofibrancy preservers}~(i). 
Combining these two facts shows that the 
double latching morphism \eqref{eq:double_latch} is a $(G\times\Sigma_n)$-cofibration. 

(ii) Part~(i) says that the morphism of simplicial $(G\times\Sigma_n)$-spaces 
\[  \colim_{V\subsetneq\{1,\dots,n\}} F\circ (A\sm V_+)  \ \to \ 
F\circ (A\sm n_+) \]
is a Reedy $(G\times\Sigma_n)$-cofibration.
Geometric realization is a left Quillen functor, 
so it takes Reedy $(G\times\Sigma_n)$-cofibrations
of simplicial $(G\times\Sigma_n)$-spaces to $(G\times\Sigma_n)$-cofibrations.
So the map 
\[ | \colim_{V\subsetneq\{1,\dots,n\}} F\circ (A\sm V_+) | \ \to \ 
|F\circ (A\sm n_+)| \]
is a $(G\times\Sigma_n)$-cofibration.
Realization commutes with colimits, so the source is a colimit, 
over proper subsets of $\{1,\dots,n\}$, of the spaces $|F\circ(A\sm V_+)|$. 
The Fubini isomorphism of Proposition \ref{prop:prolong versus realize}
identifies the $G$-space $|F\circ (A\sm V_+) |$ 
with $F(|A\sm V_+|)\iso F(|A|\sm V_+)=F_{|A|}(V_+)$.  
This proves the claim.

Part~(iii) is the special case of part~(i) for $n=1$.

(iv) The restriction functor from $G$-spaces to $H$-spaces
preserves colimits, and hence latching objects, 
and takes $G$-cofibrations to $H$-cofibrations
by Proposition \ref{prop:cofibrancy preservers}~(i); so the latching map
$l_n:L_n(F\circ A)\to F(A_n)$ is an $H$-cofibration by part~(iii).
The $H$-fixed point map
\[ (l_n)^H \ :\ ( L_n(F\circ A))^H\ \to \ (F(A_n))^H \]
is then a non-equivariant cofibration by 
Proposition \ref{prop:cofamily pushout property}.
Taking $H$-fixed points commutes with taking latching objects, 
by Proposition \ref{prop:G-fix preserves pushouts}~(iv).
So the $n$-th latching map for the simplicial space $(F\circ A)^H$
is a cofibration.
\end{proof}

The following proposition provides a way to reduce certain questions
about $\bGamma$-$G$-spaces for compact Lie groups to the special case
of finite groups.

\begin{prop}\label{prop:Gamma fixed points} 
  Let $G$ be a connected compact Lie group and $F$ a $\bGamma$-$G$-space.
  Then for every based $G$-space $K$ the map
    \[ (F^G)(K^G)\ \to \ (F(K))^G \]
    induced by the fixed point inclusions $F^G\to F$ and $K^G\to K$ is a homeomorphism.
\end{prop}
\begin{proof}
  Since $F(n_+)^G$ is closed inside $F(n_+)$ and $K^G$ is a closed subset of $K$, 
  Proposition \ref{prop:prolongation in bT}~(ii)
  shows that the inclusions induce a closed embedding $(F^G)(K^G)\to F(K)$.
  The image of this map is contained in $F(K)^G$, 
  so it only remains to show that every $G$-fixed point of $F(K)$
  is the image of a point in $(F^G)(K^G)$.
  We consider a point of $F(K)$ represented by
  a tuple $(x;\,k_1,\dots,k_n)$ in $F(n_+)\times K^n$.
  We assume that the number $n$ has been chosen minimally, so that
  $x$ is non-degenerate and the entries $k_i$ are pairwise distinct 
  and different from the basepoint of $K$.
  If the point  $[x;\,k_1,\dots,k_n]$ of $F(K)$ is $G$-fixed, 
  then for every group element $g$
  the tuple $(g x;\,g k_1,\dots,g k_n)$ is equivalent to the original tuple.
  By Proposition \ref{prop:Gamma minimal representative}~(iii)
  there is a unique permutation $\sigma(g)\in \Sigma_n$ such that
  \[ (g x;\, g k_1,\dots,g k_n)\ = \ 
  (F(\sigma(g)^{-1})(x);\, k_{\sigma(g)(1)},\dots,k_{\sigma(g)(n)} )\ . \]
  The map $\sigma:G\to \Sigma_n$ is a homomorphism,
  and it is continuous since $G$ acts continuously on $K$.
  Since $G$ is connected, the homomorphism must be trivial, i.e.,
  $\sigma(g)=1$ for all $g\in G$. Thus the points $x$ 
  and $k_1,\dots,k_n$ are all $G$-fixed.
\end{proof}

For us, the main purpose of equivariant $\bGamma$-spaces is to 
construct equivariant spectra by evaluation on spheres.
In more detail, we let $G$ be a compact Lie group
and $F$ a $\bGamma$-$G$-space.
We define an orthogonal $G$-spectrum $F(\mS)$ by
  \[  F(\mS)(V) \ = \ F(S^V)  \ . \]
The structure map $\sigma_{V,W}:S^V\sm F(\mS)(W)\to F(\mS)(V\oplus W)$
is the assembly map for $K=S^V$ and $L=S^W$,
followed by the effect of $F$ on the canonical homeomorphism 
$S^V\sm S^W\iso S^{V\oplus W}$.
The $O(V)$-action on $F(\mS)(V)$ is via the action on $S^V$ and the continuous
functoriality of $F$.
If $V$ is a $G$-representation, then $F(S^V)$ has the diagonal $G$-action,
from its actions on $F$ and on $V$.\index{subject}{Gamma-space@$\bGamma$-space!evaluation on spheres}

\begin{prop}\label{prop:Gamma on S^V}
  Let $G$ be a compact Lie group and $F$ a $G$-cofibrant $\bGamma$-$G$-space.
  \begin{enumerate}[\em (i)]
  \item 
  For every $G$-representation $V$,
  the fixed point space $F(S^V)^G$ is $(\dim(V^G)-1)$-connected.
\item  The orthogonal $G$-spectrum $F(\mS)$ is equivariantly connective.
  \end{enumerate}
\end{prop}
\begin{proof}
(i) 
  We start with the special case where $G$ is a finite group.
  We set $d=\dim(V^G)$.
  A choice of linear isometry $V^G\iso\mR^d$ induces a $G$-equivariant homeomorphism
  between $F(S^V)$ and $F(S^d\sm S^{V^\perp})$, 
  where $V^\perp$ is the orthogonal complement of $V^G$ inside $V$. 
  So $F(S^V)^G$ is homeomorphic to $F(S^d\sm S^{V^\perp})^G$.

  The sphere $S^d$ is homeomorphic to the geometric realization of the 
  based simplicial set $\Delta[d]/\partial\Delta[d]$, the represented simplicial set
  $\Delta[d]=\bDelta(-,[d])$ modulo its simplicial boundary. 
  The representation sphere $S^{V^\perp}$ is homeomorphic to the geometric realization 
  of a finite based $G$-simplicial set. Indeed, $S^{V^\perp}$ 
  admits the structure of a smooth $G$-manifold; 
  Illman's triangulation theorem \cite{illman-finite G} then provides
  a $G$-equivariant triangulation. Passing to the barycentric subdivision 
  provides a finite $G$-simplicial set $B$ that realizes to $S^{V^\perp}$.
  Altogether, $S^V$ is $G$-homeomorphic to the geometric realization
  of the $G$-simplicial set $A=(\Delta[d]/\partial\Delta[d])\sm B$.

  Proposition \ref{prop:prolong versus realize}
  provides a $G$-equivariant homeomorphism
  \[  F(S^V)\ \iso \ F(|A|) \ \iso \ |F\circ A | \ .\]
  Taking $G$-fixed points commutes with realization 
  by Proposition \ref{prop:G-fix preserves pushouts}~(iv),
  so $F(S^V)^G$ is homeomorphic to the realization of the simplicial space
  $(F\circ A)^G$ that takes $[m]\in\bDelta^{\op}$
  to the space $F(A_m)^G$.
  In dimensions below $d$, the simplicial set $\Delta[d]/\partial\Delta[d]$,
  and hence also the simplicial set $A=(\Delta[d]/\partial\Delta[d])\sm B$,
  consists only of the base point. So for $m<d$, the space $F( A_m)^G$ is a single point.
  Moreover, the simplicial space $(F\circ A)^G$ 
  is Reedy cofibrant by Proposition \ref{prop:G-cof Reedy cof}~(iv).
  So its realization is $(d-1)$-connected by Proposition \ref{prop:real Sing X}~(ii).
  Altogether this establishes the claim that the space $F(S^V)^G$ 
  is $(d-1)$-connected.

  Now we treat the case of a general compact Lie group.
  We let $G^\circ$ denote the connected component of the identity and
  $\bar G= G/G^\circ$ the finite group of components of $G$.
  The $\bGamma$-$\bar G$-space $F^{G^\circ}$ is $\bar G$-cofibrant 
  by Proposition \ref{prop:fix preserves Gamma-cofibrant}~(ii).
  Proposition \ref{prop:Gamma fixed points}  provides a homeomorphism
  \[ F(S^V)^G \ = \ \big( F(S^V)^{G^\circ} \big)^{\bar G} \ \iso \ 
  \big( (F^{G^\circ})(S^{V^{G^\circ}})\big)^{\bar G}\ . \]
  Because $(V^{G^\circ})^{\bar G}=V^G$, the right hand side is $(\dim(V^G)-1)$-connected 
  by the special case above, for the finite group $\bar G$, 
  the $\bGamma$-$\bar G$-space $F^{G^\circ}$
  and the $\bar G$-representation $V^{G^\circ}$.  

  (ii)
  We let $H$ be any closed subgroup of $G$ and show that the group $\pi_{-k}^H(F(\mS))$
  is trivial for all $k\geq 1$.
  We let $V$ be any $H$-representation, $K$ a closed subgroup of $H$, 
  and we set $d_K=\dim(V^K)$.
  The underlying $\bGamma$-$K$-space of $F$ is $K$-cofibrant
  by Proposition \ref{prop:fix preserves Gamma-cofibrant}~(i).
  So $F(S^{\mR^k\oplus V})^K$ is $(k+d_K-1)$-connected by part~(i).

  On the other hand, the cellular dimension of $S^V$ at $K$,
  in the sense of \cite[II.2, p.\,106]{tomDieck-transformation},
  is at most $d_K$. 
  Because $k$ is positive, the cellular dimension of $S^V$ at $K$ does not exceed
  the connectivity of $F(S^{\mR^k\oplus V})^K$.
  So every based continuous $H$-map $S^V\to F(S^{\mR^k\oplus V})$ is equivariantly
  null-homotopic by \cite[II Prop.\,2.7]{tomDieck-transformation},
  and the set $[S^V,F(S^{\mR^k\oplus V})]^H$ has only one element.
  Passage to the colimit over $V\in s(\Uc_H)$ proves the claim.  
\end{proof}

We can also show that prolongation of $G$-cofibrant $\bGamma$-$G$-spaces
is homotopical in the $\bGamma$-space variable, as long as
we evaluate on finite based $G$-CW-complexes.

\begin{defn}
Let $G$ be a compact Lie group. 
A morphism $\psi:E\to F$ of $\bGamma$-$G$-spaces
is a {\em strict equivalence}\index{subject}{strict equivalence!of $\bGamma$-$G$-spaces}
if for every $n\geq 1$ the map of $(G\times \Sigma_n)$-spaces 
$\psi(n_+):E(n_+)\to F(n_+)$ is an
$\Fc(G;\Sigma_n)$-weak equivalence, where $\Fc(G;\Sigma_n)$
is the family of graph subgroups of $G\times\Sigma_n$.\index{subject}{graph subgroup}
\end{defn}

\begin{rk}\label{rk:alternative strict level} 
A strict equivalence $\psi:E\to F$ of $\bGamma$-$G$-spaces
also satisfies the following condition:
for every closed subgroup $H$ of $G$ and every finite $H$-set $S$, the map
$\psi(S_+):E(S_+)\to F(S_+)$ is an $H$-weak equivalence with respect
to the diagonal $H$-actions. 
Indeed, if $S$ has $n$ elements, we may suppose that $S=\{1,\dots,n\}$
with $H$-action specified by a continuous homomorphism
$\rho:H\to\Sigma_n$.
We let $K$ be a closed subgroup of $H$ and $\Gamma$
the graph of $\rho|_K:K\to\Sigma_n$.
The map
$\psi(S_+)^\Gamma:E(S_+)^\Gamma:\to F(S_+)^\Gamma$
is a weak equivalence by the hypothesis on $\psi$.
Moreover,
\[ F(S_+)^K \ = \  F(S_+)^\Gamma\ , \]
where the fixed points on the left hand side are with respect to the diagonal $K$-action.
This proves the claim.
\end{rk}

The following result will allow us to extend equivariant homotopical properties
from the class of finite $G$-simplicial sets to the class of finite $G$-CW-complexes.
The non-equivariant case is treated for example in \cite[Thm.\,2C.5]{hatcher}.
I am sure that also the equivariant result is well-known; 
however, I don't know a reference, so I provide a proof.

\begin{prop}\label{prop:equivariant simplicial realization}
Let $G$ be a finite group.
\begin{enumerate}[\em (i)]
\item 
Let $X$ be a finite $G$-CW-complex, $A$ a finite $G$-simplicial set
and $f:|A|\to X$ a continuous $G$-map. 
Then there is a finite $G$-simplicial set $B$, a monomorphism of
$G$-simplicial sets $i:A\to B$ and a $G$-homotopy equivalence 
$h:|B|\to X$ such that
\[ h\circ |i| \ = \ f \ : \ |A|\ \to \ X\ .\]
\item
Every finite based $G$-CW-complex is based $G$-homotopy equivalent
to the realization of a finite based $G$-simplicial set.
\end{enumerate}
\end{prop}
\begin{proof}
(i)
We let $c Y= (Y\times[0,1]) / (Y\times\{1\})$ denote 
the unreduced cone of a space $Y$.
We start with a very special case, namely when there is 
a pushout square of $G$-spaces
\[ \xymatrix{ 
G/H\times|\partial\Delta[k]|\ar[r]^-{\text{incl}}\ar[d]_{\alpha} &
G/H\times c|\partial\Delta[k]| \ar[d]\\
|A|\ar[r]_-f & X} \]
for some subgroup $H$ of $G$ and some $k\geq 0$;
in other words, we suppose that $X$ is obtained from $|A|$ 
by attaching one equivariant cell. 

The continuous map
\[ \alpha(e H,-)\ : \ |\partial\Delta[k]|\ \to\ |A| \]
lands in the $H$-fixed points.
Since fixed points commute with realization
(Proposition \ref{prop:G-fix preserves pushouts}~(iv)), 
we may view it as a continuous map to $|A^H|$.
Now we use a `simplicial approximation', by which we mean the following data: 
\begin{itemize}
  \item a finite simplicial set $D$, 
  \item a morphism of simplicial sets
    $\Phi : D\to A^H$, and
  \item a homotopy equivalence $\varphi :|D|\to |\partial\Delta[k]|$
    such that $\alpha(e H,-)\circ \varphi$ is homotopic to the realization of $\Phi$.
\end{itemize}
For example, we can use  $D=\text{Sd}^m(\partial\Delta[k])$ for a suitably
large $m\geq 0$, where Sd is Kan's subdivision functor \cite[Sec.\,7]{kan-on css}; 
the remaining data is then provided 
by Lemma~7.5 and Theorem~8.5 of \cite{kan-on css}.
Alternatively, we can take $D=SD^m(\partial\Delta[k])$ for a suitably
large $m\geq 0$, where $SD$ is the `double simplicial subdivision'
of \cite[Def.\,(12.5)]{curtis}; the other data is provided by \cite[Thm.\,(12.7)]{curtis}.

We let $\tilde\Phi:G/H\times D\to A$ be the $G$-equivariant extension of $\Phi$,
i.e., $\tilde\Phi_n(g H,x)=g\cdot \Phi(x)$ for $g\in G$ and $x\in D_n$.
The geometric realization of $\tilde\Phi$ 
is then $G$-equivariantly homotopic to the composite
\[ G/H\times |D| \ \xra{G/H\times \varphi}\ 
G/H\times |\partial\Delta[k]|
\ \xra{\ \alpha\ } \ |A|\ .\]
We choose a $G$-equivariant homotopy
\[ K \ : \ G/H\times|D|\times[0,1] \ \to \ |A| \]
from $|\tilde\Phi|$ to the map $\alpha\circ(G/H\times \varphi)$.
We absorb the homotopy into the mapping cylinder of the map 
$|\tilde\Phi|:G/H\times|D|\to|A|$, and obtain a commutative diagram of $G$-spaces:
\[ \xymatrix{
  |A|\cup_{ |\tilde\Phi|} (G/H\times|D|\times[0,1])
 \ar[d]_{\Id_{|A|}\cup K} &  G/H\times |D| \ar[r]^-{\text{incl}} 
\ar[l]_-{(-,1)} \ar[d]^{G/H\times  \varphi} & 
G/H\times c| D|\ar[d]^{G/H\times \bar \varphi}\\
  |A| & G/H\times |\partial\Delta[k]| \ar[r]_-{\text{incl} } 
\ar[l]^-\alpha & G/H\times c|\partial\Delta[k]|
} \]
All vertical maps in the diagram are $G$-homotopy equivalences,
and the right horizontal maps are $G$-cofibrations.
The gluing lemma (Proposition \ref{prop:gluing lemma G-spaces})
then shows that the induced map on pushouts
\begin{align}  \label{eq:map_on_pushouts}
 |A|\cup_{|\tilde\Phi|} (G/H\times |D|\times [0,1]) &\cup_{G/H\times |D|} (G/H\times c|D|)\\ 
&\to \ 
|A|\cup_{\alpha} (G/H\times c|\partial\Delta[k]|) \ \iso \ X  \nonumber
\end{align}
is a $G$-weak equivalence;
moreover, its restriction to $|A|$ is the original map $f:|A|\to X$.
The source of \eqref{eq:map_on_pushouts}
 is homeomorphic to the unreduced mapping cone of the map
$|\tilde\Phi|$.
Mapping cones can also be formed in the category of $G$-simplicial sets,
so the source of \eqref{eq:map_on_pushouts} is equivariantly homeomorphic
to the realization of the unreduced mapping cone 
of $\tilde\Phi:G/H\times D\to A$.
This simplicial mapping cone is the desired $G$-simplicial set $B$.

The rest of the argument is now straightforward.
Induction over the number of relative equivariant cells
proves the case where $(X,|A|)$ is a finite relative $G$-CW-pair 
and where $f:|A|\to X$ is the inclusion.
In the general case, the filtration of $A$ by simplicial skeleta induces a filtration
on the geometric realization that gives $|A|$ the structure of a 
finite $G$-CW-complex.
If the map $f$ is cellular for this structure,
then the mapping cylinder $Z=|A|\times [0,1]\cup_f X$
inherits a finite $G$-CW-structure in which $|A|$ is an equivariant subcomplex. 
So the previous case provides a finite $G$-simplicial set $B$, a monomorphism of
$G$-simplicial sets $i:A\to B$ and a $G$-homotopy equivalence 
$h:|B|\to Z$ such that the composite
\[ h\circ |i| \ : \ |A|\ \to \ Z \]
is the `front inclusion' $(-,0):|A|\to Z$.
The projection $Z\to X$ is a $G$-equivariant homotopy equivalence
and $f=p\circ(-,0):|A|\to X$.
So the triple $(B,i,p\circ h)$ has the desired properties.
If the map $f$ is arbitrary, we use the equivariant cellular approximation theorem 
(see for example \cite[II Prop.\,5.6]{bredon-cohomology}, \cite[Thm.\,4.4]{matumoto}
or \cite[Ch.\,II, Thm.\,2.1]{tomDieck-transformation})
and the equivariant homotopy extension property of the $G$-map $|i|:|A|\to|B|$
to reduce to the cellular case.

(ii) We let $X$ be a finite based $G$-CW-complex. 
We apply part~(i) to the terminal simplicial $G$-set $A=\ast$
and the inclusion of the basepoint $|\ast|\to X$; 
part~(i) provides a finite based $G$-simplicial set $B$ 
and a based $G$-homotopy equivalence $h:|B|\to X$.
\end{proof}

Now we have all tools ready to show that prolongation
preserves strict equivalences between cofibrant $\bGamma$-$G$-spaces,
at least on finite $G$-CW-complexes.

\begin{prop}\label{prop:cofibrant Gamma-G on G-CW}
Let $G$ be a finite group and $\psi:E\to F$ a strict equivalence between   
cofibrant $\bGamma$-$G$-spaces.
Then for every finite based $G$-CW-complex $X$,
the map $\psi(X)^G:E(X)^G\to F(X)^G$ is a weak equivalence.
\end{prop}
\begin{proof}
We start with the special case where $X=|B|$ is 
the geometric realization of a finite based $G$-simplicial set $B$.
In this situation, the morphism 
of simplicial spaces $(\psi\circ B)^G:(E\circ B)^G\to (F\circ B)^G$
is levelwise a weak equivalence by Remark \ref{rk:alternative strict level}.
Moreover, source and target are Reedy cofibrant 
by Proposition \ref{prop:G-cof Reedy cof}~(iv).
Geometric realization takes levelwise weak equivalences
between Reedy cofibrant simplicial spaces to weak equivalences
(Proposition \ref{prop:realization invariance in Top}).
This proves the special case.

In the general case we choose a finite based $G$-simplicial set $B$ 
and a based $G$-homotopy equivalence $h:|B|\to X$,
as provided by Proposition \ref{prop:equivariant simplicial realization}~(ii).
Prolonged $\bGamma$-spaces are continuous functors,
so they preserve equivariant based  homotopies.
In the commutative square
\[ \xymatrix@C=12mm{ 
E(|B|)^G \ar[r]^-{\psi(|B|)^G} \ar[d]_{E(h)^G} &
F(|B|)^G \ar[d]^{F(h)^G} \\
E(X)^G \ar[r]_-{\psi(X)^G} & F(X)^G } \]
both vertical maps are thus homotopy equivalences.
The upper map is a weak equivalence by the previous paragraph,
hence so is the lower map.
\end{proof}

Now we move on to the analysis of special and very special $\bGamma$-$G$-spaces.
We will eventually assume that the $\bGamma$-$G$-space is $G$-cofibrant,
in order to have homotopical control over its prolongation.
The final aim is to show that for finite $G$ and very special (respectively special) $F$,
the evaluation on spheres $F(\mS)$ is a $G$-$\Omega$-spectrum
(respectively `positive' $G$-$\Omega$-spectrum),
see Theorem \ref{thm:prolonged delooping} respectively
Theorem \ref{thm:special cofibrant Gamma positive Omega} below.

If $F$ is any $\bGamma$-space and $S$ a finite set, then we define the map
\[  P_S \ : \  F(S_+)\ \to \ \map(S, F(1_+))   \]
by $P_S(x)(s)= F(p_s)(x)$, where $p_s:S_+\to 1_+$ sends $s$ to $1$ and all other elements
of $S_+$ to the basepoint. If a group $G$ acts on $F$ and $S$,
then the map $P_S$ is $G$-equivariant for the diagonal $G$-action on the source
and the conjugation action on the target.

\begin{defn}\label{def:special G-Gamma}
Let $G$ be a compact Lie group. 
A $\bGamma$-$G$-space $F$ is {\em special}\index{subject}{Gamma-space@$\bGamma$-space!special}
if for every closed subgroup $H$ of $G$ and every finite $H$-set $S$ the map
\[ (P_S)^H \ : \ F(S_+)^H\ \to \ \map^H(S,F(1_+)) \]
is a weak equivalence. 
\end{defn}

We showed in Proposition \ref{prop:fix preserves Gamma-cofibrant}
that the notion of cofibrancy for equivariant $\bGamma$-spaces
is stable under passage to closed subgroups and fixed points by normal subgroups.
Now we show the analogous statement for specialness.

\begin{prop}\label{prop:fix preserves special}
Let $H$ be a closed subgroup of a compact Lie group $G$ 
and $F$ a special $\bGamma$-$G$-space.
\begin{enumerate}[\em (i)]
\item The underlying $\bGamma$-$H$-space of $F$ is special.
\item If $H$ is normal, then the $\bGamma$-$G/H$-space $F^H$ is special.
\end{enumerate}
\end{prop}
\begin{proof}
Part~(i) is clear by definition.
For part~(ii) we consider a closed subgroup of $G/H$, which must be of
the form $\Delta/H$ for a closed subgroup $\Delta$ of $G$ with $H\leq \Delta$.
We let $S$ be a finite $\Delta/H$-set,
which we can also view as a finite $\Delta$-set by restriction along
the projection $\Delta\to\Delta/H$.
The claim then follows from the hypothesis that $F$ is special and
the relations
\[ ( F^H(S) )^{\Delta/H} \ = \ F(S)^\Delta\text{\quad and\quad}
\map^{\Delta/H}(S, F^H(1_+))\ = \ \map^\Delta(S, F(1_+)) \ .\qedhere\]
\end{proof}

The functor obtained by prolonging a $\bGamma$-$G$-space 
comes with `Wirthm{\"u}ller type' maps, defined as follows.
We let $H$ be a closed subgroup of finite index in $G$ and $Z$ a based $H$-space.
We recall that $l_Z:G\ltimes_H Z\to Z$ denotes the $H$-equivariant
projection to the wedge summand indexed by the preferred coset $e H$, i.e.,
\[ l_Z[g,z]\ = \
\begin{cases}
  g\cdot z & \text{ if $g\in H$, and}\\
\ \ast & \text{else.} 
\end{cases}\]
The $H$-equivariant map $F(l_Z): F(G\ltimes_H Z) \to F(Z)$
is then adjoint to a $G$-equivariant Wirthm{\"u}ller map\index{subject}{Wirthm{\"u}ller map!of an equivariant $\bGamma$-space}
\begin{equation}  \label{eq:define_omega_A}
  \omega_Z\ : \ F(G\ltimes_H Z) \ \to \ \map^H(G,F(Z)) \ .  
\end{equation}

We recall some equivalent characterizations of special $\bGamma$-$G$-spaces.
The equivalence of the first two conditions 
was first noted in \cite[Cor.]{shimakawa-note}. 

\begin{prop}\label{prop:special characterizations}
  Let $G$ be a compact Lie group and $F$ a $\bGamma$-$G$-space. 
  Then the following conditions are equivalent.
  \begin{enumerate}[\em (i)]
  \item The $\bGamma$-$G$-space $F$ is special. 
  \item For every $n\geq 1$ the map
    \[ P_n \ : \ F(n_+)\ \to \ \map(\{1,\dots,n\},F(1_+)) \ = \ F(1_+)^n\]
    is an $\Fc(G;\Sigma_n)$-equivalence, where
    $\Fc(G;\Sigma_n)$ is the family of graph subgroups.\index{subject}{graph subgroup}
  \item For all pairs of finite $G$-sets $T$ and $U$, the map
    \[  (F(p_T),F(p_U))\ : \ F((T\amalg U)_+)\ \to \ F(T_+)\times F(U_+)\]
    is a $G$-weak equivalence, and for every finite index subgroup $H$ of $G$
    and every finite $H$-set $S$, the Wirthm{\"u}ller map
    $\omega_{S_+}: F( (G\times_H S)_+)\to \map^H(G,F(S_+))$
    is a $G$-weak equivalence. 
  \end{enumerate}
\end{prop}
\begin{proof}
  (i)$\Longleftrightarrow$ (ii)
  We let $H$ be a closed subgroup of $G$.
  Every finite $H$-set is isomorphic to $\alpha^*\{1,\dots,n\}$
  for some continuous homomorphism $\alpha:H\to \Sigma_n$.
  If $\Gamma\leq G\times\Sigma_n$ denotes the graph of $\alpha$,
  then $F(n_+)^\Gamma=F(\alpha^*\{1,\dots,n\})^H$ and
  $(F(1_+)^n)^\Gamma= \map^H(\alpha^*\{1,\dots,n\},F(1_+))$.
  So the map $(P_n)^\Gamma$ is a weak equivalence if and only if
  the map $(P_{\alpha^*\{1,\dots,n\}})^H$  is a weak equivalence.

  (i)$\Longrightarrow$ (iii) 
  In the commutative square
  \begin{equation} \begin{aligned}\label{eq:summand_reduction}
 \xymatrix{
    F( (T\amalg U)_+)\ar[r]^-{(F(p_T),F(p_U))}\ar[d]_{P_{T\amalg U}}^\simeq & 
    F(T_+)\times F(U_+) \ar[d]^{p_T\times p_U}_\simeq\\
    \map( T\amalg U,F(1_+))\ar[r]_-\iso & \map(T,F(1_+))\times \map(U,F(1_+)) }       
    \end{aligned}  \end{equation}
  both vertical maps are $G$-weak equivalences because $F$ is special,
  and the lower horizontal map is a homeomorphism.
  So the upper horizontal map is a $G$-weak equivalence.

  Now we suppose that $H$ has finite index in $G$ and we let $S$ be a finite $H$-set.
  Then the Wirthm{\"u}ller map participates in the commutative diagram:
\[ \xymatrix{ 
  F( (G\times_H S)_+)\ar[r]^-{\omega_{S_+}}\ar[d]_{P_{G\times_H S}} &
  \map^H(G,F(S_+)) \ar[d]^{\map^H(G,P_S)} \\ 
  \map( G\times_H S, F(1_+)) \ar[r]_-\iso & \map^H(G,\map(S,F(1_+)))
} \]
The left and right vertical maps are $G$-weak equivalences by specialness, 
applied to the $G$-set $G\times_H S$ respectively the $H$-set $S$.
The lower horizontal homeomorphism sends $f:G\times_H S\to F(1_+)$
to the adjoint $H$-map $\hat f:G\to \map(S,F(1_+))$ 
with $\hat f(g)(s)=f[g,s]$.
Altogether this shows that the Wirthm{\"u}ller map $\omega_{S_+}$
is a $G$-weak equivalence.

(iii)$\Longrightarrow$ (i) 
We start by showing that for every finite $G$-set $S$ the map
$ P_S : F(S_+)\to \map(S,F(1_+))$ is a $G$-weak equivalence. 
If $S=T\amalg U$ is the disjoint union of two finite $G$-sets $T$ and $U$,
we contemplate the commutative square \eqref{eq:summand_reduction}.
The upper horizontal map is a $G$-weak equivalence by hypothesis.
The lower horizontal map is an equivariant homeomorphism.
So if the claim holds for the $G$-sets $T$ and $U$, then it holds for
their disjoint union. This reduces the claim to the case of transitive $G$-sets.
For every finite index subgroup $H$ of $G$, there is a commutative square
of $G$-maps
\[ \xymatrix{ 
F( (G/H)_+)\ar[r]^-{P_{G/H}}\ar[d]_\iso & \map(G/H,F(1_+)) \ar[d]^\iso\\
F(G\ltimes_H 1_+)\ar[r]_-{\omega_{1_+}} &\map^H(G,F(1_+))
} \]
where the right vertical map is adjoint to the $H$-map
\[ \map(G/H,F(1_+)) \ \to \ F(1_+)\ ,\quad f \ \longmapsto \ f(e H) \ . \]
The two vertical maps are homeomorphisms.
The lower horizontal map is the Wirthm{\"u}ller map 
for the $H$-set $\{1\}$;
this map is a $G$-weak equivalence by hypothesis.
So the map $P_S$ is a $G$-weak equivalence for transitive $G$-sets;
this completes the proof of the first claim.

Now we let $H$ be a finite index subgroup of $G$ and $S$ a finite $H$-set.
The $H$-maps $[1,-]:S_+\to (G\times_H S)_+$ and $l_S:(G\times_H S)_+\to S_+$
express $S_+$ as an $H$-equivariant retract of the based $G$-set $(G\times_H S)_+$.
The middle vertical map in the commutative diagram
\[ \xymatrix@C=12mm{ 
F( S_+)^H\ar[d]_{ (P_S)^H}\ar[r]^-{F([1,-])^H} & 
F( (G\times_H S)_+)^H\ar[r]^-{F(l_S)^H} \ar[d]_{(P_{G\times_H S})^H} &
F(S_+)^H\ar[d]^{ (P_S)^H} \\
\map^H(S,F(1_+)) \ar[r]_-{(l_S)^*} &
\map^H(G\times_H S,F(1_+)) \ar[r]_-{[1,-]^*} &
\map^H(S,F(1_+)) } \]
is a weak equivalence by the previous paragraph,
and both horizontal composites are the identity maps.
Since weak equivalences are closed under retracts, the map $(P_S)^H$
is a weak equivalence. 
This proves that the $\bGamma$-$G$-space $F$ is special.
\end{proof}

\begin{prop}\label{prop:prolonged Wirthmuller}
  Let $G$ be a compact Lie group and $F$ a $G$-cofibrant special $\bGamma$-$G$-space. 
  \begin{enumerate}[\em (i)]
  \item 
    For all finite based $G$-CW-complexes $X$ and $Y$ the map
    \[ (F(p_X),F(p_Y))\ : \ F(X\vee Y)\ \to \ F(X)\times F(Y) \]
    is a $G$-weak equivalence.
  \item 
        For every finite index subgroup $H$ of $G$ and every
    finite based $H$-CW-complex $Z$, the Wirthm{\"u}ller map\index{subject}{Wirthm{\"u}ller map!of an equivariant $\bGamma$-space}
    \[  \omega_Z\ : \ F(G\ltimes_H Z) \ \to \ \map^H(G,F(Z)) \]
    is a $G$-weak equivalence. 
  \item For every finite based $G$-CW-complex $X$,
    the shifted $\bGamma$-$G$-space $F_X$ is special.
  \end{enumerate}
\end{prop}
\begin{proof}
(i)
We wish to show that for every closed subgroup $H$ of $G$ the map
    \[ (F(p_X)^H,F(p_Y)^H)\ : \ F(X\vee Y)^H\ \to \ F(X)^H\times F(Y)^H \]
is a weak equivalence.
The underlying $\bGamma$-$H$-space of $F$ is $H$-cofibrant 
by Proposition \ref{prop:fix preserves Gamma-cofibrant}~(i)
and special by Proposition \ref{prop:fix preserves special}~(i).
Moreover, the underlying $H$-spaces of $X$ and $Y$ are $H$-homotopy
equivalent to finite $H$-CW-complexes, 
by \cite[Cor.\,B]{illman-restricting equivariance}.
So we may assume without loss of generality that $H=G$.

We start with the special case where the group $G$ is finite and
$X=|A|$ and $Y=|B|$ are realizations of two finite based $G$-simplicial sets $A$ and $B$.
For every $n\geq 0$, the map
\[ (F(p_{A_n})^G,F(p_{B_n})^G)\ : \ F(A_n\vee B_n)^G\ \to\ F(A_n)^G\times F(B_n)^G \]
is a weak equivalence by Proposition \ref{prop:special characterizations}~(iii).
The map of simplicial spaces
\[
 ( (F\circ p_A)^G, (F\circ p_B)^G)\ : \ 
 ( F\circ(A\vee B))^G\ \to \ (F\circ A)^G\times (F\circ B)^G  \]
is thus levelwise a weak equivalence.
The simplicial spaces 
$( F\circ(A\vee B))^G$, $(F\circ A)^G$ and $(F\circ B)^G$ are Reedy cofibrant by 
Proposition \ref{prop:G-cof Reedy cof}~(iv).
The product simplicial space $ (F\circ A)^G\times (F\circ B)^G$ is then 
Reedy cofibrant by Proposition \ref{prop:box preserves Reedy cofibrant}.
So the map induced on geometric realizations 
\[ | (F\circ(A\vee B))^G|\ \to \ |(F\circ A)^G\times (F\circ B)^G|  \]
is a weak equivalence by Proposition \ref{prop:realization invariance in Top}.
Realization commutes with fixed points 
(Proposition \ref{prop:G-fix preserves pushouts}~(iv))
and products (Proposition \ref{prop:iterated geometric realization}~(ii)), so the map 
\[
 ( |F \circ p_A|, |F\circ p_B|)^G\ : \ 
 | F\circ(A\vee B)|^G \ \to \ ( |F\circ A|\times |F\circ B| )^G  \]
is a weak equivalence.
The Fubini isomorphism $|F\circ A|^G\iso F(|A|)^G$ 
of Proposition \ref{prop:prolong versus realize} then translates this into the claim 
for the geometric realizations of $A$ and $B$.

Now we continue to assume that $G$ is finite, but $X$ and $Y$
are arbitrary finite based $G$-CW-complexes.
We choose based $G$-homotopy equivalences 
$|A|\simeq X$ and $|B|\simeq Y$ for suitable finite based $G$-simplicial sets
$A$ and $B$ as provided by Proposition \ref{prop:equivariant simplicial realization} (ii).
Prolonged $\bGamma$-$G$-spaces preserve equivariant homotopy equivalences,
so the general case follows from the special case.

It remains to treat the general case of a compact Lie group.
We let $G^\circ$ be the identity component of $G$
and we write $\bar G=G/G^\circ$ for the finite group of path components.
The $\bGamma$-$\bar G$-space $F^{G^\circ}$ is $\bar G$-cofibrant 
by Proposition \ref{prop:fix preserves Gamma-cofibrant}~(ii)
and special by Proposition \ref{prop:fix preserves special}~(ii).
Since $X$ and $Y$ are finite $G$-CW-complexes, $X^{G^\circ}$ 
and $Y^{G^\circ}$ are finite $\bar G$-CW-complexes.
So the map
\[ 
( F^{G^\circ}(p_{X^{G^\circ}}),F^{G^\circ}(p_{Y^{G^\circ}}))\ : \ 
 F^{G^\circ}(X^{G^\circ}\vee Y^{G^\circ})\ \to \ 
F^{G^\circ}(X^{G^\circ})\times F^{G\circ}(Y^{G^\circ})\]
induces a weak equivalence on $\bar G$-fixed points by the previous paragraph.
Moreover,
\[
(F^{G^\circ}(X^{G^\circ}))^{\bar G}\ \iso \ (F(X)^{G^\circ})^{\bar G} \ = \   F(X)^G 
 \]
by Proposition \ref{prop:Gamma fixed points}, and similarly for the
$G$-fixed points of $F(Y)$ and $F(X\vee Y)$.
This completes the proof.

(ii) 
We start by showing that the $G$-fixed point map
\[  (\omega_Z)^G\ : \ F(G\ltimes_H Z)^G \ \to \ ( \map^H(G,F(Z)) )^G \]
is a weak equivalence. 
Evaluation at $1\in G$ identifies the target of $(\omega_Z)^G$
with $F(Z)^H$, so we may show that the composite map
\begin{equation}  \label{eq:Wirth_rewritten}
  F(G\ltimes_H Z)^G \ \xra{\text{incl}} \ F(G\ltimes_H Z)^H \ \xra{F(l_Z)^H} \ F(Z)^H   
\end{equation}
is a weak equivalence.

We start with the special where $G$ is finite and 
$Z=|B|$ is the geometric realization 
of a finite based $H$-simplicial set $B$;
then $G\ltimes_H Z$ is $G$-homeomorphic to the geometric realization 
of the finite based $G$-simplicial set $G\ltimes_H B$.
By Proposition \ref{prop:prolong versus realize}
and \ref{prop:iterated geometric realization}~(ii), 
the space $F(|G\ltimes_H B|)^G$ is homeomorphic to $| (F\circ (G\ltimes_H B))^G|$. 
Similarly,  $F(|B|)^H$ is homeomorphic to $| (F\circ B)^H|$.
Under these homeomorphisms, the map \eqref{eq:Wirth_rewritten} 
becomes the geometric realization
of the morphism of simplicial spaces
\begin{equation}  \label{eq:simplicial_wirthmuller}
(F\circ l_B)^H\circ \text{incl}\ : \ (F\circ (G\ltimes_H B))^G \ \to \ (F\circ B)^H   
\end{equation}
whose value at the object $[n]$ of $\bDelta^{\op}$ is the map
$F(l_{B_n})^H\circ \text{incl}$.
The morphism of simplicial spaces \eqref{eq:simplicial_wirthmuller}
is a weak equivalence in every simplicial dimension
by Proposition \ref{prop:special characterizations}~(iii);
moreover, source and target are Reedy cofibrant 
by Proposition \ref{prop:G-cof Reedy cof}~(iv).
As a levelwise weak equivalence between Reedy cofibrant simplicial spaces,
the morphism \eqref{eq:simplicial_wirthmuller} induces
a weak equivalence on geometric realizations,
by Proposition \ref{prop:realization invariance in Top}. 
This completes the proof that the map \eqref{eq:Wirth_rewritten}
is a weak equivalence in the special case $Z=|B|$.

Now we treat the case where $G$ is finite, but $Z$ is any finite based $H$-CW-complex.
Proposition \ref{prop:equivariant simplicial realization}~(ii)
provides a finite based $H$-simplicial set $B$ 
and a based $H$-homotopy equivalence $h:|B|\to Z$.
The map $G\ltimes_H h:G\ltimes_H |B|\to G\ltimes_H Z$
is then a based $G$-homotopy equivalence.
In the commutative square
\[ \xymatrix@C=18mm{ 
F(G\ltimes_H |B|)^G\ar[r]^-{F(l_{|B|})^H\circ\text{incl}}
\ar[d]_{F(G\ltimes_H h)^G} &
F(|B|)^H\ar[d]^{F(h)^H} \\
F(G\ltimes_H Z)^G \ar[r]_-{F(l_Z)^H\circ\text{incl}} & F(Z)^H } \]
both vertical maps are then homotopy equivalences.
The upper map is a weak equivalence by the previous paragraph,
hence so is the lower map.
This completes the proof that the map \eqref{eq:Wirth_rewritten}
is a weak equivalence in the special case when $G$ is finite.

It remains to treat the case of a general compact Lie group.
We let $G^\circ$ be the identity component of $G$
and we write $\bar G=G/G^\circ$ for the finite group of path components.
The $\bGamma$-$\bar G$-space $F^{G^\circ}$ is $\bar G$-cofibrant 
by Proposition \ref{prop:fix preserves Gamma-cofibrant}~(ii)
and special by Proposition \ref{prop:fix preserves special}~(ii).
Since $H$ has finite index in $G$, the identity components 
of $H$ and $G$ are the same, i.e., $H^\circ=G^\circ$.
Moreover, $\bar H=H/H^\circ$ is a subgroup of $\bar G$.

Since $Z$ is a finite $H$-CW-complex, $X^{H^\circ}$ is a finite $\bar H$-CW-complex.
So the map
\[ 
( F^{G^\circ}(l_{Z^{H^\circ}}) )^{\bar H}\circ \text{incl}\ : \ 
( F^{G^\circ}( \bar G\ltimes_{\bar H} Z^{H^\circ}) )^{\bar G}\ \to \ 
( F^{G^\circ}(Z^{H^\circ}))^{\bar H} \]
is a weak equivalence by the previous paragraph.
Proposition \ref{prop:Gamma fixed points} lets us rewrite 
source and target of this map as
\begin{align*}
( F^{G^\circ}( \bar G\ltimes_{\bar H} Z^{H^\circ}))^{\bar G} 
\ &\iso\ ( F^{G^\circ}( (G\ltimes_H Z)^{G^\circ}))^{\bar G} \\ 
\ &\iso\ ( F( G\ltimes_H Z)^{G^\circ})^{\bar G} \ = \  F(G\ltimes_H Z)^G
\end{align*}
respectively 
\begin{align*}
( F^{G^\circ}(Z^{H^\circ}))^{\bar H} \ &= \ 
( F^{H^\circ}(Z^{H^\circ}))^{\bar H} \ \iso\
  (F(Z)^{H^\circ})^{\bar H} \ \iso \ F(Z)^H\ .
\end{align*}
Under these homeomorphisms, the map 
$( F^{G^\circ}(l_{Z^{H^\circ}}) )^{\bar H}\circ \text{incl}$ becomes 
the map \eqref{eq:Wirth_rewritten}.
This completes the proof that the map \eqref{eq:Wirth_rewritten}
is a weak equivalence for all compact Lie groups $G$ and all finite based
$H$-CW-complexes $Z$.
Hence the map $(\omega_Z)^G$ is a weak equivalence.

Now we let $K$ be any closed subgroup of $G$. Since $H$ has finite index in $G$,
the set $G/H$ is finite, and there is a finite set $\{\gamma_i\}_{i=1,\dots,n}$
of $K$-$H$-double coset representatives.
These determine a $K$-equivariant wedge decomposition
\[ D\ : \ {\bigvee}_{i=1,\dots,n}\ K\ltimes_{K\cap {^{\gamma_i} H}} c_{\gamma_i}^*(Z)
\ \iso \ G\ltimes_H Z\ , \quad
[k,z]_i \ \longmapsto \ [k\gamma_i, z]\ .\]
We obtain a commutative square:
\[  \xymatrix@C=12mm{ 
 F(G\ltimes_H Z)^K\ar[r]^-{(\omega_Z)^K}\ar[d]_\simeq &  ( \map^H(G,F(Z)) )^K \ar[d]^\iso \\
\prod_{i=1}^n F(K\ltimes_{K\cap {^{\gamma_i} H}} c_{\gamma_i}^*(Z))^K
\ar[r]_-{\prod \omega_{c_{\gamma_i}^*(Z)}} & \prod_{i=1}^n F(Z)^{K^{\gamma_i}\cap H}}\]
The right vertical map is evaluation at the chosen coset representatives,
and it is a homeomorphism.
The left map is $F( D^{-1})^K$ 
followed by $F(-)^K$ applied to the projections to the wedge summands.
The latter is a weak equivalence by the additivity property of part~(i),
applied to the underlying $\bGamma$-$K$-space of $F$,
and the spaces $K\ltimes_{K\cap {^{\gamma_i} H}} c_{\gamma_i}^*(Z)$;
there is a slight caveat, namely that the underlying $(K^{\gamma_i}\cap H)$-space
of $Z$ need not admit the structure of an equivariant CW-complex. However, it is always
$(K^{\gamma_i}\cap H)$-homotopy equivalent to a finite 
$(K^{\gamma_i}\cap H)$-CW-complex, by \cite[Cor.\,B]{illman-restricting equivariance};
this is enough to run the argument because prolonged
$\bGamma$-$K$-spaces preserve equivariant homotopy equivalences.
The lower horizontal map in the above square is a weak equivalence
by the previous paragraph,
applied to the underlying $\bGamma$-$K$-space of $F$,
and the subgroups $K\cap {^{\gamma_i} H}$.
So we can conclude that the map $(\omega_Z)^K$ is a weak equivalence. 
Since $K$ was an arbitrary closed subgroup of $G$,
this completes the proof that the Wirthm{\"u}ller
map $\omega_Z$ is a $G$-weak equivalence.

(iii)
We verify the criterion of Proposition \ref{prop:special characterizations}~(iii)
for the shifted $\bGamma$-$G$-space $F_X$.
Given two finite $G$-sets $T$ and $U$, part~(i) shows that 
the map
\begin{align*}
  (F_X(p_T),F_X(p_U))\ : \ F_X((T\amalg U)_+)\ &= \ F(X\sm (T\amalg U)_+) \\ 
\to \ F(X\sm T_+)&\times F(X\sm U_+) \ = \ F_X(T_+)\times F_X(U_+) 
\end{align*}
is a $G$-weak equivalence.
For every finite index subgroup $H$ of $G$ and every finite $H$-set $S$,
the effect of the shearing isomorphism
\[ X\sm (G\times_H S)_+ \ \iso \ G\ltimes_H(X\sm S_+) \]
identifies the Wirthm{\"u}ller map
\begin{align*}
\omega_{S_+}\ : \ F_X((G\times_H S)_+)\ &= \ F(X\sm (G\times_H S)_+)\\ 
\ &\to \ \map^H(G,F(X\sm S_+))\ = \ \map^H(G,F_X(S_+))    
\end{align*}
for the shifted $\bGamma$-$G$-space and the $H$-space $S_+$
with the Wirthm{\"u}ller map
\[ \omega_{X\sm S_+}\ : \ F(G\ltimes_H (X\sm S_+))\ \ \to \ \map^H(G,F(X\sm S_+)) \]
of the original $\bGamma$-$G$-space and the $H$-space $X\sm S_+$.
The latter map is a $G$-weak equivalence by part~(ii),
hence so is the former.
Proposition \ref{prop:special characterizations} now applies and shows that
the shifted $\bGamma$-$G$-space $F_X$ is special.
\end{proof}

We still have to recall the notion of a `very special' $\bGamma$-$G$-space.
We let $G$ be a compact Lie group and $F$ a special $\bGamma$-$G$-space.
We let $p_1,p_2:2_+\to 1_+$ denote the two projections.
The map
\[ (F(p_1),F(p_2))\ : \ F(2_+)\ \to \ F(1_+)\times F(1_+) \]
is a $G$-weak equivalence by specialness.
We let $\nabla:2_+\to 1_+$ denote the fold map.
We obtain a diagram of set maps
\[ \pi_0( F(1_+)^G)\times \pi_0(F(1_+)^G)\ \xla[\iso]{ \left(\pi_0(F(p_1)^G),\pi_0(F(p_2)^G)\right) } 
 \pi_0( F( 2_+)^G)\ \xra{ \pi_0(F(\nabla)^G)} 
 \pi_0( F(1_+)^G) \]
the left of which is bijective.
So the map
\begin{align*}
 \pi_0(F(\nabla)^G)\circ (\pi_0(F(p_1)^G),&\ \pi_0(F(p_2)^G))^{-1} \ : \\ 
& \pi_0( F(1_+)^G)\times \pi_0(F(1_+)^G)\ \to\  \pi_0( F(1_+)^G)   
\end{align*}
is a binary operation on the set $\pi_0(F(1_+)^G)$.

If $\tau:2_+\to 2_+$ is the involution that interchanges~1 and~2,
then composition with $\tau$ 
fixes $\nabla$ and interchanges $p_1$ and $p_2$;
this implies that the operation $+$ is commutative.
Contemplating the different ways to fold and project from the based set $3_+$
leads to the proof that the operation is also associative,
and hence an abelian monoid structure on the set $\pi_0(F(1_+)^G)$.

For every closed subgroup $H$ of $G$, 
the underlying $\bGamma$-$H$-space is again special. 
So the same argument provides an abelian monoid structure on $\pi_0(F(1_+)^H)$.

\begin{defn}\label{def:very special G-Gamma}
Let $G$ be a compact Lie group. 
A $\bGamma$-$G$-space $F$ is {\em very special}\index{subject}{Gamma-space@$\bGamma$-space!very special}\index{subject}{Gamma-space@$\bGamma$-space!equivariant}  
if it is special and for every closed subgroup $H$ of $G$ 
the abelian monoid $\pi_0(F(1_+)^H)$ is an abelian group.
\end{defn}

\begin{rk}\label{rk:label very special gives sAb}
We let $G$ be a compact Lie group and $F$ a special $\bGamma$-$G$-space. 
Then the construction of the abelian monoid structure
on $\pi_0(F(1_+)^H)$ generalizes as follows. 
We let $H$ be a closed subgroup of $G$ and $S$ a finite $G$-set.
We let $p_1,p_2:(S\amalg S)_+\to S_+$ denote the two projections.
The map
\[ (\pi_0(F(p_1)^H),\pi_0(F(p_2)^H))\ : \ \pi_0(F((S\amalg S)_+)^H)
\ \to \ \pi_0(F(S_+)^H)\times \pi_0(F(S_+)^H) \]
is bijective by Proposition \ref{prop:special characterizations}~(iii);
inverting this map and composing with the effect of 
the fold map $\nabla:(S\amalg S)_+\to S_+$ on $\pi_0(F(-)^H)$ 
yields a binary operation on the set $\pi_0( F(S_+)^H)$.
The same arguments as in the special case $S=\{1\}$ show that this operation
defines an abelian monoid structure on $\pi_0(F(S_+)^H)$, which is moreover natural
for $G$-maps in $S$. 
If $K$ is a closed subgroup of $H$, then the inclusion $F(S_+)^H\subset F(S_+)^K$ 
induces  a monoid homomorphism $\pi_0(F(S_+)^H)\to\pi_0(F(S_+)^K)$. 

We claim that for very special $F$, all the abelian monoids $\pi_0(F(S_+)^H)$
are in fact groups.
To see this, we choose a set of representatives $s_1,\dots,s_n$
for the $H$-orbits of $S$, and we let $K_i$ be the stabilizer group
of $s_i$ inside $H$. We consider the composite
\[ F(S_+)^H\ \xra{(P_S)^H} \ \map^H(S,F(1_+))\ \iso \ {\prod}_{i=1}^n F(1_+)^{K_i}\ , \]
where the first map is a weak equivalence by specialness, 
and the homeomorphism is evaluation at the elements $s_1,\dots,s_n$.
The induced bijection
\[ \pi_0(F(S_+)^H)\ \iso \ {\prod}_{i=1}^n \, \pi_0( F(1_+)^{K_i}) \]
is a homomorphism of abelian monoids by the previous paragraph. 
The target is an abelian group since $F$ is very special.
So the abelian monoid $\pi_0(F(S_+)^H)$ is also an abelian group.
\end{rk}

\begin{eg}\label{eg:Gamma is pi_*-Kan} 
We recall how very special $\bGamma$-$G$-spaces
give rise to simplicial spaces that satisfy the $\pi_*$-Kan condition
(see \cite[B.3]{bousfield-friedlander} or Definition \ref{def:pi_*-Kan}).
We let $G$ be a finite group and $F:\bGamma\to G\bT_*$ 
a very special $\bGamma$-$G$-space. 
We let $A$ be a simplicial finite based $G$-set. 
Then for every subgroup $H$ of $G$ we obtain a simplicial space
$(F\circ A)^H$ whose space of $n$-simplices is $F(A_n)^H$
(where are usual, $F(A_n)$ has the diagonal $H$-action by restriction of
the $G$-actions on $F$ and $A_n$). We claim that this simplicial space
$(F\circ A)^H$ satisfies the $\pi_*$-Kan condition.
To see this, we borrow the argument from Bousfield 
and Friedlander \cite[Proof of Lemma 4.3]{bousfield-friedlander}:
since $F$ is special, the space $F(A_n)^H$ comes with
a homotopy-commutative H-space structure, so each of its path components is
a simple space. Moreover, the map of simplicial sets
$\beta:\pi_t((F\circ A)^H)_{\text{free}}\to \pi_0( (F\circ A)^H)$
discussed in B.3 of \cite{bousfield-friedlander}
is underlying a surjective morphism of simplicial groups,
and is thus a Kan fibration by \cite[Lemma 3.2]{curtis}.
So \cite[(B.3.1)]{bousfield-friedlander} proves that
$(F\circ A)^H$ satisfies the $\pi_*$-Kan condition.
\end{eg}

\begin{prop}\label{prop:fix preserves very special}
Let $H$ be a closed subgroup of a compact Lie group $G$ 
and $F$ a very special $\bGamma$-$G$-space.
\begin{enumerate}[\em (i)]
\item The underlying $\bGamma$-$H$-space of $F$ is very special.
\item If $H$ is normal, then the $\bGamma$-$G/H$-space $F^H$ is very special.
\end{enumerate}
\end{prop}
\begin{proof}
Part~(i) is clear by definition.
For part~(ii) we recall from Proposition \ref{prop:fix preserves special}~(ii)
that the $\bGamma$-$G/H$-space $F^H$ is special.
Any closed subgroup of $G/H$ is of the form $\Delta/H$ 
for a closed subgroup $\Delta$ of $G$ with $H\leq \Delta$.
Then $(F^H(1_+))^{\Delta/H}=F(1_+)^\Delta$.
Since the monoid $\pi_0(F(1_+)^\Delta)$ is a group by hypothesis,
we conclude that the monoid $\pi_0(F^H(1_+)^{\Delta/H})$ is a group.
\end{proof}

The following Theorem \ref{thm:prolonged delooping}
is essentially a reformulation of the equivariant delooping results 
of Segal \cite{segal-some_equivariant} respectively Shimakawa \cite{shimakawa};
however, I am not aware that the result has been established in this form. 
The difference is that we formulate the result for the
prolongation (i.e., categorical Kan extension), whereas Segal and Shimakawa
work with a bar construction (also known as a homotopy coend
or homotopy Kan extension) instead.
We also give a partial extension of the machinery to compact Lie groups. 

A minor technical difference is that in Shimakawa's 
paper \cite{shimakawa}, an equivariant $\bGamma$-space is assumed to take
values in non-degenerately based $G$-spaces with the homotopy type
of a based $G$-CW-complex; similarly, Segal \cite{segal-some_equivariant} 
is implicitly assuming that the values of equivariant $\bGamma$-spaces are $G$-ANRs.
Our proof includes the observation that the key arguments work
in the larger category of compactly generated $G$-spaces and $G$-weak equivalences,
a long as we restrict to $G$-cofibrant $\bGamma$-$G$-spaces.

Our proof verifies that the prolonged cofibrant $\bGamma$-$G$-space 
satisfies the delooping criterion of Blumberg \cite[Thm.\,1.2]{blumberg}.
An advantage of this approach is that Blumberg does not quote
Segal or Shimakawa; 
rather, he works in the category of compactly generated $G$-spaces throughout,
and adapts the relevant arguments to this context.

The paper \cite{mmo} by May, Merling and Osorno also contains a modern
perspective on the work of Segal and Shimakawa, 
and there is some overlap with our present discussion.
In fact, our exposition of equivariant $\bGamma$-spaces
partially arose from discussions with Mona Merling and Peter May,
to whom I am indepted for valuable feedback. 
In \cite{mmo} the method of producing an equivariant spectrum 
by evaluating a $\bGamma$-space
on spheres is referred to as the {\em conceptual Segal machine}.
In various respects, \cite{mmo} goes further than we do here;
for example, it contains detailed comparisons of prolongation 
(the `conceptual Segal machine'),
bar construction (the `homotopical Segal machine')
and the operadic approach to equivariant delooping (via a `generalized Segal machine').
Moreover, for special $\bGamma$-$G$-spaces,
the value at $S^1$ is identified as an equivariant group completion;
this is particularly relevant since many $\bGamma$-$G$-spaces
that arise in practice are special, but {\em not} very special.

\begin{theorem}\label{thm:prolonged delooping}
  Let $G$ be a compact Lie group and $F$ a $G$-cofibrant 
  very special $\bGamma$-$G$-space.
  \begin{enumerate}[\em (i)]
  \item 
    The prolonged functor $F$ takes $G$-homotopy pushout squares 
    of finite based $G$-CW-complexes to $G$-homotopy pullback squares.
  \item
    For every finite based $G$-CW-complex $X$ and every $G$-representation $V$
    on which the identity component $G^\circ$ acts trivially, the adjoint assembly map
    \[ \tilde\alpha\ : \ F(X)\ \to \ \map_*(S^V,F(X\sm S^V)) \]
    is a $G$-weak equivalence.
    \item If $G$ is finite, then the orthogonal $G$-spectrum $F(\mS)$
      is a $G$-$\Omega$-spectrum.\index{subject}{G-Omega-spectrum@$G$-$\Omega$-spectrum}
  \end{enumerate}
\end{theorem}
\begin{proof}
(i) We start with the special case where $G$ is finite, and we consider a monomorphism
of finite based $G$-simplicial sets $i:A\to B$. We show that then the sequence
\begin{equation}  \label{eq:F_homotopy_fiber_sequence}
 F(|A|) \ \xra{F(|i|)}\ F(|B|) \ \xra{F(|q|)}\ F(|B/A|)   
\end{equation}
is a $G$-homotopy fiber sequence, where $q:B\to B/A$ is the projection.
We adapt an argument that Bousfield and Friedlander
use for non-equivariant $\bGamma$-spaces of simplicial sets 
in \cite[Lemma 4.3]{bousfield-friedlander}. 
The space $F(|A|)$ is homeomorphic to the geometric realization
of the simplicial space $F\circ A$, by Proposition \ref{prop:prolong versus realize}.
Moreover, for $H\leq G$, taking $H$-fixed points commutes with geometric realization by 
Proposition \ref{prop:G-fix preserves pushouts}~(iv),
so $F(|A|)^H$ is homeomorphic to the geometric realization of the
simplicial space $(F\circ A)^H$, sending $[n]$ to $F(A_n)^H$.
The same applies to $B$ and $B/A$,
so we may show that the sequence
\begin{equation}\label{eq:fix inside realization}
|(F\circ A)^H|\ \xra{|(F\circ i)^H|} \
|(F\circ B)^H|\ \xra{|(F\circ q)^H|} \ |(F\circ (B/A))^H|   
\end{equation}
is a homotopy fiber sequence.
We let $r:B_n\to A_n$ be the retraction to $i_n$ that sends the complement of $i_n(A_n)$
to the basepoint.
For fixed $n\geq 0$ the lower row in the commutative diagram
\[ \xymatrix@C=15mm{ 
F(A_n)^H\ar[r]^-{F(i_n)^H} \ar@{=}[d] &
F(B_n)^H\ar[r]^-{F(q_n)^H} \ar[d]^{(F(r)^H,F(q_n)^H)}_\simeq & 
F(B_n/A_n)^H  \ar@{=}[d]\\ 
F(A_n)^H \ar[r] & F(A_n)^H\times F( B_n/A_n)^H\ar[r]_-{\text{proj}}&  F( B_n/A_n )^H} \]
is a homotopy fiber sequence.
The middle vertical map is a weak equivalence by 
Proposition \ref{prop:special characterizations}~(iii).
So the upper row in the diagram is a homotopy fiber sequence.

Because $F$ is $G$-cofibrant, 
the three simplicial spaces $(F\circ A)^H$,  $(F\circ B)^H$
and $(F\circ (B/A))^H$ are Reedy cofibrant by
Proposition \ref{prop:G-cof Reedy cof}~(iv).
So we are considering a sequence of simplicial spaces that is a 
homotopy fiber sequence in every simplicial dimension.
Moreover, the simplicial spaces $(F\circ B)^H$
and $(F\circ (B/A))^H$ satisfy the $\pi_*$-Kan condition 
by Example \ref{eg:Gamma is pi_*-Kan}.
The morphism of simplicial sets 
$\pi_0((F\circ q)^H):\pi_0( (F\circ B)^H)\to\pi_0( (F\circ (B/A))^H)$
underlies a homomorphism of simplicial abelian groups
by Remark \ref{rk:label very special gives sAb};
this morphism is split surjective in every simplicial dimension.
Every surjective homomorphism of simplicial groups is a Kan fibration,
see for example \cite[Lemma 3.2]{curtis} or \cite[V Lemma 2.6]{goerss-jardine}.
So the sequence $(F\circ A)^H\to (F\circ B)^H\to (F\circ (B/A))^H$
satisfies the hypotheses of Bousfield and Friedlander's 
Theorem \ref{thm:realization preserves cartesian}; 
we conclude that the sequence \eqref{eq:fix inside realization}
is a homotopy fiber sequence.
This completes the proof that the sequence \eqref{eq:F_homotopy_fiber_sequence}
is a $G$-homotopy fiber sequence.

Now we consider a $G$-cofibration $j:X\to Y$ 
between finite based $G$-CW-complexes, where $G$ is still a finite group.
We claim that then the sequence
\begin{equation}  \label{eq:X2Y2Y/X}
 F(X) \ \xra{F(j)}\ F(Y) \ \xra{F(p)}\ F(Y/X)      
\end{equation}
is a $G$-homotopy fiber sequence.
Proposition \ref{prop:equivariant simplicial realization}
provides a finite $G$-simplicial set $A$ 
and a $G$-homotopy equivalence $k:|A|\to X$.
Another application of Proposition \ref{prop:equivariant simplicial realization}
provides a finite $G$-simplicial set $B$, 
a monomorphism of $G$-simplicial sets $i:A\to B$ and a 
$G$-homotopy equivalence $h:|B|\to Y$ such that $h\circ|i|=j\circ k:|A|\to Y$.
The induced map on cokernels
\[ h/k \ : \ |B / A| \ \iso \ |B| /  |A| \ \to \  Y/X \]
is then a $G$-weak equivalence by the gluing lemma 
(Proposition \ref{prop:gluing lemma G-spaces}).
Source and target of $h/k$ are cofibrant as $G$-spaces,
so this $G$-weak equivalence is even a $G$-homotopy equivalence.
The upper row in the commutative diagram of $G$-spaces
\[ \xymatrix@C=18mm{ 
F(|A|)\ar[r]^-{F(|i|)} \ar[d]_{F(k)}^\simeq &
F(|B|)\ar[r]^-{F(|q|)} \ar[d]_{F(h)}^\simeq & 
F(|B/A|)  \ar[d]^{F(h/k)}_\simeq\\ 
F(X) \ar[r]_-{F(j)} & F(Y)\ar[r]_-{F(p)}&  F(Y/X)} \]
is a $G$-homotopy fiber sequence by the previous paragraph.
The prolongation of every $\bGamma$-$G$-space preserves
based $G$-homotopy equivalences; so all vertical maps are $G$-homotopy equivalences.
Hence the lower row in the diagram is a $G$-homotopy fiber sequence.

Now we let $G$ be an arbitrary compact Lie group, and $j:X\to Y$ 
a $G$-cofibration between finite based $G$-CW-complexes.
We claim that then the sequence \eqref{eq:X2Y2Y/X}
is a $G$-homotopy fiber sequence.
We let $H$ be a closed subgroup of $G$ with identity component $H^\circ$ 
and group of path components $\bar H=H/H^\circ$.
The $\bGamma$-$\bar H$-space $F^{H^\circ}$ is $\bar H$-cofibrant 
by Proposition \ref{prop:fix preserves Gamma-cofibrant}
and very special by Proposition \ref{prop:fix preserves very special}.
The map $j^{H^\circ}:X^{H^\circ}\to Y^{H^\circ}$ is an $\bar H$-cofibration
by Propositions \ref{prop:cofibrancy preservers}~(i) 
and \ref{prop:cofamily pushout property}.
Since $j$ is in particular a closed embedding, 
the canonical map $X^{H^\circ}/Y^{H^\circ}\to (X/Y)^{H^\circ}$ is a homeomorphism
by Proposition \ref{prop:G-fix preserves pushouts}~(i).
So by the previous paragraph, the sequence
\[  F^{H^\circ}(X^{H^\circ}) \ \xra{F^{H^\circ}(j^{H^\circ})}\ F^{H^\circ}(Y^{H^\circ}) \ 
\xra{F^{H^\circ}(p^{H^\circ})}\ F^{H^\circ}( (Y/X)^{H^\circ})   \]
is an $\bar H$-homotopy fiber sequence.
We may thus take $\bar H$-fixed points and obtain a non-equivariant
homotopy fiber sequence.
Proposition \ref{prop:Gamma fixed points} provides a homeomorphism
$(F^{H^\circ}(X^{H^\circ}))^{\bar H}\iso (F(X)^{H^\circ})^{\bar H}=F(X)^H$,
and similarly for the other two terms.
So we conclude that the sequence
\[  F(X)^H \ \xra{F(j)^H}\ F(Y)^H \ \xra{F(p)^H}\ F(Y/X)^H     \]
is a non-equivariant homotopy fiber sequence.
This completes the proof that the sequence \eqref{eq:X2Y2Y/X}
is a $G$-homotopy fiber sequence.

Now we treat a general homotopy cocartesian square 
of finite based $G$-CW-complexes:
\[ \xymatrix{ X\ar[r]^-j\ar[d]_k & Y\ar[d]^h \\ Z \ar[r]_i & P}     \]
By replacing $Y$ and $P$ by the reduced mapping cylinders
of $j$ respectively $i$, 
we may assume that the horizontal maps are $G$-cofibrations.
Because the square is homotopy cocartesian, the induced map
$h/k:Y/X\to P/Z$ is a $G$-weak equivalence, and hence a $G$-homotopy equivalence.
We obtain a commutative diagram of based $G$-spaces:
\[ \xymatrix@C=12mm{ 
F(X)\ar[r]^-{F(j)} \ar[d]_{F(k)} & F(Y) \ar[d]_{F(h)}\ar[r]& 
F(Y/X)  \ar[d]^{F(h/k)}_\simeq\\ 
F(Z) \ar[r]_-{F(i)} & F(P)\ar[r] &  F(P/Z)} \]
Because $F$ preserves equivariant homotopy equivalences, 
the right vertical map  $F(h/k)$ is a $G$-homotopy equivalence.
Both rows are $G$-homotopy fiber sequences by the previous paragraph;
so the left square is $G$-homotopy cartesian.
This completes the proof that the prolonged functor 
takes $G$-homotopy cocartesian squares to $G$-homotopy cartesian squares.

(ii)
As usual we start with the special case where $G$ is finite.
Restricting the prolonged functor $F$ to finite $G$-CW-complexes
gives a $\mathscr W_G$-space in the sense of \cite{blumberg}.
This $\mathscr W_G$-space satisfies conditions (1) and~(2') 
of \cite[Thm.\,1.2]{blumberg}, by part~(i) of this theorem
respectively Proposition \ref{prop:prolonged Wirthmuller}~(ii).
So \cite[Thm.\,1.2]{blumberg} applies to the prolonged functor,
and gives the desired conclusion.

Now we let $G$ be a general compact Lie group.
We want to show that for every closed subgroup $H$ of $G$ the map
    \[ \tilde\alpha^H\ : \ F(X)^H\ \to \ \map^H_*(S^V,F(X\sm S^V)) \]
    is a weak equivalence.
The underlying $\bGamma$-$H$-space of $F$ is $H$-cofibrant 
by Proposition \ref{prop:fix preserves Gamma-cofibrant}~(i)
and very special by Proposition \ref{prop:fix preserves very special}~(i).
Moreover, the underlying $H$-space of $X$ is $H$-homotopy
equivalent to a finite $H$-CW-complex, 
by \cite[Cor.\,B]{illman-restricting equivariance}.
Finally, the identity component of $H$ is contained in $G^\circ$,
so it acts trivially on $V$.
So we may assume without loss of generality that $H=G$.

As usual we write $\bar G=G/G^\circ$ for the finite group of path components.
The $\bGamma$-$\bar G$-space $F^{G^\circ}$ is $\bar G$-cofibrant 
by Proposition \ref{prop:fix preserves Gamma-cofibrant}~(ii)
and very special by Proposition \ref{prop:fix preserves very special}~(ii).
Since $X$ is a finite $G$-CW-complex, $X^{G^\circ}$ is a finite $\bar G$-CW-complex.
Since $G^\circ$ acts trivially on $V$, we can view $V$ as a representation
of the factor group $\bar G$. 
The previous paragraph shows that the map
\[ \tilde\alpha\ : \ F^{G^\circ}(X^{G^\circ})\ \to \ \map_*(S^V,F^{G^\circ}(X^{G^\circ}\sm S^V)) \]
induces a weak equivalence on $\bar G$-fixed points.
Proposition \ref{prop:Gamma fixed points} provides a homeomorphism
\[   
( F^{G^\circ}( X^{G^\circ}) )^{\bar G}\ \iso\ (F(X)^{G^\circ})^{\bar G} \ = \  F(X)^G\ ; \]
similarly,
\begin{align*}
  \map^{\bar G}_*(S^V,F^{G^\circ}(X^{G^\circ}\sm S^V)) \ &\iso \ 
  \map^{\bar G}_*(S^V, F(X\sm S^V)^{G^\circ}) \\ 
  &= \ \map^G_*(S^V, F(X\sm S^V) ) \ .
\end{align*}
Under these identifications, the map $\tilde\alpha^{\bar G}$
corresponds to the $G$-fixed points of the adjoint assembly map
for the $\bGamma$-$G$-space $F$. This completes the proof.

Part~(iii) is simply the special case of~(ii) with $X=S^W$
for some $G$-re\-presen\-tation $W$. Finite groups have a trivial identity
component, which must necessarily act trivially on $V$.
\end{proof}

Now we let $F$ be a cofibrant $\bGamma$-$G$-space that is special
(but not necessarily very special). Then Theorem \ref{thm:prolonged delooping}
need not hold for $F$, as the adjoint assembly map $F(1_+)\to \Omega F(S^1)$
need not be a weak equivalence.
However, there is a standard argument to deduce 
from Theorem \ref{thm:prolonged delooping} that the adjoint assembly map
is still a $G$-weak equivalence on suspensions;
this argument uses the shift of the given $\bGamma$-$G$-space by $S^1$.
Since we work with prolongations throughout, there is a new ingredient, 
not already present in the classical references, 
namely that the $S^1$-shift of a $G$-cofibrant $\bGamma$-$G$-space 
is again $G$-cofibrant.

\begin{theorem}\label{thm:special cofibrant Gamma positive Omega}
  Let $G$ be a compact Lie group and $F$ a $G$-cofibrant, special  $\bGamma$-$G$-space. 
  Then for every pair of $G$-representations $V$ and $W$ 
  such that $W^G\ne 0$  and $G^\circ$ acts trivially on $V$,
  the adjoint assembly map
  \[ \tilde\alpha\ : \ F(S^W)\ \to \ \map_*(S^V,F(S^{V\oplus W})) \]
  is a $G$-weak equivalence.
\end{theorem}
\begin{proof}
One more time we start with the special case when $G$ is finite.
We consider the shifted $\bGamma$-$G$-space $F_{S^1}=F\circ(S^1\sm -)$,
which is $G$-cofibrant by Proposition \ref{prop:G-cof Reedy cof}~(ii)
and special by Proposition \ref{prop:prolonged Wirthmuller}~(iii).
Moreover, for every subgroup $H$ of $G$ the space
\[ ( F_{S^1}(1_+))^H\ \iso \ F(S^1)^H \]
is path connected by Proposition \ref{prop:Gamma on S^V}~(i).
So the abelian monoid $\pi_0( F_{S^1}(1_+)^H)$ has only one element,
and is thus an abelian group. Hence $F_{S^1}$ is very special.

So Theorem \ref{thm:prolonged delooping} applies to the 
$G$-cofibrant and very special $\bGamma$-$G$-space $F_{S^1}$.
Since $W$ has a non-zero $G$-fixed point, there is 
a $G$-representation $U$ and a $G$-equivariant linear isometry
$\varphi:\mR\oplus U\iso W$.
We let $\phi:F_{S^1}(L)\to F(S^1\sm L)$ denote the natural homeomorphism established in
Proposition \ref{prop:Fubini}. 
In the commutative diagram
\[ \xymatrix@C=15mm@R=7mm{ 
F_{S^1}(S^{U})\ar[d]^\iso_{\phi} \ar[r]^-{\tilde\alpha}&  
\map_*(S^V,F_{S^1}(S^{V\oplus U}))\ar[d]^{\map_*(S^V,\phi)}_\iso \\
F(S^1\sm S^{U}) \ar[r]^-{\tilde\alpha}\ar[d]_{F(S^\varphi)}^\iso&  
\map_*(S^V,F(S^1\sm S^{V\oplus U})) \ar[d]_\iso^{\map_*(S^V,F(S^{V\oplus\varphi}))}
\\
F(S^W)  \ar[r]_-{\tilde\alpha}& \map_*(S^V,F(S^{V\oplus W})) \\
} \]
all vertical maps are then homeomorphisms,
and the upper horizontal map is a $G$-weak equivalence
by Theorem \ref{thm:prolonged delooping}. 
So the lower horizontal map is a $G$-weak equivalence.
This completes the proof in the case of finite groups.

The case of general compact Lie groups can be reduced to the case of finite groups
by taking fixed points with respect to the identity component $G^\circ$,
and using the residual action of the finite group $G/G^\circ$;
the argument is almost literally the same as in the last paragraph of the proof of 
Theorem \ref{thm:prolonged delooping}~(ii), and we omit it.
\end{proof}

\begin{rk}\label{rk:Gamma limitations}
The previous Theorems \ref{thm:prolonged delooping} 
and \ref{thm:special cofibrant Gamma positive Omega}
indicate why equivariant $\bGamma$-spaces are usually considered
for finite groups only -- 
this is where the equivariant delooping machine reveals its full power.
We want to briefly explain why the spectrum $F(\mS)$ cannot be 
a $G$-$\Omega$-spectrum when $G$ has positive dimension,
except in degenerate cases.
The reason is that when $H$ is a closed subgroup of strictly smaller
dimension than $G$, the `dimension shifting Wirthm{\"u}ller map'
cannot be a $G$-weak equivalence, except in degenerate cases.
Given a based $H$-space $Z$, the $H$-equivariant collapse map
\[ 
 l_Z \ : \ G\ltimes_H Z\ \to  \ Z\sm S^L
 \]
was introduced in \eqref{eq:define_l_A}; here $L=T_{e H}(G/H)$
is the tangent $H$-representation.
The dimension shifting Wirthm{\"u}ller map 
\[ \omega_Z  \ : \ F(G\ltimes_H Z)\ \to \ \map^H(G,F(Z\sm S^L) ) \]
is the adjoint of the $H$-map $F(l_Z):F(G\ltimes_H Z)\to F(Z\sm S^L)$.

Blumberg's recognition theorem \cite[Thm.\,1.2]{blumberg} 
says that his conditions~(1) and~(2') are not only sufficient, 
but also necessary for $F(\mS)$ to be a $G$-$\Omega$-spectrum.
If $F$ is $G$-cofibrant and very special, then the prolonged functor
always satisfies conditions (1), by Theorem \ref{thm:prolonged delooping}~(i).
But condition~(2') requires that the dimension shifting Wirthm{\"u}ller map $\omega_Z$
is a $G$-weak equivalence; 
Proposition \ref{prop:prolonged Wirthmuller}~(ii) shows 
that this is true for compact Lie groups, 
provided $H$ has finite index in $G$ (in which case $L=0$).
However, the map $\omega_Z$ is typically {\em not} a $G$-weak equivalence
when $H$ has strictly smaller dimension than $G$. To illustrate this
we take the extreme case $H=e$ and $Z=1_+$,
in which case $G\ltimes_H Z$ is simply $G_+$.
If the identity component $G^\circ$ of $G$ is non-trivial, then
$F(G_+)^{G^\circ}$ is a single point by Proposition \ref{prop:Gamma fixed points}.
On the other hand, evaluation at a set of representatives
of the path components of $G$ is a homeomorphism
\[  ( \map(G,F(S^L)))^{G^0} \ \iso \   \map(G/G^\circ,F(S^L)) \ \iso \ F(S^L)^m \ ,\]
where $m=[G:G^\circ]$.
So if the Wirthm{\"u}ller map $\omega_{1_+}$ is a $G$-weak equivalence,
then the space $F(S^L)$ is non-equivariantly weakly contractible.
Since $F(\mS)$ is a connective non-equivariant $\Omega$-spectrum,
it must be altogether stably trivial.
\end{rk}
\index{subject}{Gamma-space@$\bGamma$-space!equivariant|)}

\chapter{Enriched functor categories}
\label{app:enriched functors}

In this final appendix we review definitions, properties and constructions
involving categories of enriched functors.
The general setup consists of:
\begin{itemize}
\item a complete and cocomplete closed symmetric monoidal category $\Vc$ 
(the `base category').
We denote the monoidal product in $\Vc$ by $\tensor$;
\item a skeletally small $\Vc$-category $\Dc$ (the `index category').
\end{itemize}
We denote by $\Dc^\ast$ the category of covariant $\Vc$-functors
$X:\Dc\to\Vc$ from the index category to the base category. 
We are mostly interested in the following special cases of such functor categories:
\begin{itemize}
\item orthogonal spaces 
(where $\Vc=\bT$ and $\Dc=\bL$, see Definition \ref{def:orthogonal space});
\item orthogonal spectra 
(where $\Vc=\bT_*$ and $\Dc=\bO$, see Definition \ref{def:orthogonal spectrum});
\item global functors 
(where $\Vc=\Ab$ and $\Dc=\bA$, see Definition \ref{def:global functor}).
\end{itemize}
We discuss the Day convolution product (Definition \ref{def:general box product}) 
that turns the functor category into a symmetric monoidal category, 
compare Theorem \ref{thm:symmetric monoidal}.
Construction \ref{con:skeleta in functor categories}
introduces the skeleton filtration,
which enters into the level model structure on
certain functor categories, see Proposition \ref{prop:general level model structure}. 
For background on symmetric monoidal categories
and enriched category theory we refer the reader to \cite{kelly-enriched category}.

\begin{construction}[$\Vc$-enrichment of functor categories]
The objects of the functor category $\Dc^\ast$  
are the $\Vc$-functors $X:\Dc\to \Vc$, and the morphisms
in $\Dc^\ast$ are the $\Vc$-natural transformations between such functors.
The functor category $\Dc^\ast$ is itself enriched in the base category $\Vc$, as follows. 
Given two objects $X,Y$ in $\Dc^\ast$, a $\Vc$-object of $\Vc$-natural transformations   
$\underline{\Dc}^\ast(X,Y)$ is defined 
as an enriched end \cite[Sec.\,2.1]{kelly-enriched category}, 
i.e., an equalizer in $\Vc$ of two morphisms
\[\xymatrix{
 {\prod}_c \ \underline{\Vc}(X(c),Y(c)) \quad \ar@<.4ex>[r]^-{\alpha}  \ar@<-.4ex>[r]_-{\beta} &
\quad {\prod}_{d,e} \ \underline{\Vc}(\Dc(d,e)\tensor X(d),Y(e))\ . }\]
Here the products run over a set of representatives for the
isomorphism classes of $\Dc$-objects and $\underline{\Vc}(-,-)$ denotes 
the internal function objects in the closed monoidal structure.
The $(d,e)$-component of the  morphism $\alpha$ is the composite
\begin{align*}
 {\prod}_c \ \underline{\Vc}(X(c),Y(c)) \ \xra{\text{proj}_d} \
 \underline{\Vc}(X(d),Y(d)) \ \xra{\quad} \  \underline{\Vc}(\Dc(d,e)\tensor X(d),Y(e)) 
\end{align*}
where the second morphism is adjoint to 
\[  \underline{\Vc}(X(d),Y(d))\tensor \Dc(d,e)\tensor X(d) \ \xra{\text{eval}}\
 \Dc(d,e)\tensor Y(d) \ \xra{\ Y\ } \  Y(e) \ . \]
The $(d,e)$-component of the morphism $\beta$ is the composite 
\begin{align*}
 {\prod}_c \,\underline{\Vc}(X(c),Y(c)) \ \xra{\text{proj}_e} \
 \underline{\Vc}(X(e),Y(e)) \ \xra{\quad} \  \underline{\Vc}(\Dc(d,e)\tensor X(d),Y(e)) 
\end{align*}
where the second morphism is adjoint to 
\[  \underline{\Vc}(X(e),Y(e))\tensor \Dc(d,e)\tensor X(d) \ \xra{\Id\tensor X} \
 \underline{\Vc}(X(e),Y(e))\tensor X(e) \ \xra{\text{eval}}\ Y(e) \ . \]
The {\em underlying set} of a $\Vc$-object $V$ is 
\[ u V \ = \ \Vc(I,V)\ , \]
the set of morphisms from the monoidal unit. 
The functor $\Vc(I,-):\Vc\to\text{(sets)}$ takes equalizers in $\Vc$ to 
equalizers of sets. When applied to the defining equalizer
for $\underline{\Dc}^\ast(X,Y)$ this shows that 
the set of morphisms $\Dc^\ast(X,Y)$ in the functor category $\Dc^\ast$
(i.e., the set of $\Vc$-natural transformations)
can be recovered as the underlying set of the $\Vc$-object $\underline{\Dc}^\ast(X,Y)$,
via a natural bijection
\[ \Dc^{\ast}(X,Y)\ \iso \ \Vc(I,\underline{\Dc}^\ast(X,Y))\ 
= \ u(\underline{\Dc}^\ast(X,Y)) \ .\]
The $\Vc$-objects $\underline{\Dc}^\ast(X,Y)$ even assemble into
an entire $\Vc$-category $\underline{\Dc}^\ast$ of $\Vc$-functors
whose underlying category is $\Dc^{\ast}$.
We refer to \cite[Sec.\,2.1]{kelly-enriched category} for details.

The category $\Dc^\ast$ is also tensored over the base category $\Vc$.
Given a functor $Y\in\Dc^\ast$ and an object $A$ of $\Vc$,
we define a new functor $Y\tensor A\in\Dc^\ast$
by objectwise product in $\Vc$, i.e., $Y\tensor A$ is the composite $\Vc$-functor
\[ \Dc\ \xra{\ Y \ }\ \Vc \ \xra{- \tensor A}\ \Vc\ . \]
This construction is an action of $\Vc$ on $\Dc^\ast$, which means
preferred natural associativity and unit isomorphisms 
\[ Y\tensor (A\tensor B)\ \iso\  (Y\tensor A)\tensor B \text{\qquad and\qquad}
Y\tensor I\ \iso\  Y\ . \]
Indeed, the value of the above isomorphisms
at an object $d\in \Dc$ is simply the associativity (respectively unit)
isomorphism in the monoidal structure of $\Vc$, for the triple
of objects $(Y(d), A,B)$. 
\end{construction}

\begin{rk}[Enriched Yoneda lemma]\label{rk:enriched Yoneda}\index{subject}{Yoneda lemma!enriched}
The Yoneda lemma has enriched versions that we use frequently
to identify morphisms out of representable $\Vc$-functors.  
For every object $d\in\Dc$ the covariant hom functor
\[ d^*\ = \ \Dc(d,-) \ : \ \Dc\ \to \ \Vc \]
is a $\Vc$-functor, hence an object in the functor category $\Dc^\ast$.
The enriched Yoneda lemma comes in two flavors,  a $\Vc$-valued version
and a set-valued version.
For a $\Vc$-functor $Y$ in $\Dc^{\ast}$,
the evaluation morphism is the composite
\begin{align*}
  \underline{\Dc}^{\ast}(d^\ast,Y)\ \xra{\text{proj}_d} \ 
\underline{\Vc}(d^\ast(d),Y(d)) \ &= \ 
\underline{\Vc}(\Dc(d,d),Y(d)) \\ 
&\xra{\underline{\Vc}(\Id_d,Y(d))} \ \underline{\Vc}(I,Y(d)) \ \iso \ Y(d) 
\end{align*}
where $\Id_d:I\to\Dc(d,d)$ is the identity of $d$.
This composite is an isomorphism in $\Vc$, see for example
\cite[Sec.\,2.4]{kelly-enriched category}.
Passing to underlying sets produces a bijection
\[ \Dc^{\ast}(d^\ast,Y)\ \xra{\ \iso \ } \ u(Y(d)) \]
between the set of $\Vc$-natural transformations and the underlying set of $Y(d)$.
This is the weak form of the enriched Yoneda lemma, 
compare \cite[Sec.\,1.9]{kelly-enriched category}.
\end{rk}

\medskip

{\bf Convolution product.}\index{subject}{convolution product!of enriched functors}
We review a general method, due to B.\,Day \cite{day-closed}, 
for constructing symmetric monoidal structures on certain functor categories.
We now make the additional assumption that the index category $\Dc$
is a {\em symmetric monoidal} $\Vc$-category.
We denote the monoidal product on $\Dc$ by $\oplus$.

We get the main cases of interest by specializing as follows:
\begin{enumerate}[(i)]
\item For $\Vc=\bT$ the category of spaces under cartesian product
and $\Dc=\bL$ the topological category of inner product spaces
under orthogonal direct sum, 
$\Dc^\ast=\spc$ is the category of orthogonal spaces.
This yields the box product of orthogonal spaces.
\item For $\Vc=\bT_*$ the category of based spaces under smash product
and $\Dc=\bO$ the Thom space category of the orthogonal complement bundles,
under orthogonal direct sum,
$\Dc^\ast=\spec$ is the category of orthogonal spectra.
This yields the smash product of orthogonal spectra.
\item For $\Vc=\Ab$ the category of abelian groups under tensor product
and $\Dc=\bA$ the pre-additive Burnside category
with the monoidal structure of Theorem \ref{thm:bA monoidal},
$\Dc^\ast=\GF$ is the category of global functors.
This yields the box product of global functors.
\end{enumerate}

\medskip

We denote by $\Dc\tensor\Dc$ the $\Vc$-category whose objects are
pairs of $\Dc$-objects, and with morphism $\Vc$-objects
\[ (\Dc\tensor\Dc)((d,d'),(e,e'))\ = \ 
\Dc(d,e)\tensor\Dc(d',e')\  . \]
Given two $\Vc$-functors $X,Y:\Dc\to\Vc$, the composite
\[ X\tensor Y \ : \ \Dc \tensor \Dc \ \xra{X\tensor Y} \ 
\Vc \tensor\Vc \ \xra{\tensor} \ \Vc  \]
is a $\Vc$-functor.

\begin{defn}
A {\em bimorphism} $b:(X,Y)\to Z$\index{subject}{bimorphism}
from a pair of objects $(X,Y)$ of $\Dc^\ast$
to another object $Z$ of $\Dc^\ast$ is a $\Vc$-natural 
transformation $b:X\tensor Y \to Z\circ\oplus$ 
of $\Vc$-functors $\Dc\tensor\Dc\to\Vc$.
\end{defn}

So a bimorphism  $b:(X,Y)\to Z$ consists of  $\Vc$-morphisms
\[ b(d,e) \ : \ X(d) \tensor Y(e) \ \to \ Z(d\oplus e) \]
for all objects $d,e$ of $\Dc$, that form a $\Vc$-natural transformation.
We can then define the box product\index{subject}{box product} 
as a universal recipient of a bimorphism from $X$ and $Y$.

\begin{defn}\label{def:general box product}
A {\em box product} for objects $X$ and $Y$ of $\Dc^\ast$ 
is a pair $(X\Box Y,i)$ consisting of a $\Dc^\ast$-object $X\Box Y$
and a universal bimorphism $i:(X,Y)\to X\Box Y$,
i.e., a bimorphism 
such that for every $\Dc^\ast$-object $Z$ the map
\[ \Dc^\ast(X\Box Y,Z) \ \to \ \text{Bimor}((X,Y),Z) \ , \quad
f\longmapsto f\circ i \]
is bijective.\index{subject}{universal property!of a convolution product} 
\end{defn}

\begin{prop}\label{prop:box exists} 
Every pair of objects of $\Dc^\ast$ has a box product.
\end{prop}
\begin{proof}
The universal property of the box product precisely means that  
$(X\Box Y,i)$ is an enriched left Kan extension of the $\Vc$-functor
$X\tensor Y:\Dc\tensor \Dc\to\Vc$ along the $\Vc$-functor
$\oplus:\Dc\tensor\Dc\to\Dc$. Such a Kan extension exists because $\Vc$ is cocomplete
and $\Dc$ is skeletally small, see \cite[Prop.\,4.33]{kelly-enriched category}.
\end{proof}

\begin{rk}\label{rk:box as enriched coend}
As we saw in the proof of Proposition \ref{prop:box exists},
the box product is an enriched Kan extension
along the functor $\oplus:\Dc\tensor\Dc\to\Dc$. We can make this a little
more explicit. Given $X,Y\in\Dc^\ast$, a box product $X\Box Y$
is a coequalizer of the two morphisms of $\Vc$-functors
\begin{align*}
   \coprod_{e,e',d,d'} 
\Dc(e\oplus e',-)\tensor\Dc(d,e)&\tensor\Dc(d',e')\tensor X(d)\tensor Y(d') \\
&\xymatrix{\ar@<.4ex>[r]\ar@<-.4ex>[r] &}
\coprod_{d,d'} \Dc(d\oplus d',-)\tensor X(d)\tensor Y(d')
\ .\end{align*}
The left coproduct is indexed over all quadruples $(e,e',d,d')$ in a set
of representatives of the isomorphism classes of $\Dc$-objects.
The right coproduct is indexed over pairs $(d,d')$ of such representatives.
One of the two morphisms to be coequalized is the coproduct
of the monoidal products of 
\[ \Dc(e\oplus e',-)\tensor\Dc(d,e)\tensor\Dc(d',e')\ \to \ \Dc(d\oplus d',-)\ 
, \quad  (\varphi,\tau,\tau') \ \longmapsto \ \varphi\circ(\tau\oplus\tau') \]
and the identity on $X(d)\tensor Y(d')$.
The other morphism is the coproduct of the monoidal products of the identity on 
$\Dc(e\oplus e',-)$
and the action maps $\Dc(d,e)\tensor X(d)\to X(e)$
respectively  $\Dc(d',e')\tensor Y(d')\to Y(e')$.
\end{rk}

\begin{rk}[Uniqueness of box products]\label{rk:box unique}
The universal property makes box products of two functors unique up to preferred
isomorphism. Indeed, if $(X\Box Y,i)$ and $(X\Box' Y,i')$ 
are two box products, then the universal properties provide unique morphisms
$f:X\Box Y\to X\Box' Y$ and $g:X\Box' Y\to X\Box Y$ that satisfy
\[ f\circ i \ = \ i' \text{\qquad respectively\qquad} g\circ i'\ = \ i\ . \]
Then $(g\circ f)\circ i=g\circ i' = i = \Id_{X\Box Y}\circ i$,
so $g\circ f=\Id_{X\Box Y}$ by the uniqueness part of the universal property. 
Reversing the roles gives $f\circ g=\Id_{X\Box' Y}$.
\end{rk}

In the following we {\em choose} a box product $(X\Box Y,i)$
for every pair of objects $X$ and $Y$ of $\Dc^\ast$.
The universal bimorphism $i$ is often omitted from the notation, but one should
remember that the pair $(X\Box Y,i)$, and not the object $X\Box Y$ alone,
has a universal property.
While the box products are choices, the rest of
the structure is canonically determined by these.

\begin{construction}[Functoriality]\label{con:box functoriality}
The box product $X\Box Y$ automatically becomes a functor in both variables:
if $f:X\to X' $ and $g:Y\to Y'$ are morphisms in $\Dc^\ast$, 
then the $\Vc$-morphisms
\[ \Big\lbrace\ X(d) \tensor Y(e) \ \xrightarrow{f(d)\tensor g(e)} \
X'(d) \tensor  Y'(e) \ \xrightarrow{i'_{d,e}} \ 
(X'\Box Y')(d\oplus e)\ \Big\rbrace_{d,e} \]
form a bimorphism $(X,Y)\to X'\Box Y'$. 
So there is a unique morphism $f\Box g:X\Box Y\to X'\Box Y'$
such that $(f\Box g)(d\oplus e)\circ i_{d,e}=i'_{d,e}\circ (f(d)\tensor g(e))$
for all objects $d,e$ of $\Dc$.
The uniqueness part of the universal property implies that 
this is compatible with identities and composition in both variables.
\end{construction}

Now that we have constructed a box product functor
$ \Box : \Dc^\ast \times \Dc^\ast \to \Dc^\ast$,
we recall that it `is' automatically symmetric monoidal.
Since `symmetric monoidal' is extra data, and not a property,
we are obliged to construct unit, associativity and symmetry isomorphisms.

\begin{construction}[Coherence isomorphisms]\label{con:box coherence}
It will be convenient to make the $\Vc$-functor $0^\ast=\Dc(0,-)$
represented by the unit object~0 of $\Dc$
into a strict unit for the box product (as opposed to a unit
up to coherent isomorphisms). So we make the following conventions:
\begin{itemize}
\item (Right unit) We choose $X\Box 0^\ast=X$ with universal bimorphism
$i:(X,0^\ast)\to X$ given by the maps
\begin{align*}
  X(d)\tensor 0^\ast(e)\ = \ X(d)\tensor \Dc(0,e)\  \xra{X(d)\tensor (d\oplus-)} \ 
&X(d)\tensor \Dc(d\oplus 0,d\oplus e)\\ 
\iso \  X(d)&\tensor \Dc(d,d\oplus e)\ \xra{\ \circ \ }\  X(d\oplus e) \ .   
\end{align*}
\item (Left unit) We choose $0^\ast\Box Y=Y$ with universal bimorphism
$i:(0^\ast,Y)\to Y$ given by the composite
\begin{align*}
0^\ast(d)\tensor Y(e) \ = \ \Dc(0,d)\tensor Y(e) \ \xra{(-\oplus e)\tensor Y(e)}\
& \Dc(0\oplus e,d\oplus e)\tensor Y(e) \\
\iso\  \Dc&(e,d\oplus e)\tensor Y(e) \ \xra{\ \circ \ } Y(d\oplus e) \ . 
\end{align*}
\end{itemize}
The unnamed isomorphisms are induced by the unit isomorphisms in $\Dc$,
and `$\circ$' is the functor structure of $X$ respectively $Y$.
For $X=Y=0^\ast$ these two bimorphisms $(0^\ast,0^\ast)\to 0^\ast$
are equal, so there is no ambiguity.

We also obtain the associativity and symmetry isomorphisms
from the universal property of the box products.
For the associativity isomorphism we notice that the family
\[ \Big\lbrace i_{d,e\oplus f}\circ (X(d)\tensor i_{e,f})\, : \ 
X(d) \tensor \ Y(e) \tensor Z(f) \ \to \
 (X\Box (Y\Box Z))(d\oplus e\oplus f) \Big\rbrace_{d,e,f} \]
and the family
\[ \Big\lbrace i_{d\oplus e,f}\circ(i_{d,e}\tensor Z(f))\, : \  
 X(d) \tensor \ Y(e) \tensor Z(f) \ \to \
 ((X\Box Y)\Box Z)(d\oplus e\oplus f)\Big\rbrace_{d,e,f} \]
both have the universal property of a {\em tri}\,morphism
(whose definition is hopefully clear) out of $X$, $Y$ and $Z$.
The uniqueness of representing objects gives a unique isomorphism
of $\Vc$-functors
\[ \alpha_{X,Y,Z} \ : \ X\Box (Y\Box Z)\ \iso\   (X\Box Y)\Box Z  \]
such that $(\alpha_{X,Y,Z})_{d\oplus e\oplus f}\circ i_{d,e\oplus f}\circ (X(d)\tensor i_{e,f})
=i_{d\oplus e,f}\circ (i_{d,e}\tensor Z(f))$. 

The symmetry isomorphism $\tau_{X,Y}:X\Box Y\to Y\Box X$ corresponds to the
bimorphism from $(X,Y)$ to $Y\Box X$ with components
\[  X(d) \tensor  Y(e)\ \xrightarrow{\tau_{X(d),Y(e)}} \
Y(e) \tensor  X(d)\ \xrightarrow{i_{e,d}} \ (Y\Box X)(e\oplus d)
\xrightarrow{(Y\Box X)(\tau_{e,d})} \ (Y\Box X)(d\oplus e) \ . \]
\end{construction}

The following theorem is due to B.\,Day. 
In fact, Theorems~3.3 and~4.1 of \cite{day-closed} 
together imply the `monoidal'  part, and
Theorem~3.6 of \cite{day-closed} deals with the symmetries.

\begin{theorem}\label{thm:symmetric monoidal}
The associativity and symmetry isomorphisms make the box product
into a symmetric monoidal product on the category $\Dc^\ast$ with 
the functor $0^\ast$ as a strict unit object. 
\end{theorem}

The box product also commutes with $\Vc$-tensors up to a natural isomorphism
that we now construct.
The family of morphisms
\[ X(d) \tensor ( Y(e) \tensor  A )
 \ \xra{\alpha_{X(d),Y(e),A}} \
  (X(d) \tensor Y(e))  \tensor A \  \xrightarrow{i_{d,e}\tensor A} \
(X\Box Y)(d\oplus e) \tensor A \]
forms a bimorphism from $(X,Y\tensor A)$ to $(X\Box Y)\tensor A$.
So the universal property provides a distinguished morphism
\[ \alpha_{X,Y,A}\ : \ X\Box(Y\tensor A)\ \to \ ( X\Box Y)\tensor A  \ .\]
Since the monoidal structure on $\Vc$ is closed, the functor $-\tensor A$
is a left adjoint, so it preserves coends. So the morphism
$\alpha_{X,Y,A}$ is an isomorphism.
The two types of associativity isomorphisms 
$X\Box (Y\tensor A)\iso (X\Box Y)\tensor A$ and 
$Y\tensor (A\tensor B)\iso (Y\tensor A)\tensor B$
are compatible in the sense of a commuting pentagon:
\[\xymatrix@C=-5mm{
&& X\Box (Y\tensor (A\tensor B)) \ar[dll]_{X\Box\alpha_{Y,A,B}}
\ar[drr]^{\alpha_{X,Y,A\tensor B}}\\
X\Box ((Y\tensor A)\tensor B) \ar[dr]_{\alpha_{X,Y\tensor A,B}}
&&&& (X\Box Y)\tensor (A\tensor B)\ar[dl]^{\alpha_{X\Box Y,A,B}}\\
& (X\Box(Y\tensor A))\tensor B \ar[rr]_{\alpha_{X,Y,A}\tensor B} &&
((X\Box Y)\tensor A)\tensor B
}\]
A minimally more complicated construction that also involves
the symmetry constraint in $\Vc$ to move the object $A$ past the functor $Y$
provides an isomorphism
\[ (X\tensor A)\Box Y\ \iso \   (X\Box Y)\tensor A \ . \]

\begin{rk}[Convolution product of representable functors]\label{rk:box representable}
The box product of the functors represented by two objects $a,b$ of $\Dc$ 
is related by a preferred isomorphism to the functor represented by $a\oplus b$.
Indeed, the $\Vc$-morphisms
\[ a^\ast(d)\tensor b^\ast(e) \ = \ \Dc(a,d)\tensor\Dc(b,e) \ \xra{\ \oplus \ } \
\Dc(a\oplus b,d\oplus e)\ = \ (a\oplus b)^\ast(d\oplus e) \]
form a bimorphism as $d$ and $e$ vary. So there is a unique morphism
\[ j_{a,b}\ : \  a^\ast\Box b^\ast \ \to \ (a\oplus b)^\ast \]
such that $j_{a,b}\circ i=\oplus$, where $i$ is 
the universal bimorphism from $(a^\ast,b^\ast)$ to $a^*\Box b^*$. 
On the other hand, the morphism
\[ I \ \iso \ I\tensor I \ \xra{\Id_a\tensor\Id_b} \
\Dc(a,a)\tensor\Dc(b,b) \ = \ a^*(a)\tensor b^*(b) \ \xra{i_{a,b}}\
(a^*\Box b^*)(a\oplus b)  \]
is an element of the underlying set of $(a^*\Box b^*)(a\oplus b)$.
The enriched Yoneda lemma thus provides a unique morphism
of $\Vc$-functors $(a\oplus b)^*\to a^*\Box b^*$
that evaluates to the element at $a\oplus b$. 
The universal properties of the convolution product and
of the representable functors show that these two morphisms
are inverse to each other. 
\end{rk}

\begin{rk}[Internal function objects]
The box product is a {\em closed} monoi\-dal product
in the sense that for all objects $Y$ and $Z$ of $\Dc^\ast$ the functor
\[ \Dc^\ast(-\Box Y, Z) \ : \ \Dc^\ast \ \to\ \text{(sets)} \ , \quad
X\ \longmapsto \ \Dc^\ast(X\Box Y,Z) \]
is representable.
In this book we don't make any substantial use of the internal function objects,
so we will not elaborate on the construction. The most relevant consequence
for us is that the box product preserves colimits
in both variables because the functors $X\Box -$ and $-\Box Y$
are left adjoints. This consequence can also be shown 
without reference to function objects by observing that
the original monoidal product in $\Vc$ preserves
colimits in both variables (because the monoidal structure of $\Vc$ is closed),
because $X\Box Y$ can be constructed as an enriched coend,
and because enriched coends commute with colimits.
\end{rk}

\begin{construction}[Skeleton filtration]\label{con:skeleta in functor categories}
For the rest of this appendix we assume that the index category $\Dc$ 
is equipped with a dimension function $\dim:\text{ob}(\Dc)\to\mN$ 
to the natural numbers, satisfying the following conditions for all
objects $d,e$ of $\Dc$:
\begin{itemize}
\item If $\dim(e)<\dim(d)$, then $\Dc(d,e)$ 
is an initial object of the base category $\Vc$.
\item If $\dim(e)=\dim(d)$, then $d$ and $e$ are isomorphic.
\end{itemize}
The two examples that we care about in this book are
the category of orthogonal spaces and the category of orthogonal spectra.
In both cases the dimension function is the dimension as an $\mR$-vector space.

We denote by $\Dc_{\leq m}$ the full $\Vc$-subcategory of $\Dc$
spanned by all objects of dimension at most $m$.
We denote by $\Dc^\ast_{\leq m}=(\Dc_{\leq m})^\ast$ 
the category of enriched functors from $\Dc_{\leq m}$
to $\Vc$. The restriction functor 
\[  \Dc^\ast \ \to \ \Dc^\ast_{\leq m}\ , \quad 
Y \ \longmapsto \ Y^{\leq m}\ = \ Y|_{\Dc_{\leq m}}\]
has a left adjoint 
\[ l_m \ : \ \ \Dc^\ast_{\leq m} \ \to \ \Dc^\ast \]
given by an enriched Kan extension as follows. 
For every $k\geq 0$ we choose an object $\mathbf k$ of $\Dc$ of dimension $k$.
The extension $l_m(Z)$ of an enriched functor $Z:\Dc_{\leq m}\to\Vc$ 
is a coequalizer of the two morphisms in $\Dc^\ast$:
\begin{equation}  \label{eq:skeleton_coequalizer}
 \xymatrix@C=8mm{ 
\coprod_{0\leq j\leq k\leq m} \Dc(\mathbf k,-)\tensor\Dc(\mathbf j,\mathbf k)\tensor Z(\mathbf j) 
 \ar@<-.4ex>[r] \ar@<.4ex>[r]  &
\coprod_{0\leq i\leq m} \Dc(\mathbf i,-)\tensor Z(\mathbf i) 
 }   
\end{equation}
One morphism arises from the composition morphisms
\[  \Dc(\mathbf k,-)\tensor\Dc(\mathbf j,\mathbf k) \ \to \ \Dc(\mathbf j,-) \]
and the identity on $Z(\mathbf j)$;
the other morphism arises from the action morphisms
\[  \Dc(\mathbf j,\mathbf k)\tensor Z(\mathbf j) \ \to \ Z(\mathbf k) \]
and the identity on the represented functor $\Dc(\mathbf k,-)$.
Colimits in the functor category $\Dc^\ast$ are created objectwise, 
so the value $l_m(Z)(d)$ at an object $d$ can be calculated 
by plugging $d$ into the variable slot 
in the coequalizer diagram \eqref{eq:skeleton_coequalizer}.
Kan extensions along a fully faithful functor do not change
the values on the given subcategory, 
see for example \cite[Prop.\,4.23]{kelly-enriched category}. 
For the inclusion $\Dc_{\leq m}\to\Dc$ this means that the adjunction unit
\[ Z \ \to \ (l_m Z)^{\leq m} \]
is an isomorphism for every functor $Z:\Dc_{\leq m}\to\Vc$.
\end{construction}

\begin{defn}
The {\em $m$-skeleton},\index{subject}{skeleton}\index{symbol}{$\sk^m X$ - {$m$-skeleton of $X$}}
for $m\geq 0$, of an enriched functor $Y:\Dc\to \Vc$ is the functor
\[ \sk^m Y\ = \ l_m(Y^{\leq m}) \ ,\]
the extension of the restriction of $Y$ to $\Dc_{\leq m}$.
It comes with a natural morphism $i_m:\sk^m Y\to Y$, the counit of
the adjunction $(l_m,(-)^{\leq m})$.
We set $\Dc(m)=\Dc(\mathbf m,\mathbf m)$, a $\Vc$-monoid under composition.
The {\em $m$-th latching object}\index{subject}{latching object!of an enriched functor}\index{symbol}{$L_m Y$ - {$m$-th latching object of $Y$}}
of $Y$ is the $\Dc(m)$-object
\[ L_m Y \ = \ (\sk^{m-1} Y)(\mathbf m) \ ;\]
it comes with a natural $\Dc(m)$-equivariant morphism
\[  \nu_m=i_{m-1}(\mathbf m)\ :\ L_m Y\ \to \ Y(\mathbf m) \ , \]
the {\em $m$-th latching morphism}.
\end{defn}

We agree to set $\sk^{-1} Y=\emptyset$, an initial object of $\Dc^\ast$;
then $L_0 Y=\emptyset$, an initial $\Vc$-object.
The two morphisms $i_{m-1}:\sk^{m-1} Y\to Y$ and
$i_m:\sk^m Y\to Y$ both restrict to isomorphisms on $\Dc^{\leq m-1}$, so there is a 
unique morphism $j_m:\sk^{m-1} Y\to \sk^m Y$ such that $i_m\circ j_m=i_{m-1}$.
The sequence of skeleta stabilizes to $Y$ in a 
very strong sense. For every object $d$, the morphisms
$j_m(d)$ and $i_m(d)$ are isomorphisms as soon as $m > \dim(d)$.
In particular, $Y(d)$ is a colimit, with respect to the morphisms $i_m(d)$, 
of the sequence of morphisms $j_m(d)$.
Since colimits in the functor category $\Dc^\ast$ are created objectwise,
we deduce that $Y$ is a colimit, 
with respect to the morphisms $i_m$, of the sequence of morphisms $j_m$.

\begin{eg}[Latching objects of represented functors]\label{eg:latching of represented}
  Let $d$ be an object of $\Dc$ of dimension $m=\dim(d)$, 
  and $A$ an object of $\Vc$.
  Then the functor $\Dc(d,-)\tensor A$
  is `purely $m$-dimensional' in the following sense.
  The evaluation functor
  \[ \ev_d\ : \ \Dc^\ast \ \to \ \Vc \]
  factors through the category $\Dc_{\leq m}^*$ as the composite
  \[  \Dc^\ast \ \to \  \Dc_{\leq m}^* \ \xra{\ev_d} \ \Vc \ .\]
  So the left adjoint free functor $A\mapsto \Dc(d,-)\tensor A$ 
  can be chosen as the composite
  of the two individual left adjoints
  \[ \Dc(d,-)\tensor - \ = \ l_m \circ (\Dc_{\leq m}(d,-)\tensor -) \ .\]
  The object $\Dc(d, e)\tensor A$ is initial for $\dim(e) < m$, by hypothesis
  on the dimension function, and hence the 
  latching object $L_n (\Dc(d,-)\tensor A)$ is initial for $n\leq m$.
  For $n>m$ the latching morphism 
  $\nu_n:L_n(\Dc(d,-)\tensor  A) \to \Dc(d,\mathbf n)\tensor A$ 
  is an isomorphism. So the skeleton $\sk^n (\Dc(d,-)\tensor A)$ is initial for $n<m$
  and $\sk^n (\Dc(d,-)\tensor A)=\Dc(d,-)\tensor A$ is the entire functor for $n \geq m$.
\end{eg}

We denote by
\[ G_m \ : \ \Dc(m)\Vc \ \to \ \Dc^\ast \ , \quad 
A \ \longmapsto \ \Dc(\mathbf m,-)\tensor_{\Dc(m)} A \]
the left adjoint to the evaluation functor $Y\mapsto Y(\mathbf m)$.
As a consequence of the previous example, the enriched functor $G_m A$
is purely $m$-dimensional for every $\Dc(m)$-object $A$.

\begin{prop}\label{prop:basic skeleton pushout square}
For every enriched functor $Y:\Dc\to\Vc$ and every $m\geq 0$ the commutative square
\begin{equation}\begin{aligned}\label{eq:filtration pushout D^ast}
\xymatrix@C=15mm{ G_m L_m Y \ar[r]^{G_m\nu_m} \ar[d] & G_m Y(\mathbf m) \ar[d]\\
\sk^{m-1} Y \ar[r]_-{j_m} & \sk^m Y}  
\end{aligned}\end{equation}
is a pushout in the category $\Dc^\ast$.
The two vertical morphisms are adjoint to the identity of $L_m Y$ 
respectively $Y(\mathbf m)$.
\end{prop}
\begin{proof}
All four functors are `$m$-dimensional', i.e., isomorphic
to the extensions of their restrictions to $\Dc_{\leq m}$;
for $G_m L_m Y$ and $G_m Y(\mathbf m)$ 
this follows from Example \ref{eg:latching of represented}.
So it suffices to check that the square is a pushout in $\Vc$
when evaluated at objects $d$ of dimension at most $m$.
When the dimension of $d$ is strictly less than $m$, 
then $(G_m L_m Y)(d)$ and $(G_m Y_m)(d)$ are initial objects and the morphism
$j_m(d):(\sk^{m-1} Y)(d)\to(\sk^m Y)(d)$ is an isomorphism.
So for $\dim(d) < m$ both horizontal morphisms are isomorphisms at $d$.
When $\dim(d)=m$, then $d$ is isomorphic to $\mathbf m$, 
so we may suppose that $d=\mathbf m$.
Then the square \eqref{eq:filtration pushout D^ast} evaluates to
\[ \xymatrix@C=25mm{ 
\Dc(\mathbf m,\mathbf m)\tensor_{\Dc(m)} L_m Y
 \ar[r]^{\Dc(\mathbf m,\mathbf m)\tensor_{\Dc(m)}\nu_m} \ar[d] & 
\Dc(\mathbf m,\mathbf m)\tensor_{\Dc(m)} Y(\mathbf m) \ar[d]\\
L_m Y \ar[r]_-{j_m(\mathbf m)} &  Y(\mathbf m) }   \]
So for $\dim(d)= m$, both vertical maps are isomorphisms at $d$.
\end{proof}

Given any morphism $f:A\to B$ of enriched functors in $\Dc^\ast$, we
define a relative skeleton using the commutative square of enriched functors:
\[\xymatrix@C=12mm{ \sk^m A\ar[r]^{\sk^m f}\ar[d]_{i_m} &  \sk^m B \ar[d]^{i_m}\\
A\ar[r]_f & B } \]
The {\em relative $m$-skeleton} of $f$ is the pushout
\begin{equation}\label{eq:abstract_relative_skeleton}
 \sk^m [f] \ = \ A\cup_{\sk^m A}\sk^m B  \ ;
\end{equation}
it comes with a unique morphism $i_m:\sk^m [f]\to B$
which restricts to $f:A\to B$ respectively to $i_m:\sk^m B\to B$.
Since $(\sk^{m-1} A)(\mathbf m)=L_m A$ we have
\[ (\sk^{m-1} [f])(\mathbf m) \ = \ A(\mathbf m)\cup_{L_m A}L_m B\ , \]
the $m$-th relative latching object.
A morphism $j_m[f]:\sk^{m-1} [f]\to \sk^m [f]$ is obtained from the commutative
diagram
\[\xymatrix@C=15mm{
A \ar@{=}[d] & \sk^{m-1} A \ar[r]^-{\sk^{m-1} f} \ar[l]_(.6){i_{m-1}} \ar[d]_{j_m^A} & 
\sk^{m-1} B\ar[d]^{j_m^B}\\
A  & \sk^m A \ar[r]_-{\sk^m f} \ar[l]^{i_m} & \sk^m B }\]
by taking pushouts.
The original morphism $f:A\to B$ 
factors as the composite of the countable sequence
\[ A=\sk^{-1}[f]\ \xra{j_0[f]}\ \sk^0[f] \ \xra{j_1[f]}\  \sk^1[f] \ \to\
\cdots \ \xra{j_m[f]}\  \sk^m [f] \ \to\ \cdots \ .     \]
If $d$ has dimension $n$, then the sequence stabilizes to 
the identity map of $B(d)$ from $(\sk^n[f])(d)$ onward; in particular,
the compatible morphisms $j_m$ exhibit $B$ as a colimit of the sequence.
The pushout square \eqref{eq:filtration pushout D^ast}
also has a relative version.
We contemplate the following commutative diagram:
\begin{equation}\begin{aligned}\label{eq:abstract_relative_skeleton_pushout}    
\xymatrix@C=15mm{ 
G_m(A(\mathbf m)\cup_{L_m A}L_m B)  \ar[r]^-{G_m(\nu_m f)}\ar[d] & 
G_m B(\mathbf m)\ar[d] \\
\sk^m A\cup_{\sk^{m-1} A}\sk^{m-1}B  \ar[r] \ar[d] & 
\sk^m B\ar[d] \\
\sk^{m-1} [f] \ar[r]_-{j_m[f]} & \sk^m [f] }
\end{aligned}\end{equation}
All three horizontal morphisms are isomorphisms at objects of dimension less
than $m$. At the object $\mathbf m$, the two upper vertical morphisms
are isomorphisms.
All four objects in the upper square are `$m$-dimensional', i.e., isomorphic
to the extensions of their restrictions to $\Dc_{\leq m}$;
so the upper square is a pushout of enriched functors.
The lower square is also a pushout, by definition of the relative skeleta.
We conclude that the outer composite square 
in the diagram \eqref{eq:abstract_relative_skeleton_pushout}    
is a pushout in the category $\Dc^\ast$.

\medskip

The following proposition is an immediate application of the
relative skeleton filtration. It is the key ingredient to
the lifting properties of the various level model structures 
in this book.
We recall that a pair $(i:A\to B,f:X\to Y)$ of morphisms in some category
has the {\em lifting property}\index{subject}{lifting property!of a pair of morphisms}
if for all morphisms $\alpha:A\to X$ and $\beta:B\to Y$ such that
$f\alpha=\beta i$ there exists a {\em lifting}, i.e.,
a morphism $\lambda:B\to X$ such that $\lambda i=\alpha$ and $f\lambda=\beta$.
Instead of saying that the pair $(i,f)$ has the lifting property
we also say `$i$ has the left lifting property
with respect to $f$' or `$f$ has the right lifting property with respect to $i$'.
\index{subject}{lifting property!left}
\index{subject}{lifting property!right}

\begin{prop}\label{prop:lifting property in spaces} 
Let $i:A\to B$ and $f:X\to Y$ be morphisms of enriched functors in $\Dc^\ast$.
If the pair 
\[ (\nu_m i=i(\mathbf m)\cup \nu_m^B:A(\mathbf m)\cup_{L_m A}L_m B\to B(\mathbf m),\quad 
f(\mathbf m):X(\mathbf m)\to Y(\mathbf m)) \]
has the lifting property in the category of $\Dc(m)$-objects for
every $m\geq 0$, then the pair $(i,f)$ has the lifting property 
in the functor category $\Dc^\ast$.
\end{prop}
\begin{proof}
We consider the class $f$-cof of all morphisms in $\Dc^\ast$
that have the left lifting property with respect to $f$; this class is closed
under cobase change and countable composition.
Since the pair $(\nu_m i,f_m)$
has the lifting property in the category of $\Dc(m)$-objects,
the morphism $G_m(\nu_m i)$ belongs to the class $f$-cof by adjointness.
The relative skeleton filtration \eqref{eq:abstract_relative_skeleton} 
shows that $i$ is a countable composite of cobase changes 
of the morphisms $G_m(\nu_m i)$, so $i$ belongs to the class $f$-cof.
\end{proof}

Now we discuss  a general recipe for constructing 
level model structures on the functor category $\Dc^\ast$.
As input we need, for every $m\geq 0$, a model structure $\Cc(m)$
on the category of $\Dc(m)$-objects.
We call a morphism $f:X\to Y$ in $\Dc^\ast$
\begin{itemize}
\item a {\em level equivalence} if  $f(\mathbf m):X(\mathbf m)\to Y(\mathbf m)$ 
is a weak equivalence in the model structure $\Cc(m)$ for all $m\geq 0$;
\item a {\em level fibration} if the  morphism
$f(\mathbf m):X(\mathbf m)\to Y(\mathbf m)$ is a fibration in the model structure 
$\Cc(m)$ for all $m\geq 0$;
\item a {\em cofibration} if the latching morphism
$\nu_m f:X(\mathbf m)\cup_{L_m X}L_m Y\to Y(\mathbf m)$ 
is a cofibration in the model structure 
$\Cc(m)$ for all $m\geq 0$.
\end{itemize}
Proposition \ref{prop:general level model structure} below shows 
that if the various model structures $\Cc(m)$
satisfy the following `consistency condition', then the 
level equivalences,  level fibrations and cofibrations define a model structure
on the functor category $\Dc^\ast$.

\begin{defn}[Consistency condition]\label{def:consistency condition}
\index{subject}{consistency condition}
For all $m,n\geq 0$ and every acyclic cofibration $i:A\to B$ in 
the model structure $\Cc(m)$ on $\Dc(m)$-objects,
every cobase change, in the category of $\Dc(m+n)$-objects, of the morphism
\[  \Dc(\mathbf m,\mathbf{m+n})\tensor_{\Dc(m)} i\ : \
 \Dc(\mathbf m,\mathbf {m+n})\tensor_{\Dc(m)} A\ \to\ 
\Dc(\mathbf m,\mathbf {m+n})\tensor_{\Dc(m)} B  \]
is a weak equivalence in the model structure $\Cc(m+n)$.
\end{defn}

For example, the consistency conditions holds if the functor
\[ \Dc(\mathbf m,\mathbf{m+n})\tensor_{\Dc(m)} -\ :\ \Dc(m)\Vc\ \to\ \Dc(m+n)\Vc \]
is a left Quillen functor in the given model structures, for all $m,n\geq 0$.

\begin{prop}\label{prop:general level model structure} 
  Let $\Cc(m)$ be a model structure 
  on the category of $\Dc(m)$-objects, for $m\geq 0$,
  such that the consistency condition holds.
   \begin{enumerate}[\em (i)]
  \item 
    The classes of level equivalences, 
    level fibrations and cofibrations define a model structure
    on the functor category $\Dc^\ast$.
  \item 
    A morphism $i:A\to B$ in $\Dc^\ast$ is simultaneously a 
    cofibration and a level equivalence if and only if 
    for all $m\geq 0$ the latching morphism 
    $\nu_m i:A(\mathbf m)\cup_{L_m A}L_m B\to B(\mathbf m)$ is an acyclic cofibration 
    in the model structure $\Cc(m)$.
  \item Suppose that the fibrations in the model structure $\Cc(m)$ 
    are detected by a set of morphisms $J(m)$; 
    then the level fibrations are detected by the set of morphisms 
    \[  \{G_m j \ | \ m\geq 0,\,  j\in J(m)\} \ .\]
    Similarly, if the acyclic fibrations in the model structure $\Cc(m)$ 
    are detected by a set of morphisms $I(m)$, 
    then the level acyclic fibrations
    are detected by the set of morphisms 
    \[  \{G_m i \ | \ m\geq 0,\,  i\in I(m)\} \ .\]
\end{enumerate}
\end{prop}
\begin{proof} 
We start by showing one of the directions of part~(ii):
we let $i:A\to B$ be a morphism 
such that the latching morphism $\nu_m i:A(\mathbf m)\cup_{L_m A}L_m B\to B(\mathbf m)$ 
is an acyclic cofibration in the model structure 
$\Cc(m)$ for all $m\geq 0$; we show that then $i$ is a level equivalence.

The map $i(\mathbf n):A(\mathbf n)\to B(\mathbf n)$ is the finite composite 
\[ A(\mathbf n)=(\sk^{-1}[i])(\mathbf n)  \xra{(j_0[i])(\mathbf n)} 
(\sk^0[i])(\mathbf n)\ \xra{(j_1[i])(\mathbf n)} \dots
\xra{(j_n[i])(\mathbf n)} (\sk^n[i])(\mathbf n)=B(\mathbf n) \ ,\]
so it suffices to show that $j_k[i]$
is a level equivalence for all $k\geq 0$.
The pushout square \eqref{eq:abstract_relative_skeleton_pushout}
in level $m+n$ is a pushout of $\Dc(m+n)$-objects
\[ \xymatrix@C=20mm{ 
\Dc(\mathbf m,\mathbf {m+n})\tensor_{\Dc(m)}(A(\mathbf m)\cup_{L_m A}L_m B)
\ar[r]^-{\Dc(\mathbf m,\mathbf {m+n})\tensor_{\Dc(m)} (\nu_m i)} \ar[d]
 & \Dc(\mathbf m,\mathbf {m+n})\tensor_{\Dc(m)} B(\mathbf m) \ar[d] \\
(\sk^{m-1} [i])(\mathbf {m+n}) \ar[r]_-{(j_m[i])(\mathbf{m+n})} & (\sk^m [i])(\mathbf{m+n}) } \]
The consistency condition guarantees that
the lower horizontal morphism is a $\Cc(m+n)$-weak equivalence.

(i)
  Several of the axioms are straightforward:
  the functor category $\Dc^\ast$ inherits all small limits and colimits
  from the base category $\Vc$; 
  the level equivalences satisfy the 2-out-of-3 property;
  the classes of level equivalences, cofibrations and 
  fibrations are closed under retracts. 
  
  Now we prove the factorization axiom, i.e., we show that every morphism 
  $f:A\to X$ in $\Dc^*$ can be factored as $f=q i$ 
  where $q$ is a level acyclic fibration and $i$ a cofibration;
  and it can be factored as $f=p j$ 
  where $p$ is a level fibration and $j$ a cofibration and level equivalence.
  We start with the first factorization and construct 
  an enriched functor $B$ and morphisms
  $i:A\to B$ and $q:B\to X$ by induction over the dimension.
  Since the objects $\mathbf m$ for $m\geq 0$ form a skeleton
  of the category $\Dc$, it suffices to construct all the relevant data on these 
  objects. 
  In level $0$ we choose a factorization
  \[ A(\mathbf 0) \ \xra{\ i(\mathbf 0)\ } \ B(\mathbf 0) \ \xra{\ q(\mathbf 0)\ } 
  \ X(\mathbf 0)\]
  of $f(\mathbf 0)$ in the category of $\Dc(0)$-objects such that $i(\mathbf 0)$ 
  is a cofibration
  and $q(\mathbf 0)$ is an acyclic fibration in the model structure $\Cc(0)$.
  Now we suppose that the enriched functor $B$ and the morphisms
  $i$ and $q$ have already been constructed for all objects of dimension less than $m$.
  Then we have all the data necessary to define the
  $m$-th latching object $L_m B$; moreover, the `partial morphism'
  $q:B\to X$ provides a $\Dc(m)$-morphism $L_m B\to X(\mathbf m)$
  such that the square
  \[ \xymatrix{L_m A \ar[r]^-{L_m i} \ar[d]_{\nu_m^A} & L_m B \ar[d] \\
    A(\mathbf m) \ar[r]_-{f(\mathbf m)} & X(\mathbf m)}\]
  commutes. We factor the resulting morphism 
  $A(\mathbf m)\cup_{L_m A} L_m B\to X(\mathbf m)$ 
  in the category of $\Dc(m)$-objects
  \begin{equation}\label{eq:mth factor spaces}
    A(\mathbf m)\cup_{L_m A} L_m B\ \xra{\ \nu_m i\ }\ B(\mathbf m) 
    \ \xra{\ q(\mathbf m)\ }\ X(\mathbf m)    
  \end{equation}
  such that $\nu_m i$ is a cofibration
  and $q(\mathbf m)$ is an acyclic fibration in the model structure $\Cc(m)$.
  The intermediate $\Dc(m)$-object $B(\mathbf m)$ defines the value of $B$
  at the object $\mathbf m$, 
  and $q(\mathbf m)$ defines the $m$-th level of the morphism $q$. 
  For $0\leq i<m$, the structure morphism
  is the composite
  \[ \Dc(\mathbf i,\mathbf m) \tensor B(\mathbf i)\ \to\ L_m B \ \to \  
  A(\mathbf m)\cup_{L_m A} L_m B \ 
  \xra{\nu_m i} \ B(\mathbf m) \]
  and the composite of $\nu_m i$ with the canonical morphism
  $A(\mathbf m)\to A(\mathbf m)\cup_{L_m A} L_m B$ is the $m$-th level of the morphism $i$.

  At the end of the day we have indeed factored $f=q i$
  in the category $\Dc^\ast$ such that $q$ is a level
  equivalence and level fibration. Moreover, the $m$-th
  latching morphism comes out to be the morphism
  $\nu_m i:A(\mathbf m)\cup_{L_m A} L_m B\to B(\mathbf m)$
  in the factorization \eqref{eq:mth factor spaces}, which is a cofibration
  in the model structure $\Cc(m)$. So the morphism $i$ is indeed a
  cofibration.

  The second factorization $f=p j$ as a cofibration and level equivalence $j$
  followed by a level fibration $p$ is similar, but instead of the 
  factorization \eqref{eq:mth factor spaces} we use a factorization,
  in the model category $\Cc(m)$, as an acyclic cofibration followed by 
  a fibration. Then the resulting morphism $p$ is a level fibration and
  the morphism $j$ has the property that all its latching morphisms $\nu_m j$ 
  are acyclic cofibrations. So $j$ is a cofibration
  (by definition) and a level equivalence 
  (by the part of~(ii) established above).

  It remains to show the lifting axioms. 
  In each of the model structures $\Cc(m)$ 
  the cofibrations have the left lifting property with respect 
  to the acyclic fibrations; so by 
  Proposition \ref{prop:lifting property in spaces} 
  the cofibrations in $\Dc^\ast$ have the left lifting property 
  with respect to level equivalences which are also level fibrations. 

  We postpone the proof of the other lifting property and 
  prove the remaining direction of~(ii) next.
  We let $i:A\to B$ be a cofibration and a level equivalence. 
  The second factorization axiom proved above provides a factorization
  $i=p j$ where $j:A\to D$ is a level equivalence such that 
  each latching morphism $\nu_m j$ is an acyclic cofibration
  in the model structure $\Cc(m)$,
  and $p:D\to B$ is a level fibration.
  Since $i$ and $j$ are level equivalences, so is $p$. 
  So the cofibration $i$ has the left lifting property 
  with respect to the level equivalence and level fibration $p$ 
  by the previous paragraph.
  In particular, a lift $\lambda:B\to D$ in the square
 \[     \xymatrix{ A\ar[r]^-j\ar[d]_i &  D\ar[d]^p\\
    B\ar@{=}[r] \ar@{..>}[ur]^(.4)\lambda & B} \]
  shows that the morphism $i$ is a retract of the morphism $j$.
  So the latching morphism $\nu_m i$ is a retract 
  of the latching morphism $\nu_m j$, hence also an acyclic cofibration
  in the model structure $\Cc(m)$. This proves~(ii).

  Now we prove the remaining lifting property.
  We let $i:A\to B$ be a cofibration that is also a level equivalence. 
  By~(ii), which has just been shown,
  each latching morphism $\nu_mi$ is an acyclic cofibration
  in the model structure $\Cc(m)$.
  So $i$ has the left lifting property with respect to all level fibrations
  by Proposition \ref{prop:lifting property in spaces}.   

  Property~(iii) is a straightforward consequence of
  the fact that the functor $G_m$ is left adjoint to evaluation at $\mathbf m$.
\end{proof}

\endappendix

\renewcommand{\refname}{Bibliography}

\printindex{symbol}{List of symbols}
\printindex{subject}{Index}

\end{document}